\documentclass[10pt,reqno]{amsart}

\usepackage[margin=
1
in]{geometry}

\usepackage[x11names]{xcolor} 
\usepackage{amssymb}
\usepackage{amsmath}
\usepackage{mathrsfs}
\usepackage{amsfonts, dsfont}
\usepackage{kotex}
\counterwithin{figure}{section}
\usepackage[colorlinks=true, pdfstartview=FitV, linkcolor=black, citecolor=black, urlcolor=black]{hyperref} 
\definecolor{grey}{gray}{0.5}
\usepackage{setspace}

\newif\ifEDreview
\EDreviewfalse 
\newif\ifLBshowred
\LBshowredfalse

\newif\ifCNreview
\CNreviewfalse

\makeatletter
\ifEDreview
\AtBeginDocument{\@ifundefined{r@eq:dyn-decay-lambda-delta}{%
  \expandafter\gdef\csname r@eq:dyn-decay-lambda-delta\endcsname
    {{13.69}{330}{}{equation.13.69}{}}}{} }
\fi
\makeatother

\let\hide\iffalse

\usepackage{amsthm} 
\usepackage[utf8]{inputenc} 
\usepackage{enumitem} 
\allowdisplaybreaks

\usepackage{graphicx}

\usepackage{mathtools}
\usepackage{esint}

\usepackage[utf8]{inputenc}
\usepackage{anyfontsize}

\usepackage{tikz}
\usetikzlibrary{
	calc,intersections,arrows.meta,patterns,
	decorations.pathmorphing,
	decorations.fractals
}

\DeclareGraphicsExtensions{.pdf,.png,.jpg}

\usepackage{color}

\theoremstyle{plain}

\newtheorem{thm}{Theorem}[section]
\newtheorem*{thm*}{Theorem}
\newtheorem{corollary}[thm]{Corollary}

\newtheorem{lemma}[thm]{Lemma}
\newtheorem*{lemma*}{Lemma}
\newtheorem{prop}[thm]{Proposition}
\newtheorem*{prop*}{Proposition}

\newtheorem*{statement*}{Statement}

\newtheorem*{assumption*}{Assumption}

\theoremstyle{definition}   
\newtheorem{defn}[thm]{Definition}

\theoremstyle{remark}  
\newtheorem{remark}[thm]{Remark}

\newtheorem{innercustomgeneric}{\customgenericname}
\providecommand{\customgenericname}{}

\newcommand{\newcustomtheorem}[2]{\newenvironment{#1}[1]
	{\renewcommand\customgenericname{#2}
		\renewcommand\theinnercustomgeneric{##1}\innercustomgeneric}{\endinnercustomgeneric}}

\newcustomtheorem{customthm}{Theorem}

\newcommand\bL{\mathbb{L}}

\newcommand\bR{\mathbb{R}}

\newcommand\bZ{\mathbb{Z}}
\newcommand\bW{\mathbb{W}}

\newcommand\bS{\mathbb{S}}

\newcommand\bE{\mathbb{E}}
\newcommand\bN{\mathbb{N}}

\newcommand\bP{\mathbb{P}}

\newcommand\cA{\mathcal{A}}
\newcommand\cB{\mathcal{B}}
\newcommand\cC{\mathcal{C}}
\newcommand\cD{\mathcal{D}}
\newcommand\cF{\mathcal{F}}
\newcommand\cG{\mathcal{G}}
\newcommand\cH{\mathcal{H}}
\newcommand\cI{\mathcal{I}}
\newcommand\cJ{\mathcal{J}}
\newcommand\cK{\mathcal{K}}
\newcommand\cL{\mathcal{L}}
\newcommand\cP{\mathcal{P}}
\newcommand\cR{\mathcal{R}}
\newcommand\cS{\mathcal{S}}
\newcommand\cM{\mathcal{M}}
\newcommand\cN{\mathcal{N}}
\newcommand\cT{\mathcal{T}}
\newcommand\cO{\mathcal{O}}
\newcommand\cU{\mathcal{U}}
\newcommand\cV{\mathcal{V}}

\newcommand\cZ{\mathcal{Z}}
\newcommand\cQ{\mathcal{Q}}

\newcommand\st{\mathrm{st}}
\newcommand\dy{\mathrm{dy}}
\newcommand\pt{\mathrm{pt}}
\newcommand\cE{\mathcal{E}}

\newcommand\la{\langle}
\newcommand\ra{\rangle}
\newcommand\dd{\,\mathrm{d}}
\newcommand\ee{\mathrm{e}}

\newcommand\p{\partial}
\newcommand\supp{\mathrm{supp}}

\newcommand\kin{\mathrm{kin}}

\renewcommand\liminf{\mathop{\lim\,\inf}\limits}%
\renewcommand\limsup{\mathop{\lim\,\sup}\displaylimits}%
\newcommand\sgn{\mathrm{sgn}}%

\numberwithin{equation}{section} 
\numberwithin{thm}{section}

\providecommand{\norm}[1]{\left\|#1\right\|}
\providecommand{\abs}[1]{\left|#1\right|}

\providecommand{\one}{\mathbf{1}}
\providecommand{\cX}{\mathcal{X}}
\providecommand{\cY}{\mathcal{Y}}

\theoremstyle{plain} 
\newtheorem{theorem}{Theorem}[section]

\begin{document}
	\makeatletter
	\newcommand{\CWaliaslabel}[2]{%
		\@ifundefined{r@#2}{}{%
			\begingroup
			\edef\CWtemp{\endgroup
				\noexpand\protected@write\noexpand\@auxout{}%
				{\string\newlabel{#1}{\expandafter\unexpanded\expandafter{\csname r@#2\endcsname}}}}%
			\CWtemp
		}%
	}
	\makeatother

	\title[Boundary-Driven Nonequilibrium for the Coulomb Landau Equation]{Boundary-Driven Nonequilibrium for the Coulomb Landau Equation}

	
	\thanks{\textit{Authors' comment.} This arXiv preprint provides a detailed companion to the more concise journal-submission version of the same title. This version is not intended for separate journal publication.}

	\author[Chanwoo Kim]{Chanwoo Kim}
	\address[Chanwoo Kim]{Department of Mathematics, University of Wisconsin-Madison, Madison, WI, 53717, USA} 
	\email{chanwookim.math@gmail.com; chanwoo.kim@wisc.edu}
	
	\author[Jinsol Seo]{Jinsol Seo}
	\address[Jinsol Seo]{School of Mathematics, Korea Institute for Advanced Study, 85 Hoegiro, Dongdaemun-gu, Seoul 02455, Republic of Korea}
	\email{seo9401@kias.re.kr}
	
	

	\begin{abstract}
		We establish a nonlinear theory of the Coulomb Landau
		equation in smooth bounded domains with non-isothermal
		diffuse reflection. We construct a stationary state, prove its perturbative uniqueness
		and global nonlinear stability, and establish exponential
		convergence when the wall temperature approaches its
		stationary profile exponentially. We also establish stationary and global nonlinear
		well-posedness for prescribed inflow, together with
		asymptotic stability.
		The analysis combines a critical kinetic integrability
		gain with finite-exponent velocity-gradient estimates
		near grazing, allowing the nonlinear Coulomb estimates
		to close without $L_x^\infty$ control of velocity gradients.
	\end{abstract}

	\maketitle
	
	\setcounter{tocdepth}{1} 
	\tableofcontents

\section{Introduction}\label{new:sec:introduction}

A central problem in nonequilibrium statistical mechanics
is to understand the existence and asymptotic stability
of nonequilibrium steady states (NESS)
\cite{BertiniDeSoleGabrielliJonaLasinioLandim2015}.
Kinetic equations provide a concrete PDE framework for
studying such states, particularly through boundary
problems. In collisional gases, interactions with
physical boundaries can sustain NESS through the
interplay of transport and collisions. Unlike Maxwellian
equilibria, NESS are not characterized by vanishing
entropy production, and Boltzmann's $H$-theorem does
not directly provide a Lyapunov mechanism for their
asymptotic stability.

In the case of the Boltzmann equation for hard-sphere
molecules, Esposito, Guo, Kim, and Marra
\cite{EspositoGuoKimMarra} established the existence
and asymptotic stability of NESS driven by non-isothermal
boundary conditions. For the Coulomb Landau equation,
however, the corresponding problem has remained open.
The Coulomb collision operator exhibits degenerating
coercivity at large velocities, while the singular
geometry of grazing velocities creates substantial
difficulties for boundary regularity. Moreover, the
natural energy estimates are insufficient to close
the nonlinear Coulomb interaction in the presence
of these boundary singularities.

{In this paper, we resolve the non-isothermal
	diffuse-reflection problem for the Coulomb Landau
	equation in smooth bounded domains. We construct
	a unique stationary state and establish
	its global nonlinear stability
	and exponential relaxation.
	We also develop a theory for prescribed inflow,
	providing stationary and global nonlinear
	well-posedness and asymptotic stability.} A key analytic advance
is a critical kinetic integrability gain that,
combined with finite-exponent control of grazing
singularities, allows us to close the nonlinear
Coulomb estimates beyond the natural energy framework.

\subsection{The Coulomb Landau equation and boundary conditions}
Let $\mathcal D\subset\mathbb R^3$ be a bounded connected domain with outward
unit normal $n$.  We consider
\begin{equation}\label{new:eq:intro-landau}
	\partial_tF+v\cdot\nabla_xF=\mathcal N(F,F),
	\qquad
	(t,x,v)\in(0,\infty)\times\mathcal D\times\mathbb R^3,
\end{equation}
where
\begin{equation}\label{new:eq:intro-landau-operator}
	\mathcal N(G,F)(v)
	:=
	\partial_{v_i}\!\left[
	\int_{\mathbb R^3}\phi^{ij}(v-v')
	\bigl(
	G(v')\partial_{v_j}F(v)
	-F(v)\partial_{v'_j}G(v')
	\bigr)\,dv'
	\right].
\end{equation}
Here and below repeated indices are summed, and the Coulomb kernel is
\begin{equation}\label{new:eq:intro-coulomb-kernel}
	\phi^{ij}(z)
	=
	\partial_{z_i}\partial_{z_j}|z|
	=
	\frac{\delta_{ij}}{|z|}
	-\frac{z_i z_j}{|z|^3},
	\qquad z\neq0.
\end{equation}
The phase-space boundary $\partial\mathcal D\times\mathbb R^3$ can be decomposed into the outgoing, incoming, and grazing boundaries:
\[
\gamma^\pm
:=
\bigl\{(x,v)\in\partial\mathcal D\times\mathbb R^3:
\ \pm n(x)\cdot v>0\bigr\},
\qquad
\gamma^0
:=
\bigl\{(x,v)\in\partial\mathcal D\times\mathbb R^3:
\ n(x)\cdot v=0\bigr\}.
\]

We consider two boundary conditions. The first is prescribed inflow, in which the incoming distribution is specified independently of the outgoing trace:
\begin{equation}\label{new:eq:intro-inflow-law}
	F(t,x,v)=H(t,x,v)
	\qquad\text{on }\gamma^-,
\end{equation}
where $H$ is a given incoming distribution. 

The second is thermal diffuse reflection,
\begin{equation}\label{new:eq:intro-diffuse-law}
	F(t,x,v)=M_\Theta(t,x,v)J[F](t,x)
	\qquad\text{on }\gamma^-,
\end{equation}
where $M_\Theta$ is the flux-normalized wall Maxwellian
at the prescribed temperature $\Theta(t,x)>0$, and
$J[F]$ is the outgoing mass flux:
\begin{equation}\label{new:eq:intro-wall-maxwellian}
	M_\Theta(t,x,v)
	:=
	\frac{1}{2\pi\Theta(t,x)^2}
	\exp\!\left(-\frac{|v|^2}{2\Theta(t,x)}\right),
	\qquad
	J[F](t,x)
	:=
	\int_{n(x)\cdot v>0}
	F(t,x,v)(n(x)\cdot v)\,dv.
\end{equation}
The normalization
\[
\int_{n(x)\cdot v<0}|n(x)\cdot v|M_\Theta(t,x,v)\,dv=1
\]
ensures that the incoming and outgoing mass fluxes balance pointwise at the boundary. Thus, the diffuse reflection conserves the total mass while allowing energy exchange with the wall.

For the diffuse problem, we write
\begin{equation}\label{new:eq:intro-wall-temperature}
	\Theta(t,x)
	=\Theta_0\bigl(1+\delta\vartheta(t,x)\bigr),
	\qquad
	\Theta_{\st}(x)
	=\Theta_0\bigl(1+\delta\vartheta_{\st}(x)\bigr),
\end{equation}
where $\Theta_0>0$ is a reference temperature,
$\delta>0$ measures the amplitude of the temperature variation,
and $\vartheta(t,x)$ and $\vartheta_{\st}(x)$ describe
the time-dependent and stationary variations of the wall
temperature, respectively.
For asymptotic stability, we assume that the time-dependent
wall temperature converges exponentially to its stationary
profile.

Since diffuse reflection conserves the total mass,
we fix $\mathsf M>0$ and introduce the reference
global Maxwellian
\begin{equation}\label{new:eq:intro-reference-maxwellian}
	\mu(v)
	:=
	\frac{\rho_0}{(2\pi\Theta_0)^{3/2}}
	\exp\!\left(-\frac{|v|^2}{2\Theta_0}\right),
	\qquad
	\rho_0:=\frac{\mathsf M}{|\mathcal D|}.
\end{equation}
For the prescribed-inflow problem, $\rho_0>0$ instead
denotes a fixed reference density, without a prescribed
total-mass constraint.

\begin{spacing}{1.02}

\subsection{Main results and scope}\label{new:subsec:intro-main-results}

The precise hypotheses and norms are stated in
Section~\ref{sec:main-results}.

\begingroup
\renewcommand{\thetheorem}{A}
\begin{theorem}[Non-isothermal diffuse reflection]\label{thm}
	Let $\mathcal D\subset\mathbb R^3$ be a sufficiently smooth
	bounded connected domain, and let $\Theta_{\st}$ be as in
	\eqref{new:eq:intro-wall-temperature}, with $\delta>0$
	sufficiently small. Then the Coulomb Landau equation
	admits a unique nonnegative, mass-normalized stationary state
	$F_{\st}\geq 0$ close to $\mu$, satisfying
	the thermal diffuse reflection boundary condition
	\[
	F_{\st}=M_{\Theta_{\st}}J[F_{\st}]
	\qquad\text{on }\gamma^-.
	\]
	
	Suppose further that the time-dependent wall temperature $\Theta(t,x)$
	converges exponentially to $\Theta_{\st}(x)$ as $t\to\infty$.
	Then every compatible nonnegative initial distribution
		sufficiently close to $F_{\st}$, with the same total
	mass, generates a unique global nonnegative solution
	$F(t,x,v)\geq 0$ satisfying
	\[
	F=M_\Theta J[F]
	\qquad\text{on }\gamma^-.
	\]
	The solution $F(t,x,v)$ 
	converges exponentially to $F_{\st}(x,v)$ as
	$t\to\infty$.
\end{theorem}
\endgroup

\begingroup
\renewcommand{\thetheorem}{B}
\begin{theorem}[Prescribed inflow]\label{thm:intro-B}
	Let $\mathcal D\subset\mathbb R^3$ be a bounded
	$C^{1,1}$ domain. Suppose that the prescribed stationary
	inflow $H_{\st}\geq0 $ is a sufficiently small admissible
	perturbation of a global Maxwellian $\mu$ on $\gamma^-$.
	Then the Coulomb Landau equation admits a unique stationary
	solution $F_{\st} \geq 0 $ close to $\mu$, satisfying
	\[
	F_{\st}=H_{\st}
	\qquad\text{on }\gamma^-.
	\] 
	
	For compatible initial data $F_0\geq 0$
		and time-dependent inflow $H(t)\geq 0$, both sufficiently close to $\mu$,
	there exists a unique global nonnegative solution
	$F(t,x,v)\geq 0$ satisfying
	\[
	F=H
	\qquad\text{on }\gamma^-.
	\] 
	
	Both the stationary and
	evolutionary solutions enjoy Lipschitz stability
	in the physical Coulomb graph norm. Moreover, for
	every finite $p>1$, we obtain
	\begin{equation}\label{new:eq:intro-all-finite-p}
		\nabla_v
		\left(\frac{F_{\st}-\mu}{\sqrt{\mu}}\right)
		\in L_x^pL_v^2,
		\qquad
		\nabla_v
		\left(\frac{F-\mu}{\sqrt{\mu}}\right)
		\in L_t^\infty L_x^pL_v^2,
	\end{equation}
	with the velocity weights specified in
	Section~\ref{sec:main-results}.
	
	If, in addition, the time-dependent inflow $H(t)$
	converges exponentially to $H_{\st}$ in the prescribed
	boundary norm, then the global solution $F(t,x,v)$
	converges exponentially to the corresponding
	stationary state $F_{\st}(x,v)$ as $t\to\infty$.
\end{theorem}
\endgroup

\noindent\textbf{Remarks.} Theorems A and B establish stationary states near $\mu$, their uniqueness in the stated classes, and nonlinear asymptotic stability. A nonconstant wall temperature in Theorem A produces a stationary state that is not a local Maxwellian. In Theorem B, the prescribed stationary inflow need not have Maxwellian form, and the resulting stationary state is, in general, not a local Maxwellian. Both theorems establish exponential relaxation toward their stationary states through the interplay of transport and boundary dissipation, in contrast to the periodic setting, where the Coulomb Landau equation lacks a spectral gap in the natural energy space. Finally, the nonnegativity of the stationary states in both Theorem A and Theorem B is established directly, independently of the asymptotic stability of the evolution.

\subsection{One-sided trace theory and stationary nonnegativity}
\label{new:subsec:intro-physical-graph}

For prescribed inflow, the incoming trace is given,
	whereas the outgoing trace and the kinetic Green identity
	must be recovered in an appropriate graph space.
	The bulk kinetic graph alone need not determine the full
	physical boundary trace on a curved bounded domain
\cite{NiVa}.

We encode this asymmetry directly in the topology.
Let $\mathcal H_{\kin}$ be the completion of
$C_c^\infty([0,T]\times\overline{\mathcal D}
\times\mathbb R_v^3)$ under
\begin{equation}\label{new:eq:one-sided-graph}
	\begin{aligned}
		\|f\|_{\mathcal H_{\kin}}^2
		:=
		\|f\|_{L^2_{t,x}H_v^1}^2
		+\|(\partial_t+v\cdot\nabla_x)f\|_{L^2_{t,x}H_v^{-1}}^2 +\|f(0)\|_{L^2_{x,v}}^2
		+\|f_{\gamma_T^-}\|_{L^2(\gamma_T^-)}^2.
	\end{aligned}
\end{equation}
Thus only the initial and incoming traces enter
the topology. The stationary space is defined
analogously, without the time trace.
The corresponding Coulomb graph
$\mathcal H_{\kin,\sigma}$ is obtained by
replacing $H_v^{\pm1}$ with $H_\sigma^{\pm1}$.

This completion admits an intrinsic characterization.
Every element admits a canonical
$C([0,T];L^2_{x,v})$ representative, a full
physical boundary trace, and the kinetic Green
identity. Conversely, every intrinsic weak
transport state satisfying the corresponding
bulk regularity, physical trace conditions,
and Green identity belongs to the completion.
Consequently,
\begin{equation}\label{new:eq:one-sided-identification}
	\text{one-sided physical completion}
	=
	\text{intrinsic full physical graph}.
\end{equation}
The converse inclusion follows from the surjectivity
	of the kinetic Fokker--Planck problem with prescribed
	inflow and time reversal, which together eliminate
	the orthogonal complement of the completion.
The underlying solvability theory is built
using stopped Kolmogorov flows and
Feynman--Kac representations.

This physical Green structure also allows us
to establish nonnegativity of the stationary
state directly, without using asymptotic
stability. Such positivity does not follow
merely from perturbative smallness around
a Maxwellian, particularly at large velocities.
A renormalized $L^1$ negative-part identity,
combined with the diffuse boundary law and
perturbative positivity on a bounded velocity
set, first eliminates negative boundary traces.
A further strictly dissipative renormalization
then eliminates the negative part in the interior.
Nonnegativity of the evolution is proved
separately through a time-dependent
negative-part inequality.

The recovered Green identity provides the favorable
	boundary sign for convex renormalization and yields
	the critical kinetic integrability gain.
In the diffuse problem, the weighted graph
characterization also recovers the physical
traces of complete tangential derivatives,
justifying the differentiated boundary laws.

\subsection{Critical kinetic integrability gain and nonlinear Coulomb closure}
\label{new:subsec:intro-coulomb-closure}

A central difficulty in closing the nonlinear Coulomb
estimate is a mismatch between the spatial integrability
of its two factors. Writing
$N[g,f]:=\mu^{-1/2}\mathcal N(g\sqrt{\mu},f\sqrt{\mu})$
for the perturbative bilinear Landau form around $\mu$,
we distinguish the rough factor $g$, controlled without
velocity derivatives, from the differentiated factor $f$.
Suppressing fixed velocity weights, the corresponding
mixed estimate takes the schematic form
\begin{equation}\label{new:eq:intro-coulomb-mixed}
	\|N[g,f]\|_{L_x^2H_\sigma^{-1}}
	\lesssim
	\|\mu^{1/4}g\|_{L_x^qL_v^2}
	\|f\|_{L_x^pH_\sigma^1},
	\qquad
	1/q+1/p=1/2.
\end{equation}
At the natural energy exponent $q=2$, this forces
$p=\infty$, precisely the endpoint excluded
	in general by grazing singularities.
Thus, the natural Coulomb energy alone does not close
the nonlinear difference estimate.

Convex renormalization and the recovered Green identity
	yield a critical integrability gain for the rough factor.
More precisely, we extend nonnegative
convex renormalizations of a kinetic graph function
by zero to the whole phase space. The resulting
distributional inequality retains the initial and
incoming boundary contributions as a nonnegative
measure controlled by the corresponding trace data,
while the terminal and outgoing contributions have
the favorable sign and can be discarded. Whole-space Kolmogorov potential estimates then yield
the following critical weak-type bounds for $h$
	in the respective kinetic graph spaces:
\begin{equation}\label{new:eq:intro-critical-phase-space-gain}
	\begin{aligned}
		\|h\|_{L^{14/6,\infty}_{t,x,v}}
		\lesssim \|h\|_{\mathcal H_{\kin}} \ \ 
		\text{(evolution)},\ \ \ \  \  \  \ \  \  \  
		\|h\|_{L^{12/5,\infty}_{x,v}}
		\lesssim \|h\|_{\mathcal H_{\kin}} \ \ 
		\text{(stationary)}.
	\end{aligned}
\end{equation}
The distinct critical exponents reflect the kinetic
homogeneous dimensions $14$ and $12$, respectively.
Importantly, this argument requires neither a reflection
reduction nor division by $n(x)\cdot v$, and hence
remains valid without excluding grazing velocities.

After kinetic localization and the introduction of
Maxwellian weights, these estimates improve the spatial
integrability of the rough factor beyond $L_x^2$:
one obtains $q<7/3$ in evolution and $q<12/5$
stationarily, with $q>2$. Through
$1/q+1/p=1/2$, the forbidden endpoint $p=\infty$
is thereby replaced by the finite requirements
$p>14$ in evolution and $p>12$ stationarily.
The nonlinear Coulomb estimate can therefore be
closed using finite-exponent regularity of the
differentiated factor, rather than an
$L_x^\infty$ velocity-gradient bound.

\subsection{Grazing regularity and finite-$p$ gradient estimates}
\label{new:subsec:intro-grazing}

The integrability gain of the preceding subsection reduces
nonlinear Coulomb closure to finite-$p$ velocity-gradient
estimates near grazing.

We first establish boundary $L^2$--$L^\infty$ and H\"older
estimates for kinetic Fokker--Planck equations. The
renormalized Green identity provides the favorable boundary
sign in Caccioppoli inequalities, allowing kinetic Sobolev
estimates and Moser iteration to yield local boundedness
\cite{PascucciPolidoro2004,2022AR}.
For H\"older regularity, we adapt the kinetic De Giorgi
method \cite{2019GIMV,SilvestreBoundaryKFP} to the same
weak-solution class. Away from grazing, incoming patches
with homogeneous data admit zero extension across the
boundary. At outgoing patches, where the trace is not
prescribed, we instead solve an auxiliary kinetic problem
on the opposite side, using the recovered outgoing trace
as its incoming datum, and glue the two solutions through
their Green identities.

At grazing, neither extension applies uniformly. The positive
part of the zero-extended solution satisfies a whole-space
kinetic inequality with a nonnegative defect measure.
Hypoelliptic compactness and a boundary De Giorgi
isoperimetric argument then yield oscillation decay and
H\"older regularity.

The quantitative grazing decay needed for velocity gradients
comes from a barrier comparison. The flat normal Kolmogorov
operator $\mathcal K_0:=\zeta_3\partial_{y_3}-\partial_{\zeta_3\zeta_3}$
on $\{y_3>0\}$ has a classical positive harmonic profile $H$,
vanishing on the incoming boundary and homogeneous of order
$1/2$ under kinetic scaling, which underlies the sharp grazing
regularity \cite{KimWeidner2026}.
At the critical homogeneity, $\mathcal K_0H=0$ leaves
no margin for variable-coefficient errors.
We therefore use the subcritical power $H^\beta$ with
$\beta=1-2\varepsilon\in(0,1)$, whose defect
$\mathcal K_0(H^\beta)$ controls both $|(H^\beta)_{\zeta_3\zeta_3}|$
and $d_{\rm graz}^{-1}|(H^\beta)_{\zeta_3}|$.
The resulting margin absorbs the coefficient errors, yielding
almost-critical $d_{\rm graz}^{1/2-\varepsilon}$ decay for
homogeneous incoming data, with controlled source terms.

Transferring this barrier to the Coulomb Landau operator
requires an anisotropic normalization. Its radial and
tangential diffusion eigenvalues satisfy
$\lambda_\parallel(v)\simeq\langle v\rangle^{-3}$
and $\lambda_\perp(v)\simeq\langle v\rangle^{-1}$,
so scalar rescaling is insufficient. We combine
transport-preserving boundary flattening with a matrix
normalization that preserves the sign of the normal velocity.
Kinetic H\"older estimates and Landau convolution bounds
provide the coefficient control needed to transfer the
barrier to the variable-coefficient operator.

We next convert grazing decay into a velocity-gradient
estimate. For regular coefficients, divergence-form kinetic
$L^p$ estimates first control $\nabla_v f$ in $L^p$.
Nondivergence-form estimates and kinetic Sobolev--Morrey
embedding then give an interior pointwise gradient bound
\cite{dong2024global,dong2022global,polidoro1998sobolev},
with leading scale $r^{-1}\|f\|_{L^\infty}$ on a kinetic
cylinder of radius $r$, together with controlled coefficient
and source terms. At kinetic distance $d_{\rm graz}$ from
grazing, we choose a non-grazing cylinder of comparable
radius. Combining the barrier decay with the interior
gradient estimate gives the key scale conversion
\[
|f|\lesssim d_{\rm graz}^{1/2-\varepsilon}
\quad\Longrightarrow\quad
|\nabla_v f|
\lesssim
d_{\rm graz}^{-1}
d_{\rm graz}^{1/2-\varepsilon}
=
d_{\rm graz}^{-1/2-\varepsilon},
\]
up to controlled source terms and velocity weights.
Anisotropic normalization and localization then yield
the corresponding weighted pointwise gradient estimate.

The grazing singularity persists, but its precise order
makes it integrable in the mixed norms required by the
Coulomb estimate. In boundary-flattening coordinates,
$d_{\rm graz}\asymp|\widehat x_3|^{1/3}
+|\widehat v_3|$. At the critical grazing exponent,
integration in the normal velocity produces only
logarithmic growth. The arbitrarily small power loss
remains integrable in $L_x^pL_v^2$ for every finite
$p$, by choosing $\varepsilon>0$ sufficiently small
depending on $p$. The velocity weights in the Landau
estimates control the large-velocity region.

Nonhomogeneous incoming data are handled through weighted
inflow liftings. We thus obtain, under the corresponding
weighted assumptions, $\nabla_v f\in L_t^\infty L_x^pL_v^2$
for every finite $p>1$, with the analogous stationary
estimate. These bounds cover the exponents required by
the Coulomb bilinear estimate, namely $p>14$ in evolution
and $p>12$ stationarily. Combined with the rough-factor
integrability gain of the preceding subsection, they
close the nonlinear Coulomb difference estimate without
an $L_x^\infty$ velocity-gradient bound.

\subsection{Self-consistent diffuse feedback and centered dynamics}
\label{new:subsec:intro-diffuse-feedback}

Applying the prescribed-inflow theory to diffuse reflection
requires spatial Lipschitz regularity of the self-generated
outgoing mass flux $J[F_{\st}]$.
Such regularity does not follow from the basic energy estimates
and must be established to close the self-consistent
boundary feedback.

To differentiate the diffuse boundary condition while
preserving the grazing geometry, we introduce complete
phase-space tangential derivatives.
For a vector field $X$ tangent to the level sets of
the signed distance, define
\begin{equation}\label{new:eq:intro-complete-derivative}
	\mathfrak D_X
	:=X(x)\cdot\nabla_x
	+\bigl((\nabla_xX(x))v\bigr)\cdot\nabla_v,
	\qquad
	\mathfrak D_X(n\cdot v)=0
	\quad\text{on }\partial\mathcal D\times\mathbb R^3.
\end{equation}
The associated flow preserves the incoming, outgoing,
and grazing sectors.
The diffuse boundary operator has the rank-one
velocity-moment structure familiar from
\cite{Guo2010,EspositoGuoKimMarra}.
Under complete differentiation, the principal diffuse
operator remains unchanged, while the boundary
commutators involve only lower complete derivatives.
These derivatives are constructed through finite complete
flows in the physical Coulomb graph.
The graph characterization of
Subsection~\ref{new:subsec:intro-physical-graph}
then recovers their physical traces and Green identity,
justifying the differentiated boundary conditions
without formally differentiating rough traces across
grazing.

Combining one-dimensional Landau--Kolmogorov interpolation
along the complete flows with the differentiated boundary
hierarchy and the energy and trace estimates in the
	physical graph, we control the high-order outgoing moments.
Boundary spectral regularization closes the hierarchy
without presupposing the desired flux regularity.
This yields
\begin{equation}\label{new:eq:intro-H3-outgoing}
	\|J[F_{\st}]-J[\mu]\|_{H^{\mathrm N}(\partial\mathcal D)}
	\lesssim \delta,
	\qquad \mathrm N\ge4,
\end{equation}
under the corresponding smoothness assumptions.
Sobolev embedding provides the spatial Lipschitz
regularity needed for the prescribed-inflow estimates.
The stationary state is then constructed by a fixed-point
argument in a strong invariant class, with contraction
in the basic physical graph.

For the nonlinear evolution, we center the solution
at the stationary state corresponding to the limiting
wall temperature.
The mass constraint eliminates the neutral mode,
while diffuse boundary dissipation, transport-tail
estimates, and macroscopic recovery, drawing on the
boundary rigidity framework of
\cite{DesvillettesVillani2002,EspositoGuoKimMarra2018},
yield energy estimates uniform over arbitrary time horizons.
These estimates accommodate the forcing generated
by a wall temperature converging exponentially
to its stationary profile.
A delayed-dissipation argument then gives exponential
convergence to the stationary state.
The complete-derivative hierarchy is propagated
along the evolution, completing the nonlinear
stability theory of Theorem~A.

\subsection{Relation to previous work}
\label{new:subsec:intro-previous-work}

Several approaches to physical boundary traces
	for kinetic equations have been developed, including
the variational theory of
Albritton--Armstrong--Mourrat--Novack, the operator-theoretic
approach of Nier, and the Dirichlet theory of
Litsg{\aa}rd--Nystr\"om
\cite{AlbrittonArmstrongMourratNovack,Nier2018,
	LitsgardNystrom2021}.
More recently, Valentini established stationary Kolmogorov
well-posedness with transport-defined traces and one-sided
flux estimates \cite{Valentini2026}.
Niebel and Valentini established sharp kinetic trace
estimates and counterexamples, including the failure
of the natural unrestricted Gaussian trace estimate
on smooth bounded domains \cite{NiVa}.
In the weighted Coulomb setting considered here,
	the smooth completion defined by the incoming trace is
	identified with an intrinsic graph admitting both physical
	boundary traces and the Green identity. This identification
	also underlies the stationary nonnegativity argument, the
	critical integrability gain, and the differentiated diffuse
	boundary conditions.

 Boundary regularity depends sharply on the boundary law.
	Grazing singularities already occur for transport equations,
	where regularity may depend on the geometry of the domain
	\cite{KimDiscontinuity,GuoKimTononTrescases,GuoKimTononTrescasesBV,KimLeeHolder,ChenKim2022,ChenKim2024Gradient,ChenKimLi2020LVPB,CaoKimLee2019GVPB}.
	For absorbing (zero-inflow) boundary conditions,
	Hwang, J. Jang, and Vel\'azquez established a pioneering boundary
	regularity result for kinetic Fokker--Planck equations, proving
	H\"older continuity up to the grazing set in an interval
	\cite{HwangJangVelazquez2014}.
	 Silvestre established H\"older estimates up to the boundary for
		kinetic Fokker--Planck equations with rough coefficients and prescribed
		inflow, together with higher-order asymptotic estimates near the
		incoming boundary \cite{SilvestreBoundaryKFP}.
			Further boundary regularity results include
			\cite{HendersonLucertiniWang2024} (see also \cite{ZhuBoundaryKFP}). 
	For prescribed inflow and diffuse reflection,
	K. Kim and Weidner established the sharp $C^{1/2}$
	grazing regularity for the corresponding linear kinetic
	Fokker--Planck equations \cite{KimWeidner2026}.
	Boundary regularity under Maxwell accommodation, which
	interpolates between diffuse and specular reflection,
	was also studied in \cite{KimWeidnerMaxwell2026}.
	
	Pure specular reflection has a different boundary structure. It couples incoming and outgoing velocities
		by reflection, allowing reflection-based extensions across the boundary \cite{DongGuoYastrzhembskiy2022}.
	Unlike inflow, it does not prescribe an independent incoming trace: the two boundary branches are linked by symmetry. This permits higher regularity at grazing, although the grazing geometry still limits smoothness.
	Ros-Oton and Weidner established optimal $C_{\rm kin}^{4,1}$
	regularity up to grazing for linear kinetic Fokker--Planck
	equations with specular reflection \cite{RosOtonWeidner2025}.

For the nonlinear Landau equation, J. Kim, Guo, and Hwang
developed an $L^2$--$L^\infty$ framework on the torus,
combining De Giorgi regularity and velocity estimates
to establish uniqueness \cite{2020KGH}.
For boundary problems, a global mild-solution theory with specular reflection in a finite channel was established in
\cite{DuanLiuSakamotoStrain}.\footnote{
	In the proof of the prescribed-inflow result
		(Theorem~2.3), formula~(7.4) replaces the boundary values
		of the nonlocal collision operators by expressions involving
		only the prescribed incoming data $g_\pm$. This substitution
	is not justified by the inflow condition, since the
	collision operators are nonlocal in velocity and
	depend on the unknown outgoing distribution.
	The boundary estimate~(7.8) explicitly relies on
	(7.4), as do the subsequent estimates (7.9)--(7.10).
}
Guo, Hwang, J.W. Jang, and Ouyang subsequently established
global nonlinear stability for the Coulomb Landau
equation with specular reflection in a general smooth
bounded domain
\cite{GuoHwangJangOuyang,GuoHwangJangOuyangCorrection}.
The reflection symmetry underlying these theories
is unavailable for prescribed inflow.

Carrapatoso and Mischler developed a near-equilibrium
theory under general Maxwell reflection, stating global
existence and constructive asymptotic decay for the
full range of interaction potentials, including Coulomb
\cite{CarrapatosoMischler}.\footnote{
	For the Coulomb case, the coefficient bound used in
	Step~2 of Proposition~5.5 does not hold.
	Writing $\widetilde a=a*(\mu+g)$, equation~(5.6)
	introduces
	$S_1=(\widetilde a-I)\nabla_v\bar h$.
	The subsequent bound for $S_1$ treats
	$\widetilde a-I$ as decaying like
	$\langle v\rangle^{\gamma+2}$.
	For $\gamma=-3$, however, this fails even for $g=0$:
	$(a*\mu)(v)\to0$ as $|v|\to\infty$, whereas
	$(a*\mu)(v)-I\to-I$.
	Thus the displayed bound for $S_1$ does not follow
	when $\gamma=-3$, and the deduction of (5.7) in
	that case requires another argument.
}
Their nonlinear existence theory does not establish
uniqueness in the corresponding solution class.
For prescribed inflow, we establish uniqueness near
	a Maxwellian in the physical Coulomb graph. Under
	non-isothermal diffuse reflection, we construct the
	mass-normalized stationary state and prove its uniqueness
	and nonlinear asymptotic stability.

\end{spacing}

\newpage
\section{Main Results}\label{sec:main-results}

We state the results for the Coulomb Landau equation
\eqref{new:eq:intro-landau} with prescribed inflow
\eqref{new:eq:intro-inflow-law} or diffuse reflection
\eqref{new:eq:intro-diffuse-law}. The boundary notation and the
reference Maxwellian are those of Section~\ref{new:sec:introduction},
with arbitrary fixed $\rho_0,\Theta_0>0$.
In the diffuse problem,
$\mathsf M=\rho_0|\cD|$ is the prescribed mass; for inflow,
$\rho_0$ is a reference density.
We present the diffuse-reflection results first,
	followed by the prescribed-inflow results used in their proof.

\subsection{Notation for the main results}\label{subsec:main-common-framework}

Write
\begin{align}
	F_{\st}&=\mu+\sqrt\mu\,f_{\st},
	\label{eq:stationary-perturbation-decomposition-main}\\
	F&=F_{\st}+\sqrt\mu\,f_{\dy}
	=\mu+\sqrt\mu\,f_{\pt},
	\qquad f_{\pt}=f_{\st}+f_{\dy}.
	\label{eq:dynamic-perturbation-decomposition-main}
\end{align}
Set $f_0=f_{\pt}(0)$.

Set $w(v)=(1+|v|^2/\Theta_0)^{1/2}$,
$\cD_T=(0,T)\times\cD$, and $\gamma_T^\pm=(0,T)\times\gamma^\pm$.
Boundary $L^2$ norms are taken with respect to
\begin{equation}\label{def:dgamma-main}
	\dd\gamma=|n(x)\cdot v|\,\dd S_x\dd v.
\end{equation}
The Coulomb velocity spaces $H_\sigma^{\pm1}$ and the stationary and
evolutionary graph spaces $\cH_{\kin,\sigma}$ are defined in
Definitions~\ref{def:H_sigma} and~\ref{def:H_kin_sigma}.
The kinetic H\"older norm $C^\alpha_{\kin}$ uses the distance
$\mathbf d$ of \eqref{kineticdistance}; for stationary functions we use
\begin{equation}\label{eq:stationary-kinetic-distance-main}
	\mathbf d_{\st}((x,v),(y,v'))
	:=\inf_{\tau\in\bR}\mathbf d((\tau,x,v),(0,y,v')).
\end{equation}
In both settings, $\|g\|_{C^\alpha_{\kin}}$ denotes the sum of the
supremum norm and the corresponding $\alpha$-H\"older seminorm.
We write
$L^{p,2}_{x,v}=L_x^pL_v^2$ and
$L^{\infty,p,2}_{t,x,v}=L_t^\infty L_x^pL_v^2$.
\subsection{Non-isothermal diffuse reflection}\label{subsec:main-diffuse}

The precise initial-data regularity, smallness and compatibility assumptions,
		together with the full weighted estimates, are given in
		Theorems~\ref{thm:stationary-diffuse} and~\ref{thm:dynamic-global}.

Fix integers $\mathrm N\ge4$ and $K_\cD\ge2\mathrm N$, and assume
that $\cD\subset\bR^3$ is a bounded connected $C^{K_\cD+2}$ domain.
We impose the diffuse reflection law \eqref{new:eq:intro-diffuse-law},
with $J$ and $M_\Theta$ from \eqref{new:eq:intro-wall-maxwellian}.
The stationary wall temperature $\Theta_\st$ is given by
\eqref{new:eq:intro-wall-temperature}, where
$\vartheta_\st\in C^{K_\cD}(\partial\cD)$ and
$\|\vartheta_\st\|_{C^{K_\cD}}\le1$.
\begin{theorem}[Stationary diffuse reflection]\label{thm:main}
	There are constants $\delta_\st,\varepsilon_\st>0$ and
	$\alpha_\st\in(0,1)$ such that, for $0<\delta\le\delta_\st$,
	the Coulomb Landau equation admits a stationary solution
	$F_\st$, decomposed as in
	\eqref{eq:stationary-perturbation-decomposition-main}, satisfying
	\eqref{new:eq:intro-diffuse-law} at the temperature $\Theta_\st$, with
	\[
	F_\st\ge0,\qquad
	\int_{\cD\times\bR^3}F_\st\,\dd x\dd v=\mathsf M.
	\]
	It is unique among mass-normalized stationary solutions
	$\mu+\sqrt\mu\,f$ with
	$f\in\cH_{\kin,\sigma}(\cD\times\bR^3)$ and
	$\|\mu^{1/4}f\|_{L^\infty}\le\varepsilon_\st$.
	
	For every $\ell\ge0$ and $1<p<\infty$, the stationary perturbation satisfies
	\begin{equation}\label{eq:main-diffuse-stationary-readable-bound}
		\begin{aligned}
			&\|w^\ell f_\st\|_{\cH_{\kin,\sigma}}
			+\|w^\ell f_\st\|_{L^\infty}
			+\|w^\ell f_\st\|_{C^{\alpha_\st}_{\kin}}\\
			&\qquad+\|w^\ell f_\st\|_{L_x^pH_\sigma^1}
			+\|w^\ell\nabla_v f_\st\|_{L_x^pL_v^2}
			\le C_{\ell,p}\delta.
		\end{aligned}
	\end{equation}
	The incoming trace satisfies weighted boundary estimates. The outgoing flux satisfies
	\[
	\|J[F_\st]-J[\mu]\|_{H^{\mathrm N}(\partial\cD)}\le C\delta.
	\]
	The stationary perturbation also satisfies weighted graph estimates for complete tangential derivatives through order $\mathrm N$.  In particular, the flux difference belongs to
	$C^{1,\beta}(\partial\cD)$ for every $0<\beta<1$.
	If $\vartheta_\st$ is nonconstant, $F_\st$ cannot be a global Maxwellian.
\end{theorem}

For the evolutionary problem, assume in addition that
$\vartheta\in C^{\mathrm N+1}([0,\infty);C^{K_\cD}(\partial\cD))$
and, for some $\lambda_w>0$,
\begin{equation}\label{eq:main-diffuse-readable-wall-assumptions}
	\begin{aligned}
		\max_{0\le j\le\mathrm N+1}\sup_{t\ge0}
		\|\partial_t^j\vartheta(t)\|_{C^{K_\cD}}&\le1,\\
		\max_{0\le j\le\mathrm N+1}
		\|\partial_t^j(\vartheta(t)-\vartheta_\st)\|_{C^{K_\cD}}
		&\le e^{-\lambda_wt}\qquad(t\ge0).
	\end{aligned}
\end{equation}
In particular, $\Theta(t,\cdot)$ converges exponentially to
$\Theta_\st$.

\begin{theorem}[Global diffuse-reflection dynamics]
	\label{thm:dynamic-global-centered-solution}
	Assume
		\eqref{eq:main-diffuse-readable-wall-assumptions} and let $F_\st$
		be the stationary state above. For each finite target weight
		$\ell\ge0$, there are constants
		$\varepsilon_\dy,C,\lambda_{\rm exp}>0$, depending on
		$\ell$ and the fixed structural parameters, including the wall
		relaxation rate $\lambda_w$, such that the following holds.
	
	Let $F_0=\mu+\sqrt\mu\,f_0$ have mass $\mathsf M$.
		Assume that the weighted data size $\mathsf I$ of
		\eqref{eq:dyn-data-size} satisfies $\delta+\mathsf I\le\delta_\star$,
		where $\delta_\star$ is the threshold in
		Theorem~\ref{thm:dynamic-global}, and that the generated
		initial time jets satisfy
		\eqref{eq:dynamic-initial-u-compatibility} and the zero-mass
		condition \eqref{eq:dynamic-initial-zero-perturbation-mass}. 
	
	Then \eqref{new:eq:intro-landau}, subject to the time-dependent
	diffuse law \eqref{new:eq:intro-diffuse-law}, has a global solution
	$F=\mu+\sqrt\mu\,f_\pt$ of mass $\mathsf M$.
	It is unique in the class
	$f_\pt\in\bigcap_{T>0}\cH_{\kin,\sigma}(\cD_T\times\bR^3)$,
	$\|\mu^{1/4}f_\pt\|_{L^\infty}\le\varepsilon_\dy$.
	If $F_0\ge0$, then $F(t)\ge0$ for all $t\ge0$.
	
	The solution converges exponentially to $F_\st$.  With
	$\lambda=\min\{\lambda_{\rm exp},\lambda_w\}>0$,
	\begin{equation}\label{eq:main-diffuse-readable-decay}
		\begin{aligned}
			&\|f_\dy(t)\|_{L^2_{x,v}}^2
			+\int_t^\infty\|f_\dy(s)\|_{L_x^2H_\sigma^1}^2\,\dd s
			+\|(f_\dy)_{\gamma^+}\|_{L^2((t,\infty)\times\gamma^+)}^2\\
			&\qquad\le C\bigl(\|f_0-f_\st\|_{L^2_{x,v}}^2+\delta^2\bigr)
			e^{-\lambda t},\qquad t\ge0.
		\end{aligned}
	\end{equation}
	For a time-independent wall, the term $\delta^2$ may be omitted
	and one may take $\lambda=\lambda_{\rm exp}$.
	
	At the chosen weight, $f_\pt$ satisfies global weighted
		$L^\infty$ and kinetic H\"older bounds, as well as
		$\nabla_v f_\pt\in L_t^\infty L_x^pL_v^2$ for some $p>3$.
		The perturbation $f_\dy$ satisfies weighted graph estimates for
		complete space-time derivatives through order $\mathrm N-1$.
\end{theorem}

\noindent\textit{Initial-data regularity.}
	The regularity hypotheses on $f_0$ permit the construction
	of the initial time jets and the negative-time lift in Section~13.

\subsection{Prescribed inflow}\label{subsec:main-inflow}

For inflow data we write $H_{\st}=\mu+\sqrt\mu\,h_{\st}$ and
$H=\mu+\sqrt\mu\,h_{\pt}$.\par

For the inflow results we use separate Lipschitz control in space and
velocity.  For $A\subseteq\overline\cD$, put
\[
\begin{aligned}
	\|g\|_{(C_x^{0,1}\cap C_v^{1,1})(A\times\bR^3)}
	:={}&\sup|g|
	+\sup_v\sup_{a\ne b}\frac{|g(a,v)-g(b,v)|}{|a-b|}\\
	&+\sup|\nabla_vg|
	+\sup_a\sup_{v\ne z}\frac{|\nabla_vg(a,v)-\nabla_vg(a,z)|}{|v-z|}.
\end{aligned}
\]
At the incoming boundary define the velocity and fixed-velocity spatial
fibers by
$\gamma^-_{(x,*)}=\{z:n(x)\cdot z<0\}$ and
$\gamma^-_{(*,v)}=\{y\in\partial\cD:n(y)\cdot v<0\}$.
The corresponding boundary norm is
\begin{equation}\label{eq:stationary-inflow-C11-norm-main}
	\begin{aligned}
		\|h\|_{(C_x^{0,1}\cap C_v^{1,1})(\gamma^-)}
		:={}&\sup_{\gamma^-}|h|
		+\sup_v\sup_{\substack{x,y\in\gamma^-_{(*,v)}\\x\ne y}}
		\frac{|h(x,v)-h(y,v)|}{|x-y|}\\
		&+\sup_{\gamma^-}|\nabla_vh|
		+\sup_x\sup_{\substack{v,z\in\gamma^-_{(x,*)}\\v\ne z}}
		\frac{|\nabla_vh(x,v)-\nabla_vh(x,z)|}{|v-z|}.
	\end{aligned}
\end{equation}
Empty fiberwise suprema are zero, and velocity derivatives of incoming
data are understood by continuous extension to the closed fibers.

For prescribed inflow \eqref{new:eq:intro-inflow-law}, let
$\cD\subset\bR^3$ be a bounded $C^{1,1}$ domain.
The stationary data are measured by
\eqref{eq:stationary-inflow-C11-norm-main}.  For $0<T\le\infty$, set
\[
[h]_{t,x;\gamma_T^-}
:=\sup_v\sup_{\substack{(t,x)\ne(s,y)\in(0,T)\times\gamma^-_{(*,v)}}}
\frac{|h(t,x,v)-h(s,y,v)|}{|t-s|+|x-y|},
\]
\[
[\nabla_vh]_{v;\gamma_T^-}
:=\sup_{(t,x)\in(0,T)\times\partial\cD}\sup_{\substack{v\ne z\in\gamma^-_{(x,*)}}}
\frac{|\nabla_vh(t,x,v)-\nabla_vh(t,x,z)|}{|v-z|},
\]
and define
\begin{equation}\label{eq:dynamic-inflow-C11-norm-main}
	\|h\|_{(C_{t,x}^{0,1}\cap C_v^{1,1})(\gamma_T^-)}
	:=\sup_{\gamma_T^-}|h|+[h]_{t,x;\gamma_T^-}
	+\sup_{\gamma_T^-}|\nabla_vh|+[\nabla_vh]_{v;\gamma_T^-}.
\end{equation}
\begin{theorem}[Stationary prescribed-inflow problem]\label{Landau_st_exi}
	There exist $\varepsilon_0>0$, $\alpha_0\in(0,1)$, an integer
	$A_{\st}\ge0$, and $\delta_{\st}>0$, depending on $\cD$ and the fixed
	reference parameters, such that the following holds.  If
	\begin{equation}\label{eq:main-inflow-st-smallness}
		\|w^{A_{\st}}h_{\st}\|_{(C_x^{0,1}\cap C_v^{1,1})(\gamma^-)}
		\le\delta_{\st},
	\end{equation}
	then the stationary Coulomb Landau equation
	\eqref{new:eq:intro-landau}, with incoming datum $H_{\st}$ as in
	\eqref{new:eq:intro-inflow-law}, has a solution
	$F_{\st}=\mu+\sqrt\mu\,f_{\st}$ with
	$f_{\st}\in\cH_{\kin,\sigma}(\cD\times\bR^3)$ satisfying
	\begin{equation}\label{eq:main-inflow-st-Linf}
		\|f_{\st}\|_{L^\infty}\le\varepsilon_0.
	\end{equation}
	If $H_{\st}\ge0$ almost everywhere on $\gamma^-$, then
	$F_{\st}\ge0$ almost everywhere in $\cD\times\bR^3$.
	
	For every $\ell\ge0$, if the boundary norm below is finite, then
	\begin{align}\label{eq:main-inflow-st-regularity}
		&\|w^\ell f_{\st}\|_{\cH_{\kin,\sigma}}
		+\|w^\ell f_{\st}\|_{L^\infty}
		+\|w^\ell f_{\st}\|_{C^{\alpha_0}_{\kin}}
		\le C_{\ell,\cD}
		\|w^{\ell+A_{\st}}h_{\st}\|_{(C_x^{0,1}\cap C_v^{1,1})(\gamma^-)},
	\end{align}
	and, for every $1<p<\infty$,
	\begin{equation}\label{eq:main-inflow-st-gradient}
		\|w^\ell \nabla_vf_{\st}\|_{L_x^pL_v^2}
		\le C_{\ell,p,\cD}
		\|w^{\ell+A_{\st}}h_{\st}\|_{(C_x^{0,1}\cap C_v^{1,1})(\gamma^-)}.
	\end{equation}
	
	Let $h_{\st}^1,h_{\st}^2$ satisfy
	\eqref{eq:main-inflow-st-smallness}, and let $f_{\st}^1,f_{\st}^2$ be
	corresponding stationary solutions.  If
	\[
	\|\mu^{1/4}f_{\st}^i\|_{L^\infty}\le\varepsilon_0,
	\qquad i=1,2,
	\]
	then
	\begin{equation}\label{eq:main-inflow-st-stability}
		\|f_{\st}^1-f_{\st}^2\|_{\cH_{\kin,\sigma}}
		\le C_\cD\|h_{\st}^1-h_{\st}^2\|_{L^2(\gamma^-)}.
	\end{equation}
	Consequently, the solution is unique among stationary graph solutions
	satisfying $\|\mu^{1/4}f\|_{L^\infty}\le\varepsilon_0$.
\end{theorem}
The proof is given in Subsection~\ref{subsec:inflow-stationary-proof}.

For $m\ge0$ define
\begin{align}\label{eq:dynamic-nonlinear-data-quantity-main}
	\mathsf D_m:={}&
	\|w^mf_0\|_{(C_x^{0,1}\cap C_v^{1,1})(\cD\times\bR^3)}
	+\|w^{m+1}h_{\pt}\|_{(C_{t,x}^{0,1}\cap C_v^{1,1})(\bR_+\times\gamma^-)}\notag\\
	&+\|w^mh_{\st}\|_{(C_x^{0,1}\cap C_v^{1,1})(\gamma^-)}
	+\|w^m(h_{\pt}-h_{\st})\|_{L^2(\bR_+\times\gamma^-)}.
\end{align}

\begin{theorem}[Global prescribed-inflow dynamics]\label{thm:global-dynamic-landau}
	With $\varepsilon_0$ and $\alpha_0$ as in
	Theorem~\ref{Landau_st_exi}, there exist an integer $A_{\dy}\ge0$ and
	$\delta_{\dy}>0$ such that the following holds.  Assume
	\begin{equation}\label{eq:nonlinear-dynamic-compatibility-main}
		h_{\pt}(0,x,v)=f_0(x,v),\qquad (x,v)\in\gamma^-
	\end{equation}
	and
	\begin{equation}\label{eq:dynamic-nonlinear-smallness-main}
		\mathsf D_{A_{\dy}}\le\delta_{\dy}.
	\end{equation}
	Let $f_{\st}$ be the stationary solution corresponding to $h_{\st}$ from
	Theorem~\ref{Landau_st_exi}.  Then the Coulomb Landau equation \eqref{new:eq:intro-landau},
	with initial datum $F_0=\mu+\sqrt\mu\,f_0$ and incoming datum $H$
	under \eqref{new:eq:intro-inflow-law}, has a global solution
	$F=\mu+\sqrt\mu\,f_{\pt}$, where
	$f_{\pt}\in\bigcap_{T>0}\cH_{\kin,\sigma}(\cD_T\times\bR^3)$ satisfies
	\begin{equation}\label{eq:main-inflow-dynamic-small-Linf}
		\|\mu^{1/4}f_{\pt}\|_{L^\infty(\bR_+\times\cD\times\bR^3)}
		\le\varepsilon_0.
	\end{equation}
	If $F_0\ge0$ almost everywhere in $\cD\times\bR^3$ and
	$H\ge0$ almost everywhere on $\bR_+\times\gamma^-$, then
	$F(t)\ge0$ almost everywhere for every $t\ge0$.
	
	For every $\ell\ge0$, if $\mathsf D_{\ell+A_{\dy}}<\infty$, then
	\begin{align}\label{eq:main-inflow-dynamic-regularity}
		&\sup_{T>0}\|w^\ell(f_{\pt}-f_{\st})\|_{\cH_{\kin,\sigma}(\cD_T\times\bR^3)}
		+\|w^\ell f_{\pt}\|_{L^\infty(\bR_+\times\cD\times\bR^3)}
		+\|w^\ell f_{\pt}\|_{C^{\alpha_0}_{\kin}(\bR_+\times\cD\times\bR^3)}\notag\\
		&\qquad\le C_{\ell,\cD}\mathsf D_{\ell+A_{\dy}},
	\end{align}
	and, for every $1<p<\infty$,
	\begin{equation}\label{eq:global-dynamic-mixed-gradient-main}
		\|w^\ell \nabla_vf_{\pt}\|_{L_t^\infty L_x^pL_v^2}
		\le C_{\ell,p,\cD}\mathsf D_{\ell+A_{\dy}}.
	\end{equation}
	
	For stability, let $(f_0^i,h_{\pt}^i,h_{\st}^i)$, $i=1,2$, satisfy
	\eqref{eq:nonlinear-dynamic-compatibility-main}, let $\mathsf D_m^i$ be the
	corresponding data sizes, and assume
	$\mathsf D_{A_{\dy}}^i\le\delta_{\dy}$.  If $f_{\pt}^i$ solve the
	corresponding problems on $[0,T]$ and
	\[
	\|\mu^{1/4}f_{\pt}^i\|_{L^\infty(\cD_T\times\bR^3)}\le\varepsilon_0,
	\qquad i=1,2,
	\]
	then
	\begin{align}\label{eq:main-inflow-dynamic-stability}
		\|f_{\pt}^1-f_{\pt}^2\|_{\cH_{\kin,\sigma}(\cD_T\times\bR^3)}
		\le C_\cD\Bigl(
		\|f_0^1-f_0^2\|_{L^2(\cD\times\bR^3)}
		+\|h_{\pt}^1-h_{\pt}^2\|_{L^2(\gamma_T^-)}
		\Bigr),
	\end{align}
	with $C_\cD$ independent of $T$.  In particular, the global solution is
	unique in the class
	\begin{equation}\label{eq:main-inflow-dynamic-uniqueness-class}
		f\in\bigcap_{T>0}\cH_{\kin,\sigma}(\cD_T\times\bR^3),
		\qquad
		\|\mu^{1/4}f\|_{L^\infty(\bR_+\times\cD\times\bR^3)}\le\varepsilon_0.
	\end{equation}
\end{theorem}
The proof is given in Subsection~\ref{subsec:inflow-global-proof}.

\begin{theorem}[Exponential relaxation under convergent inflow]
	\label{thm:asymptotic-stability}
	Let $\lambda_0>0$.  There exist
	$\delta_{\rm as}\in(0,\delta_{\dy}]$ and
	$\lambda\in(0,\lambda_0]$, depending only on $\cD$, the fixed reference
	parameters, and $\lambda_0$, with the following property.  Assume
	\eqref{eq:nonlinear-dynamic-compatibility-main},
	\begin{equation}\label{eq:asymptotic-weighted-data-assumption-main}
		\mathsf D_{A_{\dy}}\le\delta_{\rm as},
	\end{equation}
	and
	\begin{equation}\label{eq:asymptotic-inflow-decay-main}
		K_0:=\sup_{s\ge0}e^{\lambda_0s}
		\|h_{\pt}-h_{\st}\|_{L^2((s,\infty)\times\gamma^-)}^2<\infty.
	\end{equation}
	Let $f_{\st}$ and $f_{\pt}$ be the solutions above and set
	$f_{\dy}:=f_{\pt}-f_{\st}$.  Then
	\begin{align}\label{eq:asymptotic-stability-estimate-main}
		&\|f_{\dy}(t)\|_{L^2_{x,v}}^2
		+\int_t^\infty\|f_{\dy}(r)\|_{L_x^2H_\sigma^1}^2\,\dd r
		+\int_t^\infty
		\|(\partial_r+v\cdot\nabla_x)f_{\dy}(r)\|_{L_x^2H_\sigma^{-1}}^2\,\dd r\notag\\
		&\qquad
		+\|(f_{\dy})_{\gamma^+}\|_{L^2((t,\infty)\times\gamma^+)}^2
		\le C_\cD\bigl(\mathsf D_{A_{\dy}}^2+K_0\bigr)e^{-\lambda t},
		\qquad t\ge0.
	\end{align}
	Moreover, for every $2\le p<\infty$,
	\begin{equation}\label{eq:asymptotic-Lp-stability-main}
		\|f_{\dy}(t)\|_{L^p(\cD\times\bR^3)}
		\le C_\cD\bigl(\mathsf D_{A_{\dy}}^2+K_0\bigr)^{1/2}e^{-\lambda t/p},
		\qquad t\ge0,
	\end{equation}
	and hence
	\begin{equation}\label{eq:physical-asymptotic-stability-main}
		\|F(t)-F_{\st}\|_{L^p(\cD\times\bR^3)}
		\le C_\cD\bigl(\mathsf D_{A_{\dy}}^2+K_0\bigr)^{1/2}e^{-\lambda t/p},
		\qquad t\ge0.
	\end{equation}
\end{theorem}
The proof is given in Subsection~\ref{sec:asymptotic-behavior}.

\medskip
\noindent\textbf{Normalization for the subsequent analysis.}
	The main results above use the general Maxwellian
	\eqref{new:eq:intro-reference-maxwellian}, with arbitrary fixed
	$\rho_0,\Theta_0>0$. In the Landau analysis of Sections~5--14,
	we rescale space, time, velocity, and the distribution and fix
	$\mu(v)=e^{-|v|^2}$. We continue to denote the rescaled variables
	and domain by $(t,x,v)$ and $\cD$. The statements of
	Section~\ref{sec:main-results} in the original variables follow by
	reversing the rescaling.

\subsection{\texorpdfstring{Organization of the paper}{Organization of the paper}}\label{new:subsec:intro-organization}
Sections~3--8 develop the kinetic Fokker--Planck
		and Landau theory for prescribed inflow, including boundary regularity,
		integrability gains, and nonlinear well-posedness. Sections~9--14
	treat diffuse reflection, the stationary state, and its nonlinear stability.
	The appendices contain auxiliary results.

\section{One-Sided Kinetic Fokker–Planck Theory}\label{Sec_FP}

Let $T>0$, and let $\cD\subset\bR^d$ be a bounded open set of class \(C^{1,1}\) in the sense of Definition~\ref{domain_assumption}; 
in particular \(\partial\cD\) carries a unit outward normal \(n_x\) depending continuously on
\(x\in\partial\cD\).
We denote \(\cD_T:=(0,T)\times \cD\).

In this section, we study the dynamic inflow problem for a linear
kinetic Fokker--Planck operator in divergence form:
\begin{align}\label{FP_dy}
	\begin{cases}
		\big(\partial_t +v\cdot\nabla_x\big) f + L f=g\quad &\text{in}\quad \cD_T\times \bR^d_v\,;\\
		f(0,\,\cdot\,)= f_0\quad &\text{in}\quad \cD\times \bR^d_v\,;\\
		f=h\quad &\text{on}\quad \gamma^-_T\,.
	\end{cases}
\end{align}
We also study its stationary counterpart:
\begin{align}\label{FP_st}
	\begin{cases}
		v\cdot\nabla_x f+  L f+ \lambda f=g\quad &\text{in}\quad \cD\times \bR^d_v\,;\\
		f=h\quad &\text{on}\quad \gamma^-\,,
	\end{cases}
\end{align}
for $\lambda>0$.
Here \(L\) is given in divergence form by
\begin{alignat*}{2}
	L f\,
	:=-D_{v^i}\!\left(a^{ij} D_{v^j}f\right)+b^i D_{v^i}f+D_{v^i}\big(\bar{b}^i  f\big)+c f\,.
\end{alignat*}
We assume that the coefficients \(a^{ij}\), \(b^i\), \(\bar{b}^i\), and \(c\)
are measurable on the corresponding domain, namely
\(\cD_T\times\bR^d_v\) in the dynamic case and
\(\cD\times\bR^d_v\) in the stationary case, and that there exist constants
\(0<\nu\leq 1\) and \(K>0\) such that 
\begin{align}\label{coeff_ellip}
	\nu|\xi|^2 \le a^{ij}\xi_i\xi_j
	\quad \text{for all } \xi=(\xi_1,\ldots,\xi_d)\in\bR^d,\qquad \sup_{i,j}|a^{ij}|\le \nu^{-1}\,,
\end{align}
and
\begin{align}\label{coeff_bdd}
	|b^i|,\, |\bar{b}^i|,\, |c|\le K,
	\qquad i,j\in \{1,\ldots,d\},
\end{align}
on the corresponding domain.
Moreover, \eqref{coeff_ellip} and \eqref{coeff_bdd} give
\begin{align*}
	\|Lf\|_{L^2(\cD_T;H^{-1}_v)}
	\le
	C(d,\nu,K)\|f\|_{L^2(\cD_T;H^{1}_v)},
\end{align*}
and similarly in the stationary case.
In particular, the left-hand side of \eqref{260615144} below is a well-defined
element of \(L^2(\cD_T;H^{-1}_v)\).

For the boundary traces, let \(\dd S_x\) denote the surface measure on
\(\partial\cD\), and, for general \(d\), set
\begin{equation}\label{def:gamma}
	\begin{gathered}
		\gamma:=\partial\cD\times \bR^d_v,\qquad
		\gamma^\pm:=\{(x,v)\in\gamma:\ \pm n_x\cdot v>0\},\qquad
		\gamma^0:=\{(x,v)\in\gamma:\ n_x\cdot v=0\},\\
		\gamma_T:=(0,T)\times\gamma,\qquad
		\gamma_T^\pm:=(0,T)\times\gamma^\pm,\qquad
		\gamma_T^0:=(0,T)\times\gamma^0 .
	\end{gathered}
\end{equation}
We call \(\gamma^+\), \(\gamma^0\) and \(\gamma^-\) the outgoing, grazing and
incoming boundary, respectively.
The associated phase-space boundary measure is defined by
\begin{equation}\label{def:dgamma}
	\dd\gamma:=|n_x\cdot v|\dd S_x\dd v .
\end{equation}
We then define the corresponding \(L^2\)-spaces by
\begin{alignat*}{2}
	L^2(\gamma)
	&:=L^2(\gamma,\dd\gamma),
	&\qquad
	L^2(\gamma_T)
	&:=L^2((0,T)\times\gamma,\dd \gamma\dd t),\\
	L^2(\gamma^\pm)
	&:=L^2(\gamma^\pm,\dd\gamma),
	&\qquad
	L^2(\gamma_T^\pm)
	&:=L^2((0,T)\times\gamma^\pm,\dd \gamma\dd t).
\end{alignat*}
In particular, for any nonnegative measurable function \(\varphi\) on \(\gamma\),
\[
\int_\gamma \varphi\dd\gamma
=
\int_{\partial\cD\times\bR^d}
\varphi(x,v)\,|n_x\cdot v|\dd S_x\dd v .
\]

Throughout, \(\la\cdot,\cdot\ra_v\) denotes the duality pairing between
\(H^{-1}_v:=H^{-1}(\bR^d_v)\) and \(H^1_v:=H^1(\bR^d_v)\), normalized so that
\begin{align*}
	\la u,\varphi\ra_v=\int_{\bR^d_v}u\varphi\dd v
	\qquad\text{for }u\in L^2(\bR^d_v),\ \varphi\in H^1_v ;
\end{align*}
we use the same symbol for the pairing in either order. For a Banach space
\(B\), we write \(L^2(\cD_T;B)\) for the corresponding Bochner space. When
\(B\) is a Hilbert space, we write \(\la F,G\ra_{L^2(\cD_T;B)}\) for the
associated inner product.

\subsection{Solution space $\cH_\kin$ and the well-posedness theorem}

We introduce the dynamic and stationary \(\cH_{\kin}\)-spaces, their trace
structure, and the corresponding notion of solution before stating the main
well-posedness theorem.

\begin{defn}[$\cH_{\kin}$-spaces]\label{def:H_kin}\,
	
	\begin{enumerate}
		\item For \(f\in C_c^\infty\big([0,T]\times \overline{\cD}\times \bR^d_v\big)\), we define
		\begin{align}\label{innerprod}
			\begin{aligned}
				\|f\|_{\cH_{\kin}(\cD_T \times \bR^d_v)}^2:=\,&
				\|f\|_{L^2(\cD_T;H_v^1)}^2
				+\|(\partial_t+v\cdot \nabla_x )f\|_{L^2(\cD_T;H_v^{-1})}^2\\
				&\quad
				+\|f(0,\,\cdot\,)\|_{L^2(\cD\times \bR^d_v)}^2
				+\big\|f|_{\gamma_T^-}\big\|_{L^2(\gamma^-_T)}^2\,.
			\end{aligned}
		\end{align}
		We define \(\cH_{\kin}(\cD_T\times \bR^d_v)\) as the completion of
		\(C_c^\infty\big([0,T]\times \overline{\cD}\times \bR^d_v\big)\) with respect to the norm \eqref{innerprod}. 
		Equivalently (Proposition~\ref{260128509}), we regard \(f\) as an element of
		\(\cH_{\kin}(\cD_T\times \bR^d_v)\) if there exists a sequence
		\(f_n\in C_c^\infty\big([0,T]\times \overline{\cD}\times \bR^d_v\big)\) such that
		\begin{equation}\label{def:Hkin-approx-L2}
			f_n\longrightarrow f\quad \text{in}\quad L^2(\cD_T\times\bR^d_v),
		\end{equation}
		and
		\begin{equation}\label{def:Hkin-approx-Cauchy}
			\|f_n-f_m\|_{\cH_{\kin}(\cD_T\times \bR^d_v)}
			\longrightarrow 0\quad \text{as}\quad n,m\rightarrow \infty\,.
		\end{equation}
		For \(f\in \cH_{\kin}(\cD_T\times \bR^d_v)\), we define
		\((\p_t+v\cdot \nabla_x)f\), \(f(0,\,\cdot\,)\), and \(f_{\gamma^-}\) 
		as the corresponding limits
		\begin{align}\label{lim_transp}
			\begin{alignedat}{2}
				(\p_t+v\cdot\nabla_x)f_n
				&\longrightarrow (\p_t+v\cdot \nabla_x)f
				\quad &&\text{in}\quad L^2(\cD_T;H^{-1}_v),\\
				f_n(0,\cdot)
				&\longrightarrow f(0,\cdot)
				\quad &&\text{in}\quad L^2(\cD\times \bR^d_v),\\
				f_n|_{\gamma^-_T}
				&\longrightarrow f_{\gamma^-}
				\quad &&\text{in}\quad L^2(\gamma^-_T),
			\end{alignedat}
		\end{align}
		where \((f_n)\) is any sequence converging to \(f\) in
		\(\cH_{\kin}(\cD_T\times \bR^d_v)\).
		
		\item For \(f\in C_c^\infty\big(\overline{\cD}\times \bR^d_v\big)\), we define
		\begin{align}\label{innerprod_st}
			\begin{aligned}
				\|f\|_{ {\cH}_{\kin}(\cD \times \bR^d_v)}^2:=\,&
				\|f\|_{L^2(\cD ;H_v^1)}^2
				+\|v\cdot \nabla_x f\|_{L^2(\cD ;H_v^{-1})}^2
				+\big\|f|_{\gamma^-}\big\|_{L^2({\gamma^-})}^2\,.
			\end{aligned}
		\end{align}
		We define \(\cH_{\kin}(\cD\times \bR^d_v)\) as the completion of
		\(C_c^\infty\big(\overline{\cD}\times \bR^d_v\big)\) with respect to the norm \eqref{innerprod_st}.
		For \(f\in\cH_{\kin}(\cD\times \bR^d_v)\), we define
		\(v\cdot \nabla_x f\) and \(f_{\gamma^-}\) as the corresponding limits
		\begin{align}\label{lim_transp_st}
			\begin{alignedat}{2}
				v\cdot\nabla_x f_n
				&\longrightarrow v\cdot \nabla_x f
				\quad &&\text{in}\quad L^2(\cD;H^{-1}_v),\\
				f_n|_{\gamma^-}
				&\longrightarrow f_{\gamma^-}
				\quad &&\text{in}\quad L^2(\gamma^-),
			\end{alignedat}
		\end{align}
		where \((f_n)\) is any sequence converging to \(f\) in
		\(\cH_{\kin}(\cD\times \bR^d_v)\). 
	\end{enumerate}
\end{defn}

\begin{remark}[Time translation]\label{rem:Hkin-time-translation}
	The spaces \(\cH_{\kin}\) have been defined on \((0,T)\). For a bounded
	interval \(I=(t_1,t_2)\), we define \(\cH_{\kin}(I\times\cD\times\bR^d_v)\)
	by \eqref{innerprod} with \((0,T)\) replaced by \(I\) and with the initial
	trace taken at \(t=t_1\). Both the norm \eqref{innerprod} and the class
	\(C_c^\infty\big([0,T]\times\overline{\cD}\times\bR^d_v\big)\) are invariant
	under translation in \(t\). Hence every statement of this section holds
	verbatim on \(I\), with \(f(0,\cdot)\) replaced by \(f(t_1,\cdot)\) and with
	the constants depending on \(t_2-t_1\) only through their dependence on
	\(T\).
\end{remark}

\begin{defn}[$\cH_{\kin}$-solutions]\label{def:Hkin-solution}
	Given
	\[
	g\in L^2(\cD_T;H^{-1}_v),\qquad
	f_0\in L^2(\cD\times \bR^d_v),\qquad
	h\in L^2(\gamma^-_T),
	\]
	we say that \(f\in \cH_\kin(\cD_T\times \bR^d_v)\) is an
	\(\cH_{\kin}\)-solution of \eqref{FP_dy} with initial datum
	\(f_0\) and incoming boundary datum \(h\) if
	\begin{align}\label{260615144}
		(\p_t+v\cdot \nabla_x)f+Lf=g
		\quad \text{in}\quad L^2(\cD_T;H^{-1}_v),
	\end{align}
	and
	\begin{align}\label{260615145}
		f(0,\cdot)=f_0
		\quad \text{in}\quad L^2(\cD\times \bR^d_v),
		\qquad
		f_{\gamma^-}=h
		\quad \text{in}\quad L^2(\gamma^-_T).
	\end{align}
	
	In the stationary case, given
	\[
	g\in L^2(\cD;H^{-1}_v),\qquad
	h\in L^2(\gamma^-),
	\]
	we say that \(f\in \cH_\kin(\cD\times \bR^d_v)\) is an
	\(\cH_{\kin}\)-solution of \eqref{FP_st} with incoming boundary
	datum \(h\) if
	\begin{equation}\label{260615146}
		v\cdot \nabla_x f+Lf+\lambda f=g
		\quad \text{in}\quad L^2(\cD;H^{-1}_v),
		\qquad
		f_{\gamma^-}=h
		\quad \text{in}\quad L^2(\gamma^-).
	\end{equation}
\end{defn}

In the weak formulations below, we use test functions
\(\phi\in C_c^\infty([0,T]\times \overline{\cD}\times \bR_v^d)\) subject to the
support condition
\begin{equation}\label{260616633}
	\supp\,\phi\subset
	[0,T)\times \left((\cD\times \bR_v^d)\cup \gamma^-\right),
\end{equation}
that is, \(\phi\) vanishes near \(\{T\}\times\overline\cD\times\bR^d_v\) and
near \([0,T]\times(\gamma^+\cup\gamma^0)\).

\begin{remark}[Weak formulation]\label{weak_formulation}
	The following weak formulation is equivalent to
	Definition~\ref{def:Hkin-solution}.
	
	Every \(g\in L^2(\cD_T;H^{-1}_v)\) can be written as
	$$
	g=g^0+D_i g^i
	\quad\text{in }L^2(\cD_T;H^{-1}_v),
	\qquad
	g^0,g^1,\ldots,g^d\in L^2(\cD_T\times\bR^d_v),
	$$
	with 
	$$\sum_{i=0}^d\|g^i\|_{L^2}\le C(d)\|g\|_{L^2(\cD_T;H^{-1}_v)}.
	$$
	
	Fix one such decomposition, and let \(f_0\in L^2(\cD\times\bR^d_v)\) and
	\(h\in L^2(\gamma_T^-)\). Then \(f\in \cH_{\kin}(\cD_T\times\bR^d_v)\) is an
	\(\cH_{\kin}\)-solution of \eqref{FP_dy} with data \((g,f_0,h)\) if and only
	if, for every \(\phi\in C_c^\infty([0,T]\times\overline{\cD}\times\bR^d_v)\)
	satisfying \eqref{260616633},
	\begin{align}\label{260615125}
		\begin{aligned}
			&\iiint_{\cD_T\times \bR^d_v}
			\left(
			a^{ij}D_j fD_i\phi
			+b^i(D_i f)\phi
			-\bar b^i fD_i\phi
			+cf\phi
			-f(\p_t\phi+v\cdot\nabla_x\phi)
			\right)\dd v\dd x\dd t
			\\
			=\,&
			\iiint_{\cD_T\times \bR^d_v}
			\left(g^0\phi-g^iD_i\phi\right)\dd v\dd x\dd t
			+\iint_{\cD\times\bR^d_v}f_0\,\phi(0,\cdot)\dd v\dd x
			+\iint_{\gamma_T^-}h\phi\dd\gamma\dd t .
		\end{aligned}
	\end{align}
	Indeed, by \eqref{ibp_dy}, the two sides of \eqref{260615125} differ by the
	pairings of \((\p_t+v\cdot\nabla_x)f+Lf-g\) with \(\phi\), of
	\(f(0,\cdot)-f_0\) with \(\phi(0,\cdot)\), and of \(f_{\gamma^-}-h\) with
	\(\phi|_{\gamma_T^-}\). Arguing as in the paragraph following
	\eqref{2606171247}, the vanishing of their sum for every such \(\phi\) is
	equivalent to
	\[
	(\p_t+v\cdot\nabla_x)f+Lf=g,\qquad
	f(0,\cdot)=f_0,\qquad
	f_{\gamma^-}=h
	\]
	in the corresponding spaces.
	
	The stationary case is identical: for \(g=g^0+D_i g^i\) in
	\(L^2(\cD;H^{-1}_v)\) with \(g^0,\ldots,g^d\in L^2(\cD\times\bR^d_v)\), and
	for \(h\in L^2(\gamma^-)\), a function \(f\in \cH_{\kin}(\cD\times\bR^d_v)\)
	is an \(\cH_{\kin}\)-solution of \eqref{FP_st} with incoming datum \(h\) if
	and only if \eqref{260615125} holds with the time integrations and the
	initial term removed, with \(-f\p_t\phi\) dropped and \(\lambda f\phi\)
	added on the left-hand side, for every
	\(\phi\in C_c^\infty(\overline{\cD}\times\bR^d_v)\) supported in
	\((\cD\times \bR^d_v)\cup \gamma^-\).
\end{remark}

The main well-posedness result for \eqref{FP_dy} and \eqref{FP_st} is proved in
Subsection~\ref{260616916}.

\begin{thm}[Well-posedness in $\cH_{\kin}$]\label{solv}
	\,
	
	\begin{enumerate}
		\item For every
		\begin{align}\label{260617519}
			g\in L^2\big(\cD_T;H_v^{-1}\big),\qquad
			f_0\in L^2(\cD\times \bR_v^d),\qquad
			h\in L^2(\gamma_T^-),
		\end{align}
		equation \eqref{FP_dy} admits a unique solution \(f\) in
		\(\cH_{\kin}(\cD_T \times \bR^d_v)\).
		Moreover, there exists a constant \(C=C(d,\nu,K,T)>0\) such that
		\begin{align}\label{dyHbdd}
			\|f\|_{\cH_{\kin}(\cD_T\times\bR^d_v)}
			\le
			C\left(
			\|g\|_{L^2(\cD_T;H_v^{-1})}
			+\|f_0\|_{L^2(\cD\times \bR_v^d)}
			+\|h\|_{L^2(\gamma_T^-)}
			\right).
		\end{align}
		
		\item There exists \(\lambda_0=\lambda_0(d,\nu,K)\ge 0\) such that, for every
		\(\lambda\ge\lambda_0\) and
		\begin{align}\label{260617520}
			g\in L^2\big(\cD;H_v^{-1}\big),\qquad
			h\in L^2(\gamma^-),
		\end{align}
		equation \eqref{FP_st} admits a unique solution \(f\) in
		\(\cH_{\kin}(\cD\times \bR^d_v)\).
		Moreover, there exists a constant \(C=C(d,\nu,K,\lambda)>0\) such that
		\begin{align}\label{stHbdd}
			\|f\|_{\cH_{\kin}(\cD\times\bR^d_v)}
			\le
			C\left(
			\|g\|_{L^2(\cD;H_v^{-1})}
			+\|h\|_{L^2(\gamma^-)}
			\right).
		\end{align}
	\end{enumerate}
\end{thm}

The spaces \(\cH_{\kin}\) were defined above by completion of smooth
functions. In Subsection~\ref{sec:weakform}, we use
Theorem~\ref{solv} to show that membership in these spaces can be
characterized intrinsically by velocity Sobolev regularity, the transport
derivative, the appropriate time and boundary traces, and the corresponding
Green identity. This characterization connects the completion-based
definition with the usual weak transport framework. It is logically
subsequent to the well-posedness theorem and is not used in its proof.

\subsection{Properties of the \(\cH_{\kin}\)-solution spaces}\label{Prop_H_kin}

We establish the basic analytic properties of the \(\cH_{\kin}\)-spaces:
their Hilbert-space structure, the existence of time and phase-space traces
and the corresponding Green identities, stability under cut-offs and time
reversal, renormalized identities, and preservation of essential bounds
under traces.

\begin{prop}[Well-definedness and Hilbert-space structure of $\cH_{\kin}$]
	\label{260128509}
	\,
	
	\begin{enumerate}
		\item For each \(f\in \cH_{\kin}(\cD_T \times \bR^d_v)\),
		\((\p_t+v\cdot \nabla_x)f\), \(f(0,\,\cdot\,)\), and \(f_{\gamma^-}\)
		are uniquely determined. Moreover,
		\(\cH_{\kin}(\cD_T\times \bR^d_v)\) is a Hilbert space equipped
		with the inner product
		\begin{align*}
			(f,g)_{\cH_{\kin}(\cD_T\times\bR^d_v)}:=\,&
			(f,g)_{L^2(\cD_T;H^1_v)}
			+\big((\p_t+v\cdot \nabla_x)f,(\p_t+v\cdot \nabla_x)g\big)_{L^2(\cD_T;H^{-1}_v)}\\
			&\quad
			+\big(f(0,\cdot),g(0,\cdot)\big)_{L^2(\cD\times \bR^d_v)}
			+\big(f_{\gamma^-},g_{\gamma^-}\big)_{L^2(\gamma^-_T)}.
		\end{align*}
		
		\item For each \(f\in \cH_{\kin}(\cD \times \bR^d_v)\),
		\(v\cdot \nabla_x f\) and \(f_{\gamma^-}\) are uniquely determined.
		Moreover, \(\cH_{\kin}(\cD\times \bR^d_v)\) is a Hilbert space
		equipped with the inner product
		\begin{align*}
			(f,g)_{\cH_{\kin}(\cD\times\bR^d_v)}:=\,&
			(f,g)_{L^2(\cD;H^1_v)}
			+\big(v\cdot \nabla_x f,v\cdot \nabla_x g\big)_{L^2(\cD;H^{-1}_v)}
			+\big(f_{\gamma^-},g_{\gamma^-}\big)_{L^2(\gamma^-)}.
		\end{align*}
	\end{enumerate}
\end{prop}

The proof of Proposition~\ref{260128509} uses the following density result for test functions and trace pairs, which we establish first.

\begin{lemma}\label{260418814}\,
	\begin{enumerate}
		\item
		\(C_c^\infty(\cD_T\times \bR_v^d)\) is dense in
		\(L^2(\cD_T;H_v^{-1})\) and in \(L^2(\cD_T;H_v^{1})\). In addition, let \(\cE\) denote the set of all
		pairs of the form
		\[
		\big(\phi(0,\cdot,\cdot),\phi|_{\gamma_T^-}\big),
		\]
		where \(\phi\in C_c^\infty([0,T]\times \overline{\cD}\times \bR_v^d)\)
		satisfies \eqref{260616633}.
		Then \(\cE\) is dense in
		\(L^2(\cD\times \bR_v^d)\times L^2(\gamma_T^-)\).
		
		\item
		\(C_c^\infty(\cD\times \bR_v^d)\) is dense in
		\(L^2(\cD;H_v^{-1})\) and in \(L^2(\cD;H_v^{1})\). In addition, the set
		\[
		\left\{
		\phi|_{\gamma^-}:
		\phi\in C_c^\infty(\overline{\cD}\times \bR_v^d),
		\quad
		\supp\,\phi\subset (\cD\times \bR_v^d)\cup\gamma^-
		\right\}
		\]
		is dense in \(L^2(\gamma^-)\).
	\end{enumerate}
\end{lemma}

\begin{proof}
	It suffices to prove \((1)\); the proof of \((2)\) is the same, and simpler,
	since only one of the two pieces occurs in Step 2 below.
	
	\medskip
	\noindent
	\emph{Step 1. Density in \(L^2(\cD_T;H_v^{\mp1})\).}
	Let \(B\) be either \(H^{-1}_v\) or \(H^1_v\). Since \(C_c^\infty(\bR^d_v)\)
	is dense in \(B\), the algebraic tensor product
	\(C_c^\infty(\cD_T)\otimes C_c^\infty(\bR_v^d)\) is dense in
	\(L^2(\cD_T;B)\) by the standard density theorem for Bochner spaces; it is
	contained in \(C_c^\infty(\cD_T\times\bR_v^d)\), which proves both density
	assertions.
	
	\medskip
	\noindent
	\emph{Step 2. Density of the trace pairs.}
	Fix
	$(f_0,h)\in L^2(\cD\times \bR_v^d)\times L^2(\gamma_T^-)$.
	
	Since \(\cD\times\bR_v^d\) and \(\gamma_T^-\) are locally compact Hausdorff
	spaces endowed with Radon measures, \(C_c\) is dense in the corresponding
	\(L^2\)-spaces. Hence there exist $F_n\in C_c(\cD\times\bR_v^d)$ and $H_n\in C_c(\gamma_T^-)$
	such that
	\begin{align}\label{2606171201}
		\|F_n-f_0\|_{L^2(\cD\times\bR_v^d)}
		+
		\|H_n-h\|_{L^2(\gamma_T^-)}
		\le \frac1n .
	\end{align}
	
	Set
	\[
	S:=\overline{\big(\{0\}\times\cD\times\bR_v^d\big)\cup\gamma_T^-}.
	\]
	Define \(\Phi_n\) on \(S\) by
	\[
	\Phi_n(0,x,v):=F_n(x,v),
	\qquad
	\Phi_n|_{\gamma_T^-}:=H_n,
	\]
	and by zero on the remaining part of \(S\). Since \(F_n\) and \(H_n\) are
	compactly supported in \(\cD\times\bR_v^d\) and \(\gamma_T^-\), respectively,
	\(\Phi_n\) is continuous on \(S\), and \(\supp\,\Phi_n\) is compact and
	disjoint from
	\[
	A:=
	\big(\{T\}\times\overline{\cD}\times\bR_v^d\big)
	\cup \big([0,T]\times(\gamma^+\cup\gamma^0)\big).
	\]
	
	Choose \(\eta_n\in C_c^\infty(\bR^{1+2d})\) such that
	\[
	0\le\eta_n\le1,
	\qquad
	\eta_n\equiv1
	\quad\text{on } \supp\,\Phi_n,
	\qquad
	\supp\,\eta_n\cap A=\emptyset .
	\]
	By the Tietze extension theorem, \(\Phi_n\) has a bounded continuous extension
	to \(\bR^{1+2d}\). The measures \(\dd x\dd v\) and \(\dd\gamma\dd t\) are
	finite on the compact set \(\supp\,\eta_n\), so a uniform bound there yields
	an \(L^2\) bound. Approximating this extension uniformly on
	\(\supp\,\eta_n\) by smooth functions, we obtain
	\(\psi_{n,m}\in C^\infty(\bR^{1+2d})\) such that, setting
	\[
	\phi_{n,m}:=
	(\eta_n\psi_{n,m})|_{[0,T]\times\overline{\cD}\times\bR_v^d},
	\]
	we have
	\[
	\phi_{n,m}\in C_c^\infty([0,T]\times\overline{\cD}\times\bR_v^d),
	\qquad
	\supp\,\phi_{n,m}\cap A=\emptyset,
	\]
	and
	\begin{align}\label{2606171202}
		\|(\Phi_n-\phi_{n,m})(0,\cdot)\|_{L^2(\cD\times\bR_v^d)}
		+
		\|(\Phi_n-\phi_{n,m})|_{\gamma_T^-}\|_{L^2(\gamma_T^-)}
		\le \frac1m .
	\end{align}
	
	Since \(\Phi_n(0,\cdot)=F_n\) and \(\Phi_n|_{\gamma_T^-}=H_n\), combining
	\eqref{2606171201} and \eqref{2606171202} gives
	\[
	\|f_0-\phi_{n,m}(0,\cdot)\|_{L^2(\cD\times\bR_v^d)}
	+
	\|h-\phi_{n,m}|_{\gamma_T^-}\|_{L^2(\gamma_T^-)}
	\le \frac1n+\frac1m .
	\]
	Letting \(m\to\infty\) and then \(n\to\infty\), we obtain the desired density
	of \(\cE\).
\end{proof}

\begin{proof}[\textbf{Proof of Proposition~\ref{260128509}}]
	Since the proof of (2) is similar, we prove only (1).
	
	Let $(f_n)_{n\in\bN}$ and $(\bar f_n)_{n\in\bN}$ be two approximating sequences
	for $f$ in the definition of
	$\cH_{\kin}(\cD_T\times \bR^d_v)$. Let $(F,G,H)$ and
	$(\bar F,\bar G,\bar H)$ denote the respective limits of
	\[
	\big((\p_t +v\cdot \nabla_x) f_n,\ f_n(0,\cdot),\ f_n|_{\gamma_T^-}\big)
	\quad \text{and}\quad
	\big((\p_t +v\cdot \nabla_x) \bar f_n,\ \bar f_n(0,\cdot),\ \bar f_n|_{\gamma_T^-}\big)
	\]
	in the corresponding function spaces.
	For any $\varphi\in C_c^\infty([0,T]\times \overline{\cD}\times \bR^d_v)$
	satisfying \eqref{260616633}, the classical Green formula for the transport
	operator gives
	\begin{align*}
		&\iint_{\cD_T}
		\Big(
		\big\la (\p_t+v\cdot \nabla_x)f_n,\varphi\big\ra_{v}
		+\big\la f_n,(\partial_t+v\cdot \nabla_x)\varphi\big\ra_{v}
		\Big)\dd x\dd t \\
		=\,&
		-\iint_{\cD\times \bR^d_v}
		f_n(0,\cdot)\varphi(0,\cdot)\dd x\dd v
		+\iint_{\gamma_T^-}
		f_n|_{\gamma_T^-}(t,x,v)\,
		\varphi|_{\gamma_T^-}(t,x,v)\,
		\sgn(n_x\cdot v)\dd \gamma \dd t .
	\end{align*}
	Passing to the limit as $n\to\infty$, and using
	$\sgn(n_x\cdot v)=-1$ on $\gamma_T^-$, we obtain
	\begin{align*}
		\begin{aligned}
			&\iint_{\cD_T}
			\Big(
			\la F,\varphi\ra_{v}
			+\la f,(\partial_t+v\cdot \nabla_x)\varphi\ra_{v}
			\Big)\dd x\dd t\\
			=\,&
			-\iint_{\cD\times \bR^d_v}
			G\,\varphi(0,\cdot)\dd x\dd v
			-\iint_{\gamma_T^-}
			H\,\varphi|_{\gamma_T^-}\dd \gamma \dd t .
		\end{aligned}
	\end{align*}
	The same identity holds with $(F,G,H)$ replaced by $(\bar F,\bar G,\bar H)$.
	Subtracting the two identities gives, for every
	$\varphi\in C_c^\infty([0,T]\times \overline{\cD}\times \bR^d_v)$ satisfying
	\eqref{260616633},
	\begin{align}\label{2606171247}
		\iint_{\cD_T}\la F-\bar F,\varphi\ra_{v}\dd x\dd t
		&+\iint_{\cD\times \bR^d_v}(G-\bar G)\,\varphi(0,\cdot)\dd x\dd v \\
		&+\iint_{\gamma_T^-}(H-\bar H)\,\varphi|_{\gamma_T^-}\dd \gamma \dd t
		=0 .\nonumber
	\end{align}
	Taking $\varphi\in C_c^\infty(\cD_T\times \bR^d_v)$, the trace terms drop out.
	Since $L^2(\cD_T;H^1_v)$ is the predual of $L^2(\cD_T;H^{-1}_v)$, and since
	$C_c^\infty(\cD_T\times \bR^d_v)$ is dense in $L^2(\cD_T;H^1_v)$ by
	Lemma~\ref{260418814}.(1), we obtain $F=\bar F$. Therefore
	\[
	\iint_{\cD\times \bR^d_v}(G-\bar G)\,\varphi(0,\cdot)\dd x\dd v
	+
	\iint_{\gamma_T^-}(H-\bar H)\,\varphi|_{\gamma_T^-}\dd \gamma \dd t
	=0
	\]
	for every admissible $\varphi$. By the trace-pair density in
	Lemma~\ref{260418814}.(1), it follows that
	\(
	(G,H)=(\bar G,\bar H)\).
	Therefore
	$(\p_t+v\cdot \nabla_x)f$, $f(0,\cdot)$, and $f_{\gamma^-}$ are well defined.
	
	This also justifies the equivalent description in Definition~\ref{def:H_kin}(1).
	If $(f_n)$ is Cauchy in $\|\cdot\|_{\cH_\kin(\cD_T\times \bR^d_v)}$, then it is
	Cauchy in $L^2(\cD_T;H^1_v)$, hence converges there, and its limit is
	determined by its limit in $L^2(\cD_T\times \bR^d_v)$. Together with the
	uniqueness of $(F,G,H)$ just proved, this shows that the map from the abstract
	completion to $L^2(\cD_T;H^1_v)$ is injective.
	
	It remains to verify the Hilbert space structure of
	$\cH_\kin(\cD_T\times \bR^d_v)$.
	By the uniqueness proved above, the inner product is well defined.
	Its positive definiteness is immediate from the first term
	$\|f\|_{L^2(\cD_T;H^1_v)}$.
	
	Let $(f^m)_{m\in\mathbb N}$ be a Cauchy sequence in
	$\cH_\kin(\cD_T\times \bR^d_v)$.
	Then $(f^m)$ is Cauchy in $L^2(\cD_T;H_v^1)$, and the sequences
	\[
	\big((\partial_t+v\cdot\nabla_x)f^m\big)_{m\in\mathbb N},
	\qquad
	\big(f^m(0,\cdot)\big)_{m\in\mathbb N},
	\qquad
	\big((f^m)_{\gamma^-}\big)_{m\in\mathbb N}
	\]
	are Cauchy in their respective Hilbert spaces. Hence, for some
	\[
	f\in L^2(\cD_T;H_v^1),\qquad
	F\in L^2(\cD_T;H_v^{-1}),\qquad
	G\in L^2(\cD\times\bR_v^d),\qquad
	H\in L^2(\gamma_T^-),
	\]
	we have
	\[
	f^m\to f,\qquad
	(\partial_t+v\cdot\nabla_x)f^m\to F,\qquad
	f^m(0,\cdot)\to G,\qquad
	(f^m)_{\gamma^-}\to H
	\]
	in the corresponding spaces.
	
	For each $m$, choose an approximating sequence
	$(f_n^m)_n\subset C_c^\infty([0,T]\times\overline{\cD}\times\bR_v^d)$
	for $f^m$. Then choose $n(m)$ sufficiently large so that
	\[
	\|f_{n(m)}^m-f^m\|_{\cH_\kin(\cD_T\times \bR^d_v)}<2^{-m}.
	\]
	Set $\phi_m:=f_{n(m)}^m$. Then the diagonal sequence $(\phi_m)_{m\in\mathbb N}$
	satisfies
	\[
	\phi_m\to f \quad \text{in } L^2(\cD_T;H_v^1),
	\qquad
	\lim_{m,\ell\rightarrow \infty}
	\|\phi_m-\phi_\ell\|_{\cH_\kin(\cD_T\times \bR^d_v)}=0.
	\]
	Moreover,
	\[
	(\p_t+v\cdot \nabla_x)\phi_m\to F,\qquad
	\phi_m(0,\cdot)\to G,\qquad
	\phi_m|_{\gamma_T^-}\to H
	\]
	in the corresponding spaces. Therefore, by Definition~\ref{def:H_kin}, we have
	$f\in \cH_\kin(\cD_T\times \bR^d_v)$, and
	\[
	f^m\to f
	\quad \text{in } \cH_\kin(\cD_T\times \bR^d_v).
	\]
	Thus $\cH_\kin(\cD_T\times \bR^d_v)$ is complete.
\end{proof}

\begin{prop}[Trace theorem and Green's identities for \(\cH_{\kin}\)]
	\label{trace_theorem}\,
	
	\begin{enumerate}
		\item
		The space \(\cH_{\kin}(\cD_T \times \bR^d_v)\) continuously embeds into
		\[
		C\big([0,T];L^2(\cD\times \bR^d_v)\big).
		\]
		Moreover, there exists a phase-space boundary trace operator
		\[
		T_\gamma:\cH_{\kin}(\cD_T\times \bR^d_v)\to L^2(\gamma_T),
		\qquad
		T_\gamma f:=f_\gamma,
		\]
		such that, for all \(f,g\in \cH_{\kin}(\cD_T\times \bR^d_v)\) and all
		\(0\le s\le t\le T\),
		\begin{align}\label{ibp_dy}
			\begin{aligned}
				&\int_s^t\int_{\cD}
				\Big(
				\la (\partial_r +v\cdot \nabla_x) f,g\ra_{v}
				+\la f,(\partial_r+v\cdot \nabla_x) g\ra_{v}
				\Big)\dd x\dd r\\
				=\,&
				\iint_{\cD\times \bR^d_v} f(t,\cdot)g(t,\cdot)\dd x\dd v
				-\iint_{\cD\times \bR^d_v} f(s,\cdot)g(s,\cdot)\dd x\dd v\\
				&\quad
				+\int_s^t\int_{\gamma}
				f_{\gamma} g_{\gamma}\,\sgn(n_x\cdot v)\dd \gamma\dd r .
			\end{aligned}
		\end{align}
		In addition, for every \(f\in \cH_{\kin}(\cD_T\times \bR^d_v)\),
		\begin{align*}
			\sup_{0\le t\le T}\|f(t,\cdot)\|_{L^2(\cD\times \bR^d_v)}
			+\|f_\gamma\|_{L^2(\gamma_T)}
			\le C \|f\|_{\cH_{\kin}(\cD_T\times\bR^d_v)},
		\end{align*}
		where \(C\) is a universal constant.
		
		\item
		There exists a phase-space boundary trace operator
		\[
		T_\gamma:\cH_{\kin}(\cD\times \bR^d_v)\to L^2(\gamma),
		\qquad
		T_\gamma f:=f_\gamma,
		\]
		such that, for all \(f,g\in \cH_{\kin}(\cD\times \bR^d_v)\),
		\begin{align}\label{ibp_st}
			\int_{\cD}
			\Big(
			\la v\cdot \nabla_x f,g\ra_{v}
			+\la f,v\cdot \nabla_x g\ra_{v}
			\Big)\dd x
			=
			\int_{\gamma}f_{\gamma} g_{\gamma}\,\sgn(n_x\cdot v)\dd \gamma .
		\end{align}
		In addition, for every \(f\in \cH_{\kin}(\cD\times \bR^d_v)\),
		\begin{align*}
			\|f_\gamma\|_{L^2(\gamma)}
			\le C \|f\|_{\cH_{\kin}(\cD\times\bR^d_v)},
		\end{align*}
		where \(C\) is a universal constant.
	\end{enumerate}
	Here and below we use that \(\gamma^0\) is \(\dd\gamma\)-null, by
	\eqref{def:dgamma}, so that
	\(L^2(\gamma_T)=L^2(\gamma^+_T)\oplus L^2(\gamma^-_T)\).
\end{prop}

\begin{proof}
	Since the proof of (2) is similar, we only provide the proof of (1).
	
	\medskip
	\noindent
	\emph{Step 1.}
	Let \(f,g\in C_c^\infty([0,T]\times \overline{\cD}\times \bR^d_v)\).
	Then \eqref{ibp_dy} follows from the fundamental theorem of calculus in time
	and the divergence theorem in the \(x\)-variable.
	
	For \(f\in C_c^\infty([0,T]\times \overline{\cD}\times \bR^d_v)\), taking
	\(g=f\) and \(s=0\) in \eqref{ibp_dy}, we obtain
	\begin{align*}
		&\|f(t,\cdot)\|_{L^2(\cD\times \bR^d_v)}^2
		+\int_0^t\int_{\gamma^+}|f_\gamma|^2 \dd \gamma \dd r\\
		=\,&
		\|f(0,\cdot)\|_{L^2(\cD\times \bR^d_v)}^2
		+\int_0^t\int_{\gamma^-}|f_\gamma|^2 \dd \gamma \dd r
		+2\int_0^t\int_{\cD}
		\la (\p_t+v\cdot \nabla_x)f,f\ra_v\dd x\dd r\\
		\le\,&
		\|f(0,\cdot)\|_{L^2(\cD\times \bR^d_v)}^2
		+\|f_{\gamma^-}\|_{L^2(\gamma^-_T)}^2
		+\|(\p_t+v\cdot \nabla_x)f\|_{L^2(\cD_T;H^{-1}_v)}^2
		+\|f\|_{L^2(\cD_T;H^1_v)}^2,
	\end{align*}
	where the last inequality follows from the duality between \(H^{-1}_v\) and
	\(H^1_v\), together with Young's inequality. Consequently, there exists a universal
	constant \(C>0\) such that, for every
	\(f\in C_c^\infty([0,T]\times \overline{\cD}\times \bR^d_v)\),
	\begin{align}\label{2604201108}
		\|f\|_{C([0,T];L^2(\cD\times \bR^d_v))}
		+\|f_\gamma\|_{L^2(\gamma_T)}
		\le
		C\|f\|_{\cH_{\kin}(\cD_T\times\bR^d_v)}.
	\end{align}
	Since \(C_c^\infty([0,T]\times \overline{\cD}\times \bR^d_v)\) is dense in
	\(\cH_\kin(\cD_T\times \bR^d_v)\), \eqref{2604201108} yields the continuous
	embedding
	\[
	\cH_\kin(\cD_T\times \bR^d_v)
	\hookrightarrow
	C\big([0,T];L^2(\cD\times \bR^d_v)\big),
	\]
	as well as the existence of the trace operator \(T_\gamma\) stated in the
	proposition.
	
	\medskip
	\noindent
	\emph{Step 2.}
	Fix \(f,g\in \cH_\kin(\cD_T\times \bR^d_v)\), and take approximating sequences
	\[
	(f_n)_n,\ (g_n)_n\subset C_c^\infty([0,T]\times \overline{\cD}\times \bR^d_v)
	\]
	for \(f\) and \(g\), respectively. Then
	\[
	f_n\to f,\qquad g_n\to g
	\quad \text{in } \cH_\kin(\cD_T\times \bR^d_v).
	\]
	In particular, by the continuous embedding proved in Step 1,
	\[
	f_n\to f,\qquad g_n\to g
	\quad \text{in } C\big([0,T];L^2(\cD\times \bR^d_v)\big).
	\]
	Moreover,
	\[
	T_\gamma f_n=f_n|_{\gamma_T}\to f_\gamma:=T_\gamma f,
	\qquad
	T_\gamma g_n=g_n|_{\gamma_T}\to g_\gamma:=T_\gamma g
	\quad \text{in } L^2(\gamma_T).
	\]
	
	For each \(n\), the identity \eqref{ibp_dy} holds with \(f_n\) and \(g_n\)
	in place of \(f\) and \(g\). Passing to the limit as \(n\to\infty\) in
	\eqref{ibp_dy}, using the above convergences, we conclude that
	\eqref{ibp_dy} holds for \(f\) and \(g\).
	
	In particular, a strongly convergent sequence in $\cH_\kin$
	cannot concentrate its boundary flux at grazing velocities.  Indeed,
	if $f_j\to f$ in $\cH_\kin$, write
	$E_\eta:=\{(t,x,v)\in\gamma_T:|n_x\cdot v|<\eta\}$.
	The trace estimate just proved gives
	\[
	\int_{E_\eta}|(f_j)_\gamma|^2\dd\gamma\dd t
	\le 2C^2\|f_j-f\|_{\cH_\kin}^2
	+2\int_{E_\eta}|f_\gamma|^2\dd\gamma\dd t.
	\]
	The first term tends to zero uniformly for $j$ sufficiently large;
	the second tends to zero as $\eta\downarrow0$ by dominated
	convergence.  The finitely many remaining traces satisfy the same
	conclusion individually.  Thus
	\[
	\lim_{\eta\downarrow0}\sup_j
	\int_{E_\eta}|(f_j)_\gamma|^2\dd\gamma\dd t=0.
	\]
	The stationary assertion follows from the same argument.
\end{proof}

\begin{remark}[Compatibility of traces]\label{rem:trace-compatibility}
	The traces defined by completion in Definition~\ref{def:H_kin} are compatible
	with the full trace given by Proposition~\ref{trace_theorem}.
	
	Let \(f\in\cH_{\kin}(\cD_T\times\bR^d_v)\), and denote by the same
	symbol its continuous representative in
	\(C\big([0,T];L^2(\cD\times\bR^d_v)\big)\). Then \(f(0,\cdot)\) coincides with the
	initial trace defined in Definition~\ref{def:H_kin}(1). Moreover, if
	\(f_\gamma\in L^2(\gamma_T)\) is the full phase-space trace given by
	Proposition~\ref{trace_theorem}, then
	\[
	f_\gamma|_{\gamma_T^-}=f_{\gamma^-}
	\quad\text{in }L^2(\gamma_T^-).
	\]
	Indeed, for any approximating sequence
	\((f_n)_n\subset C_c^\infty([0,T]\times\overline{\cD}\times\bR^d_v)\),
	Proposition~\ref{trace_theorem} gives
	\[
	f_n\to f
	\quad\text{in } C([0,T];L^2(\cD\times\bR^d_v)),
	\qquad
	(f_n)_\gamma\to f_\gamma
	\quad\text{in }L^2(\gamma_T),
	\]
	whereas Definition~\ref{def:H_kin}(1) gives
	\[
	f_n(0,\cdot)\to f(0,\cdot),
	\qquad
	f_n|_{\gamma_T^-}\to f_{\gamma^-}.
	\]
	
	The stationary case is identical: if
	\(f\in\cH_{\kin}(\cD\times\bR^d_v)\) and
	\(f_\gamma\in L^2(\gamma)\) is the full trace, then
	\[
	f_\gamma|_{\gamma^-}=f_{\gamma^-}
	\quad\text{in }L^2(\gamma^-).
	\]
\end{remark}

\begin{lemma}
	\label{lem:Hkin-operations}
	\,
	
	\begin{enumerate}
		\item[\textup{(i)}] \(C_c^2([0,T]\times\overline\cD\times\bR^d_v)\subset
		\cH_{\kin}(\cD_T\times\bR^d_v)\), and for such functions the objects
		of Definition~\ref{def:H_kin} and the trace of
		Proposition~\ref{trace_theorem} coincide with the classical ones.
		\item[\textup{(ii)}] For \(\zeta\in C_c^2([0,T]\times\overline\cD\times\bR^d_v)\)
		and \(f\in\cH_{\kin}(\cD_T\times\bR^d_v)\) one has
		\(\zeta f\in\cH_{\kin}(\cD_T\times\bR^d_v)\),
		\[
		\|\zeta f\|_{\cH_{\kin}(\cD_T\times\bR^d_v)}
		\le
		C(\zeta)\|f\|_{\cH_{\kin}(\cD_T\times\bR^d_v)},
		\]
		and \((\zeta f)_\gamma=\zeta|_{\gamma_T}f_\gamma\),
		\((\zeta f)(t,\cdot)=\zeta(t,\cdot)f(t,\cdot)\) for every \(t\in[0,T]\).
		\item[\textup{(iii)}] For \(f\in\cH_{\kin}(\cD_T\times\bR^d_v)\), the
		reflection \(f^\ast(t,x,v):=f(T-t,x,-v)\) again belongs to
		\(\cH_{\kin}(\cD_T\times\bR^d_v)\), with
		\[
		(\p_t+v\cdot\nabla_x)f^\ast=-\big((\p_t+v\cdot\nabla_x)f\big)^\ast,
		\qquad
		f^\ast(t,\cdot)=f(T-t,\cdot,-\,\cdot\,),
		\qquad
		(f^\ast)_\gamma=(f_\gamma)^\ast ,
		\]
		and \(\gamma^\pm\) is exchanged with \(\gamma^\mp\).
		Moreover, there exists a universal constant \(C>0\) such that
		\[\|f^\ast\|_{\cH_{\kin}(\cD_T\times\bR^d_v)}
		\le C\|f\|_{\cH_{\kin}(\cD_T\times\bR^d_v)}.\]
	\end{enumerate}
	The stationary analogues hold as well, with \(v\cdot\nabla_x\) in place of
	\(\p_t+v\cdot\nabla_x\) and with \(f^\ast(x,v):=f(x,-v)\) in (iii).
\end{lemma}

\begin{proof}
	(i) Mollifying \(\zeta\in C_c^2\) in all variables produces
	\(\zeta_\varepsilon\in C_c^\infty(\bR^{1+2d})\) with
	\(\zeta_\varepsilon\to\zeta\) in \(C^1\),
	\(\sup_\varepsilon\|\nabla^2\zeta_\varepsilon\|_\infty<\infty\), and all
	supports contained in one compact set. Restricting to
	\([0,T]\times\overline\cD\times\bR^d_v\), we obtain a sequence which is
	Cauchy in \(\|\cdot\|_{\cH_{\kin}}\) and converges to \(\zeta\) in
	\(L^2(\cD_T\times\bR^d_v)\).
	
	(ii) For \(f\in C_c^\infty([0,T]\times\overline\cD\times\bR^d_v)\),
	\[
	(\p_t+v\cdot\nabla_x)(\zeta f)
	=
	\zeta\,(\p_t+v\cdot\nabla_x)f+f\,(\p_t\zeta+v\cdot\nabla_x\zeta).
	\]
	Multiplication by \(\zeta\) is bounded on \(H^1_v\), hence, by duality, on
	\(H^{-1}_v\), with norm at most \(\|\zeta\|_\infty+\|\nabla_v\zeta\|_\infty\).
	Moreover, \(\p_t\zeta+v\cdot\nabla_x\zeta\) is bounded, since \(\zeta\) has
	compact support in \(v\). Together with the corresponding bounds for the two
	trace terms, this gives
	\(\|\zeta f\|_{\cH_{\kin}}\le C(\zeta)\|f\|_{\cH_{\kin}}\)
	for smooth \(f\). For general \(f\in\cH_{\kin}\), pick
	\(f_n\in C_c^\infty([0,T]\times\overline\cD\times\bR^d_v)\) with
	\(f_n\to f\) in \(\cH_{\kin}\). Then
	\(\zeta f_n\in C_c^2([0,T]\times\overline\cD\times\bR^d_v)
	\subset\cH_{\kin}\) by (i), and the bound just proved,
	applied to \(f_n-f_{n'}\), shows that \((\zeta f_n)_n\) is Cauchy in
	\(\cH_{\kin}\). Since \(\cH_{\kin}\) is complete by
	Proposition~\ref{260128509} and \(\zeta f_n\to\zeta f\) in
	\(L^2(\cD_T\times\bR^d_v)\), we get \(\zeta f\in\cH_{\kin}\)
	together with the stated bound. The identification of the traces follows from
	the same passage to the limit.
	
	(iii) The map \(f\mapsto f^\ast\) preserves
	\(C_c^\infty([0,T]\times\overline\cD\times\bR^d_v)\), the measures
	\(\dd t\dd x\dd v\) and \(\dd\gamma\dd t\) are invariant under
	\((t,v)\mapsto(T-t,-v)\), and \(v\mapsto-v\) is an isometry of \(H^{\pm1}_v\).
	Hence, for smooth \(f\),
	\[
	\|f^\ast\|^2_{\cH_{\kin}}
	=
	\|f\|^2_{L^2(\cD_T;H^1_v)}
	+\|(\p_t+v\cdot\nabla_x)f\|^2_{L^2(\cD_T;H^{-1}_v)}
	+\|f(T,\cdot)\|^2_{L^2(\cD\times\bR^d_v)}
	+\|f_\gamma\|^2_{L^2(\gamma^+_T)},
	\]
	which is bounded by \(C\|f\|^2_{\cH_{\kin}}\) by
	Proposition~\ref{trace_theorem}. Applying this to a defining sequence for
	\(f\) gives \(f^\ast\in\cH_{\kin}(\cD_T\times\bR^d_v)\) together with
	the stated identities. The stationary case is identical, \(\dd\gamma\) being
	invariant under \(v\mapsto-v\).
\end{proof}

\begin{prop}[Renormalized identities in \(\cH_{\kin}\)]
	\label{prop:renormalized-identity}
	Let \(\beta:\bR\to\bR\) be a twice continuously differentiable function satisfying
	\[
	\beta(0)=\beta'(0)=0,
	\qquad
	\|\beta''\|_{L^\infty(\bR)}<\infty.
	\]
	\begin{enumerate}
		\item
		Let \(f\in \cH_{\kin}(\cD_T\times \bR^d_v)\), and let
		\(g^0,g^1,\ldots,g^d\in L^2(\cD_T\times \bR^d_v)\) satisfy
		\[
		(\partial_t+v\cdot \nabla_x)f
		=
		g^0+D_i g^i
		\quad\text{in}\quad L^2(\cD_T;H_v^{-1})\,.
		\]
		Then, for every
		\(\varphi\in C_c^2([0,T]\times \overline{\cD}\times \bR^d_v)\)
		and every \(0\le s<t\le T\), we have
		\begin{align}\label{reno_dyn}
			\begin{aligned}
				&\iint_{\cD\times \bR^d_v}
				\Big(
				\beta\big(f(t,\cdot)\big)\varphi(t,\cdot)
				-\beta\big(f(s,\cdot)\big)\varphi(s,\cdot)
				\Big)\dd x\dd v \\
				&\quad
				+\iint_{(s,t)\times \gamma}
				\beta(f_{\gamma})\,\varphi \,\operatorname{sgn}(n_x\cdot v)\dd\gamma\dd r \\
				&\quad
				-\iiint_{(s,t)\times \cD\times \bR^d_v}
				\beta(f)\big(\partial_t\varphi +v\cdot\nabla_x \varphi\big)\dd v\dd x\dd r \\
				=\,&
				\iiint_{(s,t)\times \cD\times \bR^d_v}
				\left[
				\beta'(f)g^0\varphi
				-\Big(
				\beta''(f)(D_i f)g^i\varphi
				+\beta'(f)g^iD_i\varphi
				\Big)
				\right]\dd v\dd x\dd r\, .
			\end{aligned}
		\end{align}
		
		\item
		Let \(f\in\cH_{\kin}(\cD\times \bR^d_v)\), and let
		\(g^0,g^1,\ldots,g^d\in L^2(\cD\times \bR^d_v)\) satisfy
		\begin{align}\label{2604301138}
			v\cdot \nabla_x f
			=
			g^0+D_i g^i
			\quad\text{in}\quad L^2(\cD;H_v^{-1})\,.
		\end{align}
		Then, for every \(\varphi\in C_c^2(\overline{\cD}\times \bR^d_v)\), we have
		\begin{align}\label{reno_std}
			\begin{aligned}
				&\iint_{\gamma}
				\beta(f_\gamma)\varphi \,\operatorname{sgn}(n_x\cdot v)\dd\gamma
				-\iint_{\cD\times \bR^d_v}
				\beta(f)\big(v\cdot\nabla_x \varphi\big)\dd x\dd v\\
				=\,&
				\iint_{\cD\times \bR^d_v}
				\left[
				\beta'(f)g^0\varphi
				-\Big(
				\beta''(f)(D_i f)g^i\varphi
				+\beta'(f)g^iD_i\varphi
				\Big)
				\right]\dd x\dd v .
			\end{aligned}
		\end{align}
	\end{enumerate}
\end{prop}

\begin{proof}
	We prove \((1)\). The proof of \((2)\) is the same, with the time variable
	omitted.
	
	Let $f_n\in C_c^\infty([0,T]\times \overline{\cD}\times \bR^d_v)$ be a defining sequence for \(f\in \cH_{\kin}(\cD_T\times\bR^d_v)\).
	Then \(f_n\to f\) in \(\cH_{\kin}(\cD_T\times\bR^d_v)\). In particular,
	\[
	f_n\to f\quad\text{in }L^2(\cD_T;H^1_v),
	\qquad
	(\partial_t+v\cdot\nabla_x)f_n
	\to
	(\partial_t+v\cdot\nabla_x)f
	\quad\text{in }L^2(\cD_T;H^{-1}_v).
	\]
	Moreover, by Proposition~\ref{trace_theorem},
	\begin{align}\label{261020233}
		f_n\to f
		\quad\text{in }C([0,T];L^2(\cD\times\bR^d_v)),
		\qquad
		(f_n)_\gamma\to f_\gamma
		\quad\text{in }L^2(\gamma_T).
	\end{align}
	
	By $\beta'(0)=0$ and the bound on $\beta''$,
	\[
	|\beta(a)-\beta(b)|
	\le \frac12\|\beta''\|_\infty(|a|+|b|)|a-b|,
	\qquad a,b\in\bR.
	\]
	Combining this with \eqref{261020233}, we have
	\[
	\beta(f_n)\to\beta(f)
	\quad\text{in }C([0,T];L^1(\cD\times\bR_v^d)),
	\qquad
	\beta((f_n)_\gamma)\to\beta(f_\gamma)
	\quad\text{in }L^1(\gamma_T);
	\]
	no velocity cutoff at the grazing set is needed, since the trace measure
	is $\dd\gamma=|n_x\cdot v|\dd S_x\dd v$, and the
	last convergence holds in this full measure.
	In addition, since \(\beta'(0)=0\) and \(\beta''\in L^\infty(\bR)\), the Sobolev
	chain rule in the \(v\)-variable gives \(\beta'(u)\in L^2(\cD_T;H^1_v)\) with
	\begin{align}\label{2604301150}
		\nabla_v\beta'(u)=\beta''(u)\nabla_v u
		\qquad\text{for every }u\in L^2(\cD_T;H^1_v);
	\end{align}
	see, for instance, the proof of \cite[Theorem~4.4(ii)]{evans2025measure}.
	We claim that, in addition,
	\begin{align}\label{2608061020}
		\beta'(f_n)\to\beta'(f)
		\quad\text{in }L^2(\cD_T;H^1_v).
	\end{align}
	Since \(\beta''\in L^\infty(\bR)\), the function \(\beta'\) is Lipschitz, so
	that
	$$
	\|\beta'(f_n)-\beta'(f)\|_{L^2(\cD_T\times\bR^d_v)}
	\le\|\beta''\|_\infty\|f_n-f\|_{L^2(\cD_T\times\bR^d_v)}\to0.
	$$
	For the gradients, write
	\[
	\beta''(f_n)\nabla_vf_n-\beta''(f)\nabla_vf
	=
	\beta''(f_n)\big(\nabla_vf_n-\nabla_vf\big)
	+\big(\beta''(f_n)-\beta''(f)\big)\nabla_vf .
	\]
	The first term is bounded in \(L^2\) by
	\(\|\beta''\|_\infty\|\nabla_vf_n-\nabla_vf\|_{L^2}\to0\). For the second,
	let \((n_j)_j\) be an arbitrary subsequence and extract a further subsequence
	along which \(f_{n_j}\to f\) almost everywhere; since \(\beta''\) is
	continuous and bounded, dominated convergence with dominating function
	\(2\|\beta''\|_\infty|\nabla_vf|\in L^2\) gives convergence to \(0\) along
	that further subsequence. As every subsequence thus has a further subsequence
	converging to the same limit, the full sequence converges, which proves
	\eqref{2608061020}.
	
	We now fix $\varphi\in C_c^2([0,T]\times\overline{\cD}\times\bR^d_v)$ and $0\le s<t\le T$.
	For smooth \(f_n\), the classical chain rule and integration by parts yield
	\begin{align*}
		&\iint_{\cD\times \bR^d_v}
		\Big(
		\beta(f_n(t ,\cdot))\varphi(t ,\cdot)
		-\beta(f_n(s ,\cdot))\varphi(s ,\cdot)
		\Big)\dd x\dd v \\
		&\quad
		+\iint_{(s,t)\times \gamma}
		\beta\big((f_n)_\gamma\big)\varphi\,\operatorname{sgn}(n_x\cdot v)\dd\gamma\dd r
		-\iiint_{(s ,t )\times\cD\times\bR^d_v}
		\beta(f_n)(\partial_t\varphi+v\cdot\nabla_x\varphi)\dd v\dd x\dd r \\
		=\,&
		\iint_{(s,t)\times \cD}
		\left\langle
		(\partial_t+v\cdot\nabla_x)f_n,\beta'(f_n)\varphi
		\right\rangle_v
		\dd x\dd r .
	\end{align*}
	Passing to the limit as \(n\to\infty\), using the above convergences and the
	duality between \(H^{-1}_v\) and \(H^1_v\), we obtain the same identity with
	\(f\) in place of \(f_n\), namely with right-hand side
	\[
	\iint_{(s,t)\times \cD}
	\left\langle
	(\partial_t+v\cdot\nabla_x)f,\beta'(f)\varphi
	\right\rangle_v
	\dd x\dd r .
	\]
	Finally, using
	\[
	(\partial_t+v\cdot\nabla_x)f
	=
	g^0+D_i g^i
	\quad\text{in }L^2(\cD_T;H^{-1}_v),
	\]
	together with \eqref{2604301150}, we rewrite this term as
	\begin{align*}
		&\iiint_{(s,t)\times\cD\times\bR^d_v}
		\left(
		g^0\beta'(f)\varphi
		- g^iD_i(\beta'(f)\varphi)
		\right)\dd v\dd x\dd r \\
		=\,&
		\iiint_{(s,t)\times\cD\times\bR^d_v}
		\left(
		\beta'(f)g^0\varphi
		-\Big(
		\beta''(f)(D_i f)g^i\varphi
		+\beta'(f)g^iD_i\varphi
		\Big)
		\right)\dd v\dd x\dd r.
	\end{align*}
	     This proves \eqref{reno_dyn}.
\end{proof}

The trace bounds below are obtained from the preceding renormalized identities by applying them to suitable truncations.

\begin{prop}\label{2606291053}
	Let \(\cO\neq\emptyset\) be a relatively open subset of
	\([0,T]\times\overline{\cD}\times\bR^d_v\). For \(t\in[0,T]\), set
	\[
	\cO_t:=\{(x,v)\in\cD\times\bR^d_v:(t,x,v)\in\cO\}.
	\]
	Then, for any \(f\in \cH_{\kin}(\cD_T\times\bR^d_v)\), if
	\[
	\operatorname*{ess\,sup}_{\cO} f<+\infty,
	\]
	then
	\begin{alignat}{2}\label{trace-int}
		&\operatorname*{ess\,sup}_{\gamma_T\cap\cO,\dd\gamma\dd t} f_\gamma
		\le
		\operatorname*{ess\,sup}_{\cO} f,
		\qquad
		&&
		\sup_{0\le t\le T}
		\left(
		\operatorname*{ess\,sup}_{\cO_t} f(t,\cdot)
		\right)
		\le
		\operatorname*{ess\,sup}_{\cO} f.
	\end{alignat}
	Similarly, if
	\[
	\operatorname*{ess\,inf}_{\cO} f>-\infty,
	\]
	then
	\[
	\operatorname*{ess\,inf}_{\gamma_T\cap\cO,\dd\gamma\dd t} f_\gamma
	\ge
	\operatorname*{ess\,inf}_{\cO} f,
	\qquad 
	\inf_{0\le t\le T}
	\left(
	\operatorname*{ess\,inf}_{\cO_t} f(t,\cdot)
	\right)
	\ge
	\operatorname*{ess\,inf}_{\cO} f.
	\]
	Here \(f_\gamma\) denotes the trace, and \(f(t,\cdot)\) is understood through the
	representative of \(f\) in
	\(C\big([0,T];L^2(\cD\times\bR^d_v)\big)\),
	which is provided by Proposition~\ref{trace_theorem}.
	
	In particular,
	\[
	\|f_\gamma\|_{L^\infty(\gamma_T\cap\cO,\dd\gamma\dd t)}
	\le
	\|f\|_{L^\infty(\cO)}
	\]
	and
	\[
	\sup_{0\le t\le T}
	\left\|f(t,\cdot)1_{\cO_t}\right\|_{L^\infty(\cD\times\bR^d_v)}
	\le
	\|f\|_{L^\infty(\cO)}.
	\]
	
	The stationary case is analogous: if \(\widetilde{\cO}\neq\emptyset\) is a relatively
	open subset of \(\overline\cD\times\bR^d_v\) and
	\(f\in\cH_{\kin}(\cD\times\bR^d_v)\), then
	\begin{align}\label{trace-int-st}
		\operatorname*{ess\,sup}_{\gamma\cap \widetilde{\cO},\dd\gamma} f_\gamma
		\le
		\operatorname*{ess\,sup}_{\widetilde{\cO}} f,
		\qquad
		\operatorname*{ess\,inf}_{\gamma\cap\widetilde{\cO},\dd\gamma} f_\gamma
		\ge
		\operatorname*{ess\,inf}_{\widetilde{\cO}} f .
	\end{align}
\end{prop}
\begin{proof}
	By a standard exhaustion argument, using
	\begin{align}\label{2606301225}
		E_n\uparrow E
		\quad\Longrightarrow\quad
		\operatorname*{ess\,sup}_{E_n}F\uparrow
		\operatorname*{ess\,sup}_{E}F,
	\end{align}
	we may assume throughout the proof that \(\cO\) is bounded.
	
	\medskip
	\noindent
	\emph{Step 1. An auxiliary estimate.}
	We first prove that
	\begin{align}\label{260317433}
		\begin{aligned}
			&\operatorname*{ess\,sup}_{\gamma_T^+\cap \cO,\dd\gamma\dd t} f_\gamma
			\vee
			\sup_{0\le t\le T}
			\left(\operatorname*{ess\,sup}_{\cO_t} f(t,\cdot)\right)
			\\
			&\qquad\le
			\operatorname*{ess\,sup}_{\cO} f
			\vee
			\operatorname*{ess\,sup}_{\gamma_T^-\cap \cO,\dd\gamma\dd t} f_\gamma
			\vee
			\operatorname*{ess\,sup}_{\cO_0} f(0,\cdot).
		\end{aligned}
	\end{align}
	
	Let
	\[
	M:=\operatorname*{ess\,sup}_{\cO}f.
	\]
	If \(M=+\infty\), there is nothing to prove. Otherwise, choose
	\(\psi_\infty\in C_c^\infty([0,T]\times\overline\cD\times\bR_v^d)\) such that
	\(\psi_\infty\equiv1\) on a neighborhood of \(\overline\cO\), and set
	\[
	h:=f-M\psi_\infty,
	\]
	which lies in \(\cH_{\kin}(\cD_T\times\bR^d_v)\) by Lemma~\ref{lem:Hkin-operations}(i).
	Then
	\[
	h=f-M \quad\text{a.e. on }\cO,
	\qquad
	h_\gamma=f_\gamma-M \quad\text{a.e. on }\gamma_T\cap\cO,\qquad \operatorname*{ess\,sup}_{\cO}h=0.
	\]
	Thus \eqref{260317433} is reduced to proving
	\begin{align}\label{260317430}
		\begin{aligned}
			&\operatorname*{ess\,sup}_{\gamma_T^+\cap \cO,\dd\gamma\dd t} h_\gamma
			\vee
			\sup_{0\le t\le T}
			\left(\operatorname*{ess\,sup}_{\cO_t} h(t,\cdot)\right)
			\vee 0
			\\
			&\qquad\le
			\operatorname*{ess\,sup}_{\gamma_T^-\cap \cO,\dd\gamma\dd t} h_\gamma
			\vee
			\operatorname*{ess\,sup}_{\cO_0} h(0,\cdot)
			\vee 0.
		\end{aligned}
	\end{align}
	
	By the decomposition recorded in Remark~\ref{weak_formulation}, choose
	\(g_0,\ldots,g_d\in L^2(\cD_T\times\bR_v^d)\) such that
	\[
	(\partial_t+v\cdot\nabla_x)f
	=
	g_0+\sum_{i=1}^d D_i g_i
	\qquad\text{in }L^2(\cD_T;H_v^{-1}).
	\]
	Then
	\[
	(\partial_t+v\cdot\nabla_x)h
	=
	\tilde g_0+\sum_{i=1}^d D_i\tilde g_i
	\qquad\text{in }L^2(\cD_T;H_v^{-1}),
	\]
	where
	\[
	\tilde g_0:=g_0-M(\partial_t+v\cdot\nabla_x)\psi_\infty,
	\qquad
	\tilde g_i:=g_i,\quad i=1,\dots,d.
	\]
	
	Fix \(p>2\). Let \(\beta_n\in C^2(\bR)\) be convex functions such that
	\[
	\beta_n(0)=\beta_n'(0)=0,
	\qquad
	\beta_n''\in C_c^\infty(\bR_+),
	\]
	and
	\[
	0\le \beta_n(r)\uparrow r_+^p,
	\qquad
	0\le \beta_n'(r)\uparrow p r_+^{p-1},
	\qquad
	0\le \beta_n''(r)\uparrow p(p-1)r_+^{p-2}
	\]
	for every \(r\in\bR\).
	
	Let \(\{\cO_m\}_{m\ge1}\) be a bounded relatively open exhaustion of \(\cO\)
	such that \(\overline{\cO_m}\subset\cO\), and choose
	\(\psi_m\in C_c^\infty([0,T]\times\overline\cD\times\bR_v^d)\) satisfying
	\[
	0\le \psi_m\le1,
	\qquad
	\psi_m\equiv1\quad\text{on }\overline{\cO_m},
	\qquad
	\supp\,\psi_m\subset\cO.
	\]
	Applying the renormalized identity \eqref{reno_dyn} to \(h\) with
	\((\beta,\varphi)=(\beta_n,\psi_m)\) on \((0,t)\), we obtain, for every
	\(t\in(0,T]\),
	\begin{align*}
		&\iint_{\cD\times\bR_v^d}\beta_n(h(t,\cdot))\psi_m(t,\cdot)\dd x\dd v
		+
		\int_0^t\int_{\gamma^+}\beta_n(h_\gamma)\psi_m\dd\gamma\dd \tau
		\\
		&\qquad\le
		\iint_{\cD\times\bR_v^d}\beta_n(h(0,\cdot))\psi_m(0,\cdot)\dd x\dd v
		+
		\int_0^t\int_{\gamma^-}\beta_n(h_\gamma)\psi_m\dd\gamma\dd \tau .
	\end{align*}
	Indeed, all interior terms vanish because
	$h\le0$ almost everywhere on $\supp\psi_m\subset\cO$.
	
	Using \(\psi_m\equiv1\) on \(\cO_m\) and letting \(n\to\infty\), monotone
	convergence gives, for every \(t\in(0,T]\),
	\begin{align}\label{260317920}
		\begin{aligned}
			&\left(
			\iint_{(\cO_m)_t}|h_+(t,\cdot)|^p\dd x\dd v\right)^{1/p}
			\vee
			\left(\iint_{\gamma_t^+\cap\cO_m}|(h_\gamma)_+|^p\dd\gamma\dd r\right)^{1/p}
			\\
			&\qquad\le
			\left(\iint_{\cO_0}|h_+(0,\cdot)|^p\dd x\dd v
			+
			\iint_{\gamma_T^-\cap\cO}|(h_\gamma)_+|^p\dd\gamma\dd t \right)^{1/p}.
		\end{aligned}
	\end{align}
	For \(t=0\) the same inequality holds trivially, since
    \((\cO_m)_0\subset\cO_0\).
    Since \(\cO\) is bounded and \(|\partial\cD|<\infty\), all measures involved
    are finite. Letting \(p\to\infty\) in
    \eqref{260317920} at fixed \(t\), and then taking the supremum over
    \(t\in[0,T]\) in the first term and \(t=T\) in the second, we obtain
	\begin{align*}
		\begin{aligned}
			&\sup_{0\le t\le T}
			\left(
			\operatorname*{ess\,sup}_{(\cO_m)_t} h(t,\cdot)
			\right)
			\vee
			\operatorname*{ess\,sup}_{\gamma_T^+\cap\cO_m,\dd\gamma\dd t} h_\gamma
			\vee 0\le
			\operatorname*{ess\,sup}_{\cO_0} h(0,\cdot)
			\vee
			\operatorname*{ess\,sup}_{\gamma_T^-\cap\cO,\dd\gamma\dd t} h_\gamma
			\vee 0.
		\end{aligned}
	\end{align*}
	Letting \(m\to\infty\), we obtain \eqref{260317430}, and hence
	\eqref{260317433}.
	
	\medskip
	\noindent
	\emph{Step 2. The outgoing trace and the time slices.}
	We now derive the desired upper bounds from \eqref{260317433}. Put
	\[
	S:=\operatorname*{ess\,sup}_{\cO}f.
	\]
	If \(S=+\infty\), there is nothing to prove. Otherwise, after replacing \(f\)
	by
	$f-S\psi_\infty$
	with \(\psi_\infty\equiv1\) near \(\overline\cO\), we may assume that
	\begin{align}\label{2606291114}
		\operatorname*{ess\,sup}_{\cO}f=0.
	\end{align}
	Let
	\[
	\cO^\ast:=
	\cO\setminus
	\Big(
	(\{0\}\times\overline\cD\times\bR_v^d
	)\cup
	([0,T]\times(\gamma^-\cup\gamma^0))
	\Big).
	\]
	Choose compact sets \(K_n\uparrow\cO^\ast\) and functions
	\(\zeta_n\in C_c^\infty([0,T]\times\overline\cD\times\bR_v^d)\) with
	\[
	0\le \zeta_n\le1,
	\qquad
	\zeta_n\equiv1\quad\text{on }K_n,
	\qquad
	\supp\zeta_n\subset\cO^\ast.
	\]
	By Lemma~\ref{lem:Hkin-operations}(ii), \(f\zeta_n\in\cH_{\kin}(\cD_T\times\bR^d_v)\),
	with \((f\zeta_n)_\gamma=\zeta_nf_\gamma\) and
	\((f\zeta_n)(t,\cdot)=\zeta_n(t,\cdot)f(t,\cdot)\). We apply
	\eqref{260317433} to \(f\zeta_n\). The right-hand side of
	\eqref{260317433} is then nonpositive. The
	first term is \(\le0\), because \(0\le\zeta_n\le1\) and
	\(f\le0\) a.e.\ on \(\cO\) by \eqref{2606291114}. The second and the third
	vanish, because \(\supp\,\zeta_n\) avoids both
	\(\{0\}\times\overline\cD\times\bR_v^d\)
	and \([0,T]\times(\gamma^-\cup\gamma^0)\). Hence
	\[
	\operatorname*{ess\,sup}_{\gamma_T^+\cap K_n,\dd\gamma\dd t} f_\gamma
	\vee
	\sup_{0<t\le T}
	\left(
	\operatorname*{ess\,sup}_{(K_n)_t}f(t,\cdot)
	\right)
	\le0.
	\]
    	Letting \(n\to\infty\), and using that \((\cO^\ast)_t=\cO_t\) for every
	\(t>0\), since \(\cO_t\subset\cD\times\bR^d_v\) is disjoint from
	\(\gamma^-\cup\gamma^0\), the last inequality holds with \(\cO\) in place of \(K_n\).
	Undoing the normalization \eqref{2606291114}, we obtain
	\begin{align}\label{260317523}
		\operatorname*{ess\,sup}_{\gamma_T^+\cap\cO,\dd\gamma\dd t} f_\gamma
		\vee
		\sup_{0<t\le T}
		\left(
		\operatorname*{ess\,sup}_{\cO_t}f(t,\cdot)
		\right)
		\le
		\operatorname*{ess\,sup}_{\cO}f.
	\end{align}

	\medskip
	\noindent
	\emph{Step 3. Time reversal and conclusion.}
	We apply \eqref{260317523} to
	\[
	\widehat f(t,x,v):=f(T-t,x,-v),
	\qquad
	\widehat\cO:=\{(t,x,v):(T-t,x,-v)\in\cO\}.
	\]
	By Lemma~\ref{lem:Hkin-operations}(iii), \(\widehat f=f^\ast\) belongs to
	\(\cH_{\kin}(\cD_T\times\bR^d_v)\), its kinetic trace is
	\[
	\widehat f_\gamma(t,x,v)=f_\gamma(T-t,x,-v),
	\]
	and \(\gamma^+\) is mapped onto \(\gamma^-\). Moreover, \(\dd t\dd x\dd v\)
	and \(\dd\gamma\dd t\) are invariant under \((t,v)\mapsto(T-t,-v)\), so that
	\(\operatorname*{ess\,sup}_{\widehat\cO}\widehat f
	=\operatorname*{ess\,sup}_{\cO}f\). Therefore
		\[
	\operatorname*{ess\,sup}_{\gamma_T^-\cap\cO,\dd\gamma\dd t} f_\gamma
	\vee
	\sup_{0\le t<T}
	\left(
	\operatorname*{ess\,sup}_{\cO_t}f(t,\cdot)
	\right)
	\le
	\operatorname*{ess\,sup}_{\cO}f.
	\]
	Together with \eqref{260317523}, this proves \eqref{trace-int}.
	
	Applying the same result to \(-f\), we obtain the
	\(\operatorname*{ess\,inf}\) assertions.
	The stated \(L^\infty\) estimates follow immediately.
	
	\medskip
	\noindent
	\emph{Step 4. The stationary case.}
	It remains to prove \eqref{trace-int-st}. Let
	\(f\in\cH_{\kin}(\cD\times\bR^d_v)\), and define
	\[
	F(t,x,v):=f(x,v)
	\qquad
	\text{for }(t,x,v)\in[0,T]\times\overline\cD\times\bR_v^d.
	\]
	If \((f_n)\subset C_c^\infty(\overline\cD\times\bR^d_v)\) is a defining
	sequence for \(f\), then \(F_n(t,x,v):=f_n(x,v)\) lies in
	\(C_c^\infty([0,T]\times\overline\cD\times\bR^d_v)\) and
	\[
	\|F_n-F_m\|_{\cH_{\kin}(\cD_T\times\bR^d_v)}^2
	\le
	(1+T)\|f_n-f_m\|_{\cH_{\kin}(\cD\times\bR^d_v)}^2 ,
	\]
	since \((\p_t+v\cdot\nabla_x)F_n=v\cdot\nabla_xf_n\), \(F_n(0,\cdot)=f_n\) and
	\(\|f_n\|_{L^2(\cD\times\bR^d_v)}\le\|f_n\|_{L^2(\cD;H^1_v)}\). Hence
	\(F\in\cH_{\kin}(\cD_T\times\bR^d_v)\) and
	\[
	F_\gamma(t,x,v)=f_\gamma(x,v)
	\qquad
	\text{for a.e. }(t,x,v)\in\gamma_T.
	\]
	Applying the result just proved with
	\(\widetilde\cO_T:=(0,T)\times \widetilde{\cO}\)
	gives \eqref{trace-int-st}, and the corresponding \(L^\infty\) bound follows
	immediately.
\end{proof}

\subsection{Exit times of the Kolmogorov diffusion and nonattainability of the grazing set}\label{subsec:grazing}

We establish almost sure finiteness of the exit time
(Theorem~\ref{thm:finiteness_exit_time}), nonattainability of \(\gamma^0\),
immediate exit from \(\gamma^+\cup\gamma^0\), and transversality at exit
(Theorem~\ref{thm:exit_time_properties}), continuity in
the initial point (Corollary~\ref{cor:exit_time_indicator_continuity}), and
stability under vanishing elliptic regularization
(Corollary~\ref{cor:exit_time_eps_stability}).
Versions of these results appear in \cite{2022LRR}; we give self-contained proofs in the present setting.
These properties are used in Subsection~\ref{subsec:prob-solutions}; in particular, the two corollaries
justify the limiting passage in the Feynman--Kac representation.
The proofs are given in \S\ref{subsubsec:grazing-proofs}.

Recall that \(\cD\) is a bounded \(C^{1,1}\) domain. Throughout this subsection, let
\((X_t(x,v),V_t(x,v))\) be the Kolmogorov diffusion in
\(\bR^d_x\times\bR^d_v\) given by
\begin{equation}\label{eq:kolmogorov_XV}
	X_t(x,v)=x+\int_0^t V_s(x,v)\,ds,\qquad
	V_t(x,v)=v+\sqrt{2}\,B_t,
\end{equation}
where \(B_t\) is a fixed standard \(d\)-dimensional Brownian motion.
The backward generator acting on smooth functions \(\varphi=\varphi(x,v)\) is
\[
\mathscr{L}\varphi
=
v\cdot\nabla_x\varphi+\Delta_v\varphi .
\]

We also define the exit time of \(X_t(x,v)\) from \(\cD\) by
\begin{align}\label{260619759}
	\tau_\cD(x,v):=\inf\{t>0:\ X_t(x,v)\notin\cD\}.
\end{align}
When the initial point \((x,v)\) is fixed, we simply write \(X_t\), \(V_t\), and
\(\tau_\cD\) for \(X_t(x,v)\), \(V_t(x,v)\), and \(\tau_\cD(x,v)\), respectively.
By Lemma~\ref{lem:hitting_time}, \(\tau_\cD(x,v)\) is a stopping time.

\begin{thm}\label{thm:finiteness_exit_time}
	For any $(x,v)\in\overline\cD\times \bR^d_v$, $\tau_\cD(x,v)<\infty$ almost surely.
\end{thm}

\begin{thm}\label{thm:exit_time_properties}
	For any \((x,v)\in \overline{\cD}\times\bR^d_v\), the following hold.
	
	\begin{enumerate}
		\item[(1)] (\emph{Non-attainability of \(\gamma^0\).})
		If \((x,v)\notin\gamma^0\), then
		\[
		\mathbb{P}\big(\exists t>0:\ (X_t,V_t)\in\gamma^0\big)=0.
		\]
		
		\item[(2)] (\emph{Immediate exit from \(\gamma^+\cup\gamma^0\).})
		If \((x,v)\in\gamma^+\cup\gamma^0\), then \(\tau_\cD=0\) almost surely.
		
		\item[(3)] (\emph{Attainability of the boundary from the interior.})
		If \(x\in\cD\) or \((x,v)\in\gamma^-\), then \(\tau_\cD>0\) almost surely.
		Moreover, on \(\{\tau_\cD<\infty\}\),
		\[
		V_{\tau_\cD}\cdot n(X_{\tau_\cD})>0
		\quad\text{almost surely.}
		\]
	\end{enumerate}
\end{thm}

\begin{corollary}\label{cor:exit_time_indicator_continuity}\,
	
	\begin{enumerate}
		\item[\textup{(i)}] Let \((x,v)\in\overline{\cD}\times\bR^d\), and let
		\((x_n,v_n)\in\overline{\cD}\times\bR^d\) be any sequence converging
		to \((x,v)\). Then
		\begin{equation}\label{eq:tau_pathwise_conv}
			\tau_{\cD}(x_n,v_n)
			\longrightarrow
			\tau_{\cD}(x,v)
			\qquad\text{almost surely.}
		\end{equation}
		
		\item[\textup{(ii)}] Let \((x,v)\) and \((x_n,v_n)\) be as in (i),
		and assume in addition that
		\[
		(t,x,v)\in
		\big([0,\infty)\times \overline{\cD}\times\bR^d\big)
		\setminus
		\big(\{0\}\times(\gamma^+\cup\gamma^0)\big),
		\]
		and let \(t_n\in[0,\infty)\) satisfy \(t_n\to t\). Then
		\begin{equation}\label{eq:indicator_pathwise_conv}
			\mathbf 1_{\{\tau_{\cD}(x_n,v_n)>t_n\}}
			\longrightarrow
			\mathbf 1_{\{\tau_{\cD}(x,v)>t\}}
			\qquad\text{almost surely.}
		\end{equation}
	\end{enumerate}
\end{corollary}

\begin{corollary}\label{cor:exit_time_eps_stability}
	Let \((x,v)\in \cD\times\bR^d\), and let \((X_t,V_t)\) be as in
	\eqref{eq:kolmogorov_XV}. Let \(B^x\) be a \(d\)-dimensional Brownian motion
	independent of the Brownian motion \(B\) in \eqref{eq:kolmogorov_XV}.
	For \(\varepsilon>0\), define
	\[
	X_t^{\varepsilon}:=X_t+\sqrt{2\varepsilon}\,B_t^x,
	\qquad
	\tau_{\cD}^{\varepsilon}:=\inf\{t>0:\ X_t^{\varepsilon}\notin\cD\}.
	\]
	Then $\tau_{\cD}^{\varepsilon}
	\longrightarrow
	\tau_{\cD}$
	almost surely, as $\varepsilon \downarrow 0$.
	Moreover, for every fixed \(t>0\),
	\begin{align}\label{260621441}
		{\bf 1}_{\{\tau_{\cD}^{\varepsilon}>t\}}
		\longrightarrow
		{\bf 1}_{\{\tau_{\cD}>t\}}
		\qquad\text{almost surely as }\varepsilon \downarrow 0.
	\end{align}
\end{corollary}

\subsubsection{The 1D Kolmogorov operator and a homogeneous gauge}\label{subsubsec:gauge}

The kernel and resolvent estimates below produce the homogeneous barrier
used to prove the bounded-drift polarity statements and exclude grazing
hits.

In $\bR^2$ with coordinates $z=(y,u)$, let
\[
\cL:=u\,\partial_y+\partial_{uu}.
\]
Define the anisotropic dilation and homogeneous quasi-norm
\[
\delta_\mu(y,u):=(\mu^3 y,\mu u)\quad(\mu>0),\qquad
N(y,u):=\big(u^4+|y|^{4/3}\big)^{1/4},
\]
so that $N(\delta_\mu z)=\mu N(z)$.

We record the following standard covariance computation for integrated
Brownian motion; see, for example, \cite{Kolmogorov1934}.

\begin{lemma}\label{lem:integrated_BM}
	Let $(B_t)_{t\ge0}$ be a standard $d$-dimensional Brownian motion, and set
	$J_t:=\int_0^tB_s\dd s$. Then, for every $t>0$, the pair $(J_t,B_t)$ is a
	centered Gaussian vector in $\bR^d\times\bR^d$ with
	\[
	\mathrm{Var}(B_t)=t\,I_d,
	\qquad
	\mathrm{Cov}(J_t,B_t)=\frac{t^2}{2}\,I_d,
	\qquad
	\mathrm{Var}(J_t)=\frac{t^3}{3}\,I_d .
	\]
	Its covariance matrix $\Sigma_t$ is nondegenerate, with
	$\det\Sigma_t=(t^4/12)^d$; consequently $(J_t,B_t)$ has the Gaussian density
	$(2\pi)^{-d}(\det\Sigma_t)^{-1/2}
	\exp\bigl(-\tfrac12 z\cdot\Sigma_t^{-1}z\bigr)$.
\end{lemma}

\begin{proof}
	Let \((\pi_n)_{n\ge1}\) be a sequence of partitions of \([0,t]\)
	whose mesh tends to zero, and let \(\kappa_n(s)\) denote the left
	endpoint of the interval of \(\pi_n\) containing \(s\). Define
	\[
	J_t^{(n)}
	:=
	\int_0^t B_{\kappa_n(s)}\dd s,
	\]
	so that \(J_t^{(n)}\) is a Riemann sum for \(J_t\).
	Since
	\[
	\mathbb E\bigl[|B_s-B_r|^2\bigr]
	=
	d|s-r|,
	\]
	the Cauchy--Schwarz inequality gives
	\[
	\begin{aligned}
		\mathbb E\bigl[|J_t^{(n)}-J_t|^2\bigr]
		&\le
		t\int_0^t
		\mathbb E\bigl[
		|B_{\kappa_n(s)}-B_s|^2
		\bigr]\dd s \le
		dt^2|\pi_n|
		\longrightarrow 0.
	\end{aligned}
	\]
	Thus \(J_t^{(n)}\to J_t\) in
	\(L^2(\Omega;\bR^d)\).
	
	Let \(\cH_B\) be the Gaussian Hilbert space generated by
	\[
	\{B_s^i:0\le s\le t,\ 1\le i\le d\}.
	\]
	Each component of \(J_t^{(n)}\) belongs to the linear span of this
	family. Hence each \(J_t^i\) belongs to \(\cH_B\) by the
	\(L^2\)-convergence above. It follows that every linear combination
	of the components of \(J_t\) and \(B_t\) is a centered Gaussian
	random variable. Therefore \((J_t,B_t)\) is a centered Gaussian
	vector; see
	\cite[Chapter~1, especially Example~1.10]{Janson1997}.
	
	The covariance identities follow from
	$\mathbb E[B_s^iB_r^j]
	=
	\delta_{ij}(s\wedge r)$
	and Fubini's theorem:
	\[
	\mathbb E[J_t^iB_t^j]
	=
	\delta_{ij}\int_0^t s\dd s
	=
	\delta_{ij}\frac{t^2}{2},
	\qquad
	\mathbb E[J_t^iJ_t^j]
	=
	\delta_{ij}
	\int_0^t\!\int_0^t(s\wedge r)\dd s\dd r
	=
	\delta_{ij}\frac{t^3}{3}.
	\]
	After reordering the coordinates into the \(d\) pairs
	\((J_t^i,B_t^i)\), the covariance matrix \(\Sigma_t\) is block
	diagonal with identical blocks
	\[
	\begin{pmatrix}
		t^3/3 & t^2/2\\
		t^2/2 & t
	\end{pmatrix},
	\]
	each having determinant \(t^4/12>0\). Consequently,
	\(\det\Sigma_t=(t^4/12)^d\), and the asserted Gaussian density
	follows.
\end{proof}

\begin{lemma}
	\label{lem:intBM_osc}
	Let \((\beta_t)_{t\ge0}\) be a standard one-dimensional Brownian motion, and set
	\[
	I_t:=\int_0^t \beta_s\dd s.
	\]
	Then, almost surely,
	\[
	\limsup_{t\downarrow0}\frac{I_t}{t^{3/2}}=+\infty
	\qquad\text{and}\qquad
	\liminf_{t\downarrow0}\frac{I_t}{t^{3/2}}=-\infty.
	\]
	Equivalently, almost surely, for every \(T>0\),
	\[
	\sup_{0<t\le T}\frac{I_t}{t^{3/2}}=+\infty
	\qquad\text{and}\qquad
	\inf_{0<t\le T}\frac{I_t}{t^{3/2}}=-\infty.
	\]
\end{lemma}

A sharper description follows from Strassen's functional law of the iterated logarithm \cite{strassen1964}. 
The version needed here, however, has a short elementary proof; for completeness, we include it.

\begin{proof}[\textbf{Proof of Lemma~\ref{lem:intBM_osc}}]
	We first prove that
	\[
	\limsup_{t\downarrow0}\frac{I_t}{t^{3/2}}=+\infty
	\qquad\text{almost surely.}
	\]
	
	Fix \(t_1\in(0,1/4)\), and define inductively
	\[
	t_{n+1}:=t_n^2,\qquad n\ge1.
	\]
	Then \(t_n\downarrow0\) super-exponentially. Set
	\[
	\Delta_n:=t_n-t_{n+1}=t_n-t_n^2.
	\]
	Since \(t_n\le 1/4\) for all \(n\), we have
	\[
	\frac12 t_n\le \Delta_n\le t_n.
	\]
	
	For each \(n\ge1\), decompose
	\[
	I_{t_n}
	=
	I_{t_{n+1}}+\int_{t_{n+1}}^{t_n}\beta_s\,ds
	=
	I_{t_{n+1}}+\beta_{t_{n+1}}\Delta_n+J_n,
	\]
	where
	\[
	J_n:=\int_{t_{n+1}}^{t_n}
	\bigl(\beta_s-\beta_{t_{n+1}}\bigr)\dd s.
	\]
	If we define
	\[
	\widetilde\beta^{(n)}_u
	:=
	\beta_{t_{n+1}+u}-\beta_{t_{n+1}},
	\qquad 0\le u\le \Delta_n,
	\]
	then, by Lemma~\ref{lem:brownian_shift}, \(\widetilde\beta^{(n)}\) is a
	standard Brownian motion independent of $\beta_{t_{n+1}}$ and $I_{t_{n+1}}$.
	In addition,
	\[
	J_n=\int_0^{\Delta_n}\widetilde\beta^{(n)}_u\dd u.
	\]
	Hence, by Lemma~\ref{lem:integrated_BM}, \(J_n\) is centered Gaussian with
	$\operatorname{Var}(J_n)=\frac{\Delta_n^3}{3}$.
	Moreover, the family \(\{J_n\}_{n\ge1}\) is independent, since each \(J_n\)
	depends only on the Brownian increment over the disjoint interval
	\((t_{n+1},t_n]\).
	
	Set
	\[
	\sigma_n^2:=\operatorname{Var}(J_n)=\frac{\Delta_n^3}{3}\geq \frac{t_n^3}{24},
	\]
	where the last inequality follows from \(\Delta_n\ge t_n/2\).
	
	\medskip
	\noindent
	\emph{Step 1.}
	Fix \(M>0\), and define
	\[
	E_n^+(M):=\{J_n\ge M\,t_n^{3/2}\}.
	\]
	If \(Z\sim N(0,1)\), then
	\[
	\mathbb P(E_n^+(M))
	=
	\mathbb P\!\left(Z\ge \frac{M t_n^{3/2}}{\sigma_n}\right)
	\ge
	\mathbb P\!\left(Z\ge \sqrt{24}\,M\right)
	=:p_M>0.
	\]
	Since the events \(\{E_n^+(M)\}_{n\ge1}\) are independent and
	\[
	\sum_{n=1}^\infty \mathbb P(E_n^+(M))=\infty,
	\]
	Lemma~\ref{lem:BC} (II) implies that $J_n\ge M\,t_n^{3/2}$ for infinitely many $n$, almost surely.
	
	\medskip
	\noindent
	\emph{Step 2.}
	Define
	\[
	F_n(M):=\left\{|\beta_{t_{n+1}}|\le \frac M4\,t_n^{1/2}\right\},
	\qquad
	G_n(M):=\left\{|I_{t_{n+1}}|\le \frac M4\,t_n^{3/2}\right\}.
	\]
	
	First,
	\[
	\frac{\beta_{t_{n+1}}}{t_n}
	=
	\frac{\beta_{t_n^2}}{t_n}
	\sim N(0,1),
	\]
	and hence
	\[
	\mathbb P(F_n(M)^c)
	=
	\mathbb P\!\left(|N(0,1)|>\frac M4\,t_n^{-1/2}\right).
	\]
	This is summable in \(n\), by the Gaussian tail bound and the
	super-exponential decay of \(t_n\). Hence, by Lemma~\ref{lem:BC} (I),
	applied to the events \(F_n(M)^c\), we obtain that $|\beta_{t_{n+1}}|\le \frac M4\,t_n^{1/2}$ for all sufficiently large \(n\), almost surely.
	
	Next, \(I_{t_{n+1}}\) is centered Gaussian with variance
	\[
	\operatorname{Var}(I_{t_{n+1}})
	=
	\frac{t_{n+1}^3}{3}
	=
	\frac{t_n^6}{3}.
	\]
	Thus
	\[
	\frac{I_{t_{n+1}}}{t_n^{3/2}}
	\sim
	N\!\left(0,\frac{t_n^3}{3}\right),
	\]
	and so
	\[
	\mathbb P(G_n(M)^c)
	=
	\mathbb P\!\left(
	|N(0,1)|>\frac M4\,\frac{\sqrt3}{t_n^{3/2}}
	\right).
	\]
	Again, this is summable in \(n\), by the Gaussian tail bound and the
	super-exponential decay of \(t_n\). Therefore, by Lemma~\ref{lem:BC} (I),
	applied to the events \(G_n(M)^c\), we obtain that $|I_{t_{n+1}}|\le \frac M4\,t_n^{3/2}$ for all sufficiently large $n$, almost surely.
	
	\medskip
	\noindent
	\emph{Step 3.}
	On the event \(E_n^+(M)\cap F_n(M)\cap G_n(M)\), using \(\Delta_n\le t_n\), 
	we have
	\[
	I_{t_n}
	=
	I_{t_{n+1}}+\beta_{t_{n+1}}\Delta_n+J_n
	\ge
	-\frac M4 t_n^{3/2}
	-\frac M4 t_n^{3/2}
	+M t_n^{3/2}
	=
	\frac M2 t_n^{3/2}.
	\]
	For each fixed \(M>0\), the event \(E_n^+(M)\) occurs infinitely often
	almost surely, while \(F_n(M)\) and \(G_n(M)\) hold for all sufficiently large
	\(n\) almost surely. Hence
	\[
	\limsup_{t\downarrow0}\frac{I_t}{t^{3/2}}\ge \frac M2
	\qquad\text{almost surely.}
	\]
	Taking \(M=1,2,3,\dots\) and intersecting the corresponding probability-one
	events, we obtain
	\[
	\limsup_{t\downarrow0}\frac{I_t}{t^{3/2}}=+\infty
	\qquad\text{almost surely.}
	\]
	
	Applying the same argument to \(-\beta\), we obtain
	\[
	\liminf_{t\downarrow0}\frac{I_t}{t^{3/2}}=-\limsup_{t\downarrow0}\frac{-I_t}{t^{3/2}}=-\infty
	\qquad\text{almost surely.}
	\]
	
	The equivalent formulation with \(\sup_{0<t\le T}\) and \(\inf_{0<t\le T}\)
	is immediate, since any sequence tending to \(0\) is eventually contained in
	\((0,T]\) for every fixed \(T>0\). This completes the proof.
\end{proof}

\begin{lemma}\label{lem:kernel_scaling}
	Let $W_t$ be a one-dimensional Brownian motion, and set
	\[
	Y_t:=\int_0^tU_s\dd s,\qquad U_t:=\sqrt 2\,W_t\,.
	\]
	Then for any $t>0$, $(Y_t,U_t)$ admits a Lebesgue density
	\begin{equation}\label{eq:pt_explicit}
		p_t(y,u):=\frac{\sqrt{3}}{2\pi}\,t^{-2}\exp\!\Big(-\frac{u^2}{4t}-3\Big(\frac{y}{t^{3/2}}-\frac{u}{2t^{1/2}}\Big)^{\!2}\Big),
		\qquad t>0.
	\end{equation}
	Moreover, for all $\mu>0$,
	\begin{equation}\label{eq:pt_scaling}
		p_t(\delta_\mu z)=\mu^{-4}p_{\mu^{-2}t}(z).
	\end{equation}
\end{lemma}

\begin{proof}
	By Lemma~\ref{lem:integrated_BM}, $(Y_t,U_t)=\big(\sqrt2 J_t,\sqrt2 W_t\big)$, where $J_t:=\int_0^tW_s\dd s$, is a centered Gaussian vector with covariance
	\[
	\mathrm{Var}(U_t)=2t,\qquad \mathrm{Cov}(Y_t,U_t)=t^2,\qquad \mathrm{Var}(Y_t)=\frac{2}{3}t^3.
	\]
	A direct inversion of the covariance matrix $\Sigma_t$ yields \eqref{eq:pt_explicit} with prefactor
	$\frac{1}{2\pi\sqrt{\det\Sigma_t}}=\frac{\sqrt{3}}{2\pi}t^{-2}$ (see \cite{Kolmogorov1934,Hormander1967}).
	The scaling \eqref{eq:pt_scaling} follows directly from \eqref{eq:pt_explicit}.
\end{proof}

\begin{remark}\label{2606211211}
	Similarly, for the Kolmogorov process in \eqref{eq:kolmogorov_XV} with
	initial condition \((X_0,V_0)=(x_0,v_0)\), we have
	\[
	V_t=v_0+\sqrt2 B_t,
	\qquad
	X_t=x_0+tv_0+\sqrt2\int_0^tB_s\dd s.
	\]
	Hence the centered increment
	\[
	(X_t-x_0-tv_0,\;V_t-v_0)
	\]
	is, by Lemma~\ref{lem:integrated_BM}, a centered Gaussian vector with covariance
	\[
	\mathrm{Var}(V_t-v_0)=2t\,I_d,
	\qquad
	\mathrm{Cov}(X_t-x_0-tv_0,V_t-v_0)=t^2 I_d,
	\qquad
	\mathrm{Var}(X_t-x_0-tv_0)=\frac{2}{3}t^3 I_d.
	\]
	Therefore the transition density of \((X_t,V_t)\) from \((x_0,v_0)\) is
	\begin{align*}
		P_t^{(d)}(x-x_0-tv_0,\ v-v_0),
	\end{align*}
	where
	\begin{align}\label{2606211207}
		P_t^{(d)}(x,v)
		:=
		\left(\frac{\sqrt{3}}{2\pi}\right)^d
		t^{-2d}
		\exp\!\left(
		-\left|\frac{v}{2t^{1/2}}\right|^2
		-
		3\left|
		\frac{x}{t^{3/2}}
		-
		\frac{v}{2t^{1/2}}
		\right|^{2}
		\right).
	\end{align}
	In particular, for every \(t>0\), the law of \(X_t\) is absolutely continuous
	with respect to the Lebesgue measure.
\end{remark}

\begin{lemma}\label{lem:kernel_bounds_away}
	Let $K$ be a compact subset of $\bR^2\setminus \{0\}$, and write
	$k=(k_y,k_u)\in\mathbb{N}_0^2$ for a multi-index. There exist constants
	$c_K>0$ and $C_{k,K}<\infty$ such that, for all $z\in K$ and all $t>0$,
	\begin{align}
		\label{eq:kernel_small}
		|\partial_y^{k_y}\partial_u^{k_u} p_t(z)|
		&\le C_{k,K}\, t^{-2-\frac{k_u}{2}-\frac{3k_y}{2}}\,\ee^{-c_K/t}\,.
	\end{align}
	In particular, for each $k$,
	\begin{align}\label{260619835}
		\sup_{z\in K}\int_0^\infty \,|\partial_y^{k_y}\partial_u^{k_u}p_t(z)|\dd t<\infty.
	\end{align}
\end{lemma}

\begin{proof}
	Fix a compact set \(K\subset\bR^2\setminus\{0\}\).
	From \eqref{eq:pt_explicit}, we have
	\[
	\partial_y^{k_y}\partial_u^{k_u} p_t(y,u)
	=
	t^{-2-\frac{3k_y}{2}-\frac{k_u}{2}}\,
	P_{k}\!\Big(\frac{y}{t^{3/2}},\frac{u}{t^{1/2}}\Big)\,
	\exp\big(-Q_t(y,u)\big),
	\]
	where \(P_{k}\) is a polynomial depending on \(k\), and
	\[
	Q_t(y,u):=
	\frac{u^2}{4t}
	+
	3\left(\frac{y}{t^{3/2}}-\frac{u}{2t^{1/2}}\right)^2 .
	\]
	Hence, for some \(m_k\ge0\),
	\[
	\big|\partial_y^{k_y}\partial_u^{k_u} p_t(y,u)\big|
	\le
	C_{k}\,
	t^{-2-\frac{3k_y}{2}-\frac{k_u}{2}}\,
	\left(
	1+\left|\frac{y}{t^{3/2}}\right|
	+\left|\frac{u}{t^{1/2}}\right|
	\right)^{m_k}
	\exp\big(-Q_t(y,u)\big).
	\]
	
	Observe that, in terms of \(a:=u\,t^{-1/2}\) and \(b:=y\,t^{-3/2}\),
	\[
	Q_t(y,u)
	=
	\frac{u^2}{t}
	-3\frac{uy}{t^2}
	+3\frac{y^2}{t^3}
	=
	a^2-3ab+3b^2 ,
	\]
	and
	\[
	a^2-3ab+3b^2-\frac18\bigl(a^2+b^2\bigr)
	=
	\frac78a^2-3ab+\frac{23}{8}b^2
	\ge0,
	\]
	since the discriminant \(9-4\cdot\frac78\cdot\frac{23}{8}=-\frac{17}{16}\) is
	negative. Hence
	\[
	Q_t(y,u)
	\ge
	\frac{1}{8}\left(\frac{y^2}{t^3}+\frac{u^2}{t}\right).
	\]
	Therefore,
	\[
	\left(
	1+\left|\frac{y}{t^{3/2}}\right|
	+\left|\frac{u}{t^{1/2}}\right|
	\right)^{m_k}
	\exp\big(-Q_t(y,u)\big)
	\le
	C_k
	\exp\left(
	-\frac{1}{16}
	\left(\frac{y^2}{t^3}+\frac{u^2}{t}\right)
	\right).
	\]
	Since \(K\) is compact and \(0\notin K\), there exist
	\(0<d_K\le D_K<\infty\) such that
	\[
	d_K\le |z|\le D_K
	\qquad\text{for all } z=(y,u)\in K.
	\]
	Moreover, for all \(t>0\),
	\[
	\frac{y^2}{t^3}+\frac{u^2}{t}
	\ge
	y^2\left(\frac{1}{t}-1\right)+\frac{u^2}{t}
	\ge
	-D_K^2+\frac{d_K^2}{t}.
	\]
	Thus there exist $C_{k,K}>0$ and $c_K>0$ such that
	\[
	\big|\partial_y^{k_y}\partial_u^{k_u} p_t(z)\big|
	\le
	C_{k,K}\,
	t^{-2-\frac{k_u}{2}-\frac{3k_y}{2}}
	\ee^{-c_K/t},
	\qquad t>0,\ z\in K.
	\]
	
	Finally, the integrability statement \eqref{260619835} follows directly from \eqref{eq:kernel_small},
	since \(\ee^{-c_K/t}\) controls the singularity as \(t\downarrow0\), while
	\(t^{-2-\frac{k_u}{2}-\frac{3k_y}{2}}\) is integrable as
	\(t\to\infty\).
\end{proof}

\begin{lemma}[Resolvent kernel]\label{lem:resolvent}
	For \(\lambda\ge0\), define \(\tilde p_t(y,u):=p_t(y,-u)\) and the resolvent kernel
	\[
	G_\lambda(z):=\int_0^\infty \ee^{-\lambda t}\tilde p_t(z)\,dt,
	\qquad z\neq 0.
	\]
	Then for any fixed $\lambda\ge 0$, the following assertions hold.
	\begin{enumerate}
		\item[(i)]
		\(G_\lambda\in C^\infty(\bR^2\setminus\{0\})\),
		\(G_\lambda>0\), and
		\begin{equation}\label{eq:pde_Glambda}
			(\lambda-\cL)G_\lambda=0
			\qquad\text{on }\bR^2\setminus\{0\}
		\end{equation}
		in the classical sense.
		
		\item[(ii)]
		For all \(\mu>0\) and \(z\neq 0\),
		\begin{equation}\label{eq:G_scaling}
			G_\lambda(\delta_\mu z)=\mu^{-2}G_{\mu^2\lambda}(z).
		\end{equation}
		
		\item[(iii)]
		There exists \(C_\lambda<\infty\) such that, for all \(z\neq 0\),
		\begin{equation}\label{eq:G_global_upper}
			G_\lambda(z)\le C_\lambda\,N(z)^{-2}.
		\end{equation}
		
		\item[(iv)]
		For any $R>0$, there exist constants \(0<c_{\lambda,R}\le C'_{\lambda,R}<\infty\)
		such that, for all \(z\neq 0\) with \(N(z)\le R\),
		\begin{equation}\label{eq:G_local_two_sided}
			c_{\lambda,R}\,N(z)^{-2}
			\le
			G_\lambda(z)
			\le
			C'_{\lambda,R}\,N(z)^{-2}.
		\end{equation}
		In particular, \(G_\lambda(z)\to\infty\) as \(z\to 0\).
	\end{enumerate}
\end{lemma}

\begin{remark}\label{rem:resolvent_distribution}
	Interpreted as a distribution on $\bR^2$, $G_\lambda$ is the (positive) resolvent kernel and satisfies
	\[
	(\lambda-\cL)G_\lambda=\delta_0
	\quad\text{in }\mathscr{D}'(\bR^2),
	\]
	while \eqref{eq:pde_Glambda} is the classical equation used on $\bR^2\setminus\{0\}$
	(see \cite{DynkinMarkov,BlumenthalGetoor}).
\end{remark}

\begin{proof}[\textbf{Proof of Lemma~\ref{lem:resolvent}}]
	(i) Fix a compact set \(K\Subset \bR^2\setminus\{0\}\). Since the same estimates
	as in Lemma~\ref{lem:kernel_bounds_away} hold for \(\tilde p_t\), we have, for
	each multi-index \(k\),
	\[
	\sup_{z\in K}
	\int_0^\infty \ee^{-\lambda t}
	|\partial_y^{k_y}\partial_u^{k_u}\tilde p_t(z)|\,dt
	<\infty.
	\]
	It follows that \(G_\lambda\in C^\infty(\bR^2\setminus\{0\})\), and
	differentiation under the integral is justified. Positivity is immediate from
	\(\tilde p_t(z)>0\) for \(t>0\).
	
	Next, we record the evolution equation satisfied by \(\tilde p_t\). A direct
	computation from \eqref{eq:pt_explicit} gives, for \(t>0\),
	\[
	\partial_t p_t=\cL^*p_t,
	\qquad
	\cL^*:=-u\partial_y+\partial_{uu}.
	\]
	For \(\tilde p_t(y,u):=p_t(y,-u)\), a direct computation gives
	\[
	(\partial_t-\cL)\tilde p_t(y,u)
	=
	\bigl((\partial_t-\cL^*)p_t\bigr)(y,-u)
	=
	0,
	\qquad t>0.
	\]
	Therefore, for \(z\neq0\),
	\[
	(\lambda-\cL)G_\lambda(z)
	=
	\int_0^\infty \ee^{-\lambda t}
	(\lambda \tilde p_t-\cL\tilde p_t)(z)\,dt
	=
	\int_0^\infty \ee^{-\lambda t}
	(\lambda \tilde p_t-\partial_t\tilde p_t)(z)\,dt.
	\]
	Fix \(z\neq0\). By Lemma~\ref{lem:kernel_bounds_away}, applied with
	\(K=\{z\}\) and with the multi-indices \(k\) satisfying \(|k|\le2\), both
	\(t\mapsto\tilde p_t(z)\) and
	\(t\mapsto\partial_t\tilde p_t(z)=\cL\tilde p_t(z)\) are integrable on
	\((0,\infty)\); moreover \(\tilde p_t(z)\to0\) as \(t\to\infty\) and as
	\(t\downarrow0\), by \eqref{eq:kernel_small}. Hence \(t\mapsto \ee^{-\lambda t}\tilde p_t(z)\) is
	absolutely continuous on \((0,\infty)\) with limits \(0\) at both endpoints, and
	\[
	\int_0^\infty \ee^{-\lambda t}
	(\lambda \tilde p_t-\partial_t\tilde p_t)(z)\,dt
	=
	\Big[-\ee^{-\lambda t}\tilde p_t(z)\Big]_{t=0}^{t=\infty}
	=
	0.
	\]
	This proves \eqref{eq:pde_Glambda}.
	
	(ii) Using \eqref{eq:pt_scaling} and \(\tilde p_t(y,u)=p_t(y,-u)\), we have
	\[
	\tilde p_t(\delta_\mu z)=\mu^{-4}\tilde p_{\mu^{-2}t}(z).
	\]
	With the change of variables \(t=\mu^2s\),
	\[
	G_\lambda(\delta_\mu z)
	=
	\int_0^\infty \ee^{-\lambda t}\tilde p_t(\delta_\mu z)\,dt
	=
	\int_0^\infty \ee^{-\lambda\mu^2s}
	\mu^{-4}\tilde p_s(z)\,\mu^2\,ds
	=
	\mu^{-2}G_{\mu^2\lambda}(z).
	\]
	This proves \eqref{eq:G_scaling}.
	
	(iv) Set
	\[
	\mathbb S:=\{z:N(z)=1\}
	=
	\{(y,u):u^4+|y|^{4/3}=1\}.
	\]
	Then \(\mathbb S\) is a compact subset of \(\bR^2\setminus\{0\}\).
	
	By Lemma~\ref{lem:kernel_bounds_away}, applied to the reflected compact set
	corresponding to \(\mathbb S\), there exist \(c_{\mathbb S}>0\) and
	\(C_{\mathbb S}<\infty\) such that
	\[
	\tilde p_t(\zeta)
	\le
	C_{\mathbb S} t^{-2}\ee^{-c_{\mathbb S}/t}
	\qquad \text{for all } \zeta\in\mathbb S,\ t>0.
	\]
	The right-hand side is integrable on \((0,\infty)\). Hence, by dominated
	convergence, the map
	$(\theta,\zeta)\mapsto G_\theta(\zeta)$
	is continuous on \([0,\infty)\times\mathbb S\). Moreover,
	\(G_\theta(\zeta)>0\) for all \(\theta\ge0\) and \(\zeta\in\mathbb S\).
	
	By compactness, there exist \(0<m_{\lambda,R}\le M_{\lambda,R}<\infty\) such that
	\[
	m_{\lambda,R}
	\le
	G_\theta(\zeta)
	\le
	M_{\lambda,R}
	\qquad
	\text{for all }(\theta,\zeta)\in[0,R^2\lambda]\times\mathbb S.
	\]
	
	Now let \(z\neq0\) with \(N(z)\le R\), and set
	\[
	\mu:=N(z)\in(0,R],
	\qquad
	\zeta:=\delta_{\mu^{-1}}z\in\mathbb S,
	\qquad
	\theta:=\mu^2\lambda\in [0,R^2\lambda].
	\]
	By the scaling relation \eqref{eq:G_scaling},
	$G_\lambda(z)=\mu^{-2}G_\theta(\zeta)$.
	Therefore,
	\[
	m_{\lambda,R} N(z)^{-2}
	\le
	G_\lambda(z)
	\le
	M_{\lambda,R} N(z)^{-2},
	\qquad 0<N(z)\le R.
	\]
	Setting \(c_{\lambda,R}:=m_{\lambda,R}\) and \(C'_{\lambda,R}:=M_{\lambda,R}\), this proves \eqref{eq:G_local_two_sided}.
	
	(iii) For \(N(z)\le1\), the desired upper bound follows from (iv). Now suppose
	\(N(z)\ge1\). Set
	\[
	\mu:=N(z)\ge1,
	\qquad
	\zeta:=\delta_{\mu^{-1}}z\in\mathbb S.
	\]
	By \eqref{eq:G_scaling},
	$G_\lambda(z)=\mu^{-2}G_{\mu^2\lambda}(\zeta)$.
	Since \(\tilde p_t\ge0\), the map \(\theta\mapsto G_\theta(\zeta)\) is
	nonincreasing. Hence, with 
	$$\widetilde M_\lambda:=\max_{\zeta\in\mathbb S}G_\lambda(\zeta)<\infty
	$$
	and \(\mu^2\lambda\ge\lambda\), we have
	$G_{\mu^2\lambda}(\zeta)
	\le
	G_\lambda(\zeta)
	\le
	\widetilde M_\lambda$.
	Thus
	$G_\lambda(z)
	\le
	\widetilde M_\lambda N(z)^{-2}$.
	Combining this with the case \(N(z)\le 1\) and setting
	\(C_\lambda:=\max\{C'_{\lambda,1},\widetilde M_\lambda\}\)
	proves \eqref{eq:G_global_upper}.
\end{proof}

\begin{lemma}[A polar barrier]\label{lem:drift_robust_barrier}
	Fix \(\alpha\in(0,1)\) and \(\lambda>0\). Let
	\(\Phi:=G_\lambda^\alpha\) on \(\bR^2\setminus\{0\}\).
	Then \(\Phi\in C^\infty(\bR^2\setminus\{0\})\), \(\Phi>0\), and
	\begin{equation}\label{eq:chain_rule_phi}
		(\lambda-\cL)\Phi
		=(1-\alpha)\lambda\,\Phi
		+\alpha(1-\alpha)\,G_\lambda^{\alpha-2}\,(\partial_u G_\lambda)^2
		\qquad\text{on }\bR^2\setminus\{0\}.
	\end{equation}
	Moreover, for every \(A\ge 0\), we have the following pointwise inequality on
	\(\bR^2\setminus\{0\}\):
	\begin{equation}\label{eq:drift_absorption_strong}
		(\lambda-\cL)\Phi
		\ge
		A\,|\partial_u \Phi|
		+\Big((1-\alpha)\lambda-\frac{A^2\alpha}{2(1-\alpha)}\Big)\Phi.
	\end{equation}
	In particular, if
	\begin{equation}\label{eq:lambda_condition}
		\lambda \ge \frac{A^2\alpha}{2(1-\alpha)^2},
	\end{equation}
	then
	\begin{equation}\label{eq:drift_robust_pde}
		(\lambda-\cL)\Phi \ge A\,|\partial_u\Phi|
		\qquad\text{on }\bR^2\setminus\{0\}.
	\end{equation}
\end{lemma}

\begin{proof}
	Since \(G_\lambda\) is smooth and positive on \(\bR^2\setminus\{0\}\),
	so is \(\Phi=G_\lambda^\alpha\). Since the second-order part of
	\(\cL\) acts only in the \(u\)-variable, the chain rule yields
	\[
	\cL(G_\lambda^\alpha)
	=
	\alpha G_\lambda^{\alpha-1}\cL G_\lambda
	+\alpha(\alpha-1)G_\lambda^{\alpha-2}(\partial_u G_\lambda)^2.
	\]
	Using \(\cL G_\lambda=\lambda G_\lambda\), which follows from
	\eqref{eq:pde_Glambda}, we obtain \eqref{eq:chain_rule_phi}.
	
	Next, \(\partial_u\Phi=\alpha G_\lambda^{\alpha-1}\partial_u G_\lambda\).
	Applying Young's inequality \(ab\le \frac12(a^2+b^2)\) with
	\[
	a:=\sqrt{\alpha(1-\alpha)}\,G_\lambda^{\frac{\alpha-2}{2}}
	|\partial_u G_\lambda|,
	\qquad
	b:=A\sqrt{\frac{\alpha}{1-\alpha}}\,G_\lambda^{\frac{\alpha}{2}},
	\]
	we get
	\[
	A|\partial_u\Phi|
	=
	A\alpha G_\lambda^{\alpha-1}|\partial_u G_\lambda|
	\le
	\frac{\alpha(1-\alpha)}{2}
	G_\lambda^{\alpha-2}(\partial_u G_\lambda)^2
	+
	\frac{A^2\alpha}{2(1-\alpha)}\Phi.
	\]
	Combining this bound with \eqref{eq:chain_rule_phi}, we obtain
	\[
	(\lambda-\cL)\Phi
	\ge
	A|\partial_u\Phi|
	+(1-\alpha)\lambda\,\Phi
	-\frac{A^2\alpha}{2(1-\alpha)}\Phi
	+\frac{\alpha(1-\alpha)}{2}
	G_\lambda^{\alpha-2}(\partial_u G_\lambda)^2.
	\]
	Dropping the last nonnegative term yields \eqref{eq:drift_absorption_strong}.
	If \eqref{eq:lambda_condition} holds, then the coefficient of \(\Phi\) in
	\eqref{eq:drift_absorption_strong} is nonnegative, and hence
	\eqref{eq:drift_robust_pde} follows.
\end{proof}

\begin{lemma}\label{lem:localized_bounded_drift_polar}
	Let $Z_t=(Y_t,U_t)$ be a continuous adapted process on a filtered probability
	space $(\Omega,\cF,(\cF_t)_{t\ge0},\mathbb P)$. Assume that there is a stopping
	time $\zeta$ such that
	\[
	Y_t=Y_0+\int_0^tU_s \dd s,\qquad U_t=U_0+\int_0^ta_s\dd s+\sqrt{2}W_t\quad \text{for all }0\le t\le \zeta\,,
	\]
	where $W$ is a one-dimensional $(\cF_t)$-Brownian motion, and $a_t$ is
	progressively measurable with $|a_t|\le A$ on $[0,\zeta)$ for some
	deterministic constant $A\ge0$.
	Then
	\[
	\mathbb{P}\left(
	\{Z_0\neq (0,0)\}\cap\{\exists t\in(0,\zeta]: Z_t=(0,0)\}
	\right)=0.
	\]
\end{lemma}

\begin{proof}
	Set
	\[
	E:=\{Z_0\neq(0,0)\}\in\cF_0 .
	\]
	Fix \(\alpha\in(0,1)\) and choose \(\lambda>0\) satisfying
	\eqref{eq:lambda_condition}. Let
	\begin{equation*}
		\Phi:=G_\lambda^\alpha .
	\end{equation*}
	
	For \(r>0\) and \(R>1\), define
	\[
	\sigma_r:=\inf\left\{t>0:\ N(Z_t)\le r\right\}\wedge \zeta,
	\qquad
	\tau_R:=\inf\left\{t>0:\ N(Z_t)\ge R\right\}\wedge \zeta,
	\]
	with the convention \(\inf\emptyset=\infty\); these are stopping times by
	Lemma~\ref{lem:hitting_time}.
	For \(K\in\mathbb{N}\), set
	\begin{equation*}
		E_K:=
		E\cap\left\{\frac1K\le N(Z_0)\le K\right\}.
	\end{equation*}
	Since \(Z_0\neq(0,0)\) on \(E\), we have
	$E=\bigcup_{K\in \bN}E_K$.
	
	We first prove the assertion on each \(E_K\) and before each finite time
	\(T_0\). Fix \(K,T_0\in\mathbb N\). For
	$0<r<\frac1K$ and $R\ge K$, set
	\[
	S:=\zeta\wedge\sigma_r\wedge\tau_R\wedge T_0.
	\]
	On \(E_K\), the path \((Z_t)_{0\le t\le S}\) remains in the closed annulus
	\begin{equation*}
		\left\{z\in\bR^2:\ r\le N(z)\le R\right\}
		\subset \bR^2\setminus\{0\}.
	\end{equation*}
	Moreover, if \(N(Z_0)\notin[r,R]\), then \(S=0\) by the definitions of
	\(\sigma_r\) and \(\tau_R\). Since
	\(\Phi\in C^\infty(\bR^2\setminus\{0\})\) and the closed annulus
	$\{z:\ r\le N(z)\le R\}$
	is compactly contained in the open annulus
	\(\{z:\ r/2<N(z)<2R\}\), choose
	\(\psi\in C_c^\infty(\bR^2\setminus\{0\})\) such that
	\(\psi\equiv1\) on an open neighborhood of
	\(\{z:\ r\le N(z)\le R\}\) and
	\(\operatorname{supp}\psi\subset\{z:\ r/2<N(z)<2R\}\).
	Set \(\Phi_{r,R}:=\psi\Phi\) on \(\bR^2\setminus\{0\}\), extended by
	\(0\) to all of \(\bR^2\). Then
	\(\Phi_{r,R}\in C_c^\infty(\bR^2)\), so its first and second derivatives
	are bounded, and \(\Phi_{r,R}\) coincides with \(\Phi\), together with all
	its derivatives, on an open neighborhood of
	\(\{z:\ r\le N(z)\le R\}\).
	
	Since the hypotheses only give the equations for \(Z\) on \([0,\zeta]\), we
	apply It\^o's formula to the stopped process \(Z^\zeta_t:=Z_{t\wedge\zeta}\).
	This is a globally defined It\^o process,
	\[
	Y_{t\wedge\zeta}
	=
	Y_0+\int_0^t \mathbf 1_{\{s<\zeta\}}U_{s\wedge\zeta}\dd s,
	\qquad
	U_{t\wedge\zeta}
	=
	U_0+\int_0^t \mathbf 1_{\{s<\zeta\}}a_s\dd s
	+\sqrt2\int_0^t \mathbf 1_{\{s<\zeta\}}\dd W_s ,
	\]
	with progressively measurable, almost surely locally bounded integrands.
	Applying It\^o's formula (Lemma~\ref{lem:ito}) to
	\(\ee^{-\lambda t}\Phi_{r,R}(Z^\zeta_t)\) up to \(S\), and using \(S\le\zeta\)
	together with \(\Phi_{r,R}=\Phi\) along the stopped path on \(E_K\), we get, on
	\(E_K\),
	\[
	\ee^{-\lambda S}\Phi(Z_S)
	=
	\Phi(Z_0)
	+\int_0^S \ee^{-\lambda t}
	\bigl(-\lambda\Phi+a_t\partial_u\Phi+\cL\Phi\bigr)(Z_t)\,dt
	+\sqrt2\int_0^S \ee^{-\lambda t}\partial_u\Phi_{r,R}(Z_t)\,dW_t .
	\]
	By Lemma~\ref{lem:drift_robust_barrier},
	\[
	\lambda\Phi-\cL\Phi-a_t\partial_u\Phi
	\ge
	\lambda\Phi-\cL\Phi-A|\partial_u\Phi|
	\ge0.
	\]
	Hence, on \(E_K\),
	\begin{equation*}
		\ee^{-\lambda S}\Phi(Z_S)
		\le
		\Phi(Z_0)
		+\sqrt2\int_0^S \ee^{-\lambda t}\partial_u\Phi_{r,R}(Z_t)\,dW_t .
	\end{equation*}
	
	Since \(\partial_u\Phi_{r,R}\) is bounded on \(\bR^2\),
	\[
	M_t:=\sqrt2\int_0^{S\wedge t}\ee^{-\lambda s}\partial_u\Phi_{r,R}(Z_s)\,dW_s
	\]
	is a square-integrable martingale by Lemma~\ref{lem:stochastic_integral}(ii).
	Since \(S\le T_0\) and
	\(E_K\in\cF_0\), conditioning on \(\cF_0\) gives
	\[
	\mathbb E\left[
	\mathbf 1_{E_K}\int_0^S \ee^{-\lambda s}\partial_u\Phi_{r,R}(Z_s)\,dW_s
	\right]
	=
	\mathbb E\bigl[\mathbf 1_{E_K}\,\mathbb E[M_{T_0}\mid\cF_0]\bigr]
	=0.
	\]
	Therefore,
	\[
	\mathbb E\left[
	\mathbf 1_{E_K}\ee^{-\lambda S}\Phi(Z_S)
	\right]
	\le
	\mathbb E\left[
	\mathbf 1_{E_K}\Phi(Z_0)
	\right].
	\]
	Since \(N(Z_0)\in[1/K,K]\) on \(E_K\), and \(\Phi\) is continuous on the
	compact annulus \(\{1/K\le N(z)\le K\}\), there exists a constant
	\(C_K<\infty\) such that
	$\mathbb E\left[\mathbf 1_{E_K}\Phi(Z_0)\right]\le C_K$.
	
	On the event
	\[
	E_K\cap
	\left\{\exists t\in(0,\zeta\wedge\tau_R\wedge T_0]:
	N(Z_t)=r\right\},
	\]
	we have \(S=\sigma_r\). By the lower bound in Lemma~\ref{lem:resolvent}(iv), applied with radius $1$ (recall $r<1/K\le1$),
	\[
	\ee^{-\lambda S}\Phi(Z_S)
	\ge
	\ee^{-\lambda T_0}c_{\lambda,1}^\alpha r^{-2\alpha}.
	\]
	Consequently,
	\[
	\mathbb P\left(
	E_K\cap
	\left\{\exists t\in(0,\zeta\wedge\tau_R\wedge T_0]:
	N(Z_t)=r\right\}
	\right)
	\le
	C(K,T_0,\alpha,\lambda)r^{2\alpha}.
	\]
	
	For \(R\ge K\), set
	\[
	A_R(r):=
	E_K\cap
	\left\{\exists t\in(0,\zeta\wedge\tau_R\wedge T_0]:
	N(Z_t)=r\right\}.
	\]
	Then \(A_R(r)\) is increasing in \(R\). Moreover, by the continuity of \(Z\),
	\[
	\bigcup_{R\ge K,\ R\in\mathbb N} A_R(r)
	=
	E_K\cap
	\left\{\exists t\in(0,\zeta\wedge T_0]:
	N(Z_t)=r\right\}.
	\]
	Letting \(R\to\infty\) through integer values, we obtain
	\[
	\mathbb P\left(
	E_K\cap
	\left\{\exists t\in(0,\zeta\wedge T_0]:
	N(Z_t)=r\right\}
	\right)
	\le
	C(K,T_0,\alpha,\lambda)r^{2\alpha}.
	\]
	
	If \(Z_t\) hits \((0,0)\) at some time \(t\in(0,\zeta\wedge T_0]\), then,
	since \(N(Z_0)\ge1/K\) on \(E_K\) and \(Z\) is continuous, the path must hit
	the level set \(\{N=r\}\) before \(\zeta\wedge T_0\), for every sufficiently
	small \(r>0\). Hence, for every sufficiently small \(r<1/K\),
	\[
	\mathbb P\left(
	E_K\cap
	\left\{\exists t\in(0,\zeta\wedge T_0]:
	Z_t=(0,0)\right\}
	\right)
	\le
	C(K,T_0,\alpha,\lambda)r^{2\alpha}.
	\]
	Letting \(r\downarrow0\) gives
	\[
	\mathbb P\left(
	E_K\cap
	\left\{\exists t\in(0,\zeta\wedge T_0]:
	Z_t=(0,0)\right\}
	\right)=0.
	\]
	
	Finally, taking the countable union over \(K,T_0\in\mathbb N\), we conclude that
	\[
	\mathbb P\left(
	E\cap
	\left\{\exists t\in(0,\zeta]:
	Z_t=(0,0)\right\}
	\right)=0.
	\]
\end{proof}

\subsubsection{Proof of the Results}\label{subsubsec:grazing-proofs}

\begin{proof}[\textbf{Proof of Theorem~\ref{thm:finiteness_exit_time}}]
	Since \(\cD\) is bounded, choose \(R>0\) such that \(\overline{\cD}\subset B(0,R)\).
	For each fixed \(t>0\), we have
	\[
	X_t
	=
	(x+tv)+\sqrt2\int_0^t B_s\dd s.
	\]
	Thus \(X_t\) is a nondegenerate Gaussian random vector with mean \(x+tv\) and covariance
	$\mathrm{Cov}(X_t)
	=
	\frac{2}{3}t^3 I_d$
	(see Remark~\ref{2606211211}). Hence the density of \(X_t\) is given by
	\begin{align}\label{260621406}
		f_t(y)
		=
		\left(\frac{3}{4\pi t^3}\right)^{d/2}
		\exp\left(
		-\frac34\left|\frac{y-x-tv}{t^{3/2}}\right|^2
		\right).
	\end{align}
	Since \(\overline{\cD}\subset B(0,R)\), we have
	\begin{align}\label{2606211219}
		\mathbb P(|X_t|\le R)
		=
		\int_{B(0,R)} f_t(y)\dd y\le
		|B(0,R)|\sup_{y\in\bR^d}f_t(y) \le
		C_R t^{-3d/2}.
	\end{align}
	
	Taking \(t_k:=2^k\), we obtain from \eqref{2606211219} that
	\[
	\sum_{k=1}^\infty
	\mathbb P\bigl(X_{t_k}\in\overline{\cD}\bigr)
	\le
	\sum_{k=1}^\infty
	\mathbb P\bigl(|X_{t_k}|\le R\bigr)
	\le
	C_R\sum_{k=1}^\infty 2^{-3dk/2}
	<\infty.
	\]
	By the first Borel--Cantelli lemma (Lemma~\ref{lem:BC}\textup{(I)}), almost surely, \(X_{t_k}\notin\overline{\cD}\) for all sufficiently large \(k\). In particular, \(X_{t_k}\notin\overline{\cD}\) for some \(k\), and therefore
	$\tau_\cD(x,v)\le t_k<\infty$.
	Thus \(\tau_\cD(x,v)<\infty\) almost surely.
\end{proof}

\begin{proof}[\textbf{Proof of Theorem~\ref{thm:exit_time_properties}}]
	Let \(\rho\) be the signed distance to \(\partial\cD\) defined by
	\begin{align}\label{2606211233}
		\rho(x)
		=
		-d(x,\partial\cD){\bf 1}_{x\in\cD}
		+
		d(x,\partial\cD){\bf 1}_{x\notin\cD}.
	\end{align}
	Since \(\partial\cD\) is \(C^{1,1}\), there exists \(r_0>0\) such that the following holds:
	\begin{itemize}
		\item \(\rho\in C^{1,1}\) on the tubular neighborhood
		\[
		\cT_2:=\{|\rho|<2r_0\}.
		\]
		
		\item For every \(x\in \cT_2\), there exists a unique \(\pi(x)\in\partial\cD\) such that
		\[
		d(x,\partial\cD)=|x-\pi(x)|,
		\qquad
		\nabla_x\rho(x)=n(\pi(x)).
		\]
	\end{itemize}
	Set the smaller tube
	\begin{align}\label{2606211239}
		\cT:=\{|\rho|<r_0\}\subset \cT_2,
	\end{align}
	and take \(\chi\in C_c^\infty(\bR^d_x)\) such that
	$
	\chi\equiv1\quad\text{in }\cT$ and 
	$\supp\chi\subset  \cT_2$.
	We set
	\[
	\tilde\rho:=\rho\chi.
	\]
	Then \(\tilde\rho\in C_c^{1,1}(\bR^d_x)\), and
	$\nabla_x\tilde\rho(x)=n(\pi(x))$
	for all $x\in \cT$.
	Since \(\tilde\rho\in C_c^{1,1}(\bR^d_x)\), its second-order weak derivatives \(\nabla_x^2\tilde\rho\) exist and satisfy
	$\nabla_x^2\tilde\rho\in L^\infty(\bR^d_x)$.
	Once and for all, we fix a Borel representative of \(\nabla_x^2\tilde\rho\)
	satisfying the bound
	\begin{align}\label{2608051320}
		\bigl|\nabla_x^2\tilde\rho(x)\bigr|
		\le
		\|\nabla_x^2\tilde\rho\|_{L^\infty(\bR^d_x)}
		\qquad\text{for \emph{every} }x\in\bR^d_x,
	\end{align}
	obtained by replacing an arbitrary representative by \(0\) on the null set
	where the bound fails. All pointwise statements about
	\(\nabla_x^2\tilde\rho\) below refer to this representative.
	
	Define
	\[
	Y_t:=\tilde\rho(X_t),
	\qquad
	U_t:=V_t\cdot\nabla_x\tilde\rho(X_t).
	\]
	We claim that, with probability one, for all \(t>0\),
	\begin{align}\label{260620916}
		Y_t
		=
		Y_0+\int_0^t U_s\dd s,
	\end{align}
	and
	\begin{align}\label{260620917}
		U_t
		=
		U_0
		+
		\int_0^t
		V_s^\top\nabla_x^2\tilde\rho(X_s)V_s\dd s
		+
		\sqrt2\int_0^t\nabla_x\tilde\rho(X_s)\cdot\dd B_s.
	\end{align}
	\eqref{260620916} follows from the chain rule.
	Therefore, we only prove \eqref{260620917}.
	
	Let \(\eta_\varepsilon\) be a standard mollifier in \(\bR^d\), and define
	\[
	\tilde\rho_\varepsilon:=\tilde\rho*\eta_\varepsilon.
	\]
	Then \(\tilde\rho_\varepsilon\in C_c^\infty(\bR^d)\) and
	\begin{align*}
		\tilde\rho_\varepsilon&\to\tilde\rho
		\quad\text{in }C^1(\bR^d_x)
		\quad\text{as }\varepsilon\to0,\qquad
		\sup_{\varepsilon\in(0,1)}
		\|\nabla_x^2\tilde\rho_\varepsilon\|_{L^\infty(\bR^d_x)}
		\le
		\|\nabla_x^2\tilde\rho\|_{L^\infty(\bR^d_x)}.
	\end{align*}
	Moreover, by the Lebesgue differentiation theorem, there exists a null set
	\(\mathscr N\subset\bR^d_x\) such that
	\begin{align}\label{260620835}
		\nabla_x^2\tilde\rho_\varepsilon(x)\to\nabla_x^2\tilde\rho(x)
		\qquad
		\text{for every }x\in\bR^d_x\setminus\mathscr N.
	\end{align}
	Set
	\[
	Y_t^\varepsilon:=\tilde\rho_\varepsilon(X_t),
	\qquad
	U_t^\varepsilon:=V_t\cdot\nabla_x\tilde\rho_\varepsilon(X_t).
	\]
	Since the function
	\[
	F_\varepsilon(x,v):=v\cdot\nabla_x\tilde\rho_\varepsilon(x)
	\]
	is smooth, It\^o's formula, together with
	\(dX_t=V_t\,dt\) and \(dV_t=\sqrt2\,dB_t\), gives
	\begin{equation}\label{eq:Uprime_final}
		U_t^\varepsilon
		=
		U_0^\varepsilon
		+
		\int_0^t
		V_s^\top
		(\nabla_x^2\tilde\rho_\varepsilon)(X_s)
		V_s\dd s
		+
		\sqrt2
		\int_0^t
		\nabla_x\tilde\rho_\varepsilon(X_s)\cdot\dd B_s.
	\end{equation}
	
	Fix \(S>0\) and \(R\ge1\), and define the localization time
	\[
	\theta_R:=\inf\{t\ge0: |V_t|\ge R\},
	\]
	which is a stopping time by Lemma~\ref{lem:hitting_time}.
	
	The convergence of \(\tilde\rho_\varepsilon\) to \(\tilde\rho\) in
	\(C^1(\bR^d_x)\) implies
	\begin{align}\label{260620901}
		\sup_{0\le t\le S}|U_t^\varepsilon-U_t|
		&\le
		\left(\sup_{0\le t\le S}|V_t|\right)
		\|\nabla_x\tilde\rho_\varepsilon-\nabla_x\tilde\rho\|_{L^\infty(\bR^d_x)}
		\to0
	\end{align}
	almost surely. Moreover, by Lemma~\ref{lem:stochastic_integral}(iii),
	\begin{align}\label{260620902}
		\begin{aligned}
			\mathbb E\left[
			\sup_{0\le t\le S}
			\left|
			\sqrt2\int_0^t
			\bigl(\nabla_x\tilde\rho_\varepsilon-\nabla_x\tilde\rho\bigr)(X_s)
			\cdot\dd B_s
			\right|^2
			\right]\,&\le
			8
			\mathbb E\int_0^S
			\left|
			\bigl(\nabla_x\tilde\rho_\varepsilon-\nabla_x\tilde\rho\bigr)(X_s)
			\right|^2\dd s\\
			&
			\le
			8S
			\|\nabla_x\tilde\rho_\varepsilon-\nabla_x\tilde\rho\|_{L^\infty(\bR^d_x)}^2
			\to0.
		\end{aligned}
	\end{align}
	
	In addition,
	\begin{align}\label{260620903-pre}
		&\mathbb E\left[
		\sup_{0\le t\le S\wedge\theta_R}
		\left|
		\int_0^t
		V_s^\top
		\bigl[
		\nabla_x^2\tilde\rho_\varepsilon(X_s)
		-
		\nabla_x^2\tilde\rho(X_s)
		\bigr]
		V_s\dd s
		\right|
		\right]
		\le
		R^2
		\mathbb E\int_0^S
		\left|
		\bigl(\nabla_x^2\tilde\rho_\varepsilon
		-
		\nabla_x^2\tilde\rho\bigr)(X_s)
		\right|\dd s.
	\end{align}
	Let \(\mathscr N\subset\bR^d_x\) be the null set in \eqref{260620835}.
	Since \((X_0,V_0)=(x,v)\), the transition density of \((X_t,V_t)\) is
	given by $P_t^{(d)}(y-x-tv,w-v)$, where \(P_t^{(d)}\) is defined in \eqref{2606211207}. Hence, for every
	\(t>0\),
	\[
	\mathbb P(X_t\in\mathscr N)
	=
	\int_{\mathscr N}\int_{\bR^d}
	P_t^{(d)}(y-x-tv,w-v)\dd w\dd y
	=
	0.
	\]
	Hence
	\[
	\mathbb E\int_0^S {\bf 1}_{\mathscr N}(X_t)\dd t
	=
	\int_0^S\mathbb P(X_t\in\mathscr N)\dd t
	=0.
	\]
	Thus \(X_t(\omega)\notin\mathscr N\) for a.e.
	\((\omega,t)\in\Omega\times[0,S]\). Therefore,
	\[
	\nabla_x^2\tilde\rho_\varepsilon(X_t(\omega))
	\to
	\nabla_x^2\tilde\rho(X_t(\omega))
	\]
	for a.e. \((\omega,t)\in\Omega\times[0,S]\). Moreover,
	\[
	\left|
	\nabla_x^2\tilde\rho_\varepsilon(X_t)
	-
	\nabla_x^2\tilde\rho(X_t)
	\right|
	\le
	2\|\nabla_x^2\tilde\rho\|_{L^\infty(\bR^d_x)}.
	\]
	By the dominated convergence theorem,
	\[
	\mathbb E\int_0^S
	\left|
	\bigl(\nabla_x^2\tilde\rho_\varepsilon
	-
	\nabla_x^2\tilde\rho\bigr)(X_s)
	\right|\dd s
	\to0.
	\]
	Together with \eqref{260620903-pre}, this gives
	\begin{align}\label{260620903}
		&\mathbb E\left[
		\sup_{0\le t\le S\wedge\theta_R}
		\left|
		\int_0^t
		V_s^\top
		\bigl[
		\nabla_x^2\tilde\rho_\varepsilon(X_s)
		-
		\nabla_x^2\tilde\rho(X_s)
		\bigr]
		V_s\dd s
		\right|
		\right]
		\to0.
	\end{align}
	
	Combining \eqref{260620901}, \eqref{260620902}, and \eqref{260620903}, we
	may choose a sequence \(\varepsilon_k=\varepsilon_k(S,R)\downarrow0\) along
	which all terms in \eqref{eq:Uprime_final} converge almost surely,
	uniformly on \([0,S\wedge\theta_R]\). Passing to the limit yields
	\eqref{260620917} for all \(0\le t\le S\wedge\theta_R\), almost surely.
	
	Taking a countable intersection over \(S,R\in\mathbb N\), and using
	\(\theta_R\to\infty\) almost surely as \(R\to\infty\), we obtain that
	\eqref{260620917} holds for all \(t>0\), almost surely.
	
	\smallskip
	(1) Define
	\[
	\tau_{\gamma^0}:=\inf\{t>0:(X_t,V_t)\in\gamma^0\}.
	\]
	For each \(q\in\mathbb Q_+:=\mathbb Q\cap(0,\infty)\), define
	\[
	\zeta_q:=\inf\{t\ge q:X_t\notin \cT\},
	\]
	with the convention that \(\inf\emptyset=\infty\). By
	Lemma~\ref{lem:hitting_time}, \(\zeta_q\) is a stopping time; moreover
	\(\{X_q\in\cT\}\in\cF_q\).
	
	Since \(\partial\cD\) is compact and the outer unit normal \(n\) is
	continuous on \(\partial\cD\), the set \(\gamma^0\) is closed in
	\(\bR^d_x\times\bR^d_v\). Hence, since
	\((X_0,V_0)=(x,v)\notin\gamma^0\), the continuity of
	\((X_t,V_t)\) yields
	\[
	\tau_{\gamma^0}>0
	\qquad\text{almost surely,}
	\]
	and, by Lemma~\ref{lem:hitting_time},
	\((X_{\tau_{\gamma^0}},V_{\tau_{\gamma^0}})\in\gamma^0\) on
	\(\{\tau_{\gamma^0}<\infty\}\).
	
	On the event that there exists \(t_0>0\) such that
	\((X_{t_0},V_{t_0})\in\gamma^0\), we have
	\(0<\tau_{\gamma^0}<\infty\), and
	\[
	X_{\tau_{\gamma^0}}\in\partial\cD\subset \cT,
	\qquad
	Y_{\tau_{\gamma^0}}=0,
	\qquad
	U_{\tau_{\gamma^0}}
	=
	V_{\tau_{\gamma^0}}\cdot n(X_{\tau_{\gamma^0}})=0.
	\]
	Since \(X_{\tau_{\gamma^0}}\in \cT\) and \(\cT\) is open, the continuity of
	\(X_t\) implies that there exist \(q_1,q_2\in\mathbb Q_+\) such that
	\[
	q_1<\tau_{\gamma^0}<q_2,
	\qquad
	X_t\in \cT
	\quad\text{for all }t\in[q_1,q_2].
	\]
	For this choice of \(q_1\), we have
	\[
	X_{q_1}\in\cT,
	\qquad
	q_1<\tau_{\gamma^0}<\zeta_{q_1}.
	\]
	Consequently,
	\begin{align}\label{2606201036}
		&\mathbb P\left(
		\left\{\exists\,t>0:\ (X_t,V_t)\in\gamma^0\right\}
		\right) \notag\\
		&\qquad\le
		\sum_{q\in\mathbb Q_+}
		\mathbb P\left(
		\left\{X_q\in\cT\right\}\cap\left\{q<\tau_{\gamma^0}\right\}
		\cap
		\left\{
		\exists\,t\in[q,\zeta_q):\ (X_t,V_t)\in\gamma^0
		\right\}
		\right).
	\end{align}
	We claim that, for each \(q\in\mathbb Q_+\),
	\begin{align}\label{260620953}
		\mathbb P\left(
		\left\{X_q\in\cT\right\}\cap\left\{q<\tau_{\gamma^0}\right\}
		\cap
		\left\{
		\exists\,t\in[q,\zeta_q):\ (X_t,V_t)\in\gamma^0
		\right\}
		\right)
		=0.
	\end{align}
	
	Fix \(q\in\mathbb Q_+\), and let
	\[
	\hat{\cF}_t:=\cF_{t+q},
	\qquad
	\hat B_t:=B_{t+q}-B_q,
	\qquad t\ge0.
	\]
	By Lemma~\ref{lem:brownian_shift}, \(\hat B\) is a \(d\)-dimensional Brownian
	motion with respect to \((\hat{\cF}_t)_{t\ge0}\). Put
	\[
	\hat Y_t:=Y_{t+q},
	\qquad
	\hat U_t:=U_{t+q}.
	\]
	Then
	\[
	\hat Y_t=\hat Y_0+\int_0^t \hat U_s\dd s
	\]
	and
	\[
	\hat U_t
	=
	\hat U_0
	+
	\int_0^t
	V_{s+q}^{\top}\nabla_x^2\tilde\rho(X_{s+q})V_{s+q}\dd s
	+
	\sqrt2
	\int_0^t
	\nabla_x\tilde\rho(X_{s+q})\cdot\dd\hat B_s,
	\]
	where the stochastic integral has been rewritten in terms of \(\hat B\) by
	Lemma~\ref{lem:brownian_shift}.
	
	Fix a unit vector \(e_0\in\bR^d\), write \(\Lambda_q:=\{X_q\in\cT\}\),
	and define
	\[
	\hat\nu_t
	:=
	\begin{cases}
		\nabla_x\tilde\rho(X_{t+q}),
		& \text{if } \Lambda_q\text{ occurs and } 0\le t\le \zeta_q-q,\\
		e_0,
		& \text{otherwise}.
	\end{cases}
	\]
	Since \(\Lambda_q\in\cF_q=\hat{\cF}_0\) and \(\zeta_q-q\) is an
	\((\hat{\cF}_t)\)-stopping time, the first alternative is a
	progressively measurable subset of \(\Omega\times[0,\infty)\).
	If \(\Lambda_q\) occurs and \(0\le t\le\zeta_q-q\), then
	\(X_{t+q}\in\overline{\cT}\) by the definition of \(\zeta_q\) and the
	continuity of \(X\). Since \(\chi\equiv1\) on \(\cT\), we have \(\nabla_x\tilde\rho=\nabla_x\rho\) and
	\(|\nabla_x\tilde\rho|=1\) on \(\overline{\cT}\); hence \(\hat\nu\) is
	progressively measurable and unit-valued. Set
	\[
	\hat W_t
	:=
	\int_0^t \hat\nu_s\cdot\dd\hat B_s,
	\qquad t\ge0.
	\]
	Then, by Lemma~\ref{lem:stochastic_integral}(iv), \(\hat W\) is a real-valued
	continuous local martingale with respect to
	\((\hat{\cF}_t)_{t\ge0}\), and
	\[
	\langle \hat W\rangle_t
	=
	\int_0^t |\hat\nu_s|^2\dd s
	=
	t.
	\]
	By Lemma~\ref{lem:Levy_characterization}, \(\hat W\) is a standard
	one-dimensional Brownian motion with respect to
	\((\hat{\cF}_t)_{t\ge0}\).
	
	Define
	\[
	\widetilde U_t
	:=
	\hat U_0
	+
	\int_0^t
	V_{s+q}^{\top}\nabla_x^2\tilde\rho(X_{s+q})V_{s+q}\dd s
	+
	\sqrt2\,\hat W_t,
	\]
	and
	\[
	\widetilde Y_t
	:=
	\hat Y_0+\int_0^t \widetilde U_s\dd s.
	\]
	On \(\Lambda_q\), the two progressively measurable integrands \(\hat\nu_s\)
	and \(\nabla_x\tilde\rho(X_{s+q})\) agree for \(0\le s\le\zeta_q-q\). Since
	\(\Lambda_q\in\hat{\cF}_0\) and \(\zeta_q-q\) is an
	\((\hat{\cF}_t)\)-stopping time, Lemma~\ref{lem:stochastic_integral}(i)
	gives, on \(\Lambda_q\),
	\[
	\hat W_t
	=
	\int_0^{t}\nabla_x\tilde\rho(X_{s+q})\cdot\dd\hat B_s
	\qquad\text{for all }0\le t\le\zeta_q-q.
	\]
	Therefore, again on \(\Lambda_q\),
	\[
	(\widetilde Y_t,\widetilde U_t)=(\hat Y_t,\hat U_t)
	\qquad
	\text{for all }0\le t\le \zeta_q-q .
	\]
	For \(R\ge1\), put
	\[
	\hat\theta_R:=\inf\{t\ge0: |V_{t+q}|\ge R\},
	\]
	again a stopping time by Lemma~\ref{lem:hitting_time}.
	On \([0,\hat\theta_R]\), the drift satisfies, by \eqref{2608051320},
	\[
	\left|V_{t+q}^{\top}\nabla_x^2\tilde\rho(X_{t+q})V_{t+q}\right|\le R^2\|\nabla_x^2\tilde\rho\|_{L^\infty(\bR^d_x)},
	\]
	where \(\|\cdot\|_{L^\infty}\) denotes the essential supremum of the operator
	norm. Hence, applying Lemma~\ref{lem:localized_bounded_drift_polar} to
	\((\widetilde Y,\widetilde U)\) with \(\zeta=\hat\theta_R\), we obtain
	\[
	\mathbb P\left(
	(\widetilde Y_0,\widetilde U_0)\neq(0,0),\
	\exists t\in(0,\hat\theta_R]:
	(\widetilde Y_t,\widetilde U_t)=(0,0)
	\right)=0.
	\]
	
	Set
	\[
	A_q
	:=
	\left\{X_q\in\cT\right\}\cap\left\{q<\tau_{\gamma^0}\right\}
	\cap
	\left\{
	\exists\,t\in[q,\zeta_q):
	(X_t,V_t)\in\gamma^0
	\right\}.
	\]
	On \(A_q\), the hitting time inside the last event is strictly larger
	than \(q\), since \(q<\tau_{\gamma^0}\). Moreover, \(X_q\in \cT\) on \(A_q\), hence
	\[
	(Y_q,U_q)=\big(\rho(X_q),V_q\cdot n(\pi(X_q))\big)\neq(0,0).
	\]
	Indeed, if \(\rho(X_q)=0\) and \(U_q=0\), then \((X_q,V_q)\in\gamma^0\),
	which contradicts \(q<\tau_{\gamma^0}\).
	
	If \((X_{q+t},V_{q+t})\in\gamma^0\) for some
	\(t\in(0,\zeta_q-q)\), then, since \(X_{q+t}\in \cT\),
	\[
	(Y_{q+t},U_{q+t})=(0,0).
	\]
	By the continuity of \(V\), for such a \(t\) there exists \(R\in\mathbb N\)
	with \(t\le\hat\theta_R\). Therefore
	\[
	A_q
	=
	\bigcup_{R=1}^{\infty}
	\left(
	A_q\cap
	\left\{
	\exists t\in(0,(\zeta_q-q)\wedge\hat\theta_R]:
	(Y_{q+t},U_{q+t})=(0,0)
	\right\}
	\right).
	\]
	Since \(A_q\subset\{X_q\in\cT\}\), we have on \(A_q\)
	\[
	(\widetilde Y_t,\widetilde U_t)=(Y_{q+t},U_{q+t})
	\qquad
	\text{for }0\le t\le \zeta_q-q,
	\]
	and \((\widetilde Y_0,\widetilde U_0)=(Y_q,U_q)\neq(0,0)\). Hence each set in
	the union above is contained in
	\[
	\left\{
	(\widetilde Y_0,\widetilde U_0)\neq(0,0),\
	\exists t\in(0,\hat\theta_R]:
	(\widetilde Y_t,\widetilde U_t)=(0,0)
	\right\},
	\]
	which has probability zero. Consequently,
	\[
	\mathbb P(A_q)=0.
	\]
	This proves \eqref{260620953}. Together with \eqref{2606201036}, this proves
	assertion~(1).
	
	\smallskip
	(2) Recall the processes \(Y\) and \(U\) defined at the beginning of the proof, which satisfy
	\eqref{260620916} and \eqref{260620917}.
	If \((x,v)\in\gamma^+\), then \(x\in\partial\cD\) and \(U_0=v\cdot n(x)>0\).
	For a fixed \(\omega\), by continuity of \(X_t\) and \(U_t\), there exists \(\delta_1(\omega)>0\) such that
	\[
	X_t\in \cT\quad \text{and}\quad U_t(\omega)\ge \frac12U_0\quad \text{for all }t\in[0,\delta_1(\omega)].
	\]
	Since \(Y_0=\rho(x)=0\) and \(U_0>0\), for every \(0< t\le \delta_1(\omega)\),
	\[
	Y_t(\omega)=\int_0^t U_s(\omega)\dd s \ \ge\ \frac 12 U_0\,t>0.
	\]
	Since \(X_t\in \cT\) and \(\tilde\rho=\rho\) on \(\cT\), we get \(\rho(X_t(\omega))>0\),
	that is, \(X_t(\omega)\notin\overline{\cD}\) for arbitrarily small \(t>0\). Hence
	\(\tau_{\cD}(\omega)=0\).
	Therefore \(\tau_{\cD}=0\) almost surely when \((x,v)\in\gamma^+\).
	
	If \((x,v)\in\gamma^0\), let \(n_0:=n(x)=\nabla\rho(x)\). Define the one-dimensional Brownian motion
	\[
	\beta_t:=n_0\cdot B_t,\qquad J_t:=\int_0^t \beta_s\dd s.
	\]
	From \eqref{eq:kolmogorov_XV},
	\[
	X_t-x = vt+\sqrt2\int_0^t B_s\dd s,
	\quad\text{hence}\quad
	n_0\cdot (X_t-x)=\sqrt2\,J_t,
	\]
	because \(v\cdot n_0=0\) on \(\gamma^0\).
	Since \(\tilde \rho\in C^{1,1}(\bR^d)\) with \(\nabla\tilde \rho(x)=n_0\),
	\[
	\big|\tilde \rho(x+h)-n_0\cdot h\big|=O(|h|^2)\quad \text{as } |h|\downarrow 0.
	\]
	Applying this with \(h=X_t-x\) gives
	\[
	\tilde \rho(X_t)=\sqrt2\,J_t+O(|X_t-x|^2)\quad \text{as }t\downarrow 0.
	\]
	
	Since \(V_\cdot\) has continuous paths, for each \(\omega\)
	\[
	|X_t(\omega)-x|
	=\left|\int_0^t V_s(\omega)\dd s\right|
	\le t\,\sup_{0\le s\le 1}|V_s(\omega)|
	=O(t)\quad\text{as }t\downarrow0.
	\]
	Hence \(|X_t-x|^2=O(t^2)\) pathwise.
	Since \(X_t\in \cT\) for all small enough \(t>0\), we conclude that
	\[
	|\rho(X_t)-\sqrt{2}J_t|=|\tilde \rho(X_t)-\sqrt{2}J_t|=O(t^2)\quad \text{as }t\downarrow 0.
	\]
	By Lemma~\ref{lem:intBM_osc}, almost surely
	\[
	\limsup_{t\downarrow0}\frac{J_t}{t^{3/2}}=+\infty.
	\]
	Since
	\[
	\frac{\rho(X_t)}{t^{3/2}}
	=
	\sqrt2\,\frac{J_t}{t^{3/2}}+e_t,
	\qquad
	|e_t|=\frac{|\rho(X_t)-\sqrt2\,J_t|}{t^{3/2}}=O(t^{1/2})\xrightarrow[t\downarrow0]{}0 ,
	\]
	we get, almost surely,
	\begin{align}\label{260621102}
		\limsup_{t\downarrow 0}\frac{\rho(X_t)}{t^{3/2}}
		=
		\sqrt2\,\limsup_{t\downarrow 0}\frac{J_t}{t^{3/2}}
		=+\infty.
	\end{align}
	Hence, for every \(\varepsilon>0\) there exists \(t\in(0,\varepsilon)\) such that \(\rho(X_t)>0\), equivalently, \(X_t\notin \overline\cD\).
	Therefore \(\tau_{\cD}=0\) almost surely when \((x,v)\in\gamma^0\).
	
	\smallskip
	(3) If \(x\in\cD\), then the continuity of \(X_t\) implies that
	\(X_t\in\cD\) for all sufficiently small \(t>0\). Hence
	\(\tau_{\cD}>0\).
	
	If \((x,v)\in\gamma^-\), then \(x\in\partial\cD\) and
	$U_0=v\cdot n(x)<0$.
	By the continuity of \(X_t\) and \(U_t\), for each \(\omega\), there exists
	\(\delta(\omega)>0\) such that
	\[
	X_t(\omega)\in \cT,
	\qquad
	U_t(\omega)\le \frac12 U_0<0
	\qquad
	\text{for all }t\in[0,\delta(\omega)].
	\]
	Since \(Y_0=\rho(x)=0\) and \(\tilde\rho=\rho\) on \(\cT\), we have
	\[
	\rho(X_t(\omega))
	=
	Y_t(\omega)
	=
	\int_0^t U_s(\omega)\dd s
	<0
	\qquad
	\text{for all }t\in(0,\delta(\omega)].
	\]
	Thus \(X_t(\omega)\in\cD\) for all sufficiently small \(t>0\).
	Therefore, $\tau_{\cD}>0$ almost surely.
	
	We next prove transversality. Assume that \(x\in\cD\) or
	\((x,v)\in\gamma^-\). Then \(\tau_{\cD}>0\) almost surely.
	Almost surely on the event \(\{\tau_{\cD}<\infty\}\) we have
	\(\tau_{\cD}>0\) and \(X_t\in\cD\) for all \(0< t<\tau_{\cD}\),
	so that \(X_{\tau_{\cD}}\in\overline{\cD}\). On the other hand,
	there are \(t_j\downarrow \tau_{\cD}\) with \(X_{t_j}\notin\cD\),
	and \(\cD\) is open, so that \(X_{\tau_{\cD}}\notin\cD\).
	Therefore
	\[
	X_{\tau_{\cD}}\in\partial\cD,
	\qquad
	Y_{\tau_{\cD}}
	=
	\rho(X_{\tau_{\cD}})
	=
	0.
	\]
	By the continuity of \(X_t\) and the openness of \(\cT\), for each
	\(\omega\in\{\tau_{\cD}<\infty\}\), there exists
	\(\delta(\omega)>0\), which we may take smaller than
	\(\tau_{\cD}(\omega)\), such that
	\[
	X_t(\omega)\in \cT
	\qquad
	\text{for all }
	t\in
	\big(\tau_{\cD}(\omega)-\delta(\omega),
	\tau_{\cD}(\omega)+\delta(\omega)\big).
	\]
	In particular, \(Y_t=\rho(X_t)\) on this time interval.
	
	We first show that
	\begin{align}\label{261006414}
		U_{\tau_{\cD}}\ge0
		\qquad
		\text{on }\{\tau_{\cD}<\infty\}.
	\end{align}
	Fix \(\omega\in\{\tau_{\cD}<\infty\}\). If
	\(U_{\tau_{\cD}}(\omega)<0\), then, by continuity, there exists
	\(\delta_1(\omega)\in(0,\delta(\omega))\) such that
	\[
	U_t(\omega)<0
	\qquad
	\text{for all }
	t\in
	(\tau_{\cD}(\omega)-\delta_1(\omega),
	\tau_{\cD}(\omega)].
	\]
	For \(t\in
	(\tau_{\cD}(\omega)-\delta_1(\omega),
	\tau_{\cD}(\omega))\), we then have
	\[
	Y_t(\omega)
	=
	Y_{\tau_{\cD}}(\omega)
	-
	\int_t^{\tau_{\cD}} U_s(\omega)\dd s
	=
	-
	\int_t^{\tau_{\cD}} U_s(\omega)\dd s
	>0.
	\]
	Hence \(\rho(X_t(\omega))>0\), so \(X_t(\omega)\notin\cD\) for some
	\(t<\tau_{\cD}(\omega)\). This contradicts the definition of
	\(\tau_{\cD}\) as the first exit time from \(\cD\). Therefore \eqref{261006414} holds.
	
	If \(U_{\tau_{\cD}}=0\), then $(X_{\tau_{\cD}},V_{\tau_{\cD}})\in\gamma^0$.
	Since the initial point satisfies \((x,v)\notin\gamma^0\), part (1) of this proof implies
	\[
	\mathbb P\left(
	\tau_{\cD}<\infty,\,
	(X_{\tau_{\cD}},V_{\tau_{\cD}})\in\gamma^0
	\right)=0.
	\]
	Hence
	$U_{\tau_{\cD}}>0$ almost surely on $\{\tau_{\cD}<\infty\}$.
	Equivalently,
	\[
	V_{\tau_{\cD}}\cdot n(X_{\tau_{\cD}})>0
	\qquad
	\text{almost surely on }\{\tau_{\cD}<\infty\}.
	\]
	This proves the transversality assertion.
\end{proof}

\begin{proof}[\textbf{Proof of Corollary~\ref{cor:exit_time_indicator_continuity}}]
	\noindent\emph{Step 1. An explicit coupling and uniform convergence of trajectories.}
	For all \(t\ge0\),
	\[
	V_t^{(n)}=v_n+\sqrt{2}\,B_t,\qquad
	X_t^{(n)}=x_n+t v_n+\sqrt{2}\int_0^t B_s\dd s,
	\]
	and similarly for \((X_t,V_t)\) started at \((x,v)\).
	Since the same Brownian motion \(B\) drives both processes, the stochastic terms
	cancel: \(X_s^{(n)}-X_s=(x_n-x)+s(v_n-v)\) and \(V_s^{(n)}-V_s=v_n-v\) for every
	\(s\ge0\). Hence, for every \(T>0\),
	\begin{equation}\label{eq:uniform_X_convergence}
		\sup_{0\le s\le T}|X_s^{(n)}-X_s|
		\le |x_n-x|+T|v_n-v|\ \xrightarrow[n\to\infty]{}\ 0,
	\end{equation}
	and likewise \(\sup_{0\le s\le T}|V_s^{(n)}-V_s|\to0\).
	
	\smallskip
	\noindent\emph{Step 2. Pathwise continuity of the exit time.}
	Let \(\rho\) be the signed distance to \(\partial\cD\) (see \eqref{2606211233}).
	Then \([\rho]_{C^{0,1}(\bR^d_x)}\le 1\). Indeed, if \(x,y\) lie on the same side of
	\(\partial\cD\), then \(|\rho(x)-\rho(y)|=|d(x,\partial\cD)-d(y,\partial\cD)|\le|x-y|\).
	If they lie on opposite sides, the segment \([x,y]\) meets \(\partial\cD\) at
	some \(z\), and \(|\rho(x)-\rho(y)|=d(x,\partial\cD)+d(y,\partial\cD)\le|x-z|+|z-y|=|x-y|\).
	Define \(Y_t:=\rho(X_t)\) and \(Y_t^{(n)}:=\rho(X_t^{(n)})\).
	By Lipschitzness of \(\rho\) and \eqref{eq:uniform_X_convergence},
	\begin{equation}\label{eq:uniform_Y_convergence}
		\sup_{0\le s\le T}|Y_s^{(n)}-Y_s|
		\le \sup_{0\le s\le T}|X_s^{(n)}-X_s|\ \xrightarrow[n\to\infty]{}\ 0
		\qquad\text{for every }T>0.
	\end{equation}
	
	Write
	\[
	\tau:=\tau_\cD(x,v)\quad \text{and}\quad \tau_n:=\tau_\cD(x_n,v_n).
	\]
	We prove \(\tau_n\to\tau\) almost surely.
	
	\smallskip
	\noindent{Case 1: \(x\in\cD\) or \((x,v)\in\gamma^-\).}
	By Theorem~\ref{thm:finiteness_exit_time} and
	Theorem~\ref{thm:exit_time_properties}(3),
	$0<\tau<\infty$ and $V_\tau\cdot n(X_\tau)>0$, almost surely.
	Fix \(\omega\) in the full-probability event where
	\[
	0<\tau(\omega)<\infty,
	\qquad
	V_\tau(\omega)\cdot n(X_\tau(\omega))>0,
	\]
	and where the uniform convergences
	\[
	\sup_{0\le s\le T}|X_s^{(n)}(\omega)-X_s(\omega)|\to0,
	\qquad
	\sup_{0\le s\le T}|V_s^{(n)}(\omega)-V_s(\omega)|\to0
	\]
	hold for every finite \(T>0\).
	
	We prove separately the lower and upper bounds for \(\tau_n(\omega)\).
	
	\medskip
	\noindent\emph{Lower bound.}
	Fix \(\eta\in(0,\tau(\omega))\). We prove that
	\[
	\tau_n(\omega)>\tau(\omega)-\eta
	\]
	for all sufficiently large \(n\).
	
	First assume \(x\in\cD\). Since \(X_s(\omega)\in\cD\) for every
	\(0\le s\le \tau(\omega)-\eta\), compactness gives
	\[
	m_\eta:=
	\min_{0\le s\le \tau(\omega)-\eta}\big(-\rho(X_s(\omega))\big)>0.
	\]
	By the uniform convergence \eqref{eq:uniform_Y_convergence}, for all large \(n\),
	\[
	\min_{0\le s\le \tau(\omega)-\eta}
	\big(-\rho(X_s^{(n)}(\omega))\big)
	\ge \frac12 m_\eta>0.
	\]
	Hence \(X_s^{(n)}(\omega)\in\cD\) for all
	\(0\le s\le \tau(\omega)-\eta\), and therefore
	$\tau_n(\omega)>\tau(\omega)-\eta$.
	
	Next assume \((x,v)\in\gamma^-\). Set
	\[
	u_0:=v\cdot n(x)<0.
	\]
	Choose \(\ell\in(0,\tau(\omega)-\eta)\) so small that, for \(0\le s\le\ell\),
	the path \(X_s(\omega)\) stays in the tubular neighborhood \(\cT\) in
	\eqref{2606211239}, and
	\[
	V_s(\omega)\cdot\nabla\rho(X_s(\omega))
	\le \frac34 u_0<0,
	\qquad 0\le s\le\ell.
	\]
	Since \(X_s(\omega)\in \cT\) for \(0\le s\le\ell\) and \(\cT\) is open, the compact set
	$\{X_s(\omega):0\le s\le\ell\}$
	has positive distance from \(\cT^c\). Hence, by the uniform convergence of
	\(X^{(n)}\) to \(X\), we have
	$X_s^{(n)}(\omega)\in \cT$ for all $0\le s\le\ell$, for all sufficiently large \(n\). Moreover, by the uniform convergence of
	\((X^{(n)},V^{(n)})\) to \((X,V)\) and the continuity of \(\nabla\rho\) in
	\(\cT\), for all sufficiently large \(n\),
	\[
	V_s^{(n)}(\omega)\cdot\nabla\rho(X_s^{(n)}(\omega))
	\le \frac12 u_0<0,
	\qquad 0\le s\le\ell.
	\]
	Since \(x_n\in\overline{\cD}\), we have \(\rho(x_n)\le0\). Therefore, for
	\(0<s\le\ell\),
	\[
	\rho(X_s^{(n)}(\omega))
	=
	\rho(x_n)
	+
	\int_0^s
	V_r^{(n)}(\omega)\cdot\nabla\rho(X_r^{(n)}(\omega))\dd r
	<0.
	\]
	Thus
	$\tau_n(\omega)>\ell$
	for all sufficiently large \(n\).
	Since \(0<\ell<\tau(\omega)-\eta\), the interval
	\([\ell,\tau(\omega)-\eta]\) is a compact subset of \((0,\tau(\omega))\), on
	which \(X_s(\omega)\in\cD\), and therefore
	\[
	m_{\eta,\ell}
	:=
	\min_{\ell\le s\le \tau(\omega)-\eta}\big(-\rho(X_s(\omega))\big)
	>0.
	\]
	By the uniform convergence \eqref{eq:uniform_Y_convergence}, for all
	sufficiently large \(n\),
	\[
	\min_{\ell\le s\le \tau(\omega)-\eta}\big(-\rho(X_s^{(n)}(\omega))\big)
	\ge \frac12 m_{\eta,\ell}>0.
	\]
	Combining the estimates on \([0,\ell]\) and
	\([\ell,\tau(\omega)-\eta]\), we obtain
	\[
	\rho(X_s^{(n)}(\omega))<0
	\qquad
	\text{for all }0<s\le\tau(\omega)-\eta.
	\]
	Consequently,
	$\tau_n(\omega)>\tau(\omega)-\eta$.
	
	\medskip
	\noindent\emph{Upper bound.}
	Set
	\[
	U_\tau(\omega):=V_\tau(\omega)\cdot n(X_\tau(\omega))>0.
	\]
	Since \(\nabla\rho(X_\tau)=n(X_\tau)\) and \((X_s,V_s)\) is continuous,
	there exists \(h\in(0,\eta)\) such that \(X_s(\omega)\) stays in the tubular
	neighborhood \(\cT\) for \(\tau(\omega)\le s\le\tau(\omega)+h\), and
	\[
	V_s(\omega)\cdot\nabla\rho(X_s(\omega))
	\ge \frac12 U_\tau(\omega)>0,
	\qquad
	\tau(\omega)\le s\le\tau(\omega)+h.
	\]
	Since \(\rho(X_{\tau(\omega)}(\omega))=0\), we get
	\[
	\rho\bigl(X_{\tau(\omega)+h}(\omega)\bigr)
	=
	\int_{\tau(\omega)}^{\tau(\omega)+h}
	V_s(\omega)\cdot\nabla\rho(X_s(\omega))\dd s
	>0.
	\]
	Let
	\[
	a:=\rho\bigl(X_{\tau(\omega)+h}(\omega)\bigr)>0.
	\]
	By the uniform convergence \eqref{eq:uniform_Y_convergence}, for all sufficiently
	large \(n\),
	$\rho\bigl(X_{\tau(\omega)+h}^{(n)}(\omega)\bigr)
	\ge \frac{a}{2}>0$.
	Thus \(X_{\tau(\omega)+h}^{(n)}(\omega)\notin\cD\), and hence
	\[
	\tau_n(\omega)\le\tau(\omega)+h<\tau(\omega)+\eta.
	\]
	
	Combining the lower and upper bounds, for every
	\(\eta\in(0,\tau(\omega))\), we have, for all sufficiently large \(n\),
	\[
	\tau(\omega)-\eta<\tau_n(\omega)<\tau(\omega)+\eta.
	\]
	Therefore,
	$\tau_n(\omega)\to\tau(\omega)$.
	
	\medskip
	\noindent{Case 2: \((x,v)\in\gamma^+\cup\gamma^0\).}
	Since \(\gamma^-\), \(\gamma^0\), \(\gamma^+\) partition
	\(\partial\cD\times\bR^d_v\), Cases~1 and~2 exhaust all
	\((x,v)\in\overline{\cD}\times\bR^d_v\). By
	Theorem~\ref{thm:exit_time_properties}(2), \(\tau=0\) almost surely. Fix \(\omega\) such that
	\(B_t(\omega)\to0\) as \(t\downarrow0\) and
	\((X^{(n)},V^{(n)})\to(X,V)\) uniformly on compact time intervals.
	
	First assume \((x,v)\in\gamma^+\), and set
	\[
	u_0:=v\cdot n(x)>0.
	\]
	Let \(\varepsilon>0\). Choose \(s=s(\omega,\varepsilon)\in(0,\varepsilon)\)
	so small that
	\[
	\sqrt{2}\sup_{0\le r\le s}|B_r(\omega)|
	\le
	\frac{u_0}{4},
	\]
	and such that \(X_r(\omega)\in \cT\) for all \(0\le r\le s\). Then, for
	all \(0\le r\le s\),
	\[
	V_r(\omega)\cdot n(x)
	\ge
	u_0-\sqrt{2}\sup_{0\le q\le s}|B_q(\omega)|
	\ge
	\frac{3u_0}{4}.
	\]
	Shrinking \(s\) further if necessary, and using the continuity of
	\(\nabla\rho\) on \(\cT\), we may also assume that
	\[
	V_r(\omega)\cdot\nabla\rho(X_r(\omega))
	\ge
	\frac{2u_0}{3},
	\qquad
	0\le r\le s.
	\]
	By the uniform convergence of \((X^{(n)},V^{(n)})\) to \((X,V)\) on
	\([0,s]\), for all sufficiently large \(n\), we have
	\[
	X_r^{(n)}(\omega)\in \cT,
	\qquad
	V_r^{(n)}(\omega)\cdot\nabla\rho(X_r^{(n)}(\omega))
	\ge
	\frac{u_0}{2},
	\qquad
	0\le r\le s.
	\]
	Moreover, since \(x_n\to x\) and \(\rho(x)=0\), for all sufficiently large
	\(n\),
	$|\rho(x_n)|\le \frac{u_0}{4}s$.
	Hence, for such \(n\),
	\[
	\rho(X_s^{(n)}(\omega))
	=
	\rho(x_n)
	+
	\int_0^s
	V_r^{(n)}(\omega)\cdot
	\nabla\rho(X_r^{(n)}(\omega))\dd r
	\ge
	-\frac{u_0}{4}s+\frac{u_0}{2}s
	>0.
	\]
	Thus \(X_s^{(n)}(\omega)\notin\cD\), and consequently
	$\tau_n(\omega)\le s<\varepsilon$
	for all sufficiently large \(n\). Since \(\varepsilon>0\) was arbitrary, we
	obtain $\tau_n(\omega)\to0=\tau(\omega)$.
	
	Next assume that \((x,v)\in\gamma^0\). Fix \(\omega\) in a full-probability
	event on which the following hold:
	\[
	\tau(\omega)=0,
	\]
	the convergence in \eqref{eq:uniform_Y_convergence} holds on every compact
	time interval, and there exists a sequence \(t_m=t_m(\omega)\downarrow0\) such that
	\begin{equation*}
		\rho\big(X_{t_m}(\omega)\big)>0
		\qquad\text{for all }m.
	\end{equation*}
	The last property follows from \eqref{260621102}, which in particular
	implies \(\tau(\omega)=0\).
	
	Fix \(s>0\) and choose \(m\) so large that \(t_m<s\).
	By the uniform convergence \eqref{eq:uniform_Y_convergence} on \([0,s]\),
	\[
	Y^{(n)}_{t_m}(\omega)
	\longrightarrow
	Y_{t_m}(\omega)
	=
	\rho\big(X_{t_m}(\omega)\big)
	>0.
	\]
	Therefore, for all sufficiently large \(n\),
	\[
	\rho\big(X^{(n)}_{t_m}(\omega)\big)>0.
	\]
	Thus \(X^{(n)}_{t_m}(\omega)\notin\cD\), and consequently
	$\tau_n(\omega)\le t_m<s$
	for all sufficiently large \(n\). Since \(s>0\) was arbitrary and
	\(\tau_n(\omega)\ge0\), we obtain
	$\tau_n(\omega)\to0=\tau(\omega)$
	in the \(\gamma^0\) case.
	
	Combining the two cases yields \eqref{eq:tau_pathwise_conv}, that is, part~(i).
	
	\smallskip
	\noindent\emph{Step 3. Almost sure convergence of the indicators.}
	We now prove part~(ii). Since
	$\tau_n\to \tau$ almost surely, and \(t_n\to t\), we have
	\[
	{\bf 1}_{\{\tau>t\}}
	\le
	\liminf_{n\to\infty}
	{\bf 1}_{\{\tau_n>t_n\}}
	\le
	\limsup_{n\to\infty}
	{\bf 1}_{\{\tau_n>t_n\}}
	\le
	{\bf 1}_{\{\tau\ge t\}}
	\]
	almost surely. Therefore it remains to prove that
	\begin{align}\label{260621404}
		{\bf 1}_{\{\tau>t\}}
		=
		{\bf 1}_{\{\tau\ge t\}},
		\quad\textit{i.e.,}\quad
		\mathbb P\big(\tau=t\big)=\mathbb P\big(\tau_{\cD}(x,v)=t\big)=0.
	\end{align}
	
	First assume \(t>0\). Then \eqref{260621404} holds. Indeed, if
	\(\tau_{\cD}(x,v)=t>0\), then \(X_s(x,v)\in\cD\) for all \(0<s<t\),
	while \(X_t(x,v)\notin\cD\) because \(\cD^c\) is closed and
	\(X\) is continuous; hence
	\[
	\{\tau_{\cD}(x,v)=t\}
	\subset
	\{X_t(x,v)\in\partial\cD\}.
	\]
	Since \(X_t(x,v)\) has the density \(f_t\) in \eqref{260621406} and
	\(\partial\cD\) has Lebesgue measure zero, we have
	\begin{align}\label{2610161048}
		\mathbb P\big(\tau_{\cD}(x,v)=t\big)
		\le
		\mathbb P\big(X_t(x,v)\in\partial\cD\big)
		=
		\int_{\bR^d}
		{\bf 1}_{\partial\cD}(y)f_t(y)\dd y
		=0.
	\end{align}
	
	It remains to consider the case \(t=0\). Since
	\[
	(t,x,v)\notin \{0\}\times(\gamma^+\cup\gamma^0),
	\]
	we have either \(x\in\cD\) or \((x,v)\in\gamma^-\). In this case,
	Theorem~\ref{thm:exit_time_properties}(3) gives \(\tau_{\cD}(x,v)>0\) almost surely.
	Hence
	$\mathbb P\big(\tau_{\cD}(x,v)=0\big)=0$,
	and \eqref{260621404} also holds for \(t=0\).
	
	Therefore,
	${\bf 1}_{\{\tau=t\}}=0$ almost surely,
	and the above two-sided bound yields \eqref{eq:indicator_pathwise_conv}.
\end{proof}

\begin{proof}[\textbf{Proof of Corollary~\ref{cor:exit_time_eps_stability}}]
	Since \(x\in \cD\), Theorems~\ref{thm:finiteness_exit_time} and
	\ref{thm:exit_time_properties}(3) give
	\(0<\tau_\cD<\infty\) almost surely.
	
	Fix \(\omega\) in the full-probability event on which
	\[
	0<\tau_\cD(\omega)<\infty,
	\qquad
	V_{\tau_\cD}(\omega)\cdot n(X_{\tau_\cD}(\omega))>0,
	\]
	and on which, for every finite \(T>0\),
	\begin{equation}\label{eq:uniform_eps_perturb}
		\sup_{0\le s\le T}|X_s^\varepsilon(\omega)-X_s(\omega)|
		\le
		\sqrt{2\varepsilon}\,\sup_{0\le s\le T}|B_s^x(\omega)|
		\longrightarrow0
		\qquad\text{as }\varepsilon\downarrow0.
	\end{equation}
	We prove
	$\tau_\cD^\varepsilon(\omega)\to\tau_\cD(\omega)$.
	
	Fix \(\eta\in(0,\tau_\cD(\omega))\).
	
	\smallskip
	\noindent\emph{Lower bound.}
	Since \(X_s(\omega)\in\cD\) for every
	\(0\le s\le\tau_\cD(\omega)-\eta\), compactness gives
	\[
	m_\eta:=
	\min_{0\le s\le\tau_\cD(\omega)-\eta}\big(-\rho(X_s(\omega))\big)>0.
	\]
	By \eqref{eq:uniform_eps_perturb} and the Lipschitz continuity of \(\rho\), for
	all sufficiently small \(\varepsilon>0\),
	\[
	\min_{0\le s\le\tau_\cD(\omega)-\eta}
	\big(-\rho(X_s^\varepsilon(\omega))\big)
	\ge \frac12 m_\eta>0.
	\]
	Hence
	$X_s^\varepsilon(\omega)\in\cD$ for all $0\le s\le\tau_\cD(\omega)-\eta$,
	and therefore
	$\tau_\cD^\varepsilon(\omega)>\tau_\cD(\omega)-\eta$.
	
	\smallskip
	\noindent\emph{Upper bound.}
	Set
	\[
	U_{\tau_\cD}(\omega):=
	V_{\tau_\cD}(\omega)\cdot n(X_{\tau_\cD}(\omega))>0.
	\]
	Since \(\nabla\rho(X_{\tau_\cD}(\omega))=n(X_{\tau_\cD}(\omega))\) and
	\((X_s,V_s)\) is continuous, there exists \(h\in(0,\eta)\) such that
	\[
	X_s(\omega)\in \cT
	\qquad
	\text{for all }\tau_\cD(\omega)\le s\le \tau_\cD(\omega)+h,
	\]
	where \(\cT\) is the tubular neighborhood of \(\partial\cD\) defined in
	\eqref{2606211239}. Since \(\nabla\rho\) is continuous on
	\(\cT_2\supset \overline{\cT}\), shrinking \(h\) if necessary, we also have
	\[
	V_s(\omega)\cdot\nabla\rho(X_s(\omega))
	\ge
	\frac12 U_{\tau_\cD}(\omega)>0,
	\qquad
	\tau_\cD(\omega)\le s\le\tau_\cD(\omega)+h.
	\]
	Since \(\rho(X_{\tau_\cD(\omega)}(\omega))=0\), we have
	\[
	\rho(X_{\tau_\cD(\omega)+h}(\omega))
	=
	\int_{\tau_\cD(\omega)}^{\tau_\cD(\omega)+h}
	V_s(\omega)\cdot\nabla\rho(X_s(\omega))\dd s
	>0.
	\]
	Set
	\[
	a:=\rho(X_{\tau_\cD(\omega)+h}(\omega))>0.
	\]
	By \eqref{eq:uniform_eps_perturb} and the Lipschitz continuity of \(\rho\), for
	all sufficiently small \(\varepsilon>0\),
	\[
	\rho(X_{\tau_\cD(\omega)+h}^\varepsilon(\omega))
	\ge \frac a2>0.
	\]
	Thus
	$X_{\tau_\cD(\omega)+h}^\varepsilon(\omega)\notin\cD$,
	and hence
	\[
	\tau_\cD^\varepsilon(\omega)
	\le
	\tau_\cD(\omega)+h
	<
	\tau_\cD(\omega)+\eta.
	\]
	
	Combining the lower and upper bounds, for every
	\(\eta\in(0,\tau_\cD(\omega))\), we have, for all sufficiently small
	\(\varepsilon>0\),
	\[
	\tau_\cD(\omega)-\eta
	<
	\tau_\cD^\varepsilon(\omega)
	<
	\tau_\cD(\omega)+\eta.
	\]
	Therefore,
	$\tau_\cD^\varepsilon(\omega)\to\tau_\cD(\omega)$.
	
	Finally, fix \(t>0\). Since $
	\tau_\cD^{\varepsilon}\to\tau_\cD$ almost surely,
	we have
	\[
	{\bf 1}_{\{\tau_\cD>t\}}
	\le
	\liminf_{\varepsilon\downarrow0}
	{\bf 1}_{\{\tau_\cD^{\varepsilon}>t\}}
	\le
	\limsup_{\varepsilon\downarrow0}
	{\bf 1}_{\{\tau_\cD^{\varepsilon}>t\}}
	\le
	{\bf 1}_{\{\tau_\cD\ge t\}}
	\]
	almost surely. Since
	$\mathbb P(\tau_\cD=t)=0$
	for every \(t>0\) (see \eqref{2610161048}), we obtain
	\eqref{260621441}.
\end{proof}

\subsection{Probabilistic solutions for the kinetic Fokker--Planck equation}\label{subsec:prob-solutions} 

In this subsection, we construct classical solutions to the dynamic and
stationary Fokker--Planck problems:
\begin{align}\label{FP_dy_classic}
	\begin{cases}
		\big(\partial_t +v\cdot\nabla_x\big) f -\Delta_v f=g
		\quad &\text{in } \cD_T\times \bR^d_v,\\
		f(0,\,\cdot\,)= f_0
		\quad &\text{in } \cD\times \bR^d_v,\\
		f=h
		\quad &\text{on } \gamma^-_T,
	\end{cases}
\end{align}
and
\begin{align}\label{FP_st_classic}
	\begin{cases}
		v\cdot\nabla_x f-\Delta_v f+\lambda f=g
		\quad &\text{in } \cD\times \bR^d_v,\\
		f=h
		\quad &\text{on } \gamma^-.
	\end{cases}
\end{align}

Since the Kolmogorov diffusion \eqref{eq:kolmogorov_XV} has backward generator
\(v\cdot\nabla_x+\Delta_v\), its probabilistic representation naturally gives
the backward equation with data on \(\gamma^+\). The desired formulation is
recovered by the reflection \(v\mapsto-v\): for
\(\check f(t,x,v):=f(t,x,-v)\) and \(\check g(t,x,v):=g(t,x,-v)\),
\begin{align}\label{2608051010}
	\big(\partial_t+v\cdot\nabla_x-\Delta_v\big)f=g
	\iff
	\big(\partial_t-v\cdot\nabla_x-\Delta_v\big)\check f=\check g ,
\end{align}
and \((x,v)\in\gamma^\mp\) if and only if \((x,-v)\in\gamma^\pm\), while
\(\gamma^0\) is invariant. We therefore prove the backward representation in
Proposition~\ref{prop:probabilistic-representation} and obtain
\eqref{FP_dy_classic} and \eqref{FP_st_classic} in
Corollary~\ref{prop_sol}. The general stochastic input, including the stopped
Feynman--Kac formula of Lemma~\ref{lem:feynman_kac}, is collected in
Appendix~\ref{sec:stochastic-appendix}; only the boundary-continuity,
approximation, regularity, and decay steps specific to the present problem are
carried out here.

For $(x,v)\in \overline{\cD}\times \bR^d_v$, we write
\[
Z_t(x,v):=\big(X_t(x,v),V_t(x,v)\big),
\qquad
\tau(x,v):=\tau_\cD(x,v),
\]
where $\big(X_t(x,v),V_t(x,v)\big)$ and $\tau_\cD(x,v)$ are defined in
\eqref{eq:kolmogorov_XV} and \eqref{260619759}, respectively.

\begin{prop}\label{prop:probabilistic-representation}
	\,
	\begin{enumerate}
		\item 
		Let
		\[
		g\in C_b([0,\infty)\times \overline{\cD}\times \bR^d_v)
		\cap C^\infty(\bR_+\times \cD\times \bR^d_v),
		\qquad
		f_0\in C_b(\overline{\cD}\times \bR^d_v),
		\qquad
		h\in C_b\big([0,\infty)\times (\gamma^+\cup \gamma^0)\big),
		\]
		where \(C_b\) denotes the space of bounded continuous functions.
		For \((t,x,v)\in [0,\infty)\times \overline{\cD}\times \bR^d_v\), define
		\begin{align}\label{prob-rep-dynamic}
			\begin{aligned}
				f(t,x,v)
				:=
				\bE\Bigg[
				&\int_0^{\tau(x,v)\wedge t}
				g\bigl(t-s,Z_s(x,v)\bigr)\dd s
				\\
				&\quad
				+h\bigl(t-\tau(x,v),Z_{\tau(x,v)}(x,v)\bigr)
				\,1_{\{\tau(x,v)\le t\}}
				+f_0\bigl(Z_t(x,v)\bigr)\,1_{\{\tau(x,v)>t\}}
				\Bigg].
			\end{aligned}
		\end{align}
		Then \(f\) is well-defined, bounded on
		\([0,T]\times \overline{\cD}\times \bR^d_v\) for every \(T>0\), belongs to
		\[
		C\bigl(([0,\infty)\times \overline{\cD}\times \bR^d_v)
		\setminus (\{0\}\times (\gamma^+\cup \gamma^0))\bigr)
		\cap C^\infty(\bR_+\times \cD\times \bR^d_v),
		\]
		and satisfies
		\begin{align}\label{prob-pde-dynamic}
			\partial_t f-v\cdot\nabla_x f
			=
			\Delta_v f + g
			\qquad\text{in }\bR_+\times \cD\times \bR^d_v
		\end{align}
		in the classical sense.
		Moreover,
		\[
		f=f_0
		\qquad\text{on }\{0\}\times
		\bigl((\overline{\cD}\times \bR^d_v)\setminus (\gamma^+\cup \gamma^0)\bigr),
		\qquad
		f=h
		\qquad\text{on }[0,\infty)\times (\gamma^+\cup \gamma^0),
		\]
		in the classical sense. If, in addition,
		\[
		f_0(x,v)=h(0,x,v)
		\qquad\text{for every }(x,v)\in \gamma^+,
		\]
		then \(f\in C([0,\infty)\times \overline{\cD}\times \bR^d_v)\).
		
		\item Let
		\[
		g\in C_b(\overline{\cD}\times \bR^d_v)
		\cap C^\infty(\cD\times \bR^d_v),
		\qquad
		h\in C_b(\gamma^+\cup\gamma^0).
		\]
		Fix \(\lambda>0\). For \((x,v)\in \overline{\cD}\times \bR^d_v\), define
		\begin{align}\label{prob-rep-stationary}
			\begin{aligned}
				f(x,v)
				:=
				\bE\Bigg[
				&\int_0^{\tau(x,v)}\ee^{-\lambda t}
				g\bigl(Z_t(x,v)\bigr)\dd t
				+\ee^{-\lambda \tau(x,v)}h\bigl(Z_{\tau(x,v)}(x,v)\bigr)
				\,1_{\{\tau(x,v)<\infty\}}
				\Bigg].
			\end{aligned}
		\end{align}
		Then \(f\) is well defined and bounded on
		\(\overline{\cD}\times \bR^d_v\). Moreover,
		\[
		f\in
		C\bigl(\overline{\cD}\times \bR^d_v\bigr)
		\cap C^\infty(\cD\times \bR^d_v),
		\]
		and \(f\) satisfies
		\[
		\lambda f-v\cdot\nabla_x f
		=
		\Delta_v f + g
		\qquad\text{in }\cD\times \bR^d_v
		\]
		in the classical sense.
		Moreover,
		\[
		f=h
		\qquad\text{on }\gamma^+\cup\gamma^0,
		\]
		in the classical sense.
	\end{enumerate}
\end{prop}

\begin{proof}
	We only prove (1), since the proof of (2) is similar.
	
	\medskip
	\noindent
	\emph{Step 1. Well-definedness and continuity.}
	It follows directly from the definition of \(f\) that
	\[
	|f(t,x,v)|
	\le
	t\|g\|_{L^\infty([0,t]\times\overline\cD\times\bR^d_v)}
	+\|h\|_{L^\infty([0,t]\times(\gamma^+\cup\gamma^0))}
	+\|f_0\|_{L^\infty(\overline\cD\times\bR^d_v)}.
	\]
	Moreover, the explicit formulas for \(X_t(x,v)\) and \(V_t(x,v)\)
	imply their continuity on
	\([0,\infty)\times \bR^d_x\times \bR^d_v\).
	On \(\{\tau(x,v)<\infty\}\) we also have
	\(Z_{\tau(x,v)}(x,v)\in\gamma^+\cup\gamma^0\) almost surely.
	Indeed, this holds by Theorem~\ref{thm:exit_time_properties}(3) if
	\(x\in\cD\) or \((x,v)\in\gamma^-\), and trivially if
	\((x,v)\in\gamma^+\cup\gamma^0\), since then \(\tau(x,v)=0\) almost surely
	by Theorem~\ref{thm:exit_time_properties}(2). Hence the \(h\)-term in
	\eqref{prob-rep-dynamic} is well defined.
	
	We first read off the values of \(f\) at \(t=0\) and on
	\(\gamma^+\cup\gamma^0\) directly from \eqref{prob-rep-dynamic}.
	If \((x,v)\in(\cD\times\bR^d_v)\cup\gamma^-\), then \(\tau(x,v)>0\) almost
	surely by Theorem~\ref{thm:exit_time_properties}(3), so that taking \(t=0\)
	in \eqref{prob-rep-dynamic} gives
	\[
	f(0,x,v)=\bE\bigl[f_0\bigl(Z_0(x,v)\bigr)\bigr]=f_0(x,v).
	\]
	If \((x,v)\in\gamma^+\cup\gamma^0\), then \(\tau(x,v)=0\) almost surely, so
	that for every \(t\ge0\)
	\[
	f(t,x,v)=\bE\bigl[h\bigl(t,Z_0(x,v)\bigr)\bigr]=h(t,x,v).
	\]
	Hence
	\[
	f=f_0
	\quad\text{on }\{0\}\times\bigl((\overline\cD\times \bR^d_v)\setminus(\gamma^+\cup\gamma^0)\bigr),
	\qquad
	f=h
	\quad\text{on }[0,\infty)\times(\gamma^+\cup\gamma^0).
	\]
	
	We prove the continuity of \(f\). Let
	\[
	(t,x,v)\in
	\big([0,\infty)\times \overline{\cD}\times\bR^d_v\big)
	\setminus
	\big(\{0\}\times(\gamma^+\cup\gamma^0)\big),
	\]
	and let
	\((t_n,x_n,v_n)\in [0,\infty)\times \overline\cD\times \bR^d_v\)
	be a sequence converging to \((t,x,v)\). Set
	\[
	\tau:=\tau(x,v),
	\qquad
	\tau_n:=\tau(x_n,v_n),
	\qquad
	Z_s:=Z_s(x,v),
	\qquad
	Z_s^n:=Z_s(x_n,v_n).
	\]
	By Corollary~\ref{cor:exit_time_indicator_continuity}, we have
	\[
	\tau_n\to\tau,
	\qquad
	\tau_n\wedge t_n\to \tau\wedge t,
	\qquad
	1_{\{\tau_n>t_n\}}
	\to
	1_{\{\tau>t\}}
	\qquad \text{almost surely.}
	\]
	Moreover, by Step~1 of the proof of
	Corollary~\ref{cor:exit_time_indicator_continuity},
	\(Z^n_\cdot\to Z_\cdot\) uniformly on compact time intervals, for every
	\(\omega\).
	We claim that each of the three terms appearing in the expectation in
	\eqref{prob-rep-dynamic} converges almost surely.
	For the first one, writing \(s_n:=(\tau_n\wedge t_n)\wedge(\tau\wedge t)\),
	\begin{align*}
		&\left|
		\int_0^{\tau_n\wedge t_n}\!\! g(t_n-s,Z^n_s)\dd s
		-\int_0^{\tau\wedge t}\!\! g(t-s,Z_s)\dd s
		\right|
		\\
		\le\,&
		\int_0^{s_n}\!\!\bigl|g(t_n-s,Z^n_s)-g(t-s,Z_s)\bigr|\dd s
		+\|g\|_{L^\infty}\bigl|\tau_n\wedge t_n-\tau\wedge t\bigr| ,
	\end{align*}
	and the right-hand side tends to \(0\). Indeed, for each fixed \(\omega\) the
	points \((t_n-s,Z^n_s)\) and \((t-s,Z_s)\), \(0\le s\le s_n\), \(n\ge1\), all
	lie in one compact subset \(\cK=\cK(\omega)\) of
	\([0,\infty)\times\overline\cD\times\bR^d_v\). The time components lie in
	\([0,\sup_n t_n]\); we have \(X^n_s\in\overline\cD\) for \(s\le\tau_n\); and
	\(V^n_s=v_n+\sqrt2 B_s\) is bounded uniformly in \(n\) on
	\([0,\sup_n t_n]\). Since \(g\) is uniformly continuous on \(\cK\),
	\(t_n\to t\), and \(Z^n\to Z\) uniformly on \([0,\sup_n t_n]\), the first
	summand tends to \(0\); the second tends to \(0\) by
	\(\tau_n\wedge t_n\to\tau\wedge t\).
	For the second and the third term, we use
	\(Z^n_{\tau_n}\to Z_\tau\) and \(Z^n_{t_n}\to Z_t\) almost surely together
	with the convergence of the indicators and the continuity of \(h\) on
	\([0,\infty)\times(\gamma^+\cup\gamma^0)\) and of \(f_0\) on
	\(\overline\cD\times\bR^d_v\).
	Since \(g\), \(h\) and \(f_0\) are bounded, these terms are uniformly dominated
	by the constant 
	$$
	(\sup_n t_n)\|g\|_{L^\infty}+\|h\|_{L^\infty}+\|f_0\|_{L^\infty}<\infty,
	$$
	so the dominated convergence theorem yields
	\(
	f(t_n,x_n,v_n)\to f(t,x,v).
	\)
	Thus
	\[
	f\in C\bigl(([0,\infty)\times \overline\cD\times \bR^d_v)
	\setminus (\{0\}\times (\gamma^+\cup \gamma^0)) \bigr).
	\]
	
	It remains to prove continuity at points of
	\(\{0\}\times(\gamma^+\cup\gamma^0)\) under the compatibility condition.
	Although the compatibility condition is assumed only on \(\gamma^+\), it
	extends to \(\gamma^+\cup\gamma^0\) by the continuity of \(f_0\) and \(h\):
	given \((x_*,v_*)\in\gamma^0\), the points \((x_*,v_*+k^{-1}n(x_*))\),
	\(k\in\bN\), belong to \(\gamma^+\) and converge to \((x_*,v_*)\). Hence
	$f_0=h(0,\cdot)
	\qquad\text{on }\gamma^+\cup\gamma^0$.
	
	Let \((x,v)\in\gamma^+\cup\gamma^0\), and let
	\((t_n,x_n,v_n)\to(0,x,v)\) with
	\((t_n,x_n,v_n)\in [0,\infty)\times\overline\cD\times\bR^d_v\). Set
	\[
	\tau_n:=\tau(x_n,v_n),
	\qquad
	Z_s^n:=Z_s(x_n,v_n).
	\]
	Since
	$0\le \tau_n\wedge t_n\le t_n\to0$,
	we have
	\[
	\int_0^{\tau_n\wedge t_n}
	g(t_n-s,Z_s^n)\dd s\to0
	\qquad\text{almost surely.}
	\]
	Moreover, by the compatibility condition,
	\[
	f_0(x,v)
	=
	h(0,x,v)1_{\{\tau_n\le t_n\}}
	+
	f_0(x,v)1_{\{\tau_n>t_n\}}.
	\]
	Therefore,
	\begin{align*}
		&\left|f(t_n,x_n,v_n)-f_0(x,v)\right|
		\\
		\le\;&
		\bE\left[
		\left|\int_0^{\tau_n\wedge t_n}
		g(t_n-s,Z_s^n)\dd s\right|
		\right]
		\\
		&+
		\bE\left[
		\big|h(t_n-\tau_n,Z_{\tau_n}^n)-h(0,x,v)\big|
		1_{\{\tau_n\le t_n\}}
		\right]
		+
		\bE\left[
		\big|f_0(Z_{t_n}^n)-f_0(x,v)\big|
		1_{\{\tau_n>t_n\}}
		\right].
	\end{align*}
	The first term tends to \(0\).
	
	Since \(0\le \tau_n\le t_n\to0\) on \(\{\tau_n\le t_n\}\), we have
	\[
	\big|h(t_n-\tau_n,Z_{\tau_n}^n)-h(0,x,v)\big|
	1_{\{\tau_n\le t_n\}}
	\to0
	\qquad\text{almost surely.}
	\]
	Similarly, since \(t_n\to0\),
	\[
	\big|f_0(Z_{t_n}^n)-f_0(x,v)\big|
	1_{\{\tau_n>t_n\}}
	\to0
	\qquad\text{almost surely.}
	\]
	By dominated convergence, the last two terms also tend to \(0\). Consequently,
	\[
	f(t_n,x_n,v_n)\to f_0(x,v)=h(0,x,v).
	\]
	Hence \(f\) is continuous at every point of
	\(\{0\}\times(\gamma^+\cup\gamma^0)\). Therefore, under the compatibility
	condition,
	$f\in C([0,\infty)\times \overline\cD\times \bR^d_v)$.

	\medskip
	\noindent
	\emph{Step 2. The distributional equation in
		\((0,\infty)\times \cD\times \bR^d_v\).}
	Let \(\{\cO_n\}_{n\ge1}\) be an exhaustion of
	\(\cD\times\bR^d_v\) by bounded \(C^2\)-domains such that
	\[
	\overline{\cO_n}\subset \cO_{n+1},
	\qquad
	\bigcup_{n=1}^\infty \cO_n=\cD\times\bR^d_v.
	\]
	Fix \(\varepsilon>0\). For \(n\in\bN\), consider the initial--boundary value
	problem
	\begin{align}\label{2608051105}
		\begin{cases}
			(\partial_t-v\cdot\nabla_x) f^{\varepsilon,n}
			-\Delta_v f^{\varepsilon,n}
			-\varepsilon\Delta_x f^{\varepsilon,n}=g
			\quad &\text{in } (0,\infty)\times \cO_n,\\
			f^{\varepsilon,n}(0,\,\cdot\,)= f_0
			\quad &\text{in } \cO_n,\\
			f^{\varepsilon,n}=f
			\quad &\text{on } (0,\infty)\times\partial\cO_n.
		\end{cases}
	\end{align}
	The source \(g\) is bounded and continuous on
	\([0,\infty)\times\overline{\cO_n}\) and smooth in
	\((0,\infty)\times\cO_n\). Since \(\overline{\cO_n}\subset\cD\times\bR^d_v\)
	and, by Step~1, \(f\) is continuous on
	\([0,\infty)\times(\cD\times\bR^d_v)\), the data \(f_0|_{\overline{\cO_n}}\)
	and \(f|_{[0,\infty)\times\partial\cO_n}\) are continuous and bounded on
	bounded time intervals. Finally, the zeroth-order compatibility condition
	holds, since \(f(0,\cdot)=f_0\) on \(\partial\cO_n\) by Step~1.
	
	Let $f^{\varepsilon,n}\in
	C\bigl([0,\infty)\times\overline{\cO_n}\bigr)
	\cap
	C^{1,2}\bigl((0,\infty)\times\cO_n\bigr)$ be the classical solution of \eqref{2608051105}.
	Let \(B^x\) be a \(d\)-dimensional Brownian motion independent of the Brownian
	motion driving \(V_t\), and define
	\[
	X_t^\varepsilon:=X_t+\sqrt{2\varepsilon}\,B_t^x,
	\qquad
	Z_t^\varepsilon:=(X_t^\varepsilon,V_t).
	\]
	For \((x,v)\in\cO_n\), set
	\[
	\tau^{\varepsilon,n}
	:=
	\inf\{s>0:\,Z_s^\varepsilon(x,v)\notin\cO_n\}.
	\]
	This is a stopping time by Lemma~\ref{lem:hitting_time}.
	The pair \((b,\sigma)\) associated with \(Z^\varepsilon\) through		\[
	\dd Z_t^\varepsilon
	=
	b(Z_t^\varepsilon)\dd t
	+\sqrt2\,\sigma\dd (B^x,B)_t,
	\qquad
	b(x,v)=(v,0),
	\qquad
	\sigma=\mathrm{diag}\bigl(\sqrt\varepsilon\,I_d,\,I_d\bigr),
	\]
	is globally Lipschitz, and \((B^x,B)\) is a \(2d\)-dimensional Brownian motion
	because \(B^x\) and \(B\) are independent. Lemma~\ref{lem:feynman_kac} is
	thus applicable with diffusion matrix \(\sqrt2\,\sigma\).
	Since \(\overline{\cO_n}\) is compact,
	\(f^{\varepsilon,n}\in C_b([0,T]\times\overline{\cO_n})\) for every finite
	\(T\). Applying the Feynman--Kac formula
	(Lemma~\ref{lem:feynman_kac}(i)) with \(\lambda=0\) on the time interval
	\([0,T]\), for an arbitrary finite \(T\), and then letting
	\(T\uparrow\infty\), we
	obtain, for
	\((t,x,v)\in(0,\infty)\times\cO_n\),
	\begin{align}\label{FK-eps-n}
		\begin{aligned}
			f^{\varepsilon,n}(t,x,v)
			=
			\bE\Bigg[
			&\int_0^{\tau^{\varepsilon,n}\wedge t}
			g\bigl(t-s,Z_s^\varepsilon(x,v)\bigr)\dd s
			\\
			&\quad
			+f\bigl(t-\tau^{\varepsilon,n},
			Z_{\tau^{\varepsilon,n}}^\varepsilon(x,v)\bigr)
			1_{\{\tau^{\varepsilon,n}\le t\}}
			+f_0\bigl(Z_t^\varepsilon(x,v)\bigr)
			1_{\{\tau^{\varepsilon,n}>t\}}
			\Bigg].
		\end{aligned}
	\end{align}
	In particular, for each fixed \(T>0\),
	\begin{align}\label{260622106}
		\sup_{(t,x,v)\in(0,T]\times\cO_n}
		|f^{\varepsilon,n}(t,x,v)|
		\le
		T\|g\|_{L^\infty([0,T]\times\overline\cD\times\bR^d_v)}+M_T,
	\end{align}
	where
	\[
	M_T:=
	\|f\|_{L^\infty([0,T]\times\overline\cD\times\bR^d_v)}
	+
	\|f_0\|_{L^\infty(\overline\cD\times\bR^d_v)} .
	\]
	The right-hand side of \eqref{260622106} does not depend on
	\(\varepsilon\) or \(n\).
	
	We now let \(n\to\infty\), along \(n\) large enough that \((x,v)\in\cO_n\). Set
	\[
	\tau^\varepsilon
	:=
	\inf\{s>0:\,Z_s^\varepsilon(x,v)\notin\cD\times\bR^d_v\}
	=
	\inf\{s>0:\,X_s^\varepsilon(x,v)\notin\cD\}.
	\]
	Since the domains \(\cO_n\) are increasing, \(\tau^{\varepsilon,n}\) is
	nondecreasing in \(n\) and \(\tau^{\varepsilon,n}\le\tau^\varepsilon\); write
	\(\tau^*:=\lim_{n\to\infty}\tau^{\varepsilon,n}\). Fix \(\omega\) with
	\(\tau^*(\omega)<\tau^\varepsilon(\omega)\). Then \(\tau^*(\omega)<\infty\) and
	\(\{Z_s^\varepsilon(x,v)(\omega):0\le s\le \tau^*(\omega)\}\) is a compact
	subset of the open set \(\cD\times\bR^d_v\); since
	\(\overline{\cO_n}\subset\cO_{n+1}\) and \(\bigcup_n\cO_n=\cD\times\bR^d_v\),
	it is contained in \(\cO_{n}\) for some \(n\), which
	forces \(\tau^{\varepsilon,n}(\omega)>\tau^*(\omega)\). This contradicts the
	definition of \(\tau^*\).
	Therefore
	$\tau^{\varepsilon,n}\uparrow \tau^\varepsilon$ almost surely.
	
	Since \(x\in\cD\) and \(X^\varepsilon\) has continuous paths, we have
	\(\tau^\varepsilon>0\) and \(X^\varepsilon_s(x,v)\in\cD\) for
	\(0\le s<\tau^\varepsilon\), while
	\(X^\varepsilon_{\tau^\varepsilon}(x,v)\notin\cD\) because \(\cD\) is open.
	Hence, on the event \(\{\tau^\varepsilon<\infty\}\),
	$X_{\tau^\varepsilon}^\varepsilon(x,v)\in\partial\cD$ almost surely.
	The same argument shows that \(\{\tau^\varepsilon=t\}\subset\{X_t^\varepsilon(x,v)\in\partial\cD\}\)
	for every \(t>0\). Moreover, for every fixed \(t>0\),
	$\bP(\tau^\varepsilon=t)=0$.
	Indeed, conditionally on \(X_t(x,v)\), the random variable
	\[
	X_t^\varepsilon(x,v)
	=
	X_t(x,v)+\sqrt{2\varepsilon}\,B_t^x
	\]
	is Gaussian with mean \(X_t(x,v)\) and covariance matrix
	\(2\varepsilon t I_d\). Hence its conditional law is absolutely continuous
	with respect to the Lebesgue measure. Since \(\partial\cD\) has Lebesgue
	measure zero,
	$\bP\bigl(X_t^\varepsilon(x,v)\in\partial\cD\bigr)=0$,
	and therefore \(\bP(\tau^\varepsilon=t)=0\).
	
	We claim that, for every \((t,x,v)\in(0,\infty)\times\cD\times\bR^d_v\),
	\[
	f^{\varepsilon,n}(t,x,v)\to f^\varepsilon(t,x,v),
	\]
	where
	\begin{align}\label{eps-representation}
		\begin{aligned}
			f^\varepsilon(t,x,v)
			:=
			\bE\Bigg[
			&\int_0^{\tau^\varepsilon\wedge t}
			g\bigl(t-s,Z_s^\varepsilon(x,v)\bigr)\dd s
			\\
			&\quad
			+f\bigl(t-\tau^\varepsilon,
			Z_{\tau^\varepsilon}^\varepsilon(x,v)\bigr)
			1_{\{\tau^\varepsilon\le t\}}
			+f_0\bigl(Z_t^\varepsilon(x,v)\bigr)
			1_{\{\tau^\varepsilon>t\}}
			\Bigg].
		\end{aligned}
	\end{align}
	Indeed, since \(\tau^{\varepsilon,n}\wedge t\uparrow\tau^\varepsilon\wedge t\) and
	\(g\) is bounded and continuous, the integral term converges. We next treat
	the terminal terms. On
	\(\{\tau^\varepsilon<t\}\) we have
	\(\tau^{\varepsilon,n}\to\tau^\varepsilon\),
	\(Z^\varepsilon_{\tau^{\varepsilon,n}}\to Z^\varepsilon_{\tau^\varepsilon}\)
	and \(\tau^{\varepsilon,n}\le t\) for all sufficiently large \(n\). Since
	\(t-\tau^\varepsilon>0\) and \(f\) is continuous on
	\((0,\infty)\times\overline\cD\times\bR^d_v\), the boundary
	term converges. On
	\(\{\tau^\varepsilon>t\}\), we have
	\(\tau^{\varepsilon,n}>t\) for all sufficiently large \(n\), and the initial
	term is unchanged. The remaining event \(\{\tau^\varepsilon=t\}\) has
	probability zero. Thus the random variable inside the expectation in
	\eqref{FK-eps-n} converges almost
	surely to the one in \eqref{eps-representation}; both are bounded in absolute
	value by the constant
	\(t\|g\|_{L^\infty}+\|f\|_{L^\infty([0,t]\times\overline\cD\times\bR^d_v)}+\|f_0\|_{L^\infty}\),
	uniformly in \(n\), so dominated convergence proves the claim.
	
	We next show that \(f^\varepsilon\) solves the regularized equation in the
	distributional sense. Let
	\[
	\phi\in C_c^\infty((0,\infty)\times\cD\times\bR^d_v).
	\]
	Choose \(n_0\) such that
	\[
	\operatorname{supp}\phi
	\subset (0,\infty)\times \cO_{n_0}.
	\]
	For every \(n\ge n_0\), since \(f^{\varepsilon,n}\) is classical and
	\(\phi\) is compactly supported in \((0,\infty)\times\cO_n\), integration by
	parts gives
	\begin{align*}
		\iiint_{(0,\infty)\times \cO_n} f^{\varepsilon,n}
		\bigl(
		-\partial_t\phi
		+v\cdot\nabla_x\phi
		-\Delta_v\phi
		-\varepsilon\Delta_x\phi
		\bigr)
		\dd v \dd x \dd t
		=
		\iiint_{(0,\infty)\times \cO_n} g\phi\dd v \dd x \dd t
		.
	\end{align*}
	Letting \(n\to\infty\), using the pointwise convergence
	\(f^{\varepsilon,n}\to f^\varepsilon\) on \(\operatorname{supp}\phi\) and the
	bound \eqref{260622106}, we obtain by dominated convergence
	\begin{align}\label{eps-distributional-equation}
		\iiint_{(0,\infty)\times \cD\times \bR^d_v} f^\varepsilon
		\bigl(
		-\partial_t\phi
		+v\cdot\nabla_x\phi
		-\Delta_v\phi
		-\varepsilon\Delta_x\phi
		\bigr)
		\dd v \dd x \dd t
		=
		\iiint_{(0,\infty)\times \cD\times \bR^d_v} g\phi\dd v \dd x \dd t
		.
	\end{align}
	
	It remains to let \(\varepsilon\to0\); we do so along an arbitrary sequence
	\(\varepsilon_k\downarrow0\). Since
	\begin{align}\label{260622104}
		X_s^\varepsilon=X_s+\sqrt{2\varepsilon}\,B_s^x
		\to X_s
		\qquad\text{locally uniformly in }s\text{ almost surely,}
	\end{align}
	Corollary~\ref{cor:exit_time_eps_stability} gives, for every
	\((t,x,v)\in(0,\infty)\times\cD\times\bR^d_v\),
	\begin{align}\label{260622105}
		\tau^\varepsilon\wedge t\to \tau\wedge t,
		\qquad
		1_{\{\tau^\varepsilon>t\}}\to1_{\{\tau>t\}}
		\qquad\text{almost surely,}
	\end{align}
	where
	\[
	\tau:=\inf\{s>0:\,Z_s(x,v)\notin\cD\times\bR^d_v\}.
	\]
	We also use that, for every fixed \(t>0\),
	$\bP(\tau=t)=0$.
	Indeed, \(\{\tau=t\}\subset\{X_t(x,v)\in\partial\cD\}\) by the continuity of
	\(X\), while by Remark~\ref{2606211211} the law of \(X_t(x,v)\) is absolutely
	continuous and \(|\partial\cD|=0\).
	
	We rewrite the original representation of \(f\) in a form compatible with
	\eqref{eps-representation}. Since, by Step~1, \(f=h\) on
	\([0,\infty)\times(\gamma^+\cup\gamma^0)\) and, almost surely on
	\(\{\tau<\infty\}\), \(Z_\tau(x,v)\in\gamma^+\cup\gamma^0\), we may replace
	\(h\) by \(f\) in \eqref{prob-rep-dynamic} and obtain
	\begin{align}\label{original-terminal-form}
		\begin{aligned}
			f(t,x,v)
			=
			\bE\Bigg[
			&\int_0^{\tau\wedge t}
			g\bigl(t-s,Z_s(x,v)\bigr)\dd s
			\\
			&\quad
			+f\bigl(t-\tau,Z_\tau(x,v)\bigr)
			1_{\{\tau\le t\}}
			+
			f_0\bigl(Z_t(x,v)\bigr)
			1_{\{\tau>t\}}
			\Bigg].
		\end{aligned}
	\end{align}
	We now compare \eqref{eps-representation} with \eqref{original-terminal-form}
	and argue exactly as in the passage \(n\to\infty\), with
	\eqref{260622104}, \eqref{260622105} and \(\bP(\tau=t)=0\) in place of
	\(\tau^{\varepsilon,n}\uparrow\tau^\varepsilon\) and
	\(\bP(\tau^\varepsilon=t)=0\). We also use that, on \(\{\tau<t\}\),
	\[
	\bigl|Z^\varepsilon_{\tau^\varepsilon}-Z_\tau\bigr|
	\le
	\sqrt{2\varepsilon}\sup_{0\le s\le t}|B^x_s|
	+\bigl|Z_{\tau^\varepsilon}-Z_\tau\bigr|
	\longrightarrow0
	\]
	almost surely, by \eqref{260622105} and the continuity of \(Z\). Therefore
	\[
	f^\varepsilon(t,x,v)\to f(t,x,v)
	\qquad\text{for every }(t,x,v)\in(0,\infty)\times\cD\times\bR^d_v.
	\]
	Moreover, for each fixed \(T>0\), \eqref{260622106} and the convergence
	\(f^{\varepsilon,n}\to f^\varepsilon\) give
	\[
	\sup_{(t,x,v)\in(0,T]\times\cD\times\bR^d_v}
	|f^\varepsilon(t,x,v)|
	\le
	T\|g\|_{L^\infty([0,T]\times\overline\cD\times\bR^d_v)}+M_T.
	\]
	Thus, by dominated convergence,
	$f^\varepsilon\to f$
	in $L^1_{\mathrm{loc}}
	\bigl((0,\infty)\times\cD\times\bR^d_v\bigr)$.
	
	Since \(|f^\varepsilon|\) is bounded on \(\supp\phi\) uniformly in
	\(\varepsilon\), the term
	\(\varepsilon\iiint f^\varepsilon\Delta_x\phi\) vanishes in the limit.
	Letting \(\varepsilon\to0\) in \eqref{eps-distributional-equation}, we obtain
	\begin{align*}
		\iiint_{(0,\infty)\times \cD\times \bR^d_v} f
		\bigl(
		-\partial_t\phi
		+v\cdot\nabla_x\phi
		-\Delta_v\phi
		\bigr)
		\dd v \dd x \dd t
		=
		\iiint_{(0,\infty)\times \cD\times \bR^d_v} g\phi \dd v \dd x \dd t.
	\end{align*}
	Since \(\phi\in C_c^\infty((0,\infty)\times\cD\times\bR^d_v)\) was arbitrary,
	\(f\) satisfies
	\[
	\partial_t f-v\cdot\nabla_x f=\Delta_v f+g
	\]
	in \((0,\infty)\times\cD\times\bR^d_v\) in the distributional sense.

	\medskip
	\noindent
	\emph{Step 3. Smoothness in \(\bR_+\times \cD\times \bR^d_v\).}
	Set
	\[
	\mathsf X_0:=\partial_t-v\cdot\nabla_x,
	\qquad
	\mathsf X_i:=\partial_{v_i},
	\qquad 1\le i\le d.
	\]
	Then
	$[\mathsf X_0,\mathsf X_i]=\partial_{x_i}$, $1\le i\le d$.
	Hence
	\[
	\mathsf X_0+\sum_{i=1}^d v^i[\mathsf X_0,\mathsf X_i]
	=
	\partial_t.
	\]
	Therefore the Lie algebra generated by
	$\mathsf X_0,\mathsf X_1,\dots,\mathsf X_d$
	spans all directions
	$\partial_t$, $\partial_{x_1}$, $\dots$, $\partial_{x_d}$,
	$\partial_{v_1}$, $\dots$, $\partial_{v_d}$
	at every point of \(\bR_+\times \cD\times \bR^d_v\).
	Thus the operator
	\[
	\partial_t-v\cdot\nabla_x-\Delta_v
	=
	\mathsf X_0-\sum_{i=1}^d \mathsf X_i^2
	\]
	is hypoelliptic in the sense of H\"ormander \cite{Hormander1967}; indeed
	\(\mathsf X_0-\sum_i\mathsf X_i^2=-\bigl(\sum_i\mathsf X_i^2+\mathsf Y_0\bigr)\) with \(\mathsf Y_0:=-\mathsf X_0\), and
	\(\mathsf Y_0,\mathsf X_1,\dots,\mathsf X_d\) generate the same Lie algebra as \(\mathsf X_0,\mathsf X_1,\dots,\mathsf X_d\).
	
	Since \(g\in C^\infty(\bR_+\times \cD\times \bR^d_v)\) and \(f\) is a
	distributional solution of
	\[
	(\partial_t-v\cdot\nabla_x-\Delta_v)f=g
	\]
	in \(\bR_+\times \cD\times \bR^d_v\), H\"ormander's hypoellipticity theorem \cite{Hormander1967}
	yields
	\[
	f\in C^\infty(\bR_+\times \cD\times \bR^d_v).
	\]
	Consequently, \eqref{prob-pde-dynamic} holds pointwise in
	\(\bR_+\times \cD\times \bR^d_v\).
	This completes the proof.
\end{proof}

\begin{corollary}\label{prop_sol}
	\,
	
	\begin{enumerate}
		\item 
		Let
		\[
		g\in C_c\left([0,\infty)\times \overline{\cD}\times \bR^d_v\right)
		\cap C^\infty\left((0,\infty)\times \cD\times \bR^d_v\right),
		\]
		and assume that there exists \(\varepsilon_0>0\) such that
		\begin{align}\label{260622943}
			g\equiv 0
			\quad\text{on}\quad
			[0,\varepsilon_0]\times \overline\cD\times \bR^d_v.
		\end{align}
		Then there exists
		\[
		f\in C\big([0,\infty)\times \overline{\cD}\times \bR^d_v\big)
		\cap C^\infty\big((0,\infty)\times \cD\times \bR^d_v\big)
		\]
		such that
		\begin{align}\label{260622911}
			\begin{cases}
				\big(\partial_t+v\cdot\nabla_x\big)f-\Delta_v f=g
				&\quad \text{in } (0,\infty)\times \cD\times \bR^d_v,\\
				f(0,\cdot)=0
				&\quad \text{in } \cD\times \bR^d_v,\\
				f=0
				&\quad \text{on } (0,\infty)\times \gamma^-
			\end{cases}
		\end{align}
		holds pointwise.
		Moreover, for every \(T>0\), there exists \(C_T>0\), independent of
		\((t,x,v)\), such that
		\begin{align}\label{260418620}
			|f(t,x,v)|
			\le
			C_T\,1_{[\varepsilon_0\wedge T,T]}(t)
			\exp\left(
			-\frac{|x|^{2}+|v|^2}{C_T}
			\right)
		\end{align}
		for all \((t,x,v)\in [0,T]\times \overline{\cD}\times \bR^d_v\).
		
		\item 
		Let
		\[
		g\in C_c\left(\overline{\cD}\times \bR^d_v\right)
		\cap C^\infty\left(\cD\times \bR^d_v\right).
		\]
		For every \(\lambda>0\), there exists
		\[
		f\in C\big(\overline{\cD}\times \bR^d_v\big)
		\cap C^{\infty}\big(\cD\times \bR^d_v\big)
		\]
		such that
		\begin{align*}
			\begin{cases}
				v\cdot\nabla_x f -\Delta_v f+\lambda f=g
				&\quad \text{in } \cD\times \bR^d_v,\\
				f=0
				&\quad \text{on } \gamma^-,
			\end{cases}
		\end{align*}
		holds pointwise. Moreover, there exists \(C>0\), independent of \((x,v)\),
		such that
		\begin{align}\label{260418621}
			|f(x,v)|
			\le
			C\exp\left(-\frac{|x|^{1/2}+|v|}{C}\right)
		\end{align}
		for all \((x,v)\in \overline{\cD}\times \bR^d_v\).
	\end{enumerate}
\end{corollary}

\begin{proof}
	We first prove (1). Let \(\varepsilon_0>0\) be the constant in
	\eqref{260622943}. By \eqref{2608051010}, the substitution
	\(f\mapsto \check f(t,x,v):=f(t,x,-v)\) turns \eqref{260622911} into the
	problem \eqref{260622915} below with source \(\check g(t,x,v):=g(t,x,-v)\).
	The source \(\check g\) again lies in
	\(C_c([0,\infty)\times\overline\cD\times\bR^d_v)\cap C^\infty((0,\infty)\times\cD\times\bR^d_v)\),
	again vanishes on \([0,\varepsilon_0]\times\overline\cD\times\bR^d_v\), and has
	its support contained in the same ball \(K\) used below. Renaming the
	transformed source as \(g\), it therefore suffices to consider the problem
	\begin{align}\label{260622915}
		\begin{cases}
			\big(\partial_t-v\cdot\nabla_x\big)f-\Delta_v f=g
			&\quad \text{in } (0,\infty)\times \cD\times \bR^d_v,\\
			f(0,\cdot)=0
			&\quad \text{in } \cD\times \bR^d_v,\\
			f=0
			&\quad \text{on } (0,\infty)\times \gamma^+.
		\end{cases}
	\end{align}
	By Proposition~\ref{prop:probabilistic-representation},
	\begin{align}\label{260622917}
		f(t,x,v)
		:=
		\bE\left[
		\int_0^{\tau(x,v)\wedge t}
		g\bigl(t-s,Z_s(x,v)\bigr)\dd s
		\right]
	\end{align}
	belongs to
	\[
	C\big([0,\infty)\times \overline{\cD}\times \bR^d_v\big)
	\cap C^\infty\big((0,\infty)\times \cD\times \bR^d_v\big)
	\]
	and satisfies \eqref{260622915} pointwise. It remains to prove the decay
	estimate.
	
	By \eqref{260622943} and \eqref{260622917},
	$f(t,x,v)=0$
	for all 
	$(t,x,v)\in [0,\varepsilon_0]\times\overline\cD\times\bR^d_v$.
	
	Fix \(T>0\). Since
	$g\in C_c\big([0,\infty)\times \overline\cD\times\bR^d_v\big)$,
	there exists \(R>0\) such that
	\[
	\operatorname{supp} g
	\subset
	[0,\infty)\times K,
	\qquad
	K:=\{(y,w)\in \bR^d_y\times \bR^d_w:\ |y|+|w|\le R\}.
	\]
	Let \(\bar g\) denote the extension of \(|g|\) by \(0\) to
	\([0,\infty)\times\bR^d_y\times\bR^d_w\); then
	\(\operatorname{supp}\bar g\subset[0,\infty)\times K\), and since
	\(Z_s(x,v)\in\overline\cD\times\bR^d_v\) for \(s\le\tau(x,v)\), we have
	\(1_{\{s\le\tau(x,v)\}}|g(t-s,Z_s(x,v))|\le\bar g(t-s,Z_s(x,v))\).
	Hence, for \(t\in[0,T]\), using \eqref{260622917} and
	Remark~\ref{2606211211}, we obtain, with \(P_s^{(d)}\) denoting the kernel
	in \eqref{2606211207},
	\begin{align}\label{260622922}
		\begin{aligned}
			|f(t,x,v)|
			&\le
			\int_0^t
			\iint_{\bR^d\times\bR^d}
			P_s^{(d)}(y-x-sv,w-v)
			\bar g(t-s,y,w)\dd w\dd y\dd s
			\\
			&\le
			\|g\|_{L^\infty}
			\int_0^t
			\iint_K
			P_s^{(d)}(y-x-sv,w-v)\dd w\dd y\dd s.
		\end{aligned}
	\end{align}
	By the explicit form of \(P_s^{(d)}\), there exists \(C>0\), depending only
	on \(d\), such that
	\[
	P_s^{(d)}(y-x-sv,w-v)
	\le
	C s^{-2d}\ee^{-E_s(x,v;y,w)},
	\]
	where
	\[
	E_s(x,v;y,w)
	:=
	\frac{|v-w|^2}{4s}
	+
	\frac{3\left|y-x-\frac{s}{2}(v+w)\right|^2}{s^3}.
	\]
	
	We claim that there exist constants \(c>0\) and \(C>0\), depending only on
	\(T\), \(R\), and \(d\), such that for every \(s\in(0,T]\) and every
	\((y,w)\in K\),
	\begin{align}\label{260622923}
		E_s(x,v;y,w)
		\ge
		c\bigl(|x|^2+|v|^2\bigr)-C.
	\end{align}
	Indeed, since \(s\le T\), \(|w|\le R\), and
	\(|v-w|^2\ge\frac12|v|^2-|w|^2\),
	\[
	E_s(x,v;y,w)
	\ge
	\frac{|v-w|^2}{4s}
	\ge
	\frac{|v-w|^2}{4T}
	\ge
	\frac{|v|^2}{8T}-\frac{R^2}{4T}.
	\]
	Hence
	\begin{align}\label{260622924}
		|v|^2
		\le
		C\bigl(1+E_s(x,v;y,w)\bigr).
	\end{align}
	Moreover, using \(|y|\le R\), \(|w|\le R\), and \(s\le T\), we have
	\[
	|x|
	\le
	\left|y-x-\frac{s}{2}(v+w)\right|
	+C(1+|v|).
	\]
	Therefore, using \(s\le T\) once more (so that \(1\le T^3/s^3\)),
	\begin{align*}
		|x|^2
		&\le
		C\left(
		1+
		\left|y-x-\frac{s}{2}(v+w)\right|^2
		+
		|v|^2
		\right)
		\\
		&\le
		C\left(
		1+
		\frac{\left|y-x-\frac{s}{2}(v+w)\right|^2}{s^3}
		+
		|v|^2
		\right)\le
		C\bigl(1+E_s(x,v;y,w)\bigr),
	\end{align*}
	where the constants may depend on \(T\). Combining this with
	\eqref{260622924}, we prove \eqref{260622923}.
	
	It follows from \eqref{260622923} that
	\[
	\ee^{-E_s(x,v;y,w)}
	\le
	C
	\exp\left(
	-\frac{|x|^2+|v|^2}{C}
	\right)
	\ee^{-\frac12 E_s(x,v;y,w)}.
	\]
	Consequently,
	\begin{align*}
		&\iint_K
		P_s^{(d)}(y-x-sv,w-v)\dd w\dd y
		\\
		&\qquad\le
		C
		\exp\left(
		-\frac{|x|^2+|v|^2}{C}
		\right)
		\iint_{\bR^d\times\bR^d}
		s^{-2d}\ee^{-\frac12 E_s(x,v;y,w)}
		\dd w\dd y.
	\end{align*}
	The last integral is uniformly bounded for \(s\in(0,T]\). Indeed, under the
	change of variables
	\[
	\eta:=w-v,
	\qquad
	\xi:=y-x-\frac{s}{2}(v+w),
	\]
	whose Jacobian is equal to one, we have
	\begin{align*}
		\begin{aligned}
			&\iint_{\bR^d\times\bR^d}
			s^{-2d}\ee^{-\frac12 E_s(x,v;y,w)}
			\dd w\dd y
			\\
			&\qquad=
			\iint_{\bR^d\times\bR^d}
			s^{-2d}
			\exp\left(
			-\frac{|\eta|^2}{8s}
			-\frac{3|\xi|^2}{2s^3}
			\right)
			\dd\eta\dd\xi
			\le C(d).
		\end{aligned}
	\end{align*}
	Therefore, for every \(s\in(0,T]\),
	\[
	\iint_K
	P_s^{(d)}(y-x-sv,w-v)\dd w\dd y
	\le
	C
	\exp\left(
	-\frac{|x|^2+|v|^2}{C}
	\right),
	\]
	where \(C>0\) may depend on \(T\), \(R\), and \(d\), but is independent of
	\((s,x,v)\). Substituting this into \eqref{260622922}, we obtain
	\[
	|f(t,x,v)|
	\le
	C
	\exp\left(
	-\frac{|x|^2+|v|^2}{C}
	\right)
	\]
	for all \((t,x,v)\in[0,T]\times\overline\cD\times\bR^d_v\), where \(C>0\)
	is independent of \((t,x,v)\). Since \(f(t,x,v)=0\) for \(0\le t\le
	\varepsilon_0\), we conclude that
	\[
	|f(t,x,v)|
	\le
	C\,1_{[\varepsilon_0\wedge T,T]}(t)
	\exp\left(
	-\frac{|x|^2+|v|^2}{C}
	\right)
	\]
	for all \((t,x,v)\in[0,T]\times\overline\cD\times\bR^d_v\). Reversing the sign
	of \(v\) back leaves \(|x|^2+|v|^2\) invariant and turns the
	boundary condition on \(\gamma^+\) into one on \(\gamma^-\). We thus obtain a
	solution of \eqref{260622911} satisfying \eqref{260418620}.
	
	We next prove (2). As in part~(1), reversing the sign of \(v\) preserves
	\(C_c(\overline\cD\times\bR^d_v)\cap C^\infty(\cD\times\bR^d_v)\) and the
	support ball \(K\), and exchanges \(\gamma^-\) with \(\gamma^+\). Renaming the
	transformed source as \(g\), it therefore suffices to consider the problem
	\begin{align}\label{260622930}
		\begin{cases}
			-v\cdot\nabla_x f-\Delta_v f+\lambda f=g
			&\quad \text{in } \cD\times \bR^d_v,\\
			f=0
			&\quad \text{on } \gamma^+.
		\end{cases}
	\end{align}
	By the stationary part of Proposition~\ref{prop:probabilistic-representation},
	for every \(\lambda>0\), the function
	\begin{align}\label{260622931}
		f(x,v)
		:=
		\bE\left[
		\int_0^{\tau(x,v)}
		\ee^{-\lambda s}g\bigl(Z_s(x,v)\bigr)\dd s
		\right]
	\end{align}
	belongs to
	$C\big(\overline{\cD}\times \bR^d_v\big)
	\cap C^\infty\big(\cD\times \bR^d_v\big)$
	and satisfies \eqref{260622930} pointwise. It remains to prove the decay estimate.
	
	Since
	$g\in C_c\big(\overline\cD\times\bR^d_v\big)$,
	there exists \(R>0\) such that
	\[
	\operatorname{supp}g\subset K,
	\qquad
	K:=\{(y,w)\in\bR^d_y\times\bR^d_w:\ |y|+|w|\le R\}.
	\]
	Extending \(|g|\) by \(0\) outside \(\overline\cD\times\bR^d_v\) as in part~(1),
	and using \eqref{260622931} together with Remark~\ref{2606211211}, we obtain,
	with \(P_s^{(d)}\) denoting the kernel in \eqref{2606211207},
	\begin{align}\label{260622936}
		\begin{aligned}
			|f(x,v)|
			&\le
			\|g\|_{L^\infty}
			\int_0^\infty \ee^{-\lambda s}
			\iint_K
			P_s^{(d)}(y-x-sv,w-v)\dd w\dd y\dd s.
		\end{aligned}
	\end{align}
	By the explicit form of \(P_s^{(d)}\), there exists \(C>0\), depending only
	on \(d\), such that
	\[
	P_s^{(d)}(y-x-sv,w-v)
	\le
	C s^{-2d}\ee^{-E_s(x,v;y,w)},
	\]
	where
	\[
	E_s(x,v;y,w)
	:=
	\frac{|v-w|^2}{4s}
	+
	\frac{3\left|y-x-\frac{s}{2}(v+w)\right|^2}{s^3}.
	\]
	Set
	\[
	A_s(x,v;y,w):=\lambda s+E_s(x,v;y,w).
	\]
	We claim that there exist constants \(c>0\) and \(C>0\), depending only on
	\(\lambda\) and \(R\), such that for every \(s>0\) and every \((y,w)\in K\),
	\begin{align}\label{260622937}
		A_s(x,v;y,w)
		\ge
		c\bigl(|x|^{1/2}+|v|\bigr)-C.
	\end{align}
	Indeed, since \(|w|\le R\), we have
	\[
	\lambda s+\frac{|v-w|^2}{4s}
	\ge
	\sqrt{\lambda}\,|v-w|
	\ge
	\sqrt{\lambda}\,|v|-\sqrt{\lambda}R.
	\]
	Thus
	\begin{align}\label{260622938}
		A_s(x,v;y,w)\ge c|v|-C.
	\end{align}
	To estimate the \(x\)-variable, put
	\[
	\eta:=w-v,
	\qquad
	\xi:=y-x-\frac{s}{2}(v+w).
	\]
	Then
	\[
	x
	=
	y-\frac{s}{2}(v+w)-\xi
	=
	y-sw+\frac{s}{2}\eta-\xi.
	\]
	Since \((y,w)\in K\), it follows that
	\[
	|x|
	\le
	|\xi|+\frac{s}{2}|\eta|+R(1+s).
	\]
	On the other hand,
	$A_s\ge \lambda s$, 
	$A_s\ge \frac{|\eta|^2}{4s}$, 
	and $A_s\ge \frac{3|\xi|^2}{s^3}$.
	Therefore,
	\[
	s\le \lambda^{-1}A_s,
	\qquad
	|\eta|\le 2s^{1/2}A_s^{1/2},
	\qquad
	|\xi|\le \frac{1}{\sqrt 3}s^{3/2}A_s^{1/2}.
	\]
	Consequently,
	\[
	|\xi|+s|\eta|
	\le
	C s^{3/2}A_s^{1/2}
	\le
	C A_s^2,
	\]
	where \(C>0\) depends only on \(\lambda\). Since \(s\le C A_s\), we get
	$|x|\le C(1+A_s^2)$.
	Hence
	\begin{align}\label{260622939}
		A_s(x,v;y,w)\ge c|x|^{1/2}-C.
	\end{align}
	Combining \eqref{260622938} and \eqref{260622939}, and adjusting the
	constants, we obtain \eqref{260622937}.
	
	It follows from \eqref{260622937} that
	\[
	\ee^{-A_s(x,v;y,w)}
	\le
	C
	\exp\left(
	-\frac{|x|^{1/2}+|v|}{C}
	\right)
	\ee^{-\frac12 A_s(x,v;y,w)}.
	\]
	Therefore,
	\begin{align*}
		&\int_0^\infty \ee^{-\lambda s}
		\iint_K
		P_s^{(d)}(y-x-sv,w-v)\dd w\dd y\dd s
		\\
		&\qquad\le
		C
		\exp\left(
		-\frac{|x|^{1/2}+|v|}{C}
		\right)
		\int_0^\infty \ee^{-\frac{\lambda s}{2}}
		\iint_{\bR^d\times\bR^d}
		s^{-2d}\ee^{-\frac12 E_s(x,v;y,w)}
		\dd w\dd y\dd s.
	\end{align*}
	The last integral is bounded uniformly in \((x,v)\). Indeed, for each
	\(s>0\), under the change of variables
	\[
	\eta:=w-v,
	\qquad
	\xi:=y-x-\frac{s}{2}(v+w),
	\]
	whose Jacobian is equal to one, we have
	\[
	\iint_{\bR^d\times\bR^d}
	s^{-2d}\ee^{-\frac12 E_s(x,v;y,w)}
	\dd w\dd y
	=
	\iint_{\bR^d\times\bR^d}
	s^{-2d}
	\exp\left(
	-\frac{|\eta|^2}{8s}
	-\frac{3|\xi|^2}{2s^3}
	\right)
	\dd\eta\dd\xi
	\le C(d).
	\]
	Hence
	\[
	\int_0^\infty \ee^{-\frac{\lambda s}{2}}
	\iint_{\bR^d\times\bR^d}
	s^{-2d}\ee^{-\frac12 E_s(x,v;y,w)}
	\dd w\dd y\dd s
	\le C(d,\lambda).
	\]
	Substituting this into \eqref{260622936}, we obtain
	\[
	|f(x,v)|
	\le
	C
	\exp\left(
	-\frac{|x|^{1/2}+|v|}{C}
	\right)
	\]
	for all \((x,v)\in\overline\cD\times\bR^d_v\), where \(C>0\) depends on
	\(\lambda\), \(R\), \(d\), and \(\|g\|_{L^\infty}\), but is independent of
	\((x,v)\). Reversing the sign of \(v\) back, we obtain a solution to the
	original stationary problem with boundary condition on \(\gamma^-\), and the
	same estimate holds. This proves \eqref{260418621}.
	
\end{proof}

\subsection{Proof of Theorem~\ref{solv}: Well-posedness for the kinetic Fokker--Planck equation}
\label{260616916}

We prove Theorem~\ref{solv} in two stages. We first obtain well-posedness for the constant-diffusion operator from the classical solutions constructed above, using energy estimates and density arguments. Uniform a priori estimates and the method of continuity then extend the result to the general measurable-coefficient operator.

\begin{prop}[Well-posedness for the classical kinetic Fokker--Planck equation]
	\label{solv_Del}\,
	
	\begin{enumerate}
		\item For every \(g\), \(f_0\), and \(h\) satisfying \eqref{260617519},
		\eqref{FP_dy_classic} admits a unique solution in
		\(\cH_{\kin}(\cD_T\times \bR^d_v)\), satisfying estimate
		\eqref{dyHbdd}.
		
		\item For any \(\lambda_0>0\) and \(\lambda\ge \lambda_0\), and for every
		\(g\) and \(h\) satisfying \eqref{260617520}, \eqref{FP_st_classic} admits a
		unique solution in
		\(\cH_{\kin}(\cD\times \bR^d_v)\), satisfying estimate
		\eqref{stHbdd}.
	\end{enumerate}
\end{prop}

\begin{proof}
	(1) \noindent\emph{Step 1. Estimate \eqref{dyHbdd} and uniqueness.}
	
	We first prove estimate \eqref{dyHbdd}. Let
	\(f\in \cH_{\kin}(\cD_T\times\bR_v^d)\) be a solution of
	\eqref{FP_dy_classic}. Then, by the integration-by-parts formula
	\eqref{ibp_dy}, for every \(t\in[0,T]\) we have
	\begin{align*}
		&2\iint_{(0,t)\times \cD}
		\Big\la
		\Big(\p_s+v\cdot \nabla_x-\Delta_v\Big) f,f
		\Big\ra_v\dd x\dd s
		\\
		=\,&
		-2\iint_{(0,t)\times \cD}\la \Delta_v f,f\ra_v\dd x\dd s
		+\|f(t,\cdot)\|_{L^2(\cD\times \bR^d_v)}^2
		-\|f(0,\cdot)\|_{L^2(\cD\times \bR^d_v)}^2
		\\
		&\qquad
		+\iint_{(0,t)\times \gamma}
		|f_\gamma|^2\sgn(n_x\cdot v)\dd \gamma \dd s .
	\end{align*}
	Since \(f_\gamma=h\) on \(\gamma_T^-\), it follows that
	\begin{align}\label{260617551}
		\begin{aligned}
			&\|\nabla_v f\|_{L^2(\cD_T\times \bR^d_v)}^2
			+\sup_{0\le t\le T}\|f(t,\cdot)\|_{L^2(\cD\times \bR^d_v)}^2
			+\|f_{\gamma}\|_{L^2(\gamma_T)}^2
			\\
			\le\,&
			C\left(
			\|
			\la
			g,f
			\ra_v\|_{L^1(\cD_T)}
			+\|f_0\|_{L^2(\cD\times \bR^d_v)}^2
			+\|h\|_{L^2(\gamma^-_T)}^2
			\right)
			\\
			\le\,&
			C\left(
			\|g\|_{L^2(\cD_T;H^{-1}_v)}
			\|f\|_{L^2(\cD_T;H^1_v)}
			+\|f_0\|_{L^2(\cD\times \bR^d_v)}^2
			+\|h\|_{L^2(\gamma^-_T)}^2
			\right),
		\end{aligned}
	\end{align}
	where \(C\) is a universal constant. Moreover,
	\begin{align}\label{260617552}
		\|f\|_{L^2(\cD_T;H^1_v)}^2
		\le
		(1+T)\left(
		\|\nabla_v f\|_{L^2(\cD_T\times \bR^d_v)}^2
		+\sup_{0\le t\le T}
		\|f(t,\cdot)\|_{L^2(\cD\times \bR^d_v)}^2
		\right).
	\end{align}
	Combining \eqref{260617551} and \eqref{260617552} and applying Young's
	inequality, we obtain
	\[
	\|f\|_{L^2(\cD_T;H^1_v)}
	\le
	C(T+1)\left(
	\|g\|_{L^2(\cD_T;H_v^{-1})}
	+\|f_0\|_{L^2(\cD\times \bR_v^d)}
	+\|h\|_{L^2(\gamma_T^-)}
	\right),
	\]
	where \(C\) is a universal constant. Using the equation, we also have
	\begin{align*}
		\begin{aligned}
			\left\|(\p_t+v\cdot \nabla_x)f\right\|_{L^2(\cD_T;H^{-1}_v)}
			&\le
			\left\|\Delta_v f\right\|_{L^2(\cD_T;H^{-1}_v)}
			+\left\|g\right\|_{L^2(\cD_T;H^{-1}_v)}
			\\
			&\le
			\|\nabla_v f\|_{L^2(\cD_T\times \bR^d_v)}
			+\|g\|_{L^2(\cD_T;H^{-1}_v)}.
		\end{aligned}
	\end{align*}
	Since \(f(0,\cdot)=f_0\) and \(f_{\gamma^-}=h\), the above estimates yield
	\eqref{dyHbdd}.
	
	In particular, if \(g\), \(f_0\), and \(h\) vanish identically, then
	\(f=0\). Hence uniqueness follows immediately.

	\noindent\emph{Step 2. Existence.}
	
	\noindent\emph{Step 2.1.}
	We first treat the case
	\(g\in C_c^\infty(\cD_T\times \bR^d_v)\), \(f_0\equiv0\), and \(h\equiv0\).
	Fix such a \(g\) and extend it by zero to
	\([0,\infty)\times\overline\cD\times\bR^d_v\). Since \(\supp g\) is a compact
	subset of \(\cD_T\times\bR^d_v\), there is \(\delta\in(0,T)\) such that
	\eqref{260622943} holds with \(\delta\) in place of \(\varepsilon_0\). Let
	\[
	f\in C([0,T]\times\overline{\cD}\times\bR_v^d)
	\cap C^\infty((0,T]\times\cD\times\bR_v^d)
	\]
	be the function given by Corollary~\ref{prop_sol}.(1). Then \(f\) satisfies
	\eqref{FP_dy_classic} with \(f_0\equiv0\) and \(h\equiv0\) in the classical
	sense, and \eqref{260418620} holds. In particular \(f\) vanishes identically
	on \([0,\delta)\times\overline\cD\times\bR^d_v\), so that
	\[
	f\in C^\infty\big([0,T]\times\cD\times\bR^d_v\big).
	\]
	
	Choose \(\eta\in C_c^\infty(\bR_v^d)\) such that
	$\eta\equiv1$ on $B_1$ and 
	$\supp\,\eta\subset B_2$,
	and define
	\[
	f_m(t,x,v):=f(t,x,v)\eta_m(v),
	\qquad
	\eta_m(v):=\eta(v/m).
	\]
	Then
	\begin{align}\label{260418802}
		\Big(\p_t+v\cdot \nabla_x -\Delta_v\Big)f_m
		=g_m,
		\qquad
		f_m(0,\cdot)\equiv0,\qquad f_m|_{\gamma_T^-}=0 ,
	\end{align}
	where
	\[
	g_m:=g\eta_m-2\nabla_v\cdot(f\nabla_v\eta_m)
	+f\Delta_v \eta_m.
	\]
	By \eqref{260418620},
	\begin{align*}
		\begin{aligned}
			&\left\|
			2\nabla_v\cdot(f\nabla_v\eta_m)-f\Delta_v \eta_m
			\right\|_{L^2(\cD_T;H^{-1}_v)}
			\\
			\lesssim\,&
			\|f\nabla_v \eta_m\|_{L^2(\cD_T\times \bR^d_v)}
			+\|f\Delta_v \eta_m\|_{L^2(\cD_T\times \bR^d_v)}
			=O\left(\frac1m\right).
		\end{aligned}
	\end{align*}
	Thus \(g_m\in L^2(\cD_T;H^{-1}_v)\) and
	\(g_m\to g\) in \(L^2(\cD_T;H^{-1}_v)\).
	
	Since \(\partial\cD\) is compact and of class \(C^{1,1}\), we may take
	\(x_1,\ldots,x_K\in \partial\cD\), \(r_0>0\), and
	\(\varepsilon_0\in(0,r_0/2]\) such that
	\[
	\partial\cD\subset \bigcup_{k=1}^K B_{r_0/4}(x_k),
	\]
	and, after decreasing \(r_0\) and \(\varepsilon_0\) if necessary,
	\[
	(\overline{\cD}\cap B_{r_0}(x_k))-\varepsilon n(x_k)
	\subset \cD
	\qquad
	\text{for }0<\varepsilon<\varepsilon_0,\quad k=1,\ldots,K.
	\]
	Choose \(\zeta_0,\zeta_1,\ldots,\zeta_K\in C_c^\infty(\bR^d_x)\) such that
	\[
	\supp\,\zeta_0\subset \cD,
	\qquad
	\supp\,\zeta_k\subset B_{r_0/2}(x_k)\quad\forall\,k=1,\ldots,K,
	\qquad
	\sum_{k=0}^K\zeta_k=1\quad\text{on }\overline{\cD}.
	\]
	Define
	\[
	f_{m,k}(t,x,v):=f_m(t,x,v)\zeta_k(x).
	\]
	
	By \eqref{260418802}, for every \(0\le k\le K\),
	\begin{align}\label{2608061631}
		\Big(\p_t+v\cdot \nabla_x-\Delta_v\Big)f_{m,k}
		=Q_{m,k},
		\qquad
		Q_{m,k}:=g_m\zeta_k+f_m(v\cdot \nabla_x \zeta_k),
	\end{align}
	pointwise on \((0,T]\times\cD\times\bR^d_v\). Since
	\(\sum_{k=0}^K\zeta_k\equiv1\) on the open set \(\cD\), we have
	\(\sum_{k=0}^K\nabla_x\zeta_k\equiv0\) there, and hence
	\begin{align}\label{2608061632}
		\sum_{k=0}^K Q_{m,k}=g_m .
	\end{align}
	
	For \(k=0\), we have
	\[
	f_{m,0}\in C_c^\infty([0,T]\times \cD\times \bR^d_v)
	\subset \cH_{\kin}(\cD_T\times\bR_v^d),
	\qquad
	f_{m,0}(0,\cdot)\equiv0,\qquad
	(f_{m,0})_{\gamma^-}\equiv0.
	\]
	
	For \(1\le k\le K\), we regard \(f_{m,k}\) and \(Q_{m,k}\) as defined on all
	of \([0,T]\times\bR^d_x\times\bR^d_v\) by extending them by zero outside
	\([0,T]\times\overline\cD\times\bR^d_v\); this is consistent with
	\(\supp\,\zeta_k\subset B_{r_0/2}(x_k)\). Define
	\[
	T_{k,\varepsilon}F(t,x,v):=
	F\big(t,x-\varepsilon n(x_k),v\big),
	\qquad 0<\varepsilon<\varepsilon_0.
	\]
	Let \(x\in\overline\cD\). If \(x\in B_{r_0}(x_k)\), then
	\(x-\varepsilon n(x_k)\in\cD\) by the choice of \(r_0\) and
	\(\varepsilon_0\); if \(x\notin B_{r_0}(x_k)\), then
	\(|x-\varepsilon n(x_k)-x_k|\ge r_0-\varepsilon_0\ge r_0/2\), so that
	\(T_{k,\varepsilon}f_{m,k}\) and \(T_{k,\varepsilon}Q_{m,k}\) vanish on a
	neighborhood of \(x\). Since \(f\in C^\infty([0,T]\times\cD\times\bR^d_v)\),
	it follows that
	\[
	T_{k,\varepsilon}f_{m,k}
	\in C_c^\infty([0,T]\times\overline{\cD}\times \bR^d_v),
	\]
	and, by \eqref{2608061631},
	\[
	\Big(\p_t+v\cdot \nabla_x-\Delta_v\Big)T_{k,\varepsilon}f_{m,k}
	=
	T_{k,\varepsilon}Q_{m,k}
	\qquad\text{on }[0,T]\times\overline\cD\times\bR^d_v .
	\]
	
	We now show that \(f_{m,k}\in \cH_{\kin}(\cD_T\times\bR_v^d)\).
	Observe that
	\begin{align*}
		&\|T_{k,\varepsilon}Q_{m,k}-Q_{m,k}\|_{L^2(\cD_T;H^{-1}_v)}^2\\
		=\,&
		\int_{\cD_T}
		\left\|
		Q_{m,k}\big(t,x-\varepsilon n(x_k),\cdot\big)
		-
		Q_{m,k}(t,x,\cdot)
		\right\|_{H^{-1}_v}^2\dd x\dd t
		\to 0
	\end{align*}
	as \(\varepsilon\to0\), by the continuity of translations in
	\(L^2(\bR^d_x;H^{-1}_v)\). Moreover, since
	\(f_{m,k}\in C_c([0,T]\times \overline{\cD}\times \bR^d_v)\) and
	\(f_{m,k}=0\) on \(\gamma_T^-\), we obtain
	\begin{align}\label{260418716}
		\lim_{\varepsilon\to0}
		\big\|(T_{k,\varepsilon}f_{m,k})|_{\gamma_T^-}\big\|_{L^2(\gamma_T^-)}
		=0.
	\end{align}
	
	For \(0<\varepsilon,\varepsilon'<\varepsilon_0\), applying estimate
	\eqref{dyHbdd} to
	\(T_{k,\varepsilon}f_{m,k}-T_{k,\varepsilon'}f_{m,k}\) gives
	\begin{align*}
		\|T_{k,\varepsilon}f_{m,k}
		-T_{k,\varepsilon'}f_{m,k}\|_{\cH_{\kin}(\cD_T\times \bR^d_v)}
		\lesssim\,&
		\|T_{k,\varepsilon}Q_{m,k}
		-T_{k,\varepsilon'}Q_{m,k}\|_{L^2(\cD_T;H_v^{-1})}
		\\
		&\quad
		+\|(T_{k,\varepsilon}f_{m,k})|_{\gamma_T^-}
		-(T_{k,\varepsilon'}f_{m,k})|_{\gamma_T^-}\|_{L^2(\gamma_T^-)}
		\to0
	\end{align*}
	as \(\varepsilon,\varepsilon'\to0\). Hence
	\((T_{k,\varepsilon}f_{m,k})_\varepsilon\) is Cauchy in
	\(\cH_{\kin}(\cD_T\times\bR_v^d)\). Since
	\[
	T_{k,\varepsilon}f_{m,k}\to f_{m,k}
	\quad\text{in }C_c([0,T]\times \overline{\cD}\times \bR^d_v)
	\subset L^2(\cD_T\times\bR_v^d),
	\]
	we conclude that
	$f_{m,k}\in\cH_{\kin}(\cD_T\times\bR_v^d)$.
	Moreover, by Lemma~\ref{lem:Hkin-operations}(i) the objects attached to
	\(T_{k,\varepsilon}f_{m,k}\) by Definition~\ref{def:H_kin} are the classical
	ones; letting \(\varepsilon\to0\) in
	\[
	\big(\p_t+v\cdot\nabla_x\big)T_{k,\varepsilon}f_{m,k}
	=T_{k,\varepsilon}Q_{m,k}+\Delta_v T_{k,\varepsilon}f_{m,k},
	\qquad
	\big(T_{k,\varepsilon}f_{m,k}\big)(0,\cdot)\equiv0,
	\]
	and using \eqref{260418716} together with the convergence in
	\(\cH_\kin(\cD_T\times\bR^d_v)\), we obtain
	\[
	\Big(\p_t+v\cdot\nabla_x-\Delta_v\Big)f_{m,k}=Q_{m,k},
	\qquad
	f_{m,k}(0,\cdot)\equiv0,
	\qquad
	(f_{m,k})_{\gamma^-}\equiv0,
	\]
	the first identity holding in \(L^2(\cD_T;H^{-1}_v)\).
	
	Consequently, by \eqref{2608061632},
	\[
	f_m=\sum_{k=0}^K f_{m,k}\in \cH_{\kin}(\cD_T\times\bR_v^d),
	\]
	and \eqref{260418802} holds in \(\cH_\kin(\cD_T\times\bR^d_v)\), that is,
	with the transport derivative and the traces of Definition~\ref{def:H_kin}.
	
	Since \(g_m\to g\) in \(L^2(\cD_T;H_v^{-1})\), applying estimate
	\eqref{dyHbdd} to \(f_m-f_{m'}\), we obtain that \((f_m)_m\) is Cauchy in
	\(\cH_\kin(\cD_T\times \bR^d_v)\). Since \(f_m\to f\) in
	\(L^2(\cD_T\times\bR_v^d)\), it follows that
	\(f\in \cH_\kin(\cD_T\times \bR^d_v)\). Letting \(m\to\infty\) in
	\eqref{260418802}, we conclude that
	\[
	\Big(\p_t+v\cdot\nabla_x-\Delta_v\Big)f=g,
	\qquad
	f(0,\cdot)=0,
	\qquad
	f_{\gamma^-}=0.
	\]
	This proves existence for \eqref{FP_dy_classic} when
	\(g\in C_c^\infty(\cD_T\times \bR^d_v)\), \(f_0\equiv0\), and \(h\equiv0\).
	
	\noindent\emph{Step 2.2.}
	We next treat \eqref{FP_dy_classic} for
	\(g\in L^2(\cD_T;H^{-1}_v)\), \(f_0\equiv 0\), and \(h\equiv 0\).
	By Lemma~\ref{260418814}.(1), choose a sequence
	\(g_n\in C_c^\infty(\cD_T\times \bR^d_v)\) such that
	$g_n\to g\quad \text{in }L^2(\cD_T;H^{-1}_v)$.
	Then Step 2.1 yields a solution
	\(f_n\in \cH_\kin(\cD_T\times \bR^d_v)\) to \eqref{FP_dy_classic} with data
	\((g_n,0,0)\).
	Applying estimate \eqref{dyHbdd} to \(f_n-f_{n'}\), we see that \((f_n)_n\) is
	Cauchy in \(\cH_\kin(\cD_T\times \bR^d_v)\); let \(f\) be its limit. Since the
	maps 
	$$
	f\mapsto(\p_t+v\cdot\nabla_x)f,\qquad f\mapsto f(0,\cdot),\qquad f\mapsto f_{\gamma^-}
	$$
	are bounded on \(\cH_\kin(\cD_T\times \bR^d_v)\), and
	since \(f_n\to f\) in \(L^2(\cD_T;H^1_v)\), passing to the limit gives a
	solution \(f\) to \eqref{FP_dy_classic} with \(f_0\equiv 0\) and
	\(h\equiv 0\).
	
	\noindent\emph{Step 2.3. General data.}
	Finally, we treat \eqref{FP_dy_classic} for
	\(g\in L^2(\cD_T;H^{-1}_v)\),
	\(f_0\in L^2(\cD\times \bR^d_v)\), and \(h\in L^2(\gamma^-_T)\).
	By the trace-pair density in Lemma~\ref{260418814}.(1), choose
	\(H_n\in C_c^\infty([0,T]\times\overline{\cD}\times \bR^d_v)\) such that
	\[
	\big(H_n(0,\cdot),H_n|_{\gamma_T^-}\big)\to (f_0,h)
	\qquad\text{in }
	L^2(\cD\times \bR^d_v)\times L^2(\gamma_T^-).
	\]
	For each \(n\), define
	\[
	\widetilde g_n
	:=g-\Big(\p_t+v\cdot \nabla_x- \Delta_v\Big)H_n .
	\]
	Then \(\widetilde g_n\in L^2(\cD_T;H_v^{-1})\), and by Step 2.2 there exists
	\(\widetilde f_n\in \cH_{\kin}(\cD_T\times\bR_v^d)\) such that
	\[
	\Big(\p_t+v\cdot \nabla_x -\Delta_v\Big)\widetilde{f}_n
	=\widetilde g_n,
	\qquad
	\widetilde f_n(0,\cdot)\equiv0,\qquad
	(\widetilde f_n)_{\gamma^-}\equiv 0 .
	\]
	Set
	\[
	f_n:=\widetilde f_n+H_n.
	\]
	Then \(f_n\in \cH_{\kin}(\cD_T\times\bR_v^d)\) and
	\[
	\Big(\p_t+v\cdot \nabla_x -\Delta_v\Big)f_n
	=g,
	\qquad
	f_n(0,\cdot)=H_n(0,\cdot),\qquad
	(f_n)_{\gamma^-}=H_n|_{\gamma_T^-}.
	\]
	The difference \(f_n-f_{n'}\) solves \eqref{FP_dy_classic} with data
	\(\big(0,\,H_n(0,\cdot)-H_{n'}(0,\cdot),\,(H_n-H_{n'})|_{\gamma_T^-}\big)\), so
	\eqref{dyHbdd} shows that \((f_n)_n\) is Cauchy in
	\(\cH_\kin(\cD_T\times \bR^d_v)\). Passing to the limit as in Step 2.2, we
	obtain a solution \(f\) to \eqref{FP_dy_classic} with data \((g,f_0,h)\).
	
	This completes the proof of (1).
	
	(2) \noindent\emph{Step 1. Estimate \eqref{stHbdd} and uniqueness.}
	
	We first prove estimate \eqref{stHbdd}. Let
	\(f\in \cH_{\kin}(\cD\times\bR_v^d)\) be a solution of
	\eqref{FP_st_classic}. Then, by the integration-by-parts formula
	\eqref{ibp_st}, we have
	\begin{align*}
		&2\int_{\cD}
		\la
		(v\cdot \nabla_x-\Delta_v+\lambda) f,f
		\ra_v\dd x
		\\
		=\,&
		-2\int_{\cD}\la \Delta_v f,f\ra_v\dd x
		+2\lambda \|f\|_{L^2(\cD\times \bR^d_v)}^2
		+\int_{\gamma}
		|f_\gamma|^2\sgn(n_x\cdot v)\dd \gamma .
	\end{align*}
	Since \(f_\gamma=h\) on \(\gamma^-\), it follows that
	\begin{align*}
		\begin{aligned}
			&\|\nabla_v f\|_{L^2(\cD\times \bR^d_v)}^2
			+\lambda \|f\|_{L^2(\cD\times \bR^d_v)}^2
			+\|f_{\gamma}\|_{L^2(\gamma^+)}^2
			\\
			\le\,&C \left(
			\|g\|_{L^2(\cD;H^{-1}_v)}
			\|f\|_{L^2(\cD;H^1_v)}
			+\|h\|_{L^2(\gamma^-)}^2 \right),
		\end{aligned}
	\end{align*}
	where \(C\) is a universal constant.
	Applying Young's inequality, we obtain
	\[
	\|f\|_{L^2(\cD;H^1_v)}
	\le
	C(\lambda \wedge 1)^{-1}\left(
	\|g\|_{L^2(\cD;H_v^{-1})}
	+\|h\|_{L^2(\gamma^-)}
	\right),
	\]
	where \(C\) is a universal constant. Using the equation, we also have
	\begin{align*}
		\begin{aligned}
			\left\|v\cdot \nabla_xf\right\|_{L^2(\cD;H^{-1}_v)}
			&\le
			\left\|\Delta_v f\right\|_{L^2(\cD;H^{-1}_v)}
			+\lambda \|f\|_{L^2(\cD;H^{-1}_v)}
			+\left\|g\right\|_{L^2(\cD;H^{-1}_v)}
			\\
			&\lesssim
			(1+\lambda)\|f\|_{L^2(\cD;H^1_v)}
			+\|g\|_{L^2(\cD;H^{-1}_v)} .
		\end{aligned}
	\end{align*}
	Since \(f_{\gamma^-}=h\), the above estimates yield
	\eqref{stHbdd}.
	
	In particular, if \(g\) and \(h\) vanish identically, then \(f=0\). Hence
	uniqueness follows immediately.
	
	\noindent\emph{Step 2. Existence.}
	The existence result in (2) is proved in the same way as in the proof of (1),
	using Corollary~\ref{prop_sol}.(2) instead of Corollary~\ref{prop_sol}.(1).
\end{proof}

\begin{proof}[\textbf{Proof of Theorem~\ref{solv}}]
	By Proposition~\ref{solv_Del} and the method of continuity, it suffices to prove the a priori estimates \eqref{dyHbdd} and \eqref{stHbdd}.
	
	\smallskip
	\noindent
	\textbf{(1)}
	Let \(f\in \cH_{\kin}(\cD_T\times \bR^d_v)\), and put
	\begin{align}\label{260618457}
		g:=(\p_t+v\cdot \nabla_x+L)f,\qquad
		f_0:=f(0,\cdot),\qquad h:=f_{\gamma^-}.
	\end{align}
	Then \(f\) satisfies
	\begin{align*}
		\begin{cases}
			\big(\partial_t +v\cdot\nabla_x\big) f -D_i(a^{ij} D_jf)
			=g-\big(b^iD_if+D_i(\bar{b}^i f)+cf\big)
			\quad &\text{in}\quad \cD_T\times \bR^d_v\,;\\
			f(0,\,\cdot\,)= f_0
			\quad &\text{in}\quad \cD\times \bR^d_v\,;\\
			f=h
			\quad &\text{on}\quad \gamma^-_T\,.
		\end{cases}
	\end{align*}
	For every \(t\in[0,T]\), observe that
	\begin{align}\label{260618504}
		\iint_{(0,t)\times \cD}
		\la  -D_i(a^{ij} D_jf),f\ra_v \dd x \dd s
		=
		\iiint_{\cD_t\times \bR^d_v}
		a^{ij}D_ifD_jf\dd v\dd x\dd s
		\ge
		\nu \|\nabla_v f\|_{L^2(\cD_t\times \bR^d_v)}^2,
	\end{align}
	and
	\begin{align}\label{260618505}
		\iint_{(0,t)\times \cD}
		\left|\la  b^iD_if+D_i(\bar{b}^i f)+cf,f\ra_v \right|
		\dd x \dd s
		\le
		C(d)K
		\|f\|_{L^2(\cD_t\times \bR^d_v)}
		\|f\|_{L^2(\cD_t;H^1_v)}.
	\end{align}
	Thus, arguing as in the derivation of \eqref{260617551} and using Young's
	inequality, we obtain
	\begin{align}\label{260618439}
		\begin{aligned}
			&\|\nabla_v f\|_{L^2(\cD_t\times \bR^d_v)}^2
			+\sup_{0\le s\le t}\|f(s,\cdot)\|_{L^2(\cD\times \bR^d_v)}^2
			\\
			\le\,&
			C(d,\nu,K)\left(
			\|g\|_{L^2(\cD_t;H^{-1}_v)}^2
			+\|f_0\|_{L^2(\cD\times \bR^d_v)}^2
			+\|h\|_{L^2(\gamma_t^-)}^2
			+\int_0^t\|f(s,\cdot)\|_{L^2(\cD\times \bR^d_v)}^2\dd s
			\right).
		\end{aligned}
	\end{align}
	By Gronwall's inequality, we have
	\begin{align}\label{260418912-fixed2}
		\sup_{0\le t\le T}\|f(t,\cdot)\|_{L^2(\cD\times \bR^d_v)}^2
		\le
		C
		\left(
		\|g\|_{L^2(\cD_T;H^{-1}_v)}^2
		+\|f_0\|_{L^2(\cD\times \bR^d_v)}^2
		+\|h\|_{L^2(\gamma_T^-)}^2
		\right),
	\end{align}
	where \(C=C(d,\nu,K,T)\). Taking \(t=T\) in
	\eqref{260618439} and using \eqref{260418912-fixed2}, we obtain
	\begin{align*}
		\|f\|_{L^2(\cD_T;H^1_v)}^2
		\le C\left(
		\|g\|_{L^2(\cD_T;H^{-1}_v)}^2
		+\|f_0\|_{L^2(\cD\times \bR^d_v)}^2
		+\|h\|_{L^2(\gamma_T^-)}^2
		\right).
	\end{align*}
	Moreover, by \eqref{260618457},
	\[
	\|(\partial_t+v\cdot\nabla_x)f\|_{L^2(\cD_T;H^{-1}_v)}
	\le
	\|g\|_{L^2(\cD_T;H^{-1}_v)}
	+C(d,\nu,K)\|f\|_{L^2(\cD_T;H^1_v)}.
	\]
	Together with \(f(0,\cdot)=f_0\) and \(f_{\gamma^-}=h\), the preceding
	estimates yield \eqref{dyHbdd}.
	
	\medskip
	\noindent
	\textbf{(2)}
	Let \(f\in \cH_{\kin}(\cD\times \bR^d_v)\), and put
	\begin{align}\label{260618457_st}
		g:=(v\cdot \nabla_x+L+\lambda)f,\qquad h:=f_{\gamma^-}.
	\end{align}
	Then \(f\) satisfies
	\begin{align*}
		\begin{cases}
			v\cdot\nabla_x f -D_i(a^{ij} D_jf)+\lambda f
			=g-\big(b^iD_if+D_i(\bar{b}^i f)+cf\big)
			\quad &\text{in}\quad \cD\times \bR^d_v\,;\\
			f=h
			\quad &\text{on}\quad \gamma^-\,.
		\end{cases}
	\end{align*}
	Using the stationary analogues of \eqref{260618504} and \eqref{260618505},
	and arguing as in Step 1 of Proposition~\ref{solv_Del}.(2), we obtain
	\begin{align}\label{260618439_st}
		\begin{aligned}
			\|\nabla_v f\|_{L^2(\cD\times \bR^d_v)}
			+\sqrt{\lambda}\|f\|_{L^2(\cD\times \bR^d_v)}
			\le
			C_0\left(
			\|g\|_{L^2(\cD;H^{-1}_v)}
			+\|h\|_{L^2(\gamma^-)}
			+\|f\|_{L^2(\cD\times \bR^d_v)}
			\right),
		\end{aligned}
	\end{align}
	where \(C_0=C(d,\nu,K)\). In particular, \(C_0\) is independent of
	\(\lambda\).
	
	Set \(\lambda_0:=(C_0+1)^2\). Then, for every \(\lambda\ge\lambda_0\),
	\eqref{260618439_st} yields
	\begin{align*}
		\|\nabla_v f\|_{L^2(\cD\times \bR^d_v)}
		+\|f\|_{L^2(\cD\times \bR^d_v)}
		&\le
		\|\nabla_v f\|_{L^2(\cD\times \bR^d_v)}
		+(\sqrt{\lambda}-C_0)\|f\|_{L^2(\cD\times \bR^d_v)}
		\\
		&\le
		C_0\left(
		\|g\|_{L^2(\cD;H^{-1}_v)}
		+\|h\|_{L^2(\gamma^-)}
		\right).
	\end{align*}
	Moreover, by \eqref{260618457_st},
	$v\cdot\nabla_xf=g-Lf-\lambda f$.
	Hence
	\begin{align*}
		\|v\cdot\nabla_xf\|_{L^2(\cD;H^{-1}_v)}
		&\le
		\|g\|_{L^2(\cD;H^{-1}_v)}
		+C(d,\nu,K)\|f\|_{L^2(\cD;H^1_v)}
		+\lambda\|f\|_{L^2(\cD\times\bR_v^d)}
		\\
		&\le
		C(d,\nu,K,\lambda)
		\left(
		\|g\|_{L^2(\cD;H^{-1}_v)}
		+\|h\|_{L^2(\gamma^-)}
		\right).
	\end{align*}
	Together with \(f_{\gamma^-}=h\), the preceding estimates yield
	\eqref{stHbdd}.
\end{proof}

\subsection{A weak characterization of $\cH_{\kin}$}\label{sec:weakform}

We identify the completion-based \(\cH_{\kin}\)-spaces with the weak
transport spaces characterized by velocity Sobolev regularity, the
appropriate derivatives and traces, and Green's identities.

\begin{prop}[Weak characterization of $\cH_{\kin}$]
	\label{H_kin=weak}\,
	
	\begin{enumerate}
		\item The space \(\cH_{\kin}(\cD_T\times\bR^d_v)\) coincides with
		the set of all \(f\in L^2(\cD_T;H^1_v)\) for which there exist
		\[
		(\partial_t+v\cdot\nabla_x)f\in L^2(\cD_T;H^{-1}_v),
		\qquad
		f(0,\cdot),\,f(T,\cdot)\in L^2(\cD\times \bR^d_v),
		\qquad
		f_\gamma\in L^2(\gamma_T),
		\]
		such that, for every
		\(\phi\in C_c^\infty([0,T]\times \overline{\cD}\times \bR^d_v)\),
		\begin{align}\label{260618611}
			\begin{aligned}
				&\iint_{\cD_T}
				\Big(
				\big\la (\partial_t+v\cdot\nabla_x)f,\phi\big\ra_{v}
				+\big\la f,(\partial_t+v\cdot \nabla_x) \phi\big\ra_{v}
				\Big)\dd x\dd t\\
				=\,&
				\iint_{\cD\times \bR^d_v}f(T,\cdot)\,\phi(T,\cdot)\dd x\dd v
				-\iint_{\cD\times \bR^d_v}f(0,\cdot)\,\phi(0,\cdot)\dd x\dd v
				+\iint_{\gamma_T}
				f_{\gamma}\,\phi|_{\gamma}\,\operatorname{sgn}(n_x\cdot v)\dd \gamma\dd t .
			\end{aligned}
		\end{align}
		
		\item The space \(\cH_{\kin}(\cD\times\bR^d_v)\) coincides with
		the set of all \(f\in L^2(\cD;H^1_v)\) for which there exist
		\[
		v\cdot\nabla_x f\in L^2(\cD;H^{-1}_v),
		\qquad
		f_\gamma\in L^2(\gamma),
		\]
		such that, for every
		\(\phi\in C_c^\infty(\overline{\cD}\times \bR^d_v)\),
		\begin{align*}
			\int_{\cD}
			\Big(
			\la v\cdot\nabla_x f,\phi\ra_{v}
			+\la f,v\cdot \nabla_x \phi\ra_{v}
			\Big)\dd x
			=
			\int_{\gamma}
			f_{\gamma}\,\phi|_{\gamma}\,\operatorname{sgn}(n_x\cdot v)\dd \gamma .
		\end{align*}
	\end{enumerate}
\end{prop}

\begin{proof}
	We prove \((1)\) in detail. The proof of \((2)\) is the same, with
	\(\p_t+v\cdot\nabla_x\) replaced by \(v\cdot \nabla_x\), the two time-trace
	terms omitted, and \(f^\ast(x,v):=f(x,-v)\). The only point that requires
	attention is that the solvability invoked in Step 1 must be applied with
	\(\lambda=1\). Such solvability is supplied by Proposition~\ref{solv_Del}.(2)
	rather than by Theorem~\ref{solv}.(2), the latter being available only for
	\(\lambda\ge\lambda_0\).
	
	For notational simplicity, we write
	\[
	\cH_\kin:=\cH_\kin(\cD_T\times \bR^d_v),
	\qquad
	\mathfrak T:=\p_t+v\cdot \nabla_x,
	\]
	and let \(X\) denote the class described in the statement of the proposition.
	Let \(F\in X\). Applying the argument in the proof of
	Proposition~\ref{260128509} to \(F\) and to \(F(T-t,x,-v)\), we see that
	\[
	\mathfrak T F\in L^2(\cD_T;H^{-1}_v),\qquad
	F(0,\cdot),F(T,\cdot)\in L^2(\cD\times \bR^d_v),\qquad
	F_\gamma\in L^2(\gamma_T)
	\]
	are uniquely determined.
	
	It follows that \(X\) is an inner product space under the inner
	product
	\begin{align*}
		\la F,G\ra_X:=\,&
		\la F,G\ra_{L^2(\cD_T;H^1_v)}
		+\la \mathfrak TF,\mathfrak T G\ra_{L^2(\cD_T;H^{-1}_v)}
		\\
		&+\la F(0,\cdot),G(0,\cdot)\ra_{L^2(\cD\times \bR^d_v)}
		+\la F(T,\cdot),G(T,\cdot)\ra_{L^2(\cD\times \bR^d_v)}
		+\la F_\gamma,G_\gamma\ra_{L^2(\gamma_T)}.
	\end{align*}
	The space $X$ is complete.  Let $(F_j)$ be a Cauchy sequence in
	$X$.
	Its bulk functions, transport derivatives, initial and
	terminal traces, and full boundary traces converge respectively to
	\[
	F\in L^2(\cD_T;H_v^1),\quad
	G\in L^2(\cD_T;H_v^{-1}),\quad
	F_0,F_T\in L^2(\cD\times\bR_v^d),\quad
	H\in L^2(\gamma_T).
	\]
	For every compactly supported smooth test function $\phi$,
	passage to the limit in \eqref{260618611} is legitimate by
	$H_v^{-1}$--$H_v^1$ duality, Cauchy--Schwarz for the two time
	traces, and
	\[
	\left|\int_{\gamma_T}(F_j)_\gamma\phi\,
	\operatorname{sgn}(n_x\cdot v)\dd\gamma\dd t
	-\int_{\gamma_T}H\phi\,
	\operatorname{sgn}(n_x\cdot v)\dd\gamma\dd t\right|
	\le\|(F_j)_\gamma-H\|_{L^2(\gamma_T)}
	\|\phi\|_{L^2(\gamma_T)}.
	\]
	Testing first in the interior identifies
	$G=(\partial_t+v\cdot\nabla_x)F$ in distributions.  The
	limiting Green identity then identifies $F_0,F_T,H$ as the
	traces of $F$ in the sense of the definition of $X$.
	Consequently $F\in X$ and $F_j\to F$ in $X$.  In the stationary case a Cauchy sequence has limits
	$F\in L^2(\cD;H_v^1)$, $G\in L^2(\cD;H_v^{-1})$ and
	$H\in L^2(\gamma)$.  Passing to the limit gives
	\[
	\int_\cD\bigl(\langle G,\phi\rangle_v
	+\langle F,v\cdot\nabla_x\phi\rangle_v\bigr)\dd x
	=\int_\gamma H\phi\,\operatorname{sgn}(n_x\cdot v)\dd\gamma.
	\]
	Interior tests identify $G=v\cdot\nabla_xF$, so the limit
	lies in the stationary intrinsic class and the sequence converges
	in its norm.
	
	Moreover, \(\cH_\kin\subset X\) by Proposition~\ref{trace_theorem} and
	\eqref{ibp_dy}, and
	\[
	\|f\|_{\cH_\kin}
	\simeq
	\|f\|_X
	\quad \text{for all }f\in \cH_\kin;
	\]
	here \(\|f\|_{\cH_\kin}\le \|f\|_X\) follows from
	Remark~\ref{rem:trace-compatibility}, and the reverse inequality from
	Proposition~\ref{trace_theorem}. Thus \(\cH_\kin\) is a closed subspace
	of the Hilbert space $X$, and orthogonal projection onto
	$\cH_\kin$ is well defined.
	
	Suppose, to the contrary, that
	\[
	\cH_\kin\subsetneq X,
	\]
	and pick \(G\in X\setminus \cH_\kin\). Writing \(P\) for the orthogonal
	projection onto \(\cH_\kin\) and setting \(F_0:=G-PG\), we obtain
	\begin{align}\label{260618549}
		F_0\neq 0,\qquad
		\la F_0,f\ra_X=0
		\quad \text{for all }f\in \cH_\kin.
	\end{align}
	
	\smallskip
	\noindent
	\emph{Step 1.}
	We use the relations
	\[
	\la u,v\ra_{H^1_v}
	=
	\big\la u,(1-\Delta_v)v\big\ra_v,
	\qquad
	\la u,v\ra_{H^{-1}_v}
	=
	\big\la (1-\Delta_v)^{-1}u,v\big\ra_v.
	\]
	For any \(f\in C_c^\infty([0,T]\times \overline{\cD}\times \bR^d_v)\), set
	\[
	g:=(\mathfrak T-\Delta_v+1)f .
	\]
	Then \(\mathfrak T f=g-(1-\Delta_v)f\). By \eqref{260618611}, applied to
	\(F_0\), we have
	\begin{align*}
		\la \mathfrak T F_0,\mathfrak T f\ra_{L^2(\cD_T;H^{-1}_v)}
		&=
		-\la \mathfrak T F_0,(1-\Delta_v)f\ra_{L^2(\cD_T;H^{-1}_v)}
		+\la \mathfrak T F_0,g \ra_{L^2(\cD_T;H^{-1}_v)}
		\\
		&=
		-\iint_{\cD_T}\la \mathfrak T F_0,f\ra_v\dd x\dd t
		+\iint_{\cD_T}\la (1-\Delta_v)^{-1}\mathfrak T F_0,g \ra_v\dd x\dd t
		\\
		&=
		\iint_{\cD_T}\la F_0,\mathfrak Tf\ra_v\dd x\dd t
		-\la F_0(T,\cdot),f(T,\cdot)\ra_{L^2(\cD\times \bR^d_v)}
		\\
		&\qquad
		+\la F_0(0,\cdot),f(0,\cdot)\ra_{L^2(\cD\times \bR^d_v)}
		-\la (F_0)_\gamma,f_\gamma\ra_{L^2(\gamma^+_T)}
		\\
		&\qquad
		+\la (F_0)_\gamma,f_\gamma\ra_{L^2(\gamma_T^-)}
		+\iint_{\cD_T}\la (1-\Delta_v)^{-1}\mathfrak T F_0,g \ra_v\dd x\dd t .
	\end{align*}
	Therefore, \eqref{260618549} gives
	\begin{align}\label{260618621}
		\begin{aligned}
			0=\la F_0,f\ra_{X}
			=\,&
			\iint_{\cD_T}
			\la F_0+(1-\Delta_v)^{-1}\mathfrak T F_0,g\ra_v\dd x\dd t
			\\
			&\quad
			+2\la F_0(0,\cdot),f(0,\cdot)\ra_{L^2(\cD\times \bR^d_v)}
			+2\la (F_0)_\gamma,f_\gamma\ra_{L^2(\gamma_T^-)}\\
			=\,&
			\la (1-\Delta_v)F_0+\mathfrak T F_0,g\ra_{L^2(\cD_T;H^{-1}_v)}\\
			&+2\la F_0(0,\cdot),f(0,\cdot)\ra_{L^2(\cD\times \bR^d_v)}
			+2\la (F_0)_\gamma,f_\gamma\ra_{L^2(\gamma_T^-)} .
		\end{aligned}
	\end{align}
	Since \(C_c^\infty([0,T]\times \overline{\cD}\times \bR^d_v)\) is dense in
	\(\cH_\kin\), \eqref{260618621} holds for every \(f\in\cH_\kin\).
	
	The operator \(L=-\Delta_v+1\) satisfies \eqref{coeff_ellip} and
	\eqref{coeff_bdd} with \(\nu=K=1\). Hence, by Theorem~\ref{solv}.(1), there
	exists \(f\in\cH_\kin\) such that
	\[
	\Big((\mathfrak T-\Delta_v+1)f,\ f(0,\cdot),\ f_{\gamma^-}\Big)
	=
	\Big(
	(1-\Delta_v)F_0+\mathfrak T F_0,\,
	F_0(0,\cdot),\,
	(F_0)_\gamma|_{\gamma_T^-}
	\Big)
	\]
	in
	\[
	L^2(\cD_T;H^{-1}_v)
	\times L^2(\cD\times \bR^d_v)
	\times L^2(\gamma_T^-).
	\]
	Substituting this \(f\) into \eqref{260618621}, we obtain
	\[
	\left\|(1-\Delta_v)F_0+\mathfrak T F_0\right\|_{L^2(\cD_T;H^{-1}_v)}^2
	+2\|F_0(0,\cdot)\|_{L^2(\cD\times\bR^d_v)}^2
	+2\|(F_0)_\gamma\|_{L^2(\gamma_T^-)}^2
	=0.
	\]
	Hence
	\begin{align}\label{260618706}
		(1-\Delta_v)F_0+\mathfrak T F_0=0
		\quad \text{in }L^2(\cD_T;H^{-1}_v),
	\end{align}
	and
	\begin{align}\label{260618707}
		\Big(F_0(0,\cdot),(F_0)_\gamma|_{\gamma_T^-}\Big)=(0,0)
		\quad \text{in } L^2(\cD\times \bR^d_v)\times L^2(\gamma_T^-).
	\end{align}
	
	\smallskip
	\noindent
	\emph{Step 2.}
	For \(\phi(t,x,v)\), denote
	\[
	\phi^\ast(t,x,v):=\phi(T-t,x,-v).
	\]
	By Lemma~\ref{lem:Hkin-operations}(iii), if \(f\in \cH_\kin\), then
	\(f^\ast\in \cH_\kin\), with
	\[
	\mathfrak T f^\ast =-(\mathfrak T f)^\ast,
	\qquad
	\Delta_v f^\ast=(\Delta_v f)^\ast,
	\qquad
	(f^\ast)_\gamma=(f_\gamma)^\ast ,
	\]
	the last identity exchanging \(\gamma^\pm_T\) with \(\gamma^\mp_T\). The class
	\(X\) is invariant under \(F\mapsto F^\ast\) as well: if \(F\in X\), then
	\eqref{260618611} for \(F^\ast\) and a test function \(\phi\) is
	\eqref{260618611} for \(F\) and the test function \(\phi^\ast\), both sides
	changing sign. Since \(\dd t\dd x\dd v\) and \(\dd\gamma\dd t\) are invariant
	under \((t,v)\mapsto (T-t,-v)\) and \(v\mapsto-v\) is an isometry of
	\(H^{\pm1}_v\), we also have
	\(\la F^\ast,G^\ast\ra_X=\la F,G\ra_X\) for all \(F,G\in X\). Since
	\((f^\ast)^\ast=f\), it follows from \eqref{260618549} that
	\[
	\la F_0^\ast,f\ra_X
	=
	\la F_0,f^\ast\ra_X
	=0
	\qquad \text{for all }f\in\cH_\kin .
	\]
	Applying Step 1 to \(F_0^\ast\), and using that \((1-\Delta_v)^{-1}\)
	commutes with the transformation \(u\mapsto u^\ast\), we obtain
	\begin{align}\label{26061870611}
		(1-\Delta_v)F_0-\mathfrak T F_0
		=
		\big((1-\Delta_v)F_0^\ast+\mathfrak T F_0^\ast\big)^\ast
		=0
		\quad \text{in }L^2(\cD_T;H^{-1}_v),
	\end{align}
	and
	\begin{align}\label{26061870711}
		\big(F_0^\ast(0,\cdot),(F_0^\ast)_{\gamma}|_{\gamma_T^-}\big)=(0,0)
		\quad \text{in }
		L^2(\cD\times \bR^d_v)\times L^2(\gamma_T^-).
	\end{align}
	Since \(F_0^\ast(0,x,v)=F_0(T,x,-v)\) and
	\((F_0^\ast)_\gamma=(F_0)_\gamma^\ast\) carries \(\gamma_T^-\) onto
	\(\gamma_T^+\), and since \(v\mapsto-v\) preserves the relevant measures,
	\eqref{26061870711} is equivalent to
	\begin{align}\label{2608061711}
		\big(F_0(T,\cdot),(F_0)_{\gamma}|_{\gamma_T^+}\big)=(0,0)
		\quad \text{in }
		L^2(\cD\times \bR^d_v)\times L^2(\gamma_T^+).
	\end{align}
	Adding \eqref{260618706} and \eqref{26061870611} gives
	\((1-\Delta_v)F_0=0\) in \(L^2(\cD_T;H^{-1}_v)\). Since
	\(1-\Delta_v:H^1_v\to H^{-1}_v\) is an isomorphism, we get
	\begin{align*}
		(F_0,\mathfrak T F_0)=(0,0)
		\quad \text{in}\quad
		L^2(\cD_T;H^1_v)\times L^2(\cD_T;H^{-1}_v).
	\end{align*}
	Together with \eqref{260618707} and \eqref{2608061711}, this implies
	\(\|F_0\|_X=0\), contradicting \eqref{260618549}.
\end{proof}

\section{Boundary Regularity for Kinetic Fokker–Planck Equations}\label{2606261120}

This section develops local boundary regularity in stages. We
first set up kinetic cylinders and weak solutions and prove an
\(L^2\)--\(L^\infty\) estimate; we then pass to half-space geometry, boundary
H\"older regularity, grazing decay, and finally velocity-gradient bounds under
regular coefficients.

Throughout this section we use \(\cH_{\kin}\) on arbitrary bounded time
intervals instead of $(0,T)$, as permitted by Remark~\ref{rem:Hkin-time-translation}, and we
denote by \(f(t_1,\cdot)\), \(f(t_2,\cdot)\), \(f_\gamma\) and
\(f_{\gamma^-}\) the corresponding traces, the terminal one being furnished by
Proposition~\ref{trace_theorem}.

\subsection{Kinetic cylinders and local weak solutions}

\begin{defn}\label{def:KC}
	Let \(Z_0:=(t_0,x_0,v_0)\in\bR\times\bR_x^d\times\bR_v^d\) and \(r>0\).
	
	\begin{enumerate}
		\item We define the \emph{kinetic cylinder} centered at \(Z_0\) with radius
		\(r\) by
		\begin{align}\label{KC}
			\begin{split}
				Q_r(Z_0)
				:=
				\big\{(t,x,v)\in \bR\times\bR_x^d\times\bR_v^d:\;&
				t\in(t_0-r^2,t_0],
				\\
				& |x-x_0-(t-t_0)v_0|<r^3,\quad |v-v_0|<r
				\big\}.
			\end{split}
		\end{align}
		For a set \(E\subset\bR_x^d\), we write
		\begin{align*}
			Q_r^E(Z_0)
			&:=
			Q_r(Z_0)\cap\bigl(\bR\times E\times\bR_v^d\bigr).
		\end{align*}
		
		\item We define the \emph{kinetic distance} by
		\begin{align}\label{kineticdistance}
			\mathbf d(Z,\bar Z)
			:=
			\min_{w\in\bR_v^d}
			\max\left(
			|t-\bar t|^{1/2},
			|x-\bar x-(t-\bar t)w|^{1/3},
			|v-w|,
			|\bar v-w|
			\right),
		\end{align}
		where \(Z=(t,x,v)\) and \(\bar Z=(\bar t,\bar x,\bar v)\).
	\end{enumerate}
\end{defn}

For \(Z:=(t,x,v)\) and \(\bar Z:=(\bar t,\bar x,\bar v)\) in
\(\bR^{1+2d}\), we use the group law
\[
Z\circ\bar Z
:=
(t+\bar t,\ x+\bar x+\bar t v,\ v+\bar v).
\]
The inverse of \(Z\) is then given by
\[
Z^{-1}=(-t,-x+tv,-v).
\]
With this notation,
\[
Q_r(Z)=Z\circ Q_r(0),
\qquad
Q_r(0)=(-r^2,0]\times B_{r^3}\times B_r.
\]
Moreover,
\[
\mathbf d(Z,\bar Z)
=
\mathbf d(\bar Z^{-1}\circ Z,0)
=
\mathbf d(Z^{-1}\circ \bar Z,0),
\]
and the distance is left-invariant:
\[
\mathbf d(Z,\bar Z)
=
\mathbf d(W\circ Z,W\circ \bar Z)
\qquad
\text{for all }W\in\bR^{1+2d}.
\]

We also use the homogeneous quasi-norm
\begin{align*}
	|Z|_{\kin}
	:=
	\max\bigl(|t|^{1/2},|x|^{1/3},|v|\bigr).
\end{align*}
It is comparable to \(Z\mapsto \mathbf d(Z,0)\). More precisely,
\[
\mathbf d(Z,0)\le |Z|_{\kin}\le 2\mathbf d(Z,0).
\]
We also mention that $\mathbf{d}$ is a metric on $\bR\times \bR^d_x\times \bR^d_v$ (see, \textit{e.g.}, \cite[Proposition~2.2]{imbert2021schauder}).

\begin{lemma}\label{260127226}
	Let \(r>0\), \(Z=(t,x,v)\in\bR^{1+2d}\), and
	\(\bar Z=(\bar t,\bar x,\bar v)\in\bR^{1+2d}\) with \(\bar t\le t\). Then the
	following assertions hold.
	
	\begin{enumerate}
		\item If \(\bar Z\in Q_r(Z)\), then \(\mathbf d(Z,\bar Z)<r\).
		
		\item If \(\mathbf d(Z,\bar Z)<r\), then \(\bar Z\in Q_{2r}(Z)\).
		
		\item For all \(\sigma,\tau>0\) and \(Z_1\in Q_\sigma(Z)\), one has
		\(Q_\tau(Z_1)\subset Q_{\sigma+\tau}(Z)\).
	\end{enumerate}
\end{lemma}

\begin{proof}
	The first assertion follows by taking \(w=v\) in the definition of
	\(\mathbf d\). For the second assertion, choose \(w\in\bR_v^d\) such that all
	four quantities in the definition of \(\mathbf d(Z,\bar Z)\) are less than
	\(r\). Then
	\[
	|\bar v-v|<2r,
	\qquad
	|\bar x-x-(\bar t-t)v|
	\le r^3+r^2r=2r^3<(2r)^3.
	\]
	Together with \(\bar t\le t\), this gives \(\bar Z\in Q_{2r}(Z)\).
	
	For the third assertion, write \(Z=(t,x,v)\), \(Z_1=(t_1,x_1,v_1)\), and let
	\(\bar Z=(\bar t,\bar x,\bar v)\in Q_\tau(Z_1)\).  Then
	\(t-\sigma^2-\tau^2<\bar t\le t\) and
	\(|\bar v-v|\le|\bar v-v_1|+|v_1-v|<\tau+\sigma\), while
	\begin{align*}
		|\bar x-x-(\bar t-t)v|
		&\le
		|\bar x-x_1-(\bar t-t_1)v_1|
		+|x_1-x-(t_1-t)v|
		+|\bar t-t_1|\,|v_1-v|
		\\
		&<
		\tau^3+\sigma^3+\tau^2\sigma
		\le
		(\sigma+\tau)^3 .
	\end{align*}
	Since \(\sigma^2+\tau^2\le(\sigma+\tau)^2\), this gives
	\(\bar Z\in Q_{\sigma+\tau}(Z)\).
\end{proof}

\begin{lemma}\label{lem:cylinder-inclusion}
	Let \(Z_0\in\bR^{1+2d}\), \(r,s>0\), and \(Z\in Q_r(Z_0)\).  Then
	\(Q_s(Z)\subset Q_{r+s}(Z_0)\).
\end{lemma}

\begin{proof}
	Write \(Z_0=(t_0,x_0,v_0)\), \(Z=(t_1,x_1,v_1)\), and let
	\(Y=(\tau,y,w)\in Q_s(Z)\).  Then
	\(\tau\le t_1\le t_0\) and
	\(\tau>t_1-s^2>t_0-r^2-s^2\ge t_0-(r+s)^2\), while
	\(|w-v_0|\le|w-v_1|+|v_1-v_0|<s+r\).  Moreover
	\begin{align*}
		y-x_0-(\tau-t_0)v_0
		&=\bigl(y-x_1-(\tau-t_1)v_1\bigr)
		+\bigl(x_1-x_0-(t_1-t_0)v_0\bigr)
		+(\tau-t_1)(v_1-v_0),
	\end{align*}
	so that \(|y-x_0-(\tau-t_0)v_0|<s^3+r^3+s^2r\le(r+s)^3\).
\end{proof}

Throughout, for a set \(A\subset\bR^{1+2d}\) which need not be open, we write
\(u\in L^2_{\mathrm{loc}}(A)\) if \(u\in L^2(K\cap A)\) for every compact
\(K\subset A\); the same convention is used for the class
\(C_c^\infty(A)\), whose elements are the smooth functions whose support is a
compact subset of \(A\).  This is relevant because the sets \(Q_R(Z_0)\) and
\(Q_R^{\overline\cD}(Z_0)\) contain part of their topological boundary.

\begin{defn}[Weak solution with kinetic boundary trace]
	\label{def:L2sol}
	Let \(Z_0:=(t_0,x_0,v_0)\in\bR\times\bR_x^d\times\bR_v^d\) and \(R>0\).
	Let
	\[
	a^{ij},\,b^i,\,c\in L^\infty\bigl(Q_R^{\cD}(Z_0)\bigr),
	\qquad i,j=1,\ldots,d,
	\]
	and
	\[
	g^0,\,g^1,\ldots,g^d\in L^2\bigl(Q_R^{\cD}(Z_0)\bigr),
	\qquad
	h\in L^2((\bR\times \gamma^-)\cap Q_R(Z_0),\mathrm{d}\gamma\dd t).
	\]
	We say that \(f\in L^2\bigl(Q_R^{\cD}(Z_0)\bigr)\) is a
	\emph{weak solution} of
	\begin{align}\label{lineqtn}
		\begin{cases}
			(\partial_t+v\cdot\nabla_x)f
			=
			D_i(a^{ij}D_j f)+b^iD_i f+cf+g^0+D_i g^i
			&\quad \text{in } Q_R^{\cD}(Z_0),\\
			f=h
			&\quad \text{on } (\bR\times\gamma^-)\cap Q_R(Z_0),
		\end{cases}
	\end{align}
	if the following condition holds.
	
	For every
	$\varphi\in C_c^\infty\bigl(Q_R^{\overline\cD}(Z_0)\bigr)$,
	the function \(f\varphi\), extended by zero outside \(Q_R^{\cD}(Z_0)\),
	belongs to
	$\cH_{\kin}\bigl((t_0-R^2,t_0)\times\cD\times\bR_v^d\bigr)
	$
	and satisfies
	\begin{align}\label{eq:localized-weak-solution}
		\begin{aligned}
			(\partial_t+v\cdot\nabla_x)(f\varphi)
			=\;&
			D_i\bigl(a^{ij}D_j(f\varphi)\bigr)
			+b^iD_i(f\varphi)
			+cf\varphi
			\\
			&\quad
			+g^0\varphi
			+D_i(g^i\varphi)
			+\cR_\varphi[f]
		\end{aligned}
	\end{align}
	in \(L^2((t_0-R^2,t_0)\times\cD;H_v^{-1})\), where
	\begin{align}\label{2608052010}
		\cR_\varphi[f]
		:=
		\,&f(\partial_t\varphi+v\cdot\nabla_x\varphi)
		\nonumber\\
		&-\left(D_i(a^{ij}fD_j\varphi)
		+a^{ij}D_j(fD_i\varphi)
		+fb^iD_i\varphi
		+g^iD_i\varphi
		-a^{ij}fD_{ij}\varphi\right).
	\end{align}
	Moreover,
	\begin{equation}\label{def:L2sol-boundary-conditions}
		(f\varphi)(t_0-R^2,\cdot)\equiv0
		\quad\text{in }L^2(\cD\times\bR_v^d),\qquad (f\varphi)_{\gamma^-}=h\,\varphi
		\quad\text{in }L^2\bigl((t_0-R^2,t_0)\times\gamma^-\bigr).
	\end{equation}
	
	In \eqref{eq:localized-weak-solution} and \eqref{2608052010}, every
	coefficient and every datum is multiplied by \(\varphi\) or by one of its
	derivatives; since \(\operatorname{supp}\varphi\) is a compact subset of
	\(Q_R^{\overline\cD}(Z_0)\), all the terms are well defined on
	\(Q_R^{\cD}(Z_0)\) and are understood to be extended by zero outside this
	set.  The same convention applies to \(h\varphi\).
\end{defn}

\begin{corollary}[Local Green formula and trace bounds]
	\label{cor:trace-bound-scaled-cylinder}
	Let \(f\in L^2(Q_R^\cD(Z_0))\).  Then \(f\) is a weak solution of
	\eqref{lineqtn}, in the sense of Definition~\ref{def:L2sol}, if and only if
	\(\nabla_v f\in L^2_{\mathrm{loc}}(Q_R^{\overline\cD}(Z_0))\) and there exist
	\[
	f(t_0,\cdot)
	\in
	L^2_{\mathrm{loc}}
	\big((\overline\cD\cap B_{R^3}(x_0))\times B_R(v_0)\big),
	\qquad
	f_\gamma
	\in
	L^2_{\mathrm{loc}}
	\big((\bR\times\gamma)\cap Q_R(Z_0),\dd\gamma\dd t\big),
	\]
	such that
	\[
	f_\gamma=h
	\qquad
	\text{on }(\bR\times\gamma^-)\cap Q_R(Z_0),
	\]
	and, for every
	$\varphi\in C_c^\infty(Q_R^{\overline\cD}(Z_0))$,
	one has
	\begin{align}\label{260701105}
		\begin{aligned}
			&\iint_{(\cD\cap B_{R^3}(x_0))\times B_R(v_0)}
			f(t_0,x,v)\varphi(t_0,x,v)\dd x\dd v\\
			&\quad +
			\iint_{(\bR\times\gamma)\cap Q_R(Z_0)}
			f_\gamma\,\varphi|_\gamma\,\operatorname{sgn}(n_x\cdot v)
			\dd\gamma\dd t\\
			&=
			\iiint_{Q_R^\cD(Z_0)}
			f(\partial_t+v\cdot\nabla_x)\varphi
			\dd Z\\
			&\quad
			+
			\iiint_{Q_R^\cD(Z_0)}
			\left(
			-a^{ij}D_jfD_i\varphi
			+b^iD_if\,\varphi
			+cf\varphi
			+g^0\varphi
			-g^iD_i\varphi
			\right)
			\dd Z .
		\end{aligned}
	\end{align}
	
	Furthermore, if \(f\in L^\infty(Q_R^\cD(Z_0))\), then
	\begin{alignat}{2}
		&
		\operatorname*{ess\,sup}_{(\bR\times\gamma)\cap Q_R(Z_0)}
		f_\gamma
		\le
		\operatorname*{ess\,sup}_{Q_R^\cD(Z_0)} f,
		\qquad
		&&
		\operatorname*{ess\,sup}_{(\cD\cap B_{R^3}(x_0))\times B_R(v_0)}
		f(t_0,\cdot)
		\le
		\operatorname*{ess\,sup}_{Q_R^\cD(Z_0)} f,
		\label{2606301239}
		\\
		&
		\operatorname*{ess\,inf}_{(\bR\times\gamma)\cap Q_R(Z_0)}
		f_\gamma
		\ge
		\operatorname*{ess\,inf}_{Q_R^\cD(Z_0)} f,
		\qquad
		&&
		\operatorname*{ess\,inf}_{(\cD\cap B_{R^3}(x_0))\times B_R(v_0)}
		f(t_0,\cdot)
		\ge
		\operatorname*{ess\,inf}_{Q_R^\cD(Z_0)} f .
		\label{2606301240}
	\end{alignat}
\end{corollary}

We use the convention that, for a null set \(E\),
\[
\operatorname*{ess\,sup}_{E} f:=-\infty,
\qquad
\operatorname*{ess\,inf}_{E} f:=+\infty .
\]

Since \(C_c^\infty(Q_{R'}^{\overline\cD}(Z_0))\subset
C_c^\infty(Q_R^{\overline\cD}(Z_0))\) for \(0<R'\le R\),
Corollary~\ref{cor:trace-bound-scaled-cylinder} shows in particular that the
restriction of a weak solution of \eqref{lineqtn} in \(Q_R^\cD(Z_0)\) to
\(Q_{R'}^\cD(Z_0)\) is again a weak solution.  We use this without further
comment.

\begin{proof}[\textbf{Proof of Corollary~\ref{cor:trace-bound-scaled-cylinder}}]
	Set
	\[
	I_R:=(t_0-R^2,t_0).
	\]
	
	We first prove the equivalence.  Assume that \(f\) is a weak solution in the
	sense of Definition~\ref{def:L2sol}.  For every
	\(\eta\in C_c^\infty(Q_R^{\overline\cD}(Z_0))\), Definition~\ref{def:L2sol}
	gives \(f\eta\in \cH_{\kin}(I_R\times\cD\times\bR^d_v)\).
	In particular,
	$\nabla_v f\in L^2_{\mathrm{loc}}(Q_R^{\overline\cD}(Z_0))$.
	By Proposition~\ref{H_kin=weak}, \(f\eta\) has terminal and kinetic boundary
	traces and satisfies the Green formula.  The uniqueness of traces implies that
	the traces obtained from different cutoffs agree on overlaps.  Hence they patch
	together to define
	\[
	f(t_0,\cdot)
	\in
	L^2_{\mathrm{loc}}
	\big((\overline\cD\cap B_{R^3}(x_0))\times B_R(v_0)\big),
	\qquad
	f_\gamma
	\in
	L^2_{\mathrm{loc}}
	\big((\bR\times\gamma)\cap Q_R(Z_0),\dd\gamma\dd t\big),
	\]
	characterized by \((f\eta)(t_0,\cdot)=f(t_0,\cdot)\eta(t_0,\cdot)\) and
	\((f\eta)_\gamma=f_\gamma\,\eta|_\gamma\).  Moreover,
	\[
	f_\gamma=h
	\qquad
	\text{on }(\bR\times\gamma^-)\cap Q_R(Z_0),
	\]
	by the boundary condition in Definition~\ref{def:L2sol}.
	
	Now fix \(\varphi\in C_c^\infty(Q_R^{\overline\cD}(Z_0))\) and choose
	\(\eta\in C_c^\infty(Q_R^{\overline\cD}(Z_0))\) with \(\eta\equiv1\) on a
	relative neighborhood of \(\operatorname{supp}\varphi\).
	Applying Proposition~\ref{H_kin=weak} to \(f\eta\) and using
	\eqref{eq:localized-weak-solution}, we obtain \eqref{260701105}.
	
	Conversely, suppose that the stated local traces exist, that
	\[
	\nabla_v f\in L^2_{\mathrm{loc}}(Q_R^{\overline\cD}(Z_0)),
	\qquad
	f_\gamma=h
	\quad
	\text{on }(\bR\times\gamma^-)\cap Q_R(Z_0),
	\]
	and that \eqref{260701105} holds for every
	\(\varphi\in C_c^\infty(Q_R^{\overline\cD}(Z_0))\).
	Fix \(\eta\in C_c^\infty(Q_R^{\overline\cD}(Z_0))\), and set \(u:=f\eta\),
	extended by zero.  Then
	$u\in L^2(I_R\times\cD;H_v^1)$,
	because \(\operatorname{supp}\eta\) is a compact subset of
	\(Q_R^{\overline\cD}(Z_0)\) and \(\nabla_vf\in
	L^2_{\mathrm{loc}}(Q_R^{\overline\cD}(Z_0))\).
	Define
	\[
	u(t_0,\cdot):=f(t_0,\cdot)\eta(t_0,\cdot)\in L^2(\cD\times \bR^d_v),
	\qquad
	u_\gamma:=f_\gamma\,\eta|_\gamma\in L^2(I_R\times\gamma,\dd\gamma\dd t),
	\]
	again extended by zero; both are square integrable because \(\eta\) has
	compact support in \(Q_R^{\overline\cD}(Z_0)\).
	
	Let \(\phi\in C_c^\infty([t_0-R^2,t_0]\times\overline\cD\times\bR^d_v)\).
	Then \(\eta\phi\in C_c^\infty(Q_R^{\overline\cD}(Z_0))\), so \eqref{260701105}
	may be applied with \(\varphi=\eta\phi\).  Expanding the derivatives of
	\(\eta\phi\) and rearranging the resulting terms exactly as in the derivation
	of \eqref{2608052010}, we obtain
	\begin{align}\label{2608051032}
		\begin{aligned}
			&\iint_{\cD\times\bR^d_v}u(t_0,\cdot)\,\phi(t_0,\cdot)\dd x\dd v
			+\iint_{I_R\times\gamma}
			u_\gamma\,\phi|_\gamma\,\sgn(n_x\cdot v)\dd\gamma\dd t
			\\
			=\;&
			\iiint_{I_R\times\cD\times\bR^d_v}
			u\,(\partial_t+v\cdot\nabla_x)\phi\dd Z
			\\
			&\quad
			+\iiint_{I_R\times\cD\times\bR^d_v}
			\bigl(-a^{ij}D_juD_i\phi+b^iD_iu\,\phi+cu\phi
			+G^0\phi-G^iD_i\phi\bigr)\dd Z,
		\end{aligned}
	\end{align}
	where \(G^0,\ldots,G^d\in L^2(I_R\times\cD\times\bR^d_v)\) are the
	coefficients of \(g^0\eta+D_i(g^i\eta)+\cR_\eta[f]\) written in
	divergence form, namely
	\[
	G^i:=g^i\eta-a^{ij}fD_j\eta,
	\qquad
	G^0:=g^0\eta+f(\partial_t+v\cdot\nabla_x)\eta
	-a^{ij}D_jfD_i\eta-b^ifD_i\eta-g^iD_i\eta .
	\]
	Taking \(\phi\) supported in \(I_R\times\cD\times\bR^d_v\) in
	\eqref{2608051032} shows that
	\[
	(\partial_t+v\cdot\nabla_x)u
	=
	D_i(a^{ij}D_ju)+b^iD_iu+cu+G^0+D_iG^i
	\]
	holds in the sense of distributions, and the right-hand side belongs to
	\(L^2(I_R\times\cD;H_v^{-1})\); hence
	\((\partial_t+v\cdot\nabla_x)u\in L^2(I_R\times\cD;H_v^{-1})\).  Since
	\(u(t_0-R^2,\cdot)=0\) and \(\eta\) vanishes near \(t=t_0-R^2\),
	\eqref{2608051032} is precisely the identity \eqref{260618611} for \(u\),
	with terminal and initial data \(u(t_0,\cdot)\) and \(0\).  By
	Proposition~\ref{H_kin=weak},
	$u\in \cH_{\kin}(I_R\times\cD\times\bR^d_v)$,
	and its terminal value and trace are the functions \(u(t_0,\cdot)\) and
	\(u_\gamma=f_\gamma\,\eta|_\gamma\) defined above.
	In particular \(u(t_0-R^2,\cdot)=0\) and
	$u_{\gamma^-}=f_\gamma\eta=h\eta$ on $I_R\times\gamma^-$.
	Therefore Definition~\ref{def:L2sol} holds, and \(f\) is a weak solution of
	\eqref{lineqtn}.
	
	It remains to prove the \(L^\infty\)-trace bounds.  Assume
	\(f\in L^\infty(Q_R^\cD(Z_0))\).  Fix \(0<\tau<R\), and choose
	$\eta\in C_c^\infty(Q_R^{\overline\cD}(Z_0))$
	such that
	$0\le \eta\le 1$, and 
	$\eta\equiv 1$
	on $Q_\tau^{\overline\cD}(Z_0)$.
	Set
	\[
	M:=\operatorname*{ess\,sup}_{Q_R^\cD(Z_0)} f,
	\qquad
	u:=(f-M)\eta .
	\]
	Since \(\eta\in C_c^\infty(Q_R^{\overline\cD}(Z_0))\) satisfies
	\(0\le\eta\le1\) and \(f\le M\) a.e. in \(Q_R^\cD(Z_0)\), we have
	\(u\le0\) a.e. in \(I_R\times\cD\times\bR^d_v\); moreover
	\[
	u=f\eta-M\eta\in \cH_{\kin}(I_R\times\cD\times\bR^d_v),
	\]
	since \(f\eta\in\cH_\kin\) by Definition~\ref{def:L2sol} and
	\(\eta\in C_c^\infty([t_0-R^2,t_0]\times\overline\cD\times\bR^d_v)\).
	By the first part of the proof, the traces of \(u\) are
	\[
	u(t_0,\cdot)
	=
	(f(t_0,\cdot)-M)\eta(t_0,\cdot),
	\qquad
	u_\gamma
	=
	(f_\gamma-M)\eta|_\gamma .
	\]
	Applying Proposition~\ref{2606291053} with
	\(\cO=[t_0-R^2,t_0]\times\overline\cD\times\bR^d_v\), we obtain
	\(u_\gamma\le0\) \(\dd\gamma\dd t\)-a.e. and \(u(t,\cdot)\le0\) for a.e.
	\(t\in I_R\); since \(u\in C([t_0-R^2,t_0];L^2(\cD\times\bR^d_v))\) by
	Proposition~\ref{trace_theorem}, letting \(t\uparrow t_0\) also gives
	\(u(t_0,\cdot)\le0\).  Since
	\(\eta\equiv1\) on \(Q_\tau^{\overline\cD}(Z_0)\), we get
	\[
	\operatorname*{ess\,sup}_{(\bR\times\gamma)\cap Q_\tau(Z_0)}
	f_\gamma
	\le M,
	\]
	and
	\[
	\operatorname*{ess\,sup}_{(\cD\cap B_{\tau^3}(x_0))\times B_\tau(v_0)}
	f(t_0,\cdot)
	\le M.
	\]
	Letting \(\tau\uparrow R\) and using \eqref{2606301225}, we obtain
	\eqref{2606301239}.  Applying the same argument to \(-f\) gives
	\eqref{2606301240}.  This completes the proof.
\end{proof}

\begin{lemma}\label{260312141}
	Let \(f\) be a weak solution of \eqref{lineqtn} in the sense of
	Definition~\ref{def:L2sol}. Let \(\beta\in C^2(\bR)\) satisfy
	\[
	\beta(0)=\beta'(0)=0,
	\qquad
	\sup_{\bR}\bigl(|\beta'|+|\beta''|\bigr)<\infty.
	\]
	Then, for every \(\varphi\in C_c^\infty(Q_R^{\overline\cD}(Z_0))\), the following
	identity holds:
	\begin{align}
		\begin{split}\label{2601161211}
			&\iint_{(\cD\cap B_{R^3}(x_0))\times B_R(v_0)}
			\beta\bigl(f(t_0,\cdot)\bigr)\varphi(t_0,\cdot)\dd x\dd v
			+
			\iint_{(\bR\times\gamma)\cap Q_R(Z_0)}
			\beta(f_\gamma)\varphi\,\operatorname{sgn}(n_x\cdot v)\dd\gamma\dd t
			\\
			=\;&
			\iiint_{Q_R^\cD(Z_0)}
			\beta(f)(\partial_t\varphi+v\cdot\nabla_x\varphi)\dd Z
			+
			\iiint_{Q_R^\cD(Z_0)}
			\beta'(f)(b^iD_i f+cf+g^0)\varphi\dd Z
			\\
			&\quad-
			\iiint_{Q_R^\cD(Z_0)}
			\Bigl(
			\beta''(f)a^{ij}D_i fD_j f\,\varphi
			+
			\beta'(f)a^{ij}D_jfD_i\varphi
			\Bigr)\dd Z
			\\
			&\quad-
			\iiint_{Q_R^\cD(Z_0)}
			\Bigl(
			\beta''(f)D_i f\,g^i\varphi
			+
			\beta'(f)g^iD_i\varphi
			\Bigr)\dd Z.
		\end{split}
	\end{align}
	
	In particular, assume in addition that \(\beta\ge0\) on \(\bR\) and
	\[
	\beta(f_{\gamma^-})=0
	\quad
	\dd\gamma\dd t\text{-a.e. on }(\bR\times\gamma^-)\cap Q_R(Z_0).
	\]
	Then, for every nonnegative
	\(\varphi\in C_c^\infty(Q_R^{\overline\cD}(Z_0))\),
	\begin{align}\label{260113330}
		\begin{split}
			&\iiint_{Q_R^\cD(Z_0)}
			\beta(f)(-\partial_t\varphi-v\cdot\nabla_x\varphi)\dd Z
			\\
			&\quad+
			\iiint_{Q_R^\cD(Z_0)}
			\Bigl(
			\beta''(f)a^{ij}D_i fD_j f\,\varphi
			+
			\beta'(f)a^{ij}D_jfD_i\varphi
			\Bigr)\dd Z
			\\
			\le\;&
			-\iiint_{Q_R^\cD(Z_0)}
			\Bigl(
			\beta''(f)D_i f\,g^i\varphi
			+
			\beta'(f)g^iD_i\varphi
			\Bigr)\dd Z
			\\
			&\quad+
			\iiint_{Q_R^\cD(Z_0)}
			\beta'(f)(b^iD_i f+cf+g^0)\varphi\dd Z.
		\end{split}
	\end{align}
\end{lemma}

\begin{proof}
	Set
	\[
	I_R:=(t_0-R^2,t_0),
	\qquad
	\cQ_R:=Q_R^\cD(Z_0).
	\]
	Fix \(\varphi\in C_c^\infty(Q_R^{\overline\cD}(Z_0))\). Choose
	\(\eta\in C_c^\infty(Q_R^{\overline\cD}(Z_0))\) such that
	$\eta\equiv1$ in a neighborhood of $\operatorname{supp}\varphi$
	relative to $\bR\times\overline\cD\times\bR^d_v$.
	Set $u:=f\eta$.
	By Definition~\ref{def:L2sol}, we have
	\[
	u\in \cH_{\kin}(I_R\times\cD\times\bR_v^d),
	\qquad
	u(t_0-R^2,\cdot)=0,
	\qquad
	u_\gamma=f_\gamma\eta
	\quad\text{on }I_R\times\gamma.
	\]
	
	We first compute the equation satisfied by \(u\). Since
	\[
	(\partial_t+v\cdot\nabla_x)f
	=
	D_i(a^{ij}D_j f)+b^iD_i f+cf+g^0+D_i g^i
	\quad\text{in }\cQ_R,
	\]
	we get
	\[
	(\partial_t+v\cdot\nabla_x)u
	=
	D_i(a^{ij}D_j u)+b^iD_i u+cu+G_\eta^0+D_iG_\eta^i
	\quad\text{in }I_R\times\cD\times\bR_v^d,
	\]
	where
	\[
	G_\eta^i
	:=
	g^i\eta-a^{ij}fD_j\eta,
	\]
	and
	\[
	G_\eta^0
	:=
	g^0\eta
	+
	f(\partial_t\eta+v\cdot\nabla_x\eta)
	-
	a^{ij}D_jfD_i\eta
	-
	b^ifD_i\eta
	-
	g^iD_i\eta .
	\]
	This is \eqref{2608052010} with \(\varphi\) replaced by \(\eta\), after
	using \(a^{ij}D_j(fD_i\eta)-a^{ij}fD_{ij}\eta=a^{ij}D_jf\,D_i\eta\).
	
	We now apply the renormalized identity \eqref{reno_dyn} to \(u\) on
	\(I_R\times\cD\times\bR^d_v\), with the test function \(\varphi\) and with
	the divergence-form data
	\[
	\tilde g^i:=a^{ij}D_ju+G_\eta^i,
	\qquad
	\tilde g^0:=b^iD_iu+cu+G_\eta^0,
	\]
	which belong to \(L^2(I_R\times\cD\times\bR^d_v)\).  This gives
	\begin{align*}
		&\iint_{\cD\times\bR^d_v}
		\beta\bigl(u(t_0,\cdot)\bigr)\varphi(t_0,\cdot)\dd x\dd v
		+
		\iint_{I_R\times\gamma}
		\beta(u_\gamma)\varphi\,\operatorname{sgn}(n_x\cdot v)\dd\gamma\dd t
		\\
		=\;&
		\iiint_{I_R\times\cD\times\bR_v^d}
		\beta(u)(\partial_t\varphi+v\cdot\nabla_x\varphi)\dd Z
		+
		\iiint_{I_R\times\cD\times\bR_v^d}
		\beta'(u)(b^iD_i u+cu+G_\eta^0)\varphi\dd Z
		\\
		&\quad-
		\iiint_{I_R\times\cD\times\bR_v^d}
		\Bigl(
		\beta''(u)a^{ij}D_i uD_j u\,\varphi
		+
		\beta'(u)a^{ij}D_juD_i\varphi
		\Bigr)\dd Z
		\\
		&\quad-
		\iiint_{I_R\times\cD\times\bR_v^d}
		\Bigl(
		\beta''(u)D_i u\,G_\eta^i\varphi
		+
		\beta'(u)G_\eta^iD_i\varphi
		\Bigr)\dd Z .
	\end{align*}
	There is no initial-time term because \(u(t_0-R^2,\cdot)=0\) and
	\(\beta(0)=0\).
	
	Since \(\eta\equiv1\) in a relative neighborhood of
	\(\operatorname{supp}\varphi\), we have \(D\eta=0\) and
	\((\partial_t+v\cdot\nabla_x)\eta=0\) there, whence, on
	\(\operatorname{supp}\varphi\),
	\[
	u=f,\qquad
	\nabla_vu=\nabla_vf,\qquad
	G_\eta^0=g^0,\qquad
	G_\eta^i=g^i,
	\qquad
	u_\gamma=f_\gamma .
	\]
	Since every integrand in the identity above carries a factor \(\varphi\),
	\(D_i\varphi\) or \(\partial_t\varphi+v\cdot\nabla_x\varphi\), they all vanish
	outside \(\operatorname{supp}\varphi\subset Q_R^{\overline\cD}(Z_0)\), and the
	identity reduces exactly to \eqref{2601161211}.
	
	It remains to prove \eqref{260113330}. Assume that \(\beta\ge0\),
	\(\varphi\ge0\), and
	\[
	\beta(f_{\gamma^-})=0
	\quad
	\dd\gamma\dd t\text{-a.e. on }((t_0-R^2,t_0)\times\gamma^-)\cap Q_R(Z_0).
	\]
	Then the final-time term in \eqref{2601161211} is nonnegative. Moreover, the
	boundary term is also nonnegative, since its contribution on \(\gamma^-\)
	vanishes and its contribution on \(\gamma^+\) is nonnegative. Dropping these
	two terms from \eqref{2601161211} gives \eqref{260113330}.
\end{proof}

\begin{lemma}
	\label{2605311047}
	Let \(Z_0\in \bR\times\overline{\cD}\times\bR^d\) and \(R>0\). Assume that
	\[
	a^{ij},\,b^i,\,c\in L^\infty(Q_R^\cD(Z_0)),
	\qquad i,j=1,\dots,d,
	\]
	and that there exist constants \(0<\lambda\le \Lambda\) and \(M\ge0\)
	such that, for a.e. \(Z\in Q_R^\cD(Z_0)\) and every \(\xi\in\bR^d\),
	\begin{align}\label{260702750}
		\lambda |\xi|^2
		\le
		a^{ij}(Z)\xi_i\xi_j,
		\qquad
		|a^{ij}(Z)|\le \Lambda,
		\qquad
		|b(Z)|\le M,
		\qquad
		|c(Z)|\le M^2.
	\end{align}
	Let
	\[
	g^0,\ldots,g^d\in L^2(Q_R^\cD(Z_0)),
	\qquad
	h\in L^2((\bR\times\gamma^-)\cap Q_R(Z_0),\mathrm d\gamma\dd t),
	\]
	and let \(f\) be a weak solution of
	\[
	\begin{cases}
		(\partial_t+v\cdot\nabla_x)f
		=
		D_i(a^{ij}D_j f)+b^iD_i f+cf+g^0+D_ig^i
		&\quad\text{in }Q_R^\cD(Z_0),\\
		f=h
		&\quad\text{on }(\bR\times\gamma^-)\cap Q_R(Z_0).
	\end{cases}
	\]
	
	Let
	\[
	a_k^{ij},\,b_k^i,\,c_k\in L^\infty(Q_R^\cD(Z_0)),
	\qquad i,j=1,\dots,d,
	\]
	satisfy the same bounds as in \eqref{260702750}, with
	\((a,b,c)\) replaced by \((a_k,b_k,c_k)\). Assume that
	\[
	a_k^{ij}\to a^{ij},
	\qquad
	b_k^i\to b^i,
	\qquad
	c_k\to c
	\quad \text{a.e. in } Q_R^\cD(Z_0).
	\]
	Let
	\[
	g_k^i\in L^2(Q_R^\cD(Z_0)),
	\qquad i=0,\dots,d,
	\]
	and
	\[
	h_k\in L^2((\bR\times\gamma^-)\cap Q_R(Z_0),\mathrm d\gamma\dd t)
	\]
	satisfy
	\[
	g_k^i\to g^i
	\quad \text{in }L^2(Q_R^\cD(Z_0)),
	\qquad i=0,\dots,d,
	\]
	and
	\[
	h_k\to h
	\quad \text{in }
	L^2((\bR\times\gamma^-)\cap Q_R(Z_0),\mathrm d\gamma\dd t).
	\]
	
	Then, for each \(\kappa\in(0,1)\), there exists a sequence \(\{f_k\}\) such
	that \(f_k\) is a weak solution of
	\[
	\begin{cases}
		(\partial_t+v\cdot\nabla_x)f_k
		=
		D_i(a_k^{ij}D_j f_k)+b_k^iD_i f_k+c_kf_k+g_k^0+D_ig_k^i
		&\quad\text{in }Q_{\kappa R}^\cD(Z_0),\\
		f_k=h_k
		&\quad\text{on }(\bR\times\gamma^-)\cap Q_{\kappa R}(Z_0),
	\end{cases}
	\]
	and
	\[
	f_k\to f,
	\qquad
	D_v f_k\to D_v f
	\quad \text{in }L^2(Q_{\kappa R}^\cD(Z_0)).
	\]
\end{lemma}

\begin{proof}
	Set
	\[
	I_R:=(t_0-R^2,t_0),
	\qquad
	\cQ_R:=Q_R^\cD(Z_0),
	\qquad
	\cQ_\kappa:=Q_{\kappa R}^\cD(Z_0),
	\]
	and
	\[
	\Gamma_R^-:=(\bR\times\gamma^-)\cap Q_R(Z_0).
	\]
	We extend all the data to \(I_R\times\cD\times\bR_v^d\) by setting, outside
	\(\cQ_R\),
	\(a^{ij}=a_k^{ij}:=\lambda\delta^{ij}\), \(b^i=b_k^i:=0\), \(c=c_k:=0\) and
	\(g^i=g_k^i:=0\) for \(i=0,\dots,d\).  The extended
	coefficients still satisfy \eqref{260702750} and converge a.e., and the
	extended source terms still converge in
	\(L^2(I_R\times\cD\times\bR^d_v)\).
	
	Choose a smooth cutoff
	$\eta\in C_c^\infty(Q_R^{\overline\cD}(Z_0))$ such that
	$\eta\equiv1
	\quad\text{on }Q_{\kappa R}^{\overline\cD}(Z_0)$.
	Set
	\[
	u:=f\eta.
	\]
	By Definition~\ref{def:L2sol},
	\[
	u\in \cH_{\kin}(I_R\times\cD\times\bR_v^d),
	\qquad
	u(t_0-R^2,\cdot)=0,
	\qquad
	u_{\gamma^-}=h\eta
	\quad\text{on }I_R\times\gamma^-.
	\]
	
	Define
	\[
	\cL u
	:=
	D_i(a^{ij}D_j u)+b^iD_i u+cu,
	\qquad
	\cL_k u
	:=
	D_i(a_k^{ij}D_j u)+b_k^iD_i u+c_ku.
	\]
	Then the localized weak formulation in Definition~\ref{def:L2sol} gives
	\[
	(\partial_t+v\cdot\nabla_x)u
	=
	\cL u+D_iF^i+F^0
	\quad\text{in }I_R\times\cD\times\bR_v^d,
	\]
	where
	\begin{align*}
		\begin{aligned}
			F^i
			&:=
			g^i\eta-a^{ij}fD_j\eta,
			\\
			F^0
			&:=
			g^0\eta
			+
			f(\partial_t+v\cdot\nabla_x)\eta
			-
			a^{ij}D_j fD_i\eta
			-
			b^i fD_i\eta
			-
			g^iD_i\eta .
		\end{aligned}
	\end{align*}
	Here \(F=(F^i)_{i=1}^d\) and \(F^0\) belong to
	\(L^2(I_R\times\cD\times\bR_v^d)\).
	
	For each \(k\), set
	\begin{align*}
		F_k^i
		&:=
		g_k^i\eta-a_k^{ij}fD_j\eta,
		\\
		F_k^0
		&:=
		g_k^0\eta
		+
		f(\partial_t+v\cdot\nabla_x)\eta
		-
		a_k^{ij}D_j fD_i\eta
		-
		b_k^i fD_i\eta
		-
		g_k^iD_i\eta .
	\end{align*}
	Since the extended coefficients satisfy \eqref{coeff_ellip} and
	\eqref{coeff_bdd} with \(\nu=\min\{\lambda,\Lambda^{-1}\}\) and
	\(K=\max\{M,M^2\}\), and since \(\bar b^i=0\), Theorem~\ref{solv}\,(1)
	provides a unique
	\[
	u_k\in \cH_{\kin}(I_R\times\cD\times\bR_v^d)
	\]
	such that
	$(\partial_t+v\cdot\nabla_x)u_k
	=
	\cL_ku_k+D_iF_k^i+F_k^0
	\quad\text{in }I_R\times\cD\times\bR_v^d$,
	with
	\[
	u_k(t_0-R^2,\cdot)=0,
	\qquad
	(u_k)_{\gamma^-}=h_k\eta
	\quad\text{on }I_R\times\gamma^-.
	\]
	
	Since \(\eta\equiv1\) on \(Q_{\kappa R}^{\overline\cD}(Z_0)\), we have
	\[
	D_v\eta=0,
	\qquad
	(\partial_t+v\cdot\nabla_x)\eta=0
	\qquad
	\text{on }\cQ_\kappa.
	\]
	Thus
	\[
	F_k^i=g_k^i,
	\qquad
	F_k^0=g_k^0
	\quad\text{on }\cQ_\kappa.
	\]
	Since \(u_k\in\cH_\kin(I_R\times\cD\times\bR^d_v)\), Propositions~\ref{trace_theorem}
	and \ref{H_kin=weak} give the Green identity \eqref{260701105}
	for \(u_k\) on \(\cQ_\kappa\).  Hence, denoting by \(f_k\) the
	restriction of \(u_k\) to \(\cQ_\kappa\),
	Corollary~\ref{cor:trace-bound-scaled-cylinder} shows that \(f_k\) is a weak
	solution of
	\[
	(\partial_t+v\cdot\nabla_x)f_k
	=
	D_i(a_k^{ij}D_j f_k)+b_k^iD_i f_k+c_kf_k+g_k^0+D_ig_k^i
	\quad\text{in }\cQ_\kappa,
	\]
	with
	\[
	f_k=h_k
	\quad\text{on }(\bR\times\gamma^-)\cap Q_{\kappa R}(Z_0).
	\]
	
	It remains to prove the convergence. Since
	\[
	(\partial_t+v\cdot\nabla_x)u
	=
	\cL_ku
	+
	D_i\bigl((a^{ij}-a_k^{ij})D_j u\bigr)
	+
	(b^i-b_k^i)D_i u
	+
	(c-c_k)u
	+
	D_iF^i+F^0,
	\]
	we have
	\begin{align*}
		(\partial_t+v\cdot\nabla_x)(u-u_k)
		&=
		\cL_k(u-u_k)
		+
		D_iE_k^i
		+
		E_k^0,
	\end{align*}
	where
	\[
	E_k^i
	:=
	(a^{ij}-a_k^{ij})D_j u+F^i-F_k^i
	=
	(a^{ij}-a_k^{ij})\eta D_jf+(g^i-g_k^i)\eta,
	\]
	and
	\begin{align*}
		E_k^0
		:={}&
		(b^i-b_k^i)D_i u
		+
		(c-c_k)u
		+
		F^0-F_k^0
		\\
		={}&
		(b^i-b_k^i)\eta D_i f
		+
		(c-c_k)u
		+
		(g^0-g_k^0)\eta
		\\
		&\quad
		-
		(a^{ij}-a_k^{ij})D_j fD_i\eta
		-
		(g^i-g_k^i)D_i\eta .
	\end{align*}
	The initial trace of \(u-u_k\) is zero, while its inflow trace is
	$(h-h_k)\eta$.
	Therefore the energy estimate \eqref{dyHbdd} gives
	\begin{align}\label{260531516}
		\|u-u_k\|_{\cH_{\kin}(I_R\times\cD\times\bR_v^d)}
		\lesssim\;&
		\|(h-h_k)\eta\|_{L^2(I_R\times\gamma^-,\,\mathrm d\gamma\dd t)}
		\nonumber\\
		&+
		\sum_{i=1}^d\|E_k^i\|_{L^2(I_R\times\cD\times\bR_v^d)}
		+
		\|E_k^0\|_{L^2(I_R\times\cD\times\bR_v^d)}.
	\end{align}
	The implicit constant is independent of \(k\), by the uniform ellipticity and
	the uniform \(L^\infty\)-bounds of \(a_k,b_k,c_k\).
	
	By assumption,
	\[
	h_k\to h
	\quad\text{in }L^2(\Gamma_R^-,\mathrm d\gamma\dd t),
	\qquad
	g_k^i\to g^i
	\quad\text{in }L^2(\cQ_R),
	\quad i=0,\dots,d.
	\]
	Since \(\eta\) is supported in \(Q_R(Z_0)\), the boundary term in
	\eqref{260531516} tends to zero, and the terms involving \(g_k^i-g^i\) also
	tend to zero. Moreover, \(u,D_vu\in L^2\), and the uniform bounds together
	with the a.e. convergence of \(a_k,b_k,c_k\) imply, by the dominated
	convergence theorem, that
	\[
	\sum_{i=1}^d\|E_k^i\|_{L^2(I_R\times\cD\times\bR_v^d)}
	+
	\|E_k^0\|_{L^2(I_R\times\cD\times\bR_v^d)}
	\to0.
	\]
	Hence \eqref{260531516} yields
	$u_k\to u$ and 
	$D_vu_k\to D_vu$ in $L^2(\cQ_\kappa)$.
	Finally, \(\eta\equiv1\) on \(\cQ_\kappa\), so \(u=f\) and
	\(D_vu=D_vf\) there. Therefore,
	$f_k\to f$ and 
	$D_vf_k\to D_vf$ in $L^2(Q_{\kappa R}^\cD(Z_0))$.
	This completes the proof.
\end{proof}

\subsection{$L^2$--$L^\infty$ estimate}\label{L2_Linfty_estimate}

The local \(L^2\)--\(L^\infty\) estimate is standard in the regularity
theory for kinetic equations with rough coefficients.  The whole-space
estimate can be found in \cite{2022AR}, following the Sobolev--Moser
approach of \cite{PascucciPolidoro2004}, and the corresponding
homogeneous-inflow estimate is given in
\cite[Theorem~1.1]{SilvestreBoundaryKFP}.
We revisit the standard Sobolev--Moser argument in the weak-solution framework used here.

\begin{thm}[\(L^2\)--\(L^\infty\) estimate]\label{cor:linfty-zero-inflow}
	Let \(Z_0\in \bR^{1+2d}\) and \(R>0\). Assume that
	$a^{ij},\,b^i,\,c\in L^\infty(Q_R^\cD(Z_0))$, $i,j=1,\dots,d$,
	and that there exist constants \(0<\lambda\le \Lambda\) and \(M\ge0\) such that,
	for a.e. \(Z\in Q_R^\cD(Z_0)\) and every \(\xi\in\bR^d\),
	\begin{equation}\label{con_abc}
		\lambda |\xi|^2
		\le
		a^{ij}(Z)\xi_i\xi_j,
		\qquad
		|a^{ij}(Z)|\le \Lambda,
		\qquad
		|b(Z)|\le M,
		\qquad
		|c(Z)|\le M^2.
	\end{equation}
	Let \(f\) be a weak solution of
	\[
	\begin{cases}
		(\partial_t +v\cdot \nabla_x)f
		=
		D_i(a^{ij}D_j f)+b^iD_i f+cf+g
		&\quad\text{in }Q_R^{\cD}(Z_0),\\
		f=h\le 0
		&\quad\text{on }(\bR\times\gamma^-)\cap Q_R(Z_0),
	\end{cases}
	\]
	where \(g\in L^\infty(Q_R^\cD(Z_0))\). Then, for every \(0<\kappa<1\),
	\begin{align}\label{260628242}
		\|f_+\|_{L^\infty(Q_{\kappa R}^\cD(Z_0))}
		\le
		C(d,\lambda,\Lambda,\kappa)\,(MR+1)^{4d+2}R^{-1-2d}
		\left(
		\|f_+\|_{L^2(Q_R^\cD(Z_0))}
		+
		R^{3+2d}\|g\|_{L^\infty(Q_R^\cD(Z_0))}
		\right).
	\end{align}
	In particular, if \(h\equiv0\) on \((\bR\times\gamma^-)\cap Q_R(Z_0)\),
	then
	\begin{align}\label{260628241}
		\|f\|_{L^\infty(Q_{\kappa R}^\cD(Z_0))}
		\le
		C(d,\lambda,\Lambda,\kappa)\,(MR+1)^{4d+2}R^{-1-2d}
		\left(
		\|f\|_{L^2(Q_R^\cD(Z_0))}
		+
		R^{3+2d}\|g\|_{L^\infty(Q_R^\cD(Z_0))}
		\right).
	\end{align}
\end{thm}

After the auxiliary lemmas and the Moser gain of
Proposition~\ref{cor:boundary-moser-gain}, we prove
Theorem~\ref{cor:linfty-zero-inflow}.

\begin{lemma}
	\label{lem:global-potential-AR}
	Let
	\begin{align}\label{260424252}
		\mathfrak G(t,x,v):=
		\left(\frac{3}{4\pi^2 t^4}\right)^{d/2}
		\exp\!\left(
		-\frac14\left|\frac{v}{t^{1/2}}\right|^2
		-3\left|\frac{x}{t^{3/2}}-\frac{v}{2t^{1/2}}\right|^2
		\right)\mathbf 1_{\{t>0\}},
	\end{align}
	and define, for \(Z=(t,x,v)\) and \(Y=(s,y,w)\) in \(\bR^{1+2d}\),
	\begin{equation}\label{Gf}
		\begin{alignedat}{2}
			&\mathfrak G[f](Z)
			\,&&:=
			\int_{\bR^{1+2d}}
			\mathfrak G(Y^{-1}\circ Z)\,f(Y)\dd Y,
			\\
			&\mathfrak G_i[f](Z)
			\,&&:=-
			\int_{\bR^{1+2d}}
			\p_{w^i}\big(\mathfrak G(Y^{-1}\circ Z)\big)\,f(Y)\dd Y.
		\end{alignedat}
	\end{equation}
	Set
	\begin{equation}\label{eq:QpRp}
		Q_p:=\frac{(2d+1)p}{(2d+1)-p}
		\quad(p<2d+1),
		\qquad
		R_p:=\frac{(4d+2)p}{(4d+2)-p}
		\quad(p<4d+2).
	\end{equation}
	All norms below are taken over \(\bR^{1+2d}\).
	\begin{enumerate}
		\item (Weak \(L^1\) estimate) For every \(f\in L^1\),
		\begin{align}
			\label{eq:weak-L1-endpoint}
			\|\mathfrak G[f]\|_{L^{Q_1,\infty}}
			+
			\|\mathfrak G_i[f]\|_{L^{R_1,\infty}}
			\le
			C(d)\|f\|_{L^1}.
		\end{align}
		
		\item (Weak \(L^p\) estimate) For every \(p>1\) and every
		\(f\in L^{p,\infty}\),
		\begin{align}
			\label{eq:weak-Lp-endpoint}
			\|\mathfrak G[f]\|_{L^{Q_p,\infty}}
			\one_{\{p<2d+1\}}
			+
			\|\mathfrak G_i[f]\|_{L^{R_p,\infty}}
			\one_{\{p<4d+2\}}
			\le
			C(d,p)\|f\|_{L^{p,\infty}}.
		\end{align}
		
		\item (Strong \(L^p\) estimate) For every \(p>1\) and every
		\(f\in L^p\),
		\begin{align}
			\label{eq:strong-Lp-estimate}
			\|\mathfrak G[f]\|_{L^{Q_p}}
			\one_{\{p<2d+1\}}
			+
			\|\mathfrak G_i[f]\|_{L^{R_p}}
			\one_{\{p<4d+2\}}
			\le
			C(d,p)\|f\|_{L^{p}}.
		\end{align}
	\end{enumerate}
\end{lemma}

\begin{proof}
	A direct computation gives
	\[
	\p_{w^i}\bigl(\mathfrak G(Y^{-1}\circ Z)\bigr)
	=
	\mathfrak K_i(Y^{-1}\circ Z),
	\]
	where
	\[
	\mathfrak K_i(t,x,v)
	:=
	\left(
	\frac{3x_i}{t^2}
	-
	\frac{v_i}{t}
	\right)\mathfrak G(t,x,v).
	\]
	Hence
	\[
	\mathfrak G_i[f](Z)
	=
	-\int_{\bR^{1+2d}}
	\mathfrak K_i(Y^{-1}\circ Z)f(Y)\dd Y.
	\]
	
	Set
	\[
	\delta_\rho(t,x,v):=(\rho^2t,\rho^3x,\rho v).
	\]
	The explicit formulas give
	\[
	\mathfrak G(\delta_\rho Z)
	=
	\rho^{-4d}\mathfrak G(Z),
	\qquad
	\mathfrak K_i(\delta_\rho Z)
	=
	\rho^{-(4d+1)}\mathfrak K_i(Z).
	\]
	Both \(\mathfrak G\) and \(\mathfrak K_i\) are continuous on
	\(\bR^{1+2d}\setminus\{0\}\) and vanish on \(\{t\le0\}\setminus\{0\}\), so
	they are bounded on the compact set \(\{|Z|_{\kin}=1\}\).  Since
	\(|\delta_\rho Z|_{\kin}=\rho|Z|_{\kin}\), the homogeneity relations above
	give, for every \(Z\neq0\),
	\[
	|\mathfrak G(Z)|
	\le
	C|Z|_{\kin}^{-4d},
	\qquad
	|\mathfrak K_i(Z)|
	\le
	C|Z|_{\kin}^{-(4d+1)}.
	\]
	
	Write \(Q=4d+2\).  For \(0<\alpha<Q\) define
	\[
	I_\alpha h(Z):=\int_{\bR^{1+2d}}
	|Y^{-1}\circ Z|_{\kin}^{\alpha-Q}|h(Y)|\dd Y.
	\]
	The balls
	\(\mathrm{B}_r(Z)=\{Y:|Y^{-1}\circ Z|_{\kin}<r\}\)
	are left translates of the inverse of \(\mathrm{B}_r(0)\), and
	\(|\mathrm{B}_r(Z)|=c_dr^Q\): left translations and the inversion \(Z\mapsto Z^{-1}\) preserve
	Lebesgue measure, and \(\delta_r\) has Jacobian \(r^Q\).
	Directly from the group law, the homogeneous gauge satisfies
	\(|Z\circ Y|_{\kin}\le |Z|_{\kin}+|Y|_{\kin}\) and
	\(|Z^{-1}|_{\kin}\le2^{1/3}|Z|_{\kin}\).
	Thus these balls satisfy a uniform doubling condition and
	the standard disjoint-ball covering property.  Their centered Hardy--Littlewood maximal operator
	\[
	\mathcal Mh(Z):=\sup_{r>0}\frac1{|\mathrm{B}_r(Z)|}
	\int_{\mathrm{B}_r(Z)}|h(Y)|\dd Y
	\]
	is therefore of weak type \((1,1)\), and is bounded on
	\(L^{p,\infty}\) for \(1<p<\infty\).  These statements follow
	from the disjoint-ball covering argument and interpolation;
	the constants depend only on the homogeneous group.
	
	We give the potential estimate without appealing to a
	Euclidean convolution theorem.  Decomposing \(\mathrm{B}_r(Z)\)
	into dyadic annuli yields, for every \(\alpha>0\),
	\[
	\int_{\mathrm{B}_r(Z)}|Y^{-1}\circ Z|_{\kin}^{\alpha-Q}
	|h(Y)|\dd Y\le Cr^\alpha\mathcal Mh(Z).
	\]
	If \(1<p<Q/\alpha\), the Lorentz estimate
	\(\int_E|h|\le C_p\|h\|_{L^{p,\infty}}|E|^{1-1/p}\)
	on each annulus gives
	\[
	\int_{\bR^{1+2d}\setminus \mathrm{B}_r(Z)}
	|Y^{-1}\circ Z|_{\kin}^{\alpha-Q}|h(Y)|\dd Y
	\le C_{p,\alpha}\|h\|_{L^{p,\infty}}
	r^{\alpha-Q/p}.
	\]
	Taking
	\(r=(\|h\|_{L^{p,\infty}}/\mathcal Mh(Z))^{p/Q}\)
	when \(0<\mathcal Mh(Z)<\infty\), we find
	\[
	I_\alpha h(Z)\le C_{p,\alpha}
	\|h\|_{L^{p,\infty}}^{\alpha p/Q}
	(\mathcal Mh(Z))^{1-\alpha p/Q}.
	\]
	This also holds when \(\mathcal Mh(Z)=\infty\), trivially, and
	when \(\mathcal Mh(Z)=0\), since then \(h=0\) almost everywhere.
	Since \(\mathcal M:L^{p,\infty}\to L^{p,\infty}\),
	\begin{equation}\label{eq:kinetic-young}
		\|I_\alpha h\|_{L^{q,\infty}}
		\le C_{p,\alpha}\|h\|_{L^{p,\infty}},
		\qquad \frac1q=\frac1p-\frac\alpha Q,
		\quad 1<p<\frac Q\alpha.
	\end{equation}
	For \(p=1\), the integral over the complement of \(\mathrm{B}_r(Z)\)
	is at most \(Cr^{\alpha-Q}\|h\|_{L^1}\).  Optimizing the same
	near-field bound and using the weak \((1,1)\) maximal estimate yields
	\[
	\|I_\alpha h\|_{L^{Q/(Q-\alpha),\infty}}
	\le C_\alpha\|h\|_{L^1}.
	\]
	For \(1<p<Q/\alpha\), real interpolation between two
	neighboring weak-type estimates in \eqref{eq:kinetic-young}
	also yields
	\(\|I_\alpha h\|_{L^{q,p}}\le C\|h\|_{L^p}\), hence
	\(\|I_\alpha h\|_{L^q}\le C\|h\|_{L^p}\).
	Because
	\[
	|\mathfrak G[f](Z)|\le CI_2f(Z),\qquad
	|\mathfrak G_i[f](Z)|\le CI_1f(Z),
	\]
	the preceding estimates with \(\alpha=2\) and \(\alpha=1\)
	give \eqref{eq:weak-L1-endpoint}, \eqref{eq:weak-Lp-endpoint},
	and \eqref{eq:strong-Lp-estimate}, with the exponent ranges
	stated in \eqref{eq:QpRp}.
\end{proof}

\begin{lemma}\label{thm:beta-sobolev-subsolution-zero-inflow}
	Let \(Z_0\in\bR^{1+2d}\) and \(R>0\). Assume that
	\(a^{ij},b^i,c\in L^\infty(Q_R^\cD(Z_0))\), \(i,j=1,\dots,d\), satisfy
	\eqref{con_abc} for some constants \(0<\lambda\le\Lambda\) and \(M\ge0\).
	Let \(f\) be a weak solution of
	\[
	\begin{cases}
		(\partial_t+v\cdot\nabla_x)f
		=
		D_i(a^{ij}D_j f)+b^iD_i f+cf+g
		&\quad\text{in }Q_R^\cD(Z_0),\\
		f=h\le0
		&\quad\text{on }(\bR\times\gamma^-)\cap Q_R(Z_0),
	\end{cases}
	\]
	where \(g\in L^2(Q_R^\cD(Z_0))\). Let
	\(\beta:\bR\to[0,\infty)\) be a \(C^2\) function satisfying
	\[
	\beta(t)=0\quad \text{for all } t\le 0,
	\qquad
	\beta''\in C_c^\infty(\bR),
	\qquad
	\beta''\ge 0.
	\]
	Then there exists a constant \(C=C(d)>0\) such that, for every \(0<r<R\),
	\begin{align}
		\begin{split}\label{L2+1/d}
			\|\beta(f)\|_{L^{2+1/d}(Q_r^\cD(Z_0))}
			\le\,&C\bigg(
			r\|g_\beta\|_{L^2(Q_R^\cD(Z_0))}
			+\frac{1}{R-r}\|\beta(f)\|_{L^2(Q_R^\cD(Z_0))}
			\\
			&\qquad\qquad
			+\Bigl(1+\frac{R\Lambda}{R-r}\Bigr)
			\|\beta'(f)\nabla_v f\|_{L^2(Q_R^\cD(Z_0))}
			\bigg),
		\end{split}
	\end{align}
	where $g_\beta:=\bigl(b^iD_i f+cf+g\bigr)\beta'(f)$.
\end{lemma}

The above lemma will mainly be used in the case of zero inflow data.
However, for the broader applications in Subsection~\ref{Subsec:holderconti}, we state it under the slightly more general assumption that the inflow trace is nonpositive.

\begin{proof}[\textbf{Proof of Lemma~\ref{thm:beta-sobolev-subsolution-zero-inflow}}]
	Set
	\[
	q_*:=2+\frac1d,
	\qquad
	q^*:=2+\frac{4}{2d-1}.
	\]
	Since $q_*\leq q^*$ and $|Q_r(Z_0)|=C(d) r^{4d+2}$,
	H\"older's inequality yields
	\begin{equation}\label{eq:local-embed-qstar}
		\|u\|_{L^{q_*}(E)}
		\le C(d)r\|u\|_{L^{q^*}(E)}
		\qquad\text{for every measurable }E\subset Q_r(Z_0).
	\end{equation}
	
	Recall that \(Q_\rho(0)=(-\rho^2,0]\times B_{\rho^3}\times B_\rho\).  Choose
	\(\psi_0\in C_c^\infty(Q_R(0))\) of the product form
	\(\psi_0(t,x,v)=\theta_1(t)\theta_2(x)\theta_3(v)\) such that
	\[
	0\le \psi_0\le 1,
	\qquad
	\psi_0\equiv 1 \quad\text{on }Q_r(0),
	\]
	with \(\theta_1,\theta_2,\theta_3\) making their transitions on the scales
	\(R^2-r^2\), \(R^3-r^3\) and \(R-r\), respectively.  Since
	\(R^2-r^2\ge R(R-r)\), \(R^3-r^3\ge R^2(R-r)\), and \(|v|<R\) on
	\(Q_R(0)\), such a choice satisfies
	\[
	\|\nabla_v\psi_0\|_{\infty}\le \frac{C}{R-r},
	\qquad
	\|(\partial_t+v\cdot\nabla_x)\psi_0\|_{\infty}
	\le \frac{C}{R(R-r)}.
	\]
	Put \(\psi(Z):=\psi_0(Z_0^{-1}\circ Z)\).  Since the vector fields
	\(\p_t+v\cdot\nabla_x\) and \(\nabla_v\) are left-invariant, we obtain
	\begin{align}\label{2604221112}
		\psi\in C_c^\infty(Q_R(Z_0)),
		\qquad
		0\le \psi\le 1,
		\qquad
		\psi\equiv 1 \quad\text{on }Q_r(Z_0),
	\end{align}
	and
	\begin{align}\label{2604221113}
		\|\nabla_v\psi\|_{\infty}
		\le \frac{C}{R-r},
		\qquad
		\|(\partial_t+v\cdot\nabla_x)\psi\|_{\infty}
		\le \frac{C}{R(R-r)}.
	\end{align}
	
	Let
	\[
	\zeta\in C_c^\infty\big((0,\infty)\times\bR_x^d\times\bR_v^d\big),
	\qquad
	\zeta\ge0,
	\qquad
	\iiint_{\bR^{1+2d}}\zeta\dd Z=1,
	\]
	and define
	\[
	\zeta_\varepsilon(t,x,v)
	:=
	\frac{1}{\varepsilon^{2d+1}}
	\zeta\!\left(\frac{t}{\varepsilon},\frac{x}{\varepsilon},\frac{v}{\varepsilon}\right),
	\qquad
	u^{(\varepsilon)}:=u*\zeta_\varepsilon ,
	\]
	the convolution being the Euclidean one on \(\bR^{1+2d}\).
	Set \(U:=\beta(f)\psi\), extended by zero outside \(Q_R^\cD(Z_0)\), and define
	\[
	F_\varepsilon:=U^{(\varepsilon)} .
	\]
	Then \(F_\varepsilon\in C_c^\infty(\bR^{1+2d})\).  Since \(\mathfrak G\) is the
	fundamental solution of \(\p_t+v\cdot\nabla_x-\Delta_v\) and \(F_\varepsilon\)
	is compactly supported,
	\begin{equation}\label{eq:Fe-representation}
		F_\varepsilon(Z)
		=
		\mathfrak G\!\left[
		\bigl(\partial_t+v\cdot\nabla_x-\Delta_v\bigr)F_\varepsilon
		\right](Z).
	\end{equation}
	
	Let us fix \(\varepsilon>0\). For \(Z=(t,x,v)\) and
	\(Y=(s,y,w)\) in \(\bR^{1+2d}\), define
	\[
	\widetilde\zeta_Z(Y):=\zeta_\varepsilon(Z-Y).
	\]
	A direct computation gives
	\begin{equation}\label{eq:Ik-decomposition-sub}
		\bigl(\partial_t+v\cdot\nabla_x-\Delta_v\bigr)F_\varepsilon(Z)
		=
		I_1(Z)+\cdots+I_7(Z),
	\end{equation}
	where
	\begin{align}
		I_1(Z)
		:={}&
		\iiint
		\Bigl(
		\beta(f)\big((-\partial_s-w\cdot\nabla_y)(\psi\widetilde\zeta_Z)\big)
		+
		a^{ij}\partial_{w^j}\beta(f)\,\partial_{w^i}(\psi\widetilde\zeta_Z)
		\Bigr)\dd Y,
		\notag\\
		I_2(Z)
		:={}&
		\iiint
		\beta(f)\,\widetilde\zeta_Z\,(\partial_s+w\cdot\nabla_y)\psi\dd Y,
		\label{eq:I2-def}
		\\
		I_3(Z)
		:={}&
		\iiint
		\partial_{w^j}\beta(f)\,\psi\,\partial_{w^j}\widetilde\zeta_Z\dd Y,
		\label{eq:I3-def}
		\\
		I_4(Z)
		:={}&
		\iiint
		\beta(f)\,\partial_{w^j}\psi\,\partial_{w^j}\widetilde\zeta_Z\dd Y,
		\label{eq:I4-def}
		\\
		I_5(Z)
		:={}&
		-\iiint
		a^{ij}\partial_{w^j}\beta(f)\,\psi\,\partial_{w^i}\widetilde\zeta_Z\dd Y,
		\label{eq:I5-def}
		\\
		I_6(Z)
		:={}&
		-\iiint
		a^{ij}\partial_{w^j}\beta(f)\,\partial_{w^i}\psi\,\widetilde\zeta_Z\dd Y,
		\label{eq:I6-def}
		\\
		I_7(Z)
		:={}&
		\iiint
		U(Y)\,(v-w)\cdot\nabla_x\zeta_\varepsilon(Z-Y)\dd Y.
		\label{eq:I7-def}
	\end{align}
	Here \(I_1\) is obtained by adding and subtracting the terms \(I_5\) and
	\(I_6\).  The integration by parts in \(w\) that produces \(I_3\) and
	\(I_4\) is justified because \(\beta(f)\psi\in L^2(\bR^{1+2d};H^1_w)\).
	Since \(\mathfrak G\ge0\), \eqref{eq:Fe-representation} and
	\eqref{eq:Ik-decomposition-sub} imply
	\begin{equation}\label{eq:Fe-upper-bound}
		0\le F_\varepsilon(Z)\le \sum_{k=1}^7 \mathfrak G[I_k](Z),
	\end{equation}
	the lower bound because \(\beta\ge0\), \(\psi\ge0\) and \(\zeta\ge0\).
	
	\medskip
	\noindent\emph{Estimate of \(I_1\).}
	Let
	\[
	\varphi_Z(Y):=\psi(Y)\widetilde\zeta_Z(Y).
	\]
	Then \(\varphi_Z\in C_c^\infty(Q_R^{\overline{\cD}}(Z_0))\) and
	\(\varphi_Z\ge0\). Moreover, since
	\[
	f|_{(\bR\times\gamma^-)\cap Q_R(Z_0)}\le 0,
	\qquad
	\beta(t)=0\quad\text{for }t\le 0,
	\qquad
	\beta\ge0,
	\]
	we have \(\beta(f_{\gamma^-})=0\), so the hypotheses of the second part of
	Lemma~\ref{260312141} are met.  Applying \eqref{260113330} with the test
	function \(\varphi_Z\) and with \(g^0=g\), \(g^1=\dots=g^d=0\), and
	discarding on the left-hand side the nonnegative term
	\(\beta''(f)a^{ij}D_ifD_jf\,\varphi_Z\), we obtain
	\[
	I_1(Z)
	\le
	\iiint g_\beta(Y)\,\psi(Y)\widetilde\zeta_Z(Y)\dd Y
	=
	\bigl(g_\beta\psi\bigr)^{(\varepsilon)}(Z),
	\]
	where \(g_\beta\psi\) is extended by zero outside \(Q_R^\cD(Z_0)\).
	By the positivity of \(\mathfrak G\), we therefore have
	\[
	\mathfrak G[I_1](Z)
	\le
	\mathfrak G\big[(g_\beta\psi)^{(\varepsilon)}\big](Z).
	\]
	Consequently, Lemma~\ref{lem:global-potential-AR} and
	\eqref{eq:local-embed-qstar} give
	\begin{equation}\label{eq:I1-est}
		\left\|
		\mathfrak G\big[(g_\beta\psi)^{(\varepsilon)}\big]
		\right\|_{L^{q_*}(Q_r(Z_0))}
		\le
		C(d)\,r\,\|g_\beta\|_{L^2(Q_R^\cD(Z_0))}.
	\end{equation}
	
	\medskip
	\noindent\emph{Estimates of \(I_2\) and \(I_4\).}
	Set
	\[
	F_2:=\beta(f)(\partial_t+v\cdot\nabla_x)\psi,
	\]
	extended by zero outside \(Q_R^\cD(Z_0)\). Then \(I_2=(F_2)^{(\varepsilon)}\) by \eqref{eq:I2-def}. Hence Lemma~\ref{lem:global-potential-AR} and \eqref{eq:local-embed-qstar} give
	\begin{align*}
		\begin{split}
			\|\mathfrak G[I_2]\|_{L^{q_*}(Q_r(Z_0))}
			&\le
			C(d)r\|F_2\|_{L^2(Q_R^\cD(Z_0))}
			\\
			&\le
			\frac{C(d)r}{R(R-r)}
			\|\beta(f)\|_{L^2(Q_R^\cD(Z_0))}
			\le
			\frac{C(d)}{R-r}
			\|\beta(f)\|_{L^2(Q_R^\cD(Z_0))}.
		\end{split}
	\end{align*}
	
	Next define
	\[
	H_j:=\beta(f)D_j\psi,
	\qquad j=1,\dots,d,
	\]
	extended by zero outside \(Q_R^\cD(Z_0)\). By \eqref{eq:I4-def},
	\[
	I_4
	=
	-\sum_{j=1}^d D_j(H_j)^{(\varepsilon)}
	\quad\text{in }\cD'(\bR^{1+2d}).
	\]
	Therefore Lemma~\ref{lem:global-potential-AR} gives
	\begin{align*}
		\begin{split}
			\|\mathfrak G[I_4]\|_{L^{q_*}(Q_r(Z_0))}
			&\le
			\sum_{j=1}^d
			\|\mathfrak G_j[(H_j)^{(\varepsilon)}]\|_{L^{q_*}(Q_r(Z_0))}
			\\
			&\le
			C(d)\sum_{j=1}^d
			\|H_j\|_{L^2(Q_R^\cD(Z_0))}
			\le
			\frac{C(d)}{R-r}
			\|\beta(f)\|_{L^2(Q_R^\cD(Z_0))}.
		\end{split}
	\end{align*}
	
	\medskip
	\noindent\emph{Estimate of \(I_3+I_5\).}
	Define
	\[
	G_j
	:=
	\sum_{i=1}^d(\delta_{ij}-a^{ji})D_i\beta(f)\,\psi,
	\qquad j=1,\dots,d,
	\]
	extended by zero outside \(Q_R^\cD(Z_0)\). Combining \eqref{eq:I3-def} and \eqref{eq:I5-def}, we obtain
	\[
	I_3+I_5
	=
	-\sum_{j=1}^d D_j(G_j)^{(\varepsilon)}
	\quad\text{in }\cD'(\bR^{1+2d}).
	\]
	Consequently, Lemma~\ref{lem:global-potential-AR} yields
	\begin{align*}
		\begin{split}
			\|\mathfrak G[I_3+I_5]\|_{L^{q_*}(Q_r(Z_0))}
			&\le
			C(d)\sum_{j=1}^d
			\|G_j\|_{L^2(Q_R^\cD(Z_0))}
			\\
			&\le
			C(d)(1+\Lambda)
			\|\nabla_v\beta(f)\|_{L^2(Q_R^\cD(Z_0))}.
		\end{split}
	\end{align*}
	
	\medskip
	\noindent\emph{Estimate of \(I_6\).}
	Set
	\[
	F_6:=-a^{ij}D_j\beta(f)\,D_i\psi,
	\]
	extended by zero outside \(Q_R^\cD(Z_0)\). Then \(I_6=(F_6)^{(\varepsilon)}\) by \eqref{eq:I6-def}. Hence Lemma~\ref{lem:global-potential-AR} and \eqref{eq:local-embed-qstar} give
	\begin{align*}
		\begin{split}
			\|\mathfrak G[I_6]\|_{L^{q_*}(Q_r(Z_0))}
			&\le
			C(d)\,r\,\|F_6\|_{L^2(Q_R^\cD(Z_0))}
			\\
			&\le
			C(d)\frac{r\Lambda}{R-r}
			\|\nabla_v\beta(f)\|_{L^2(Q_R^\cD(Z_0))}.
		\end{split}
	\end{align*}
	
	\medskip
	\noindent\emph{Estimate of \(I_7\).}
	For \((s,y,w)\in\bR^{1+2d}\), set
	\[
	K(s,y,w):=w\cdot\nabla_y\zeta(s,y,w),
	\]
	and define
	\[
	K_\varepsilon(s,y,w)
	:=
	w\cdot\nabla_y\zeta_\varepsilon(s,y,w)=
	\frac{1}{\varepsilon^{2d+1}}
	K\!\left(\frac{s}{\varepsilon},\frac{y}{\varepsilon},\frac{w}{\varepsilon}\right).
	\]
	Then \(I_7=U*K_\varepsilon\) by \eqref{eq:I7-def}, the convolution being the
	Euclidean one on \(\bR^{1+2d}\).  Moreover, \(K\in C_c^\infty(\bR^{1+2d})\)
	and
	\[
	\iiint_{\bR^{1+2d}} K(s,y,w)\dd s\dd y\dd w
	=
	\iint\left(\int_{\bR^d_y}\nabla_y\zeta\dd y\right)\cdot w\dd s\dd w
	=0.
	\]
	Since \(U\in L^2(\bR^{1+2d})\) and \(\{K_\varepsilon\}\) is the
	\(L^1\)-dilation family generated by \(K\), the standard approximate-identity
	argument gives \(U*K_\varepsilon\to \left(\iiint K\right)U=0\) in
	\(L^2(\bR^{1+2d})\), that is,
	\[
	\|I_7\|_{L^2(\bR^{1+2d})}
	=
	\|U*K_\varepsilon\|_{L^2(\bR^{1+2d})}
	\to 0
	\qquad\text{as }\varepsilon\to0.
	\]
	Consequently, Lemma~\ref{lem:global-potential-AR} and \eqref{eq:local-embed-qstar} yield
	\begin{equation}\label{eq:I7-est}
		\|\mathfrak G[I_7]\|_{L^{q_*}(Q_r(Z_0))}
		\le
		C(d)\,r\,\|I_7\|_{L^2(\bR^{1+2d})}
		\to 0\qquad \text{as }\varepsilon\to0.
	\end{equation}
	
	\medskip
	\noindent\textit{Conclusion.}
	By the standard properties of mollification,
	\(F_\varepsilon=(\beta(f)\psi)^{(\varepsilon)}\to \beta(f)\psi\) almost everywhere.
	Fatou's lemma, together with \(\psi\equiv1\) on \(Q_r(Z_0)\), yields
	\[
	\|\beta(f)\|_{L^{q_*}(Q_r^\cD(Z_0))}
	\le
	\liminf_{\varepsilon\to0}
	\|F_\varepsilon\|_{L^{q_*}(Q_r(Z_0))}.
	\]
	Therefore, combining \eqref{eq:Fe-upper-bound} with
	\eqref{eq:I1-est}--\eqref{eq:I7-est} and letting \(\varepsilon\to0\), we obtain
	\begin{align*}
		\begin{split}
			\|\beta(f)\|_{L^{q_*}(Q_r^\cD(Z_0))}
			\le\,
			C(d)\Bigg(
			&r\|g_\beta\|_{L^2(Q_R^\cD(Z_0))}
			+\frac{1}{R-r}\|\beta(f)\|_{L^2(Q_R^\cD(Z_0))}
			\\
			&\qquad
			+\Bigl(1+\frac{R\Lambda}{R-r}\Bigr)
			\|\nabla_v\beta(f)\|_{L^2(Q_R^\cD(Z_0))}
			\Bigg).
		\end{split}
	\end{align*}
	Since \(q_*=2+\frac1d\) and
	\(\nabla_v\beta(f)=\beta'(f)\nabla_v f\), this proves \eqref{L2+1/d}.
\end{proof}

\begin{lemma}\label{2601141041}
	Suppose that \(a^{ij}\), \(b^i\), \(c\), \(f\), and \(g\) satisfy the hypotheses of
	Lemma~\ref{thm:beta-sobolev-subsolution-zero-inflow}.
	Then, for every
	\(p\ge 1\) and every \(0<r<R\),
	\begin{align}
		\begin{split}\label{eq:lp-sobolev-cor}
			\left\|f_+^p\right\|_{L^{2+1/d}(Q_r^\cD(Z_0))}
			\le C(d)p\Bigg(
			r\left\||g|f_+^{p-1}\right\|_{L^2(Q_R^\cD(Z_0))}+\left(\frac{1}{R-r}+M^2r\right)
			\left\|f_+^p\right\|_{L^2(Q_R^\cD(Z_0))}\quad &
			\\
			+\left(1+\frac{R\Lambda}{R-r}+Mr\right)
			\left\|f_+^{p-1}\nabla_v f\right\|_{L^2(Q_R^\cD(Z_0))}
			&\Bigg).
		\end{split}
	\end{align}
\end{lemma}

\begin{proof}
	Let \(\{\beta_n\}_{n\ge1}\subset C^\infty(\bR)\) be a sequence of nonnegative
	functions such that
	\[
	\beta_n(t)=0\quad\text{for }t\le 0,
	\qquad
	\beta_n''\in C_c^\infty(\bR),
	\qquad
	\beta_n''\ge 0,
	\]
	and, for every \(t\in\bR\),
	\[
	0\le \beta_n(t)\uparrow t_+^p,
	\qquad
	0\le \beta_n'(t)\le p\,t_+^{p-1}.
	\]
	
	Applying Lemma~\ref{thm:beta-sobolev-subsolution-zero-inflow} with
	\(\beta=\beta_n\), we obtain
	\begin{align}
		\begin{split}\label{eq:beta-n-est}
			\|\beta_n(f)\|_{L^{2+1/d}(Q_r^\cD(Z_0))}
			\le C(d)\Bigg(
			r\|g_{\beta_n}\|_{L^2(Q_R^\cD(Z_0))}
			+\frac{1}{R-r}\|\beta_n(f)\|_{L^2(Q_R^\cD(Z_0))}\quad &
			\\
			+\Bigl(1+\frac{R\Lambda}{R-r}\Bigr)
			\|\beta_n'(f)\nabla_v f\|_{L^2(Q_R^\cD(Z_0))}
			\,&\Bigg),
		\end{split}
	\end{align}
	where
	\[
	g_{\beta_n}
	:=
	\beta_n'(f)\bigl(b^iD_i f+cf+g\bigr).
	\]
	
	We estimate the terms on the right-hand side of \eqref{eq:beta-n-est}.
	Since \(0\le \beta_n(t)\le t_+^p\), we have
	\[
	\|\beta_n(f)\|_{L^2(Q_R^\cD(Z_0))}
	\le
	\|f_+^p\|_{L^2(Q_R^\cD(Z_0))}.
	\]
	Also, since \(0\le \beta_n'(t)\le p\,t_+^{p-1}\), we have
	\[
	\|\beta_n'(f)\nabla_v f\|_{L^2(Q_R^\cD(Z_0))}
	\le
	p\|f_+^{p-1}\nabla_v f\|_{L^2(Q_R^\cD(Z_0))}.
	\]
	Moreover, using \eqref{con_abc},
	\begin{align*}
		\|g_{\beta_n}\|_{L^2(Q_R^\cD(Z_0))}
		&\le
		pM\|f_+^{p-1}\nabla_v f\|_{L^2(Q_R^\cD(Z_0))}
		\\
		&\quad
		+pM^2\|f_+^p\|_{L^2(Q_R^\cD(Z_0))}
		+p\left\||g|f_+^{p-1}\right\|_{L^2(Q_R^\cD(Z_0))}.
	\end{align*}
	Substituting these bounds into \eqref{eq:beta-n-est}, we get
	\begin{align}
		\begin{split}\label{eq:beta-n-est-final}
			\|\beta_n(f)\|_{L^{2+1/d}(Q_r^\cD(Z_0))}
			\le C(d)p\Bigg(
			r\left\||g|f_+^{p-1}\right\|_{L^2(Q_R^\cD(Z_0))}
			+\left(\frac{1}{R-r}+M^2r\right)
			\left\|f_+^p\right\|_{L^2(Q_R^\cD(Z_0))}\quad &
			\\
			+\left(1+\frac{R\Lambda}{R-r}+Mr\right)
			\left\|f_+^{p-1}\nabla_v f\right\|_{L^2(Q_R^\cD(Z_0))}
			&\Bigg),
		\end{split}
	\end{align}
	where the constant is independent of \(n\).
	
	Finally, since \(\beta_n(f)\uparrow f_+^p\) almost everywhere, the monotone
	convergence theorem gives
	\[
	\|\beta_n(f)\|_{L^{2+1/d}(Q_r^\cD(Z_0))}
	\longrightarrow
	\|f_+^p\|_{L^{2+1/d}(Q_r^\cD(Z_0))}.
	\]
	Passing to the limit in \eqref{eq:beta-n-est-final}, we obtain
	\eqref{eq:lp-sobolev-cor}.
\end{proof}

\begin{lemma}\label{thm:caccioppoli}
	Suppose that \(a^{ij}\), \(b^i\), \(c\), \(f\), and \(g\) satisfy the hypotheses of
	Lemma~\ref{thm:beta-sobolev-subsolution-zero-inflow}.
	Then, for every
	\(p\ge1\) and every \(0<r<R\),
	\begin{align}
		\begin{split}\label{est:caccioppoli}
			\left\|f_+^{p-1}\nabla_v f\right\|_{L^2(Q_r^\cD(Z_0))}^2
			\le
			\frac{C(d)}{\lambda p}
			\Bigg(\,&
			\iiint_{Q_R^\cD(Z_0)} f_+^{2p-1}|g|\dd Z
			\\
			&\quad
			+\frac{(\Lambda+1)^2}{\lambda}
			\left(M^2+\frac{1}{(R-r)^2}\right)
			\left\|f_+^p\right\|_{L^2(Q_R^\cD(Z_0))}^2
			\Bigg).
		\end{split}
	\end{align}
\end{lemma}

\begin{proof}
	Let \(\psi\) be as in \eqref{2604221112}--\eqref{2604221113}; its restriction
	to \(\bR\times\overline\cD\times\bR^d_v\) belongs to
	\(C_c^\infty(Q_R^{\overline{\cD}}(Z_0))\).
	
	Let \(\{\alpha_n\}_{n\ge1}\subset C^\infty(\bR)\) be a sequence of
	nonnegative functions such that
	\[
	\alpha_n(t)=0\quad\text{for }t\le 0,
	\qquad
	\alpha_n''\in C_c^\infty(\bR),\qquad \alpha_n''\ge0,
	\]
	and, for every \(t\in\bR\),
	\[
	0\le \alpha_n(t)\uparrow t_+^p,
	\qquad
	0\le \alpha_n'(t)\uparrow p\,t^{p-1}\mathbf 1_{t>0}.
	\]
	Set
	\[
	\beta_n(t):=\alpha_n(t)^2.
	\]
	Then
	\[
	\beta_n(t)=0\quad\text{for }t\le0,
	\qquad
	\beta_n'(t)=2\alpha_n(t)\alpha_n'(t),
	\]
	and
	\[
	\beta_n''(t)
	=
	2(\alpha_n'(t))^2+2\alpha_n(t)\alpha_n''(t)
	\ge 2(\alpha_n'(t))^2,
	\]
	where \(\beta_n''\) is bounded because \(\alpha_n'\) is bounded and
	\(\alpha_n''\) has compact support.
	Since \(\beta_n'=2\alpha_n\alpha_n'\) is unbounded, we truncate.  For \(k>0\), choose \(\chi_k\in C^\infty(\bR)\) with
	\[
	0\le\chi_k\le1,
	\qquad
	\chi_k\equiv1\ \text{on }(-\infty,k],
	\qquad
	\chi_k\equiv0\ \text{on }[k+1,\infty),
	\]
	and let \(\beta_{n,k}\in C^2(\bR)\) be defined by
	\[
	\beta_{n,k}''=\chi_k\beta_n'' ,
	\qquad
	\beta_{n,k}(0)=\beta_{n,k}'(0)=0 .
	\]
	Then \(\beta_{n,k}\) satisfies all the hypotheses of Lemma~\ref{260312141}:
	indeed \(\beta_{n,k}''\) is bounded and compactly supported in
	\([0,k+1]\), so \(\beta_{n,k}'\) is bounded, and \(\beta_{n,k}=0\) on
	\((-\infty,0]\).  Moreover \(\beta_{n,k}\ge0\), and, since
	\(0\le\beta_{n,k}''\le\beta_n''\),
	\begin{align}\label{2608052120}
		0\le\beta_{n,k}\le\beta_n,
		\qquad
		0\le\beta_{n,k}'\le\beta_n',
		\qquad
		\beta_{n,k}''\ge 2(\alpha_n')^2\mathbf 1_{(-\infty,k]} .
	\end{align}
	
	We apply Lemma~\ref{260312141} with \(\beta=\beta_{n,k}\) and
	\(\varphi=\psi^2\). Since \(f\le0\) on the inflow boundary and
	\(\beta_{n,k}=0\) on \((-\infty,0]\), \eqref{260113330} gives
	\begin{align*}
		&\iiint_{Q_R^\cD(Z_0)}
		\beta_{n,k}''(f)\,a^{ij}D_i fD_j f\,\psi^2\dd Z
		\\
		\le\;&
		\iiint_{Q_R^\cD(Z_0)}
		\beta_{n,k}'(f)(cf+g)\psi^2\dd Z
		+
		\iiint_{Q_R^\cD(Z_0)}
		\beta_{n,k}(f)\,(\partial_t+v\cdot\nabla_x)(\psi^2)\dd Z
		\\
		&\quad
		-\iiint_{Q_R^\cD(Z_0)}
		\beta_{n,k}'(f)\,a^{ij}D_j fD_i(\psi^2)\dd Z
		+\iiint_{Q_R^\cD(Z_0)}
		\beta_{n,k}'(f)\,b^iD_i f\,\psi^2\dd Z .
	\end{align*}
	We bound the left-hand side from below using the last relation in
	\eqref{2608052120} and \(a^{ij}\xi_i\xi_j\ge\lambda|\xi|^2\), and the
	right-hand side from above using
	\(0\le\beta_{n,k}\le\alpha_n^2\) and
	\(0\le\beta_{n,k}'\le2\alpha_n\alpha_n'\).  Since the resulting right-hand side
	does not depend on \(k\), we may let \(k\to\infty\) and use monotone
	convergence on the left-hand side to obtain
	\begin{align*}
		&2\lambda
		\|\alpha_n'(f)\nabla_v f\,\psi\|_{L^2(Q_R^\cD(Z_0))}^2\\
		\le\;&
		2\iiint_{Q_R^\cD(Z_0)}
		|\alpha_n(f)\alpha_n'(f)|\,|cf+g|\,\psi^2\dd Z
		+\iiint_{Q_R^\cD(Z_0)}
		|\alpha_n(f)|^2\,|(\partial_t+v\cdot\nabla_x)(\psi^2)|\dd Z
		\\
		&\quad
		+C(d)\Lambda\iiint_{Q_R^\cD(Z_0)}
		|\alpha_n(f)\alpha_n'(f)|\,|\nabla_v f|\,\psi\,|\nabla_v\psi|\dd Z
		+C(d)\iiint_{Q_R^\cD(Z_0)}
		|\alpha_n(f)\alpha_n'(f)|\,|b|\,|\nabla_v f|\,\psi^2\dd Z .
	\end{align*}
	
	Set
	\[
	A_n:=\|\alpha_n'(f)\nabla_v f\,\psi\|_{L^2(Q_R^\cD(Z_0))},
	\qquad
	B:=\|f_+^p\|_{L^2(Q_R^\cD(Z_0))}.
	\]
	Since
	\[
	0\le \alpha_n(f)\le f_+^p,
	\qquad
	0\le \alpha_n'(f)\le p f_+^{p-1},
	\qquad
	|b|\le M,
	\qquad
	|c|\le M^2,
	\]
	we estimate the four terms on the right-hand side as follows:
	\begin{align*}
		2\iiint_{Q_R^\cD(Z_0)}
		|\alpha_n(f)\alpha_n'(f)|\,|g|\,\psi^2\dd Z
		&\le
		2p\iiint_{Q_R^\cD(Z_0)} f_+^{2p-1}|g|\dd Z,
		\\
		2\iiint_{Q_R^\cD(Z_0)}
		|\alpha_n(f)\alpha_n'(f)|\,|c|\,|f|\,\psi^2\dd Z
		&\le
		2pM^2\iiint_{Q_R^\cD(Z_0)} f_+^{2p}\dd Z
		=
		2pM^2B^2,
		\\
		\iiint_{Q_R^\cD(Z_0)}
		|\alpha_n(f)|^2\,|(\partial_t+v\cdot\nabla_x)(\psi^2)|\dd Z
		&\le
		\frac{C(d)}{R(R-r)}B^2
		\le
		\frac{C(d)}{(R-r)^2}B^2,
		\\
		C(d)\Lambda\iiint_{Q_R^\cD(Z_0)}
		|\alpha_n(f)\alpha_n'(f)|\,|\nabla_v f|\,\psi\,|\nabla_v\psi|\dd Z
		&\le
		\frac{C(d)\Lambda}{R-r}A_nB\le
		\frac{\lambda}{4}A_n^2
		+
		\frac{C(d)\Lambda^2}{\lambda(R-r)^2}B^2,
		\\
		C(d)\iiint_{Q_R^\cD(Z_0)}
		|\alpha_n(f)\alpha_n'(f)|\,|b|\,|\nabla_v f|\,\psi^2\dd Z
		&\le
		C(d)M A_nB\le
		\frac{\lambda}{4}A_n^2
		+
		\frac{C(d)}{\lambda}M^2B^2.
	\end{align*}
	Absorbing the two \(\lambda A_n^2/4\) terms into the left-hand side, we obtain
	\begin{align*}
		A_n^2
		\le\;&
		C(d)\lambda^{-1}
		\bigg(
		p\iiint_{Q_R^\cD(Z_0)} f_+^{2p-1}|g|\dd Z
		+pM^2B^2
		+\frac{1}{(R-r)^2}B^2
		+\frac{\Lambda^2}{\lambda(R-r)^2}B^2
		+\frac{M^2}{\lambda}B^2
		\bigg)
		\\
		\le\;&
		C(d)\lambda^{-1}p
		\bigg(
		\iiint_{Q_R^\cD(Z_0)} f_+^{2p-1}|g|\dd Z
		+
		\frac{(\Lambda+1)^2}{\lambda}
		\left(M^2+\frac{1}{(R-r)^2}\right)B^2
		\bigg),
	\end{align*}
	where we used \(p\ge1\) in the last inequality.
	
	Letting \(n\to\infty\), we have
	$\alpha_n'(f)\to p f_+^{p-1}$ almost everywhere.
	Since \(\psi\equiv1\) on \(Q_r(Z_0)\), Fatou's lemma yields
	\begin{align*}
		&p^2\left\|f_+^{p-1}\nabla_v f\right\|_{L^2(Q_r^\cD(Z_0))}^2
		\\
		\le\;&
		\liminf_{n\to\infty} A_n^2
		\\
		\le\;&
		C(d)\lambda^{-1}p
		\bigg(
		\iiint_{Q_R^\cD(Z_0)} f_+^{2p-1}|g|\dd Z
		+
		\frac{(\Lambda+1)^2}{\lambda}
		\left(M^2+\frac{1}{(R-r)^2}\right)
		\left\|f_+^p\right\|_{L^2(Q_R^\cD(Z_0))}^2
		\bigg).
	\end{align*}
	Dividing by \(p^2\), this proves \eqref{est:caccioppoli}.
\end{proof}

\begin{prop}\label{cor:boundary-moser-gain}
	Suppose that \(a^{ij}\), \(b^i\), \(c\), \(f\), and \(g\) satisfy the hypotheses of
	Lemma~\ref{thm:beta-sobolev-subsolution-zero-inflow}.
	Let \(p\ge1\), and assume in addition that
	$g\in L^{2p}(Q_R^\cD(Z_0))$ and $f_+^p\in L^2(Q_R^\cD(Z_0))$.
	Then, for every \(0<r<R\),
	\begin{align}
		\begin{split}\label{260428431}
			\|f_+^p\|_{L^{2+1/d}(Q_r^\cD(Z_0))}
			\le C(d)p\Bigg(
			&R\left\||g|f_+^{p-1}\right\|_{L^2(Q_R^\cD(Z_0))}
			\\
			&\quad +\frac{(1+\Lambda)^2}{\lambda}
			\left(M^2R+\frac{R}{(R-r)^2}\right)
			\left\|f_+^p\right\|_{L^2(Q_R^\cD(Z_0))}
			\Bigg).
		\end{split}
	\end{align}
	Moreover,
	\begin{equation}\label{eq:L2p-source-product-bound}
		\left\||g|f_+^{p-1}\right\|_{L^2(Q_R^\cD(Z_0))}
		\le
		\|g\|_{L^{2p}(Q_R^\cD(Z_0))}
		\|f_+^p\|_{L^2(Q_R^\cD(Z_0))}^{1-1/p},
	\end{equation}
	where the last factor is interpreted as $1$ when $p=1$.
\end{prop}

\begin{proof}
	Estimate \eqref{eq:L2p-source-product-bound} follows directly from
	H\"older's inequality.  In particular, the source terms in
	Lemmas~\ref{2601141041} and \ref{thm:caccioppoli} are finite.
	Put $	\widetilde R:=\frac{R+r}{2}$ and set
	\[
	\cX:=\left\||g|f_+^{p-1}\right\|_{L^2(Q_R^\cD(Z_0))},
	\qquad
	\cY:=\left\|f_+^p\right\|_{L^2(Q_R^\cD(Z_0))},\qquad 	\cZ:=\left\|f_+^{p-1}\nabla_v f\right\|_{L^2(Q_{\widetilde R}^\cD(Z_0))}.
	\]
	Applying Lemma~\ref{2601141041} with inner radius \(r\) and outer radius
	\(\widetilde R\), and using
	\[
	\widetilde R-r=\frac{R-r}{2},
	\qquad
	\widetilde R\le R,
	\qquad
	Q_{\widetilde R}^\cD(Z_0)\subset Q_R^\cD(Z_0),
	\]
	we obtain
	\[
	\|f_+^p\|_{L^{2+1/d}(Q_r^\cD(Z_0))}
	\le
	C(d)p\left(R\cX+A\cY+B\cZ\right),
	\]
	where
	\[
	A:=M^2R+\frac{1}{R-r},
	\qquad
	B:=(1+\Lambda)\left(MR+\frac{R}{R-r}\right).
	\]
	Next, applying Lemma~\ref{thm:caccioppoli} with inner radius
	\(\widetilde R\) and outer radius \(R\), we get
	\begin{align*}
		\cZ^2
		&\le
		\frac{C(d)}{\lambda p}
		\bigg(
		\iiint_{Q_R^\cD(Z_0)} f_+^{2p-1}|g|\dd Z
		+
		\frac{(1+\Lambda)^2}{\lambda}
		\left(M^2+\frac{1}{(R-\widetilde R)^2}\right)\cY^2
		\bigg)
		\\
		&\le
		\frac{C(d)}{\lambda p}
		\bigg(
		\cX\cY
		+
		\frac{(1+\Lambda)^2}{\lambda}
		\left(M^2+\frac{1}{(R-r)^2}\right)\cY^2
		\bigg).
	\end{align*}
	Here we used H\"older's inequality:
	\[
	\iiint_{Q_R^\cD(Z_0)} f_+^{2p-1}|g|\dd Z
	\le
	\left\||g|f_+^{p-1}\right\|_{L^2(Q_R^\cD(Z_0))}
	\left\|f_+^p\right\|_{L^2(Q_R^\cD(Z_0))}
	=
	\cX\cY.
	\]
	Since \(p\ge1\), it follows that
	\[
	\cZ
	\le
	C(d)\lambda^{-1/2}\sqrt{\cX\cY}
	+
	C(d)\frac{1+\Lambda}{\lambda}
	\left(M+\frac{1}{R-r}\right)\cY.
	\]
	Therefore,
	\begin{align*}
		\|f_+^p\|_{L^{2+1/d}(Q_r^\cD(Z_0))}
		&\le
		C(d)p
		\left(
		R\cX
		+\lambda^{-1/2}B\sqrt{\cX\cY}
		+
		\left(A+\frac{B^2}{\lambda R}\right)\cY
		\right)
		\\
		&\le
		C(d)p
		\left(
		R\cX
		+
		\left(A+\frac{B^2}{\lambda R}\right)\cY
		\right),
	\end{align*}
	where Young's inequality was used in the last step.
	
	It remains to estimate the coefficient of \(\cY\). Since
	\[
	A=M^2R+\frac{1}{R-r}
	\le
	M^2R+\frac{R}{(R-r)^2},
	\]
	and
	\begin{align*}
		\frac{B^2}{\lambda R}
		&=
		\frac{(1+\Lambda)^2}{\lambda R}
		\left(MR+\frac{R}{R-r}\right)^2\le
		\frac{C(1+\Lambda)^2}{\lambda}
		\left(M^2R+\frac{R}{(R-r)^2}\right),
	\end{align*}
	we obtain
	\[
	A+\frac{B^2}{\lambda R}
	\le
	\frac{C(1+\Lambda)^2}{\lambda}
	\left(M^2R+\frac{R}{(R-r)^2}\right).
	\]
	Substituting this into the previous estimate gives \eqref{260428431}.
\end{proof}

\begin{proof}[\textbf{Proof of Theorem~\ref{cor:linfty-zero-inflow}}]
	Set
	\[
	q:=1+\frac{1}{2d},
	\qquad
	p_n:=q^n=\left(1+\frac{1}{2d}\right)^n,
	\]
	and define
	\[
	R_n:=\left(\kappa+\frac{1-\kappa}{2^n}\right)R.
	\]
	Then
	\[
	R_0=R,
	\qquad
	R_n\searrow \kappa R,
	\qquad
	R_n-R_{n+1}=\frac{1-\kappa}{2^{n+1}}R.
	\]
	
	Let
	\[
	a_n:=\|f_+\|_{L^{2p_n}(Q_{R_n}^\cD(Z_0))},
	\qquad
	\cG:=\|g\|_{L^\infty(Q_R^\cD(Z_0))}.
	\]
	We first show that \(a_n<\infty\) for every \(n\).  Indeed, \(a_0<\infty\)
	because \(f\in L^2(Q_R^\cD(Z_0))\), and \(g\in L^{2p}(Q_R^\cD(Z_0))\) for
	every \(p\ge1\) because \(g\in L^\infty\) and \(Q_R^\cD(Z_0)\) has finite
	measure.  Hence, if \(a_n<\infty\), the qualitative hypotheses of
	Proposition~\ref{cor:boundary-moser-gain} are satisfied with \(p=p_n\) on
	\(Q_{R_n}^\cD(Z_0)\), and the right-hand side of \eqref{260428431}, applied
	with \((p,R,r)=(p_n,R_n,R_{n+1})\), is finite.  Therefore \(a_{n+1}<\infty\).
	
	By H\"older's inequality, we have
	\begin{align*}
		\left\||g| f_+^{p_n-1}\right\|_{L^2(Q_{R_n}^\cD(Z_0))}
		&\le
		\left\||g|^{p_n}\right\|_{L^2(Q_{R_n}^\cD(Z_0))}^{1/p_n}
		\left\|f_+^{p_n}\right\|_{L^2(Q_{R_n}^\cD(Z_0))}^{1-1/p_n}
		\\
		&\le
		C(d)R^{(1+2d)/p_n}a_n^{p_n-1}\cG.
	\end{align*}
	
	\medskip
	\noindent
	Applying \eqref{260428431} with \((p,R,r)=(p_n,R_n,R_{n+1})\), and using
	$2p_{n+1}=\left(2+\frac1d\right)p_n$,
	we obtain
	\begin{align}
		a_{n+1}^{p_n}
		&\le
		C_0p_n\Big(
		R^{1+(1+2d)/p_n}a_n^{p_n-1}\cG
		+\bigl(M^2R+2^{2n}R^{-1}\bigr)a_n^{p_n}
		\Big)
		\notag
		\\
		&\le
		C_0(MR+1)^2R^{-1}p_n2^{2n}a_n^{p_n-1}
		\Big(
		a_n+R^{2+(1+2d)/p_n}\cG
		\Big),
		\label{eq:linfty-step1}
	\end{align}
	where \(C_0=C_0(d,\lambda,\Lambda,\kappa)\). Here we used
	\[
	\frac{R_n}{(R_n-R_{n+1})^2}
	\le
	C(\kappa)\,2^{2n}R^{-1}.
	\]
	
	Define
	\[
	K_n:=
	\Bigl(
	C_0\,p_n\,2^{2n}(MR+1)^2R^{-1}
	\Bigr)^{1/p_n}.
	\]
	Then \eqref{eq:linfty-step1} gives
	\[
	a_{n+1}
	\le
	K_n\,a_n^{1-1/p_n}
	\Bigl(a_n+R^{2+(1+2d)/p_n}\cG\Bigr)^{1/p_n}.
	\]

	Now set
	\[
	b_n:=a_n+R^{2+(1+2d)/p_n}\cG .
	\]
	Since
	$p_{n+1}=\left(1+\frac1{2d}\right)p_n$,
	we have
	\[
	2+\frac{1+2d}{p_{n+1}}
	=
	2+\frac{1+2d}{p_n}-\frac1{p_n}.
	\]
	Moreover, since \(K_n\) contains a factor \(R^{-1/p_n}\), by enlarging
	\(C_0\) if necessary, we may assume that
	\[
	R^{2+(1+2d)/p_{n+1}}\cG
	\le
	K_n R^{2+(1+2d)/p_n}\cG .
	\]
	Therefore,
	\[
	b_{n+1}
	\le
	K_n\left(
	a_n^{1-1/p_n}b_n^{1/p_n}
	+
	R^{2+(1+2d)/p_n}\cG
	\right).
	\]
	Put
	\[
	A:=a_n,
	\qquad
	B:=R^{2+(1+2d)/p_n}\cG,
	\qquad
	\varepsilon:=\frac1{p_n}.
	\]
	Then
	\begin{align}\label{e-inequality}
		\begin{aligned}
			A^{1-\varepsilon}(A+B)^\varepsilon+B
			&\le
			A^{1-\varepsilon}(A+B)^\varepsilon
			+
			B^{1-\varepsilon}(A+B)^\varepsilon  \\
			&=
			(A+B)^\varepsilon
			\left(A^{1-\varepsilon}+B^{1-\varepsilon}\right) \\
			&\le
			2^\varepsilon(A+B),
		\end{aligned}
	\end{align}
	where the last inequality follows from the concavity of
	\(s\mapsto s^{1-\varepsilon}\):
	\[
	\frac{A^{1-\varepsilon}+B^{1-\varepsilon}}{2}
	\le
	\left(\frac{A+B}{2}\right)^{1-\varepsilon}.
	\]
	Consequently,
	$b_{n+1}\le K_n2^{1/p_n}b_n$.
	Iterating this inequality gives
	\[
	b_n
	\le
	b_0\prod_{j=0}^{n-1}\bigl(K_j2^{1/p_j}\bigr).
	\]

	Since \(Q_{\kappa R}^\cD(Z_0)\subset Q_{R_n}^\cD(Z_0)\) for every \(n\),
	\[
	\|f_+\|_{L^{2p_n}(Q_{\kappa R}^\cD(Z_0))}\le a_n.
	\]
	Because \(Q_{\kappa R}^\cD(Z_0)\) has finite measure and \(p_n\to\infty\),
	the left-hand side converges to \(\|f_+\|_{L^\infty(Q_{\kappa R}^\cD(Z_0))}\),
	so that
	\[
	\|f_+\|_{L^\infty(Q_{\kappa R}^\cD(Z_0))}
	\le
	\liminf_{n\to\infty}a_n
	\le
	\liminf_{n\to\infty}b_n.
	\]
	Moreover,
	\[
	b_0
	=
	a_0+R^{2+(1+2d)}\cG
	=
	\|f_+\|_{L^2(Q_R^\cD(Z_0))}+R^{3+2d}\cG.
	\]
	
	Finally,
	\[
	K_n2^{1/p_n}
	=
	\Bigl(
	C_0\,p_n\,2^{2n}(MR+1)^2R^{-1}
	\Bigr)^{1/p_n}2^{1/p_n},
	\]
	and
	\[
	\sum_{n=0}^\infty \frac1{p_n}=2d+1,
	\qquad
	\sum_{n=0}^\infty \frac{\log p_n}{p_n}<\infty,
	\qquad
	\sum_{n=0}^\infty \frac{n}{p_n}<\infty.
	\]
	Hence
	\[
	\prod_{n=0}^\infty \bigl(K_n2^{1/p_n}\bigr)
	\le
	C(d,\lambda,\Lambda,\kappa)\,(MR+1)^{4d+2}R^{-1-2d}.
	\]
	Combining the above estimates, we obtain \eqref{260628242}.
	
	If \(h\equiv0\) on \((\bR\times \gamma^-)\cap Q_R(Z_0)\), then the same argument applies to both \(f\) and \(-f\), with source terms \(g\) and \(-g\), respectively. Hence \eqref{260628241} follows from \eqref{260628242}.
\end{proof}

\subsection{Half-space geometry and extension lemmas}
\label{Subsec:3preparing}

The remaining boundary estimates are formulated in \(\bR^d_+\).
We collect here the geometry and extension tools used in Subsections~\ref{Subsec:holderconti}--\ref{Subsec:Dvest}. 
For
\(Z_0=(t_0,x_0,v_0)\in\bR\times \bR^d_x\times\bR^d_v\),
we write
\[
Q_r^\pm(Z_0):=Q_r(Z_0)\cap\{\pm x^d>0\},
\]
and we use, for the half-space, the notation
\begin{align}\label{2608051140}
	\gamma:=\{(x,v):x^d=0\},\quad
	\gamma^\pm:=\{(x,v):x^d=0,\ \pm v^d<0\},\quad
	\gamma^0:=\{(x,v):x^d=0,\ v^d=0\},
\end{align}
together with \(\dd\gamma:=|v^d|\dd S_x\dd v\).  This is consistent with
\eqref{def:gamma} and \eqref{def:dgamma}, since the outward unit normal of
\(\bR^d_+\) along \(\{x^d=0\}\) is \(-e_d\).

We first realize half-space cylinders inside bounded \(C^{1,1}\)-domains and
record the kinetic scaling, which makes the preceding estimates available in
\(\bR^d_+\). We then compare
\(\mathbf d(\cdot,\bR\times\gamma^0)\) with
$(x^d)^{1/3}+|v^d|$,
to localize near and away from grazing. The extension lemmas then reduce
non-grazing boundary estimates to interior ones, while
Corollary~\ref{cor:scaled-smallness-implies-half} provides the smallness
criterion used in Appendix~\ref{proofofholder}.

\subsubsection{Half-space cylinders, localization, and scaling}

\begin{lemma}
	\label{2606301252}
	Let \(R>0\) and
	$Z_0=(t_0,x_0,v_0)
	\in
	\bR\times\bR^d_x\times\bR^d_v$.
	Then there exists a bounded \(C^{1,1}\)-domain \(U\subset\bR^d_x\) such that
	\[
	Q_R^+(Z_0)=Q_R^U(Z_0).
	\]
	Moreover, for every \((t,x,v)\in Q_R(Z_0)\), denoting $x=(x',x^d)\in \bR^{d-1}\times \bR$,
	\[
	x^d=0
	\quad\Longleftrightarrow\quad
	x\in\partial U,
	\]
	and, in this case,
	$n_U(x)=-e_d$.
\end{lemma}	
\begin{proof}
	Set
	\[
	R_0:=R^3+|x_0|+R^2|v_0|.
	\]
	If \((t,x,v)\in Q_R(Z_0)\), then
	\[
	x=x_0+(t-t_0)v_0+y,
	\qquad
	|t-t_0|<R^2,
	\qquad
	|y|<R^3,
	\]
	and hence \(|x|<R_0\).  Thus the \(x\)-projection of \(Q_R(Z_0)\) is contained
	in \(B_{R_0}\).
	
	Choose \(L>R_0+2\), set
	\[
	K:=[-L,L]^{d-1}\times[1,L],
	\qquad
	U:=\{x\in\bR^d:\operatorname{dist}(x,K)<1\}.
	\]
	Then \(U\) is bounded and of class \(C^{1,1}\), being the open unit
	parallel neighborhood of the compact convex set \(K\).  Indeed, for a convex
	body the distance function is \(C^{1,1}\) with unit gradient on
	\(\{\operatorname{dist}(\cdot,K)=1\}\).
	
	We claim that
	\begin{align}\label{260706001}
		U\cap B_{R_0}=\bR^d_+\cap B_{R_0},
		\qquad
		\partial U\cap B_{R_0}=\{x^d=0\}\cap B_{R_0}.
	\end{align}
	Indeed, if \(x\in B_{R_0}\), then \(x'\in(-L,L)^{d-1}\) and \(x^d<L\). Hence
	\[
	\operatorname{dist}(x,K)
	=
	\begin{cases}
		1-x^d, & x^d<1,\\
		0, & 1\le x^d<L.
	\end{cases}
	\]
	Therefore, inside \(B_{R_0}\),
	\[
	\operatorname{dist}(x,K)<1
	\quad\Longleftrightarrow\quad
	x^d>0,
	\qquad
	\operatorname{dist}(x,K)=1
	\quad\Longleftrightarrow\quad
	x^d=0,
	\]
	which proves \eqref{260706001}.  On this flat portion of the boundary, the
	outward unit normal is \(-e_d\).
	
	Since the \(x\)-projection of \(Q_R(Z_0)\) is contained in \(B_{R_0}\),
	\eqref{260706001} gives
	\[
	Q_R^U(Z_0)
	=
	Q_R(Z_0)\cap(\bR\times U\times\bR^d_v)
	=
	Q_R(Z_0)\cap\{x^d>0\}
	=
	Q_R^+(Z_0).
	\]
	The boundary characterization and the formula for \(n_U\) follow from the
	second identity in \eqref{260706001}.
\end{proof}

\begin{remark}\label{260630607}
	Definition~\ref{def:L2sol} and all the results of the two preceding
	subsections are stated for a bounded \(C^{1,1}\)-domain \(\cD\).  By
	Lemma~\ref{2606301252} they apply verbatim on a half-space cylinder
	\(Q_R^+(Z_0)\): one chooses a bounded \(C^{1,1}\)-domain \(U\) with
	\(Q_R^+(Z_0)=Q_R^U(Z_0)\), and inside \(Q_R(Z_0)\) the boundary
	\(\p U\) coincides with \(\{x^d=0\}\) and \(n_U=-e_d\).  Consequently, the
	boundary portion \((\bR\times\gamma_U)\cap Q_R(Z_0)\), the measure
	\(\dd\gamma_U\) and the splitting \(\gamma_U^\pm\) coincide with the
	half-space objects \eqref{2608051140} throughout \(Q_R(Z_0)\), so that the
	notions of a weak solution in \(Q_R^U(Z_0)\) and in \(Q_R^+(Z_0)\)
	coincide.  We use this identification without further comment, and we say
	that \(f\) is a weak solution in \(Q_R^+(Z_0)\) accordingly.
	
	The same device applies to a full cylinder \(Q_R(Z_0)\): choosing for
	\(U\) a ball containing the \(x\)-projection of \(Q_R(Z_0)\) gives
	\(Q_R^U(Z_0)=Q_R(Z_0)\) and \((\bR\times\gamma_U)\cap Q_R(Z_0)=\emptyset\),
	so that Definition~\ref{def:L2sol} and
	Corollary~\ref{cor:trace-bound-scaled-cylinder} also provide the notion of a
	weak solution in \(Q_R(Z_0)\), with no boundary term in \eqref{260701105}.
\end{remark}

\begin{lemma}[Kinetic scaling]
	\label{lem:kinetic-scaling}
	Let \(R>0\), \(K>0\), and
	$Z_0=(t_0,x_0,v_0)
	\in
	\bR\times\overline{\bR^d_+}\times\bR^d$.
	Let \(f\) be a weak solution of
	\[
	\partial_t f+v\cdot\nabla_x f
	=
	D_i(a^{ij}D_j f)+b^iD_i f+cf+g^0+D_i g^i
	\qquad\text{in }Q_R^+(Z_0),
	\]
	with
	$f=h$ on $(\bR\times\gamma^-)\cap Q_R(Z_0)$
	in the sense of traces.
	
	For \(Z=(t,x,v)\), define
	\[
	Z_K:=(K^2t,K^3x,Kv),
	\qquad
	Z_{0,K}:=(K^2t_0,K^3x_0,Kv_0).
	\]
	For any function \(\phi\), set
	\[
	\phi_K(Z):=\phi(Z_{K^{-1}}).
	\]
	Then \(f_K\) is a weak solution of
	\[
	\partial_t f_K+v\cdot\nabla_x f_K
	=
	D_i(a_K^{ij}D_j f_K)
	+
	\frac{b_K^i}{K}D_i f_K
	+
	\frac{c_K}{K^2}f_K
	+
	\frac{g_K^0}{K^2}
	+
	D_i\left(\frac{g_K^i}{K}\right)\qquad \text{in } Q_{KR}^+(Z_{0,K}), 
	\]
	with
	$f_K=h_K$ on $(\bR\times\gamma^-)\cap Q_{KR}(Z_{0,K})$
	in the sense of traces.
\end{lemma}

\begin{proof}
	The map \(Z\mapsto Z_K\) sends \(Q_R^+(Z_0)\) onto
	\(Q_{KR}^+(Z_{0,K})\).  Since \(K>0\), it also preserves
	\(\bR\times\gamma^-\).
	
	In \(Q_{KR}^+(Z_{0,K})\), the following identities hold in the sense of
	distributions:
	\[
	(\partial_t+v\cdot\nabla_x)f_K
	=
	K^{-2}
	\bigl((\partial_t+v\cdot\nabla_x)f\bigr)_K,
	\]
	\[
	D_i f_K=K^{-1}(D_i f)_K,
	\qquad
	D_i(a_K^{ij}D_j f_K)
	=
	K^{-2}\bigl(D_i(a^{ij}D_j f)\bigr)_K.
	\]
	Moreover,
	\[
	\frac{b_K^i}{K}D_i f_K
	=
	K^{-2}(b^iD_i f)_K,
	\qquad
	\frac{c_K}{K^2}f_K
	=
	K^{-2}(cf)_K,
	\]
	and
	\[
	\frac{g_K^0}{K^2}
	=
	K^{-2}(g^0)_K,
	\qquad
	D_i\left(\frac{g_K^i}{K}\right)
	=
	K^{-2}(D_i g^i)_K .
	\]
	Composing the original equation with \(Z_{K^{-1}}\) and multiplying by
	\(K^{-2}\), we obtain the scaled equation.
	
	To see that \(f_K\) is a weak solution in the sense of
	Definition~\ref{def:L2sol}, we use the characterization of
	Corollary~\ref{cor:trace-bound-scaled-cylinder}, which is available in the
	half-space by Remark~\ref{260630607}.  The map \(\phi\mapsto\phi_{K^{-1}}\)
	is a bijection from
	\(C_c^\infty(Q_{KR}^{\overline{\bR^d_+}}(Z_{0,K}))\) onto
	\(C_c^\infty(Q_R^{\overline{\bR^d_+}}(Z_0))\), and under \(Z\mapsto Z_K\)
	\[
	\dd Z\mapsto K^{4d+2}\dd Z,
	\qquad
	\dd x\dd v\mapsto K^{4d}\dd x\dd v,
	\qquad
	\dd\gamma\dd t\mapsto K^{4d}\dd\gamma\dd t .
	\]
	Since every volume integrand in \eqref{260701105} carries the extra factor
	\(K^{-2}\) coming from the identities above, both sides of \eqref{260701105}
	are multiplied by \(K^{4d}\).  Hence \eqref{260701105} for \(f\) on
	\(Q_R^+(Z_0)\) is equivalent to \eqref{260701105} for \(f_K\) on
	\(Q_{KR}^+(Z_{0,K})\), with terminal trace \((f(t_0,\cdot))_K\) and kinetic
	trace \((f_\gamma)_K\).  In particular,
	\[
	Z\in(\bR\times\gamma^-)\cap Q_{KR}(Z_{0,K})
	\quad\Longleftrightarrow\quad
	Z_{K^{-1}}\in(\bR\times\gamma^-)\cap Q_R(Z_0),
	\]
	so the trace condition \(f=h\) scales to \(f_K=h_K\) on
	\((\bR\times\gamma^-)\cap Q_{KR}(Z_{0,K})\).
\end{proof}

\subsubsection{Kinetic cylinders and the grazing set}

Recall the definition of the kinetic distance $\mathbf d$ in \eqref{kineticdistance}.

\begin{lemma}[Distance to the grazing set]
	\label{260522259}
	For every \(Z=(t,x,v)\in\bR\times\bR_x^d\times\bR_v^d\),
	\begin{align}\label{260702722}
		\frac{1}{4}\left(|x^d|^{1/3}+|v^d|\right)
		\le
		\mathbf d\bigl(Z,\bR\times\gamma^0\bigr)
		\le
		|x^d|^{1/3}+|v^d|.
	\end{align}
\end{lemma}

\begin{proof}
	For the lower bound, let \(Y=(s,y,u)\in\bR\times\gamma^0\) and \(\xi\in\bR_v^d\), and put
	\[
	\rho:=
	\max\left\{
	|t-s|^{1/2},
	|x-y-(t-s)\xi|^{1/3},
	|v-\xi|,
	|u-\xi|
	\right\}.
	\]
	Since \(y^d=u^d=0\),
	\[
	|v^d|
	\le |v-\xi|+|u-\xi|
	\le 2\rho,
	\]
	and
	\[
	|x^d|
	\le
	|x^d-y^d-(t-s)\xi^d|
	+
	|t-s|\,|\xi^d-u^d|
	\le
	\rho^3+\rho^2\rho
	=
	2\rho^3.
	\]
	Hence
	\[
	|x^d|^{1/3}+|v^d|
	\le
	(2^{1/3}+2)\rho
	\le
	4\rho.
	\]
	Taking the infimum over all \(Y\in\bR\times\gamma^0\) and
	\(\xi\in\bR_v^d\) gives the lower bound.
	
	For the upper bound, write \(x=(x',x^d)\) and \(v=(v',v^d)\) in $\bR^{d-1}\times \bR$.  Choose
	\[
	Y:=\bigl(t,(x',0),(v',0)\bigr)\in\bR\times\gamma^0,
	\qquad
	\xi:=(v',0).
	\]
	This choice gives
	\[
	\mathbf d\bigl(Z,\bR\times\gamma^0\bigr)
	\le
	\max\left\{|x^d|^{1/3},|v^d|\right\}
	\le
	|x^d|^{1/3}+|v^d|.
	\]
	This completes the proof.
\end{proof}

\begin{lemma}[Recentering cylinders at the grazing set]
	\label{260515154}
	Let \(r>0\) and \(Z_0\in\bR\times\bR_x^d\times\bR_v^d\).  Assume that
	\[
	Q_r(Z_0)\cap(\bR\times\gamma^0)\neq\emptyset .
	\]
	Then there exists \(\overline Z_0\in\bR\times\gamma^0\) such that
	\[
	Q_r(Z_0)\subset Q_{2r}(\overline Z_0),
	\qquad
	Q_R(\overline Z_0)\subset Q_{R+2^{1/3}r}(Z_0)
	\quad\text{for every }R>0.
	\]
\end{lemma}

\begin{proof}
	Write \(Z_0=(t_0,x_0,v_0)\), and choose
	\[
	Z_1=(t_1,x_1,v_1)\in Q_r(Z_0)\cap(\bR\times\gamma^0).
	\]
	Then \(x_1^d=v_1^d=0\) and
	\[
	t_0-r^2<t_1\le t_0,\qquad
	|v_1-v_0|<r,\qquad
	|x_1-x_0-(t_1-t_0)v_0|<r^3 .
	\]
	Set
	\[
	\overline Z_0
	:=
	\bigl(t_0,\overline x_0,\overline v_0\bigr)
	:=
	\bigl(t_0,\ x_1-(t_1-t_0)v_1,\ v_1\bigr).
	\]
	Since \(\overline x_0^d=x_1^d-(t_1-t_0)v_1^d=0\) and
	\(\overline v_0^d=v_1^d=0\), we have
	\(\overline Z_0\in\bR\times\gamma^0\).
	
	Let \(Z=(s,y,w)\in Q_r(Z_0)\).  Then
	\[
	t_0-(2r)^2<s\le t_0,
	\qquad
	|w-\overline v_0|
	\le |w-v_0|+|v_0-v_1|
	<2r .
	\]
	Moreover,
	\begin{align*}
		|y-\overline x_0-(s-t_0)\overline v_0|
		&=
		|y-x_1-(s-t_1)v_1|                                      \\
		&\le
		|y-x_0-(s-t_0)v_0|
		+
		|x_1-x_0-(t_1-t_0)v_0|
		+
		|s-t_1|\,|v_0-v_1|                                      \\
		&<
		r^3+r^3+r^3
		<
		(2r)^3 .
	\end{align*}
	Thus
	$Q_r(Z_0)\subset Q_{2r}(\overline Z_0)$.
	
	Next fix \(R>0\) and let \(Z=(s,y,w)\in Q_R(\overline Z_0)\).  Then
	\[
	t_0-(R+2^{1/3}r)^2<s\le t_0,
	\qquad
	|w-v_0|
	\le |w-v_1|+|v_1-v_0|
	<R+r
	<R+2^{1/3}r .
	\]
	For the spatial component, using
	\(|s-t_1|<R^2+r^2\), we get
	\begin{align*}
		|y-x_0-(s-t_0)v_0|
		&\le
		|y-\overline x_0-(s-t_0)\overline v_0|
		+
		|x_1-x_0-(t_1-t_0)v_0|
		+
		|s-t_1|\,|v_1-v_0|                                      \\
		&<
		R^3+r^3+(R^2+r^2)r                                      \\
		&=
		R^3+R^2r+2r^3 .
	\end{align*}
	Since \(1\le 3\cdot 2^{1/3}\),
	\[
	R^3+R^2r+2r^3
	\le
	R^3+3\cdot 2^{1/3}R^2r+3\cdot 2^{2/3}Rr^2+2r^3
	=
	(R+2^{1/3}r)^3 .
	\]
	Therefore
	$Q_R(\overline Z_0)\subset Q_{R+2^{1/3}r}(Z_0)$.
\end{proof}

\begin{lemma}
	\label{lem:grazing-cylinder-geometry}
	Let \(r>0\) and
	$Z_0=(t_0,x_0,v_0)
	\in
	\bR\times\overline{\bR^d_+}\times\bR^d$.
	Assume that
	$Q_r(Z_0)\cap(\bR\times\gamma^0)\neq\emptyset$.
	Then
	\begin{align}\label{260304848}
		|x_0^d|<2r^3,
		\qquad
		|v_0^d|<r.
	\end{align}
	Moreover,
	\begin{align}\label{260304851}
		|Q_r^+(Z_0)|
		\ge
		\varepsilon_0 |Q_r(Z_0)|,
	\end{align}
	where
	\[
	\varepsilon_0
	:=
	\frac{
		|\{y\in B_1:y^d>1/2\}|
	}{
		2|B_1|
	}
	>0 .
	\]
	In particular, \(\varepsilon_0\) depends only on \(d\).
\end{lemma}

\begin{proof}
	Choose
	\[
	\widehat Z=(\widehat t,\widehat x,\widehat v)
	\in Q_r(Z_0)\cap(\bR\times\gamma^0),
	\]
	and set \(\sigma:=\widehat t-t_0\).  Then
	\[
	-r^2<\sigma\le 0,
	\qquad
	|\widehat v-v_0|<r,
	\qquad
	|\widehat x-x_0-\sigma v_0|<r^3 .
	\]
	Since \(\widehat x^d=\widehat v^d=0\),
	\[
	|v_0^d|<r,
	\qquad
	|x_0^d+\sigma v_0^d|<r^3 .
	\]
	Hence
	\[
	|x_0^d|
	\le
	|x_0^d+\sigma v_0^d|
	+
	|\sigma|\,|v_0^d|
	<
	2r^3,
	\]
	which proves \eqref{260304848}.
	
	We prove \eqref{260304851}.  Use the coordinates
	\[
	y:=x-x_0-(t-t_0)v_0
	\]
	inside \(Q_r(Z_0)\), and define
	\[
	E:=
	\left\{
	(t,x,v):
	t_0-\frac{r^2}{2}<t\le t_0,\ 
	|y|<r^3,\ 
	y^d>\frac{r^3}{2},\ 
	|v-v_0|<r
	\right\}.
	\]
	Then \(E\subset Q_r(Z_0)\).  Since \(x_0^d\ge0\), \(|v_0^d|<r\), and
	\(|t-t_0|<r^2/2\), every point of \(E\) satisfies
	\[
	x^d
	=
	y^d+x_0^d+(t-t_0)v_0^d
	>
	\frac{r^3}{2}-\frac{r^2}{2}r
	=
	0.
	\]
	Thus \(E\subset Q_r^+(Z_0)\).
	
	The map \((t,x,v)\mapsto(t,y,v)\) has Jacobian one.  Therefore
	\begin{align*}
		|Q_r^+(Z_0)|
		\ge |E|                                                        
		&=
		\frac{r^2}{2}
		\left|\{y\in B_{r^3}:y^d>r^3/2\}\right|
		|B_r|                                                          \\
		&=
		\frac12
		|\{y\in B_1:y^d>1/2\}|
		|B_1|
		r^{2+4d}.
	\end{align*}
	On the other hand,
	\[
	|Q_r(Z_0)|
	=
	|B_1|^2 r^{2+4d}.
	\]
	Hence
	\[
	|Q_r^+(Z_0)|
	\ge
	\frac{
		|\{y\in B_1:y^d>1/2\}|
	}{
		2|B_1|
	}
	|Q_r(Z_0)|,
	\]
	which proves \eqref{260304851}.
\end{proof}

\begin{lemma}
	\label{lem:grazing-shifted-minus-volume}
	Let
	\[
	Z_0=(t_0,x_0,v_0)
	\in
	\bR\times\overline{\bR^d_+}\times\bR^d,
	\qquad
	\widehat Z_0:=(t_0-1,x_0-v_0,v_0).
	\]
	Let \(\omega\in(0,1)\) and \(r\in(0,1)\) satisfy
	\[
	r\le \frac{\omega^3}{8},
	\]
	and assume that
	\[
	Q_{r}(Z_0)\cap(\bR\times\gamma^0)\neq\emptyset .
	\]
	Then
	\[
	|Q_\omega^-(\widehat Z_0)|
	\ge
	\varepsilon_1 |Q_\omega(\widehat Z_0)|,\qquad \varepsilon_1
	:=
	\frac{
		|\{z\in B_1:z^d<-1/2\}|
	}{
		|B_1|
	}
	>0 .
	\]
	In particular, \(\varepsilon_1\) depends only on \(d\).
\end{lemma}

\begin{proof}
	By Lemma~\ref{lem:grazing-cylinder-geometry},
	\begin{align}\label{260312458}
		|x_0^d|<2r^3,
		\qquad
		|v_0^d|<r .
	\end{align}
	Write
	\[
	\widehat t_0:=t_0-1,
	\qquad
	\widehat x_0:=x_0-v_0.
	\]
	For \(t\in(\widehat t_0-\omega^2,\widehat t_0]\), set
	\[
	s:=\widehat t_0-t\in[0,\omega^2).
	\]
	The \(x\)-center of \(Q_\omega(\widehat Z_0)\) at time \(t\) is
	\[
	\widehat x_0+(t-\widehat t_0)v_0
	=
	x_0-(1+s)v_0.
	\]
	Let
	\[
	a_s:=x_0^d-(1+s)v_0^d .
	\]
	Since \(0\le s<\omega^2<1\), \eqref{260312458} gives
	\[
	|a_s|
	\le
	|x_0^d|+(1+s)|v_0^d|
	<
	2r^3+2r
	\le
	4r
	\le
	\frac{\omega^3}{2}.
	\]
	Therefore, with
	\[
	c_d:=|\{z\in B_1:z^d<-1/2\}|>0,
	\]
	we have, for every such \(s\),
	\[
	\left|
	B_{\omega^3}\bigl(x_0-(1+s)v_0\bigr)
	\cap\{x^d<0\}
	\right|
	\ge
	c_d\,\omega^{3d}.
	\]
	Indeed, after translating by \(x_0-(1+s)v_0\) and scaling by
	\(\omega^3\), the left-hand side becomes
	\[
	\omega^{3d}
	\left|
	\left\{
	z\in B_1:
	z^d<-\frac{a_s}{\omega^3}
	\right\}
	\right|,
	\]
	and \(|a_s|/\omega^3\le 1/2\).
	
	It follows that
	\begin{align*}
		|Q_\omega^-(\widehat Z_0)|
		&=
		|B_\omega|
		\int_{\widehat t_0-\omega^2}^{\widehat t_0}
		\left|
		B_{\omega^3}
		\bigl(\widehat x_0+(t-\widehat t_0)v_0\bigr)
		\cap\{x^d<0\}
		\right|
		\dd t\ge
		|B_\omega|\,\omega^2\,c_d\omega^{3d}.
	\end{align*}
	Since
	\[
	|Q_\omega(\widehat Z_0)|
	=
	|B_\omega|\,|B_{\omega^3}|\,\omega^2
	=
	|B_1|^2\omega^{2+4d},
	\]
	we obtain
	\[
	|Q_\omega^-(\widehat Z_0)|
	\ge
	\frac{c_d}{|B_1|}
	|Q_\omega(\widehat Z_0)|,
	\]
	which is the assertion, since \(\varepsilon_1=c_d/|B_1|\).
\end{proof}

\subsubsection{Extensions away from the grazing set}

Let
\[
Z_0=(t_0,x_0,v_0)
\in
\bR\times\overline{\bR^d_+}\times\bR^d,
\qquad
R>0,
\]
and assume the non-grazing condition
\[
Q_R(Z_0)\cap(\bR\times\gamma^0)=\emptyset .
\]
We extend weak solutions from \(Q_R^+(Z_0)\) to the full cylinder
\(Q_R(Z_0)\).  Accordingly, the coefficients \(a^{ij}\) and \(b^i\) are
assumed to be defined on all of \(Q_R(Z_0)\), whereas the source term is only
given on \(Q_R^+(Z_0)\).

Let
\[
a^{ij},\, b^i\in L^\infty(Q_R(Z_0)),
\qquad
g\in L^\infty(Q_R^+(Z_0)).
\]
Assume that there exist constants \(0<\lambda\le \Lambda\) and \(M\ge0\) such
that, for a.e. \(Z\in Q_R(Z_0)\) and every \(\xi\in\bR^d\),
\begin{equation}\label{260702120}
	\lambda |\xi|^2
	\le
	a^{ij}(Z)\xi_i\xi_j,
	\qquad
	|a^{ij}(Z)|\le \Lambda,
	\qquad
	|b(Z)|\le M .
\end{equation}
Let \(f\in L^\infty(Q_R^+(Z_0))\) be a weak solution of
\begin{align}\label{260706944}
	\begin{cases}
		(\partial_t+v\cdot\nabla_x)f
		=
		D_i(a^{ij}D_j f)+b^iD_i f+g
		&\quad\text{in } Q_R^+(Z_0),\\
		f=0
		&\quad\text{on } (\bR\times\gamma^-)\cap Q_R(Z_0),
	\end{cases}
\end{align}
where the boundary condition is understood in the trace sense.

The two lemmas below treat the cases in which the boundary portion of
\(Q_R(Z_0)\) is entirely incoming or entirely outgoing.  The following
observation shows that, under the non-grazing condition, these two cases are
exhaustive.

\begin{lemma}
	\label{2608051530}
	Let \(Z_0=(t_0,x_0,v_0)\in\bR\times\overline{\bR^d_+}\times\bR^d\) and
	\(R>0\) satisfy \(Q_R(Z_0)\cap(\bR\times\gamma^0)=\emptyset\).  Then either
	\[
	Q_R(Z_0)\cap(\bR\times\gamma)\subset\bR\times\gamma^-
	\qquad\text{or}\qquad
	Q_R(Z_0)\cap(\bR\times\gamma)\subset\bR\times\gamma^+ .
	\]
\end{lemma}

\begin{proof}
	Set \(\Sigma:=Q_R(Z_0)\cap(\bR\times\gamma)\).  If \(\Sigma=\emptyset\), both
	inclusions hold.  Otherwise, pick \((t,x,v)\in\Sigma\).  The constraints
	defining \(Q_R(Z_0)\) decouple in \(x\) and \(v\), so
	\(\{(t,x)\}\times B_R(v_0)\subset Q_R(Z_0)\); since \(x^d=0\), this gives
	\(\{(t,x)\}\times B_R(v_0)\subset\Sigma\).  The non-grazing condition then
	forces \(B_R(v_0)\cap\{w^d=0\}=\emptyset\), that is, \(|v_0^d|\ge R\).  If
	\(v_0^d\ge R\), then \(w^d>0\) for every \(w\in B_R(v_0)\), and
	\(\Sigma\subset\bR\times\gamma^-\); if \(v_0^d\le-R\), then \(w^d<0\) on
	\(B_R(v_0)\), and \(\Sigma\subset\bR\times\gamma^+\).
\end{proof}

\begin{lemma}[Zero extension across an inflow boundary]
	\label{lem:inflow-zero-extension}
	Assume that
	\[
	Q_R(Z_0)\cap(\bR\times\gamma)
	\subset
	\bR\times\gamma^- .
	\]
	Under the setting of \eqref{260702120} and \eqref{260706944}, define
	\begin{equation}\label{260309612}
		\tilde f:=f\mathbf 1_{Q_R^+(Z_0)},
		\qquad
		\tilde g:=g\mathbf 1_{Q_R^+(Z_0)}
		\qquad\text{in }Q_R(Z_0).
	\end{equation}
	Then \(\tilde f\) is a weak solution of
	\[
	\partial_t \tilde f+v\cdot\nabla_x\tilde f
	=
	D_i\bigl(a^{ij}D_j\tilde f\bigr)
	+
	b^iD_i\tilde f
	+
	\tilde g
	\qquad\text{in }Q_R(Z_0).
	\]
	Moreover,
	\begin{equation}\label{26312551}
		\operatorname*{osc}_{Q_R(Z_0)}\tilde f
		=
		\operatorname*{osc}_{Q_R^+(Z_0)}f ,
		\qquad
		\|\tilde f\|_{L^\infty(Q_R(Z_0))}
		=
		\|f\|_{L^\infty(Q_R^+(Z_0))} .
	\end{equation}
\end{lemma}

\begin{proof}
	We verify \eqref{260701105} for \(\tilde f\) on the full cylinder
	\(Q_R(Z_0)\), which is legitimate by the second part of
	Remark~\ref{260630607}.  Set
	\(\Sigma:=Q_R(Z_0)\cap(\bR\times\gamma)\) and let
	\(\varphi\in C_c^\infty(Q_R(Z_0))\).  Since \(\Sigma\subset\bR\times\gamma^-\)
	and \(f_\gamma=0\) on \(\Sigma\), the boundary term in \eqref{260701105} for
	\(f\) on \(Q_R^+(Z_0)\) vanishes, and we obtain
	\begin{align*}
		&\iint_{\{x^d>0\}\cap B_{R^3}(x_0)\times B_R(v_0)}\!\!\!
		f(t_0,\cdot)\varphi(t_0,\cdot)\dd x\dd v \\
		&\qquad =
		\iiint_{Q_R^+(Z_0)}
		\Big(
		f(\partial_t\varphi+v\cdot\nabla_x\varphi)
		-
		a^{ij}D_j fD_i\varphi
		+
		b^iD_i f\,\varphi
		+g\varphi
		\Big)\dd Z .
	\end{align*}
	Since the extension is made across the \(x\)-boundary only,
	\[
	D_v\tilde f=(D_v f)\mathbf 1_{Q_R^+(Z_0)}
	\qquad\text{in }\cD'(Q_R(Z_0)),
	\]
	so, by \eqref{260309612}, the preceding identity is exactly
	\eqref{260701105} for \(\tilde f\) on \(Q_R(Z_0)\), with terminal trace
	\(f(t_0,\cdot)\mathbf 1_{\{x^d>0\}}\), source \(\tilde g\), and no boundary
	term.  By Corollary~\ref{cor:trace-bound-scaled-cylinder}, \(\tilde f\) is a
	weak solution in \(Q_R(Z_0)\).
	
	We prove \eqref{26312551}.  If \(\Sigma=\emptyset\), then
	\(Q_R(Z_0)=Q_R^+(Z_0)\) and both identities are immediate.  Otherwise
	\(\Sigma\) is a nonempty relatively open subset of
	\(Q_R(Z_0)\cap(\bR\times\gamma^-)\), hence has positive
	\(\dd\gamma\dd t\)-measure.  If we had
	\(\operatorname*{ess\,sup}_{Q_R^+(Z_0)}f<0\), then
	\eqref{2606301239} would force \(f_\gamma<0\) \(\dd\gamma\dd t\)-a.e. on
	\(\Sigma\), contradicting \(f_\gamma=0\) there;
	applying the same argument to \(-f\), we conclude
	\[
	\operatorname*{ess\,inf}_{Q_R^+(Z_0)} f
	\le 0
	\le
	\operatorname*{ess\,sup}_{Q_R^+(Z_0)} f .
	\]
	Since \(\tilde f=0\) on \(Q_R(Z_0)\setminus Q_R^+(Z_0)\), this gives
	\[
	\operatorname*{ess\,sup}_{Q_R(Z_0)}\tilde f
	=
	\operatorname*{ess\,sup}_{Q_R^+(Z_0)} f,
	\qquad
	\operatorname*{ess\,inf}_{Q_R(Z_0)}\tilde f
	=
	\operatorname*{ess\,inf}_{Q_R^+(Z_0)} f,
	\]
	and \eqref{26312551} follows.
\end{proof}

\begin{lemma}[Extension across an outgoing boundary]
	\label{lem:boundary-extension-away-from-grazing}
	Assume that
	\[
	Q_R(Z_0)\cap(\bR\times\gamma)
	\subset
	\bR\times\gamma^+ .
	\]
	Under the setting of \eqref{260702120} and \eqref{260706944}, there exists an
	extension \(\bar f\) of \(f\) to \(Q_R(Z_0)\) such that:
	\begin{enumerate}
		\item \(\bar f=f\) a.e. in \(Q_R^+(Z_0)\);
		
		\item \(\bar f\) is a weak solution of
		\begin{align}\label{260702332}
			(\partial_t+v\cdot\nabla_x)\bar f
			=
			D_i(a^{ij}D_j\bar f)
			+
			b^iD_i\bar f
			+
			g\,\mathbf 1_{Q_R^+(Z_0)}
			\qquad\text{in }Q_R(Z_0);
		\end{align}
		
		\item
		\[
		\|\bar f\|_{L^\infty(Q_R(Z_0))}
		=
		\|f\|_{L^\infty(Q_R^+(Z_0))}.
		\]
	\end{enumerate}
\end{lemma}

\begin{proof}
	Set
	\[
	\Sigma:=Q_R(Z_0)\cap(\bR\times\gamma).
	\]
	If \(\Sigma=\emptyset\), then \(Q_R(Z_0)=Q_R^+(Z_0)\), and there is nothing to
	prove.  Hence we assume \(\Sigma\neq\emptyset\).
	
	By Corollary~\ref{cor:trace-bound-scaled-cylinder}, \(f\) has a kinetic trace
	\(f_\gamma\) on \(\Sigma\), and
	\begin{align}\label{260702342}
		\|f_\gamma\|_{L^\infty(\Sigma,\dd\gamma\dd t)}
		\le
		\|f\|_{L^\infty(Q_R^+(Z_0))}.
	\end{align}
	
	By Lemma~\ref{2606301252}, applied after the reflection
	$(x^d,v^d)\mapsto(-x^d,-v^d)$,
	there exists a bounded \(C^{1,1}\)-domain \(U\subset\bR^d\) such that
	\[
	Q_R^U(Z_0)=Q_R^-(Z_0),\qquad Q_R(Z_0)\cap \{x^d=0\}=Q_R(Z_0)\cap \{x\in \p U\}.
	\]
	Moreover, for every \((t,x,v)\in Q_R(Z_0)\cap \{x^d=0\}\), one has
	\[
	n_U(x)=e_d .
	\]
	We use the convention
	\[
	\gamma_U^\pm
	:=
	\{(x,v)\in\partial U\times\bR^d_v:\ \pm n_U(x)\cdot v>0\},
	\]
	where
	\[
	\dd\gamma_U:=|n_U(x)\cdot v|\dd S_{\p U}(x)\dd v,
	\]
	and \(\dd S_{\p U}\) is the surface measure on \(\p U\).
	Since \(n_U=e_d\) on \(\Sigma\), the assumption
	\(\Sigma\subset\bR\times\gamma^+=\bR\times\{x^d=0,\ v^d<0\}\) gives
	$\Sigma\subset \bR\times\gamma_U^-$.
	As \(\Sigma\) is bounded, it has finite \(\dd\gamma_U\dd t\)-measure, so
	\eqref{260702342} gives
	\[
	f_\gamma\mathbf 1_\Sigma
	\in
	L^2(\gamma_U^-,\dd\gamma_U\dd t).
	\]
	
	The coefficients \(a^{ij},b^i\) are defined only on \(Q_R(Z_0)\); we extend
	them to \((t_0-R^2,t_0)\times U\times\bR^d_v\) by setting
	\(a^{ij}:=\lambda\delta^{ij}\) and \(b^i:=0\) outside \(Q_R(Z_0)\), which
	preserves \eqref{260702120}.  Since the extended coefficients satisfy
	\eqref{coeff_ellip} and \eqref{coeff_bdd} with
	\(\nu=\min\{\lambda,\Lambda^{-1}\}\), \(K=M\) and \(\bar b^i=c=0\),
	Theorem~\ref{solv}\,(1) provides a unique
	\[
	f^\ast\in
	\cH_{\kin}\big((t_0-R^2,t_0)\times U\times\bR^d_v\big)
	\]
	satisfying
	\[
	\partial_t f^\ast+v\cdot\nabla_x f^\ast
	=
	D_i(a^{ij}D_jf^\ast)+b^iD_i f^\ast
	\qquad\text{in }(t_0-R^2,t_0)\times U\times\bR^d_v,
	\]
	with
	\[
	f^\ast(t_0-R^2,\cdot)=0,
	\qquad
	(f^\ast)_{\gamma_U^-}=f_\gamma\mathbf 1_\Sigma
	\quad\text{on }(t_0-R^2,t_0)\times\gamma_U^- .
	\]
	The extension of the coefficients is irrelevant for the conclusion:
	on \(Q_R^-(Z_0)\subset Q_R(Z_0)\) the extended coefficients coincide with the
	original ones.
	
	Define
	\[
	\bar f(Z):=
	\begin{cases}
		f(Z), & Z\in Q_R^+(Z_0),\\
		f^\ast(Z), & Z\in Q_R^-(Z_0).
	\end{cases}
	\]
	Then \(\bar f=f\) a.e. in \(Q_R^+(Z_0)\).
	
	We prove \eqref{260702332} by verifying the identity \eqref{260701105} for
	\(\bar f\) on the full cylinder \(Q_R(Z_0)\), which is legitimate by the
	second part of Remark~\ref{260630607}.  Let
	\(\varphi\in C_c^\infty(Q_R(Z_0))\), and set
	\begin{align*}
		I:=\,\,&
		\iiint_{Q_R(Z_0)}
		\Big(
		\bar f(-\partial_t\varphi-v\cdot\nabla_x\varphi)
		+
		a^{ij}D_j\bar fD_i\varphi
		-
		b^iD_i\bar f\,\varphi
		\Big)\dd Z\\
		=\,\,&\iiint_{Q_R^+(Z_0)}(\cdots) \dd Z + \iiint_{Q_R^-(Z_0)}(\cdots) \dd Z
		=:I_f+I_{f^\ast}\,.
	\end{align*}
	On \(\Sigma\) the outward unit normal of \(\bR^d_+\) is \(-e_d\), so
	\(\sgn(n_x\cdot v)=\sgn(-v^d)=1\) and \(\dd\gamma=(-v^d)\dd x'\dd v\).
	Hence \eqref{260701105} for \(f\) on \(Q_R^+(Z_0)\) gives
	\[
	I_f
	=
	-\iint_{\left(\{x^d>0\}\cap B_{R^3}(x_0)\right)\times B_R(v_0)}\!\!\!
	f(t_0,\cdot)\varphi(t_0,\cdot)\dd x\dd v
	+\iiint_{Q_R^+(Z_0)}g\varphi\dd Z
	-
	\iiint_{\Sigma}
	f_\gamma\,\varphi\,(-v^d)\dd x'\dd v\dd t .
	\]
	Similarly, \(n_U=e_d\) on \(\Sigma\), so \(\sgn(n_U\cdot v)=-1\) and
	\(\dd\gamma_U=(-v^d)\dd x'\dd v\) there; applying \eqref{ibp_dy} to \(f^\ast\) and
	\(\varphi\) on \((t_0-R^2,t_0)\times U\times\bR^d_v\), together with the
	equation for \(f^\ast\) and \(f^\ast(t_0-R^2,\cdot)=0\), we obtain
	\[
	I_{f^\ast}
	=
	-\iint_{\left(\{x^d<0\}\cap B_{R^3}(x_0)\right)\times B_R(v_0)}\!\!\!
	f^\ast(t_0,\cdot)\varphi(t_0,\cdot)\dd x\dd v
	+
	\iiint_{\Sigma}
	(f^\ast)_{\gamma_U^-}\,\varphi\,(-v^d)\dd x'\dd v\dd t .
	\]
	Since \((f^\ast)_{\gamma_U^-}=f_\gamma\) on \(\Sigma\), the two interface terms cancel.
	Adding the two identities and moving the terminal terms to the left-hand
	side, we arrive at
	\[
	\iint_{B_{R^3}(x_0)\times B_R(v_0)}\!\!\!\bar f(t_0,\cdot)\varphi(t_0,\cdot)\dd x\dd v
	+I
	=
	\iiint_{Q_R(Z_0)}
	g\,\mathbf 1_{Q_R^+(Z_0)}\varphi\dd Z ,
	\]
	which is exactly \eqref{260701105} for \(\bar f\) on \(Q_R(Z_0)\), with
	terminal trace
	\(\bar f(t_0,\cdot)=f(t_0,\cdot)\mathbf 1_{\{x^d>0\}}+f^\ast(t_0,\cdot)\mathbf 1_{\{x^d<0\}}\)
	and no boundary term.  Since \(\nabla_v f\in
	L^2_{\mathrm{loc}}(Q_R^{\overline{\bR^d_+}}(Z_0))\) and
	\(\nabla_vf^\ast\in L^2\), we have
	\(\nabla_v\bar f\in L^2_{\mathrm{loc}}(Q_R(Z_0))\), so
	Corollary~\ref{cor:trace-bound-scaled-cylinder} proves \eqref{260702332}.
	
	It remains to prove the \(L^\infty\)-bound.
	We show that
	\begin{align}\label{260702343}
		\|f^\ast\|_{L^\infty(Q_R^-(Z_0))}\le \|f_\gamma\|_{L^\infty(\Sigma,\dd\gamma\dd t)}=:K .
	\end{align}
	Let
	\[
	\beta(s):=\frac12(|s|-K)_+^2.
	\]
	Choose smooth convex functions \(\beta_m\in C^2(\bR)\) such that
	\[
	0\le \beta_m\uparrow \beta,
	\qquad
	|\beta_m'(s)|\le (|s|-K)_+\le |s|,
	\qquad
	0\le \beta_m''\uparrow \mathbf 1_{\{|s|>K\}},\qquad \beta_m''\in C_c^\infty(\bR).
	\]
	Let \(\eta_N(v):=\eta(v/N)\), where
	\(\eta\in C_c^\infty(\bR^d_v)\), \(0\le\eta\le1\), and
	\(\eta\equiv1\) on \(B_1\).  We apply the renormalized identity
	\eqref{reno_dyn} to \(f^\ast\), with renormalization \(\beta_m\), test function
	\(\eta_N\), and divergence-form data \(g^i=a^{ij}D_jf^\ast\),
	\(g^0=b^iD_if^\ast\).  Since \((\p_t+v\cdot\nabla_x)\eta_N=0\), letting
	\(N\to\infty\) gives, for every \(t\in(t_0-R^2,t_0]\),
	\begin{align}\label{260702344}
		&\iint_{U\times\bR^d_v}
		\beta_m(f^\ast(t))\dd x\dd v
		+
		\iiint_{(t_0-R^2,t)\times U\times\bR^d_v}
		\beta_m''(f^\ast)a^{ij}D_i f^\ast D_j f^\ast\dd Z
		\nonumber\\
		&\qquad
		+
		\iint_{(t_0-R^2,t)\times\gamma_U^+}
		\beta_m((f^\ast)_{\gamma_U^+})\dd\gamma_U\dd s
		\nonumber\\
		&\le
		\iiint_{(t_0-R^2,t)\times U\times\bR^d_v}
		\beta_m'(f^\ast)b^iD_i f^\ast\dd Z
		\nonumber\\
		&\qquad
		+
		\iint_{U\times\bR^d_v}
		\beta_m(f^\ast(t_0-R^2))\dd x\dd v
		+
		\iint_{(t_0-R^2,t)\times\gamma_U^-}
		\beta_m((f^\ast)_{\gamma_U^-})\dd\gamma_U\dd s .
	\end{align}
	Indeed, the only cutoff error is bounded by
	\[
	\frac{C\Lambda}{N}
	\left(
	\iiint |f^\ast|^2\dd Z
	\right)^{1/2}
	\left(
	\iiint |D_vf^\ast|^2\dd Z
	\right)^{1/2},
	\]
	and hence vanishes as \(N\to\infty\).  Letting \(m\to\infty\) in
	\eqref{260702344}, we obtain
	\begin{align*}
		&\iint_{U\times\bR^d_v}
		\beta(f^\ast(t))\dd x\dd v
		+
		\iiint_{(t_0-R^2,t)\times U\times\bR^d_v}
		\mathbf 1_{\{|f^\ast|>K\}}a^{ij}D_i f^\ast D_j f^\ast\dd Z
		\nonumber\\
		&\qquad
		+
		\iint_{(t_0-R^2,t)\times\gamma_U^+}
		\beta((f^\ast)_{\gamma_U^+})\dd\gamma_U\dd s
		\nonumber\\
		&\le
		\iiint_{(t_0-R^2,t)\times U\times\bR^d_v}
		\beta'(f^\ast)b^iD_i f^\ast\dd Z
		\nonumber\\
		&\qquad
		+
		\iint_{U\times\bR^d_v}
		\beta(f^\ast(t_0-R^2))\dd x\dd v
		+
		\iint_{(t_0-R^2,t)\times\gamma_U^-}
		\beta((f^\ast)_{\gamma_U^-})\dd\gamma_U\dd s .
	\end{align*}
	The initial term vanishes because \(f^\ast(t_0-R^2,\cdot)=0\).  The incoming
	boundary term also vanishes since
	\[
	\big|(f^\ast)_{\gamma_U^-}\big|=|f_\gamma\mathbf 1_\Sigma|\le K
	\quad\text{a.e. on }(t_0-R^2, t_0)\times \gamma_U^-.
	\]
	Since $\beta\ge 0$, discarding the nonnegative \(\gamma_U^+\)-term and setting
	\[
	u:=(|f^\ast|-K)_+,
	\]
	we get
	\[
	\frac12\iint_{U\times\bR^d_v}u^2(t)\dd x\dd v
	+
	\lambda
	\iiint_{(t_0-R^2,t)\times U\times\bR^d_v}
	|D_vu|^2\dd Z
	\le
	\iiint_{(t_0-R^2,t)\times U\times\bR^d_v}
	|b|\,u\,|D_vu|\dd Z .
	\]
	Using \eqref{260702120} and Young's inequality,
	\[
	\iint_{U\times\bR^d_v}u^2(t)\dd x\dd v
	\le
	\frac{M^2}{\lambda}
	\int_{t_0-R^2}^t
	\iint_{U\times\bR^d_v}u^2(s)\dd x\dd v\dd s \qquad \text{for all }t\in (t_0-R^2,t_0].
	\]
	Gronwall's inequality gives \(u\equiv0\).  Hence \(|f^\ast|\le K\) a.e. in
	$(t_0-R^2,t_0)\times U\times\bR^d_v$,
	and \eqref{260702343} follows.
	
	Combining \eqref{260702343} with \eqref{260702342}, we obtain
	\[
	\|\bar f\|_{L^\infty(Q_R(Z_0))}
	\le
	\max\left\{
	\|f\|_{L^\infty(Q_R^+(Z_0))},
	\|f_\gamma\|_{L^\infty(\Sigma,\dd\gamma\dd t)}
	\right\}
	\le
	\|f\|_{L^\infty(Q_R^+(Z_0))}.
	\]
	This completes the proof.
\end{proof}

\subsubsection{A smallness criterion in the half-space}

We record the following corollary of Theorem~\ref{cor:linfty-zero-inflow}, which
will be used in Appendix~\ref{proofofholder}.

\begin{corollary}
	\label{cor:scaled-smallness-implies-half}
	Let \(Z_0\in\bR\times\overline{\bR^d_+}\times\bR^d\), \(R>0\), and \(K>0\).
	Assume that \(a^{ij}\), \(b^i\), and \(c\) satisfy \eqref{con_abc} in
	\(Q_R^+(Z_0)\), with constants \(\lambda,\Lambda\), and \(M\).  Then there
	exists \(\varepsilon>0\), depending only on \(d\), \(\lambda\), \(\Lambda\),
	\(M\), \(K\), and \(R\), such that the following holds.
	
	Let \(g\in L^\infty(Q_R^+(Z_0))\), and let \(f\) be a weak solution of
	\[
	(\partial_t+v\cdot\nabla_x)f
	=
	D_i(a^{ij}D_jf)+b^iD_if+cf+g
	\qquad\text{in }Q_R^+(Z_0)
	\]
	such that
	\[
	f=h\le 0
	\qquad\text{on}\qquad
	(\bR\times\gamma^-)\cap Q_R(Z_0)
	\]
	in the sense of traces. If
	\[
	\|f_+\|_{L^2(Q_R^+(Z_0))}\le \varepsilon,
	\qquad
	\|g\|_{L^\infty(Q_R^+(Z_0))}\le K,
	\]
	then
	\[
	f\le \frac12
	\qquad\text{in}\qquad
	Q_{R/2}^+(Z_0).
	\]
\end{corollary}

\begin{proof}
	Let
	$C_1:=C(d,\lambda,\Lambda,1/2)$ be the constant in
	Theorem~\ref{cor:linfty-zero-inflow} with $\kappa:=\frac 12$. Set
	\[
	C_0:=
	C_1
	\left(\frac M2+1\right)^{4d+2}2^{1+2d}.
	\]
	Choose \(0<\rho<1\) such that
	\begin{equation}\label{260629001}
		\rho R\le 1
		\qquad\text{and}\qquad
		C_0\,2^{-(3+2d)}(\rho R)^2K\le \frac14.
	\end{equation}
	Such a \(\rho\) exists, and it depends only on
	\(d,\lambda,\Lambda,M,K\), and \(R\).
	
	Set
	\[
	A:=\frac1{\rho R}.
	\]
	Then \(A\ge1\). For \(Z=(t,x,v)\), define
	\[
	Z_A:=\delta_A Z=(A^2t,A^3x,Av),
	\qquad
	Z_{0,A}:=\delta_A Z_0.
	\]
	For any function \(\phi\), write
	\[
	\phi_A(Z):=\phi(\delta_{A^{-1}}Z).
	\]
	In particular, define
	\[
	\widetilde f:=f_A,
	\qquad
	\widetilde a^{ij}:=a_A^{ij},
	\qquad
	\widetilde b^i:=A^{-1}b_A^i,
	\qquad
	\widetilde c:=A^{-2}c_A,
	\qquad
	\widetilde g:=A^{-2}g_A.
	\]
	By Lemma~\ref{lem:kinetic-scaling}, applied with \(K=A\), \(\widetilde f\) is
	defined in
	\[
	Q_{AR}^+(Z_{0,A})=Q_{1/\rho}^+(Z_{0,A})
	\]
	and satisfies
	\[
	(\partial_t+v\cdot\nabla_x)\widetilde f
	=
	D_i(\widetilde a^{ij}D_j\widetilde f)
	+\widetilde b^iD_i\widetilde f
	+\widetilde c\,\widetilde f
	+\widetilde g
	\qquad\text{in }Q_{1/\rho}^+(Z_{0,A}).
	\]
	Moreover, since \(A^{-1}=\rho R\le1\), the rescaled coefficients satisfy
	\[
	\lambda|\xi|^2
	\le
	\widetilde a^{ij}(Z)\xi_i\xi_j,
	\qquad
	|\widetilde a^{ij}(Z)|\le\Lambda,
	\qquad
	|\widetilde b(Z)|\le M,
	\qquad
	|\widetilde c(Z)|\le M^2
	\]
	for a.e. \(Z\in Q_{1/\rho}^+(Z_{0,A})\) and every \(\xi\in\bR^d\). Also,
	\[
	\|\widetilde g\|_{L^\infty(Q_{1/\rho}^+(Z_{0,A}))}
	\le
	A^{-2}\|g\|_{L^\infty(Q_R^+(Z_0))}
	\le
	(\rho R)^2K.
	\]
	
	The inflow boundary condition is preserved by the scaling. Hence
	$\widetilde f\le0$ on $
	(\bR\times\gamma^-)\cap Q_{1/\rho}(Z_{0,A})$
	in the sense of traces.
	
	Now fix
	\[
	Z_1\in Q_{1/(2\rho)}^+(Z_{0,A}).
	\]
	By Lemma~\ref{260127226}\,(3) and \(\rho<1\),
	\[
	Q_{1/2}(Z_1)\subset Q_{\frac{1}{2\rho}+\frac12}(Z_{0,A})
	\subset Q_{1/\rho}(Z_{0,A}),
	\qquad\text{hence}\qquad
	Q_{1/2}^+(Z_1)\subset Q_{1/\rho}^+(Z_{0,A}).
	\]
	Applying Theorem~\ref{cor:linfty-zero-inflow} to \(\widetilde f\) in
	\(Q_{1/2}^+(Z_1)\), with \(\kappa=1/2\), which is legitimate by
	Remark~\ref{260630607}, gives
	\begin{align*}
		\|\widetilde f_+\|_{L^\infty(Q_{1/4}^+(Z_1))}
		\le{}&
		C_1
		\left(\frac M2+1\right)^{4d+2}
		2^{1+2d}
		\\
		&\quad\times
		\left(
		\|\widetilde f_+\|_{L^2(Q_{1/2}^+(Z_1))}
		+
		2^{-(3+2d)}
		\|\widetilde g\|_{L^\infty(Q_{1/2}^+(Z_1))}
		\right)
		\\
		\le{}&
		C_0
		\|\widetilde f_+\|_{L^2(Q_{1/\rho}^+(Z_{0,A}))}
		+
		C_0\,2^{-(3+2d)}(\rho R)^2K.
	\end{align*}
	By \eqref{260629001}, the second term is at most \(1/4\).
	
	Under the scaling \(Z\mapsto Z_A=:\widetilde Z\), the kinetic measure scales as
	\[
	\dd \widetilde Z=A^{4d+2}\dd Z.
	\]
	Therefore
	\[
	\|\widetilde f_+\|_{L^2(Q_{1/\rho}^+(Z_{0,A}))}
	=
	A^{2d+1}
	\|f_+\|_{L^2(Q_R^+(Z_0))}.
	\]
	Thus, if
	\[
	\|f_+\|_{L^2(Q_R^+(Z_0))}
	\le
	\frac{1}{4C_0}A^{-(2d+1)},
	\]
	then
	\[
	C_0
	\|\widetilde f_+\|_{L^2(Q_{1/\rho}^+(Z_{0,A}))}
	\le
	\frac14.
	\]
	Consequently,
	\[
	\|\widetilde f_+\|_{L^\infty(Q_{1/4}^+(Z_1))}
	\le
	\frac12.
	\]
	
	Since \(Z_1\in Q_{1/(2\rho)}^+(Z_{0,A})\) was arbitrary and the cylinders
	\(Q_{1/4}^+(Z_1)\) cover \(Q_{1/(2\rho)}^+(Z_{0,A})\), we obtain
	$\widetilde f\le \frac12$ in $Q_{1/(2\rho)}^+(Z_{0,A})$.
	Scaling back gives
	$f\le \frac12$ in $Q_{R/2}^+(Z_0)$.
	
	Finally, define
	\[
	\varepsilon
	:=
	\frac{1}{4C_0}A^{-(2d+1)}
	=
	\frac{1}{4C_0}(\rho R)^{2d+1}.
	\]
	Then \(\varepsilon\) depends only on \(d,\lambda,\Lambda,M,K\), and \(R\).
	This completes the proof.
\end{proof}

\subsection{H\"older regularity up to the boundary for rough coefficients}
\label{Subsec:holderconti}

Boundary H\"older regularity for kinetic equations with rough
coefficients is established in \cite[Theorem~1.2]{SilvestreBoundaryKFP}; see also \cite{2019GIMV} for the corresponding whole-space theory.

\begin{thm}[Scaled boundary H\"older estimate]
	\label{cor:scaled-boundary-holder-estimate}
	Let
	$Z_0\in\bR\times\overline{\bR^d_+}\times\bR^d$, $r>0$.
	Assume that
	\[
	a^{ij},\,b^i\in L^\infty(Q_{2r}^+(Z_0)),
	\qquad i,j=1,\ldots,d,
	\]
	and that there exist constants \(0<\lambda\le \Lambda\) and \(M\ge0\) such
	that, for a.e. \(Z\in Q_{2r}^+(Z_0)\) and every \(\xi\in\bR^d\),
	\begin{align}\label{con_ab}
		\lambda |\xi|^2
		\le
		a^{ij}(Z)\xi_i\xi_j,
		\qquad
		|a^{ij}(Z)|\le \Lambda,
		\qquad
		|b(Z)|\le M .
	\end{align}
	Let \(g\in L^\infty(Q_{2r}^+(Z_0))\), and let
	\(f\in L^\infty(Q_{2r}^+(Z_0))\) be a weak solution of
	\[
	(\partial_t+v\cdot\nabla_x)f
	=
	D_i(a^{ij}D_j f)+b^iD_i f+g
	\qquad\text{in }Q_{2r}^+(Z_0),
	\]
	with
	\begin{align}\label{2606291020}
		f=0
		\qquad
		\text{on }(\bR\times\gamma^-)\cap Q_{2r}(Z_0)
	\end{align}
	in the sense of traces.
	
	There exist \(\alpha=\alpha(d,\lambda,\Lambda)\in(0,1)\) and
	\(C=C(d,\lambda,\Lambda)>0\), independent of \(Z_0\), \(r\), \(M\),
	\(f\), and \(g\), such that the following holds.  Under the above
	assumptions, \(f\) admits a representative, still denoted by \(f\), for
	which
	\begin{align}\label{eq:scaled-boundary-holder-estimate}
		\frac{|f(Z_1)-f(Z_2)|}{\mathbf d(Z_1,Z_2)^\alpha}
		\le
		C\left(M+r^{-1}\right)^\alpha
		\left(
		\|f\|_{L^\infty(Q_{2r}^+(Z_0))}
		+
		\left(M+r^{-1}\right)^{-2}
		\|g\|_{L^\infty(Q_{2r}^+(Z_0))}
		\right)
	\end{align}
	for all distinct \(Z_1,Z_2\in Q_r^+(Z_0)\).
\end{thm}

\begin{prop}\label{2601271146}
	Theorem~\ref{cor:scaled-boundary-holder-estimate} holds when \(r=1\).
\end{prop}

A version of Proposition \ref{2601271146} is proved in \cite[Theorem~1.2]{SilvestreBoundaryKFP}.
For completeness, Appendix~\ref{proofofholder} gives a proof for weak solutions in the half-space, based on the standard De Giorgi oscillation-decay argument.

\begin{proof}[\textbf{Proof of Theorem~\ref{cor:scaled-boundary-holder-estimate}}]
	For \(Z=(t,x,v)\) and \(r>0\), write
	\[
	Z_r:=(r^2t,r^3x,rv),
	\qquad
	Z_{0,r^{-1}}:=(r^{-2}t_0,r^{-3}x_0,r^{-1}v_0).
	\]
	Define
	\[
	f_r(Z):=f(Z_r),
	\qquad
	a_r^{ij}(Z):=a^{ij}(Z_r),
	\qquad
	b_r^i(Z):=b^i(Z_r),
	\qquad
	g_r(Z):=g(Z_r).
	\]
	By Lemma~\ref{lem:kinetic-scaling}, applied with \(K=r^{-1}\), the function
	\(f_r\) is a weak solution of
	\[
	(\partial_t+v\cdot\nabla_x)f_r
	=
	D_i(a_r^{ij}D_j f_r)
	+
	r b_r^iD_i f_r
	+
	r^2 g_r
	\qquad\text{in }Q_2^+(Z_{0,r^{-1}}),
	\]
	with $f_r=0$ on $(\bR\times\gamma^-)\cap Q_2(Z_{0,r^{-1}})$
	in the sense of traces. 
	Moreover,
	\[
	|r b_r|\le rM
	\qquad\text{a.e. in }Q_2^+(Z_{0,r^{-1}}).
	\]
	
	Applying Proposition~\ref{2601271146} to \(f_r\), with drift bound \(rM\) and
	source \(r^2g_r\), we obtain, for all distinct
	\(\widetilde Z_1,\widetilde Z_2\in Q_1^+(Z_{0,r^{-1}})\),
	\[
	\frac{|f_r(\widetilde Z_1)-f_r(\widetilde Z_2)|}
	{\mathbf d(\widetilde Z_1,\widetilde Z_2)^\alpha}
	\le
	C(rM+1)^\alpha
	\left(
	\|f_r\|_{L^\infty(Q_2^+(Z_{0,r^{-1}}))}
	+
	(rM+1)^{-2}
	\|r^2g_r\|_{L^\infty(Q_2^+(Z_{0,r^{-1}}))}
	\right).
	\]
	
	Now let \(Z_1,Z_2\in Q_r^+(Z_0)\) be distinct and set
	\[
	\widetilde Z_i:=Z_{i,r^{-1}},
	\qquad i=1,2.
	\]
	Then
	\[
	\widetilde Z_1,\widetilde Z_2\in Q_1^+(Z_{0,r^{-1}}),
	\qquad
	f_r(\widetilde Z_i)=f(Z_i),\qquad \mathbf d(\widetilde Z_1,\widetilde Z_2)
	=
	r^{-1}\mathbf d(Z_1,Z_2).
	\]
	Hence
	\[
	\frac{|f(Z_1)-f(Z_2)|}{\mathbf d(Z_1,Z_2)^\alpha}
	=
	r^{-\alpha}
	\frac{|f_r(\widetilde Z_1)-f_r(\widetilde Z_2)|}
	{\mathbf d(\widetilde Z_1,\widetilde Z_2)^\alpha}.
	\]
	Using
	\[
	\|f_r\|_{L^\infty(Q_2^+(Z_{0,r^{-1}}))}
	=
	\|f\|_{L^\infty(Q_{2r}^+(Z_0))}\quad \text{and}\quad \|r^2g_r\|_{L^\infty(Q_2^+(Z_{0,r^{-1}}))}
	=
	r^2\|g\|_{L^\infty(Q_{2r}^+(Z_0))},
	\]
	we get
	\[
	\frac{|f(Z_1)-f(Z_2)|}{\mathbf d(Z_1,Z_2)^\alpha}
	\le
	C r^{-\alpha}(rM+1)^\alpha
	\left(
	\|f\|_{L^\infty(Q_{2r}^+(Z_0))}
	+
	(rM+1)^{-2}r^2
	\|g\|_{L^\infty(Q_{2r}^+(Z_0))}
	\right),
	\]
	which proves the scaled estimate.
\end{proof}

\subsection{Decay near the grazing set for regular coefficients}

We prove the weighted grazing decay estimate by constructing a local barrier,
transferring the comparison bound to weak solutions through smooth
approximation, and then scaling the result.

For a set \(E\subset \bR^{1+2d}\) and \(\alpha\in(0,1)\), we denote by
\(C^\alpha_{\kin}(E)\) the space of all \(\varphi\in C(E)\) such that
\[
[\varphi]_{C^\alpha_{\kin}(E)}
:=
\sup_{\substack{Z,Y\in E\\ Z\neq Y}}
\frac{|\varphi(Z)-\varphi(Y)|}{\mathbf d(Z,Y)^\alpha}
<\infty.
\]
For bounded \(\varphi\in C^\alpha_{\kin}(E)\), we write
\begin{equation}\label{eq:Ckin-norm}
	\|\varphi\|_{C^\alpha_{\kin}(E)}
	:=\|\varphi\|_{L^\infty(E)}+[\varphi]_{C^\alpha_{\kin}(E)},
\end{equation}
and drop \(E\) when it is clear from the context.

\begin{thm}[Weighted decay estimate near the grazing set]
	\label{260129331}
	Let $
	Z_0=(t_0,x_0,v_0)\in \bR\times\gamma^0$
	and $r>0$.
	Assume that there exists $\alpha\in(0,1)$ such that 
	\[
	a^{ij}\in C^\alpha_{\kin}(Q_{2r}(Z_0)),\quad
	\nabla_v a^{ij},\,b^i\in L^\infty(Q_{2r} (Z_0)),
	\]
	where \(\nabla_v a^{ij}\) is understood in the weak sense. Assume moreover
	that there exist constants \(0<\lambda\le \Lambda\) and \(M_1,M_2\ge0\) such
	that, for a.e. \(Z\in Q_{2r} (Z_0)\) and every \(\xi\in\bR^d\),
	\begin{align}\label{260605507}
		\lambda|\xi|^2
		\le
		a^{ij}(Z)\xi_i\xi_j,
		\qquad
		|a^{ij}(Z)|\le \Lambda,
	\end{align}
	and
	\begin{align}\label{2607041007}
		[a]_{C^\alpha_{\kin}(Q_{2r} (Z_0))}\le M_1,
		\qquad
		\|\nabla_v a\|_{L^\infty(Q_{2r} (Z_0))}
		+
		\|b\|_{L^\infty(Q_{2r} (Z_0))}
		\le M_2.
	\end{align}
	Set
	\[
	M:=\max\{M_1^{1/\alpha},\,M_2\}.
	\]
	
	Let \(f\in L^\infty(Q_{2r}^+(Z_0))\) be a weak solution of
	\begin{align}\label{2605311027}
		\partial_t f+v\cdot\nabla_x f
		=
		D_i(a^{ij}D_{j}f)+b^iD_i f+g
		\qquad\text{in }Q_{2r}^+(Z_0),
	\end{align}
	with $f=0$ on $(\bR\times\gamma^-)\cap Q_{2r}(Z_0)$
	in the sense of traces, where \(g\in L^\infty(Q_{2r}^+(Z_0))\).
	
	Then for any $0<\varepsilon<\frac 12$, there exists a constant
	$C=C(d,\lambda,\Lambda,\alpha,\varepsilon)>0$
	such that
	\begin{align}\label{eq:grazing-weighted-zero}
		\begin{aligned}
			&\left\|
			\frac{f}{\mathbf{d} (\,\cdot\,,\bR\times\gamma^0)^{1/2-\varepsilon}}
			\right\|_{L^\infty(Q_{r/2}^+(Z_0))}
			\\
			\le\,&C\left(M+r^{-1}\right)^{1/2-\varepsilon}
			\left(
			\|f\|_{L^\infty(Q_r^+(Z_0))}
			+
			\left(M+r^{-1}\right)^{-2}
			\|g\|_{L^\infty(Q_r^+(Z_0))}
			\right).
		\end{aligned}
	\end{align}
\end{thm}
Theorem~\ref{260129331} is proved at the end of this subsection.

\begin{lemma}[Local barrier near the grazing set]
	\label{260605518}
	For every \(0<\lambda\le\Lambda\), \(\alpha\in(0,1)\),
	\(M_1,M_2\ge0\), and \(\varepsilon\in(0,1/2)\), there exist
	\(r_0\in(0,1)\) and \(C\ge1\), depending only on
	$d,\,\lambda,\,\Lambda,\,\alpha,\,M_1,\,M_2,\,\varepsilon$,
	such that the following holds.
	
	Let
	\[
	a^{ij}\in C^\alpha_{\kin}(Q_{r_0} (0)),
	\qquad
	\bar b^i\in L^\infty(Q_{r_0} (0)),
	\qquad i,j=1,\dots,d.
	\]
	Assume that, for a.e. \(Z\in Q_{r_0} (0)\) and every
	\(\xi\in\bR^d\),
	\[
	\lambda|\xi|^2
	\le
	a^{ij}(Z)\xi_i\xi_j,
	\qquad
	|a^{ij}(Z)|\le \Lambda,
	\]
	and
	\[
	[a]_{C^\alpha_{\kin}(Q_{r_0} (0))}\le M_1,
	\qquad
	\|\bar b\|_{L^\infty(Q_{r_0} (0))}\le M_2.
	\]
	Then there exists a function
	\[
	B\in C(\overline{Q_{r_0}^+(0)})
	\cap
	C^\infty\Bigl(Q_{r_0}^+(0)\cup\bigl((\bR\times\gamma^+)\cap Q_{r_0}(0)\bigr)\Bigr)
	\]
	such that
	\[
	C^{-1}\bigl(|t|+|x|^2+|v|^2\bigr)
	\le
	B(t,x,v)
	\le
	C\bigl(|t|^{1/2}+|x|^{1/3}+|v|\bigr)^{1/2-\varepsilon}
	\quad\text{in }Q_{r_0}^+(0),
	\]
	and
	\[
	(\partial_t+v\cdot\nabla_x)B
	-
	a^{ij}D_{ij}B
	-
	\bar b^iD_iB
	\ge 1
	\quad\text{a.e. in }Q_{r_0}^+(0).
	\]
\end{lemma}

\begin{proof}
	Let
	\[
	\beta:=1-2\varepsilon\in(0,1),
	\qquad
	\mathfrak T:=\partial_t+v\cdot\nabla_x,
	\qquad
	L:=a^{ij}D_{ij}+\bar b^iD_i.
	\]
	Only the symmetric part of \(a\) enters \(a^{ij}D_{ij}\), so we may assume
	that \(a\) is symmetric.
	
	Fix \(r\in(0,1)\) for the moment, and suppose that the hypotheses hold in
	\(Q_r^+(0)\). We shall impose smallness conditions on \(r\), depending only on
	the prescribed constants, and then set \(r_0:=r\). Since
	\(a\in C^\alpha_{\kin}(Q_r^+(0))\), it has a continuous extension to
	\(\overline{Q_r^+(0)}\). We still denote this extension by \(a\). In
	particular, \(a(0)\) is well-defined and satisfies the same ellipticity
	bounds.
	
	\medskip
	\noindent
	\emph{Step 1. The frozen barrier.}
	Set
	\[
	A_{\mathrm{fr}}:=(a^{ij}(0))_{i,j=1}^d.
	\]
	Choose \(P\in SO(d)\) such that
	\[
	Pe_d=\frac{A_{\mathrm{fr}}^{1/2}e_d}{|A_{\mathrm{fr}}^{1/2}e_d|},
	\]
	and define
	\[
	\cS:=(A_{\mathrm{fr}}^{1/2}P)^{-1}.
	\]
	Then
	\[
	\cS A_{\mathrm{fr}}\cS^\top=I_d,
	\qquad
	\Lambda^{-1/2}|\xi|\le |\cS\xi|\le \lambda^{-1/2}|\xi|,\qquad 		(\cS\xi)^d=\frac{\xi^d}{|A_{\mathrm{fr}}^{1/2}e_d|}.
	\]
	In particular, \(\cS\) preserves the sign of the \(d\)-th coordinate.
	
	Let \(H\) be the function from Lemma~\ref{260603514}, and define
	\[
	B_0(x,v):=
	\Big(H\bigl((\cS x)^d,(\cS v)^d\bigr)\Big)^\beta.
	\]
	Put
	\[
	L_0:=a^{ij}(0)D_{ij},
	\qquad
	y:=(\cS x)^d,
	\qquad
	w:=(\cS v)^d.
	\]
	By the choice of \(\cS\),
	\[
	(\mathfrak T-L_0)B_0
	=
	\left[
	(w\partial_y-\partial_{ww})H^\beta
	\right](y,w).
	\]
	Since \(w\partial_yH-\partial_{ww}H=0\),
	\[
	(w\partial_y-\partial_{ww})(H^\beta)
	=
	\beta(1-\beta)H^\beta\frac{|\partial_wH|^2}{H^2}.
	\]
	Hence
	\[
	(\mathfrak T-L_0)B_0>0.
	\]
	
	Set
	\[
	q(y,w):=
	\frac{\eta(y^{-1/3}w)}{y^{1/3}+|w|},
	\]
	where \(\eta\) is as in Lemma~\ref{260603514}. By the estimates for \(H\),
	\[
	\frac{|\partial_wH|}{H}\simeq q,
	\qquad
	\frac{|\partial_{ww}H|}{H}\lesssim q^2.
	\]
	Therefore
	\[
	(\mathfrak T-L_0)B_0\simeq H^\beta q^2.
	\]
	Moreover, by the chain rule and the bounds on \(\cS\),
	\[
	|D_vB_0|\lesssim H^\beta q,
	\qquad
	|D_{vv}B_0|\lesssim H^\beta q^2.
	\]
	Consequently,
	\[
	|D_{vv}B_0|
	\lesssim
	(\mathfrak T-L_0)B_0.
	\]
	
	Since \(\eta\ge1\),
	\[
	q^{-1}\le y^{1/3}+|w|.
	\]
	By the special choice of \(\cS\),
	\[
	y=\frac{x^d}{|A_{\mathrm{fr}}^{1/2}e_d|},
	\qquad
	w=\frac{v^d}{|A_{\mathrm{fr}}^{1/2}e_d|}.
	\]
	Hence
	\[
	q^{-1}
	\lesssim
	|x^d|^{1/3}+|v^d|
	\lesssim
	\mathbf d((t,x,v),0),
	\]
	where the last inequality follows from \eqref{260702722}. Thus
	\[
	|D_vB_0|
	\lesssim
	\mathbf d(\cdot,0)(\mathfrak T-L_0)B_0,
	\qquad
	|D_{vv}B_0|
	\lesssim
	(\mathfrak T-L_0)B_0
	\]
	in \(Q_r^+(0)\).
	
	\medskip
	\noindent
	\emph{Step 2. Perturbing the coefficients.}
	Since
	\[
	(L-L_0)B_0
	=
	(a^{ij}-a^{ij}(0))D_{ij}B_0+\bar b^iD_iB_0,
	\]
	the estimates from Step~1 and the assumptions on \(a\) and \(\bar b\) give
	\[
	|(L-L_0)B_0|
	\le
	C_2
	\bigl(
	M_1\mathbf d(\cdot,0)^\alpha
	+
	M_2\mathbf d(\cdot,0)
	\bigr)
	(\mathfrak T-L_0)B_0
	\]
	in \(Q_r^+(0)\), where
	$C_2=C_2(d,\lambda,\Lambda,\varepsilon)>0$.
	Since \(\mathbf d(\cdot,0)\le r\) in \(Q_r^+(0)\), choose
	\(r_1\in(0,1)\), depending only on
	\(d,\lambda,\Lambda,\alpha,M_1,M_2,\varepsilon\), such that
	\[
	C_2(M_1r_1^\alpha+M_2r_1)\le\frac12.
	\]
	Then, for every \(r\le r_1\),
	\[
	(\mathfrak T-L)B_0
	\ge
	\frac12(\mathfrak T-L_0)B_0
	>0
	\quad\text{in }Q_r^+(0).
	\]
	
	\medskip
	\noindent
	\emph{Step 3. Correction near the outgoing side.}
	Set
	\[
	\kappa:=\frac{2(M_2+1)}{\lambda},
	\qquad
	r_2:=r_1\wedge \frac{\lambda}{4(M_2+1)}.
	\]
	Then \(\kappa r_2\le1/2\). Hence, for \(|v^d|\le r_2\),
	\[
	v^d-\kappa |v^d|^2\ge0
	\quad\text{if }v^d\ge0,
	\]
	and
	\[
	v^d-\kappa |v^d|^2
	\ge
	-2|v^d|
	\quad\text{if }v^d<0.
	\]
	Moreover, since \(|\bar b|\le M_2\) and \(a^{dd}\ge\lambda\),
	\[
	\begin{aligned}
		-L\bigl(v^d-\kappa |v^d|^2\bigr)
		&=
		-\bar b^d+2\kappa a^{dd}+2\kappa\bar b^d v^d \ge
		2\kappa(\lambda-M_2r_2)-M_2
		\ge 2
	\end{aligned}
	\]
	in \(Q_{r_2}^+(0)\).
	
	By Lemma~\ref{260603514},
	\[
	H(y,w)\gtrsim \sqrt{|w|}\,\mathbf 1_{\{w<0\}}.
	\]
	Since \(\cS\) preserves the sign of the \(d\)-th coordinate, there
	exists
	\[
	c_\ast=c_\ast(d,\lambda,\Lambda,\varepsilon)>0
	\]
	such that
	\[
	B_0(x,v)
	\ge
	c_\ast |v^d|^{\beta/2}\mathbf 1_{\{v^d<0\}}.
	\]
	Define
	\[
	B_1:=B_0+\frac{c_\ast}{2}\left(v^d-\kappa |v^d|^2\right).
	\]
	Using the positivity of \((\mathfrak T-L)B_0\) in \(Q_{r_2}^+(0)\) and the
	preceding computation, we get
	\[
	(\mathfrak T-L)B_1
	\ge c_\ast
	\quad\text{in }Q_{r_2}^+(0).
	\]
	Also \(B_1\ge0\) in \(Q_{r_2}^+(0)\). Indeed, this is clear when
	\(v^d\ge0\), while for \(v^d<0\),
	\[
	B_1
	\ge
	c_\ast\left(|v^d|^{\beta/2}-|v^d|\right)\ge0,
	\]
	because \(|v^d|\le r_2<1\) and \(\beta/2<1\).
	
	\medskip
	\noindent
	\emph{Step 4. Adding a coercive correction.}
	Set
	\[
	\phi(t,x,v):=-t+|x|^2+|v|^2.
	\]
	Since \(Q_{r_2}^+(0)\subset Q_1^+(0)\),
	\[
	|(\mathfrak T-L)\phi|
	\le
	1+2|x||v|+2\Lambda d+2M_2|v|
	\le
	C_3
	\quad\text{in }Q_{r_2}^+(0),
	\]
	where $C_3=C_3(d,\Lambda,M_2)\ge1$.
	Put
	\[
	c_1:=\frac{c_\ast}{2C_3},
	\qquad
	\widetilde B:=B_1+c_1\phi.
	\]
	Then
	\[
	(\mathfrak T-L)\widetilde B
	\ge
	c_\ast-c_1C_3
	=
	\frac{c_\ast}{2}
	\quad\text{in }Q_{r_2}^+(0).
	\]
	Since \(B_1\ge0\) and \(t\le0\) in \(Q_{r_2}^+(0)\),
	\[
	\widetilde B
	\ge
	c_1\phi
	=
	c_1\bigl(|t|+|x|^2+|v|^2\bigr)
	\quad\text{in }Q_{r_2}^+(0).
	\]
	
	It remains to prove the upper bound. By Lemma~\ref{260603514},
	\[
	H(y,w)\lesssim y^{1/6}+|w|^{1/2}.
	\]
	Thus
	\[
	B_0=H(y,w)^\beta
	\lesssim
	y^{\beta/6}+|w|^{\beta/2}.
	\]
	Using the bounds on \(\cS\), the definition of \(B_1\), and
	\(Q_{r_2}^+(0)\subset Q_1^+(0)\), we obtain
	\[
	\widetilde B
	\le
	C
	\left(
	|x^d|^{\beta/6}
	+
	|v^d|^{\beta/2}
	+
	|v^d|
	+
	|v^d|^2
	+
	|t|
	+
	|x|^2
	+
	|v|^2
	\right).
	\]
	Since \(\beta/2\in(0,1)\) and \(Q_{r_2}^+(0)\subset Q_1^+(0)\), the right-hand
	side is bounded by
	\[
	C\bigl(|t|^{1/2}+|x|^{1/3}+|v|\bigr)^{\beta/2};
	\]
	recall that
	$\frac{\beta}{2}=\frac12-\varepsilon$.
	
	Finally set
	\[
	r_0:=r_2,
	\qquad
	B:=\frac{2}{c_\ast}\widetilde B.
	\]
	Then
	\[
	(\mathfrak T-L)B\ge1
	\quad\text{a.e. in }Q_{r_0}^+(0),
	\]
	and, after changing the constant \(C\),
	\[
	C^{-1}\bigl(|t|+|x|^2+|v|^2\bigr)
	\le
	B(t,x,v)
	\le
	C\bigl(|t|^{1/2}+|x|^{1/3}+|v|\bigr)^{1/2-\varepsilon}
	\]
	in \(Q_{r_0}^+(0)\).
	
	The asserted continuity and smoothness follow from Lemma~\ref{260603514} and
	from the polynomial correction terms. This completes the proof.
\end{proof}

\begin{lemma}[Smooth approximation near the half-space boundary]
	\label{lem:smooth-approximation-half-space}
	Let
	$Z_0\in \bR\times\overline{\bR^d_+}\times\bR^d$.
	Assume that \(a^{ij}\), \(b^i\), \(f\), and \(g\) satisfy, in
	\(Q_2(Z_0)\) and \(Q_2^+(Z_0)\), the equation and structural assumptions appearing in
	Theorem~\ref{260129331} with \(r=1\).
	
	Then there exist
	\[
	a_n^{ij},\, b_n^i,\, g_n\in C^\infty(\bR^{1+2d})
	\]
	and
	\[
	f_n\in
	C\big(\overline{Q_{3/2}^+(Z_0)}\big)
	\cap
	C^\infty\big(Q_{3/2}^+(Z_0)\big)
	\]
	such that \(a_n^{ij}\) and \(b_n^i\) satisfy the structural bounds
	\eqref{260605507} and \eqref{2607041007} in \(Q_{3/2}(Z_0)\) with the same
	constants, and
	\[
	(\partial_t+v\cdot\nabla_x)f_n
	=
	D_i(a_n^{ij}D_jf_n)+b_n^iD_if_n+g_n
	\quad\text{in }Q_{3/2}^+(Z_0).
	\]
	Moreover,
	\[
	f_n=0
	\qquad\text{on }(\bR\times\gamma^-)\cap Q_{3/2}(Z_0)
	\]
	in the sense of traces,
	\[ f_n\to f, \qquad \nabla_v f_n\to \nabla_v f \quad\text{in }L^2(Q_{3/2}^+(Z_0)), \]
	and
	\[
	\|g_n\|_{L^\infty(Q_{3/2}^+(Z_0))}
	\le
	\|g\|_{L^\infty(Q_2^+(Z_0))}.
	\]
\end{lemma}

\begin{proof}
	Choose \(\zeta\in C_c^\infty(\bR^{1+2d})\) such that
	\[
	\zeta\ge0,
	\qquad
	\iiint_{\bR^{1+2d}}\zeta\dd Z=1,\qquad \operatorname{supp}\zeta
	\subset
	\{(t,x,v):0<t<1,\ |x|<1,\ |v|<1\}.
	\]
	
	We extend the coefficients to \(\bR^{1+2d}\) by setting
	\[
	a^{ij}:=\Lambda\delta^{ij},
	\qquad
	b^i:=0
	\quad\text{outside }Q_2(Z_0),
	\]
	and extend \(g\) by zero outside \(Q_2^+(Z_0)\). We keep the same notation
	for the extended functions. We do not claim any global regularity for the extended coefficient \(a\).
	
	For \(n\in\mathbb N\), set
	\[
	\varepsilon_n:=\min\left\{\frac1{8n},\ (1+|v_0|)^{-1/2},\ (1+|t_0|)^{-1}\right\},
	\]
	and define
	\[
	\widehat a_n^{ij}:=(a^{ij})^{(\varepsilon_n),\star},
	\qquad
	\widehat b_n^i:=(b^i)^{(\varepsilon_n),\star},
	\qquad
	\widehat g_n:=g^{(\varepsilon_n),\star},
	\]
	where the mollification is understood in the kinetic sense defined in
	\eqref{2606061226}. By the support condition on \(\zeta\), the choice
	of \(\varepsilon_n\), and the elementary cylinder inclusion, we have
	\[
	Y^{-1}\circ Z\in Q_2(Z_0)
	\]
	whenever
	\[
	Z\in Q_{7/4}(Z_0),
	\qquad
	Y\in\operatorname{supp}\zeta_{\varepsilon_n,\star}.
	\]
	Consequently, on \(Q_{7/4}(Z_0)\), the mollifications
	\(\widehat a_n^{ij}\), \(\widehat b_n^i\), and \(\widehat g_n\) only use the
	values of the extended functions in \(Q_2(Z_0)\).
	
	By Lemma~\ref{2606061211},
	\[
	\widehat a_n^{ij},\,\widehat b_n^i,\,\widehat g_n
	\in C^\infty(\bR^{1+2d}),
	\]
	and
	\[
	\widehat a_n^{ij}\to a^{ij},
	\qquad
	\widehat b_n^i\to b^i,
	\qquad
	\widehat g_n\to g
	\quad\text{a.e. in }Q_{7/4}^+(Z_0).
	\]
	Moreover,
	\[
	\widehat g_n\to g
	\quad\text{in }L^2(Q_{7/4}^+(Z_0)).
	\]
	Since \(g\) is extended by zero outside \(Q_2^+(Z_0)\), the
	\(L^\infty\)-contraction property gives
	\[
	\|\widehat g_n\|_{L^\infty(Q_{7/4}^+(Z_0))}
	\le
	\|g\|_{L^\infty(Q_2^+(Z_0))}.
	\]
	
	We next verify the structural bounds in \(Q_{7/4}(Z_0)\). Since the
	mollification only samples \(Q_2(Z_0)\), the positivity and unit mass of
	\(\zeta_{\varepsilon_n,\star}\) imply that, for every
	\(Z\in Q_{7/4}(Z_0)\) and every \(\xi\in\bR^d\),
	\[
	\lambda|\xi|^2
	\le
	\widehat a_n^{ij}(Z)\xi_i\xi_j,
	\qquad
	|\widehat a_n^{ij}(Z)|\le\Lambda.
	\]
	Since \(D_v\) commutes with the kinetic mollification in the weak sense, and
	since the mollification only samples \(Q_2(Z_0)\), we have $D_v\widehat a_n=(D_va)^{(\varepsilon_n),\star}$ in
	\(Q_{7/4}(Z_0)\).
	Hence the \(L^\infty\)-contraction property gives
	\[
	\|\nabla_v\widehat a_n\|_{L^\infty(Q_{7/4}(Z_0))}
	+
	\|\widehat b_n\|_{L^\infty(Q_{7/4}(Z_0))}
	\le M_2.
	\]
	
	We also check the kinetic H\"older bound. Let
	\(Z_1,Z_2\in Q_{7/4}(Z_0)\). If
	\(Y\in\operatorname{supp}\zeta_{\varepsilon_n,\star}\), then
	\[
	Y^{-1}\circ Z_1,\,
	Y^{-1}\circ Z_2
	\in Q_2(Z_0).
	\]
	By the left-invariance of \(\mathbf d\),
	\[
	\mathbf d(Y^{-1}\circ Z_1,Y^{-1}\circ Z_2)
	=
	\mathbf d(Z_1,Z_2).
	\]
	Therefore
	\[
	\begin{aligned}
		|\widehat a_n^{ij}(Z_1)-\widehat a_n^{ij}(Z_2)|
		&\le
		\iiint_{\bR^{1+2d}}
		\zeta_{\varepsilon_n,\star}(Y)
		\left|
		a^{ij}(Y^{-1}\circ Z_1)
		-
		a^{ij}(Y^{-1}\circ Z_2)
		\right|\dd Y \\
		&\le
		[a^{ij}]_{C^\alpha_{\kin}(Q_2(Z_0))}
		\mathbf d(Z_1,Z_2)^\alpha \\
		&\le
		M_1\mathbf d(Z_1,Z_2)^\alpha.
	\end{aligned}
	\]
	Thus
	$[\widehat a_n^{ij}]_{C^\alpha_{\kin}(Q_{7/4}(Z_0))}
	\le M_1$.
	
	By Lemma~\ref{2605311047}, applied with \(\kappa=7/8\), there exists a weak
	solution \(\widehat f_n\) of
	\[
	\begin{cases}
		(\partial_t+v\cdot\nabla_x)\widehat f_n
		=
		D_i(\widehat a_n^{ij}D_j\widehat f_n)
		+\widehat b_n^iD_i\widehat f_n+\widehat g_n
		&\quad\text{in }Q_{7/4}^+(Z_0),\\
		\widehat f_n=0
		&\quad\text{on }(\bR\times\gamma^-)\cap Q_{7/4}(Z_0),
	\end{cases}
	\]
	such that
	\[
	\widehat f_n\to f,
	\qquad
	\nabla_v\widehat f_n\to\nabla_v f
	\quad\text{in }L^2(Q_{7/4}^+(Z_0)).
	\]
	
	We next record the interior smoothness of \(\widehat f_n\). Let
	\[
	A_n:=\frac{\widehat a_n+\widehat a_n^\top}{2}.
	\]
	Then \(A_n\) is smooth and uniformly positive definite. Let
	\(\sigma_n:=A_n^{1/2}\), and set
	\[
	X_{n,k}:=\sigma_n^{ik}D_i,
	\qquad k=1,\dots,d.
	\]
	Since the antisymmetric part of \(\widehat a_n\) contributes only to
	first-order terms, the equation can be written as
	\[
	(\partial_t+v\cdot\nabla_x)\widehat f_n
	=
	\sum_{k=1}^d X_{n,k}^2\widehat f_n
	+
	\widetilde b_n^iD_i\widehat f_n
	+
	\widehat g_n
	\quad\text{in }\cD'\big(Q_{7/4}^+(Z_0)\cap\{t<t_0\}\big),
	\]
	for some smooth vector field \(\widetilde b_n\). The vector fields
	\[
	X_0:=\partial_t+v\cdot\nabla_x-\widetilde b_n^iD_i,
	\qquad
	X_{n,1},\dots,X_{n,d}
	\]
	satisfy H\"ormander's bracket condition in the interior. Hence, by
	H\"ormander's hypoellipticity theorem \cite{Hormander1967},
	\[
	\widehat f_n\in
	C^\infty_{\mathrm{loc}}\big(Q_{7/4}^+(Z_0)\cap\{t<t_0\}\big).
	\]
	
	We now perform the final time shift. Choose \(\tau_n\downarrow0\) so small
	that
	\[
	Q_{3/2}(Z_0)+(-\tau_n,0,0)\subset Q_{7/4}(Z_0)
	\]
	and
	\[
	\begin{aligned}
		&\|f(t-\tau_n,x,v)-f(t,x,v)\|_{L^2(Q_{3/2}^+(Z_0))} \\
		&\quad
		+
		\|\nabla_v f(t-\tau_n,x,v)-\nabla_v f(t,x,v)\|_{L^2(Q_{3/2}^+(Z_0))}
		\le 2^{-n}.
	\end{aligned}
	\]
	This is possible by the translation continuity of \(L^2\), since
	\(Q_{3/2}(Z_0)+(-\tau_n,0,0)\subset Q_{7/4}(Z_0)\).
	
	Define
	\begin{align*}
		\begin{gathered}
			f_n(t,x,v):=\widehat f_n(t-\tau_n,x,v),\qquad g_n(t,x,v):=\widehat g_n(t-\tau_n,x,v),\\
			a_n^{ij}(t,x,v):=\widehat a_n^{ij}(t-\tau_n,x,v),
			\qquad
			b_n^i(t,x,v):=\widehat b_n^i(t-\tau_n,x,v).
		\end{gathered}
	\end{align*}
	Then
	$a_n^{ij},\,b_n^i,\,g_n\in C^\infty(\bR^{1+2d})$,
	and the structural bounds \eqref{260605507} and \eqref{2607041007} hold in
	\(Q_{3/2}(Z_0)\) with the same constants. Moreover,
	\[
	\|g_n\|_{L^\infty(Q_{3/2}^+(Z_0))}
	\le
	\|g\|_{L^\infty(Q_2^+(Z_0))}.
	\]
	
	Since the shifted cylinder lies strictly below the terminal time \(t=t_0\),
	the preceding interior regularity gives
	$f_n\in C^\infty(Q_{3/2}^+(Z_0))$.
	The shifted equation is
	$(\partial_t+v\cdot\nabla_x)f_n
	=
	D_i(a_n^{ij}D_jf_n)+b_n^iD_if_n+g_n
	\quad\text{in }Q_{3/2}^+(Z_0)$.
	Since \(\tau_n>0\), the equation holds on all of \(Q_{3/2}^+(Z_0)\),
	including the terminal slice \(t=t_0\), as required by the comparison
	principle of Lemma~\ref{lem:classical-kinetic-comparison}.
	Also, since \(\widehat f_n=0\) on
	\((\bR\times\gamma^-)\cap Q_{7/4}(Z_0)\), the inclusion above implies
	$f_n=0$ on $(\bR\times\gamma^-)\cap Q_{3/2}(Z_0)$
	in the sense of traces.
	
	It remains to check the convergence. We have
	\[
	\begin{aligned}
		\|f_n-f\|_{L^2(Q_{3/2}^+(Z_0))}
		&\le
		\|\widehat f_n-f\|_{L^2(Q_{7/4}^+(Z_0))}
		+
		\|f(t-\tau_n,x,v)-f(t,x,v)\|_{L^2(Q_{3/2}^+(Z_0))},
	\end{aligned}
	\]
	and similarly,
	\[
	\begin{aligned}
		\|\nabla_vf_n-\nabla_vf\|_{L^2(Q_{3/2}^+(Z_0))}
		\le\,&
		\|\nabla_v\widehat f_n-\nabla_vf\|_{L^2(Q_{7/4}^+(Z_0))}
		\\
		&\quad+
		\|\nabla_v f(t-\tau_n,x,v)-\nabla_v f(t,x,v)\|_{L^2(Q_{3/2}^+(Z_0))}.
	\end{aligned}
	\]
	Therefore,
	\[
	f_n\to f,
	\qquad
	\nabla_vf_n\to\nabla_vf
	\quad\text{in }L^2(Q_{3/2}^+(Z_0)).
	\]
	
	Finally, we prove continuity up to the kinetic boundary. Note that
	\[
	\overline{Q_{3/2}^+(Z_0)}+(-\tau_n,0,0)
	\subset \overline{Q_{7/4}^+(Z_0)}\cap\{t<t_0\}.
	\]
	Since \(\widehat f_n=0\) on \((\bR\times\gamma^-)\cap Q_{7/4}(Z_0)\)
	and \(\widehat g_n\in L^\infty\), Theorem~\ref{cor:linfty-zero-inflow}
	gives \(\widehat f_n\in L^\infty(Q_{13/8}^+(Z_0))\).  Hence
	Theorem~\ref{cor:scaled-boundary-holder-estimate} applies, with uniform local
	cylinders, to \(\widehat f_n\) in a neighborhood of the shifted set. A finite
	covering argument yields that \(\widehat f_n\) is continuous up to the kinetic
	boundary on
	$\overline{Q_{3/2}^+(Z_0)}+(-\tau_n,0,0)$.
	Therefore
	$f_n\in C\big(\overline{Q_{3/2}^+(Z_0)}\big)$.
	This completes the proof.
\end{proof}

\begin{lemma}[Smooth approximation, interior version]
	\label{lem:smooth-approximation-interior}
	Let \(Z_0\in\bR^{1+2d}\), and assume that \(a^{ij},b^i,f,g\) satisfy in
	\(Q_2(Z_0)\) the structural assumptions and the equation of
	Lemma~\ref{lem:smooth-approximation-half-space}, the equation now being
	imposed on the full cylinder and no boundary condition being imposed.
	Then the conclusions of Lemma~\ref{lem:smooth-approximation-half-space}
	hold with \(Q_{3/2}^+(Z_0)\), \(Q_{7/4}^+(Z_0)\), \(Q_2^+(Z_0)\)
	replaced throughout by \(Q_{3/2}(Z_0)\), \(Q_{7/4}(Z_0)\),
	\(Q_2(Z_0)\), and \(f_n\in C^\infty(Q_{3/2}(Z_0))\) solves the
	approximate equation on all of \(Q_{3/2}(Z_0)\).
\end{lemma}

\begin{proof}
	Apply the proof of Lemma~\ref{lem:smooth-approximation-half-space} with
	\(\cD\) a ball containing the \(x\)-projection of \(Q_2(Z_0)\),
	which is admissible by Remark~\ref{260630607}; then
	\((\bR\times\gamma)\cap Q_{7/4}(Z_0)=\emptyset\), so the trace
	condition is vacuous and every step is carried out on full cylinders.  As
	there, the time shift \(\tau_n\downarrow0\) yields smoothness up to the
	terminal slice \(t=t_0\).
\end{proof}

\begin{lemma}\label{260705834}
	Let \(\varepsilon\in(0,1/2)\), and let \(f_n\) be as in
	Lemma~\ref{lem:smooth-approximation-half-space}. Then
	\begin{align}\label{2607041155}
		\left\|
		\frac{f_n}{\mathbf d(\,\cdot\,,\bR\times\gamma^0)^{1/2-\varepsilon}}
		\right\|_{L^\infty(Q_1^+(Z_0))}
		\le
		C
		\left(
		\|f_n\|_{L^2(Q_{3/2}^+(Z_0))}
		+
		\|g_n\|_{L^\infty(Q_{3/2}^+(Z_0))}
		\right),
	\end{align}
	where
	$C=C(d,\lambda,\Lambda,\alpha,M,\varepsilon)>0$.
\end{lemma}

\begin{proof}
	Put
	\[
	\theta:=\frac12-\varepsilon.
	\]
	Fix a dimensional constant \(C_d\ge1\) such that, whenever
	$\bar b^j=D_i a^{ij}+b^j$,
	one has
	\[
	|\bar b|
	\le
	C_d\bigl(|\nabla_v a|+|b|\bigr).
	\]
	We apply Lemma~\ref{260605518} with the parameters
	$M_1:=M^\alpha$ and $M_2:=C_dM$.
	Let \(r_0\in(0,1/16]\) be the corresponding radius in
	Lemma~\ref{260605518}. We may impose the additional restriction
	\(r_0\le1/16\) by shrinking the radius in that lemma if necessary. Thus
	\(r_0\) depends only on
	$d,\,\lambda,\,\Lambda,\,\alpha,\,M,\,\varepsilon$.
	Choose \(r_\ast\in(0,1)\), depending only on \(r_0\), such that
	\[
	16r_\ast\le \frac{r_0}{2},
	\qquad
	r_0+2^{10/3}r_\ast\le \frac18.
	\]
	All constants below may change from line to line, but depend only on
	$d,\,\lambda,\,\Lambda,\,\alpha,\,M,\,\varepsilon$.
	
	\medskip
	\noindent
	\emph{Step 1. Spatially shifted smooth solutions.}
	For \(\delta>0\), define
	\[
	f_{n,\delta}(t,x,v):=f_n(t,x+\delta e_d,v),
	\]
	and define \(a_{n,\delta}^{ij}\), \(b_{n,\delta}^i\), and
	\(g_{n,\delta}\) analogously. For all sufficiently small \(\delta>0\) such
	that
	\[
	\overline{Q_{9/8}^+(Z_0)}+(0,\delta e_d,0)
	\subset Q_{5/4}^+(Z_0),
	\]
	the shifted functions are well-defined on \(Q_{9/8}^+(Z_0)\), and
	\(f_{n,\delta}\) is smooth in a neighborhood of
	\(\overline{Q_{9/8}^+(Z_0)}\). Since the shift is only in the \(x\)-variable
	by a constant vector, the kinetic transport operator is unchanged. Hence
	\[
	(\partial_t+v\cdot\nabla_x)f_{n,\delta}
	=
	D_i(a_{n,\delta}^{ij}D_jf_{n,\delta})
	+
	b_{n,\delta}^iD_if_{n,\delta}
	+
	g_{n,\delta}
	\quad\text{in }Q_{9/8}^+(Z_0).
	\]
	Equivalently, in the interior,
	\[
	(\partial_t+v\cdot\nabla_x)f_{n,\delta}
	-
	a_{n,\delta}^{ij}D_{ij}f_{n,\delta}
	-
	\bar b_{n,\delta}^jD_jf_{n,\delta}
	=
	g_{n,\delta},
	\]
	where
	\[
	\bar b_{n,\delta}^j:=D_i a_{n,\delta}^{ij}+b_{n,\delta}^j .
	\]
	The coefficient \(a_{n,\delta}\) satisfies the same ellipticity and
	boundedness assumptions as \(a_n\), and
	\[
	[a_{n,\delta}]_{C^\alpha_{\kin}(Q_{9/8}(Z_0))}\le M^\alpha.
	\]
	Furthermore, by the choice of \(C_d\),
	\[
	\|\bar b_{n,\delta}\|_{L^\infty}
	\le
	C_d\left(
	\|\nabla_v a_{n,\delta}\|_{L^\infty}
	+
	\|b_{n,\delta}\|_{L^\infty}
	\right)
	\le C_dM.
	\]
	Thus \(a_{n,\delta}\) and \(\bar b_{n,\delta}\) satisfy the assumptions of
	Lemma~\ref{260605518} with the parameters $M_1:=M^\alpha$ and $M_2:=C_dM$.
	
	Set
	\[
	m_\delta:=
	\sup_{(\bR\times(\gamma^-\cup\gamma^0))\cap Q_{9/8}(Z_0)}
	|f_{n,\delta}|.
	\]
	We claim that
	\begin{equation}\label{eq:shifted-grazing-bound}
		|f_{n,\delta}(Z)|
		\le
		m_\delta
		+
		C\,\mathbf d(Z,\bR\times\gamma^0)^\theta
		\left(
		\|f_{n,\delta}\|_{L^\infty(Q_{9/8}^+(Z_0))}
		+
		\|g_{n,\delta}\|_{L^\infty(Q_{9/8}^+(Z_0))}
		\right)
	\end{equation}
	for every \(Z\in Q_1^+(Z_0)\).
	
	\medskip
	\noindent
	\emph{Step 2. Localization near the grazing set.}
	Fix \(Z\in Q_1^+(Z_0)\), and set
	\[
	\rho_Z:=\mathbf d(Z,\bR\times\gamma^0).
	\]
	If \(\rho_Z\ge r_\ast\), then
	\[
	|f_{n,\delta}(Z)|
	\le
	\|f_{n,\delta}\|_{L^\infty(Q_{9/8}^+(Z_0))}
	\le
	r_\ast^{-\theta}\rho_Z^\theta
	\|f_{n,\delta}\|_{L^\infty(Q_{9/8}^+(Z_0))}.
	\]
	Thus \eqref{eq:shifted-grazing-bound} follows in this case.
	
	It remains to consider the case
	$\rho_Z<r_\ast$.
	Write \(Z=(t,x,v)\).  The point \((t,(x',0),(v',0))\) belongs to
	\(\bR\times\gamma^0\) and, by \eqref{260702722}, lies at kinetic
	distance at most \(|x^d|^{1/3}+|v^d|\le4\rho_Z\) from \(Z\).  Since it
	has the same time coordinate as \(Z\), Lemma~\ref{260127226}(2) gives
	\[
	Q_{8\rho_Z}(Z)\cap(\bR\times\gamma^0)\neq\emptyset.
	\]
	The smaller radius \(2\rho_Z\) would not do: as \(x^d\downarrow0\), one
	has \(\rho_Z\to|v^d|/2\), whereas any
	\(\bar Z\in Q_{2\rho_Z}(Z)\cap(\bR\times\gamma^0)\) would have
	to satisfy \(|v^d|\le|\bar v-v|<2\rho_Z\).
	Hence Lemma~\ref{260515154}, applied with radius \(8\rho_Z\), gives a point
	\[
	Y_Z\in\bR\times\gamma^0
	\]
	such that
	\[
	Q_{8\rho_Z}(Z)\subset Q_{16\rho_Z}(Y_Z),
	\qquad
	Q_R(Y_Z)\subset Q_{R+2^{10/3}\rho_Z}(Z)
	\quad\text{for every }R>0.
	\]
	In particular,
	\[
	Z\in Q_{16\rho_Z}(Y_Z).
	\]
	Since
	$16\rho_Z<16r_\ast\le \frac{r_0}{2}$,
	we have $Z\in Q_{r_0/2}(Y_Z)$.
	Moreover,
	\[
	Q_{r_0}(Y_Z)
	\subset
	Q_{r_0+2^{10/3}\rho_Z}(Z)
	\subset
	Q_{1/8}(Z).
	\]
	Since \(Z\in Q_1(Z_0)\), Lemma~\ref{lem:cylinder-inclusion} gives
	$Q_{1/8}(Z)\subset Q_{9/8}(Z_0)$.
	Therefore
	\[
	Q_{r_0}^+(Y_Z)\subset Q_{9/8}^+(Z_0).
	\]
	
	Set
	\[
	\Xi_Z:=Y_Z^{-1}\circ Z.
	\]
	Then
	$\Xi_Z\in Q_{16\rho_Z}(0)\subset Q_{r_0/2}(0)$,
	and the remaining part of the proof of
	\eqref{eq:shifted-grazing-bound} is reduced to a barrier comparison in
	\(Q_{r_0}^+(0)\).
	
	\medskip
	\noindent
	\emph{Step 3. Barrier comparison in \(Q_{r_0}^+(0)\).}
	We continue to work in the case \(\rho_Z<r_\ast\). For
	\(\Xi\in Q_{r_0}^+(0)\), define
	\[
	\tilde f(\Xi):=f_{n,\delta}(Y_Z\circ \Xi),
	\qquad
	\tilde g(\Xi):=g_{n,\delta}(Y_Z\circ \Xi),
	\]
	and
	\[
	\tilde a^{ij}(\Xi):=a_{n,\delta}^{ij}(Y_Z\circ \Xi),
	\qquad
	\tilde b^i(\Xi):=\bar b_{n,\delta}^i(Y_Z\circ \Xi).
	\]
	Since \(Y_Z\in\bR\times\gamma^0\), left translation by \(Y_Z\)
	preserves the half-space boundary:
	\[
	Y_Z\circ Q_{r_0}^+(0)=Q_{r_0}^+(Y_Z).
	\]
	By Step~2,
	$Q_{r_0}^+(Y_Z)\subset Q_{9/8}^+(Z_0)$.
	Hence the coefficients \(\tilde a\) and \(\tilde b\) are defined on
	\(Q_{r_0}(0)\) and satisfy the assumptions of Lemma~\ref{260605518} with
	the same parameters
	$M_1:=M^\alpha$ and $M_2:=C_dM$.
	
	The kinetic transport operator is invariant under left translation. Hence
	\[
	(\partial_t+v\cdot\nabla_x)\tilde f
	-
	\tilde a^{ij}D_{ij}\tilde f
	-
	\tilde b^iD_i\tilde f
	=
	\tilde g
	\quad\text{in }Q_{r_0}^+(0).
	\]
	
	By Lemma~\ref{260605518}, applied to \(\tilde a\) and \(\tilde b\), there
	exists a barrier
	\[
	B\in C(\overline{Q_{r_0}^+(0)})
	\cap
	C^\infty\Bigl(Q_{r_0}^+(0)\cup\bigl((\bR\times\gamma^+)\cap Q_{r_0}(0)\bigr)\Bigr)
	\]
	such that
	\[
	(\partial_t+v\cdot\nabla_x)B
	-
	\tilde a^{ij}D_{ij}B
	-
	\tilde b^iD_iB
	\ge 1
	\quad\text{in }Q_{r_0}^+(0),
	\]
	and
	\[
	0\le B(\Xi)\le C\,\mathbf d(\Xi,0)^\theta
	\quad\text{in }Q_{r_0}^+(0).
	\]
	
	Define
	\[
	\cA_{r_0}
	:=
	\overline{Q_{r_0}^+(0)}
	\cap
	\left(
	\{t=-r_0^2\}\cup\{|x|=r_0^3\}\cup\{|v|=r_0\}
	\right).
	\]
	The lower bound in Lemma~\ref{260605518} gives, on \(\cA_{r_0}\),
	\[
	B\ge c_0r_0^6,
	\]
	where $c_0=c_0(d,\lambda,\Lambda,\alpha,M,\varepsilon)>0$.
	
	Let
	\[
	K_\delta
	:=
	C_*
	\left(
	\|f_{n,\delta}\|_{L^\infty(Q_{9/8}^+(Z_0))}
	+
	\|g_{n,\delta}\|_{L^\infty(Q_{9/8}^+(Z_0))}
	\right),\qquad C_*:=1+(c_0r_0^6)^{-1}.
	\]
	Since
	\[
	Y_Z\circ Q_{r_0}^+(0)=Q_{r_0}^+(Y_Z)
	\subset Q_{9/8}^+(Z_0),
	\]
	we have
	\[
	K_\delta\ge \|\tilde g\|_{L^\infty(Q_{r_0}^+(0))},
	\qquad
	K_\delta c_0r_0^6
	\ge
	\|\tilde f\|_{L^\infty(Q_{r_0}^+(0))}.
	\]
	
	Define
	\[
	W_\pm:=K_\delta B+m_\delta\pm\tilde f
	\quad\text{in }Q_{r_0}^+(0).
	\]
	Since
	\[
	\left|
	(\partial_t+v\cdot\nabla_x)\tilde f
	-
	\tilde a^{ij}D_{ij}\tilde f
	-
	\tilde b^iD_i\tilde f
	\right|
	=
	|\tilde g|,
	\]
	the choice of \(K_\delta\) implies
	\[
	(\partial_t+v\cdot\nabla_x)W_\pm
	-
	\tilde a^{ij}D_{ij}W_\pm
	-
	\tilde b^iD_iW_\pm
	\ge0
	\quad\text{in }Q_{r_0}^+(0).
	\]
	
	We check the comparison boundary. If
	\[
	\Xi\in(\bR\times(\gamma^-\cup\gamma^0))\cap Q_{r_0}(0),
	\]
	then, since \(Y_Z\in\bR\times\gamma^0\),
	\[
	Y_Z\circ\Xi
	\in
	(\bR\times(\gamma^-\cup\gamma^0))\cap Q_{r_0}(Y_Z)
	\subset
	(\bR\times(\gamma^-\cup\gamma^0))\cap Q_{9/8}(Z_0).
	\]
	Therefore, by the definition of \(m_\delta\),
	\[
	W_\pm(\Xi)\ge0
	\qquad
	\text{on }
	(\bR\times(\gamma^-\cup\gamma^0))\cap Q_{r_0}(0).
	\]
	On the artificial boundary \(\cA_{r_0}\), the lower bound
	\(B\ge c_0r_0^6\) and the definition of \(K_\delta\) give
	\[
	W_\pm
	\ge
	K_\delta c_0r_0^6
	-
	\|\tilde f\|_{L^\infty(Q_{r_0}^+(0))}
	\ge0.
	\]
	Thus, by the classical kinetic comparison principle in \(Q_{r_0}^+(0)\)
	(Lemma~\ref{lem:classical-kinetic-comparison}), with comparison boundary
	\[
	\cA_{r_0}
	\cup
	\bigl((\bR\times(\gamma^-\cup\gamma^0))\cap Q_{r_0}(0)\bigr),
	\]
	we obtain
	\[
	W_\pm\ge0
	\quad\text{in }Q_{r_0}^+(0).
	\]
	
	Evaluating at
	\[
	\Xi_Z=Y_Z^{-1}\circ Z\in Q_{r_0/2}^+(0),
	\]
	we get
	\[
	|f_{n,\delta}(Z)|
	=
	|\tilde f(\Xi_Z)|
	\le
	m_\delta+K_\delta B(\Xi_Z).
	\]
	Since \(\Xi_Z\in Q_{16\rho_Z}(0)\), Lemma~\ref{260127226}(1) gives
	\[
	\mathbf d(\Xi_Z,0)
	=
	\mathbf d(Z,Y_Z)
	\le 16\rho_Z.
	\]
	Therefore
	\[
	B(\Xi_Z)
	\le
	C\,\mathbf d(\Xi_Z,0)^\theta
	\le
	C\,\rho_Z^\theta
	=
	C\,\mathbf d(Z,\bR\times\gamma^0)^\theta.
	\]
	Consequently,
	\[
	|f_{n,\delta}(Z)|
	\le
	m_\delta
	+
	C\,\mathbf d(Z,\bR\times\gamma^0)^\theta
	\left(
	\|f_{n,\delta}\|_{L^\infty(Q_{9/8}^+(Z_0))}
	+
	\|g_{n,\delta}\|_{L^\infty(Q_{9/8}^+(Z_0))}
	\right).
	\]
	This proves \eqref{eq:shifted-grazing-bound} in the case
	\(\rho_Z<r_\ast\), and therefore completes the proof of the shifted estimate.
	
	\medskip
	\noindent
	\emph{Step 4. Letting the shift parameter vanish.}
	We now let \(\delta\downarrow0\) in \eqref{eq:shifted-grazing-bound}.
	
	Since $f_n\in C(\overline{Q_{3/2}^+(Z_0)})$
	and
	$f_n=0$ on $(\bR\times\gamma^-)\cap Q_{3/2}(Z_0)$,
	we also have
	\[
	f_n=0
	\qquad\text{on }(\bR\times\gamma^0)\cap Q_{3/2}(Z_0),
	\]
	because \(\gamma^0\subset\overline{\gamma^-}\).
	
	Moreover, \(f_n\) is uniformly continuous on
	\(\overline{Q_{5/4}^+(Z_0)}\). If
	\[
	P=(t,x,v)\in(\bR\times(\gamma^-\cup\gamma^0))\cap Q_{9/8}(Z_0),
	\]
	then
	\[
	P_\delta:=(t,x+\delta e_d,v)
	\in Q_{5/4}^+(Z_0)
	\]
	for all sufficiently small \(\delta>0\). Hence
	\[
	|f_{n,\delta}(P)|
	=
	|f_n(P_\delta)|
	=
	|f_n(P_\delta)-f_n(P)|.
	\]
	Taking the supremum over such \(P\), the uniform continuity of \(f_n\) gives
	\[
	m_\delta\to0
	\qquad\text{as }\delta\downarrow0.
	\]
	
	Next, by the choice of \(\delta\),
	\[
	\|f_{n,\delta}\|_{L^\infty(Q_{9/8}^+(Z_0))}
	\le
	\|f_n\|_{L^\infty(Q_{5/4}^+(Z_0))}
	\]
	and
	\[
	\|g_{n,\delta}\|_{L^\infty(Q_{9/8}^+(Z_0))}
	\le
	\|g_n\|_{L^\infty(Q_{5/4}^+(Z_0))}
	\le
	\|g_n\|_{L^\infty(Q_{3/2}^+(Z_0))}.
	\]
	Since \(f_n\in C(\overline{Q_{3/2}^+(Z_0)})\),
	$f_{n,\delta}(Z)\to f_n(Z)$ for every $Z\in Q_1^+(Z_0)$.
	Letting \(\delta\downarrow0\) in \eqref{eq:shifted-grazing-bound}, we obtain
	\[
	|f_n(Z)|
	\le
	C\,\mathbf d(Z,\bR\times\gamma^0)^\theta
	\left(
	\|f_n\|_{L^\infty(Q_{5/4}^+(Z_0))}
	+
	\|g_n\|_{L^\infty(Q_{3/2}^+(Z_0))}
	\right)
	\]
	for every \(Z\in Q_1^+(Z_0)\). Therefore,
	\[
	\left\|
	\frac{f_n}{\mathbf d(\,\cdot\,,\bR\times\gamma^0)^\theta}
	\right\|_{L^\infty(Q_1^+(Z_0))}
	\le
	C
	\left(
	\|f_n\|_{L^\infty(Q_{5/4}^+(Z_0))}
	+
	\|g_n\|_{L^\infty(Q_{3/2}^+(Z_0))}
	\right).
	\]
	
	\medskip
	\noindent
	\emph{Step 5. Replacing the \(L^\infty\)-norm by an \(L^2\)-norm.}
	By the local \(L^2\)-to-\(L^\infty\) estimate
	Theorem~\ref{cor:linfty-zero-inflow},
	\[
	\|f_n\|_{L^\infty(Q_{5/4}^+(Z_0))}
	\le
	C
	\left(
	\|f_n\|_{L^2(Q_{3/2}^+(Z_0))}
	+
	\|g_n\|_{L^\infty(Q_{3/2}^+(Z_0))}
	\right).
	\]
	Combining this with the estimate from Step~4 gives
	\[
	\left\|
	\frac{f_n}{\mathbf d(\,\cdot\,,\bR\times\gamma^0)^\theta}
	\right\|_{L^\infty(Q_1^+(Z_0))}
	\le
	C
	\left(
	\|f_n\|_{L^2(Q_{3/2}^+(Z_0))}
	+
	\|g_n\|_{L^\infty(Q_{3/2}^+(Z_0))}
	\right).
	\]
	Since \(\theta=1/2-\varepsilon\), this proves \eqref{2607041155}.
\end{proof}

\begin{proof}[\textbf{Proof of Theorem~\ref{260129331}}]
	Since \(Z_0\in \bR\times\gamma^0\), the left translation by
	\(Z_0^{-1}\) preserves the half-space geometry:
	\[
	Z_0^{-1}\circ Q_R^+(Z_0)=Q_R^+(0),
	\qquad
	Z_0^{-1}\circ(\bR\times\gamma^-)=\bR\times\gamma^-,\qquad Z_0^{-1}\circ(\bR\times\gamma^0)=\bR\times\gamma^0.
	\]
	Thus it suffices to prove the estimate for \(Z_0=0\). By the kinetic scaling
	in Lemma~\ref{lem:kinetic-scaling}, we may also assume that \(r=1\).
	
	Put
	\[
	\theta:=\frac12-\varepsilon.
	\]
	
	\medskip
	\noindent
	\emph{Step 1. A normalized estimate.}
	We first prove the following estimate under the normalized assumptions
	\[
	[a]_{C^\alpha_{\kin}(Q_2(0))}\le1,
	\qquad
	\|\nabla_v a\|_{L^\infty(Q_2(0))}
	+
	\|b\|_{L^\infty(Q_2(0))}
	\le1.
	\]
	By Lemma~\ref{lem:smooth-approximation-half-space}, there exist approximating
	solutions \(f_n\) and sources \(g_n\) such that
	\[
	f_n\to f
	\quad\text{in }L^2(Q_{3/2}^+(0)),
	\qquad
	\|g_n\|_{L^\infty(Q_{3/2}^+(0))}
	\le
	\|g\|_{L^\infty(Q_2^+(0))}.
	\]
	Applying Lemma~\ref{260705834} to \(f_n\), we obtain
	\begin{align}\label{260705838}
		\begin{aligned}
			\left\|
			\frac{f_n}{\mathbf d(\,\cdot\,,\bR\times\gamma^0)^\theta}
			\right\|_{L^\infty(Q_1^+(0))}
			&\le
			C
			\left(
			\|f_n\|_{L^2(Q_{3/2}^+(0))}
			+
			\|g_n\|_{L^\infty(Q_{3/2}^+(0))}
			\right) \\
			&\le
			C
			\left(
			\|f_n\|_{L^2(Q_{3/2}^+(0))}
			+
			\|g\|_{L^\infty(Q_2^+(0))}
			\right),
		\end{aligned}
	\end{align}
	where
	$C=C(d,\lambda,\Lambda,\alpha,\varepsilon)>0$.
	
	Passing to a subsequence, still denoted by \(f_n\), we may assume that
	\[
	f_n\to f
	\quad\text{a.e. in }Q_{3/2}^+(0).
	\]
	Since $\mathbf d(Z,\bR\times\gamma^0)>0$ for every $Z\in Q_1^+(0)$,
	we may pass to the pointwise a.e. limit in \eqref{260705838}. Using also
	$\|f_n\|_{L^2(Q_{3/2}^+(0))}
	\to
	\|f\|_{L^2(Q_{3/2}^+(0))}$,
	we get
	\begin{align}\label{260705839}
		\begin{aligned}
			\left\|
			\frac{f}{\mathbf d(\,\cdot\,,\bR\times\gamma^0)^\theta}
			\right\|_{L^\infty(Q_1^+(0))}
			\,&\le
			C
			\left(
			\|f\|_{L^2(Q_{3/2}^+(0))}
			+
			\|g\|_{L^\infty(Q_2^+(0))}
			\right)\\
			&\le
			C
			\left(
			\|f\|_{L^\infty(Q_{3/2}^+(0))}
			+
			\|g\|_{L^\infty(Q_2^+(0))}
			\right).
		\end{aligned}
	\end{align}
	
	By the same argument after left translation, \eqref{260705839} also holds in
	any cylinder centered at \(Y\in\bR\times\gamma^0\). Namely, if the
	normalized assumptions hold in \(Q_2(Y)\), then
	\begin{align}\label{eq:normalized-translated}
		\left\|
		\frac{f}{\mathbf d(\,\cdot\,,\bR\times\gamma^0)^\theta}
		\right\|_{L^\infty(Q_1^+(Y))}
		\le
		C
		\left(
		\|f\|_{L^\infty(Q_{3/2}^+(Y))}
		+
		\|g\|_{L^\infty(Q_2^+(Y))}
		\right).
	\end{align}
	
	\medskip
	\noindent
	\emph{Step 2. Scaling and removal of the normalization.}
	We now recover the dependence on \(M\). Set
	\[
	R_M:=16(M+1),
	\qquad
	\delta_M:=\frac1{64R_M}.
	\]
	For \(Z=(t,x,v)\), write
	\[
	Z_{R_M}:=(R_M^2t,R_M^3x,R_Mv),
	\qquad
	Z_{R_M^{-1}}:=(R_M^{-2}t,R_M^{-3}x,R_M^{-1}v),
	\]
	and define
	\[
	f_{R_M}(Z):=f(Z_{R_M^{-1}}),
	\qquad
	a_{R_M}^{ij}(Z):=a^{ij}(Z_{R_M^{-1}}),
	\]
	\[
	b_{R_M}^i(Z):=b^i(Z_{R_M^{-1}}),
	\qquad
	g_{R_M}(Z):=g(Z_{R_M^{-1}}).
	\]
	Then \(f_{R_M}\) satisfies
	\[
	(\partial_t+v\cdot\nabla_x)f_{R_M}
	=
	D_i\bigl(a_{R_M}^{ij}D_jf_{R_M}\bigr)
	+
	\frac{b_{R_M}^i}{R_M}D_if_{R_M}
	+
	\frac{g_{R_M}}{R_M^2}
	\quad\text{in }Q_{2R_M}^+(0),
	\]
	and
	\[
	f_{R_M}=0
	\quad\text{on }(\bR\times\gamma^-)\cap Q_{2R_M}(0).
	\]
	Moreover,
	\[
	[a_{R_M}]_{C^\alpha_{\kin}(Q_{2R_M}(0))}
	=
	R_M^{-\alpha}[a]_{C^\alpha_{\kin}(Q_2(0))}
	\le
	R_M^{-\alpha}M_1
	\le1,
	\]
	and
	\[
	\|D_v a_{R_M}\|_{L^\infty(Q_{2R_M}(0))}
	+
	\left\|
	\frac{b_{R_M}}{R_M}
	\right\|_{L^\infty(Q_{2R_M}(0))}
	\le
	\frac{2M_2}{R_M}
	\le1.
	\]
	Thus the rescaled equation satisfies the normalized assumptions in
	\(Q_{2R_M}(0)\).
	
	We estimate the weighted norm on \(Q_{1/2}^+(0)\). Decompose
	$Q_{1/2}^+(0)=I_1\cup I_2$,
	where
	\[
	I_1:=
	\left\{
	Z\in Q_{1/2}^+(0):
	\mathbf d(Z,\bR\times\gamma^0)<\delta_M
	\right\},\qquad 
	I_2:=
	\left\{
	Z\in Q_{1/2}^+(0):
	\mathbf d(Z,\bR\times\gamma^0)\ge \delta_M
	\right\}.
	\]
	
	Let \(Z\in I_1\), and put
	\[
	\rho_Z:=\mathbf d(Z,\bR\times\gamma^0).
	\]
	Writing \(Z=(t,(x',x^d),(v',v^d))\), take the explicit grazing point
	\[
	Y_Z:=\bigl(t,(x',0),(v',0)\bigr)\in\bR\times\gamma^0 ,
	\]
	for which \eqref{260702722} gives
	\(\mathbf d(Z,Y_Z)\le|x^d|^{1/3}+|v^d|\le4\rho_Z\).  Then
	\[
	\mathbf d(Z,Y_Z)<4\delta_M=\frac1{16R_M}.
	\]
	Since \(Y_Z\) has the same time coordinate as \(Z\),
	Lemma~\ref{260127226}(2) gives \(Y_Z\in Q_{1/(8R_M)}(Z)\); since
	\(Z\in Q_{1/2}(0)\), Lemma~\ref{lem:cylinder-inclusion} then gives
	\[
	Y_Z\in Q_{1/2+1/(8R_M)}(0).
	\]
	Applying the same inclusion once more,
	\[
	Q_{2/R_M}(Y_Z)
	\subset
	Q_{1/2+1/(8R_M)+2/R_M}(0).
	\]
	Since \(R_M\ge16\), we have
	$\frac12+\frac1{8R_M}+\frac2{R_M}<1$,
	and hence
	\[
	Q_{2/R_M}(Y_Z)\subset Q_1(0).
	\]
	In particular,
	\[
	Q_{3/(2R_M)}^+(Y_Z)\subset Q_1^+(0),
	\qquad
	Q_{2/R_M}^+(Y_Z)\subset Q_1^+(0).
	\]
	After the scaling \(Z\mapsto Z_{R_M}\), these inclusions become
	\[
	Q_{3/2}^+((Y_Z)_{R_M})\subset Q_{R_M}^+(0),
	\qquad
	Q_2^+((Y_Z)_{R_M})\subset Q_{R_M}^+(0).
	\]
	Therefore the translated normalized estimate \eqref{eq:normalized-translated}
	applies to \(f_{R_M}\) in the cylinder centered at \((Y_Z)_{R_M}\).
	
	Since \(Z_{R_M}\in Q_1^+((Y_Z)_{R_M})\) and
	\[
	\mathbf d(Z_{R_M},\bR\times\gamma^0)
	=
	R_M\,\mathbf d(Z,\bR\times\gamma^0),
	\]
	we obtain
	\[
	\frac{|f(Z)|}{\mathbf d(Z,\bR\times\gamma^0)^\theta}
	\le
	C R_M^\theta
	\left(
	\|f_{R_M}\|_{L^\infty(Q_{3/2}^+((Y_Z)_{R_M}))}
	+
	R_M^{-2}
	\|g_{R_M}\|_{L^\infty(Q_2^+((Y_Z)_{R_M}))}
	\right).
	\]
	The inclusions above imply
	\[
	\|f_{R_M}\|_{L^\infty(Q_{3/2}^+((Y_Z)_{R_M}))}
	\le
	\|f\|_{L^\infty(Q_1^+(0))},\qquad 
	\|g_{R_M}\|_{L^\infty(Q_2^+((Y_Z)_{R_M}))}
	\le
	\|g\|_{L^\infty(Q_1^+(0))}.
	\]
	Thus
	\[
	\sup_{I_1}
	\frac{|f|}{\mathbf d(\,\cdot\,,\bR\times\gamma^0)^\theta}
	\le
	C R_M^\theta
	\left(
	\|f\|_{L^\infty(Q_1^+(0))}
	+
	R_M^{-2}\|g\|_{L^\infty(Q_1^+(0))}
	\right).
	\]
	
	If \(Z\in I_2\), then
	$\mathbf d(Z,\bR\times\gamma^0)\ge\delta_M=\frac1{64R_M}$.
	Therefore
	\[
	\frac{|f(Z)|}{\mathbf d(Z,\bR\times\gamma^0)^\theta}
	\le
	C R_M^\theta
	\|f\|_{L^\infty(Q_{1/2}^+(0))}
	\le
	C R_M^\theta
	\|f\|_{L^\infty(Q_1^+(0))}.
	\]
	Hence
	\[
	\sup_{I_2}
	\frac{|f|}{\mathbf d(\,\cdot\,,\bR\times\gamma^0)^\theta}
	\le
	C R_M^\theta
	\|f\|_{L^\infty(Q_1^+(0))}.
	\]
	
	Combining the estimates on \(I_1\) and \(I_2\), and recalling that
	\(R_M=16(M+1)\), we obtain
	\[
	\left\|
	\frac{f}{\mathbf d(\,\cdot\,,\bR\times\gamma^0)^\theta}
	\right\|_{L^\infty(Q_{1/2}^+(0))}
	\le
	C(M+1)^\theta
	\left(
	\|f\|_{L^\infty(Q_1^+(0))}
	+
	(M+1)^{-2}\|g\|_{L^\infty(Q_1^+(0))}
	\right).
	\]
	Since \(\theta=1/2-\varepsilon\), this proves the estimate for \(Z_0=0\) and
	\(r=1\). Undoing the initial left translation and kinetic scaling gives
	\eqref{eq:grazing-weighted-zero}.
\end{proof}

\subsection{Velocity-gradient bounds for regular coefficients}\label{Subsec:Dvest}

Under kinetic H\"older regularity of the leading coefficients, the extension
lemmas reduce the non-grazing boundary estimate to the interior case. 
Combining local \(S_p\)-estimates at the scale of the distance to \(\gamma^0\) with the grazing decay estimate, which is obtained from a supersolution barrier, we obtain the weighted \(D_v\)-bound below;
the proof is given in \S\ref{260706939}. 
For a related use of superharmonic functions and boundary distance to describe boundary behavior of solutions and their derivatives in elliptic \(L^p\)-theory, see \cite{Seo2026}.

\begin{thm}[Scaled \(D_v\)-estimate]
	\label{prop:scaled-Dv-estimate}
	Let
	$Z_0\in \bR\times\overline{\bR^d_+}\times\bR^d$
	and $r>0$.
	Assume that there exists \(\alpha\in(0,1)\) such that
	\[
	a^{ij}\in C^\alpha_{\kin}(Q_{2r}(Z_0)),
	\qquad
	\nabla_v a^{ij},\,b^i\in L^\infty(Q_{2r}(Z_0)),
	\]
	where \(\nabla_v a^{ij}\) is understood in the weak sense.
	Assume moreover that
	there exist constants \(0<\lambda\le \Lambda\) and \(M_1,M_2\ge0\) such that,
	for a.e. \(Z\in Q_{2r}(Z_0)\) and every \(\xi\in\bR^d\),
	\[
	\lambda|\xi|^2
	\le
	a^{ij}(Z)\xi_i\xi_j,
	\qquad
	|a^{ij}(Z)|\le \Lambda,
	\]
	and
	\begin{align}\label{2607041007_11}
		[a]_{C^\alpha_{\kin}(Q_{2r}(Z_0))}\le M_1,
		\qquad
		\|\nabla_v a\|_{L^\infty(Q_{2r}(Z_0))}
		+
		\|b\|_{L^\infty(Q_{2r}(Z_0))}
		\le M_2.
	\end{align}
	Set
	\[
	M:=\max\{M_1^{1/\alpha},\,M_2\}.
	\]
	
	Let \(f\in L^\infty(Q_{2r}^+(Z_0))\) be a weak solution of
	\begin{align}\label{2605311025}
		\partial_t f+v\cdot \nabla_x f
		=
		D_i(a^{ij}D_j f)+b^iD_i f+g
		\qquad\text{in }Q_{2r}^+(Z_0),
	\end{align}
	with
	$f=0$ on $(\bR\times\gamma^-)\cap Q_{2r}(Z_0)$
	in the sense of traces, where \(g\in L^\infty(Q_{2r}^+(Z_0))\). Then the following hold.
	
	\begin{enumerate}
		\item
		If
		\begin{align}\label{2607261203}
			Q_{2r}(Z_0)\cap(\bR\times\gamma^0)=\emptyset,
		\end{align}
		then there exists
		$C=C(d,\lambda,\Lambda,\alpha)>0$
		such that
		\begin{align}\label{eq:scaled-nongrazing-Dv}
			\|\nabla_v f\|_{L^\infty(Q_r^+(Z_0))}
			\le
			C(M+r^{-1})
			\left(
			\|f\|_{L^\infty(Q_{2r}^+(Z_0))}
			+
			(M+r^{-1})^{-2}
			\|g\|_{L^\infty(Q_{2r}^+(Z_0))}
			\right).
		\end{align}
		
		\item
		If
		\[
		Z_0\in\bR\times\gamma^0,
		\]
		then for every \(\varepsilon\in(0,1/2)\), there exists
		$C=C(d,\lambda,\Lambda,\alpha,\varepsilon)>0$
		such that, for a.e. \(Z\in Q_r^+(Z_0)\),
		\begin{align}\label{eq:scaled-grazing-Dv}
			|\nabla_v f(Z)|
			\le
			C\left(
			\mathbf d(Z,\bR\times\gamma^0)^{-1/2-\varepsilon}I_1+I_2
			\right),
		\end{align}
		where
		\begin{align*}
			I_1\,&:=
			(M+r^{-1})^{1/2-\varepsilon}
			\|f\|_{L^\infty(Q_{2r}^+(Z_0))}
			+
			(M+r^{-1})^{-3/2-\varepsilon}
			\|g\|_{L^\infty(Q_{2r}^+(Z_0))},\\
			I_2\,&:=
			Mr^{1/2-\varepsilon}I_1
			+
			(M+r^{-1})^{-1}
			\|g\|_{L^\infty(Q_{2r}^+(Z_0))}.
		\end{align*}
	\end{enumerate}
\end{thm}

\subsubsection{Interior $D_vf$ estimate for regular coefficients}

We first combine global divergence- and nondivergence-form \(L^p\)-estimates
with the kinetic Sobolev--Morrey embedding to obtain an interior estimate for
\(D_vf\). We then combine this estimate with the extension procedures of
Lemmas~\ref{lem:inflow-zero-extension} and
\ref{lem:boundary-extension-away-from-grazing} to derive the corresponding
boundary estimates.	

\begin{prop}\label{prop:local-Sp}
	Let \(Z_0=(t_0,x_0,v_0)\in\bR^{1+2d}\)
	and \(r>0\). Assume that there exists \(\alpha\in(0,1)\) such that
	\[
	a^{ij}\in C^\alpha_{\kin}(Q_{2r}(Z_0)),
	\qquad
	\nabla_v a^{ij},\,b^i\in L^\infty(Q_{2r}(Z_0)),
	\]
	where \(\nabla_v a^{ij}\) is understood in the weak sense. Assume moreover
	that there exist constants \(0<\lambda\le\Lambda\) and \(M_1,M_2\ge0\) such
	that, for a.e. \(Z\in Q_{2r}(Z_0)\) and every \(\xi\in\bR^d\),
	\begin{align}\label{260605507_11}
		\lambda|\xi|^2
		\le
		a^{ij}(Z)\xi_i\xi_j,
		\qquad
		|a^{ij}(Z)|\le \Lambda,
	\end{align}
	and
	\begin{align*}
		[a]_{C^\alpha_{\kin}(Q_{2r}(Z_0))}\le M_1,
		\qquad
		\|\nabla_v a\|_{L^\infty(Q_{2r}(Z_0))}
		+
		\|b\|_{L^\infty(Q_{2r}(Z_0))}
		\le M_2.
	\end{align*}
	Set
	\[
	M:=\max\{M_1^{1/\alpha},\,M_2\}.
	\]
	
	Let \(f\in L^\infty(Q_{2r}(Z_0))\) be a weak solution of
	\begin{align}\label{2605311027_11}
		\partial_t f+v\cdot\nabla_x f
		=
		D_i(a^{ij}D_jf)+b^iD_if+g
		\qquad\text{in }Q_{2r}(Z_0),
	\end{align}
	where \(g\in L^\infty(Q_{2r}(Z_0))\). Then there exists a constant
	$C=C(d,\lambda,\Lambda,\alpha)>0$
	such that
	\begin{align}\label{eq:Dvestimate}
		\|\nabla_v f\|_{L^\infty(Q_r(Z_0))}
		\le
		C\left(M+r^{-1}\right)
		\left(
		\|f\|_{L^\infty(Q_{2r}(Z_0))}
		+
		\left(M+r^{-1}\right)^{-2}
		\|g\|_{L^\infty(Q_{2r}(Z_0))}
		\right).
	\end{align}
\end{prop}

The proof of Proposition~\ref{prop:local-Sp} relies on global
\(L^p\)-estimates for the kinetic Fokker--Planck equation in divergence form
\cite{dong2024global} and in nondivergence form \cite{dong2022global},
together with a kinetic Sobolev--Morrey embedding. We first recall the
coefficient assumption and the solution spaces used in
\cite{dong2022global,dong2024global}.

\begin{remark}
	The leading-coefficient assumption in
	\cite{dong2022global,dong2024global} is a small-oscillation condition in the
	\((x,v)\)-variables. Namely, for sufficiently small \(\gamma_0>0\), there
	exists \(R_0>0\) such that, for every
	\(Z_0=(t_0,x_0,v_0)\in\bR^{1+2d}\) and every \(R\in(0,R_0]\),
	\[
	\operatorname{osc}_{x,v}(a,Q_R(Z_0))\le \gamma_0,
	\]
	where
	\[
	\operatorname{osc}_{x,v}(a,Q_R(Z_0))
	:=
	\fint_{t_0-R^2}^{t_0}
	\fint_{D_R(Z_0,t)\times D_R(Z_0,t)}
	|a(t,x_1,v_1)-a(t,x_2,v_2)|
	\dd x_1\dd v_1\dd x_2\dd v_2\dd t,
	\]
	and
	\[
	D_R(Z_0,t)
	:=
	\left\{
	(x,v):
	|x-x_0-(t-t_0)v_0|<R^3,\quad |v-v_0|<R
	\right\}.
	\]
	
	We shall use this assumption after a cutoff. Let
	\(a^{ij}\in C^\alpha_{\kin}(Q_2(0))\) satisfy the structural assumptions of
	Proposition~\ref{prop:local-Sp} with \(r=1\). Choose
	\(\psi_0\in C_c^\infty(Q_1(0))\), \(0\le\psi_0\le1\), and define
	\[
	\hat a^{ij}(t,x,v)
	:=
	\begin{cases}
		\psi_0(t,x,v)a^{ij}(t,x,v)
		+
		(1-\psi_0(t,x,v))\Lambda\delta^{ij},
		& t\le0,\\[2mm]
		\Lambda\delta^{ij},
		& t>0.
	\end{cases}
	\]
	Then \(\hat a\) satisfies the above small-oscillation assumption for some
	\(R_0>0\) depending only on
	$d,\,\Lambda,\,\alpha,\,M_1,\,\psi_0,\,\gamma_0$.
	Indeed, by the cutoff construction,
	\[
	[\hat a]_{C^\alpha_{\kin}((-\infty,0]\times\bR^{2d})}
	\le
	C(d,\Lambda,\alpha,M_1,\psi_0)=:C_0,
	\]
	while \(\hat a\equiv \Lambda I_d\) for \(t>0\). 
	Hence, for
	\((x_1,v_1),(x_2,v_2)\in D_R(Z_0,t)\), denoting
	\(Z_i:=(t,x_i,v_i)\), we have
	\[
	|\hat a(Z_1)-\hat a(Z_2)|=0
	\qquad\text{if }t>0,\qquad |\hat a(Z_1)-\hat a(Z_2)|
	\le
	C_0\mathbf d(Z_1,Z_2)^\alpha
	\le
	CR^\alpha
	\qquad\text{if }t\le0,
	\]
	where $C=C(d,\Lambda,\alpha,M_1,\psi_0)$.
	Consequently,
	\[
	\sup_{Z_0\in\bR^{1+2d}}
	\operatorname{osc}_{x,v}(\hat a,Q_R(Z_0))
	\le
	C(d,\Lambda,\alpha,M_1,\psi_0)R^\alpha.
	\]
	Choosing \(R_0>0\) sufficiently small, depending only on
	\(d,\Lambda,\alpha,M_1,\psi_0\), and \(\gamma_0\), gives the desired condition.
\end{remark}

\begin{prop}\label{prop:global-Sobolev-estimates}
	Let \(1<p<\infty\). For \(k\in\mathbb Z\), set
	\[
	\bW^{k,p}(-4,0)
	:=
	L^p\big((-4,0)\times\bR_x^d;W^{k,p}(\bR_v^d)\big),
	\qquad
	\bL^p(-4,0):=\bW^{0,p}(-4,0).
	\]
	Let
	$u\in C_c^\infty((-4,0]\times\bR_x^d\times\bR_v^d)$,
	and assume that \(\hat a^{ij}\) and \(\hat b^i\) satisfy, for a.e.
	\(Z\in(-\infty,0]\times\bR^{2d}\) and every \(\xi\in\bR^d\),
	\[
	\lambda|\xi|^2
	\le
	\hat a^{ij}(Z)\xi_i\xi_j,
	\qquad
	|\hat a^{ij}(Z)|\le \Lambda,
	\]
	and
	\[
	[\hat a]_{C^\alpha_{\kin}((-\infty,0]\times\bR^{2d})}
	\le M_1,
	\qquad
	\|\nabla_v\hat a\|_{L^\infty((-\infty,0]\times\bR^{2d})}
	+
	\|\hat b\|_{L^\infty((-\infty,0]\times\bR^{2d})}
	\le M_2.
	\]
	Denote
	\[
	M:=\max(M_1^{1/\alpha}, M_2).
	\]
	Then the following estimates hold.
	
	\begin{enumerate}
		\item Suppose that there exist
		\(g^0,g^1,\ldots,g^d\in\bL^p(-4,0)\) such that
		\[
		(\partial_t+v\cdot\nabla_x)u
		=
		D_i(\hat a^{ij}D_j u)+\hat b^iD_i u+g^0+D_ig^i
		\qquad
		\text{in }\bW^{-1,p}(-4,0).
		\]
		Then
		\begin{align}\label{Spdiv}
			\|u\|_{\bS_p(-4,0)}
			:=
			\|u\|_{\bW^{1,p}(-4,0)}
			+
			\|(\partial_t+v\cdot\nabla_x)u\|_{\bW^{-1,p}(-4,0)}
			\le
			C\sum_{i=0}^d\|g^i\|_{\bL^p(-4,0)},
		\end{align}
		where $C=C(d,p,\lambda,\Lambda,\alpha,M)$.
		
		\item Suppose that there exists \(g\in\bL^p(-4,0)\) such that
		\[
		(\partial_t+v\cdot\nabla_x)u
		=
		\hat a^{ij}D_{ij}u
		+
		(D_j\hat a^{ji}+\hat b^i)D_i u
		+
		g
		\qquad
		\text{in }\bL^p(-4,0).
		\]
		Then
		\begin{align}\label{Spnondiv}
			\|u\|_{\cS_p(-4,0)}
			:=
			\|u\|_{\bW^{2,p}(-4,0)}
			+
			\|(\partial_t+v\cdot\nabla_x)u\|_{\bL^p(-4,0)}
			\le
			C\|g\|_{\bL^p(-4,0)},
		\end{align}
		where
		$C=C(d,p,\lambda,\Lambda,\alpha,M)$.
		
		\item If \(p>4d+2\), then
		\begin{align}\label{S_M_indq}
			\|D_vu\|_{L^\infty((-4,0]\times\bR^{2d})}
			\le
			C(d,p)\|u\|_{\cS_p(-4,0)}.
		\end{align}
	\end{enumerate}
\end{prop}

Estimates \((1)\) and \((2)\) follow from
\cite[Theorem~1.15]{dong2024global} and \cite[Theorem~2.6]{dong2022global}, respectively, after extending
\(\hat a=\Lambda I_d\) and \(\hat b=0\) to \(\{t>0\}\).
Estimate \((3)\) is the kinetic Sobolev embedding given in
\cite[Theorem~2.1(4)]{polidoro1998sobolev}.

\begin{proof}[\textbf{Proof of Proposition~\ref{prop:local-Sp}}]
	\noindent\emph{Step 1. Normalization.}
	By left translation and the kinetic scaling
	\[
	\delta_r(t,x,v):=(r^2t,r^3x,rv),
	\]
	it suffices to prove the case \(Z_0=0\) and \(r=1\). Indeed, if
	\[
	f_r(Z):=f(Z_0\circ\delta_r Z),
	\]
	then
	\[
	(\partial_t+v\cdot\nabla_x)f_r
	=
	D_i(a_r^{ij}D_jf_r)+b_r^iD_if_r+g_r
	\quad\text{in }Q_2(0),
	\]
	where
	\[
	a_r^{ij}(Z):=a^{ij}(Z_0\circ\delta_r Z),
	\qquad
	b_r^i(Z):=r\,b^i(Z_0\circ\delta_r Z),
	\qquad
	g_r(Z):=r^2g(Z_0\circ\delta_r Z).
	\]
	Moreover,
	\[
	[a_r]_{C^\alpha_{\kin}(Q_2(0))}\le r^\alpha M_1,
	\qquad
	\|\nabla_v a_r\|_{L^\infty(Q_2(0))}
	+
	\|b_r\|_{L^\infty(Q_2(0))}
	\le rM_2,
	\]
	and hence
	\[
	M_r=\max\{(r^\alpha M_1)^{1/\alpha},rM_2\}=rM.
	\]
	Thus the normalized estimate scales back to \eqref{eq:Dvestimate}.
	
	We therefore assume from now on that
	\[
	Z_0=0,
	\qquad
	r=1.
	\]
	It remains to prove
	\begin{align}\label{eq:normalized-Dv-estimate-local-Sp}
		\|D_vf\|_{L^\infty(Q_1(0))}
		\le
		C(M+1)
		\left(
		\|f\|_{L^\infty(Q_2(0))}
		+
		(M+1)^{-2}\|g\|_{L^\infty(Q_2(0))}
		\right),
	\end{align}
	where \(C=C(d,\lambda,\Lambda,\alpha)\).
	
	\medskip
	\noindent
	\emph{Step 2. Reduction to the unit-coefficient case.}
	We first prove that \eqref{eq:normalized-Dv-estimate-local-Sp} follows from
	the following unit-coefficient estimate. Assume that \eqref{260605507_11}
	holds in \(Q_2(0)\) and that
	\begin{align}\label{eq:unit-coeff-bound-local-Sp}
		[a]_{C^\alpha_{\kin}(Q_2(0))}\le1,
		\qquad
		\|\nabla_v a\|_{L^\infty(Q_2(0))}
		+
		\|b\|_{L^\infty(Q_2(0))}
		\le1.
	\end{align}
	Then
	\begin{align}\label{eq:unit-coeff-Dv-estimate-local-Sp}
		\|D_vf\|_{L^\infty(Q_1(0))}
		\le
		C
		\left(
		\|f\|_{L^\infty(Q_2(0))}
		+
		\|g\|_{L^\infty(Q_2(0))}
		\right),
	\end{align}
	where \(C=C(d,\lambda,\Lambda,\alpha)\).
	
	Assume \eqref{eq:unit-coeff-Dv-estimate-local-Sp}. Set
	\[
	\mu:=M+1,
	\qquad
	\rho:=\frac1{8\mu}.
	\]
	For each \(\bar Z\in Q_1(0)\),
	\begin{align}\label{eq:scaled-cylinder-inclusion-Q2-local-Sp}
		\bar Z\circ\delta_\rho Q_2(0)
		=
		Q_{2\rho}(\bar Z)
		\subset Q_2(0).
	\end{align}
	Define
	\[
	F(Z):=f(\bar Z\circ\delta_\rho Z).
	\]
	Then
	\[
	(\partial_t+v\cdot\nabla_x)F
	=
	D_i(A^{ij}D_jF)+B^iD_iF+G
	\quad\text{in }Q_2(0),
	\]
	where
	\[
	A^{ij}(Z):=a^{ij}(\bar Z\circ\delta_\rho Z),
	\qquad
	B^i(Z):=\rho b^i(\bar Z\circ\delta_\rho Z),
	\qquad
	G(Z):=\rho^2 g(\bar Z\circ\delta_\rho Z).
	\]
	By \eqref{260605507_11}, \eqref{2607041007_11}, and
	\eqref{eq:scaled-cylinder-inclusion-Q2-local-Sp}, the coefficients \(A,B\)
	satisfy \eqref{260605507_11} and \eqref{eq:unit-coeff-bound-local-Sp} in
	\(Q_2(0)\). Indeed,
	\[
	[A]_{C^\alpha_{\kin}(Q_2(0))}
	\le
	\rho^\alpha M_1
	\le1,
	\qquad
	\|\nabla_vA\|_{L^\infty(Q_2(0))}
	+
	\|B\|_{L^\infty(Q_2(0))}
	\le
	\rho M_2
	\le1.
	\]
	Therefore \eqref{eq:unit-coeff-Dv-estimate-local-Sp} applies to \(F\). Since
	\[
	D_vF(Z)=\rho D_vf(\bar Z\circ\delta_\rho Z),
	\]
	we obtain
	\[
	\rho |D_vf(\bar Z)|
	\le
	C
	\left(
	\|F\|_{L^\infty(Q_2(0))}
	+
	\|G\|_{L^\infty(Q_2(0))}
	\right).
	\]
	Using \eqref{eq:scaled-cylinder-inclusion-Q2-local-Sp},
	\[
	\|F\|_{L^\infty(Q_2(0))}
	\le
	\|f\|_{L^\infty(Q_2(0))},
	\qquad
	\|G\|_{L^\infty(Q_2(0))}
	\le
	\rho^2\|g\|_{L^\infty(Q_2(0))}.
	\]
	Hence
	\[
	|D_vf(\bar Z)|
	\le
	C
	\left(
	\rho^{-1}\|f\|_{L^\infty(Q_2(0))}
	+
	\rho\|g\|_{L^\infty(Q_2(0))}
	\right).
	\]
	Since \(\rho^{-1}=8(M+1)\) and \(\rho\le(M+1)^{-1}\), taking the supremum
	over \(\bar Z\in Q_1(0)\) gives \eqref{eq:normalized-Dv-estimate-local-Sp}.
	Thus it remains to prove \eqref{eq:unit-coeff-Dv-estimate-local-Sp}.
	
	\medskip
	\noindent
	\emph{Step 3. Proof of the unit-coefficient estimate.}
	Assume \eqref{260605507_11} and \eqref{eq:unit-coeff-bound-local-Sp}. By the
	interior version of the approximation lemma
	(Lemma~\ref{lem:smooth-approximation-interior}), there
	exist
	\[
	a_n^{ij},\,b_n^i,\,g_n\in C^\infty(\bR^{1+2d}),
	\qquad
	f_n\in C(\overline{Q_{3/2}(0)})\cap C^\infty(Q_{3/2}(0)),
	\]
	such that \(a_n^{ij}\) and \(b_n^i\) satisfy \eqref{260605507_11} in
	\(Q_{3/2}(0)\) and
	\begin{align*}
		[a_n]_{C^\alpha_{\kin}(Q_{3/2}(0))}\le1,
		\qquad
		\|\nabla_v a_n\|_{L^\infty(Q_{3/2}(0))}
		+
		\|b_n\|_{L^\infty(Q_{3/2}(0))}
		\le1.
	\end{align*}
	Moreover,
	\[
	(\partial_t+v\cdot\nabla_x)f_n
	=
	D_i(a_n^{ij}D_jf_n)+b_n^iD_if_n+g_n
	\quad\text{in }Q_{3/2}(0),
	\]
	\[
	f_n\to f,
	\qquad
	D_vf_n\to D_vf
	\quad\text{in }L^2(Q_{3/2}(0)),
	\]
	and
	\[
	\|g_n\|_{L^\infty(Q_{3/2}(0))}
	\le
	\|g\|_{L^\infty(Q_2(0))}.
	\]
	
	We first prove the smooth estimate
	\begin{align}\label{eq:smooth-unit-Dv-estimate-local-Sp}
		\|D_vf\|_{L^\infty(Q_1(0))}
		\le
		C
		\left(
		\|f\|_{L^\infty(Q_{5/4}(0))}
		+
		\|g\|_{L^\infty(Q_{5/4}(0))}
		\right)
	\end{align}
	for smooth \(f,a,b,g\) satisfying \eqref{260605507_11} in \(Q_{3/2}(0)\) and
	\[
	[a]_{C^\alpha_{\kin}(Q_{3/2}(0))}\le1,
	\qquad
	\|\nabla_v a\|_{L^\infty(Q_{3/2}(0))}
	+
	\|b\|_{L^\infty(Q_{3/2}(0))}
	\le1.
	\]
	Assuming this for the moment, applying
	\eqref{eq:smooth-unit-Dv-estimate-local-Sp} to \(f_n\) and using the
	unit-coefficient \(L^2\)-\(L^\infty\) estimate
	\begin{align*}
		\|f_n\|_{L^\infty(Q_{5/4}(0))}
		\le
		C
		\left(
		\|f_n\|_{L^2(Q_{3/2}(0))}
		+
		\|g_n\|_{L^\infty(Q_{3/2}(0))}
		\right)
	\end{align*}
	gives
	\[
	\|D_vf_n\|_{L^\infty(Q_1(0))}
	\le
	C
	\left(
	\|f_n\|_{L^2(Q_{3/2}(0))}
	+
	\|g_n\|_{L^\infty(Q_{3/2}(0))}
	\right).
	\]
	After passing to a subsequence, \(D_vf_n\to D_vf\) a.e. in \(Q_1(0)\). Hence,
	for a.e. \(Z\in Q_1(0)\),
	\[
	|D_vf(Z)|
	\le
	\limsup_{n\to\infty}\|D_vf_n\|_{L^\infty(Q_1(0))}
	\le
	C
	\left(
	\|f\|_{L^2(Q_{3/2}(0))}
	+
	\|g\|_{L^\infty(Q_2(0))}
	\right).
	\]
	Therefore
	\[
	\|D_vf\|_{L^\infty(Q_1(0))}
	\le
	C
	\left(
	\|f\|_{L^2(Q_{3/2}(0))}
	+
	\|g\|_{L^\infty(Q_2(0))}
	\right).
	\]
	Since
	\[
	\|f\|_{L^2(Q_{3/2}(0))}
	\le
	C(d)\|f\|_{L^\infty(Q_2(0))},
	\]
	we obtain \eqref{eq:unit-coeff-Dv-estimate-local-Sp}.
	
	It remains to prove \eqref{eq:smooth-unit-Dv-estimate-local-Sp}.
	
	\smallskip
	\noindent
	\emph{Proof of \eqref{eq:smooth-unit-Dv-estimate-local-Sp}.}
	Let $p:=4d+3$.
	Choose
	\begin{alignat*}{3}
		&\zeta\in C_c^\infty(Q_{9/8}(0)),
		\qquad
		&&0\le \zeta\le 1\quad\text{in }Q_{9/8}(0),
		\qquad
		&&\zeta\equiv1\quad\text{on }Q_1(0),\\
		&\eta\in C_c^\infty(Q_{5/4}(0)),
		\qquad
		&&0\le \eta\le 1\quad\text{in }Q_{5/4}(0),
		\qquad
		&&\eta\equiv1\quad\text{on }Q_{9/8}(0),\\
		&\chi\in C_c^\infty(Q_{4/3}(0)),
		\qquad
		&&0\le \chi\le 1\quad\text{in }Q_{4/3}(0),
		\qquad
		&&\chi\equiv1\quad\text{on }Q_{5/4}(0).			
	\end{alignat*}
	The cutoffs are fixed once and for all.
	
	Set \(u:=f\eta\). Define
	\[
	\hat a_\eta^{ij}:=\chi a^{ij}+(1-\chi)\Lambda\delta^{ij},
	\qquad
	\hat b_\eta^i:=\chi b^i .
	\]
	By \eqref{260605507_11}, the unit bounds in \(Q_{3/2}(0)\), and the cutoff
	construction, \(\hat a_\eta,\hat b_\eta\) satisfy the hypotheses of
	Proposition~\ref{prop:global-Sobolev-estimates}, with coefficient parameter
	bounded by a constant depending only on \(d,\lambda,\Lambda,\alpha\). Since
	\(\chi\equiv1\) on \(\operatorname{supp}\eta\),
	\begin{align}\label{eq:cutoff-div-equation}
		(\partial_t+v\cdot\nabla_x)u
		=
		D_i(\hat a_\eta^{ij}D_ju)+\hat b_\eta^iD_iu
		+
		G_\eta^0+D_iG_\eta^i
		\quad\text{in }\bW^{-1,p}(-4,0),
	\end{align}
	where
	\[
	G_\eta^i
	:=
	-\big(a^{ij}D_j\eta+a^{ji}D_j\eta\big)f,
	\]
	and
	\[
	G_\eta^0
	:=
	\eta g
	+
	\left[
	(\partial_t+v\cdot\nabla_x)\eta
	-
	b^iD_i\eta
	+
	D_i a^{ji}D_j\eta
	+
	a^{ji}D_{ij}\eta
	\right]f .
	\]
	Using \eqref{260605507_11}, the unit bounds in \(Q_{3/2}(0)\), and the fixed
	choice of \(\eta\), we have
	\begin{align}\label{eq:cutoff-div-source-bound}
		\|G_\eta^0\|_{\bL^p(-4,0)}
		+
		\sum_{i=1}^d\|G_\eta^i\|_{\bL^p(-4,0)}
		\le
		C
		\left(
		\|f\|_{L^\infty(Q_{5/4}(0))}
		+
		\|g\|_{L^\infty(Q_{5/4}(0))}
		\right).
	\end{align}
	Applying \eqref{Spdiv} to \eqref{eq:cutoff-div-equation} and using
	\eqref{eq:cutoff-div-source-bound}, we obtain
	\begin{align*}
		\|f\eta\|_{\bS_p(-4,0)}
		\le
		C
		\left(
		\|f\|_{L^\infty(Q_{5/4}(0))}
		+
		\|g\|_{L^\infty(Q_{5/4}(0))}
		\right).
	\end{align*}
	Since \(\eta\equiv1\) on \(\operatorname{supp}\zeta\),
	\begin{align}\label{eq:Dv-Lp-on-zeta-support}
		\|D_vf\|_{L^p(\operatorname{supp}\zeta)}
		\le
		C
		\left(
		\|f\|_{L^\infty(Q_{5/4}(0))}
		+
		\|g\|_{L^\infty(Q_{5/4}(0))}
		\right).
	\end{align}
	
	Next set \(w:=f\zeta\). Define
	\[
	\hat a_\zeta^{ij}:=\eta a^{ij}+(1-\eta)\Lambda\delta^{ij},
	\qquad
	\hat b_\zeta^i:=\eta b^i .
	\]
	The same argument shows that \(\hat a_\zeta,\hat b_\zeta\) satisfy the
	hypotheses of Proposition~\ref{prop:global-Sobolev-estimates}, with coefficient
	parameter bounded only in terms of \(d,\lambda,\Lambda,\alpha\). Since
	\(\eta\equiv1\) on \(\operatorname{supp}\zeta\),
	\begin{align}\label{eq:cutoff-nondiv-equation}
		(\partial_t+v\cdot\nabla_x)w
		=
		\hat a_\zeta^{ij}D_{ij}w
		+
		(D_j\hat a_\zeta^{ji}+\hat b_\zeta^i)D_iw
		+
		H_\zeta
		\quad\text{in }\bL^p(-4,0),
	\end{align}
	where
	\[
	\begin{aligned}
		H_\zeta
		:={}&
		\zeta g
		+
		\big((\partial_t+v\cdot\nabla_x)\zeta\big)f  \\
		&-
		a^{ij}\big(D_ifD_j\zeta+D_jfD_i\zeta\big)
		-
		a^{ij}fD_{ij}\zeta
		-
		(D_ja^{ji}+b^i)fD_i\zeta .
	\end{aligned}
	\]
	Using \eqref{260605507_11}, the unit bounds in \(Q_{3/2}(0)\),
	\eqref{eq:Dv-Lp-on-zeta-support}, and the fixed choice of \(\zeta\), we get
	\begin{align}\label{eq:nondiv-source-bound}
		\|H_\zeta\|_{\bL^p(-4,0)}
		\le
		C
		\left(
		\|f\|_{L^\infty(Q_{5/4}(0))}
		+
		\|g\|_{L^\infty(Q_{5/4}(0))}
		\right).
	\end{align}
	Applying \eqref{Spnondiv} to \eqref{eq:cutoff-nondiv-equation} and using
	\eqref{eq:nondiv-source-bound}, we obtain
	\begin{align}\label{eq:local-Sp-bound}
		\|f\zeta\|_{\cS_p(-4,0)}
		\le
		C
		\left(
		\|f\|_{L^\infty(Q_{5/4}(0))}
		+
		\|g\|_{L^\infty(Q_{5/4}(0))}
		\right).
	\end{align}
	Finally, by \eqref{S_M_indq},
	\[
	\|D_v(f\zeta)\|_{L^\infty((-4,0]\times\bR^{2d})}
	\le
	C(d,p)\|f\zeta\|_{\cS_p(-4,0)}.
	\]
	Since \(\zeta\equiv1\) on \(Q_1(0)\), \eqref{eq:local-Sp-bound} implies
	\[
	\|D_vf\|_{L^\infty(Q_1(0))}
	\le
	C
	\left(
	\|f\|_{L^\infty(Q_{5/4}(0))}
	+
	\|g\|_{L^\infty(Q_{5/4}(0))}
	\right).
	\]
	This proves \eqref{eq:smooth-unit-Dv-estimate-local-Sp}.
\end{proof}

\subsubsection{Proof of Theorem~\ref{prop:scaled-Dv-estimate}}\label{260706939}

\begin{proof}
	(1)
	By assumption \eqref{2607261203} and Lemma~\ref{2608051530}, applied
	with \(R=2r\), either
	\[
	Q_{2r}(Z_0)\cap (\bR\times \gamma)
	\subset \bR\times \gamma^-
	\qquad\text{or}\qquad
	Q_{2r}(Z_0)\cap (\bR\times \gamma)
	\subset \bR\times \gamma^+.
	\]
	Applying Lemma~\ref{lem:inflow-zero-extension} in the first case and
	Lemma~\ref{lem:boundary-extension-away-from-grazing} in the second, we obtain
	an extension $\bar f$ defined on $Q_{2r}(Z_0)$ that is a weak solution of
	\begin{align*}
		(\partial_t+v\cdot\nabla_x)\bar f
		={}&
		D_i\bigl(a^{ij}D_j\bar f\bigr)
		+b^iD_i\bar f
		+g\,\mathbf 1_{Q_{2r}^+(Z_0)}
		\qquad\text{in }Q_{2r}(Z_0).
	\end{align*}
	Moreover,
	\begin{align}\label{2607261215}
		\bar f=f
		\quad\text{in }Q_{2r}^+(Z_0),
		\qquad
		\|\bar f\|_{L^\infty(Q_{2r}(Z_0))}
		=
		\|f\|_{L^\infty(Q_{2r}^+(Z_0))}.
	\end{align}
	Since $a^{ij}$ and $b^i$ are defined on $Q_{2r}(Z_0)$ and satisfy
	\eqref{2607041007_11}, Proposition~\ref{prop:local-Sp} yields
	\[
	\|\nabla_v\bar f\|_{L^\infty(Q_r(Z_0))}
	\le
	C(M+r^{-1})
	\left(
	\|\bar f\|_{L^\infty(Q_{2r}(Z_0))}
	+
	(M+r^{-1})^{-2}
	\|g\|_{L^\infty(Q_{2r}^+(Z_0))}
	\right).
	\]
	Combining this estimate with \eqref{2607261215} gives
	\eqref{eq:scaled-nongrazing-Dv}.

	(2)
	Fix \(\varepsilon\in(0,1/2)\). We prove the estimate for a.e.
	\(Z\in Q_r^+(Z_0)\). Fix such a point \(Z\), and set
	\[
	d_Z:=\mathbf d(Z,\bR\times\gamma^0).
	\]
	Since \(Z\in Q_r^+(Z_0)\), we have \(x^d>0\), so \(d_Z>0\) by
	Lemma~\ref{260522259}. Since
	\(Z_0\in\bR\times\gamma^0\) and \(Z\in Q_r^+(Z_0)\), Lemma~\ref{260522259}
	gives
	\begin{align}\label{eq:dZ-upper}
		d_Z\le C r .
	\end{align}
	
	Set \(C_0:=10\).  Write \(Z=(t,x,v)\) and
	\(Z_*:=(t,(x',0),(v',0))\in\bR\times\gamma^0\).  The proof of
	Lemma~\ref{260522259} gives
	\(\mathbf d(Z,Z_*)\le|x^d|^{1/3}+|v^d|\le4d_Z\).  Since \(Z_*\) has
	the same time coordinate as \(Z\), Lemma~\ref{260127226}(2) applies and
	yields \(Z_*\in Q_{10d_Z}(Z)\).  Hence
	\[
	Q_{C_0d_Z}(Z)\cap(\bR\times\gamma^0)\neq\emptyset,
	\]
	and, whenever \(0<s\le C_0^{-1}d_Z\), Lemma~\ref{260127226}(1) and the
	triangle inequality for \(\mathbf d\) give
	\begin{align}\label{eq:distance-comparison-small-cylinder}
		\mathbf d(Y,\bR\times\gamma^0)\simeq d_Z
		\qquad
		\text{for every }Y\in Q_s(Z).
	\end{align}
	Choose \(0<\theta<\kappa<1\), depending only on \(d\), such that
	\begin{align}\label{eq:theta-kappa-choice}
		2\kappa+2C_0\theta\le 1,
		\qquad
		8C_0\theta\le\kappa,
		\qquad
		4\theta\le C_0^{-1}.
	\end{align}
	
	\smallskip
	\noindent
	\emph{Case 1: \(d_Z<\theta r\).}
	By Lemma~\ref{260515154}, applied to \(Q_{C_0d_Z}(Z)\), there exists
	\[
	\widehat Z\in\bR\times\gamma^0
	\]
	such that, using \(2^{1/3}\le2\) in Lemma~\ref{260515154},
	\begin{align}\label{eq:touching-cylinder}
		Q_{C_0d_Z}(Z)\subset Q_{2C_0d_Z}(\widehat Z),
		\qquad
		Q_R(\widehat Z)\subset Q_{R+2C_0d_Z}(Z)
		\quad\text{for every }R>0.
	\end{align}
	Set
	\[
	R:=\kappa r.
	\]
	Using \(Z\in Q_r^+(Z_0)\), \(d_Z<\theta r\), and
	\eqref{eq:theta-kappa-choice}, we have
	\begin{align}\label{eq:grazing-local-cylinder}
		Q_{2R}^+(\widehat Z)\subset Q_{2r}^+(Z_0),
		\qquad
		Q_{2\theta d_Z}^+(Z)\subset Q_{R/2}^+(\widehat Z).
	\end{align}
	Indeed, by \eqref{eq:touching-cylinder} and
	Lemma~\ref{lem:cylinder-inclusion},
	\[
	Q_{2R}(\widehat Z)
	\subset Q_{2R+2C_0d_Z}(Z)
	\subset Q_{r+2\kappa r+2C_0\theta r}(Z_0)
	\subset Q_{2r}(Z_0),
	\]
	and
	\[
	Q_{2\theta d_Z}(Z)
	\subset Q_{C_0d_Z}(Z)
	\subset Q_{2C_0d_Z}(\widehat Z)
	\subset Q_{R/2}(\widehat Z),
	\]
	because \(2C_0d_Z\le 2C_0\theta r\le R/4\). Moreover, by
	\eqref{eq:distance-comparison-small-cylinder} and \(4\theta\le C_0^{-1}\),
	\begin{align}\label{eq:d-comparable}
		\mathbf d(Y,\bR\times\gamma^0)\simeq d_Z
		\qquad
		\text{for every }Y\in Q_{2\theta d_Z}(Z).
	\end{align}
	In particular,
	\[
	Q_{2\theta d_Z}(Z)\cap(\bR\times\gamma^0)=\emptyset.
	\]
	
	Apply Theorem~\ref{260129331} in \(Q_{2R}^+(\widehat Z)\). Since
	\(Q_R^+(\widehat Z)\subset Q_{2r}^+(Z_0)\) and
	$M+R^{-1}\simeq M+r^{-1}$,
	we obtain
	\begin{align}\label{eq:local-grazing-decay}
		\left\|
		\frac{f}{\mathbf d(\cdot,\bR\times\gamma^0)^{1/2-\varepsilon}}
		\right\|_{L^\infty(Q_{R/2}^+(\widehat Z))}
		\le
		C I_1 ,
	\end{align}
	where \(C=C(d,\lambda,\Lambda,\alpha,\varepsilon)\). Therefore
	\eqref{eq:grazing-local-cylinder}, \eqref{eq:d-comparable}, and
	\eqref{eq:local-grazing-decay} imply
	\begin{align}\label{eq:f-local-bound}
		\|f\|_{L^\infty(Q_{2\theta d_Z}^+(Z))}
		\le
		C d_Z^{1/2-\varepsilon}I_1 .
	\end{align}
	
	Applying \eqref{eq:scaled-nongrazing-Dv} with center \(Z\) and radius
	\(\theta d_Z\), and absorbing powers of \(\theta\) into the constant, gives
	\begin{align}\label{eq:local-Dv-near}
		|\nabla_v f(Z)|
		\le
		C\left[
		(M+d_Z^{-1})
		\|f\|_{L^\infty(Q_{2\theta d_Z}^+(Z))}
		+
		(M+d_Z^{-1})^{-1}
		\|g\|_{L^\infty(Q_{2r}^+(Z_0))}
		\right].
	\end{align}
	Combining \eqref{eq:local-Dv-near} with \eqref{eq:f-local-bound}, we get
	\[
	|\nabla_v f(Z)|
	\le
	C\left[
	\left(Md_Z^{1/2-\varepsilon}
	+d_Z^{-1/2-\varepsilon}\right)I_1
	+
	(M+d_Z^{-1})^{-1}
	\|g\|_{L^\infty(Q_{2r}^+(Z_0))}
	\right].
	\]
	Using \eqref{eq:dZ-upper} and \(d_Z<\theta r\),
	\[
	Md_Z^{1/2-\varepsilon}I_1
	\le
	CMr^{1/2-\varepsilon}I_1,
	\qquad
	(M+d_Z^{-1})^{-1}
	\le
	C(M+r^{-1})^{-1}.
	\]
	Therefore
	\[
	|\nabla_v f(Z)|
	\le
	C\left(
	d_Z^{-1/2-\varepsilon}I_1
	+
	Mr^{1/2-\varepsilon}I_1
	+
	(M+r^{-1})^{-1}
	\|g\|_{L^\infty(Q_{2r}^+(Z_0))}
	\right),
	\]
	which gives \eqref{eq:scaled-grazing-Dv} in Case 1.
	
	\smallskip
	\noindent
	\emph{Case 2: \(d_Z\ge\theta r\).}
	Set
	\[
	\rho:=\frac{\theta r}{16C_0}.
	\]
	Then
	\begin{align*}
		Q_{2\rho}^+(Z)\subset Q_{2r}^+(Z_0),
		\qquad
		Q_{2\rho}(Z)\cap(\bR\times\gamma^0)=\emptyset .
	\end{align*}
	Indeed, by Lemma~\ref{lem:cylinder-inclusion},
	\[
	Q_{2\rho}(Z)
	\subset Q_{r+2\rho}(Z_0)
	\subset Q_{2r}(Z_0),
	\]
	and since \(2\rho\le C_0^{-1}d_Z\), \eqref{eq:distance-comparison-small-cylinder}
	gives
	\[
	\mathbf d(Y,\bR\times\gamma^0)\simeq d_Z>0
	\qquad
	\text{for every }Y\in Q_{2\rho}(Z).
	\]
	Applying \eqref{eq:scaled-nongrazing-Dv} with center \(Z\) and radius \(\rho\),
	and using
	$M+\rho^{-1}\simeq M+r^{-1}$,
	we obtain
	\begin{align}\label{eq:far-nongrazing-Dv}
		|\nabla_v f(Z)|
		\le
		C\left[
		(M+r^{-1})\|f\|_{L^\infty(Q_{2r}^+(Z_0))}
		+
		(M+r^{-1})^{-1}\|g\|_{L^\infty(Q_{2r}^+(Z_0))}
		\right].
	\end{align}
	By \eqref{eq:dZ-upper},
	$d_Z^{-1/2-\varepsilon}\ge cr^{-1/2-\varepsilon}$.
	Thus
	\[
	\begin{aligned}
		d_Z^{-1/2-\varepsilon}I_1+I_2
		&\ge
		c\left(r^{-1/2-\varepsilon}
		+
		Mr^{1/2-\varepsilon}\right)
		(M+r^{-1})^{1/2-\varepsilon}
		\|f\|_{L^\infty(Q_{2r}^+(Z_0))}
		\\
		&\quad
		+
		(M+r^{-1})^{-1}
		\|g\|_{L^\infty(Q_{2r}^+(Z_0))}.
	\end{aligned}
	\]
	Since
	\[
	\left(r^{-1/2-\varepsilon}
	+
	Mr^{1/2-\varepsilon}\right)
	(M+r^{-1})^{1/2-\varepsilon}
	=
	(M+r^{-1})(1+Mr)^{1/2-\varepsilon}
	\ge
	M+r^{-1},
	\]
	the right-hand side of \eqref{eq:far-nongrazing-Dv} is bounded by
	$C\left(d_Z^{-1/2-\varepsilon}I_1+I_2\right)$.
	This proves \eqref{eq:scaled-grazing-Dv} in Case 2, and hence completes the
	proof of (2).
\end{proof}

\section{Coulomb Landau Operators and One-Sided Kinetic Estimates}\label{sec:landau-framework}

This section collects the material on the Coulomb Landau operator. 
From this section on, the dimension is fixed to be $d=3$.

Recall the Landau equation \eqref{new:eq:intro-landau}, with the
collision operator and Coulomb kernel defined in
\eqref{new:eq:intro-landau-operator} and
\eqref{new:eq:intro-coulomb-kernel}.
More precisely, for functions $G$ and $F$ for which the expressions below
are locally integrable, we interpret $\cN(G,F)$ distributionally by
\[
\begin{aligned}
	\big\langle \cN(G,F),\varphi\big\rangle
	={}&
	-\int_{\bR^3}
	(\phi^{ij}\ast G)\,
	D_jF\,D_i\varphi
	\dd v+
	\int_{\bR^3}
	D_j(\phi^{ij}\ast G)\,
	F\,D_i\varphi
	\dd v
\end{aligned}
\]
for every $\varphi\in C_c^\infty(\bR^3)$.

With the normalization fixed at the end of
	Section~\ref{sec:main-results}, we write
\begin{equation}\label{def:reference-maxwellian}
	F(t,x,v)
	:=
	\mu(v)+\sqrt{\mu(v)}\,f(t,x,v),
	\qquad
	\mu(v):=\ee^{-|v|^2}.
\end{equation}
Since $\cN(\mu,\mu)=0$, the perturbation $f$ satisfies
\[
\begin{aligned}
	\p_t f+v\cdot\nabla_x f
	=
	\frac{1}{\sqrt{\mu}}
	\Big(
	\cN(\mu,\sqrt{\mu}f)
	+
	\cN(\sqrt{\mu}f,\mu)
	+
	\cN(\sqrt{\mu}f,\sqrt{\mu}f)
	\Big).
\end{aligned}
\]

Define
\begin{equation}\label{def:N}
	N[g,f]
	:=
	\frac{1}{\sqrt{\mu}}\,
	\cN(\sqrt{\mu}g,\sqrt{\mu}f),
\end{equation}
and set
\begin{equation}\label{def:LAK}
	Af:=N[\sqrt{\mu},f],
	\qquad
	Kf:=N[f,\sqrt{\mu}],
	\qquad
	L:=-A-K.
\end{equation}
Then the perturbative Landau equation takes the form
\[
\p_t f+v\cdot\nabla_x f+Lf=N[f,f].
\]

\begin{remark}[Well-definedness of the collision term]\label{rem:well-defined-N}
	The formula above defines $\cN(G,F)$ as a distribution on $\bR_v^3$ for
	$G\in L^1(\bR_v^3)$ and $F\in H^1(\bR_v^3)$. Indeed, since
	$|\phi^{ij}(z)|\lesssim|z|^{-1}$ and $|\nabla\phi^{ij}(z)|\lesssim|z|^{-2}$,
	we have $\phi^{ij}\in L^{3,\infty}(\bR_v^3)$ and
	$\p_j\phi^{ij}\in L^{3/2,\infty}(\bR_v^3)$, so that Young's convolution
	inequality in Lorentz spaces gives
	$\phi^{ij}\ast G\in L^{3,\infty}(\bR_v^3)$ and
	$(\p_j\phi^{ij})\ast G\in L^{3/2,\infty}(\bR_v^3)$, with norms controlled by
	$\|G\|_{L^1(\bR_v^3)}$.
	Together with $\nabla F\in L^2(\bR_v^3)$ and
	$F\in H^1(\bR_v^3)\hookrightarrow L^6(\bR_v^3)$, Hölder's inequality on the
	support of $\varphi$ yields
	\begin{equation}\label{eq:N-distributional-bound}
		\left|\big\langle\cN(G,F),\varphi\big\rangle\right|
		\le
		C_\varphi\|G\|_{L^1(\bR_v^3)}\|F\|_{H^1(\bR_v^3)},
		\qquad
		\varphi\in C_c^\infty(\bR_v^3).
	\end{equation}
	In particular, $N[g,f]=\mu^{-1/2}\cN(\sqrt{\mu}g,\sqrt{\mu}f)$ is well defined
	whenever $\sqrt{\mu}\,g\in L^1(\bR_v^3)$ and $\sqrt{\mu}\,f\in H^1(\bR_v^3)$.
\end{remark}

\subsection{The linearized Landau operator and the anisotropic Sobolev spaces $H^1_\sigma$}

\subsubsection{Landau operators and diffusion coefficients}

Following \cite{GuoCMP2002}, for any function $G$ for which the
convolution is well defined, we set
\[
\sigma_G^{ij}:=\phi^{ij}\ast G,
\qquad
\sigma_G^i:=v_j\sigma_G^{ij}. 
\]
We also write
\[
\sigma^{ij}:=\sigma_\mu^{ij},
\qquad
\sigma^i:=\sigma_\mu^i.
\]

\begin{prop}\label{prop:explicit-landau-operators}
	Let $f,g\in \cS(\bR_v^3)$.
	Then
	\begin{align}
		N[g,f]
		={}&
		\p_i\left(\sigma_{\sqrt{\mu}g}^{ij}\p_jf\right)
		-
		\left(
		\p_j\sigma_{\sqrt{\mu}g}^{ij}
		+
		2v_j\sigma_{\sqrt{\mu}g}^{ij}
		\right)\p_i f+
		\left(
		8\pi\sqrt{\mu}\,g
		+
		(v_iv_j-\delta_{ij})\sigma_{\sqrt{\mu}g}^{ij}
		\right)f\label{eq:explicit-landau-nonlinearity}
		\\
		={}&
		\p_i\left(\sigma_{\sqrt{\mu}g}^{ij}\p_jf\right)
		-
		2v_j\sigma_{\sqrt{\mu}g}^{ij}
		\p_i f-\p_i\left(
		\p_j\sigma_{\sqrt{\mu}g}^{ij}\,f\right)+
		(v_iv_j-\delta_{ij})\sigma_{\sqrt \mu g}^{ij}f.\label{eq:explicit-landau-nonlinearity_2}
	\end{align}
	Moreover,
	\begin{align}
		Af
		={}&
		\p_i\left(\sigma^{ij}\p_jf\right)
		-
		v_iv_j\sigma^{ij}f
		+
		4\pi\mu f,
		\label{eq:explicit-landau-A}
		\\
		Kf
		={}&
		\sqrt{\mu}
		\left(
		4v_iv_j-2\delta_{ij}
		\right)
		\sigma_{\sqrt{\mu}f}^{ij}
		+
		8\pi\mu f.
		\label{eq:explicit-landau-K}
	\end{align}
\end{prop}

\begin{proof}
	Set $G:=\sqrt{\mu}\,g$, so that
	$\sigma_G^{ij}=\sigma_{\sqrt{\mu}g}^{ij}$.
	Since
	\[
	\phi^{ij}=D_{ij}|v|\qquad \text{and}\qquad \Delta_v^2|v|=\Delta_v\left(\frac{2}{|v|}\right)=-8\pi\delta_0
	\]
	in the sense of distributions, for every $H\in\cS(\bR_v^3)$, we have
	\begin{equation}\label{eq:double-divergence-sigma}
		D_{ij}\sigma_H^{ij}
		=
		(D_{ij}\phi^{ij})\ast H
		=
		\left(\Delta_v^2|v|\right)\ast H
		=
		-8\pi H.
	\end{equation}
	
	Since
	\[
	D_i\sqrt{\mu}=-v_i\sqrt{\mu},
	\qquad
	\phi^{ij}\ast D_jG=D_j\sigma_G^{ij},
	\]
	the definition of $N[g,f]$ gives
	\begin{align*}
		N[g,f]
		&=
		\frac{1}{\sqrt{\mu}}
		D_i\left(
		\sigma_G^{ij}D_j(\sqrt{\mu}f)
		-
		(D_j\sigma_G^{ij})\sqrt{\mu}f
		\right)=
		(D_i-v_i)
		\left[
		\sigma_G^{ij}D_jf
		-
		\left(
		v_j\sigma_G^{ij}
		+
		D_j\sigma_G^{ij}
		\right)f
		\right].
	\end{align*}
	Expanding this expression and using the symmetry
	$\sigma_G^{ij}=\sigma_G^{ji}$, we obtain
	\begin{align*}
		N[g,f]
		={}&
		D_i\left(\sigma_G^{ij}D_jf\right)
		-
		\left(
		D_j\sigma_G^{ij}
		+
		2v_j\sigma_G^{ij}
		\right)D_i f+
		\left[
		(v_iv_j-\delta_{ij})\sigma_G^{ij}
		-
		D_{ij}\sigma_G^{ij}
		\right]f.
	\end{align*}
	Formula \eqref{eq:explicit-landau-nonlinearity} follows from
	\eqref{eq:double-divergence-sigma}.
	
	We next take $g=\sqrt{\mu}$, so that $G=\mu$ and
	$Af=N[\sqrt{\mu},f]$. Since
	\[
	D_j\mu=-2v_j\mu,
	\qquad
	\phi^{ij}(z)z_j=0,
	\]
	we have
	\begin{equation}\label{eq:divergence-maxwellian-sigma}
		D_j\sigma^{ij}
		=
		-2v_j\sigma^{ij}
		=
		-2\sigma^i.
	\end{equation}
	Together with \eqref{eq:double-divergence-sigma}, this yields
	\[
	D_i\sigma^i
	=
	-\frac12D_{ij}\sigma^{ij}
	=
	4\pi\mu.
	\]
	On the other hand,
	\begin{align*}
		D_i\sigma^i
		=
		D_i(v_j\sigma^{ij})
		=
		\delta_{ij}\sigma^{ij}
		+
		v_jD_i\sigma^{ij}
		=
		\delta_{ij}\sigma^{ij}
		-
		2v_iv_j\sigma^{ij}.
	\end{align*}
	Consequently,
	\begin{equation}\label{eq:trace-maxwellian-sigma}
		\delta_{ij}\sigma^{ij}
		=
		2v_iv_j\sigma^{ij}
		+
		4\pi\mu.
	\end{equation}
	Applying \eqref{eq:explicit-landau-nonlinearity} with
	$g=\sqrt{\mu}$ and using
	\eqref{eq:divergence-maxwellian-sigma} and
	\eqref{eq:trace-maxwellian-sigma}, we obtain
	\begin{align*}
		Af
		&=
		D_i\left(\sigma^{ij}D_jf\right)
		+
		\left[
		8\pi\mu
		+
		(v_iv_j-\delta_{ij})\sigma^{ij}
		\right]f
		=
		D_i\left(\sigma^{ij}D_jf\right)
		-
		v_iv_j\sigma^{ij}f
		+
		4\pi\mu f.
	\end{align*}
	
	Finally, applying \eqref{eq:explicit-landau-nonlinearity} to
	$N[f,\sqrt{\mu}]=Kf$ gives
	\begin{align*}
		Kf
		={}&
		D_i\left(
		\sigma_{\sqrt{\mu}f}^{ij}D_j\sqrt{\mu}
		\right)
		-
		\left(
		D_j\sigma_{\sqrt{\mu}f}^{ij}
		+
		2v_j\sigma_{\sqrt{\mu}f}^{ij}
		\right)D_i\sqrt{\mu}
		+
		\left(
		8\pi\sqrt{\mu}\,f
		+
		(v_iv_j-\delta_{ij})
		\sigma_{\sqrt{\mu}f}^{ij}
		\right)\sqrt{\mu}.
	\end{align*}
	Using $D_i\sqrt{\mu}=-v_i\sqrt{\mu}$, we find
	\begin{align*}
		Kf
		={}&
		\sqrt{\mu}
		\left[
		(4v_iv_j-2\delta_{ij})
		\sigma_{\sqrt{\mu}f}^{ij}
		+
		v_iD_j\sigma_{\sqrt{\mu}f}^{ij}
		-
		v_jD_i\sigma_{\sqrt{\mu}f}^{ij}
		\right]
		+
		8\pi\mu f.
	\end{align*}
	The last two terms in the bracket cancel after exchanging the dummy
	indices $i$ and $j$. This proves
	\eqref{eq:explicit-landau-K}.
\end{proof}

\begin{corollary}[Divergence-form expansion of $N$]\label{cor:expansion-N}
	For any $g$ with $\mu^{1/4}g\in L^\infty(\bR_v^3)$, the term
	$N[g,\cdot]$ admits the divergence-form expansion
	\begin{equation}\label{eq:N-divergence-expansion}
		N[g,f]
		=
		\p_i\!\left(\sigma_{\sqrt{\mu}g}^{ij}\,\p_jf\right)
		+
		\mathfrak b_g^i\,\p_i f
		+
		\mathfrak c_g\,f,
	\end{equation}
	where
	\[
	\mathfrak b_g^i
	:=
	-\p_j\sigma_{\sqrt{\mu}g}^{ij}
	-2\sigma_{\sqrt{\mu}g}^{i},
	\qquad
	\mathfrak c_g
	:=
	(v_iv_j-\delta_{ij})\sigma_{\sqrt{\mu}g}^{ij}
	+8\pi\sqrt{\mu}\,g,
	\]
	and the coefficients satisfy, for every $v\in\bR^3$, the pointwise
	bounds
	\begin{equation}\label{eq:N-coefficient-bounds}
		|\sigma_{\sqrt{\mu}g}^{ij}(v)|
		\le
		C\|\mu^{1/4}g\|_{L^\infty}\la v\ra^{-1},
		\quad
		|\mathfrak b_g^i(v)|
		\le
		C\|\mu^{1/4}g\|_{L^\infty}\la v\ra^{-1},
		\quad
		|\mathfrak c_g(v)|
		\le
		C\|\mu^{1/4}g\|_{L^\infty}.
	\end{equation}
\end{corollary}

\begin{proof}
	Identity \eqref{eq:N-divergence-expansion} is
	\eqref{eq:explicit-landau-nonlinearity} rewritten with
	$v_j\sigma_{\sqrt{\mu}g}^{ij}=\sigma_{\sqrt{\mu}g}^{i}$, and the
	coefficient bounds \eqref{eq:N-coefficient-bounds} follow from
	Lemma~\ref{lem:perturbed-landau-ellipticity}: the bound for
	$\sigma_{\sqrt{\mu}g}^{ij}$ from \eqref{eq:perturbation-sigma-bound} and
	$\mathfrak q(v,\xi)\le\la v\ra^{-1}|\xi|^2$, the bound for
	$\mathfrak b_g^i$ from \eqref{eq:perturbed-sigma-i-bound} and
	\eqref{eq:derivative-perturbed-sigma-bound}, and the bound for
	$\mathfrak c_g$ from \eqref{eq:perturbed-sigma-i-bound} for the term
	$v_iv_j\sigma^{ij}_{\sqrt\mu g}$, from \eqref{eq:perturbation-sigma-bound}
	with $\xi=\bar\xi=e_i$ for the trace term
	$\delta_{ij}\sigma^{ij}_{\sqrt\mu g}$, and from
	$|\sqrt{\mu}\,g|\le\|\mu^{1/4}g\|_{L^\infty}$ for the last term.
\end{proof}

For $v\in\bR^3$, let $P_v$ denote the orthogonal projection onto the
radial direction and $Q_v$ the complementary projection onto $v^\perp$:
\begin{equation}
	\label{eq:radial-tangential-projections}
	P_v:=\frac{v\otimes v}{|v|^2}
	\quad\text{if }v\neq0,
	\qquad
	P_0:=0,
	\qquad
	Q_v:=I_3-P_v,
\end{equation}
so that $P_v\xi=\dfrac{\xi\cdot v}{|v|^2}\,v$.
Define
\begin{align}\label{26720259}
	\mathfrak q(v,\xi)
	:=
	(1+|v|)^{-3}|P_v\xi|^2
	+
	(1+|v|)^{-1}|Q_v\xi|^2.
\end{align}

By \eqref{new:eq:intro-coulomb-kernel} one has $\phi^{ij}(Rz)=R\,\phi^{ij}(z)\,R^T$ for
every $R\in SO(3)$, and $\mu$ is radial, so the substitution $w\mapsto Rw$ in
$\sigma^{ij}=\phi^{ij}\ast\mu$ gives $\sigma(Rv)=R\,\sigma(v)\,R^T$.  Choosing
$R$ in the stabilizer of $v$ shows that $\sigma(v)$ commutes with every
rotation about $v$, hence is diagonal in the frame $(P_v,Q_v)$.  Writing
$r=|v|$ and
\begin{equation}
	\label{eq:coulomb-eigenvalue-definition}
	\lambda_\parallel(r):=e_1^T\sigma(re_1)e_1,
	\qquad
	\lambda_\perp(r):=e_2^T\sigma(re_1)e_2,
\end{equation}
we obtain the radial decomposition
\begin{equation}
	\label{eq:coulomb-eigen-decomposition}
	\sigma(v)=\lambda_\parallel(|v|)P_v+\lambda_\perp(|v|)Q_v,
	\qquad v\in\bR^3,
\end{equation}
the parallel eigenvalue being simple and the perpendicular one double; at
$v=0$ the conventions $P_0=0$, $Q_0=I_3$ and
$\lambda_\parallel(0)=\lambda_\perp(0)$ are consistent with $\sigma(0)$ being
a multiple of the identity.

\begin{lemma}[Coulomb eigenvalues and logarithmic symbol bounds]
	\label{lem:coulomb-symbol-bounds}
	The matrix $\sigma$ belongs to $C^\infty(\bR^3)$ and is positive definite at
	every $v\in\bR^3$; in particular $\lambda_\parallel$ and $\lambda_\perp$ are
	smooth and strictly positive on $[0,\infty)$.  There are constants $c,C>0$
	and $C_k>0$, $k\in\bN$, such that for all $r\ge0$,
	\begin{align}
		c\la r\ra^{-3}\le\lambda_\parallel(r)\le C\la r\ra^{-3},
		\qquad
		c\la r\ra^{-1}\le\lambda_\perp(r)\le C\la r\ra^{-1},
		\label{eq:eigenvalue-size}
	\end{align}
	and
	\begin{align}
		\abs{(r\p_r)^k\lambda_\parallel(r)}
		\le C_k\lambda_\parallel(r),
		\qquad
		\abs{(r\p_r)^k\lambda_\perp(r)}
		\le C_k\lambda_\perp(r).
		\label{eq:eigenvalue-log-derivatives}
	\end{align}
\end{lemma}

\begin{proof}
	By rotation invariance it suffices to work at $v=re_1$, and the size bounds
	\eqref{eq:eigenvalue-size} are proved in
	\cite[Propositions 2.3 and 2.4]{degond1997dispersion}; it therefore remains
	to prove the smoothness and the positivity of $\sigma$, and to establish
	\eqref{eq:eigenvalue-log-derivatives}.
	
	\emph{Smoothness and positivity.}
	The smoothness of $\sigma$
	follows from the smoothness of the Gaussian $\mu$ and its exponential decay
	at infinity.
	To prove the positivity of $\sigma$, let $\xi\ne0$.
	Then, by \eqref{new:eq:intro-coulomb-kernel},
	\begin{equation*}
		\xi^T\sigma(v)\xi
		=\int_{\bR^3}\xi_i\xi_j\phi^{ij}(z)\mu(v-z)\dd z=\int_{\bR^3}\frac{|Q_z \xi|^2}{|z|}\mu(v-z)\dd z.
	\end{equation*}
	The integrand is nonnegative and $\mu$ is strictly positive, while
	$|Q_z\xi|^2$ vanishes only on the line $\bR\xi$, a set of
	three-dimensional measure zero; hence the integral is strictly positive.
	Thus $\sigma(v)$ is positive definite for every
	$v$, and by \eqref{eq:coulomb-eigenvalue-definition} the functions
	$\lambda_\parallel,\lambda_\perp$ are smooth and strictly positive on
	$[0,\infty)$.  At $v=0$ rotation invariance gives $\sigma(0)=\lambda_0I_3$
	with $\lambda_0>0$, so that $\lambda_\parallel(0)=\lambda_\perp(0)=\lambda_0$.
	
	\emph{Bounded velocities.}
	Throughout $D:=r\p_r$, and $\lambda$ denotes either one of the eigenvalues
	$\lambda_\parallel,\lambda_\perp$.  By induction,
	\begin{equation}\label{eq:log-derivative-expansion}
		D^k=\sum_{j=1}^kS(k,j)\,r^j\p_r^j,
		\qquad k\in\bN,
	\end{equation}
	with $S(k,j)$ the Stirling numbers of the second kind.  Since $\lambda$ is
	smooth, each term of \eqref{eq:log-derivative-expansion} applied to
	$\lambda$ is bounded on $0\le r\le4$, while $\lambda\ge c>0$ there by
	positivity and continuity; this is
	\eqref{eq:eigenvalue-log-derivatives} for $r\le4$.
	
	\emph{Large velocities: bulk part.}
	Let $r\ge4$.  Directly from the kernel,
	\begin{align}
		\lambda_\parallel(r)
		&=\int_{\bR^3}\mu(w)
		\frac{|w_\perp|^2}{|re_1-w|^3}\dd w,
		\qquad w_\perp:=(w_2,w_3),
		\label{eq:lambda-parallel-integral}\\
		\lambda_\perp(r)
		&=\int_{\bR^3}\mu(w)
		\left\{
		\frac1{|re_1-w|}-\frac{w_2^2}{|re_1-w|^3}
		\right\}\dd w.
		\label{eq:lambda-perp-integral}
	\end{align}
	Choose $\chi\in C_c^\infty(\bR^3)$ with $\chi=1$ on $|w|\le1/4$ and $\chi=0$
	on $|w|\ge1/3$, and put $\chi_r(w):=\chi(w/r)$.  Split each of the integrals
	\eqref{eq:lambda-parallel-integral}--\eqref{eq:lambda-perp-integral} into its
	bulk part and its tail part, defined by
	$$
	\int_{\bR^3}(\cdots)\chi_r(w)\dd w\qquad \text{and}\qquad
	\int_{\bR^3}(\cdots)\left(1-\chi_r(w)\right)\dd w,
	$$
	respectively.
	
	On $\supp\chi_r$ one has $|w|\le r/3$, hence $|re_1-w|\simeq r$.  Induction
	from
	\[
	D|re_1-w|^{-\alpha}
	=-\alpha\frac{r(r-w_1)}{|re_1-w|^2}|re_1-w|^{-\alpha}
	\]
	gives
	\[
	D^k|re_1-w|^{-\alpha}
	=P_{k,\alpha}\!\left(
	\frac{r(r-w_1)}{|re_1-w|^2},
	\frac{r^2}{|re_1-w|^2}
	\right)|re_1-w|^{-\alpha},
	\]
	where $P_{k,\alpha}$ is a polynomial.  Applying this with
	$\alpha=1,3$ gives, for every
	$k\in\bN$ and every $|w|\le r/3$,
	\begin{align}
		\abs{D^k
			\frac{|w_\perp|^2}{|re_1-w|^3}}
		&\le C_k r^{-3}|w_\perp|^2,
		\label{eq:near-par-derivative}\\
		\abs{D^k
			\left(
			\frac1{|re_1-w|}-\frac{w_2^2}{|re_1-w|^3}
			\right)}
		&\le C_k r^{-1}(1+|w|^2r^{-2}).
		\label{eq:near-perp-derivative}
	\end{align}
	The cutoff derivatives produce no power loss.
	Indeed, putting $y=w/r$, one has
	$$
	D^j\chi_r(w)=(-y\cdot\nabla_y)^j\chi(y)\big|_{y=w/r},
	$$
	so that $|D^j\chi_r|\le C_j$ with support in $r/4\le|w|\le r/3$.
	Expanding $D^k$ of the product by the Leibniz rule and integrating against
	$\mu$, the terms carrying a cutoff derivative are supported in
	$r/4\le|w|\le r/3$, where $\mu(w)\le\ee^{-r^2/16}$, hence are bounded by a
	polynomial in $r$ times $\ee^{-r^2/16}$, and the remaining ones
	obey \eqref{eq:near-par-derivative}--\eqref{eq:near-perp-derivative}.  Hence
	the bulk parts and all their logarithmic derivatives are bounded by
	$C_kr^{-3}$ in the parallel case and by $C_kr^{-1}$ in the perpendicular
	case.
	
	\emph{Large velocities: tail part.}
	Write the tail contribution to the matrix as
	\[
	I^{ij}(r)
	:=\int_{\bR^3}\phi^{ij}(re_1-w)\bigl(1-\chi_r(w)\bigr)\mu(w)\dd w
	=\int_{\bR^3}\phi^{ij}(z)\,g_r(z)\dd z,
	\]
	where
	\[
	g_r(z):=\Bigl(1-\chi\bigl(\tfrac{re_1-z}{r}\bigr)\Bigr)\mu(re_1-z).
	\]
	Since $\chi\equiv 1$ on $|w|\le\frac14$, the factor
	$1-\chi\bigl((re_1-z)/r\bigr)$ vanishes unless $|re_1-z|\ge\frac r4$, so
	that $g_r$ is supported in $\{|re_1-z|\ge r/4\}$; moreover, since
	$\chi\in C_c^\infty(\bR^3)$, every $r$-derivative of
	$\chi\bigl((re_1-z)/r\bigr)$ is bounded uniformly for $r\ge4$.
	Consequently, since $\mu(y)=\ee^{-|y|^2}$,
	\[
	\abs{\p_r^mg_r(z)}
	\le C_m\la re_1-z\ra^{m}\one_{\{|re_1-z|\ge r/4\}}\ee^{-|re_1-z|^2}
	\le C_m\ee^{-r^2/64}\ee^{-|re_1-z|^2/2},
	\qquad m\in\bN.
	\]
	Differentiating under the integral sign, which is legitimate by the
	domination just obtained, and using $|\phi^{ij}(z)|\le2|z|^{-1}$ together
	with $\int_{\bR^3}|z|^{-1}\ee^{-|re_1-z|^2/2}\dd z\le C$ uniformly in $r$, we get
	$\abs{\p_r^kI^{ij}(r)}\le C_k\ee^{-r^2/64}$ for every $k\in\bN$.  By
	\eqref{eq:log-derivative-expansion}, every polynomial factor is absorbed by a
	smaller Gaussian exponent, so that
	\begin{equation}\label{eq:far-exponential}
		\abs{D^kI^{ij}(r)}
		\le C_k\ee^{-c_kr^2},
		\qquad k\in\bN.
	\end{equation}
	The tail contributions to $\lambda_\parallel$ and $\lambda_\perp$ are
	$e_1^TI(r)e_1$ and $e_2^TI(r)e_2$.  Combining the bulk estimates,
	\eqref{eq:far-exponential}, and the lower bounds of
	\eqref{eq:eigenvalue-size} proves
	\eqref{eq:eigenvalue-log-derivatives} for $r\ge4$, which completes the
	proof.
\end{proof}

\begin{lemma}[Eigenvalues of the Landau diffusion matrix]\label{lem:landau-diffusion-eigenvalues}
	For every $\varepsilon>0$, there exists a constant
	$C_\varepsilon\ge1$, depending only on $\varepsilon$, such that
	\begin{equation}\label{eq:coulomb-maxwellian-anisotropy}
		C_\varepsilon^{-1}\mathfrak q(v,\xi)
		\le
		\sigma_{\mu^\varepsilon}^{ij}(v)\xi_i\xi_j
		\le
		C_\varepsilon\mathfrak q(v,\xi)
	\end{equation}
	for every $v,\xi\in\bR^3$.  In particular, taking $\varepsilon=1$ and
	$\xi=v$, so that $P_vv=v$ and $Q_vv=0$, we obtain
	\begin{equation}\label{eq:coulomb-sigmavv-comparison}
		C^{-1}\la v\ra^{-1}
		\le v_i\sigma^{ij}(v)v_j+\la v\ra^{-3}
		\le C\la v\ra^{-1},
	\end{equation}
	where the added $\la v\ra^{-3}$ covers bounded velocities.
\end{lemma}

\begin{proof}
	For $\varepsilon=1$, \eqref{eq:coulomb-eigen-decomposition} gives
	$\sigma^{ij}(v)\xi_i\xi_j
	=\lambda_\parallel(|v|)|P_v\xi|^2+\lambda_\perp(|v|)|Q_v\xi|^2$, which is
	comparable to $\mathfrak q(v,\xi)$ by \eqref{eq:eigenvalue-size} and
	$\la v\ra\simeq1+|v|$.  For general $\varepsilon>0$, a direct change of
	variables and the homogeneity of $\phi^{ij}$ give
	\[
	\sigma_{\mu^\varepsilon}^{ij}(v)
	=
	\varepsilon^{-1}
	\sigma^{ij}(\sqrt{\varepsilon}\,v),
	\]
	and since
	\[
	P_{\sqrt{\varepsilon}v}=P_v,
	\qquad
	1+\sqrt{\varepsilon}|v|
	\simeq_\varepsilon
	1+|v|,
	\]
	the result follows from the case $\varepsilon=1$.  Finally, $\xi=v$ in
	\eqref{26720259} gives $\mathfrak q(v,v)=(1+|v|)^{-3}|v|^2$, so
	\eqref{eq:coulomb-maxwellian-anisotropy} with $\varepsilon=1$ yields
	$v_i\sigma^{ij}(v)v_j\simeq\la v\ra^{-3}|v|^2$; adding $\la v\ra^{-3}$ and
	using $|v|^2+1=\la v\ra^2$ gives
	\eqref{eq:coulomb-sigmavv-comparison}.
\end{proof}

\begin{lemma}[Perturbative bounds for the Landau diffusion matrix]
	\label{lem:perturbed-landau-ellipticity}
	Let $g$ be measurable on $\bR_v^3$ and assume that
	\[
	\mu^{1/4}g\in L^\infty(\bR_v^3).
	\]
	Then there exists a universal constant $C>0$ such that the following
	statements hold.
	
	\begin{enumerate}
		\item
		For every $v,\xi,\bar\xi\in\bR^3$,
		\begin{equation}\label{eq:perturbation-sigma-bound}
			\left|
			\sigma_{\sqrt\mu g}^{ij}(v)\xi_i\bar\xi_j
			\right|
			\le
			C
			\|\mu^{1/4}g\|_{L^\infty(\bR_v^3)}
			\mathfrak q(v,\xi)^{1/2}
			\mathfrak q(v,\bar\xi)^{1/2}.
		\end{equation}
		Consequently, there exists a universal constant
		$\varepsilon_{\mathrm{ell}}>0$ such that, if
		\begin{equation}\label{eq:smallness-perturbed-sigma}
			\|\mu^{1/4}g\|_{L^\infty(\bR_v^3)}
			\le
			\varepsilon_{\mathrm{ell}},
		\end{equation}
		then     for every $v,\xi\in\bR^3$,
		\begin{equation}\label{260502449}
			\frac12
			\sigma^{ij}(v)\xi_i\xi_j
			\le
			\sigma_{\mu+\sqrt\mu g}^{ij}(v)\xi_i\xi_j
			\le
			2
			\sigma^{ij}(v)\xi_i\xi_j.
		\end{equation}
		
		\item
		For every $v,\xi\in\bR^3$,
		\begin{equation}\label{eq:perturbed-sigma-i-bound}
			\left|
			v_j\sigma_{\sqrt\mu g}^{ij}(v)\xi_i
			\right|
			\le
			C
			\|\mu^{1/4}g\|_{L^\infty(\bR_v^3)}
			\langle v\rangle^{-1}|\xi|,\qquad \left|
			v_iv_j\sigma_{\sqrt\mu g}^{ij}(v)
			\right|
			\le
			C
			\|\mu^{1/4}g\|_{L^\infty(\bR_v^3)}
			\langle v\rangle^{-1}.
		\end{equation}
		
		\item
		For every $1\le i,j,k\le3$ and $v\in\bR^3$,
		\begin{equation}\label{eq:derivative-perturbed-sigma-bound}
			\left|
			\p_k\sigma_{\sqrt\mu g}^{ij}(v)
			\right|
			\le
			C
			\|\mu^{1/4}g\|_{L^\infty(\bR_v^3)}
			\langle v\rangle^{-2}.
		\end{equation}
	\end{enumerate}
\end{lemma}

\begin{proof}
	Set
	\[
	M:=\|\mu^{1/4}g\|_{L^\infty(\bR_v^3)}.
	\]
	Since $|\sqrt\mu g|\le M\mu^{1/4}$ and $\phi(z)$ is nonnegative
	definite, the matrix Cauchy--Schwarz inequality followed by the
	integral Cauchy--Schwarz inequality gives
	\begin{align*}
		\left|
		\sigma_{\sqrt\mu g}^{ij}(v)\xi_i\bar\xi_j
		\right|
		&\le
		M
		\int_{\bR^3}
		\mu(v_*)^{1/4}
		\left(
		\phi^{ij}(v-v_*)\xi_i\xi_j
		\right)^{1/2}
		\left(
		\phi^{ij}(v-v_*)\bar\xi_i\bar\xi_j
		\right)^{1/2}
		\dd v_*
		\\
		&\le
		M
		\left(
		\sigma_{\mu^{1/4}}^{ij}(v)\xi_i\xi_j
		\right)^{1/2}
		\left(
		\sigma_{\mu^{1/4}}^{ij}(v)\bar\xi_i\bar\xi_j
		\right)^{1/2}.
	\end{align*}
	Lemma~\ref{lem:landau-diffusion-eigenvalues}, applied with
	$\varepsilon=1/4$, now proves
	\eqref{eq:perturbation-sigma-bound}.
	
	The same lemma, applied with $\varepsilon=1$, yields
	$\mathfrak q(v,\xi)
	\le
	C\sigma^{ij}(v)\xi_i\xi_j$.
	Thus
	\[
	\left|
	\sigma_{\sqrt\mu g}^{ij}(v)\xi_i\xi_j
	\right|
	\le
	CM\sigma^{ij}(v)\xi_i\xi_j.
	\]
	Choosing $\varepsilon_{\mathrm{ell}}>0$ so that $C\varepsilon_{\mathrm{ell}}\le1/2$, we
	obtain
	\[
	\frac12\sigma^{ij}\xi_i\xi_j
	\le
	\sigma_{\mu+\sqrt\mu g}^{ij}\xi_i\xi_j
	\le
	\frac32\sigma^{ij}\xi_i\xi_j,
	\]
	which implies \eqref{260502449}.
	
	\eqref{eq:perturbed-sigma-i-bound} directly follows from \eqref{eq:perturbation-sigma-bound} and the definition of $\mathfrak q(v,\xi)$ in \eqref{26720259}.
	
	Finally, since $|\p_k\phi^{ij}(z)|\le C|z|^{-2}$, we have
	\begin{align*}
		\left|
		\p_k\sigma_{\sqrt\mu g}^{ij}(v)
		\right|
		&\le
		CM
		\int_{\bR^3}
		\frac{\mu(v_*)^{1/4}}{|v-v_*|^2}
		\dd v_*.
	\end{align*}
	The standard Gaussian convolution estimate
	\[
	\int_{\bR^3}
	\frac{\mu(v_*)^{1/4}}{|v-v_*|^2}
	\dd v_*
	\le
	C\langle v\rangle^{-2}
	\]
	follows by splitting into $|v_*|\le|v|/2$ and
	$|v_*|>|v|/2$, with a further split at $|v-v_*|=1$ in the latter
	region. This proves \eqref{eq:derivative-perturbed-sigma-bound}.
\end{proof}

\begin{lemma}\label{260506438}
	There exists a universal constant $c\in(0,1)$ such that, for every
	$v_0,v,\xi\in\bR^3$ satisfying
	\[
	|v-v_0|
	<
	\frac{c}{1+|v_0|},
	\]
	we have
	\begin{equation}\label{eq:local-comparability-sigma}
		\frac12
		\sigma^{ij}(v_0)\xi_i\xi_j
		\le
		\sigma^{ij}(v)\xi_i\xi_j
		\le
		2
		\sigma^{ij}(v_0)\xi_i\xi_j.
	\end{equation}
\end{lemma}

\begin{proof}
	Taking $g=\sqrt{\mu}$ in
	\eqref{eq:derivative-perturbed-sigma-bound}, we obtain
	\begin{equation}\label{eq:derivative-sigma-bound}
		|\p_k\sigma^{ij}(w)|
		\le
		C(1+|w|)^{-2}.
	\end{equation}
	Let
	\[
	v_t:=v_0+t(v-v_0),
	\qquad
	0\le t\le1.
	\]
	If $c\le1/2$, then
	\[
	1+|v_t|
	\ge
	\frac12(1+|v_0|).
	\]
	Therefore, by the fundamental theorem of calculus and
	\eqref{eq:derivative-sigma-bound},
	\begin{align*}
		\left|
		\big(
		\sigma^{ij}(v)-\sigma^{ij}(v_0)
		\big)\xi_i\xi_j
		\right|
		&\le
		|v-v_0|
		\int_0^1
		|\nabla_v\sigma^{ij}(v_t)|
		\dd t\,|\xi|^2
		\\
		&\le
		C|v-v_0|(1+|v_0|)^{-2}|\xi|^2
		\le
		Cc(1+|v_0|)^{-3}|\xi|^2.
	\end{align*}
	On the other hand, Lemma~\ref{lem:landau-diffusion-eigenvalues} with
	$\varepsilon=1$ gives
	\[
	\sigma^{ij}(v_0)\xi_i\xi_j
	\ge
	C^{-1}\mathfrak q(v_0,\xi)
	\ge
	C^{-1}(1+|v_0|)^{-3}|\xi|^2.
	\]
	Consequently,
	\[
	\left|
	\big(
	\sigma^{ij}(v)-\sigma^{ij}(v_0)
	\big)\xi_i\xi_j
	\right|
	\le
	Cc\,
	\sigma^{ij}(v_0)\xi_i\xi_j.
	\]
	Choosing the universal constant $c>0$ sufficiently small, we obtain
	\[
	\left|
	\sigma^{ij}(v)\xi_i\xi_j
	-
	\sigma^{ij}(v_0)\xi_i\xi_j
	\right|
	\le
	\frac12
	\sigma^{ij}(v_0)\xi_i\xi_j,
	\]
	which proves \eqref{eq:local-comparability-sigma}.
\end{proof}

\subsubsection{Anisotropic Sobolev spaces associated with the Landau operator}

\begin{defn}\label{def:H_sigma}
	We define
	\begin{equation}\label{def:H1sigma-space}
		H_\sigma^1(\bR_v^3)
		:=
		\left\{
		f\in H_{\mathrm{loc}}^1(\bR_v^3)
		:
		\|f\|_{H_\sigma^1}<\infty
		\right\},
	\end{equation}
	where
	\begin{equation}\label{eq:H_sigma_norm}
		\|f\|_{H_\sigma^1}^2
		:=
		\int_{\bR^3}
		\left(
		\sigma^{ij}\p_i f\,\p_j f
		+
		\sigma^{ij}v_iv_j f^2
		\right)
		\dd v.
	\end{equation}
	We denote by
	\begin{equation}\label{def:Hminus1sigma-space}
		H_\sigma^{-1}(\bR_v^3)
		:=
		\big(H_\sigma^1(\bR_v^3)\big)^*
	\end{equation}
	the dual space of $H_\sigma^1(\bR_v^3)$, equipped with the norm
	\begin{equation}\label{def:Hminus1sigma-norm}
		\|F\|_{H_\sigma^{-1}}
		:=
		\sup_{\substack{\varphi\in H_\sigma^1\\
				\|\varphi\|_{H_\sigma^1}\le1}}
		\left|
		\langle F,\varphi\rangle_{H_\sigma^{-1},H_\sigma^1}
		\right|.
	\end{equation}
\end{defn}

\begin{lemma}[Weighted coercive norm equivalence]
	\label{lem:weighted-Hsigma-equivalence}
	Write $w(v):=\la v\ra:=(1+|v|^2)^{1/2}$.  For each fixed $\ell\ge0$, there is $C_\ell\ge1$ such
	that, for every $f\in H_\sigma^1(\bR_v^3)$,
	\begin{equation}\label{eq:weighted-Hsigma-equivalence}
		C_\ell^{-1}\norm{w^\ell f}_{H^1_\sigma}^2
		\le\int_{\bR^3}w^{2\ell}\left(
		\sigma^{ij}\p_i f\,\p_j f
		+
		\sigma^{ij}v_iv_j f^2
		\right)\dd v
		\le C_\ell\norm{w^\ell f}_{H^1_\sigma}^2.
	\end{equation}
\end{lemma}

\begin{proof}
	Write $\abs\xi_\sigma^2:=\sigma^{ij}(v)\xi_i\xi_j$.  By \eqref{eq:H_sigma_norm}
	the integrand in \eqref{eq:weighted-Hsigma-equivalence} is
	\[
	D_f:=\abs{\nabla f}_\sigma^2+\sigma^{ij}v_iv_jf^2,
	\qquad\text{while}\qquad
	\norm{w^\ell f}_{H^1_\sigma}^2
	=\int_{\bR^3}\left(
	\abs{\nabla(w^\ell f)}_\sigma^2+\sigma^{ij}v_iv_j\,w^{2\ell}f^2
	\right)\dd v .
	\]
	The zeroth-order terms of $\norm{w^\ell f}_{H^1_\sigma}^2$ and of
	$\int w^{2\ell}D_f\dd v$ coincide, so only the gradient terms have to be compared.
	Since $\nabla\log w^\ell=\ell v/\la v\ra^2$,
	\begin{equation}\label{eq:weight-gradient-sigma}
		\abs{\nabla\log w^\ell}_\sigma^2
		=\frac{\ell^2}{\la v\ra^4}\,\sigma^{ij}v_iv_j
		\le\ell^2\,\sigma^{ij}v_iv_j.
	\end{equation}
	The product rule
	gives
	\[
	\nabla(w^\ell f)=w^\ell\left(\nabla f+f\nabla\log w^\ell\right),
	\qquad
	w^\ell\nabla f=\nabla(w^\ell f)-(w^\ell f)\nabla\log w^\ell ,
	\]
	so that $\abs{A+B}_\sigma^2\le2\abs A_\sigma^2+2\abs B_\sigma^2$ together with
	\eqref{eq:weight-gradient-sigma} yields the pointwise bounds
	\[
	\abs{\nabla(w^\ell f)}_\sigma^2
	\le2w^{2\ell}\abs{\nabla f}_\sigma^2+2\ell^2w^{2\ell}\sigma^{ij}v_iv_jf^2,
	\qquad
	w^{2\ell}\abs{\nabla f}_\sigma^2
	\le2\abs{\nabla(w^\ell f)}_\sigma^2+2\ell^2w^{2\ell}\sigma^{ij}v_iv_jf^2 .
	\]
	Adding the common zeroth-order term to each and integrating gives
	\eqref{eq:weighted-Hsigma-equivalence} with $C_\ell=2\ell^2+2$.
\end{proof}

\begin{lemma}\label{lem:local-embedding-H_sigma}
	For every $f\in H_\sigma^1(\bR_v^3)$,
	\begin{equation}\label{eq:local-H1-H_sigma}
		\|\la v\ra^{-1/2}f\|_{L^2(\bR_v^3)}+\|\la v\ra^{-3/2}\nabla_vf\|_{L^2(\bR_v^3)}
		\le
		C\|f\|_{H_\sigma^1},
	\end{equation}
	where $\la v\ra:=(1+|v|^2)^{1/2}$.
	Consequently,
	\begin{equation}\label{eq:local-L6-H_sigma}
		\|\la v\ra^{-3/2}f\|_{L^6(\bR_v^3)}
		\le
		C\|f\|_{H_\sigma^1}.
	\end{equation}
	Dually, every measurable $F$ with $\la v\ra^{1/2}F\in L^2(\bR_v^3)$
	defines an element of $H_\sigma^{-1}(\bR_v^3)$, and
	\begin{equation}\label{eq:Hsigma-dual-weighted}
		\|F\|_{H_\sigma^{-1}}
		\le
		C\|\la v\ra^{1/2}F\|_{L^2(\bR_v^3)} .
	\end{equation}
\end{lemma}

\begin{proof}
	Fix $R\ge1$. By Lemma~\ref{lem:landau-diffusion-eigenvalues}, there exists
	$c_R>0$ such that for any $v\in B_{2R}$ and $\xi\in\bR^3$,
	$\sigma^{ij}(v)\xi_i\xi_j\ge c_R|\xi|^2$.
	Therefore,
	\begin{equation}\label{eq:local-gradient-H_sigma}
		\|\nabla_v f\|_{L^2(B_{2R})}
		\le
		C_R\|f\|_{H_\sigma^1}.
	\end{equation}
	
	Set
	\[
	A_R:=B_{2R}\setminus\overline{B_R}.
	\]
	Since $|v|\ge R$ on $A_R$, Lemma~\ref{lem:landau-diffusion-eigenvalues}
	also gives
	\[
	\sigma^{ij}(v)v_iv_j
	\ge
	c_R
	\qquad
	\text{on }A_R.
	\]
	Hence
	\begin{equation}\label{eq:annular-L2-H_sigma}
		\|f\|_{L^2(A_R)}
		\le
		C_R\|f\|_{H_\sigma^1}.
	\end{equation}
	
	The Poincar\'e inequality with the average of $f$ over $A_R$ yields
	\[
	\|f\|_{L^2(B_R)}
	\le
	C_R
	\left(
	\|\nabla_v f\|_{L^2(B_{2R})}
	+
	\|f\|_{L^2(A_R)}
	\right).
	\]
	Combining this estimate with
	\eqref{eq:local-gradient-H_sigma} and
	\eqref{eq:annular-L2-H_sigma} proves
	\eqref{eq:local-H1-H_sigma}.
	
	Estimate \eqref{eq:Hsigma-dual-weighted} is immediate from
	\eqref{eq:local-H1-H_sigma}. For $\varphi\in H^1_\sigma$,
	\[
	|\la F,\varphi\ra|\le\|\la v\ra^{1/2}F\|_{L^2}\|\la v\ra^{-1/2}\varphi\|_{L^2}
	\le C\|\la v\ra^{1/2}F\|_{L^2}\|\varphi\|_{H^1_\sigma}.
	\]
	
	Estimate \eqref{eq:local-L6-H_sigma} follows by applying the Sobolev
	embedding $\dot H^1(\bR^3)\hookrightarrow L^6(\bR^3)$ to
	$u:=\la v\ra^{-3/2}f$. Indeed,
	\[
	\nabla_vu=\la v\ra^{-3/2}\nabla_vf-\tfrac32\la v\ra^{-7/2}v\,f,
	\]
	so that
	\[|\nabla_vu|\le\la v\ra^{-3/2}|\nabla_vf|+\tfrac32\la v\ra^{-5/2}|f|,
	\]
	and both terms are bounded in $L^2(\bR_v^3)$ by $C\|f\|_{H^1_\sigma}$ on account
	of \eqref{eq:local-H1-H_sigma}.
\end{proof}

\begin{lemma}[Hilbert structure of $H^1_\sigma$]\label{lem:H_sigma-Hilbert}
	\begin{enumerate}
		\item The space $H_\sigma^1(\bR_v^3)$ is a Hilbert space with inner product
		\[
		(f,h)_{H_\sigma^1}
		:=
		\int_{\bR^3}
		\left(
		\sigma^{ij}\p_i f\,\p_j h
		+
		\sigma^{ij}v_iv_jfh
		\right)
		\dd v,
		\]
		and $C_c^\infty(\bR_v^3)$ is dense in $H_\sigma^1(\bR_v^3)$.  More
		precisely, let $\chi\in C_c^\infty(\bR^3)$ be radial with
		\begin{equation}\label{eq:radial-cutoff}
			0\le\chi\le1,
			\qquad
			\chi=1\quad\text{on }B_1,
			\qquad
			\chi=0\quad\text{outside }B_2,
		\end{equation}
		and set $\chi_R(v):=\chi(v/R)$ for $R\ge1$.  Then
		\begin{equation}\label{eq:radial-cutoff-convergence}
			\chi_Rf\longrightarrow f
			\quad\text{in }H_\sigma^1(\bR_v^3)
			\qquad\text{for every }f\in H_\sigma^1(\bR_v^3).
		\end{equation}
		
		\item For every $f\in H_v^1(\bR_v^3)$ and every
		$F\in H^{-1}_\sigma(\bR_v^3)$,
		\begin{align}\label{eq:Hv-into-Hsigma}
			\|f\|_{H^1_\sigma}\le C\|f\|_{H^1_v},
			\qquad
			\|F\|_{H^{-1}_v}\le C\|F\|_{H^{-1}_\sigma} ;
		\end{align}
		in particular $H^1_v(\bR_v^3)\subset H^1_\sigma(\bR_v^3)$ and
		$H^{-1}_\sigma(\bR_v^3)\subset H^{-1}_v(\bR_v^3)$.
		
		\item For every $\zeta\in C^1(\bR_v^3)$ with
		$\|\zeta\|_{C^1(\bR_v^3)}<\infty$, and with
		$\la\zeta F,\varphi\ra:=\la F,\zeta\varphi\ra$,
		\begin{align}\label{eq:Hsigma-multiplier}
			\|\zeta f\|_{H_\sigma^1}\le C\|\zeta\|_{C^1(\bR_v^3)}\|f\|_{H_\sigma^1},
			\qquad
			\|\zeta F\|_{H_\sigma^{-1}}\le C\|\zeta\|_{C^1(\bR_v^3)}\|F\|_{H_\sigma^{-1}}.
		\end{align}
		
		\item Let $R>0$. If $f\in H^1_\sigma(\bR_v^3)$ and
		$F\in H^{-1}_v(\bR_v^3)$ are supported in $B_R$, then
		\begin{align}\label{eq:Hv-Hsigma-compact}
			\|f\|_{H^1_v}\le C_R\|f\|_{H^1_\sigma},
			\qquad
			\|F\|_{H^{-1}_\sigma}\le C_R\|F\|_{H^{-1}_v}.
		\end{align}
		Consequently, for every
		$\zeta\in C_c^\infty(\bR_v^3)$ there is $C_\zeta>0$ with
		\begin{align}\label{eq:local-coulomb-sobolev-comparison}
			\|\zeta f\|_{H_v^1}\le C_\zeta\|f\|_{H_\sigma^1},
			\qquad
			\|\zeta F\|_{H_v^{-1}}\le C_\zeta\|F\|_{H_\sigma^{-1}}.
		\end{align}
	\end{enumerate}
\end{lemma}

\begin{proof}
	\emph{(1).} Positive definiteness of the inner product and completeness
	follow from the continuous embedding
	$H^1_\sigma(\bR_v^3)\hookrightarrow H^1_{\mathrm{loc}}(\bR_v^3)$ of
	Lemma~\ref{lem:local-embedding-H_sigma}.
	
	For \eqref{eq:radial-cutoff-convergence}, put $f_R:=\chi_Rf$.  Since
	$\nabla_v(f-f_R)=(1-\chi_R)\nabla_vf-f\nabla_v\chi_R$, the terms carrying
	$(1-\chi_R)f$ and $(1-\chi_R)\nabla_vf$ tend to zero in the corresponding
	weighted $L^2$ spaces by dominated convergence.  As $\chi_R$ is radial,
	$\nabla_v\chi_R=a_R(v)v$ with $|a_R(v)|\le CR^{-2}$ on the support of
	$\nabla_v\chi_R$, which lies in $\{R\le|v|\le2R\}$, whence
	\[
	\int_{\bR^3}
	\sigma^{ij}\p_i\chi_R\p_j\chi_R f^2
	\dd v
	\le
	\frac{C}{R^4}
	\int_{\{R\le|v|\le2R\}}
	\sigma^{ij}v_iv_jf^2
	\dd v\longrightarrow0 .
	\]
	By the Cauchy--Schwarz inequality for the nonnegative matrix $\sigma$,
	this proves \eqref{eq:radial-cutoff-convergence}.  For the density, each
	$f_R$ lies in $H^1(\bR_v^3)$ and has
	compact support, and $\rho_\delta\ast f_R\to f_R$ in $H^1(\bR_v^3)$ for a
	standard mollifier $\rho_\delta$; since $\sigma^{ij}$ and
	$\sigma^{ij}v_iv_j$ are bounded on compact sets, the convergence also
	holds in $H^1_\sigma$.  A diagonal argument completes the proof.
	
	\emph{(2).} By Lemma~\ref{lem:landau-diffusion-eigenvalues} and
	$\mathfrak q(v,\xi)\le C\la v\ra^{-1}|\xi|^2$ we have
	$\sigma^{ij}\xi_i\xi_j\le C|\xi|^2$ and $\sigma^{ij}v_iv_j\le C$, hence
	the first bound in \eqref{eq:Hv-into-Hsigma}.  The second follows by
	restricting $F$ to $H^1_v(\bR_v^3)$, the restriction being injective
	because $C_c^\infty(\bR_v^3)\subset H^1_v(\bR_v^3)$ is dense in
	$H^1_\sigma(\bR_v^3)$ by \textup{(1)}.
	
	\emph{(3).} Writing $\p_i(\zeta f)=\zeta\p_if+f\p_i\zeta$ and using the
	Cauchy--Schwarz inequality for $\sigma$,
	\[
	\sigma^{ij}\p_i(\zeta f)\p_j(\zeta f)
	\le
	2\zeta^2\sigma^{ij}\p_if\p_jf
	+
	2f^2\sigma^{ij}\p_i\zeta\p_j\zeta .
	\]
	Since $\sigma^{ij}\p_i\zeta\p_j\zeta\le C\|\zeta\|_{C^1}^2\la v\ra^{-1}$,
	the second term is bounded by
	$C\|\zeta\|_{C^1}^2\|\la v\ra^{-1/2}f\|_{L^2}^2$, hence by
	$C\|\zeta\|_{C^1}^2\|f\|_{H^1_\sigma}^2$ on account of
	\eqref{eq:local-H1-H_sigma}.  Together with
	$\int\sigma^{ij}v_iv_j\zeta^2f^2\le\|\zeta\|_{C^1}^2\|f\|_{H^1_\sigma}^2$
	this gives the first bound in \eqref{eq:Hsigma-multiplier}; the second
	follows from $\la\zeta F,\varphi\ra=\la F,\zeta\varphi\ra$.
	
	\emph{(4).} Since $\la v\ra\le\la R\ra$ on $B_R$,
	\eqref{eq:local-H1-H_sigma} gives
	\[
	\|u\|_{L^2}\le\la R\ra^{1/2}\|\la v\ra^{-1/2}u\|_{L^2}\le C_R\|u\|_{H^1_\sigma},
	\qquad
	\|\nabla_vu\|_{L^2}\le\la R\ra^{3/2}\|\la v\ra^{-3/2}\nabla_vu\|_{L^2}\le C_R\|u\|_{H^1_\sigma},
	\]
	which is the first bound in \eqref{eq:Hv-Hsigma-compact}.  For the
	second, fix $\zeta\in C_c^\infty(B_{2R})$ with $\zeta=1$ on $B_R$, so
	that $\la F,\varphi\ra=\la F,\zeta\varphi\ra$ for every $\varphi$.  As
	$\zeta\varphi$ is supported in $B_{2R}$, the first bound and
	\eqref{eq:Hsigma-multiplier} give
	\[
	\|F\|_{H^{-1}_\sigma}
	\le
	\|F\|_{H^{-1}_v}\sup_{\|\varphi\|_{H^1_\sigma}\le1}\|\zeta\varphi\|_{H^1_v}
	\le C_R\|F\|_{H^{-1}_v} .
	\]
	Finally, let $\vartheta\in C_c^\infty(\bR_v^3)$ with
	$\supp\vartheta\subset B_R$.  Then $\vartheta f$ is supported in $B_R$,
	so \eqref{eq:Hv-Hsigma-compact} and \eqref{eq:Hsigma-multiplier} give the
	first bound in \eqref{eq:local-coulomb-sobolev-comparison}, and the
	second follows from \eqref{eq:Hsigma-multiplier} and
	\eqref{eq:Hv-into-Hsigma}.
\end{proof}

\begin{remark}\label{rem:Riesz-H_sigma}
	For every $F\in H_\sigma^{-1}$, there exists a unique
	$g\in H_\sigma^1$ such that
	\[
	\langle F,\varphi\rangle_{H_\sigma^{-1},H_\sigma^1}
	=
	\int_{\bR^3}
	\left(
	\sigma^{ij}\p_i g\,\p_j\varphi
	+
	\sigma^{ij}v_iv_jg\varphi
	\right)
	\dd v
	\]
	for every $\varphi\in H_\sigma^1$. Equivalently, in the sense of
	distributions,
	\[
	F
	=
	\sigma^{ij}v_iv_jg
	-
	\p_i\left(\sigma^{ij}\p_jg\right).
	\]
	Thus
	\[
	H_\sigma^{-1}
	=
	\left\{
	\sigma^{ij}v_iv_jg
	-
	\p_i\left(\sigma^{ij}\p_jg\right)
	\,:\,
	g\in H_\sigma^1
	\right\},
	\]
	and the representing function $g$ satisfies
	\[
	\|F\|_{H_\sigma^{-1}}
	=
	\|g\|_{H_\sigma^1}.
	\]
\end{remark}

We use the standard macro--micro decomposition
\[
f
=
\mathbf P f
+
(\mathbf I-\mathbf P)f,
\]
where $\mathbf P$ denotes the orthogonal projection in
$L^2(\bR_v^3)$ onto
\[
\ker L
=
\operatorname{span}
\left\{
\sqrt{\mu},
\,
v_1\sqrt{\mu},
\,
v_2\sqrt{\mu},
\,
v_3\sqrt{\mu},
\,
\left(|v|^2-\frac32\right)\sqrt{\mu}
\right\}.
\]
Since $\mu(v)=\ee^{-|v|^2}$, these functions are mutually orthogonal
in $L^2(\bR_v^3)$ and satisfy
\begin{equation*}
	\|\sqrt{\mu}\|_{L^2_v}^2
	=
	\pi^{3/2},
	\qquad
	\|v_i\sqrt{\mu}\|_{L^2_v}^2
	=
	\frac12\pi^{3/2},
	\qquad
	\left\|
	\left(|v|^2-\frac32\right)\sqrt{\mu}
	\right\|_{L^2_v}^2
	=
	\frac32\pi^{3/2}.
\end{equation*}
It follows that
\hypertarget{comp:macro-projection}{}
\begin{equation*}
	\mathbf P f
	=
	\left[
	\mathbf a
	+
	\mathbf b\cdot v
	+
	\mathbf c
	\left(|v|^2-\frac32\right)
	\right]
	\sqrt{\mu},
\end{equation*}
\hypertarget{comp:macro-moments}{}
where
\begin{align*}
	\begin{aligned}
		\mathbf a
		&:=
		\pi^{-3/2}
		\int_{\bR^3}
		f(v)\sqrt{\mu(v)}
		\dd v,
		\\
		\mathbf b_i
		&:=
		2\pi^{-3/2}
		\int_{\bR^3}
		f(v)v_i\sqrt{\mu(v)}
		\dd v,
		\qquad i=1,2,3,
		\\
		\mathbf c
		&:=
		\frac23\pi^{-3/2}
		\int_{\bR^3}
		f(v)
		\left(|v|^2-\frac32\right)
		\sqrt{\mu(v)}
		\dd v.
	\end{aligned}
\end{align*}

By Corollary~\ref{cor:expansion-N}, the operator obtained by adding the
perturbation $N[g,\cdot]$ to $A$ is again in divergence form. We define
\begin{equation}\label{def_Ag}
	\cA_gf
	:=
	Af+8\pi\mu f+N[g,f]
	=
	D_i\left(a_g^{ij}D_jf\right)
	+
	b_g^iD_if
	+
	c_gf,
\end{equation}
where
\begin{equation}\label{260507202}
	\begin{aligned}
		a_g^{ij}
		&:=
		\sigma^{ij}
		+
		\sigma_{\sqrt\mu g}^{ij},
		\\
		b_g^i
		&:=
		-D_j\sigma_{\sqrt\mu g}^{ij}
		-
		2\sigma_{\sqrt\mu g}^{i},
		\\
		c_g
		&:=
		-v_iv_j\sigma^{ij}
		+
		(v_iv_j-\delta_{ij})
		\sigma_{\sqrt\mu g}^{ij}
		+
		12\pi\mu
		+
		8\pi\sqrt\mu g.
	\end{aligned}
\end{equation}
We also define
\begin{equation}\label{def_Kg}
	\cK f
	:=
	\sqrt\mu
	\sum_{i,j=1}^3
	\left(
	4v_iv_j-2\delta_{ij}
	\right)
	\sigma_{\sqrt\mu f}^{ij}.
\end{equation}
Then
\[
Af+Kf+N[g,f]
=
\cA_gf+\cK f.
\]
Accordingly we set
\begin{equation}\label{def:L_g}
	L_gf:=Lf-N[g,f],
\end{equation}
so that
\[
L_g=-\cA_g-\cK.
\]

\begin{lemma}[Weighted estimates for the compact part]
	\label{2607161158}
	There exists a universal constant $C>0$ such that
	\begin{equation}\label{eq:compact-part-pointwise}
		|\cK f(v)|
		\le
		C\mu(v)^{1/4}\|\sigma_{\sqrt \mu f}\|_{L^\infty(\bR_v^3)}\leq C\mu(v)^{1/4}
		\|\mu^{1/4}f\|_{L^2(\bR_v^3)}
	\end{equation}
	whenever $\mu^{1/4}f\in L^2(\bR_v^3)$. Consequently, for every
	$\alpha\ge0$,
	\begin{equation}\label{eq:compact-part-weighted-Hminus}
		\|w^\alpha\cK f\|_{H_\sigma^{-1}}\le 
		C\|w^{\alpha+1/2}\cK f\|_{L^2(\bR_v^3)} \le 
		C_\alpha
		\|\mu^{1/4}f\|_{L^2(\bR_v^3)},
	\end{equation}
	where $w(v):=(1+|v|^2)^{1/2}$, $C$ is a universal constant, and $C_\alpha>0$ depends only on $\alpha$.
\end{lemma}

\begin{proof}
	Since $|\phi^{ij}(z)|\le C|z|^{-1}$, the Cauchy--Schwarz
	inequality gives
	\begin{align*}
		|\sigma_{\sqrt\mu f}^{ij}(v)|
		&\le
		C
		\left(
		\int_{\bR^3}
		\frac{\mu(v_*)^{1/2}}{|v-v_*|^2}
		\dd v_*
		\right)^{1/2}
		\|\mu^{1/4}f\|_{L^2(\bR_v^3)}.
	\end{align*}
	Moreover,
	\begin{align*}
		\sup_{v\in\bR^3}
		\int_{\bR^3}
		\frac{\mu(v_*)^{1/2}}{|v-v_*|^2}
		\dd v_*
		&\le
		\|\mu^{1/2}\|_{L^\infty}
		\int_{|z|\le1}|z|^{-2}\dd z
		+
		\|\mu^{1/2}\|_{L^1}
		\le C.
	\end{align*}
	Thus
	\[
	|\sigma_{\sqrt\mu f}^{ij}(v)|
	\le
	C\|\mu^{1/4}f\|_{L^2(\bR_v^3)}.
	\]
	Since
	\[
	|4v_iv_j-2\delta_{ij}|
	\le
	Cw(v)^2,
	\qquad
	w(v)^2\sqrt{\mu(v)}
	\le
	C\mu(v)^{1/4},
	\]
	\eqref{eq:compact-part-pointwise} follows from \eqref{def_Kg}.
	
	By \eqref{eq:Hsigma-dual-weighted},
	$\|F\|_{H_\sigma^{-1}}\le C\|w^{1/2}F\|_{L^2(\bR_v^3)}$.
	Applying this to $F=w^\alpha\cK f$ and using
	\eqref{eq:compact-part-pointwise}, we obtain
	\begin{align*}
		\|w^\alpha\cK f\|_{H_\sigma^{-1}}
		&\le
		C
		\|w^{\alpha+1/2}\mu^{1/4}\|_{L^2(\bR_v^3)}
		\|\mu^{1/4}f\|_{L^2(\bR_v^3)},
	\end{align*}
	which proves \eqref{eq:compact-part-weighted-Hminus}.
\end{proof}

The first three of the following properties of the linearized collision
operator are elementary consequences of
Proposition~\ref{prop:explicit-landau-operators}. The coercivity
\eqref{eq:L-coercivity} is the standard spectral-gap estimate for the
Landau operator with Coulomb interaction, quoted from the literature.

\begin{lemma}[Basic properties of $L$ and {$N[g,\cdot]$}]\label{260715201}
	Fix $1/4\le\varpi<1/2$ and suppose
	$\mu^\varpi g\in L^2(\bR_v^3)$.  The following assertions hold
	for $f,h\in H^1_\sigma(\bR_v^3)$, with an absolute constant $C$
	for $L$ and a constant $C_\varpi$ depending only on $\varpi$
	for the estimate of $N[g,\cdot]$.
	
	\begin{enumerate}
		\item
		$L$ and $N[g,\cdot]$ are bounded from $H_\sigma^1(\bR_v^3)$ to $H_\sigma^{-1}(\bR_v^3)$:
		\begin{align}\label{boundedness of L}
			\|Lf\|_{H^{-1}_\sigma}\le C\|f\|_{H^1_\sigma}\,,
			\qquad
			\|N[g,f]\|_{H^{-1}_\sigma}\le C_\varpi \|\mu^\varpi g\|_{L^2(\bR_v^3)}\|f\|_{H^1_\sigma}.
		\end{align}
		In particular $L_g=L-N[g,\cdot]$ of \eqref{def:L_g} is bounded
		from $H^1_\sigma(\bR_v^3)$ to $H^{-1}_\sigma(\bR_v^3)$.
		Moreover, if $\mu^{1/4}g\in L^\infty(\bR^3_v)$, then
		\begin{align}\label{260911445}
			\|N[g,f]\|_{H^{-1}_\sigma}\leq C \|\mu^{1/4}g\|_{L^\infty(\bR^3_v)}\|f\|_{H^1_\sigma}.
		\end{align}
		
		\item
		$L$ is symmetric:
		\begin{align}\label{eq:L-symmetry}
			\la Lf,h\ra_v=\la Lh,f\ra_v.
		\end{align}
		
		\item
		The collision invariants lie in the kernel of $L$, and the
		mass invariant $\sqrt\mu$ is orthogonal to the range of $N[g,\cdot]$:
		\begin{align}\label{eq:L-kernel}
			\la Lf,\zeta\ra_v=0
			\quad\text{for every }\zeta\in\operatorname{ran}\mathbf P,
			\qquad
			L\mathbf Pf=0
			\quad\text{in }H^{-1}_\sigma(\bR_v^3),
		\end{align}
		and
		\begin{align}\label{eq:N-mass-kernel}
			\la N[g,f],\sqrt\mu\ra_v=0 .
		\end{align}
		
		\item
		$L$ is coercive on the orthogonal complement of its kernel:
		\begin{align}\label{eq:L-coercivity}
			\|(\mathbf{I}-\mathbf{P})f\|_{H^1_\sigma}^2\le C\la Lf,f\ra_v\,.
		\end{align}
	\end{enumerate}
\end{lemma}

\begin{proof}
	By \eqref{eq:local-H1-H_sigma} and
	$\|\sqrt\mu\,\psi\|_{L^1}\le\|\la v\ra^{1/2}\sqrt\mu\|_{L^2}
	\|\la v\ra^{-1/2}\psi\|_{L^2}$,
	\[
	\|\sqrt\mu\,\psi\|_{H^1(\bR_v^3)}+\|\sqrt\mu\,\psi\|_{L^1(\bR_v^3)}
	\le C\|\psi\|_{H^1_\sigma},
	\qquad \psi\in H^1_\sigma(\bR_v^3),
	\]
	so that Remark~\ref{rem:well-defined-N} defines $L$ and $N[g,\cdot]$ on
	$H^1_\sigma(\bR_v^3)$, and \eqref{eq:N-distributional-bound} makes them
	continuous from $H^1_\sigma(\bR_v^3)$ to $\cD'(\bR_v^3)$.  Since
	$C_c^\infty(\bR_v^3)$ is dense in $H^1_\sigma(\bR_v^3)$ by
	Lemma~\ref{lem:H_sigma-Hilbert}.(1), it suffices to prove (1) for
	$f\in C_c^\infty(\bR_v^3)$; (2) and (3) are then identities between
	continuous pairings and follow from their validity on
	$C_c^\infty(\bR_v^3)$.  The functions in $\operatorname{ran}\mathbf P$
	lie in $H^1_v(\bR_v^3)\subset H^1_\sigma(\bR_v^3)$ by
	Lemma~\ref{lem:H_sigma-Hilbert}.(2).  By
	Lemma~\ref{lem:landau-diffusion-eigenvalues} with $\varepsilon=1$ and
	\eqref{eq:H_sigma_norm},
	\begin{equation*}
		\int_{\bR^3}\mathfrak q(v,\nabla_v\psi)\dd v
		+\int_{\bR^3}\mathfrak q(v,v)\,\psi^2\dd v
		\le C\|\psi\|_{H^1_\sigma}^2,
		\qquad \psi\in H^1_\sigma .
	\end{equation*}
	
	\emph{(1).}
	Since $Af=N[\sqrt\mu,f]$ by \eqref{def:LAK}, only the second inequality in
	\eqref{boundedness of L} has to be proved.  Write
	\(G:=\sqrt\mu g=\mu^{1/2-\varpi}\,\mu^\varpi g\) and
	\(M_\varpi:=\|\mu^\varpi g\|_{L^2_v}\); since \(\varpi<1/2\), every
	polynomial moment of \(G\) is bounded by \(C_\varpi M_\varpi\).  Besides
	the diffusion coefficient \(\sigma_G^{ij}\) of \eqref{def:N} we use its
	divergence and the two conjugated derivatives
	\begin{equation}\label{eq:velalg-AF-BF-definition}
		b^i:=\p_j\sigma_G^{ij}=(\p_j\phi^{ij})\ast G,
		\qquad
		D_i^\pm:=\p_i\pm v_i .
	\end{equation}
	Integrating the \(v_*\)-derivative in the first Landau slot of
	\eqref{new:eq:intro-landau-operator} by parts, and using
	\(\p_i(\sqrt\mu f)=\sqrt\mu\,D_i^-f\) and
	\(\p_i(\mu^{-1/2}\varphi)=\mu^{-1/2}D_i^+\varphi\), gives the exact
	conjugated weak form
	\begin{equation}\label{eq:velalg-exact-weak-form}
		\langle N[g,f],\varphi\rangle_v
		=-\int_{\bR^3}
		\bigl\{\sigma_G^{ij}D_j^-f-b^i f\bigr\}D_i^+\varphi\dd v ,
		\qquad \varphi\in C_c^\infty(\bR_v^3);
	\end{equation}
	no derivative of \(g\) and no further coefficient occurs.  By
	\eqref{eq:coulomb-sigmavv-comparison},
	\begin{equation}\label{eq:velalg-Dpm}
		\|\sigma^{1/2}D^\pm q\|_{L^2_v}\le C\|q\|_{H^1_\sigma},
		\qquad
		\|\,(\sigma^{ij}v_iv_j)^{1/2}q\|_{L^2_v}\le C\|q\|_{H^1_\sigma} .
	\end{equation}
	
	\emph{The diffusion coefficient.}
	Let \(r=|v|\ge4\), \(e:=v/r\), and use the frame \(P_v,Q_v\) of
	\eqref{eq:radial-tangential-projections}.  On \(\{|v_*|\le r/2\}\) the
	Coulomb kernel \eqref{new:eq:intro-coulomb-kernel} satisfies
	\[
	e^T\phi(v-v_*)e=|v-v_*|^{-3}|Q_vv_*|^2,
	\qquad
	\|P_v\phi(v-v_*)Q_v\|\le Cr^{-2}|v_*|,
	\qquad
	\|Q_v\phi(v-v_*)Q_v\|\le Cr^{-1}.
	\]
	Integrating against \(G\), whose zeroth and second moments are bounded by
	\(C_\varpi M_\varpi\), gives
	\begin{equation}\label{eq:AF-anisotropy}
		\|P_v\sigma_GP_v\|\le C_\varpi M_\varpi\la v\ra^{-3},
		\quad
		\|P_v\sigma_GQ_v\|\le C_\varpi M_\varpi\la v\ra^{-2},
		\quad
		\|Q_v\sigma_GQ_v\|\le C_\varpi M_\varpi\la v\ra^{-1}.
	\end{equation}
	On the complementary region, Cauchy--Schwarz gives
	\begin{equation}\label{eq:AF-complementary-CS}
		\int_{|v_*|>r/2}|v-v_*|^{-1}|G(v_*)|\dd v_*
		\le C_\varpi M_\varpi
		\left(
		\int_{|v_*|>r/2}|v-v_*|^{-2}
		\mu(v_*)^{1-2\varpi}\dd v_*
		\right)^{1/2},
	\end{equation}
	and splitting the last integral at \(|v-v_*|=1\), with
	\(\int_{|z|<1}|z|^{-2}\dd z=4\pi\) on the first part and the Gaussian tail
	on the second,
	\begin{equation}\label{eq:AF-complementary-exponential}
		\int_{|v_*|>r/2}|v-v_*|^{-1}|G(v_*)|\dd v_*
		\le C_\varpi\ee^{-c_\varpi r^2}M_\varpi ,
	\end{equation}
	which is dominated by each right-hand side of \eqref{eq:AF-anisotropy}.
	Since \(\sigma\) is uniformly positive definite on \(\{r\le4\}\),
	\eqref{eq:AF-anisotropy} holds for every \(v\), and by
	\eqref{eq:coulomb-eigen-decomposition} it is the matrix estimate
	\begin{equation}\label{eq:AF-sigma-form}
		\left|\sigma_G^{ij}(v)\xi_i\bar\xi_j\right|
		\le C_\varpi M_\varpi
		\left(\sigma^{ij}(v)\xi_i\xi_j\right)^{1/2}
		\left(\sigma^{ij}(v)\bar\xi_i\bar\xi_j\right)^{1/2}.
	\end{equation}
	With \eqref{eq:velalg-Dpm} this bounds the first term of
	\eqref{eq:velalg-exact-weak-form} by
	\(C_\varpi M_\varpi\|f\|_{H^1_\sigma}\|\varphi\|_{H^1_\sigma}\).
	
	\emph{The drift coefficient.}
	Here \((\p_j\phi^{ij})(z)=-2z_i|z|^{-3}\) is not square integrable near the
	origin, so we split it by a radial cutoff \(\chi\in C_c^\infty(\bR^3)\)
	with \(\chi=1\) on \(|z|\le1\) and \(\chi=0\) on \(|z|\ge2\):
	\begin{equation}\label{eq:velalg-BF-split}
		b=b^{\rm loc}+b^{\rm far},
		\qquad
		b^{i,\rm loc}:=\bigl(\chi\,\p_j\phi^{ij}\bigr)\ast G,
		\qquad
		b^{i,\rm far}:=\bigl((1-\chi)\,\p_j\phi^{ij}\bigr)\ast G .
	\end{equation}
	On \(|v_*|\le r/2\) one has \(|P_v\p_j\phi^{ij}(v-v_*)|\le Cr^{-2}\) and
	\(|Q_v\p_j\phi^{ij}(v-v_*)|\le Cr^{-3}|v_*|\); integrating against \(G\),
	and treating \(|v_*|>r/2\) by
	\eqref{eq:AF-complementary-CS}--\eqref{eq:AF-complementary-exponential}
	with the bounded kernel \(|(1-\chi)\p_j\phi^{ij}|\le C\), we obtain
	\begin{equation}\label{eq:BF-far-component-bounds}
		|P_vb^{\rm far}|\le C_\varpi M_\varpi\la v\ra^{-2},
		\qquad
		|Q_vb^{\rm far}|\le C_\varpi M_\varpi\la v\ra^{-3}.
	\end{equation}
	For $|v|\ge4$, \eqref{eq:BF-far-component-bounds} and the
	radial--tangential eigenvalue bounds for $\sigma$ imply
	$b^{\rm far}\cdot\sigma^{-1}b^{\rm far}
	\lesssim M_\varpi^2\langle v\rangle^{-1}
	\lesssim M_\varpi^2\sigma^{ij}v_iv_j$.
	On $|v|\le4$, the far kernel is bounded, so
	$\|b^{\rm far}\|_{L^\infty_v}\lesssim\|G\|_{L^1_v}\lesssim M_\varpi$;
	there $\sigma^{-1}$ is bounded.  Consequently, for all $v$,
	\begin{equation*}
		b^{\rm far}\cdot\sigma^{-1}b^{\rm far}
		\le C_\varpi M_\varpi^2
		\bigl(\mathbf1_{\{|v|\le4\}}+\sigma^{ij}v_iv_j\bigr).
	\end{equation*}
	By the local embedding \eqref{eq:local-H1-H_sigma},
	\[
	\int_{\bR^3}
	\bigl(b^{\rm far}\cdot\sigma^{-1}b^{\rm far}\bigr)|f|^2\dd v
	\le C_\varpi M_\varpi^2\|f\|_{H^1_\sigma}^2.
	\]
	Cauchy--Schwarz in the metric $\sigma$ and \eqref{eq:velalg-Dpm} now give
	\[
	\left|\int_{\bR^3}b^{i,\rm far}f\,D_i^+\varphi\dd v\right|
	\le C_\varpi M_\varpi
	\|f\|_{H^1_\sigma}\|\varphi\|_{H^1_\sigma}.
	\]
	
	For \(b^{\rm loc}\) we use the order-one potential estimate
	\begin{equation}\label{eq:velalg-Euclidean-I1-L2-L6}
		\bigl\||z|^{-2}\ast G\bigr\|_{L^6_v}\le C\|G\|_{L^2_v}
		\le C_\varpi M_\varpi ,
	\end{equation}
	which follows from the pointwise bound
	\(\bigl||z|^{-2}\ast G\bigr|\le C\|G\|_{L^2_v}^{2/3}(\cM G)^{1/3}\)
	and the \(L^2\) bound for the Hardy--Littlewood maximal function \(\cM\).
	The pointwise bound is obtained by splitting the convolution at \(|z|=R\),
	bounding the near part by \(CR\,\cM G\) through dyadic annuli and the far
	part by \(CR^{-1/2}\|G\|_{L^2_v}\) by Cauchy--Schwarz, and optimizing in
	\(R\).
	
	Cover a dyadic shell \(S_R:=\{R\le|v|\le2R\}\), \(R\ge4\), by unit cubes
	\(Q\) whose fixed enlargements \(Q^*\subset Q^{**}\) have uniformly finite
	overlap.  If \(v\in Q\) and \(|v-v_*|\le2\) then \(v_*\in Q^{**}\) and
	\(|v_*|\simeq R\), so that
	\begin{equation*}
		\|G\|_{L^2(Q^{**})}\le C_\varpi\ee^{-cR^2}M_\varpi .
	\end{equation*}
	Applying \eqref{eq:velalg-Euclidean-I1-L2-L6} to \(G\one_{Q^{**}}\),
	\begin{equation}\label{eq:velalg-local-B-L6}
		\bigl\|\one_Qb^{\rm loc}\bigr\|_{L^6_v}
		\le C\|G\|_{L^2(Q^{**})}
		\le C_\varpi\ee^{-cR^2}M_\varpi .
	\end{equation}
	On \(Q^*\), \eqref{eq:local-H1-H_sigma} and
	\eqref{eq:coulomb-sigmavv-comparison} give
	\(\|q\|_{L^2(Q^*)}\le CR^{1/2}\|q\|_{H^1_\sigma(Q^*)}\) and
	\(\|\nabla_vq\|_{L^2(Q^*)}\le CR^{3/2}\|q\|_{H^1_\sigma(Q^*)}\), hence, by
	Euclidean Sobolev interpolation on a set of fixed size,
	\begin{equation}\label{eq:velalg-cube-sobolev}
		\|q\|_{L^3(Q^*)}\le CR\|q\|_{H^1_\sigma(Q^*)},
		\qquad
		\|D^\pm q\|_{L^2(Q^*)}\le CR^{3/2}\|q\|_{H^1_\sigma(Q^*)} .
	\end{equation}
	H\"older with \(1/6+1/3+1/2=1\), \eqref{eq:velalg-local-B-L6} and
	\eqref{eq:velalg-cube-sobolev} therefore give
	\[
	\left|\int_Qb^{i,\rm loc}f\,D_i^+\varphi\dd v\right|
	\le C_\varpi\ee^{-cR^2}R^{5/2}M_\varpi
	\|f\|_{H^1_\sigma(Q^*)}\|\varphi\|_{H^1_\sigma(Q^*)} .
	\]
	Summing over \(Q\subset S_R\) by Cauchy--Schwarz in the cube index, the
	finite overlap of the \(Q^*\) costs no factor equal to the number of
	cubes, and \(R^{5/2}\ee^{-cR^2}\) is summable over the dyadic shells; on
	\(\{|v|\le8\}\) the two norms are equivalent by
	Lemma~\ref{lem:H_sigma-Hilbert}.(4).  Hence
	\[
	\left|\int_{\bR^3}b^{i,\rm loc}f\,D_i^+\varphi\dd v\right|
	\le C_\varpi M_\varpi\|f\|_{H^1_\sigma}\|\varphi\|_{H^1_\sigma}.
	\]
	
	Adding the three bounds, taking the supremum over
	\(\|\varphi\|_{H^1_\sigma}\le1\), and using the density of
	\(C_c^\infty(\bR_v^3)\) in \(H^1_\sigma(\bR_v^3)\) give the second
	inequality in \eqref{boundedness of L}; the boundedness of
	$L_g=L-N[g,\cdot]$ follows.  Finally,
	\(\mu^{3/8}g=\mu^{1/8}\,\mu^{1/4}g\) gives
	\(\|\mu^{3/8}g\|_{L^2_v}\le\|\mu^{1/8}\|_{L^2_v}\|\mu^{1/4}g\|_{L^\infty}\),
	so that \eqref{260911445} follows from the second inequality in
	\eqref{boundedness of L} with \(\varpi=3/8\) and from
	\eqref{eq:Hv-into-Hsigma}.

	\emph{(2).} The symmetry
	$\la A f,h\ra_v=\la f,Ah\ra_v$
	follows from \eqref{eq:explicit-landau-A}.
	For $K$, \eqref{eq:explicit-landau-K} gives
	\[
	\la Kf,h\ra_v
	=
	\iint_{\bR^3\times\bR^3}
	\mathfrak k(v,v_*)f(v_*)h(v)
	\dd v_*\dd v
	+
	8\pi\int_{\bR^3}\mu fh\dd v,
	\]
	where
	\[
	\mathfrak k(v,v_*)
	:=
	\left(4v_iv_j-2\delta_{ij}\right)
	\phi^{ij}(v-v_*)\sqrt{\mu(v)}\sqrt{\mu(v_*)},
	\]
	the double integral being absolutely convergent since, by
	\eqref{eq:compact-part-pointwise} applied to $|f|$ and by
	$|4v_iv_j-2\delta_{ij}|\sqrt\mu\le C\mu^{1/4}$,
	\[
	\iint_{\bR^3\times\bR^3}
	|\mathfrak k(v,v_*)||f(v_*)||h(v)|\dd v_*\dd v
	\le
	C\|\mu^{1/4}f\|_{L^2(\bR_v^3)}
	\|\mu^{1/4}h\|_{L^1(\bR_v^3)}
	<\infty .
	\]
	Since $\phi^{ij}(z)z_j=0$, we have
	$\phi^{ij}(v-v_*)v_j=\phi^{ij}(v-v_*)v_{*j}$ and hence
	$\phi^{ij}(v-v_*)v_iv_j=\phi^{ij}(v-v_*)v_{*i}v_{*j}$.  Since $\phi^{ij}$
	is even, $\mathfrak k(v,v_*)=\mathfrak k(v_*,v)$.  Since $L=-A-K$, this
	proves \eqref{eq:L-symmetry}.
	
	\emph{(3).} Let $F,G\in\cS(\bR_v^3)$ and write $G_*:=G(v_*)$.  Adding the
	distributional definitions of $\cN(G,F)$ and $\cN(F,G)$ and exchanging $v$
	and $v_*$ in the second one, all the integrals converging absolutely for
	Schwartz $F,G$, we obtain from the evenness of $\phi^{ij}$, for every
	polynomial $\psi$ of degree at most two,
	\begin{equation}\label{eq:landau-weak-symmetrized}
		\begin{aligned}
			&\la\cN(G,F)+\cN(F,G),\psi\ra
			\\
			={}&
			-\iint_{\bR^3\times\bR^3}
			\phi^{ij}(v-v_*)
			\left[
			G_*\p_jF(v)-F(v)(\p_jG)_*
			\right]
			\left(
			\p_i\psi(v)-\p_i\psi(v_*)
			\right)
			\dd v_*\dd v .
		\end{aligned}
	\end{equation}
	For $\psi\in\{1,v_1,v_2,v_3\}$ the gradient $\nabla\psi$ is constant, so
	the last factor vanishes.  For $\psi=|v|^2$ it equals $2(v-v_*)_i$, and
	$\phi^{ij}(z)z_i=0$.  Hence the right-hand side of
	\eqref{eq:landau-weak-symmetrized} vanishes for
	$\psi\in\{1,v_1,v_2,v_3,|v|^2\}$.  Applying this with $G=\mu$ and
	$F=\sqrt\mu f$ and recalling \eqref{def:N} and \eqref{def:LAK}, so that
	\[
	\cN(\mu,\sqrt\mu f)+\cN(\sqrt\mu f,\mu)
	=
	\sqrt\mu\left(Af+Kf\right)
	=
	-\sqrt\mu\,Lf,
	\]
	we obtain $\la Lf,\psi\sqrt\mu\ra_v=0$ for every such $\psi$.  Since every
	$\zeta\in\operatorname{ran}\mathbf P$ is a linear combination of these
	functions $\psi\sqrt\mu$, the first assertion in \eqref{eq:L-kernel}
	follows.  The second follows from the
	first and \eqref{eq:L-symmetry}: for every $h\in H^1_\sigma$,
	$\la L\mathbf Pf,h\ra_v=\la Lh,\mathbf Pf\ra_v=0$, since
	$\mathbf Pf\in\operatorname{ran}\mathbf P$.
	
	By \eqref{def:N} we have $\sqrt\mu\,N[g,f]=\cN(\sqrt\mu g,\sqrt\mu f)$,
	and the distributional definition of $\cN$ pairs it with $\nabla\varphi$
	only, so that testing with $\varphi\equiv1$ gives
	\[
	\la N[g,f],\sqrt\mu\ra_v
	=\int_{\bR^3}\cN(\sqrt\mu g,\sqrt\mu f)\dd v
	=\la\cN(\sqrt\mu g,\sqrt\mu f),1\ra=0
	\]
	for $f\in C_c^\infty(\bR_v^3)$, all the integrals converging absolutely;
	the identity extends to $f\in H^1_\sigma(\bR_v^3)$ by (1).  This is
	\eqref{eq:N-mass-kernel}.
	
	\emph{(4).} Estimate \eqref{eq:L-coercivity} is the microscopic coercivity
	of the linearized Coulomb Landau operator, proved in \cite{GuoCMP2002}.
\end{proof}

\begin{lemma}\label{lem:poly-maxwellian}
	For every polynomial $p$ on $\bR_v^3$ there exists $C_p>0$ such that
	\begin{equation}\label{eq:poly-maxwellian}
		\big\|\mu^{-1/4}L\big(p\sqrt\mu\big)\big\|_{L^2(\bR_v^3)}
		\le
		C_p .
	\end{equation}
\end{lemma}

\begin{proof}
	Write $\psi:=p\sqrt\mu$ and recall $L=-A-K$ together with the explicit
	formulas \eqref{eq:explicit-landau-A} and \eqref{eq:explicit-landau-K}.
	All the coefficients occurring there are bounded: by
	Lemma~\ref{lem:landau-diffusion-eigenvalues} and
	\eqref{eq:derivative-sigma-bound},
	\[
	|\sigma^{ij}|+|v_iv_j\sigma^{ij}|\le C\la v\ra^{-1},
	\qquad
	|\p_k\sigma^{ij}|\le C\la v\ra^{-2}.
	\]
	Since every derivative of $\psi$ is again a polynomial times $\sqrt\mu$,
	\eqref{eq:explicit-landau-A} gives $|A\psi|\le C_p\la v\ra^{N}\sqrt\mu$
	for some $N=N(p)$. Therefore,
	\[
	\big|\mu^{-1/4}A\psi\big|\le C_p\la v\ra^{N}\mu^{1/4}\in L^2(\bR_v^3).
	\]
	For the compact part, $\|\mu^{1/4}\psi\|_{L^2(\bR_v^3)}\le C_p$, so
	\eqref{eq:compact-part-pointwise} gives
	$\|\sigma^{ij}_{\sqrt\mu\psi}\|_{L^\infty}\le C_p$. Hence
	\eqref{eq:explicit-landau-K} gives
	$|K\psi|\le C_p\la v\ra^{2}\sqrt\mu$ and therefore
	\[
	\big|\mu^{-1/4}K\psi\big|\le C_p\la v\ra^{2}\mu^{1/4}\in L^2(\bR_v^3).
	\]
	Adding the two bounds proves \eqref{eq:poly-maxwellian}.
\end{proof}

\subsection{The solution space $\cH_{\kin,\sigma}$}

We introduce the dynamic and stationary $\cH_{\kin,\sigma}$-spaces, the
corresponding notion of solution, and the analytic properties of these spaces.

\begin{defn}[$\cH_{\kin,\sigma}$-spaces]\label{def:H_kin_sigma}\,
	
	\begin{enumerate}
		\item For \(f\in C_c^\infty\big([0,T]\times \overline{\cD}\times \bR_v^3\big)\), we define
		\begin{align}\label{innerprod_sigma}
			\begin{aligned}
				\|f\|_{\cH_{\kin,\sigma}(\cD_T \times \bR_v^3)}^2:=\,&
				\|f\|_{L^2(\cD_T;H_\sigma^1)}^2
				+\|(\partial_t+v\cdot \nabla_x )f\|_{L^2(\cD_T;H_\sigma^{-1})}^2\\
				&\quad
				+\|f(0,\,\cdot\,)\|_{L^2(\cD\times \bR_v^3)}^2
				+\big\|f|_{\gamma_T^-}\big\|_{L^2(\gamma^-_T)}^2\,.
			\end{aligned}
		\end{align}
		We define \(\cH_{\kin,\sigma}(\cD_T\times \bR_v^3)\) as the completion of
		$C_c^\infty\big([0,T]\times \overline{\cD}\times \bR_v^3\big)$ with respect to the norm \eqref{innerprod_sigma}. 
		For \(f\in \cH_{\kin,\sigma}(\cD_T\times \bR_v^3)\), we define
		\((\p_t+v\cdot \nabla_x)f\), \(f(0,\,\cdot\,)\), and \(f_{\gamma^-}\) 
		as the corresponding limits
		\begin{align}\label{lim_transp_sigma}
			\begin{alignedat}{2}
				(\p_t+v\cdot\nabla_x)f_n
				&\longrightarrow (\p_t+v\cdot \nabla_x)f
				\quad &&\text{in}\quad L^2(\cD_T;H^{-1}_\sigma),\\
				f_n(0,\cdot)
				&\longrightarrow f(0,\cdot)
				\quad &&\text{in}\quad L^2(\cD\times \bR_v^3),\\
				f_n|_{\gamma^-_T}
				&\longrightarrow f_{\gamma^-}
				\quad &&\text{in}\quad L^2(\gamma^-_T),
			\end{alignedat}
		\end{align}
		where \((f_n)\) is any sequence converging to \(f\) in
		\(\cH_{\kin,\sigma}(\cD_T\times \bR_v^3)\).
		
		\item For \(f\in C_c^\infty\big(\overline{\cD}\times \bR_v^3\big)\), we define
		\begin{align}\label{innerprod_st_sigma}
			\begin{aligned}
				\|f\|_{ {\cH}_{\mathrm{kin,\sigma}}(\cD \times \bR_v^3)}^2:=\,&
				\|f\|_{L^2(\cD ;H_\sigma^1)}^2
				+\|v\cdot \nabla_x f\|_{L^2(\cD ;H_\sigma^{-1})}^2
				+\big\|f|_{\gamma^-}\big\|_{L^2({\gamma^-})}^2\,.
			\end{aligned}
		\end{align}
		We define \(\cH_{\kin,\sigma}(\cD\times \bR_v^3)\) as the completion of
		\(C_c^\infty\big(\overline{\cD}\times \bR_v^3\big)\) with respect to the norm \eqref{innerprod_st_sigma}.
		For \(f\in\cH_{\kin,\sigma}(\cD\times \bR_v^3)\), we define
		\(v\cdot \nabla_x f\) and \(f_{\gamma^-}\) as the corresponding limits
		\begin{align}\label{lim_transp_st_sigma}
			\begin{alignedat}{2}
				v\cdot\nabla_x f_n
				&\longrightarrow v\cdot \nabla_x f
				\quad &&\text{in}\quad L^2(\cD;H^{-1}_\sigma),\\
				f_n|_{\gamma^-}
				&\longrightarrow f_{\gamma^-}
				\quad &&\text{in}\quad L^2(\gamma^-),
			\end{alignedat}
		\end{align}
		where \((f_n)\) is any sequence converging to \(f\) in
		\(\cH_{\kin,\sigma}(\cD\times \bR_v^3)\). 
	\end{enumerate}
\end{defn}

\begin{prop}\label{2607151218}
	The limits \eqref{lim_transp_sigma} exist and do not depend on the
	approximating sequence, and
	$\cH_{\kin,\sigma}(\cD_T\times\bR_v^3)$ and
	$\cH_{\kin,\sigma}(\cD\times\bR_v^3)$ are Hilbert spaces
	with the norms \eqref{innerprod_sigma} and
	\eqref{innerprod_st_sigma}, respectively.
	
	In particular, $\cH_{\kin,\sigma}(\cD_T\times\bR_v^3)$ continuously embeds into
	$C\big([0,T];L^2(\cD\times\bR_v^3)\big)$. Moreover, there exists a continuous
	phase-space boundary trace operator
	\[
	T_\gamma:
	\cH_{\kin,\sigma}(\cD_T\times\bR_v^3)
	\longrightarrow L^2(\gamma_T),
	\qquad
	T_\gamma f=:f_\gamma,
	\]
	such that \eqref{ibp_dy}, with \(\la\cdot,\cdot\ra_v\) the \(H_\sigma^{-1}\)--\(H_\sigma^1\) pairing, holds for all
	$f,g\in\cH_{\kin,\sigma}(\cD_T\times\bR_v^3)$. In addition,
	\begin{align}\label{260730346}
		\sup_{0\le t\le T}
		\|f(t,\cdot)\|_{L^2(\cD\times\bR_v^3)}
		+\|f_\gamma\|_{L^2(\gamma_T)}
		\le
		C\|f\|_{\cH_{\kin,\sigma}(\cD_T\times\bR_v^3)},
	\end{align}
	where $C>0$ is a universal constant.
	
	Similarly, there exists a continuous phase-space boundary trace operator
	\[
	T_\gamma:
	\cH_{\kin,\sigma}(\cD\times\bR_v^3)
	\longrightarrow L^2(\gamma),
	\qquad
	T_\gamma f=:f_\gamma,
	\]
	such that \eqref{ibp_st} holds for all
	$f,g\in\cH_{\kin,\sigma}(\cD\times\bR_v^3)$. In addition,
	\[
	\|f_\gamma\|_{L^2(\gamma)}
	\le
	C\|f\|_{\cH_{\kin,\sigma}(\cD\times\bR_v^3)},
	\]
	where $C>0$ is a universal constant.
	
	The two full traces agree with the incoming traces in
	\eqref{lim_transp_sigma} and \eqref{lim_transp_st_sigma}, respectively:
	\[
	f_\gamma|_{\gamma_T^-}=f_{\gamma^-}
	\quad\text{in }L^2(\gamma_T^-),
	\qquad
	f_\gamma|_{\gamma^-}=f_{\gamma^-}
	\quad\text{in }L^2(\gamma^-).
	\]
\end{prop}

\begin{proof}
	We give the evolutionary argument.  The stationary case follows from
	the same Green identity without its time endpoints.  Write
	$\mathfrak T=\partial_t+v\cdot\nabla_x$ and let $f$ be smooth with compact velocity
	support.  Integration in $(t,x)$ gives, for $0\le s\le t\le T$,
	\begin{align*}
		\|f(t)\|_2^2+\|f_{\gamma^+}\|_{L^2((s,t)\times\gamma^+)}^2
		={}&\|f(s)\|_2^2+\|f_{\gamma^-}\|_{L^2((s,t)\times\gamma^-)}^2
		+2\int_s^t\!\int_\cD\langle \mathfrak Tf,f\rangle_{H_\sigma^{-1},H_\sigma^1}\,\dd x\dd r.
	\end{align*}
	At $s=0$, Cauchy--Schwarz and Young's inequality imply
	\begin{align*}
		\sup_{0\le t\le T}\|f(t)\|_2^2
		+\|f_{\gamma^+}\|_{L^2(\gamma_T^+)}^2
		\le C\big(&\|f(0)\|_2^2
		+\|f_{\gamma^-}\|_{L^2(\gamma_T^-)}^2
		+\|\mathfrak T f\|_{L^2(\cD_T;H_\sigma^{-1})}^2
		+\|f\|_{L^2(\cD_T;H_\sigma^1)}^2\big),
	\end{align*}
	with $C$ independent of $T$.  In particular, the graph norm $\|f\|_{\cH_{\kin,\sigma}}$ controls
	the full boundary trace $f_\gamma$ and the supremum of the time slices.  Applying
	this estimate to the difference of two smooth functions extends both
	coordinate maps continuously to the completion in
	Definition~\ref{def:H_kin_sigma}.  Uniform convergence of the time
	slices of an approximating sequence gives a representative in
	$C([0,T];L^2_{x,v})$.  The bilinear Green identity follows by applying
	the smooth identity to pairs of approximants and passing to the limit:
	the bulk terms converge in the $H_\sigma^{-1}$--$H_\sigma^1$ duality,
	and both traces converge in the flux $L^2$ space.  The resulting
	outgoing and incoming traces agree with the limits in the definition.
	
	For completeness, the abstract graph completion embeds injectively
	into its bulk $L^2_tH_\sigma^1$ coordinate.  Indeed, convergence in this
	space implies convergence in $L^2_{\mathrm{loc}}$ in velocity, by
	\eqref{eq:local-H1-H_sigma}.  If the bulk coordinate vanishes, its
	transport derivative vanishes as a distribution.  The Green identity,
	tested successively with functions supported near either time endpoint
	and near the non-grazing parts of the phase boundary, then forces the
	initial value and both boundary traces to vanish.  Thus an element
	with zero bulk coordinate has zero graph norm.  As a completion of an
	inner-product space, $\cH_{\kin,\sigma}$ is therefore Hilbert.
	
	In the stationary case the same calculation gives
	\[
	\|f_{\gamma^+}\|_{L^2(\gamma^+)}^2
	\le \|f_{\gamma^-}\|_{L^2(\gamma^-)}^2
	+2\|v\cdot\nabla_xf\|_{L^2_xH_\sigma^{-1}}
	\|f\|_{L^2_xH_\sigma^1}.
	\]
	Density yields the full trace, the stationary Green identity, and the
	corresponding Hilbert-space assertion.  None of these estimates uses
	coercivity of $H_\sigma^1$ over all velocities in the unweighted
	$L^2_v$ norm.
\end{proof}

\begin{lemma}\label{lem:compact-velocity-support}\,
	
	\begin{enumerate}
		\item For every $\zeta\in C^1(\bR_v^3)$ with
		$\|\zeta\|_{C^1(\bR_v^3)}<\infty$, multiplication by $\zeta$ maps each of
		$\cH_{\kin}(\cD_T\times\bR_v^3)$ and
		$\cH_{\kin,\sigma}(\cD_T\times\bR_v^3)$ boundedly into itself, with
		\[
		\|\zeta f\|_{\cH_{\kin,\sigma}(\cD_T\times\bR_v^3)}
		\le
		C\|\zeta\|_{C^1(\bR_v^3)}\|f\|_{\cH_{\kin,\sigma}(\cD_T\times\bR_v^3)}
		\]
		and
		\[
		(\p_t+v\cdot\nabla_x)(\zeta f)=\zeta\,(\p_t+v\cdot\nabla_x)f,
		\qquad
		(\zeta f)(0,\cdot)=\zeta f(0,\cdot),
		\qquad
		(\zeta f)_{\gamma^-}=\zeta f_{\gamma^-} .
		\]
		
		\item Let $R>0$ and let $f$ belong to
		$\cH_{\kin}(\cD_T\times\bR_v^3)$ or to
		$\cH_{\kin,\sigma}(\cD_T\times\bR_v^3)$ and satisfy
		\[
		\supp\big(f(t,x,\cdot)\big)\subset B_R
		\qquad\text{for a.e. }(t,x)\in\cD_T.
		\]
		Then $f$ belongs to both spaces, and
		\[
		\|f\|_{\cH_{\kin,\sigma}(\cD_T\times \bR_v^3)}\simeq_R\|f\|_{\cH_{\kin}(\cD_T\times \bR_v^3)}.
		\]
	\end{enumerate}
	The same assertions hold, with the same proofs, for the stationary spaces
	$\cH_\kin(\cD\times\bR_v^3)$ and $\cH_{\kin,\sigma}(\cD\times\bR_v^3)$,
	with
	\[
	v\cdot\nabla_x(\zeta f)=\zeta\,v\cdot\nabla_xf,
	\qquad
	(\zeta f)_{\gamma^\pm}=\zeta f_{\gamma^\pm}.
	\]
\end{lemma}

\begin{proof}
	A function of $v$ commutes with $\p_t+v\cdot\nabla_x$ and multiplies the
	initial and inflow traces in \eqref{innerprod_sigma} pointwise, so both
	assertions reduce to their velocity counterparts
	\eqref{eq:Hsigma-multiplier} and
	\eqref{eq:Hv-into-Hsigma}--\eqref{eq:Hv-Hsigma-compact} of
	Lemma~\ref{lem:H_sigma-Hilbert}; the unweighted analogues are immediate
	from the Leibniz rule.
	
	Let $\phi\in C_c^\infty([0,T]\times\overline\cD\times\bR_v^3)$.  Since
	$|\zeta|\le\|\zeta\|_{C^1}$ on the traces,
	\eqref{eq:Hsigma-multiplier} applied pointwise in $(t,x)$ bounds each of
	the four terms of \eqref{innerprod_sigma} for $\zeta\phi$ by
	$C\|\zeta\|_{C^1}$ times the corresponding term for $\phi$.  When
	$\zeta\in C_c^\infty(\bR_v^3)$ the product $\zeta\phi$ is again an
	admissible test function; for general $\zeta\in C^1$ we apply the bound
	to a mollification $\zeta_\delta$ of $\zeta$ in $v$, for which
	$\|\zeta_\delta\|_{C^1}\le\|\zeta\|_{C^1}$ and $\zeta_\delta$,
	$\nabla\zeta_\delta$ converge uniformly to $\zeta$, $\nabla\zeta$ on the
	compact velocity support of $\phi$, so that
	$\zeta_\delta\phi\to\zeta\phi$ in
	$\cH_{\kin,\sigma}(\cD_T\times\bR_v^3)$.  Since $\cH_{\kin,\sigma}(\cD_T\times\bR_v^3)$ is by
	definition the completion of
	$C_c^\infty([0,T]\times\overline\cD\times\bR_v^3)$ under this norm, the
	multiplier extends by density.  This proves~(1), and the same argument
	applies to $\cH_\kin$.
	
	For (2), suppose in addition that the velocity support of $\phi$ lies in
	$B_R$; then so does that of $(\p_t+v\cdot\nabla_x)\phi(t,x,\cdot)$.
	Hence \eqref{eq:Hv-into-Hsigma}--\eqref{eq:Hv-Hsigma-compact}, which
	together are two-sided on functions supported in $B_R$, applied
	pointwise in $(t,x)$, together with the fact that the two trace terms
	are literally the same, yield
	\begin{equation}\label{eq:Hkin-Hkinsigma-compact}
		\|\phi\|_{\cH_{\kin,\sigma}(\cD_T\times\bR_v^3)}
		\simeq_R
		\|\phi\|_{\cH_{\kin}(\cD_T\times\bR_v^3)} .
	\end{equation}
	Now let $f$ be as in (2) and assume, say,
	$f\in\cH_\kin(\cD_T\times\bR_v^3)$.  Fix $\zeta\in C_c^\infty(B_{2R})$
	with $\zeta=1$ on $B_R$, so that $\zeta f=f$, and choose
	$\phi_n\in C_c^\infty([0,T]\times\overline\cD\times\bR_v^3)$ with
	$\phi_n\to f$ in $\cH_\kin$.  By (1), $\zeta\phi_n\to f$ in $\cH_\kin$,
	and each $\zeta\phi_n$ has velocity support in $B_{2R}$; by
	\eqref{eq:Hkin-Hkinsigma-compact} applied to
	$\zeta\phi_n-\zeta\phi_m$, the sequence $(\zeta\phi_n)$ is Cauchy in the
	complete space $\cH_{\kin,\sigma}(\cD_T\times\bR_v^3)$.  Its limit there
	coincides with $f$, since convergence in either of the two kinetic norms
	implies convergence in $L^2_{\mathrm{loc}}(\cD_T\times\bR_v^3)$ by
	\eqref{eq:local-H1-H_sigma}.  Letting $n\to\infty$ in
	\eqref{eq:Hkin-Hkinsigma-compact} gives the asserted equivalence, and
	exchanging the roles of the two spaces proves the converse implication.
\end{proof}

\begin{prop}\label{H_kin=weak_landau}
	\,
	
	\begin{enumerate}
		
		\item
		\(\cH_{\kin,\sigma}(\cD_T\times\bR_v^3)\) coincides with the set of all
		\(f\in L^2(\cD_T;H^1_\sigma)\) for which there exist
		\begin{align}\label{2607151229}
			(\partial_t+v\cdot \nabla_x)f
			&\in L^2(\cD_T;H^{-1}_\sigma),
			\qquad
			f(0,\cdot),\,f(T,\cdot)
			\in L^2(\cD\times \bR_v^3),
			\qquad 
			f_\gamma\in L^2(\gamma_T),
		\end{align}
		such that, for every
		$\phi\in
		C_c^\infty
		([0,T]\times \overline{\cD}\times \bR_v^3)$,
		\begin{align}\label{2607151226}
			\begin{aligned}
				&\iint_{\cD_T}
				\Big(
				\la (\partial_t+v\cdot\nabla_x)f,\phi\ra_v
				+
				\la f,(\partial_t+v\cdot\nabla_x)\phi\ra_v
				\Big)
				\dd x\dd t
				\\
				={}&
				\iint_{\cD\times\bR_v^3}
				(f\phi)(T,\cdot)
				\dd x\dd v
				-
				\iint_{\cD\times\bR_v^3}
				(f\phi)(0,\cdot)
				\dd x\dd v
				+
				\iint_{\gamma_T}
				f_\gamma\,
				\phi|_\gamma\,
				\operatorname{sgn}(n_x\cdot v)
				\dd\gamma\dd t.
			\end{aligned}
		\end{align}
		
		\item
		\(\cH_{\kin,\sigma}(\cD\times\bR_v^3)\) coincides with the set of all
		\(f\in L^2(\cD;H^1_\sigma)\) for which there exist
		\[
		v\cdot\nabla_x f\in L^2(\cD;H^{-1}_\sigma),
		\qquad
		f_\gamma\in L^2(\gamma),
		\]
		such that, for every
		$\phi\in
		C_c^\infty
		(\overline{\cD}\times \bR_v^3)$,
		\begin{align}\label{eq:H_kin_sigma-weak-stationary}
			&\int_{\cD}
			\Big(
			\la v\cdot\nabla_x f,\phi\ra_v
			+
			\la f,v\cdot\nabla_x\phi\ra_v
			\Big)
			\dd x=
			\int_\gamma
			f_\gamma\,
			\phi|_\gamma\,
			\operatorname{sgn}(n_x\cdot v)
			\dd\gamma.
		\end{align}
		
	\end{enumerate}
	In either case the objects occurring in these conditions are
	uniquely determined by $f$, and they coincide with the transport term and
	the traces attached to $f$ by Definition~\ref{def:H_kin_sigma} and
	Proposition~\ref{2607151218}. Uniqueness follows by taking the difference
	of two admissible choices and testing with $\phi$ supported successively
	away from $\{t=0,T\}\cup\gamma_T$, near $\{t=0\}$, near $\{t=T\}$ and near
	$\gamma_T$; the coincidence then follows because $f\in\cH_{\kin,\sigma}$
	satisfies the same identity with those traces, by
	Proposition~\ref{2607151218}.
\end{prop}

\begin{proof}
	We prove only (1), since the proof of (2) is identical after omitting
	the temporal trace terms.
	
	If $f\in\cH_{\kin,\sigma}(\cD_T\times\bR_v^3)$,
	then \eqref{2607151229} and \eqref{2607151226} follow from
	Proposition~\ref{2607151218}. This proves one inclusion.
	
	Conversely, suppose that
	$f\in L^2(\cD_T;H_\sigma^1)$ satisfies
	\eqref{2607151229} and \eqref{2607151226}. Let
	$\chi\in C_c^\infty(\bR_v^3)$ be radial and satisfy \eqref{eq:radial-cutoff}.
	Set
	\[
	\chi_R(v):=\chi(v/R),
	\qquad
	f_R:=\chi_Rf.
	\]
	Since $\chi_R$ depends only on $v$, applying \eqref{2607151226} with the test function
	$\chi_R\varphi$, we see that $f_R$ satisfies the same identity with the transport term
	\[
	(\p_t+v\cdot\nabla_x)f_R
	:=
	\chi_R(\p_t+v\cdot\nabla_x)f,
	\]
	and the traces
	\[
	f_R(0)=\chi_Rf(0),
	\qquad
	f_R(T)=\chi_Rf(T),
	\qquad
	(f_R)_\gamma=\chi_Rf_\gamma.
	\]
	
	Since $f_R$ and
	$(\p_t+v\cdot\nabla_x)f_R$ have velocity support in $B_{2R}$,
	Lemma~\ref{lem:H_sigma-Hilbert}.(4), namely
	\eqref{eq:Hv-Hsigma-compact}, gives
	\[
	f_R\in L^2(\cD_T;H_v^1),
	\qquad
	(\p_t+v\cdot\nabla_x)f_R
	\in L^2(\cD_T;H_v^{-1}).
	\]
	Hence Proposition~\ref{H_kin=weak} yields
	\[
	f_R\in\cH_\kin(\cD_T\times\bR_v^3).
	\]
	
	Choose
	\[
	f_{R,k}\in
	C_c^\infty
	\big(
	[0,T]\times\overline{\cD}\times\bR_v^3
	\big)
	\]
	such that
	\[
	f_{R,k}\to f_R
	\qquad
	\text{in }\cH_\kin.
	\]
	By Lemma~\ref{lem:compact-velocity-support}.(1) and $\chi_{2R}f_R=f_R$, we
	may replace $f_{R,k}$ by $\chi_{2R}f_{R,k}$, and so assume that the
	functions $f_{R,k}$ have velocity support in $B_{4R}$. On functions and
	distributions with velocity support in a fixed compact set, the weighted
	and unweighted velocity norms are equivalent by
	Lemma~\ref{lem:compact-velocity-support}.(2). Therefore,
	$f_{R,k}\to f_R$ in $\cH_{\kin,\sigma}$,
	and hence
	\[
	f_R\in
	\cH_{\kin,\sigma}(\cD_T\times\bR_v^3).
	\]
	
	Finally, we let $R\to\infty$. Estimate
	\eqref{eq:radial-cutoff-convergence} gives
	\begin{equation*}
		f_R=\chi_Rf\longrightarrow f
		\qquad
		\text{in }L^2(\cD_T;H_\sigma^1).
	\end{equation*}
	For the transport term, write
	$F:=(\p_t+v\cdot\nabla_x)f\in L^2(\cD_T;H^{-1}_\sigma)$ and let
	$u(t,x,\cdot)\in H^1_\sigma$ be its Riesz representative in the sense of
	Remark~\ref{rem:Riesz-H_sigma}, so that
	$\|F(t,x,\cdot)\|_{H^{-1}_\sigma}=\|u(t,x,\cdot)\|_{H^1_\sigma}$. For
	$\varphi\in H^1_\sigma$,
	\[
	\la \chi_RF-F,\varphi\ra
	=
	\big(u,(\chi_R-1)\varphi\big)_{H^1_\sigma}
	=
	\int_{\bR^3}\Big(\sigma^{ij}\p_iu\,\p_j\big[(\chi_R-1)\varphi\big]
	+\sigma^{ij}v_iv_j\,u\,(\chi_R-1)\varphi\Big)\dd v .
	\]
	For the two terms carrying the factor $(\chi_R-1)$ we first apply the
	Cauchy--Schwarz inequality for the nonnegative matrix $\sigma$, which
	bounds their sum by
	\[
	\left(\int_{\{|v|\ge R\}}
	\big(\sigma^{ij}\p_iu\,\p_ju+\sigma^{ij}v_iv_j\,u^2\big)\dd v\right)^{1/2}
	\|\varphi\|_{H^1_\sigma},
	\]
	and the first factor tends to $0$ as $R\to\infty$ by dominated
	convergence, uniformly in $\varphi$. The remaining commutator term
	contains
	$\nabla_v\chi_R=a_R(v)v$ with $|a_R|\le CR^{-2}$ supported in
	$\{R\le|v|\le2R\}$, so that, by the Cauchy--Schwarz inequality for the
	nonnegative matrix $\sigma$,
	\[
	\left|\int_{\bR^3}\sigma^{ij}\p_iu\,\varphi\,\p_j\chi_R\dd v\right|
	\le
	\frac{C}{R^{2}}
	\left(\int_{\bR^3}\sigma^{ij}\p_iu\,\p_ju\dd v\right)^{1/2}
	\left(\int_{\{R\le|v|\le2R\}}\sigma^{ij}v_iv_j\varphi^2\dd v\right)^{1/2}
	\longrightarrow 0 .
	\]
	Hence $\chi_RF\to F$ in $H^{-1}_\sigma$ for a.e. $(t,x)$, and since
	$\|\chi_RF\|_{H^{-1}_\sigma}\le C\|F\|_{H^{-1}_\sigma}$ uniformly in $R$
	by Lemma~\ref{lem:H_sigma-Hilbert}.(3) and duality, dominated convergence
	gives
	\begin{equation*}
		(\p_t+v\cdot\nabla_x)f_R=\chi_R(\p_t+v\cdot\nabla_x)f
		\longrightarrow (\p_t+v\cdot\nabla_x)f
		\qquad
		\text{in }L^2(\cD_T;H_\sigma^{-1}).
	\end{equation*}
	Together with $f_R(0)=\chi_Rf(0)\to f(0)$ in
	$L^2(\cD\times\bR_v^3)$ and $(f_R)_\gamma=\chi_Rf_\gamma\to f_\gamma$ in
	$L^2(\gamma_T)$, both by dominated convergence, this shows that
	$(f_R)_R$ is a Cauchy family in $\cH_{\kin,\sigma}(\cD_T\times\bR_v^3)$.
	That space being complete, $f_R$ converges in it to some limit. This
	limit must equal $f$, since $f_R\to f$ already in $L^2(\cD_T;H^1_\sigma)$.
	Therefore
	$f\in\cH_{\kin,\sigma}(\cD_T\times\bR_v^3)$.
\end{proof}

\begin{lemma}
	\label{lem:transport-tail-L2}\,
	
	\begin{enumerate}
		\item
		Let $f\in\cH_{\kin,\sigma}(\cD_T\times\bR^3_v)$.  Then there exists
		$\eta_f\in C([0,T])$ such that
		\begin{equation*}
			|\eta_f(t)|
			\le
			C_\cD\|f(t)\|_{L^2(\cD\times\bR^3_v)}^2
			\qquad
			\text{for every }t\in[0,T],
		\end{equation*}
		and such that, for every $0\le s<t\le T$,
		\begin{align}\label{eq:transport-tail-L2-evolutionary}
			\begin{aligned}
				&\|\la v\ra^{1/2}f\|_{L^2((s,t)\times \cD\times\bR^3_v)}^2\\
				\leq\,&C_{\cD}\Big(
				\|f\|_{L^2((s,t)\times \cD;H_\sigma^1)}^2
				+
				\|(\partial_t+v\cdot\nabla_x)f\|_
				{L^2((s,t)\times \cD;H_\sigma^{-1})}^2
				+
				\|f_\gamma\|_
				{L^2((s,t)\times \gamma)}^2\Big)
				\\
				&+\eta_f(t)-\eta_f(s).
			\end{aligned}
		\end{align}
		In particular,
		\begin{equation}\label{260829706}
			\|\la v\ra^{1/2}f\|_{L^2(\cD_T\times\bR^3_v)}
			\leq
			C_{\cD}
			\|f\|_{\cH_{\kin,\sigma}
				(\cD_T\times\bR^3_v)}.
		\end{equation}
		\item
		For every $f\in\cH_{\kin,\sigma}(\cD\times\bR^3_v)$,
		\begin{equation}\label{eq:transport-tail-L2-stationary}
			\|\la v\ra^{1/2}f\|_{L^2(\cD\times\bR^3_v)}
			\leq C_{\cD}\Big(
			\|f\|_{L^2(\cD;H_\sigma^1)}
			+
			\|v\cdot\nabla_xf\|_
			{L^2(\cD;H_\sigma^{-1})}
			+
			\|f_\gamma\|_{L^2(\gamma,\dd\gamma)}
			\Big).
		\end{equation}
		In particular,
		\begin{equation*}
			\|\la v\ra^{1/2}f\|_{L^2(\cD\times\bR^3_v)}
			\leq
			C_{\cD}
			\|f\|_{\cH_{\kin,\sigma}
				(\cD\times\bR^3_v)}.
		\end{equation*}
	\end{enumerate}
\end{lemma}

\begin{proof}
	We prove only the evolutionary estimate, since the stationary one follows
	from the same argument once the endpoint terms in time are discarded.
	
	Fix $x_0\in\cD$, let
	\begin{equation*}
		\psi(x,v)
		:=
		(x-x_0)\cdot\frac{v}{\la v\ra},
	\end{equation*}
	choose a radial function $\chi\in C_c^\infty(\bR^3)$ satisfying
	\eqref{eq:radial-cutoff}, and set
	\begin{equation*}
		\chi_R(v):=\chi(v/R),
		\qquad
		\psi_R(x,v):=\chi_R(v)\psi(x,v),
		\qquad R\geq1.
	\end{equation*}
	Since $\cD$ is bounded,
	\begin{equation*}
		|\psi_R|\leq C_{\cD},
		\qquad
		v\cdot\nabla_x\psi_R
		=
		\chi_R(v)\frac{|v|^2}{\la v\ra},
		\qquad
		\nabla_v\psi_R
		=
		\chi_R\nabla_v\psi
		+
		\psi\nabla_v\chi_R .
	\end{equation*}
	As $\nabla_v\chi_R$ is radial and supported in $\{R\leq|v|\leq2R\}$,
	\begin{equation*}
		(\nabla_v\psi_R)^T
		\sigma(v)\nabla_v\psi_R
		\leq
		C_{\cD}\la v\ra^{-3}
		\qquad\text{uniformly in }R\ge1,
	\end{equation*}
	and hence the product rule and \eqref{eq:local-H1-H_sigma} yield
	\begin{equation*}
		\|\psi_Rh\|_{H_\sigma^1}
		\leq
		C_{\cD}\|h\|_{H_\sigma^1},
		\qquad R\geq1.
	\end{equation*}
	
	Let $f\in\cH_{\kin,\sigma}(\cD_T\times\bR^3_v)$ and put $g_R:=\psi_Rf$.
	Since $\chi_R$ has compact velocity support,
	Lemma~\ref{lem:compact-velocity-support} gives
	$g_R\in\cH_{\kin,\sigma}(\cD_T\times\bR^3_v)$, and
	\begin{equation*}
		(\partial_t+v\cdot\nabla_x)g_R
		=
		\psi_R(\partial_t+v\cdot\nabla_x)f
		+
		\chi_R\frac{|v|^2}{\la v\ra}f
	\end{equation*}
	in the sense of distributions.  We may therefore apply the evolutionary
	Green identity to the pair $(f,g_R)$.  Since
	\begin{equation*}
		\left\langle
		\psi_R(\partial_t+v\cdot\nabla_x)f,f
		\right\rangle_v
		=
		\left\langle
		(\partial_t+v\cdot\nabla_x)f,\psi_Rf
		\right\rangle_v ,
	\end{equation*}
	it follows that
	\begin{align*}
		&\iiint_{(s,t)\times\cD\times\bR^3_v}
		\chi_R(v)\frac{|v|^2}{\la v\ra}|f(\tau,x,v)|^2\dd v\dd x\dd \tau
		\\
		={}&
		\iint_{\cD\times\bR^3_v}
		\psi_R|f(t)|^2\dd v\dd x
		-
		\iint_{\cD\times\bR^3_v}
		\psi_R|f(s)|^2\dd v\dd x
		+
		\iint_{(s,t)\times \gamma}
		\psi_R|f_{\gamma}|^2\sgn(n_x\cdot v)\dd \gamma\dd \tau\\
		&-
		2\iint_{(s,t)\times\cD}\left\langle
		(\partial_t+v\cdot\nabla_x)f,\psi_Rf
		\right\rangle_v\dd x\dd \tau .
	\end{align*}
	Consequently, Young's inequality and \eqref{eq:Hsigma-multiplier} give
	\begin{align*}
		&\iiint_{(s,t)\times\cD\times\bR^3_v}
		\chi_R(v)\frac{|v|^2}{\la v\ra}|f|^2\dd v\dd x\dd \tau\\
		\leq\,&\iint_{\cD\times\bR^3_v}
		\psi_R|f(t)|^2\dd v\dd x
		-
		\iint_{\cD\times\bR^3_v}
		\psi_R|f(s)|^2\dd v\dd x
		+
		\iint_{(s,t)\times \gamma}
		\psi_R|f_{\gamma}|^2\sgn(n_x\cdot v)\dd \gamma\dd \tau\\
		&\qquad +
		C\left(\|(\partial_t+v\cdot\nabla_x)f\|_
		{L^2(\cD_T;H_\sigma^{-1})}^2+
		\|f\|_{L^2(\cD_T;H_\sigma^1)}^2
		\right).
	\end{align*}
	Now set
	\begin{equation*}
		\eta_f(t)
		:=
		2\iint_{\cD\times \bR^3_v}\psi(x,v)|f(t,x,v)|^2\dd v\dd x .
	\end{equation*}
	Since $0\leq \chi_R\leq 1$ and $\chi_R\to1$ pointwise, Fatou's lemma
	together with the dominated convergence theorem yields
	\begin{align*}
		&\iiint_{(s,t)\times\cD\times\bR^3_v}
		\frac{|v|^2}{\la v\ra}|f|^2\dd v\dd x\dd \tau
		\\
		\leq\,&C\left(
		\|f\|_{L^2(\cD_T;H_\sigma^1)}^2
		+
		\|(\partial_t+v\cdot\nabla_x)f\|_
		{L^2(\cD_T;H_\sigma^{-1})}^2
		+
		\|f_\gamma\|_
		{L^2(\gamma_T,\dd t\dd\gamma)}^2\right)
		\\
		&\quad +\frac 12\big(\eta_f(t)-\eta_f(s)\big).
	\end{align*}
	On $\{|v|\leq 1\}$, the local estimate \eqref{eq:local-H1-H_sigma} gives
	\begin{equation*}
		\|f\|_{L^2_{t,x,v}(|v|\leq1)}
		\leq
		C\|f\|_{L^2(\cD_T;H_\sigma^1)},
	\end{equation*}
	whereas $\frac{|v|^2}{\la v\ra}\geq \frac 12 \la v\ra$ on $\{|v|\geq 1\}$.
	Splitting the velocity space into these two regions proves
	\eqref{eq:transport-tail-L2-evolutionary}, and \eqref{260829706} then
	follows from Proposition~\ref{2607151218}.
	The continuity of $\eta_f$ on $[0,T]$ follows from $f\in C([0,T];L^2(\cD\times \bR^3_v))$ (see Proposition \ref{2607151218}).
	Therefore the proof is completed.
\end{proof}

\subsection{Nonlinear estimates in the weighted kinetic spaces}\label{ss:nonlinear-weighted-estimates}

For a Banach space \(X\), set
\[
L_a^p(A;X):=L^p(A,\dd a;X),
\]
and use the mixed-norm conventions
\begin{align}
	\|h\|_{L^{p,q}_{a,b}(A\times B)}
	:=\|h\|_{L_a^p(A;L_b^q(B))},\qquad 
	\|h\|_{L^{p,q,r}_{a,b,c}(A\times B\times C)}
	:=\|h\|_{L_a^p(A;L_b^q(B;L_c^r(C)))}.
	\label{eq:mixed-norm-convention}
\end{align}
In \eqref{eq:mixed-norm-convention} the exponents always carry
variable subscripts.  A pair of exponents without them, such as
\(L^{7/3,\infty}\), denotes a Lorentz space.
Theorem~\ref{260529130}, proved at the end of the subsection, gives the required
$H^{-1}_\sigma$ estimate for the nonlinear Landau operator.

\begin{thm}[Nonlinear estimate in kinetic energy spaces]
	\label{260529130}
	Let \(\ell\ge0\) and \(14<p<\infty\).  There exists
	\(C_{\ell,p}>0\), independent of \(T\) and of \(\cD\), such that the
	following assertions hold.
	\begin{enumerate}
		\item
		Let \(0<T<\infty\).  If
		\(\mu^{1/8}g\in\cH_{\kin,\sigma}(\cD_T\times\bR_v^3)\), then
		\begin{align}
			\label{eq:N-estimate-evolutionary-sharp}
			\|w^\ell N[g,f]\|_{L^2(\cD_T;H^{-1}_\sigma)}
			\le{}&
			C_{\ell,p}
			\|\mu^{1/8}g\|_{\cH_{\kin,\sigma}(\cD_T\times\bR_v^3)}
			\|w^\ell f\|_{L_t^\infty L_x^pH^1_\sigma(\cD_T\times\bR_v^3)} .
		\end{align}
		
		\item
		If \(\mu^{1/8}g\in\cH_{\kin,\sigma}(\cD\times\bR_v^3)\), then
		\begin{align}
			\label{eq:N-estimate-stationary-sharp}
			\|w^\ell N[g,f]\|_{L^2(\cD;H^{-1}_\sigma)}
			\le{}&
			C_{\ell,p}
			\|\mu^{1/8}g\|_{\cH_{\kin,\sigma}(\cD\times\bR_v^3)}
			\|w^\ell f\|_{L_x^pH^1_\sigma(\cD\times\bR_v^3)} .
		\end{align}
	\end{enumerate}
	In particular, if \(w^{\ell+2}f\in L^\infty\) and
	\(w^{\ell+2}D_vf\in L^{\infty,p,2}_{t,x,v}\), respectively
	\(L^{p,2}_{x,v}\), then
	\begin{align}
		\label{eq:N-estimate-evolutionary}
		\|w^\ell N[g,f]\|_{L^2(\cD_T;H^{-1}_\sigma)}
		\le{}&
		C_{\ell,p}(1+|\cD|^{1/p})
		\|\mu^{1/8}g\|_{\cH_{\kin,\sigma}(\cD_T\times\bR_v^3)}
		\notag\\
		&\times
		\left(
		\|w^{\ell+2}f\|_{L^\infty(\cD_T\times\bR_v^3)}
		+
		\|w^{\ell+2}D_vf\|_{L^{\infty,p,2}_{t,x,v}(\cD_T\times\bR_v^3)}
		\right),
		\\
		\label{eq:N-estimate-stationary}
		\|w^\ell N[g,f]\|_{L^2(\cD;H^{-1}_\sigma)}
		\le{}&
		C_{\ell,p}(1+|\cD|^{1/p})
		\|\mu^{1/8}g\|_{\cH_{\kin,\sigma}(\cD\times\bR_v^3)}
		\notag\\
		&\times
		\left(
		\|w^{\ell+2}f\|_{L^\infty(\cD\times\bR_v^3)}
		+
		\|w^{\ell+2}D_vf\|_{L^{p,2}_{x,v}(\cD\times\bR_v^3)}
		\right).
	\end{align}
\end{thm}

\begin{remark}\label{rem:stationary-nonlinear-p12}
	The stationary estimates \eqref{eq:N-estimate-stationary-sharp} and
	\eqref{eq:N-estimate-stationary} hold for $12<p<\infty$, with the
	remaining hypotheses of Theorem~\ref{260529130} unchanged. Indeed,
	$q:=2p/(p-2)$ satisfies $2<q<12/5$ and $1/p+1/q=1/2$.
	The stationary integrability gain \eqref{IntGain_std}, the stationary
	version of \eqref{eq:N-Maxwellian-multiplier}, and the two-slot bound
	\eqref{estq:N_st} give \eqref{eq:N-estimate-stationary-sharp}.
	The weighted bound \eqref{eq:N-estimate-stationary} follows as in the
	proof of Theorem~\ref{260529130}. The evolutionary restriction
	$p>14$ is unchanged.
\end{remark}

The theorem combines two ingredients, established next.  The first is a mixed-norm estimate for the nonlinear operator, Lemma~\ref{lem:velalg}, and the second is an integrability gain for $\cH_\kin$-functions.

\begin{lemma}[Gaussian-\(L^2\) first slot]
	\label{lem:velalg}
	Fix \(\ell\ge0\) and \(1/4\le\varpi<1/2\).  For every \(g\) with
	\(\mu^\varpi g\in L^2(\bR_v^3)\) and every \(f\in H^1_\sigma(\bR_v^3)\),
	\begin{equation}\label{eq:velalg}
		\norm{w^\ell N[g,f]}_{H^{-1}_\sigma(\bR_v^3)}
		\le C_{\ell,\varpi}
		\norm{\mu^\varpi g}_{L_v^2}
		\norm{w^\ell f}_{H^1_\sigma(\bR_v^3)} .
	\end{equation}
\end{lemma}

\begin{proof}
	Put \(u:=w^\ell f\), let \(\varphi\in C_c^\infty(\bR_v^3)\), and write
	\(M_\varpi:=\|\mu^\varpi g\|_{L^2_v}\), the conjugated derivatives
	\(D_i^\pm\) being those of \eqref{eq:velalg-AF-BF-definition}.  Since
	\(\la w^\ell N[g,f],\varphi\ra_v=\la N[g,f],w^\ell\varphi\ra_v\), inserting
	\(f=w^{-\ell}u\) into \eqref{eq:velalg-exact-weak-form} and using
	\[
	D^-f=w^{-\ell}\left(D^-u-\ell w^{-2}v\,u\right),
	\qquad
	D^+\left(w^\ell\varphi\right)
	=w^{\ell}\left(D^+\varphi+\ell w^{-2}v\,\varphi\right),
	\]
	we obtain
	\begin{equation*}
		\la w^\ell N[g,f],\varphi\ra_v
		=
		\la N[g,u],\varphi\ra_v
		+\cR_\ell ,
	\end{equation*}
	where
	\[
	\cR_\ell
	:=
	\ell\int_{\bR^3}w^{-2}
	\left\{
	u\,\sigma_G^{ij}v_j\,D_i^+\varphi
	-\varphi\,\bigl(\sigma_G^{ij}D_j^-u-b^iu\bigr)v_i
	\right\}\dd v
	+\ell^2\int_{\bR^3}w^{-4}\sigma_G^{ij}v_iv_j\,u\varphi\dd v .
	\]
	The main term is bounded by Lemma~\ref{260715201}.(1), which gives
	\[
	\left|\la N[g,u],\varphi\ra_v\right|
	\le C_\varpi M_\varpi\|u\|_{H^1_\sigma}\|\varphi\|_{H^1_\sigma}.
	\]
	Every term of \(\cR_\ell\) carries an extra factor \(w^{-2}v\), which only
	improves the estimates just obtained.  Indeed
	\eqref{eq:AF-sigma-form} and
	\(\bigl(\sigma^{ij}v_iv_j\bigr)^{1/2}\le C\la v\ra^{-1/2}\) give
	\[
	\left|\ell\int_{\bR^3}w^{-2}u\,\sigma_G^{ij}v_j\,D_i^+\varphi\dd v\right|
	\le
	C_\ell M_\varpi
	\|\la v\ra^{-5/2}u\|_{L^2_v}
	\|\sigma^{1/2}D^+\varphi\|_{L^2_v}
	\le
	C_{\ell,\varpi}M_\varpi\|u\|_{H^1_\sigma}\|\varphi\|_{H^1_\sigma}
	\]
	by \eqref{eq:local-H1-H_sigma} and \eqref{eq:velalg-Dpm}; the term
	carrying \(\sigma_G^{ij}D_j^-u\) is symmetric in \(u\) and \(\varphi\), and
	\(|w^{-4}\sigma_G^{ij}v_iv_j|\le C_\varpi M_\varpi\la v\ra^{-5}\) handles
	the last one.
	For the term carrying \(b\), the splitting \eqref{eq:velalg-BF-split} and
	\eqref{eq:BF-far-component-bounds} give
	\(w^{-2}|v||b^{\rm far}|\le C_\varpi M_\varpi\la v\ra^{-3}\), while the
	local part is estimated exactly as above through \eqref{eq:velalg-local-B-L6} and
	\eqref{eq:velalg-cube-sobolev}, the additional factor \(w^{-1}\) only
	improving the dyadic sum.  Taking the supremum over
	\(\|\varphi\|_{H^1_\sigma}\le1\) and using the density of
	\(C_c^\infty(\bR_v^3)\) in \(H^1_\sigma(\bR_v^3)\) proves \eqref{eq:velalg}.
\end{proof}

\begin{lemma}[Velocity and mixed-norm estimates]
	\label{lemma:estN}
	Let \(\ell\ge0\) and \(1/4\le\varpi<1/2\), and let \(p_x,\,p_t,\,q_x,\,q_t\in [2,\infty]\) satisfy
	\begin{align}\label{260909908}
		\frac 1{p_x}+\frac 1{q_x}=\frac 1{p_t}+\frac 1{q_t}=\frac 12.
	\end{align}
	Then
	\begin{align}
		\label{estq:N_st}
		\|w^\ell N[g,f]\|_{L^2(\cD;H^{-1}_\sigma)}
		\le{}&
		C_{\ell,\varpi}
		\|\mu^\varpi g\|_{L^{q_x,2}_{x,v}}
		\|w^\ell f\|_{L_x^{p_x}H^1_\sigma},
	\end{align}
	and, for every \(T>0\),
	\begin{align}
		\label{estq:N_dy}
		\|w^\ell N[g,f]\|_{L^2(\cD_T;H^{-1}_\sigma)}
		\le{}&
		C_{\ell,\varpi}
		\|\mu^\varpi g\|_{L^{q_t,q_x,2}_{t,x,v}}
		\|w^\ell f\|_{L_{t,x}^{p_t,p_x}H^1_\sigma}.
	\end{align}
	The estimates hold whenever the right-hand sides are finite.
	All mixed norms are taken over the displayed phase-space domain, and
	\(C_{\ell,\varpi}\) is independent of \(p_x\), \(p_t\), \(T\), and \(\cD\).
\end{lemma}

\begin{proof}
	Apply \eqref{eq:velalg} pointwise in \(x\), and in \(t,x\).
	Using H\"older's
	inequality with \eqref{260909908}, we obtain \eqref{estq:N_st} and \eqref{estq:N_dy}.
\end{proof}

\begin{prop}[Integrability gain for evolutionary kinetic functions]
	\label{lem:int_gain}
	Let \(0<T<\infty\).  There exists a universal constant \(C>0\)
	such that, for every $f\in\cH_\kin(\cD_T\times\bR_v^3)$, we have
	\begin{equation}
		\label{2603291118}
		\|f\|_
		{L^{7/3,\infty}(\cD_T\times\bR_v^3)}
		\le
		C\|f\|_{\cH_\kin(\cD_T\times\bR_v^3)}.
	\end{equation}
	Consequently, for every \(2\le r<7/3\),
	\begin{equation}
		\label{2603291119}
		\|f\|_{L^r(\cD_T\times\bR_v^3)}
		\le
		C_r\|f\|_{\cH_\kin(\cD_T\times\bR_v^3)},
	\end{equation}
	where \(C_r\) is independent of \(T\) and \(\cD\).
\end{prop}

\begin{proof}
	Set \(\mathfrak T:=\partial_t+v\cdot\nabla_x\).  There exist
	\(g_0,g_1,g_2,g_3\in L^2(\cD_T\times\bR_v^3)\) such that
	\begin{equation}
		\label{eq:local-gain-source-decomposition}
		(\mathfrak T-\Delta_v)f
		=g_0+\sum_{i=1}^3D_ig_i
	\end{equation}
	and
	\begin{equation}
		\label{eq:local-gain-source-bound}
		\sum_{i=0}^3\|g_i\|_{L^2}
		\le
		C\|f\|_{\cH_\kin(\cD_T\times\bR_v^3)}.
	\end{equation}
	Choose \(\alpha_n\in C_c^\infty(\bR)\) such that
	\[
	0\le\alpha_n\le1,
	\qquad
	\alpha_n\uparrow1
	\quad\text{pointwise},
	\]
	and define
	\[
	\beta_n(r):=\int_0^r\int_0^s\alpha_n(\tau)\dd\tau\dd s.
	\]
	Then \(\beta_n\) is smooth and convex, with
	\begin{equation}
		\label{eq:local-gain-renormalizations}
		\begin{gathered}
			\beta_n(0)=\beta_n'(0)=0,
			\qquad
			\|\beta_n'\|_{L^\infty(\bR)}<\infty,
			\\
			0\le\beta_n(r)\uparrow\frac{r^2}{2},
			\qquad
			|\beta_n'(r)|\le|r|,
			\qquad
			0\le\beta_n''(r)\le1,\qquad 
			|\beta_n'(r)|^2\le2\beta_n(r).
		\end{gathered}
	\end{equation}
	The last inequality follows from \(0\le\beta_n''\le1\): for \(r>0\) one
	has \(\beta_n'(s)\ge(\beta_n'(r)-r+s)_+\) on \([0,r]\), whence
	\(\beta_n(r)\ge|\beta_n'(r)|^2/2\), and \(r<0\) follows by reflection.
	Extend
	\[
	U_n:=\beta_n(f)\mathbf1_{\cD_T\times\bR_v^3}
	\]
	by zero to \(\bR^7\).  
	The distributional derivatives of the
	characteristic functions are
	\[
	\partial_t\mathbf1_{(0,T)}=\delta_{t=0}-\delta_{t=T},
	\qquad
	\nabla_x\mathbf1_\cD=-n_x\delta_{\partial\cD}.
	\]
	Thus the boundary part of $\mathfrak T U_n$ is exactly
	\[
	\beta_n(f(0))\delta_{t=0}
	-\beta_n(f(T))\delta_{t=T}
	-(n_x\cdot v)\beta_n(f_\gamma)
	\mathbf1_{(0,T)}\,\delta_{\partial\cD}.
	\]
	Each term is a finite measure: $\beta_n(f(0))$ and
	$\beta_n(f(T))$ are integrable, and
	$\beta_n(f_\gamma)\in L^1(\gamma_T)$.
	In particular, restriction to $|n_x\cdot v|<\eta$ has total
	variation tending to zero as $\eta\downarrow0$; no additional
	measure is supported on grazing velocities.  We apply
	\eqref{reno_dyn} to
	\eqref{eq:local-gain-source-decomposition} written as
	\(\mathfrak Tf=g_0+D_i(g_i+D_if)\).  The field \(D_if\) contributes
	\(-\beta_n''(f)|D_vf|^2\), which is nonpositive and is discarded, and
	\(-\beta_n'(f)D_if\,D_i\varphi=-D_i\beta_n(f)\,D_i\varphi\), which is
	the weak form of \(\Delta_v\beta_n(f)\) and is moved to the left-hand
	side.  Hence, in \(\cD'(\bR^7)\),
	\begin{align}
		\label{eq:local-gain-zero-extension-inequality}
		(\mathfrak T-\Delta_v)U_n
		\le{}&
		\mathfrak m_n
		+\beta_n'(f)g_0
		-\beta_n''(f)D_if\,g_i
		+D_i\bigl(\beta_n'(f)g_i\bigr).
	\end{align}
	Here the repeated index \(i\) is summed, all interior terms are
	extended by zero, and \(\mathfrak m_n\) is the positive measure
	supported on the initial and inflow sets:
	\begin{align*}
		\mathfrak m_n(E)
		:={}&
		\iint_{\cD\times\bR_v^3}
		\mathbf1_E(0,x,v)\beta_n(f(0,x,v))\dd x\dd v+
		\int_{\gamma_T^-}
		\mathbf1_E(t,x,v)\beta_n(f_{\gamma^-}(t,x,v))\dd\gamma.
	\end{align*}
	The nonpositive dissipation, terminal, and outgoing-boundary terms
	have been discarded.  By
	\eqref{eq:local-gain-renormalizations},
	\eqref{eq:local-gain-source-bound}, and Cauchy--Schwarz,
	\begin{align}
		\label{eq:local-gain-L1-source-bound}
		&\|\mathfrak m_n\|_{\cM(\bR^7)}
		+\|\beta_n'(f)g_0\|_{L^1}
		+\sum_{i=1}^3\|\beta_n''(f)D_if\,g_i\|_{L^1}
		+\sum_{i=1}^3\|\beta_n'(f)g_i\|_{L^1}
		\le
		C\|f\|_{\cH_\kin(\cD_T\times\bR_v^3)}^2.
	\end{align}
	
	Let \(\zeta\in C_c^\infty((0,\infty)\times\bR_x^3\times\bR_v^3)\)
	be nonnegative with \(\int_{\bR^7}\zeta=1\), let
	\(\zeta_\varepsilon(Z):=\varepsilon^{-7}\zeta(Z/\varepsilon)\), and
	set \(F_{n,\varepsilon}:=U_n*\zeta_\varepsilon\).  Since \(\zeta\) is
	supported in positive times and \(U_n\) vanishes for \(t<0\), so does
	\(F_{n,\varepsilon}\).  The commutator between
	ordinary convolution and \(\mathfrak T\) is
	\[
	R_{n,\varepsilon}(Z)
	:=\iiint_{\bR^7}
	U_n(Y)(v-w)\cdot\nabla_x
	\zeta_\varepsilon(Z-Y)\dd Y,
	\qquad Z=(t,x,v),\quad Y=(s,y,w).
	\]
	It has the form \(R_{n,\varepsilon}=U_n*K_\varepsilon\), where
	\(K\in C_c^\infty(\bR^7)\) and \(\int K=0\).  Hence
	\begin{equation}
		\label{eq:local-gain-commutator}
		\|R_{n,\varepsilon}\|_{L^1(\bR^7)}
		\longrightarrow0
		\qquad\text{as }\varepsilon\to0.
	\end{equation}
	Set
	\begin{align*}
		A_{n,\varepsilon}
		&:={}
		\bigl(
		\mathfrak m_n+\beta_n'(f)g_0
		-\beta_n''(f)D_if\,g_i
		\bigr)*\zeta_\varepsilon,\qquad 
		B_{i,n,\varepsilon}
		:=\bigl(\beta_n'(f)g_i\bigr)*\zeta_\varepsilon.
	\end{align*}
	Convolving \eqref{eq:local-gain-zero-extension-inequality} gives
	\begin{equation*}
		(\mathfrak T-\Delta_v)F_{n,\varepsilon}
		\le
		A_{n,\varepsilon}+D_iB_{i,n,\varepsilon}
		+R_{n,\varepsilon}.
	\end{equation*}
	By \eqref{eq:local-gain-L1-source-bound} and Young's inequality,
	\[
	\|A_{n,\varepsilon}\|_{L^1}
	+\sum_{i=1}^3\|B_{i,n,\varepsilon}\|_{L^1}
	\le
	C\|f\|_{\cH_\kin(\cD_T\times\bR_v^3)}^2,
	\]
	uniformly in \(n\) and \(\varepsilon\); moreover,
	\eqref{eq:local-gain-commutator} gives
	\(\|R_{n,\varepsilon}\|_{L^1}\to0\) for each fixed \(n\).
	
	Let \(\mathfrak G\) denote the whole-space Kolmogorov potential of
	\eqref{Gf} and \(\mathfrak G_i[H]:=\mathfrak G[D_iH]\).  For
	\(H\in L^1(\bR^7)\) the latter is defined by density: on
	\(C_c^\infty(\bR^7)\), \eqref{eq:weak-L1-endpoint} bounds
	\(\mathfrak G_i[H]\) in \(L^{14/13,\infty}(\bR^7)\) by \(\|H\|_{L^1}\),
	so the map extends continuously to \(L^1\), and the extension equals
	\(\mathfrak G[D_iH]\) in \(\cD'\).
	
	By \eqref{eq:weak-L1-endpoint} with \(d=3\), so that \(Q_1=7/6\),
	\begin{equation}
		\label{eq:local-gain-potential-endpoints}
		\|\mathfrak G[H]\|_{L^{7/6,\infty}(\bR^7)}
		\le C\|H\|_{L^1},
	\end{equation}
	while \eqref{eq:weak-Lp-endpoint} with \(d=3\) and \(p=14/13\), for which
	\(R_p=7/6\), gives
	\begin{equation}
		\label{eq:local-gain-critical-lorentz}
		\|\mathfrak G_i[H]\|_{L^{7/6,\infty}(\bR^7)}
		\le C\|H\|_{L^{14/13,\infty}(\bR^7)} .
	\end{equation}
	Set
	\begin{equation*}
		X_n:=\|U_n\|_{L^{7/6,\infty}(\bR^7)} ,
	\end{equation*}
	which is finite for each fixed \(n\).  Indeed, \(0\le U_n\le|f|^2/2\)
	gives \(U_n\in L^1\), and \(U_n\le\|\beta_n'\|_{L^\infty}|f|\) gives
	\(U_n\in L^2\).  Since
	\(|\beta_n'(f)g_i|\le\sqrt2\,U_n^{1/2}|g_i|\) by
	\eqref{eq:local-gain-renormalizations}, and since
	\(\|U_n^{1/2}\|_{L^{7/3,\infty}}^2=X_n\), H\"older's inequality in
	Lorentz spaces with \(\tfrac{13}{14}=\tfrac37+\tfrac12\) and Young's
	inequality for the mollification give
	\begin{equation}
		\label{eq:local-gain-B-lorentz}
		\|B_{i,n,\varepsilon}\|_{L^{14/13,\infty}(\bR^7)}
		\le C\|U_n^{1/2}\|_{L^{7/3,\infty}}\|g_i\|_{L^2}
		\le C\|f\|_{\cH_\kin(\cD_T\times\bR_v^3)}X_n^{1/2},
	\end{equation}
	uniformly in \(\varepsilon\).
	
	The function \(F_{n,\varepsilon}\) is smooth, nonnegative, bounded,
	integrable, and vanishes for \(t<0\).  Hence
	\(F_{n,\varepsilon}
	=\mathfrak G\bigl[(\mathfrak T-\Delta_v)F_{n,\varepsilon}\bigr]\)
	by the Duhamel formula for the Kolmogorov operator.  Since its
	kernel is nonnegative, monotonicity of \(\mathfrak G\) gives the
	pointwise bound
	\[
	0\le F_{n,\varepsilon}
	\le
	|\mathfrak G[A_{n,\varepsilon}]|
	+\sum_{i=1}^3|\mathfrak G_i[B_{i,n,\varepsilon}]|
	+|\mathfrak G[R_{n,\varepsilon}]|.
	\]
	Consequently, \eqref{eq:local-gain-L1-source-bound},
	\eqref{eq:local-gain-commutator},
	\eqref{eq:local-gain-potential-endpoints},
	\eqref{eq:local-gain-critical-lorentz}, and
	\eqref{eq:local-gain-B-lorentz} give
	\[
	\limsup_{\varepsilon\to0}
	\|F_{n,\varepsilon}\|_{L^{7/6,\infty}(\bR^7)}
	\le
	C\|f\|_{\cH_\kin(\cD_T\times\bR_v^3)}^2
	+C\|f\|_{\cH_\kin(\cD_T\times\bR_v^3)}X_n^{1/2}.
	\]
	Since \(U_n\in L^1(\bR^7)\), \(F_{n,\varepsilon}\to U_n\) at every
	Lebesgue point of \(U_n\) as \(\varepsilon\to0\); by the Fatou property
	of Lorentz spaces,
	\[
	X_n
	\le
	C\|f\|_{\cH_\kin(\cD_T\times\bR_v^3)}^2
	+C\|f\|_{\cH_\kin(\cD_T\times\bR_v^3)}X_n^{1/2} .
	\]
	Since \(X_n<\infty\), Young's inequality absorbs the last term and leaves
	\[
	X_n\le C\|f\|_{\cH_\kin(\cD_T\times\bR_v^3)}^2
	\]
	uniformly in \(n\).  Letting \(n\to\infty\) and using
	\(U_n\uparrow|f\mathbf1_{\cD_T\times\bR_v^3}|^2/2\), we obtain
	\eqref{2603291118} after taking the square root.  Finally,
	\eqref{2603291119} follows by interpolating \eqref{2603291118} with
	\(\|f\|_{L^2}\le C\|f\|_{\cH_\kin(\cD_T\times\bR_v^3)}\).
\end{proof}

\begin{lemma}[Stationary Kolmogorov potentials]
	\label{lem:stationary-kolmogorov-potential-sharp}
	Let $\mathfrak G$ be the Kolmogorov fundamental solution
	\eqref{260424252} with $d=3$.
	For $s>0$ and
	$Z=(x,v)$, $Y=(y,w)$ in $\bR^3\times\bR^3$ put
	\[
	p_s(Z;Y):=\mathfrak G(s,x-y-sw,v-w),
	\]
	and define
	\begin{equation}\label{eq:stationary-green-kernels}
		\mathfrak G_\st(Z;Y):=\int_0^\infty p_s(Z;Y)\dd s,
		\qquad
		\mathfrak K_{\st,i}(Z;Y):=-\int_0^\infty\p_{w^i}p_s(Z;Y)\dd s ,
	\end{equation}
	together with the operators
	\begin{equation}\label{eq:stationary-Gf}
		\mathfrak G_\st[H](Z):=\int_{\bR^3_x\times \bR^3_v}\mathfrak G_\st(Z;Y)H(Y)\dd Y,
		\qquad
		\mathfrak G_{\st,i}[H](Z):=\int_{\bR^3_x\times \bR^3_v}\mathfrak K_{\st,i}(Z;Y)H(Y)\dd Y .
	\end{equation}
	Then,
	\begin{align}
		\|\mathfrak G_\st[H]\|_{L^{6/5,\infty}(\bR^3_x\times \bR^3_v)}
		+\|\mathfrak G_{\st,i}[H]\|_{L^{12/11,\infty}(\bR^3_x\times \bR^3_v)}
		&\le C\|H\|_{L^1(\bR^3_x\times \bR^3_v)},
		\label{eq:stationary-weak-L1-endpoint}
		\\
		\|\mathfrak G_{\st,i}[H]\|_{L^{6/5,\infty}(\bR^3_x\times \bR^3_v)}
		&\le C\|H\|_{L^{12/11,\infty}(\bR^3_x\times \bR^3_v)}.
		\label{eq:stationary-weak-Lp-endpoint}
	\end{align}
\end{lemma}

\begin{proof}
	Write $Q:=4d=12$ and, for fixed $s>0$,
	\[
	\xi:=x-y-sw,
	\qquad
	\eta:=v-w,
	\]
	so that $p_s(Z;Y)=\mathfrak G(s,\xi,\eta)$.
	
	\medskip
	\noindent\emph{Step 1. Bounds on $p_s$ and $\p_{w^i}p_s$.}
	By \eqref{260424252},
	\begin{equation}\label{eq:stationary-kernel-scaling}
		\mathfrak G(s,\xi,\eta)
		=s^{-Q/2}\,
		\mathfrak G\bigl(1,s^{-3/2}\xi,s^{-1/2}\eta\bigr),
	\end{equation}
	and $\mathfrak G(1,\cdot,\cdot)$ and its first derivatives lie in $L^1\cap L^\infty(\bR^6)$.
	Differentiating \eqref{eq:stationary-kernel-scaling} and changing
	variables, we have
	\begin{equation}\label{eq:stationary-scaled-norms}
		\begin{gathered}
			\|\mathfrak G(s,\cdot)\|_{L^1}=1,
			\qquad
			\|\mathfrak G(s,\cdot)\|_{L^\infty}\le Cs^{-\frac Q2},
			\\
			\bigl\|\bigl(s\p_{\xi_i}+\p_{\eta_i}\bigr)\mathfrak G(s,\cdot)\bigr\|_{L^1}
			\le Cs^{-\frac12},
			\qquad
			\bigl\|\bigl(s\p_{\xi_i}+\p_{\eta_i}\bigr)\mathfrak G(s,\cdot)\bigr\|_{L^\infty}
			\le Cs^{-\frac{Q+1}2}.
		\end{gathered}
	\end{equation}
	Since $\xi=x-y-sw$ and $\eta=v-w$, the chain rule gives
	\begin{equation}\label{eq:stationary-chain-rule}
		\p_{w^i}p_s(Z;Y)
		=-\bigl(s\p_{\xi_i}+\p_{\eta_i}\bigr)\mathfrak G(s,\xi,\eta).
	\end{equation}
	Furthermore, for fixed $Z$ the map $Y=(y,w)\mapsto(\xi,\eta)$ is affine
	with Jacobian determinant of modulus one and range
	$\bR^3\times\bR^3$, and so is the map $Z=(x,v)\mapsto(\xi,\eta)$ for
	fixed $Y$.  Hence, for every measurable $F$ on $\bR^3\times\bR^3$,
	\[
	\int_{\bR^3_y\times \bR^3_w}|F(\xi,\eta)|\dd Y
	=\int_{\bR^3_x\times \bR^3_v}|F(\xi,\eta)|\dd Z
	=\|F\|_{L^1(\bR^3\times\bR^3)},
	\qquad
	\sup_{Z,Y}|F(\xi,\eta)|=\|F\|_{L^\infty(\bR^3\times\bR^3)} ,
	\]
	the first two integrals being taken at any fixed value of the other
	variable.  Applying this to $F=\mathfrak G(s,\cdot)$ and, by
	\eqref{eq:stationary-chain-rule}, to
	$F=(s\p_{\xi_i}+\p_{\eta_i})\mathfrak G(s,\cdot)$,
	\eqref{eq:stationary-scaled-norms} becomes
	\begin{equation}\label{eq:stationary-kernel-sections}
		\begin{gathered}
			\sup_Y\|p_s(\cdot;Y)\|_{L^1}
			=\sup_Z\|p_s(Z;\cdot)\|_{L^1}=1,
			\qquad
			\sup_{Z,Y}p_s(Z;Y)\le Cs^{-\frac Q2},
			\\
			\sup_Y\|\p_{w^i}p_s(\cdot;Y)\|_{L^1}
			=\sup_Z\|\p_{w^i}p_s(Z;\cdot)\|_{L^1}
			\le Cs^{-\frac12},
			\qquad
			\sup_{Z,Y}|\p_{w^i}p_s(Z;Y)|\le Cs^{-\frac{Q+1}2}.
		\end{gathered}
	\end{equation}
	
	\medskip
	\noindent\emph{Step 2. Weak-type bounds on the kernel sections.}
	Fix $Y$ and $R>0$.  By \eqref{eq:stationary-kernel-sections},
	\[
	\Bigl\|\int_0^Rp_s(\cdot;Y)\dd s\Bigr\|_{L^1(\bR^3_x\times \bR^3_v)}\le R,
	\qquad
	\Bigl\|\int_R^\infty p_s(\cdot;Y)\dd s\Bigr\|_{L^\infty(\bR^3_x\times \bR^3_v)}
	\le C\int_R^\infty s^{-\frac Q2}\dd s
	=C_0R^{-\frac{Q-2}2}.
	\]
	Given $\lambda>0$, choose $R:=(2C_0/\lambda)^{2/(Q-2)}$, so that $C_0R^{-\frac{Q-2}2}=\lambda/2$.  Chebyshev's inequality applied
	to the first gives
	\[
	\bigl|\{Z:\mathfrak G_\st(Z;Y)>\lambda\}\bigr|
	\le  \left|\left\{Z:\int_0^Rp_s(Z;Y)\dd s>\frac \lambda 2\right\}\right|
	\le\frac2\lambda R
	=C\lambda^{-\frac Q{Q-2}}.
	\]
	Since $|\mathfrak K_{\st,i}|\le\int_0^\infty|\p_{w^i}p_s|\dd s$, the same
	argument, choosing 
	$R\simeq\lambda^{-2/(Q-1)}$, gives
	\[
	\bigl|\{Z:|\mathfrak K_{\st,i}(Z;Y)|>\lambda\}\bigr|
	\le C\lambda^{-Q/(Q-1)}.
	\]
	The bounds
	\eqref{eq:stationary-kernel-sections} are symmetric in $Z$ and $Y$, so
	both computations hold in either variable.  Since $Q/(Q-2)=6/5$ and
	$Q/(Q-1)=12/11$,
	\begin{equation}\label{eq:stationary-kernel-weak}
		\begin{aligned}
			\sup_Y\|\mathfrak G_\st(\cdot;Y)\|_{L^{6/5,\infty}}
			+\sup_Z\|\mathfrak G_\st(Z;\cdot)\|_{L^{6/5,\infty}}
			&\le C,
			\\
			\sup_Y\|\mathfrak K_{\st,i}(\cdot;Y)\|_{L^{12/11,\infty}}
			+\sup_Z\|\mathfrak K_{\st,i}(Z;\cdot)\|_{L^{12/11,\infty}}
			&\le C.
		\end{aligned}
	\end{equation}
	
	\medskip
	\noindent\emph{Step 3. Endpoints and interpolation.}
	Now we use the equivalence of the $L^{q,\infty}$ norm, for $1<q<\infty$,
	\[
	\|g\|_{L^{q,\infty}}
	\simeq_q
	\sup_{0<|E|<\infty}|E|^{-1/q'}\int_E|g|
	\]
	(see, \textit{e.g.}, \cite[Exercise~1.1.12]{grafakos2014classical}).
	For a measurable kernel $K(Z;Y)$ on $\bR^6\times \bR^6$ and $H\in L^1(\bR^6)$, define
	$TH(Z):=\int K(Z;Y)H(Y)\dd Y$.
	Then for every measurable $E$ with $0<|E|<\infty$,
	\[
	\int_E|TH|
	\le
	\int|H(Y)|\Bigl(\int_E|K(Z;Y)|\dd Z\Bigr)\dd Y
	\le
	C_q\Bigl(\sup_Y\|K(\cdot;Y)\|_{L^{q,\infty}}\Bigr)
	|E|^{1-1/q}\|H\|_{L^1},
	\]
	so that $T:L^1\to L^{q,\infty}$.
	Taking $(T,K,q)=(\mathfrak G_\st[\cdot],\mathfrak G_\st(\cdot;\cdot),6/5)$ and
	$(\mathfrak G_{\st,i}[\cdot],\mathfrak K_{\st,i}(\cdot;\cdot),12/11)$ in
	\eqref{eq:stationary-kernel-weak} gives
	\eqref{eq:stationary-weak-L1-endpoint}.
	
	The duality of Lorentz spaces \cite[Theorem 1.4.16]{grafakos2014classical} gives
	\[
	|\mathfrak G_{\st,i}[H](Z)|
	\le C\|\mathfrak K_{\st,i}(Z;\cdot)\|_{L^{12/11,\infty}(\bR^6)}\|H\|_{L^{12,1}(\bR^6)},
	\]
	hence $\mathfrak G_{\st,i}:L^{12,1}\to L^\infty$.
	Real interpolation of
	Lorentz spaces \cite[Section~5.3]{BerghLofstrom1976} between
	$\mathfrak G_{\st,i}:L^1\to L^{12/11,\infty}$ and
	$\mathfrak G_{\st,i}:L^{12,1}\to L^\infty$, with index $(\theta,p):=(\tfrac1{11},\infty)$, gives
	\[
	\|\mathfrak G_{\st,i}[H]\|_{L^{6/5,\infty}(\bR^6)}\leq C\|H\|_{L^{12/11,\infty}(\bR^6)},
	\]
	which is \eqref{eq:stationary-weak-Lp-endpoint}.
\end{proof}

\begin{prop}[Integrability gain for stationary kinetic functions]
	\label{prop:stationary-int-gain-critical-sharp}
	There is a universal constant $C$ such that every
	$f\in\cH_\kin(\cD\times\bR_v^3)$ satisfies
	\begin{equation}
		\|f\|_{L^{12/5,\infty}(\cD\times\bR_v^3)}
		\le C\|f\|_{\cH_\kin(\cD\times\bR_v^3)}.
		\label{eq:stationary-critical-weak-sharp}
	\end{equation}
	Consequently, for every $2\le r<12/5$,
	\begin{equation}
		\|f\|_{L^r(\cD\times\bR_v^3)}
		\le C_r\|f\|_{\cH_\kin(\cD\times\bR_v^3)}.
		\label{eq:stationary-strong-subcritical-sharp}
	\end{equation}
	The constants are independent of $\cD$.
\end{prop}

\begin{proof}
	Write $M:=\|f\|_{\cH_\kin(\cD\times\bR_v^3)}$ and $P:=v\cdot\nabla_x-\Delta_v$.
	The proof of Proposition~\ref{lem:int_gain} applies on $\bR^6$ in place of
	$\bR^7$, with the following substitutions.
	
	The renormalized identity is \eqref{reno_std}; the measure $\mathfrak m_n$ in
	\eqref{eq:local-gain-zero-extension-inequality} is carried by the inflow set
	alone, and \eqref{eq:local-gain-L1-source-bound} holds with $M^2$ on the right.
	The mollification is in $(x,v)$, and the Duhamel formula is replaced by
	$F=\mathfrak G_\st[PF]$, valid for smooth nonnegative bounded integrable $F$ by
	\eqref{eq:stationary-kernel-sections}; the kernel of $\mathfrak G_\st$ is
	nonnegative.  The homogeneous dimension is $12$ in place of $14$:
	$X_n:=\|U_n\|_{L^{6/5,\infty}(\bR^6)}=\|U_n^{1/2}\|_{L^{12/5,\infty}(\bR^6)}^2$,
	the H\"older exponents are $\tfrac{11}{12}=\tfrac5{12}+\tfrac12$, and
	\eqref{eq:stationary-weak-L1-endpoint}--\eqref{eq:stationary-weak-Lp-endpoint}
	replace \eqref{eq:local-gain-potential-endpoints}--%
	\eqref{eq:local-gain-critical-lorentz}.
	
	This gives $X_n\le CM^2$ uniformly in $n$, and
	\eqref{eq:stationary-critical-weak-sharp} follows on letting $n\to\infty$ and
	taking square roots.  Interpolating it with $\|f\|_{L^2}\le M$ gives
	\eqref{eq:stationary-strong-subcritical-sharp}.
\end{proof}

\begin{remark}[Stationary versus evolutionary critical exponents]
	\label{rem:stationary-critical-scaling-sharp}
	The stationary kinetic homogeneous dimension is $12$, whereas the
	evolutionary one is $14$.  Accordingly the renormalized-square endpoints are
	$12/5$ and $7/3$, respectively.  The sharper stationary endpoint is recorded here, but the later nonlinear and
	diffuse arguments deliberately retain the more conservative evolutionary range
	$q<7/3$, and hence $p>14$, in order to work with a single set of exponents
	throughout the rest of the paper; when a fixed exponent is needed, they use the value $p=30$.
\end{remark}

\begin{prop}[Mixed-norm gain of integrability]
	\label{prop:IntGain}
	Let $0<T<\infty$ and $f\in \cH_\kin(\cD_T\times \bR^3_v)$, and set
	\begin{equation}
		\label{eq:IntGain-windows}
		I_n:=[n,n+1]\cap[0,T],
		\qquad 0\le n\le\lfloor T\rfloor ,
	\end{equation}
	so that the $I_n$ cover $[0,T]$ and meet only at endpoints.  For every
	$2\le r\le q<\tfrac73$ there is a constant $C_q$, depending only on $q$, such
	that
	\begin{equation}
		\label{IntGain_dyn_11}
		\left(\sum_{n=0}^{\lfloor T\rfloor}
		\|f\|_{L^{r,q,q}_{t,x,v}(I_n\times \cD\times\bR_v^3)}^2\right)^{1/2}
		\le
		C_q\|f\|_{\cH_\kin(\cD_T\times\bR_v^3)}.
	\end{equation}
	In particular, 
	\begin{equation}
		\label{eq:dynamic-global-mixed-integrability}
		\|f\|_{L_{t,x,v}^{r,q,q}(\cD_T\times\bR_v^3)}
		\le C_q\|f\|_{\cH_\kin(\cD_T\times\bR_v^3)}.
	\end{equation}
\end{prop}

\begin{proof}
	Choose \(\eta\in C_c^\infty((-2,2))\) with $\eta\equiv1$ on $[-1,1]$, and
	put $\eta_n(t):=\eta(t-n)$ for $n\in\bN_0$.
	
	Since $|I_n|\le1$ and $r\le q$, H\"older's inequality in $t$ gives
	\[
	\|f\|_{L^{r,q,q}_{t,x,v}(I_n\times\cD\times\bR_v^3)}
	\le|I_n|^{\frac1r-\frac1q}
	\|f\|_{L^{q}(I_n\times\cD\times\bR_v^3)}
	\le\|f\|_{L^{q}(I_n\times\cD\times\bR_v^3)} ,
	\]
	and $\eta_n\equiv1$ on $I_n\subset[n-1,n+1]$, so that
	\begin{equation}
		\label{eq:IntGain-local-Lq-reduction}
		\sum_{n=0}^{\lfloor T\rfloor}
		\|f\|_{L^{r,q,q}_{t,x,v}(I_n\times \cD\times\bR_v^3)}^2
		\le
		\sum_{n=0}^{\lfloor T\rfloor}
		\|f\eta_n\|_{L^q(\cD_T\times \bR_v^3)}^2 .
	\end{equation}
	
	Since $\eta$ is supported in $(-2,2)$,
	\begin{equation}
		\label{eq:IntGain-cutoff-bounds}
		|\eta_n|
		+|\p_t\eta_n|
		\le
		C\mathbf1_{(n-2,n+2)}.
	\end{equation}
	Since $\eta_n$ depends only on $t$, we have
	\begin{equation}
		\label{eq:IntGain-H1-bound}
		\sum_{n=0}^{\lfloor T\rfloor}
		\|f\eta_{n}\|_{L^2(\cD_T;H^1_v)}^2
		\le
		C\|f\|_{L^2(\cD_T;H^1_v)}^2.
	\end{equation}
	
	Denote $\mathfrak T:=\p_t+v\cdot \nabla_x$, and choose \(g_0,g_1,g_2,g_3\in L^2(\cD_T\times\bR_v^3)\) so that
	\[
	\mathfrak Tf=g_0+D_ig_i,
	\qquad
	\sum_{i=0}^3\|g_i\|_{L^2}^2
	\le C\|\mathfrak Tf\|_{L^2(\cD_T;H^{-1}_v)}^2.
	\]
	Then
	\[
	\mathfrak T(f\eta_n)
	=
	f\p_t\eta_n
	+g_0\eta_n
	-g_iD_i\eta_n
	+D_i(g_i\eta_n).
	\]
	The \(H^{-1}_v\) characterization and \eqref{eq:IntGain-cutoff-bounds} imply
	\begin{equation}
		\label{eq:IntGain-transport-bound}
		\sum_{n=0}^{\lfloor T\rfloor}
		\|\mathfrak T(f\eta_n)\|_
		{L^2(\cD_T;H^{-1}_v)}^2
		\le
		C\left(
		\|f\|_{L^2(\cD_T\times \bR_v^3)}^2
		+\|\mathfrak Tf\|_{L^2(\cD_T;H^{-1}_v)}^2
		\right).
	\end{equation}
	
	By \eqref{eq:IntGain-cutoff-bounds},
	\begin{equation}
		\label{eq:IntGain-initial-trace-bound}
		\begin{aligned}
			\sum_{n=0}^{\lfloor T\rfloor}
			\|(f\eta_n)(0)\|_{L^2(\cD\times \bR_v^3)}^2
			\,&\le
			C\|f(0)\|_{L^2(\cD\times \bR_v^3)}^2,
			\\
			\sum_{n=0}^{\lfloor T\rfloor}
			\|(f\eta_n)_{\gamma^-}\|_
			{L^2(\gamma_T^-)}^2
			\,&\le
			C\|f_{\gamma^-}\|_{L^2(\gamma_T^-)}^2.
		\end{aligned}
	\end{equation}
	Combining \eqref{eq:IntGain-H1-bound},
	\eqref{eq:IntGain-transport-bound} and
	\eqref{eq:IntGain-initial-trace-bound}, and applying
	\eqref{2603291119} to each $f\eta_n$, we obtain
	\begin{align}\label{260911726}
		\sum_{n=0}^{\lfloor T\rfloor}
		\|f\eta_n\|_{L^q(\cD_T\times \bR_v^3)}^2
		\le C\|f\|_{\cH_\kin(\cD_T\times \bR^3_v)}^2
	\end{align}
	where $C>0$ is a constant depending only on $\eta$ and $q$.
	Combining \eqref{eq:IntGain-local-Lq-reduction} and \eqref{260911726}, we have \eqref{IntGain_dyn_11}.
	
	\eqref{eq:dynamic-global-mixed-integrability} follows from \eqref{IntGain_dyn_11} and the embedding $\ell^2(\bN_0)\hookrightarrow\ell^r(\bN_0)$, $r\ge2$.
\end{proof}

\begin{corollary}\label{2609201144}
		For any $g\in \cH_{\kin,\sigma}(\cD_T\times \bR^3_v)$ and $2\leq r\leq q<\frac 73$, 
		\begin{align}\label{IntGain_dyn}
			\|\mu^{1/8}g\|_{L^{r,q,2}_{t,x,v}(\cD_T\times \bR^3_v)}\leq C_q\|g\|_{\cH_{\kin,\sigma}(\cD_T\times \bR^3_v)}.
		\end{align}
		Similarly, for any $g\in \cH_{\kin,\sigma}(\cD\times \bR^3_v)$, $2\leq q<\frac {12}5$, 
		\begin{align}\label{IntGain_std}
			\|\mu^{1/8}g\|_{L^{q,2}_{x,v}(\cD\times \bR^3_v)}\leq C_q\|g\|_{\cH_{\kin,\sigma}(\cD\times \bR^3_v)}.
		\end{align}
	\end{corollary}
	
\begin{proof}
	Let $I_n$, $0\le n\le\lfloor T\rfloor$, be the windows
	\eqref{eq:IntGain-windows}. Since $r\ge2$,
	\begin{align*}
	\|\mu^{1/16}g\|_{L^{r,q,q}_{t,x,v}(\cD_T\times\bR_v^3)}
	&=\Big(\sum_{n=0}^{\lfloor T\rfloor}
	\|\mu^{1/16}g\|_{L^{r,q,q}_{t,x,v}(I_n\times\cD\times\bR_v^3)}^r\Big)^{1/r}\\
	&\le\Big(\sum_{n=0}^{\lfloor T\rfloor}
	\|\mu^{1/16}g\|_{L^{r,q,q}_{t,x,v}(I_n\times\cD\times\bR_v^3)}^2\Big)^{1/2}.
	\end{align*}
	Applying \eqref{IntGain_dyn_11} to $\mu^{1/16}g$ yields
	\begin{align}\label{2609201138}
	\|\mu^{1/16}g\|_{L^{r,q,q}_{t,x,v}(\cD_T\times\bR_v^3)}
	\le C_q\|\mu^{1/16}g\|_{\cH_{\kin}(\cD_T\times\bR_v^3)}.
	\end{align}
	H\"older in velocity, with $1/2=1/s+1/q$, gives
	\begin{equation}\label{eq:IntGain-velocity-exponent}
	\|\mu^{1/8}g\|_{L_v^2}
	\le\|\mu^{1/16}\|_{L_v^s}\|\mu^{1/16}g\|_{L_v^q},
	\end{equation}
	where the Maxwellian norm is finite for $2\le q<\infty$.
	Consequently \eqref{2609201138} controls the claimed
	$L^{r,q,2}_{t,x,v}$ norm. Moreover, Maxwellian multiplication gives
	\begin{equation}\label{eq:N-Maxwellian-multiplier}
	\|\mu^{1/16}g\|_{\cH_\kin(\cD_T\times\bR_v^3)}
	\le C\|g\|_{\cH_{\kin,\sigma}(\cD_T\times\bR_v^3)}.
	\end{equation}
	Indeed, $\|\mu^{1/16}u\|_{H_v^1}\le C\|u\|_{H_\sigma^1}$ and
	$\|\mu^{1/16}U\|_{H_v^{-1}}\le C\|U\|_{H_\sigma^{-1}}$;
	the initial and inflow traces are controlled directly. This proves
	\eqref{IntGain_dyn}.
	For \eqref{IntGain_std}, apply \eqref{eq:IntGain-velocity-exponent}
	at each $x$ and \eqref{eq:stationary-strong-subcritical-sharp}
	to $\mu^{1/16}g$ to obtain
	\[
	\|\mu^{1/8}g\|_{L^{q,2}_{x,v}(\cD\times\bR_v^3)}
	\le C_q\|\mu^{1/16}g\|_{\cH_\kin(\cD\times\bR_v^3)}.
	\]
	The stationary form of \eqref{eq:N-Maxwellian-multiplier} concludes.
\end{proof}

\begin{proof}[\textbf{Proof of Theorem~\ref{260529130}}]
	Only the evolutionary case \eqref{eq:N-estimate-evolutionary-sharp} is
	proved; the stationary one is identical.
	
	Set \(q:=2p/(p-2)\).  The condition \(p>14\) gives precisely
	\[
	2<q<\frac73,
	\qquad
	\frac1p+\frac1q=\frac12.
	\]
	Corollary~\ref{2609201144}, applied to $\mu^{1/8}g$, gives
	\begin{equation}
		\label{eq:N-estimate-g-bound}
		\|\mu^{1/4}g\|_{L^{2,q,2}_{t,x,v}(\cD_T\times \bR_v^3)}\le 
		C_p\|\mu^{1/8}g\|_{\cH_{\kin,\sigma}(\cD_T\times \bR_v^3)}.
	\end{equation}
	
	Apply \eqref{estq:N_dy} with \(\varpi=1/4\), followed by
	\eqref{eq:N-estimate-g-bound}, to obtain
	\eqref{eq:N-estimate-evolutionary-sharp}.
	
	For the last two estimates, \eqref{eq:H_sigma_norm} and
	Lemma~\ref{lem:landau-diffusion-eigenvalues} give
	\[
	\|w^\ell f(t,x,\cdot)\|_{H^1_\sigma}
	\le
	C_\ell\left(
	\|w^{\ell-1/2}f(t,x,\cdot)\|_{L^2_v}
	+\|w^{\ell-1/2}D_vf(t,x,\cdot)\|_{L^2_v}
	\right),
	\]
	while
	\begin{align*}
		\|w^{\ell-1/2}f\|_{L^{\infty,p,2}_{t,x,v}(\cD_T\times \bR_v^3)}
		&\le
		|\cD|^{1/p}\|w^{-5/2}\|_{L^2(\bR_v^3)}
		\|w^{\ell+2}f\|_{L^\infty(\cD_T\times \bR_v^3)},
		\\
		\|w^{\ell-1/2}D_vf\|_{L^{\infty,p,2}_{t,x,v}(\cD_T\times \bR_v^3)}
		&\le
		\|w^{\ell+2}D_vf\|_{L^{\infty,p,2}_{t,x,v}(\cD_T\times \bR_v^3)}.
	\end{align*}
	Inserting these into
	\eqref{eq:N-estimate-evolutionary-sharp} and
	\eqref{eq:N-estimate-stationary-sharp} gives
	\eqref{eq:N-estimate-evolutionary} and
	\eqref{eq:N-estimate-stationary}.
\end{proof}

\section{Linearized Landau Inflow Problems: Coercivity and Solvability}\label{sec:landau-regularity}

We work in the framework of \S\ref{sec:landau-framework}.  For sufficiently
small $g$, Theorem~\ref{260322442} establishes unique
solvability in the $\cH_{\kin,\sigma}$ spaces of
Definition~\ref{def:H_kin_sigma} for the time-dependent equation
\begin{equation}\label{solv:linland_dy}
	\begin{cases}
		(\partial_t+v\cdot\nabla_x)f+L_gf=S
		&\text{in }\cD_T\times\bR_v^3,\\
		f(0,\cdot)=f_0
		&\text{in }\cD\times\bR_v^3,\\
		f=h
		&\text{on }\gamma_T^-,
	\end{cases}
\end{equation}
and its stationary counterpart
\begin{equation}\label{solv:linland_st}
	\begin{cases}
		v\cdot\nabla_xf+L_gf=S
		&\text{in }\cD\times\bR_v^3,\\
		f=h
		&\text{on }\gamma^-,
	\end{cases}
\end{equation}
where $L_g$ is defined in \eqref{def:L_g}.

The theorem, together with the velocity-weighted estimates that accompany
it, is proved in Subsection~\ref{solv_linlandau_dy}.

\begin{defn}[$\cH_{\kin,\sigma}$-solutions]\label{def:Hkinsigma-solution}
	Let $\cL$ be a bounded operator from $L^2(\cD_T;H_\sigma^1)$ to
	$L^2(\cD_T;H_\sigma^{-1})$. Given
	\[
	S\in L^2(\cD_T;H_\sigma^{-1}),
	\qquad
	f_0\in L^2(\cD\times\bR_v^3),
	\qquad
	h\in L^2(\gamma_T^-),
	\]
	we say that $f\in\cH_{\kin,\sigma}(\cD_T\times\bR_v^3)$ is an
	$\cH_{\kin,\sigma}$-solution of the inflow problem
	\begin{equation}\label{def:Hkinsigma-inflow-dy}
		\begin{cases}
			(\partial_t+v\cdot\nabla_x)f+\cL f=S
			&\text{in }\cD_T\times\bR_v^3,\\
			f(0,\cdot)=f_0
			&\text{in }\cD\times\bR_v^3,\\
			f=h
			&\text{on }\gamma_T^-
		\end{cases}
	\end{equation}
	if
	\begin{align*}
		\begin{gathered}
			(\partial_t+v\cdot\nabla_x)f+\cL f=S
			\quad\text{in }L^2(\cD_T;H_\sigma^{-1}),\\
			f(0,\cdot)=f_0
			\quad\text{in }L^2(\cD\times\bR_v^3),
			\qquad
			f_{\gamma^-}=h
			\quad\text{in }L^2(\gamma_T^-).		
		\end{gathered}
	\end{align*}
	
	In the stationary case, let $\cL$ be a bounded operator from
	$L^2(\cD;H_\sigma^1)$ to $L^2(\cD;H_\sigma^{-1})$. Given
	\[
	S\in L^2(\cD;H_\sigma^{-1}),
	\qquad
	h\in L^2(\gamma^-),
	\]
	we say that $f\in\cH_{\kin,\sigma}(\cD\times\bR_v^3)$ is an
	$\cH_{\kin,\sigma}$-solution of the inflow problem
	\begin{equation}\label{def:Hkinsigma-inflow-st}
		\begin{cases}
			v\cdot\nabla_xf+\cL f=S
			&\text{in }\cD\times\bR_v^3,\\
			f=h
			&\text{on }\gamma^-
		\end{cases}
	\end{equation}
	if
	\begin{equation}\label{def:Hkinsigma-solution-st}
		v\cdot\nabla_xf+\cL f=S
		\quad\text{in }L^2(\cD;H_\sigma^{-1}),
		\qquad
		f_{\gamma^-}=h
		\quad\text{in }L^2(\gamma^-).
	\end{equation}
\end{defn}

\begin{thm}[Perturbed linearized Landau inflow problems]
	\label{260322442}
	Here $L_g$ is as in \eqref{def:L_g}, and all inflow problems below are understood in the sense of
	Definition~\ref{def:Hkinsigma-solution}.
	There exist constants
	\begin{align*}
		\varepsilon_0=\varepsilon_0(\cD)>0,
		\qquad
		C=C(\cD)>0
	\end{align*}
	such that the following assertions hold.
	
	\begin{enumerate}
		\item
		Let $0<T<\infty$, and assume that
		\begin{align*}
			\|\mu^{1/4}g\|_
			{L^\infty(\cD_T\times\bR_v^3)}
			\le
			\varepsilon_0.
		\end{align*}
		Then, for every
		\[
		S\in L^2(\cD_T;H_\sigma^{-1}),
		\qquad
		f_0\in L^2(\cD\times\bR_v^3),
		\qquad
		h\in L^2(\gamma_T^-),
		\]
		the problem
		\begin{equation}\label{260322508}
			\begin{cases}
				(\partial_t+v\cdot\nabla_x)f+L_gf=S
				&\text{in }\cD_T\times\bR_v^3,\\
				f(0,\cdot)=f_0
				&\text{in }\cD\times\bR_v^3,\\
				f=h
				&\text{on }\gamma_T^-
			\end{cases}
		\end{equation}
		admits a unique solution
		$
		f\in\cH_{\kin,\sigma}(\cD_T\times\bR_v^3)
		$.
		Moreover,
		\begin{equation}\label{eq:perturbed-landau-dynamic-estimate}
			\|f\|_{\cH_{\kin,\sigma}(\cD_T\times\bR_v^3)}
			\le
			C\left(
			\|S\|_{L^2(\cD_T;H_\sigma^{-1})}
			+
			\|f_0\|_{L^2(\cD\times\bR_v^3)}
			+
			\|h\|_{L^2(\gamma_T^-)}
			\right).
		\end{equation}
		
		\item
		Assume that
		\[
		\|\mu^{1/4}g\|_
		{L^\infty(\cD\times\bR_v^3)}
		\le
		\varepsilon_0.
		\]
		Then, for every
		\[
		S\in L^2(\cD;H_\sigma^{-1}),
		\qquad
		h\in L^2(\gamma^-),
		\]
		the problem
		\begin{equation}\label{eq:perturbed-landau-stationary-inflow}
			\begin{cases}
				v\cdot\nabla_xf+L_gf=S
				&\text{in }\cD\times\bR_v^3,\\
				f=h
				&\text{on }\gamma^-
			\end{cases}
		\end{equation}
		admits a unique solution
		$
		f\in\cH_{\kin,\sigma}(\cD\times\bR_v^3)$.
		Moreover,
		\begin{equation}\label{eq:perturbed-landau-stationary-estimate}
			\|f\|_{\cH_{\kin,\sigma}(\cD\times\bR_v^3)}
			\le
			C\left(
			\|S\|_{L^2(\cD;H_\sigma^{-1})}
			+
			\|h\|_{L^2(\gamma^-)}
			\right).
		\end{equation}
	\end{enumerate}
	
	Here, the constant $C$ depends only on $\cD$. In particular, it is independent of $T$, $g$, and the data.
	
\end{thm}

\subsection{The macroscopic estimate}\label{subsec:macro-landau}

Recall the macro--micro decomposition
\[
f=\mathbf P f+(\mathbf I-\mathbf P)f,
\]
where $\mathbf P f$ and its coefficients $\mathbf a$, $\mathbf b$, $\mathbf c$
are as in the earlier \hyperlink{comp:macro-projection}{projection formula} and
\hyperlink{comp:macro-moments}{Gaussian moment formulas} of
the operator section. Here the coefficient functions depend on
$(t,x)$ (or only $x$ in the stationary case), and the earlier formulas
are applied pointwise in these variables. The macroscopic quantities $\mathbf a$, $\mathbf b$ and $\mathbf c$ are not controlled by the dissipation of $L$, which by \eqref{eq:L-coercivity} is coercive only on $(\mathbf I-\mathbf P)f$.
We recover them from the transport operator by testing the equation against velocity moments of solutions of auxiliary Poisson problems in $x$, following the test-function method of \cite{Guo2010,EspositoGuoKimMarra}.

Let
\begin{equation}\label{260715226}
	\cZ
	:=
	\left\{
	v^\alpha\sqrt\mu
	\ :\
	\alpha\in\bN_0^3,\ |\alpha|\le4
	\right\},
\end{equation}
where $v^\alpha:=v_1^{\alpha_1}v_2^{\alpha_2}v_3^{\alpha_3}$. Thus
$p\sqrt\mu\in\operatorname{span}\cZ$ for every polynomial $p$ of degree
at most $4$, hence for every velocity profile of the test functions used
below, each of which is a polynomial of degree at most $4$ times $\sqrt\mu$.
For an $\cS'(\bR_v^3)$-valued function $S$ on $\cD$,
define
\begin{equation*}
	\|S\|_{\cM}^2
	:=
	\sum_{\zeta\in\cZ}
	\int_{\cD}
	\left|
	\langle S(x,\cdot),\zeta\rangle_v
	\right|^2
	\dd x,
\end{equation*}
whenever the right-hand side is finite. For an
$\cS'(\bR_v^3)$-valued function $S$ on
$[0,T]\times\cD$, set
\begin{equation*}
	\|S\|_{\cM(s,t)}^2
	:=
	\int_s^t
	\|S(r)\|_{\cM}^2
	\dd r,
	\qquad
	\|S\|_{\cM_T}
	:=
	\|S\|_{\cM(0,T)}.
\end{equation*}
Here $\langle\cdot,\cdot\rangle_v$ denotes the
$\cS'(\bR_v^3)$--$\cS(\bR_v^3)$ duality
pairing.

Since $\cZ$ is a finite subset of
$\cS(\bR_v^3)\subset H_\sigma^1$, duality gives
\begin{equation*}
	\sum_{\zeta\in\cZ}
	\left|
	\langle S,\zeta\rangle_v
	\right|^2
	\le
	C\|S\|_{H_\sigma^{-1}}^2
	\qquad
	\text{for every }S\in H_\sigma^{-1}.
\end{equation*}
Similarly,
\begin{equation*}
	\sum_{\zeta\in\cZ}
	\left|
	(f,\zeta)_{L^2_v}
	\right|^2
	\le
	C\|f\|_{L^2(\bR_v^3)}^2
	\qquad
	\text{for every }f\in L^2(\bR_v^3).
\end{equation*}
Consequently,
\begin{equation*}
	\begin{aligned}
		\|S\|_{\cM}
		&\le
		C\|S\|_{L^2(\cD;H_\sigma^{-1})},
		\\
		\|f\|_{\cM}
		&\le
		C\|f\|_{L^2(\cD\times\bR_v^3)},
	\end{aligned}
\end{equation*}
and, for every $0\le s<t\le T$,
\begin{equation}\label{eq:appendix-source-M-time-Hminus1}
	\|S\|_{\cM(s,t)}
	\le
	C\|S\|_{L^2((s,t)\times\cD;H_\sigma^{-1})}.
\end{equation}

\begin{lemma}\label{lem:appendix-boundary-estimate}
	Let $\phi\in H^2(\cD)\cap H_0^1(\cD)$, and define
	\[
	\Phi(x,v)
	:=
	\sum_{k=1}^3
	p_k(v)D_{x_k}\phi(x)\sqrt{\mu(v)},
	\]
	where $p_1,p_2,p_3$ are fixed polynomials. Then
	\begin{equation}\label{2607161211}
		\int_\gamma|\Phi|^2\dd\gamma
		\le
		C\|\partial_n\phi\|_{L^2(\partial\cD)}^2
		\le
		C\|\phi\|_{H^2(\cD)}^2.
	\end{equation}
	In particular, if
	\begin{equation}\label{260716118}
		\begin{cases}
			-\Delta\phi=r
			&\text{in }\cD,\\
			\phi=0
			&\text{on }\partial\cD,
		\end{cases}
	\end{equation}
	for some $r\in L^2(\cD)$, then
	\begin{equation}\label{26071612112}
		\int_\gamma|\Phi|^2\dd\gamma
		\le
		C\|r\|_{L^2(\cD)}^2.
	\end{equation}
	Here, $C>0$ depends only on $p_1,p_2,p_3$ and $\cD$.
\end{lemma}

\begin{proof}
	By the trace theorem and the tangential decomposition of the boundary
	gradient on $C^{1,1}$ domains
	(see \cite[Theorem~1.5.1.2 and Section~3.1.1]{Grisvard1985}),
	the trace of $\nabla_x\phi$ is well defined on $\partial\cD$ and
	\begin{equation}\label{eq:normal-trace-estimate}
		\|\partial_n\phi\|_{L^2(\partial\cD)}
		\le
		C\|\phi\|_{H^2(\cD)}.
	\end{equation}
	Since $\phi\in H_0^1(\cD)$, its trace vanishes on $\partial\cD$.
	Hence
	\[
	(I_3-n\otimes n)\nabla_x\phi
	=
	\nabla_\tau(\phi|_{\partial\cD})
	=
	0
	\qquad
	\text{a.e. on }\partial\cD,
	\]
	and therefore
	\[
	\nabla_x\phi
	=
	\partial_n\phi\,n
	\qquad
	\text{a.e. on }\partial\cD.
	\]
	It follows that
	\[
	\Phi(x,v)
	=
	\partial_n\phi(x)
	\sum_{k=1}^3
	p_k(v)n_k(x)\sqrt{\mu(v)}
	\qquad
	\text{for a.e. }(x,v)\in\gamma.
	\]
	Using
	\[
	\dd\gamma
	=
	|n(x)\cdot v|\,\dd v\,\dd S_x,
	\]
	together with $|n(x)\cdot v|\le |v|$ and
	\[
	\left|
	\sum_{k=1}^3 p_k(v)n_k(x)
	\right|^2
	\le
	\sum_{k=1}^3 |p_k(v)|^2,
	\]
	we obtain
	\begin{align*}
		\int_\gamma |\Phi|^2\dd\gamma
		&\le
		\|\partial_n\phi\|_{L^2(\partial\cD)}^2
		\int_{\bR^3}
		|v|
		\sum_{k=1}^3 |p_k(v)|^2
		\mu(v)\,\dd v\le
		C\|\partial_n\phi\|_{L^2(\partial\cD)}^2
		\le
		C\|\phi\|_{H^2(\cD)}^2,
	\end{align*}
	where the last inequality follows from
	\eqref{eq:normal-trace-estimate}. This proves
	\eqref{2607161211}.
	
	If $\phi$ satisfies \eqref{260716118}, the standard $H^2$-regularity
	estimate for the Dirichlet problem on $C^{1,1}$ domains
	(see, e.g., \cite[Chapter~2]{Grisvard1985}) yields
	\[
	\|\phi\|_{H^2(\cD)}
	\le
	C\|r\|_{L^2(\cD)}.
	\]
	Combining this with \eqref{2607161211} gives
	\eqref{26071612112}.
\end{proof}

\begin{lemma}\label{lem:admissible-test-functions}
	Let $p_1,\dots,p_M$ be fixed polynomials on $\bR_v^3$ and let
	\[
	\Phi(t,x,v):=\sum_{k=1}^M\psi_k(t,x)\,p_k(v)\sqrt{\mu(v)} .
	\]
	If $\psi_1,\dots,\psi_M\in L^2\big(0,T;H^1(\cD)\big)$ with
	$\p_t\psi_k\in L^2(\cD_T)$, then
	\[
	\Phi\in\cH_{\kin,\sigma}(\cD_T\times\bR_v^3),
	\qquad
	T_\gamma\Phi=\Phi|_\gamma,
	\]
	and
	\[
	\|\Phi\|_{\cH_{\kin,\sigma}(\cD_T\times\bR_v^3)}
	\le
	C\sum_{k=1}^M
	\Big(
	\|\psi_k\|_{L^2(0,T;H^1(\cD))}+\|\p_t\psi_k\|_{L^2(\cD_T)}
	+\|\psi_k(0,\cdot)\|_{L^2(\cD)}
	\Big),
	\]
	with $C$ depending only on $p_1,\dots,p_M$ and $\cD$. The analogous
	statement holds in the stationary space, for
	$\psi_k\in H^1(\cD)$.
\end{lemma}

\begin{proof}
	It suffices to treat a single summand $\Phi=\psi\,p\sqrt\mu$. Since
	$p\sqrt\mu$ and $v_ip\sqrt\mu$ are Schwartz functions, they belong to
	$H^1_v(\bR_v^3)\subset H^1_\sigma(\bR_v^3)$ by
	Lemma~\ref{lem:H_sigma-Hilbert}.(2). Consequently
	\[
	\|\Phi(t,x,\cdot)\|_{H^1_\sigma}\le C_p|\psi(t,x)|,
	\qquad
	(\p_t+v\cdot\nabla_x)\Phi
	=
	\big(\p_t\psi\big)p\sqrt\mu+\big(D_{x_i}\psi\big)v_ip\sqrt\mu,
	\]
	so that, by \eqref{eq:Hsigma-dual-weighted} applied to the
	Schwartz functions $p\sqrt\mu$ and $v_ip\sqrt\mu$,
	\[
	\|(\p_t+v\cdot\nabla_x)\Phi(t,x,\cdot)\|_{H^{-1}_\sigma}
	\le C_p\big(|\p_t\psi(t,x)|+|\nabla_x\psi(t,x)|\big).
	\]
	Hence
	$\Phi\in L^2(\cD_T;H^1_\sigma)$ and
	$(\p_t+v\cdot\nabla_x)\Phi\in L^2(\cD_T;H^{-1}_\sigma)$, with the stated
	bounds. The initial value $\Phi(0,\cdot)=\psi(0,\cdot)p\sqrt\mu$ lies in
	$L^2(\cD\times\bR_v^3)$. Finally, since
	$\int_{\bR^3}|p\sqrt\mu|^2|n_x\cdot v|\dd v\le C_p$, the ordinary trace
	theorem on $\cD$ gives
	\[
	\|\Phi\|_{L^2(\gamma_T)}^2
	\le
	C_p\|\psi\|_{L^2((0,T)\times\p\cD)}^2
	\le
	C_p\|\psi\|_{L^2(0,T;H^1(\cD))}^2 .
	\]
	
	To conclude that $\Phi$ belongs to $\cH_{\kin,\sigma}$, which is defined
	as a completion, we exhibit an approximating sequence. Choose
	$\psi_n\in C^\infty([0,T]\times\overline\cD)$ with $\psi_n\to\psi$ in
	$L^2(0,T;H^1(\cD))$ and $\p_t\psi_n\to\p_t\psi$ in $L^2(\cD_T)$, and let
	$\chi_R$ be the radial cut-offs of \eqref{eq:radial-cutoff}. Then
	$$
	\psi_n\chi_Rp\sqrt\mu\in C_c^\infty([0,T]\times\overline\cD\times\bR_v^3),
	$$
	and \eqref{eq:radial-cutoff-convergence} gives
	$$
	\chi_Rp\sqrt\mu\to p\sqrt\mu\quad \text{and}\quad \chi_Rv_ip\sqrt\mu\to v_ip\sqrt\mu\quad \text{in }H^1_\sigma,
	$$
	as $R\to\infty$.
	The corresponding convergence in $H^{-1}_\sigma$ follows from
	\eqref{eq:Hsigma-dual-weighted}, since
	$\|\la v\ra^{1/2}(1-\chi_R)q\sqrt\mu\|_{L^2}\to0$ for every polynomial $q$.
	Combining the two approximations, a diagonal argument produces a sequence
	in $C_c^\infty([0,T]\times\overline\cD\times\bR_v^3)$ converging to $\Phi$
	in the norm \eqref{innerprod_sigma}. Hence
	$\Phi\in\cH_{\kin,\sigma}(\cD_T\times\bR_v^3)$. Since the same sequence
	converges to $\Phi|_\gamma$ in $L^2(\gamma_T)$, the trace operator
	$T_\gamma$ of Proposition~\ref{2607151218} returns $\Phi|_\gamma$. The
	stationary case is identical, with the temporal terms omitted.
\end{proof}

\begin{thm}[Macroscopic estimate]
	\label{thm:macro-estimate-landau}
	There exists a constant $C=C(\cD)>0$ such that the following
	statements hold.
	
	\begin{enumerate}
		\item
		Let
		$
		f\in \cH_{\kin,\sigma}(\cD_T\times\bR_v^3),
		$ and set
		\begin{equation}\label{eq:appendix-linearized-landau-inflow}
			S
			:=
			(\partial_t+v\cdot\nabla_x)f+Lf ,
		\end{equation}
		which belongs to $L^2(\cD_T;H_\sigma^{-1})$ automatically, by
		Definition~\ref{def:H_kin_sigma} and \eqref{boundedness of L}. Here and
		below, $f_\gamma:=T_\gamma f$ denotes the phase-space trace supplied by
		Proposition~\ref{2607151218}.
		We identify $f$ with its representative in $C\bigl([0,T];L^2(\cD\times\bR_v^3)\bigr)$.
		Then there exists a function $\eta_f\in C([0,T])$ such that
		\begin{equation}\label{eq:macro-eta-bound}
			|\eta_f(t)|
			\le
			C\|f(t)\|_{\cM}^2
			\le
			C\|f(t)\|_{L^2(\cD\times\bR_v^3)}^2
			\qquad
			\text{for every }t\in[0,T],
		\end{equation}
		and
		\begin{equation}\label{eq:macro-eta-difference-bound}
			|\eta_f(t)-\eta_f(s)|
			\le
			C
			\left(\|f(t)\|_{\cM}+\|f(s)\|_{\cM}\right)\|f(t)-f(s)\|_{\cM}
			\qquad
			\text{for every }0\le s<t\le T.
		\end{equation}
		Moreover, for every $0\le s<t\le T$,
		\begin{align}
			\int_s^t
			\|\mathbf P f(r)\|_{L^2(\cD\times\bR_v^3)}^2
			\dd r
			\le{}&
			C\Big(
			\|\mu^{1/4}(\mathbf I-\mathbf P)f\|_
			{L^2((s,t)\times\cD\times\bR_v^3)}^2
			+
			\|S\|_{\cM(s,t)}^2
			\notag\\
			&\qquad
			+
			\|f_\gamma\|_{L^2((s,t)\times\gamma)}^2
			\Big)
			+
			\eta_f(t)-\eta_f(s).
			\label{eq:appendix-macro-estimate-endpoint}
		\end{align}
		
		\item
		Let
		$f\in\cH_{\kin,\sigma}(\cD\times\bR_v^3)$, and suppose that
		\begin{equation}\label{eq:appendix-linearized-landau-stationary}
			S
			:=
			v\cdot\nabla_x f+Lf
			\in
			L^2(\cD;H_\sigma^{-1}).
		\end{equation}
		Then
		\begin{equation}\label{eq:appendix-macro-estimate-endpoint_st}
			\|\mathbf P f\|_{L^2(\cD\times\bR_v^3)}^2
			\le
			C\left(
			\|\mu^{1/4}(\mathbf I-\mathbf P)f\|_
			{L^2(\cD\times\bR_v^3)}^2
			+
			\|S\|_{\cM}^2
			+
			\|f_\gamma\|_{L^2(\gamma)}^2
			\right).
		\end{equation}
	\end{enumerate}
\end{thm}

\begin{proof}
	We first prove (1). 
	Since
	$C_c^\infty([0,T]\times\overline{\cD}\times\bR^3_v)$
	is dense in
	$\cH_{\kin,\sigma}(\cD_T\times\bR^3_v)$
	and all the terms in \eqref{eq:appendix-macro-estimate-endpoint} are stable under
	this approximation, it suffices to prove the result for
	\[
	f\in
	C_c^\infty([0,T]\times\overline{\cD}\times\bR^3_v).
	\]
	
	\medskip
	\noindent\emph{Macroscopic moment equations.}
	Set
	\[
	Q:=(\mathbf I-\mathbf P)f,
	\qquad
	\mathbf P f
	=
	\left[
	\mathbf a
	+
	\mathbf b\cdot v
	+
	\mathbf c
	\left(|v|^2-\frac32\right)
	\right]\sqrt\mu,
	\]
	where $\mathbf a$, $\mathbf b$, and $\mathbf c$ are defined in
	\hyperlink{comp:macro-moments}{Gaussian moment formulas}. Taking the velocity moments of
	\eqref{eq:appendix-linearized-landau-inflow} against
	\[
	\pi^{-3/2}\sqrt\mu,
	\qquad
	2\pi^{-3/2}v_i\sqrt\mu,
	\qquad
	\frac23\pi^{-3/2}
	\left(|v|^2-\frac32\right)\sqrt\mu,
	\]
	and using \eqref{eq:L-kernel}, we obtain
	\begin{align}
		\partial_t\mathbf a
		+
		\frac12\nabla_x\cdot\mathbf b
		&=
		S_{\mathbf a},
		\label{eq:appendix-moment-a}
		\\
		\partial_t\mathbf b
		+
		\nabla_x(\mathbf a+\mathbf c)
		+
		\nabla_x\cdot\Theta
		&=
		S_{\mathbf b},
		\label{eq:appendix-moment-b}
		\\
		\partial_t\mathbf c
		+
		\frac13\nabla_x\cdot\mathbf b
		+
		\nabla_x\cdot\Lambda
		&=
		S_{\mathbf c}
		\label{eq:appendix-moment-c}
	\end{align}
	in the sense of distributions on $\cD_T$, where
	\begin{align*}
		S_{\mathbf a}
		&:=
		\pi^{-3/2}
		\langle S,\sqrt\mu\rangle_v,
		\\
		S_{\mathbf b,i}
		&:=
		2\pi^{-3/2}
		\langle S,v_i\sqrt\mu\rangle_v,
		\\
		S_{\mathbf c}
		&:=
		\frac23\pi^{-3/2}
		\left\langle
		S,
		\left(|v|^2-\frac32\right)\sqrt\mu
		\right\rangle_v,
	\end{align*}
	and
	\begin{align*}
		\Theta_{ij}
		&:=
		2\pi^{-3/2}
		\left(
		Q,
		\left(v_iv_j-\frac12\delta_{ij}\right)\sqrt\mu
		\right)_{L^2_v},
		\\
		\Lambda_j
		&:=
		\frac23\pi^{-3/2}
		\left(
		Q,
		\left(|v|^2-\frac32\right)v_j\sqrt\mu
		\right)_{L^2_v}.
	\end{align*}
	Here repeated indices are summed, and
	\[
	(\nabla_x\cdot\Theta)_i
	:=
	\partial_{x_j}\Theta_{ij},
	\qquad
	\nabla_x\cdot\Lambda
	:=
	\partial_{x_j}\Lambda_j.
	\]
	
	By the definition of $\cM$, for a.e. $r\in(0,T)$,
	\begin{equation}\label{eq:appendix-source-moment-control}
		\|S_{\mathbf a}(r)\|_{L^2(\cD)}^2
		+
		\|S_{\mathbf b}(r)\|_{L^2(\cD)}^2
		+
		\|S_{\mathbf c}(r)\|_{L^2(\cD)}^2
		\le
		C\|S(r)\|_{\cM}^2.
	\end{equation}
	Moreover, \eqref{eq:appendix-source-M-time-Hminus1} gives
	\[
	\|S\|_{\cM(s,t)}
	\le
	C
	\|S\|_{L^2((s,t)\times\cD;H_\sigma^{-1})}.
	\]
	
	For every polynomial $p$,
	\[
	\left|
	(Q,p\sqrt\mu)_{L^2_v}
	\right|
	\le
	\|\mu^{1/4}Q\|_{L^2_v}
	\|p\mu^{1/4}\|_{L^2_v}.
	\]
	Hence, for a.e. $r\in(0,T)$,
	\begin{equation}\label{eq:appendix-higher-moment-bound}
		\|\Theta(r)\|_{L^2(\cD)}^2
		+
		\|\Lambda(r)\|_{L^2(\cD)}^2
		\le
		C
		\|\mu^{1/4}Q(r)\|_{L^2(\cD\times\bR_v^3)}^2.
	\end{equation}
	
	Since $L\mathbf P f=0$ by \eqref{eq:L-kernel}, we have $Lf=LQ$. Let
	\[
	\Phi(x,v)
	=
	\sum_{k=1}^M
	\psi_k(x)p_k(v)\sqrt\mu(v),
	\]
	where $p_1,\ldots,p_M$ are fixed polynomials and
	$\psi_1,\ldots,\psi_M\in L^2(\cD)$. By the symmetry \eqref{eq:L-symmetry} of $L$,
	\begin{align}
		\left|
		\int_{\cD}
		\langle LQ,\Phi\rangle_v
		\dd x
		\right|
		&=
		\left|
		(Q,L\Phi)_{L^2_{x,v}}
		\right|
		\notag\\
		&\le
		\|\mu^{1/4}Q\|_{L^2(\cD\times\bR_v^3)}
		\|\mu^{-1/4}L\Phi\|_{L^2(\cD\times\bR_v^3)}
		\notag\\
		&\le
		C
		\|\mu^{1/4}Q\|_{L^2(\cD\times\bR_v^3)}
		\left(
		\sum_{k=1}^M
		\|\psi_k\|_{L^2(\cD)}^2
		\right)^{1/2}.
		\label{eq:appendix-collision-Q-bound}
	\end{align}
	The last step is Lemma~\ref{lem:poly-maxwellian} applied to each
	$p_k$, which in particular guarantees that $L\Phi(x,\cdot)$ is a function
	with $\mu^{-1/4}L\Phi(x,\cdot)\in L^2(\bR_v^3)$. Hence the pairing
	$\la LQ,\Phi\ra_v$ may be written as the $L^2$ inner product
	$(Q,L\Phi)_{L^2_v}$. The constant $C$ depends only on the fixed
	polynomials $p_1,\ldots,p_M$.
	
	For $0\le s<t\le T$, let
	\[
	\cD_{s,t}:=(s,t)\times \cD,\qquad   \gamma_{s,t}:=(s,t)\times\gamma
	\]
	and define
	\begin{equation*}
		\cR(s,t)
		:=
		\|\mu^{1/4}Q\|_
		{L^2(\cD_{s,t}\times\bR_v^3)}^2
		+
		\|S\|_{\cM(s,t)}^2
		+
		\|f_\gamma\|_{L^2(\gamma_{s,t})}^2.
	\end{equation*}

	\medskip
	\noindent\emph{Estimate of $\mathbf c$.}
	For each $r\in[0,T]$, let $\phi_{\mathbf c}(r)$ solve
	\begin{equation}\label{eq:appendix-phi-c}
		\begin{cases}
			-\Delta_x\phi_{\mathbf c}(r,\cdot)=\mathbf c(r,\cdot)
			&\text{in }\cD,\\
			\phi_{\mathbf c}(r,\cdot)=0
			&\text{on }\partial\cD.
		\end{cases}
	\end{equation}
	The Dirichlet elliptic estimate gives
	\begin{equation}\label{eq:appendix-elliptic-c}
		\|\phi_{\mathbf c}(r)\|_{H^2(\cD)}
		\le
		C\|\mathbf c(r)\|_{L^2(\cD)}.
	\end{equation}
	Under the present smoothness assumption on $f$,
	\[
	\mathbf c\in C^\infty([0,T];L^2(\cD)).
	\]
	Since the inverse Dirichlet Laplacian is linear and bounded from
	$L^2(\cD)$ to $H^2(\cD)\cap H_0^1(\cD)$, it follows that
	\[
	\phi_{\mathbf c}
	\in
	C^\infty\bigl([0,T];H^2(\cD)\cap H_0^1(\cD)\bigr).
	\]
	
	Define
	\begin{equation*}
		\Phi_{\mathbf c}
		:=
		\frac{4}{5\pi^{3/2}}
		\left(|v|^2-\frac52\right)
		v\cdot\nabla_x\phi_{\mathbf c}\sqrt\mu.
	\end{equation*}
	By the Gaussian identities
	\begin{align}\label{260717408}
		\begin{gathered}
			\qquad\int_{\bR^3}
			\left(|v|^2-\frac52\right)
			v_iv_j\mu
			\dd v
			=0,\qquad \int_{\bR^3}
			\left(|v|^2-\frac32\right)
			\left(|v|^2-\frac52\right)
			v_iv_j\mu
			\dd v
			=
			\frac54\pi^{3/2}\delta_{ij}    
		\end{gathered}
	\end{align}
	for all $i,\,j\in \{1,2,3\}$, and the oddness of the relevant integrands, we obtain
	\[
	\left(
	\mathbf P f,
	v\cdot\nabla_x\Phi_{\mathbf c}
	\right)_{x,v}
	=
	\int_{\cD}
	\mathbf c\,\Delta_x\phi_{\mathbf c}
	\dd x.
	\]
	Therefore,
	\begin{equation}\label{eq:appendix-main-c}
		-
		\int_s^t
		\left(
		\mathbf P f,
		v\cdot\nabla_x\Phi_{\mathbf c}
		\right)_{x,v}
		\dd r
		=
		\|\mathbf c\|_{L^2(s,t;L^2(\cD))}^2.
	\end{equation}

	Set
	\begin{equation*}
		\mathscr E_{\mathbf c}(\tau)
		:=
		\left(
		f(\tau),\Phi_{\mathbf c}(\tau)
		\right)_{x,v}.
	\end{equation*}
	By \eqref{eq:appendix-elliptic-c} and the definition of
	$\cM$,
	\begin{equation}\label{eq:appendix-endpoint-c-bound}
		|\mathscr E_{\mathbf c}(\tau)|
		\le
		C\|f(\tau)\|_{\cM}
		\|\mathbf c(\tau)\|_{L^2(\cD)}
		\le
		C\|f(\tau)\|_{\cM}^2
		\le
		C\|f(\tau)\|_{L^2(\cD\times\bR_v^3)}^2.
	\end{equation}
	Moreover, since the maps
	\[
	f\mapsto\mathbf c
	\qquad\text{and}\qquad
	\mathbf c\mapsto\Phi_{\mathbf c}
	\]
	are linear, the same argument yields, for every $0\le s<t\le T$,
	\begin{align}
		|\mathscr E_{\mathbf c}(t)-\mathscr E_{\mathbf c}(s)|
		&\le
		\left|
		\left(
		f(t)-f(s),\Phi_{\mathbf c}(t)
		\right)_{x,v}
		\right|
		+
		\left|
		\left(
		f(s),
		\Phi_{\mathbf c}(t)-\Phi_{\mathbf c}(s)
		\right)_{x,v}
		\right|
		\notag\\
		&\le
		C\|f(t)-f(s)\|_{\cM}
		\|\mathbf c(t)\|_{L^2(\cD)}
		+
		C\|f(s)\|_{\cM}
		\|\mathbf c(t)-\mathbf c(s)\|_{L^2(\cD)}
		\notag\\
		&\le
		C
		\left(
		\|f(t)\|_{\cM}
		+
		\|f(s)\|_{\cM}
		\right)
		\|f(t)-f(s)\|_{\cM}.
		\label{eq:appendix-interaction-c-difference}
	\end{align}
	
	Since $Lf=LQ$, the Green identity \eqref{ibp_dy}, applicable by
	Lemma~\ref{lem:admissible-test-functions}, together with the oddness of
	the relevant integrands and
	\eqref{eq:appendix-main-c} yield
	\begin{align}
		\|\mathbf c\|_{L^2(s,t;L^2(\cD))}^2
		={}&
		\mathscr E_{\mathbf c}(s)
		-
		\mathscr E_{\mathbf c}(t)
		+
		\int_s^t
		\left(
		f,\partial_r\Phi_{\mathbf c}
		\right)_{x,v}
		\dd r
		\notag\\
		&+
		\int_s^t
		\left(
		Q,v\cdot\nabla_x\Phi_{\mathbf c}
		\right)_{x,v}
		\dd r
		-
		\int_s^t
		\langle LQ,\Phi_{\mathbf c}\rangle_{x,v}
		\dd r
		\notag\\
		&-
		\int_s^t\int_\gamma
		f_\gamma\Phi_{\mathbf c}
		\operatorname{sgn}(n_x\cdot v)
		\dd\gamma\dd r
		+
		\int_s^t
		\langle S,\Phi_{\mathbf c}\rangle_{x,v}
		\dd r.
		\label{eq:appendix-c-identity}
	\end{align}
	
	Differentiating \eqref{eq:appendix-phi-c} in time and using
	\eqref{eq:appendix-moment-c}, we obtain
	\[
	-\Delta_x\partial_r\phi_{\mathbf c}
	=
	-\frac13\nabla_x\cdot\mathbf b
	-
	\nabla_x\cdot\Lambda
	+
	S_{\mathbf c}.
	\]
	Hence
	\begin{equation*}
		\|\nabla_x\partial_r\phi_{\mathbf c}\|_{L^2(\cD)}
		\le
		C\left(
		\|\mathbf b\|_{L^2(\cD)}
		+
		\|\Lambda\|_{L^2(\cD)}
		+
		\|S_{\mathbf c}\|_{L^2(\cD)}
		\right).
	\end{equation*}
	Moreover, by the Gaussian identity \eqref{260717408},
	\begin{equation*}
		\left(
		\mathbf P f,\partial_r\Phi_{\mathbf c}
		\right)_{x,v}
		=
		0.
	\end{equation*}
	Therefore 
	\eqref{eq:appendix-source-moment-control},
	\eqref{eq:appendix-higher-moment-bound}, and Young's inequality give,
	for every $\delta>0$,
	\begin{align}
		\left|
		\int_s^t
		\left(
		f,\partial_r\Phi_{\mathbf c}
		\right)_{x,v}
		\dd r
		\right|
		\le{}&
		\frac{\delta}{2}
		\|\mathbf b\|_{L^2(\cD_{s,t})}^2
		+
		C\left(1+\delta^{-1}\right)
		\left(
		\|\mu^{1/4}Q\|_
		{L^2(\cD_{s,t}\times\bR_v^3)}^2
		+
		\|S\|_{\cM(s,t)}^2
		\right).
		\label{eq:appendix-time-c-bound}
	\end{align}
	
	Lemma~\ref{lem:poly-maxwellian},
	\eqref{eq:appendix-elliptic-c}, and
	\eqref{eq:appendix-collision-Q-bound} imply
	\begin{align*}
		&
		\left|
		\left(
		Q,v\cdot\nabla_x\Phi_{\mathbf c}
		\right)_{x,v}
		\right|
		+
		\left|
		\langle LQ,\Phi_{\mathbf c}\rangle_{x,v}
		\right|
		\le
		C
		\|\mu^{1/4}Q\|_{L^2(\cD\times\bR_v^3)}
		\|\mathbf c\|_{L^2(\cD)}.
	\end{align*}
	Moreover, since the velocity profile of $\Phi_{\mathbf c}$ lies in
	$\operatorname{span}\cZ$ by \eqref{260715226},
	\[
	\left|
	\langle S,\Phi_{\mathbf c}\rangle_{x,v}
	\right|
	\le
	C
	\|S\|_{\cM}
	\|\mathbf c\|_{L^2(\cD)},
	\]
	and Lemma~\ref{lem:appendix-boundary-estimate} gives
	\[
	\|\Phi_{\mathbf c}\|_{L^2(\gamma)}
	\le
	C\|\mathbf c\|_{L^2(\cD)}.
	\]
	Integrating in time and applying Young's inequality, we obtain
	\begin{align}
		&\left|
		\int_s^t
		\left(
		Q,v\cdot\nabla_x\Phi_{\mathbf c}
		\right)_{x,v}
		\dd r
		\right|
		+
		\left|
		\int_s^t
		\langle LQ,\Phi_{\mathbf c}\rangle_{x,v}
		\dd r
		\right|
		\notag\\
		&\quad+
		\left|
		\int_s^t
		\langle S,\Phi_{\mathbf c}\rangle_{x,v}
		\dd r
		\right|
		+
		\left|
		\int_s^t\int_\gamma
		f_\gamma\Phi_{\mathbf c}
		\operatorname{sgn}(n_x\cdot v)
		\dd\gamma\dd r
		\right|\le
		\frac12
		\|\mathbf c\|_{L^2(s,t;L^2(\cD))}^2
		+
		C\cR(s,t).
		\label{eq:appendix-remainder-c-bound}
	\end{align}
	
	Combining \eqref{eq:appendix-c-identity},
	\eqref{eq:appendix-time-c-bound}, and
	\eqref{eq:appendix-remainder-c-bound}, and absorbing the
	$\mathbf c$-term, we conclude that, for every $\delta>0$,
	\begin{equation}\label{eq:appendix-c-estimate}
		\begin{aligned}
			\|\mathbf c\|_{L^2(s,t;L^2(\cD))}^2
			\le{}&
			\delta
			\|\mathbf b\|_{L^2(s,t;L^2(\cD))}^2
			+
			C\left(1+\delta^{-1}\right)\cR(s,t)
			+
			2\left(
			\mathscr E_{\mathbf c}(s)
			-
			\mathscr E_{\mathbf c}(t)
			\right).
		\end{aligned}
	\end{equation}

	\medskip
	\noindent\emph{Estimate of $\mathbf b$.}
	For $j\in\{1,2,3\}$ and $r\in[0,T]$, let $\phi_j(r)$ solve
	\begin{equation}\label{eq:appendix-phi-b}
		\begin{cases}
			-\Delta_x\phi_j=\mathbf b_j
			&\text{in }\cD,\\
			\phi_j=0
			&\text{on }\partial\cD.
		\end{cases}
	\end{equation}
	Then $\mathbf b_j\in C^\infty([0,T];L^2(\cD))$ and
	\begin{equation}\label{eq:appendix-elliptic-b}
		\|\phi_j(r)\|_{H^2(\cD)}
		\le
		C\|\mathbf b_j(r)\|_{L^2(\cD)}.
	\end{equation}
	
	Set
	\[
	\Phi_{\mathbf b}^j
	:=
	\sum_{m=1}^3\Phi_{\mathbf b}^{j,m},
	\]
	where, for $m\neq j$,
	\begin{equation*}
		\Phi_{\mathbf b}^{j,m}
		:=
		\left[
		\frac{8}{7\pi^{3/2}}
		|v|^2v_mv_j\,\partial_{x_m}\phi_j
		-
		\frac{2}{\pi^{3/2}}
		\left(v_m^2-\frac12\right)
		\partial_{x_j}\phi_j
		\right]\sqrt\mu,
	\end{equation*}
	and
	\begin{equation*}
		\Phi_{\mathbf b}^{j,j}
		:=
		\frac{2}{\pi^{3/2}}
		\left(v_j^2-\frac12\right)
		\partial_{x_j}\phi_j\sqrt\mu.
	\end{equation*}
	
	For fixed indices $i,\ell,m,j$, with no summation,
	\begin{align*}
		\int_{\bR^3}
		|v|^2v_iv_\ell v_mv_j\mu
		\dd v
		&=
		\frac78\pi^{3/2}
		\left(
		\delta_{i\ell}\delta_{mj}
		+
		\delta_{im}\delta_{\ell j}
		+
		\delta_{ij}\delta_{\ell m}
		\right),
		\\
		\int_{\bR^3}
		v_iv_\ell
		\left(v_m^2-\frac12\right)
		\mu
		\dd v
		&=
		\frac12\pi^{3/2}
		\delta_{im}\delta_{\ell m}.
	\end{align*}
	Using these identities and the oddness of the relevant integrands, we obtain, for all $j,\,m\in\{1,2,3\}$,
	\begin{align*}
		\left(
		\mathbf P f,
		v\cdot\nabla_x\Phi_{\mathbf b}^{j,m}
		\right)_{x,v}
		=
		\int_{\cD}
		\mathbf b_j\partial_{x_m}^2\phi_j
		\dd x.
	\end{align*}
	Therefore,
	\[
	\left(
	\mathbf P f,
	v\cdot\nabla_x\Phi_{\mathbf b}^j
	\right)_{x,v}
	=
	\int_{\cD}
	\mathbf b_j\Delta_x\phi_j
	\dd x,
	\]
	and hence
	\begin{equation}\label{eq:appendix-main-b}
		-
		\int_s^t
		\left(
		\mathbf P f,
		v\cdot\nabla_x\Phi_{\mathbf b}^j
		\right)_{x,v}
		\dd r
		=
		\|\mathbf b_j\|_{L^2(s,t;L^2(\cD))}^2.
	\end{equation}

	Define
	\begin{equation*}
		\mathscr E_{\mathbf b}^j(\tau)
		:=
		\left(
		f(\tau),\Phi_{\mathbf b}^j(\tau)
		\right)_{x,v},
		\qquad
		\mathscr E_{\mathbf b}(\tau)
		:=
		\sum_{j=1}^3\mathscr E_{\mathbf b}^j(\tau).
	\end{equation*}
	By the same argument as in
	\eqref{eq:appendix-endpoint-c-bound} and
	\eqref{eq:appendix-interaction-c-difference}, we obtain
	\begin{equation*}
		|\mathscr E_{\mathbf b}(\tau)|
		\le
		C\|f(\tau)\|_{\cM}^2
		\le
		C\|f(\tau)\|_{L^2(\cD\times\bR_v^3)}^2
	\end{equation*}
	for every $\tau\in[0,T]$, and
	\begin{equation*}
		|\mathscr E_{\mathbf b}(t)-\mathscr E_{\mathbf b}(s)|
		\le
		C
		\left(
		\|f(t)\|_{\cM}
		+
		\|f(s)\|_{\cM}
		\right)
		\|f(t)-f(s)\|_{\cM}
	\end{equation*}
	for every $0\le s<t\le T$.
	
	Since $Lf=LQ$, the Green identity \eqref{ibp_dy}, applicable by
	Lemma~\ref{lem:admissible-test-functions}, and
	\eqref{eq:appendix-main-b} give
	\begin{align}
		\|\mathbf b_j\|_{L^2(s,t;L^2(\cD))}^2
		={}&
		\mathscr E_{\mathbf b}^j(s)
		-
		\mathscr E_{\mathbf b}^j(t)
		+
		\int_s^t
		\left(
		f,\partial_r\Phi_{\mathbf b}^j
		\right)_{x,v}
		\dd r
		\notag\\
		&+
		\int_s^t
		\left(
		Q,v\cdot\nabla_x\Phi_{\mathbf b}^j
		\right)_{x,v}
		\dd r
		-
		\int_s^t
		\langle LQ,\Phi_{\mathbf b}^j\rangle_{x,v}
		\dd r
		\notag\\
		&-
		\int_s^t\int_\gamma
		f_\gamma\Phi_{\mathbf b}^j
		\operatorname{sgn}(n_x\cdot v)
		\dd\gamma\dd r
		+
		\int_s^t
		\langle S,\Phi_{\mathbf b}^j\rangle_{x,v}
		\dd r.
		\label{eq:appendix-b-identity}
	\end{align}
	
	The Gaussian identity
	\[
	\int_{\bR^3}
	\left(|v|^2-\frac32\right)
	\left(v_m^2-\frac12\right)
	\mu
	\dd v
	=
	\frac12\pi^{3/2}
	\]
	implies
	\begin{equation*}
		\left(
		\mathbf P f,
		\partial_r\Phi_{\mathbf b}^j
		\right)_{x,v}
		=
		-
		\int_{\cD}
		\mathbf c\,
		\partial_{x_j}\partial_r\phi_j
		\dd x.
	\end{equation*}
	Differentiating \eqref{eq:appendix-phi-b} in time and using
	\eqref{eq:appendix-moment-b}, we obtain
	\begin{equation*}
		-\Delta_x\partial_r\phi_j
		=
		-\partial_{x_j}(\mathbf a+\mathbf c)
		-
		\partial_{x_m}\Theta_{jm}
		+
		S_{\mathbf b,j}
	\end{equation*}
	in the sense of distributions. Hence
	\begin{equation*}
		\|\nabla_x\partial_r\phi_j\|_{L^2(\cD)}
		\le
		C\left(
		\|\mathbf a\|_{L^2(\cD)}
		+
		\|\mathbf c\|_{L^2(\cD)}
		+
		\|\Theta\|_{L^2(\cD)}
		+
		\|S_{\mathbf b}\|_{L^2(\cD)}
		\right).
	\end{equation*}
	Moreover,
	\[
	\left|
	\left(
	Q,\partial_r\Phi_{\mathbf b}^j
	\right)_{x,v}
	\right|
	\le
	C
	\|\mu^{1/4}Q\|_{L^2(\cD\times\bR_v^3)}
	\|\nabla_x\partial_r\phi_j\|_{L^2(\cD)}.
	\]
	Combining these estimates with
	\eqref{eq:appendix-source-moment-control} and
	\eqref{eq:appendix-higher-moment-bound}, we obtain, for every
	$\delta>0$,
	\begin{align}
		&\sum_{j=1}^3
		\left|
		\int_s^t
		\left(
		f,\partial_r\Phi_{\mathbf b}^j
		\right)_{x,v}
		\dd r
		\right|
		\notag\\
		&\qquad\le
		\frac{\delta}{2}
		\|\mathbf a\|_{L^2(s,t;L^2(\cD))}^2
		\notag\\
		&\qquad\quad+
		C\left(1+\delta^{-1}\right)
		\left(
		\|\mathbf c\|_{L^2(s,t;L^2(\cD))}^2
		+
		\|\mu^{1/4}Q\|_
		{L^2(\cD_{s,t}\times\bR_v^3)}^2
		+
		\|S\|_{\cM(s,t)}^2
		\right).
		\label{eq:appendix-time-b-bound}
	\end{align}
	
	By Lemma~\ref{lem:poly-maxwellian},
	\eqref{eq:appendix-collision-Q-bound}, and
	\eqref{eq:appendix-elliptic-b},
	\begin{align*}
		&
		\left|
		\left(
		Q,v\cdot\nabla_x\Phi_{\mathbf b}^j
		\right)_{x,v}
		\right|
		+
		\left|
		\langle LQ,\Phi_{\mathbf b}^j\rangle_{x,v}
		\right|
		\\
		&\qquad\le
		C
		\|\mu^{1/4}Q\|_{L^2(\cD\times\bR_v^3)}
		\|\mathbf b_j\|_{L^2(\cD)}.
	\end{align*}
	Likewise, the velocity profile of $\Phi_{\mathbf b}^j$ lies in
	$\operatorname{span}\cZ$ by \eqref{260715226}, so
	\[
	\left|
	\langle S,\Phi_{\mathbf b}^j\rangle_{x,v}
	\right|
	\le
	C
	\|S\|_{\cM}
	\|\mathbf b_j\|_{L^2(\cD)},
	\]
	while Lemma~\ref{lem:appendix-boundary-estimate} gives
	\[
	\|\Phi_{\mathbf b}^j\|_{L^2(\gamma)}
	\le
	C\|\mathbf b_j\|_{L^2(\cD)}.
	\]
	Integrating in time, summing over $j$, and applying Young's
	inequality, we find
	\begin{align}
		&\sum_{j=1}^3
		\Bigg[
		\left|
		\int_s^t
		\left(
		Q,v\cdot\nabla_x\Phi_{\mathbf b}^j
		\right)_{x,v}
		\dd r
		\right|
		+
		\left|
		\int_s^t
		\langle LQ,\Phi_{\mathbf b}^j\rangle_{x,v}
		\dd r
		\right|
		\notag\\
		&\qquad\qquad+
		\left|
		\int_s^t
		\langle S,\Phi_{\mathbf b}^j\rangle_{x,v}
		\dd r
		\right|
		+
		\left|
		\int_s^t\int_\gamma
		f_\gamma\Phi_{\mathbf b}^j
		\operatorname{sgn}(n_x\cdot v)
		\dd\gamma\dd r
		\right|
		\Bigg]
		\notag\\
		&\qquad\le
		\frac12
		\|\mathbf b\|_{L^2(s,t;L^2(\cD))}^2
		+
		C\cR(s,t).
		\label{eq:appendix-remainder-b-bound}
	\end{align}
	
	Summing \eqref{eq:appendix-b-identity} over $j$, using
	\eqref{eq:appendix-time-b-bound} and
	\eqref{eq:appendix-remainder-b-bound}, and absorbing the
	$\mathbf b$-term, we conclude that, for every $\delta>0$,
	\begin{equation}\label{eq:appendix-b-estimate}
		\begin{aligned}
			\|\mathbf b\|_{L^2(s,t;L^2(\cD))}^2
			\le{}&
			\delta
			\|\mathbf a\|_{L^2(s,t;L^2(\cD))}^2
			\\
			&+
			C\left(1+\delta^{-1}\right)
			\left(
			\|\mathbf c\|_{L^2(s,t;L^2(\cD))}^2
			+
			\cR(s,t)
			\right)
			\\
			&+
			2\left(
			\mathscr E_{\mathbf b}(s)
			-
			\mathscr E_{\mathbf b}(t)
			\right).
		\end{aligned}
	\end{equation}

	\medskip
	\noindent\emph{Estimate of $\mathbf a$.}
	For each $r\in[0,T]$, let $\phi_{\mathbf a}(r)$ solve
	\begin{equation}\label{eq:appendix-phi-a}
		\begin{cases}
			-\Delta_x\phi_{\mathbf a}=\mathbf a
			&\text{in }\cD,\\
			\phi_{\mathbf a}=0
			&\text{on }\partial\cD.
		\end{cases}
	\end{equation}
	Then
	\begin{equation}\label{eq:appendix-elliptic-a}
		\|\phi_{\mathbf a}(r)\|_{H^2(\cD)}
		\le
		C\|\mathbf a(r)\|_{L^2(\cD)}.
	\end{equation}
	Define
	\begin{equation*}
		\Phi_{\mathbf a}
		:=
		-
		\frac{4}{5\pi^{3/2}}
		\left(|v|^2-5\right)
		v\cdot\nabla_x\phi_{\mathbf a}\sqrt\mu.
	\end{equation*}
	The Gaussian identities
	\begin{align*}
		\int_{\bR^3}
		\left(|v|^2-5\right)v_iv_j\mu
		\dd v
		&=
		-\frac54\pi^{3/2}\delta_{ij},
		\\
		\int_{\bR^3}
		\left(|v|^2-\frac32\right)
		\left(|v|^2-5\right)v_iv_j\mu
		\dd v
		&=
		0,
	\end{align*}
	together with parity, imply
	\[
	\left(
	\mathbf P f,
	v\cdot\nabla_x\Phi_{\mathbf a}
	\right)_{x,v}
	=
	\int_{\cD}
	\mathbf a\Delta_x\phi_{\mathbf a}
	\dd x.
	\]
	Hence
	\begin{equation}\label{eq:appendix-main-a}
		-
		\int_s^t
		\left(
		\mathbf P f,
		v\cdot\nabla_x\Phi_{\mathbf a}
		\right)_{x,v}
		\dd r
		=
		\|\mathbf a\|_{L^2(s,t;L^2(\cD))}^2.
	\end{equation}

	Set
	\begin{equation*}
		\mathscr E_{\mathbf a}(\tau)
		:=
		\left(
		f(\tau),\Phi_{\mathbf a}(\tau)
		\right)_{x,v}.
	\end{equation*}
	By the same argument as in
	\eqref{eq:appendix-endpoint-c-bound} and
	\eqref{eq:appendix-interaction-c-difference}, we obtain
	\begin{equation*}
		|\mathscr E_{\mathbf a}(\tau)|
		\le
		C\|f(\tau)\|_{\cM}^2
		\le
		C\|f(\tau)\|_{L^2(\cD\times\bR_v^3)}^2
	\end{equation*}
	for every $\tau\in[0,T]$, and
	\begin{equation*}
		|\mathscr E_{\mathbf a}(t)-\mathscr E_{\mathbf a}(s)|
		\le
		C
		\left(
		\|f(t)\|_{\cM}
		+
		\|f(s)\|_{\cM}
		\right)
		\|f(t)-f(s)\|_{\cM}
	\end{equation*}
	for every $0\le s<t\le T$.

	Since $Lf=LQ$, the Green identity \eqref{ibp_dy}, applicable by
	Lemma~\ref{lem:admissible-test-functions}, and
	\eqref{eq:appendix-main-a} give
	\begin{align}
		\|\mathbf a\|_{L^2(s,t;L^2(\cD))}^2
		={}&
		\mathscr E_{\mathbf a}(s)
		-
		\mathscr E_{\mathbf a}(t)
		+
		\int_s^t
		\left(
		f,\partial_r\Phi_{\mathbf a}
		\right)_{x,v}
		\dd r
		\notag\\
		&+
		\int_s^t
		\left(
		Q,v\cdot\nabla_x\Phi_{\mathbf a}
		\right)_{x,v}
		\dd r
		-
		\int_s^t
		\langle LQ,\Phi_{\mathbf a}\rangle_{x,v}
		\dd r
		\notag\\
		&-
		\int_s^t\int_\gamma
		f_\gamma\Phi_{\mathbf a}
		\operatorname{sgn}(n_x\cdot v)
		\dd\gamma\dd r
		+
		\int_s^t
		\langle S,\Phi_{\mathbf a}\rangle_{x,v}
		\dd r.
		\label{eq:appendix-a-identity}
	\end{align}
	
	A direct Gaussian computation gives
	\begin{equation*}
		\left(
		\mathbf P f,
		\partial_r\Phi_{\mathbf a}
		\right)_{x,v}
		=
		\int_{\cD}
		\mathbf b\cdot
		\nabla_x\partial_r\phi_{\mathbf a}
		\dd x.
	\end{equation*}
	Differentiating \eqref{eq:appendix-phi-a} in time and using
	\eqref{eq:appendix-moment-a}, we obtain
	\begin{equation*}
		-\Delta_x\partial_r\phi_{\mathbf a}
		=
		-\frac12\nabla_x\cdot\mathbf b
		+
		S_{\mathbf a}
	\end{equation*}
	in the sense of distributions. Therefore,
	\begin{equation*}
		\|\nabla_x\partial_r\phi_{\mathbf a}\|_{L^2(\cD)}
		\le
		C\left(
		\|\mathbf b\|_{L^2(\cD)}
		+
		\|S_{\mathbf a}\|_{L^2(\cD)}
		\right).
	\end{equation*}
	Moreover,
	\[
	\left|
	\left(
	Q,\partial_r\Phi_{\mathbf a}
	\right)_{x,v}
	\right|
	\le
	C
	\|\mu^{1/4}Q\|_{L^2(\cD\times\bR_v^3)}
	\|\nabla_x\partial_r\phi_{\mathbf a}\|_{L^2(\cD)}.
	\]
	Thus, by \eqref{eq:appendix-source-moment-control},
	\begin{equation}\label{eq:appendix-time-a-bound}
		\left|
		\int_s^t
		\left(
		f,\partial_r\Phi_{\mathbf a}
		\right)_{x,v}
		\dd r
		\right|
		\le
		C\|\mathbf b\|_{L^2(s,t;L^2(\cD))}^2
		+
		C\cR(s,t).
	\end{equation}
	
	Lemma~\ref{lem:poly-maxwellian},
	\eqref{eq:appendix-collision-Q-bound}, and
	\eqref{eq:appendix-elliptic-a} imply
	\begin{align*}
		&
		\left|
		\left(
		Q,v\cdot\nabla_x\Phi_{\mathbf a}
		\right)_{x,v}
		\right|
		+
		\left|
		\langle LQ,\Phi_{\mathbf a}\rangle_{x,v}
		\right|
		\le
		C
		\|\mu^{1/4}Q\|_{L^2(\cD\times\bR_v^3)}
		\|\mathbf a\|_{L^2(\cD)}.
	\end{align*}
	Since the velocity profile of $\Phi_{\mathbf a}$ lies in
	$\operatorname{span}\cZ$ by \eqref{260715226},
	\[
	\left|
	\langle S,\Phi_{\mathbf a}\rangle_{x,v}
	\right|
	\le
	C\|S\|_{\cM}
	\|\mathbf a\|_{L^2(\cD)}.
	\]
	Finally, Lemma~\ref{lem:appendix-boundary-estimate} gives
	\[
	\|\Phi_{\mathbf a}\|_{L^2(\gamma)}
	\le
	C\|\mathbf a\|_{L^2(\cD)}.
	\]
	Integrating in time and applying Young's inequality, we obtain
	\begin{align}
		&\left|
		\int_s^t
		\left(
		Q,v\cdot\nabla_x\Phi_{\mathbf a}
		\right)_{x,v}
		\dd r
		\right|
		+
		\left|
		\int_s^t
		\langle LQ,\Phi_{\mathbf a}\rangle_{x,v}
		\dd r
		\right|
		\notag\\
		&\quad+
		\left|
		\int_s^t
		\langle S,\Phi_{\mathbf a}\rangle_{x,v}
		\dd r
		\right|
		+
		\left|
		\int_s^t\int_\gamma
		f_\gamma\Phi_{\mathbf a}
		\operatorname{sgn}(n_x\cdot v)
		\dd\gamma\dd r
		\right|
		\notag\\
		&\qquad\le
		\frac12
		\|\mathbf a\|_{L^2(s,t;L^2(\cD))}^2
		+
		C\cR(s,t).
		\label{eq:appendix-remainder-a-bound}
	\end{align}
	
	Combining \eqref{eq:appendix-a-identity},
	\eqref{eq:appendix-time-a-bound}, and
	\eqref{eq:appendix-remainder-a-bound}, and absorbing the
	$\mathbf a$-term, we conclude that
	\begin{equation}\label{eq:appendix-a-estimate}
		\begin{aligned}
			\|\mathbf a\|_{L^2(s,t;L^2(\cD))}^2
			\le{}&
			C
			\|\mathbf b\|_{L^2(s,t;L^2(\cD))}^2
			+
			C\cR(s,t)
			+
			2\left(
			\mathscr E_{\mathbf a}(s)
			-
			\mathscr E_{\mathbf a}(t)
			\right).
		\end{aligned}
	\end{equation}

	\medskip
	\noindent\emph{Closing the estimates.}
	Set
	\[
	\mathsf A
	:=
	\|\mathbf a\|_{L^2(s,t;L^2(\cD))}^2,
	\qquad
	\mathsf B
	:=
	\|\mathbf b\|_{L^2(s,t;L^2(\cD))}^2,
	\qquad
	\mathsf C
	:=
	\|\mathbf c\|_{L^2(s,t;L^2(\cD))}^2,
	\]
	and
	\[
	\mathbf D[\mathscr E_\alpha]
	:=
	\mathscr E_\alpha(s)-\mathscr E_\alpha(t),
	\qquad
	\alpha\in\{\mathbf a,\mathbf b,\mathbf c\}.
	\]
	After enlarging $C\ge1$ if necessary,
	\eqref{eq:appendix-c-estimate},
	\eqref{eq:appendix-b-estimate}, and
	\eqref{eq:appendix-a-estimate} give
	\begin{align}
		\mathsf C
		&\le
		\varepsilon_{\mathbf c}\mathsf B
		+
		C\left(1+\varepsilon_{\mathbf c}^{-1}\right)
		\cR(s,t)
		+
		2\mathbf D[\mathscr E_{\mathbf c}],
		\label{eq:appendix-closing-c}
		\\
		\mathsf B
		&\le
		\varepsilon_{\mathbf b}\mathsf A
		+
		C\left(1+\varepsilon_{\mathbf b}^{-1}\right)
		\left(
		\mathsf C+\cR(s,t)
		\right)
		+
		2\mathbf D[\mathscr E_{\mathbf b}],
		\label{eq:appendix-closing-b}
		\\
		\mathsf A
		&\le
		C\left(
		\mathsf B+\cR(s,t)
		\right)
		+
		2\mathbf D[\mathscr E_{\mathbf a}].
		\label{eq:appendix-closing-a}
	\end{align}
	
	Choose
	\[
	\varepsilon_{\mathbf b}
	:=
	\frac{1}{8C},
	\qquad
	K
	:=
	C\left(1+\varepsilon_{\mathbf b}^{-1}\right),
	\qquad
	\varepsilon_{\mathbf c}
	:=
	\frac{1}{8K}.
	\]
	Multiplying \eqref{eq:appendix-closing-a} by
	$2\varepsilon_{\mathbf b}$ and
	\eqref{eq:appendix-closing-c} by $2K$, and adding the resulting
	inequalities to \eqref{eq:appendix-closing-b}, we obtain
	\begin{align*}
		2\varepsilon_{\mathbf b}\mathsf A
		+
		\mathsf B
		+
		2K\mathsf C
		\le{}&
		\varepsilon_{\mathbf b}\mathsf A
		+
		\left(
		2C\varepsilon_{\mathbf b}
		+
		2K\varepsilon_{\mathbf c}
		\right)\mathsf B
		+
		K\mathsf C
		\\
		&+
		C\cR(s,t)
		+
		4\varepsilon_{\mathbf b}\mathbf D[\mathscr E_{\mathbf a}]
		+
		2\mathbf D[\mathscr E_{\mathbf b}]
		+
		4K\mathbf D[\mathscr E_{\mathbf c}].
	\end{align*}
	Here $C$ is still the constant fixed above, so that
	\[
	2C\varepsilon_{\mathbf b}
	=
	2K\varepsilon_{\mathbf c}
	=
	\frac14,
	\]
	and therefore
	\begin{align}
		\varepsilon_{\mathbf b}\mathsf A
		+
		\frac12\mathsf B
		+
		K\mathsf C
		\le{}&
		C\cR(s,t)
		+
		4\varepsilon_{\mathbf b}\mathbf D[\mathscr E_{\mathbf a}]
		+
		2\mathbf D[\mathscr E_{\mathbf b}]
		+
		4K\mathbf D[\mathscr E_{\mathbf c}].
		\label{eq:appendix-closing-weighted}
	\end{align}
	Here, and only in the coefficient of $\cR(s,t)$, the
	letter $C$ now denotes a possibly larger constant depending only on
	$\cD$. The values of $\varepsilon_{\mathbf b}$, $K$ and
	$\varepsilon_{\mathbf c}$ remain those fixed above.
	Since $C\ge1$,
	\[
	0<\varepsilon_{\mathbf b}\le\frac12,
	\qquad
	\varepsilon_{\mathbf b}\le K,
	\]
	and hence
	\[
	\varepsilon_{\mathbf b}
	\left(
	\mathsf A+\mathsf B+\mathsf C
	\right)
	\le
	\varepsilon_{\mathbf b}\mathsf A
	+
	\frac12\mathsf B
	+
	K\mathsf C.
	\]
	Dividing \eqref{eq:appendix-closing-weighted} by
	$\varepsilon_{\mathbf b}$, we obtain
	\begin{equation}\label{eq:appendix-closing-abc}
		\mathsf A+\mathsf B+\mathsf C
		\le
		C\cR(s,t)
		+
		\kappa_{\mathbf a}\mathbf D[\mathscr E_{\mathbf a}]
		+
		\kappa_{\mathbf b}\mathbf D[\mathscr E_{\mathbf b}]
		+
		\kappa_{\mathbf c}\mathbf D[\mathscr E_{\mathbf c}],
	\end{equation}
	where
	\[
	\kappa_{\mathbf a}:=4,
	\qquad
	\kappa_{\mathbf b}:=2\varepsilon_{\mathbf b}^{-1},
	\qquad
	\kappa_{\mathbf c}:=4K\varepsilon_{\mathbf b}^{-1}.
	\]

	Define
	\begin{equation*}
		\eta_f(\tau)
		:=
		-
		\frac32\pi^{3/2}
		\left(
		\kappa_{\mathbf a}\mathscr E_{\mathbf a}(\tau)
		+
		\kappa_{\mathbf b}\mathscr E_{\mathbf b}(\tau)
		+
		\kappa_{\mathbf c}\mathscr E_{\mathbf c}(\tau)
		\right).
	\end{equation*}
	The bounds \eqref{eq:macro-eta-bound} and
	\eqref{eq:macro-eta-difference-bound} for $\eta_f$ are proved, for arbitrary
	$f\in C\big([0,T];L^2(\cD\times\bR_v^3)\big)$, in the final step of this
	proof. They are not needed for what follows.
	
	By the orthogonality of the collision invariants,
	\begin{equation}\label{eq:appendix-Pf-macro-norm}
		\|\mathbf P f\|_{L^2(\cD\times\bR_v^3)}^2
		=
		\pi^{3/2}
		\left(
		\|\mathbf a\|_{L^2(\cD)}^2
		+
		\frac12\|\mathbf b\|_{L^2(\cD)}^2
		+
		\frac32\|\mathbf c\|_{L^2(\cD)}^2
		\right).
	\end{equation}
	Therefore, by \eqref{eq:appendix-closing-abc},
	\begin{align*}
		\int_s^t
		\|\mathbf P f(r)\|_{L^2(\cD\times\bR_v^3)}^2
		\dd r
		&\le
		\frac32\pi^{3/2}
		\left(
		\mathsf A+\mathsf B+\mathsf C
		\right)
		\\
		&\le
		C\cR(s,t)
		+
		\eta_f(t)-\eta_f(s).
	\end{align*}
	This proves \eqref{eq:appendix-macro-estimate-endpoint}.
	
	\medskip
	\noindent\emph{The stationary case.}
	As in the time-dependent case, it suffices to
	argue for $f\in C_c^\infty(\overline{\cD}\times\bR_v^3)$.
	For the stationary problem, set
	\[
	Q:=(\mathbf I-\mathbf P)f\quad \text{and}\quad 
	\cR_{\mathrm{st}}
	:=
	\|\mu^{1/4}Q\|_{L^2(\cD\times\bR_v^3)}^2
	+
	\|S\|_{\cM}^2
	+
	\|f_\gamma\|_{L^2(\gamma)}^2.
	\]
	Using the same test functions as above, but without the time
	derivative and endpoint terms, we obtain the following three estimates:
	\begin{align*}
		\|\mathbf c\|_{L^2(\cD)}^2
		&\le
		\frac12\|\mathbf c\|_{L^2(\cD)}^2
		+
		C\cR_{\mathrm{st}},
		\\
		\|\mathbf b\|_{L^2(\cD)}^2
		&\le
		\frac12\|\mathbf b\|_{L^2(\cD)}^2
		+
		C\cR_{\mathrm{st}},
		\\
		\|\mathbf a\|_{L^2(\cD)}^2
		&\le
		\frac12\|\mathbf a\|_{L^2(\cD)}^2
		+
		C\cR_{\mathrm{st}}.
	\end{align*}
	After absorption,
	\[
	\|\mathbf a\|_{L^2(\cD)}^2
	+
	\|\mathbf b\|_{L^2(\cD)}^2
	+
	\|\mathbf c\|_{L^2(\cD)}^2
	\le
	C\cR_{\mathrm{st}}.
	\]
	Together with \eqref{eq:appendix-Pf-macro-norm}, this yields
	\[
	\|\mathbf P f\|_{L^2(\cD\times\bR_v^3)}^2
	\le
	C\left(
	\|\mu^{1/4}(\mathbf I-\mathbf P)f\|_
	{L^2(\cD\times\bR_v^3)}^2
	+
	\|S\|_{\cM}^2
	+
	\|f_\gamma\|_{L^2(\gamma)}^2
	\right),
	\]
	which proves
	\eqref{eq:appendix-macro-estimate-endpoint_st}.
\end{proof}

\subsection{Proof of Theorem~\ref{260322442} (well-posedness theorem) and velocity-weighted estimates}
\label{solv_linlandau_dy}

We present the time-dependent proof. In the stationary case, the time
derivative and endpoint terms are omitted, the stationary Green identity
\eqref{ibp_st} supplies the outflow term
$\|f_\gamma\|_{L^2(\gamma^+)}^2$, and
\eqref{eq:appendix-macro-estimate-endpoint_st} replaces the dynamic macroscopic estimate.
Only the base point of the continuity argument requires a separate treatment.

We introduce the coercive model Landau operator
\begin{equation}\label{def:A_0}
	A_0f
	:=
	-D_i\left(\sigma^{ij}D_jf\right)
	+
	\sigma^{ij}v_iv_jf.
\end{equation}
More precisely, $A_0:H_\sigma^1\to H_\sigma^{-1}$ is defined by
\begin{equation}\label{eq:model-landau-operator}
	\langle A_0f,\varphi\rangle_{H_\sigma^{-1},H_\sigma^1}
	:=
	\int_{\bR^3}
	\left(
	\sigma^{ij}D_jf\,D_i\varphi
	+
	\sigma^{ij}v_iv_jf\varphi
	\right)
	\dd v.
\end{equation}
By the Cauchy--Schwarz inequality,
\[
\left|
\langle A_0f,\varphi\rangle_{H_\sigma^{-1},H_\sigma^1}
\right|
\le
\|f\|_{H_\sigma^1}
\|\varphi\|_{H_\sigma^1}.
\]
Taking the supremum over
$\|\varphi\|_{H_\sigma^1}=1$ and then choosing
$\varphi=f/\|f\|_{H_\sigma^1}$ when $f\neq0$, we obtain
\begin{equation}\label{eq:model-landau-isometry}
	\|A_0f\|_{H_\sigma^{-1}}
	=
	\|f\|_{H_\sigma^1}.
\end{equation}
The bilinear form in \eqref{eq:model-landau-operator} is symmetric,
and
\begin{equation}\label{260715202}
	\langle A_0f,f\rangle_{H_\sigma^{-1},H_\sigma^1}
	=
	\|f\|_{H_\sigma^1}^2.
\end{equation}
In particular, $A_0$ is coercive with coercivity constant $1$.

\begin{lemma}
	\label{260502320}
	Let $0<T<\infty$. For every
	\[
	S\in L^2(\cD_T;H_\sigma^{-1}),
	\qquad
	f_0\in L^2(\cD\times\bR_v^3),
	\qquad
	h\in L^2(\gamma_T^-),
	\]
	the problem (recall $A_0$ defined in \eqref{def:A_0})
	\begin{equation}\label{FP_dy_sigma}
		\begin{cases}
			(\partial_t+v\cdot\nabla_x)f+A_0f=S
			&\text{in }\cD_T\times\bR_v^3,\\
			f(0,\cdot)=f_0
			&\text{in }\cD\times\bR_v^3,\\
			f=h
			&\text{on }\gamma_T^-
		\end{cases}
	\end{equation}
	admits a unique solution
	$
	f\in\cH_{\kin,\sigma}(\cD_T\times\bR_v^3)
	$, and
	\begin{equation}\label{eq:model-landau-inflow-estimate}
		\begin{aligned}
			&\|f\|_{\cH_{\kin,\sigma}(\cD_T\times\bR_v^3)}
			+\|f(T,\cdot)\|_{L^2(\cD\times\bR_v^3)}
			+\|f_\gamma\|_{L^2(\gamma_T^+)}\\
			\le\,&
			C\left(
			\|S\|_{L^2(\cD_T;H_\sigma^{-1})}
			+
			\|f_0\|_{L^2(\cD\times\bR_v^3)}
			+
			\|h\|_{L^2(\gamma_T^-)}
			\right),
		\end{aligned}
	\end{equation}
	where $C$ is a universal constant.
\end{lemma}

\begin{proof}
	\medskip
	\noindent\emph{A priori estimate and uniqueness.}
	Let $f\in\cH_{\kin,\sigma}(\cD_T\times\bR_v^3)$ be a solution of
	\eqref{FP_dy_sigma}. Applying the Green identity from
	Proposition~\ref{2607151218}, and using
	$f_\gamma=h$ on $\gamma_T^-$ and \eqref{260715202}, we obtain
	\begin{align*}
		&\frac12
		\|f(T)\|_{L^2(\cD\times\bR_v^3)}^2
		+
		\frac12
		\|f_\gamma\|_{L^2(\gamma_T^+)}^2
		+
		\|f\|_{L^2(\cD_T;H_\sigma^1)}^2
		\notag\\
		&\qquad=
		\frac12
		\|f_0\|_{L^2(\cD\times\bR_v^3)}^2
		+
		\frac12
		\|h\|_{L^2(\gamma_T^-)}^2
		+
		\iint_{\cD_T}
		\langle S,f\rangle_{H_\sigma^{-1},H_\sigma^1}
		\dd x\dd t.
	\end{align*}
	By duality and Young's inequality,
	\begin{align*}
		\left|
		\iint_{\cD_T}
		\langle S,f\rangle_{H_\sigma^{-1},H_\sigma^1}
		\dd x\dd t
		\right|
		&\le
		\|S\|_{L^2(\cD_T;H_\sigma^{-1})}
		\|f\|_{L^2(\cD_T;H_\sigma^1)}
		\\
		&\le
		\frac12
		\|S\|_{L^2(\cD_T;H_\sigma^{-1})}^2
		+
		\frac12
		\|f\|_{L^2(\cD_T;H_\sigma^1)}^2.
	\end{align*}
	Therefore,
	\begin{align}
		&\|f(T)\|_{L^2(\cD\times\bR_v^3)}^2
		+
		\|f_\gamma\|_{L^2(\gamma_T^+)}^2
		+
		\|f\|_{L^2(\cD_T;H_\sigma^1)}^2
		\notag\\
		&\qquad\le
		\|f_0\|_{L^2(\cD\times\bR_v^3)}^2
		+
		\|h\|_{L^2(\gamma_T^-)}^2
		+
		\|S\|_{L^2(\cD_T;H_\sigma^{-1})}^2.
		\label{eq:model-landau-energy-estimate}
	\end{align}
	
	On the other hand, by \eqref{FP_dy_sigma} and
	\eqref{eq:model-landau-isometry},
	\begin{align}
		\|(\partial_t+v\cdot\nabla_x)f\|_
		{L^2(\cD_T;H_\sigma^{-1})}
		&\le
		\|S\|_{L^2(\cD_T;H_\sigma^{-1})}
		+
		\|A_0f\|_{L^2(\cD_T;H_\sigma^{-1})}
		\notag\\
		&=
		\|S\|_{L^2(\cD_T;H_\sigma^{-1})}
		+
		\|f\|_{L^2(\cD_T;H_\sigma^1)}.
		\label{eq:model-landau-transport-estimate}
	\end{align}
	Hence, by
	\eqref{eq:model-landau-energy-estimate},
	\eqref{eq:model-landau-transport-estimate}, and the definition of
	$\cH_{\kin,\sigma}(\cD_T\times\bR_v^3)$,
	\begin{align*}
		\|f\|_{\cH_{\kin,\sigma}(\cD_T\times\bR_v^3)}\le
		C\left(
		\|S\|_{L^2(\cD_T;H_\sigma^{-1})}
		+
		\|f_0\|_{L^2(\cD\times\bR_v^3)}
		+
		\|h\|_{L^2(\gamma_T^-)}
		\right),
	\end{align*}
	where $C>0$ is universal. In particular,
	\eqref{eq:model-landau-inflow-estimate} holds.

	Uniqueness follows immediately by applying
	\eqref{eq:model-landau-energy-estimate} to the difference of two
	solutions with the same data.

	\medskip
	\noindent\emph{Existence.}
	By Lemma~\ref{lem:H_sigma-Hilbert},
	\begin{equation*}
		H_\sigma^{-1}(\bR_v^3)
		\hookrightarrow
		H_v^{-1}(\bR_v^3),
		\qquad
		\|F\|_{H_v^{-1}}
		\le
		C\|F\|_{H_\sigma^{-1}}.
	\end{equation*}
	In particular,
	\[
	\|S\|_{L^2(\cD_T;H_v^{-1})}
	\le
	C\|S\|_{L^2(\cD_T;H_\sigma^{-1})}.
	\]
	
	For $\varepsilon\in[0,1]$, define
	\begin{equation*}
		A_0^\varepsilon
		:=
		\varepsilon A_0
		+
		(1-\varepsilon)(-\Delta_v+1).
	\end{equation*}
	Equivalently,
	\[
	A_0^\varepsilon f
	=
	-D_i\left[
	\left(
	\varepsilon\sigma^{ij}
	+
	(1-\varepsilon)\delta_{ij}
	\right)D_jf
	\right]
	+
	\left(
	\varepsilon\sigma^{ij}v_iv_j
	+
	1-\varepsilon
	\right)f.
	\]
	
	Let $\varepsilon\in[0,1)$, and suppose that
	$\widetilde f\in\cH_\kin(\cD_T\times\bR_v^3)$ solves
	\begin{equation}\label{eq:model-epsilon-problem}
		\begin{cases}
			(\partial_t+v\cdot\nabla_x)\widetilde f
			+
			A_0^\varepsilon\widetilde f
			=
			\widetilde S
			&\text{in }\cD_T\times\bR_v^3,\\
			\widetilde f(0,\cdot)=\widetilde f_0
			&\text{in }\cD\times\bR_v^3,\\
			\widetilde f=\widetilde h
			&\text{on }\gamma_T^-,
		\end{cases}
	\end{equation}
	where
	\[
	\widetilde S\in L^2(\cD_T;H_v^{-1}),
	\qquad
	\widetilde f_0\in L^2(\cD\times\bR_v^3),
	\qquad
	\widetilde h\in L^2(\gamma_T^-).
	\]
	Applying the Green identity from
	Proposition~\ref{trace_theorem}, we obtain
	\begin{align}
		&\frac12
		\|\widetilde f(T)\|_{L^2(\cD\times\bR_v^3)}^2
		+
		\frac12
		\|\widetilde f_\gamma\|_{L^2(\gamma_T^+)}^2
		+
		(1-\varepsilon)
		\|\widetilde f\|_{L^2(\cD_T;H_v^1)}^2
		+
		\varepsilon
		\|\widetilde f\|_{L^2(\cD_T;H_\sigma^1)}^2
		\notag\\
		&\qquad=
		\frac12
		\|\widetilde f_0\|_{L^2(\cD\times\bR_v^3)}^2
		+
		\frac12
		\|\widetilde h\|_{L^2(\gamma_T^-)}^2
		+
		\iint_{\cD_T}
		\langle
		\widetilde S,
		\widetilde f
		\rangle_{H_v^{-1},H_v^1}
		\dd x\dd t.
		\label{260718829}
	\end{align}
	By duality and Young's inequality,
	\begin{align*}
		\left|
		\iint_{\cD_T}
		\langle
		\widetilde S,
		\widetilde f
		\rangle_{H_v^{-1},H_v^1}
		\dd x\dd t
		\right|
		&\le
		\|\widetilde S\|_{L^2(\cD_T;H_v^{-1})}
		\|\widetilde f\|_{L^2(\cD_T;H_v^1)}
		\\
		&\le
		\frac{1-\varepsilon}{2}
		\|\widetilde f\|_{L^2(\cD_T;H_v^1)}^2
		+
		\frac{1}{2(1-\varepsilon)}
		\|\widetilde S\|_{L^2(\cD_T;H_v^{-1})}^2.
	\end{align*}
	Therefore,
	\begin{align}
		(1-\varepsilon)
		\|\widetilde f\|_{L^2(\cD_T;H_v^1)}^2
		\le
		\|\widetilde f_0\|_{L^2(\cD\times\bR_v^3)}^2
		+
		\|\widetilde h\|_{L^2(\gamma_T^-)}^2
		+
		\frac{1}{1-\varepsilon}
		\|\widetilde S\|_{L^2(\cD_T;H_v^{-1})}^2.
		\label{eq:model-epsilon-energy-estimate}
	\end{align}
	
	Fix $\varepsilon_*\in[0,1)$. Since
	$1-\varepsilon\ge1-\varepsilon_*$ for
	$\varepsilon\in[0,\varepsilon_*]$,
	\eqref{eq:model-epsilon-energy-estimate} gives
	\begin{equation}\label{eq:model-epsilon-H1-estimate}
		\|\widetilde f\|_{L^2(\cD_T;H_v^1)}
		\le
		\frac{C}{1-\varepsilon_*}
		\left(
		\|\widetilde f_0\|_{L^2(\cD\times\bR_v^3)}
		+
		\|\widetilde h\|_{L^2(\gamma_T^-)}
		+
		\|\widetilde S\|_{L^2(\cD_T;H_v^{-1})}
		\right).
	\end{equation}
	
	Since the coefficients $\sigma^{ij}$ and
	$\sigma^{ij}v_iv_j$ are bounded, $A_0$ is bounded from $H_v^1$ to
	$H_v^{-1}$:
	\begin{equation*}
		\|A_0g\|_{H_v^{-1}}
		\le
		C\|g\|_{H_v^1}.
	\end{equation*}
	Together with the boundedness of
	$-\Delta_v+1:H_v^1\to H_v^{-1}$, this yields
	\begin{equation*}
		\|A_0^\varepsilon g\|_{H_v^{-1}}
		\le
		C\|g\|_{H_v^1},
		\qquad
		\varepsilon\in[0,1],
	\end{equation*}
	where $C>0$ is independent of $\varepsilon$. Hence
	\eqref{eq:model-epsilon-problem} implies
	\begin{equation}\label{eq:model-epsilon-transport-estimate}
		\|(\partial_t+v\cdot\nabla_x)\widetilde f\|_
		{L^2(\cD_T;H_v^{-1})}
		\le
		\|\widetilde S\|_{L^2(\cD_T;H_v^{-1})}
		+
		C\|\widetilde f\|_{L^2(\cD_T;H_v^1)}.
	\end{equation}
	Combining
	\eqref{eq:model-epsilon-H1-estimate} and
	\eqref{eq:model-epsilon-transport-estimate}, and using the initial
	and inflow boundary conditions, we obtain
	\begin{equation}\label{eq:model-epsilon-Hkin-estimate}
		\|\widetilde f\|_{\cH_\kin(\cD_T\times\bR_v^3)}
		\le
		\frac{C}{1-\varepsilon_*}
		\left(
		\|\widetilde f_0\|_{L^2(\cD\times\bR_v^3)}
		+
		\|\widetilde h\|_{L^2(\gamma_T^-)}
		+
		\|\widetilde S\|_{L^2(\cD_T;H_v^{-1})}
		\right)
	\end{equation}
	for every $\varepsilon\in[0,\varepsilon_*]$.
	
	For $\varepsilon=0$, we have
	\[
	A_0^0=-\Delta_v+1,
	\]
	and the corresponding inflow problem is uniquely solvable in
	$\cH_\kin(\cD_T\times\bR_v^3)$ by Theorem~\ref{solv}.
	(The base point of the stationary counterpart is \eqref{FP_st_classic} with	$\lambda=1$, uniquely solvable in $\cH_\kin(\cD\times\bR_v^3)$ by
	Proposition~\ref{solv_Del}.(2).)
	
	For each fixed $\varepsilon_*\in(0,1)$,
	\eqref{eq:model-epsilon-Hkin-estimate}
	allows us to apply the method of continuity along the family
	$\{A_0^{\theta\varepsilon_*}\}_{0\le\theta\le1}$.
	Hence the inflow problem with $A_0^{\varepsilon_*}$ is uniquely
	solvable in $\cH_\kin(\cD_T\times\bR_v^3)$.
	
	Since $\varepsilon_*\in(0,1)$ is arbitrary, for every $\varepsilon\in[0,1)$, the problem
	\begin{equation}\label{eq:model-epsilon-inflow-problem}
		\begin{cases}
			(\partial_t+v\cdot\nabla_x)f^\varepsilon
			+
			A_0^\varepsilon f^\varepsilon
			=
			S
			&\text{in }\cD_T\times\bR_v^3,\\
			f^\varepsilon(0,\cdot)=f_0
			&\text{in }\cD\times\bR_v^3,\\
			f^\varepsilon=h
			&\text{on }\gamma_T^-,
		\end{cases}
	\end{equation}
	admits a unique solution
	\[
	f^\varepsilon
	\in
	\cH_\kin(\cD_T\times\bR_v^3).
	\]

	\medskip
	\noindent\emph{Passage to the limit as $\varepsilon\uparrow1$.}
	Consider the family
	$\{f^\varepsilon\}_{1/2\le\varepsilon<1}$.
	Applying \eqref{260718829} to $f^\varepsilon$ and using
	\begin{align*}
		\left|
		\iint_{\cD_T}
		\left\langle
		S,
		f^\varepsilon
		\right\rangle_{H_\sigma^{-1},H_\sigma^1}
		\dd x\dd t
		\right|
		&\le
		\|S\|_{L^2(\cD_T;H_\sigma^{-1})}
		\|f^\varepsilon\|_{L^2(\cD_T;H_\sigma^1)}
		\\
		&\le
		\frac{\varepsilon}{2}
		\|f^\varepsilon\|_{L^2(\cD_T;H_\sigma^1)}^2
		+
		\frac{1}{2\varepsilon}
		\|S\|_{L^2(\cD_T;H_\sigma^{-1})}^2,
	\end{align*}
	we obtain
	\begin{align*}
		&\|f^\varepsilon(T)\|_{L^2(\cD\times\bR_v^3)}^2
		+
		\|f_\gamma^\varepsilon\|_{L^2(\gamma_T^+)}^2
		+
		\varepsilon
		\|f^\varepsilon\|_{L^2(\cD_T;H_\sigma^1)}^2
		\notag\\
		&\qquad\le
		\|f_0\|_{L^2(\cD\times\bR_v^3)}^2
		+
		\|h\|_{L^2(\gamma_T^-)}^2
		+
		\frac{1}{\varepsilon}
		\|S\|_{L^2(\cD_T;H_\sigma^{-1})}^2.
	\end{align*}
	Since $\varepsilon\ge1/2$,
	\begin{align}
		\sup_{\frac12\le\varepsilon<1}
		\Big(
		&\|f^\varepsilon(T)\|_{L^2(\cD\times\bR_v^3)}^2
		+
		\|f_\gamma^\varepsilon\|_{L^2(\gamma_T^+)}^2
		+
		\|f^\varepsilon\|_{L^2(\cD_T;H_\sigma^1)}^2
		\Big)
		\notag\\
		&\qquad\le
		2\|f_0\|_{L^2(\cD\times\bR_v^3)}^2
		+
		2\|h\|_{L^2(\gamma_T^-)}^2
		+
		4\|S\|_{L^2(\cD_T;H_\sigma^{-1})}^2.
		\label{eq:model-epsilon-uniform-Hsigma-bound}
	\end{align}
	
	Since $f_\gamma^\varepsilon=h$ on $\gamma_T^-$, the family
	$\{f_\gamma^\varepsilon\}_{1/2\le\varepsilon<1}$ is bounded in
	$L^2(\gamma_T)$. Hence, by weak compactness, there exist a sequence
	$\varepsilon_k\uparrow1$ and
	\[
	f\in L^2(\cD_T;H_\sigma^1),
	\qquad
	f_{[T]}\in L^2(\cD\times\bR_v^3),
	\qquad
	f_{[\gamma]}\in L^2(\gamma_T)
	\]
	such that
	\begin{alignat}{2}
		f^{\varepsilon_k}
		&\rightharpoonup
		f\qquad 
		&&\text{weakly in }
		L^2(\cD_T;H_\sigma^1),
		\label{eq:model-epsilon-weak-f}
		\\
		f^{\varepsilon_k}(T)
		&\rightharpoonup
		f_{[T]}\qquad 
		&&\text{weakly in }
		L^2(\cD\times\bR_v^3),
		\notag\\
		f_\gamma^{\varepsilon_k}
		&\rightharpoonup
		f_{[\gamma]}\qquad 
		&&\text{weakly in }
		L^2(\gamma_T).
		\label{eq:model-epsilon-weak-trace}
	\end{alignat}
	Since the restriction map from $L^2(\gamma_T)$ to
	$L^2(\gamma_T^-)$ is bounded,
	\eqref{eq:model-epsilon-weak-trace} implies
	\begin{equation}\label{eq:model-limit-inflow-trace}
		f_{[\gamma]}=h
		\qquad
		\text{on }\gamma_T^-.
	\end{equation}
	
	Fix
	$\phi\in C_c^\infty\bigl([0,T]\times\overline\cD\times\bR_v^3\bigr)$.
	The kinetic Green identity for $f^{\varepsilon_k}$ reads
	\begin{align}
		&\iint_{\cD_T}
		\left\langle
		(\partial_t+v\cdot\nabla_x)f^{\varepsilon_k},
		\phi
		\right\rangle_v
		\dd x\dd t
		+
		\iint_{\cD_T}
		\left\langle
		f^{\varepsilon_k},
		(\partial_t+v\cdot\nabla_x)\phi
		\right\rangle_v
		\dd x\dd t
		\notag\\
		&\qquad=
		\iint_{\cD\times\bR_v^3}
		f^{\varepsilon_k}(T,x,v)\phi(T,x,v)
		\dd x\dd v
		-
		\iint_{\cD\times\bR_v^3}
		f_0(x,v)\phi(0,x,v)
		\dd x\dd v
		\notag\\
		&\qquad\quad
		+
		\int_0^T\int_\gamma
		f_\gamma^{\varepsilon_k}
		\phi|_\gamma
		\operatorname{sgn}(n_x\cdot v)
		\dd\gamma\dd t.
		\label{260718922}
	\end{align}
	By \eqref{eq:model-epsilon-inflow-problem},
	\begin{align*}
		&\iint_{\cD_T}
		\left\langle
		(\partial_t+v\cdot\nabla_x)f^{\varepsilon_k},
		\phi
		\right\rangle_v
		\dd x\dd t
		\notag\\
		&\qquad=
		\iint_{\cD_T}
		\left\langle
		S,
		\phi
		\right\rangle_{H_\sigma^{-1},H_\sigma^1}
		\dd x\dd t
		-
		\varepsilon_k
		\iint_{\cD_T}
		\left\langle
		A_0f^{\varepsilon_k},
		\phi
		\right\rangle_{H_\sigma^{-1},H_\sigma^1}
		\dd x\dd t
		\notag\\
		&\qquad\quad
		+
		(1-\varepsilon_k)
		\iiint_{\cD_T\times\bR_v^3}
		f^{\varepsilon_k}
		(\Delta_v-1)\phi
		\dd v\dd x\dd t.
	\end{align*}
	
	Since $A_0:H_\sigma^1\to H_\sigma^{-1}$ is bounded,
	\eqref{eq:model-epsilon-weak-f} gives
	\begin{equation*}
		A_0f^{\varepsilon_k}
		\rightharpoonup
		A_0f
		\qquad
		\text{weakly in }
		L^2(\cD_T;H_\sigma^{-1}).
	\end{equation*}
	Consequently,
	\begin{align*}
		&\varepsilon_k
		\iint_{\cD_T}
		\left\langle
		A_0f^{\varepsilon_k},
		\phi
		\right\rangle_{H_\sigma^{-1},H_\sigma^1}
		\dd x\dd t\longrightarrow
		\iint_{\cD_T}
		\left\langle
		A_0f,
		\phi
		\right\rangle_{H_\sigma^{-1},H_\sigma^1}
		\dd x\dd t.
	\end{align*}
	Moreover, since
	\[
	(\Delta_v-1)\phi
	\in
	L^2(\cD_T;H_\sigma^{-1}),
	\]
	\eqref{eq:model-epsilon-uniform-Hsigma-bound} yields
	\begin{align*}
		&(1-\varepsilon_k)
		\left|
		\iiint_{\cD_T\times\bR_v^3}
		f^{\varepsilon_k}
		(\Delta_v-1)\phi
		\dd v\dd x\dd t
		\right|
		\\
		&\qquad\le
		(1-\varepsilon_k)
		\|f^{\varepsilon_k}\|_{L^2(\cD_T;H_\sigma^1)}
		\|(\Delta_v-1)\phi\|_{L^2(\cD_T;H_\sigma^{-1})}
		\longrightarrow0.
	\end{align*}
	
	Passing to the limit in \eqref{260718922}, using
	\eqref{eq:model-epsilon-weak-f}--\eqref{eq:model-epsilon-weak-trace},
	we obtain
	\begin{align}
		&\iint_{\cD_T}
		\left\langle
		S-A_0f,
		\phi
		\right\rangle_{H_\sigma^{-1},H_\sigma^1}
		\dd x\dd t
		+
		\iint_{\cD_T}
		\left\langle
		f,
		(\partial_t+v\cdot\nabla_x)\phi
		\right\rangle_v
		\dd x\dd t
		\notag\\
		&\qquad=
		\iint_{\cD\times\bR_v^3}
		f_{[T]}(x,v)\phi(T,x,v)
		\dd x\dd v
		-
		\iint_{\cD\times\bR_v^3}
		f_0(x,v)\phi(0,x,v)
		\dd x\dd v
		\notag\\
		&\qquad\quad
		+
		\int_0^T\int_\gamma
		f_{[\gamma]}
		\phi|_\gamma
		\operatorname{sgn}(n_x\cdot v)
		\dd\gamma\dd t.
		\label{eq:model-limit-green-identity}
	\end{align}
	
	Since
	\[
	S-A_0f
	\in
	L^2(\cD_T;H_\sigma^{-1}),
	\]
	Proposition~\ref{H_kin=weak_landau}, applied to
	\eqref{eq:model-limit-green-identity}, gives
	\[
	f\in\cH_{\kin,\sigma}(\cD_T\times\bR_v^3)
	\]
	and
	\[
	(\partial_t+v\cdot\nabla_x)f=S-A_0f,
	\qquad
	f(0,\cdot)=f_0,
	\qquad
	f(T,\cdot)=f_{[T]},
	\qquad
	f_\gamma=f_{[\gamma]}.
	\]
	Together with \eqref{eq:model-limit-inflow-trace}, this proves that
	$f$ solves \eqref{FP_dy_sigma}. The terminal and outflow terms in
	\eqref{eq:model-landau-inflow-estimate} are those already produced by
	\eqref{eq:model-landau-energy-estimate}.
\end{proof}

\begin{lemma}[Time-dependent linearized Landau inflow problem]
	\label{260502321}
	Let $0<T<\infty$. For every
	\[
	S\in L^2(\cD_T;H_\sigma^{-1}),
	\qquad
	f_0\in L^2(\cD\times\bR_v^3),
	\qquad
	h\in L^2(\gamma_T^-),
	\]
	the problem (recall $L$ defined in \eqref{def:LAK})
	\begin{equation}\label{260322508_11}
		\begin{cases}
			(\partial_t+v\cdot\nabla_x)f+Lf=S
			&\text{in }\cD_T\times\bR_v^3,\\
			f(0,\cdot)=f_0
			&\text{in }\cD\times\bR_v^3,\\
			f=h
			&\text{on }\gamma_T^-
		\end{cases}
	\end{equation}
	admits a unique solution $
	f\in\cH_{\kin,\sigma}(\cD_T\times\bR_v^3)$, where $\cH_{\kin,\sigma}(\cD_T\times\bR_v^3)$ is defined in Definition~\ref{def:H_kin_sigma}.
	Moreover there exists a constant
	$C=C(\cD)>0$ such that
	\begin{align}
		&\|f\|_{\cH_{\kin,\sigma}(\cD_T\times\bR_v^3)}
		+
		\|f_\gamma\|_{L^2(\gamma_T^+)}
		\le
		C\left(
		\|S\|_{L^2(\cD_T;H_\sigma^{-1})}
		+
		\|f_0\|_{L^2(\cD\times\bR_v^3)}
		+
		\|h\|_{L^2(\gamma_T^-)}
		\right).
		\label{eq:linearized-landau-inflow-estimate}
	\end{align}
\end{lemma}

\begin{proof}
	\noindent\emph{Reduction to a uniform a priori estimate.}
	For $\theta\in[0,1]$, define
	\[
	L_\theta
	:=
	(1-\theta)L+\theta A_0.
	\]
	Since $L_1=A_0$, Lemma~\ref{260502320} gives unique solvability for
	$\theta=1$. Moreover, $\theta\mapsto L_\theta$ is continuous as a
	family of bounded operators from $H_\sigma^1$ to
	$H_\sigma^{-1}$. Thus, by the method of continuity, it suffices to
	establish an a priori estimate, uniform in $\theta\in[0,1]$, for
	solutions
	\[
	f\in\cH_{\kin,\sigma}(\cD_T\times\bR_v^3)
	\]
	of
	\begin{equation}\label{260715210}
		\begin{cases}
			(\partial_t+v\cdot\nabla_x)f+L_\theta f=S
			&\text{in }\cD_T\times\bR_v^3,\\
			f(0,\cdot)=f_0
			&\text{in }\cD\times\bR_v^3,\\
			f_\gamma=h
			&\text{on }\gamma_T^-.
		\end{cases}
	\end{equation}
	More precisely, we shall prove
	\begin{align*}
		\|f\|_{\cH_{\kin,\sigma}(\cD_T\times\bR_v^3)}
		\le
		C\Big(
		&\|S\|_{L^2(\cD_T;H_\sigma^{-1})}
		+
		\|f_0\|_{L^2(\cD\times\bR_v^3)}
		+
		\|h\|_{L^2(\gamma_T^-)}
		\Big),
	\end{align*}
	where $C=C(\cD)>0$ is independent of $\theta$, $T$, and the data.
	
	Applying the Green identity from
	Proposition~\ref{2607151218} and using
	$f_\gamma=h$ on $\gamma_T^-$, we obtain
	\begin{align}
		&\frac12
		\|f(T)\|_{L^2(\cD\times\bR_v^3)}^2
		-
		\frac12
		\|f_0\|_{L^2(\cD\times\bR_v^3)}^2
		+
		\frac12
		\iint_{\gamma_T}
		|f_\gamma|^2
		\operatorname{sgn}(n(x)\cdot v)
		\dd\gamma\dd t
		\notag\\
		&\qquad
		+
		\iint_{\cD_T}
		\langle L_\theta f,f\rangle_
		{H_\sigma^{-1},H_\sigma^1}
		\dd x\dd t
		=
		\iint_{\cD_T}
		\langle S,f\rangle_
		{H_\sigma^{-1},H_\sigma^1}
		\dd x\dd t.
		\label{260715203_1}
	\end{align}
	The boundary term is given by
	\begin{equation*}
		\iint_{\gamma_T}
		|f_\gamma|^2
		\operatorname{sgn}(n(x)\cdot v)
		\dd\gamma\dd t
		=
		\|f_\gamma\|_{L^2(\gamma_T^+)}^2
		-
		\|h\|_{L^2(\gamma_T^-)}^2.
	\end{equation*}
	
	\medskip
	\noindent\emph{A preliminary estimate for $\theta>0$.}
	Let $\theta\in(0,1]$. By \eqref{eq:L-coercivity} and
	\eqref{260715202},
	\begin{align}\label{260715203_3}
		\begin{aligned}
			&\iint_{\cD_T}
			\langle L_\theta f,f\rangle_
			{H_\sigma^{-1},H_\sigma^1}
			\dd x\dd t\\
			=\,&
			(1-\theta)
			\iint_{\cD_T}
			\langle Lf,f\rangle_
			{H_\sigma^{-1},H_\sigma^1}
			\dd x\dd t
			+
			\theta
			\iint_{\cD_T}
			\langle A_0f,f\rangle_
			{H_\sigma^{-1},H_\sigma^1}
			\dd x\dd t
			\\
			\ge\,&
			\theta
			\|f\|_{L^2(\cD_T;H_\sigma^1)}^2.
		\end{aligned}
	\end{align}
	Furthermore,
	\begin{align*}
		\left|
		\iint_{\cD_T}
		\langle S,f\rangle_
		{H_\sigma^{-1},H_\sigma^1}
		\dd x\dd t
		\right|
		&\le
		\|S\|_{L^2(\cD_T;H_\sigma^{-1})}
		\|f\|_{L^2(\cD_T;H_\sigma^1)}
		\\
		&\le
		\frac{1}{2\theta}
		\|S\|_{L^2(\cD_T;H_\sigma^{-1})}^2
		+
		\frac{\theta}{2}
		\|f\|_{L^2(\cD_T;H_\sigma^1)}^2.
	\end{align*}
	Combining
	\eqref{260715203_1}--\eqref{260715203_3}, we obtain
	\begin{align*}
		&\frac12
		\|f(T)\|_{L^2(\cD\times\bR_v^3)}^2
		+
		\frac12
		\|f_\gamma\|_{L^2(\gamma_T^+)}^2
		+
		\frac{\theta}{2}
		\|f\|_{L^2(\cD_T;H_\sigma^1)}^2
		\notag\\
		&\qquad\le
		\frac12
		\|f_0\|_{L^2(\cD\times\bR_v^3)}^2
		+
		\frac12
		\|h\|_{L^2(\gamma_T^-)}^2
		+
		\frac{1}{2\theta}
		\|S\|_{L^2(\cD_T;H_\sigma^{-1})}^2.
	\end{align*}
	Consequently,
	\begin{align}
		\|f\|_{L^2(\cD_T;H_\sigma^1)}
		&\le
		\frac{1}{\theta}
		\|S\|_{L^2(\cD_T;H_\sigma^{-1})} +
		\frac{1}{\sqrt{\theta}}
		\left(
		\|f_0\|_{L^2(\cD\times\bR_v^3)}^2
		+
		\|h\|_{L^2(\gamma_T^-)}^2
		\right)^{1/2}
		\notag\\
		&\le
		\frac{1}{\theta}
		\left(
		\|S\|_{L^2(\cD_T;H_\sigma^{-1})}
		+
		\|f_0\|_{L^2(\cD\times\bR_v^3)}
		+
		\|h\|_{L^2(\gamma_T^-)}
		\right).
		\label{260715324}
	\end{align}

	\medskip
	\noindent\emph{Uniform estimate for small $\theta$.}
	Throughout this part, $C_\cD>0$ denotes a constant depending only
	on $\cD$ and may vary from line to line.
	
	Equation \eqref{260715210} can be rewritten as
	\begin{equation}\label{260715302}
		\begin{cases}
			(\partial_t+v\cdot\nabla_x)f+Lf=S_\theta:=S+\theta(L-A_0)f
			&\text{in }\cD_T\times\bR_v^3,\\
			f(0,\cdot)=f_0
			&\text{in }\cD\times\bR_v^3,\\
			f_\gamma=h
			&\text{on }\gamma_T^-.
		\end{cases}
	\end{equation}
	
	Applying the Green identity from
	Proposition~\ref{2607151218} to \eqref{260715302}, and using the
	microscopic coercivity \eqref{eq:L-coercivity}, we obtain
	\begin{align}
		&\|f(T)\|_{L^2(\cD\times\bR_v^3)}^2
		+
		\|f_\gamma\|_{L^2(\gamma_T^+)}^2
		+
		\|(\mathbf I-\mathbf P)f\|_
		{L^2(\cD_T;H_\sigma^1)}^2
		\notag\\
		&\qquad\le
		C\left(
		\|f_0\|_{L^2(\cD\times\bR_v^3)}^2
		+
		\|h\|_{L^2(\gamma_T^-)}^2
		+
		\left|
		\iint_{\cD_T}
		\langle S_\theta,f\rangle_
		{H_\sigma^{-1},H_\sigma^1}
		\dd x\dd t
		\right|
		\right),
		\label{260715248}
	\end{align}
	where $C>0$ is universal.
	
	Applying
	\eqref{eq:appendix-macro-estimate-endpoint} with $s=0$ and $t=T$,
	and using
	\eqref{eq:macro-eta-bound} and
	\eqref{eq:appendix-source-M-time-Hminus1}, we obtain
	\begin{align}
		\|\mathbf P f\|_{L^2(\cD_T\times\bR_v^3)}^2
		\le
		C_\cD\Big(
		&\|\mu^{1/4}(\mathbf I-\mathbf P)f\|_
		{L^2(\cD_T\times\bR_v^3)}^2
		+
		\|S_\theta\|_{L^2(\cD_T;H_\sigma^{-1})}^2
		\notag\\
		&+
		\|f_\gamma\|_{L^2(\gamma_T^+)}^2
		+
		\|h\|_{L^2(\gamma_T^-)}^2
		+
		\|f(T)\|_{L^2(\cD\times\bR_v^3)}^2
		+
		\|f_0\|_{L^2(\cD\times\bR_v^3)}^2
		\Big).
		\label{260715249}
	\end{align}
	Indeed,
	\begin{align*}
		|\eta_f(T)-\eta_f(0)|
		&\le
		C\left(
		\|f(T)\|_{\cM}
		+
		\|f_0\|_{\cM}
		\right)
		\|f(T)-f_0\|_{\cM}
		\\
		&\le
		C\left(
		\|f(T)\|_{L^2(\cD\times\bR_v^3)}^2
		+
		\|f_0\|_{L^2(\cD\times\bR_v^3)}^2
		\right).
	\end{align*}
	Moreover,
	\begin{equation}\label{eq:micro-Maxwellian-control}
		\|\mu^{1/4}(\mathbf I-\mathbf P)f\|_
		{L^2(\cD_T\times\bR_v^3)}
		\le
		C
		\|(\mathbf I-\mathbf P)f\|_
		{L^2(\cD_T;H_\sigma^1)}.
	\end{equation}
	
	Since the range of $\mathbf P$ is finite-dimensional, the
	$L_v^2$ and $H_\sigma^1$ norms are equivalent on
	the range of $\mathbf P$. Hence
	\begin{equation}\label{eq:macro-Hsigma-equivalence}
		\|\mathbf P f\|_{L^2(\cD_T;H_\sigma^1)}
		\le
		C
		\|\mathbf P f\|_{L^2(\cD_T\times\bR_v^3)}.
	\end{equation}
	Combining
	\eqref{260715248}, \eqref{260715249},
	\eqref{eq:micro-Maxwellian-control}, and
	\eqref{eq:macro-Hsigma-equivalence}, we obtain
	\begin{align}
		\|f\|_{L^2(\cD_T;H_\sigma^1)}^2
		\le
		C_\cD\Big(
		&\|f_0\|_{L^2(\cD\times\bR_v^3)}^2
		+
		\|h\|_{L^2(\gamma_T^-)}^2
		+
		\|S_\theta\|_{L^2(\cD_T;H_\sigma^{-1})}^2
		\notag\\
		&+
		\left|
		\iint_{\cD_T}
		\langle S_\theta,f\rangle_
		{H_\sigma^{-1},H_\sigma^1}
		\dd x\dd t
		\right|
		\Big).
		\label{260715248_11}
	\end{align}
	
	For every $\delta>0$,
	\begin{align*}
		\left|
		\iint_{\cD_T}
		\langle S_\theta,f\rangle_
		{H_\sigma^{-1},H_\sigma^1}
		\dd x\dd t
		\right|
		&\le
		\|S_\theta\|_{L^2(\cD_T;H_\sigma^{-1})}
		\|f\|_{L^2(\cD_T;H_\sigma^1)}
		\\
		&\le
		\delta
		\|f\|_{L^2(\cD_T;H_\sigma^1)}^2
		+
		\frac{1}{4\delta}
		\|S_\theta\|_{L^2(\cD_T;H_\sigma^{-1})}^2.
	\end{align*}
	Choosing $\delta>0$ so that $C_\cD\delta\le1/2$ for $C_\cD$ in 
	\eqref{260715248_11}, we obtain
	\begin{equation}\label{260715321}
		\|f\|_{L^2(\cD_T;H_\sigma^1)}
		\le
		C_2\left(
		\|S_\theta\|_{L^2(\cD_T;H_\sigma^{-1})}
		+
		\|f_0\|_{L^2(\cD\times\bR_v^3)}
		+
		\|h\|_{L^2(\gamma_T^-)}
		\right),
	\end{equation}
	where $C_2=C_2(\cD)>0$ is independent of $\theta$, $T$, and the
	data.
	
	By the definition of $S_\theta$,
	\begin{align*}
		\|S_\theta\|_{L^2(\cD_T;H_\sigma^{-1})}
		&\le
		\|S\|_{L^2(\cD_T;H_\sigma^{-1})}
		+
		\theta
		\|(L-A_0)f\|_{L^2(\cD_T;H_\sigma^{-1})}.
	\end{align*}
	By \eqref{boundedness of L} and
	\eqref{eq:model-landau-isometry}, there exists a universal constant
	$C_3>0$ such that
	\begin{equation*}
		\|(L-A_0)f\|_{L^2(\cD_T;H_\sigma^{-1})}
		\le
		C_3
		\|f\|_{L^2(\cD_T;H_\sigma^1)}.
	\end{equation*}
	It follows from \eqref{260715321} that
	\begin{align*}
		\|f\|_{L^2(\cD_T;H_\sigma^1)}
		\le{}&
		C_2\left(
		\|S\|_{L^2(\cD_T;H_\sigma^{-1})}
		+
		\|f_0\|_{L^2(\cD\times\bR_v^3)}
		+
		\|h\|_{L^2(\gamma_T^-)}
		\right)
		+
		C_2C_3\theta
		\|f\|_{L^2(\cD_T;H_\sigma^1)}.
	\end{align*}
	Set
	\[
	\theta_*
	:=
	\min\left\{
	1,\frac{1}{2C_2C_3}
	\right\}.
	\]
	Then $\theta_*=\theta_*(\cD)>0$, and for every
	$\theta\in[0,\theta_*]$,
	\begin{equation}\label{260715325}
		\|f\|_{L^2(\cD_T;H_\sigma^1)}
		\le
		2C_2\left(
		\|S\|_{L^2(\cD_T;H_\sigma^{-1})}
		+
		\|f_0\|_{L^2(\cD\times\bR_v^3)}
		+
		\|h\|_{L^2(\gamma_T^-)}
		\right).
	\end{equation}

	\medskip
	\noindent\emph{Conclusion.}
	For $\theta\in(\theta_*,1]$, estimate \eqref{260715324} gives
	\begin{align*}
		\|f\|_{L^2(\cD_T;H_\sigma^1)}
		&\le
		\frac{1}{\theta_*}
		\left(
		\|S\|_{L^2(\cD_T;H_\sigma^{-1})}
		+
		\|f_0\|_{L^2(\cD\times\bR_v^3)}
		+
		\|h\|_{L^2(\gamma_T^-)}
		\right).
	\end{align*}
	Combining this with \eqref{260715325}, we obtain
	\begin{equation}\label{260715328}
		\|f\|_{L^2(\cD_T;H_\sigma^1)}
		\le
		C_4
		\left(
		\|S\|_{L^2(\cD_T;H_\sigma^{-1})}
		+
		\|f_0\|_{L^2(\cD\times\bR_v^3)}
		+
		\|h\|_{L^2(\gamma_T^-)}
		\right)
	\end{equation}
	for every $\theta\in[0,1]$, where
	\[
	C_4
	:=
	\max\left\{
	2C_2,\frac{1}{\theta_*}
	\right\}.
	\]
	In particular, $C_4$ depends only on $\cD$.
	
	By \eqref{260715210}, \eqref{boundedness of L}, and
	\eqref{eq:model-landau-isometry},
	\begin{align*}			
		\begin{aligned}
			&\|(\partial_t+v\cdot\nabla_x)f\|_
			{L^2(\cD_T;H_\sigma^{-1})}\\
			\le\,&
			\|S\|_{L^2(\cD_T;H_\sigma^{-1})}
			+
			(1-\theta)
			\|Lf\|_{L^2(\cD_T;H_\sigma^{-1})}
			+
			\theta
			\|A_0f\|_{L^2(\cD_T;H_\sigma^{-1})}
			\\
			\le\,&
			\|S\|_{L^2(\cD_T;H_\sigma^{-1})}
			+
			C
			\|f\|_{L^2(\cD_T;H_\sigma^1)},				
		\end{aligned}
	\end{align*}
	where $C>0$ is independent of $\theta$, $T$, and the data.
	Therefore, \eqref{260715328}, together with the initial and inflow
	boundary conditions, yields
	\begin{equation}\label{eq:Ltheta-uniform-estimate-final}
		\|f\|_{\cH_{\kin,\sigma}(\cD_T\times\bR_v^3)}
		\le
		C_\cD
		\left(
		\|S\|_{L^2(\cD_T;H_\sigma^{-1})}
		+
		\|f_0\|_{L^2(\cD\times\bR_v^3)}
		+
		\|h\|_{L^2(\gamma_T^-)}
		\right),
	\end{equation}
	where $C_\cD>0$ depends only on $\cD$.
	
	Thus the required a priori estimate is uniform in
	$\theta\in[0,1]$. The method of continuity yields the unique
	solvability of \eqref{260322508_11}, and
	\eqref{eq:Ltheta-uniform-estimate-final} proves the claimed estimate. The
	outflow term in \eqref{eq:linearized-landau-inflow-estimate} is supplied by
	\eqref{260730346}.
\end{proof}

\begin{proof}[\textbf{Proof of Theorem~\ref{260322442}}]
	As mentioned at the beginning of this subsection, we prove only (1).
	
	For $\theta\in[0,1]$, set
	\[
	L_{g,\theta}
	:=
	L-\theta N[g,\cdot].
	\]
	Then
	\[
	L_{g,0}=L,
	\qquad
	L_{g,1}=L_g.
	\]
	
	Let $C_1=C_1(\cD)>0$ be the constant in Lemma~\ref{260502321},
	and let $C_2>0$ be the constant in \eqref{260911445}. Set
	\begin{equation*}
		\varepsilon_0
		:=
		\min\left\{\varepsilon_{\mathrm{ell}},\
		\frac{1}{2C_1C_2}\right\},
	\end{equation*}
	which is the constant $\varepsilon_0=\varepsilon_0(\cD)$ announced in
	Theorem~\ref{260322442}, and assume that
	\begin{equation}\label{eq:Ng-continuity-smallness}
		\|\mu^{1/4}g\|_
		{L^\infty(\cD_T\times\bR_v^3)}
		\le
		\varepsilon_0.
	\end{equation}
	The minimum with $\varepsilon_{\mathrm{ell}}$ makes
	\eqref{eq:smallness-perturbed-sigma}, hence \eqref{260502449}, available
	under \eqref{eq:Ng-continuity-smallness}.
	
	Let $\theta\in[0,1]$, and suppose that
	$f\in\cH_{\kin,\sigma}(\cD_T\times\bR_v^3)$ satisfies
	\begin{equation*}
		\begin{cases}
			(\partial_t+v\cdot\nabla_x)f+L_{g,\theta}f=S
			&\text{in }\cD_T\times\bR_v^3,\\
			f(0,\cdot)=f_0
			&\text{in }\cD\times\bR_v^3,\\
			f=h
			&\text{on }\gamma_T^-.
		\end{cases}
	\end{equation*}
	Equivalently,
	\[
	(\partial_t+v\cdot\nabla_x)f+Lf
	=
	S+\theta N[g,f].
	\]
	By \eqref{260911445},
	\begin{equation}\label{eq:Ng-spacetime-bound}
		\|N[g,f]\|_{L^2(\cD_T;H_\sigma^{-1})}
		\le
		C_2
		\|\mu^{1/4}g\|_
		{L^\infty(\cD_T\times\bR_v^3)}
		\|f\|_{L^2(\cD_T;H_\sigma^1)}.
	\end{equation}
	Applying Lemma~\ref{260502321} and using \eqref{eq:Ng-spacetime-bound}, we
	obtain
	\begin{align*}
		\|f\|_{\cH_{\kin,\sigma}(\cD_T\times\bR_v^3)}
		\le{}&
		C_1\left(
		\|S\|_{L^2(\cD_T;H_\sigma^{-1})}
		+
		\|f_0\|_{L^2(\cD\times\bR_v^3)}
		+
		\|h\|_{L^2(\gamma_T^-)}
		\right)
		\\
		&\qquad +
		\theta C_1C_2
		\|\mu^{1/4}g\|_
		{L^\infty(\cD_T\times\bR_v^3)}
		\|f\|_{L^2(\cD_T;H_\sigma^1)}.
	\end{align*}
	Since
	\[
	\|f\|_{L^2(\cD_T;H_\sigma^1)}
	\le
	\|f\|_{\cH_{\kin,\sigma}(\cD_T\times\bR_v^3)},
	\]
	and since \eqref{eq:Ng-continuity-smallness} gives
	$\theta C_1C_2\|\mu^{1/4}g\|_{L^\infty(\cD_T\times\bR_v^3)}\le\frac12$,
	the last term is absorbed into the left-hand side. Hence
	\begin{equation}\label{eq:Lg-inflow-final-estimate}
		\|f\|_{\cH_{\kin,\sigma}(\cD_T\times\bR_v^3)}
		\le
		2C_1\left(
		\|S\|_{L^2(\cD_T;H_\sigma^{-1})}
		+
		\|f_0\|_{L^2(\cD\times\bR_v^3)}
		+
		\|h\|_{L^2(\gamma_T^-)}
		\right),
	\end{equation}
	uniformly in $\theta\in[0,1]$.
	
	Since $L_{g,0}=L$, the problem at $\theta=0$ is uniquely solvable by
	Lemma~\ref{260502321}. The method of continuity and
	\eqref{eq:Lg-inflow-final-estimate} therefore yield unique solvability at
	$\theta=1$, together with
	\eqref{eq:perturbed-landau-dynamic-estimate}. This proves (1).
\end{proof}

We conclude this subsection by extending the preceding a priori estimates to
polynomial velocity weights for both the time-dependent and stationary
perturbed Landau inflow problems.

\begin{thm}[Weighted estimates for perturbed Landau inflow problems]
	\label{thm:weighted-steady}
	Let $\varepsilon_0=\varepsilon_0(\cD)>0$ be the smallness constant
	in Theorem~\ref{260322442}.
	For each real $\ell \ge0$, there exists a constant
	$C_\ell=C(\cD,\ell)>0$ such that the following assertions hold.
	
	\begin{enumerate}
		\item
		Let $0<T<\infty$, and assume that
		\begin{equation}
			\label{eq:weighted-dynamic-coefficient}
			\|\mu^{1/4}g\|_
			{L^\infty(\cD_T\times\bR_v^3)}
			\le
			\varepsilon_0.
		\end{equation}
		Let
		\[
		w^\ell S\in L^2(\cD_T;H_\sigma^{-1}),
		\qquad
		w^\ell f_0\in L^2(\cD\times\bR_v^3),
		\qquad
		w^\ell h\in L^2(\gamma_T^-),
		\]
		and let $f\in
		\cH_{\kin,\sigma}(\cD_T\times\bR_v^3)$
		satisfy \eqref{solv:linland_dy}.
		Then
		\[
		w^\ell f
		\in
		\cH_{\kin,\sigma}(\cD_T\times\bR_v^3),
		\]
		and
		\begin{align}
			\|w^\ell f\|_
			{\cH_{\kin,\sigma}(\cD_T\times\bR_v^3)}
			\le
			C_\ell \Big(
			&\|w^\ell S\|_{L^2(\cD_T;H_\sigma^{-1})}
			+
			\|w^\ell f_0\|_{L^2(\cD\times\bR_v^3)}
			+
			\|w^\ell h\|_{L^2(\gamma_T^-)}
			\Big).
			\label{eq:weighted-dynamic-estimate}
		\end{align}
		In particular,
		\begin{align}
			\|w^{\ell}f\|_
			{L^2(\cD_T\times\bR_v^3)}
			\le
			C_\ell \Big(
			&\|w^\ell S\|_{L^2(\cD_T;H_\sigma^{-1})}
			+
			\|w^\ell f_0\|_{L^2(\cD\times\bR_v^3)}
			+
			\|w^\ell h\|_{L^2(\gamma_T^-)}
			\Big).
			\label{eq:weighted-dynamic-L2-estimate}
		\end{align}
		
		\item
		Assume that
		\[
		\|\mu^{1/4}g\|_
		{L^\infty(\cD\times\bR_v^3)}
		\le
		\varepsilon_0.
		\]
		Let
		\[
		w^\ell S\in L^2(\cD;H_\sigma^{-1}),
		\qquad
		w^\ell h\in L^2(\gamma^-),
		\]
		and let 
		$f\in
		\cH_{\kin,\sigma}(\cD\times\bR_v^3)$
		satisfy \eqref{solv:linland_st}.
		Then
		\[
		w^\ell f
		\in
		\cH_{\kin,\sigma}(\cD\times\bR_v^3),
		\]
		and
		\begin{equation}\label{eq:weighted-steady-estimate}
			\|w^\ell f\|_
			{\cH_{\kin,\sigma}(\cD\times\bR_v^3)}
			\le
			C_\ell \left(
			\|w^\ell S\|_{L^2(\cD;H_\sigma^{-1})}
			+
			\|w^\ell h\|_{L^2(\gamma^-)}
			\right).
		\end{equation}
		In particular,
		\begin{equation}\label{eq:weighted-steady-L2-estimate}
			\|w^{\ell}f\|_
			{L^2(\cD\times\bR_v^3)}
			\le
			C_\ell \left(
			\|w^\ell S\|_{L^2(\cD;H_\sigma^{-1})}
			+
			\|w^\ell h\|_{L^2(\gamma^-)}
			\right).
		\end{equation}
	\end{enumerate}
	
	The constant $C_\ell$ is independent of $T$, $g$, and the data
	under the stated smallness condition on $g$.
\end{thm}

\begin{proof}
	We prove the time-dependent assertion. The stationary assertion
	follows from the same argument, using the stationary part of
	Theorem~\ref{260322442}.
	
	For real $k\ge0$, set
	\[
	\mathfrak S_k
	:=
	\|w^kS\|_{L^2(\cD_T;H_\sigma^{-1})}
	+
	\|w^kf_0\|_{L^2(\cD\times\bR_v^3)}
	+
	\|w^kh\|_{L^2(\gamma_T^-)}.
	\]
	We argue by induction on $\lceil\ell/2\rceil$. The statement at
	level $\ell$ is deduced from the statement at level
	$k_\ell:=(\ell-2)_+<\ell$, and $\lceil k_\ell/2\rceil<\lceil\ell/2\rceil$
	for every $\ell>0$. Nothing below uses that $\ell$ is an integer.
	
	For $\ell=0$, Theorem~\ref{260322442} and
	\eqref{eq:local-H1-H_sigma} give
	\begin{equation*}
		\|f\|_{\cH_{\kin,\sigma}(\cD_T\times\bR_v^3)}
		+
		\|w^{-1/2}f\|_{L^2(\cD_T\times\bR_v^3)}
		\le
		C_0\mathfrak S_0.
	\end{equation*}
	
	Let $\ell>0$, and suppose that
	\eqref{eq:weighted-dynamic-estimate} and
	\eqref{eq:weighted-dynamic-L2-estimate} hold at the level
	$k_\ell$ defined by
	\[
	k_\ell:=(\ell-2)_+.
	\]
	For $0\le k\le\ell$, the derivative bounds
	\[
	|D^\alpha w^{k-\ell}|
	\le
	C_\ell w^{k-\ell-|\alpha|},
	\qquad
	|\alpha|\le1,
	\]
	show that multiplication by $w^{k-\ell}$ is bounded on
	$H_\sigma^1$ and, by duality, on $H_\sigma^{-1}$. Hence
	\begin{equation*}
		\mathfrak S_k
		\le
		C_\ell \mathfrak S_\ell,
		\qquad
		0\le k\le\ell.
	\end{equation*}
	Since
	\[
	\ell-\frac52
	\le
	k_\ell-\frac12,
	\]
	the induction hypothesis at level $k_\ell$ yields
	\begin{align}
		\|w^{\ell-5/2}f\|_
		{L^2(\cD_T\times\bR_v^3)}
		&\le
		\|w^{k_\ell-1/2}f\|_
		{L^2(\cD_T\times\bR_v^3)}
		\notag\\
		&\le
		C_\ell \mathfrak S_{k_\ell}
		\le
		C_\ell \mathfrak S_\ell.
		\label{eq:weighted-induction-L2-control}
	\end{align}
	
	Choose a radial function $\chi\in C_c^\infty(\bR^3)$ such that
	\[
	0\le\chi\le1,
	\qquad
	\chi=1
	\quad\text{on }B_1,
	\qquad
	\supp\chi\subset B_2.
	\]
	For $R\ge1$, set
	\[
	\chi_R(v):=\chi(v/R),
	\qquad
	W_R(v):=w(v)^\ell \chi_R(v),
	\qquad
	f_R:=W_Rf.
	\]
	The radial derivative bounds for $\chi_R$ imply
	\begin{equation}\label{eq:dynamic-cutoff-multiplier-bound}
		\sup_{R\ge1}
		\left(
		\|\chi_R\|_{\cL(H_\sigma^1)}
		+
		\|\chi_R\|_{\cL(H_\sigma^{-1})}
		\right)
		<\infty.
	\end{equation}
	Consequently,
	\begin{align}
		\|W_RS\|_{L^2(\cD_T;H_\sigma^{-1})}
		&\le
		C
		\|w^\ell S\|_{L^2(\cD_T;H_\sigma^{-1})},
		\notag\\
		\|W_Rf_0\|_{L^2(\cD\times\bR_v^3)}
		&\le
		\|w^\ell f_0\|_{L^2(\cD\times\bR_v^3)},
		\notag\\
		\|W_Rh\|_{L^2(\gamma_T^-)}
		&\le
		\|w^\ell h\|_{L^2(\gamma_T^-)}.
		\label{eq:dynamic-truncated-data-bound}
	\end{align}
	
	We use the convention
	\[
	[B,W]f:=B(Wf)-W(Bf).
	\]
	For a radial multiplier $W$, the symmetry of $a_g^{ij}$ gives
	\begin{equation*}
		[\cA_g,W]f
		=
		2D_j\left(a_g^{ij}(D_iW)f\right)
		+
		Q_{g,W}f,
	\end{equation*}
	where
	\[
	Q_{g,W}
	:=
	b_g^iD_iW
	-
	(D_ja_g^{ij})D_iW
	-
	a_g^{ij}D_{ij}W.
	\]
	By Corollary~\ref{cor:expansion-N},
	Lemma~\ref{lem:perturbed-landau-ellipticity}, and
	\eqref{260507202},
	\begin{align*}
		|a_g^{ij}\xi_i\overline{\zeta_j}|
		&\le
		C
		\mathfrak q(v,\xi)^{1/2}
		\mathfrak q(v,\zeta)^{1/2},
		\notag\\
		|Da_g|
		&\le
		Cw^{-2},
		\notag\\
		|v_ib_g^i|
		+
		|\operatorname{tr}a_g|
		+
		|a_g^{ij}v_iv_j|
		&\le
		Cw^{-1},
	\end{align*}
	where $C>0$ is independent of $g$ under the stated smallness
	condition.
	
	On $\{R\le|v|\le2R\}$,
	\[
	|D^m\chi_R|
	\le
	C_mR^{-m},
	\qquad
	m=1,2.
	\]
	Since $W_R$ is radial, it follows that
	\begin{align*}
		\mathfrak q(v,DW_R)
		&\le
		C_\ell w^{2\ell-5},
		\notag\\
		|Q_{g,W_R}|
		&\le
		C_\ell w^{\ell-3},
	\end{align*}
	uniformly in $R\ge1$ and in $g$ satisfying the stated smallness
	condition.
	
	Therefore, for every $\varphi\in C_c^\infty(\bR_v^3)$,
	\begin{align*}
		\left|
		\left\langle
		[\cA_g,W_R]f,\varphi
		\right\rangle_{H_\sigma^{-1},H_\sigma^1}
		\right|
		&\le
		2\left|
		\int_{\bR^3}
		a_g^{ij}(D_iW_R)fD_j\varphi
		\dd v
		\right|
		+
		\left|
		\int_{\bR^3}
		Q_{g,W_R}f\varphi
		\dd v
		\right|
		\\
		&\le
		C_\ell
		\|w^{\ell-5/2}f\|_{L^2(\bR_v^3)}
		\|\varphi\|_{H_\sigma^1},
	\end{align*}
	where \eqref{eq:local-H1-H_sigma} is used for the zero-order term.
	Hence
	\begin{equation*}
		\|[\cA_g,W_R]f\|_{H_\sigma^{-1}}
		\le
		C_\ell
		\|w^{\ell-5/2}f\|_{L^2(\bR_v^3)}.
	\end{equation*}
	
	For the compact part, Lemma~\ref{2607161158} and
	\eqref{eq:dynamic-cutoff-multiplier-bound} give
	\begin{align*}
		\|[\cK,W_R]f\|_{H_\sigma^{-1}}
		&\le
		\|\cK(W_Rf)\|_{H_\sigma^{-1}}
		+
		\|\chi_R(w^\ell \cK f)\|_{H_\sigma^{-1}}
		\\
		&\le
		C_\ell
		\left(
		\|\mu^{1/4}W_Rf\|_{L^2(\bR_v^3)}
		+
		\|\mu^{1/4}f\|_{L^2(\bR_v^3)}
		\right)
		\\
		&\le
		C_\ell
		\|w^{\ell-5/2}f\|_{L^2(\bR_v^3)}.
	\end{align*}
	Here we used
	\[
	\mu^{1/4}(w^\ell+1)
	\le
	C_\ell w^{\ell-5/2},
	\qquad
	\ell \ge0,
	\]
	which holds since
	$\mu^{1/4}(w^\ell+1)w^{5/2-\ell}=\mu^{1/4}\big(w^{5/2}+w^{5/2-\ell}\big)$
	is bounded for each fixed $\ell \ge0$.
	Since $L_g=-\cA_g-\cK$, we conclude that
	\begin{equation}\label{eq:dynamic-truncated-commutator}
		\|[L_g,W_R]f\|_{H_\sigma^{-1}}
		\le
		C_\ell
		\|w^{\ell-5/2}f\|_{L^2(\bR_v^3)},
	\end{equation}
	uniformly in $R\ge1$.
	
	Since $W_R$ is independent of $(t,x)$, $f_R$ satisfies
	\begin{equation}\label{eq:dynamic-truncated-problem}
		\begin{cases}
			(\partial_t+v\cdot\nabla_x)f_R+L_gf_R=S_R
			&\text{in }\cD_T\times\bR_v^3,\\
			f_R(0,\cdot)=W_Rf_0
			&\text{in }\cD\times\bR_v^3,\\
			f_R=W_Rh
			&\text{on }\gamma_T^-,
		\end{cases}
	\end{equation}
	where
	\[
	S_R:=W_RS+[L_g,W_R]f.
	\]
	By
	\eqref{eq:dynamic-truncated-data-bound},
	\eqref{eq:dynamic-truncated-commutator}, and
	\eqref{eq:weighted-induction-L2-control},
	\begin{equation*}
		\|S_R\|_{L^2(\cD_T;H_\sigma^{-1})}
		\le
		C_\ell \mathfrak S_\ell.
	\end{equation*}
	Since $W_R\in C_c^\infty(\bR_v^3)$, Lemma~\ref{lem:compact-velocity-support}.(1)
	gives $f_R=W_Rf\in\cH_{\kin,\sigma}(\cD_T\times\bR_v^3)$. Hence $f_R$ is
	\emph{the} solution of \eqref{eq:dynamic-truncated-problem} provided by
	Theorem~\ref{260322442}.
	Applying Theorem~\ref{260322442} to
	\eqref{eq:dynamic-truncated-problem}, we obtain
	\begin{equation}\label{eq:dynamic-truncated-uniform-bound}
		\sup_{R\ge1}
		\|W_Rf\|_{\cH_{\kin,\sigma}(\cD_T\times\bR_v^3)}
		\le
		C_\ell \mathfrak S_\ell.
	\end{equation}
	
	By weak compactness, there exist $R_j\uparrow\infty$ and
	\[
	F\in\cH_{\kin,\sigma}(\cD_T\times\bR_v^3)
	\]
	such that
	\[
	W_{R_j}f
	\rightharpoonup
	F
	\quad\text{weakly in }
	\cH_{\kin,\sigma}(\cD_T\times\bR_v^3).
	\]
	Since $W_{R_j}f\to w^\ell f$ in
	$\cD'(\cD_T\times\bR_v^3)$, we have
	\[
	F=w^\ell f.
	\]
	Moreover,
	\[
	W_{R_j}f_0
	\longrightarrow
	w^\ell f_0
	\quad\text{in }L^2(\cD\times\bR_v^3),
	\]
	and
	\[
	W_{R_j}h
	\longrightarrow
	w^\ell h
	\quad\text{in }L^2(\gamma_T^-).
	\]
	Thus the initial and inflow traces of $F=w^\ell f$ are
	$w^\ell f_0$ and $w^\ell h$, respectively. By weak lower
	semicontinuity and \eqref{eq:dynamic-truncated-uniform-bound},
	\begin{equation*}
		\|w^\ell f\|_
		{\cH_{\kin,\sigma}(\cD_T\times\bR_v^3)}
		\le
		C_\ell \mathfrak S_\ell.
	\end{equation*}
	This proves \eqref{eq:weighted-dynamic-estimate}.
	
	Finally,
	\eqref{260829706} gives
	\[
	\|w^\ell f\|_{L^2(\cD_T\times\bR_v^3)}
	\le C\|w^\ell f\|_{\cH_{\kin,\sigma}(\cD_T\times\bR_v^3)}
	\le C_\ell\mathfrak S_\ell.
	\]
	This proves \eqref{eq:weighted-dynamic-L2-estimate} and completes
	the induction.
\end{proof}

\section{Grazing Regularity for Linearized Landau Equations}

We apply \S\ref{2606261120} to the linearized Landau equation, obtaining
a weighted local-in-time $L^2$--$L^\infty$ estimate and, after boundary
flattening and diffusion normalization, kinetic H\"older and velocity-gradient
estimates.

\subsection{Weighted $L^2$--$L^\infty$ estimate}

Set $w:=\langle v\rangle$.  Throughout this subsection $g$ is measurable
in $(t,x,v)$ and satisfies \eqref{eq:smallness-perturbed-sigma} for almost
every $(t,x)$; the variables $(t,x)$ enter only as parameters, and every
bound below is pointwise in $(t,x)$ and depends on $g$ only through the
bound 
$$
\varepsilon_{\mathrm{ell}}
$$
in \eqref{eq:smallness-perturbed-sigma}.
For $\ell \ge0$, conjugation by $w^\ell$
gives the operators
\begin{align}\label{eq:weighted-Ag-expansion}
	\begin{aligned}
		\cA_{g,\ell}u
		&:=
		w^\ell \cA_g(w^{-\ell}u)
		=
		D_i\left(a_{g,\ell}^{ij}D_ju\right)
		+
		b_{g,\ell}^iD_iu
		+
		c_{g,\ell}u,\qquad \cK_\ell u
		:=
		w^\ell \cK(w^{-\ell}u),
	\end{aligned}
\end{align}
and
\begin{align}
	a_{g,\ell}^{ij}
	&:=
	a_g^{ij},
	\label{eq:weighted-ag-coefficient}
	\\
	b_{g,\ell}^i
	&:=
	b_g^i
	-
	2\ell a_g^{ij}\frac{v_j}{w^2},
	\notag\\
	c_{g,\ell}
	&:=
	c_g
	-
	\ell \frac{(D_ia_g^{ij})v_j}{w^2}
	-
	\ell \frac{\operatorname{tr}a_g}{w^2}
	-
	\ell \frac{b_g^iv_i}{w^2}
	+
	\ell(\ell+2)
	\frac{a_g^{ij}v_iv_j}{w^4}.
	\label{eq:weighted-cg-expanded}
\end{align}

By \eqref{eq:weighted-ag-coefficient} and the definition of $a_g$ in
\eqref{260507202}, $a_{g,\ell}=\sigma_{\mu+\sqrt\mu g}$ is symmetric,
and \eqref{260502449} gives $a_{g,\ell}\ge\frac12\sigma>0$ under
\eqref{eq:smallness-perturbed-sigma}.  The Cauchy--Schwarz inequality for
the nonnegative form $a_{g,\ell}$ is therefore available, and combining
it with \eqref{260502449} and
Lemma~\ref{lem:landau-diffusion-eigenvalues} (applied with
$\varepsilon=1$), we obtain universal constants $c,C>0$
such that, for all $v,\xi,\zeta\in\bR^3$,
\begin{equation}\label{eq:weighted-ag-ellipticity}
	\begin{aligned}
		c\,\mathfrak q(v,\xi)
		&\le
		a_{g,\ell}^{ij}\xi_i\xi_j
		\le
		C\mathfrak q(v,\xi),
		\\
		\left|a_{g,\ell}^{ij}\xi_i\zeta_j\right|
		&\le
		C\mathfrak q(v,\xi)^{1/2}
		\mathfrak q(v,\zeta)^{1/2}.
	\end{aligned}
\end{equation}

Writing $(a_{g,\ell}v)_i:=a_{g,\ell}^{ij}v_j$, the coefficient
definitions \eqref{260507202} and
\eqref{eq:weighted-ag-coefficient}--\eqref{eq:weighted-cg-expanded},
together with \eqref{eq:weighted-ag-ellipticity},
\eqref{eq:derivative-sigma-bound},
\eqref{eq:perturbed-sigma-i-bound},
\eqref{eq:derivative-perturbed-sigma-bound},
\eqref{eq:perturbation-sigma-bound}, and
\eqref{eq:smallness-perturbed-sigma}, give
\begin{equation}\label{eq:weighted-coefficient-pointwise-bounds}
	\begin{aligned}
		|\nabla_va_{g,\ell}|
		\le Cw^{-2},\qquad
		|a_{g,\ell}v|\le Cw^{-1},\qquad
		|b_{g,\ell}|+|c_{g,\ell}|
		\le C_\ell w^{-1}.
	\end{aligned}
\end{equation}
Furthermore, \eqref{eq:compact-part-pointwise} yields
\begin{equation}\label{eq:weighted-compact-part-pointwise}
	|\cK_\ell u(v)|
	\le
	C_\ell \mu(v)^{1/8}
	\|\mu^{1/8}u\|_{L^2(\bR_v^3)}.
\end{equation}
In \eqref{eq:weighted-coefficient-pointwise-bounds} and
\eqref{eq:weighted-compact-part-pointwise}, $C>0$ is universal and
$C_\ell>0$ depends only on $\ell$; both are uniform in $g$ under
\eqref{eq:smallness-perturbed-sigma}.
By \eqref{eq:weighted-compact-part-pointwise}, $\cK_\ell f$ has the
form of the first summand of \eqref{eq:L2-Linfty-source-bound} with
$\varepsilon=1/8$, and is absorbed into the source $G$ in the
applications below.

\begin{thm}[Local $L^2$--$L^\infty$ estimate]
	\label{260504554}
	Let $\ell \ge0$, and assume that $g$ is defined on
	$\bR_t\times\bR_x^3\times\bR_v^3$ and satisfies
	\[
	\mathop{\rm ess\,sup}_{(t,x)\in(-4,0]\times\bR_x^3}
	\|\mu^{1/4}g(t,x,\cdot)\|_{L^\infty(\bR_v^3)}
	\le\varepsilon_{\mathrm{ell}},
	\]
	where $\varepsilon_{\mathrm{ell}}$ is the constant in
	Lemma~\ref{lem:perturbed-landau-ellipticity}.
	Let $f\in\cH_{\kin,\sigma}((-4,0)\times \cD\times \bR_v^3)$ satisfy
	\begin{equation}\label{eq:weighted-Lg-equation}
		\begin{cases}
			(\partial_t+v\cdot\nabla_x)f
			=
			\cA_{g,\ell}f+G
			&\text{in }(-4,0]\times\cD\times\bR_v^3,\\
			f=0
			&\text{on }(-4,0]\times\gamma^-,
		\end{cases}
	\end{equation}
	where $\cA_{g,\ell}$ is defined by
	\eqref{eq:weighted-Ag-expansion}--\eqref{eq:weighted-cg-expanded}; the
	first line is understood as an identity in
	$L^2((-4,0)\times\cD;H_\sigma^{-1})$ and the second as
	$f_{\gamma^-}=0$.
	Assume that
	\[
	w^{21}f
	\in
	L^2((-4,0]\times\cD\times\bR_v^3)
	\]
	and that, for some $C_G,K\ge0$ and $\varepsilon>0$,
	\begin{equation}\label{eq:L2-Linfty-source-bound}
		|G(t,x,v)|
		\le
		C_G\mu(v)^\varepsilon
		\|\mu^\varepsilon f(t,x,\cdot)\|_{L^2(\bR_v^3)}
		+
		Kw(v)^{-20}
	\end{equation}
	for almost every
	$(t,x,v)\in(-4,0]\times\cD\times\bR_v^3$.
	
	Then
	\begin{equation}\label{eq:L2-Linfty-estimate}
		\|f\|_{L^\infty((-1,0]\times\cD\times\bR_v^3)}
		\le
		C\left(
		\|w^{21}f\|_
		{L^2((-4,0]\times\cD\times\bR_v^3)}
		+
		K
		\right),
	\end{equation}
	where
	\[
	C=C(\cD,\ell,C_G,\varepsilon)>0.
	\]
	In particular, $C$ depends on $G$ only through the constants
	$C_G,\varepsilon$ of \eqref{eq:L2-Linfty-source-bound}, and is
	independent of $g$, $f$ and $K$ under the stated smallness
	condition on $g$.
\end{thm}

\begin{proof}
	Throughout the proof, $C>0$ depends only on
	$\cD$, $\ell$, $C_G$, and $\varepsilon$, and may vary from line to
	line.  Set
	\[
	p_n:=\left(\frac76\right)^n,
	\qquad
	R_n:=1+2^{-n},
	\qquad
	\cO_n:=(-R_n^2,0]\times\cD\times\bR_v^3,
	\]
	and write
	\[
	\|u\|_{L_A^p(E)}
	:=\left(\iiint_E |u|^p w^A\dd v\dd x\dd t\right)^{1/p}.
	\]
	
	Choose the lattices
	\[
	x_k\in\frac13\bZ^3,
	\qquad
	v_l\in\frac13\bZ^3,
	\qquad
	X_{k,l}:=(0,x_k,v_l).
	\]
	
	For every $R\in[1,2]$ the cylinders $Q_R^\cD(X_{k,l})$ cover
	$(-R^2,0]\times\cD\times\bR_v^3$, and
	\begin{equation}\label{eq:lattice-weight-comparability}
		w(v)\simeq w(v_l)
		\qquad\text{on }Q_2^\cD(X_{k,l}).
	\end{equation}
	
	\medskip
	\noindent\emph{Step 1. Global Moser gain.}
	On $Q_2^\cD(X_{k,l})$, estimates
	\eqref{eq:weighted-ag-ellipticity} and
	\eqref{eq:weighted-coefficient-pointwise-bounds} give admissible constants
	\[
	\lambda_l=cw(v_l)^{-3},
	\qquad
	\Lambda_l=Cw(v_l)^{-1},
	\qquad
	M_l:=\max\{C_\ell,1\}
	\]
	in Proposition~\ref{cor:boundary-moser-gain}.
	The bounds $c|\xi|^2w^{-3}\le\mathfrak q(v,\xi)\le C|\xi|^2w^{-1}$ are
	immediate from \eqref{26720259}, and we used them together with
	\eqref{eq:lattice-weight-comparability}.  
	Applying Proposition~\ref{cor:boundary-moser-gain} to
	$f$ and $-f$ therefore yields
	\begin{align}
		&\big\||f|^{p_n}\big\|_
		{L^{7/3}(Q_{R_{n+1}}^\cD(X_{k,l}))}
		\notag\\
		&\qquad\le
		Cp_n\Big(
		\big\||G||f|^{p_n-1}\big\|_
		{L^2(Q_{R_n}^\cD(X_{k,l}))}
		+
		2^{2n}w(v_l)^3
		\big\||f|^{p_n}\big\|_
		{L^2(Q_{R_n}^\cD(X_{k,l}))}
		\Big).
		\label{eq:local-moser-iteration}
	\end{align}
	Indeed, $R_n-R_{n+1}=2^{-n-1}$ and hence the coefficient factor in
	\eqref{260428431} is at most $C2^{2n}w(v_l)^3$.
	
	Multiply \eqref{eq:local-moser-iteration} by $w(v_l)^{18}$, square,
	and sum over $(k,l)$.
	The bounded overlap and $\ell^2\hookrightarrow\ell^{7/3}$ give
	\begin{equation}\label{eq:global-moser-gain}
		\big\||f|^{p_n}\big\|_
		{L_{42}^{7/3}(\cO_{n+1})}
		\le
		Cp_n\left(
		\big\||G||f|^{p_n-1}\big\|_
		{L_{36}^2(\cO_n)}
		+
		2^{2n}
		\big\||f|^{p_n}\big\|_
		{L_{42}^2(\cO_n)}
		\right).
	\end{equation}
	
	\medskip
	\noindent\emph{Step 2. Source estimate.}
	Fix $n\ge0$ and set
	\[
	p:=p_n,
	\qquad
	Z_n:=
	\big\||f|^p\big\|_{L_{42}^2(\cO_n)}.
	\]
	H\"older's inequality gives
	\begin{align}
		\big\||G||f|^{p-1}\big\|_{L_{36}^2(\cO_n)}^2
		&\le
		\left(
		\iiint_{\cO_n}
		|G|^{2p}w^{-6p+42}
		\dd v\dd x\dd t
		\right)^{1/p}
		Z_n^{2-2/p}.
		\label{260429515}
	\end{align}
	
	Set
	\[
	A(t,x)
	:=
	\|\mu^\varepsilon f(t,x,\cdot)\|_{L^2(\bR_v^3)}.
	\]
	Let $p'=p/(p-1)$, with $p'=\infty$ when $p=1$.
	
	Set
	\[
	c_\varepsilon
	:=\max\left\{
	\int_{\bR_v^3}\mu^{2\varepsilon}w^{36}\dd v,\ 1\right\}<\infty .
	\]
	Since $\mu\le1$, $w\ge1$ and $p,p'\ge1$, we have
	$\mu^{2\varepsilon p}w^{-6p+42}\le\mu^{2\varepsilon}w^{36}$ and
	$\mu^{2\varepsilon p'}w^{-42p'/p}\le\mu^{2\varepsilon}$.  The first
	two summands below are therefore at most $c_\varepsilon^{1/p}$ and
	$c_\varepsilon^{1/p'}$ respectively, hence at most $c_\varepsilon$.
	In the third summand the $(t,x)$-integration
	contributes the finite factor $R_n^2|\cD|\le4|\cD|$, and since
	$-46p+42\le-4$ for $p\ge1$, the velocity integration is dominated by
	the convergent integral $c_*:=\int_{\bR_v^3}w^{-4}\dd v$; that
	summand is therefore at most $\max\{1,4c_*|\cD|\}$.  This is where the
	$\cD$-dependence of the constant in \eqref{eq:L2-Linfty-estimate}
	originates.  Hence,
	uniformly for $p\ge1$,
	\begin{equation}\label{eq:uniform-velocity-moments}
		\begin{aligned}
			&\left(\int_{\bR_v^3}
			(w^{-3}\mu^\varepsilon)^{2p}w^{42}\dd v
			\right)^{1/p}
			+\left\|\mu^{2\varepsilon}w^{-42/p}\right\|_{L^{p'}(\bR_v^3)}
			+\left(\iiint_{\cO_n}w^{-46p+42}\dd v\dd x\dd t\right)^{1/p}
			\le C_{\varepsilon,\cD}.
		\end{aligned}
	\end{equation}
	Moreover, H\"older's inequality in $v$ gives
	\begin{align*}
		A(t,x)^2
		&\le
		\left(
		\int_{\bR_v^3}
		|f(t,x,v)|^{2p}w^{42}
		\dd v
		\right)^{1/p}
		\left\|\mu^{2\varepsilon}w^{-42/p}\right\|_{L^{p'}(\bR_v^3)}.
	\end{align*}
	Consequently, \eqref{eq:L2-Linfty-source-bound},
	$(a+b)^{2p}\le2^{2p}(a^{2p}+b^{2p})$, and
	\eqref{eq:uniform-velocity-moments} yield
	\begin{equation}\label{eq:source-moment-bound}
		\left(
		\iiint_{\cO_n}
		|G|^{2p}w^{-6p+42}
		\dd v\dd x\dd t
		\right)^{1/p}
		\le
		CZ_n^{2/p}
		+
		CK^2.
	\end{equation}
	Combining \eqref{260429515} and \eqref{eq:source-moment-bound}, we obtain
	\begin{equation*}
		\big\||G||f|^{p-1}\big\|_{L_{36}^2(\cO_n)}
		\le
		C\left(
		Z_n
		+
		KZ_n^{1-1/p}
		\right),
	\end{equation*}
	where $Z_n^{1-1/p}$ is interpreted as $1$ when $p=1$.
	
	Substituting this into \eqref{eq:global-moser-gain} and using
	$2^{2n}\ge1$, we obtain
	\begin{equation}\label{eq:global-moser-recursion}
		\big\||f|^{p_n}\big\|_
		{L_{42}^{7/3}(\cO_{n+1})}
		\le
		Cp_n2^{2n}
		\left(
		Z_n
		+
		KZ_n^{1-1/p_n}
		\right).
	\end{equation}
	
	\medskip
	\noindent\emph{Step 3. Iteration.}
	Set
	\[
	Y_n:=\|f\|_{L_{42}^{2p_n}(\cO_n)}.
	\]
	Since
	\[
	2p_{n+1}=\frac73p_n,
	\qquad
	Z_n=Y_n^{p_n},
	\]
	\eqref{eq:global-moser-recursion} gives
	\begin{equation*}
		Y_{n+1}^{p_n}
		\le
		Cp_n2^{2n}
		\left(
		Y_n^{p_n}
		+
		KY_n^{p_n-1}
		\right).
	\end{equation*}
	Put $A_n:=Cp_n2^{2n}\ge1$.
	
	Factoring $Y_n^{p_n}+KY_n^{p_n-1}=Y_n^{p_n-1}(Y_n+K)$ and taking
	the $p_n$-th root gives
	\[
	Y_{n+1}\le A_n^{1/p_n}Y_n^{1-1/p_n}(Y_n+K)^{1/p_n}.
	\]
	Adding $K$, using $A_n^{1/p_n}\ge1$, and applying the elementary
	inequality \eqref{e-inequality}, valid for all $A,B\ge0$ and every
	exponent in $(0,1]$, with $A=Y_n$, $B=K$ and exponent
	$1/p_n$, we obtain
	\begin{equation}\label{eq:Yn-linear-recursion}
		Y_{n+1}+K
		\le
		A_n^{1/p_n}
		\Bigl(
		Y_n^{1-1/p_n}(Y_n+K)^{1/p_n}+K
		\Bigr)
		\le
		(2A_n)^{1/p_n}(Y_n+K).
	\end{equation}
	
	The cruder bound $Y_n^{p_n}+KY_n^{p_n-1}\le(Y_n+K)^{p_n}$
	would instead give the factor $A_n^{1/p_n}+1\ge2$, whose infinite
	product diverges.
	Because $p_n=(7/6)^n$,
	\[
	\sum_{n=0}^\infty\frac{1+n}{p_n}<\infty.
	\]
	Thus iteration of \eqref{eq:Yn-linear-recursion} gives
	$\sup_nY_n\le C(Y_0+K)$.  Since $w^{42}\ge1$,
	$\cO_n\downarrow(-1,0]\times\cD\times\bR_v^3$, and
	$p_n\uparrow\infty$, we conclude that
	\begin{align*}
		\|f\|_
		{L^\infty((-1,0]\times\cD\times\bR_v^3)}
		&\le
		C(Y_0+K)
		=
		C\left(
		\|w^{21}f\|_
		{L^2((-4,0]\times\cD\times\bR_v^3)}
		+
		K
		\right).
	\end{align*}
\end{proof}

\subsection{Boundary flattening and normalization of the diffusion}
\label{flat_and_norm}

Two grazing sets occur below, and we distinguish them by a subscript:
\(\gamma_{\cD}^0\) is the set \(\{(x,v):x\in\partial\cD,\ n_x\cdot v=0\}\)
of \eqref{def:gamma}, while
\(\gamma_{\bR_+^3}^0=\{(x,v):x_3=0,\ v_3=0\}\) is its half-space
counterpart, introduced in Subsection~\ref{Subsec:3preparing}; the same
convention is used for \(\gamma^\pm\) when the ambient domain is not clear from the context.

\subsubsection{Boundary charts and their phase-space lifts}

We use the following quantitative definition to fix the geometric constants	\(R_{\cD}\) and \(\Lambda_{\cD}\).

\begin{defn}[\(C^{1,1}\) domain]
	\label{domain_assumption}
	A bounded domain \(\cD\subset\bR^d\), $d\ge 2$, is said to be of class \(C^{1,1}\)
	if there exist constants \(0<R_{\cD}\le1\) and \(\Lambda_{\cD}>0\)
	such that, for every \(x_0\in\partial\cD\), there exist
	\(P_{x_0}\in SO(d)\) and
	\(\eta_{x_0}\in C^{1,1}(B_{R_{\cD}}^{d-1})\) satisfying
	\begin{equation}\label{def:C11-graph-bounds}
		\eta_{x_0}(0)=0,
		\qquad
		\nabla\eta_{x_0}(0)=0,
		\qquad
		[\nabla\eta_{x_0}]_{C^{0,1}(B_{R_{\cD}}^{d-1})}
		\le \Lambda_{\cD},
	\end{equation}
	and, writing \(y=(y',y_d)\in\bR^{d-1}\times\bR\),
	\[
	\begin{aligned}
		\cD
		&=
		\bigl\{
		x_0+P_{x_0}y:
		y_d>\eta_{x_0}(y')
		\bigr\},\\
		\partial\cD
		&=
		\bigl\{
		x_0+P_{x_0}y:
		y_d=\eta_{x_0}(y')
		\bigr\}
	\end{aligned}
	\qquad\text{in}\qquad
	\bigl\{
	x_0+P_{x_0}y:
	|y'|<R_{\cD},\
	|y_d|<R_{\cD}
	\bigr\}.
	\]
\end{defn}

Although Definition~\ref{domain_assumption} is stated for arbitrary
\(d\ge2\), from \S\ref{sec:landau-framework} onward we restrict our
attention to the case \(d=3\), as announced there.

Fix \(x_0\in\partial\cD\), abbreviate \(P=P_{x_0}\) and
\(\eta=\eta_{x_0}\), and set \(U=B_{R_{\cD}}(x_0)\).  Define
\[
\psi\bigl(x_0+P(y',y_3)\bigr)
:=(y',y_3-\eta(y')),\qquad V:=\psi(U).
\]
Then \(\psi:U\to V\) is a \(C^{1,1}\)-diffeomorphism with inverse
\[
\phi(\widehat x)
=x_0+P\bigl(\widehat x',
\widehat x_3+\eta(\widehat x')\bigr).
\]
This chart sends \(\cD\), \(\partial\cD\), and \(\cD^c\) locally to
\(\{\widehat x_3>0\}\), \(\{\widehat x_3=0\}\), and
\(\{\widehat x_3\le0\}\), respectively.  Moreover,
\begin{equation}\label{260509301}
	\psi(x_0)=0,\qquad D\psi(x_0)=P^\top,\qquad
	\nabla\psi_3(x)=-|\nabla\psi_3(x)|\,n(x)
	\quad(x\in U\cap\partial\cD).
\end{equation}
There is \(r_{\cD}>0\) such that
\begin{equation}\label{eq:boundary-chart-uniform-bounds}
	B_{r_{\cD}}(0)\subset V,\qquad
	\|D\psi\|_{C^{0,1}(U)}+\|D\phi\|_{C^{0,1}(V)}
	\le C_\cD.
\end{equation}
For example, we may take \(r_{\cD}=R_{\cD}/(\Lambda_{\cD}+1)\); the inclusion
then follows from \(|\eta(y')|\le \Lambda_{\cD}|y'|^2/2\).

Write
\begin{equation*}
	A_x:=D\psi(x),\qquad B_{\widehat x}:=D\phi(\widehat x).
\end{equation*}
The phase-space lifts
\begin{equation}\label{260508457}
	\Psi(x,v):=(\psi(x),A_xv),\qquad
	\Phi(\widehat x,\widehat v)
	:=(\phi(\widehat x),B_{\widehat x}\widehat v)
\end{equation}
are mutually inverse locally Lipschitz maps.  In particular,
\(A_{\phi(\widehat x)}B_{\widehat x}=I_3\) and
\(\det A_x=\det B_{\widehat x}=1\).  At the boundary,
\begin{equation}\label{2605131213}
	\widehat v_3=(A_xv)_3
	=-|\nabla\psi_3(x)|\,n(x)\cdot v.
\end{equation}
Thus the incoming set \(n(x)\cdot v<0\) is mapped to
\(\{\widehat x_3=0,\widehat v_3>0\}\).

\subsubsection{Distance to the grazing set under flattening}

\begin{lemma}[Kinetic distance under boundary flattening]
	\label{260518305}
	Let $x_0\in\partial\cD$, and let $\psi$ and $A_x=D\psi(x)$ be the
	boundary chart defined above around $x_0$.  Set
	\[
	\Theta(t,x,v):=(t,\psi(x),A_xv).
	\]
	There exist $0<c_0<1$ and
	$C_\cD\ge1$, depending only
	on the geometric constants of $\cD$, such that the following hold.
	
	\begin{enumerate}
		\item\label{item:flattening-distance-comparison}
		If
		\[
		Z_i=(t_i,x_i,v_i)\in
		\bR\times B_{R_{\cD}}(x_0)\times\bR_v^3,\qquad i=1,2,
		\]
		and
		\begin{equation}\label{260515443}
			\mathbf d(Z_1,Z_2)\le\langle v_1\rangle^{-2},
		\end{equation}
		then
		\begin{equation}\label{260515229}
			C_\cD^{-1}\mathbf d(Z_1,Z_2)
			\le\mathbf d\bigl(\Theta(Z_1),\Theta(Z_2)\bigr)
			\le C_\cD\mathbf d(Z_1,Z_2).
		\end{equation}
		
		\item\label{item:flattening-grazing-distance}
		If \(Z=(t,x,v)\in\bR\times\overline\cD\times\bR_v^3\) satisfies
		\begin{equation}\label{eq:flattening-grazing-range}
			x\in B_{c_0R_{\cD}}(x_0),
			\qquad
			\mathbf d(Z,\bR\times\gamma_{\cD}^0)
			<c_0\langle v\rangle^{-2},
		\end{equation}
		then
		\begin{equation}\label{260515320}
			C_\cD^{-1}\mathbf d(Z,\bR\times\gamma_{\cD}^0)
			\le
			\mathbf d\bigl(\Theta(Z),\bR\times\gamma_{\bR_+^3}^0\bigr)
			\le C_\cD\,\mathbf d(Z,\bR\times\gamma_{\cD}^0).
		\end{equation}
	\end{enumerate}
\end{lemma}

\begin{proof}
	We first prove (\ref{item:flattening-distance-comparison}).
	Put $d=\mathbf d(Z_1,Z_2)$ and choose an intermediate velocity
	$\bar v$ in \eqref{kineticdistance} such that all four defining
	quantities are at most $2d$.  Then $|\bar v|\le|v_1|+2d$, so
	\[
	|x_1-x_2|
	\le
	|x_1-x_2-(t_1-t_2)\bar v|+|t_1-t_2|\,|\bar v|
	\le C\bigl(d^3+\langle v_1\rangle d^2\bigr).
	\]
	Consequently, for every $K\ge1$, the hypothesis
	\begin{equation}\label{260515443K}
		\mathbf d(Z_1,Z_2)\le K\langle v_1\rangle^{-2}
	\end{equation}
	gives $\langle v_1\rangle^2d\le K$ and $d\le K$, hence
	$\langle v_1\rangle^2d^4\le Kd^3$, $d^6\le K^3d^3$,
	$\langle v_1\rangle^2d^2\le Kd$ and $\langle v_1\rangle d^3\le K^2d$,
	and therefore
	\begin{equation}\label{eq:flattening-basic-errors}
		|x_1-x_2|^2\le C_Kd^3,
		\qquad
		|x_1-x_2|\,|v_1|\le C_Kd .
	\end{equation}
	In particular \eqref{eq:flattening-basic-errors} holds with $K=1$ under
	\eqref{260515443}.
	
	Since $B_{R_{\cD}}(x_0)$ is convex and $\psi\in C^{1,1}$ there with
	$D\psi=A$, Taylor's formula and the Lipschitz continuity of $A$ give
	\[
	\begin{aligned}
		|\psi(x_1)-\psi(x_2)-(t_1-t_2)A_{x_2}\bar v|
		&\le C_\cD\bigl(|x_1-x_2|^2+d^3\bigr),\\
		|A_{x_1}v_1-A_{x_2}\bar v|
		+|A_{x_2}v_2-A_{x_2}\bar v|
		&\le C_\cD\bigl(|x_1-x_2|\,|v_1|+d\bigr).
	\end{aligned}
	\]
	Using $A_{x_2}\bar v$ as intermediate velocity in
	\eqref{kineticdistance} and inserting
	\eqref{eq:flattening-basic-errors}, we obtain
	\begin{equation}\label{eq:flattening-upper-bound}
		\mathbf d\bigl(\Theta(Z_1),\Theta(Z_2)\bigr)
		\le C_\cD d.
	\end{equation}
	
	For the reverse inequality, write $\widehat Z_i=\Theta(Z_i)
	=(t_i,\widehat x_i,\widehat v_i)$, $\phi=\psi^{-1}$ and
	$B_{\widehat x}=A_{\phi(\widehat x)}^{-1}$.  Here $\phi$ is understood as
	the explicit map
	$\phi(\widehat x)=x_0+P(\widehat x',\widehat x_3+\eta(\widehat x'))$,
	defined and $C^{1,1}$ on the convex slab
	$\{|\widehat x'|<R_{\cD}\}\supset V$, whereas $V$ itself need not be
	convex.  Since $A_x$ and $A_x^{-1}$ are bounded by $C_\cD$, we have
	$\langle\widehat v_i\rangle\asymp_{\cD}\langle v_i\rangle$, and therefore
	\eqref{eq:flattening-upper-bound} and \eqref{260515443} give
	\[
	\widehat d:=\mathbf d\bigl(\widehat Z_1,\widehat Z_2\bigr)
	\le C_\cD\langle\widehat v_1\rangle^{-2}.
	\]
	Thus the pair $(\widehat Z_1,\widehat Z_2)$ satisfies
	\eqref{260515443K} with $K=C_\cD$, and the argument of the previous
	two paragraphs applies verbatim in the hatted variables, with
	$(\psi,A)$ replaced by $(\phi,B)$.  Since
	$\Theta^{-1}(\widehat Z_i)=Z_i$, that argument yields
	$d\le C_\cD\widehat d$, which together with
	\eqref{eq:flattening-upper-bound} proves \eqref{260515229}.
	
	We now prove (\ref{item:flattening-grazing-distance}) and write
	$r:=\mathbf d(Z,\bR\times\gamma_{\cD}^0)$ and
	$\widehat r:=\mathbf d\bigl(\Theta(Z),\bR\times\gamma_{\bR_+^3}^0\bigr)$.
	Both $\bR\times\gamma_{\cD}^0$ and $\bR\times\gamma_{\bR_+^3}^0$ are
	closed, and by \eqref{260509301} and \eqref{2605131213} the map $\Theta$
	carries the first into the second on
	$\bR\times B_{R_{\cD}}(x_0)\times\bR_v^3$, while $\Theta^{-1}$ carries
	the second into the first on $\bR\times V\times\bR_v^3$.
	Hence $r=0$ forces
	$\widehat r=0$ and conversely, and \eqref{260515320} is trivial in that
	case; we may therefore assume $r>0$.
	
	Choose $\bar Z=(\bar t,\bar x,\bar v)\in\bR\times\gamma_{\cD}^0$ with
	$\mathbf d(Z,\bar Z)<2r<2c_0\langle v\rangle^{-2}$.  As in
	the proof of (\ref{item:flattening-distance-comparison}),
	\[
	|x-\bar x|\le C\bigl(r^3+\langle v\rangle r^2\bigr)
	\le Cc_0^2 .
	\]
	After decreasing $c_0$, both $x$ and $\bar x$ lie in
	$B_{R_{\cD}}(x_0)$, and \eqref{260515443} holds for the pair
	$(Z,\bar Z)$; consequently, \eqref{260515229} and $\Theta(\bar Z)
	\in\bR\times\gamma_{\bR_+^3}^0$ yield
	\[
	\widehat r\le\mathbf d\bigl(\Theta(Z),\Theta(\bar Z)\bigr)
	\le C_\cD r .
	\]
	In particular $0<\widehat r\le C_\cD c_0
	\langle A_xv\rangle^{-2}$ by the previous paragraph.
	Choose $\bar{\widehat Z}\in\bR\times\gamma_{\bR_+^3}^0$ with
	$\mathbf d(\Theta(Z),\bar{\widehat Z})<2\widehat r$.  Its spatial
	displacement from $\Theta(Z)$ is bounded by
	$C_\cD\bigl(\widehat r^3+\langle A_xv\rangle\widehat r^2\bigr)
	\le C_\cD c_0^2$, while
	$|\psi(x)|=|\psi(x)-\psi(x_0)|\le C_\cD c_0R_{\cD}$;
	hence, after decreasing $c_0$ once more, the spatial
	component of $\bar{\widehat Z}$ lies in $B_{r_{\cD}}(0)\subset V$ by
	\eqref{eq:boundary-chart-uniform-bounds}, and the image of
	$\bar{\widehat Z}$ under $\Theta^{-1}$ belongs to
	$\bR\times\gamma_{\cD}^0$.  Since
	$\mathbf d(\Theta(Z),\bar{\widehat Z})\le2\widehat r
	\le C_\cD c_0\langle A_xv\rangle^{-2}$ is of the form
	\eqref{260515443K}, the reverse half of the proof of
	(\ref{item:flattening-distance-comparison}) applies to the pair
	$(\Theta(Z),\bar{\widehat Z})$ and gives
	\[
	r\le\mathbf d\bigl(Z,\Theta^{-1}(\bar{\widehat Z})\bigr)
	\le C_\cD\widehat r .
	\]
	This proves \eqref{260515320}.
\end{proof}

\subsubsection{Boundary flattening and normalization}

The boundary chart constructed above reduces the local geometry
to the half-space. We now compose it with a linear transformation in velocity
chosen to normalize the diffusion at the center while preserving the boundary
orientation.

\begin{lemma}[Boundary localization and normalization]
	\label{landau_local_bdry}
	Let \(Z_0=(t_0,x_0,v_0)\in\bR\times\partial\cD\times\bR^3\), and
	let \(A_*=A_{x_0}=D\psi(x_0)\) and \(\sigma_*=\sigma(v_0)\).  Choose
	\(\cR\in SO(3)\) so that
	\[
	\cR e_3
	=\frac{\sigma_*^{1/2}A_*^{\top}e_3}
	{|\sigma_*^{1/2}A_*^{\top}e_3|},
	\]
	and set
	\begin{equation}\label{260508503}
		\cS:=A_*\sigma_*^{1/2}\cR,\qquad
		\cT:=\cS^{-1}
		=\cR^\top\sigma_*^{-1/2}A_*^\top,
	\end{equation}
	where the last identity uses \(A_*=P^\top\in SO(3)\) from
	\eqref{260509301}.
	Then
	\begin{equation}\label{eq:boundary-normal-direction}
		\cS^{\top}e_3=\kappa_*e_3,\qquad
		\kappa_*:=|\sigma_*^{1/2}A_*^{\top}e_3|>0.
	\end{equation}
	
	There exist \(0<\varepsilon_{\rm b}<1\) and \(C_\cD\ge1\), depending
	only on \(\cD\), with the following properties.  Put
	\begin{equation}\label{260512414}
		r_b:=\varepsilon_{\rm b}^2\langle v_0\rangle^{-4},\qquad
		Y_0:=(t_0,0,\cT A_*v_0),
	\end{equation}
	and define
	\[
	\widetilde\Phi_b(t,y,w)
	:=\bigl(t,\Phi(\cS y,\cS w)\bigr).
	\]
	Then, for every \(0<r\le r_b\), the following hold.
	
	\begin{enumerate}
		\item\label{item:boundary-cylinder-comparison}
		The map \(\widetilde\Phi_b\) is well defined on \(Q_r(Y_0)\), and
		\begin{equation}\label{260512642}
			Q_{\varepsilon_{\rm b}\langle v_0\rangle^{-3/2}r}(Z_0)
			\subset \widetilde\Phi_b(Q_r(Y_0))
			\subset Q_{\varepsilon_{\rm b}^{-1}r}(Z_0)
			\subset Q_1(Z_0)\cap(\bR\times U\times\bR^3).
		\end{equation}
		It maps \(Q_r^+(Y_0)\) into
		\(Q_{\varepsilon_{\rm b}^{-1}r}^{\cD}(Z_0)\), maps
		\(Q_r(Y_0)\cap\{y_3=0,\,w_3>0\}\) into the incoming boundary, and
		maps \(Q_r(Y_0)\cap\{y_3=0,\,w_3=0\}\) into
		\(\bR\times\gamma_{\cD}^0\).
		
		\item\label{item:boundary-transformed-equation}
		Let \(a:Q_1(Z_0)\to\bR_{\rm sym}^{3\times3}\) and
		\(b:Q_1(Z_0)\to\bR^3\) be defined on the full cylinder and, for some
		\(M_0\ge0\), satisfy a.e.
		\begin{equation}\label{2605131249}
			\frac12\,\sigma(v)\le a(t,x,v)\le2\,\sigma(v),\qquad
			|\nabla_va(t,x,v)|\le M_0\langle v\rangle^{-2},\qquad
			|b(t,x,v)|\le M_0\langle v\rangle^{-1}.
		\end{equation}
		The matrix inequalities are understood as inequalities of quadratic forms.
		The derivative of \(a\) is understood in the weak sense.
		Let \(G\in L^2\bigl(Q_{\varepsilon_{\rm b}^{-1}r}^{\cD}(Z_0)\bigr)\) and
		suppose that \(f\) is a weak solution, in the sense of
		Definition~\ref{def:L2sol}, of
		\begin{equation}\label{eq:boundary-original-equation}
			\begin{cases}
				(\partial_t+v\cdot\nabla_x)f
				=D_{v_i}(a^{ij}D_{v_j}f)+b^iD_{v_i}f+G
				&\text{in }Q_{\varepsilon_{\rm b}^{-1}r}^{\cD}(Z_0),\\
				f=0
				&\text{on }Q_{\varepsilon_{\rm b}^{-1}r}(Z_0)
				\cap(\bR\times\gamma^-).
			\end{cases}
		\end{equation}
		For \(Y=(t,y,w)\in Q_r(Y_0)\), write
		\[
		(x,v):=\Phi(\cS y,\cS w),\qquad Z:=(t,x,v),
		\]
		and define
		\begin{align}
			\widetilde a(Y)
			&:=\cT A_xa(Z)A_x^\top\cT^\top,
			\label{eq:boundary-transformed-a}\\
			\widetilde b(Y)
			&:=\cT\bigl(A_xb(Z)-\mathfrak b_A(x,v)\bigr),
			\qquad
			\mathfrak b_A(x,v):=(v\cdot\nabla_xA_x)v.
			\label{eq:boundary-transformed-b}
		\end{align}
		Here \(\mathfrak b_A\) is defined a.e. because \(A\) is Lipschitz.  Then
		\(F:=f\circ\widetilde\Phi_b\) is a weak solution, in the same sense with
		\(\cD\) replaced by the half-space
		\(\bR^3_+:=\{y\in\bR^3:y_3>0\}\), of
		\begin{equation}\label{260512656}
			\begin{cases}
				(\partial_t+w\cdot\nabla_y)F
				=D_{w_i}(\widetilde a^{ij}D_{w_j}F)+\widetilde b^iD_{w_i}F
				+G\circ\widetilde\Phi_b
				&\text{in }Q_r^+(Y_0),\\
				F=0&\text{on }Q_r(Y_0)\cap\{y_3=0,w_3>0\}.
			\end{cases}
		\end{equation}
		
		\item\label{item:boundary-transformed-coefficients}
		Under \eqref{2605131249},
		\begin{align}
			C_\cD^{-1}|\xi|^2
			&\le \widetilde a^{ij}(Y)\xi_i\xi_j
			\le C_\cD|\xi|^2,
			\label{eq:boundary-normalized-ellipticity}\\
			|\nabla_w\widetilde a(Y)|+|\widetilde b(Y)|
			&\le C_\cD(M_0+1)\langle v_0\rangle^4
			\label{eq:boundary-normalized-coefficient-bound}
		\end{align}
		for a.e. \(Y\in Q_r(Y_0)\) and every \(\xi\in\bR^3\).
		
		\item\label{item:boundary-transformed-holder}
		If, in addition, \(0<\alpha_0<1\), \(M_1\ge0\), and
		\begin{equation}\label{eq:boundary-original-holder-assumption}
			[a]_{C_{\kin}^{\alpha_0}(Q_1(Z_0))}\le M_1,
		\end{equation}
		then
		\begin{equation}\label{eq:boundary-normalized-holder-bound}
			[\widetilde a]_{C_{\kin}^{\alpha_0}(Q_r(Y_0))}
			\le C_\cD(M_1+1)\langle v_0\rangle^3.
		\end{equation}
	\end{enumerate}
\end{lemma}

\begin{proof}
	Since \(A_*=P^\top\in SO(3)\), the spectral bounds for \(\sigma(v_0)\)
	give
	\begin{equation}\label{eq:boundary-ST-bounds}
		|\cS|\le C_\cD\langle v_0\rangle^{-1/2},\qquad
		|\cT|\le C_\cD\langle v_0\rangle^{3/2},\qquad
		|\cT A_*v_0|\le C_\cD\langle v_0\rangle^{5/2}.
	\end{equation}
	Set \(w_*=\cT A_*v_0\).  If \(Y\in Q_r(Y_0)\), then
	\[
	|w-w_*|<r,\qquad
	|y-(t-t_0)w_*|<r^3.
	\]
	Hence, by \eqref{eq:boundary-ST-bounds},
	\[
	\begin{aligned}
		|\cS y|
		&\le |\cS|
		\bigl(r^3+r^2|w_*|\bigr)\le C_\cD\bigl(
		\langle v_0\rangle^{-1/2}r^3
		+\langle v_0\rangle^2r^2\bigr)
		\le C_\cD\varepsilon_{\rm b}^4.
	\end{aligned}
	\]
	We choose \(\varepsilon_{\rm b}\) so that the last quantity is smaller
	than \(r_{\cD}\).  Thus \(\cS y\in V\), and
	\(\widetilde\Phi_b\) is well defined on \(Q_r(Y_0)\).
	
	Using \(\cS w_*=A_*v_0\), the \(C^{1,1}\) bounds in
	\eqref{eq:boundary-chart-uniform-bounds}, and Taylor expansion at the
	origin, we obtain
	\begin{align}
		|v-v_0|
		&\le C_\cD\bigl(
		r\langle v_0\rangle^{-1/2}
		+r^2\langle v_0\rangle^3\bigr),
		\label{eq:boundary-forward-velocity}\\
		|x-x_0-(t-t_0)v_0|
		&\le C_\cD\bigl(
		r^3\langle v_0\rangle^{-1/2}
		+r^4\langle v_0\rangle^4\bigr).
		\label{eq:boundary-forward-position}
	\end{align}
	In particular, after decreasing \(\varepsilon_{\rm b}\) once more,
	\begin{equation}\label{eq:boundary-local-sigma-condition}
		|v-v_0|
		\le C_\cD\varepsilon_{\rm b}^2
		\langle v_0\rangle^{-9/2}
		<\frac{c}{1+|v_0|},
	\end{equation}
	where \(c\) is the constant in Lemma~\ref{260506438}, and
	\eqref{eq:local-comparability-sigma} is therefore valid throughout the
	image cylinder.
	
	For \(r\le r_b\), the bounds \eqref{eq:boundary-forward-velocity} and
	\eqref{eq:boundary-forward-position} imply the middle and right inclusions
	in \eqref{260512642}, after choosing \(\varepsilon_{\rm b}\) sufficiently
	small.  Conversely, set
	\(q=\varepsilon_{\rm b}\langle v_0\rangle^{-3/2}r\).
	If \(Z=(t,x,v)\in Q_q(Z_0)\), then
	\[
	|x-x_0|
	\le q^3+q^2|v_0|
	\le C\varepsilon_{\rm b}^6<R_{\cD},
	\]
	after reducing \(\varepsilon_{\rm b}\) if necessary.  Hence \(x\in U\),
	so the inverse map is well defined.  Taylor expansion of \(\psi\) at
	\(x_0\) now gives
	\[
	\begin{aligned}
		|\cT A_xv-w_*|
		&\le C_\cD\bigl(
		\langle v_0\rangle^{3/2}q
		+\langle v_0\rangle^{7/2}q^2
		+\langle v_0\rangle^{5/2}q^3\bigr)<r,\\
		|\cT\psi(x)-(t-t_0)w_*|
		&\le C_\cD\bigl(
		\langle v_0\rangle^{3/2}q^3
		+\langle v_0\rangle^{7/2}q^4\bigr)<r^3.
	\end{aligned}
	\]
	This proves the left inclusion; the same choice of
	\(\varepsilon_{\rm b}\) ensures
	\(Q_{\varepsilon_{\rm b}^{-1}r}(Z_0)\subset Q_1(Z_0)\cap
	(\bR\times U\times\bR^3)\).  Finally,
	\eqref{eq:boundary-normal-direction} and \eqref{2605131213} give the
	claims concerning the half-space and the incoming boundary.
	
	To verify \eqref{260512656}, first set
	\(\widehat F(t,\widehat x,\widehat v)
	:=f(t,\Phi(\widehat x,\widehat v))\).  The chain rule gives
	\[
	\begin{aligned}
		(\partial_t+\widehat v\cdot\nabla_{\widehat x})\widehat F
		&=
		D_{\widehat v_i}\!\left(
		(A_xaA_x^\top)^{ij}D_{\widehat v_j}\widehat F\right)+(A_xb-\mathfrak b_A)\cdot\nabla_{\widehat v}\widehat F
		+G(t,\Phi(\widehat x,\widehat v)),
	\end{aligned}
	\]
	where \((x,v)=\Phi(\widehat x,\widehat v)\).  Substituting
	\((\widehat x,\widehat v)=(\cS y,\cS w)\) yields
	\eqref{260512656}.
	
	By \eqref{260508457}, $\Phi$ is locally bi-Lipschitz with
	$\det D\Phi=1$ a.e., and the linear substitution
	$(y,w)\mapsto(\cS y,\cS w)$ contributes the constant Jacobian
	$(\det\cS)^2$, which cancels from both sides.  The pullback
	$\varphi:=\widetilde\varphi\circ\widetilde\Phi_b^{-1}$ of a test function
	$\widetilde\varphi$ for the transformed problem is compactly supported and
	Lipschitz, and \eqref{260701105} therefore applies to $\varphi$, since it
	involves only $\varphi$, $D_i\varphi$ and
	$(\partial_t+v\cdot\nabla_x)\varphi$ and so extends from $C_c^\infty$ to
	that class by mollification.  At fixed $x$ the velocity substitution is
	linear with an $\widehat x$-dependent but $\widehat v$-independent matrix,
	so the divergence structure in $v$ is preserved.
	
	The only term that is not immediate is the transport error: changing
	variables in
	$\iiint f\,(\partial_t+v\cdot\nabla_x)\varphi$ produces the extra term
	$\iiint\widehat F\,\mathfrak b_A\cdot\nabla_{\widehat v}\widehat\varphi$,
	where $\widehat F$ and $\widehat\varphi$ denote $f$ and $\varphi$ in the
	hatted variables, and integrating by parts in $\widehat v$ turns it into
	$-\iiint(\mathfrak b_A\cdot\nabla_{\widehat v}\widehat F)\widehat\varphi
	-\iiint\widehat F(\operatorname{div}_{\widehat v}\mathfrak b_A)
	\widehat\varphi$.  The first of these is the drift
	\eqref{eq:boundary-transformed-b}.  The second is a zeroth-order term that
	\eqref{260512656} does not carry, and it must therefore be checked that
	$\operatorname{div}_{\widehat v}\mathfrak b_A$ vanishes.  It does: writing
	\(y=P^\top(x-x_0)\), we have \(A_x=M(y')P^\top\), where
	\[
	M(y')=\begin{pmatrix}I_2&0\\-\nabla\eta(y')^\top&1\end{pmatrix},
	\]
	\(I_2\) is the \(2\times2\) identity, and \(y'=\widehat x'\) when
	\(x=\phi(\widehat x)\).  Only the third row of
	\(M\) varies, and it varies only in \(y'\).  Since
	\((P^{\top}v)'=\widehat v'\) for \(v=B_{\widehat x}\widehat v\), we obtain
	\begin{equation}\label{eq:transport-error-chart-form}
		\mathfrak b_A\bigl(\phi(\widehat x),B_{\widehat x}\widehat v\bigr)
		=\Bigl(0,\,0,\,
		-\sum_{k,l=1}^2\partial_k\partial_l\eta(\widehat x')
		\widehat v_k\widehat v_l\Bigr)
		\qquad\text{for a.e. }\widehat x.
	\end{equation}
	This expression is independent of \(\widehat v_3\) and has vanishing
	first two components.  Hence
	$\operatorname{div}_{\widehat v}\mathfrak b_A=0$, and no zeroth-order term
	arises.
	
	Finally, by
	\eqref{2605131213} and
	\eqref{eq:boundary-normal-direction}, $\widetilde\Phi_b$ maps
	$Q_r(Y_0)\cap\{y_3=0,\,w_3>0\}$ bi-Lipschitz onto a subset of
	$(\bR\times\gamma^-)\cap Q_{\varepsilon_{\rm b}^{-1}r}(Z_0)$, and the
	two kinetic boundary measures are mutually absolutely continuous there,
	with density bounded above and below by positive constants depending on
	$Z_0$.  The vanishing of the kinetic trace of
	Corollary~\ref{cor:trace-bound-scaled-cylinder} is therefore carried
	along, and this gives the second line of \eqref{260512656}.  No constant
	from this step enters the estimates below.
	
	At the center,
	\[
	\cT A_*\sigma_*A_*^\top\cT^\top=I_3.
	\]
	The anisotropic effect of replacing \(A_*\) by \(A_x\) is controlled by
	\[
	E_x:=
	\sigma_*^{1/2}(A_x^{\top}A_*-I_3)\sigma_*^{-1/2},
	\qquad
	\cT A_x\sigma_*A_x^\top\cT^\top
	=\cR^\top(I_3+E_x)^\top(I_3+E_x)\cR,
	\]
	and
	\[
	\begin{aligned}
		\|E_x\|
		&\le
		C_\cD\langle v_0\rangle |x-x_0|
		\le C_\cD\langle v_0\rangle^2r^2
		\le C_\cD\varepsilon_{\rm b}^4.
	\end{aligned}
	\]
	Indeed, the middle inequality follows from
	\eqref{eq:boundary-forward-position} and
	\(|t-t_0||v_0|\le r^2|v_0|\); for
	\(\varepsilon_{\rm b}\) sufficiently small, this estimate,
	\eqref{eq:boundary-local-sigma-condition}, and the first bound in
	\eqref{2605131249} yield
	\eqref{eq:boundary-normalized-ellipticity}.  Moreover,
	\(\langle v\rangle\asymp\langle v_0\rangle\) on the image cylinder,
	while \eqref{eq:transport-error-chart-form} and
	\([\nabla\eta]_{C^{0,1}}\le \Lambda_{\cD}\) give
	\[
	|\mathfrak b_A(x,v)|\le \Lambda_{\cD}|\widehat v|^2
	\le C_\cD\langle v_0\rangle^2.
	\]
	Consequently,
	\[
	|\nabla_w\widetilde a|
	\le C_\cD M_0\langle v_0\rangle^{1/2},\qquad
	|\widetilde b|
	\le C_\cD\bigl(
	M_0\langle v_0\rangle^{1/2}
	+\langle v_0\rangle^{7/2}\bigr),
	\]
	which proves \eqref{eq:boundary-normalized-coefficient-bound}, the
	exponent $4$ being a convenient bound for the $7/2$ just obtained.
	
	It remains to prove (\ref{item:boundary-transformed-holder}).
	Let \(Y,Y'\in Q_r(Y_0)\), put
	\(Z=\widetilde\Phi_b(Y)\), \(Z'=\widetilde\Phi_b(Y')\) and
	\(\delta:=\mathbf d(Y,Y')\le2r\le2\varepsilon_{\rm b}^2
	\langle v_0\rangle^{-4}\), and let \(\bar w\) be an almost minimizer for
	\(\delta\) in \eqref{kineticdistance}.  Since \(\cS w_*=A_*v_0\)
	and \(|\cS|\le C_\cD\langle v_0\rangle^{-1/2}\), we have
	\(|\cS\bar w|\le C_\cD\langle v_0\rangle\), and therefore
	\[
	|\cS(y-y')|
	\le C_\cD\bigl(
	\langle v_0\rangle^{-1/2}\delta^3+\langle v_0\rangle\delta^2\bigr).
	\]
	Expanding \(\phi\) to second order at \(\cS y'\), with error
	\(C_\cD|\cS(y-y')|^2\le C_\cD\varepsilon_{\rm b}^2\delta^3\),
	and using \(B_{\cS y'}\cS\bar w\) as intermediate velocity in
	\eqref{kineticdistance}, we obtain
	\[
	\mathbf d(Z,Z')\le C_\cD\delta,
	\qquad
	|x-x'|\le C_\cD|\cS(y-y')|
	\le C_\cD\langle v_0\rangle\delta^{\alpha_0},
	\]
	the last step using \(\delta\le2\varepsilon_{\rm b}^2
	\langle v_0\rangle^{-4}\le1\) and \(\alpha_0<1\).  Writing
	\[
	\widetilde a(Y)-\widetilde a(Y')
	=\cT\Bigl(
	A_x\bigl(a(Z)-a(Z')\bigr)A_x^\top
	+(A_x-A_{x'})a(Z')A_x^\top
	+A_{x'}a(Z')(A_x-A_{x'})^\top\Bigr)\cT^\top
	\]
	and combining this identity with the two bounds obtained above,
	\(|a|\le C\langle v_0\rangle^{-1}\),
	\(|\cT|\le C_\cD\langle v_0\rangle^{3/2}\),
	\eqref{eq:boundary-original-holder-assumption}, and the Lipschitz bound
	for \(A\), we obtain
	\[
	[\widetilde a]_{C_{\kin}^{\alpha_0}(Q_r(Y_0))}
	\le C_\cD(M_1+1)\langle v_0\rangle^3,
	\]
	the first summand contributing \(C_\cD M_1\langle v_0\rangle^3\) and the
	remaining two \(C_\cD\langle v_0\rangle^3\); this
	is \eqref{eq:boundary-normalized-holder-bound}.
\end{proof}

\subsubsection{Interior localization and normalization}

The next two lemmas select an admissible interior cylinder and normalize the diffusion.

\begin{lemma}
	\label{lem:interior-kinetic-cylinder}
	If \(Z_0=(t_0,x_0,v_0)\in\bR\times\cD\times\bR^3\),
	\(0<R\le1\), and
	\(\operatorname{dist}(x_0,\partial\cD)\ge R\), then
	\[
	Q_{\sqrt{R/(2\langle v_0\rangle)}}(Z_0)
	\subset\bR\times\cD\times\bR^3.
	\]
\end{lemma}

\begin{proof}
	For \(Z=(t,x,v)\) in the indicated cylinder, with
	\(r=\sqrt{R/(2\langle v_0\rangle)}\),
	\[
	|x-x_0|
	\le |x-x_0-(t-t_0)v_0|+|t-t_0||v_0|
	<r^3+r^2|v_0|\le R.
	\]
\end{proof}

\begin{lemma}
	\label{landau_local_interior}
	Let \(Z_0=(t_0,x_0,v_0)\in\bR\times\cD\times\bR^3\), and set
	\[
	\cS:=\sigma(v_0)^{1/2},\qquad
	\cT:=\cS^{-1}.
	\]
	There are universal constants \(0<\varepsilon_{\rm i}<1\) and \(C\ge1\)
	with the following properties.  Define
	\begin{equation}\label{eq:interior-radius}
		R_0:=\min\{1,\operatorname{dist}(x_0,\partial\cD)\},\qquad
		r_i:=\varepsilon_{\rm i}\sqrt{\frac{R_0}{2\langle v_0\rangle}},
	\end{equation}
	\begin{equation}\label{eq:interior-normalized-center}
		Y_0:=(t_0,\cT x_0,\cT v_0),
		\qquad
		\widetilde\Phi_i(t,y,w):=(t,\cS y,\cS w).
	\end{equation}
	For every \(0<r\le r_i\):
	
	\begin{enumerate}
		\item\label{item:interior-cylinder-comparison}
		\begin{equation}\label{260512642_11}
			Q_{\varepsilon_{\rm i}\langle v_0\rangle^{-3/2}r}(Z_0)
			\subset\widetilde\Phi_i(Q_r(Y_0))
			\subset Q_{\varepsilon_{\rm i}^{-1}r}(Z_0)
			\subset Q_1(Z_0)\cap(\bR\times\cD\times\bR^3).
		\end{equation}
		
		\item\label{item:interior-transformed-equation}
		Let \(a:Q_1(Z_0)\to\bR_{\rm sym}^{3\times3}\) and
		\(b:Q_1(Z_0)\to\bR^3\) be defined on the full cylinder and, for some
		\(M_0\ge0\), satisfy a.e.
		\begin{equation}\label{eq:interior-coefficient-assumption}
			\frac12\,\sigma(v)\le a(t,x,v)\le2\,\sigma(v),\qquad
			|\nabla_va(t,x,v)|\le M_0\langle v\rangle^{-2},\qquad
			|b(t,x,v)|\le M_0\langle v\rangle^{-1}.
		\end{equation}
		The derivative of \(a\) is understood in the weak sense.
		Let \(G\in L^2\bigl(Q_{\varepsilon_{\rm i}^{-1}r}(Z_0)\bigr)\) and suppose
		that \(f\) is a weak solution, in the sense of
		Definition~\ref{def:L2sol}, of
		\begin{equation}\label{eq:interior-original-equation}
			(\partial_t+v\cdot\nabla_x)f
			=D_{v_i}(a^{ij}D_{v_j}f)+b^iD_{v_i}f+G
			\quad\text{in }Q_{\varepsilon_{\rm i}^{-1}r}(Z_0).
		\end{equation}
		For \(Y=(t,y,w)\in Q_r(Y_0)\), put
		\(Z=\widetilde\Phi_i(Y)\) and
		\begin{equation}\label{eq:interior-transformed-coefficients}
			\widetilde a(Y):=\cT a(Z)\cT^\top,\qquad
			\widetilde b(Y):=\cT b(Z).
		\end{equation}
		Then \(F:=f\circ\widetilde\Phi_i\) satisfies
		\begin{equation}\label{eq:interior-transformed-equation}
			(\partial_t+w\cdot\nabla_y)F
			=D_{w_i}(\widetilde a^{ij}D_{w_j}F)+\widetilde b^iD_{w_i}F
			+G\circ\widetilde\Phi_i
			\quad\text{in }Q_r(Y_0).
		\end{equation}
		
		\item\label{item:interior-transformed-coefficients}
		Under \eqref{eq:interior-coefficient-assumption},
		\begin{align}
			&C^{-1}|\xi|^2
			\le\widetilde a^{ij}(Y)\xi_i\xi_j\le C|\xi|^2,
			\label{eq:interior-normalized-ellipticity}\\
			&|\nabla_w\widetilde a(Y)|+|\widetilde b(Y)|
			\le CM_0\langle v_0\rangle^{1/2}.
			\label{eq:interior-normalized-coefficient-bound}
		\end{align}
		
		\item\label{item:interior-transformed-holder}
		If, in addition, \(M_1\ge0\) and
		\[
		0<\alpha_0<1,\qquad
		[a]_{C_{\kin}^{\alpha_0}(Q_1(Z_0))}\le M_1,
		\]
		then
		\begin{equation}\label{eq:interior-normalized-holder-bound}
			[\widetilde a]_{C_{\kin}^{\alpha_0}(Q_r(Y_0))}
			\le CM_1\langle v_0\rangle^3.
		\end{equation}
	\end{enumerate}
\end{lemma}

\begin{proof}
	The spectral bounds for \(\sigma(v_0)\) give
	\begin{equation}\label{eq:interior-ST-bounds}
		|\cS|\le C\langle v_0\rangle^{-1/2},\qquad
		|\cT|\le C\langle v_0\rangle^{3/2},\qquad
		|\cT v_0|\le C\langle v_0\rangle^{5/2}.
	\end{equation}
	Because \(\widetilde\Phi_i\) is linear and applies the same matrix to
	\(y\) and \(w\), it preserves the form of
	\(\partial_t+w\cdot\nabla_y\).  The first two bounds in
	\eqref{eq:interior-ST-bounds} give the first two inclusions in
	\eqref{260512642_11}; the last inclusion follows from
	Lemma~\ref{lem:interior-kinetic-cylinder} and the definition of \(r_i\).
	
	The transformed equation follows directly from the chain rule.  On the
	image cylinder,
	\[
	|v-v_0|
	\le|\cS|\,|w-\cT v_0|
	\le C\varepsilon_{\rm i}\langle v_0\rangle^{-1}
	<\frac{c}{1+|v_0|},
	\]
	where \(c\) is the constant in Lemma~\ref{260506438} and
	\(\varepsilon_{\rm i}\) has been chosen sufficiently small.  Therefore,
	\(\langle v\rangle\asymp\langle v_0\rangle\), and
	\eqref{eq:local-comparability-sigma} applies throughout the image
	cylinder.
	The identity
	\(\cT\sigma(v_0)\cT^\top=I_3\), together with
	\eqref{eq:local-comparability-sigma} and the first bound in
	\eqref{eq:interior-coefficient-assumption}, proves
	\eqref{eq:interior-normalized-ellipticity}.  Finally,
	\[
	|\nabla_w\widetilde a|
	\le|\cT|^2|\cS|\,|\nabla_va|
	\le CM_0\langle v_0\rangle^{1/2},
	\qquad
	|\widetilde b|\le|\cT|\,|b|
	\le CM_0\langle v_0\rangle^{1/2},
	\]
	while
	\[
	\mathbf d(\widetilde\Phi_i(Y),\widetilde\Phi_i(Y'))
	\le C\mathbf d(Y,Y')
	\]
	follows from \eqref{kineticdistance} and
	\(|\cS|\le C\); multiplying \(a\) by \(\cT\) on both
	sides, we obtain \eqref{eq:interior-normalized-holder-bound}.
\end{proof}

\subsubsection{Boundary and interior kinetic covering}

For \(M\ge1\), set
\begin{equation*}
	\rho_M(v):=M^{-1}\langle v\rangle^{-6},\qquad
	\varsigma_M(v):=M^{-3}\langle v\rangle^{-12}.
\end{equation*}
The first scale is used in the boundary layer, and the smaller second
scale at interior points.

\begin{lemma}[Boundary--interior kinetic covering]
	\label{covering_landau}
	For every \(M\ge1\), there are sequences
	\[
	Z_{j,k}=(t_{j,k},x_{j,k},v_{j,k})
	\in(-1,0]\times\overline\cD\times\bR^3,
	\qquad j=1,2,\quad k\in\bN,
	\]
	such that
	\begin{equation}\label{eq:covering-interior-distance}
		x_{1,k}\in\partial\cD,\qquad
		\operatorname{dist}(x_{2,k},\partial\cD)
		\ge\frac{\langle v_{2,k}\rangle^{-18}}{2M^3},
	\end{equation}
	and
	\begin{equation}\label{eq:boundary-interior-covering}
		(-1,0]\times\overline\cD\times\bR^3
		\subset
		\bigcup_{k\in\bN}Q_{\rho_M(v_{1,k})}(Z_{1,k})
		\,\cup\,
		\bigcup_{k\in\bN}Q_{\varsigma_M(v_{2,k})}(Z_{2,k}).
	\end{equation}
	In particular,
	\begin{equation}\label{260522351}
		\sqrt{\frac{\operatorname{dist}(x_{2,k},\partial\cD)}
			{2\langle v_{2,k}\rangle}}
		\ge\frac{\langle v_{2,k}\rangle^{-19/2}}{2M^{3/2}}
		\ge\frac{\langle v_{2,k}\rangle^{-10}}{2M^2}.
	\end{equation}
\end{lemma}

\begin{proof}
	The balls \(B(v,\varsigma_M(v)/4)\), \(v\in\bR^3\), cover \(\bR^3\), so by
	the Lindel\"of property there is a countable subcover
	\begin{equation*}
		\bR^3=\bigcup_{m\in\bN}B(v_m,\varsigma_m/4),\qquad
		\varsigma_m:=M^{-3}\langle v_m\rangle^{-12},
	\end{equation*}
	and we write \(\rho_m:=M^{-1}\langle v_m\rangle^{-6}\).  Moreover,
	\(\varsigma_m\le\rho_m\).  For each \(m\), split
	\[
	E_m:=\{x\in\overline\cD:\operatorname{dist}(x,\partial\cD)\le\rho_m^3/2\},
	\qquad
	F_m:=\{x\in\overline\cD:\operatorname{dist}(x,\partial\cD)\ge\rho_m^3/2\}.
	\]
	Both \(\partial\cD\) and \(F_m\) are compact, so there are finite families
	\[
	x_{1,m,\ell}\in\partial\cD,\qquad
	x_{2,m,\ell}\in F_m
	\]
	with
	\(\partial\cD\subset\bigcup_\ell B(x_{1,m,\ell},\rho_m^3/4)\) and
	\(F_m\subset\bigcup_\ell B(x_{2,m,\ell},\varsigma_m^3/4)\).  A nearest
	boundary point of \(x\in E_m\) lies within \(\rho_m^3/2\) of \(x\), and
	therefore
	\begin{equation}\label{eq:covering-spatial-subcovers}
		E_m\subset\bigcup_\ell B(x_{1,m,\ell},3\rho_m^3/4),
		\qquad
		F_m\subset\bigcup_\ell B(x_{2,m,\ell},\varsigma_m^3/4).
	\end{equation}
	Empty families are allowed.
	
	Let \((t_n)\) enumerate \(\mathbb Q\cap(-1,0]\), and form the countable
	families
	\[
	(t_n,x_{1,m,\ell},v_m),\qquad
	(t_n,x_{2,m,\ell},v_m).
	\]
	After relabeling, these are the required \(Z_{1,k}\) and \(Z_{2,k}\).
	The second family satisfies
	\[
	\operatorname{dist}(x_{2,m,\ell},\partial\cD)
	\ge\rho_m^3/2
	=\frac{\langle v_m\rangle^{-18}}{2M^3},
	\]
	which proves \eqref{eq:covering-interior-distance}; taking square roots
	gives \eqref{260522351}.
	
	It remains to check the covering.  Fix
	\(Z=(t,x,v)\in(-1,0]\times\overline\cD\times\bR^3\), and choose \(m\)
	so that \(|v-v_m|<\varsigma_m/4\).  If \(x\in E_m\), select
	\(x_{1,m,\ell}\) from \eqref{eq:covering-spatial-subcovers} and choose
	\(t_n\in\mathbb Q\cap(-1,0]\) with
	\[
	t\le t_n,\qquad
	t_n-t<\min\{\rho_m^2,\rho_m^3/[4(1+|v_m|)]\}.
	\]
	Then \(t\in(t_n-\rho_m^2,t_n]\), \(|v-v_m|<\varsigma_m/4\le\rho_m\), and
	\[
	|x-x_{1,m,\ell}-(t-t_n)v_m|
	\le\frac{3\rho_m^3}{4}+\frac{\rho_m^3|v_m|}{4(1+|v_m|)}
	<\rho_m^3,
	\]
	so \(Z\in Q_{\rho_m}(t_n,x_{1,m,\ell},v_m)\).  If \(x\in F_m\), choose
	\(x_{2,m,\ell}\) and \(t_n\) in the same way with \(\rho_m\) replaced by
	\(\varsigma_m\) and with the time bound
	\(t_n-t<\min\{\varsigma_m^2,3\varsigma_m^3/[4(1+|v_m|)]\}\).  This gives
	\(Z\in Q_{\varsigma_m}(t_n,x_{2,m,\ell},v_m)\).  Therefore,
	\eqref{eq:boundary-interior-covering} holds.
\end{proof}

\subsection{Weighted H\"older regularity}

We first record the coefficient input.  If \(g\) is defined on
\(\bR\times\bR_x^3\times\bR_v^3\) and satisfies
\eqref{eq:smallness-perturbed-sigma} uniformly in \((t,x)\), then
Lemma~\ref{lem:perturbed-landau-ellipticity},
\eqref{eq:derivative-sigma-bound},
\eqref{eq:perturbed-sigma-i-bound},
\eqref{eq:derivative-perturbed-sigma-bound}, and
\eqref{260507202} show that the globally defined coefficients
\(a_g,b_g\) satisfy, for a universal constant \(C\),
\begin{equation}\label{eq:global-holder-landau-coefficient-bounds}
	\begin{gathered}
		\frac12\,\sigma^{ij}(v)\xi_i\xi_j
		\le a_g^{ij}(t,x,v)\xi_i\xi_j
		\le2\,\sigma^{ij}(v)\xi_i\xi_j,\\
		|\nabla_va_g(t,x,v)|
		\le C\langle v\rangle^{-2},
		\qquad
		|b_g(t,x,v)|
		\le C\langle v\rangle^{-1}
	\end{gathered}
\end{equation}
for a.e. \((t,x,v)\) and every \(\xi\in\bR^3\).

\begin{thm}[Global weighted kinetic H\"older estimate]
	\label{260526316}
	Let \(g\) be defined on \(\bR\times\bR_x^3\times\bR_v^3\) and satisfy
	\eqref{eq:smallness-perturbed-sigma} uniformly in \((t,x)\), so that
	\(a_g\) and \(b_g\) satisfy
	\eqref{eq:global-holder-landau-coefficient-bounds}.  Let
	\[
	f\in\cH_{\kin,\sigma}
	\bigl((-4,0)\times\cD\times\bR_v^3\bigr)
	\]
	satisfy
	\begin{equation}\label{eq:global-holder-original-equation}
		\begin{cases}
			(\partial_t+v\cdot\nabla_x)f
			=D_{v_i}(a_g^{ij}D_{v_j}f)+b_g^iD_{v_i}f+G
			&\text{in }(-4,0]\times\cD\times\bR_v^3,\\
			f=0&\text{on }(-4,0]\times\gamma^-,
		\end{cases}
	\end{equation}
	where the boundary condition is understood in the trace sense of
	Proposition~\ref{2607151218}.
	
	There exists \(\alpha=\alpha(\cD)\in(0,1)\) such that the following
	holds.  If \(L\ge0\) and
	\[
	\|\langle v\rangle^{L+12}f\|_
	{L^\infty((-4,0]\times\cD\times\bR_v^3)}
	+
	\|\langle v\rangle^LG\|_
	{L^\infty((-4,0]\times\cD\times\bR_v^3)}<\infty,
	\]
	then \(f\) has a representative which is continuous on
	\((-1,0]\times\overline\cD\times\bR^3\) and satisfies
	\begin{align}\label{eq:global-weighted-holder-estimate}
		\begin{aligned}
			&\bigl(\langle v\rangle^L+\langle v'\rangle^L\bigr)
			\frac{|f(Z)-f(Z')|}{\mathbf d(Z,Z')^\alpha}
			\\
			\le\,&C_{L,\cD}\left(\|\langle v\rangle^{L+12}f\|_
			{L^\infty((-4,0]\times\cD\times\bR_v^3)}
			+
			\|\langle v\rangle^LG\|_
			{L^\infty((-4,0]\times\cD\times\bR_v^3)}\right)
		\end{aligned}
	\end{align}
	whenever
	\(Z=(t,x,v),Z'=(t',x',v')
	\in(-1,0]\times\overline\cD\times\bR^3\) and
	\(0<\mathbf d(Z,Z')\le1\).
\end{thm}

\begin{proof}	
	We write
	\begin{equation*}
		\mathfrak R_L(f,G)
		:=
		\|\langle v\rangle^{L+12}f\|_
		{L^\infty((-4,0]\times\cD\times\bR_v^3)}
		+
		\|\langle v\rangle^LG\|_
		{L^\infty((-4,0]\times\cD\times\bR_v^3)}
		<\infty.
	\end{equation*}
	
	We first record the time-ordered engulfing property.  There is a universal
	\(C_0\ge1\) such that
	\begin{equation}\label{eq:kinetic-cylinder-engulfing}
		\left.
		\begin{gathered}
			Z=(t,x,v)\in Q_s(Z_*),\qquad
			\bar Z=(\bar t,\bar x,\bar v),\quad \bar t\le t,\\
			\mathbf d(Z,\bar Z)\le C_0^{-1}s
		\end{gathered}
		\right\}
		\quad\Longrightarrow\quad
		Z,\bar Z\in Q_{C_0s}(Z_*).
	\end{equation}
	Indeed, choose an intermediate velocity \(u\) in
	\eqref{kineticdistance} so that each quantity defining
	\(\mathbf d(Z,\bar Z)\) is at most \(2\mathbf d(Z,\bar Z)\).  Then
	\[
	\begin{gathered}
		0\le t_*-\bar t
		=(t_*-t)+(t-\bar t)
		\le s^2+4C_0^{-2}s^2,\\
		|\bar v-v_*|
		\le|\bar v-u|+|u-v|+|v-v_*|
		\le(1+4C_0^{-1})s,
	\end{gathered}
	\]
	while
	\[
	\begin{aligned}
		|\bar x-x_*-(\bar t-t_*)v_*|
		&\le |x-x_*-(t-t_*)v_*|
		+|\bar x-x-(\bar t-t)u|
		+|t-\bar t|\,|u-v_*|\\
		&\le s^3+8C_0^{-3}s^3
		+4C_0^{-2}(1+2C_0^{-1})s^3
		\le Cs^3.
	\end{aligned}
	\]
	Taking the universal constant \(C_0\) sufficiently large proves
	\eqref{eq:kinetic-cylinder-engulfing}; the condition \(\bar t\le t\)
	is used only in the first of these bounds.
	
	Let \(\varepsilon_{\rm b}\) and \(\varepsilon_{\rm i}\) be the constants
	in Lemmas~\ref{landau_local_bdry} and
	\ref{landau_local_interior}.  Fix \(K\ge2C_0\), and choose
	\(m_{\cD}\in\bN\) so that
	\begin{equation}\label{eq:global-holder-choice-m}
		m_{\cD}\ge
		\max\{2K\varepsilon_{\rm b}^{-3},
		4K\varepsilon_{\rm i}^{-2},1\}.
	\end{equation}
	Set \(C_1:=6\) and
	\begin{equation*}
		\delta:=
		\frac{1}{2C_0C_1^{12}m_{\cD}^3}.
	\end{equation*}
	These constants depend only on the geometric constants of \(\cD\).
	Let \(\alpha\in(0,1)\) be the exponent furnished by
	Theorem~\ref{cor:scaled-boundary-holder-estimate} for the uniform
	ellipticity constants in Lemmas~\ref{landau_local_bdry} and
	\ref{landau_local_interior}; thus \(\alpha\) depends only on \(\cD\).
	
	It is enough first to prove the estimate for Lebesgue points of \(f\) in
	the physical interior.
	Fix two Lebesgue points \(Z,Z'\).  By symmetry, assume
	\(t'\le t\), and put \(d=\mathbf d(Z,Z')\).  Since
	\(|v-v'|\le2d\le2\),
	\begin{equation}\label{eq:global-holder-weight-comparison}
		C_1^{-1}\langle v\rangle
		\le\langle v'\rangle
		\le C_1\langle v\rangle.
	\end{equation}
	It suffices to prove
	\begin{equation}\label{eq:global-holder-one-sided}
		\langle v\rangle^L
		\frac{|f(Z)-f(Z')|}{d^\alpha}
		\le C_{L,\cD}\mathfrak R_L(f,G).
	\end{equation}
	
	If
	\begin{equation*}
		d\ge\delta
		\bigl(\langle v\rangle\wedge\langle v'\rangle\bigr)^{-12},
	\end{equation*}
	then \eqref{eq:global-holder-weight-comparison} and \(\alpha<1\) give
	\[
	\langle v\rangle^L
	\frac{|f(Z)-f(Z')|}{d^\alpha}
	\le C_{L,\cD}
	\|\langle v\rangle^{L+12}f\|_{L^\infty}.
	\]
	We henceforth assume
	\begin{equation}\label{eq:global-holder-local-case}
		d<\delta
		\bigl(\langle v\rangle\wedge\langle v'\rangle\bigr)^{-12}.
	\end{equation}
	
	Apply Lemma~\ref{covering_landau} with \(M=m_{\cD}\) to the later point
	\(Z\).
	For some \(j\in\{1,2\}\) and \(k\in\bN\),
	\[
	Z\in Q_{s_{j,k}}(Z_{j,k}),\qquad
	s_{1,k}:=m_{\cD}^{-1}\langle v_{1,k}\rangle^{-6},\qquad
	s_{2,k}:=m_{\cD}^{-3}\langle v_{2,k}\rangle^{-12}.
	\]
	Since $s_{j,k}\le1$ we have $|v-v_{j,k}|<1$, hence
	$\tfrac12\langle v_{j,k}\rangle\le\langle v\rangle
	\le2\langle v_{j,k}\rangle$, while $|v-v'|\le2$ gives
	$\tfrac13\langle v\rangle\le\langle v'\rangle\le3\langle v\rangle$;
	together with the choice $C_1=6$, these two bounds give
	\begin{equation}\label{eq:global-holder-center-weight-comparison}
		C_1^{-1}\langle v_{j,k}\rangle
		\le\langle v\rangle,\langle v'\rangle
		\le C_1\langle v_{j,k}\rangle.
	\end{equation}
	By \eqref{eq:global-holder-local-case} and the choice of \(\delta\),
	\(d\le C_0^{-1}s_{j,k}\).  Since \(t'\le t\),
	\eqref{eq:kinetic-cylinder-engulfing} yields
	\begin{equation}\label{eq:global-holder-common-cylinder}
		Z,Z'\in Q_{C_0s_{j,k}}(Z_{j,k}).
	\end{equation}
	
	Suppose first that \(j=1\).  Let
	\(\psi_k\), \(A_{k,x}=D\psi_k(x)\), \(\cS_k\),
	\(\cT_k\), and \(\widetilde\Phi_{{\rm b},k}\) be the boundary
	chart, normalization matrices, and normalized map associated with the
	center \(Z_{1,k}\) in Lemma~\ref{landau_local_bdry}.  Define
	\[
	Y_{1,k}:=
	(t_{1,k},0,\cT_kA_{k,x_{1,k}}v_{1,k}),
	\qquad
	R_{\rm b}:=
	K\varepsilon_{\rm b}^{-1}m_{\cD}^{-1}\langle v_{1,k}\rangle^{-4}.
	\]
	Then
	\[
	C_0s_{1,k}
	\le\varepsilon_{\rm b}\langle v_{1,k}\rangle^{-3/2}R_{\rm b},
	\qquad
	2R_{\rm b}\le\varepsilon_{\rm b}^2\langle v_{1,k}\rangle^{-4},
	\]
	where the second inequality follows from
	\eqref{eq:global-holder-choice-m}.  Thus
	Lemma~\ref{landau_local_bdry} applies at radius \(2R_{\rm b}\), and by its
	items \ref{item:boundary-transformed-equation} and
	\ref{item:boundary-transformed-coefficients}, it turns
	\eqref{eq:global-holder-original-equation} into an equation of the form
	treated in Theorem~\ref{cor:scaled-boundary-holder-estimate} on
	\(Q_{2R_{\rm b}}^+(Y_{1,k})\), with homogeneous inflow trace and with
	ellipticity constants depending only on \(\cD\).  Moreover, the first
	inclusion in \eqref{260512642}, at radius \(R_{\rm b}\), together with
	\eqref{eq:global-holder-common-cylinder} and the first of the two
	displayed inequalities, places the normalized images \(Y,Y'\) of
	\(Z,Z'\) in \(Q_{R_{\rm b}}^+(Y_{1,k})\).
	Let \(\widetilde b_{{\rm b},k}\) denote the transformed drift in
	\eqref{eq:boundary-transformed-b}, and set
	\[
	F_k:=f\circ\widetilde\Phi_{{\rm b},k},\qquad
	\widetilde G_k:=G\circ\widetilde\Phi_{{\rm b},k}.
	\]
	Explicitly,
	\[
	Y=(t,\cT_k\psi_k(x),\cT_kA_{k,x}v),
	\qquad
	Y'=(t',\cT_k\psi_k(x'),
	\cT_kA_{k,x'}v').
	\]
	
	By \eqref{eq:global-holder-local-case} and
	\eqref{eq:global-holder-weight-comparison},
	$d\le\langle v\rangle^{-2}$.
	Moreover, $x,x'\in B_{R_{\cD}}(x_{1,k})$, by the bound
	$|x-x_{1,k}|\le C\varepsilon_{\rm b}^6<R_{\cD}$ in the proof of
	Lemma~\ref{landau_local_bdry}.  Hence
	Lemma~\ref{260518305}(\ref{item:flattening-distance-comparison}),
	applied to
	\[
	\Theta_k(t,x,v):=(t,\psi_k(x),A_{k,x}v),
	\]
	gives
	\[
	\mathbf d\bigl(\Theta_k(Z),\Theta_k(Z')\bigr)
	\le C_\cD d.
	\]
	Since $Y$ and $Y'$ are obtained by applying $\cT_k$ to both
	the spatial and velocity components of $\Theta_k(Z)$ and
	$\Theta_k(Z')$, respectively, \eqref{kineticdistance} and
	\eqref{eq:boundary-ST-bounds} yield
	\begin{equation}\label{eq:global-boundary-distance-comparison}
		\mathbf d(Y,Y')
		\le\max\{1,|\cT_k|\}
		\mathbf d\bigl(\Theta_k(Z),\Theta_k(Z')\bigr)
		\le C_\cD\langle v_{1,k}\rangle^{3/2}d.
	\end{equation}
	
	The normalized drift bound in Lemma~\ref{landau_local_bdry} and the
	definition of \(R_{\rm b}\) imply
	\[
	c_{\cD}\langle v_{1,k}\rangle^4
	\le
	\|\widetilde b_{{\rm b},k}\|_{L^\infty(Q_{2R_{\rm b}}(Y_{1,k}))}
	+R_{\rm b}^{-1}
	\le C_\cD\langle v_{1,k}\rangle^4.
	\]
	Theorem~\ref{cor:scaled-boundary-holder-estimate}, followed by continuous
	extension to the flat boundary, therefore gives
	\begin{align*}
		[F_k]_{C_{\kin}^\alpha(Q_{R_{\rm b}}^+(Y_{1,k}))}
		\le C_\cD\Bigl(
		\langle v_{1,k}\rangle^{4\alpha}
		\|F_k\|_{L^\infty(Q_{2R_{\rm b}}^+(Y_{1,k}))}
		+\langle v_{1,k}\rangle^{-4(2-\alpha)}
		\|\widetilde G_k\|_
		{L^\infty(Q_{2R_{\rm b}}^+(Y_{1,k}))}
		\Bigr).
	\end{align*}
	Using \eqref{eq:global-boundary-distance-comparison},
	\eqref{eq:global-holder-center-weight-comparison}, and
	\(0<\alpha<1\), we obtain
	\[
	\langle v\rangle^L
	\frac{|f(Z)-f(Z')|}{d^\alpha}
	\le C_{L,\cD}\Bigl(
	\langle v_{1,k}\rangle^{L+6}
	\|F_k\|_{L^\infty(Q_{2R_{\rm b}}^+(Y_{1,k}))}
	+\langle v_{1,k}\rangle^L
	\|\widetilde G_k\|_
	{L^\infty(Q_{2R_{\rm b}}^+(Y_{1,k}))}
	\Bigr)
	\le C_{L,\cD}\mathfrak R_L(f,G).
	\]
	The last inequality uses the outer-cylinder inclusion in
	Lemma~\ref{landau_local_bdry} and
	\(\langle v\rangle\asymp\langle v_{1,k}\rangle\) on its image.
	
	Suppose now that \(j=2\), and let
	\[
	\cS_k:=\sigma(v_{2,k})^{1/2},\qquad
	\cT_k:=\cS_k^{-1},\qquad
	\widetilde\Phi_{{\rm i},k}(s,y,w)
	:=(s,\cS_ky,\cS_kw).
	\]
	The normalized center and the two relevant radii are
	\[
	Y_{2,k}:=(t_{2,k},\cT_kx_{2,k},\cT_kv_{2,k}),
	\qquad
	R_{\rm i}:=
	K\varepsilon_{\rm i}^{-1}m_{\cD}^{-3}\langle v_{2,k}\rangle^{-10},\qquad
	R_{0,k}:=\min\{1,\operatorname{dist}(x_{2,k},\partial\cD)\}.
	\]
	By \eqref{260522351}
	(with $\operatorname{dist}$ replaced by $R_{0,k}$, the case
	$\operatorname{dist}\ge1$ being immediate),
	\[
	\sqrt{\frac{R_{0,k}}{2\langle v_{2,k}\rangle}}
	\ge\frac{\langle v_{2,k}\rangle^{-19/2}}{2m_{\cD}^{3/2}}
	\ge\frac{\langle v_{2,k}\rangle^{-10}}{2m_{\cD}^2}.
	\]
	Consequently,
	\[
	C_0s_{2,k}
	\le\varepsilon_{\rm i}\langle v_{2,k}\rangle^{-3/2}R_{\rm i},
	\qquad
	2R_{\rm i}
	\le\varepsilon_{\rm i}\sqrt{\frac{R_{0,k}}{2\langle v_{2,k}\rangle}},
	\]
	so Lemma~\ref{landau_local_interior} applies at radius \(2R_{\rm i}\);
	by its items \ref{item:interior-transformed-equation} and
	\ref{item:interior-transformed-coefficients}, the transformed equation has
	universal ellipticity constants, and the first inclusion in
	\eqref{260512642_11} gives \(Y,Y'\in Q_{R_{\rm i}}(Y_{2,k})\).
	Let \(\widetilde b_{{\rm i},k}\) denote the transformed drift in
	\eqref{eq:interior-transformed-coefficients}, and set
	\[
	F_k:=f\circ\widetilde\Phi_{{\rm i},k},\qquad
	\widetilde G_k:=G\circ\widetilde\Phi_{{\rm i},k},
	\]
	\[
	Y=(t,\cT_kx,\cT_kv),\qquad
	Y'=(t',\cT_kx',\cT_kv').
	\]
	Moreover,
	\begin{equation}\label{eq:global-interior-distance-comparison}
		\mathbf d(Y,Y')\le C\langle v_{2,k}\rangle^{3/2}d,
	\end{equation}
	which follows directly from \eqref{kineticdistance} and
	\(|\cT_k|\le C\langle v_{2,k}\rangle^{3/2}\).
	
	Here
	\[
	c_{\cD}\langle v_{2,k}\rangle^{10}
	\le
	\|\widetilde b_{{\rm i},k}\|_{L^\infty(Q_{2R_{\rm i}}(Y_{2,k}))}
	+R_{\rm i}^{-1}
	\le C_\cD\langle v_{2,k}\rangle^{10}.
	\]
	After a constant translation in \(y\), the full cylinder lies in
	\(\bR_+^3\), and such a translation preserves
	\(\partial_t+w\cdot\nabla_y\) and the velocity-divergence operator.
	Hence Theorem~\ref{cor:scaled-boundary-holder-estimate} applies with an
	empty boundary portion and yields
	\begin{align*}
		[F_k]_{C_{\kin}^\alpha(Q_{R_{\rm i}}(Y_{2,k}))}
		\le C_\cD\Bigl(
		\langle v_{2,k}\rangle^{10\alpha}
		\|F_k\|_{L^\infty(Q_{2R_{\rm i}}(Y_{2,k}))}
		+\langle v_{2,k}\rangle^{-10(2-\alpha)}
		\|\widetilde G_k\|_
		{L^\infty(Q_{2R_{\rm i}}(Y_{2,k}))}
		\Bigr).
	\end{align*}
	Together with \eqref{eq:global-interior-distance-comparison}, this gives
	\[
	\langle v\rangle^L
	\frac{|f(Z)-f(Z')|}{d^\alpha}
	\le C_{L,\cD}\Bigl(
	\langle v_{2,k}\rangle^{L+12}
	\|F_k\|_{L^\infty(Q_{2R_{\rm i}}(Y_{2,k}))}
	+\langle v_{2,k}\rangle^L
	\|\widetilde G_k\|_
	{L^\infty(Q_{2R_{\rm i}}(Y_{2,k}))}
	\Bigr)
	\le C_{L,\cD}\mathfrak R_L(f,G).
	\]
	Here we used the outer-cylinder inclusion in
	Lemma~\ref{landau_local_interior} and the corresponding velocity
	comparability.
	This proves \eqref{eq:global-holder-one-sided} in all cases.
	Equation~\eqref{eq:global-holder-weight-comparison} gives
	\eqref{eq:global-weighted-holder-estimate}.
	
	Finally, \eqref{eq:global-holder-one-sided} holds, with a constant
	independent of the pair, for any two Lebesgue points of \(f\) in
	\((-1,0]\times\cD\times\bR^3\) at \(\mathbf d\)-distance at most
	\(1\).  On each velocity-bounded subset of the Lebesgue set of \(f\)
	this makes \(f\) uniformly \(\mathbf d\)-continuous, and that set is
	dense in \((-1,0]\times\overline\cD\times\bR^3\), since it has full
	measure in \((-1,0)\times\cD\times\bR^3\) and \(\mathbf d\) induces
	the Euclidean topology.  Hence \(f\) extends uniquely and continuously
	to \((-1,0]\times\overline\cD\times\bR^3\).  This extension is the
	asserted representative, and \eqref{eq:global-weighted-holder-estimate}
	passes to it because both of its sides are continuous in \((Z,Z')\).
\end{proof}

\begin{lemma}[H\"older continuity of \(\sigma_{\sqrt\mu g}\)]
	\label{260527138}
	Let \(0<\alpha<1\) and
	\[
	\mu^{1/4}g\in
	C_{\kin}^\alpha(\bR\times\bR_x^3\times\bR_v^3),
	\]
	and write
	\(\sigma_{\sqrt\mu g}^{ij}(Z)
	:=\bigl(\phi^{ij}\ast[\sqrt\mu\,g(t,x,\cdot)]\bigr)(v)\)
	for \(Z=(t,x,v)\).  Then
	\begin{equation}\label{eq:holder-sigma-f-mu}
		\left|
		\sigma_{\sqrt\mu g}^{ij}(Z)
		-\sigma_{\sqrt\mu g}^{ij}(Z')
		\right|
		\le C_\alpha
		\|\mu^{1/4}g\|_{C_{\kin}^\alpha}
		\mathbf d(Z,Z')^{2\alpha/3}
	\end{equation}
	for \(i,j=1,2,3\) and all
	\(Z,Z'\in\bR\times\bR_x^3\times\bR_v^3\); equivalently,
	\(\bigl[\sigma_{\sqrt\mu g}^{ij}\bigr]_{C_{\kin}^{2\alpha/3}}
	\le C_\alpha\|\mu^{1/4}g\|_{C_{\kin}^\alpha}\), the norm being
	\eqref{eq:Ckin-norm}.  Its \(L^\infty\)-part cannot be dropped: for
	constant \(\mu^{1/4}g\) the seminorm vanishes while
	\(\sigma_{\sqrt\mu g}^{ij}\) does not.
\end{lemma}

\begin{proof}
	Set
	\[
	H:=\mu^{1/4}g,\qquad
	M:=\|H\|_{C_{\kin}^\alpha},\qquad
	d:=\mathbf d(Z,Z').
	\]
	The case \(d=0\) is immediate.  For \(0<d\le1\), choose
	\(\overline v\in\bR^3\) such that
	\[
	\max\bigl\{
	|t-t'|^{1/2},
	|x-x'-(t-t')\overline v|^{1/3},
	|v-\overline v|,
	|v'-\overline v|
	\bigr\}\le2d.
	\]
	For \(w\in\bR^3\), put
	\[
	Z_w:=(t,x,v-w),\qquad
	Z'_w:=(t',x',v'-w),\qquad
	d_w:=\mathbf d(Z_w,Z'_w).
	\]
	Using \(\overline v-w\) as the intermediate velocity in
	\eqref{kineticdistance}, we obtain
	\begin{equation}\label{260514852}
		d_w
		\le
		\max\left\{2d,
		(8d^3+4d^2|w|)^{1/3}\right\}
		\le C d^{2/3}\langle w\rangle^{1/3}.
	\end{equation}
	
	Since \(\sqrt\mu\,g=\mu^{1/4}H\),
	\[
	\sigma_{\sqrt\mu g}^{ij}(Z)
	=
	\int_{\bR^3}\phi^{ij}(w)
	H(Z_w)\mu^{1/4}(v-w)\dd w.
	\]
	Therefore,
	\begin{align*}
		\left|
		\sigma_{\sqrt\mu g}^{ij}(Z)
		-\sigma_{\sqrt\mu g}^{ij}(Z')
		\right|
		\le{}&
		\int_{\bR^3}|\phi^{ij}(w)|
		|H(Z_w)-H(Z'_w)|\mu^{1/4}(v-w)\dd w
		\notag\\
		&+
		\int_{\bR^3}|\phi^{ij}(w)|\,|H(Z'_w)|
		\left|\mu^{1/4}(v-w)-\mu^{1/4}(v'-w)\right|
		\dd w.
	\end{align*}
	By \eqref{260514852},
	\[
	|H(Z_w)-H(Z'_w)|
	\le CMd^{2\alpha/3}\langle w\rangle^{\alpha/3},
	\]
	and \(|\phi^{ij}(w)|\le C|w|^{-1}\).  Hence the first integral is
	bounded by
	\[
	CMd^{2\alpha/3}
	\int_{\bR^3}|w|^{-1}\langle w\rangle^{\alpha/3}
	\mu^{1/4}(v-w)\dd w
	\le C_\alpha Md^{2\alpha/3}.
	\]
	Here the last convolution is uniformly bounded in \(v\), since
	\(|w|^{-1}\) is locally integrable in dimension three and
	\(\langle w\rangle^{\alpha/3}/|w|\) is bounded for \(|w|\ge1\).
	
	For the second integral, write
	\[
	\mu^{1/4}(v-w)-\mu^{1/4}(v'-w)
	=(v-v')\cdot\int_0^1\nabla\mu^{1/4}
	\bigl(v'+\theta(v-v')-w\bigr)\dd\theta ,
	\]
	so that, by \(|H|\le M\), Fubini's theorem, and \(|v-v'|\le2d\),
	\[
	\begin{aligned}
		&\int_{\bR^3}|\phi^{ij}(w)|\,|H(Z'_w)|
		\left|\mu^{1/4}(v-w)-\mu^{1/4}(v'-w)\right|\dd w\\
		&\qquad\le
		CM|v-v'|
		\sup_{z\in\bR^3}
		\int_{\bR^3}|w|^{-1}|\nabla\mu^{1/4}(z-w)|\dd w
		\le CMd.
	\end{aligned}
	\]
	The supremum is finite by the same splitting as above, since
	\(\nabla\mu^{1/4}\in L^1\cap L^\infty(\bR^3)\).
	Since \(d\le d^{2\alpha/3}\), the two bounds prove
	\eqref{eq:holder-sigma-f-mu} in the range \(0<d\le1\).
	
	For \(d>1\) it suffices to use \(|H|\le M\) and the same convolution
	bound, giving \(|\sigma_{\sqrt\mu g}^{ij}|\le CM\) pointwise and
	hence
	\(|\sigma_{\sqrt\mu g}^{ij}(Z)-\sigma_{\sqrt\mu g}^{ij}(Z')|
	\le2CM\le2CMd^{2\alpha/3}\).  This proves
	\eqref{eq:holder-sigma-f-mu} for all \(Z,Z'\).
\end{proof}

\subsection{Weighted velocity-gradient and mixed-norm estimates}\label{sec:C1v-estimate}

\begin{thm}[Global weighted velocity-gradient estimate]
	\label{260520311}
	Let $0<\alpha<1$ and $0<\delta<1/2$.  Let $a$ and $b$ be
	defined on $\bR\times\bR_x^3\times\bR_v^3$, with $a$ symmetric and
	$\nabla_va$ defined weakly.  Assume that, for some
	$K_0,K_\alpha\ge1$, a.e. $Z=(t,x,v)$ and every $\xi\in\bR^3$ satisfy
	\begin{equation}\label{eq:global-gradient-structural-assumption}
		\frac12\sigma^{ij}(v)\xi_i\xi_j
		\le a^{ij}(Z)\xi_i\xi_j
		\le2\sigma^{ij}(v)\xi_i\xi_j,
		\qquad
		|\nabla_va(Z)|\le K_0\langle v\rangle^{-2},
		\qquad
		|b(Z)|\le K_0\langle v\rangle^{-1},
	\end{equation}
	\begin{equation}\label{eq:global-gradient-holder-assumption}
		[a]_{C_{\kin}^\alpha
			(\bR\times\bR_x^3\times\bR_v^3)}\le K_\alpha.
	\end{equation}
	Set
	\begin{equation}\label{eq:global-gradient-weight-loss}
		A=A(\alpha)
		:=\left\lceil
		\max\left\{
		\frac32\max\left\{4,\frac3\alpha\right\},
		3+\max\left\{10,\frac3\alpha\right\}
		\right\}
		\right\rceil.
	\end{equation}
	Suppose that
	\[
	f\in\cH_{\kin,\sigma}
	\bigl((-4,0)\times\cD\times\bR_v^3\bigr)
	\]
	satisfies
	\begin{equation}\label{eq:global-weighted-gradient-equation}
		\begin{cases}
			(\partial_t+v\cdot\nabla_x)f
			=D_{v_i}(a^{ij}D_{v_j}f)+b^iD_{v_i}f+G
			&\text{in }(-4,0]\times\cD\times\bR_v^3,\\
			f=0&\text{on }(-4,0]\times\gamma^- ,
		\end{cases}
	\end{equation}
	where the boundary condition is understood in the trace sense of
	Proposition~\ref{2607151218}.  Then, for every $L\ge0$,
	\begin{align}\label{eq:global-weighted-gradient-estimate}
		&\left\|
		\langle v\rangle^L
		\big(\mathbf d(\cdot,\bR\times\gamma_{\cD}^0)\wedge 1\big)^{1/2+\delta}D_vf
		\right\|_{L^\infty((-1,0]\times\cD\times\bR_v^3)}
		\notag\\
		&\qquad\le C\left(
		\|\langle v\rangle^{L+A}f\|_
		{L^\infty((-4,0]\times\cD\times\bR_v^3)}
		+
		\|\langle v\rangle^LG\|_
		{L^\infty((-4,0]\times\cD\times\bR_v^3)}
		\right),
	\end{align}
	provided the right-hand side is finite.  Here
	$C=C(\alpha,\delta,L,K_0,K_\alpha,\cD)$.
\end{thm}

\begin{proof}
	For the proof, set
	\[
	q_{\rm b}:=\max\left\{4,\frac3\alpha\right\},
	\qquad
	q_{\rm i}:=\max\left\{10,\frac3\alpha\right\},
	\qquad
	\eta:=\frac12+\delta
	\]
	and denote the right-hand side of
	\eqref{eq:global-weighted-gradient-estimate} by
	\begin{align*}
		\cN_L:={}&
		\|\langle v\rangle^{L+A}f\|_
		{L^\infty((-4,0]\times\cD\times\bR_v^3)}
		+
		\|\langle v\rangle^LG\|_
		{L^\infty((-4,0]\times\cD\times\bR_v^3)}.
	\end{align*}
	
	We first fix a localization scale that works in every normalization chart.
	Let \(\varepsilon_{\rm b}\) and \(\varepsilon_{\rm i}\) be the constants in
	Lemmas~\ref{landau_local_bdry} and \ref{landau_local_interior}, and let
	\(c_0\) be the constant in \eqref{eq:flattening-grazing-range}.  Choose \(m=m_{\cD}\in\bN\) so large that
	\begin{equation}\label{eq:global-gradient-choice-m}
		\frac{20}{m}\le\varepsilon_{\rm b}^3,\qquad
		\frac4m\le\varepsilon_{\rm i}^2,\qquad
		\frac{C_\cD}{m^2}\le
		\frac{c_0}2R_{\cD},\qquad
		\frac{C_\cD(1+\varepsilon_{\rm b}^{-2})}{m}
		\le\frac{c_0}2.
	\end{equation}
	The first two inequalities in \eqref{eq:global-gradient-choice-m} make
	the normalized cylinders admissible.  The last two allow us to compare
	the physical and the flattened grazing distances.
	
	Fix
	\[
	Z=(t,x,v)\in(-1,0]\times\cD\times\bR_v^3,
	\qquad
	d(Z):=\mathbf d(Z,\bR\times\gamma_{\cD}^0),
	\qquad
	d_*(Z):=d(Z)\wedge1.
	\]
	By Lemma~\ref{covering_landau}, \(Z\) lies in either a boundary
	covering cylinder or an interior covering cylinder; in both cases we
	write \(v_*\) for the velocity at the center of that cylinder, and
	\(\widetilde v\) for the velocity coordinate of a point
	in such an enlarged cylinder.  Then
	\begin{equation}\label{eq:global-gradient-weight-comparison}
		\langle \widetilde v\rangle\asymp_{\cD}\langle v_*\rangle
	\end{equation}
	throughout every enlarged cylinder used below.  The proof therefore
	reduces to three local estimates: a non-grazing boundary estimate, a
	grazing boundary estimate, and an interior estimate.
	
	\medskip
	\noindent\emph{Step 1. Boundary cylinders.}
	Suppose that
	\[
	Z\in Q_{m^{-1}\langle v_{\rm b}\rangle^{-6}}(Z_{\rm b}),
	\qquad
	Z_{\rm b}=(t_{\rm b},x_{\rm b},v_{\rm b}),
	\qquad x_{\rm b}\in\partial\cD,
	\qquad v_*=v_{\rm b}.
	\]
	Let \(\psi_{\rm b}\), \(A_{{\rm b},x}=D\psi_{\rm b}(x)\),
	\(\cS_{\rm b}\), \(\cT_{\rm b}\), and
	\(\widetilde\Phi_{\rm b}\) be the chart, normalization matrices, and
	normalized map supplied by Lemma~\ref{landau_local_bdry} at
	\(Z_{\rm b}\), and set
	\[
	r:=\varepsilon_{\rm b}^{-1}m^{-1}\langle v_{\rm b}\rangle^{-4}.
	\]
	The choice of \(m\) gives
	\[
	m^{-1}\langle v_{\rm b}\rangle^{-6}
	\le\varepsilon_{\rm b}\langle v_{\rm b}\rangle^{-3/2}r,
	\qquad
	20r\le\varepsilon_{\rm b}^2\langle v_{\rm b}\rangle^{-4}.
	\]
	Hence, writing
	\[
	Y:=\widetilde\Phi_{\rm b}^{-1}(Z),\qquad
	Y_{\rm b}:=
	(t_{\rm b},0,\cT_{\rm b}A_{{\rm b},x_{\rm b}}v_{\rm b}),
	\]
	and
	\[
	F:=f\circ\widetilde\Phi_{\rm b},
	\qquad
	H:=G\circ\widetilde\Phi_{\rm b},
	\]
	we have
	\begin{equation*}
		Y\in Q_r^+(Y_{\rm b}),\qquad
		\widetilde\Phi_{\rm b}(Q_{20r}^+(Y_{\rm b}))
		\subset Q_1^\cD(Z_{\rm b})
		\subset(-2,0]\times\cD\times\bR_v^3.
	\end{equation*}
	The chain rule and \eqref{eq:boundary-ST-bounds} give
	\begin{equation}\label{eq:global-boundary-gradient-transform}
		|D_vf(Z)|\le C_\cD\langle v_{\rm b}\rangle^{3/2}|D_wF(Y)|.
	\end{equation}
	
	The transformed coefficient bounds in
	Lemma~\ref{landau_local_bdry}, together with
	\(q_{\rm b}\ge4\) and \(\alpha q_{\rm b}\ge3\), allow us to use
	\[
	M_1=C_\cD(K_\alpha+1)\langle v_{\rm b}\rangle^{\alpha q_{\rm b}},
	\qquad
	M_2=C_\cD(K_0+1)\langle v_{\rm b}\rangle^{q_{\rm b}}
	\]
	in Theorem~\ref{prop:scaled-Dv-estimate}.  Thus, with the notation
	used there,
	\[
	M:=\max\{M_1^{1/\alpha},M_2\},
	\qquad M+r^{-1}\asymp\langle v_{\rm b}\rangle^{q_{\rm b}},
	\]
	after enlarging the admissible bounds by fixed factors if necessary.
	
	If \(Q_{2r}(Y_{\rm b})\) does not meet the flat grazing set, the
	non-grazing part of Theorem~\ref{prop:scaled-Dv-estimate} yields
	\begin{equation*}
		\|D_wF\|_{L^\infty(Q_r^+(Y_{\rm b}))}
		\le C\left(
		\langle v_{\rm b}\rangle^{q_{\rm b}}\|F\|_{L^\infty(Q_{2r}^+(Y_{\rm b}))}
		+\langle v_{\rm b}\rangle^{-q_{\rm b}}\|H\|_{L^\infty(Q_{2r}^+(Y_{\rm b}))}
		\right).
	\end{equation*}
	
	Assume now that \(Q_{2r}(Y_{\rm b})\) meets the flat grazing set.
	Lemma~\ref{260515154} provides
	\(\overline Y\in\bR\times\gamma_{\bR_+^3}^0\) such that
	\begin{equation}\label{eq:global-grazing-cylinder-comparison}
		Q_{2r}(Y_{\rm b})\subset Q_{4r}(\overline Y),
		\qquad
		Q_{16r}(\overline Y)\subset Q_{20r}(Y_{\rm b}).
	\end{equation}
	By Lemma~\ref{landau_local_bdry}(\ref{item:boundary-cylinder-comparison}),
	the point \(\widetilde\Phi_{\rm b}(\overline Y)\) lies in
	\(\bR\times\gamma_{\cD}^0\).
	Moreover, \eqref{eq:global-grazing-cylinder-comparison} and
	\eqref{260512642} at radius \(20r\) place both \(Z\) and
	\(\widetilde\Phi_{\rm b}(\overline
	Y)\) in \(Q_{20\varepsilon_{\rm b}^{-1}r}(Z_{\rm b})\).  Since
	\(\mathbf d\) is a metric and Lemma~\ref{260127226} bounds the
	\(\mathbf d\)-radius of \(Q_s\) by \(s\), we obtain
	\[
	d(Z)\le\mathbf d\bigl(Z,\widetilde\Phi_{\rm b}(\overline Y)\bigr)
	\le40\varepsilon_{\rm b}^{-1}r=:C_\cD r .
	\]
	This bound is available only in the present sub-case; for a general
	boundary point, \(d(Z)\) need not be small.  With the covering radius and
	\eqref{eq:global-gradient-choice-m},
	\[
	d(Z)\le C_\cD r
	\le\frac{c_0}2\langle v\rangle^{-2},
	\qquad
	|x-x_{\rm b}|
	\le Cm^{-2}\langle v_{\rm b}\rangle^{-11}
	\le\frac{c_0}2R_{\cD}.
	\]
	Therefore, \(d_*(Z)=d(Z)\), \eqref{eq:flattening-grazing-range} holds, and
	Lemma~\ref{260518305}(\ref{item:flattening-grazing-distance}) applies
	at \(Z\) and gives
	\begin{equation}\label{260805433}
		d(Z)\le C_\cD
		\mathbf d\bigl(\Theta_{\rm b}(Z),
		\bR\times\gamma_{\bR_+^3}^0\bigr),
		\qquad
		\Theta_{\rm b}(t,x,v):=(t,\psi_{\rm b}(x),A_{{\rm b},x}v).
	\end{equation}
	To pass from \(\Theta_{\rm b}(Z)\) to \(Y\), set
	\(L_{\cS}(t,y,w):=(t,\cS_{\rm b}y,\cS_{\rm b}w)\),
	so that \(\Theta_{\rm b}(Z)=L_{\cS}(Y)\).  By
	\eqref{eq:boundary-normal-direction},
	\((\cS_{\rm b}y)_3=\kappa_*y_3\) and
	\((\cS_{\rm b}w)_3=\kappa_*w_3\).  Therefore, \(L_{\cS}\) maps
	\(\bR\times\gamma_{\bR_+^3}^0\) into itself.  Moreover,
	\(|\cS_{\rm b}|\le C_\cD\) by
	\eqref{eq:boundary-ST-bounds}, and \eqref{kineticdistance} then gives
	\(\mathbf d(L_{\cS}U,L_{\cS}U')\le
	C_\cD\mathbf d(U,U')\).  Hence
	\(\mathbf d(\Theta_{\rm b}(Z),\bR\times\gamma_{\bR_+^3}^0)\le
	C_\cD\mathbf d(Y,\bR\times\gamma_{\bR_+^3}^0)\), and
	\eqref{260805433} yields
	\begin{equation}\label{eq:global-normalized-grazing-comparison}
		d(Z)\le C_\cD
		\mathbf d(Y,\bR\times\gamma_{\bR_+^3}^0).
	\end{equation}
	Only this direction is available, and only it is used: the reverse
	comparison would need \(|\cT_{\rm b}|\le C_\cD\), whereas
	\(|\cT_{\rm b}|\asymp\langle v_{\rm b}\rangle^{3/2}\).
	We now apply the grazing part of
	Theorem~\ref{prop:scaled-Dv-estimate} at \(\overline Y\) with
	radius \(8r\).  Together with
	\(r\asymp_{\cD}\langle v_{\rm b}\rangle^{-4}\) and
	\(M+r^{-1}\asymp\langle v_{\rm b}\rangle^{q_{\rm b}}\), this gives
	\begin{align*}
		&\left\|
		\mathbf d(\cdot,\bR\times\gamma_{\bR_+^3}^0)^\eta D_wF
		\right\|_{L^\infty(Q_{4r}^+(\overline Y))}
		\\
		&\qquad\le C\left(
		\langle v_{\rm b}\rangle^{3q_{\rm b}/2-2-4\delta}
		\|F\|_{L^\infty(Q_{20r}^+(Y_{\rm b}))}
		+
		\langle v_{\rm b}\rangle^{-q_{\rm b}/2-2-4\delta}
		\|H\|_{L^\infty(Q_{20r}^+(Y_{\rm b}))}
		\right).
	\end{align*}
	
	\medskip
	\noindent\emph{Step 2. Interior cylinders.}
	Suppose instead that the interior alternative in
	Lemma~\ref{covering_landau} holds.  Thus there exists an interior
	covering center
	\[
	Z_{\rm i}=(t_{\rm i},x_{\rm i},v_{\rm i}),
	\qquad x_{\rm i}\in\cD,
	\]
	such that \(v_*=v_{\rm i}\) and
	\[
	Z\in Q_{m^{-3}\langle v_{\rm i}\rangle^{-12}}(Z_{\rm i}).
	\]
	Set
	\[
	r:=\varepsilon_{\rm i}^{-1}m^{-3}\langle v_{\rm i}\rangle^{-10}.
	\]
	The interior part of Lemma~\ref{covering_landau} and
	\eqref{eq:global-gradient-choice-m} give
	\[
	m^{-3}\langle v_{\rm i}\rangle^{-12}
	\le\varepsilon_{\rm i}\langle v_{\rm i}\rangle^{-3/2}r,
	\qquad
	2r\le\varepsilon_{\rm i}
	\sqrt{\frac{\min\{1,\operatorname{dist}(x_{\rm i},\partial\cD)\}}
		{2\langle v_{\rm i}\rangle}}.
	\]
	Let
	\[
	\cS_{\rm i}:=\sigma(v_{\rm i})^{1/2},\qquad
	\cT_{\rm i}:=\cS_{\rm i}^{-1},\qquad
	\widetilde\Phi_{\rm i}(s,y,w)
	:=(s,\cS_{\rm i}y,\cS_{\rm i}w)
	\]
	be the normalization in Lemma~\ref{landau_local_interior}, and define
	\[
	Y:=\widetilde\Phi_{\rm i}^{-1}(Z),\qquad
	Y_{\rm i}:=
	(t_{\rm i},\cT_{\rm i}x_{\rm i},
	\cT_{\rm i}v_{\rm i}),
	\qquad
	F:=f\circ\widetilde\Phi_{\rm i},\qquad
	H:=G\circ\widetilde\Phi_{\rm i}.
	\]
	Then
	\begin{equation*}
		Y\in Q_r(Y_{\rm i}),\qquad
		\widetilde\Phi_{\rm i}(Q_{2r}(Y_{\rm i}))
		\subset Q_1^\cD(Z_{\rm i})
		\subset(-2,0]\times\cD\times\bR_v^3,
	\end{equation*}
	and
	\begin{equation}\label{eq:global-interior-gradient-transform}
		|D_vf(Z)|\le C\langle v_{\rm i}\rangle^{3/2}|D_wF(Y)|.
	\end{equation}
	Since \(q_{\rm i}\ge10\) and \(\alpha q_{\rm i}\ge3\), we may take
	\[
	M_1=C(K_\alpha+1)\langle v_{\rm i}\rangle^{\alpha q_{\rm i}},
	\qquad
	M_2=C(K_0+1)\langle v_{\rm i}\rangle^{q_{\rm i}}
	\]
	in Theorem~\ref{prop:scaled-Dv-estimate}.  The effective scale
	\(\max\{M_1^{1/\alpha},M_2\}+r^{-1}\) is then comparable to
	\(\langle v_{\rm i}\rangle^{q_{\rm i}}\).
	The non-grazing estimate in Theorem~\ref{prop:scaled-Dv-estimate} gives
	\begin{equation*}
		\|D_wF\|_{L^\infty(Q_r(Y_{\rm i}))}
		\le C\left(
		\langle v_{\rm i}\rangle^{q_{\rm i}}\|F\|_{L^\infty(Q_{2r}(Y_{\rm i}))}
		+\langle v_{\rm i}\rangle^{-q_{\rm i}}\|H\|_{L^\infty(Q_{2r}(Y_{\rm i}))}
		\right).
	\end{equation*}
	
	\medskip
	\noindent\emph{Step 3. Return to the global weights.}
	By \eqref{eq:global-gradient-weight-comparison}, on every image
	cylinder used above,
	\[
	\|F\|_{L^\infty}
	\le C\langle v_*\rangle^{-L-A}
	\|\langle v\rangle^{L+A}f\|_{L^\infty((-4,0]\times\cD\times\bR_v^3)},
	\]
	and
	\[
	\|H\|_{L^\infty}
	\le C\langle v_*\rangle^{-L}
	\|\langle v\rangle^LG\|_{L^\infty((-4,0]\times\cD\times\bR_v^3)}.
	\]
	Combining these bounds with the chain-rule estimates
	\eqref{eq:global-boundary-gradient-transform} and
	\eqref{eq:global-interior-gradient-transform} gives,
	respectively in the boundary non-grazing, boundary grazing, and
	interior cases,
	\begin{align*}
		\langle v\rangle^Ld_*(Z)^\eta|D_vf(Z)|
		&\le C\cN_L
		\left(\langle v_{\rm b}\rangle^{q_{\rm b}+3/2-A}
		+\langle v_{\rm b}\rangle^{-q_{\rm b}+3/2}\right),\\
		\langle v\rangle^Ld_*(Z)^\eta|D_vf(Z)|
		&\le C\cN_L
		\left(
		\langle v_{\rm b}\rangle^{3q_{\rm b}/2-1/2-4\delta-A}
		+\langle v_{\rm b}\rangle^{-q_{\rm b}/2-1/2-4\delta}
		\right),\\
		\langle v\rangle^Ld_*(Z)^\eta|D_vf(Z)|
		&\le C\cN_L
		\left(\langle v_{\rm i}\rangle^{q_{\rm i}+3/2-A}
		+\langle v_{\rm i}\rangle^{-q_{\rm i}+3/2}\right).
	\end{align*}
	Here the first and third estimates use \(d_*\le1\), while the second
	uses \eqref{eq:global-normalized-grazing-comparison}.
	
	Finally, \(q_{\rm b}\le q_{\rm i}\), \(q_{\rm b}\ge4\),
	\(q_{\rm i}\ge10\), and the definition of \(A\) in
	\eqref{eq:global-gradient-weight-loss} imply
	\[
	q_{\rm b}+\frac32\le A,\qquad
	\frac32q_{\rm b}-\frac12-4\delta\le A,\qquad
	q_{\rm i}+\frac32\le A,
	\]
	as well as
	\[
	-q_{\rm b}+\frac32\le0,\qquad
	-\frac12q_{\rm b}-\frac12-4\delta\le0,\qquad
	-q_{\rm i}+\frac32\le0.
	\]
	Every power of \(\langle v_*\rangle\) in the preceding three estimates
	is therefore nonpositive, and hence
	\[
	\langle v\rangle^Ld_*(Z)^\eta|D_vf(Z)|
	\le C\cN_L.
	\]
	Taking the essential supremum over \(Z\) proves
	\eqref{eq:global-weighted-gradient-estimate}.
\end{proof}

For the Landau coefficients \(a=a_g\) and \(b=b_g\),
\eqref{eq:global-gradient-structural-assumption} follows from
\eqref{eq:global-holder-landau-coefficient-bounds}.  If
\(\mu^{1/4}g\in C_{\kin}^{\alpha_0}\), then
Lemma~\ref{260527138} and
\eqref{eq:derivative-sigma-bound} give
\eqref{eq:global-gradient-holder-assumption} with
\(\alpha=2\alpha_0/3\) and
\[
K_\alpha\le C_{\alpha_0}
\left(1+\|\mu^{1/4}g\|_
{C_{\kin}^{\alpha_0}(\bR\times\bR_x^3\times\bR_v^3)}
\right).
\]
Both gradient estimates in this subsection therefore apply to the Landau
coefficients with \(\alpha=2\alpha_0/3\).

Below, $\|\cdot\|_{L^{\infty,p,2}_{t,x,v}}$ denotes the iterated norm
$L_t^\infty\bigl(L_x^p(L_v^2)\bigr)$, the velocity integration being taken
first.

\begin{thm}[Global mixed-norm velocity-gradient estimate]
	\label{260527139}
	Let $0<\alpha<1$, and assume that $a$, $b$, $f$ and $G$ satisfy
	\eqref{eq:global-gradient-structural-assumption},
	\eqref{eq:global-gradient-holder-assumption} and
	\eqref{eq:global-weighted-gradient-equation}, as in
	Theorem~\ref{260520311}.  Let $A=A(\alpha)$ be defined by
	\eqref{eq:global-gradient-weight-loss}.  Then, for every
	$1<p<\infty$ and $L\ge0$,
	\begin{align}\label{eq:global-mixed-gradient-estimate}
		&\|\langle v\rangle^LD_vf\|_
		{L^{\infty,p,2}_{t,x,v}((-1,0]\times\cD\times\bR_v^3)}
		\notag\\
		&\qquad\le C\left(
		\|\langle v\rangle^{L+A+3}f\|_
		{L^\infty((-4,0]\times\cD\times\bR_v^3)}
		+
		\|\langle v\rangle^{L+3}G\|_
		{L^\infty((-4,0]\times\cD\times\bR_v^3)}
		\right),
	\end{align}
	provided the right-hand side is finite.  Here
	$C=C(\alpha,p,L,K_0,K_\alpha,\cD)$.
\end{thm}

\begin{proof}
	Set
	\[
	\delta_p:=\min\left\{\frac18,\frac{3}{2p}\right\},
	\qquad
	d(Z):=\mathbf d(Z,\bR\times\gamma_{\cD}^0),
	\]
	and denote the right-hand side data norm by
	\[
	\cN_{L,p}:=
	\|\langle v\rangle^{L+A+3}f\|_
	{L^\infty((-4,0]\times\cD\times\bR_v^3)}
	+
	\|\langle v\rangle^{L+3}G\|_
	{L^\infty((-4,0]\times\cD\times\bR_v^3)}.
	\]
	Then \(0<\delta_p<1/2\).
	Theorem~\ref{260520311}, applied with
	\(\delta=\delta_p\) and velocity weight \(L+3\), gives
	\begin{equation}\label{eq:mixed-gradient-pointwise-reduction}
		\langle v\rangle^L|D_vf(Z)|
		\le C\cN_{L,p}\,
		\langle v\rangle^{-3}
		\bigl(d(Z)\wedge1\bigr)^{-1/2-\delta_p}
	\end{equation}
	for a.e. \(Z\in(-1,0]\times\cD\times\bR_v^3\).  It remains to prove
	the purely geometric bound
	\begin{equation}\label{eq:mixed-gradient-geometric-integrability}
		\left\|
		\langle v\rangle^{-3}
		(d(\cdot)\wedge1)^{-1/2-\delta_p}
		\right\|_
		{L^{\infty,p,2}_{t,x,v}((-1,0]\times\cD\times\bR_v^3)}
		\le C_{p,\cD}.
	\end{equation}
	
	Choose \(0<\varepsilon<1\), depending only on \(\cD\), smaller than
	the localization constant in Lemma~\ref{260518305}; on the region
	\[
	d(Z)\ge\varepsilon\langle v\rangle^{-2},
	\]
	we have
	\[
	\langle v\rangle^{-3}
	(d(Z)\wedge1)^{-1/2-\delta_p}
	\le C_{p,\cD}
	\langle v\rangle^{-2+2\delta_p}.
	\]
	The squared velocity integral of the right-hand side is finite because
	\[
	\int_{\bR^3}\langle v\rangle^{-4+4\delta_p}\,\dd v<\infty,
	\qquad
	\delta_p\le\frac18<\frac14.
	\]
	The boundedness of \(\cD\) therefore gives the required mixed-norm
	bound on this region.  For each fixed \((t,x)\), this region and its
	complement partition \(\bR_v^3\).  The two are estimated separately and
	recombined by the triangle inequality, first in \(L^2_v\) and then in
	\(L^p_x\).
	
	It remains to consider
	\[
	d(Z)<\varepsilon\langle v\rangle^{-2}.
	\]
	An almost minimizer in \eqref{kineticdistance} yields
	\[
	\operatorname{dist}(x,\partial\cD)
	\le C\left(d(Z)^3+\langle v\rangle d(Z)^2\right)
	\le C\varepsilon^2.
	\]
	After decreasing \(\varepsilon\), this region is contained in a fixed
	tubular neighborhood covered by finitely many boundary charts.  In one
	such chart, write
	\[
	\widehat x:=\psi(x),
	\qquad
	\widehat v:=A_xv.
	\]
	Here \(d(Z)<1\), and Lemmas~\ref{260518305} and \ref{260522259} give
	\begin{equation*}
		d(Z)\asymp_{\cD}
		|\widehat x_3|^{1/3}+|\widehat v_3|,
		\qquad
		\langle v\rangle\asymp_{\cD}\langle\widehat v\rangle.
	\end{equation*}
	The coordinate Jacobians are uniformly bounded above and below, and,
	for each fixed \(\widehat x\), integration first in
	\(\widehat v'=(\widehat v_1,\widehat v_2)\) and then in
	\(s=\widehat v_3\) gives
	\begin{align*}
		&\int_{\bR^3}\langle\widehat v\rangle^{-6}
		\left(|\widehat x_3|^{1/3}+|\widehat v_3|\right)^{-1-2\delta_p}
		\,\dd\widehat v
		\\
		&\qquad\le
		C\int_{\bR}(1+s^2)^{-2}
		\left(|\widehat x_3|^{1/3}+|s|\right)^{-1-2\delta_p}
		\,\dd s
		\\
		&\qquad\le
		C_p\left(1+|\widehat x_3|^{-2\delta_p/3}\right).
	\end{align*}
	Consequently, after taking the \(L^2_{\widehat v}\) norm and its \(p\)-th
	power, the singular factor in the spatial integral is at most
	\[
	C_p\left(1+|\widehat x_3|^{-\delta_pp/3}\right).
	\]
	This is locally integrable up to \(\widehat x_3=0\), since
	\[
	\frac{\delta_pp}{3}\le\frac12<1.
	\]
	The tangential chart domains are bounded, the number of charts is
	finite, and all bounds are uniform in \(t\).  This proves
	\eqref{eq:mixed-gradient-geometric-integrability}.  Combining
	\eqref{eq:mixed-gradient-geometric-integrability} with
	\eqref{eq:mixed-gradient-pointwise-reduction} proves
	\eqref{eq:global-mixed-gradient-estimate}.
\end{proof}

\section{Nonlinear Landau Theory with Prescribed Inflow}
\label{sec:landau-wellposedness}

We prove well-posedness and stability for the nonlinear Landau inflow problem,
first in the stationary setting and then for the global evolution, and finally
establish convergence to the steady state.

We use the stationary and dynamic perturbations fixed in
\eqref{eq:stationary-perturbation-decomposition-main} and
\eqref{eq:dynamic-perturbation-decomposition-main}, including
$f_{\pt}=f_{\st}+f_{\dy}$.

The incoming trace relations, with $h_{\st}$ and $h_{\pt}$ from Subsection~\ref{subsec:main-inflow}, are
\[
F_{\st}|_{\gamma^-}
=\mu+\sqrt{\mu}\,h_{\st},
\qquad
F|_{\gamma_T^-}
=\mu+\sqrt{\mu}\,h_{\pt}.
\]
Accordingly,
\[
f_{\st}|_{\gamma^-}=h_{\st},
\qquad
f_{\pt}|_{\gamma_T^-}=h_{\pt}.
\]
Regarding \(h_{\st}\) as time independent, we set
\[
h_{\dy}:=h_{\pt}-h_{\st},
\qquad
f_{\dy}|_{\gamma_T^-}=h_{\dy}.
\]

The stationary perturbation \(f_{\st}\) satisfies
\begin{equation*}
	\begin{cases}
		v\cdot\nabla_x f_{\st}+Lf_{\st}
		=N[f_{\st},f_{\st}]
		&\text{in }\cD\times\bR_v^3,\\
		f_{\st}|_{\gamma^-}=h_{\st}
		&\text{on }\gamma^-.
	\end{cases}
\end{equation*}
The total time-dependent perturbation \(f_{\pt}\) satisfies
\begin{equation*}
	\begin{cases}
		(\partial_t+v\cdot\nabla_x)f_{\pt}+Lf_{\pt}
		=N[f_{\pt},f_{\pt}]
		&\text{in }\cD_T\times\bR_v^3,\\
		f_{\pt}(0,\cdot)=f_0
		&\text{in }\cD\times\bR_v^3,\\
		f_{\pt}|_{\gamma_T^-}=h_{\pt}
		&\text{on }\gamma_T^-.
	\end{cases}
\end{equation*}
By \eqref{eq:dynamic-perturbation-decomposition-main} and the bilinearity of \(N\),
the difference \(f_{\dy}\) satisfies
\begin{equation}
	\begin{cases}
		(\partial_t+v\cdot\nabla_x)f_{\dy}+Lf_{\dy}
		=N[f_{\pt},f_{\dy}]
		+N[f_{\dy},f_{\st}]
		&\text{in }\cD_T\times\bR_v^3,\\
		f_{\dy}(0,\cdot)=f_0-f_{\st}
		&\text{in }\cD\times\bR_v^3,\\
		f_{\dy}|_{\gamma_T^-}=h_{\dy}
		&\text{on }\gamma_T^-.
	\end{cases}
	\label{eq:dynamic-perturbation-equation}
\end{equation}

\subsection{Anisotropic extensions of the inflow data}
\label{sec:inflow-anisotropic-extensions}

Throughout this subsection, \(\cD\subset\bR^3\) is a bounded
\(C^{1,1}\) domain satisfying Definition~\ref{domain_assumption}.  The notation \(C_{\ell,\cD}\) means that
the constant depends only on \(\ell\) and the \(C^{1,1}\)-character
\((R_\cD,K_\cD)\) of the domain.
The notation \(C_\cD\) has the analogous meaning without dependence on \(\ell\).

The inflow fibers and stationary and evolutionary anisotropic norms are
those defined in Subsection~\ref{subsec:main-inflow}, in particular
\eqref{eq:stationary-inflow-C11-norm-main} and
\eqref{eq:dynamic-inflow-C11-norm-main}.
On each velocity fiber, finite $C^{1,1}$ regularity includes the
continuous extension of the function and its velocity gradient to the
closed fiber. The empty fiberwise supremum is interpreted as zero. We use the same convention on product sets.  More precisely, if
\(A\subset\bR^m\) and \(g\in C(A\times\bR_v^3)\), then
\[
\begin{aligned}
	\|g\|_{(C_\alpha^{0,1}\cap C_v^{1,1})(A\times\bR_v^3)}
	:={}&
	\sup_{A\times \bR_v^3}|g|
	+
	\sup_{v\in\bR^3}
	\sup_{\substack{\alpha,\beta\in A\\\alpha\neq\beta}}
	\frac{|g(\alpha,v)-g(\beta,v)|}{|\alpha-\beta|}
	\\
	&+
	\sup_{A\times \bR_v^3}|\nabla_vg|
	+
	\sup_{\alpha\in A}
	\sup_{\substack{v,z\in\bR^3\\v\neq z}}
	\frac{|\nabla_vg(\alpha,v)-\nabla_vg(\alpha,z)|}{|v-z|}.
\end{aligned}
\]
For \(A=I\times\cD\), with \(\alpha=(t,x)\), we use the equivalent
product distance \(|t-s|+|x-y|\).  Thus the notation covers
\(A=\partial\cD\), \(A=\cD\), and \(A=I\times\cD\).

We first complete the data from each inflow velocity half-space to all
velocities, preserving one velocity derivative across the grazing
plane.  We then extend the completed data from \(\partial\cD\) into a
fixed collar of \(\cD\).

\begin{lemma}[Velocity completion across the grazing set]
	\label{lem:anisotropic-inflow-velocity-extension}
	There exists a linear operator \(\cV\), independent of the
	weight exponent, such that
	\[
	\cV h=h\qquad\text{on }\gamma^-.
	\]
	For every \(\ell\ge0\),
	\begin{equation}
		\|w^\ell\cV h\|_{(C_x^{0,1}\cap C_v^{1,1})
			(\partial\cD\times\bR_v^3)}
		\le
		C_{\ell,\cD}
		\|w^{\ell+1}h\|_{(C_x^{0,1}\cap C_v^{1,1})(\gamma^-)}.
		\label{eq:velocity-extension-weighted-bound}
	\end{equation}
	The same rule, applied at each time, satisfies
	\begin{equation}
		\|w^\ell\cV h\|_{(C_{t,x}^{0,1}\cap C_v^{1,1})
			((0,T)\times\partial\cD\times\bR_v^3)}
		\le
		C_{\ell,\cD}
		\|w^{\ell+1}h\|_{(C_{t,x}^{0,1}\cap C_v^{1,1})(\gamma_T^-)}
		\label{eq:velocity-extension-evolutionary-bound}
	\end{equation}
	for every \(0<T\le\infty\), with a constant independent of \(T\).
\end{lemma}

\begin{proof}
	On either \(\bR_v^3\) or a half-space \(H\subset\bR_v^3\), the product
	rule gives, for \(m\ge0\),
	\begin{equation}
		\|w^m u\|_{C^{1,1}(H)}
		\simeq_m
		\sum_{|\beta|\le2}
		\|w^mD_v^\beta u\|_{L^\infty(H)},
		\label{eq:inflow-weighted-C11-equivalence}
	\end{equation}
	where second derivatives are taken almost everywhere.  Here we used
	\[
	|D_v^jw^m|\le C_mw^{m-j},
	\qquad j=1,2.
	\]
	
	Set
	\[
	(a_1,a_2):=(2,3),
	\qquad
	(c_1,c_2):=(3,-2),
	\]
	and define
	\[
	Q_j(x,v):=v-a_j(n_x\cdot v)_+n_x,
	\qquad j=1,2.
	\]
	Then
	\[
	Q_j(x,v)\in\overline{\gamma^-_{(x,\ast)}},
	\qquad
	w(v)\le w(Q_j(x,v))\le2w(v).
	\]
	Let
	\begin{equation}
		(\cV h)(x,v)
		:=
		\sum_{j=1}^2c_jh(x,Q_j(x,v)).
		\label{eq:inflow-velocity-completion}
	\end{equation}
	The values on the grazing plane are understood in the trace sense
	specified above.  Since
	\[
	\sum_{j=1}^2c_j=1,
	\qquad
	\sum_{j=1}^2c_ja_j=0,
	\]
	we have \(\cV h=h\) when \(n_x\cdot v\le0\).
	
	Fix \(x\in\partial\cD\).  On \(n_x\cdot v>0\), put
	\[
	A_{j,x}:=I_3-a_jn_x\otimes n_x.
	\]
	Then
	\[
	\nabla_v\cV h(x,v)
	=
	\sum_{j=1}^2c_jA_{j,x}\nabla_vh(x,Q_j(x,v)),
	\]
	and, almost everywhere,
	\[
	D_v^2\cV h(x,v)
	=
	\sum_{j=1}^2c_jA_{j,x}
	D_v^2h(x,Q_j(x,v))A_{j,x}.
	\]
	Moreover,
	\[
	\sum_{j=1}^2c_jA_{j,x}=I_3.
	\]
	Hence both \(\cV h\) and its velocity gradient have matching
	traces across \(n_x\cdot v=0\).  The distributional Hessian therefore
	has no singular part on that plane.  By
	\eqref{eq:inflow-weighted-C11-equivalence} and the weight comparison,
	if
	\[
	M_\ell:=
	\|w^{\ell+1}h\|_{(C_x^{0,1}\cap C_v^{1,1})(\gamma^-)},
	\]
	then
	\begin{equation}
		\sup_{x\in\partial\cD}
		\|w^\ell\cV h(x,\cdot)\|_{C^{1,1}(\bR_v^3)}
		\le C_\ell M_\ell.
		\label{eq:velocity-derivative-weighted-bound}
	\end{equation}
	
	It remains to estimate the dependence on \(x\).  Since
	\(\partial\cD\) is compact and \(C^{1,1}\), there are \(r_0>0\) and
	\(L_\cD>0\) such that, whenever \(x,y\in\partial\cD\) and
	\(|x-y|<r_0\),
	\begin{equation}
		|n_x-n_y|\le L_\cD|x-y|,
		\qquad
		n_x\cdot n_y\ge\frac12.
		\label{eq:inflow-normal-comparison}
	\end{equation}
	We claim that, for each fixed \(K\ge1\), after decreasing
	\(r_0\) by an amount depending only on \(K\) and \(\cD\),
	\begin{equation}
		w^\ell(a)|h(x,a)-h(y,b)|
		\le
		C_{\ell,K,\cD}M_\ell|x-y|
		\label{eq:cross-fiber-comparison}
	\end{equation}
	whenever \(x,y\in\partial\cD\) satisfy \(|x-y|<r_0\) and
	\[
	a\in\overline{\gamma^-_{(x,\ast)}},
	\quad
	b\in\overline{\gamma^-_{(y,\ast)}},
	\quad
	|a-b|\le Kw(a)|x-y|.
	\]
	Indeed, let \(\rho:=|x-y|\), \(\Lambda:=2(L_\cD+1)\), and set
	\[
	\xi:=a-\Lambda w(a)\rho(n_x+n_y).
	\]
	By \eqref{eq:inflow-normal-comparison},
	\[
	n_x\cdot \xi
	\le-\frac32\Lambda w(a)\rho<0,
	\qquad
	n_y\cdot \xi
	\le
	\left(L_\cD-\frac32\Lambda\right)w(a)\rho<0.
	\]
	Thus \(\xi\) belongs to both open inflow fibers.  If
	\[
	\rho\le\frac{1}{2(K+2\Lambda)},
	\]
	then every point \(z\) on either segment \([a,\xi]\) or \([b,\xi]\)
	satisfies
	\[
	\frac12w(a)\le w(z)\le\frac32w(a);
	\]
	indeed \(z\mapsto|z-a|\) is convex.  Hence on these two segments it
	is bounded by \(\max\{|\xi-a|,|b-a|\}\le(K+2\Lambda)w(a)\rho\).
	The segments lie in the appropriate closed inflow half-spaces, and
	\eqref{eq:inflow-weighted-C11-equivalence} gives
	\[
	|\nabla_vh(x,z)|+|\nabla_vh(y,z)|
	\le C_\ell M_\ell w^{-\ell-1}(a)
	\]
	on those segments, while
	\[
	|h(x,\xi)-h(y,\xi)|
	\le M_\ell w^{-\ell-1}(\xi)\rho.
	\]
	Consequently,
	\[
	\begin{aligned}
		|h(x,a)-h(y,b)|
		&\le
		\int_{[a,\xi]}|\nabla_vh(x,z)|\dd s
		+|h(x,\xi)-h(y,\xi)|
		+\int_{[\xi,b]}|\nabla_vh(y,z)|\dd s
		\\
		&\le
		C_{\ell,K,\cD}M_\ell w^{-\ell}(a)\rho,
	\end{aligned}
	\]
	which proves \eqref{eq:cross-fiber-comparison}.
	
	For \(j=1,2\), the Lipschitz continuity of \(n\) and of the positive
	part yields
	\[
	|Q_j(x,v)-Q_j(y,v)|
	\le C_\cD|v|\,|x-y|
	\le C_\cD w(Q_j(x,v))|x-y|.
	\]
	Applying \eqref{eq:cross-fiber-comparison} to the corresponding
	\(Q_j\)'s proves the required \(x\)-Lipschitz bound when \(x\) and
	\(y\) are close.  For the remaining pairs it follows from
	\eqref{eq:velocity-derivative-weighted-bound} and the positive lower
	bound on \(|x-y|\).  This proves
	\eqref{eq:velocity-extension-weighted-bound}.
	
	For time-dependent data, the same construction is applied at each
	time.  Since \(Q_j\) is independent of \(t\),
	\[
	w^\ell(v)
	|\cV h(t,x,v)-\cV h(s,x,v)|
	\le
	C_\ell
	\|w^{\ell+1}h\|_{(C_{t,x}^{0,1}\cap C_v^{1,1})(\gamma_T^-)}
	|t-s|.
	\]
	Combining this with the spatial estimate proves
	\eqref{eq:velocity-extension-evolutionary-bound}; none of the
	constants depends on \(T\).
\end{proof}

\begin{thm}[Weighted anisotropic inflow liftings]
	\label{thm:universal-anisotropic-inflow-extension}
	\begin{enumerate}
		\item[(i)]
		There exists a linear operator \(\cI_{\st}\), independent
		of \(\ell\), such that, for every \(\ell\ge0\) and every \(h\) satisfying
		\[
		w^{\ell+1}h\in(C_x^{0,1}\cap C_v^{1,1})(\gamma^-),
		\]
		we have
		\[
		w^\ell\cI_{\st}h
		\in
		(C_x^{0,1}\cap C_v^{1,1})(\cD\times\bR_v^3),
		\qquad
		\operatorname{Tr}_{\gamma^-}\cI_{\st}h=h,
		\]
		and
		\begin{equation}
			\|w^\ell\cI_{\st}h\|_{(C_x^{0,1}\cap C_v^{1,1})
				(\cD\times\bR_v^3)}
			\le
			C_{\ell,\cD}
			\|w^{\ell+1}h\|_{(C_x^{0,1}\cap C_v^{1,1})(\gamma^-)}.
			\label{eq:stationary-inflow-extension-bound}
		\end{equation}
		
		\item[(ii)]
		There exists a linear lifting rule \(\cI_{\rm ev}\),
		independent of \(\ell\) and \(T\), with the following property.  Let
		\(0<T\le\infty\), and suppose that
		\[
		w^\ell f_0
		\in
		(C_x^{0,1}\cap C_v^{1,1})(\cD\times\bR_v^3),
		\qquad
		w^{\ell+1}h
		\in
		(C_{t,x}^{0,1}\cap C_v^{1,1})(\gamma_T^-).
		\]
		Let \(h_0:=\lim_{t\downarrow0}h(t,\cdot)\), and assume
		\begin{equation}
			h_0=\operatorname{Tr}_{\gamma^-}f_0.
			\label{eq:initial-inflow-compatibility}
		\end{equation}
		Then
		\[
		w^\ell\cI_{\rm ev}(f_0,h)
		\in
		(C_{t,x}^{0,1}\cap C_v^{1,1})
		(\cD_T\times\bR_v^3),
		\]
		\[
		\cI_{\rm ev}(f_0,h)(0)=f_0,
		\qquad
		\operatorname{Tr}_{\gamma_T^-}\cI_{\rm ev}(f_0,h)=h,
		\]
		and
		\begin{equation}
			\begin{aligned}
				&\|w^\ell\cI_{\rm ev}(f_0,h)\|_
				{(C_{t,x}^{0,1}\cap C_v^{1,1})(\cD_T\times\bR_v^3)}
				\le
				C_{\ell,\cD}
				\left(
				\|w^\ell f_0\|_{(C_x^{0,1}\cap C_v^{1,1})(\cD\times\bR_v^3)}
				+
				\|w^{\ell+1}h\|_{(C_{t,x}^{0,1}\cap C_v^{1,1})(\gamma_T^-)}
				\right).
			\end{aligned}
			\label{eq:evolutionary-inflow-extension-bound}
		\end{equation}
	\end{enumerate}
\end{thm}

\begin{proof}
	By the \(C^{1,1}\) regularity of \(\partial\cD\), there is
	\(\delta_\cD>0\) such that the nearest-point projection
	\[
	\pi:\{x\in\cD:d(x)<\delta_\cD\}\longrightarrow\partial\cD,
	\qquad
	d(x):=\operatorname{dist}(x,\partial\cD),
	\]
	is well defined and Lipschitz.  Choose
	\(\zeta\in C^\infty([0,\infty))\) with
	\[
	0\le\zeta\le1,
	\qquad
	\zeta=1\quad\text{on }[0,\delta_\cD/3],
	\qquad
	\zeta=0\quad\text{on }[2\delta_\cD/3,\infty).
	\]
	For \(g:\partial\cD\times\bR_v^3\to\bR\), define
	\[
	(\cC g)(x,v)
	:=
	\begin{cases}
		\zeta(d(x))g(\pi(x),v),&d(x)<\delta_\cD,\\
		0,&d(x)\ge\delta_\cD.
	\end{cases}
	\]
	Then \(\operatorname{Tr}_{\partial\cD}\cC g=g\).  The Lipschitz
	bounds for \(d\), \(\pi\), and \(\zeta\), followed by the product rule
	in \(v\), give
	\begin{equation}
		\|w^\ell\cC g\|_{(C_x^{0,1}\cap C_v^{1,1})
			(\cD\times\bR_v^3)}
		\le
		C_\cD
		\|w^\ell g\|_{(C_x^{0,1}\cap C_v^{1,1})
			(\partial\cD\times\bR_v^3)}.
		\label{eq:spatial-collar-extension-bound}
	\end{equation}
	Applied pointwise in time, the same argument gives
	\begin{equation}
		\|w^\ell\cC g\|_{(C_{t,x}^{0,1}\cap C_v^{1,1})
			((0,T)\times\cD\times\bR_v^3)}
		\le
		C_\cD
		\|w^\ell g\|_{(C_{t,x}^{0,1}\cap C_v^{1,1})
			((0,T)\times\partial\cD\times\bR_v^3)},
		\label{eq:spatial-collar-extension-evolutionary-bound}
	\end{equation}
	with a constant independent of \(T\).
	
	Define
	\[
	\cI_{\st}:=\cC\cV,
	\]
	where \(\cV\) is given by
	Lemma~\ref{lem:anisotropic-inflow-velocity-extension}.  The trace
	identities and estimates
	\eqref{eq:velocity-extension-weighted-bound} and
	\eqref{eq:spatial-collar-extension-bound} prove
	\eqref{eq:stationary-inflow-extension-bound} and part~(i).
	
	Now we prove part~(ii).
	The time-Lipschitz estimate gives a uniform limit
	\(w^{\ell+1}h(t)\to w^{\ell+1}h_0\) as \(t\downarrow0\).  Difference
	quotients in \(v\), followed by weak-* compactness in
	\(W^{2,\infty}\) on each velocity fiber, give
	\begin{equation}
		\|w^{\ell+1}h_0\|_{(C_x^{0,1}\cap C_v^{1,1})(\gamma^-)}
		\le \|w^{\ell+1}h\|_{(C_{t,x}^{0,1}\cap C_v^{1,1})(\gamma_T^-)}.
		\label{eq:initial-trace-bound}
	\end{equation}
	Extend \(h\) to \(t=0\) by \(h(0):=h_0\), and put
	\[
	G(t):=\cC\cV h(t),
	\qquad 0\le t<T.
	\]
	By \eqref{eq:velocity-extension-evolutionary-bound},
	\eqref{eq:spatial-collar-extension-evolutionary-bound}, and
	\eqref{eq:initial-trace-bound},
	\begin{equation}
		\|w^\ell G\|_{(C_{t,x}^{0,1}\cap C_v^{1,1})
			([0,T)\times\cD\times\bR_v^3)}
		\le C_{\ell,\cD}\|w^{\ell+1}h\|_{(C_{t,x}^{0,1}\cap C_v^{1,1})(\gamma_T^-)},
		\qquad
		\operatorname{Tr}_{\gamma_T^-}G=h.
		\label{eq:evolutionary-boundary-extension-bound}
	\end{equation}
	
	Let
	\[
	B:=f_0-G(0).
	\]
	The compatibility condition \eqref{eq:initial-inflow-compatibility}
	implies \(\operatorname{Tr}_{\gamma^-}B=0\).  Define
	\begin{equation}
		\cI_{\rm ev}(f_0,h)(t,x,v)
		:=
		G(t,x,v)+\ee^{-t}B(x,v).
		\label{eq:evolutionary-inflow-lifting-formula}
	\end{equation}
	Then
	\[
	\cI_{\rm ev}(f_0,h)(0)=f_0,
	\qquad
	\operatorname{Tr}_{\gamma_T^-}\cI_{\rm ev}(f_0,h)=h.
	\]
	Moreover,
	\[
	\|w^\ell \ee^{-t}B\|_{(C_{t,x}^{0,1}\cap C_v^{1,1})
		([0,T)\times\cD\times\bR_v^3)}
	\le
	2\|w^\ell B\|_{(C_x^{0,1}\cap C_v^{1,1})
		(\cD\times\bR_v^3)},
	\]
	since \(|\ee^{-t}-\ee^{-s}|\le|t-s|\).  In addition,
	\[
	\begin{aligned}
		\|w^\ell B\|_{(C_x^{0,1}\cap C_v^{1,1})
			(\cD\times\bR_v^3)}
		&\le
		\|w^\ell f_0\|_{(C_x^{0,1}\cap C_v^{1,1})
			(\cD\times\bR_v^3)}
		+\|w^\ell G(0)\|_{(C_x^{0,1}\cap C_v^{1,1})
			(\cD\times\bR_v^3)}
		\\
		&\le
		\|w^\ell f_0\|_{(C_x^{0,1}\cap C_v^{1,1})
			(\cD\times\bR_v^3)}
		+C_{\ell,\cD}\|w^{\ell+1}h\|_{(C_{t,x}^{0,1}\cap C_v^{1,1})(\gamma_T^-)},
	\end{aligned}
	\]
	Together with \eqref{eq:evolutionary-boundary-extension-bound}, this proves
	\eqref{eq:evolutionary-inflow-extension-bound}.  The formulas
	\eqref{eq:inflow-velocity-completion} and
	\eqref{eq:evolutionary-inflow-lifting-formula} also show that the
	operators are linear and independent of \(\ell\) and \(T\).
\end{proof}

\begin{lemma}
	\label{lem:dynamic-weighted-Lipschitz-to-kinetic-holder}
	Let \(0<\alpha<1\), \(m\ge0\), and let
	\(\cQ=I\times E\times\bR_v^3\) with \(I\subset\bR\) an interval and
	\(E\subset\bR_x^3\).  Let \(H\) be defined on \(\cQ\).  If
	\begin{equation}\label{eq:dynamic-weighted-Lipschitz-size}
		\mathfrak L_{m+1}(H)
		:=\|w^{m+1}H\|_{(C_{t,x}^{0,1}\cap C_v^{1,1})(\cQ)}<\infty,
	\end{equation}
	then, for all \(Z=(t,x,v)\) and \(Z'=(t',x',v')\) in \(\cQ\),
	\begin{equation}\label{eq:dynamic-weighted-Lipschitz-kinetic-holder}
		|w(v)^mH(Z)-w(v')^mH(Z')|
		\le C_{m,\alpha}\mathfrak L_{m+1}(H)\,\mathbf d(Z,Z')^\alpha.
	\end{equation}
\end{lemma}

\begin{proof}
	Put \(d:=\mathbf d(Z,Z')\).  For \(0<d\le1\), choose an intermediate
	velocity \(\bar v\) in \eqref{kineticdistance} such that all four
	defining quantities are at most \(2d\).  Since
	\(|\bar v|\le|v|+2d\le Cw(v)\), this gives
	\begin{equation}\label{eq:dynamic-weighted-Lipschitz-metric-components}
		|t-t'|\le4d^2,
		\qquad
		|v-v'|\le4d,
		\qquad
		|x-x'|\le8d^3+4d^2|\bar v|\le Cw(v)d^2.
	\end{equation}
	Set \(K:=w^{m+1}H\), so that \(w^mH=w^{-1}K\).  The coordinate
	Lipschitz bounds in \eqref{eq:dynamic-weighted-Lipschitz-size}, the
	estimates \eqref{eq:dynamic-weighted-Lipschitz-metric-components}, and
	\(|w(v)^{-1}-w(v')^{-1}|\le|v-v'|\) give
	\[
	\begin{aligned}
		|w(v)^mH(Z)-w(v')^mH(Z')|
		&\le
		w(v)^{-1}|K(Z)-K(Z')|
		+|w(v)^{-1}-w(v')^{-1}|\,|K(Z')|
		\\
		&\le
		C\mathfrak L_{m+1}(H)\left(d+d^2\right)
		\le
		C\mathfrak L_{m+1}(H)d^\alpha.
	\end{aligned}
	\]
	For \(d\ge1\), the left-hand side is at most
	\(2\|w^mH\|_{L^\infty(\cQ)}\le2\mathfrak L_{m+1}(H)d^\alpha\).
\end{proof}

\subsection{Existence and uniqueness of the steady state}\label{subsec:inflow-stationary-proof}

In this subsection, we prove Theorem~\ref{Landau_st_exi} in the normalized variables.
Writing $F_\st=\mu+\sqrt\mu\,f_\st$ and $H_\st=\mu+\sqrt\mu\,h_\st$, the stationary problem reads
\begin{equation}
	\begin{cases}
		v\cdot\nabla_xf_{\st}+Lf_{\st}=N[f_{\st},f_{\st}]
		&\text{in }\cD\times\bR_v^3,\\
		f_{\st}|_{\gamma^-}=h_{\st}
		&\text{on }\gamma^-.
	\end{cases}
	\label{2605041257}
\end{equation}
The proof combines stationary linear H\"older and
velocity-gradient estimates with an iteration for the nonlinear problem.
The stationary kinetic metric and the extension lemma below place coefficients
defined on \(\cD\times\bR_v^3\) in the global H\"older framework required by
the regularity estimates.

Recall the definitions of $\mathbf d_{\st}$ and $C^\alpha_{\kin}(\cD\times\bR_v^3)$: for \(Z=(x,v)\) and \(Z'=(y,v')\),
\begin{align}
	\mathbf d_{\st}(Z,Z')
	:=\,&
	\inf_{\tau\in\bR}
	\mathbf d\bigl((\tau,x,v),(0,y,v')\bigr)
	\label{2607271135}\\
	=\,&\inf_{(\tau,\widehat v)\in\bR\times\bR_v^3}
	\max\left\{
	|\tau|^{1/2},
	|x-y-\tau\widehat v|^{1/3},
	|v-\widehat v|,
	|v'-\widehat v|
	\right\},\nonumber
\end{align}
which is a metric on $\bR^3_x\times \bR^3_v$, and for \(\alpha\in(0,1]\),
\[
\|u\|_{C^\alpha_{\kin}(\cD\times\bR_v^3)}
:=
\|u\|_{L^\infty(\cD\times\bR_v^3)}
+
\sup_{\substack{Z,Z'\in\cD\times\bR_v^3\\Z\neq Z'}}
\frac{|u(Z)-u(Z')|}{\mathbf d_{\st}(Z,Z')^\alpha}.
\]

Multiplication by \(\mu^{1/4}\) is bounded on both kinetic H\"older
spaces.  Indeed, we split the increment of \(\mu^{1/4}u\) as
\(\mu^{1/4}(v)\bigl(u(Z)-u(Z')\bigr)
+\bigl(\mu^{1/4}(v)-\mu^{1/4}(v')\bigr)u(Z')\).  Since
\(|v-v'|\le2\mathbf d_{\st}\) and \(|v-v'|\le2\mathbf d\), the
boundedness and the Lipschitz continuity of \(\mu^{1/4}\) bound the
second term by \(C\|u\|_{L^\infty}\) times
\(\min\{1,\mathbf d_{\st}\}\le\mathbf d_{\st}^\alpha\),
respectively \(\min\{1,\mathbf d\}\le\mathbf d^\alpha\).  Hence
there is a universal \(C_*\ge1\) such that, for every
\(\alpha\in(0,1]\),
\begin{equation}
	\|\mu^{1/4}u\|_{C^\alpha_{\kin}}
	\le
	C_*\|u\|_{C^\alpha_{\kin}},
	\label{eq:maxwellian-holder-multiplier}
\end{equation}
on \(\cD\times\bR_v^3\) and on \(\cD_T\times\bR_v^3\) alike.

The data in this subsection are prescribed only on \(\gamma^-\).
The weighted \(C^{1,1}\)-equivalence
\eqref{eq:inflow-weighted-C11-equivalence} implies that, whenever
\(0\le m\le m'\),
\begin{equation}
	\|w^m h\|_{(C_x^{0,1}\cap C_v^{1,1})(\gamma^-)}
	\le C_{m,m'}\|w^{m'} h\|_{(C_x^{0,1}\cap C_v^{1,1})(\gamma^-)}.
	\label{eq:stationary-boundary-weight-monotonicity}
\end{equation}
Throughout this subsection, constants depend only on the parameters
displayed in their subscripts and on the \(C^{1,1}\)-character of
\(\cD\).  In particular, they are independent of \(g\), \(h\), and
the solution.
\begin{lemma}[Norm-preserving H\"older extension]
	\label{lem:norm-preserving-holder-extension}
	Let \((\bS,d_\bS)\) be a metric space, let
	\(\varnothing\neq A\subset\bS\), and let \(\alpha\in(0,1]\).  Every
	bounded \(g\in C^{0,\alpha}(A)\) has an extension
	\(\mathfrak g\in C^{0,\alpha}(\bS)\) satisfying
	\begin{equation}
		\mathfrak g|_A=g,
		\qquad
		\|\mathfrak g\|_{L^\infty(\bS)}=\|g\|_{L^\infty(A)},
		\qquad
		[\mathfrak g]_{C^{0,\alpha}(\bS)}=[g]_{C^{0,\alpha}(A)}.
		\label{eq:norm-preserving-holder-extension}
	\end{equation}
\end{lemma}

\begin{proof}
	Set
	\[
	L:=[g]_{C^{0,\alpha}(A)},
	\qquad
	M:=\|g\|_{L^\infty(A)}.
	\]
	Since \(\alpha\in(0,1]\), the function
	\[
	d_{\bS,\alpha}(X,Y):=d_\bS(X,Y)^\alpha
	\]
	is a metric on \(\bS\).  Thus \(g\) is \(L\)-Lipschitz on \(A\)
	with respect to \(d_{\bS,\alpha}\).  The McShane theorem provides an
	extension \(\bar g:\bS\to\bR\) such that
	\[
	\bar g|_A=g,
	\qquad
	|\bar g(X)-\bar g(Y)|
	\le Ld_\bS(X,Y)^\alpha,
	\qquad X,Y\in\bS.
	\]
	The McShane extension need not preserve the supremum norm.  We
	therefore define
	\[
	\eta(r):=\max\{-M,\min\{r,M\}\},
	\qquad
	\mathfrak g:=\eta\circ\bar g.
	\]
	The map \(\eta\) is one-Lipschitz and fixes every value of \(g\).
	Consequently,
	\[
	\mathfrak g|_A=g,
	\qquad
	\|\mathfrak g\|_{L^\infty(\bS)}\le M,
	\qquad
	[\mathfrak g]_{C^{0,\alpha}(\bS)}\le L.
	\]
	Restriction to \(A\) gives the reverse inequalities.  This proves
	\eqref{eq:norm-preserving-holder-extension}.
\end{proof}

\begin{lemma}[Stationary linearized Landau problem]
	\label{260527559}
	Let \(\varepsilon_0=\varepsilon_0(\cD)>0\) be the smallness constant
	in Theorem~\ref{260322442}, and let
	\(\alpha_0=\alpha_0(\cD)\in(0,1)\) be the exponent in
	Theorem~\ref{260526316}.  Set
	\begin{equation}
		B^{\st}_0:=37 .
		\label{eq:stationary-B0}
	\end{equation}
	
	Assume that
	\[
	\|\mu^{1/4}g\|_{L^\infty(\cD\times\bR_v^3)}
	\le\varepsilon_0.
	\]
	For every \(\ell\ge0\), every inflow datum \(h_{\st}\) and every forcing term $S$ satisfying
	\[
	\big\|w^{\ell+B^{\st}_0} h_{\st}\big\|_{(C_x^{0,1}\cap C_v^{1,1})(\gamma^-)}+\|w^{\ell+B^{\st}_0}S\|_{L^\infty(\cD\times \bR^3_v)}<\infty,
	\]
	the problem
	\begin{equation}
		\begin{cases}
			v\cdot\nabla_x f+Lf=N[g,f]+S
			&\text{in }\cD\times\bR_v^3,\\
			f|_{\gamma^-}=h_{\st}
			&\text{on }\gamma^-
		\end{cases}
		\label{260526322}
	\end{equation}
	has a unique solution
	\(f\in\cH_{\kin,\sigma}(\cD\times\bR_v^3)\), and
	\begin{equation}
		\begin{aligned}
			&\|w^\ell f\|_{\cH_{\kin,\sigma}}
			+\|w^\ell f\|_{L^\infty}
			+\|w^\ell f\|_{C^{\alpha_0}_{\kin}}\\
			\le\,&
			C_{\ell,\cD}\left(\big\|w^{\ell+B^{\st}_0} h_{\st}\big\|_{(C_x^{0,1}\cap C_v^{1,1})(\gamma^-)}+\big\|w^{\ell+B^{\st}_0} S\big\|_{L^\infty(\cD\times \bR^3_v)}\right).
		\end{aligned}
		\label{260527600}
	\end{equation}
	All unmarked norms in this subsection are taken over
	\(\cD\times\bR_v^3\).
	
	Suppose in addition that, for some \(\alpha\in(0,1)\) and \(K>0\),
	\begin{equation}
		\|\mu^{1/4} g\|_{C^\alpha_{\kin}(\cD\times\bR_v^3)}\le K.
		\label{eq:stationary-coefficient-holder-assumption}
	\end{equation}
	Set
	\begin{equation}
		B^{\st}_1(\alpha)
		:=
		\max\left\{
		B^{\st}_0,\,A(2\alpha/3)+28
		\right\},
		\label{eq:stationary-B1}
	\end{equation}
	where \(A(\cdot)\) is defined in
	\eqref{eq:global-gradient-weight-loss}.  Then, for every
	\(p\in(1,\infty)\),
	\begin{equation}
		\|w^\ell D_vf\|_{L^{p,2}_{x,v}}
		\le
		C_{\ell,\cD,\alpha,p,K}
		\left(\big\|w^{\ell+B^{\st}_1(\alpha)} h_{\st}\big\|_{(C_x^{0,1}\cap C_v^{1,1})(\gamma^-)}+\big\|w^{\ell+B^{\st}_1(\alpha)} S\big\|_{L^\infty(\cD\times \bR^3_v)}\right),
		\label{eq:stationary-mixed-gradient-estimate}
	\end{equation}
	provided the norm on the right is finite.
\end{lemma}

\begin{proof}
	We denote
	$$
	{\rm D}_{\ell}:=\big\|w^{\ell} h_{\st}\big\|_{(C_x^{0,1}\cap C_v^{1,1})(\gamma^-)}+\big\|w^{\ell} S\big\|_{L^\infty(\cD\times \bR^3_v)}.
	$$

	\medskip
	\noindent\emph{Step 1. Energy estimate.}
	For every \(m\ge0\), the boundary norm satisfies
	\begin{equation}
		\|w^m h_{\st}\|_{L^2(\gamma^-)}
		\le C_{\cD}\big\|w^{m+3} h_{\st}\big\|_{(C_x^{0,1}\cap C_v^{1,1})(\gamma^-)},
		\label{eq:stationary-inflow-trace-bound}
	\end{equation}
	since
	\(\int_{\bR^3}w^{-6}|n_x\cdot v|\dd v<\infty\), uniformly in
	\(x\in\partial\cD\).  
	In addition, by \eqref{eq:Hsigma-dual-weighted},
	\begin{equation}\label{260915411}
		\|w^mS\|_{L^2(\cD;H^{-1}_\sigma)}\leq C_\cD \|w^{m+1/2}S\|_{L^2(\cD\times \bR^3_v)}\leq C_\cD \|w^{m+3}S\|_{L^\infty(\cD\times \bR^3_v)}.
	\end{equation}
	By
	\eqref{eq:stationary-inflow-trace-bound} and \eqref{260915411},
	Theorems~\ref{260322442} and \ref{thm:weighted-steady} yield the
	unique solution of
	\eqref{260526322} and
	\begin{align}
		\|w^\ell f\|_{\cH_{\kin,\sigma}}+\|w^\ell f\|_{L^2}
		&\le C_{\ell,\cD}{\rm D}_{\ell+3}.
		\label{260727959}
	\end{align}
	
	\medskip
	\noindent\emph{Step 2. Weighted \(L^\infty\) estimate.}
	Set
	\[
	H:=\cI_{\st}h_{\st},
	\]
	where \(\cI_{\st}\) is given by
	Theorem~\ref{thm:universal-anisotropic-inflow-extension}.  Its trace
	identity and estimate give, for every \(q\ge0\),
	\begin{equation}
		\operatorname{Tr}_{\gamma^-}H=h_{\st},
		\qquad
		\|w^qH\|_{(C_x^{0,1}\cap C_v^{1,1})(\cD\times\bR_v^3)}
		\le C_{q,\cD}\big\|w^{q+1} h_{\st}\big\|_{(C_x^{0,1}\cap C_v^{1,1})(\gamma^-)}.
		\label{eq:stationary-lift-anisotropic-bound}
	\end{equation}
	Since \(w\) is independent of \(x\),
	\eqref{eq:inflow-weighted-C11-equivalence} turns this into
	\begin{equation}
		\sum_{|\beta|\le2}
		\|w^qD_v^\beta H\|_{L^\infty}
		+\|w^q\nabla_xH\|_{L^\infty}
		\le C_{q,\cD}\big\|w^{q+1} h_{\st}\big\|_{(C_x^{0,1}\cap C_v^{1,1})(\gamma^-)}.
		\label{eq:stationary-lift-derivative-bound}
	\end{equation}
	
	Let \(\bar f:=f-H\).  By
	\eqref{eq:stationary-lift-anisotropic-bound},
	\(\bar f|_{\gamma^-}=0\), and
	\begin{equation}
		v\cdot\nabla_x\bar f=\cA_g\bar f+\bar S,
		\qquad
		\bar S:=S+\cK f-\bigl(v\cdot\nabla_xH-\cA_gH\bigr).
		\label{2607271033}
	\end{equation}
	The identity \(\cA_g=D_i(a_g^{ij}D_j)+b_g^iD_i+c_g\) puts
	\eqref{2607271033} in the form required by the regularity estimates:
	\begin{equation}
		v\cdot\nabla_x\bar f
		=D_i(a_g^{ij}D_j\bar f)+b_g^iD_i\bar f+G,
		\qquad
		G:=c_g\bar f+\bar S.
		\label{eq:stationary-bar-f-regularity-equation}
	\end{equation}
	From \eqref{260727959} and
	\eqref{eq:stationary-lift-derivative-bound},
	\begin{equation}
		\|w^m\bar f\|_{L^2}
		\le C_{m,\cD}{\rm D}_{m+4}.
		\label{eq:bar-f-weighted-L2}
	\end{equation}
	In \(\bar S\) the compact part acts on \(f=\bar f+H\).  We split
	\(\cK f=\cK\bar f+\cK H\) and absorb
	\(\cK H\) into the term carrying the data, using
	\(\|\mu^{1/8}H(x,\cdot)\|_{L^2_v}\le C\|w^{q}H\|_{L^\infty}\) for any
	\(q\ge0\).  The coefficient bounds in
	Lemma~\ref{lem:perturbed-landau-ellipticity}, the pointwise estimate
	\eqref{eq:compact-part-pointwise} for \(\cK\), and
	\eqref{eq:stationary-lift-derivative-bound} then imply
	\begin{align}
		\begin{aligned}
			w^m|\bar S(x,v)|
			\le
			C_{m,\cD}\Big[
			\,&\mu^{1/8}(v)\|\mu^{1/8}\bar f(x,\cdot)\|_{L^2_v}+{\rm D}_{m+22}w(v)^{-20}
			\Big].
			\label{2607271112}
		\end{aligned}
	\end{align}
	
	Extend \(g\) by zero in \(x\), and regard \(\bar f\) as a
	time-independent function on \((-4,0)\times\cD\times\bR_v^3\).
	For a fixed \(\ell\ge0\), set \(u:=w^\ell\bar f\).  Conjugating
	\eqref{2607271033} by \(w^\ell\) gives
	\[
	(\partial_t+v\cdot\nabla_x)u
	=\cA_{g,\ell}u+w^\ell \bar S,
	\qquad u|_{\gamma^-}=0.
	\]
	Estimate \eqref{2607271112} with \(m=\ell\) has precisely the source
	form required in Theorem~\ref{260504554} (with, for instance,
	\(\varepsilon=1/8\)).  Moreover,
	\(\|\mu^{1/8}\bar f\|_{L_v^2}
	\le\|\mu^{1/8}u\|_{L_v^2}\).
	Applying that theorem and then
	\eqref{eq:bar-f-weighted-L2} with \(m=\ell+21\) gives
	\begin{equation}
		\|w^\ell\bar f\|_{L^\infty}
		\le
		C_{\ell,\cD}{\rm D}_{\ell+25}.
		\label{2607271120}
	\end{equation}
	Using this once more in \eqref{2607271112},
	\begin{equation}
		\|w^\ell \bar S\|_{L^\infty}
		\le
		C_{\ell,\cD}{\rm D}_{\ell+25}.
		\label{2607271121}
	\end{equation}
	The coefficient bound \(|c_g(v)|\le Cw(v)^{-1}\), together with
	\eqref{2607271120} and \eqref{2607271121}, therefore gives
	\begin{equation}
		\|w^\ell G\|_{L^\infty}
		\le
		C_{\ell,\cD}{\rm D}_{\ell+25}.
		\label{eq:stationary-full-source-bound}
	\end{equation}
	Together with \eqref{eq:stationary-lift-derivative-bound}, this proves
	\begin{equation*}
		\|w^\ell f\|_{L^\infty}
		\le C_{\ell,\cD}{\rm D}_{\ell+25}.
	\end{equation*}
	
	\medskip
	\noindent\emph{Step 3. Stationary kinetic H\"older estimate.}
	The weight loss in Theorem~\ref{260526316} is \(12\).  If
	\(\rho=\mathbf d_{\st}((x,v),(y,v'))<1/4\), the definition of
	\(\mathbf d_{\st}\) supplies \(\tau\in\bR\) with
	\[
	\mathbf d\bigl((\tau,x,v),(0,y,v')\bigr)<2\rho,
	\]
	and then \(|\tau|^{1/2}<2\rho\) forces \(|\tau|<1\).  Since
	\(\mathbf d\) is symmetric and invariant under time translations,
	we may interchange the two points and assume \(\tau\in(-1,0]\).
	Since \(\bar f\) is time independent, Theorem~\ref{260526316},
	applied to \eqref{eq:stationary-bar-f-regularity-equation}, and
	\eqref{eq:stationary-full-source-bound} give
	\begin{align*}
		&\bigl(w^\ell(v)+w^\ell(v')\bigr)
		|\bar f(x,v)-\bar f(y,v')|\le
		C_{\ell,\cD}{\rm D}_{\ell+37}
		\rho^{\alpha_0}.
	\end{align*}
	Moreover, the choice of \(\tau\) above implies
	\(|v-v'|<4\rho\) and \(w(v')\simeq w(v)\).  Hence
	\[
	|w^\ell(v)-w^\ell(v')|\,|\bar f(y,v')|
	\le
	C_\ell\rho\|w^\ell\bar f\|_{L^\infty}
	\le
	C_{\ell,\cD}
	{\rm D}_{\ell+25}
	\rho^{\alpha_0},
	\]
	where we used \(\rho\le\rho^{\alpha_0}\).  Decomposing
	\[
	(w^\ell\bar f)(x,v)-(w^\ell\bar f)(y,v')
	=w^\ell(v)\bigl(\bar f(x,v)-\bar f(y,v')\bigr)
	+\bigl(w^\ell(v)-w^\ell(v')\bigr)\bar f(y,v')
	\]
	proves the required estimate when \(\rho<1/4\).
	For \(\rho\ge1/4\), the same conclusion follows from
	\eqref{2607271120}.  Hence
	\begin{equation}
		\|w^\ell\bar f\|_{C^{\alpha_0}_{\kin}}
		\le
		C_{\ell,\cD}{\rm D}_{\ell+37}.
		\label{eq:stationary-bar-f-holder}
	\end{equation}
	
	The definition of \(\mathbf d_{\st}\) also gives, when
	\(0<\rho<1\),
	\begin{equation}
		|x-y|+|v-v'|
		\le Cw(v)\rho,
		\qquad
		|v-v'|\le4\rho,
		\qquad
		w\simeq w(v)\ \text{ on }[v,v'].
		\label{eq:stationary-distance-euclidean-control}
	\end{equation}
	Indeed, choose \((\tau,\widehat v)\) within a factor two of the
	infimum in \eqref{2607271135} and use
	\(|\tau|\lesssim\rho^2\),
	\(|x-y-\tau\widehat v|\lesssim\rho^3\), and
	\(|v-\widehat v|+|v'-\widehat v|\lesssim\rho\).
	Equations \eqref{eq:stationary-lift-anisotropic-bound} and
	\eqref{eq:stationary-distance-euclidean-control} show that
	\begin{equation}
		\|w^\ell H\|_{C^{\alpha_0}_{\kin}}
		\le C_{\ell,\cD}\big\|w^{\ell+2} h_{\st}\big\|_{(C_x^{0,1}\cap C_v^{1,1})(\gamma^-)}.
		\label{eq:stationary-extension-holder}
	\end{equation}
	Indeed, writing \(w^\ell H=w^{-1}(w^{\ell+1}H)\) and using
	\(|D_vw^{-1}|\le w^{-1}\), the increment of \(w^\ell H\) between two
	points at \(\mathbf d_{\st}\)-distance \(\rho<1\) is at most
	\[
	C_\ell\sup_{[v,v']}w^{-1}\,
	\|w^{\ell+1}H\|_{(C_x^{0,1}\cap C_v^{1,1})}
	\bigl(|x-y|+|v-v'|\bigr)
	\le
	C_\ell\rho\,\|w^{\ell+1}H\|_{(C_x^{0,1}\cap C_v^{1,1})},
	\]
	by \eqref{eq:stationary-distance-euclidean-control}.  Now use
	\(\rho\le\rho^{\alpha_0}\), and the supremum norm when
	\(\rho\ge1\).
	Combining \eqref{260727959},
	\eqref{eq:stationary-bar-f-holder}, and
	\eqref{eq:stationary-extension-holder} proves \eqref{260527600}.
	
	\medskip
	\noindent\emph{Step 4. Mixed-norm estimate for \(D_vf\).}
	Assume \eqref{eq:stationary-coefficient-holder-assumption}.  Apply
	Lemma~\ref{lem:norm-preserving-holder-extension}, with the metric
	\(\mathbf d_{\st}\), to \(\mu^{1/4}g\), and define the extension of
	\(g\) to \(\bR_x^3\times\bR_v^3\) by
	\(g:=\mu^{-1/4}(\mu^{1/4}g)\) outside \(\cD\).  Only the values of
	\(g\) on \(\cD\times\bR_v^3\) enter \eqref{260526322}.  Therefore
	this extension changes neither the equation nor the solution.  It
	replaces, from here on, the extension by zero used in Step~2.  Since
	Lemma~\ref{lem:norm-preserving-holder-extension} preserves both the
	supremum norm and the H\"older seminorm, the extension satisfies
	\begin{equation}
		\|\mu^{1/4}g\|_{L^\infty(\bR^6)}\le\varepsilon_0,
		\qquad
		[\mu^{1/4}g]_{C^{\alpha}_{\kin}(\bR^6)}
		\le K.
		\label{eq:stationary-g-holder}
	\end{equation}
	Its time-independent extension has the same bound for the full kinetic
	metric since \(\mathbf d_{\st}\le\mathbf d\).  By
	\eqref{eq:stationary-g-holder},
	Lemma~\ref{260527138} gives the coefficient H\"older bounds required
	in Theorem~\ref{260527139}.
	
	Let
	\[
	A':=A(2\alpha/3),
	\]
	where \(A(\cdot)\) is defined in
	\eqref{eq:global-gradient-weight-loss}.  Applying
	Theorem~\ref{260527139} to the
	time-independent equation
	\eqref{eq:stationary-bar-f-regularity-equation}, and then using
	\eqref{2607271120} and \eqref{eq:stationary-full-source-bound}, yields
	\begin{equation}
		\|w^\ell D_v\bar f\|_{L^{p,2}_{x,v}}
		\le
		C_{\ell,\cD,\alpha,p,K}
		{\rm D}_{\ell+A'+28}.
		\label{eq:stationary-bar-f-gradient}
	\end{equation}
	On the other hand, \eqref{eq:stationary-lift-derivative-bound} and
	H\"older's inequality in \(v\) give
	\begin{equation}
		\|w^\ell D_vH\|_{L^{p,2}_{x,v}}
		\le
		C_{\ell,\cD,\alpha,p}
		\big\|w^{\ell+A'+28} h_{\st}\big\|_{(C_x^{0,1}\cap C_v^{1,1})(\gamma^-)}.
		\label{eq:stationary-extension-gradient}
	\end{equation}
	Combining \eqref{eq:stationary-bar-f-gradient} and
	\eqref{eq:stationary-extension-gradient} proves
	\eqref{eq:stationary-mixed-gradient-estimate} by the definition of
	\(B^{\st}_1(\alpha)\).
\end{proof}

\begin{proof}[\textbf{Proof of Theorem~\ref{Landau_st_exi}}]
	Let $\varepsilon_0$ be the constant of Theorem~\ref{260322442} and $\alpha_0$ the exponent of Theorem~\ref{260526316}, and fix
	\begin{equation}
		A_{\st}:=\max\{B^{\st}_0+5,B^{\st}_1(\alpha_0)+2\},
		\label{eq:nonlinear-stationary-A0-choice}
	\end{equation}
	where \(B^{\st}_0\) and \(B^{\st}_1(\alpha_0)\) are the weights \eqref{eq:stationary-B0} and
	\eqref{eq:stationary-B1} of Lemma~\ref{260527559}, the latter at \(\alpha=\alpha_0\).
	
	\medskip
	\noindent\emph{Step 1. Iteration and uniform estimates.}
	Set \(f^{(0)}:=0\).  Given \(f^{(k)}\), define \(f^{(k+1)}\) by
	\begin{equation}
		\begin{cases}
			v\cdot\nabla_xf^{(k+1)}+Lf^{(k+1)}
			=N[f^{(k)},f^{(k+1)}]
			&\text{in }\cD\times\bR_v^3,\\
			f^{(k+1)}|_{\gamma^-}=h_{\st}
			&\text{on }\gamma^-.
		\end{cases}
		\label{260527636}
	\end{equation}
	Assume inductively that
	\[
	\|\mu^{1/4}f^{(k)}\|_{L^\infty}\le\varepsilon_0,
	\qquad
	\|f^{(k)}\|_{C^{\alpha_0}_{\kin}}\le1.
	\]
	By \eqref{eq:maxwellian-holder-multiplier}, the second induction
	hypothesis gives
	\eqref{eq:stationary-coefficient-holder-assumption} for \(g=f^{(k)}\)
	with \(\alpha=\alpha_0\) and \(K=C_*\).  Therefore the constants
	below do not depend on \(k\).  Lemma~\ref{260527559}, used with
	\(\ell=0,5\), and then with \((\ell,p)=(2,30)\), gives
	\begin{equation}
		\begin{aligned}
			&\|f^{(k+1)}\|_{\cH_{\kin,\sigma}}
			+\|f^{(k+1)}\|_{L^\infty}
			+\|f^{(k+1)}\|_{C^{\alpha_0}_{\kin}}
			\\
			&\qquad
			+\|w^5f^{(k+1)}\|_{L^\infty}
			+\|w^2D_vf^{(k+1)}\|_{L^{30,2}_{x,v}}
			\le C_\cD
			\|w^{A_{\st}}h_{\st}\|_
			{(C_x^{0,1}\cap C_v^{1,1})(\gamma^-)}.
		\end{aligned}
		\label{eq:nonlinear-iteration-uniform-bound}
	\end{equation}
	Here \eqref{eq:stationary-boundary-weight-monotonicity} and the choice
	\eqref{eq:nonlinear-stationary-A0-choice} were used.  Choose
	\(\delta_{\st}\le1\) so that
	\begin{equation*}
		C_\cD\delta_{\st}\le\min\{\varepsilon_0,1\}.
	\end{equation*}
	Since the induction hypotheses hold for \(f^{(0)}=0\), the sequence is
	well defined and \eqref{eq:nonlinear-iteration-uniform-bound} holds
	uniformly in \(k\).
	
	\medskip
	\noindent\emph{Step 2. Contraction.}
	For \(k\ge0\), put
	\(d^{(k)}:=f^{(k+1)}-f^{(k)}\).  Subtracting two consecutive
	instances of \eqref{260527636} and using the bilinearity of \(N\)
	gives, for \(k\ge0\),
	\begin{equation}
		\begin{cases}
			v\cdot\nabla_xd^{(k+1)}+Ld^{(k+1)}
			=N[f^{(k+1)},d^{(k+1)}]
			+N[d^{(k)},f^{(k+1)}],\\
			d^{(k+1)}|_{\gamma^-}=0.
		\end{cases}
		\label{eq:nonlinear-iteration-difference}
	\end{equation}
	Apply Theorem~\ref{260322442}, with coefficient \(f^{(k+1)}\), to
	\eqref{eq:nonlinear-iteration-difference}, and estimate its source by
	the stationary nonlinear estimate \eqref{eq:N-estimate-stationary} of
	Theorem~\ref{260529130} with \((\ell,p)=(0,30)\), using
	\(\|\mu^{1/8}d^{(k)}\|_{\cH_{\kin,\sigma}}
	\le C\|d^{(k)}\|_{\cH_{\kin,\sigma}}\) and
	\eqref{eq:nonlinear-iteration-uniform-bound}.  This gives
	\begin{equation}
		\|d^{(k+1)}\|_{\cH_{\kin,\sigma}}
		\le
		C_\cD\delta_{\st}\|d^{(k)}\|_{\cH_{\kin,\sigma}}.
		\label{260529647}
	\end{equation}
	Decrease \(\delta_{\st}\), if necessary, so that
	\(C_\cD\delta_{\st}\le1/2\) in \eqref{260529647}.  Since
	\(\|d^{(0)}\|_{\cH_{\kin,\sigma}}=\|f^{(1)}\|_{\cH_{\kin,\sigma}}\)
	is finite by \eqref{eq:nonlinear-iteration-uniform-bound}, the
	sequence \((f^{(k)})_{k\ge0}\) converges in
	\(\cH_{\kin,\sigma}\) to some \(f_{\st}\).
	
	\medskip
	\noindent\emph{Step 3. Passage to the limit.}
	Convergence in \(\cH_{\kin,\sigma}\) and \eqref{eq:local-H1-H_sigma}
	give \(w^{-1/2}f^{(k)}\to w^{-1/2}f_{\st}\) and
	\(w^{-3/2}D_vf^{(k)}\to w^{-3/2}D_vf_{\st}\) in
	\(L^2(\cD\times\bR_v^3)\).  Therefore, after taking a subsequence,
	\(f^{(k)}\to f_{\st}\) and \(D_vf^{(k)}\to D_vf_{\st}\) almost
	everywhere.  Fatou's lemma and
	\eqref{eq:nonlinear-iteration-uniform-bound} give
	\begin{equation}
		\|w^5f_{\st}\|_{L^\infty}
		+\|w^2D_vf_{\st}\|_{L^{30,2}_{x,v}}
		\le C_\cD\delta_{\st},
		\qquad
		\|\mu^{1/4}f_{\st}\|_{L^\infty}\le\varepsilon_0,
		\label{eq:nonlinear-limit-auxiliary-bounds}
	\end{equation}
	and, since the difference quotient defining the
	\(C^{\alpha_0}_{\kin}\) seminorm is lower semicontinuous under almost
	everywhere convergence,
	\begin{equation}
		\|f_{\st}\|_{C^{\alpha_0}_{\kin}}
		\le C_\cD\delta_{\st} .
		\label{eq:nonlinear-limit-holder-bound}
	\end{equation}
	Moreover, \eqref{boundedness of L} and
	\eqref{eq:N-estimate-stationary} yield
	\begin{align*}
		\|N[f^{(k)},f^{(k+1)}]-N[f_{\st},f_{\st}]\|_{L^2_xH^{-1}_\sigma}
		\le\,&
		C\|\mu^{1/4}f^{(k)}\|_{L^\infty}
		\|f^{(k+1)}-f_{\st}\|_{\cH_{\kin,\sigma}}
		\\
		&\qquad
		+C\|f^{(k)}-f_{\st}\|_{\cH_{\kin,\sigma}}
		\left(
		\|w^5f_{\st}\|_{L^\infty}
		+\|w^2D_vf_{\st}\|_{L^{30,2}_{x,v}}
		\right),
	\end{align*}
	which tends to zero.  Passing to the limit in \eqref{260527636}, and
	using continuity of the kinetic trace, proves that \(f_{\st}\) solves
	\eqref{2605041257}.  This also proves \eqref{eq:main-inflow-st-Linf}.
	
	\medskip
	\noindent\emph{Step 4. Weighted estimates.}
	Regard \eqref{2605041257} as the linear problem
	\eqref{260526322} with \(g=f_{\st}\).  Lemma~\ref{260527559},
	\eqref{eq:nonlinear-limit-auxiliary-bounds}, and
	\eqref{eq:stationary-boundary-weight-monotonicity} give
	\eqref{eq:main-inflow-st-regularity}.
	
	By \eqref{eq:nonlinear-limit-holder-bound} and
	\eqref{eq:maxwellian-holder-multiplier}, the coefficient
	\(g=f_{\st}\) satisfies
	\eqref{eq:stationary-coefficient-holder-assumption} with
	\(\alpha=\alpha_0\) and \(K=C_*\).  Thus the gradient part of
	Lemma~\ref{260527559} applies for every \(p\in(1,\infty)\).  Since
	\(A_{\st}\ge B^{\st}_1(\alpha_0)\),
	\eqref{eq:stationary-boundary-weight-monotonicity} gives
	\eqref{eq:main-inflow-st-gradient}.
	
	\medskip
	\noindent\emph{Step 5. Stability and uniqueness.}
	Set
	\[
	f^\delta:=f_{\st}^1-f_{\st}^2,
	\qquad
	h^\delta:=h_{\st}^1-h_{\st}^2.
	\]
	Since \(\|\mu^{1/4}f_{\st}^2\|_{L^\infty}\le\varepsilon_0\), we may
	apply Lemma~\ref{260527559} to \(f_{\st}^2\), viewed as the solution
	of \eqref{260526322} with \(g=f_{\st}^2\).  The first part of the
	lemma gives
	\(\|f_{\st}^2\|_{C^{\alpha_0}_{\kin}}\le C_\cD\delta_{\st}\le1\), hence
	\eqref{eq:stationary-coefficient-holder-assumption} with
	\(\alpha=\alpha_0\) and \(K=C_*\) by
	\eqref{eq:maxwellian-holder-multiplier}.  The weighted estimates of
	the lemma with \(\ell=5\) and \((\ell,p)=(2,30)\) then give
	\begin{equation}
		\|w^5f_{\st}^2\|_{L^\infty}
		+\|w^2D_vf_{\st}^2\|_{L^{30,2}_{x,v}}
		\le C_\cD\delta_{\st}.
		\label{eq:nonlinear-second-solution-regularity}
	\end{equation}
	The difference satisfies
	\[
	(v\cdot\nabla_x+L)f^\delta
	=N[f_{\st}^1,f^\delta]+N[f^\delta,f_{\st}^2],
	\qquad
	f^\delta|_{\gamma^-}=h^\delta.
	\]
	Theorem~\ref{260322442}, with coefficient \(f_{\st}^1\), together
	with \eqref{eq:N-estimate-stationary} and
	\eqref{eq:nonlinear-second-solution-regularity}, yields
	\[
	\|f^\delta\|_{\cH_{\kin,\sigma}}
	\le
	C_\cD\delta_{\st}\|f^\delta\|_{\cH_{\kin,\sigma}}
	+C_\cD\|h^\delta\|_{L^2(\gamma^-)}.
	\]
	After one final decrease of \(\delta_{\st}\), the first term is absorbed,
	which proves \eqref{eq:main-inflow-st-stability}.  Taking
	\(h^\delta=0\) proves uniqueness in the class of $f$ for which $\|\mu^{1/4}f\|_{L^\infty}\leq \varepsilon_0$.
\end{proof}

\subsection{Global well-posedness of the evolution problem}\label{subsec:inflow-global-proof}

In this subsection, we prove Theorem~\ref{thm:global-dynamic-landau} in the normalized variables.
Writing $F=\mu+\sqrt\mu\,f_\pt$, $F_0=\mu+\sqrt\mu\,f_0$ and $H=\mu+\sqrt\mu\,h_\pt$, the evolution problem reads
\begin{equation}
	\begin{cases}
		(\partial_t+v\cdot\nabla_x)f_{\pt}+Lf_{\pt}
		=N[f_{\pt},f_{\pt}]
		&\text{in }\bR_+\times\cD\times\bR_v^3,\\
		f_{\pt}(0)=f_0
		&\text{in }\cD\times\bR_v^3,\\
		f_{\pt}|_{\bR_+\times\gamma^-}=h_{\pt}
		&\text{on }\bR_+\times\gamma^-.
	\end{cases}
	\label{2605041257_dy}
\end{equation}

We use the stationary and evolutionary boundary norms from
Subsection~\ref{sec:inflow-anisotropic-extensions}.
The limit \(h(0):=\lim_{t\downarrow0}h(t)\) is understood in the
weighted uniform topology.  Thus the corner compatibility condition is
\begin{equation}
	h(0)=\operatorname{Tr}_{\gamma^-}f_0.
	\label{eq:evolution-data-compatibility}
\end{equation}
The weighted \(C^{1,1}\)-equivalence
\eqref{eq:inflow-weighted-C11-equivalence} implies the usual weight
monotonicity for these boundary norms and for the weighted
anisotropic norm of \(f_0\).
Throughout this subsection, constants depend only on the parameters
displayed in their subscripts and on the \(C^{1,1}\)-character of
\(\cD\).  In particular, they are independent of \(T\) and of the data.

\begin{lemma}[Evolutionary linearized Landau problem]
	\label{260527559_dy}
	Let \(\varepsilon_0=\varepsilon_0(\cD)>0\) be the smallness constant
	in Theorem~\ref{260322442}, let
	\(\alpha_0=\alpha_0(\cD)\in(0,1)\) be the exponent in
	Theorem~\ref{260526316}, and let \(A_{\st}\) and
	\(\delta_{\st}\) be the constants in
	Theorem~\ref{Landau_st_exi}.
	
	There is an integer \(B^{\dy}_0=B^{\dy}_0(\cD)\) with the following property.
	Let \(0<T<\infty\), \(\ell\ge0\), and suppose that the stationary
	inflow datum satisfies
	\begin{equation}
		\|w^{A_{\st}}h_{\st}\|_
		{(C_x^{0,1}\cap C_v^{1,1})(\gamma^-)}
		\le\delta_{\st}.
		\label{eq:evolutionary-stationary-data-assumption}
	\end{equation}
	Let \(f_{\st}\) be the steady solution given by
	Theorem~\ref{Landau_st_exi}.  Assume that
	\begin{equation}
		\|\mu^{1/4}g\|_{L^\infty(\cD_T\times\bR_v^3)}
		\le\varepsilon_0,
		\label{eq:evolutionary-coefficient-smallness}
	\end{equation}
	assume that
	\(\mu^{1/8}(g-f_{\st})\in\cH_{\kin,\sigma}(\cD_T\times\bR_v^3)\),
	and set
	\begin{equation}
		\cG_T
		:=
		\|\mu^{1/8}(g-f_{\st})\|_{
			\cH_{\kin,\sigma}(\cD_T\times\bR_v^3)}.
		\label{eq:evolutionary-coefficient-difference-size}
	\end{equation}
	
	Let \(f_0\) be initial data, \(h_{\pt}\) be the inflow data satisfying
	\eqref{eq:evolution-data-compatibility}, $S$ be the forcing term, and put
	\[
	h_{\dy}:=h_{\pt}-h_{\st},
	\]
	where \(h_{\st}\) is regarded as time independent.
	For
	\begin{align}\label{260915524}
		\begin{aligned}
			\cM_{m}(T)
			:={}&
			\|w^{m}f_0\|_
			{(C_x^{0,1}\cap C_v^{1,1})(\cD\times\bR_v^3)}
			\\
			&+\|w^{m+1}h_{\pt}\|_
			{(C_{t,x}^{0,1}\cap C_v^{1,1})(\gamma_T^-)}+(1+\cG_T)
			\|w^{m}h_{\st}\|_
			{(C_x^{0,1}\cap C_v^{1,1})(\gamma^-)}
			+\|w^{m}h_{\dy}\|_{L^2(\gamma_T^-)}\\
			&+\|w^{m+1} S\|_{L^\infty(\cD_T\times \bR^3_v)}+\|w^mS\|_{L^2(\cD_T;H^{-1}_\sigma)},
		\end{aligned}
	\end{align}
	assume that
	\begin{align}
		\cM_{\ell+B^{\dy}_0}(T)<\infty.
		\label{eq:evolutionary-regularity-data-quantity}
	\end{align}
	Then
	\begin{equation}
		\begin{cases}
			(\partial_t+v\cdot\nabla_x)f+Lf=N[g,f]+S
			&\text{in }\cD_T\times\bR_v^3,\\
			f(0)=f_0
			&\text{in }\cD\times\bR_v^3,\\
			f|_{\gamma_T^-}=h_{\pt}
			&\text{on }\gamma_T^-
		\end{cases}
		\label{260526322_1111}
	\end{equation}
	has a unique solution
	\(f\in\cH_{\kin,\sigma}(\cD_T\times\bR_v^3)\).
	Moreover,
	\begin{equation}
		\|w^\ell N[g-f_{\st},f_{\st}]\|_{
			L^2(\cD_T;H^{-1}_\sigma)}
		\le
		C_{\ell,\cD}\cG_T
		\|w^{\ell+B^{\dy}_0}h_{\st}\|_
		{(C_x^{0,1}\cap C_v^{1,1})(\gamma^-)},
		\label{eq:evolutionary-nonlinear-source-bound}
	\end{equation}
	and
	\begin{align}
		\begin{aligned}
			&\|w^\ell(f-f_{\st})\|_{
				\cH_{\kin,\sigma}(\cD_T\times\bR_v^3)}+\|w^\ell f\|_{L^\infty(\cD_T\times\bR_v^3)}
			+\|w^\ell f\|_{C^{\alpha_0}_{\kin}(\cD_T\times\bR_v^3)}
			\le C_{\ell,\cD}\cM_{\ell+B^{\dy}_0}(T),
			\label{260728655_1111}
		\end{aligned}	
	\end{align}
	
	Suppose in addition that, for some \(\alpha\in(0,1)\) and \(K>0\),
	\begin{equation}
		\|\mu^{1/4}g\|_{C^\alpha_{\kin}(\cD_T\times\bR_v^3)}\le K.
		\label{eq:evolutionary-coefficient-holder-assumption}
	\end{equation}
	There is an integer
	\(B^{\dy}_1=B^{\dy}_1(\cD,\alpha)\ge B^{\dy}_0\) such that, for every
	\(p\in(1,\infty)\),
	\begin{align}
		&\|w^\ell D_vf\|_{L^{\infty,p,2}_{t,x,v}(\cD_T\times\bR_v^3)}
		\le
		C_{\ell,\cD,\alpha,p,K}\cM_{\ell+B^{\dy}_1}(T),
		\label{eq:evolutionary-mixed-gradient-estimate}
	\end{align}
	whenever the quantity on the right is finite.
	Every constant in this
	lemma is independent of \(T\).
\end{lemma}

\begin{proof}
	Set
	\begin{equation}
		B^{\dy}_0:=A_{\st}+38 .
		\label{eq:dynamic-B0}
	\end{equation}
	
	\medskip
	\noindent\emph{Step 1. Energy estimate.}
	Set
	\[
	F_{\st}:=N[g-f_{\st},f_{\st}].
	\]
	Theorem~\ref{260529130}, followed by
	\eqref{eq:main-inflow-st-regularity} and \eqref{eq:main-inflow-st-gradient}, yields
	\begin{align}
		\|w^mF_{\st}\|_{L^2(\cD_T;H^{-1}_\sigma)}
		&\le
		C_{m,\cD}\cG_T
		\left(
		\|w^{m+5}f_{\st}\|_{L^\infty}
		+\|w^{m+2}D_vf_{\st}\|_{L^{30,2}_{x,v}}
		\right)
		\notag\\
		&\le
		C_{m,\cD}\cG_T
		\|w^{m+A_{\st}+5}h_{\st}\|_
		{(C_x^{0,1}\cap C_v^{1,1})(\gamma^-)}.
		\label{eq:dynamic-stationary-source-bound}
	\end{align}
	Also, since \(w^{-2}\in L^2(\bR_v^3)\) and \(\cD\) is bounded,
	\begin{align*}
		\|w^m(f_0-f_{\st})\|_{L^2}
		&\le
		C_{m,\cD}
		\left(
		\|w^{m+2}f_0\|_{(C_x^{0,1}\cap C_v^{1,1})(\cD\times\bR_v^3)}
		+\|w^{m+A_{\st}+2}h_{\st}\|_
		{(C_x^{0,1}\cap C_v^{1,1})(\gamma^-)}
		\right),
	\end{align*}
	where the second term comes from
	\(\|w^mf_{\st}\|_{L^2}\le C\|w^{m+2}f_{\st}\|_{L^\infty}\) and
	\eqref{eq:main-inflow-st-regularity}.  The weight \(A_{\st}+5\)
	simultaneously dominates the two extra weights needed to embed the
	uniform norm of \(f_0\) into \(L^2_{x,v}\) and the stationary weight
	losses in \eqref{eq:dynamic-stationary-source-bound}.
	
	By the two estimates above and the definition of
	\(\cM_{m+A_{\st}+5}(T)\), Theorem~\ref{260322442}, with
	coefficient \(g\), gives a unique
	\(f_{\dy}\in\cH_{\kin,\sigma}(\cD_T\times\bR_v^3)\) solving
	\begin{equation}
		\begin{cases}
			(\partial_t+v\cdot\nabla_x)f_{\dy}+Lf_{\dy}
			=N[g,f_{\dy}]+F_{\st}+S,\\
			f_{\dy}(0)=f_0-f_{\st},\\
			f_{\dy}|_{\gamma_T^-}=h_{\dy}.
		\end{cases}
		\label{eq:dynamic-difference-equation}
	\end{equation}
	The weighted estimate in Theorem~\ref{thm:weighted-steady} gives
	\begin{equation}
		\|w^mf_{\dy}\|_{\cH_{\kin,\sigma}(\cD_T\times\bR_v^3)}
		\le C_{m,\cD}\cM_{m+A_{\st}+5}(T).
		\label{eq:dynamic-weighted-Hkin-Ben}
	\end{equation}
	Since \(f_{\st}\in\cH_{\kin,\sigma}(\cD\times\bR_v^3)\) and, by
	\eqref{eq:main-inflow-st-regularity} with \(\ell=2\),
	\(f_{\st}\in L^2(\cD\times\bR_v^3)\), the characterization in
	Proposition~\ref{H_kin=weak_landau} shows that the time-independent
	extension of \(f_{\st}\) belongs to
	\(\cH_{\kin,\sigma}(\cD_T\times\bR_v^3)\).  Set
	\(f:=f_{\st}+f_{\dy}\).  The stationary equation and the
	identity
	\[
	F_{\st}+N[f_{\st},f_{\st}]=N[g,f_{\st}]
	\]
	show directly that \(f\) solves \eqref{260526322_1111}.  Uniqueness
	follows from the uniqueness of \eqref{eq:dynamic-difference-equation}.
	This proves solvability, uniqueness, and the source estimate, as
	\(A_{\st}+5\le B^{\dy}_0\) by \eqref{eq:dynamic-B0}.  The last term of
	\(\cM_{m+A_{\st}+5}(T)\) is stronger than the boundary term
	needed in the weighted energy estimate.
	
	\medskip
	\noindent\emph{Step 2. The inflow lifting and the \(L^\infty\) estimate.}
	Let
	\begin{equation*}
		H:=\cI_{\rm ev}(f_0,h_{\pt}),
	\end{equation*}
	where \(\cI_{\rm ev}\) is given by
	Theorem~\ref{thm:universal-anisotropic-inflow-extension}. Then $H(0)=f_0$ and $H|_{\gamma_T^-}=h_{\pt}$, and for \(q\ge0\),
	\begin{align}
		&\sum_{|\beta|\le2}\|w^qD_v^\beta H\|_{L^\infty}
		+\|w^q\partial_tH\|_{L^\infty}
		+\|w^q\nabla_xH\|_{L^\infty}
		\notag\\
		&\qquad\le
		C_{q,\cD}
		\left(
		\|w^qf_0\|_{(C_x^{0,1}\cap C_v^{1,1})(\cD\times\bR_v^3)}
		+\|w^{q+1}h_{\pt}\|_
		{(C_{t,x}^{0,1}\cap C_v^{1,1})(\gamma_T^-)}
		\right).
		\label{eq:dynamic-lift-derivative-bound}
	\end{align}
	Set \(\bar f:=f-H\).  Then \(\bar f(0)=0\),
	\(\bar f|_{\gamma_T^-}=0\), and
	\begin{equation}
		(\partial_t+v\cdot\nabla_x)\bar f=\cA_g\bar f+\bar S,
		\qquad
		\bar S:=S+\cK f-
		\bigl[(\partial_t+v\cdot\nabla_x)H-\cA_gH\bigr].
		\label{eq:dynamic-lifting-source}
	\end{equation}
	In divergence form this is
	\begin{equation}
		(\partial_t+v\cdot\nabla_x)\bar f
		=D_i(a_g^{ij}D_j\bar f)+b_g^iD_i\bar f+G,
		\qquad G:=c_g\bar f+\bar S.
		\label{eq:dynamic-bar-f-regularity-equation}
	\end{equation}
	
	Extend \(\bar f\), \(\bar S\), and \(G\) by zero to \(t<0\).  Extend
	\(g\) by zero from \((0,T)\times\cD\times\bR_v^3\) to
	\(\bR\times\bR_x^3\times\bR_v^3\).  Since \(\bar f(0)=0\), the zero
	initial trace creates no distributional term.  For \(\tau\in(0,T]\),
	put \(I_\tau:=(\tau-4,\tau)\cap(0,T)\).  Then \(|I_\tau|\le4\).
	Write \(\bar f=f_{\dy}+f_{\st}-H\).  The time-uniform trace
	estimate \eqref{260730346}, applied to \(w^mf_{\dy}\), together
	with \eqref{eq:dynamic-weighted-Hkin-Ben}, bounds
	\(\|w^mf_{\dy}\|_{L^2(I_\tau\times\cD\times\bR_v^3)}\) by
	\(2\|w^mf_{\dy}\|_{\cH_{\kin,\sigma}(\cD_T\times\bR_v^3)}\).  The two
	remaining terms are bounded through their weighted uniform norms,
	using \eqref{eq:main-inflow-st-regularity} and
	\eqref{eq:dynamic-lift-derivative-bound} together with
	\(w^{-2}\in L^2(\bR_v^3)\).  Hence
	\begin{equation}
		\|w^m\bar f\|_{L^2(I_\tau\times\cD\times\bR_v^3)}
		\le C_{m,\cD}\cM_{m+A_{\st}+5}(T).
		\label{eq:dynamic-bar-f-local-L2}
	\end{equation}
	In \(S\) the compact part acts on \(f=\bar f+H\).  We split
	\(\cK f=\cK\bar f+\cK H\) and absorb
	\(\cK H\) into the terms carrying the data, using
	\(\|\mu^{1/8}H(t,x,\cdot)\|_{L^2_v}\le C\|w^{q}H\|_{L^\infty}\) for
	any \(q\ge0\).  The coefficient bounds in
	Lemma~\ref{lem:perturbed-landau-ellipticity}, the pointwise estimate
	\eqref{eq:compact-part-pointwise} for \(\cK\), and
	\eqref{eq:dynamic-lift-derivative-bound} then imply
	\begin{align}
		w^m|\bar S(t,x,v)|
		\le{}& C_m\mu^{1/8}(v)
		\|\mu^{1/8}\bar f(t,x,\cdot)\|_{L^2_v}
		+C_m\|w^{m+21}f_0\|_
		{(C_x^{0,1}\cap C_v^{1,1})(\cD\times\bR_v^3)}w(v)^{-20}
		\notag\\
		&+C_m\|w^{m+22}h_{\pt}\|_
		{(C_{t,x}^{0,1}\cap C_v^{1,1})(\gamma_T^-)}w(v)^{-20}
		+C_m\|w^{m+20}S\|_
		{L^\infty(\cD_T\times \bR^3_v)}w(v)^{-20}.
		\label{eq:dynamic-source-pointwise}
	\end{align}
	For fixed \(m\ge0\), set \(u_m:=w^m\bar f\).  Conjugating
	\eqref{eq:dynamic-lifting-source} by \(w^m\) gives
	\[
	(\partial_t+v\cdot\nabla_x)u_m
	=\cA_{g,m}u_m+w^m\bar S,
	\qquad u_m(0)=0,\qquad (u_m)_{\gamma_T^-}=0.
	\]
	Thus the zero extension of \(u_m\) to \(t<0\) satisfies the equation
	without an initial-time distribution.  Since \(m\ge0\) we have
	\(\|\mu^{1/8}\bar f(t,x,\cdot)\|_{L^2_v}
	\le\|\mu^{1/8}u_m(t,x,\cdot)\|_{L^2_v}\).  Hence estimate
	\eqref{eq:dynamic-source-pointwise} has the source form
	\eqref{eq:L2-Linfty-source-bound} required in
	Theorem~\ref{260504554}, again with \(\varepsilon=1/8\).  Apply that
	theorem on each translated interval \((\tau-4,\tau]\).  Its loss of
	\(21\) velocity weights and
	\eqref{eq:dynamic-bar-f-local-L2} give, uniformly in \(\tau\),
	\begin{equation}
		\|w^m\bar f\|_{L^\infty(\cD_T\times\bR_v^3)}
		\le C_{m,\cD}\cM_{m+A_{\st}+26}(T).
		\label{eq:dynamic-bar-f-Linfty-Binfty}
	\end{equation}
	Using this estimate in \eqref{eq:dynamic-source-pointwise}, and then
	using \(|c_g(v)|\le Cw(v)^{-1}\), yields
	\begin{equation}
		\|w^m\bar S\|_{L^\infty}
		+\|w^mG\|_{L^\infty}
		\le C_{m,\cD}\cM_{m+A_{\st}+26}(T).
		\label{eq:dynamic-full-source-Linfty}
	\end{equation}
	Since \(f=\bar f+H\), the same bound holds for
	\(\|w^mf\|_{L^\infty}\).
	
	\medskip
	\noindent\emph{Step 3. Kinetic H\"older estimate.}
	Fix
	\(Z=(t,x,v)\) and \(Z'=(s,y,v')\) in
	\(\cD_T\times\bR_v^3\), and let \(\rho:=\mathbf d(Z,Z')<1/2\).
	After interchanging the points, assume \(s\le t\).  Applying
	Theorem~\ref{260526316} to
	\eqref{eq:dynamic-bar-f-regularity-equation} on
	\((t-4,t]\), and using
	\eqref{eq:dynamic-bar-f-Linfty-Binfty} and
	\eqref{eq:dynamic-full-source-Linfty}, gives
	\begin{align}
		&(w^\ell(v)+w^\ell(v'))
		|\bar f(Z)-\bar f(Z')|
		\le
		C_{\ell,\cD}\cM_{\ell+B^{\dy}_0}(T)\rho^{\alpha_0}.
		\label{eq:dynamic-bar-f-unweighted-holder}
	\end{align}
	The definition of the kinetic distance also gives
	\begin{equation*}
		|v-v'|\le4\rho,
		\qquad w(v')\simeq w(v).
	\end{equation*}
	Consequently,
	\[
	|w^\ell(v)-w^\ell(v')|\,|\bar f(Z')|
	\le
	C_{\ell,\cD}\cM_{\ell+A_{\st}+26}(T)\rho
	\le
	C_{\ell,\cD}\cM_{\ell+B^{\dy}_0}(T)\rho^{\alpha_0}.
	\]
	Combining this with
	\eqref{eq:dynamic-bar-f-unweighted-holder} proves the weighted
	H\"older estimate for \(\bar f\).  For \(\rho\ge1/2\), it follows
	directly from \eqref{eq:dynamic-bar-f-Linfty-Binfty}.
	
	It remains to estimate \(H\).  When \(0<\rho<1\), the definition of
	\(\mathbf d\) gives
	\begin{equation}
		|t-s|+|x-y|+|v-v'|\le Cw(v)\rho,
		\qquad w(v')\simeq w(v).
		\label{eq:dynamic-distance-euclidean-control}
	\end{equation}
	Equations \eqref{eq:evolutionary-inflow-extension-bound} and
	\eqref{eq:dynamic-distance-euclidean-control} therefore imply
	\begin{equation*}
		\|w^\ell H\|_{C^{\alpha_0}_{\kin}(\cD_T\times\bR_v^3)}
		\le
		C_{\ell,\cD}
		\left(
		\|w^{\ell+1}f_0\|_{(C_x^{0,1}\cap C_v^{1,1})(\cD\times\bR_v^3)}
		+\|w^{\ell+2}h_{\pt}\|_
		{(C_{t,x}^{0,1}\cap C_v^{1,1})(\gamma_T^-)}
		\right).
	\end{equation*}
	\eqref{260728655_1111} follows from
	\eqref{eq:dynamic-weighted-Hkin-Ben} and \(A_{\st}+5\le B^{\dy}_0\), by \eqref{eq:dynamic-B0}.
	
	\medskip
	\noindent\emph{Step 4. Mixed-norm estimate for \(D_vf\).}
	Assume \eqref{eq:evolutionary-coefficient-holder-assumption}.
	Apply Lemma~\ref{lem:norm-preserving-holder-extension}, with the
	metric \(\mathbf d\), to \(\mu^{1/4}g\), and define the extension of
	\(g\) to \(\bR\times\bR_x^3\times\bR_v^3\) by
	\(g:=\mu^{-1/4}(\mu^{1/4}g)\) outside \(\cD_T\times\bR_v^3\).
	Only the values of \(g\) on \(\cD_T\times\bR_v^3\) enter
	\eqref{260526322_1111}.  Therefore this extension changes neither
	the equation nor the solution.  It replaces, from here on, the
	extension by zero used in Step~2.  Since
	Lemma~\ref{lem:norm-preserving-holder-extension} preserves both the
	supremum norm and the H\"older seminorm, the extension satisfies
	\begin{equation*}
		\|\mu^{1/4}g\|_{L^\infty(\bR^7)}\le\varepsilon_0,
		\qquad
		[\mu^{1/4}g]_{C^\alpha_{\kin}(\bR^7)}
		\le K.
	\end{equation*}
	Lemma~\ref{260527138} supplies the corresponding coefficient
	H\"older bounds.
	
	Let
	\[
	A':=A(2\alpha/3),
	\]
	where \(A(\cdot)\) is defined in
	\eqref{eq:global-gradient-weight-loss}, and set
	\begin{equation}
		B^{\dy}_1:=\max\{B^{\dy}_0,A_{\st}+A'+29\}.
		\label{eq:dynamic-B1}
	\end{equation}
	Applying Theorem~\ref{260527139} to
	\eqref{eq:dynamic-bar-f-regularity-equation} on the same moving
	time intervals, with the full source \(G\), gives first
	\begin{align*}
		&\|w^\ell D_v\bar f\|_{L^{\infty,p,2}_{t,x,v}}
		\quad\le
		C_{\ell,\cD,\alpha,p,K}
		\left(
		\|w^{\ell+A'+3}\bar f\|_{L^\infty}
		+\|w^{\ell+3}G\|_{L^\infty}
		\right).
	\end{align*}
	By \eqref{eq:dynamic-bar-f-Linfty-Binfty},
	\eqref{eq:dynamic-full-source-Linfty}, and the definition of \(B^{\dy}_1\),
	this becomes
	\begin{equation*}
		\|w^\ell D_v\bar f\|_{L^{\infty,p,2}_{t,x,v}}
		\le
		C_{\ell,\cD,\alpha,p,K}\cM_{\ell+B^{\dy}_1}(T).
	\end{equation*}
	Finally, since \(B^{\dy}_1\ge2\),
	\begin{align*}
		\|w^\ell D_vH\|_{L^{\infty,p,2}_{t,x,v}}
		&\le
		C_{\cD,p}\|w^{\ell+B^{\dy}_1}D_vH\|_{L^\infty}
		\notag\\
		&\le
		C_{\ell,\cD,p}\left(
		\|w^{\ell+B^{\dy}_1}f_0\|_{(C_x^{0,1}\cap C_v^{1,1})(\cD\times\bR_v^3)}
		+\|w^{\ell+B^{\dy}_1+1}h_{\pt}\|_
		{(C_{t,x}^{0,1}\cap C_v^{1,1})(\gamma_T^-)}
		\right).
	\end{align*}
	Together with \(f=\bar f+H\), this proves
	\eqref{eq:evolutionary-mixed-gradient-estimate}.
\end{proof}

\begin{proof}[\textbf{Proof of Theorem~\ref{thm:global-dynamic-landau}}]
	Recall that $\varepsilon_0$ is the constant of Theorem~\ref{260322442} and $\alpha_0$ the exponent of Theorem~\ref{260526316}, as in the proof of Theorem~\ref{Landau_st_exi}.
	Let \(B^{\dy}_0=B^{\dy}_0(\cD)\) be the weight \eqref{eq:dynamic-B0} in Lemma~\ref{260527559_dy}, and
	let \(B^{\dy}_1(\alpha_0)=B^{\dy}_1(\cD,\alpha_0)\) be the gradient weight \eqref{eq:dynamic-B1} in the
	same lemma with \(\alpha=\alpha_0\).  Fix
	\begin{equation}
		A_{\dy}:=\max\{B^{\dy}_0+5,\,B^{\dy}_1(\alpha_0)+2\}.
		\label{eq:nonlinear-dynamic-A0-choice}
	\end{equation}
	By \eqref{eq:dynamic-B0},
	\(B^{\dy}_0=A_{\st}+38\); hence
	\(A_{\dy}\ge A_{\st}+5\), and the weight monotonicity
	\eqref{eq:stationary-boundary-weight-monotonicity} for the boundary
	norms shows that every stationary datum occurring below is dominated
	by \(\mathsf D_{A_{\dy}}\).  In particular, after the first of the
	finitely many decreases of \(\delta_{\dy}\le1\) performed in this
	proof, all of which depend only on \(\cD\), we may assume
	\(C_\cD\delta_{\dy}\le\delta_{\st}\), so that
	\eqref{eq:dynamic-nonlinear-smallness-main} implies the stationary
	smallness condition
	\eqref{eq:evolutionary-stationary-data-assumption} and
	Theorem~\ref{Landau_st_exi} supplies \(f_{\st}\).
	
	\medskip
	\noindent\emph{Step 1. Iteration and uniform estimates.}
	Put \(f^{(0)}:=f_{\st}\), regarded as time independent.  Given
	\(f^{(k)}\), let \(f^{(k+1)}\) solve
	\begin{equation}
		\begin{cases}
			(\partial_t+v\cdot\nabla_x)f^{(k+1)}+Lf^{(k+1)}
			=N[f^{(k)},f^{(k+1)}],\\
			f^{(k+1)}(0)=f_0,\\
			f^{(k+1)}|_{\bR_+\times\gamma^-}=h_{\pt}.
		\end{cases}
		\label{eq:nonlinear-dynamic-iteration}
	\end{equation}
	The problem is first solved on \((0,T)\) for each finite \(T\).
	Uniqueness in Lemma~\ref{260527559_dy} makes the solutions for
	different \(T\) consistent under restriction, so that
	\eqref{eq:nonlinear-dynamic-iteration} is solved on \(\bR_+\).  Set
	\[
	X_k(T):=
	\|\mu^{1/8}(f^{(k)}-f_{\st})\|_{
		\cH_{\kin,\sigma}(\cD_T\times\bR_v^3)},
	\]
	which is exactly the coefficient-difference size
	\eqref{eq:evolutionary-coefficient-difference-size} for
	\(g=f^{(k)}\).  Assume inductively that, for every \(T>0\),
	\begin{equation}
		\|\mu^{1/4}f^{(k)}\|_{L^\infty(\cD_T\times\bR_v^3)}
		\le\varepsilon_0,
		\qquad
		\|f^{(k)}\|_{C^{\alpha_0}_{\kin}(\cD_T\times\bR_v^3)}\le1,
		\qquad
		X_k(T)\le1.
		\label{eq:nonlinear-dynamic-induction-hypothesis}
	\end{equation}
	The first hypothesis is
	\eqref{eq:evolutionary-coefficient-smallness}, and by
	\eqref{eq:maxwellian-holder-multiplier} the second gives
	\eqref{eq:evolutionary-coefficient-holder-assumption} for
	\(g=f^{(k)}\) with \(\alpha=\alpha_0\) and \(K=C_*\).  Therefore
	Lemma~\ref{260527559_dy} applies with \(g=f^{(k)}\), and the constants
	below depend neither on \(k\) nor on \(T\).
	
	When the lemma is applied, its data expression is bounded, by
	\eqref{eq:dynamic-nonlinear-data-quantity-main}, as follows: for
	\(m,B\ge0\),
	\begin{align}
		&\|w^{m+B}f_0\|_
		{(C_x^{0,1}\cap C_v^{1,1})(\cD\times\bR_v^3)}
		+\|w^{m+B+1}h_{\pt}\|_
		{(C_{t,x}^{0,1}\cap C_v^{1,1})(\gamma_T^-)}
		\notag\\
		&\quad +(1+\cG_T)
		\|w^{m+B}h_{\st}\|_
		{(C_x^{0,1}\cap C_v^{1,1})(\gamma^-)}
		+\|w^{m+B}(h_{\pt}-h_{\st})\|_{L^2(\gamma_T^-)}
		\notag\\
		&\qquad\le
		\mathsf D_{m+B}
		+\cG_T\|w^{m+B}h_{\st}\|_
		{(C_x^{0,1}\cap C_v^{1,1})(\gamma^-)},
		\label{eq:nonlinear-data-to-linear-data}
	\end{align}
	where \(\cG_T=X_k(T)\) for the present choice of coefficient.
	The energy estimate \eqref{260728655_1111} with \(\ell=0\), together
	with \(\|\mu^{1/8}u\|_{\cH_{\kin,\sigma}}
	\le C\|u\|_{\cH_{\kin,\sigma}}\),
	\eqref{eq:nonlinear-data-to-linear-data}, and
	\(A_{\dy}\ge B^{\dy}_0\), gives
	\begin{equation}
		X_{k+1}(T)\le C_\cD\mathsf D_{A_{\dy}}\bigl(1+X_k(T)\bigr).
		\label{eq:nonlinear-dynamic-energy-induction}
	\end{equation}
	Likewise \eqref{260728655_1111}, used with \(\ell=0\) and
	\(\ell=5\), and \eqref{eq:evolutionary-mixed-gradient-estimate}, used
	with \((\ell,p)=(2,30)\), give
	\begin{align}
		&\|f^{(k+1)}\|_{L^\infty}
		+\|w^5f^{(k+1)}\|_{L^\infty}
		+\|f^{(k+1)}\|_{C^{\alpha_0}_{\kin}}
		+\|w^2D_vf^{(k+1)}\|_{L^{\infty,30,2}_{t,x,v}}
		\notag\\
		&\qquad\le
		C_\cD\mathsf D_{A_{\dy}}\bigl(1+X_k(T)\bigr)
		\le2C_\cD\delta_{\dy},
		\label{eq:nonlinear-dynamic-uniform-regularity}
	\end{align}
	all norms being taken on \(\cD_T\times\bR_v^3\).  Here
	\eqref{eq:nonlinear-dynamic-A0-choice} was used to absorb the weight
	losses \(B^{\dy}_0+5\) and \(B^{\dy}_1(\alpha_0)+2\) into \(A_{\dy}\).
	
	Choose \(\delta_{\dy}\le1\) so that
	\begin{equation}
		2C_\cD\delta_{\dy}\le\min\{\varepsilon_0,1\}
		\label{eq:nonlinear-dynamic-delta-choice}
	\end{equation}
	in \eqref{eq:nonlinear-dynamic-energy-induction} and
	\eqref{eq:nonlinear-dynamic-uniform-regularity}.  Since
	\(\mu^{1/4}\le1\), these two estimates return
	\eqref{eq:nonlinear-dynamic-induction-hypothesis} for \(f^{(k+1)}\).
	The hypotheses hold for \(f^{(0)}=f_{\st}\): indeed \(X_0\equiv0\),
	while \eqref{eq:main-inflow-st-Linf} and \eqref{eq:main-inflow-st-regularity} with \(\ell=0\) give
	\(\|\mu^{1/4}f_{\st}\|_{L^\infty}\le\varepsilon_0\) and
	\(\|f_{\st}\|_{C^{\alpha_0}_{\kin}}\le C_\cD\delta_{\dy}\le1\).
	Hence the sequence is well defined and
	\eqref{eq:nonlinear-dynamic-energy-induction} and
	\eqref{eq:nonlinear-dynamic-uniform-regularity} hold uniformly in
	\(k\) and in \(T\).
	
	\medskip
	\noindent\emph{Step 2. Contraction and passage to the limit.}
	For \(k\ge0\), put \(d^{(k+1)}:=f^{(k+1)}-f^{(k)}\).  For \(k\ge1\), subtracting two
	consecutive instances of \eqref{eq:nonlinear-dynamic-iteration} and
	using the bilinearity of \(N\) gives
	\begin{equation}
		\begin{cases}
			(\partial_t+v\cdot\nabla_x)d^{(k+1)}+Ld^{(k+1)}
			=N[f^{(k)},d^{(k+1)}]+N[d^{(k)},f^{(k)}],\\
			d^{(k+1)}(0)=0,\\
			d^{(k+1)}|_{\gamma_T^-}=0.
		\end{cases}
		\label{eq:nonlinear-dynamic-iteration-difference}
	\end{equation}
	Apply Theorem~\ref{260322442}, with coefficient \(f^{(k)}\), to
	\eqref{eq:nonlinear-dynamic-iteration-difference}, and estimate its
	source by the evolutionary nonlinear estimate
	\eqref{eq:N-estimate-evolutionary} of Theorem~\ref{260529130} with
	\((\ell,p)=(0,30)\), using
	\(\|\mu^{1/8}d^{(k)}\|_{\cH_{\kin,\sigma}}
	\le C\|d^{(k)}\|_{\cH_{\kin,\sigma}}\) and
	\eqref{eq:nonlinear-dynamic-uniform-regularity}.  This gives
	\begin{equation}
		\|d^{(k+1)}\|_{\cH_{\kin,\sigma}(\cD_T\times\bR_v^3)}
		\le C_\cD\delta_{\dy}
		\|d^{(k)}\|_{\cH_{\kin,\sigma}(\cD_T\times\bR_v^3)}.
		\label{eq:nonlinear-dynamic-contraction}
	\end{equation}
	Decrease \(\delta_{\dy}\), if necessary, so that
	\(C_\cD\delta_{\dy}\le1/2\) in
	\eqref{eq:nonlinear-dynamic-contraction}.  The constant there is
	independent of \(T\), and \eqref{260728655_1111} with \(\ell=0\)
	shows that
	\[
	D_1:=\sup_{T>0}
	\|d^{(1)}\|_{\cH_{\kin,\sigma}(\cD_T\times\bR_v^3)}
	\le C_\cD\mathsf D_{A_{\dy}}<\infty.
	\]
	Taking the supremum over \(T\) in
	\eqref{eq:nonlinear-dynamic-contraction} and iterating therefore gives
	\[
	\sup_{T>0}
	\|d^{(k)}\|_{\cH_{\kin,\sigma}(\cD_T\times\bR_v^3)}
	\le2^{1-k}D_1,
	\qquad k\ge1,
	\]
	whence, for \(j>k\ge1\),
	\[
	\sup_{T>0}
	\|f^{(j)}-f^{(k)}\|_{\cH_{\kin,\sigma}(\cD_T\times\bR_v^3)}
	\le2^{1-k}D_1.
	\]
	Consequently there is
	\(f_{\pt}\in\bigcap_{T>0}\cH_{\kin,\sigma}(\cD_T\times\bR_v^3)\) such
	that
	\begin{equation}
		\sup_{T>0}
		\|f^{(k)}-f_{\pt}\|_{\cH_{\kin,\sigma}(\cD_T\times\bR_v^3)}
		\longrightarrow0,
		\qquad k\to\infty.
		\label{eq:nonlinear-dynamic-strong-convergence}
	\end{equation}
	
	By \eqref{eq:nonlinear-dynamic-strong-convergence} and
	\eqref{eq:local-H1-H_sigma}, after passing to a subsequence,
	\(f^{(k)}\to f_{\pt}\) and \(D_vf^{(k)}\to D_vf_{\pt}\) almost
	everywhere on \(\bR_+\times\cD\times\bR_v^3\).  Fatou's lemma and,
	for the H\"older seminorm, the lower semicontinuity of the difference
	quotients defining it under almost everywhere convergence, transfer
	\eqref{eq:nonlinear-dynamic-uniform-regularity} to the limit:
	\begin{equation}
		\|w^5f_{\pt}\|_{L^\infty}
		+\|f_{\pt}\|_{C^{\alpha_0}_{\kin}}
		+\|w^2D_vf_{\pt}\|_{L^{\infty,30,2}_{t,x,v}}
		\le2C_\cD\delta_{\dy},
		\label{eq:nonlinear-dynamic-limit-regularity}
	\end{equation}
	uniformly in \(T\).  Since \(\mu^{1/4}\le1\), the choice
	\eqref{eq:nonlinear-dynamic-delta-choice} turns
	\eqref{eq:nonlinear-dynamic-limit-regularity} into
	\eqref{eq:main-inflow-dynamic-small-Linf}.
	
	Finally, \eqref{boundedness of L}, applied to
	\(N[f^{(k)},f^{(k+1)}-f_{\pt}]\) with the coefficient bound
	\(\|\mu^{1/4}f^{(k)}\|_{L^\infty}\le\varepsilon_0\), and
	\eqref{eq:N-estimate-evolutionary}, applied to
	\(N[f^{(k)}-f_{\pt},f_{\pt}]\) with
	\eqref{eq:nonlinear-dynamic-limit-regularity}, show that
	\[
	N[f^{(k)},f^{(k+1)}]\longrightarrow N[f_{\pt},f_{\pt}]
	\quad\text{in }L^2(\cD_T;H^{-1}_\sigma)
	\]
	for every \(T>0\).  Continuity of the initial and kinetic trace maps
	then permits passage to the limit in
	\eqref{eq:nonlinear-dynamic-iteration}.  It follows that \(f_{\pt}\)
	solves \eqref{2605041257_dy} on every finite time interval, with
	initial trace \(f_0\) and inflow trace \(h_{\pt}\).
	
	\medskip
	\noindent\emph{Step 3. Weighted estimates.}
	Regard \eqref{2605041257_dy} as the linear problem
	\eqref{260526322_1111} with \(g=f_{\pt}\).  Keeping the data
	dependence explicit in \eqref{eq:nonlinear-data-to-linear-data}, and
	using \(\|w^{B^{\dy}_0}h_{\st}\|_{(C_x^{0,1}\cap C_v^{1,1})(\gamma^-)}
	\le C_\cD\mathsf D_{A_{\dy}}\le C_\cD\delta_{\dy}\), the energy
	estimate \eqref{260728655_1111} with \(\ell=0\) improves
	\eqref{eq:nonlinear-dynamic-energy-induction} to
	\[
	X_{k+1}(T)\le
	C_\cD\bigl(\mathsf D_{A_{\dy}}+\delta_{\dy}X_k(T)\bigr).
	\]
	Since \(C_\cD\delta_{\dy}\le1/2\) and \(X_0\equiv0\), induction gives
	\(X_k(T)\le2C_\cD\mathsf D_{A_{\dy}}\) uniformly in \(k\) and \(T\),
	and \eqref{eq:nonlinear-dynamic-strong-convergence} yields
	\begin{equation}
		\cG:=\sup_{T>0}
		\|\mu^{1/8}(f_{\pt}-f_{\st})\|_{
			\cH_{\kin,\sigma}(\cD_T\times\bR_v^3)}
		\le C_\cD\mathsf D_{A_{\dy}}.
		\label{eq:nonlinear-dynamic-coefficient-difference}
	\end{equation}
	Let \(\ell\ge0\) with \(\mathsf D_{\ell+A_{\dy}}<\infty\).  Inserting
	\eqref{eq:nonlinear-dynamic-coefficient-difference} into
	\eqref{eq:nonlinear-data-to-linear-data} with \(m=\ell\) and
	\(B=B^{\dy}_0\), and using \(\mathsf D_{A_{\dy}}\le1\) together with
	\(A_{\dy}\ge B^{\dy}_0\), bounds the data expression of
	Lemma~\ref{260527559_dy} by
	\(C_{\ell,\cD}\mathsf D_{\ell+A_{\dy}}\).  Hence
	\eqref{260728655_1111} gives
	\eqref{eq:main-inflow-dynamic-regularity}.
	
	By \eqref{eq:nonlinear-dynamic-limit-regularity} and
	\eqref{eq:maxwellian-holder-multiplier}, the coefficient
	\(g=f_{\pt}\) satisfies
	\eqref{eq:evolutionary-coefficient-holder-assumption} with
	\(\alpha=\alpha_0\) and \(K=C_*\).  Thus the gradient part of
	Lemma~\ref{260527559_dy} applies for every \(p\in(1,\infty)\).  Since
	\(A_{\dy}\ge B^{\dy}_1(\alpha_0)\), the same bound for the data expression,
	now with \(B=B^{\dy}_1(\alpha_0)\), gives
	\eqref{eq:global-dynamic-mixed-gradient-main}.
	
	\medskip
	\noindent\emph{Step 4. Stability and uniqueness.}
	Set
	\[
	d:=f_{\pt}^1-f_{\pt}^2.
	\]
	We first record the regularity of \(f_{\pt}^2\). By the smallness assumption \(\|\mu^{1/4}f_{\pt}^2\|_{L^\infty(\cD_T\times\bR_v^3)}\le\varepsilon_0\), we may apply
	Lemma~\ref{260527559_dy} to \(f_{\pt}^2\), viewed as the solution of
	\eqref{260526322_1111} with \(g=f_{\pt}^2\).  Its
	coefficient-difference term contributes
	\[
	\|\mu^{1/8}(f_{\pt}^2-f_{\st}^2)\|_
	{\cH_{\kin,\sigma}(\cD_T\times\bR_v^3)}
	\|w^{B^{\dy}_0}h_{\st}^2\|_
	{(C_x^{0,1}\cap C_v^{1,1})(\gamma^-)}
	\le
	C_\cD\delta_{\dy}
	\|f_{\pt}^2-f_{\st}^2\|_
	{\cH_{\kin,\sigma}(\cD_T\times\bR_v^3)}
	\]
	to the right-hand side of \eqref{260728655_1111} with \(\ell=0\), and
	is absorbed after the preceding choices of \(\delta_{\dy}\).  The
	pointwise and H\"older parts of the lemma then give
	\(\|w^5f_{\pt}^2\|_{L^\infty}
	+\|f_{\pt}^2\|_{C^{\alpha_0}_{\kin}}\le C_\cD\delta_{\dy}\le1\),
	hence \eqref{eq:evolutionary-coefficient-holder-assumption} with
	\(\alpha=\alpha_0\) and \(K=C_*\) by
	\eqref{eq:maxwellian-holder-multiplier}.  The weighted estimates of
	the lemma with \(\ell=5\) and \((\ell,p)=(2,30)\) therefore give
	\begin{equation}
		\|w^5f_{\pt}^2\|_{L^\infty(\cD_T\times\bR_v^3)}
		+\|w^2D_vf_{\pt}^2\|_{L^{\infty,30,2}_{t,x,v}
			(\cD_T\times\bR_v^3)}
		\le C_\cD\delta_{\dy},
		\label{eq:nonlinear-second-dynamic-solution-regularity}
	\end{equation}
	uniformly in \(T\).
	
	The difference satisfies
	\[
	\begin{cases}
		(\partial_t+v\cdot\nabla_x)d+Ld
		=N[f_{\pt}^1,d]+N[d,f_{\pt}^2],\\
		d(0)=f_0^1-f_0^2,\\
		d|_{\gamma_T^-}=h_{\pt}^1-h_{\pt}^2.
	\end{cases}
	\]
	Theorem~\ref{260322442}, with coefficient \(f_{\pt}^1\), together
	with \eqref{eq:N-estimate-evolutionary} and
	\eqref{eq:nonlinear-second-dynamic-solution-regularity}, yields
	\[
	\|d\|_{\cH_{\kin,\sigma}}
	\le
	C_\cD\delta_{\dy}\|d\|_{\cH_{\kin,\sigma}}
	+C_\cD\left(
	\|f_0^1-f_0^2\|_{L^2}
	+\|h_{\pt}^1-h_{\pt}^2\|_{L^2(\gamma_T^-)}
	\right).
	\]
	After one final decrease of \(\delta_{\dy}\), the first term is
	absorbed, which proves \eqref{eq:main-inflow-dynamic-stability} with a
	constant independent of \(T\).  Taking \(f_0^1=f_0^2\) and
	\(h_{\pt}^1=h_{\pt}^2\) proves uniqueness in the class
	\eqref{eq:main-inflow-dynamic-uniqueness-class}.
\end{proof}

\begin{prop}[Nonnegativity for prescribed inflow]
	\label{prop:inflow-nonnegativity}
	After decreasing the smallness constants in the stationary and
	nonlinear evolutionary theorems, if necessary, the following hold.
	\begin{enumerate}
		\item If the physical stationary inflow datum
		$H_{\st}=\mu+\sqrt\mu\,h_{\st}$ is nonnegative on $\gamma^-$,
		then the solution $F_{\st}=\mu+\sqrt\mu\,f_{\st}$ of
		Theorem~\ref{Landau_st_exi} is nonnegative almost everywhere.
		\item If $F_0=\mu+\sqrt\mu\,f_0\ge0$ and
		$H_{\pt}=\mu+\sqrt\mu\,h_{\pt}\ge0$ almost everywhere
		on the initial domain and incoming boundary, respectively,
		then the solution $F(t)=\mu+\sqrt\mu\,f_{\pt}(t)$ of
		Theorem~\ref{thm:global-dynamic-landau} is nonnegative
		almost everywhere for every $t\ge0$.
	\end{enumerate}
\end{prop}

\begin{proof}
	All estimates below are on $\cD\times\bR_v^3$ in the stationary
	case, and on $(0,T)\times\cD\times\bR_v^3$, with an arbitrary
	finite $T>0$, in the evolutionary case.
	Write $F=\mu+\sqrt\mu\,f$ for either solution and set
	\[
	a^{ij}:=\phi^{ij}*F,\qquad
	b^i:=\partial_{v_j}a^{ij},\qquad
	J^i:=a^{ij}\partial_{v_j}F-b^iF.
	\]
	The physical equation reads
	\[
	\mathcal TF=\partial_{v_i}J^i,
	\qquad
	\mathcal T=v\cdot\nabla_x
	\quad\text{or}\quad
	\mathcal T=\partial_t+v\cdot\nabla_x.
	\]
	The distributional kernel identity
	$\partial_{v_i}\partial_{v_j}\phi^{ij}=-8\pi\delta_0$
	(cf.\ \eqref{eq:double-divergence-sigma}) gives
	\begin{equation}\label{eq:inflow-pos-divb}
		\partial_{v_i}b^i=-8\pi F.
	\end{equation}
	We decrease the perturbative smallness constants, if necessary,
	so that $\|\mu^{1/4}f\|_{L^\infty}\le\varepsilon_{\mathrm{ell}}$
	from \eqref{eq:smallness-perturbed-sigma}.  Then
	Lemma~\ref{lem:perturbed-landau-ellipticity} gives
	\begin{equation}\label{eq:inflow-pos-ellipticity}
		\tfrac12\sigma^{ij}(v)\xi_i\xi_j
		\le a^{ij}(t,x,v)\xi_i\xi_j
		\le 2\sigma^{ij}(v)\xi_i\xi_j.
	\end{equation}
	This conclusion does not presuppose the sign of $F$.
	
	\smallskip
	\noindent\emph{A convex negative-part renormalization.}
	Let
	\[
	\beta(r)=r_- -\arctan(r_-),\qquad
	H(r)=r\beta'(r)-\beta(r),\qquad r_-:=\max\{-r,0\}.
	\]
	Then $\beta\in C^2(\bR)$ is convex,
	$\beta(0)=\beta'(0)=0$, and $\beta',\beta''$ are bounded.
	For $r=-s<0$, $s>0$, one has
	\[
	\beta''(-s)=\frac{2s}{(1+s^2)^2}>0,
	\qquad
	H(-s)=\arctan s-\frac{s}{1+s^2}>0,
	\]
	whereas $\beta''(r)=H(r)=0$ for $r\ge0$.  In particular,
	\begin{equation}\label{eq:inflow-pos-beta-sign}
		\beta\ge0,\quad
		\beta(r)=0\ \Longleftrightarrow\ r\ge0,\quad
		-rH(r)\ge0,\quad
		-rH(r)>0\ \text{when }r<0.
	\end{equation}
	Moreover $H'(r)=r\beta''(r)$ and
	\[
	0\le\beta(r)\le r_-,\qquad
	0\le H(r)\le r_-,\qquad
	|rH(r)|\le\tfrac\pi2 r_-.
	\]
	
	\smallskip
	\noindent\emph{Admissibility and velocity cutoffs.}
	The function $\sqrt\mu$ is a bounded $C^1$ velocity multiplier;
	$\mu$ belongs to both physical kinetic graph spaces.  Hence the
	Gaussian multiplier property and
	Lemma~\ref{lem:compact-velocity-support} imply that
	$F\in\cH_{\kin,\sigma}$, stationary or on $[0,T]$, respectively.
	By Lemma~\ref{lem:local-embedding-H_sigma},
	\[
	\|w^{-1/2}f\|_{L^2}
	+\|w^{-3/2}\nabla_vf\|_{L^2}
	\le C\|f\|_{L^2H_\sigma^1},
	\qquad w(v)=\langle v\rangle.
	\]
	Since $F=\mu+\sqrt\mu\,f$ and
	$\nabla_vF=\nabla_v\mu+\sqrt\mu\,\nabla_vf
	+f\nabla_v\sqrt\mu$,
	Gaussian decay absorbs any fixed polynomial weight.  Thus, for
	every $m\ge0$,
	\begin{equation}\label{eq:inflow-pos-weighted-l2}
		\|w^mF\|_{L^2}+\|w^m\nabla_vF\|_{L^2}
		\le C_{m,T}\bigl(1+\|f\|_{L^2H_\sigma^1}\bigr),
	\end{equation}
	with the dependence on $T$ omitted in the stationary case.
	Moreover, the perturbative pointwise bound yields
	$|F|\le C\mu^{1/4}$.  Since $|\phi(z)|\le C|z|^{-1}$ and
	$|\nabla\phi(z)|\le C|z|^{-2}$, both locally integrable
	in dimension three, velocity convolution with this Gaussian
	dominant gives
	\[
	\|a\|_{L^\infty}+\|b\|_{L^\infty}\le C.
	\]
	Together with $J=a\nabla_vF-bF$ and
	\eqref{eq:inflow-pos-weighted-l2}, this yields
	\begin{equation}\label{eq:inflow-pos-flux-l2}
		w^mF,\quad w^m\nabla_vF,\quad w^mJ,\quad w^mbF\in L^2
		\qquad (m\ge0).
	\end{equation}
	
	Fix $\chi_R,\zeta_R\in C_c^\infty(\bR_v^3)$, with
	$0\le\chi_R,\zeta_R\le1$, such that
	\[
	\begin{gathered}
		\chi_R=1\ \text{on }B_R,\qquad
		\operatorname{supp}\chi_R\subset B_{2R},\qquad
		|\nabla_v\chi_R|\le C/R,\\
		\zeta_R=1\ \text{on }B_{3R},\qquad
		\operatorname{supp}\zeta_R\subset B_{4R}.
	\end{gathered}
	\]
	Lemma~\ref{lem:compact-velocity-support} yields
	$u_R:=\zeta_RF\in\cH_{\kin}$; moreover
	\[
	\mathcal Tu_R=g_R^0+\partial_{v_i}g_R^i,
	\qquad
	g_R^0=-J\cdot\nabla_v\zeta_R,
	\quad g_R^i=\zeta_RJ^i.
	\]
	Both $g_R^0$ and $g_R^i$ are in $L^2$ by
	\eqref{eq:inflow-pos-flux-l2}.  On $\operatorname{supp}\chi_R$,
	we have $u_R=F$, $g_R^0=0$, and $g_R^i=J^i$.
	Therefore Proposition~\ref{prop:renormalized-identity}, with
	$u_R$ as the graph function, renormalization $\beta$, and test
	function $\chi_R$, applies to either problem.
	
	The bulk integrand in the resulting renormalized identity is
	\[
	-\beta''(F)\nabla_vF\cdot J
	=-\beta''(F)a^{ij}\partial_{v_i}F\partial_{v_j}F
	+b\cdot\nabla_vH(F).
	\]
	The chain rule is valid because $H'\in L^\infty$, $F$ is locally
	$H_v^1$, and $H(F)$ vanishes at $F=0$.  Formula
	\eqref{eq:inflow-pos-divb} and approximation of the compactly
	supported $H(F)\chi_R$ give
	\[
	\int b\cdot\nabla_vH(F)\,\chi_R
	=8\pi\int FH(F)\chi_R
	-\int bH(F)\cdot\nabla_v\chi_R.
	\]
	The other velocity-cutoff term is
	$-\int\beta'(F)J\cdot\nabla_v\chi_R$.  In view of
	$|\beta'|\le1$ and $0\le H(F)\le F_-$,
	the total error $\mathcal E_R$, the sum of these two velocity-cutoff terms, satisfies
	\[
	|\mathcal E_R|
	\le\frac{C}{R}\int_{\{R<|v|<2R\}}
	\bigl(|J|+|bF|\bigr)
	\le CR^{-1/2}
	\bigl(\|wJ\|_{L^2}+\|wbF\|_{L^2}\bigr)
	\longrightarrow0,
	\]
	where the integrals include $(t,x)$ in the evolutionary case.
	The diffusion term has the integrable majorant
	\[
	0\le\beta''(F)a^{ij}\partial_{v_i}F\partial_{v_j}F
	\le C\sigma^{ij}\partial_{v_i}F\partial_{v_j}F\in L^1
	\]
	by \eqref{eq:inflow-pos-ellipticity} and the Gaussian multiplier
	bounds.  Likewise $|FH(F)|\le(\pi/2)F_-\in L^1$, since
	$F_-\le\sqrt\mu\,|f|$.
	Finally $\beta(F_\gamma)\le(F_\gamma)_-
	\le\sqrt\mu\,|f_\gamma|$ is integrable by the full physical
	trace theorem and the Gaussian factor.
	These bounds justify letting $R\to\infty$ in all the terms
	of the renormalized identities.
	
	\smallskip
	\noindent\emph{Stationary solution.}
	The limiting stationary identity is
	\[
	\begin{aligned}
		&\int_{\gamma^+}\beta((F_{\st})_\gamma)\,\dd\gamma
		+\int_{\cD\times\bR_v^3}
		\beta''(F_{\st})a^{ij}\partial_{v_i}F_{\st}\partial_{v_j}F_{\st}
		+8\pi\int_{\cD\times\bR_v^3}(-F_{\st})H(F_{\st})\\
		&\hspace{3em}=\int_{\gamma^-}\beta(H_{\st})\,\dd\gamma.
	\end{aligned}
	\]
	If $H_{\st}\ge0$, the right-hand side vanishes.  By
	\eqref{eq:inflow-pos-ellipticity} and
	\eqref{eq:inflow-pos-beta-sign}, each term on the left is
	nonnegative.  Hence
	$(-F_{\st})H(F_{\st})=0$ almost everywhere, which by the
	strict inequality in \eqref{eq:inflow-pos-beta-sign} implies
	$F_{\st}\ge0$ almost everywhere.
	
	\smallskip
	\noindent\emph{Evolutionary solution.}
	For every $T>0$ and $0\le t\le T$, the limiting evolutionary
	identity with time endpoints $0,t$ yields
	\[
	\begin{aligned}
		&\int_{\cD\times\bR_v^3}\beta(F(t))
		+\int_0^t\int_{\gamma^+}\beta(F_\gamma)\,\dd\gamma\,\dd s
		+\int_0^t\int_{\cD\times\bR_v^3}
		\beta''(F)a^{ij}\partial_{v_i}F\partial_{v_j}F\\
		&\quad+8\pi\int_0^t\int_{\cD\times\bR_v^3}(-F)H(F)
		=\int_{\cD\times\bR_v^3}\beta(F_0)
		+\int_0^t\int_{\gamma^-}\beta(H_{\pt})\,\dd\gamma\,\dd s.
	\end{aligned}
	\]
	Under $F_0\ge0$ and $H_{\pt}\ge0$, the right-hand side is
	zero and every term on the left is nonnegative. Thus
	$\int\beta(F(t))=0$, and hence $F(t)\ge0$ almost everywhere.
	The identity is valid for each time endpoint by
	Proposition~\ref{prop:renormalized-identity}.  Alternatively,
	nonnegativity for almost every $t$ extends to every $t$ through
	the canonical $C([0,T];L^2_{x,v})$ representative of $f_{\pt}$:
	\[
	\|F(t)-F(s)\|_{L^1_{x,v}}
	\le |\cD|^{1/2}\|\sqrt\mu\|_{L^2_v}
	\|f_{\pt}(t)-f_{\pt}(s)\|_{L^2_{x,v}}.
	\]
	Since $T>0$ was arbitrary, the conclusion holds for all $t\ge0$.
\end{proof}

\subsection{Exponential convergence to the steady state}
\label{sec:asymptotic-behavior}

Throughout this subsection, \(f_{\st}\), \(f_{\pt}\), and
\(f_{\dy}\) have the meaning fixed in
\eqref{eq:stationary-perturbation-decomposition-main}--
\eqref{eq:dynamic-perturbation-decomposition-main} and in the incoming boundary relations and difference equation
in Section~\ref{sec:landau-wellposedness}.  In particular,
\[
f_{\dy}=f_{\pt}-f_{\st},
\qquad
h_{\dy}=h_{\pt}-h_{\st},
\]
and, by \eqref{eq:dynamic-perturbation-equation},
\begin{equation*}
	\begin{cases}
		(\partial_t+v\cdot\nabla_x)f_{\dy}
		+Lf_{\dy}
		=N[f_{\pt},f_{\dy}]
		+N[f_{\dy},f_{\st}]
		&\text{in }\bR_+\times\cD\times\bR_v^3,\\
		f_{\dy}(0)=f_0-f_{\st}
		&\text{in }\cD\times\bR_v^3,\\
		f_{\dy}|_{\bR_+\times\gamma^-}=h_{\dy}
		&\text{on }\bR_+\times\gamma^-.
	\end{cases}
\end{equation*}
Both the statement and the proof below use the kinetic transport operator
and the dissipation functional
\begin{equation}\label{eq:asymptotic-dissipation-functional}
	\mathfrak T:=\partial_t+v\cdot\nabla_x,
	\qquad
	\begin{aligned}
		\cD_f(s,t)
		:={}&\int_s^t\left(
		\|f(r)\|_{L^2(\cD;H_\sigma^1)}^2
		+\|\mathfrak Tf(r)\|_{L^2(\cD;H_\sigma^{-1})}^2
		\right)\dd r
		\\
		&+\|f_{\gamma^+}\|_{L^2((s,t)\times\gamma^+)}^2,
	\end{aligned}
\end{equation}
for \(0\le s\le t\le\infty\).  The following statement is formulated only
in terms of the actual initial and inflow data.  No interior lifting is
included among the hypotheses.

\begin{thm}[Exponential convergence to the steady state]
	
	Let \(A_{\dy}=A_{\dy}(\cD)\) and \(\delta_{\dy}=\delta_{\dy}(\cD)\) be
	given by Theorem~\ref{thm:global-dynamic-landau}, and let
	\(\lambda_0>0\).  There are
	\[
	\delta_{\rm as}=\delta_{\rm as}(\cD)\in(0,\delta_{\dy}],
	\qquad
	\lambda=\lambda(\cD,\lambda_0)\in(0,\lambda_0]
	\]
	with the following property.
	
	Let \(f_0\), \(h_{\pt}\), and \(h_{\st}\) be actual initial,
	evolutionary inflow, and stationary inflow data satisfying the
	compatibility condition
	\begin{equation}\label{eq:asymptotic-corner-compatibility}
		h_{\pt}(0)=\operatorname{Tr}_{\gamma^-}f_0
	\end{equation}
	and the smallness condition
	\begin{equation}\label{eq:asymptotic-weighted-data-assumption}
		\mathsf D_{A_{\dy}}\le\delta_{\rm as},
	\end{equation}
	where \(\mathsf D_m\) is the data size
	\eqref{eq:dynamic-nonlinear-data-quantity-main}.  Assume further that
	the exponentially weighted inflow-tail quantity $K_0$
	of \eqref{eq:asymptotic-inflow-decay-main} is finite,
	which holds, for instance, whenever
	\(\sup_{t\ge0}
	\ee^{\lambda_0t}\|h_{\pt}(t)-h_{\st}\|_{L^2(\gamma^-)}^2\)
	is finite.
	
	Let \(f_{\st}\) and \(f_{\pt}\) be the solutions supplied by
	Theorems~\ref{Landau_st_exi} and \ref{thm:global-dynamic-landau}, and
	set \(f_{\dy}:=f_{\pt}-f_{\st}\).  Then, with a constant
	\(C_\cD>0\) depending only on \(\cD\),
	\begin{equation}\label{eq:asymptotic-stability-estimate}
		\|f_{\dy}(t)\|_{L^2(\cD\times\bR_v^3)}^2
		+\cD_{f_{\dy}}(t,\infty)
		\le
		C_\cD\bigl(\mathsf D_{A_{\dy}}^2+K_0\bigr)\ee^{-\lambda t},
		\qquad t\ge0,
	\end{equation}
	and, for every \(p\in[2,\infty)\),
	\begin{equation}\label{eq:asymptotic-Lp-stability}
		\|f_{\dy}(t)\|_{L^p(\cD\times\bR_v^3)}
		\le
		C_\cD\bigl(\mathsf D_{A_{\dy}}^2+K_0\bigr)^{1/2}
		\ee^{-\lambda t/p},
		\qquad t\ge0 .
	\end{equation}
	In terms of the physical distributions
	\eqref{eq:stationary-perturbation-decomposition-main} and
\eqref{eq:dynamic-perturbation-decomposition-main}, the latter estimate
	reads
	\begin{equation}\label{eq:physical-asymptotic-stability}
		\|F(t)-F_{\st}\|_{L^p(\cD\times\bR_v^3)}
		\le
		C_\cD\bigl(\mathsf D_{A_{\dy}}^2+K_0\bigr)^{1/2}
		\ee^{-\lambda t/p},
		\qquad 2\le p<\infty,\ t\ge0 .
	\end{equation}
\end{thm}

\begin{lemma}[Exponential decay with a self-referential source]
	\label{lem:delayed-dissipation-decay}
	Let \(c_0>0\) and \(\lambda_0>0\), and put
	\begin{equation}\label{260904647}
		\lambda
		:=\min\left\{
		\log\left(1+\frac{c_0}2\right),\,\lambda_0
		\right\}.
	\end{equation}
	There is \(C=C(c_0)>0\) with the following property.  Let
	\(K_1,\delta\ge0\), and let \(\cE:[0,\infty)\to[0,\infty)\) be
	measurable and satisfy
	\begin{equation}\label{eq:delayed-dissipation-energy}
		\cE(t)+c_0\int_s^t\cE(r)\dd r
		\le
		\cE(s)
		+K_1\ee^{-\lambda_0s}
		+\delta\ee^{-\lambda s},
		\qquad 0\le s\le t.
	\end{equation}
	Then
	\begin{equation}\label{eq:delayed-dissipation-decay}
		\cE(t)\le C\bigl(\cE(0)+K_1+\delta\bigr)\ee^{-\lambda t},
		\qquad t\ge0.
	\end{equation}
\end{lemma}

\begin{proof}
	Write \(b:=K_1+\delta\).  Since \(\lambda\le\lambda_0\), dropping the
	integral in \eqref{eq:delayed-dissipation-energy} gives
	\begin{equation}\label{eq:delayed-dissipation-pointwise}
		\cE(t)\le\cE(s)+b\ee^{-\lambda s},
		\qquad 0\le s\le t.
	\end{equation}
	
	Fix \(n\ge1\).  Taking \(t=n\) in
	\eqref{eq:delayed-dissipation-pointwise} and using
	\(\ee^{-\lambda s}\le\ee^{-\lambda(n-1)}\) for \(s\ge n-1\) gives
	\[
	\cE(n)\le\cE(s)+b\ee^{-\lambda(n-1)},
	\qquad s\in[n-1,n].
	\]
	Integrating this over \(s\in[n-1,n]\), an interval of unit length,
	yields
	\[
	\cE(n)\le\int_{n-1}^n\cE(r)\dd r+b\ee^{-\lambda(n-1)}.
	\]
	Inserting this lower bound for the integral into
	\eqref{eq:delayed-dissipation-energy} with \((s,t)=(n-1,n)\) gives
	\[
	(1+c_0)\cE(n)
	\le\cE(n-1)+(1+c_0)b\ee^{-\lambda(n-1)}.
	\]
	Iterating this recursion from \(n\) down to \(0\), we obtain
	\[
	\cE(n)
	\le(1+c_0)^{-n}\cE(0)
	+(1+c_0)b\sum_{j=1}^{n}(1+c_0)^{-j}\ee^{-\lambda(n-j)}.
	\]
	By \eqref{260904647} we have \(\ee^\lambda\le1+c_0/2\), hence
	\((1+c_0)^{-1}\le\ee^{-\lambda}\) and
	\[
	\sum_{j=1}^{n}(1+c_0)^{-j}\ee^{-\lambda(n-j)}
	=\ee^{-\lambda n}
	\sum_{j=1}^{n}\left(\frac{\ee^\lambda}{1+c_0}\right)^{j}
	\le\ee^{-\lambda n}
	\sum_{j=1}^{\infty}\left(\frac{1+c_0/2}{1+c_0}\right)^{j}
	=\left(1+\frac2{c_0}\right)\ee^{-\lambda n}.
	\]
	Therefore
	\begin{equation}\label{eq:delayed-dissipation-integers}
		\cE(n)\le C(c_0)\bigl(\cE(0)+b\bigr)\ee^{-\lambda n},
		\qquad n\in\bN_0.
	\end{equation}
	
	For general \(t\ge0\), set \(n:=\lfloor t\rfloor\).  Then
	\eqref{eq:delayed-dissipation-pointwise} with \(s=n\),
	\eqref{eq:delayed-dissipation-integers}, and
	\(\ee^{-\lambda n}\le\ee^{\lambda}\ee^{-\lambda t}\) give
	\[
	\cE(t)\le\cE(n)+b\ee^{-\lambda n}
	\le C(c_0)\bigl(\cE(0)+b\bigr)\ee^{-\lambda t},
	\]
	which is \eqref{eq:delayed-dissipation-decay}.
\end{proof}

\begin{proof}[\textbf{Proof of Theorem~\ref{thm:asymptotic-stability}}]
	Since \(\delta_{\rm as}\le\delta_{\dy}\), Theorem~\ref{thm:global-dynamic-landau}
	applies.  Let \((f^{(k)})_{k\ge0}\) be the iteration
	\eqref{eq:nonlinear-dynamic-iteration} constructed in its proof, so that
	\(f^{(0)}=f_{\st}\) and, by
	\eqref{eq:nonlinear-dynamic-strong-convergence},
	\begin{equation}\label{eq:asymptotic-iteration-convergence}
		\sup_{T>0}
		\|f^{(k)}-f_{\pt}\|_{\cH_{\kin,\sigma}(\cD_T\times\bR_v^3)}
		\longrightarrow0,
		\qquad k\to\infty.
	\end{equation}
	Put
	\[
	f_{\dy}^{(k)}:=f^{(k)}-f_{\st}.
	\]
	Subtracting the stationary equation \eqref{2605041257} from
	\eqref{eq:nonlinear-dynamic-iteration} and using the bilinearity of
	\(N\), we find that \(f_{\dy}^{(k+1)}\) solves
	\begin{equation}
		\begin{cases}
			\mathfrak Tf_{\dy}^{(k+1)}+Lf_{\dy}^{(k+1)}=S^{(k+1)},
			\qquad
			S^{(k+1)}:=N[f^{(k)},f_{\dy}^{(k+1)}]+N[f_{\dy}^{(k)},f_{\st}],\\
			f_{\dy}^{(k+1)}(0)=f_0-f_{\st},\\
			f_{\dy}^{(k+1)}|_{\bR_+\times\gamma^-}=h_{\dy}.
		\end{cases}
		\label{260904713}
	\end{equation}
	Throughout, \(\delta_{\rm as}\le\delta_{\dy}\) is decreased finitely many
	times, always by an amount depending only on \(\cD\).  We abbreviate
	\[
	\cD_k(s,t):=\cD_{f_{\dy}^{(k)}}(s,t),
	\qquad
	\mathsf E_0:=\|f_0-f_{\st}\|_{L^2(\cD\times\bR_v^3)}^2 ,
	\]
	the latter being the common initial energy of every \(f^{(k+1)}_{\dy}\).
	By \eqref{eq:dynamic-nonlinear-data-quantity-main}, \eqref{eq:main-inflow-st-regularity} with
	\(\ell=2\), and \(w^{-2}\in L^2(\bR_v^3)\),
	\begin{equation}\label{eq:asymptotic-initial-energy-bound}
		\mathsf E_0\le C_\cD\mathsf D_{A_{\dy}}^2 .
	\end{equation}
	Finally, \eqref{eq:asymptotic-inflow-decay-main} reads
	\begin{equation}\label{eq:asymptotic-inflow-tail}
		\|h_{\dy}\|_{L^2((s,t)\times\gamma^-)}^2
		\le
		K_0\ee^{-\lambda_0s},
		\qquad 0\le s\le t\le\infty .
	\end{equation}

	\medskip
	\noindent\emph{Step 1. A dissipation inequality for the iterates.}
	Fix \(k\ge0\) and \(0\le s\le t<\infty\).  Testing
	\eqref{260904713} against \(f_{\dy}^{(k+1)}\), using the Green identity
	of Proposition~\ref{2607151218} and the coercivity
	\eqref{eq:L-coercivity} of \(L\), gives, with a universal \(c>0\),
	\begin{align}\label{eq:asymptotic-microscopic-energy}
		&\frac12\|f_{\dy}^{(k+1)}(t)\|_{L^2(\cD\times\bR_v^3)}^2
		+c\int_s^t
		\|(\mathbf I-\mathbf P)f_{\dy}^{(k+1)}(r)\|_{L^2(\cD;H_\sigma^1)}^2
		\dd r
		+\frac12\|(f_{\dy}^{(k+1)})_{\gamma^+}\|_{L^2((s,t)\times\gamma^+)}^2
		\notag\\
		&\qquad\le
		\frac12\|f_{\dy}^{(k+1)}(s)\|_{L^2(\cD\times\bR_v^3)}^2
		+\frac12\|h_{\dy}\|_{L^2((s,t)\times\gamma^-)}^2
		+\int_s^t
		\langle S^{(k+1)}(r),f_{\dy}^{(k+1)}(r)\rangle_
		{H_\sigma^{-1},H_\sigma^1}\dd r.
	\end{align}
	On the other hand, Theorem~\ref{thm:macro-estimate-landau}\,(1), applied
	to \(f_{\dy}^{(k+1)}\) on \((0,T)\) for every \(T>0\), produces
	\(\eta^{\mathbf P}_{k+1}\in C([0,\infty))\): the function it produces at
	time \(r\) is a fixed bounded quadratic form evaluated at
	\(f_{\dy}^{(k+1)}(r)\) (see the passage to general \(f\) in its proof),
	so the functions obtained for different \(T\) patch together.  By
	\eqref{eq:macro-eta-bound},
	\begin{equation}\label{eq:asymptotic-eta-bound}
		|\eta^{\mathbf P}_{k+1}(r)|
		\le
		C_\cD\|f_{\dy}^{(k+1)}(r)\|_{L^2(\cD\times\bR_v^3)}^2,
		\qquad r\ge0,
	\end{equation}
	and \eqref{eq:appendix-macro-estimate-endpoint} gives
	\begin{align}\label{eq:asymptotic-macroscopic-energy}
		&\int_s^t\|\mathbf Pf_{\dy}^{(k+1)}(r)\|_{L^2(\cD\times\bR_v^3)}^2\dd r
		\notag\\
		&\qquad\le
		C_\cD\Bigl(
		\int_s^t
		\|(\mathbf I-\mathbf P)f_{\dy}^{(k+1)}(r)\|_{L^2(\cD;H_\sigma^1)}^2
		\dd r
		+\|S^{(k+1)}\|_{L^2((s,t)\times\cD;H_\sigma^{-1})}^2
		\notag\\
		&\hspace{22mm}
		+\|(f_{\dy}^{(k+1)})_{\gamma^+}\|_{L^2((s,t)\times\gamma^+)}^2
		+\|h_{\dy}\|_{L^2((s,t)\times\gamma^-)}^2
		\Bigr)
		+\eta^{\mathbf P}_{k+1}(t)-\eta^{\mathbf P}_{k+1}(s).
	\end{align}
	Here \eqref{eq:appendix-source-M-time-Hminus1} was used to bound the
	moment norm \(\|S^{(k+1)}\|_{\cM(s,t)}\) by
	\(C\|S^{(k+1)}\|_{L^2((s,t)\times\cD;H_\sigma^{-1})}\), and
	\eqref{eq:micro-Maxwellian-control} to bound
	\(\|\mu^{1/4}(\mathbf I-\mathbf P)f_{\dy}^{(k+1)}\|
	_{L^2((s,t)\times\cD\times\bR_v^3)}\) by
	\(C\|(\mathbf I-\mathbf P)f_{\dy}^{(k+1)}\|
	_{L^2((s,t)\times\cD;H_\sigma^1)}\).  Moreover
	\(\gamma=\gamma^+\cup\gamma^-\) up to a set of measure zero and
	\(f_{\dy}^{(k+1)}|_{\gamma^-}=h_{\dy}\), whence
	\begin{equation}\label{eq:asymptotic-trace-split}
		\|(f_{\dy}^{(k+1)})_\gamma\|_{L^2((s,t)\times\gamma)}^2
		=
		\|(f_{\dy}^{(k+1)})_{\gamma^+}\|_{L^2((s,t)\times\gamma^+)}^2
		+\|h_{\dy}\|_{L^2((s,t)\times\gamma^-)}^2 .
	\end{equation}
	
	We now combine the two estimates.  Let \(C_\cD\) denote the larger of
	the constants in \eqref{eq:asymptotic-eta-bound} and
	\eqref{eq:asymptotic-macroscopic-energy}, and fix
	\(\kappa=\kappa(\cD)>0\) so small that
	\begin{equation}\label{eq:asymptotic-kappa-choice}
		\kappa C_\cD\le\min\left\{\frac c2,\frac14\right\} .
	\end{equation}
	Set
	\[
	\cE^0_{k+1}(r)
	:=\frac12\|f_{\dy}^{(k+1)}(r)\|_{L^2(\cD\times\bR_v^3)}^2
	-\kappa\eta^{\mathbf P}_{k+1}(r),
	\]
	so that \eqref{eq:asymptotic-eta-bound} and
	\eqref{eq:asymptotic-kappa-choice} give the energy equivalence
	\begin{equation}\label{eq:asymptotic-energy-equivalence}
		\frac14\|f_{\dy}^{(k+1)}(r)\|_{L^2(\cD\times\bR_v^3)}^2
		\le
		\cE^0_{k+1}(r)
		\le
		\frac34\|f_{\dy}^{(k+1)}(r)\|_{L^2(\cD\times\bR_v^3)}^2,
		\qquad r\ge0 .
	\end{equation}
	Add \(\kappa\) times \eqref{eq:asymptotic-macroscopic-energy} to
	\eqref{eq:asymptotic-microscopic-energy}.  The endpoint terms
	\(\kappa(\eta^{\mathbf P}_{k+1}(t)-\eta^{\mathbf P}_{k+1}(s))\) turn the
	two energies into \(\cE^0_{k+1}(t)\) and \(\cE^0_{k+1}(s)\).  On the
	right-hand side there appear \(\kappa C_\cD\) times the dissipation
	\(\int_s^t\|(\mathbf I-\mathbf P)f_{\dy}^{(k+1)}\|_
	{L^2(\cD;H_\sigma^1)}^2\) and \(\kappa C_\cD\) times the outgoing
	trace \(\|(f_{\dy}^{(k+1)})_{\gamma^+}\|_
	{L^2((s,t)\times\gamma^+)}^2\); by \eqref{eq:asymptotic-kappa-choice}
	both are absorbed by the left-hand side.  Hence
	\begin{align}\label{eq:asymptotic-energy-after-kappa}
		&\cE^0_{k+1}(t)
		+\frac c2\int_s^t
		\|(\mathbf I-\mathbf P)f_{\dy}^{(k+1)}(r)\|_{L^2(\cD;H_\sigma^1)}^2\dd r
		+\kappa\int_s^t
		\|\mathbf Pf_{\dy}^{(k+1)}(r)\|_{L^2(\cD\times\bR_v^3)}^2\dd r
		\notag\\
		&\qquad
		+\frac14\|(f_{\dy}^{(k+1)})_{\gamma^+}\|_{L^2((s,t)\times\gamma^+)}^2
		\notag\\
		&\quad\le
		\cE^0_{k+1}(s)
		+\|h_{\dy}\|_{L^2((s,t)\times\gamma^-)}^2
		+C_\cD\|S^{(k+1)}\|_{L^2((s,t)\times\cD;H_\sigma^{-1})}^2
		\notag\\
		&\qquad
		+\int_s^t
		\langle S^{(k+1)}(r),f_{\dy}^{(k+1)}(r)\rangle_
		{H_\sigma^{-1},H_\sigma^1}\dd r .
	\end{align}
	
	The two dissipation terms in
	\eqref{eq:asymptotic-energy-after-kappa} dominate the full norm.
	Indeed \(\|g\|_{L^2(\cD;H_\sigma^1)}^2
	\le2\|(\mathbf I-\mathbf P)g\|_{L^2(\cD;H_\sigma^1)}^2
	+2\|\mathbf Pg\|_{L^2(\cD;H_\sigma^1)}^2\), and the range of
	\(\mathbf P\) is finite dimensional, so that
	\eqref{eq:macro-Hsigma-equivalence} bounds the second summand by
	\(C\|\mathbf Pg\|_{L^2(\cD\times\bR_v^3)}^2\), where $C>0$ is a universal constant.  Consequently there
	is \(c'=c'(\cD)>0\) with
	\begin{align}\label{eq:asymptotic-full-Hsigma-lower}
		&\frac c2\int_s^t
		\|(\mathbf I-\mathbf P)f_{\dy}^{(k+1)}(r)\|_{L^2(\cD;H_\sigma^1)}^2\dd r
		+\kappa\int_s^t
		\|\mathbf Pf_{\dy}^{(k+1)}(r)\|_{L^2(\cD\times\bR_v^3)}^2\dd r
		\notag\\
		&\qquad\ge
		c'\int_s^t
		\|f_{\dy}^{(k+1)}(r)\|_{L^2(\cD;H_\sigma^1)}^2\dd r .
	\end{align}
	By Young's inequality,
	\begin{equation*}
		\int_s^t
		\langle S^{(k+1)},f_{\dy}^{(k+1)}\rangle_
		{H_\sigma^{-1},H_\sigma^1}\dd r
		\le
		\frac{c'}2\int_s^t
		\|f_{\dy}^{(k+1)}(r)\|_{L^2(\cD;H_\sigma^1)}^2\dd r
		+C_\cD
		\|S^{(k+1)}\|_{L^2((s,t)\times\cD;H_\sigma^{-1})}^2 ,
	\end{equation*}
	and the first term on the right is absorbed by
	\eqref{eq:asymptotic-full-Hsigma-lower}, leaving
	\(\frac{c'}2\int_s^t
	\|f_{\dy}^{(k+1)}\|_{L^2(\cD;H_\sigma^1)}^2\) on the left.
	
	It remains to produce the transport part of \(\cD_{k+1}\).  By
	\eqref{260904713} and \eqref{boundedness of L},
	\begin{equation*}
		\|\mathfrak Tf_{\dy}^{(k+1)}\|_{L^2((s,t)\times\cD;H_\sigma^{-1})}^2
		\le
		C_\cD\int_s^t
		\|f_{\dy}^{(k+1)}(r)\|_{L^2(\cD;H_\sigma^1)}^2\dd r
		+2\|S^{(k+1)}\|_{L^2((s,t)\times\cD;H_\sigma^{-1})}^2 .
	\end{equation*}
	Adding a sufficiently small multiple of this inequality to the left-hand
	side, at the expense of half of the surviving
	\(\frac{c'}2\int_s^t\|f_{\dy}^{(k+1)}\|_{L^2(\cD;H_\sigma^1)}^2\), we
	arrive at
	\begin{align}\label{eq:asymptotic-energy-before-source}
		&\cE^0_{k+1}(t)+c_1\cD_{k+1}(s,t)
		\notag\\
		&\qquad\le
		\cE^0_{k+1}(s)
		+\|h_{\dy}\|_{L^2((s,t)\times\gamma^-)}^2
		+C_\cD\|S^{(k+1)}\|_{L^2((s,t)\times\cD;H_\sigma^{-1})}^2,
	\end{align}
	for all \(0\le s\le t<\infty\), with some \(c_1=c_1(\cD)>0\).
	
	\medskip
	\noindent\emph{Step 2. The source term.}
	We claim that
	\begin{align}\label{eq:asymptotic-source-split}
		&\|S^{(k+1)}\|_{L^2((s,t)\times\cD;H_\sigma^{-1})}^2
		\notag\\
		&\qquad\le
		C_\cD\mathsf D_{A_{\dy}}^2
		\left(
		\cD_{k+1}(s,t)
		+\|f_{\dy}^{(k)}(s)\|_{L^2(\cD\times\bR_v^3)}^2
		+\cD_k(s,\infty)
		+K_0\ee^{-\lambda_0s}
		\right).
	\end{align}
	Indeed, by \eqref{boundedness of L} and
	\eqref{eq:nonlinear-dynamic-uniform-regularity} (for \(k=0\), by
	\eqref{eq:main-inflow-st-regularity} with \(\ell=0\)),
	\[
	\|N[f^{(k)},f_{\dy}^{(k+1)}]\|_{L^2((s,t)\times\cD;H_\sigma^{-1})}
	\le
	C\|\mu^{1/4}f^{(k)}\|_{L^\infty}
	\|f_{\dy}^{(k+1)}\|_{L^2((s,t)\times\cD;H_\sigma^1)}
	\le
	C_\cD\mathsf D_{A_{\dy}}\,\cD_{k+1}(s,t)^{1/2}.
	\]
	For the second half of \(S^{(k+1)}\) we use
	\eqref{eq:N-estimate-evolutionary} with \((\ell,p)=(0,30)\), together
	with
	\eqref{eq:main-inflow-st-regularity} with \(\ell=5\),
	\eqref{eq:main-inflow-st-gradient} with \((\ell,p)=(2,30)\), and
	\(A_{\dy}\ge A_{\st}+5\):
	\[
	\|N[f_{\dy}^{(k)},f_{\st}]\|_{L^2((s,t)\times\cD;H_\sigma^{-1})}
	\le
	C_\cD
	\|\mu^{1/8}f_{\dy}^{(k)}\|_{\cH_{\kin,\sigma}((s,t)\times\cD\times\bR_v^3)}
	\,\mathsf D_{A_{\dy}} .
	\]
	By \(\|\mu^{1/8}u\|_{\cH_{\kin,\sigma}}\le C\|u\|_{\cH_{\kin,\sigma}}\),
	the definition \eqref{innerprod_sigma} of the
	\(\cH_{\kin,\sigma}\)-norm on the slab \((s,t)\), and
	\eqref{eq:asymptotic-inflow-tail},
	\[
	\|\mu^{1/8}f_{\dy}^{(k)}\|_
	{\cH_{\kin,\sigma}((s,t)\times\cD\times\bR_v^3)}^2
	\le
	C\left(
	\|f_{\dy}^{(k)}(s)\|_{L^2(\cD\times\bR_v^3)}^2
	+\cD_k(s,\infty)
	+K_0\ee^{-\lambda_0s}
	\right),
	\]
	the last term accounting for the incoming trace.  Squaring the two
	displays, adding them and using \((a+b)^2\le2a^2+2b^2\) proves
	\eqref{eq:asymptotic-source-split}.
	
	Insert \eqref{eq:asymptotic-source-split} and
	\eqref{eq:asymptotic-inflow-tail} into
	\eqref{eq:asymptotic-energy-before-source}.  Since
	\(f_{\dy}^{(k+1)}\in\cH_{\kin,\sigma}(\cD_t\times\bR_v^3)\), the
	definition \eqref{innerprod_sigma} and \eqref{260730346} give
	\(\cD_{k+1}(s,t)<\infty\) for every \(t<\infty\).  Decreasing
	\(\delta_{\rm as}\) so that \(C_\cD\delta_{\rm as}^2\le c_1/2\)
	therefore lets the term
	\(C_\cD\mathsf D_{A_{\dy}}^2\cD_{k+1}(s,t)\) be absorbed by the
	left-hand side.  Using also \(\mathsf D_{A_{\dy}}\le1\), this gives
	\begin{equation}\label{eq:asymptotic-energy-absorbed}
		\cE^0_{k+1}(t)+\frac{c_1}2\cD_{k+1}(s,t)
		\le
		\cE^0_{k+1}(s)
		+C_\cD K_0\ee^{-\lambda_0s}
		+C_\cD\mathsf D_{A_{\dy}}^2
		\left(
		\|f_{\dy}^{(k)}(s)\|_{L^2(\cD\times\bR_v^3)}^2
		+\cD_k(s,\infty)
		\right),
	\end{equation}
	for all \(0\le s\le t<\infty\).
	
	\medskip
	\noindent\emph{Step 3. Recovery of the unweighted mass and the abstract
		lemma.}
	Let \(\bar\eta_{k+1}\) be the function supplied by Lemma~\ref{lem:transport-tail-L2}\,(1) for
	\(f_{\dy}^{(k+1)}\).  It is determined pointwise by
	\(f_{\dy}^{(k+1)}(r)\), so applying the lemma on \((0,T)\) for every
	\(T>0\) yields a single \(\bar\eta_{k+1}\in C([0,\infty))\).  By
	the pointwise bound in Lemma~\ref{lem:transport-tail-L2}\,(1) and
	\eqref{eq:transport-tail-L2-evolutionary} there is \(C_\cD>0\) with
	\begin{equation}\label{eq:asymptotic-transport-tail}
		|\bar\eta_{k+1}(r)|
		\le
		C_\cD\|f_{\dy}^{(k+1)}(r)\|_{L^2(\cD\times\bR_v^3)}^2,
		\qquad r\ge0,
	\end{equation}
	and, by \(\|u\|_{L^2}\le\|w^{1/2}u\|_{L^2}\) and
	\eqref{eq:asymptotic-trace-split},
	\begin{equation}\label{eq:asymptotic-transport-tail-slab}
		\|f_{\dy}^{(k+1)}\|_{L^2((s,t)\times\cD\times\bR_v^3)}^2
		\le
		C_\cD\left(
		\cD_{k+1}(s,t)+\|h_{\dy}\|_{L^2((s,t)\times\gamma^-)}^2
		\right)
		+\bar\eta_{k+1}(t)-\bar\eta_{k+1}(s),
	\end{equation}
	for all \(0\le s\le t<\infty\).
	Fix \(\bar\kappa=\bar\kappa(\cD)>0\) so small that
	\begin{equation}\label{eq:asymptotic-kappa-bar-choice}
		\bar\kappa C_\cD\le\min\left\{\frac{c_1}2,\frac18\right\},
	\end{equation}
	and set
	\[
	\cE_{k+1}(r)
	:=\cE^0_{k+1}(r)-\bar\kappa\bar\eta_{k+1}(r).
	\]
	By \eqref{eq:asymptotic-energy-equivalence},
	\eqref{eq:asymptotic-transport-tail}, and
	\eqref{eq:asymptotic-kappa-bar-choice},
	\begin{equation}\label{eq:asymptotic-full-energy-equivalence}
		\frac18\|f_{\dy}^{(k+1)}(r)\|_{L^2(\cD\times\bR_v^3)}^2
		\le
		\cE_{k+1}(r)
		\le
		\frac78\|f_{\dy}^{(k+1)}(r)\|_{L^2(\cD\times\bR_v^3)}^2 .
	\end{equation}
	Moreover, \eqref{eq:asymptotic-full-energy-equivalence},
	\eqref{eq:asymptotic-transport-tail-slab},
	\eqref{eq:asymptotic-kappa-bar-choice} and
	\eqref{eq:asymptotic-inflow-tail} give
	\[
	\bar\kappa\int_s^t\cE_{k+1}(r)\dd r
	\le
	\bar\kappa\|f_{\dy}^{(k+1)}\|_{L^2((s,t)\times\cD\times\bR_v^3)}^2
	\le
	\frac{c_1}2\cD_{k+1}(s,t)
	+C_\cD K_0\ee^{-\lambda_0s}
	+\bar\kappa\bigl(\bar\eta_{k+1}(t)-\bar\eta_{k+1}(s)\bigr).
	\]
	Add this to \eqref{eq:asymptotic-energy-absorbed} and cancel the common
	term \(\frac{c_1}2\cD_{k+1}(s,t)\).  Moving
	\(\bar\kappa\bar\eta_{k+1}(t)\) to the left-hand side and
	\(-\bar\kappa\bar\eta_{k+1}(s)\) into \(\cE^0_{k+1}(s)\), and recalling
	\(\cE_{k+1}=\cE^0_{k+1}-\bar\kappa\bar\eta_{k+1}\), we obtain
	\begin{equation}\label{eq:asymptotic-energy-with-transport}
		\begin{aligned}
			&\cE_{k+1}(t)+\bar\kappa\int_s^t\cE_{k+1}(r)\dd r
			\\
			\le\,&
			\cE_{k+1}(s)
			+C_\cD K_0\ee^{-\lambda_0s}
			+C_\cD\mathsf D_{A_{\dy}}^2
			\left(
			\|f_{\dy}^{(k)}(s)\|_{L^2(\cD\times\bR_v^3)}^2
			+\cD_k(s,\infty)
			\right),
		\end{aligned}
	\end{equation}
	for all \(0\le s\le t<\infty\).
	
	We may now fix the decay rate: let
	\begin{equation*}
		\lambda=\lambda(\cD,\lambda_0)
		:=\min\left\{
		\log\left(1+\frac{\bar\kappa}2\right),\,\lambda_0
		\right\}
	\end{equation*}
	be the exponent produced by Lemma~\ref{lem:delayed-dissipation-decay}
	with \(c_0=\bar\kappa\).  Note that \(\lambda\) depends only on \(\cD\)
	and \(\lambda_0\); in particular it is independent of \(k\).
	
	\medskip
	\noindent\emph{Step 4. Induction along the iteration.}
	Let \(M>0\), to be fixed at the end of this step, and consider the
	assertion
	\begin{equation}\label{eq:asymptotic-induction-hypothesis}
		\|f_{\dy}^{(k)}(t)\|_{L^2(\cD\times\bR_v^3)}^2
		+\cD_k(t,\infty)
		\le M\ee^{-\lambda t},
		\qquad t\ge0 .
	\end{equation}
	We show that \eqref{eq:asymptotic-induction-hypothesis} holds for every
	\(k\ge0\) once \(M\) is chosen suitably.  For \(k=0\) both terms
	vanish, since \(f^{(0)}_\dy=f^{(0)}-f_{\st}=0\).  Assume
	\eqref{eq:asymptotic-induction-hypothesis} for some \(k\ge0\).  Then the
	last bracket in \eqref{eq:asymptotic-energy-absorbed} and in
	\eqref{eq:asymptotic-energy-with-transport} is at most
	\(M\ee^{-\lambda s}\), so that, for all \(0\le s\le t<\infty\),
	\begin{align}
		&\cE^0_{k+1}(t)+\frac{c_1}2\cD_{k+1}(s,t)
		\le
		\cE^0_{k+1}(s)
		+C_\cD K_0\ee^{-\lambda_0s}
		+C_\cD\mathsf D_{A_{\dy}}^2M\ee^{-\lambda s},
		\label{eq:asymptotic-energy-final}
		\\
		&\cE_{k+1}(t)+\bar\kappa\int_s^t\cE_{k+1}(r)\dd r
		\le
		\cE_{k+1}(s)
		+C_\cD K_0\ee^{-\lambda_0s}
		+C_\cD\mathsf D_{A_{\dy}}^2M\ee^{-\lambda s}.
		\label{eq:asymptotic-energy-final-transport}
	\end{align}
	
	The function \(\cE_{k+1}\) is continuous on \([0,\infty)\), since
	\(f_{\dy}^{(k+1)}\in C([0,\infty);L^2(\cD\times\bR_v^3))\) by
	Proposition~\ref{2607151218} and
	\(\eta^{\mathbf P}_{k+1},\bar\eta_{k+1}\in C([0,\infty))\), and it is
	nonnegative by \eqref{eq:asymptotic-full-energy-equivalence}.  Hence
	\eqref{eq:asymptotic-energy-final-transport} is exactly the hypothesis
	\eqref{eq:delayed-dissipation-energy} of
	Lemma~\ref{lem:delayed-dissipation-decay} with \(c_0=\bar\kappa\),
	\(K_1=C_\cD K_0\), and
	\(\delta=C_\cD\mathsf D_{A_{\dy}}^2M\).  Since
	\(\cE_{k+1}(0)\le\frac78\mathsf E_0\) by
	\eqref{eq:asymptotic-full-energy-equivalence}, the lemma gives
	\begin{equation}\label{eq:asymptotic-energy-decay}
		\|f_{\dy}^{(k+1)}(t)\|_{L^2(\cD\times\bR_v^3)}^2
		\le8\cE_{k+1}(t)
		\le
		C\left(
		\mathsf E_0+K_0
		+\mathsf D_{A_{\dy}}^2M
		\right)\ee^{-\lambda t},
		\qquad t\ge0 .
	\end{equation}
	Next, apply \eqref{eq:asymptotic-energy-final} on the slab
	\((t,\tau)\), \(\tau>t\), and drop \(\cE^0_{k+1}(\tau)\ge0\):
	\[
	\frac{c_1}2\cD_{k+1}(t,\tau)
	\le
	\cE^0_{k+1}(t)
	+C_\cD K_0\ee^{-\lambda_0t}
	+C_\cD\mathsf D_{A_{\dy}}^2M\ee^{-\lambda t}.
	\]
	Letting \(\tau\to\infty\), and using
	\(\cE^0_{k+1}(t)\le\frac34\|f_{\dy}^{(k+1)}(t)\|_{L^2}^2\),
	\eqref{eq:asymptotic-energy-decay} and \(\lambda\le\lambda_0\), we
	obtain the same bound for \(\cD_{k+1}(t,\infty)\).  Adding the two,
	there is \(C_\cD\ge1\), independent of \(k\), of \(M\) and of
	\(\lambda_0\), with
	\begin{equation}\label{eq:asymptotic-induction-output}
		\|f_{\dy}^{(k+1)}(t)\|_{L^2(\cD\times\bR_v^3)}^2
		+\cD_{k+1}(t,\infty)
		\le
		C_\cD\left(
		\mathsf E_0+K_0+\mathsf D_{A_{\dy}}^2M
		\right)\ee^{-\lambda t}.
	\end{equation}
	With this constant, decrease \(\delta_{\rm as}\) once more so that
	\[C_\cD\delta_{\rm as}^2\le\frac12,\]
	and fix
	\[
	M:=2C_\cD\left(\mathsf E_0+K_0\right).
	\]
	Then
	\(C_\cD\mathsf D_{A_{\dy}}^2M\le M/2\), so the right-hand side of
	\eqref{eq:asymptotic-induction-output} is at most
	\(\frac M2\ee^{-\lambda t}+\frac M2\ee^{-\lambda t}=M\ee^{-\lambda t}\).
	This proves \eqref{eq:asymptotic-induction-hypothesis} for \(k+1\), and
	hence, by induction, for every \(k\ge0\).
	
	\medskip
	\noindent\emph{Step 5. Passage to the limit.}
	Fix \(T>0\).  By \eqref{eq:asymptotic-iteration-convergence},
	\eqref{260730346} and the definition \eqref{innerprod_sigma} of the
	\(\cH_{\kin,\sigma}\)-norm,
	\[
	f_{\dy}^{(k)}\longrightarrow f_{\dy}
	\quad\text{in }
	C\bigl([0,T];L^2(\cD\times\bR_v^3)\bigr)
	\ \text{and in }
	L^2\bigl((0,T);L^2(\cD;H_\sigma^1)\bigr),
	\]
	while \(\mathfrak Tf_{\dy}^{(k)}\to\mathfrak Tf_{\dy}\) in
	\(L^2\bigl((0,T);L^2(\cD;H_\sigma^{-1})\bigr)\) and
	\((f_{\dy}^{(k)})_\gamma\to(f_{\dy})_\gamma\) in
	\(L^2(\gamma_T)\).  Each of the three norms entering
	\eqref{eq:asymptotic-dissipation-functional} is therefore continuous
	along the sequence, so that, for \(0\le t<T\),
	\[
	\|f_{\dy}(t)\|_{L^2(\cD\times\bR_v^3)}^2
	+\cD_{f_{\dy}}(t,T)
	=\lim_{k\to\infty}
	\left(
	\|f_{\dy}^{(k)}(t)\|_{L^2(\cD\times\bR_v^3)}^2
	+\cD_k(t,T)
	\right)
	\le M\ee^{-\lambda t}
	\]
	by \eqref{eq:asymptotic-induction-hypothesis} and
	\(\cD_k(t,T)\le\cD_k(t,\infty)\).  Letting \(T\uparrow\infty\), monotone
	convergence gives
	\[
	\|f_{\dy}(t)\|_{L^2(\cD\times\bR_v^3)}^2
	+\cD_{f_{\dy}}(t,\infty)
	\le M\ee^{-\lambda t},
	\qquad t\ge0 .
	\]
	Since \(M\le C_\cD(\mathsf D_{A_{\dy}}^2+K_0)\) by
	\eqref{eq:asymptotic-initial-energy-bound}, this is
	\eqref{eq:asymptotic-stability-estimate}.
	
	\medskip
	\noindent\emph{Step 6. The \(L^p\) estimates.}
	Write \(\Lambda:=\mathsf D_{A_{\dy}}^2+K_0\), and
	enlarge \(C_\cD\) so that \(C_\cD\ge1\).  Since
	\(\mathsf D_{A_{\dy}}\le\Lambda^{1/2}\),
	\eqref{eq:main-inflow-dynamic-regularity} with \(\ell=0\) and \eqref{eq:main-inflow-st-regularity} with
	\(\ell=0\) give
	\begin{equation}\label{eq:asymptotic-Linfty-bound}
		\|f_{\dy}\|_{L^\infty(\bR_+\times\cD\times\bR_v^3)}
		\le
		\|f_{\pt}\|_{L^\infty(\bR_+\times\cD\times\bR_v^3)}
		+\|f_{\st}\|_{L^\infty(\cD\times\bR_v^3)}
		\le C_\cD\mathsf D_{A_{\dy}}
		\le C_\cD\Lambda^{1/2} .
	\end{equation}
	Let \(2\le p<\infty\) and \(t\ge0\).  From the pointwise bound
	\(|f_{\dy}(t)|^p\le
	\|f_{\dy}(t)\|_{L^\infty}^{p-2}|f_{\dy}(t)|^2\) we get
	\[
	\|f_{\dy}(t)\|_{L^p(\cD\times\bR_v^3)}
	\le
	\|f_{\dy}(t)\|_{L^\infty(\cD\times\bR_v^3)}^{1-2/p}
	\left(\|f_{\dy}(t)\|_{L^2(\cD\times\bR_v^3)}^2\right)^{1/p} ,
	\]
	so that \eqref{eq:asymptotic-Linfty-bound} and
	\eqref{eq:asymptotic-stability-estimate} give
	\[
	\|f_{\dy}(t)\|_{L^p(\cD\times\bR_v^3)}
	\le
	\left(C_\cD\Lambda^{1/2}\right)^{1-2/p}
	\left(C_\cD\Lambda\ee^{-\lambda t}\right)^{1/p}
	=
	C_\cD^{1-1/p}\Lambda^{1/2}\ee^{-\lambda t/p} .
	\]
	Since \(C_\cD\ge1\), this proves
	\eqref{eq:asymptotic-Lp-stability}.  Finally,
	\eqref{eq:stationary-perturbation-decomposition-main} and
\eqref{eq:dynamic-perturbation-decomposition-main} give
	\(F(t)-F_{\st}=\sqrt\mu\,f_{\dy}(t)\), and \(\sqrt\mu\le1\), so
	\eqref{eq:physical-asymptotic-stability} follows from
	\eqref{eq:asymptotic-Lp-stability}.
\end{proof}

\section{Diffuse Reflection: Geometric and Analytic Setup}\label{sec:diffuse-setting}

Let $\cD\subset\bR^3$ be the bounded connected $C^{K_\cD+2}$ domain
specified in \S\ref{subsec:main-diffuse}. We work with ordered complete
derivatives of order at most $\mathrm N\ge4$, where
\begin{equation}
	\label{eq:order-and-domain-regularity}
	K_\cD\ge2\mathrm N.
\end{equation}
Constants may depend on $\mathrm N$. Complete derivatives are interpreted
distributionally; every stated norm requires the corresponding Bochner-space
representative.

We use the reference Maxwellian $\mu(v)=\ee^{-|v|^2}$, of temperature $\tfrac12$, and the velocity weight
\begin{equation}
	\label{eq:diffuse-part-weight}
	w(v):=\la v\ra=(1+|v|^2)^{1/2}.
\end{equation}
For constant $\Theta>0$, we use $M_\Theta(v)$ from the first definition \eqref{new:eq:intro-wall-maxwellian} with $\Theta(t,x)\equiv\Theta$,
and set
\begin{equation}\label{eq:outgoing-diffuse-moment}
	\cA_0f_{\gamma^+}(x):=\int_{n_x\cdot v>0}f_{\gamma^+}(x,v)\sqrt{\mu(v)}\,n_x\cdot v\dd v,
	\qquad
	c_\mu:=\int_{n\cdot v>0}\mu(v)\,n\cdot v\dd v=\frac\pi2 ,
\end{equation}
the first for an outgoing trace $f_{\gamma^+}$ and $x\in\partial\cD$.
By its normalization,
$\int_{n\cdot v>0}M_\Theta(v)(n\cdot v)\,dv=1$ for every unit
vector $n$. Cauchy--Schwarz in the outgoing flux measure gives
\begin{equation}
	\label{eq:outgoing-moment-flux-bound}
	\|\cA_0f_{\gamma^+}\|_{L^2(\partial\cD)}^2
	\le c_\mu\|f_{\gamma^+}\|_{L^2(\gamma^+)}^2.
\end{equation}
The wall temperature satisfies the following conditions.
\begin{enumerate}
	\item\label{it:dyn-theta-amplitude}
	(\emph{The variation is $\delta$-small through $\mathrm N+1$ time derivatives.})
	\begin{equation}
		\label{eq:dyn-theta-amplitude}
		\Theta_w(t,x)=\tfrac12\bigl(1+\delta\vartheta(t,x)\bigr),
		\qquad
		\max_{0\le j\le\mathrm N+1}\ \sup_{t>0}\
		\|\partial_t^j\vartheta(t,\cdot)\|_{C^{K_\cD}(\partial\cD)}\le1 ,
	\end{equation}
	with $0<\delta\le\tfrac12$, so that $\tfrac14\le\Theta_w\le\tfrac34$.
	\item\label{it:dyn-theta-relaxation}
	(\emph{The wall relaxes exponentially.})  There are
	$\vartheta_\infty\in C^{K_\cD}(\partial\cD)$ and $\lambda_w>0$ with
	\begin{equation}
		\label{eq:dyn-theta-relaxation}
		\max_{0\le j\le\mathrm N+1}\
		\bigl\|\partial_t^j\bigl(\vartheta(t,\cdot)-\vartheta_\infty\bigr)
		\bigr\|_{C^{K_\cD}(\partial\cD)}
		\le\ee^{-\lambda_wt},
		\qquad t\ge0 .
	\end{equation}
	Letting $t\to\infty$ in \eqref{eq:dyn-theta-relaxation} and using
	\eqref{eq:dyn-theta-amplitude} gives
	$\|\vartheta_\infty\|_{C^{K_\cD}(\partial\cD)}\le1$.
\end{enumerate}

Set
\begin{equation}\label{eq:wall-response-profiles}
	\mathfrak m_\Theta:=\frac{M_\Theta}{\sqrt\mu}\ \ (\Theta>0),
	\qquad
	\mathfrak m_\circ:=\mathfrak m_{1/2}=\frac{\sqrt\mu}{c_\mu},
	\qquad
	d_{\delta}:=\mathfrak m_{\Theta_w}-\mathfrak m_\circ,
\end{equation}
the middle identity being the normalization $c_\mu M_{1/2}=\mu$ of
\eqref{new:eq:intro-wall-maxwellian}.  The wall operator and its datum are
\begin{equation}\label{eq:nonisothermal-diffuse-response}
	P_wf_{\gamma^+}:=\mathfrak m_{\Theta_w}\cA_0f_{\gamma^+},
	\qquad
	h_w:=c_\mu\mathfrak m_{\Theta_w}-\sqrt\mu=c_\mu d_\delta ,
\end{equation}
and the isothermal projections are
\begin{equation}\label{eq:isothermal-diffuse-projections}
	\cP_0^\pm f_{\gamma^+}:=\mathfrak m_\circ\cA_0f_{\gamma^+}
	\qquad\text{on }\gamma^\pm .
\end{equation}
The decomposition $\mathfrak m_{\Theta_w}=\mathfrak m_\circ+d_\delta$
gives
\begin{equation}\label{eq:dyn-theta-wall-split}
	P_wf_{\gamma^+}=\cP_0^-f_{\gamma^+}+r_\delta[f_{\gamma^+}],
	\qquad
	r_\delta[f_{\gamma^+}]:=d_\delta\cA_0f_{\gamma^+}.
\end{equation}
We use this decomposition of $P_w$ throughout. The limiting wall is given by
\begin{equation}
	\label{eq:limiting-wall-profiles}
	\Theta_w^\infty:=\tfrac12(1+\delta\vartheta_\infty),
	\quad
	d_\delta^\infty:=\mathfrak m_{\Theta_w^\infty}-\mathfrak m_\circ,
	\quad
	h_w^\infty:=c_\mu d_\delta^\infty,
	\quad
	P_w^\infty f_{\gamma^+}:=\mathfrak m_{\Theta_w^\infty}\cA_0f_{\gamma^+}.
\end{equation}
The stationary equation for the perturbation $f_\st$ is
\begin{equation}\label{eq:stationary-full-diffuse-problem}
	\left\{
	\begin{aligned}
		v\cdot\nabla_xf_\st+Lf_\st&=N[f_\st,f_\st]
		&&\text{in }\cD\times\bR^3,\\
		(f_\st)_{\gamma^-}&=P_w^\infty(f_\st)_{\gamma^+}+h_w^\infty
		&&\text{on }\gamma^-,
	\end{aligned}
	\right.
	\qquad
	\int_\cD\int_{\bR^3}f_\st(x,v)\sqrt\mu(v)\,\dd v\dd x=0 .
\end{equation}

On $\cH_{\kin,\sigma}$ we use the norm equivalence
\begin{align}\label{260902137}
	\begin{alignedat}{2}
		&\|f\|_{\cH_{\kin,\sigma}(\cD\times \bR^3)}\,&&\simeq
		\|f\|_{L^2(\cD\times \bR^3_v)}+\|f\|_{L^2(\cD;H^1_\sigma)}
		+\|v\cdot\nabla_xf\|_{L^2(\cD;H^{-1}_\sigma)}
		+\|f_{\gamma}\|_{L^2(\gamma)}
		\\
		&\|f\|_{\cH_{\kin,\sigma}(\cD_T\times \bR^3)}\,&&\simeq
		\|f\|_{L^2(\cD_T\times \bR^3_v)}+\|f\|_{L^2(\cD_T;H^1_\sigma)}+\|(\p_t+v\cdot\nabla_x)f\|_{L^2(\cD_T;H^{-1}_\sigma)}
		\\
		&&&\qquad +\|f_{\gamma}\|_{L^2(\gamma_T)}+\sup_{0\le t\le T}\|f(t)\|_{L^2(\cD\times\bR_v^3)}
	\end{alignedat}
\end{align}
where the constants in this equivalence are independent of $T$ (see Proposition~\ref{2607151218} and Lemma~\ref{lem:transport-tail-L2}).  We refer to \eqref{innerprod_sigma} as the graph norm and to estimates in this norm as graph estimates.

\subsection{The atlas and tangential fields}\label{ss:diffuse-fields}

For a vector field \(X\in C^1(\bR^3;\bR^3)\) we set
\begin{equation}
	\label{eq:complete-phase-space-derivative}
	\mathfrak f_X:=\bigl(X(x),(\nabla_xX(x))v\bigr),
	\qquad
	\mathfrak D_X
	:=\mathfrak f_X\cdot\nabla_{x,v}
	=X(x)\cdot\nabla_x
	+\bigl((\nabla_xX(x))v\bigr)\cdot\nabla_v
\end{equation}
on \(\bR^3\times\bR^3\), and by restriction on \(\overline\cD\times\bR^3\).
For tangential fields $X$ defined below, the flow preserves
$\overline\cD$ and exists for all time. Let
\(\Phi_s^X\) be the spatial flow of \(X\) and
\begin{equation}
	\label{eq:complete-phase-flow}
	\widehat\Phi_s^X(x,v):=\bigl(\Phi_s^X(x),D\Phi_s^X(x)v\bigr)
\end{equation}
its tangent lift, whose generator is \(\mathfrak D_X\) and whose phase
divergence is
\begin{equation}
	\label{eq:phase-divergence}
	\operatorname{div}_{x,v}\mathfrak f_X=2\operatorname{div}_xX ,
\end{equation}
since $\operatorname{div}_v((\nabla_xX)v)=\operatorname{div}_xX$. The \emph{complete pullback} of $q$ is
\begin{equation}
	\label{eq:bulk-flow-pullback}
	T_s^Xq:=q\circ\widehat\Phi_s^X .
\end{equation}

\begin{prop}[Tangential fields]
	\label{prop:complete-fields}
	Let \(n\) be the outward unit normal of \(\cD\).  There exist a collar radius \(r_\cD>0\) and a function
	\(\rho_\cD\in C^{K_\cD+2}(\bR^3)\) which agrees with the signed distance to
	\(\partial\cD\), negative in \(\cD\), on the collar
	\[
	\cN_\cD:=\{x\in\bR^3:\ |\rho_\cD(x)|<r_\cD\},
	\]
	on which the nearest-point projection \(\pi\) onto \(\partial\cD\) is well
	defined and \(C^{K_\cD+1}\), and on which the outward normal extends as
	\[
	\bar n:=\nabla\rho_\cD\ \ \text{on }\cN_\cD ,
	\qquad
	\bar n=n\circ\pi,
	\qquad
	\bar n|_{\partial\cD}=n .
	\]
	There exist further a finite index set \(\mathfrak A\) and a family
	\[
	\{\tau_\nu\}_{\nu\in\mathfrak A}
	\subset C^{K_\cD+1}(\bR^3;\bR^3),
	\qquad
	\operatorname{supp}\tau_\nu\subset\cN_\cD ,
	\]
	with the following properties, where
	\(\mathfrak D_\nu:=\mathfrak D_{\tau_\nu}\) and
	\(d_\nu:=2\operatorname{div}_x\tau_\nu\).
	
	\begin{enumerate}
		\item\label{it:cf-tangential}
		\(\tau_\nu\cdot\nabla\rho_\cD=0\) on \(\cN_\cD\) for every
		\(\nu\in\mathfrak A\), and
		\[
		\operatorname{span}\{\tau_\nu(x):\nu\in\mathfrak A\}
		=T_x\partial\cD
		\qquad\text{for every }x\in\partial\cD .
		\]
		
		\item\label{it:cf-grazing}
		\(\mathfrak D_\nu(\bar n\cdot v)=0\) on \(\cN_\cD\times\bR^3\), for every
		\(\nu\in\mathfrak A\).
		
		\item\label{it:cf-adjoint}
		\[
		\int_{\cD\times\bR^3}(\mathfrak D_\nu f)\psi\dd x\dd v
		=\int_{\cD\times\bR^3}f\,(-\mathfrak D_\nu-d_\nu)\psi\dd x\dd v
		\qquad
		\bigl(f,\psi\in C_c^\infty(\overline\cD\times\bR^3)\bigr).
		\]
		
		\item\label{it:cf-frame}
		The index set is \(\mathfrak A=\{(\kappa,a):\kappa\in\mathfrak K_{\p},\
		a\in\{1,2\}\}\) for a finite family \(\{U_\kappa\}_{\kappa\in\mathfrak
			K_{\p}}\) of boundary charts covering \(\partial\cD\), and
		\(\tau_{(\kappa,a)}\) is supported in the normal tube
		\(\cN_\kappa:=\pi^{-1}(\partial\cD\cap U_\kappa)\cap\cN_\cD\).  On
		\(\cN_\kappa\) there is a rotation field
		\(R_\kappa\in C^{K_\cD+1}(\cN_\kappa;SO(3))\) with
		\(R_\kappa(x)e_3=\bar n(x)\) and \(R_\kappa(x)=R_\kappa(\pi(x))\) such
		that, for every \(\nu\in\mathfrak A\), there is a matrix field
		\(M_{\kappa,\nu}\in C^{K_\cD}(\cN_\kappa;\bR^{3\times3})\) for which,
		writing \((R_\kappa^{*}f)(x,\xi):=f(x,R_\kappa(x)\xi)\), one has, for every
		\(f\in C^1(\cN_\kappa\times\bR^3)\), every \(x\in\cN_\kappa\) and every
		\(\xi\in\bR^3\),
		\[
		(\mathfrak D_\nu f)(x,R_\kappa\xi)
		=\tau_\nu\cdot\nabla_x(R_\kappa^{*}f)
		+(M_{\kappa,\nu}\xi)\cdot\nabla_\xi(R_\kappa^{*}f),
		\qquad
		e_3^T M_{\kappa,\nu}=0 .
		\]
	\end{enumerate}
\end{prop}

\begin{proof}
	
	Let \(d_\cD\) be the signed distance to \(\partial\cD\), negative in \(\cD\).
	By the regularity of signed distance to a $C^{K_\cD+2}$ hypersurface,
	choose $r_\cD>0$ such that $d_\cD\in C^{K_\cD+2}$ whenever
	$\operatorname{dist}(x,\partial\cD)<3r_\cD$
	\cite[Lemma~14.16]{GilbargTrudinger1983}.
	Choose a smooth odd nondecreasing \(\zeta\) with
	\[
	\zeta(s)=s\quad\text{if }|s|\le r_\cD,
	\qquad
	\zeta(s)=2r_\cD\operatorname{sgn}s\quad\text{if }|s|\ge2r_\cD,
	\]
	set \(\rho_\cD:=\zeta(d_\cD)\) on that neighborhood, and extend \(\rho_\cD\)
	by the constants \(-2r_\cD\) on the remaining interior region and \(2r_\cD\) on
	the remaining exterior region.  Since \(\zeta(d_\cD)\) is constant for
	\(|d_\cD|\ge2r_\cD\), the three definitions glue in \(C^{K_\cD+2}\), and
	\begin{equation}
		\label{eq:signed-distance}
		\rho_\cD\in C^{K_\cD+2}(\bR^3),
		\qquad
		\cD=\{\rho_\cD<0\},
		\qquad
		\partial\cD=\{\rho_\cD=0\},
		\qquad
		\nabla\rho_\cD\ne0\quad\text{on }\partial\cD.
	\end{equation}
	On the collar \(\cN_\cD=\{|\rho_\cD|<r_\cD\}\) we have \(\rho_\cD=d_\cD\); the
	nearest-point map \(\pi:\cN_\cD\to\partial\cD\) is well defined and
	\(C^{K_\cD+1}\), and for every \(x\in\cN_\cD\),
	\[
	x=\pi(x)+\rho_\cD(x)n(\pi(x)),
	\qquad
	\bar n(x):=\nabla\rho_\cD(x)=n(\pi(x)),
	\qquad
	|\nabla\rho_\cD(x)|=1.
	\]
	
	Let \(\mathfrak K_{\p}\) denote the fixed finite boundary-chart index set and
	\(e_1,e_2,e_3\) the standard basis of \(\bR^3\).  For each
	\(\kappa\in\mathfrak K_{\p}\) choose \(x_\kappa\in\partial\cD\),
	\(P_\kappa\in SO(3)\), a bounded open \(\mathsf O_\kappa'\subset\bR^2\)
	containing \(0\), a radius \(r_\kappa>0\), and
	\(\eta_\kappa\in C^{K_\cD+2}(\overline{\mathsf O_\kappa'})\) with
	\(\eta_\kappa(0)=0\) and \(\nabla\eta_\kappa(0)=0\), so that the boundary is a
	graph in the chart
	\[
	\mathsf O_\kappa
	:=\{(y',y_3):y'\in\mathsf O_\kappa',\
	|y_3-\eta_\kappa(y')|<r_\kappa\},
	\qquad
	U_\kappa:=x_\kappa+P_\kappa\mathsf O_\kappa,
	\]
	and so that
	\[
	\begin{aligned}
		\cD\cap U_\kappa
		&=\{x_\kappa+P_\kappa(y',y_3):
		y\in\mathsf O_\kappa,\ y_3>\eta_\kappa(y')\},\\
		\partial\cD\cap U_\kappa
		&=\{x_\kappa+P_\kappa(y',y_3):
		y\in\mathsf O_\kappa,\ y_3=\eta_\kappa(y')\}.
	\end{aligned}
	\]
	Such a finite atlas exists by compactness and the regularity of $\partial\cD$.  We flatten the graph by
	\begin{align}
		\label{eq:boundary-flattening-maps}
		\psi_\kappa\bigl(x_\kappa+P_\kappa(y',y_3)\bigr)
		&:=(y',y_3-\eta_\kappa(y')),\notag\\
		\phi_\kappa(y)
		&:=x_\kappa+P_\kappa
		\bigl(y',y_3+\eta_\kappa(y')\bigr),
	\end{align}
	so that
	\[
	V_\kappa:=\psi_\kappa(U_\kappa)
	=\mathsf O_\kappa'\times(-r_\kappa,r_\kappa),
	\qquad
	V_\kappa^+=\mathsf O_\kappa'\times(0,r_\kappa),
	\]
	and put \(B_{\kappa,y}:=D_y\phi_\kappa(y)\), so that, for \(y\in V_\kappa\),
	\begin{equation}
		\label{eq:flattening-matrix-formulas}
		B_{\kappa,y}
		=P_\kappa
		\begin{pmatrix}I_2&0\\\nabla\eta_\kappa(y')^T&1\end{pmatrix},
		\qquad
		\det B_{\kappa,y}=1 .
	\end{equation}
	In particular the flattening is volume preserving in \(x\).

	Complete the standing boundary charts to a finite interior/boundary atlas of
	\(\overline\cD\), with index set
	\[
	\mathfrak K=\mathfrak K_{\circ}\mathbin{\dot\cup}\mathfrak K_{\p},
	\]
	where \(\mathfrak K_{\circ}\) indexes the interior charts added now.  For
	\(\kappa\in\mathfrak K\) we write \(U_\kappa\) for the physical spatial chart,
	with the boundary charts as above.  For
	\(\kappa\in\mathfrak K_{\circ}\) choose a bounded open \(V_\kappa\subset\bR^3\),
	a point \(x_\kappa\in\cD\) and \(P_\kappa\in SO(3)\), and set
	\begin{equation}
		\label{eq:interior-affine-chart}
		U_\kappa:=x_\kappa+P_\kappa V_\kappa
		\ \ \text{with}\ \ \overline{U_\kappa}\subset\cD,
		\qquad
		\phi_\kappa(y):=x_\kappa+P_\kappa y,
		\qquad
		\psi_\kappa(x):=P_\kappa^T(x-x_\kappa) .
	\end{equation}
	For a boundary chart $\kappa\in\mathfrak K_{\p}$, the artificial
	boundary in flattened coordinates is
	\(\widehat\Sigma_\kappa:=\partial V_\kappa^+\setminus\{y_3=0\}\), while for
	\(\kappa\in\mathfrak K_\circ\) we put
	\(\widehat\Sigma_\kappa:=\partial V_\kappa\).  Each \(\psi_\kappa\) and \(\phi_\kappa\) extends with
	\(C^{K_\cD+2}\) regularity to the relevant spatial closure.  In physical
	coordinates the artificial boundary of a chart is
	\[
	\Sigma_\kappa:=\phi_\kappa(\widehat\Sigma_\kappa)\subset\partial U_\kappa .
	\]
	Finally, shrink each chart once: choose open
	\(\widetilde U_\kappa\subset U_\kappa\) with
	\begin{equation}
		\label{eq:fixed-atlas-core-cover}
		\overline\cD\subset\bigcup_{\kappa\in\mathfrak K}\widetilde U_\kappa,
		\qquad
		\overline{\widetilde U_\kappa\cap\overline\cD}\cap\Sigma_\kappa=\varnothing
		\quad(\kappa\in\mathfrak K_{\p}),
	\end{equation}
	and \(\overline{\widetilde U_\kappa}\subset U_\kappa\) for
	\(\kappa\in\mathfrak K_{\circ}\).
	
	Each boundary chart carries two coordinate tangential directions, which we index
	by
	\begin{equation*}
		\mathfrak A:=\{(\kappa,a):\kappa\in\mathfrak K_{\p},\ a\in\{1,2\}\},
		\qquad
		\sum_{\nu\in\mathfrak A}
		=\sum_{\kappa\in\mathfrak K_{\p}}\sum_{a=1}^2 .
	\end{equation*}
	For \(\kappa\in\mathfrak K_{\p}\) and \(a\in\{1,2\}\) the coordinate tangent
	fields on the physical boundary are
	\begin{equation}
		\label{eq:coordinate-tangent-fields}
		T_{\kappa,a}\bigl(\phi_\kappa(y',0)\bigr)
		:=\partial_{y_a}\phi_\kappa(y',0)
		=B_{\kappa,(y',0)}e_a ;
	\end{equation}
	for each \(x\in\partial\cD\cap U_\kappa\) they span \(T_x\partial\cD\).  Choose
	\(\chi_\kappa\in C_c^\infty(U_\kappa\cap\partial\cD)\) equal to one on
	\(\widetilde U_\kappa\cap\partial\cD\), which the shrinking of the charts
	permits, and \(\chi_\rho\in C_c^\infty((-r_\cD,r_\cD))\) equal to one near
	zero, with
	\(\operatorname{supp}\chi_\rho\) shrunk so that the closed tube
	\[
	\bigl\{x\in\cN_\cD:\ \pi(x)\in\operatorname{supp}\chi_\kappa,\
	\rho_\cD(x)\in\operatorname{supp}\chi_\rho\bigr\}
	\]
	has compact closure in \(U_\kappa\) for every \(\kappa\in\mathfrak K_{\p}\).
	The compact supports of the finitely many $\chi_\kappa$ permit this common choice of normal cutoff without changing $r_\cD$ or $\rho_\cD$.  For
	\(x\in\cN_\cD\) we set
	\begin{equation}
		\label{eq:signed-distance-tangent-extension}
		\tau_{\kappa,a}(x)
		:=\chi_\rho(\rho_\cD(x))\,
		(\chi_\kappa T_{\kappa,a})(\pi(x)),
	\end{equation}
	extended by zero outside the collar, and we write
	\(\tau_\nu:=\tau_{\kappa,a}\) for \(\nu=(\kappa,a)\in\mathfrak A\).  As $\bar n=n\circ\pi$, these fields are tangent to the level sets of $\rho_\cD$. They are supported in $U_\kappa$ and agree with $T_{\kappa,a}$ on $\widetilde U_\kappa\cap\partial\cD$. Away from the boundary, the extension need not coincide with $B_{\kappa,y}e_a$.  The $C^{K_\cD+1}$ regularity of $\pi$ and $T_{\kappa,a}$, together with the normal cutoff, gives
	\begin{equation}
		\label{eq:tangential-field-family}
		\{\tau_\nu\}_{\nu\in\mathfrak A}
		\subset C^{K_\cD+1}(\bR^3;\bR^3),
		\qquad
		\operatorname{supp}\tau_{(\kappa,a)}\subset\cN_\kappa
		:=\pi^{-1}(\partial\cD\cap U_\kappa)\cap\cN_\cD .
	\end{equation}
	Every boundary point lies in a shrunken boundary chart, since the interior charts are compactly contained in $\cD$. There $\chi_\kappa=\chi_\rho(0)=1$ and $\pi(x)=x$, so $\tau_{(\kappa,a)}(x)=T_{\kappa,a}(x)$ for $a=1,2$. Thus
	\begin{equation}
		\label{eq:tangential-field-properties}
		\tau_\nu\cdot\nabla\rho_\cD=0
		\quad\text{in the signed-distance collar},
		\qquad
		\operatorname{span}\{\tau_\nu(x):\nu\in\mathfrak A\}=T_x\partial\cD
		\quad(x\in\partial\cD).
	\end{equation}
	This is assertion~(\ref{it:cf-tangential}).
	
	Since \(\bar n=\nabla\rho_\cD\) on the collar, \(\nabla\bar n=\nabla^2\rho_\cD\)
	is symmetric; together with \(\tau_\nu\cdot\bar n=0\) this gives
	\begin{equation}
		\label{eq:complete-derivative-preserves-grazing}
		\mathfrak D_\nu(\bar n\cdot v)
		=((\nabla \bar n)\tau_\nu)\cdot v
		+((\nabla\tau_\nu)v)\cdot \bar n
		=v\cdot\nabla(\tau_\nu\cdot \bar n)=0.
	\end{equation}
	This proves (\ref{it:cf-grazing}).
	
	The phase-space divergence is $d_\nu=2\operatorname{div}\tau_\nu$.
	Integration by parts has no boundary term because $\tau_\nu\cdot n=0$;
	this proves (\ref{it:cf-adjoint}).

	We apply Gram--Schmidt to the ordered frame \((T_{\kappa,2},T_{\kappa,1},n)\)
	of \eqref{eq:coordinate-tangent-fields}, a basis at every point of
	\(\partial\cD\cap U_\kappa\).  Since \(\cD\cap U_\kappa\) lies above the
	graph of \(\eta_\kappa\), that frame is positively oriented, with
	\[
	\det(T_{\kappa,2},T_{\kappa,1},n)=(1+|\nabla\eta_\kappa|^2)^{1/2}>0.
	\]
	The normal
	\(n\) is orthogonal to \(\operatorname{span}(T_{\kappa,1},T_{\kappa,2})=T_x\partial\cD\),
	so the third column of the Gram--Schmidt output \(R_\kappa^{\p}\) is \(n\) itself.
	We extend it along the normal rays by \(R_\kappa(x):=R_\kappa^{\p}(\pi(x))\) on
	the normal tube \(\cN_\kappa=\pi^{-1}(\partial\cD\cap U_\kappa)\cap\cN_\cD\), which
	contains \(\operatorname{supp}\tau_\nu\) for \(\nu=(\kappa,a)\) by
	\eqref{eq:tangential-field-family}.  This produces
	\begin{equation}
		\label{eq:rotating-frame}
		R_\kappa\in C^{K_\cD+1}(\cN_\kappa;SO(3)),
		\qquad
		R_\kappa(x)e_3=\bar n(x),
		\qquad
		R_\kappa(x)=R_\kappa(\pi(x)) ,
	\end{equation}
	with the regularity holding because \(T_{\kappa,a},n\in C^{K_\cD+1}\) by
	\eqref{eq:flattening-matrix-formulas}, Gram--Schmidt is smooth in its arguments
	on nondegenerate frames, and \(\pi\in C^{K_\cD+1}\).
	For \(f\in C^1\) we define the rotated representative
	\begin{equation}
		\label{eq:rotating-frame-representative}
		(R_\kappa^{*}f)(x,\xi):=f(x,R_\kappa(x)\xi),
	\end{equation}
	where \(\xi\) denotes the rotated velocity variable.  The chain rule gives
	\(\nabla_\xi(R_\kappa^{*}f)=R_\kappa^T(\nabla_vf)(x,R_\kappa\xi)\) and
	\begin{align*}
		\tau_\nu\cdot\nabla_x(R_\kappa^{*}f)
		={}&(\tau_\nu\cdot\nabla_xf)(x,R_\kappa\xi)
		+((\tau_\nu\cdot\nabla_xR_\kappa)\xi)\cdot\nabla_vf(x,R_\kappa\xi),\\
		(C\xi)\cdot\nabla_vf(x,R_\kappa\xi)
		={}&(R_\kappa^TC\xi)
		\cdot\nabla_\xi (R_\kappa^{*}f)
		\qquad(C\in\bR^{3\times3}),
	\end{align*}
	the second identity because \(R_\kappa R_\kappa^T=I\).  Substituting the first
	into \eqref{eq:complete-phase-space-derivative} and applying the second twice, at
	\(C=(\nabla\tau_\nu)R_\kappa\) and at
	\(C=\tau_\nu\cdot\nabla_xR_\kappa\), we obtain
	\begin{equation}
		\label{eq:rotating-frame-complete-lift}
		(\mathfrak D_\nu f)(x,R_\kappa\xi)
		=\tau_\nu\cdot\nabla_x(R_\kappa^{*}f)
		+(M_{\kappa,\nu}\xi)\cdot\nabla_\xi (R_\kappa^{*}f),
		\qquad
		M_{\kappa,\nu}
		:=R_\kappa^T(\nabla\tau_\nu)R_\kappa
		-R_\kappa^T(\tau_\nu\cdot\nabla_xR_\kappa).
	\end{equation}
	By \eqref{eq:tangential-field-family} and \eqref{eq:rotating-frame},
	\(M_{\kappa,\nu}\in C^{K_\cD}(\cN_\kappa;\bR^{3\times3})\).  It remains to
	identify its normal row.  Differentiating in \(x\) the identity
	\(\tau_\nu\cdot\bar n=0\), valid on all of \(\cN_\cD\), and using that
	\(\nabla\bar n=\nabla^2\rho_\cD\) is symmetric, we obtain
	\[
	0=\nabla_x(\tau_\nu\cdot \bar n)^T
	=\bar n^T(\nabla\tau_\nu)
	+(\tau_\nu\cdot\nabla_x\bar n)^T ,
	\]
	while differentiating \(\bar n^T R_\kappa=e_3^T\), which is
	\eqref{eq:rotating-frame}, along \(\tau_\nu\) gives
	\[
	\bar n^T(\tau_\nu\cdot\nabla_xR_\kappa)
	=-(\tau_\nu\cdot\nabla_x\bar n)^T R_\kappa .
	\]
	Since \(e_3^T R_\kappa^T=\bar n^T\), the two identities give
	\[
	e_3^T M_{\kappa,\nu}
	=\bar n^T(\nabla\tau_\nu)R_\kappa
	-\bar n^T(\tau_\nu\cdot\nabla_xR_\kappa)
	=-(\tau_\nu\cdot\nabla_x\bar n)^T R_\kappa
	+(\tau_\nu\cdot\nabla_x\bar n)^T R_\kappa ,
	\]
	that is,
	\[
	e_3^T M_{\kappa,\nu}=0.
	\]
	Together with \eqref{eq:rotating-frame-complete-lift} and the support property
	recorded after \eqref{eq:signed-distance-tangent-extension}, this is
	assertion~(\ref{it:cf-frame}), and Proposition~\ref{prop:complete-fields} is
	proved.
\end{proof}

A vector field $X\in C^1(\bR^3;\bR^3)$ is \emph{tangential} if $X\cdot\nabla\rho_\cD=0$ on $\cN_\cD$.
By (\ref{it:cf-tangential}) of Proposition~\ref{prop:complete-fields}
every \(\tau_\nu\) is tangential, and by \eqref{eq:tangential-field-family} each
is of class \(C^{K_\cD+1}\), hence of class \(C^{2\mathrm N+1}\) by
\eqref{eq:order-and-domain-regularity}.  Since \(\overline\cD\) is compact and
\(\rho_\cD\) is constant along the flow of a tangential field on
\(\cN_\cD\), that flow exists for every \(s\), preserves \(\overline\cD\),
and preserves every level set of \(\rho_\cD\) inside the collar.

\subsection{Ordered complete words}\label{ss:diffuse-words}

For \(\nu\in\mathfrak A\), write
\(\mathfrak f_\nu=\mathfrak f_{\tau_\nu}\),
\(\Phi_s^\nu=\Phi_s^{\tau_\nu}\),
\(\widehat\Phi_s^\nu=\widehat\Phi_s^{\tau_\nu}\), and
\(T_s^\nu=T_s^{\tau_\nu}\).
The spatial flow preserves \(\rho_\cD\) by tangency, and the lifted flow
preserves \(\bar n\cdot v\) by
\eqref{eq:complete-derivative-preserves-grazing}. Hence \(\gamma^+\),
\(\gamma^-\), and the grazing set \(\gamma^0\) are invariant.

For an ordered word \(I=(\nu_1,\ldots,\nu_m)\in\mathfrak A^m\) we set
\begin{equation}
	\label{eq:ordered-complete-derivatives}
	\mathfrak D_I
	:=\mathfrak D_{\nu_1}\cdots\mathfrak D_{\nu_m},
	\qquad
	\tau_I:=(\tau_{\nu_1}\cdot\nabla_x)\cdots(\tau_{\nu_m}\cdot\nabla_x)
	\quad\text{on }\partial\cD,
	\qquad
	|I|:=m,
	\qquad
	\mathfrak D_\varnothing=\tau_\varnothing:=\operatorname{Id}.
\end{equation}
The factors are not commuted.

Adjoin to \(\mathfrak A\) one temporal letter and set
\begin{equation}
	\label{eq:dyn-mixed-alphabet}
	\mathfrak A_\bullet:=\mathfrak A\cup\{\mathrm t\},
	\qquad
	\mathfrak D_{\mathrm t}:=\partial_t,
	\qquad
	d_{\mathrm t}:=0 .
\end{equation}
For \(I=(\nu_1,\ldots,\nu_m)\in\mathfrak A_\bullet^m\), retain the indicated
order in \(\mathfrak D_I\) and put
\[
j(I):=\#\{i:\nu_i=\mathrm t\},
\qquad
I'\in\mathfrak A^{m-j(I)}.
\]
Thus \(j(I)\) counts the temporal letters, while \(I'\) is the spatial
subword in its original order. The notation \(\tau_I\) is used only for
spatial words.

For spatial words $I=(\nu_1,\ldots,\nu_m)\in\mathfrak A^m$, integration by parts in $(x,v)$ gives
\begin{alignat}{2}
	\mathfrak D_\nu^\dagger\,&:=-\mathfrak D_\nu-d_\nu,
	\qquad &\int_{\cD\times\bR^3} (\mathfrak D_\nu f)\psi
	&=\int_{\cD\times\bR^3} f\,\mathfrak D_\nu^\dagger\psi,
	\label{eq:complete-field-formal-adjoint}\\
	\mathfrak D_I^\dagger
	\,&:=\mathfrak D_{\nu_m}^\dagger\cdots\mathfrak D_{\nu_1}^\dagger,\qquad
	&\int_{\cD\times\bR^3} (\mathfrak D_If)\psi
	&=\int_{\cD\times\bR^3} f\,\mathfrak D_I^\dagger\psi,
	\label{eq:ordered-complete-formal-adjoint}
\end{alignat}
for \(f,\psi\in C_c^\infty(\overline\cD\times\bR^3)\), the first identity being
(\ref{it:cf-adjoint}) of Proposition~\ref{prop:complete-fields} and the
second its iterate. These are fixed-time spatial adjoint identities;
when a word contains temporal letters, the temporal adjoint is taken
separately on compactly supported space-time tests.

For \(m\ge1\) and \(S=\{s_1<\cdots<s_k\}\subset\{1,\ldots,m\}\) we put
\begin{equation}
	\label{eq:ordered-position-subwords}
	[m]:=\{1,\ldots,m\},
	\qquad
	I_S:=(\nu_{s_1},\ldots,\nu_{s_k}),
	\qquad
	S^c:=[m]\setminus S,
	\qquad
	I_\varnothing:=\varnothing.
\end{equation}
We write \(I_1\sqcup I_2=I\) for the ordered pair
\((I_S,I_{S^c})\), counting each subset \(S\subset[m]\) once. Thus
\begin{equation}
	\label{eq:ordered-two-coloring-convention}
	\sum_{I_1\sqcup I_2=I}F(I_1,I_2)
	:=\sum_{S\subset[m]}F(I_S,I_{S^c}).
\end{equation}
Conditions such as \(I_1\ne\varnothing\) restrict the subsets in this sum.
Equal letters in distinct positions are counted separately.

For tangential $X$ and spatial words $I\in\mathfrak A^m$, their
restrictions to smooth boundary functions are
\begin{equation}
	\label{eq:boundary-complete-derivative}
	\mathfrak D_X^{\pm}
	:=\mathfrak D_X\big|_{\gamma^\pm},
	\qquad
	\mathfrak D_I^{\pm}
	:=\mathfrak D_{\nu_1}^{\pm}\cdots\mathfrak D_{\nu_m}^{\pm},
	\qquad
	\mathfrak D_\nu^{\pm}=\mathfrak D_{\tau_\nu}^{\pm}.
\end{equation}
These restrictions are well defined because \(\mathfrak f_X\) is tangent to
the phase boundary. On functions independent of \(v\), the operator
\(\mathfrak D_I^{\pm}\) agrees with the boundary tangential word \(\tau_I\).

\subsection{Wall profiles}\label{ss:diffuse-wall}

The wall temperature appears in this part through the isothermal profile \(\mathfrak m_\circ\),
the difference \(d_\delta\) and the datum \(h_w\) of
\eqref{eq:wall-response-profiles}--\eqref{eq:nonisothermal-diffuse-response}.
The profiles have Gaussian velocity decay, with an additional factor
$\delta$ for $d_\delta$ and $h_w$. The estimates below retain these
properties under chart derivatives and complete derivatives.
The following lemma applies to any fixed
\(\vartheta\in C^{K_\cD}(\partial\cD)\) with
\(\|\vartheta\|_{C^{K_\cD}}\le1\), in particular to
\(\vartheta(t,\cdot)\) and \(\vartheta_\infty\) under
\eqref{eq:dyn-theta-amplitude}--\eqref{eq:dyn-theta-relaxation}.

\begin{lemma}[Uniform wall-profile derivatives]
	\label{lem:wall-profile-derivatives}
	Let \(\phi_\kappa(y',0)\) be the boundary
	parametrization of \eqref{eq:boundary-flattening-maps}.  For every finite
	\(L\ge0\), every integer \(K_v\ge0\) and every integer \(0\le K_x\le K_\cD\)
	there are \(c>0\) and \(C<\infty\), depending on \(L,K_v,K_x\) and on the
	fixed atlas but not on \(\delta\), such that on every boundary chart
	\begin{equation}
		\label{eq:wall-profile-chart-bound}
		\sup_{|\alpha|\le K_v,\ |\beta|\le K_x}
		\left|w(v)^L\partial_{y'}^\beta\partial_v^\alpha
		d_\delta(\phi_\kappa(y',0),v)\right|
		\le C\delta\ee^{-c|v|^2}.
	\end{equation}
	The estimate also holds for \(h_w=c_\mu d_\delta\).
	For each ordered spatial word \(I\) of length at most \(\mathrm N\),
	\eqref{eq:ordered-complete-derivatives} gives, with a constant also
	depending on \(I\),
	\begin{equation}
		\label{eq:wall-profile-complete-bound}
		|w(v)^L\mathfrak D_Ih_w(x,v)|
		\le C_{I,L}\delta\ee^{-c|v|^2}
		\qquad ((x,v)\in\gamma^-\cup\gamma^0),
	\end{equation}
	the value on \(\gamma^0\) being the incoming continuous limit.
\end{lemma}

\begin{proof}
	The normalization in \eqref{new:eq:intro-wall-maxwellian} gives
	\[
	\mathfrak m_\Theta(v)=(2\pi\Theta^2)^{-1}
	\exp\!\left[-\left(\frac1{2\Theta}-\frac12\right)|v|^2\right].
	\]
	For \(\Theta\in[1/4,3/4]\) the coefficient in the exponent is at least
	\(1/6\). Repeated differentiation in \(\Theta\) and \(v\) therefore gives
	\begin{equation}\label{eq:wall-maxwellian-ratio-bound}
		\sup_{1/4\le\Theta\le3/4}
		w(v)^L|\partial_\Theta^k\partial_v^\alpha\mathfrak m_\Theta(v)|
		\le C_{L,k,\alpha}e^{-|v|^2/12}
	\end{equation}
	for fixed \(L,k,\alpha\). Indeed, the derivatives are the same Gaussian
	multiplied by polynomials in \(v\) of degree at most \(2k+|\alpha|\),
	with coefficients bounded on this temperature interval.
	
	The mean-value formula
	\begin{equation*}
		d_\delta(x,v)=\frac\delta2\vartheta(x)
		\int_0^1\partial_\Theta\mathfrak m_{\frac12(1+s\delta\vartheta(x))}(v)\,ds
	\end{equation*}
	retains the explicit factor $\delta$ under every chart derivative.
	The derivatives of $\vartheta$ through order $K_\cD$ are bounded, while
	\eqref{eq:wall-maxwellian-ratio-bound} controls the resulting temperature
	derivatives. Thus
	\eqref{eq:wall-profile-chart-bound} holds; the estimate for \(h_w=c_\mu d_\delta\)
	follows immediately.
	
	On a boundary chart, each \(\mathfrak D_\nu\) is a tangential spatial
	derivative plus a velocity derivative with coefficient linear in \(v\).
	Expanding an ordered word of length at most \(\mathrm N\) consequently
	produces chart derivatives of total order at most \(\mathrm N\), multiplied
	by polynomials in \(v\) with bounded coefficients. Since
	\(\mathrm N\le K_\cD\), \eqref{eq:wall-profile-complete-bound} follows from
	the chart estimate after a reduction of the Gaussian exponent. The profiles
	are continuous in \(v\) across \(n_x\cdot v=0\), which also gives the
	assertion on \(\gamma^0\).
\end{proof}

\begin{lemma}[Time derivatives of $d_\delta$]
	\label{lem:dyn-theta-profile-time-derivatives}
	Assume (\ref{it:dyn-theta-amplitude}).  For every finite \(L\ge0\),
	every integer \(K_v\ge0\), every integer \(0\le K_x\le K_\cD\) and every
	\(0\le j\le\mathrm N+1\) there are \(c>0\) and \(C<\infty\), depending on
	\(L,K_v,K_x,j\) but not on \(\delta\) or on \(t\), such that on every
	boundary chart
	\begin{equation}
		\label{eq:dyn-theta-defect-time-bound}
		\sup_{|\alpha|\le K_v,\ |\beta|\le K_x}
		\bigl|w(v)^L\partial_t^j\partial_{y'}^\beta\partial_v^\alpha
		d_\delta\bigl(t,\phi_\kappa(y',0),v\bigr)\bigr|
		\le C\delta\ee^{-c|v|^2}.
	\end{equation}
	Under the additional relaxation condition
	\eqref{eq:dyn-theta-relaxation}, the same indices satisfy
	\begin{equation}
		\label{eq:dyn-theta-defect-relaxation-bound}
		\sup_{|\alpha|\le K_v,\ |\beta|\le K_x}
		\bigl|w(v)^L\partial_t^j\partial_{y'}^\beta\partial_v^\alpha
		\bigl(d_\delta(t)-d_\delta^\infty\bigr)
		\bigl(\phi_\kappa(y',0),v\bigr)\bigr|
		\le C\delta\,\ee^{-\lambda_wt}\ee^{-c|v|^2}.
	\end{equation}
	The same estimates hold for \(h_w=c_\mu d_\delta\) and
	\(h_w-h_w^\infty=c_\mu(d_\delta-d_\delta^\infty)\).
	For each ordered spatial word \(I\) with \(|I|\le\mathrm N\),
	\begin{equation}
		\label{eq:dyn-theta-defect-complete-bound}
		\bigl|w(v)^L\partial_t^j\mathfrak D_Ih_w(t,x,v)\bigr|
		\le C_{I,L}\,\delta\ee^{-c|v|^2},
		\qquad
		\bigl|w(v)^L\partial_t^j\mathfrak D_I
		\bigl(h_w(t)-h_w^\infty\bigr)(x,v)\bigr|
		\le C_{I,L}\,\delta\,\ee^{-\lambda_wt}\ee^{-c|v|^2},
	\end{equation}
	for \((x,v)\in\gamma^-\cup\gamma^0\), the value on \(\gamma^0\) being the
	incoming continuous limit.
\end{lemma}

\begin{proof}
	The range \(\Theta_w\in[1/4,3/4]\) allows the use of
	\eqref{eq:wall-maxwellian-ratio-bound} uniformly in \(t\) and \(\delta\).
	For \(j\ge1\), the chain rule expresses \(\partial_t^j d_\delta\) as a
	sum of terms
	\[
	(\partial_\Theta^k\mathfrak m_\Theta)|_{\Theta=\Theta_w}
	\prod_{i=1}^k\partial_t^{j_i}\Theta_w,
	\qquad j_i\ge1,\quad \sum_{i=1}^k j_i=j.
	\]
	Since $\partial_t^{j_i}\Theta_w=(\delta/2)\partial_t^{j_i}\vartheta$,
	every term is $O(\delta)$, including after chart and velocity
	differentiation by \eqref{eq:wall-maxwellian-ratio-bound}.
	For $j=0$ use \eqref{eq:wall-profile-chart-bound}. This proves
	\eqref{eq:dyn-theta-defect-time-bound}.
	
	For the relaxation bound, write
	\[
	d_\delta(t)-d_\delta^\infty
	=\frac\delta2(\vartheta(t,x)-\vartheta_\infty(x))
	\int_0^1\partial_\Theta\mathfrak m_{\Theta_s(t,x)}\,ds,
	\qquad
	\Theta_s=(1-s)\Theta_w^\infty+s\Theta_w(t).
	\]
	A derivative falling on the first factor retains a derivative of
	$\vartheta-\vartheta_\infty$. A positive-order time derivative falling
	on the integral yields a factor
	$\partial_t^i\Theta_s=(s\delta/2)\partial_t^i
	(\vartheta-\vartheta_\infty)$, $i\ge1$; spatial derivatives of the
	integral do not remove the first factor. Hence every resulting term
	contains a factor bounded by $e^{-\lambda_wt}$ by
	\eqref{eq:dyn-theta-relaxation}. The other factors are uniformly bounded
	by \eqref{eq:dyn-theta-amplitude} and
	\eqref{eq:wall-maxwellian-ratio-bound}. Applying further velocity
	derivatives and absorbing polynomial factors into the Gaussian proves
	\eqref{eq:dyn-theta-defect-relaxation-bound}.
	
	The statements for \(h_w\) follow by multiplication by \(c_\mu\).
	As in the proof of Lemma~\ref{lem:wall-profile-derivatives}, an ordered
	complete word is a finite sum of chart derivatives with polynomial
	velocity coefficients; it commutes with time differentiation. Thus
	\eqref{eq:dyn-theta-defect-complete-bound} follows from the two chart
	bounds, including the continuous representatives on \(\gamma^0\).
\end{proof}

In the adapted boundary coordinates \eqref{eq:rotating-frame-representative},
a boundary profile \(a(x,v)\) is represented by
\begin{equation}
	\label{eq:boundary-profile-chart-representation}
	a_\kappa(y',\xi)
	:=a\bigl(\phi_\kappa(y',0),
	R_\kappa(\phi_\kappa(y',0))\xi\bigr),
	\qquad \kappa\in\mathfrak K_{\p} .
\end{equation}

\begin{defn}[Gaussian envelope]
	\label{def:gaussian-polynomial-envelope}
	A boundary profile \(a\in C^{\mathrm N}(\partial\cD\times\bR^3)\) has a
	\emph{Gaussian envelope through order \(\mathrm N\)} if there are \(c_a>0\)
	and \(C_a<\infty\) such that
	\begin{equation}
		\label{eq:boundary-gaussian-polynomial-envelope}
		\max_{\kappa\in\mathfrak K_{\p}}
		\max_{|\alpha|+|\beta|\le\mathrm N}
		|\partial_{\xi}^{\alpha}\partial_{y'}^{\beta}a_\kappa(y',\xi)|
		\le C_a \ee^{-c_a|\xi|^2}.
	\end{equation}
\end{defn}

A derivative in \(y'\) of \eqref{eq:boundary-profile-chart-representation}
also differentiates the rotation \(R_\kappa\), producing factors linear
in \(\xi\). The chart bound \eqref{eq:wall-profile-chart-bound} with
\(K_v=K_x=\mathrm N\), together with the regularity of \(R_\kappa\) in
\eqref{eq:rotating-frame}, therefore gives Gaussian envelopes for
\(d_\delta\) and \(h_w\) with \(C_a\le C\delta\). Since
\(\mathfrak m_{\Theta_w}=\mathfrak m_\circ+d_\delta\), the same holds for
\(\mathfrak m_{\Theta_w}\), with \(C_a\le C\).

For \(|K|\le\mathrm N\), expanding \(\mathfrak D_K\) and taking chart
derivatives through order \(\mathrm N\) requires at most
\(2\mathrm N\) derivatives of the profiles. This is permitted by
\(K_\cD\ge2\mathrm N\) in \eqref{eq:order-and-domain-regularity} and
\(\tau_\nu,R_\kappa\in C^{K_\cD+1}\). Thus \(\mathfrak D_Ka\) has a
Gaussian envelope, with \(C_a\le C\delta\) for \(a=d_\delta,h_w\) and
\(C_a\le C\) for \(a=\mathfrak m_{\Theta_w}\). The polynomial
velocity coefficients are absorbed by reducing the Gaussian rate; since
the fields are tangential, these derivatives depend only on the boundary
values of \(a\). Under \eqref{eq:dyn-theta-relaxation}, the same
argument, using \eqref{eq:dyn-theta-defect-relaxation-bound}, gives
\(C_a\le C\delta\ee^{-\lambda_wt}\) for \(a=h_w(t)-h_w^\infty\).

\subsection{The collision operator and weighted estimates}

We record the decomposition of the linearized collision operator and the weighted
estimates used below, followed by equivalence and interpolation inequalities for
the coercive norm.

Throughout the diffuse part $\phi$ is the Coulomb kernel
\eqref{new:eq:intro-coulomb-kernel}.

The compact part of the Maxwellian linearization consists of a velocity
multiplier and an integral operator. Define
\begin{align}
	(\cK_{0,{\rm int}}f)(v)
	&:=\int_{\bR^3}k_0^{\rm col}(v,v')f(v')\dd v',\nonumber\\
	k_0^{\rm col}(v,v')
	&:=\sqrt{\mu(v)\mu(v')}
	\left\{2\operatorname{tr}\phi(v-v')
	-4v_iv_j\phi^{ij}(v-v')\right\}.
	\label{eq:compact-integral-kernel}
\end{align}
With $\cA_\mu=-A$ and $\cK_0=-K$ as in \eqref{def:LAK}, we have
$L=\cA_\mu+\cK_0$. The identity
$\sigma_{\sqrt\mu f}^{ij}=\phi^{ij}*(\sqrt\mu f)$ and
\eqref{def_Kg} give
\begin{equation}
	\label{eq:compact-integral-is-inflow-K}
	\cK_{0,{\rm int}}=-\cK,
	\qquad
	\cK_0=\cK_{0,{\rm int}}-8\pi\mu ,
\end{equation}
where the second identity follows from \eqref{eq:explicit-landau-K}.

For $\ell\ge0$ and $\mathrm L\in\bR$, the Gaussian first-slot bound
\eqref{eq:velalg} with $\varpi=\tfrac14$ and
$\mu^{1/4}\le C_{\mathrm L}w^{\mathrm L}$ yield
\begin{equation}
	\label{eq:velalg-poly}
	\|w^{\ell}N[g,f]\|_{H_\sigma^{-1}(\bR_v^3)}
	\le C_{\ell,\mathrm L}\,\|w^{\mathrm L}g\|_{L_v^2}\,
	\|w^{\ell}f\|_{H_\sigma^1(\bR_v^3)} .
\end{equation}

For $L\ge0$, the weighted stationary inflow estimate reads
\begin{align}
	\|w^Lu\|_{L^2(\cD\times\bR_v^3)}
	+\|w^Lu\|_{L_x^2H_\sigma^1}
	+\|w^L v\cdot\nabla_xu\|_{L_x^2H_\sigma^{-1}}
	&\notag\\
	+\|w^Lu_{\gamma^-}\|_{L^2(\gamma^-)}
	+\|w^Lu_{\gamma^+}\|_{L^2(\gamma^+)}
	&\le C_L\left(
	\|w^LS\|_{L_x^2H_\sigma^{-1}}
	+\|w^LH\|_{L^2(\gamma^-)}\right)
	\label{eq:weighted-inflow-estimate}
\end{align}
for the solution of \eqref{solv:linland_st} with source $S$ and incoming datum $H$.
This follows from Theorem~\ref{thm:weighted-steady}(2) and the trace estimate
of Proposition~\ref{2607151218}.

We use the common coefficient threshold
$\varepsilon_{\rm in}=\varepsilon_0(\cD)>0$ of
Theorem~\ref{260322442}, Lemma~\ref{260527559}, and
Theorem~\ref{thm:weighted-steady} in \eqref{eq:weighted-inflow-estimate}.

The evolutionary estimate is
\begin{align}
	&\|w^Lf\|_{L^2(\cD_T\times\bR_v^3)}
	+\|w^Lf\|_{L^2(\cD_T;H_\sigma^1)}
	+\|w^L(\partial_t+v\cdot\nabla_x)f\|_{L^2(\cD_T;H_\sigma^{-1})}\notag\\
	&+\sup_{0\le t\le T}\|w^Lf(t)\|_{L^2(\cD\times\bR_v^3)}
	+\|w^Lf_{\gamma^-}\|_{L^2(\gamma_T^-)}
	+\|w^Lf_{\gamma^+}\|_{L^2(\gamma_T^+)}\notag\\
	&\qquad \le C_L\Bigl(
	\|w^LF\|_{L^2(\cD_T;H_\sigma^{-1})}
	+\|w^Lf_0\|_{L^2(\cD\times\bR_v^3)}
	+\|w^LH\|_{L^2(\gamma_T^-)}\Bigr),
	\label{eq:weighted-graph-convention}
\end{align}
for $(\partial_t+v\cdot\nabla_x)f+L_gf=F$, $f(0)=f_0$, and
$f_{\gamma^-}=H$. Theorem~\ref{thm:weighted-steady}(1) and
Proposition~\ref{2607151218} give this bound and identify the weighted
boundary and time traces. The constant depends only on $\cD$ and $L$, not
on $T$ or on $g$ satisfying \eqref{eq:weighted-dynamic-coefficient}.
We use the same estimate on translated time intervals.

\begin{lemma}[Weighted coercive norms at real weights]
	\label{lem:real-weight-Hsigma-equivalence}
	For every $\mathrm L\in\bR$ there is $C_{\mathrm L}\ge1$ such that, for every
	$f$ with $w^{\mathrm L}f\in H^1_\sigma(\bR_v^3)$,
	\begin{equation}
		\label{eq:real-weight-Hsigma-equivalence}
		C_{\mathrm L}^{-1}\norm{w^{\mathrm L}f}_{H^1_\sigma}^2
		\le\int_{\bR^3}w^{2\mathrm L}\bigl(
		\sigma^{ij}\p_if\,\p_jf+\sigma^{ij}v_iv_jf^2\bigr)\dd v
		\le C_{\mathrm L}\norm{w^{\mathrm L}f}_{H^1_\sigma}^2 ,
	\end{equation}
	and, for every $\mathrm L'\le\mathrm L$,
	\begin{equation}
		\label{eq:real-weight-Hsigma-monotonicity}
		\norm{w^{\mathrm L'}f}_{H^1_\sigma}
		\le C_{\mathrm L',\mathrm L}\norm{w^{\mathrm L}f}_{H^1_\sigma} .
	\end{equation}
\end{lemma}

\begin{proof}
	The product-rule proof of Lemma~\ref{lem:weighted-Hsigma-equivalence}
	uses only
	\[
	|\nabla\log w^{\mathrm L}|_\sigma^2
	=\frac{\mathrm L^2}{\la v\ra^4}\,\sigma^{ij}v_iv_j
	\le\mathrm L^2\,\sigma^{ij}v_iv_j ,
	\]
	so that the constant depends on $\mathrm L$ through $|\mathrm L|$; and
	\eqref{eq:real-weight-Hsigma-monotonicity} follows from
	\eqref{eq:real-weight-Hsigma-equivalence} at the two weights together with
	$w^{2\mathrm L'}\le w^{2\mathrm L}$.
\end{proof}

\begin{prop}[Polynomial-weight interpolation]
	\label{prop:polynomial-weight-interpolation}
	Let $0\le\theta\le1$, let $\mathrm L_0,\mathrm L_1\in\bR$, let
	$1\le p_0,p_1\le\infty$, and let
	\begin{equation}
		\label{eq:polynomial-weight-slots}
		\mathrm L_\theta:=\theta\mathrm L_1+(1-\theta)\mathrm L_0,
		\qquad
		\frac1{p_\theta}:=\frac{\theta}{p_1}+\frac{1-\theta}{p_0} ,
	\end{equation}
	where the same parameter $\theta$ interpolates both weights and spatial exponents.
	Then
	\begin{equation}
		\label{eq:polynomial-weight-interpolation}
		\|w^{\mathrm L_\theta}u\|_{L_x^{p_\theta}H_\sigma^1}
		\le C_{\mathrm L_0,\mathrm L_1}
		\|w^{\mathrm L_1}u\|_{L_x^{p_1}H_\sigma^1}^{\theta}
		\|w^{\mathrm L_0}u\|_{L_x^{p_0}H_\sigma^1}^{1-\theta}.
	\end{equation}
\end{prop}

\begin{proof}
	For \(\mathrm L\in\bR\), put
	\[
	E_{\mathrm L}(u)=\int_{\bR^3}w^{2\mathrm L}
	\bigl(\sigma^{ij}\partial_i u\,\partial_j u
	+\sigma^{ij}v_iv_j|u|^2\bigr)\,dv.
	\]
	By \eqref{eq:real-weight-Hsigma-equivalence},
	\(E_{\mathrm L}(u)^{1/2}\simeq_{\mathrm L}
	\|w^{\mathrm L}u\|_{H_\sigma^1}\), uniformly for
	\(\mathrm L\) between \(\mathrm L_0\) and \(\mathrm L_1\).
	The integrand without its weight is nonnegative. Since
	\(w^{2\mathrm L_\theta}=(w^{2\mathrm L_1})^\theta
	(w^{2\mathrm L_0})^{1-\theta}\), H\"older in \(v\) gives
	\[
	E_{\mathrm L_\theta}(u(x))^{1/2}
	\le E_{\mathrm L_1}(u(x))^{\theta/2}
	E_{\mathrm L_0}(u(x))^{(1-\theta)/2}.
	\]
	H\"older in \(x\), with
	\(1/p_\theta=\theta/p_1+(1-\theta)/p_0\), completes the proof.
	The cases \(\theta=0,1\) and \(p_i=\infty\) follow by the usual
	endpoint interpretation.
\end{proof}

\begin{lemma}[Weight interpolation in the coercive slot]
	\label{lem:dyn-hier-weight-trade}
	Let $\mathsf Q$ be either $\cD\times\bR_v^3$ or $(0,T)\times\cD\times\bR_v^3$ with
	$0<T<\infty$, let $U$ be measurable on $\mathsf Q$, let $\ell\ge0$,
	let $0<\theta\le1$, let $2\le r,p\le\infty$, and let
	\begin{equation}
		\label{eq:weight-trade-output-exponents}
		\frac1{\widetilde r}:=\frac\theta2+(1-\theta)\frac1r,
		\qquad
		\frac1{\widetilde p}:=\frac\theta2+(1-\theta)\frac1p.
	\end{equation}
	Whenever the right-hand side is finite,
	\begin{equation}
		\label{eq:dyn-hier-weight-trade}
		\|w^\ell U\|_{L_t^{\widetilde r}L_x^{\widetilde p}H_\sigma^1}
		\le C_{\ell,\theta}
		\|U\|_{L_t^rL_x^pH_\sigma^1}^{1-\theta}
		\,\bigl\|w^{\ell/\theta}U\bigr\|
		_{L_t^2L_x^2H_\sigma^1}^{\theta},
	\end{equation}
	all norms over $\mathsf Q$, with the outer $L_t^{\widetilde r}$, $L_t^r$ and $L_t^2$ absent
	when $\mathsf Q=\cD\times\bR_v^3$.
	The constant depends only on $\ell$ and $\theta$.  In particular it
	does not depend on the horizon $T$, nor on any parameter introduced
	later in this part.
\end{lemma}

\begin{proof}
	Apply Proposition~\ref{prop:polynomial-weight-interpolation} with
	\((\mathrm L_0,\mathrm L_1)=(0,\ell/\theta)\) and
	\((p_0,p_1)=(p,2)\) at each fixed \(t\). It gives
	\[
	\|w^\ell U(t)\|_{L_x^{\widetilde p}H_\sigma^1}
	\le C_{\ell,\theta}
	\|U(t)\|_{L_x^pH_\sigma^1}^{1-\theta}
	\|w^{\ell/\theta}U(t)\|_{L_x^2H_\sigma^1}^\theta.
	\]
	H\"older in time with $1/\widetilde r=(1-\theta)/r+\theta/2$
	yields \eqref{eq:dyn-hier-weight-trade}. Both applications of
	H\"older have constants independent of $T$; the stationary case
	requires only the spatial estimate.
\end{proof}

\section{Linear Diffuse Problems and Macroscopic Coercivity}\label{sec:diffuse-L2}

We use the macro--micro decomposition $f=\mathbf Pf+(\mathbf I-\mathbf P)f$
of \S\ref{subsec:macro-landau}, with the coefficients $\mathbf a,\mathbf b,\mathbf c$
of \hyperlink{comp:macro-moments}{Gaussian moment formulas}. The wall temperature is $\Theta_w=\frac12(1+\delta\vartheta)$, where
$\|\vartheta\|_{C^{K_\cD}(\partial\cD)}\le1$; in the evolutionary problem the boundary
condition holds for almost every time, whereas the stationary wall is
$\Theta_w^\infty$ from \eqref{eq:limiting-wall-profiles}.
Constants are uniform over these wall profiles.

\subsection{Boundary algebra and traces}\label{ss:diffuse-L2-inputs}%

We use the diffuse projection of
Guo~\cite[(49)]{Guo2010} and the non-isothermal normalization of
Esposito--Guo--Kim--Marra~\cite{EspositoGuoKimMarra}.

Let $-\Delta_{\p\cD}\ge0$ be the Laplace--Beltrami operator of $\partial\cD$
and $A_\partial:=-\Delta_{\p\cD}$.  We use boundary Sobolev regularity through order $\mathrm N+2$: a word of length $\mathrm N$ requires $\mathrm N$ tangential derivatives, and the embedding $H^{m+2}(\partial\cD)\hookrightarrow C^{m,\alpha}(\partial\cD)$ requires two more. This regularity follows from \eqref{eq:order-and-domain-regularity}.  We use the spectral norm
\begin{equation}
	\label{eq:word-boundary-spectral-Hs}
	\|u\|_{H^s_\partial}:=
	\|(1+A_\partial)^{s/2}u\|_{L^2(\p\cD)},
	\qquad 0\le s\le2\mathrm N.
\end{equation}

For the approximate wall problem, choose
$\chi\in C_c^\infty([0,\infty))$ with $0\le\chi\le1$
and $\chi=1$ in a neighborhood of zero, and set, on $L^2(\p\cD)$,
\begin{equation}
	\label{eq:boundary-spectral-cutoff-definition}
	J_N:=\chi(-\Delta_{\p\cD}/N^2),
	\qquad E_N:=\operatorname{Ran}J_N,
	\qquad N\in\mathbb N.
\end{equation}
The cutoff $J_N$ need not be a projection. If
$-\Delta_{\p\cD}e_k=\lambda_ke_k$, then $E_N$ is spanned by the
$e_k$ for which $\chi(\lambda_k/N^2)\ne0$, and
$J_N:E_N\to E_N$ is invertible.
Since the zero eigenspace of $-\Delta_{\p\cD}$ contains the constants, we have
$J_N1=1$.  For the truncated wall $P_{w,N}$ of
\eqref{eq:boundary-spectral-diffuse-operator} below, the unit-flux
normalization of $M_\Theta$ gives, for a.e.\ $x\in\p\cD$,
\[
\int_{n_x\cdot v<0}P_{w,N}f_{\gamma^+}(x,v)\sqrt\mu(v)|n_x\cdot v|\dd v
=(J_N\cA_0f_{\gamma^+})(x),
\]
which need not agree pointwise with $\cA_0f_{\gamma^+}(x)$.
Its integral over $\p\cD$ is nevertheless unchanged: $J_N$ is self-adjoint
and $J_N1=1$ (Lemma~\ref{lem:boundary-spectral-cutoffs}), so
\begin{equation}
	\label{eq:truncated-wall-flux-integrated}
	\int_{\p\cD}J_N\cA_0f_{\gamma^+}\dd S_x
	=\langle\cA_0f_{\gamma^+},J_N1\rangle_{L^2(\p\cD)}
	=\int_{\p\cD}\cA_0f_{\gamma^+}\dd S_x ,
\end{equation}
at every $N$.  

\begin{lemma}[Spectral cutoffs]
	\label{lem:boundary-spectral-cutoffs}
	With $J_N$ and $E_N$ as in
	\eqref{eq:boundary-spectral-cutoff-definition}, the space
	$E_N$ is finite dimensional, $J_N$ is self-adjoint with $0\le J_N\le\mathbf I$
	on $L^2(\p\cD)$, hence an $L^2(\p\cD)$ contraction, $J_N1=1$, and
	\begin{equation}
		\label{eq:boundary-spectral-strong-convergence}
		J_N\psi\longrightarrow\psi
		\qquad\text{strongly in }L^2(\p\cD)
	\end{equation}
	for every fixed $\psi\in L^2(\p\cD)$.
\end{lemma}

\begin{proof}
	Let $\{e_k\}$ be an orthonormal eigenbasis of $-\Delta_{\p\cD}$, with
	eigenvalues $0\le\lambda_k\uparrow\infty$.  Then
	\[
	J_N\psi=\sum_k\chi(\lambda_k/N^2)\langle\psi,e_k\rangle e_k.
	\]
	The multiplier has finite support in the spectrum and takes values
	in $[0,1]$, so $J_N$ is finite-rank and self-adjoint with $0\le J_N\le\mathbf I$.
	As $\chi(0)=1$, it fixes constants. For each $k$,
	$\chi(\lambda_k/N^2)\to1$; dominated convergence in $\ell^2$
	gives \eqref{eq:boundary-spectral-strong-convergence}.
\end{proof}

\begin{lemma}[Boundary Sobolev norms]
	\label{lem:word-boundary-spectral-chart}
	For $0\le s\le2\mathrm N$, the spectral space in
	\eqref{eq:word-boundary-spectral-Hs} is the intrinsic chart Sobolev space
	$H^s(\p\cD)$, with equivalent norms.  For a fixed smooth partition
	$\{\varrho_\kappa\}_{\kappa\in\mathfrak K_{\p}}$ subordinate to the
	boundary atlas, set, for $s\ge0$,
	\[
	\|u\|_{H^s(\p\cD)}^2
	:=\sum_{\kappa\in\mathfrak K_{\p}}
	\|(\varrho_\kappa u)\circ\phi_\kappa(\cdot,0)\|_{H^s(\bR^2)}^2.
	\]
	Then, for every integer $0\le m\le2\mathrm N-2$,
	\begin{equation}
		\label{eq:boundary-elliptic-finite-order}
		\|u\|_{H^{m+2}(\p\cD)}
		\le C_m\bigl(\|A_\partial u\|_{H^m(\p\cD)}+\|u\|_2\bigr),
		\qquad
		\|A_\partial u\|_{H^m(\p\cD)}
		\le C_m\|u\|_{H^{m+2}(\p\cD)}.
	\end{equation}
	In particular the spectral cutoff $J_N$ of
	Lemma~\ref{lem:boundary-spectral-cutoffs} is contractive in the spectral
	$H^{\mathrm N}_\partial$ norm and uniformly bounded in the chart
	$H^{\mathrm N}$ norm.
\end{lemma}

\begin{proof}
	In boundary charts, $-\Delta_{\p\cD}$ is uniformly elliptic in divergence form with
	$C^{K_\cD+1}$ coefficients. The standard local elliptic estimates,
	patched over the compact boundary, give the first inequality in
	\eqref{eq:boundary-elliptic-finite-order}
	\cite[Theorems~8.8 and~8.10]{GilbargTrudinger1983}; the converse follows
	from the coordinate formula for $\Delta_{\p\cD}$.
	
	For $u\in C^{2\mathrm N}(\partial\cD)$, iterating \eqref{eq:boundary-elliptic-finite-order} and its
	reverse gives, for $0\le j\le\mathrm N$,
	\[
	\|u\|_{H^{2j}(\p\cD)}\simeq \|(1+A_\partial)^ju\|_2.
	\]
	The same interior estimates and \eqref{eq:order-and-domain-regularity} show that the
	eigenfunctions of $A_\partial$ belong to $C^{2\mathrm N}(\partial\cD)$, so completing
	$C^{2\mathrm N}(\partial\cD)$ identifies the corresponding Hilbert spaces
	at the even orders $0,2,\ldots,2\mathrm N$.
	Complex interpolation between
	consecutive even orders \cite[Chapter~4]{TaylorPDE1} gives the equivalence for every real
	$0\le s\le2\mathrm N$, in particular $s=\mathrm N$.  Finally
	$J_N=\chi(A_\partial/N^2)$ commutes with
	$(1+A_\partial)^{\mathrm N/2}$ and $|\chi|\le1$, so
	\[
	\|J_Nu\|_{H^{\mathrm N}_\partial}\le\|u\|_{H^{\mathrm N}_\partial}.
	\]
	The norm equivalence obtained from \eqref{eq:boundary-elliptic-finite-order}
	by complex interpolation therefore gives the chart bound.
\end{proof}

The truncated wall operator is
\begin{equation}
	\label{eq:boundary-spectral-diffuse-operator}
	P_{w,N}:=\mathfrak m_{\Theta_w}J_N\cA_0
	:L^2(\gamma^+)\longrightarrow L^2(\gamma^-).
\end{equation}

For $0<\eta<1<R<\infty$ we write
\begin{equation*}
	\gamma_{\eta,R}^\pm
	:=\{(x,v)\in\gamma^\pm:\ |n_x\cdot v|\ge\eta,\ |v|\le R\}.
\end{equation*}

\begin{lemma}[Diffuse boundary algebra]
	\label{lem:diffuse-boundary-algebra}
	Let $f_{\gamma^+}\in L^2(\gamma^+,\dd\gamma)$.
	\begin{enumerate}
		\item\label{it:diffuse-P0-projection} The operator $\cP_0^+$ is the fiberwise orthogonal projection onto
		$\operatorname{span}\{\sqrt\mu\}$, the incoming and outgoing rank-one
		norms agree, and
		\begin{equation}
			\label{eq:exact-diffuse-variance}
			\|f_{\gamma^+}\|_{L^2(\gamma^+)}^2
			-\|\cP_0^-f_{\gamma^+}\|_{L^2(\gamma^-)}^2
			=\|(\mathbf I-\cP_0^+)f_{\gamma^+}\|_{L^2(\gamma^+)}^2.
		\end{equation}
		
		\item There is $C>0$ such that
		\begin{equation}
			\label{eq:wall-defect-flux-bound}
			\sup_{x\in\partial\cD}
			\int_{n_x\cdot v<0}|d_\delta(x,v)|^2|n_x\cdot v|\dd v
			\le C\delta^2,
		\end{equation}
		\begin{equation}
			\label{eq:wall-defect-flux-orthogonality}
			\int_{n_x\cdot v<0}
			d_\delta(x,v)\sqrt\mu(v)|n_x\cdot v|\dd v=0
			\quad\text{for a.e. }x\in\partial\cD,
		\end{equation}
		and hence
		\begin{equation}
			\label{eq:wall-defect-operator-bound}
			\|r_\delta[f_{\gamma^+}]\|_{L^2(\gamma^-)}
			\le C\delta\|f_{\gamma^+}\|_{L^2(\gamma^+)},
			\qquad
			\langle\cP_0^-f_{\gamma^+},r_\delta[f_{\gamma^+}]\rangle_{L^2(\gamma^-)}=0.
		\end{equation}
		
		\item\label{it:diffuse-P0-contraction} For $0<\eta<1<R<\infty$, there is
		$\theta_{\eta,R}\ge0$, depending only on $\eta,R$, such that
		$\theta_{\eta,R}\to0$ as $\eta\downarrow0$ and $R\uparrow\infty$, and
		\begin{equation}
			\label{eq:rank-one-outgoing-tail}
			\|\cP_0^+f_{\gamma^+}\|_{L^2(\gamma^+\setminus\gamma_{\eta,R}^+)}^2
			\le\theta_{\eta,R}
			\|\cP_0^+f_{\gamma^+}\|_{L^2(\gamma^+)}^2.
		\end{equation}
	\end{enumerate}
\end{lemma}

\begin{proof}
	By rotation invariance of $\mu$,
	\[
	c_\mu^{-1}\int_{\pm n_x\cdot v>0}\mu(v)|n_x\cdot v|\dd v=1
	\qquad(x\in\partial\cD),
	\]
	so the single function $c_\mu^{-1/2}\sqrt\mu$ has unit norm in the incoming
	and in the outgoing flux space over every boundary point.  Moreover,
	\[
	\bigl\langle f_{\gamma^+},c_\mu^{-1/2}\sqrt\mu\bigr\rangle
	_{L^2(n_x\cdot v>0,\,n_x\cdot v\dd v)}
	=c_\mu^{-1/2}\cA_0f_{\gamma^+}(x).
	\]
	Hence $\cP_0^+$ is the orthogonal projection onto this line, and
	\[
	\|\cP_0^-f_{\gamma^+}(x,\cdot)\|_{L^2(n_x\cdot v<0,\,|n_x\cdot v|\dd v)}^2
	=\|\cP_0^+f_{\gamma^+}(x,\cdot)\|_{L^2(n_x\cdot v>0,\,n_x\cdot v\dd v)}^2
	=c_\mu^{-1}|\cA_0f_{\gamma^+}(x)|^2.
	\]
	The orthogonal decomposition of $f_{\gamma^+}$ proves
	\eqref{eq:exact-diffuse-variance}.
	
	For $\Theta\in[\tfrac14,\tfrac34]$,
	\[
	\frac{M_\Theta(v)}{\sqrt\mu(v)}
	=C(\Theta)\exp\!\left[-\left(\frac1{2\Theta}-\frac12\right)|v|^2\right],
	\]
	and the coefficient of $|v|^2$ is at least $1/6$.  Applying the mean-value
	formula in $\Theta$, we therefore obtain
	$|d_\delta(x,v)|\le C\delta(1+|v|^2)e^{-c|v|^2}$, which proves
	\eqref{eq:wall-defect-flux-bound}.  The unit-flux normalization gives
	\[
	\int_{n_x\cdot v<0}d_\delta\sqrt\mu|n_x\cdot v|\dd v
	=\int_{n_x\cdot v<0}M_{\Theta_w}|n_x\cdot v|\dd v
	-c_\mu^{-1}\int_{n_x\cdot v<0}\mu|n_x\cdot v|\dd v=0,
	\]
	which is \eqref{eq:wall-defect-flux-orthogonality}.
	Recall from \eqref{eq:dyn-theta-wall-split} that
	$r_\delta[f_{\gamma^+}]=d_\delta\cA_0f_{\gamma^+}$, where
	$\cA_0f_{\gamma^+}(x)$ does not depend on $v$.  Integrating first in $v$ and
	using \eqref{eq:wall-defect-flux-bound} and
	\eqref{eq:outgoing-moment-flux-bound}, we obtain
	\begin{align*}
		\|r_\delta[f_{\gamma^+}]\|_{L^2(\gamma^-)}^2
		&=\int_{\partial\cD}|\cA_0f_{\gamma^+}(x)|^2
		\int_{n_x\cdot v<0}|d_\delta(x,v)|^2|n_x\cdot v|\dd v\dd S_x\\
		&\le C\delta^2\|\cA_0f_{\gamma^+}\|_{L^2(\partial\cD)}^2
		\le Cc_\mu\delta^2\|f_{\gamma^+}\|_{L^2(\gamma^+)}^2.
	\end{align*}
	Taking square roots gives the first bound in
	\eqref{eq:wall-defect-operator-bound}.  By \eqref{eq:wall-response-profiles}
	and \eqref{eq:isothermal-diffuse-projections},
	$\cP_0^-f_{\gamma^+}=c_\mu^{-1}\sqrt\mu\,\cA_0f_{\gamma^+}$, and the same
	integration gives
	\[
	\langle\cP_0^-f_{\gamma^+},r_\delta[f_{\gamma^+}]\rangle_{L^2(\gamma^-)}
	=c_\mu^{-1}\int_{\partial\cD}|\cA_0f_{\gamma^+}(x)|^2
	\int_{n_x\cdot v<0}d_\delta(x,v)\sqrt\mu(v)|n_x\cdot v|\dd v\dd S_x
	=0
	\]
	by \eqref{eq:wall-defect-flux-orthogonality}; this is the second formula of
	\eqref{eq:wall-defect-operator-bound}.  By this orthogonality and \eqref{eq:dyn-theta-wall-split},
	\[
	\|P_wf_{\gamma^+}\|_{L^2(\gamma^-)}^2
	=\|\cP_0^-f_{\gamma^+}\|_{L^2(\gamma^-)}^2
	+\|r_\delta[f_{\gamma^+}]\|_{L^2(\gamma^-)}^2.
	\]
	
	Finally, by rotation invariance set
	\[
	J_{\eta,R}:=
	\int_{\substack{v_3>0:\ v_3<\eta\ \text{or}\ |v|>R}}
	\mu(v)v_3\dd v,
	\qquad
	\theta_{\eta,R}:=\frac{J_{\eta,R}}{c_\mu}.
	\]
	Since $\cP_0^+f_{\gamma^+}=c_\mu^{-1}\sqrt\mu\,\cA_0f_{\gamma^+}$,
	\[
	\|\cP_0^+f_{\gamma^+}\|_{L^2(\gamma^+\setminus\gamma_{\eta,R}^+)}^2
	=c_\mu^{-2}J_{\eta,R}\|\cA_0f_{\gamma^+}\|_{L^2(\partial\cD)}^2,
	\qquad
	\|\cP_0^+f_{\gamma^+}\|_{L^2(\gamma^+)}^2
	=c_\mu^{-1}\|\cA_0f_{\gamma^+}\|_{L^2(\partial\cD)}^2.
	\]
	Hence \eqref{eq:rank-one-outgoing-tail} holds, and dominated convergence
	gives the stated limit.
\end{proof}

The exact and truncated wall operators have the common form
\begin{equation}
	\label{eq:wall-contraction-operator}
	P_\cJ:=\mathfrak m_{\Theta_w}\cJ\cA_0
	:L^2(\gamma^+)\longrightarrow L^2(\gamma^-),
\end{equation}
where $\cJ$ is self-adjoint on $L^2(\p\cD)$ with $0\le\cJ\le\mathbf I$ and
$\cJ1=1$: the choice $\cJ=\mathbf I$ gives $P_w$ of
\eqref{eq:nonisothermal-diffuse-response} and $\cJ=J_N$ gives $P_{w,N}$ of
\eqref{eq:boundary-spectral-diffuse-operator}.  

\begin{lemma}[Boundary form of the reflection law]
	\label{lem:wall-boundary-form}
	Let $\cJ$ be as in \eqref{eq:wall-contraction-operator} and let
	$f_{\gamma^+}\in L^2(\gamma^+,\dd\gamma)$. 
	Write 
	$$
	\varphi:=\cA_0f_{\gamma^+},\quad 	 p_\cJ:=\mathfrak m_\circ\cJ\varphi,
	\quad
	r_\cJ:=d_\delta\cJ\varphi, 
	$$
	so that $P_\cJ f_{\gamma^+}=p_\cJ+r_\cJ$.
	Also denote
	\[
	B:=\|f_{\gamma^+}\|_{L^2(\gamma^+)},
	\quad
	M:=\|(\mathbf I-\cP_0^+)f_{\gamma^+}\|_{L^2(\gamma^+)},
	\]
	and 
	\begin{equation}
		\label{eq:wall-spectral-defect}
		\mathfrak d_\cJ^2
		:=c_\mu^{-1}\bigl(
		\|\varphi\|_{L^2(\p\cD)}^2-\|\cJ\varphi\|_{L^2(\p\cD)}^2\bigr)
		\ \ge0 .
	\end{equation}
	\begin{enumerate}
		\item\label{it:wall-J-components} The two components satisfy
		\begin{equation}
			\label{eq:wall-J-components}
			\|p_\cJ\|_{L^2(\gamma^-)}^2=B^2-M^2-\mathfrak d_\cJ^2,
			\qquad
			\langle p_\cJ,r_\cJ\rangle_{L^2(\gamma^-)}=0,
			\qquad
			\|r_\cJ\|_{L^2(\gamma^-)}\le C\delta B ;
		\end{equation}
		in particular, for every $h_{\gamma^+}\in L^2(\gamma^+)$,
		\begin{equation}
			\label{eq:wall-J-operator-bound}
			\|P_\cJ h_{\gamma^+}\|_{L^2(\gamma^-)}\le C\|h_{\gamma^+}\|_{L^2(\gamma^+)} ,
		\end{equation}
		and, for every finite $L\ge0$,
		\begin{equation}
			\label{eq:wall-J-weighted-operator-bound}
			\|w^LP_\cJ h_{\gamma^+}\|_{L^2(\gamma^-)}
			\le C_L\|h_{\gamma^+}\|_{L^2(\gamma^+)} ,
		\end{equation}
		the right-hand side carrying no weight.
		
		\item\label{it:wall-J-continuation} If $f_{\gamma^-}=\alpha P_\cJ f_{\gamma^+}+b$ with $0\le\alpha\le1$ and
		$b\in L^2(\gamma^-,\dd\gamma)$, then, with
		\begin{equation}
			\label{eq:wall-continuation-defect}
			K^2:=(1-\alpha^2)B^2+\alpha^2\mathfrak d_\cJ^2,
		\end{equation}
		the boundary quadratic form has the exact expansion
		\begin{equation}
			\label{eq:wall-J-exact}
			B^2-\|f_{\gamma^-}\|_{L^2(\gamma^-)}^2
			=K^2+\alpha^2M^2-\alpha^2\|r_\cJ\|_{L^2(\gamma^-)}^2
			-\|b\|_{L^2(\gamma^-)}^2
			-2\alpha\langle p_\cJ+r_\cJ,b\rangle_{L^2(\gamma^-)} .
		\end{equation}
		
		\item\label{it:wall-J-lower} Under the hypothesis of (\ref{it:wall-J-continuation}), for every $0<\theta<1$,
		\begin{equation}
			\label{eq:wall-J-lower}
			B^2-\|f_{\gamma^-}\|_{L^2(\gamma^-)}^2
			\ge\tfrac12K^2+\tfrac12M^2
			-(\theta+C\delta^2)B^2-C\theta^{-1}\|b\|_{L^2(\gamma^-)}^2 .
		\end{equation}
	\end{enumerate}
	The constants $C$ depend only on $\cD$, and $C_L$ also on $L$; in particular, they do not depend on $\delta$, $\vartheta$ or $\cJ$.
\end{lemma}

\begin{proof}
	Write $\gamma_x^\pm:=\{v:\pm n_x\cdot v>0\}$ with the fiber measure
	$|n_x\cdot v|\dd v$ of \eqref{def:dgamma}.  The fiber normalization in
	\eqref{eq:wall-response-profiles} and
	\eqref{eq:isothermal-diffuse-projections} gives, for a.e.
	$x\in\p\cD$,
	\begin{equation}
		\label{eq:wall-rank-one-isometry}
		\|\mathfrak m_\circ\psi\|_{L^2(\gamma_x^-)}^2=c_\mu^{-1}|\psi(x)|^2,
		\qquad
		\|\cP_0^+f_{\gamma^+}(x,\cdot)\|_{L^2(\gamma_x^+)}^2=c_\mu^{-1}|\varphi(x)|^2 .
	\end{equation}
	Integrating in $x$ and using \eqref{eq:exact-diffuse-variance} and
	\eqref{eq:wall-spectral-defect},
	\[
	\|p_\cJ\|_{L^2(\gamma^-)}^2
	=c_\mu^{-1}\|\cJ\varphi\|_{L^2(\p\cD)}^2
	=B^2-M^2-\mathfrak d_\cJ^2 .
	\]
	The orthogonality is \eqref{eq:wall-defect-flux-orthogonality} applied
	fiberwise, and \eqref{eq:wall-defect-flux-bound} with $\|\cJ\|\le1$ gives
	\[
	\|r_\cJ\|_{L^2(\gamma^-)}^2\le C\delta^2\|\cJ\varphi\|_{L^2(\p\cD)}^2\le C\delta^2B^2.
	\]
	The component bounds give \eqref{eq:wall-J-operator-bound}
	for arbitrary $h_{\gamma^+}$.
	For \eqref{eq:wall-J-weighted-operator-bound}, write
	$P_\cJ=\mathfrak m_{\Theta_w}\cJ\cA_0$ with
	$\mathfrak m_{\Theta_w}=\mathfrak m_\circ+d_\delta$.  Both profiles are Gaussian
	in $v$ after multiplication by any finite $w^L$, by
	\eqref{eq:wall-profile-chart-bound} for $d_\delta$ and by its definition for
	$\mathfrak m_\circ$, so
	\[
	\sup_{x\in\p\cD}\int_{n_x\cdot v<0}
	w^{2L}|\mathfrak m_{\Theta_w}|^2|n_x\cdot v|\,\dd v\le C_L .
	\]
	Fiberwise Cauchy--Schwarz, $\|\cJ\|\le1$ on $L^2(\p\cD)$ and the flux bound
	\eqref{eq:outgoing-moment-flux-bound} for $\cA_0$ then give
	\eqref{eq:wall-J-weighted-operator-bound}.  
	
	Expanding the boundary condition and using
	\eqref{eq:wall-J-components} gives
	\[
	\|f_{\gamma^-}\|_{L^2(\gamma^-)}^2
	=\alpha^2(B^2-M^2-\mathfrak d_\cJ^2)
	+\alpha^2\|r_\cJ\|_{L^2(\gamma^-)}^2
	+\|b\|_{L^2(\gamma^-)}^2
	+2\alpha\langle p_\cJ+r_\cJ,b\rangle_{L^2(\gamma^-)} ;
	\]
	subtracting from $B^2$ and inserting
	\eqref{eq:wall-continuation-defect} gives \eqref{eq:wall-J-exact}.
	
	For (\ref{it:wall-J-lower}), $M\le B$ and $\mathfrak d_\cJ^2\ge0$ give
	\[
	K^2+\alpha^2M^2-\tfrac12K^2-\tfrac12M^2
	=\tfrac12(1-\alpha^2)(B^2-M^2)
	+\tfrac12\alpha^2(M^2+\mathfrak d_\cJ^2)\ge0 .
	\]
	By \eqref{eq:wall-J-components}, $\|p_\cJ+r_\cJ\|_{L^2(\gamma^-)}\le CB$, so
	Young's inequality gives
	\[
	2\alpha|\langle p_\cJ+r_\cJ,b\rangle_{L^2(\gamma^-)}|\le\theta B^2+C_\theta\|b\|_{L^2(\gamma^-)}^2.
	\]
	With
	$\alpha^2\|r_\cJ\|_{L^2(\gamma^-)}^2\le C\delta^2B^2$, substitution into
	\eqref{eq:wall-J-exact} proves \eqref{eq:wall-J-lower}.
\end{proof}

Let $I=(t_0,t_1)\subset(0,T)$ and set
\begin{equation}
	\label{eq:dynamic-I-boundary-cylinders}
	\cD_I:=I\times\cD,
	\qquad
	\gamma_I^\pm:=I\times\gamma^\pm.
\end{equation}
Recall that $f\in\cH_{\kin,\sigma}(\cD_I\times\bR_v^3)$ satisfies
\begin{equation*}
	f\in C(\overline I;L^2(\cD\times\bR_v^3))\cap L^2(\cD_I;H_\sigma^1),
	\qquad
	(\partial_t+v\cdot\nabla_x)f\in L^2(\cD_I;H_\sigma^{-1}),
	\qquad
	f_{\gamma^\pm}\in L^2(\gamma_I^\pm),
\end{equation*}
and that for every $\psi\in\cH_{\kin,\sigma}(\cD_I\times\bR_v^3)$, the
dynamic Green identity holds
\begin{align}
	&\int_I\!\langle (\partial_t+v\cdot\nabla_x)f,\psi\rangle_{x,v}\dd t
	+\int_I\!\langle f,(\partial_t+v\cdot\nabla_x)\psi\rangle_{x,v}\dd t
	\notag\\
	=\,& \la f(t_1),\psi(t_1)\ra_{x,v}-\la f(t_0),\psi(t_0)\ra_{x,v}
	+\int_{\gamma_I^+}f_{\gamma^+}\psi\dd\gamma\dd t
	-\int_{\gamma_I^-}f_{\gamma^-}\psi\dd\gamma\dd t.
	\label{eq:dynamic-Green-identity}
\end{align}
Its quadratic form is
\begin{align}
	&\frac12\|f(t_1)\|_{L^2(\cD\times\bR_v^3)}^2-
	\frac12\|f(t_0)\|_{L^2(\cD\times\bR_v^3)}^2
	+\frac12\|f_{\gamma^+}\|_{L^2(\gamma_I^+)}^2-\frac12\|f_{\gamma^-}\|_{L^2(\gamma_I^-)}^2
	\notag\\
	=\,&\int_I
	\langle(\partial_t+v\cdot\nabla_x)f,f\rangle_{x,v}\dd t.
	\label{eq:dynamic-quadratic-Green-identity}
\end{align}
For time-independent tests, \eqref{eq:dynamic-Green-identity} reduces to the
stationary Green identity up to the endpoint terms.
For $0\le\alpha\le1$ and $\cJ$ as in
\eqref{eq:wall-contraction-operator}, consider
\begin{equation}
	\label{eq:dynamic-admissible-boundary-law}
	f_{\gamma^-}(t)=\alpha P_\cJ f_{\gamma^+}(t)+b(t)
	\qquad\text{for a.e. }t\in(0,T).
\end{equation}
The identities and bounds of Lemma~\ref{lem:wall-boundary-form}
extend to $\gamma_I^\pm$ by integration in time.

\begin{prop}[Bounded-velocity non-grazing trace]
	\label{prop:bounded-nongrazing-trace}
	Let $0<\eta<1<R<\infty$ and let $I=(t_0,t_1)$ be as in
	\eqref{eq:dynamic-I-boundary-cylinders}.  There is a constant $C_{\eta,R}$, depending only on
	$\eta$, $R$, and the fixed boundary geometry, with the following properties.
	\begin{enumerate}
		\item\label{it:nongrazing-stationary} Every $f\in\cH_{\kin,\sigma}(\cD\times\bR_v^3)$
		satisfies
		\begin{equation}
			\label{eq:bounded-nongrazing-trace-estimate}
			\|f_{\gamma^\pm}\|_{L^2(\gamma_{\eta,R}^\pm)}
			\le C_{\eta,R}\left(
			\|f\|_{L^2(\cD;H_\sigma^1)}
			+\|v\cdot\nabla_xf\|_{L^2(\cD;H_\sigma^{-1})}
			\right).
		\end{equation}
		
		\item\label{it:nongrazing-evolutionary} For every $f\in\cH_{\kin,\sigma}(\cD_I\times\bR_v^3)$
		there are functions $\mathfrak e_f^\pm\in C([t_0,t_1])$ with
		\begin{equation}
			\label{eq:nongrazing-endpoint-bound}
			|\mathfrak e_f^\pm(t)|
			\le C_{\eta,R}\|f(t)\|_{L^2(\cD\times\bR_v^3)}^2,
			\qquad t\in\overline I,
		\end{equation}
		such that
		\begin{align}
			\|f_{\gamma^\pm}\|_{L^2(I\times\gamma_{\eta,R}^\pm)}^2
			\le{}&C_{\eta,R}\left(
			\|f\|_{L^2(\cD_I;H_\sigma^1)}^2
			+\|(\partial_t+v\cdot\nabla_x)f\|_{L^2(\cD_I;H_\sigma^{-1})}^2
			\right)
			\notag\\
			&+\mathfrak e_f^\pm(t_0)-\mathfrak e_f^\pm(t_1).
			\label{eq:dynamic-bounded-nongrazing-trace-estimate}
		\end{align}
	\end{enumerate}
	Here $C_{\eta,R}$ does not depend on $I$.
\end{prop}

\begin{proof}
	We first prove the evolutionary assertion.
	
	\noindent\emph{Step 1. Localization at the boundary.}
	We prove the estimate on a compact boundary patch for one sign, then sum over a finite cover of $\partial\cD$, with the endpoint functions
	summed.  Choose relatively compact boundary patches
	$U_0\subset\subset U_1\subset\subset\partial\cD$ and consider
	\[
	\Gamma_0:=\{(x,v):x\in U_0,\ |v|\le R,
	\ \pm n_x\cdot v\ge\eta\}.
	\]
	After shrinking $U_1$,
	\[
	\sup_{x,y\in U_1}|n_x-n_y|\le\frac{\eta}{8R}.
	\]
	Hence every velocity occurring in $\Gamma_0$ satisfies
	$\pm n_y\cdot v\ge7\eta/8$ for all $y\in U_1$.  We may therefore choose
	$\omega\in C_c^\infty(\bR_v^3)$, $0\le\omega\le1$, such that $\omega=1$ on the
	velocity projection of $\Gamma_0$ and
	\[
	\operatorname{supp}\omega\subset\{ |v|\le2R:\ \pm n_x\cdot v\ge\eta/2
	\text{ for every }x\in U_1\}.
	\]
	Choose a tubular-collar cutoff $\zeta\in C_c^\infty$ equal to one on a collar
	of $U_0$ and supported over $U_1$, and put $\psi:=\zeta^2$, so that
	$0\le\psi\omega^2\le1$ and $\psi\omega^2=1$ on $\Gamma_0$.  Neither cutoff depends
	on $t$.
	
	\noindent\emph{Step 2. The localized Green identity.}
	First let $f\in C_c^\infty(\overline I\times\overline\cD\times\bR_v^3)$.  Then
	\[
	\partial_t|f|^2+v\cdot\nabla_x|f|^2
	=2\bigl\langle(\partial_t+v\cdot\nabla_x)f,f\bigr\rangle .
	\]
	Multiply by $\psi\omega^2$ and apply the divergence theorem in $(t,x)$ on
	$I\times\cD$:
	\begin{align}
		&\int_I\int_{\partial\cD\times\bR^3}
		\psi\omega^2|f|^2\,(v\cdot n_x)\dd S_x\dd v\dd t
		+\Bigl[\int_{\cD\times\bR^3}\psi\omega^2|f(t)|^2\dd x\dd v\Bigr]_{t=t_0}^{t=t_1}
		\notag\\
		=\,&
		\int_I\int_{\cD\times\bR^3}
		(v\cdot\nabla_x\psi)\omega^2|f|^2\dd x\dd v\dd t
		+2\int_I\int_{\cD}
		\bigl\langle(\partial_t+v\cdot\nabla_x)f,
		\psi\omega^2f\bigr\rangle_{H_\sigma^{-1},H_\sigma^1}\dd x\dd t.
		\label{eq:nongrazing-inner-green}
	\end{align}
	Every term of \eqref{eq:nongrazing-inner-green} is continuous for the norm
	\eqref{innerprod_sigma}.  For the boundary term and the endpoint terms this continuity follows from
	\eqref{260730346}; for the first term on the right we use
	$\operatorname{supp}\omega\subset\{|v|\le2R\}$ and
	\eqref{eq:local-H1-H_sigma}.  For the last term we use
	$H_\sigma^{-1}$--$H_\sigma^1$ duality together with the multiplier bound
	\[
	\|\psi\omega^2f\|_{H_\sigma^1}\le C_R\|f\|_{H_\sigma^1},
	\]
	which follows from
	$\partial_{v_i}(\omega^2f)=\omega^2\partial_{v_i}f+(\partial_{v_i}\omega^2)f$,
	Lemma~\ref{lem:landau-diffusion-eigenvalues},
	\eqref{eq:H_sigma_norm} and \eqref{eq:local-H1-H_sigma}.
	Since the smooth functions are dense by Definition~\ref{def:H_kin_sigma},
	\eqref{eq:nongrazing-inner-green} holds for every
	$f\in\cH_{\kin,\sigma}(\cD_I\times\bR_v^3)$, with $f$ replaced by
	$f_{\gamma^\pm}$ in the boundary term.
	
	\noindent\emph{Step 3. The flux bound on the patch.}
	On $\operatorname{supp}(\zeta\omega)$ we have $\pm n_x\cdot v\ge\eta/2$, so the
	boundary integrand in \eqref{eq:nongrazing-inner-green} has one sign and
	dominates $|f_{\gamma^\pm}|^2|n_x\cdot v|$ on $\Gamma_0$, where
	$\psi\omega^2=1$.  Multiplying \eqref{eq:nongrazing-inner-green} by $\pm1$
	accordingly, and setting
	\begin{equation}
		\label{eq:nongrazing-endpoint-functional}
		\mathfrak e_f^\pm(t)
		:=\pm\int_{\cD\times\bR^3}\psi\omega^2|f(t)|^2\dd x\dd v,
	\end{equation}
	which satisfies \eqref{eq:nongrazing-endpoint-bound} because
	$0\le\psi\omega^2\le1$, we obtain
	\begin{align}
		\int_I\int_{\Gamma_0}|f_{\gamma^\pm}|^2|n_x\cdot v|\dd S_x\dd v\dd t
		\le{}&\Bigl|\int_I\int_{\cD\times\bR^3}
		(v\cdot\nabla_x\psi)\omega^2|f|^2\dd x\dd v\dd t\Bigr|
		\notag\\
		&+2\|(\partial_t+v\cdot\nabla_x)f\|_{L^2(\cD_I;H_\sigma^{-1})}
		\|\psi\omega^2f\|_{L^2(\cD_I;H_\sigma^1)}
		\notag\\
		&+\mathfrak e_f^\pm(t_0)-\mathfrak e_f^\pm(t_1).
		\label{eq:nongrazing-inner-flux-bound}
	\end{align}
	By $\operatorname{supp}\omega\subset\{|v|\le2R\}$ and
	\eqref{eq:local-H1-H_sigma}, the first term on the right of
	\eqref{eq:nongrazing-inner-flux-bound} is at most
	$C_{\eta,R}\|f\|_{L^2(\cD_I;H_\sigma^1)}^2$; by the multiplier bound of
	Step 2 and Young's inequality, so is the second, up to the additional summand
	$\|(\partial_t+v\cdot\nabla_x)f\|_{L^2(\cD_I;H_\sigma^{-1})}^2$.  A finite
	cover of $\partial\cD$ by patches $U_0$, with the endpoint functionals
	\eqref{eq:nongrazing-endpoint-functional} added, proves
	\eqref{eq:dynamic-bounded-nongrazing-trace-estimate}.  The cutoffs are time independent, so the constant does not depend on the length of $I$.
	
	\noindent\emph{Proof of (1).}  Use the same $\Gamma_0$, $\psi$, $\omega$.
	There is no time variable.  For
	$f\in C_c^\infty(\overline\cD\times\bR_v^3)$ we apply the divergence theorem
	in $x$ to $v\cdot\nabla_x|f|^2=2\langle v\cdot\nabla_xf,f\rangle$.  This
	gives \eqref{eq:nongrazing-inner-green} without $\int_I$ and without the
	endpoint bracket.  It extends to every
	$f\in\cH_{\kin,\sigma}(\cD\times\bR_v^3)$ by the density in
	Definition~\ref{def:H_kin_sigma}\,(2) and the stationary trace bound of
	Proposition~\ref{2607151218}.  Step 3 then gives, with
	$\mathfrak e_f^\pm\equiv0$,
	\[
	\int_{\Gamma_0}|f_{\gamma^\pm}|^2|n_x\cdot v|\dd S_x\dd v
	\le C_{\eta,R}\bigl(\|f\|_{L^2(\cD;H_\sigma^1)}^2
	+\|v\cdot\nabla_xf\|_{L^2(\cD;H_\sigma^{-1})}^2\bigr).
	\]
	As in (2), no separate $L^2(\cD\times\bR_v^3)$ norm of $f$ is needed: the term
	with $v\cdot\nabla_x\psi$ has velocity support in $\{|v|\le2R\}$, so
	\eqref{eq:local-H1-H_sigma} bounds it by
	$C_R\|f\|_{L^2(\cD;H_\sigma^1)}^2$.  A finite cover of $\p\cD$ and a square
	root give \eqref{eq:bounded-nongrazing-trace-estimate}.
\end{proof}

\begin{prop}[Recovery of the diffuse trace]
	\label{prop:full-diffuse-trace-recovery}
	There are $0<\eta_*<1<R_*<\infty$ and $C>0$, depending only on the boundary
	geometry, with the following properties.
	\begin{enumerate}
		\item\label{it:recovery-stationary} Every $f\in\cH_{\kin,\sigma}(\cD\times\bR_v^3)$
		satisfies
		\begin{equation}
			\label{eq:full-outgoing-trace-recovery}
			\|f_{\gamma^+}\|_{L^2(\gamma^+)}
			\le C\left(
			\|f\|_{L^2(\cD;H_\sigma^1)}
			+\|v\cdot\nabla_xf\|_{L^2(\cD;H_\sigma^{-1})}
			+\|(\mathbf I-\cP_0^+)f_{\gamma^+}\|_{L^2(\gamma^+)}
			\right).
		\end{equation}
		If $f_{\gamma^-}=P_wf_{\gamma^+}+b$ with
		$b\in L^2(\gamma^-,\dd\gamma)$, then
		\begin{equation}
			\label{eq:full-incoming-trace-recovery-exact}
			\|f_{\gamma^-}\|_{L^2(\gamma^-)}
			\le C\|f_{\gamma^+}\|_{L^2(\gamma^+)}+\|b\|_{L^2(\gamma^-)}.
		\end{equation}
		
		\item\label{it:recovery-evolutionary} Let $I=(t_0,t_1)$ and
		$f\in\cH_{\kin,\sigma}(\cD_I\times\bR_v^3)$.  Then there is a function
		$\mathfrak e_f$ on $\overline I$ with
		\begin{equation}
			\label{eq:recovery-endpoint-bound}
			0\le\mathfrak e_f(t)\le C\|f(t)\|_{L^2(\cD\times\bR_v^3)}^2,
			\qquad t\in\overline I,
		\end{equation}
		such that
		\begin{align}
			\|f_{\gamma^+}\|_{L^2(\gamma_I^+)}^2
			\le{}&C\Bigl(
			\|f\|_{L^2(\cD_I;H_\sigma^1)}^2
			+\|(\partial_t+v\cdot\nabla_x)f\|_{L^2(\cD_I;H_\sigma^{-1})}^2
			\notag\\
			&\quad
			+\|(\mathbf I-\cP_0^+)f_{\gamma^+}\|_{L^2(\gamma_I^+)}^2\Bigr)
			+\mathfrak e_f(t_0)-\mathfrak e_f(t_1).
			\label{eq:dynamic-full-outgoing-trace}
		\end{align}
		If in addition \eqref{eq:dynamic-admissible-boundary-law} holds, then
		\begin{equation}
			\label{eq:dynamic-full-incoming-trace}
			\|f_{\gamma^-}\|_{L^2(\gamma_I^-)}
			\le C\|f_{\gamma^+}\|_{L^2(\gamma_I^+)}+\|b\|_{L^2(\gamma_I^-)}.
		\end{equation}
		The constants are independent of $\alpha$, of $\cJ$, and of $I$.
	\end{enumerate}
\end{prop}

\begin{proof}
	Choose $\eta_*,R_*$ in Lemma~\ref{lem:diffuse-boundary-algebra} so that
	$\theta_{\eta_*,R_*}\le1/4$, and write
	$\gamma_*^\pm:=\gamma_{\eta_*,R_*}^\pm$.  
	We first prove the time-dependent estimate.
	
	Put $Pf:=\cP_0^+f_{\gamma^+}$.  By
	Lemma~\ref{lem:diffuse-boundary-algebra}\,(\ref{it:diffuse-P0-contraction}), integrated over $I$,
	\[
	\|Pf\|_{L^2(\gamma_I^+)}^2
	\le\|Pf\|_{L^2(I\times\gamma_*^+)}^2
	+\tfrac14\|Pf\|_{L^2(\gamma_I^+)}^2,
	\qquad\text{hence}\qquad
	\|Pf\|_{L^2(\gamma_I^+)}\le2\|Pf\|_{L^2(I\times\gamma_*^+)} .
	\]
	On $I\times\gamma_*^+$ we have $Pf=f_{\gamma^+}-(\mathbf I-\cP_0^+)f_{\gamma^+}$, so
	\[
	\|Pf\|_{L^2(\gamma_I^+)}
	\le2\|f_{\gamma^+}\|_{L^2(I\times\gamma_*^+)}
	+2\|(\mathbf I-\cP_0^+)f_{\gamma^+}\|_{L^2(\gamma_I^+)} .
	\]
	Square the preceding inequality.
	Proposition~\ref{prop:bounded-nongrazing-trace}\,(\ref{it:nongrazing-evolutionary}),
	applied with $(\eta,R)=(\eta_*,R_*)$ and the outgoing sign, bounds the
	square of the first term on the right by
	\[
	C\Bigl(\|f\|_{L^2(\cD_I;H_\sigma^1)}^2
	+\|(\partial_t+v\cdot\nabla_x)f\|_{L^2(\cD_I;H_\sigma^{-1})}^2\Bigr)
	+\mathfrak e_f(t_0)-\mathfrak e_f(t_1),
	\]
	where $\mathfrak e_f$ is a fixed multiple of the functional
	\eqref{eq:nongrazing-endpoint-functional} of that proposition, taken with the
	outgoing sign.  Here $0\le\psi\omega^2\le1$, so it is nonnegative and satisfies
	\eqref{eq:nongrazing-endpoint-bound}; together, these two properties give
	\eqref{eq:recovery-endpoint-bound}.  Since $\cP_0^+$ is the fiberwise
	orthogonal projection of Lemma~\ref{lem:diffuse-boundary-algebra}\,(\ref{it:diffuse-P0-projection}),
	\[
	\|f_{\gamma^+}\|_{L^2(\gamma_I^+)}^2
	=\|Pf\|_{L^2(\gamma_I^+)}^2
	+\|(\mathbf I-\cP_0^+)f_{\gamma^+}\|_{L^2(\gamma_I^+)}^2 ,
	\]
	and this, with the squared decomposition above and the bound just stated
	for its first term, proves \eqref{eq:dynamic-full-outgoing-trace}.
	
	Under \eqref{eq:dynamic-admissible-boundary-law},
	\eqref{eq:wall-J-operator-bound} and $\alpha\le1$ give
	\eqref{eq:dynamic-full-incoming-trace}.
	
	For (1), the same argument on $\gamma^+$ has no time integral or
	endpoint term and gives
	\[
	\|Pf\|_{L^2(\gamma^+)}
	\le2\|f_{\gamma^+}\|_{L^2(\gamma_*^+)}
	+2\|(\mathbf I-\cP_0^+)f_{\gamma^+}\|_{L^2(\gamma^+)} .
	\]
	Proposition~\ref{prop:bounded-nongrazing-trace}\,(\ref{it:nongrazing-stationary})
	with $(\eta,R)=(\eta_*,R_*)$ bounds the first term on the right. The triangle inequality
	\[
	\|f_{\gamma^+}\|_{L^2(\gamma^+)}\le\|Pf\|_{L^2(\gamma^+)}+\|(\mathbf I-\cP_0^+)f_{\gamma^+}\|_{L^2(\gamma^+)}
	\]
	gives
	\eqref{eq:full-outgoing-trace-recovery}, and
	\eqref{eq:wall-J-operator-bound} with $\cJ=\mathbf I$ gives
	\eqref{eq:full-incoming-trace-recovery-exact}.
\end{proof}

\begin{lemma}[Weak closure in $\cH_{\kin,\sigma}$]
	\label{lem:dynamic-weak-graph-bulk-identification}\,
	
	\begin{enumerate}
		\item
		Let $(f_k)$ be bounded in
		$\cH_{\kin,\sigma}(\cD\times\bR_v^3)$ and assume
		\begin{equation}
			\label{eq:stationary-weak-graph-bulk-distribution-limit}
			f_k\longrightarrow f
			\quad\text{in }\cD'(\cD\times\bR_v^3).
		\end{equation}
		Then $f$ is the bulk coordinate of a unique element of
		$\cH_{\kin,\sigma}(\cD\times\bR_v^3)$, still denoted by $f$,
		$f_k\rightharpoonup f$ in that space, and
		\begin{align}
			v\cdot\nabla_xf_k&\rightharpoonup v\cdot\nabla_xf
			&&\text{in }L^2(\cD;H_\sigma^{-1}),
			\notag\\
			(f_k)_{\gamma^\pm}&\rightharpoonup f_{\gamma^\pm}
			&&\text{in }L^2(\gamma^\pm).
			\label{eq:stationary-weak-graph-coordinate-convergence}
		\end{align}
		If moreover
		$(f_k)_{\gamma^-}=P_w(f_k)_{\gamma^+}+b_k$ with
		$b_k\rightharpoonup b$ in $L^2(\gamma^-)$, then
		$f_{\gamma^-}=P_wf_{\gamma^+}+b$.
		
		\item\label{it:closure-evolutionary}
		Let $(f_k)$ be bounded in
		$\cH_{\kin,\sigma}(\cD_T\times\bR_v^3)$ and assume
		\begin{equation}
			\label{eq:dynamic-weak-graph-bulk-distribution-limit}
			f_k\longrightarrow f
			\quad\text{in }\cD'((0,T)\times\cD\times\bR_v^3).
		\end{equation}
		Then $f$ is the bulk coordinate of a unique element of
		$\cH_{\kin,\sigma}(\cD_T\times\bR_v^3)$, still denoted by $f$, and
		\begin{equation}
			\label{eq:dynamic-weak-graph-full-convergence}
			f_k\rightharpoonup f\quad\text{in }\cH_{\kin,\sigma}(\cD_T\times\bR_v^3).
		\end{equation}
		Consequently,
		\begin{align}
			(\partial_t+v\cdot\nabla_x)f_k
			&\rightharpoonup(\partial_t+v\cdot\nabla_x)f
			&&\text{in }L^2(\cD_T;H_\sigma^{-1}),
			\notag\\
			f_k(t)&\rightharpoonup f(t)
			&&\text{in }L^2(\cD\times\bR_v^3),\qquad 0\le t\le T,
			\notag\\
			(f_k)_{\gamma^\pm}&\rightharpoonup f_{\gamma^\pm}
			&&\text{in }L^2(\gamma_T^\pm).
			\label{eq:dynamic-weak-graph-coordinate-convergence}
		\end{align}
		If moreover $(f_k)_{\gamma^-}=P_w(f_k)_{\gamma^+}+b_k$ with $b_k\rightharpoonup b$ in
		$L^2(\gamma_T^-)$, then $f_{\gamma^-}=P_wf_{\gamma^+}+b$.
	\end{enumerate}
\end{lemma}

\begin{proof}
	We prove (2), then indicate the changes for (1).
	Since $\cH_{\kin,\sigma}(\cD_T\times\bR_v^3)$ is a Hilbert space, every subsequence has a further subsequence
	converging weakly in it to some $\widetilde f$.  The bulk coordinate of $\widetilde f$ equals $f$ by
	\eqref{eq:dynamic-weak-graph-bulk-distribution-limit}.  Any two weak cluster
	points therefore have the same bulk coordinate and hence coincide by
	Proposition~\ref{H_kin=weak_landau}\,(1).  This proves
	\eqref{eq:dynamic-weak-graph-full-convergence}.  The coordinate maps, bounded by
	\eqref{260902137}, and the transport coordinate of
	Definition~\ref{def:H_kin_sigma} give
	\eqref{eq:dynamic-weak-graph-coordinate-convergence}.  For the boundary law,
	$P_w$ is bounded from $L^2(\gamma_T^+)$ to $L^2(\gamma_T^-)$ by
	\eqref{eq:wall-J-operator-bound}, hence weakly continuous; by the last line of
	\eqref{eq:dynamic-weak-graph-coordinate-convergence}, the
	relation $(f_k)_{\gamma^-}=P_w(f_k)_{\gamma^+}+b_k$ therefore passes to the limit.
	
	For (1) the argument is the same in the Hilbert space
	$\cH_{\kin,\sigma}(\cD\times\bR_v^3)$.  Uniqueness of the element with a
	given bulk coordinate is Proposition~\ref{H_kin=weak_landau}\,(2).  The
	coordinate maps $f\mapsto v\cdot\nabla_xf$ and $f\mapsto f_{\gamma^\pm}$ are
	bounded by the stationary line of \eqref{260902137}.  There is no time
	slice, so the middle line of
	\eqref{eq:dynamic-weak-graph-coordinate-convergence} has no stationary
	counterpart.  The boundary law passes to the limit because $P_w$ is bounded
	from $L^2(\gamma^+)$ to $L^2(\gamma^-)$ by \eqref{eq:wall-J-operator-bound}.
\end{proof}

\subsection{Macroscopic estimate for stationary problems}\label{ss:diffuse-macro}

We use the test-function method of \S\ref{subsec:macro-landau}
and \cite{Guo2010,EspositoGuoKimMarra} to recover the macroscopic fields.
The rigid-motion reduction follows the kernel structure in
Desvillettes--Villani~\cite[(1) and Definition--Lemma~1]{DesvillettesVillani2002}.
For specular reflection, macroscopic estimates based on the test-function method and Korn's inequality were obtained by Chen and Kim \cite{ChenKimMacroscopic2024}. Under diffuse reflection, the boundary dissipation provides additional control of the rigid-motion modes, as expressed in the following lemma.

\begin{lemma}[Boundary observability]
	\label{lem:boundary-rigid-observability}
	There exist a relatively open patch $U\subset\partial\cD$ and an open ball
	$V\subset\bR_v^3$ such that $n_x\cdot v\ge1/2$ on $U\times V$ and, for every
	$a_0,c_0\in\bR$, $\mathbf b_0\in\bR^3$, skew-symmetric
	$\Omega_0\in\bR^{3\times3}$, and $\beta\in L^2(U)$,
	\begin{equation}
		\label{eq:boundary-rigid-observability}
		|\mathbf b_0|+|\Omega_0|+|c_0|
		\le C\bigl\|
		\bigl(a_0+(\mathbf b_0+\Omega_0x)\cdot v+c_0(|v|^2-\tfrac32)\bigr)\sqrt\mu
		-\mathfrak m_\circ\beta\bigr\|_{L^2(U\times V,\dd\gamma)}.
	\end{equation}
	The constant depends only on the chosen patch of $\partial\cD$.
\end{lemma}

\begin{proof}
	We fix a non-grazing patch and establish coercivity in velocity and in the rigid-motion coefficients.  
	
	\noindent\emph{Step 1. The patch and the velocity window.}
	Fix $x_*\in\partial\cD$ and rotate coordinates so that $n_{x_*}=e_3$.  After
	shrinking to a nonempty relatively open patch $U$ with compact closure in one
	boundary chart, continuity of the normal gives
	\begin{equation*}
		n_x\cdot e_3\ge\frac34
		\qquad(x\in U).
	\end{equation*}
	Set $r_v:=1/4$ and $V:=B(e_3,r_v)$.  Then
	\[
	n_x\cdot v\ge n_x\cdot e_3-|v-e_3|\ge\frac12
	\qquad((x,v)\in U\times V).
	\]
	Let $m_{U,V}>0$ be a lower bound for the flux density
	$\mu(v)(n_x\cdot v)$ on $U\times V$.
	
	\noindent\emph{Step 2. The velocity Gram form.}
	Fix $p\in\bR^3$, $c\in\bR$, and write $v=e_3+y$ with $y\in B_{r_v}=B(0,r_v)$.  Since
	\[
	\frac{1}{|B_{r_v}|}\int_{B_{r_v}}|y|^2\dd y=\frac{3r_v^2}{5},
	\]
	the minimization over the constant $\lambda$ removes the mean and gives
	\begin{align*}
		&\inf_{\lambda\in\bR}
		\int_V|p\cdot v+c(|v|^2-\tfrac32)-\lambda|^2\dd v=
		\int_{B_{r_v}}\left|
		(p+2ce_3)\cdot y
		+c\left(|y|^2-\frac{3r_v^2}{5}\right)
		\right|^2\dd y.
	\end{align*}
	\begingroup
	The odd linear part and the centered radial part are orthogonal on $B_{r_v}$,
	and both have positive $L^2$ norm unless their coefficients vanish.
	Therefore, for some $\kappa_{r_v}>0$ depending only on $r_v$,
	\begin{equation}
		\label{eq:velocity-rigid-Gram}
		\inf_{\lambda\in\bR}\int_V
		|p\cdot v+c(|v|^2-\tfrac32)-\lambda|^2\dd v
		\ge\kappa_{r_v}(|p|^2+|c|^2).
	\end{equation}
	\endgroup
	
	\noindent\emph{Step 3. Integration over the patch.}
	For $x\in U$, set $p_x:=\mathbf b_0+\Omega_0x$ and
	$\lambda_x:=c_\mu^{-1}\beta(x)-a_0$.  Since
	$\mathfrak m_\circ=c_\mu^{-1}\sqrt\mu$, the function inside the norm in
	\eqref{eq:boundary-rigid-observability} equals
	\[
	\sqrt\mu(v)\left[
	p_x\cdot v+c_0(|v|^2-\tfrac32)-\lambda_x
	\right].
	\]
	Integrating \eqref{eq:velocity-rigid-Gram} over $U$ with the flux
	weight gives
	\begin{align}
		&\bigl\|
		\bigl(a_0+(\mathbf b_0+\Omega_0x)\cdot v+c_0(|v|^2-\tfrac32)\bigr)\sqrt\mu
		-\mathfrak m_\circ\beta\bigr\|_{L^2(U\times V,\dd\gamma)}^2\notag\\
		&\qquad\ge
		m_{U,V}\kappa_{r_v}\left(
		\int_U|\mathbf b_0+\Omega_0x|^2\dd S_x+|U||c_0|^2
		\right).
		\label{eq:boundary-rigid-velocity-observability}
	\end{align}
	
	\noindent\emph{Step 4. The spatial Gram form.}
	For the spatial rigid motion, on the six-dimensional space
	\[
	\cR_{\bR}
	:=\{(b,\Omega):b\in\bR^3,
	\ \Omega\in\bR^{3\times3},\ \Omega^T=-\Omega\},
	\]
	define the quadratic form
	\[
	\cQ_U(b,\Omega):=\int_U|b+\Omega x|^2\dd S_x.
	\]
	Its kernel is trivial.  Indeed, if $(b,\Omega)\ne(0,0)$ and $\Omega=0$ then
	$b+\Omega x$ does not vanish, while if $\Omega\ne0$ then $\Omega$ has rank
	two and the solution set of $b+\Omega x=0$ is either empty or an affine line.
	In either case the zero set of $b+\Omega x$ meets the two-dimensional
	relatively open patch $U$ in a set of surface measure zero.  Hence the smallest eigenvalue of $\cQ_U$ is
	positive, and there is $\lambda_U>0$ such that
	\[
	\cQ_U(b,\Omega)
	\ge\lambda_U\bigl(|b|^2+|\Omega|^2\bigr)
	\qquad((b,\Omega)\in\cR_{\bR}).
	\]
	Combining \eqref{eq:boundary-rigid-velocity-observability} with the preceding
	lower bound, and using equivalence of Euclidean norms in the
	finite-dimensional coefficient space, proves
	\eqref{eq:boundary-rigid-observability}.
\end{proof}

For $k=1,2$ let $H_0^k(\cD)$ be the closure of $C_c^\infty(\cD)$ in
$H^k(\cD)$ and let $H^{-k}(\cD):=(H_0^k(\cD))'$, with
\[
\|q\|_{H^{-k}(\cD)}
:=\sup\bigl\{\langle q,\varphi\rangle:
\varphi\in C_c^\infty(\cD),\ \|\varphi\|_{H^k(\cD)}\le1\bigr\}.
\]

\begin{lemma}[Negative-order estimates on a bounded Lipschitz domain]
	\label{lem:negative-order-estimates}
	Let $\cD\subset\bR^3$ be bounded, Lipschitz and connected, and let
	\[
	\cR_{\rm rig}
	:=\{x\mapsto z_0+\Omega x:\ z_0\in\bR^3,
	\ \Omega^T=-\Omega\}.
	\]
	\begin{enumerate}
		\item\label{it:necas} \emph{(Ne\v{c}as)}  For $s\in\{-2,-1\}$ and
		$q\in H^{s+1}(\cD)$,
		\[
		\|q\|_{H^{s+1}(\cD)}
		\le C_\cD\bigl(
		\|q\|_{H^s(\cD)}+\|\nabla q\|_{H^s(\cD)}\bigr).
		\]
		
		\item\label{it:neg-poincare} \emph{(Mean-zero Poincar\'e)}  For every
		$q\in L^2(\cD)$,
		\[
		\|q-\fint_{\cD} q\dd x\|_{L^2(\cD)}
		\le C_\cD\|\nabla q\|_{H^{-1}(\cD)}.
		\]
		
		\item\label{it:neg-korn} \emph{(Korn rigidity)}  For every
		$z\in L^2(\cD;\bR^3)$,
		\[
		\inf_{r\in\cR_{\rm rig}}
		\|z-r\|_{L^2(\cD)}
		\le C_\cD\|\operatorname{sym}\nabla z\|_{H^{-1}(\cD)}.
		\]
	\end{enumerate}
\end{lemma}

\begin{proof}
	The first assertion is the integer-order Ne\v{c}as norm equivalence on
	bounded Lipschitz domains; see
	Amrouche--Girault~\cite[Theorem~2.3]{AmroucheGirault1994} with
	$m=-1$ and $m=0$.
	
	The remaining estimates follow by the same compactness argument.
	The kernels of $\nabla$ and $\operatorname{sym}\nabla$ are, respectively,
	the constants and the rigid motions.  Suppose
	the asserted estimate fails.  Then there are $u_j$ orthogonal to the kernel in
	$L^2(\cD)$ with $\|u_j\|_{L^2(\cD)}=1$ and the operator applied to $u_j$
	tending to zero in $H^{-1}(\cD)$.  Since
	$L^2(\cD)\hookrightarrow H^{-1}(\cD)$ is compact, after a subsequence
	$u_j\rightharpoonup u$ in $L^2$ and $u_j\to u$ in $H^{-1}$, so the operator
	annihilates $u$.  In (\ref{it:neg-poincare}) connectedness gives $u\in\bR$; in
	(\ref{it:neg-korn}) the distributional identity
	\[
	\p_j\p_kz_i
	=\p_je_{ik}(z)+\p_ke_{ij}(z)-\p_ie_{jk}(z),
	\qquad
	e_{ij}(z):=\tfrac12(\p_i z_j+\p_jz_i),
	\]
	gives $D^2u=0$, hence $u(x)=z_0+Bx$, and $e(u)=0$ forces $B^T=-B$, so
	$u\in\cR_{\rm rig}$.  In both cases orthogonality gives $u=0$, hence
	$u_j\to0$ in $H^{-1}$.  (\ref{it:necas}) of this lemma closes the argument: with
	$s=-1$ it gives
	\[
	1=\|u_j\|_{L^2(\cD)}
	\le C\bigl(\|u_j\|_{H^{-1}}+\|\nabla u_j\|_{H^{-1}}\bigr)
	\longrightarrow0
	\]
	in (\ref{it:neg-poincare}), while in (\ref{it:neg-korn}) applying it first with
	$s=-2$ to each $\p_kz_i$ and then with $s=-1$ to each $z_i$ through the
	second-derivative identity above gives
	\[
	\|z\|_{L^2}
	\le C\bigl(\|e(z)\|_{H^{-1}}+\|z\|_{H^{-1}}\bigr),
	\]
	whence $\|u_j\|_{L^2(\cD)}\to0$.  In either case this is a contradiction.
\end{proof}

The bulk mass and the centered mass mode are
\begin{equation}
	\label{eq:dynamic-mass-scalar}
	\mathfrak m_f:=
	\int_{\cD\times\bR_v^3}f(x,v)\sqrt\mu(v)\dd v\dd x,
	\qquad
	\bar{\mathbf a}:=\pi^{-3/2}|\cD|^{-1}\mathfrak m_f,
	\qquad
	\mathbf a^\circ:=\mathbf a-\bar{\mathbf a} ;
\end{equation}
for an evolutionary $f$ they are formed from $f(t)$ at each time.

\medskip
The stationary $H_\sigma^1$--$H_\sigma^{-1}$ estimate below follows
the diffuse macroscopic rigidity argument of
Esposito--Guo--Kim--Marra~\cite{EspositoGuoKimMarra2018}.  The estimate uses the outgoing variance instead of the full trace.

\begin{prop}[Macroscopic recovery]
	\label{prop:macroscopic-recovery}
	Let $f\in\cH_{\kin,\sigma}(\cD\times\bR_v^3)$ satisfy
	\[
	v\cdot\nabla_xf+Lf=S
	\qquad\text{in }L^2(\cD;H_\sigma^{-1}).
	\]
	No relation between $f_{\gamma^-}$ and $f_{\gamma^+}$ is assumed.  Then the
	centered macroscopic fields are bounded by
	\begin{equation}
		\label{eq:stationary-macro-estimate}
		\|\mathbf a^\circ\|_{L_x^2}^2
		+\|\mathbf b\|_{L_x^2}^2
		+\|\mathbf c\|_{L_x^2}^2
		\le C\bigl(\|(\mathbf I-\mathbf P)f\|_{L^2(\cD;H_\sigma^1)}^2+\|S\|_{L^2(\cD;H_\sigma^{-1})}^2+\|(\mathbf I-\cP_0^+)f_{\gamma^+}\|_{L^2(\gamma^+)}^2\bigr),
	\end{equation}
	and hence,
	\begin{equation}
		\label{eq:macroscopic-recovery-estimate}
		\|\mathbf P f\|_{L^2(\cD\times\bR_v^3)}
		\le C\bigl(X+D+M+|\mathfrak m_f|\bigr),
	\end{equation}
	where $\mathfrak m_f$ is defined in \eqref{eq:dynamic-mass-scalar} and
	$X$, $D$, $M$ denote the three norms on the right of
	\eqref{eq:stationary-macro-estimate}:
	\[
	X:=\|(\mathbf I-\mathbf P)f\|_{L^2(\cD;H_\sigma^1)},
	\qquad
	D:=\|S\|_{L^2(\cD;H_\sigma^{-1})},
	\qquad
	M:=\|(\mathbf I-\cP_0^+)f_{\gamma^+}\|_{L^2(\gamma^+)} .
	\]
	The constant depends only on $\cD$.
\end{prop}

\begin{proof}
	Use the earlier \hyperlink{comp:macro-projection}{projection formula} and \hyperlink{comp:macro-moments}{Gaussian moment formulas},
	with $X,D,M$ as in the statement.  Since $L\mathbf P f=0$, the distributional identity
	\begin{equation*}
		v\cdot\nabla_x\mathbf P f
		=S-L(\mathbf I-\mathbf P)f-v\cdot\nabla_x(\mathbf I-\mathbf P)f
	\end{equation*}
	holds.  For every fixed polynomial--Gaussian test $\psi$,
	\[
	\left\|\langle S-L(\mathbf I-\mathbf P)f,\psi\ra_v
	-\nabla_x\cdot\langle v(\mathbf I-\mathbf P)f,\psi\ra_v\right\|_{H^{-1}_x}
	\le C_\psi(X+D).
	\]
	Indeed, the first two terms are controlled by $H_\sigma^{-1}$--$H_\sigma^1$
	duality, whereas
	\[
	\|\langle v(\mathbf I-\mathbf P)f,\psi\ra_v\|_{L^2_x}
	\le
	\|\langle v\rangle^{-1/2}(\mathbf I-\mathbf P)f\|_{L^2(\cD\times\bR_v^3)}
	\|\langle v\rangle^{1/2}v\psi\|_{L^2_v}
	\le C_\psi\|(\mathbf I-\mathbf P)f\|_{L^2(\cD;H_\sigma^1)}.
	\]
	Set
	\[
	\mathfrak a:=\nabla \mathbf a-\tfrac32\nabla \mathbf c,
	\qquad \mathfrak b:=\operatorname{sym}\nabla\mathbf b,
	\qquad \mathfrak c:=\nabla \mathbf c.
	\]
	Then
	\[
	v\cdot\nabla_x\mathbf P f
	=\bigl(v\cdot\mathfrak a+v^T\mathfrak bv+|v|^2v\cdot\mathfrak c\bigr)\sqrt\mu.
	\]
	Testing with
	\begin{equation*}
		\psi_i^c=(|v|^2-\tfrac52)v_i\sqrt\mu,
		\quad
		\psi_{ik}^b=(v_iv_k-\tfrac12\delta_{ik})\sqrt\mu,
		\quad
		\psi_i^a=(\tfrac72-|v|^2)v_i\sqrt\mu
	\end{equation*}
	and using the centered Gaussian moments gives, with no cross terms,
	\begin{align*}
		\langle v\cdot\nabla_x\mathbf P f,\psi_i^c\ra_v
		&=\tfrac54\pi^{3/2}\mathfrak c_i,\\
		\langle v\cdot\nabla_x\mathbf P f,\psi_{ik}^b\ra_v
		&=\tfrac12\pi^{3/2}\mathfrak b_{ik},\\
		\langle v\cdot\nabla_x\mathbf P f,\psi_i^a\ra_v
		&=\tfrac12\pi^{3/2}\mathfrak a_i.
	\end{align*}
	The coefficients follow from
	$$
	\int v_rv_i\mu\dd v=\tfrac12\pi^{3/2}\delta_{ri}\,,\,\,
	\int |v|^2v_rv_i\mu\dd v=\tfrac54\pi^{3/2}\delta_{ri}\,,\,\,
	\int |v|^4v_rv_i\mu\dd v=\tfrac{35}8\pi^{3/2}\delta_{ri},
	$$
	and Wick's formula for
	$\int v_rv_sv_iv_k\mu\dd v$.  The Gaussian identities and the moment-remainder bound yield
	\begin{equation*}
		\|\nabla \mathbf c\|_{H^{-1}}
		+\|\operatorname{sym}\nabla\mathbf b\|_{H^{-1}}
		+\|\nabla(\mathbf a-\tfrac32\mathbf c)\|_{H^{-1}}
		\le C\left(
		\|(\mathbf I-\mathbf P)f\|_{L^2(\cD;H_\sigma^1)}
		+\|v\cdot\nabla_xf+Lf\|_{L^2(\cD;H_\sigma^{-1})}
		\right).
	\end{equation*}
	
	Write the spatial averages as
	\[
	c_0:=\fint_{\cD} \mathbf c,
	\qquad
	a_0:=\fint_{\cD} \mathbf a.
	\]
	Applying Lemma~\ref{lem:negative-order-estimates}\,(\ref{it:neg-poincare}) to $\mathbf c$ and to
	$\mathbf a-\tfrac32\mathbf c$, we obtain
	\[
	\|\mathbf c-c_0\|_{L^2_x}+\|\mathbf a-a_0\|_{L^2_x}\le C(X+D).
	\]
	Let $\mathbf b_0+\Omega_0x$ be the $L^2(\cD;\bR^3)$-orthogonal projection of
	$\mathbf b$ onto the rigid-motion space $\cR_{\rm rig}$.
	Lemma~\ref{lem:negative-order-estimates}\,(\ref{it:neg-korn}) then gives
	\[
	\|\mathbf b-(\mathbf b_0+\Omega_0x)\|_{L^2_x}\le C(X+D).
	\]
	Define the rigid hydrodynamic mode
	\begin{equation*}
		r(x,v):=
		\bigl(a_0+(\mathbf b_0+\Omega_0x)\cdot v
		+c_0(|v|^2-\tfrac32)\bigr)\sqrt\mu.
	\end{equation*}
	Because $v^T\Omega_0v=0$,
	\begin{equation*}
		v\cdot\nabla_xr=0,
		\qquad Lr=0.
	\end{equation*}
	Being a Gaussian--polynomial profile with bounded coefficients on
	$\overline\cD$ and with vanishing transport, $r$ lies in
	$\cH_{\kin,\sigma}(\cD\times\bR_v^3)$.
	The preceding bounds and finite-dimensional equivalence of the
	Maxwellian modes give
	\begin{align}\label{260915416}
		\|f-r\|_{L^2(\cD;H_\sigma^1)}
		+\|v\cdot\nabla_x(f-r)\|_{L^2(\cD;H_\sigma^{-1})}
		\le C(X+D).
	\end{align}
	Since $v\cdot\nabla_xr=0$ and $L\mathbf Pf=0$, we have
	\[
	v\cdot\nabla_x(f-r)=(v\cdot\nabla_xf+Lf)-L(\mathbf I-\mathbf P)f,
	\]
	and
	\eqref{boundedness of L} yields the bound for $\|v\cdot\nabla_x(f-r)\|_{L^2(\cD;H_\sigma^{-1})}$ in \eqref{260915416}.
	
	Let $U\times V$ be the outgoing window from
	Lemma~\ref{lem:boundary-rigid-observability}.  Proposition~\ref{prop:bounded-nongrazing-trace}\,(\ref{it:nongrazing-stationary}),
	applied to $f-r$, gives
	\[
	\|f_{\gamma^+}-r\|_{L^2(U\times V,\dd\gamma)}\le C(X+D).
	\]
	Put $\beta:=\cA_0f_{\gamma^+}$.  Since
	$\cP_0^+f_{\gamma^+}=\mathfrak m_\circ\beta$, restriction of the outgoing
	variance to $U\times V$ and the preceding trace bound imply
	\begin{equation*}
		\|r-\mathfrak m_\circ\beta\|_{L^2(U\times V,\dd\gamma)}
		\le C(X+D+M).
	\end{equation*}
	Applying Lemma~\ref{lem:boundary-rigid-observability}, we therefore
	obtain
	\[
	|\mathbf b_0|+|\Omega_0|+|c_0|\le C(X+D+M).
	\]
	Finally, the microscopic part, the momentum modes, and the centered energy
	mode have zero mass moment.  Since $\int_{\bR_v^3}\mu\dd v=\pi^{3/2}$,
	\[
	\int_{\cD}\int_{\bR_v^3}f(x,v)\sqrt\mu(v)\dd v\dd x
	=\pi^{3/2}\int_\cD \mathbf a(x)\dd x=\pi^{3/2}|\cD|a_0.
	\]
	Hence $a_0=\bar{\mathbf a}$ and $\mathbf a^\circ=\mathbf a-a_0$.  Combining the bounds for the non-rigid modes with the
	preceding rigid-coefficient estimate proves
	\eqref{eq:stationary-macro-estimate}.  Adding
	$|a_0|=\pi^{-3/2}|\cD|^{-1}|\mathfrak m_f|$ proves
	\eqref{eq:macroscopic-recovery-estimate}.
\end{proof}

\subsection{The stationary problem}\label{ss:linear}

For $\|\mu^{1/4}g\|_{L^\infty(\cD\times\bR_v^3)}\le\varepsilon_{\rm in}$,
Theorem~\ref{260322442}\,(2) supplies, for any
$F\in L^2(\cD;H_\sigma^{-1})$ and $H\in L^2(\gamma^-,\dd\gamma)$,
the unique $f\in\cH_{\kin,\sigma}$ satisfying
$v\cdot\nabla_xf+L_gf=F$ and $f_{\gamma^-}=H$, together with
\eqref{eq:perturbed-landau-stationary-estimate}.

\begin{prop}[Stationary graph estimate]
	\label{prop:inhomogeneous-full-diffuse-graph}
	There are constants
	$0<\varepsilon_{\rm dif}\le\varepsilon_{\rm in}$ and $\delta_{\rm dif}>0$
	such that the following holds.  Assume
	\begin{equation}
		\label{eq:diffuse-graph-smallness}
		\|\mu^{1/4}g\|_{L^\infty(\cD\times\bR_v^3)}\le\varepsilon_{\rm dif},
		\qquad 0<\delta\le\delta_{\rm dif}.
	\end{equation}
	Let $\cJ$ be as in \eqref{eq:wall-contraction-operator}, let
	$F\in L^2(\cD;H_\sigma^{-1})$, $b\in L^2(\gamma^-,\dd\gamma)$, and let
	$f\in\cH_{\kin,\sigma}(\cD\times\bR_v^3)$ satisfy
	\begin{equation}
		\label{eq:inhomogeneous-full-diffuse-problem}
		\begin{cases}
			v\cdot\nabla_xf+Lf-N[g,f]=F
			&\text{in }L^2(\cD;H_\sigma^{-1}),\\
			f_{\gamma^-}=P_\cJ f_{\gamma^+}+b
			&\text{in }L^2(\gamma^-,\dd\gamma).
		\end{cases}
	\end{equation}
	Then
	\begin{equation}
		\label{eq:inhomogeneous-full-diffuse-graph-estimate}
		\|f\|_{\cH_{\kin,\sigma}(\cD\times\bR_v^3)}
		\le C\left(
		\|F\|_{L^2(\cD;H_\sigma^{-1})}
		+\|b\|_{L^2(\gamma^-)}
		+\left|\int_{\cD}\int_{\bR_v^3}
		f(x,v)\sqrt\mu(v)\dd v\dd x\right|
		\right).
	\end{equation}
	The constant is independent of $\cJ$, hence of $N$ when $\cJ=J_N$.  
\end{prop}

\begin{proof}
	Set
	\begin{equation*}
		X:=\|(\mathbf I-\mathbf P)f\|_{L^2(\cD;H_\sigma^1)},
		\quad
		M:=\|(\mathbf I-\cP_0^+)f_{\gamma^+}\|_{L^2(\gamma^+)},
		\quad
		B:=\|f_{\gamma^+}\|_{L^2(\gamma^+)},
	\end{equation*}
	and
	\begin{equation*}
		Y:=\|F\|_{L^2(\cD;H_\sigma^{-1})}
		+\|b\|_{L^2(\gamma^-)}
		+\left|\int_{\cD}\int_{\bR_v^3}
		f(x,v)\sqrt\mu(v)\dd v\dd x\right|.
	\end{equation*}
	
	With $S=F+N[g,f]$, Proposition~\ref{prop:macroscopic-recovery}
	and Lemma~\ref{260715201} give
	\begin{align}
		\|\mathbf P f\|_{L^2(\cD\times\bR_v^3)}
		&\le C\left(X+M+\|F\|_{L^2(\cD;H_\sigma^{-1})}
		+\left|\int_{\cD}\int_{\bR_v^3}f(x,v)\sqrt\mu(v)\dd v\dd x\right|\right)\notag\\
		&\qquad +C\|\mu^{1/4}g\|_{L^\infty(\cD\times\bR_v^3)}\|f\|_{L^2(\cD;H_\sigma^1)}\notag\\
		&\le C\left(X+M+Y\right)
		+C\|\mu^{1/4}g\|_{L^\infty(\cD\times\bR_v^3)}\left(X+\|\mathbf P f\|_{L^2(\cD\times\bR_v^3)}\right).
		\label{eq:graph-macro-preabsorption}
	\end{align}
	Here we used the finite-dimensional estimate
	$\|\mathbf P f\|_{L^2(\cD;H_\sigma^1)}\le C\|\mathbf P f\|_{L^2(\cD\times\bR_v^3)}$.
	Choose $\varepsilon_{\rm dif}$ so that the last coefficient in
	\eqref{eq:graph-macro-preabsorption} is at most $1/4$.  Then
	\begin{equation}
		\label{eq:graph-macro-bound}
		\|\mathbf P f\|_{L^2(\cD\times\bR_v^3)}+\|f\|_{L^2(\cD;H_\sigma^1)}
		\le C_0(X+M+Y).
	\end{equation}
	Since $L\mathbf P f=0$, the equation and \eqref{boundedness of L}--\eqref{260911445}
	imply
	\begin{equation*}
		\|v\cdot\nabla_xf\|_{L^2(\cD;H_\sigma^{-1})}
		\le C_0(X+M+Y).
	\end{equation*}
	Applying
	Proposition~\ref{prop:full-diffuse-trace-recovery}\,(\ref{it:recovery-stationary}),
	we obtain
	\[
	B\le C_0(X+M+Y).
	\]
	Finally, \eqref{eq:wall-J-operator-bound} gives
	\begin{equation}
		\label{eq:graph-incoming-bound}
		\|f_{\gamma^-}\|_{L^2(\gamma^-)}
		\le C_0(X+M+Y).
	\end{equation}
	These bounds control every component of the graph norm by $C_0(X+M+Y)$.
	
	Lemma~\ref{lem:wall-boundary-form}\,(\ref{it:wall-J-lower}) at
	$\alpha=1$, with the nonnegative term $\tfrac12\mathfrak d_\cJ^2$ omitted,
	gives for every $0<\theta<1$
	\[
	B^2-\|f_{\gamma^-}\|_{L^2(\gamma^-)}^2
	\ge\tfrac12M^2-\theta B^2-C\delta^2B^2-C_\theta\|b\|_{L^2(\gamma^-)}^2,
	\]
	with constants independent of $\cJ$.
	
	Pairing \eqref{eq:inhomogeneous-full-diffuse-problem} with $f$ and using the
	quadratic Green identity \eqref{ibp_st} of Proposition~\ref{2607151218}, we
	obtain
	\begin{equation*}
		\langle Lf,f\rangle
		-\langle N[g,f],f\rangle
		+\frac12\left(B^2-\|f_{\gamma^-}\|_{L^2(\gamma^-)}^2\right)
		=\langle F,f\rangle.
	\end{equation*}
	By \eqref{eq:L-coercivity}, \eqref{boundedness of L}--\eqref{260911445}, and the preceding boundary estimate, there is a constant
	${\rm c}_0>0$ such that
	\begin{align*}
		\begin{aligned}
			&{\rm c}_0(X^2+M^2)
			\\
			\le{}& C Y\|f\|_{L^2(\cD;H_\sigma^1)}
			+C\|\mu^{1/4}g\|_{L^\infty(\cD\times\bR_v^3)}\|f\|_{L^2(\cD;H_\sigma^1)}^2
			+C(\theta+\delta^2)B^2+C_\theta Y^2.
		\end{aligned}
	\end{align*}
	We insert \eqref{eq:graph-macro-bound} and the bound $B\le C_0(X+M+Y)$ obtained above.  With
	$Z=(X^2+M^2)^{1/2}$ and $X+M\le\sqrt2Z$, there is a fixed $K_0\ge1$ such
	that
	\[
	{\rm c}_0Z^2
	\le K_0Y(\sqrt2Z+Y)
	+K_0(\|\mu^{1/4}g\|_{L^\infty(\cD\times\bR_v^3)}+\theta+\delta^2)(\sqrt2Z+Y)^2
	+C_\theta Y^2.
	\]
	\begingroup
	Choose $0<\theta<1$ sufficiently small, depending on
	$K_0$ and ${\rm c}_0$, and then decrease $\varepsilon_{\rm dif}$ and
	$\delta_{\rm dif}$ so that the coefficient
	$K_0(\|\mu^{1/4}g\|_{L^\infty}+\theta+\delta^2)$ is
	small compared with ${\rm c}_0$.  These choices apply to every wall
	law \eqref{eq:wall-contraction-operator}.  By the preceding energy
	estimate and Young's inequality, the $Z^2$ and mixed $YZ$ terms
	can be absorbed, yielding $Z^2\le CY^2$.
	\endgroup
	Hence
	\[
	X+M\le CY.
	\]
	Substitution into \eqref{eq:graph-macro-bound} and
	\eqref{eq:graph-incoming-bound} controls the velocity-coercive, transport, and
	both trace components by $CY$.  Lemma~\ref{lem:transport-tail-L2}\,(2)
	then supplies
	\[
	\|f\|_{L^2(\cD\times\bR_v^3)}\le CY.
	\]
	Together with \eqref{260902137}, these bounds prove
	\eqref{eq:inhomogeneous-full-diffuse-graph-estimate}.
\end{proof}

\begin{corollary}[Stationary homogeneous coercivity]
	\label{cor:homogeneous-diffuse-coercivity}
	Under the smallness conditions of
	Proposition~\ref{prop:inhomogeneous-full-diffuse-graph}, let
	$S\in L^2(\cD;H_\sigma^{-1})$ and let
	$w\in\cH_{\kin,\sigma}(\cD\times\bR_v^3)$ satisfy
	\[
	v\cdot\nabla_xw+Lw-N[g,w]=S,
	\qquad
	w_{\gamma^-}=P_ww_{\gamma^+},
	\qquad
	\int_{\cD}\int_{\bR_v^3}w(x,v)\sqrt\mu(v)\dd v\dd x=0.
	\]
	Then
	\begin{equation}
		\label{eq:homogeneous-diffuse-coercivity}
		\|w\|_{\cH_{\kin,\sigma}(\cD\times\bR_v^3)}
		\le C\|S\|_{L^2(\cD;H_\sigma^{-1})}.
	\end{equation}
\end{corollary}

\begin{proof}
	Apply Proposition~\ref{prop:inhomogeneous-full-diffuse-graph} with $F=S$,
	$b=0$, and $\int_{\cD}\int_{\bR_v^3}w(x,v)\sqrt\mu(v)\dd v\dd x=0$.
\end{proof}

\begin{corollary}[Uniform stationary graph estimate]
	\label{cor:spectral-uniform-estimate}
	Assume \eqref{eq:diffuse-graph-smallness}, let
	$F\in L^2(\cD;H_\sigma^{-1})$ and $b\in L^2(\gamma^-,\dd\gamma)$, and let
	$f\in\cH_{\kin,\sigma}(\cD\times\bR_v^3)$ satisfy
	\begin{equation}
		\label{eq:spectral-full-diffuse-problem}
		v\cdot\nabla_xf+Lf-N[g,f]=F,
		\qquad f_{\gamma^-}=P_{w,N}f_{\gamma^+}+b .
	\end{equation}
	Then
	\begin{equation}
		\label{eq:spectral-full-diffuse-estimate}
		\|f\|_{\cH_{\kin,\sigma}(\cD\times\bR_v^3)}
		\le C\left(
		\|F\|_{L^2(\cD;H_\sigma^{-1})}+\|b\|_{L^2(\gamma^-)}
		+\left|\int_\cD\int_{\bR_v^3}f\sqrt\mu\dd v\dd x\right|
		\right),
	\end{equation}
	with $C$ independent of $N$.
\end{corollary}

\begin{proof}
	This is Proposition~\ref{prop:inhomogeneous-full-diffuse-graph} with
	$\cJ=J_N$, whose constant does not depend on $\cJ$.
\end{proof}

\begin{corollary}[Weighted stationary graph estimate]
	\label{cor:spectral-weighted-estimate}
	Assume \eqref{eq:diffuse-graph-smallness}, let $L\ge0$ be finite, let
	$w^LF\in L^2(\cD;H_\sigma^{-1})$ and $w^Lb\in L^2(\gamma^-,\dd\gamma)$, and let
	$f\in\cH_{\kin,\sigma}(\cD\times\bR_v^3)$ satisfy
	\eqref{eq:spectral-full-diffuse-problem}.  Then
	\begin{align}
		\|w^Lf\|_{\cH_{\kin,\sigma}(\cD\times\bR_v^3)}
		+\|w^Lf_{\gamma^-}\|_{L^2(\gamma^-)}
		+\|w^Lf_{\gamma^+}\|_{L^2(\gamma^+)}
		&\notag\\
		\le C_L\Bigl(
		\|w^LF\|_{L^2(\cD;H_\sigma^{-1})}
		+\|w^Lb\|_{L^2(\gamma^-)}
		&+\Bigl|\int_\cD\int_{\bR_v^3}f\sqrt\mu\dd v\dd x\Bigr|
		\Bigr),
		\label{eq:spectral-weighted-diffuse-estimate}
	\end{align}
	with $C_L$ independent of $N$.  At $L=0$ this is
	\eqref{eq:spectral-full-diffuse-estimate}.
\end{corollary}

\begin{proof}
	Regard \eqref{eq:spectral-full-diffuse-problem} as the prescribed-inflow problem
	\eqref{solv:linland_st} with source $F$ and incoming datum
	$H:=P_{w,N}f_{\gamma^+}+b$, and apply
	\eqref{eq:weighted-inflow-estimate} at that $L$.  By
	\eqref{eq:wall-J-weighted-operator-bound} at $\cJ=J_N$ and the trace bound of
	Proposition~\ref{2607151218},
	\[
	\|w^LH\|_{L^2(\gamma^-)}
	\le C_L\|f_{\gamma^+}\|_{L^2(\gamma^+)}+\|w^Lb\|_{L^2(\gamma^-)}
	\le C_L\|f\|_{\cH_{\kin,\sigma}(\cD\times\bR_v^3)}
	+\|w^Lb\|_{L^2(\gamma^-)} ,
	\]
	and \eqref{eq:spectral-full-diffuse-estimate} bounds the graph norm on the
	right.  Since $w\ge1$, the weighted data norms also control the unweighted norms of $F$ and $b$.
\end{proof}

\begin{prop}[Stationary linear diffuse solvability]
	\label{prop:full-diffuse-graph-solvability}
	Assume \eqref{eq:diffuse-graph-smallness}, and let
	$F\in L^2(\cD;H_\sigma^{-1})$, $b\in L^2(\gamma^-,\dd\gamma)$, and
	$m\in\bR$.  Then
	\begin{equation}
		\label{eq:full-diffuse-linear-problem}
		\begin{cases}
			v\cdot\nabla_xf+Lf-N[g,f]=F,\\
			f_{\gamma^-}=P_wf_{\gamma^+}+b,\\
			\displaystyle\int_\cD\int_{\bR_v^3}f\sqrt\mu\dd v\dd x=m
		\end{cases}
	\end{equation}
	has a solution $f\in\cH_{\kin,\sigma}(\cD\times\bR_v^3)$ if and only if
	\begin{equation}
		\label{eq:linear-full-diffuse-compatibility}
		\int_\cD\langle F(x,\cdot),\sqrt\mu\rangle_{H_\sigma^{-1},H_\sigma^1}
		\dd x
		+\int_{\gamma^-}b\sqrt\mu\dd\gamma=0.
	\end{equation}
	The solution is unique and satisfies
	\begin{equation}
		\label{eq:full-diffuse-linear-bound}
		\|f\|_{\cH_{\kin,\sigma}(\cD\times\bR_v^3)}
		\le C\left(
		\|F\|_{L^2(\cD;H_\sigma^{-1})}+\|b\|_{L^2(\gamma^-)}+|m|
		\right).
	\end{equation}
	In particular, the data $F=0$, $b=h_w$, and $m=0$ are compatible.
\end{prop}

\begin{proof}
	We first derive the compatibility condition and solve the problem with a spectral wall cutoff.
	
	Apply the Green identity with the $x$-independent test function
	$\sqrt\mu$.  This is allowed: $\sqrt\mu\in\cH_{\kin,\sigma}(\cD\times\bR_v^3)$
	by Proposition~\ref{H_kin=weak_landau}\,(2).  Indeed $\cD$ is bounded and
	$\sqrt\mu\in H_\sigma^1$, so $\sqrt\mu\in L^2(\cD;H_\sigma^1)$; take
	$v\cdot\nabla_x\sqrt\mu=0$ and $(\sqrt\mu)_\gamma=\sqrt\mu\in L^2(\gamma)$;
	for $\phi\in C_c^\infty(\overline\cD\times\bR_v^3)$ the divergence theorem gives
	\[
	\int_\cD\langle\sqrt\mu,v\cdot\nabla_x\phi\rangle_v\dd x =\int_\gamma\sqrt\mu\,\phi\,\operatorname{sgn}(n_x\cdot v)\dd\gamma,
	\]
	which is \eqref{eq:H_kin_sigma-weak-stationary}.  Lemma~\ref{260715201} gives
	\begin{equation}
		\label{eq:physical-mass-flux-identity}
		\int_\cD\langle F,\sqrt\mu\rangle\dd x
		=\int_{\gamma^+}f_{\gamma^+}\sqrt\mu\dd\gamma
		-\int_{\gamma^-}f_{\gamma^-}\sqrt\mu\dd\gamma.
	\end{equation}
	For almost every $x\in\partial\cD$, the unit-flux normalization of $M_\Theta$ gives
	\[
	\int_{n_x\cdot v<0}P_wf_{\gamma^+}(x,v)\sqrt\mu(v)
	|n_x\cdot v|\dd v=\cA_0f_{\gamma^+}(x).
	\]
	Substitution of the boundary law into
	\eqref{eq:physical-mass-flux-identity} proves
	\eqref{eq:linear-full-diffuse-compatibility}.
	
	For the spectral boundary approximation, the map
	$\phi\mapsto\mathfrak m_{\Theta_w}\phi$ is bounded above and below
	from $L^2(\partial\cD)$ into the incoming flux space.  Indeed, for
	$\Theta\in[\tfrac14,\tfrac34]$, rotation of the normal to $e_3$ and direct Gaussian
	integration give
	\begin{equation*}
		\int_{v_3<0}\left|\frac{M_\Theta(v)}{\sqrt\mu(v)}\right|^2(-v_3)\dd v
		=\frac{1}{8\pi\,\Theta^2(1-\Theta)^2}
		=\frac{c_\mu^{-1}}{y^2(2-y)^2},
		\qquad y:=2\Theta.
	\end{equation*}
	Since $y\in[1/2,3/2]$,
	\begin{equation}
		\label{eq:wall-profile-injectivity}
		c_\mu^{-1}\|\phi\|_{L^2(\partial\cD)}^2
		\le\|\mathfrak m_{\Theta_w}\phi\|_{L^2(\gamma^-)}^2
		\le\frac{16}{9}c_\mu^{-1}\|\phi\|_{L^2(\partial\cD)}^2.
	\end{equation}
	
	For $\phi\in E_N$, let $U_\phi$ be the unique prescribed-inflow solution from
	Theorem~\ref{260322442}\,(2) with zero source and
	incoming datum $\mathfrak m_{\Theta_w}\phi$.  Uniqueness there makes
	$\phi\mapsto U_\phi$ linear.  Define
	\begin{equation*}
		A_N\phi:=\phi-J_N\cA_0(U_\phi)_{\gamma^+},
		\qquad
		\beta_N(\phi):=\int_\cD\int_{\bR_v^3}U_\phi\sqrt\mu\dd v\dd x.
	\end{equation*}
	Then $A_N:E_N\to E_N$ and $\beta_N:E_N\to\bR$ are linear.  Applying the
	source-free form of \eqref{eq:physical-mass-flux-identity} and the pointwise
	unit-flux normalization, we obtain
	\begin{equation*}
		\int_{\partial\cD}\cA_0(U_\phi)_{\gamma^+}\dd S_x
		=\int_{\partial\cD}\phi\dd S_x.
	\end{equation*}
	Since $J_N$ is self-adjoint and $J_N1=1$,
	\[
	\operatorname{Ran}A_N\subset E_N^0,
	\qquad
	E_N^0:=\left\{\psi\in E_N:\int_{\partial\cD}\psi\dd S_x=0\right\}.
	\]
	If $A_N\phi=0$, then $U_\phi$ satisfies the homogeneous spectral diffuse law, so
	\eqref{eq:spectral-full-diffuse-estimate} gives
	\[
	\|U_\phi\|_{\cH_{\kin,\sigma}(\cD\times\bR_v^3)}
	\le C|\beta_N(\phi)|.
	\]
	If $\dim_\bR\ker A_N\ge2$, the restriction of the linear functional
	$\beta_N$ to $\ker A_N$ has a nonzero kernel vector $\phi$.  Then the
	preceding bound gives $U_\phi=0$, and the incoming
	trace together with \eqref{eq:wall-profile-injectivity} yields $\phi=0$, a
	contradiction.  Hence $\dim\ker A_N\le1$.
	
	Since $1\in E_N$, the surface-mean functional is nonzero on $E_N$, so
	$E_N^0$ has codimension one.  The range of $A_N$ lies in $E_N^0$, so its
	rank is at most $\dim E_N-1$, and rank--nullity gives
	$\dim\ker A_N\ge1$.  Thus
	\[
	\dim\ker A_N=1,
	\qquad \operatorname{Ran}A_N=E_N^0.
	\]
	For every $0\ne\kappa_N\in\ker A_N$,
	\[
	\beta_N(\kappa_N)\ne0,
	\]
	for otherwise
	$\|U_{\kappa_N}\|_{\cH_{\kin,\sigma}(\cD\times\bR_v^3)}\le C|\beta_N(\kappa_N)|$
	and \eqref{eq:wall-profile-injectivity} would force $\kappa_N=0$.
	
	Assume \eqref{eq:linear-full-diffuse-compatibility} and let $U_{F,b}$
	be the prescribed-inflow solution with source $F$ and incoming
	datum $b$.  Equation \eqref{eq:physical-mass-flux-identity} gives
	\begin{equation*}
		\int_{\partial\cD}\cA_0(U_{F,b})_{\gamma^+}\dd S_x
		=\int_\cD\langle F,\sqrt\mu\rangle\dd x
		+\int_{\gamma^-}b\sqrt\mu\dd\gamma.
	\end{equation*}
	Under \eqref{eq:linear-full-diffuse-compatibility}, self-adjointness and
	$J_N1=1$ imply
	$J_N\cA_0(U_{F,b})_{\gamma^+}\in E_N^0$.  Since $A_N$ maps $E_N$ onto $E_N^0$,
	choose $\phi_0\in E_N$ with
	\[
	A_N\phi_0=J_N\cA_0(U_{F,b})_{\gamma^+}.
	\]
	Then $\widetilde f_N:=U_{\phi_0}+U_{F,b}$ satisfies
	\eqref{eq:spectral-full-diffuse-problem}, since the boundary law there is
	equivalent to
	\[
	\phi_0=J_N\cA_0(\widetilde f_N)_{\gamma^+}.
	\]
	The homogeneous response $U_{\kappa_N}$ satisfies the same spectral boundary
	law with zero source and residual.  Since $\beta_N(\kappa_N)\ne0$, the unique
	scalar
	\begin{equation*}
		t_N:=\frac{m-
			\displaystyle\int_\cD\int_{\bR_v^3}\widetilde f_N\sqrt\mu\dd v\dd x}
		{\beta_N(\kappa_N)}
	\end{equation*}
	makes
	$f_N:=\widetilde f_N+t_NU_{\kappa_N}$ have bulk mass $m$.
	\eqref{eq:spectral-full-diffuse-estimate} then gives the uniform
	estimate
	\[
	\|f_N\|_{\cH_{\kin,\sigma}(\cD\times\bR_v^3)}
	\le C\left(\|F\|_{L^2(\cD;H_\sigma^{-1})}
	+\|b\|_{L^2(\gamma^-)}+|m|\right).
	\]
	If $f_N^{(1)}$ and $f_N^{(2)}$ solve the same fixed-$N$ problem with the
	same source, residual, and mass, then their difference solves the homogeneous
	fixed-$N$ problem with zero incoming residual and zero mass.  Applying
	\eqref{eq:spectral-full-diffuse-estimate} to that difference gives
	$\|f_N^{(1)}-f_N^{(2)}\|_{\cH_{\kin,\sigma}(\cD\times\bR_v^3)}=0$, hence uniqueness
	at each fixed $N$.
	
	The uniform estimate makes the unweighted bulk coordinate,
	the velocity-coercive coordinate, the transport coordinate, and both flux
	traces weakly precompact.  Passing to a subsequence, denote their limits by
	\[
	f^{(0)},\quad f^{(1)},\quad G,\quad h_{\gamma^+},\quad h_{\gamma^-}.
	\]
	Compactly supported smooth tests identify $f^{(0)}=f^{(1)}=:f$
	distributionally and $G=v\cdot\nabla_xf$.  Passing the Green identity to
	the limit gives
	\begin{equation*}
		\langle v\cdot\nabla_xf,\psi\rangle
		+\int_{\cD\times\bR_v^3}f\,v\cdot\nabla_x\psi\dd v\dd x
		=\int_{\gamma^+}h_{\gamma^+}\psi\dd\gamma
		-\int_{\gamma^-}h_{\gamma^-}\psi\dd\gamma
	\end{equation*}
	for every $\psi\in C_c^\infty(\overline\cD\times\bR_v^3)$.  This identity is
	\eqref{eq:H_kin_sigma-weak-stationary} with $f_\gamma:=h_{\gamma^\pm}$ on
	$\gamma^\pm$, and $f\in L^2(\cD;H^1_\sigma)$,
	$v\cdot\nabla_xf\in L^2(\cD;H^{-1}_\sigma)$, $f_\gamma\in L^2(\gamma)$.  Hence
	Proposition~\ref{H_kin=weak_landau}\,(2) yields
	$f\in\cH_{\kin,\sigma}(\cD\times\bR_v^3)$ and identifies $h_{\gamma^\pm}$ with its
	boundary traces.  Since $L-N[g,\cdot]$ is bounded, the equation passes to the
	limit, and so does the prescribed mass, by weak convergence of the unweighted
	bulk coordinate.
	
	Set $\phi_N:=\cA_0(f_N)_{\gamma^+}$.
	Since \eqref{eq:outgoing-moment-flux-bound} bounds $\cA_0$ on the outgoing
	flux space, $\phi_N\rightharpoonup\phi:=\cA_0f_{\gamma^+}$ in
	$L^2(\partial\cD)$.  For every fixed $\psi\in L^2(\partial\cD)$,
	self-adjointness of $J_N$ gives
	\[
	\langle J_N\phi_N-\phi_N,\psi\rangle
	=\langle\phi_N,(J_N-\mathbf I)\psi\rangle\longrightarrow0
	\]
	by \eqref{eq:boundary-spectral-strong-convergence} and the uniform bound on
	$\phi_N$.  Thus $J_N\phi_N\rightharpoonup\phi$.  The upper bound in
	\eqref{eq:wall-profile-injectivity} makes
	$\psi\mapsto\mathfrak m_{\Theta_w}\psi$ a fixed bounded map from
	$L^2(\partial\cD)$ to $L^2(\gamma^-)$.  Passing to the limit in
	\[
	(f_N)_{\gamma^-}=\mathfrak m_{\Theta_w}J_N\phi_N+b
	\]
	therefore gives
	$f_{\gamma^-}=P_wf_{\gamma^+}+b$.
	Weak lower semicontinuity of each term of the graph norm proves
	\eqref{eq:full-diffuse-linear-bound}.
	
	If two solutions have the same data and mass, their difference has
	zero source, zero boundary residual, and zero bulk mass, and
	Proposition~\ref{prop:inhomogeneous-full-diffuse-graph} forces it to vanish.
	Finally, the definition
	$h_w=c_\mu\mathfrak m_{\Theta_w}-\sqrt\mu$ and the unit-flux normalization of $M_\Theta$ give
	$\int_{\gamma^-}h_w\sqrt\mu\dd\gamma=0$, proving the last assertion.
\end{proof}

The same finite-dimensional construction gives the fixed-cutoff statement.

\begin{corollary}[Stationary solvability at fixed spectral cutoff]
	\label{cor:fixedN-spectral-solvability}
	Assume \eqref{eq:diffuse-graph-smallness}.  Fix $N\in\mathbb N$ and let
	$F\in L^2(\cD;H_\sigma^{-1})$, $b\in L^2(\gamma^-,\dd\gamma)$, and
	$m\in\bR$.  Then the problem
	\begin{equation}
		\label{eq:fixedN-spectral-problem-mass}
		\begin{cases}
			v\cdot\nabla_xf+Lf-N[g,f]=F,\\
			f_{\gamma^-}=P_{w,N}f_{\gamma^+}+b,\\
			\displaystyle\int_\cD\int_{\bR_v^3}f\sqrt\mu\dd v\dd x=m
		\end{cases}
	\end{equation}
	has a solution $f\in\cH_{\kin,\sigma}(\cD\times\bR_v^3)$ if and only if
	$(F,b)$ satisfies \eqref{eq:linear-full-diffuse-compatibility}, i.e.
	\[
	\int_\cD\langle F(x,\cdot),\sqrt\mu\rangle_{H_\sigma^{-1},H_\sigma^1}
	\dd x
	+\int_{\gamma^-}b\sqrt\mu\dd\gamma=0 .
	\]
	The solution is unique and satisfies
	\eqref{eq:spectral-full-diffuse-estimate}, and hence
	\begin{equation}
		\label{eq:fixedN-spectral-solvability-bound}
		\|f\|_{\cH_{\kin,\sigma}(\cD\times\bR_v^3)}
		\le C\left(
		\|F\|_{L^2(\cD;H_\sigma^{-1})}+\|b\|_{L^2(\gamma^-)}+|m|
		\right),
	\end{equation}
	with $C$ independent of $N$.  In particular $F=0$, $b=h_w$, $m=0$ are
	compatible data for every $N$.
\end{corollary}

\begin{proof}
	\emph{Necessity.}  Let $f$ solve \eqref{eq:fixedN-spectral-problem-mass}.
	The compatibility argument in the proof of
	Proposition~\ref{prop:full-diffuse-graph-solvability} uses only the
	interior equation, so it gives \eqref{eq:physical-mass-flux-identity} for
	this $f$.  By the pointwise flux identity before
	\eqref{eq:truncated-wall-flux-integrated}, integrated over $\p\cD$, and by
	\eqref{eq:truncated-wall-flux-integrated},
	\[
	\int_{\gamma^-}P_{w,N}f_{\gamma^+}\sqrt\mu\dd\gamma
	=\int_{\p\cD}J_N\cA_0f_{\gamma^+}\dd S_x
	=\int_{\p\cD}\cA_0f_{\gamma^+}\dd S_x
	=\int_{\gamma^+}f_{\gamma^+}\sqrt\mu\dd\gamma .
	\]
	Substitution of $f_{\gamma^-}=P_{w,N}f_{\gamma^+}+b$ into
	\eqref{eq:physical-mass-flux-identity} gives
	\eqref{eq:linear-full-diffuse-compatibility}.
	
	\emph{Existence and uniqueness.}  The finite-dimensional construction in the proof of
	Proposition~\ref{prop:full-diffuse-graph-solvability} uses only
	\eqref{eq:diffuse-graph-smallness}, Theorem~\ref{260322442}\,(2)
	(note $\varepsilon_{\rm dif}\le\varepsilon_{\rm in}$),
	Corollary~\ref{cor:spectral-uniform-estimate}, and
	\eqref{eq:linear-full-diffuse-compatibility}; they do not require passage to $N\to\infty$.  Under \eqref{eq:linear-full-diffuse-compatibility} they give a
	solution $f_N$ of \eqref{eq:fixedN-spectral-problem-mass} and its uniqueness.
	
	\emph{Estimate.}  Corollary~\ref{cor:spectral-uniform-estimate} applies to
	every solution and gives \eqref{eq:spectral-full-diffuse-estimate}.  The mass
	term there equals $|m|$, which gives
	\eqref{eq:fixedN-spectral-solvability-bound}; the constant is that of
	Corollary~\ref{cor:spectral-uniform-estimate}, so it does not depend on $N$.
	Finally $\int_{\gamma^-}h_w\sqrt\mu\dd\gamma=0$, as shown at the end of the
	proof of Proposition~\ref{prop:full-diffuse-graph-solvability}.
\end{proof}

\subsection{Macroscopic estimate for evolutionary problems}\label{ss:diffuse-homotopies}

Throughout this subsection, $\cJ$ is a contraction as in \eqref{eq:wall-contraction-operator};
the \emph{physical path} is $\cJ=\mathbf I$, and the \emph{boundary-spectral
	path} is $\cJ=J_N$; in both cases $\cJ$ is self-adjoint on $L^2(\p\cD)$ with $0\le\cJ\le\mathbf I$ and $\cJ1=1$.
When $\alpha=0$, the law \eqref{eq:dynamic-admissible-boundary-law} is the prescribed inflow
$f_{\gamma^-}=b$. 

For an interval $I=(t_0,t_1)\subset(0,T)$, together with the interval notation
\eqref{eq:dynamic-I-boundary-cylinders}, we abbreviate
\begin{equation}
	\label{eq:dynamic-window-quantities}
	A_I:=\|\mathbf Pf\|_{L^2(I;L^2(\cD\times\bR_v^3))},
	\qquad
	X_I:=\|(\mathbf I-\mathbf P)f\|_{L^2(\cD_I;H_\sigma^1)},
	\qquad
	D_I:=\|S\|_{L^2(\cD_I;H_\sigma^{-1})},
\end{equation}
$S$ being the source of \eqref{eq:dynamic-linear-source-equation} below,
\begin{equation}
	\label{eq:dynamic-window-traces}
	B_I:=\|f_{\gamma^+}\|_{L^2(\gamma_I^+)},
	\qquad
	M_I:=\|(\mathbf I-\cP_0^+)f_{\gamma^+}\|_{L^2(\gamma_I^+)} .
\end{equation}
The outgoing and incoming boundary mass scalars are
\begin{equation}
	\label{eq:dynamic-boundary-mass-scalars}
	\varphi_f(t,x):=\cA_0f_{\gamma^+}(t,x),
	\qquad
	\beta_b(t,x):=
	\int_{n_x\cdot v<0}b(t,x,v)\sqrt\mu(v)|n_x\cdot v|\dd v,
\end{equation}
and Cauchy--Schwarz gives
\begin{align}
	\|\varphi_f(t)\|_{L^2(\partial\cD)}^2
	&\le c_\mu\|f_{\gamma^+}(t)\|_{L^2(\gamma^+)}^2,
	\label{eq:dynamic-varphi-bound}\\
	\|\beta_b(t)\|_{L^2(\partial\cD)}^2
	&\le c_\mu\|b(t)\|_{L^2(\gamma^-)}^2.
	\label{eq:dynamic-beta-b-bound}
\end{align}

For either path, integrate \eqref{eq:wall-spectral-defect} and
\eqref{eq:wall-continuation-defect} over $I$ and write
\[
\mathfrak d_{\cJ,I}^2
:=c_\mu^{-1}\int_I
\left(
\|\varphi_f(t)\|_{L^2(\partial\cD)}^2
-\|\cJ\varphi_f(t)\|_{L^2(\partial\cD)}^2
\right)\dd t,
\]
\begin{equation}
	\label{eq:dynamic-alpha-K}
	K_I^2:=(1-\alpha^2)B_I^2+\alpha^2\mathfrak d_{\cJ,I}^2 ,
\end{equation}
so that $\mathfrak d_{\cJ,I}=0$ and $K_I^2=(1-\alpha^2)B_I^2$ on the
physical path.  Lemma~\ref{lem:wall-boundary-form}\,(\ref{it:wall-J-lower}),
integrated over $I$, is the boundary input of this subsection: if
$f_{\gamma^\pm}\in L^2(\gamma_I^\pm)$, $b\in L^2(\gamma_I^-)$ and $0\le\alpha\le1$
satisfy \eqref{eq:dynamic-admissible-boundary-law}, then
\begin{equation}
	\label{eq:dynamic-graph-boundary-input}
	B_I^2-\|f_{\gamma^-}\|_{L^2(\gamma_I^-)}^2
	\ge\tfrac12K_I^2+\tfrac12M_I^2
	-(\theta+C\delta^2)B_I^2
	-C_\theta\|b\|_{L^2(\gamma_I^-)}^2
\end{equation}
for every $0<\theta<1$, with constants independent of $\cJ$, of $\alpha$, of
$I$, and of $T$.

Since \eqref{def:reference-maxwellian} gives
\begin{equation*}
	\int\mu\dd v=\pi^{3/2},
	\qquad
	\int v_iv_j\mu\dd v=\tfrac12\pi^{3/2}\delta_{ij},
	\qquad
	\int(|v|^2-\tfrac32)^2\mu\dd v=\tfrac32\pi^{3/2},
\end{equation*}
we have the exact identity
\begin{align}
	\|\mathbf P f(t)\|_{L^2(\cD\times\bR_v^3)}^2
	&=\pi^{3/2}\|\mathbf a^\circ(t)\|_{L^2_x}^2
	+\pi^{-3/2}|\cD|^{-1}|\mathfrak m_f(t)|^2
	+\tfrac12\pi^{3/2}\|\mathbf b(t)\|_{L^2_x}^2
	+\tfrac32\pi^{3/2}\|\mathbf c(t)\|_{L^2_x}^2.
	\label{eq:dynamic-macro-norm-exact}
\end{align}

\begin{lemma}[Macroscopic estimates along the $\alpha$-paths]
	\label{lem:dynamic-alpha-macroscopic-estimate}
	Let $I=(t_0,t_1)$, let $0\le\alpha\le1$, let
	\[
	S\in L^2(\cD_I;H_\sigma^{-1}),
	\qquad
	b\in L^2(\gamma_I^-),
	\]
	and let $f\in\cH_{\kin,\sigma}(\cD_I\times\bR_v^3)$ satisfy
	\begin{equation}
		\label{eq:dynamic-linear-source-equation}
		(\partial_t+v\cdot\nabla_x)f+Lf=S
		\qquad\text{on }I\times\cD\times\bR_v^3,
	\end{equation}
	and 
	\begin{equation}
		\label{eq:dynamic-admissible-boundary-law_111}
		f_{\gamma^-}=\alpha P_\cJ f_{\gamma^+}+b
		\qquad\text{on }\gamma^-_I.
	\end{equation}
	Then there is $\eta_f\in C([t_0,t_1])$ such that
	\begin{equation}
		\label{eq:dynamic-macro-eta-bound}
		|\eta_f(t)|\le C\|f(t)\|_{L^2(\cD\times\bR_v^3)}^2
	\end{equation}
	and
	\begin{align}
		&\|\mathbf a^\circ\|_{L^2(I;L_x^2)}^2
		+\|\mathbf b\|_{L^2(I;L_x^2)}^2
		+\|\mathbf c\|_{L^2(I;L_x^2)}^2
		\notag\\
		&\quad\le C\Bigl\{
		\|(\mathbf I-\mathbf P)f\|_{L^2(\cD_I;H_\sigma^1)}^2
		+\|S\|_{L^2(\cD_I;H_\sigma^{-1})}^2
		+\|(\mathbf I-\cP_0^+)f_{\gamma^+}\|_{L^2(\gamma_I^+)}^2
		\notag\\
		&\qquad\qquad
		+K_I^2+\delta^2\|f_{\gamma^+}\|_{L^2(\gamma_I^+)}^2
		+\|b\|_{L^2(\gamma_I^-)}^2
		+\int_I|\mathfrak m_f(t)|^2\dd t
		\Bigr\}
		+\eta_f(t_0)-\eta_f(t_1),
		\label{eq:dynamic-alpha-macro-estimate}
	\end{align}
	with $K_I$ from \eqref{eq:dynamic-alpha-K}. 
	Here $\eta_f(t)$ depends only on $f(t)$ and $\cD$, and $C$ depends only on $\cD$.
\end{lemma}

\begin{proof}
	We first argue for smooth $f$. For general graph data choose
	$f_k\in C_c^\infty(\overline I\times\overline{\cD}\times\bR^3_v)$ with
	$f_k\to f$ in $\cH_{\kin,\sigma}(\cD_I\times\bR^3_v)$. Set
	\[
	S_k:=(\partial_t+v\cdot\nabla_x)f_k+Lf_k,
	\qquad
	b_k:=(f_k)_{\gamma^-}-\alpha P_\cJ(f_k)_{\gamma^+}.
	\]
	Then $f_k$ satisfies \eqref{eq:dynamic-linear-source-equation} and
	\eqref{eq:dynamic-admissible-boundary-law_111} with $(S,b)$ replaced by
	$(S_k,b_k)$. By Proposition~\ref{2607151218}, \eqref{boundedness of L} and
	\eqref{eq:wall-J-operator-bound}, as $k\to\infty$,
	\[
	f_k\to f\ \text{in}\ C(\overline I;L^2(\cD\times\bR^3_v))\cap L^2(\cD_I;H^1_\sigma),
	\quad
	(f_k)_{\gamma^\pm}\to f_{\gamma^\pm}\ \text{in}\ L^2(\gamma_I^\pm),
	\quad
	S_k\to S,
	\quad
	b_k\to b,
	\]
	with $S_k\to S$ in $L^2(\cD_I;H^{-1}_\sigma)$ and $b_k\to b$ in $L^2(\gamma_I^-)$.
	Every term of \eqref{eq:dynamic-alpha-macro-estimate} for $(f_k,S_k,b_k)$
	converges to the same term for $(f,S,b)$; for $\eta_{f_k}(t_0)$ and
	$\eta_{f_k}(t_1)$ this holds because $\eta_f(t)$ is a bounded quadratic form of
	$f(t)\in L^2(\cD\times\bR^3_v)$ (see \eqref{eq:dynamic-macro-interaction-definition}).
	Hence the estimate for $f_k$ gives the estimate for $f$.

	Since $0\le\cJ\le\mathbf I$, $(\mathbf I-\cJ)^2\le\mathbf I-\cJ^2$, so that for any $\psi\in L^2(\p\cD)$,
	\begin{align}\label{260914420}
		\|(\mathbf I-\cJ)\psi\|_{L^2(\partial\cD)}^2
		\le\|\psi\|_{L^2(\partial\cD)}^2-\|\cJ\psi\|_{L^2(\partial\cD)}^2 ,
	\end{align}
	while $(1-\alpha)^2\le1-\alpha^2$.
	Combining
	\[
	\varphi_f-\alpha\cJ\varphi_f =(1-\alpha)\varphi_f+\alpha(\mathbf I-\cJ)\varphi_f ,
	\]
	and \eqref{eq:dynamic-varphi-bound}, and using \eqref{260914420} and \eqref{eq:wall-J-components}, we obtain
	\begin{equation}
		\label{eq:dynamic-alpha-square-payment}
		\|\varphi_f-\alpha\cJ\varphi_f\|_{L^2(I\times\partial\cD)}^2
		\le 2(1-\alpha)^2\|\varphi_f\|_{L^2(I\times\partial\cD)}^2
		+2\alpha^2\|(\mathbf I-\cJ)\varphi_f\|_{L^2(I\times\partial\cD)}^2
		\le CK_I^2 .
	\end{equation}
	Thus $K_I$ controls the boundary mass mismatch uniformly in
	$\alpha$ and $\cJ$.
	
	For $\mathbf c$ and $\mathbf b$, take the Dirichlet potentials
	\begin{equation}
		\label{eq:dynamic-macro-Dirichlet-potentials}
		\left\{
		\begin{aligned}
			-\Delta\phi_c&=\mathbf c,&\phi_c|_{\partial\cD}&=0,\\
			-\Delta\phi_j&=\mathbf b_j,&\phi_j|_{\partial\cD}&=0,
			\qquad j=1,2,3.
		\end{aligned}
		\right.
	\end{equation}
	For the mass mode solve instead
	\begin{equation}
		\label{eq:dynamic-macro-Neumann-potential}
		-\Delta\phi_a=\mathbf a^\circ,
		\qquad
		\partial_n\phi_a|_{\partial\cD}=0,
		\qquad
		\int_\cD\phi_a\dd x=0.
	\end{equation}
	The Dirichlet and mean-zero Neumann estimates give
	\begin{align}
		&\|\phi_c\|_{H^2_x}+\|\nabla\phi_c\|_{L^2(\partial\cD)}
		\le C\|\mathbf c\|_{L^2_x},
		\notag\\
		&\|\phi_j\|_{H^2_x}+\|\nabla\phi_j\|_{L^2(\partial\cD)}
		\le C\|\mathbf b_j\|_{L^2_x},
		\label{eq:dynamic-macro-Dirichlet-elliptic}\\
		&\|\phi_a\|_{H^2_x}+\|\nabla\phi_a\|_{L^2(\partial\cD)}
		\le C\|\mathbf a^\circ\|_{L^2_x}.
		\label{eq:dynamic-macro-Neumann-elliptic}
	\end{align}
	
	The Maxwellian-polynomial test functions are
	\begin{equation}
		\label{eq:dynamic-macro-c-test}
		\Psi_c:=\sqrt2\,(2|v|^2-5)\,v\cdot\nabla_x\phi_c\sqrt\mu,
		\qquad
		\Psi_a:=2\sqrt2\,(5-|v|^2)\,v\cdot\nabla_x\phi_a\sqrt\mu,
	\end{equation}
	and, for $j=1,2,3$,
	\begin{align}
		\Psi_b^j:={}&\sum_{m\ne j}
		\left\{4|v|^2v_mv_j\partial_m\phi_j
		-\frac72(2v_m^2-1)\partial_j\phi_j\right\}\sqrt\mu
		+\frac72(2v_j^2-1)\partial_j\phi_j\sqrt\mu.
		\label{eq:dynamic-macro-b-test}
	\end{align}
	Since $\mu=\ee^{-|v|^2}$ factors over the components, parity and the
	one-dimensional moments $\int_\bR s^{2k}\ee^{-s^2}\dd s
	=\tfrac12\sqrt\pi,\ \tfrac34\sqrt\pi,\ \tfrac{15}8\sqrt\pi$ for $k=1,2,3$
	give
	\begin{align}
		\int(2|v|^2-5)v_iv_k\,\mu\dd v&=0,
		&\int(|v|^2-\tfrac32)(2|v|^2-5)v_iv_k\,\mu\dd v
		&=\tfrac52\pi^{3/2}\delta_{ik},
		\notag\\
		\int(5-|v|^2)v_iv_k\,\mu\dd v&=\tfrac54\pi^{3/2}\delta_{ik},
		&\int(|v|^2-\tfrac32)(5-|v|^2)v_iv_k\,\mu\dd v&=0 .
		\label{eq:dynamic-macro-Gaussian-identities}
	\end{align}
	For \eqref{eq:dynamic-macro-b-test}, the coefficient of
	$(\mathbf b)_j$ in
	$-(\mathbf P f,v\cdot\nabla_x\Psi_b^j)_v$ is
	$-7\Delta\phi_j$, while the cross terms cancel:
	$-7\partial_{jm}\phi_j+7\partial_{mj}\phi_j=0$.  Hence
	\begin{align}
		- (\mathbf P f,v\cdot\nabla_x\Psi_c)_{x,v}
		&=\tfrac{5\sqrt2}2\pi^{3/2}\|\mathbf c\|_{L^2_x}^2,
		\label{eq:dynamic-macro-c-isolation}\\
		- (\mathbf P f,v\cdot\nabla_x\Psi_b^j)_{x,v}
		&=\tfrac72\pi^{3/2}\|(\mathbf b)_j\|_{L^2_x}^2,
		\label{eq:dynamic-macro-b-isolation}\\
		- (\mathbf P f,v\cdot\nabla_x\Psi_a)_{x,v}
		&=\tfrac{5\sqrt2}2\pi^{3/2}\|\mathbf a^\circ\|_{L^2_x}^2.
		\label{eq:dynamic-macro-a-isolation}
	\end{align}
	The last identity uses \eqref{eq:dynamic-macro-Neumann-potential}:
	\begin{equation*}
		-\int_\cD \mathbf a\Delta\phi_a\dd x
		=\int_\cD \mathbf a \mathbf a^\circ\dd x
		=\|\mathbf a^\circ\|_{L^2_x}^2.
	\end{equation*}
	
	For any fixed polynomial $P$,
	\begin{align}
		|((\mathbf I-\mathbf P)f,P\sqrt\mu)_v|
		&\le
		\|\langle v\rangle^{-1/2}(\mathbf I-\mathbf P)f\|_{L^2_v}
		\|\langle v\rangle^{1/2}P\sqrt\mu\|_{L^2_v}
		\notag\\
		&\le C_P\|(\mathbf I-\mathbf P)f\|_{H_\sigma^1}
		\label{eq:dynamic-macro-polynomial-weight-check}
	\end{align}
	by \eqref{eq:local-H1-H_sigma}.  Also
	\begin{equation}
		\label{eq:dynamic-macro-source-weight-check}
		|\langle S,P\sqrt\mu\rangle_v|
		\le \|S\|_{H_\sigma^{-1}}\|P\sqrt\mu\|_{H_\sigma^1}
		\le C_P\|S\|_{H_\sigma^{-1}}.
	\end{equation}
	Together with \eqref{boundedness of L},
	\eqref{eq:dynamic-macro-Dirichlet-elliptic}, and
	\eqref{eq:dynamic-macro-Neumann-elliptic}, these give, for every one of the
	five tests,
	\begin{align}
		&|\langle S,\Psi_\bullet\rangle_{x,v}|
		+|((\mathbf I-\mathbf P)f,v\cdot\nabla_x\Psi_\bullet)_{x,v}|
		+|\langle L(\mathbf I-\mathbf P)f,\Psi_\bullet\rangle_{x,v}|
		\notag\\
		&\qquad\le C
		\bigl(\|S\|_{L^2(\cD;H_\sigma^{-1})}
		+\|(\mathbf I-\mathbf P)f\|_{L^2(\cD;H_\sigma^1)}\bigr)
		\|\bullet\|_{L^2_x}.
		\label{eq:dynamic-macro-bulk-test-bound}
	\end{align}
	
	Taking the collision-invariant moments of
	\eqref{eq:dynamic-linear-source-equation}, using \eqref{eq:L-kernel} and the
	symmetry \eqref{eq:L-symmetry} of the Landau bilinear form, gives the system
	\eqref{eq:appendix-moment-a}--\eqref{eq:appendix-moment-c} of
	\S\ref{subsec:macro-landau}, written here for $f$:
	\begin{align}
		\partial_t \mathbf a+\tfrac12\nabla_x\cdot\mathbf b&=S_a,
		\label{eq:dynamic-moment-a}\\
		\partial_t(\mathbf b)_i
		+\partial_i(\mathbf a+\mathbf c)
		+\partial_j\Theta_{ij}[(\mathbf I-\mathbf P)f]&=S_{b,i},
		\label{eq:dynamic-moment-b}\\
		\partial_t\mathbf c+\frac13\nabla_x\cdot\mathbf b
		+\nabla_x\cdot\Lambda[(\mathbf I-\mathbf P)f]&=S_c,
		\label{eq:dynamic-moment-c}
	\end{align}
	where
	\begin{align}
		S_a&:=\pi^{-3/2}\langle S,\sqrt\mu\rangle_v,\notag\\
		S_{b,i}&:=2\pi^{-3/2}\langle S,v_i\sqrt\mu\rangle_v,
		\notag\\
		S_c&:=(\tfrac32\pi^{3/2})^{-1}
		\langle S,(|v|^2-\tfrac32)\sqrt\mu\rangle_v,
		\label{eq:dynamic-source-moments}\\
		\Theta_{ij}[(\mathbf I-\mathbf P)f]&:=2\pi^{-3/2}
		((\mathbf I-\mathbf P)f,v_iv_j\sqrt\mu)_v,
		\notag\\
		\Lambda_j[(\mathbf I-\mathbf P)f]&:=(\tfrac32\pi^{3/2})^{-1}
		((\mathbf I-\mathbf P)f,(|v|^2-\tfrac32)v_j\sqrt\mu)_v.\nonumber
	\end{align}
	Here \eqref{eq:dynamic-macro-polynomial-weight-check} and
	\eqref{eq:dynamic-macro-source-weight-check} give
	\begin{align*}
		\|\Theta[(\mathbf I-\mathbf P)f]\|_{L^2_x}+\|\Lambda[(\mathbf I-\mathbf P)f]\|_{L^2_x}
		\,&\le C\|(\mathbf I-\mathbf P)f\|_{L^2(\cD;H_\sigma^1)},
		\\
		\|S_a\|_{L^2_x}+\|S_b\|_{L^2_x}+\|S_c\|_{L^2_x}
		\,&\le C\|S\|_{L^2(\cD;H_\sigma^{-1})}.
	\end{align*}
	Differentiating \eqref{eq:dynamic-macro-Dirichlet-potentials} and using
	\eqref{eq:dynamic-moment-b}--\eqref{eq:dynamic-moment-c} with the
	preceding bounds, we obtain
	\begin{align}
		\|\nabla\partial_t\phi_c\|_{L^2_x}
		&\le C\bigl(
		\|\mathbf b\|_{L^2_x}+\|(\mathbf I-\mathbf P)f\|_{L^2(\cD;H_\sigma^1)}
		+\|S\|_{L^2(\cD;H_\sigma^{-1})}\bigr),
		\label{eq:dynamic-time-c-potential}\\
		\sum_{j=1}^3\|\nabla\partial_t\phi_j\|_{L^2_x}
		&\le C\bigl(
		\|\mathbf a^\circ\|_{L^2_x}+|\cD|^{-1/2}|\mathfrak m_f|
		+\|\mathbf c\|_{L^2_x}+\|(\mathbf I-\mathbf P)f\|_{L^2(\cD;H_\sigma^1)}
		+\|S\|_{L^2(\cD;H_\sigma^{-1})}\bigr).
		\label{eq:dynamic-time-b-potential}
	\end{align}
	
	For the Neumann potential, recall the boundary mass scalars
	\eqref{eq:dynamic-boundary-mass-scalars}.
	By \eqref{eq:outgoing-diffuse-moment},
	\eqref{eq:wall-response-profiles}, and
	\eqref{eq:wall-defect-flux-orthogonality}, the normal mass flux of
	\eqref{eq:dynamic-admissible-boundary-law_111} is
	\begin{equation}
		\label{eq:dynamic-alpha-mass-flux}
		j_f(t,x):=
		\int_{\bR_v^3}f_\gamma(t,x,v)\sqrt\mu(v)n_x\cdot v\dd v
		=\varphi_f-\alpha\cJ\varphi_f-\beta_b.
	\end{equation}
	Indeed, the isothermal part of $\alpha P_\cJ f_{\gamma^+}$ carries the incoming mass
	flux $\alpha\cJ\varphi_f$, the part $\alpha d_\delta\cJ\varphi_f$ carries none, and the residual contributes $-\beta_b$.
	Applying \eqref{eq:dynamic-Green-identity} with the time-independent test
	$\psi=\pi^{-3/2}\zeta(x)\sqrt\mu(v)$, $\zeta\in H^1(\cD)$ of mean zero, and
	using \eqref{eq:L-kernel} and \eqref{eq:dynamic-alpha-mass-flux}, we obtain
	the weak form of \eqref{eq:dynamic-moment-a},
	\begin{align}
		\langle\partial_t\mathbf a^\circ,\zeta\rangle
		&=\tfrac12\int_\cD\mathbf b\cdot\nabla\zeta\dd x
		+\int_\cD S_a\zeta\dd x
		-\pi^{-3/2}\int_{\partial\cD}(\varphi_f-\alpha\cJ\varphi_f)\zeta\dd S_x
		+\pi^{-3/2}\int_{\partial\cD}\beta_b\zeta\dd S_x.
		\label{eq:dynamic-a-circle-weak-time}
	\end{align}
	The variational form of
	\eqref{eq:dynamic-macro-Neumann-potential} gives
	\begin{align}
		\|\nabla\partial_t\phi_a\|_{L^2_x}
		&=\sup_{\substack{\zeta\in H^1(\cD),\ \int_\cD\zeta=0\\
				\|\nabla\zeta\|_{L^2_x}=1}}
		|\langle\partial_t\mathbf a^\circ,\zeta\rangle|.
		\label{eq:dynamic-time-a-variational}
	\end{align}
	For every mean-zero $\zeta\in H^1(\cD)$,
	\[
	\|\zeta\|_{L^2_x}+\|\zeta\|_{L^2(\partial\cD)}
	\le C\|\nabla\zeta\|_{L^2_x}.
	\]
	Combining \eqref{eq:dynamic-a-circle-weak-time},
	\eqref{eq:dynamic-time-a-variational}, and the preceding mean-zero
	Poincar\'e and trace inequality, we obtain
	\begin{align}
		\|\nabla\partial_t\phi_a\|_{L^2_x}
		\le C\Bigl(&\|\mathbf b\|_{L^2_x}
		+\|S\|_{L^2(\cD;H_\sigma^{-1})}
		\notag\\
		&+\|\varphi_f(t)-\alpha\cJ\varphi_f(t)\|_{L^2(\partial\cD)}
		+\|b(t)\|_{L^2(\gamma^-)}\Bigr).
		\label{eq:dynamic-time-a-potential}
	\end{align}
	
	The Gaussian identities in
	\eqref{eq:dynamic-macro-Gaussian-identities} and those used in
	\eqref{eq:dynamic-macro-b-isolation} also give
	\begin{align}
		(\mathbf P f,\partial_t\Psi_c)_{x,v}&=0,\nonumber
		\\
		(\mathbf P f,\partial_t\Psi_b^j)_{x,v}
		&=-\tfrac72\pi^{3/2}\int_\cD \mathbf c\partial_j\partial_t\phi_j\dd x,\nonumber
		\\
		(\mathbf P f,\partial_t\Psi_a)_{x,v}
		&=\tfrac{5\sqrt2}2\pi^{3/2}
		\int_\cD\mathbf b\cdot\nabla\partial_t\phi_a\dd x.
		\label{eq:dynamic-time-a-macro}
	\end{align}
	Using \eqref{eq:dynamic-time-a-macro} and
	\eqref{eq:dynamic-macro-polynomial-weight-check} for the
	$(\mathbf I-\mathbf P)f$ part, then
	\eqref{eq:dynamic-time-c-potential}--\eqref{eq:dynamic-time-a-potential},
	we obtain, for every $0<\varepsilon<1$,
	\begin{align}
		\left|\int_I(f,\partial_t\Psi_c)_{x,v}\dd t\right|
		&\le \varepsilon\|\mathbf b\|_{L^2(I;L_x^2)}^2
		+C_\varepsilon(X_I^2+D_I^2),
		\label{eq:dynamic-time-c-test-bound}\\
		\sum_j\left|\int_I(f,\partial_t\Psi_b^j)_{x,v}\dd t\right|
		&\le \varepsilon\|\mathbf a^\circ\|_{L^2(I;L_x^2)}^2
		+C_\varepsilon\Bigl(
		\|\mathbf c\|_{L^2(I;L_x^2)}^2+X_I^2+D_I^2
		+\int_I|\mathfrak m_f|^2\dd t\Bigr),
		\label{eq:dynamic-time-b-test-bound}\\
		\left|\int_I(f,\partial_t\Psi_a)_{x,v}\dd t\right|
		&\le C\|\mathbf b\|_{L^2(I;L_x^2)}^2+C(X_I^2+D_I^2)
		+CK_I^2
		+C\|b\|_{L^2(\gamma_I^-)}^2,
		\label{eq:dynamic-time-a-test-bound}
	\end{align}
	the last by \eqref{eq:dynamic-alpha-square-payment}.
	
	The diffuse boundary law separates the rank-one and microscopic contributions.
	We compute the boundary pairings in
	\eqref{eq:dynamic-Green-identity}.  Define
	\begin{equation*}
		\mathfrak B_I[\Psi]
		:=\int_{\gamma_I^+}f_{\gamma^+}\Psi\dd\gamma\dd t
		-\int_{\gamma_I^-}f_{\gamma^-}\Psi\dd\gamma\dd t.
	\end{equation*}
	By \eqref{eq:isothermal-diffuse-projections} and
	\eqref{eq:dynamic-admissible-boundary-law_111},
	\begin{align}
		\mathfrak B_I[\Psi]
		={}&\int_{\gamma_I^+}(\mathbf I-\cP_0^+)f_{\gamma^+}\Psi
		+\int_{\gamma_I^+}\cP_0^+f_{\gamma^+}\Psi
		-\alpha\int_{\gamma_I^-}p_\cJ\Psi
		\notag\\
		&-\alpha\int_{\gamma_I^-}r_\cJ\Psi
		-\int_{\gamma_I^-}b\Psi.
		\label{eq:dynamic-macro-boundary-decomposition}
	\end{align}
	All boundary integrals in \eqref{eq:dynamic-macro-boundary-decomposition}
	use $\dd\gamma\dd t$.
	
	For $\Psi_c$, the Dirichlet condition in
	\eqref{eq:dynamic-macro-Dirichlet-potentials} gives
	$\nabla\phi_c=(\partial_n\phi_c)n$ on $\partial\cD$.  Rotation to $n=e_3$ and
	\eqref{eq:dynamic-macro-Gaussian-identities} give
	\[
	\sqrt2\int_{v_3>0}(2|v|^2-5)v_3^2\,\mu(v)\dd v
	=\frac{\sqrt2}2\int_{\bR_v^3}(2|v|^2-5)v_3^2\,\mu(v)\dd v=0 .
	\]
	Replacing $v_3$ by $-v_3$ in the preceding identity gives zero on $v_3<0$
	as well.
	Thus both isothermal rank-one terms in
	\eqref{eq:dynamic-macro-boundary-decomposition} vanish for $\Psi_c$.
	
	For $\Psi_a$, \eqref{eq:dynamic-macro-Neumann-potential} gives
	$n\cdot\nabla\phi_a=0$ on $\partial\cD$. With $n=e_3$,
	$v\cdot\nabla\phi_a$ is odd in the tangential velocity variables, while
	the remaining factor is even. Hence
	\[
	\int_{v_3>0}2\sqrt2\,(5-|v|^2)(v\cdot\nabla\phi_a)
	\mu(v)v_3\,dv
	=\int_{v_3<0}2\sqrt2\,(5-|v|^2)(v\cdot\nabla\phi_a)
	\mu(v)|v_3|\,dv=0.
	\]
	Both isothermal rank-one terms vanish.
	
	For $\Psi_b^j$, the velocity polynomial in
	\eqref{eq:dynamic-macro-b-test} is even.  Applying the change $v\mapsto-v$ and
	the identity $p_\cJ=\mathfrak m_\circ\cJ\varphi_f$, we obtain
	\[
	\int_{\gamma_I^-}p_\cJ\Psi_b^j\dd\gamma\dd t
	=\int_{\gamma_I^+}\mathfrak m_\circ\cJ\varphi_f\Psi_b^j\dd\gamma\dd t,
	\]
	so the rank-one part of the $\Psi_b^j$ boundary term is
	$\int_{\gamma_I^+}\mathfrak m_\circ(\varphi_f-\alpha\cJ\varphi_f)\Psi_b^j$.
	By \eqref{eq:wall-rank-one-isometry},
	\eqref{eq:dynamic-alpha-square-payment}, and the boundary bound for the
	five tests proved below,
	\begin{align}
		\left|\int_{\gamma_I^+}
		\mathfrak m_\circ(\varphi_f-\alpha\cJ\varphi_f)\Psi_b^j\dd\gamma\dd t\right|
		&\le CK_I\|(\mathbf b)_j\|_{L^2(I;L_x^2)}
		\notag\\
		&\le \varepsilon\|(\mathbf b)_j\|_{L^2(I;L_x^2)}^2
		+C_\varepsilon K_I^2.
		\label{eq:dynamic-b-alpha-residual-payment}
	\end{align}
	The resulting error is controlled by $K_I^2$ in
	\eqref{eq:dynamic-graph-boundary-input}.
	
	For every test in \eqref{eq:dynamic-macro-c-test}--\eqref{eq:dynamic-macro-b-test},
	\eqref{eq:dynamic-macro-Dirichlet-elliptic}--\eqref{eq:dynamic-macro-Neumann-elliptic} give
	\[
	\|\Psi_\bullet\|_{L^2(\gamma_I^+)}
	+\|\Psi_\bullet\|_{L^2(\gamma_I^-)}
	\le C\|\bullet\|_{L^2(I;L_x^2)}.
	\]
	The rank-one terms vanish for $\Psi_c$ and for $\Psi_a$ by the two vanishing
	identities established above,
	and for $\Psi_b^j$ they are bounded by
	\eqref{eq:dynamic-b-alpha-residual-payment}.  Using the preceding trace
	bound and \eqref{eq:wall-J-components},
	Cauchy--Schwarz and Young's inequality then give
	\begin{align}
		|\mathfrak B_I[\Psi_c]|
		&\le \varepsilon\|\mathbf c\|_{L^2(I;L_x^2)}^2
		+C_\varepsilon\bigl(M_I^2+\delta^2B_I^2
		+\|b\|_{L^2(\gamma_I^-)}^2\bigr),
		\label{eq:dynamic-c-boundary-bound}\\
		\sum_j|\mathfrak B_I[\Psi_b^j]|
		&\le \varepsilon\|\mathbf b\|_{L^2(I;L_x^2)}^2
		+C_\varepsilon\bigl(M_I^2+K_I^2+\delta^2B_I^2
		+\|b\|_{L^2(\gamma_I^-)}^2\bigr),\nonumber
		\\
		|\mathfrak B_I[\Psi_a]|
		&\le \varepsilon\|\mathbf a^\circ\|_{L^2(I;L_x^2)}^2
		+C_\varepsilon\bigl(M_I^2+\delta^2B_I^2
		+\|b\|_{L^2(\gamma_I^-)}^2\bigr).
		\label{eq:dynamic-a-boundary-bound}
	\end{align}
	
	Collect the error terms in
	\begin{align*}
		\cR_I:={}&X_I^2+D_I^2+M_I^2
		+K_I^2+\delta^2B_I^2
		+\|b\|_{L^2(\gamma_I^-)}^2
		+\int_I|\mathfrak m_f(t)|^2\dd t.
	\end{align*}
	Put
	\begin{equation*}
		\kappa_c:=\tfrac{5\sqrt2}2\pi^{3/2},
		\qquad
		\kappa_b:=\tfrac72\pi^{3/2},
		\qquad
		\kappa_a:=\tfrac{5\sqrt2}2\pi^{3/2},
	\end{equation*}
	and define
	\begin{equation*}
		\mathscr E_c(t):=\kappa_c^{-1}(f(t),\Psi_c(t))_{x,v},
		\qquad
		\mathscr E_b(t):=\kappa_b^{-1}\sum_{j=1}^3
		(f(t),\Psi_b^j(t))_{x,v},
		\qquad
		\mathscr E_a(t):=\kappa_a^{-1}(f(t),\Psi_a(t))_{x,v}.
	\end{equation*}
	For a test $\Psi$, set
	\begin{align*}
		&\mathfrak F_I[\Psi]\\
		:={}&
		\int_I\Bigl\{
		\langle S,\Psi\rangle_{x,v}
		-\langle L(\mathbf I-\mathbf P)f,\Psi\rangle_{x,v}
		+((\mathbf I-\mathbf P)f,v\cdot\nabla_x\Psi)_{x,v}
		+(f,\partial_t\Psi)_{x,v}
		\Bigr\}\dd t
		-\mathfrak B_I[\Psi].
	\end{align*}
	Applying \eqref{eq:dynamic-Green-identity},
	\eqref{eq:dynamic-linear-source-equation}, and
	\eqref{eq:dynamic-macro-c-isolation}--\eqref{eq:dynamic-macro-a-isolation}, we obtain
	\begin{align}
		\kappa_c\|\mathbf c\|_{L^2(I;L_x^2)}^2
		&=\mathfrak F_I[\Psi_c]
		+\kappa_c\{\mathscr E_c(t_0)-\mathscr E_c(t_1)\},
		\label{eq:dynamic-c-mode-exact}\\
		\kappa_b\|\mathbf b\|_{L^2(I;L_x^2)}^2
		&=\sum_{j=1}^3\mathfrak F_I[\Psi_b^j]
		+\kappa_b\{\mathscr E_b(t_0)-\mathscr E_b(t_1)\},\label{eq:dynamic-b-mode-exact}
		\\
		\kappa_a\|\mathbf a^\circ\|_{L^2(I;L_x^2)}^2
		&=\mathfrak F_I[\Psi_a]
		+\kappa_a\{\mathscr E_a(t_0)-\mathscr E_a(t_1)\}.
		\label{eq:dynamic-a-mode-exact}
	\end{align}
	Choose the Young parameters in
	\eqref{eq:dynamic-macro-bulk-test-bound},
	\eqref{eq:dynamic-time-c-test-bound}--\eqref{eq:dynamic-time-a-test-bound}, and
	\eqref{eq:dynamic-c-boundary-bound}--\eqref{eq:dynamic-a-boundary-bound} so that
	\begin{align}
		|\mathfrak F_I[\Psi_c]|
		&\le \frac{\kappa_c}{2}\|\mathbf c\|_{L^2(I;L_x^2)}^2
		+\frac{\kappa_c\varepsilon_c}{2}
		\|\mathbf b\|_{L^2(I;L_x^2)}^2
		+C(1+\varepsilon_c^{-1})\cR_I,
		\label{eq:dynamic-c-mode-Q-bound}\\
		\left|\sum_{j=1}^3\mathfrak F_I[\Psi_b^j]\right|
		&\le \frac{\kappa_b}{2}\|\mathbf b\|_{L^2(I;L_x^2)}^2
		+\frac{\kappa_b\varepsilon_b}{2}
		\|\mathbf a^\circ\|_{L^2(I;L_x^2)}^2
		+C(1+\varepsilon_b^{-1})
		\bigl(\|\mathbf c\|_{L^2(I;L_x^2)}^2+\cR_I\bigr),\nonumber
		\\
		|\mathfrak F_I[\Psi_a]|
		&\le \frac{\kappa_a}{2}\|\mathbf a^\circ\|_{L^2(I;L_x^2)}^2
		+C\bigl(\|\mathbf b\|_{L^2(I;L_x^2)}^2+\cR_I\bigr).
		\label{eq:dynamic-a-mode-Q-bound}
	\end{align}
	Subtracting the first terms on the right
	of \eqref{eq:dynamic-c-mode-Q-bound}--\eqref{eq:dynamic-a-mode-Q-bound} from
	\eqref{eq:dynamic-c-mode-exact}--\eqref{eq:dynamic-a-mode-exact} and dividing by
	$\kappa_c/2$, $\kappa_b/2$, and $\kappa_a/2$, we obtain
	\begin{align}
		\|\mathbf c\|_{L^2(I;L_x^2)}^2
		&\le \varepsilon_c\|\mathbf b\|_{L^2(I;L_x^2)}^2
		+C(1+\varepsilon_c^{-1})\cR_I
		+2\{\mathscr E_c(t_0)-\mathscr E_c(t_1)\},
		\label{eq:dynamic-c-mode-estimate}\\
		\|\mathbf b\|_{L^2(I;L_x^2)}^2
		&\le \varepsilon_b\|\mathbf a^\circ\|_{L^2(I;L_x^2)}^2
		+C(1+\varepsilon_b^{-1})
		\bigl(\|\mathbf c\|_{L^2(I;L_x^2)}^2+\cR_I\bigr)
		+2\{\mathscr E_b(t_0)-\mathscr E_b(t_1)\},
		\label{eq:dynamic-b-mode-estimate}\\
		\|\mathbf a^\circ\|_{L^2(I;L_x^2)}^2
		&\le C\bigl(\|\mathbf b\|_{L^2(I;L_x^2)}^2+\cR_I\bigr)
		+2\{\mathscr E_a(t_0)-\mathscr E_a(t_1)\}.
		\label{eq:dynamic-a-mode-estimate}
	\end{align}
	By the elliptic estimates and Cauchy--Schwarz, we have
	\[
	|\mathscr E_a(t)|+|\mathscr E_b(t)|+|\mathscr E_c(t)|
	\le C\|f(t)\|_{L^2(\cD\times\bR_v^3)}^2.
	\]
	
	The three mode estimates are combined with fixed positive weights.
	\begingroup
	Choose $\varepsilon_b>0$ small, put
	$K:=C(1+\varepsilon_b^{-1})$, and then choose
	$\varepsilon_c>0$ so that
	$2C\varepsilon_b\le\tfrac14$ and $2K\varepsilon_c\le\tfrac14$.
	Multiply \eqref{eq:dynamic-a-mode-estimate} by $2\varepsilon_b$
	and \eqref{eq:dynamic-c-mode-estimate} by $2K$, then add
	\eqref{eq:dynamic-b-mode-estimate}.
	\endgroup
	Then
	\begin{align}
		&\varepsilon_b\|\mathbf a^\circ\|_{L^2(I;L_x^2)}^2
		+\frac12\|\mathbf b\|_{L^2(I;L_x^2)}^2
		+K\|\mathbf c\|_{L^2(I;L_x^2)}^2
		\le C\cR_I
		+\varepsilon_b\{\eta_f(t_0)-\eta_f(t_1)\},
		\label{eq:dynamic-macro-weighted-closing}
	\end{align}
	with
	\begin{equation}
		\label{eq:dynamic-macro-interaction-definition}
		\eta_f:=4\mathscr E_a+2\varepsilon_b^{-1}\mathscr E_b
		+4K\varepsilon_b^{-1}\mathscr E_c.
	\end{equation}
	The elliptic bounds for the test functions and
	\eqref{eq:dynamic-macro-interaction-definition} give
	\eqref{eq:dynamic-macro-eta-bound}.  Since
	$0<\varepsilon_b\le1/2\le K$, division of
	\eqref{eq:dynamic-macro-weighted-closing} by $\varepsilon_b$ yields
	\begin{align*}
		&\|\mathbf a^\circ\|_{L^2(I;L_x^2)}^2
		+\|\mathbf b\|_{L^2(I;L_x^2)}^2
		+\|\mathbf c\|_{L^2(I;L_x^2)}^2
		\le C\cR_I+\eta_f(t_0)-\eta_f(t_1),
	\end{align*}
	which is \eqref{eq:dynamic-alpha-macro-estimate}.
\end{proof}

\begin{lemma}[Mass scalar along the $\alpha$-paths]
	\label{lem:dynamic-alpha-mass-scalar}
	Under the hypotheses of
	Lemma~\ref{lem:dynamic-alpha-macroscopic-estimate}, on both paths the mass
	scalar \eqref{eq:dynamic-mass-scalar} satisfies 
	\begin{equation}
		\label{eq:dynamic-mass-balance-exact}
		\frac{\dd}{\dd t}\mathfrak m_f(t)
		=\pi^{3/2}\int_\cD S_a(t,x)\dd x
		-(1-\alpha)\int_{\partial\cD}\varphi_f(t,x)\dd S_x
		+\int_{\partial\cD}\beta_b(t,x)\dd S_x
	\end{equation}
	in distributions on $I$, where $\varphi_f$ is the function defined in \eqref{eq:dynamic-boundary-mass-scalars}.
	In particular, with $\tau:=t_1-t_0$,
	\begin{align}
		\int_I|\mathfrak m_f(t)|^2\dd t
		\le{}&3\tau|\mathfrak m_f(t_0)|^2
		+C\tau^2\|S\|_{L^2(\cD_I;H_\sigma^{-1})}^2
		+C\tau^2K_I^2
		+C\tau^2\|b\|_{L^2(\gamma_I^-)}^2,
		\label{eq:dynamic-mass-finite-I}
	\end{align}
	with $K_I$ from \eqref{eq:dynamic-alpha-K}, where the constant $C$ depends only on $\cD$.  At $\alpha=1$,
	\begin{equation}
		\label{eq:dynamic-physical-mass-balance}
		\mathfrak m_f(t)=\mathfrak m_f(t_0)
		+\pi^{3/2}\int_{t_0}^t\int_\cD S_a(r,x)\dd x\dd r
		+\int_{t_0}^t\int_{\partial\cD}\beta_b(r,x)\dd S_x\dd r.
	\end{equation}
\end{lemma}

\begin{proof}
	Self-adjointness of $\cJ$
	and $\cJ1=1$ give
	\begin{equation*}
		\int_{\partial\cD}(\mathbf I-\cJ)\varphi_f\dd S_x
		=\langle\varphi_f,(\mathbf I-\cJ)1\rangle_{L^2(\partial\cD)}=0.
	\end{equation*}
	Integrating \eqref{eq:dynamic-linear-source-equation} against $\sqrt\mu$ and using
	Lemma~\ref{260715201} yields \eqref{eq:dynamic-mass-balance-exact}.
	
	Integrate \eqref{eq:dynamic-mass-balance-exact} from $t_0$ to $t\in I$.
	Then \eqref{eq:dynamic-source-moments},
	\eqref{eq:dynamic-macro-source-weight-check},
	\eqref{eq:dynamic-varphi-bound}, and \eqref{eq:dynamic-beta-b-bound} give
	\begin{align}
		|\mathfrak m_f(t)|
		\le{}&|\mathfrak m_f(t_0)|
		+C\tau^{1/2}D_I
		+C(1-\alpha)\tau^{1/2}\|f_{\gamma^+}\|_{L^2(\gamma_I^+)}
		+C\tau^{1/2}\|b\|_{L^2(\gamma_I^-)}.
		\label{eq:dynamic-mass-pointwise-I}
	\end{align}
	Since $0\le\alpha\le1$,
	\begin{equation*}
		(1-\alpha)^2\|f_{\gamma^+}\|_{L^2(\gamma_I^+)}^2
		\le(1-\alpha^2)\|f_{\gamma^+}\|_{L^2(\gamma_I^+)}^2
		\le K_I^2
	\end{equation*}
	by \eqref{eq:dynamic-alpha-K}.  Squaring
	\eqref{eq:dynamic-mass-pointwise-I} and integrating over $I$ therefore gives
	\eqref{eq:dynamic-mass-finite-I}.  Setting $\alpha=1$ in
	\eqref{eq:dynamic-mass-balance-exact} and integrating gives
	\eqref{eq:dynamic-physical-mass-balance}.
\end{proof}

\subsection{The evolutionary problem}\label{ss:diffuse-dynamic-linear}

Fix $0<T<\infty$, and write
\[
\mathbb X_T:=\cH_{\kin,\sigma}(\cD_T\times\bR_v^3),
\qquad
\mathbb Y_T:=L^2(\cD_T;H_\sigma^{-1})\times L^2(\cD\times\bR_v^3)
\times L^2(\gamma_T^-).
\]

\begin{prop}[Evolutionary $\alpha$-uniform graph estimates]
	\label{prop:dynamic-alpha-uniform-graph}
	There exist $\varepsilon_{\rm dyn},\delta_{\rm dyn}>0$ with
	\begin{equation}
		\label{eq:dynamic-KS-perturbed-smallness-choice}
		\varepsilon_{\rm dyn}\le\varepsilon_0(\cD),
	\end{equation}
	where $\varepsilon_0(\cD)$ is the constant of Theorem~\ref{260322442},
	such that the following holds:  Let $0\le\alpha\le1$ and assume
	\begin{equation}
		\label{eq:dynamic-alpha-graph-smallness}
		\|\mu^{1/4}g\|_{L^\infty((0,T)\times\cD\times\bR_v^3)}
		\le\varepsilon_{\rm dyn},
		\qquad
		0<\delta\le\delta_{\rm dyn}.
	\end{equation}
	Let
	\[
	F\in L^2(\cD_T;H_\sigma^{-1}),
	\qquad
	f_0\in L^2(\cD\times\bR_v^3),
	\qquad
	b\in L^2(\gamma_T^-),
	\]
	and let $f\in\mathbb X_T$ satisfy
	\begin{equation}
		\label{eq:dynamic-alpha-frozen-problem}
		(\partial_t+v\cdot\nabla_x)f+L_gf=F,
		\qquad
		f(0)=f_0,
	\end{equation}
	and either
	\begin{align}\label{it:graph-physical}
		f_{\gamma^-}=\alpha P_wf_{\gamma^+}+b
	\end{align}
	or
	\begin{align}\label{it:graph-spectral}
		f_{\gamma^-}=\alpha P_{w,N}f_{\gamma^+}+b
	\end{align}
	for some $N\in\bN$.
	Then
	\begin{align}
		&\sup_{0\le t\le T}\|f(t)\|_{L^2(\cD\times\bR_v^3)}
		+\|f\|_{\mathbb X_T}
		+\|f_{\gamma^+}\|_{L^2(\gamma_T^+)}
		\notag\\
		&\qquad\le C_T\Bigl(
		\|f_0\|_{L^2(\cD\times\bR_v^3)}
		+\|F\|_{L^2(\cD_T;H_\sigma^{-1})}
		+\|b\|_{L^2(\gamma_T^-)}\Bigr),
		\label{eq:dynamic-alpha-uniform-graph}
	\end{align}
	where $C_T$ is independent of $\alpha$, and of $N$ in case \eqref{it:graph-spectral}.
	Moreover, $C_T$ depends only on $\cD$ and on $T$, is nondecreasing in $T$, and
	depends on $g$ only through \eqref{eq:dynamic-alpha-graph-smallness}.
	
	If in addition $\alpha=1$, then
	\begin{align}
		&\sup_{0\le t\le T}\|f(t)\|_{L^2(\cD\times\bR_v^3)}
		+\|f\|_{\mathbb X_T}
		+\|f_{\gamma^+}\|_{L^2(\gamma_T^+)}
		\notag\\
		&\qquad\le C\Bigl(
		\|f_0\|_{L^2(\cD\times\bR_v^3)}
		+\|F\|_{L^2(\cD_T;H_\sigma^{-1})}
		+\|b\|_{L^2(\gamma_T^-)}
		+\|\mathfrak m_f\|_{L^2(0,T)}\Bigr),
		\label{eq:dyn-hier-mass-graph}
	\end{align}
	with $C$ depending only on $\cD$ (not on $T$ or $\cJ$).
	Moreover, for $\alpha=1$: if
	$\iint_{\cD\times \bR^3_v}f_0(x,v)\,\sqrt \mu\dd v\dd x=0$ and
	\begin{equation}
		\label{eq:dynamic-physical-zero-mass-data}
		\int_\cD\langle F(t,x,v),\sqrt\mu\rangle_v\dd x=0,
		\quad
		\int_{\partial\cD}\beta_b(t,x)\dd S_x=0
	\end{equation}
	for a.e.\ $t$, with $\beta_b$ the incoming mass scalar
	\eqref{eq:dynamic-boundary-mass-scalars}, then $\mathfrak m_f\equiv 0$ on $[0,T]$,
	so \eqref{eq:dyn-hier-mass-graph} holds without the term
	$\|\mathfrak m_f\|_{L^2(0,T)}$.
\end{prop}

\begin{proof}
	We treat the physical and boundary-spectral laws simultaneously, writing
	$\cJ=\mathbf I$ or $\cJ=J_N$. The local estimates hold on any
	$I=(t_0,t_1)\subset(0,T)$, with constants independent of $I$, $T$,
	$\alpha$, and $\cJ$. Write $\tau=t_1-t_0$ and use the quantities
	\eqref{eq:dynamic-window-quantities}--\eqref{eq:dynamic-window-traces}. Set
	\[
	\mathsf H_I:=\|F\|_{L^2(\cD_I;H_\sigma^{-1})}^2
	+\|b\|_{L^2(\gamma_I^-)}^2,
	\qquad Q_I^2:=X_I^2+M_I^2+K_I^2.
	\]
	The finite-dimensionality of $\ker L$ gives
	\[
	\|\mathbf Pf\|_{L^2(\cD_I;H_\sigma^1)}\le CA_I,
	\qquad \|f\|_{L^2(\cD_I;H_\sigma^1)}^2\le C(A_I^2+X_I^2).
	\]
	The quadratic Green identity, microscopic coercivity, and
	\eqref{eq:dynamic-graph-boundary-input} imply, for $0<\lambda,\theta<1$,
	\begin{align}
		\|f(t_1)\|_2^2+cQ_I^2
		&\le \|f(t_0)\|_2^2
		+C(\varepsilon_{\rm dyn}+\lambda)(A_I^2+X_I^2)
		+C_\lambda\|F\|_{L^2(\cD_I;H_\sigma^{-1})}^2\notag\\
		&\quad+C(\theta+\delta^2)B_I^2
		+C_\theta\|b\|_{L^2(\gamma_I^-)}^2.
		\label{eq:r38-energy}
	\end{align}
	Here and below $\|\cdot\|_2$ at a fixed time denotes the
	$L^2(\cD\times\bR_v^3)$ norm. The positive constant $c$ and the
	coefficient of $B_I^2$ are independent of $\lambda$ and $\theta$.
	The dependence of $C_\lambda$ and $C_\theta$ is confined to the data terms.
	
	Put $S=F+N[g,f]$. Lemma~\ref{260715201} gives
	\begin{equation}
		\|S\|_{L^2(\cD_I;H_\sigma^{-1})}^2
		\le 2\|F\|_{L^2(\cD_I;H_\sigma^{-1})}^2
		+C\varepsilon_{\rm dyn}^2(A_I^2+X_I^2).
		\label{eq:r38-source}
	\end{equation}
	By Lemma~\ref{lem:dynamic-alpha-macroscopic-estimate} and
	\eqref{eq:dynamic-macro-norm-exact}, the term
	$C\varepsilon_{\rm dyn}^2 A_I^2$ can be absorbed after decreasing
	$\varepsilon_{\rm dyn}$. Multiplying the interaction functional from that
	lemma by a fixed constant, we obtain
	\begin{equation}
		A_I^2\le C\left(Q_I^2+\delta^2B_I^2+\mathsf H_I
		+\|\mathfrak m_f\|_{L^2(I)}^2\right)
		+\eta_f(t_0)-\eta_f(t_1),
		\qquad |\eta_f(t)|\le C\|f(t)\|_2^2.
		\label{eq:r38-macro}
	\end{equation}
	The collision term has zero mass moment. Thus
	Lemma~\ref{lem:dynamic-alpha-mass-scalar}, applied to $S$, yields
	\begin{equation}
		\|\mathfrak m_f\|_{L^2(I)}^2
		\le C\tau\|f(t_0)\|_2^2
		+C\tau^2\left(\mathsf H_I+K_I^2
		+\varepsilon_{\rm dyn}^2(A_I^2+X_I^2)\right).
		\label{eq:r38-mass}
	\end{equation}
	The constants in \eqref{eq:r38-energy}--\eqref{eq:r38-mass} depend only on
	the fixed domain and collision estimates.
	
	Choose $\kappa>0$ so small that
	\begin{equation}
		\mathcal E_f(t):=\|f(t)\|_2^2+\kappa\eta_f(t)
		\quad\text{satisfies}\quad
		\tfrac12\|f(t)\|_2^2\le\mathcal E_f(t)
		\le\tfrac32\|f(t)\|_2^2,
		\label{eq:r38-equivalent-energy}
	\end{equation}
	and the contribution $C\kappa Q_I^2$ in
	$\kappa\,\eqref{eq:r38-macro}$ is absorbable by the left side of
	\eqref{eq:r38-energy}. Next fix $\lambda>0$ and decrease
	$\varepsilon_{\rm dyn}$ so that the terms involving
	$(\varepsilon_{\rm dyn}+\lambda)(A_I^2+X_I^2)$ and
	$\kappa\varepsilon_{\rm dyn}^2(A_I^2+X_I^2)$ can be absorbed.
	Finally choose $0<\tau_*\le1$ so that the contribution
	$C\kappa\tau_*^2K_I^2$ from \eqref{eq:r38-mass} is also absorbable.
	To make the order of these choices explicit, one may take a fixed
	structural $C$ larger than the constants in the three estimates and their products, and require
	\begin{equation*}
		C\kappa\le \tfrac18c,
		\qquad C(\varepsilon_{\rm dyn}+\lambda)\le\tfrac18\kappa,
		\qquad C\varepsilon_{\rm dyn}^2\le\tfrac18,
		\qquad C\kappa\tau_*^2\le\tfrac18c.
	\end{equation*}
	The same $C$ is used only in these sufficient inequalities;
	none of the choices depends on $I$, $T$, $\alpha$, or $N$.
	
	For $\tau\le\tau_*$, add $\kappa$ times
	\eqref{eq:r38-macro} to \eqref{eq:r38-energy} and substitute
	\eqref{eq:r38-mass}. The $\eta_f$ terms combine into
	$\mathcal E_f(t_1)-\mathcal E_f(t_0)$.
	The choices above absorb all terms containing $A_I^2$ or $Q_I^2$ on
	the right. Consequently, with $c_*>0$ independent of $\theta$,
	\begin{equation}
		\mathcal E_f(t_1)+c_*(A_I^2+Q_I^2)
		\le C\bigl(\mathcal E_f(t_0)+\mathsf H_I\bigr)
		+C(\theta+\delta^2)B_I^2.
		\label{eq:r38-before-trace}
	\end{equation}
	
	The equation, \eqref{boundedness of L}, and Lemma~\ref{260715201}
	give
	\begin{equation}
		\|(\partial_t+v\cdot\nabla_x)f\|_{L^2(\cD_I;H_\sigma^{-1})}
		\le C\bigl(\|F\|_{L^2(\cD_I;H_\sigma^{-1})}+A_I+X_I\bigr).
		\label{eq:r38-transport}
	\end{equation}
	The evolutionary part of
	Proposition~\ref{prop:full-diffuse-trace-recovery} therefore yields
	\begin{equation}
		B_I^2\le C\bigl(\|f(t_0)\|_2^2+A_I^2+X_I^2+M_I^2
		+\|F\|_{L^2(\cD_I;H_\sigma^{-1})}^2\bigr).
		\label{eq:r38-trace-recovery}
	\end{equation}
	Fix $\theta>0$ and then $\delta_{\rm dyn}>0$ small enough that
	substitution of \eqref{eq:r38-trace-recovery} into
	\eqref{eq:r38-before-trace} absorbs its $A_I^2$, $X_I^2$, and $M_I^2$
	contributions. Adding a sufficiently small multiple of
	\eqref{eq:r38-trace-recovery} to the resulting estimate gives
	\begin{equation}
		\mathcal E_f(t_1)+c_0\bigl(A_I^2+Q_I^2+B_I^2\bigr)
		\le C\bigl(\mathcal E_f(t_0)+\mathsf H_I\bigr),
		\qquad \tau\le\tau_*.
		\label{eq:r38-local-closed}
	\end{equation}
	The positive constant $c_0$ is fixed independently of $I$, $T$,
	$\alpha$, and $N$. The boundary law and
	\eqref{eq:wall-J-operator-bound} give
	\[
	\|f_{\gamma^-}\|_{L^2(\gamma_I^-)}\le CB_I+\|b\|_{L^2(\gamma_I^-)}.
	\]
	Apply \eqref{eq:r38-local-closed} to $(t_0,t)$ for every
	$t\in[t_0,t_1]$, and use \eqref{eq:r38-transport} and
	\eqref{eq:r38-equivalent-energy}. We obtain
	\begin{align}
		&\sup_{t_0\le t\le t_1}\|f(t)\|_2^2
		+\|f\|_{L^2(\cD_I;H_\sigma^1)}^2
		+\|(\partial_t+v\cdot\nabla_x)f\|_{L^2(\cD_I;H_\sigma^{-1})}^2
		+\|f_{\gamma^+}\|_{L^2(\gamma_I^+)}^2
		+\|f_{\gamma^-}\|_{L^2(\gamma_I^-)}^2
		\notag\\
		&\hspace{10mm}\le C\bigl(\|f(t_0)\|_2^2+\mathsf H_I\bigr),
		\qquad \tau\le\tau_*.
		\label{eq:r38-local-graph}
	\end{align}
	Partition $[0,T]$ into at most $1+T/\tau_*$ consecutive intervals
	of length at most $\tau_*$. Iterating \eqref{eq:r38-local-graph} at
	successive endpoints and summing the interval estimates proves
	\eqref{eq:dynamic-alpha-uniform-graph}. The resulting constant $C_T$
	is nondecreasing in $T$ and independent of $\alpha$ and $N$.
	
	When $\alpha=1$, the mass term in \eqref{eq:r38-macro} is kept on
	the right rather than bounded by \eqref{eq:r38-mass}. The preceding
	combination of \eqref{eq:r38-energy}, \eqref{eq:r38-macro}, and
	\eqref{eq:r38-trace-recovery} may then be made directly on $(0,t)$,
	for any $t\le T$; no restriction on its length is needed. The
	parameters $\kappa$, $\lambda$, $\varepsilon_{\rm dyn}$, $\theta$,
	and $\delta_{\rm dyn}$ were fixed without reference to $T$. Hence
	\begin{align*}
		&\mathcal E_f(t)
		+c\Bigl(\|f\|_{L^2(\cD_{(0,t)};H_\sigma^1)}^2
		+\|f_{\gamma^+}\|_{L^2(\gamma_{(0,t)}^+)}^2\Bigr)
		\\
		&\qquad\le C\Bigl(\mathcal E_f(0)
		+\|F\|_{L^2(\cD_{(0,t)};H_\sigma^{-1})}^2
		+\|b\|_{L^2(\gamma_{(0,t)}^-)}^2
		+\|\mathfrak m_f\|_{L^2(0,t)}^2\Bigr),
	\end{align*}
	where $c,C$ are independent of $t$, $T$, and $N$.
	The transport estimate and the boundary law give
	\eqref{eq:dyn-hier-mass-graph}.
	
	Finally, the collision term has zero mass moment, and
	\eqref{eq:dynamic-physical-mass-balance} gives, for $0\le t\le T$,
	\[
	\mathfrak m_f(t)=\iint_{\cD\times\bR_v^3}f_0\sqrt\mu\,\dd v\dd x
	+\int_0^t\!\int_\cD\langle F(r),\sqrt\mu\rangle_v\dd x\dd r
	+\int_0^t\!\int_{\partial\cD}\beta_b(r,x)\dd S_x\dd r.
	\]
	The zero-mass assumptions of the proposition make all three terms
	vanish.
\end{proof}

\begin{prop}[Evolutionary linear diffuse solvability]
	\label{prop:dynamic-inhomogeneous-physical-solvability}
	Assume \eqref{eq:dynamic-alpha-graph-smallness} and let
	$F\in L^2(\cD_T;H_\sigma^{-1})$, $f_0\in L^2(\cD\times\bR_v^3)$, and
	$b\in L^2(\gamma_T^-)$.
	\begin{enumerate}
		\item\label{it:solv-physical}
		There is a unique $f\in\mathbb X_T$ satisfying
		\begin{equation}
			\label{eq:dynamic-inhomogeneous-physical-problem}
			(\partial_t+v\cdot\nabla_x)f+L_gf=F,
			\qquad f(0)=f_0,
			\qquad f_{\gamma^-}=P_wf_{\gamma^+}+b.
		\end{equation}
		Moreover, we have
		\begin{equation}
			\label{eq:dynamic-inhomogeneous-physical-solvability-estimate}
			\|f\|_{\cH_{\kin,\sigma}(\cD_T\times\bR_v^3)}
			\le C_T\left(
			\|F\|_{L^2(\cD_T;H_\sigma^{-1})}
			+\|f_0\|_{L^2(\cD\times\bR_v^3)}+\|b\|_{L^2(\gamma_T^-)}\right).
		\end{equation}
		
		\item\label{it:solv-spectral}
		For every $N\in\bN$ there is a unique
		$f_N\in\mathbb X_T$ satisfying
		\begin{equation}
			\label{eq:dynamic-spectral-problem}
			(\partial_t+v\cdot\nabla_x)f_N+L_gf_N=F,
			\qquad f_N(0)=f_0,
			\qquad (f_N)_{\gamma^-}=P_{w,N}(f_N)_{\gamma^+} +b.
		\end{equation}
		Moreover, we have
		\begin{equation}
			\label{eq:dynamic-spectral-solvability-estimate}
			\|f_N\|_{\cH_{\kin,\sigma}(\cD_T\times\bR_v^3)}
			\le C_T\left(
			\|F\|_{L^2(\cD_T;H_\sigma^{-1})}
			+\|f_0\|_{L^2(\cD\times\bR_v^3)}+\|b\|_{L^2(\gamma_T^-)}\right)
		\end{equation}
		with $C_T$ independent of $N$.
	\end{enumerate}
	The constants $C_T$ in
	\eqref{eq:dynamic-inhomogeneous-physical-solvability-estimate} and
	\eqref{eq:dynamic-spectral-solvability-estimate} depend only on $\cD$ and on $T$, are
	nondecreasing in $T$, and depend on $g$ only through
	\eqref{eq:dynamic-alpha-graph-smallness}.
	In addition, if
	$\iint_{\cD\times\bR^3_v}f_0\sqrt\mu\dd v\dd x=0$ and
	\eqref{eq:dynamic-physical-zero-mass-data} holds for a.e.\ $t$, then the
	constants $C_T$ in \eqref{eq:dynamic-inhomogeneous-physical-solvability-estimate}
	and \eqref{eq:dynamic-spectral-solvability-estimate} depend only on $\cD$.
\end{prop}

\begin{proof}
	For the physical problem put $P_\bullet=P_w$; for the spectral
	problem put $P_\bullet=P_{w,N}$. Define the continuation for
	$\alpha\in[0,1]$ by
	\[
	\widehat\Lambda_\alpha f
	:=\bigl(
	(\partial_t+v\cdot\nabla_x)f+L_gf,
	\ f(0),
	\ f_{\gamma^-} -\alpha P_\bullet f_{\gamma^+}
	\bigr).
	\]

	The trace estimate \eqref{eq:wall-J-operator-bound} bounds the wall
	operator uniformly in $\cJ$.  With \eqref{260902137},
	\eqref{boundedness of L} and Lemma~\ref{260715201},
	\[
	\widehat\Lambda_\alpha\in\cL(\mathbb X_T,\mathbb Y_T),
	\qquad
	\|\widehat\Lambda_\beta-\widehat\Lambda_\alpha\|
	\le C|\beta-\alpha|,
	\]
	with $C$ independent of $\alpha$ and of $N$.

	Proposition~\ref{prop:dynamic-alpha-uniform-graph}, in case
	\eqref{it:graph-physical} respectively \eqref{it:graph-spectral}, gives
	\begin{equation}
		\label{eq:dynamic-continuation-apriori}
		\|f\|_{\mathbb X_T}
		\le C_T\|\widehat\Lambda_\alpha f\|_{\mathbb Y_T},
		\qquad 0\le\alpha\le1,
	\end{equation}
	with $C_T$ independent of $\alpha$, and of $N$ in case
	(\ref{it:solv-spectral}).
	
	At $\alpha=0$, $\widehat\Lambda_0$ is the prescribed-inflow operator
	\begin{equation*}
		\cA_{\rm in}f:=
		\bigl((\partial_t+v\cdot\nabla_x)f+L_gf,\ f(0),\ f_{\gamma^-}\bigr)
	\end{equation*}
	of Theorem~\ref{260322442}\,(1).
	The choice \eqref{eq:dynamic-KS-perturbed-smallness-choice} ensures that
	\eqref{eq:dynamic-alpha-graph-smallness} meets the coefficient hypothesis of
	Theorem~\ref{260322442}\,(1).
	Then
	\begin{equation*}
		\cA_{\rm in}:\mathbb X_T\longrightarrow\mathbb Y_T
		\quad\text{is bijective},
		\qquad
		\|f\|_{\mathbb X_T}\le C_L\|\cA_{\rm in}f\|_{\mathbb Y_T},
	\end{equation*}
	with $C_L$ depending only on $\cD$.
	
	The method of continuity therefore gives bijectivity of
	$\widehat\Lambda_\alpha$ throughout $[0,1]$. At $\alpha=1$, with target
	$(F,f_0,b)$, this yields
	\eqref{eq:dynamic-inhomogeneous-physical-problem} and
	\eqref{eq:dynamic-spectral-problem}.
	The estimates \eqref{eq:dynamic-inhomogeneous-physical-solvability-estimate}
	and \eqref{eq:dynamic-spectral-solvability-estimate} follow from
	\eqref{eq:dynamic-continuation-apriori} at $\alpha=1$. 
	
	For $\alpha=1$ and either boundary operator, the zero-mass conditions
	give $\mathfrak m_f\equiv0$ by
	Proposition~\ref{prop:dynamic-alpha-uniform-graph}. The bound
	\eqref{eq:dyn-hier-mass-graph} is then independent of $T$.
\end{proof}

\begin{corollary}[Weighted evolutionary graph estimate]
	\label{cor:dynamic-weighted-estimate}
	Assume \eqref{eq:dynamic-alpha-graph-smallness}, let $L\ge0$ be finite, let
	$w^LF\in L^2(\cD_T;H_\sigma^{-1})$, $w^Lf_0\in L^2(\cD\times\bR_v^3)$ and
	$w^Lb\in L^2(\gamma_T^-)$, and let $f_N\in\mathbb X_T$ satisfy
	\eqref{eq:dynamic-spectral-problem}.  Then
	\begin{align}
		\|w^Lf_N\|_{\cH_{\kin,\sigma}(\cD_T\times\bR_v^3)}
		+\sup_{0\le t\le T}\|w^Lf_N(t)\|_{L^2(\cD\times\bR_v^3)}
		&\notag\\
		+\|w^L(f_N)_{\gamma^-}\|_{L^2(\gamma_T^-)}
		+\|w^L(f_N)_{\gamma^+}\|_{L^2(\gamma_T^+)}
		&\notag\\
		\le C_{T,L}\Bigl(
		\|w^LF\|_{L^2(\cD_T;H_\sigma^{-1})}
		+\|w^Lf_0\|_{L^2(\cD\times\bR_v^3)}
		&+\|w^Lb\|_{L^2(\gamma_T^-)}\Bigr),
		\label{eq:dynamic-weighted-diffuse-estimate}
	\end{align}
	with $C_{T,L}$ independent of $N$.  At $L=0$ this is
	\eqref{eq:dynamic-spectral-solvability-estimate}.  The same holds in case
	(\ref{it:solv-physical}), with $P_w$ in place of $P_{w,N}$ and
	\eqref{eq:dynamic-inhomogeneous-physical-solvability-estimate} in place of
	\eqref{eq:dynamic-spectral-solvability-estimate}.
	
	Moreover, there is $C'_L$, depending only on $\cD$ and $L$ (not on $T$ or $N$), such that
	the left side of \eqref{eq:dynamic-weighted-diffuse-estimate} is bounded by
	\begin{equation}
		\label{eq:dynamic-weighted-diffuse-estimate-T-free}
		C'_L\Bigl(
		\|w^LF\|_{L^2(\cD_T;H_\sigma^{-1})}
		+\|w^Lf_0\|_{L^2(\cD\times\bR_v^3)}
		+\|w^Lb\|_{L^2(\gamma_T^-)}
		+\|\mathfrak m_{f_N}\|_{L^2(0,T)}\Bigr).
	\end{equation}
	If in addition $\iint_{\cD\times\bR^3_v}f_0\sqrt\mu\dd v\dd x=0$ and
	\eqref{eq:dynamic-physical-zero-mass-data} holds for a.e.\ $t$, then
	$\mathfrak m_{f_N}\equiv0$, so \eqref{eq:dynamic-weighted-diffuse-estimate} holds
	with $C_{T,L}=C'_L$, independent of $T$.  The same holds in case
	(\ref{it:solv-physical}).
\end{corollary}

\begin{proof}
	By \eqref{eq:dynamic-KS-perturbed-smallness-choice},
	\eqref{eq:dynamic-alpha-graph-smallness} implies \eqref{eq:weighted-dynamic-coefficient}.
	Apply the prescribed-inflow estimate
	\eqref{eq:weighted-graph-convention} to
	\eqref{eq:dynamic-spectral-problem}, with source $F$, initial datum $f_0$,
	and incoming datum $H:=P_{w,N}(f_N)_{\gamma^+}+b$.  By
	\eqref{eq:wall-J-weighted-operator-bound} at $\cJ=J_N$, used at almost every
	$t$ and integrated, and by the trace bound of
	Proposition~\ref{2607151218},
	\[
	\|w^LH\|_{L^2(\gamma_T^-)}
	\le C_L\|(f_N)_{\gamma^+}\|_{L^2(\gamma_T^+)}+\|w^Lb\|_{L^2(\gamma_T^-)}
	\le C_L\|f_N\|_{\cH_{\kin,\sigma}(\cD_T\times\bR_v^3)}
	+\|w^Lb\|_{L^2(\gamma_T^-)} ,
	\]
	and \eqref{eq:dynamic-spectral-solvability-estimate} bounds the graph norm on
	the right.  Since $w\ge1$, the weighted data norms control these unweighted terms.  Case
	(\ref{it:solv-physical}) is the same with $\cJ=\mathbf I$.
	
	For the $T$-independent estimate, control the outgoing trace instead
	by \eqref{eq:dyn-hier-mass-graph}, applied with $\alpha=1$ and $\cJ=J_N$.
	The resulting constant depends only on $\cD$ and $L$:
	\[
	\|w^LH\|_{L^2(\gamma_T^-)}
	\le C\Bigl(
	\|f_0\|_{L^2(\cD\times\bR_v^3)}
	+\|F\|_{L^2(\cD_T;H_\sigma^{-1})}
	+\|b\|_{L^2(\gamma_T^-)}
	+\|\mathfrak m_{f_N}\|_{L^2(0,T)}\Bigr)
	+\|w^Lb\|_{L^2(\gamma_T^-)}.
	\]
	The constants in \eqref{eq:weighted-graph-convention} and
	\eqref{eq:dyn-hier-mass-graph} are independent of $T$ and $N$,
	which proves \eqref{eq:dynamic-weighted-diffuse-estimate-T-free}.
	Under the zero-mass conditions, Proposition~\ref{prop:dynamic-alpha-uniform-graph}
	gives $\mathfrak m_{f_N}\equiv0$.
\end{proof}

\section{Complete Tangential Derivatives and Diffuse Boundary Calculus}\label{sec:diffuse-complete-equations}

We establish the complete-derivative commutators, differentiated diffuse boundary identities, complete-flow estimates, weak closure, and interpolation used in \S\S\ref{sec:diffuse-common}--\ref{sec:dyn-realization-section}.

Unless stated otherwise, $C$ depends on $\cD$, the fixed fields $\cX_{\rm c}$ of \eqref{eq:tangential-field-family}, $\mathrm N$, and the stated parameters, and may change from line to line.  The implicit constants in $\lesssim$ and $\simeq$ are subject to the same rule.

We call $C_c^\infty(\overline\cD\times\bR_v^3)$, and on $(0,T)$ the space
$C_c^\infty([0,T]\times\overline\cD\times\bR_v^3)$, the \emph{smooth core}; the traces of its
elements on $\gamma^\pm$ form the \emph{smooth trace core}.
Since $\tau_\nu\in C^{K_\cD+1}$ by \eqref{eq:tangential-field-family}, the complete
pullbacks of \S\ref{app:diffuse-lifts} and their adjoints map $C_c^k(\cD\times\bR^3)$ and
$C_c^k(\overline\cD\times\bR^3)$ into themselves for $0\le k\le K_\cD$.  We first establish the algebraic identities on these cores.  Extensions to weak derivatives and boundary traces are justified below by the corresponding duality and approximation arguments.

\subsection{Interior algebra}\label{ss:complete-interior-calculus}

We establish the ordered Leibniz rule, transport and collision commutators, and weighted $H_\sigma^{-1}$ bounds for their lower-order terms.

\begin{lemma}[Leibniz rule]
	\label{lem:leibniz-rules}
	Let $\mathrm A$ be a linear operator mapping smooth functions to smooth functions, let $\nu\in\mathfrak A_\bullet$, and let $I=(\nu_1,\ldots,\nu_k)$ be an ordered word over $\mathfrak A_\bullet$ with $k\ge1$.  Define $\operatorname{ad}_{\mathfrak D_\nu}(\mathrm A)$ and $\operatorname{ad}_I(\mathrm A)$ by the first line of
	\begin{equation}
		\label{eq:ad-word}
		\begin{aligned}
			&\operatorname{ad}_{\mathfrak D_\nu}(\mathrm A):=[\mathfrak D_\nu,\mathrm A],
			\qquad
			\operatorname{ad}_I(\mathrm A):=\operatorname{ad}_{\mathfrak D_{\nu_1}}\cdots\operatorname{ad}_{\mathfrak D_{\nu_k}}(\mathrm A),\\
			&[\mathfrak D_I,\mathrm A]
			=\sum_{\varnothing\ne S\subset[k]}\operatorname{ad}_{I_S}(\mathrm A)\,\mathfrak D_{I_{S^c}},
		\end{aligned}
	\end{equation}
	with $[k]$, $I_S$ and $S^c$ as in \eqref{eq:ordered-position-subwords}.  Then the second line of \eqref{eq:ad-word} holds on smooth functions.
\end{lemma}

\begin{proof}
	For $k=1$ the assertion is the definition of the commutator. The general case follows by induction.  Let $k\ge2$, assume the identity for words of length $k-1$, and put $I':=(\nu_2,\ldots,\nu_k)$, so that $\mathfrak D_I=\mathfrak D_{\nu_1}\mathfrak D_{I'}$ and
	\[
	[\mathfrak D_I,\mathrm A]
	=\mathfrak D_{\nu_1}[\mathfrak D_{I'},\mathrm A]+[\mathfrak D_{\nu_1},\mathrm A]\,\mathfrak D_{I'} .
	\]
	The induction hypothesis, with positions indexed by $\{2,\ldots,k\}$, gives
	\[
	[\mathfrak D_{I'},\mathrm A]
	=\sum_{\varnothing\ne S'\subset\{2,\ldots,k\}}\operatorname{ad}_{I_{S'}}(\mathrm A)\,\mathfrak D_{I_{\{2,\ldots,k\}\setminus S'}} .
	\]
	For each such $S'$,
	\[
	\mathfrak D_{\nu_1}\operatorname{ad}_{I_{S'}}(\mathrm A)
	=\operatorname{ad}_{I_{\{1\}\cup S'}}(\mathrm A)+\operatorname{ad}_{I_{S'}}(\mathrm A)\,\mathfrak D_{\nu_1}
	\]
	by the definitions of $\operatorname{ad}_{\mathfrak D_{\nu_1}}$ and $\operatorname{ad}_I$.  The first term contributes the summand with $S=\{1\}\cup S'$, the second the summand with $S=S'$, and $[\mathfrak D_{\nu_1},\mathrm A]\,\mathfrak D_{I'}$ is the summand with $S=\{1\}$.  These are precisely the nonempty subsets of $[k]$, each occurring once.
\end{proof}

\begin{lemma}[Weighted Leibniz rule]
	\label{lem:weight-leibniz}
	Let $\ell\in\bR$.  For every ordered word $I$ over $\mathfrak A_\bullet$ with $|I|\le\mathrm N$,
	\begin{equation}
		\label{eq:weight-word-Leibniz-identity}
		w^{-\ell}\mathfrak D_I(w^{\ell}f)
		=\sum_{I_1\sqcup I_2=I}\bigl(w^{-\ell}\mathfrak D_{I_1}w^{\ell}\bigr)\mathfrak D_{I_2}f,
		\qquad
		\sup_x\bigl\|w^{-\ell}\mathfrak D_{I_1}w^{\ell}\bigr\|_{C^1(\bR_v^3)}\le C_\ell,
	\end{equation}
	with the two-coloring convention \eqref{eq:ordered-two-coloring-convention}; the bound holds for every $I_1$ in the sum, and the term with $I_1=\varnothing$ is $\mathfrak D_If$.
\end{lemma}

\begin{proof}
	For $\nu\in\mathfrak A$, \eqref{eq:complete-phase-space-derivative} gives
	$\mathfrak D_\nu w^{\ell}=\ell\,w^{\ell-2}\,v\cdot((\nabla_x\tau_\nu)v)$, and
	$\mathfrak D_{\mathrm t}w^\ell=0$.  The ordered Leibniz rule gives \eqref{eq:weight-word-Leibniz-identity}. Induction on $|I_1|$ shows that
	$w^{2|I_1|-\ell}\mathfrak D_{I_1}w^\ell$ is a polynomial in $v$ of degree at most $2|I_1|$
	with coefficients bounded on $\overline\cD$.
	The claimed bound follows.
\end{proof}

\subsubsection{The transport term}\label{ss:transport-commutator}

\begin{prop}[Ordered transport commutators]
	\label{prop:ordered-transport-weight}
	Let $\nu\in\mathfrak A$.  Then
	\begin{equation}
		\label{eq:basic-transport-commutator}
		[\mathfrak D_\nu,v\cdot\nabla_x]
		=-v_jv_k\,\p_{jk}\tau_{\nu,i}\,\p_{v_i}.
	\end{equation}
	For every ordered word $I$ over $\mathfrak A$ with
	$1\le |I|\le\mathrm N$,
	\begin{equation}
		\label{eq:ordered-transport-commutator}
		[\mathfrak D_I,v\cdot\nabla_x]
		=\sum_{|J|\le |I|-1}\cV_{I,J}\mathfrak D_J,
	\end{equation}
	where $\cV_{I,J}=b_{I,J}(x,v)\cdot\nabla_v$, and $b_{I,J}$ is a
	polynomial in $v$ of degree at most two whose $x$-coefficients involve only
	derivatives of the fields $\tau_\mu$, $\mu\in\mathfrak A$, of order at most
	$|I|+1$.  In particular, for $(x,v)\in\overline\cD\times\bR^3$,
	\begin{equation}
		\label{eq:ordered-transport-coefficient-bound}
		|b_{I,J}(x,v)|\le C\la v\ra^2 .
	\end{equation}
	For every ordered word $I$ over the mixed alphabet $\mathfrak A_\bullet$ of
	\eqref{eq:dyn-mixed-alphabet} with $1\le|I|\le\mathrm N$,
	\begin{equation}
		\label{eq:mixed-transport-commutator}
		[\mathfrak D_I,\partial_t+v\cdot\nabla_x]
		=[\mathfrak D_I,v\cdot\nabla_x]
		=\sum_{|J|\le|I|-1}\cV_{I,J}\mathfrak D_J
		=\sum_{|J'|\le|I'|-1}\cV_{I',J'}\,\partial_t^{j(I)}\mathfrak D_{J'},
	\end{equation}
	where the last sum is over spatial words $J'$, the word $J'\mathrm t^{j(I)}$ is $J'$ followed by $j(I)$ temporal letters, and
	\[
	\cV_{I,J}=\cV_{I',J'}\quad\text{if }J=J'\mathrm t^{j(I)},
	\qquad
	\cV_{I,J}=0\quad\text{otherwise};
	\]
	in particular $\cV_{\mathrm t,\varnothing}=0$.  Moreover, for every such $I$ and every $\ell\ge0$,
	\begin{equation}
		\label{eq:ordered-transport-weight}
		\|w^\ell\cV_{I,J}f\|_{H_\sigma^{-1}}
		\le C\|w^{\ell+4}f\|_{H_\sigma^1},
		\qquad |J|<|I|,
	\end{equation}
	pointwise in $x$.  Consequently,
	\begin{equation}
		\label{eq:ordered-transport-weight-integrated}
		\|w^\ell\cV_{I,J}f\|_{L_x^2H_\sigma^{-1}}
		\le C\|w^{\ell+4}f\|_{L_x^2H_\sigma^1}.
	\end{equation}
\end{prop}

\begin{proof}
	The phase vector fields of $\mathfrak D_\nu$ and $v\cdot\nabla_x$ are
	$(\tau_\nu,(\nabla\tau_\nu)v)$ and $(v,0)$.  Their Lie bracket has spatial component $0$
	and velocity component
	$-v_jv_k\p_{jk}\tau_{\nu,i}$, which proves
	\eqref{eq:basic-transport-commutator}.  If
	$\cV=b(x,v)\cdot\nabla_v$ and $b$ is polynomial in $v$ of degree at
	most two, then
	\begin{equation}
		\label{eq:complete-velocity-commutator}
		[\mathfrak D_\nu,\cV]
		=\left(
		\tau_\nu\cdot\nabla_x b
		+((\nabla\tau_\nu)v)\cdot\nabla_v b
		-(\nabla\tau_\nu)b
		\right)\cdot\nabla_v.
	\end{equation}
	The coefficient on the right again has degree at most two, since $\nabla_vb$ has degree at most one in $v$.  Write
	$I=(\nu,K)$.  If
	$[\mathfrak D_{K},v\cdot\nabla_x]=\sum_J\cV_{K,J}\mathfrak D_J$, then
	\[
	\begin{aligned}
		[\mathfrak D_\nu\mathfrak D_{K},v\cdot\nabla_x]
		&=\mathfrak D_\nu[\mathfrak D_{K},v\cdot\nabla_x]
		+[\mathfrak D_\nu,v\cdot\nabla_x]\mathfrak D_{K}\\
		&=\sum_J[\mathfrak D_\nu,\cV_{K,J}]\mathfrak D_J
		+\sum_J\cV_{K,J}\mathfrak D_\nu\mathfrak D_J
		+[\mathfrak D_\nu,v\cdot\nabla_x]\mathfrak D_{K}.
	\end{aligned}
	\]
	By \eqref{eq:complete-velocity-commutator}, every new coefficient is a
	polynomial in $v$ of degree at most two; the complete words on the right
	have length at most $|I|-1$.
	The base case \eqref{eq:basic-transport-commutator} involves second derivatives of $\tau_\nu$. Each application of \eqref{eq:complete-velocity-commutator} adds at most one spatial derivative, so $b_{I,J}$ involves derivatives of the fields of order at most $|I|+1$.
	Since $|I|+1\le\mathrm N+1\le K_\cD+1$ by \eqref{eq:order-and-domain-regularity} and $\tau_\nu\in C^{K_\cD+1}$ by \eqref{eq:tangential-field-family}, these
	$x$-coefficients are uniformly bounded, which proves
	\eqref{eq:ordered-transport-coefficient-bound}.
	
	Since $[\mathfrak D_\nu,\partial_t]=0$ for $\nu\in\mathfrak A$ and
	$[\partial_t,v\cdot\nabla_x]=0$, we have
	$$
	[\mathfrak D_I,\partial_t+v\cdot\nabla_x]=[\mathfrak D_I,v\cdot\nabla_x]\quad\text{and}\quad \mathfrak D_I=\partial_t^{j(I)}\mathfrak D_{I'},
	$$
	where $I'$ is the spatial subword of $I$ (see \eqref{eq:dyn-mixed-alphabet}).
	By \eqref{eq:ordered-transport-commutator} for $I'$,
	\[
	[\mathfrak D_I,v\cdot\nabla_x]
	=\partial_t^{j(I)}[\mathfrak D_{I'},v\cdot\nabla_x]
	=\sum_{|J'|\le|I'|-1}\cV_{I',J'}\,\partial_t^{j(I)}\mathfrak D_{J'}.
	\]  Since
	$\mathfrak D_{J'\mathrm t^{j(I)}}=\partial_t^{j(I)}\mathfrak D_{J'}$,
	the right side is $\sum_{|J'|\le|I'|-1}\cV_{I',J'}\mathfrak D_{J'\mathrm t^{j(I)}}$,
	which is \eqref{eq:mixed-transport-commutator} with the stated form of
	$\cV_{I,J}$.
	
	By \eqref{eq:mixed-transport-commutator}, it suffices to prove \eqref{eq:ordered-transport-weight} for spatial words. Then \eqref{eq:ordered-transport-coefficient-bound} gives
	\(|b_{I,J}(x,v)|\le Cw^2\).  Fix \(x\) and write \(F:=w^{\ell+4}f\).
	By \eqref{eq:Hsigma-dual-weighted},
	\[
	\bigl\|w^\ell\cV_{I,J}f\bigr\|_{H_\sigma^{-1}}
	\le C\bigl\|w^{\ell+1/2}b_{I,J}\cdot\nabla_vf\bigr\|_{L_v^2}
	\le C\bigl\|w^{\ell+5/2}\nabla_vf\bigr\|_{L_v^2} .
	\]
	The product rule
	\[
	\nabla_vF=w^{\ell+4}\nabla_vf+(\ell+4)w^{\ell+2}vf ,
	\] together with
	\(|v|\le w\) and \(w\ge1\), gives
	\[
	\bigl\|w^{\ell+5/2}\nabla_vf\bigr\|_{L_v^2}
	\le\bigl\|w^{-3/2}\nabla_vF\bigr\|_{L_v^2}
	+(\ell+4)\bigl\|w^{-1/2}F\bigr\|_{L_v^2},
	\]
	and \eqref{eq:local-H1-H_sigma} gives
	\[
	\|w^{-3/2}\nabla_vF\|_{L_v^2}+\|w^{-1/2}F\|_{L_v^2}\le C\|F\|_{H_\sigma^1}.
	\]
	These bounds prove \eqref{eq:ordered-transport-weight}; integration in $x$ gives \eqref{eq:ordered-transport-weight-integrated}.
\end{proof}

\subsubsection{The collision operators}\label{ss:deformation}

For a matrix kernel $\mathsf K=\mathsf K(x,z)$ we define, fiberwise in $x$,
\[
\mathscr Q_{\mathsf K}(G,F)
:=\partial_{v_i}\!\int_{\bR^3}\mathsf K_{ij}(x,v-v')
\{G(v')\partial_{v_j}F(v)-F(v)\partial_{v'_j}G(v')\}\,\dd v',
\]
and
\begin{equation}
	\label{eq:generalized-conjugated-landau-form}
	N_{\mathsf K}[g,f]
	:=\mu^{-1/2}\mathscr Q_{\mathsf K}(\sqrt\mu\,g,\sqrt\mu\,f).
\end{equation}
For the Coulomb kernel $\phi$, $\mathscr Q_\phi=\cN$, and hence $N_\phi=N$ by \eqref{def:N}.  For a matrix $M$, define
\[
\mathfrak D_M^{\rm ker}\mathsf K
:=(\operatorname{tr}M)\mathsf K-M\mathsf K-\mathsf KM^T
+(Mz)\cdot\nabla_z\mathsf K.
\]
For $\nu\in\mathfrak A$, with $M_\nu:=\nabla\tau_\nu$, set
\begin{equation}
	\label{eq:complete-kernel-deformation-operator}
	\mathscr D_\nu^{\rm ker}\mathsf K
	:=\tau_\nu\cdot\nabla_x\mathsf K+\mathfrak D_{M_\nu}^{\rm ker}\mathsf K,
	\qquad
	\mathfrak p_\nu(x,v):=\mathfrak D_\nu\log\sqrt\mu
	=(M_\nu v)\cdot\nabla_v\log\sqrt\mu.
\end{equation}
For $X=\tau_\nu$ we also write $M_X$, $\mathscr D_X^{\rm ker}$ and $\mathfrak p_X$ for
$M_\nu$, $\mathscr D_\nu^{\rm ker}$ and $\mathfrak p_\nu$.

\begin{prop}[Coulomb tensor deformation]
	\label{prop:complete-lift-coulomb-deformation}\,
	
	\begin{enumerate}
		\item\label{it:kernel-deformation-identity}
		Assume that $\mathsf K(x,\cdot)$ is smooth on
		$\bR^3\setminus\{0\}$ for every $x\in\overline\cD$, and that
		\begin{equation}
			\label{eq:admissible-deformation-kernel}
			|\partial_x^\alpha\partial_z^\beta \mathsf K(x,z)|
			\le C_0|z|^{-1-|\beta|},
			\qquad |\alpha|+|\beta|\le2,
		\end{equation}
		for some constant $C_0$.  Then
		\begin{equation}
			\label{eq:complete-kernel-deformation-identity}
			\mathfrak D_\nu\mathscr Q_{\mathsf K}(G,F)
			-\mathscr Q_{\mathsf K}(\mathfrak D_\nu G,F)
			-\mathscr Q_{\mathsf K}(G,\mathfrak D_\nu F)
			=\mathscr Q_{\mathscr D_\nu^{\rm ker}\mathsf K}(G,F)
		\end{equation}
		on the smooth core.  Equivalently,
		\begin{align}
			\mathfrak D_\nu N_{\mathsf K}[g,f]
			={}&N_{\mathsf K}[\mathfrak D_\nu g,f]
			+N_{\mathsf K}[g,\mathfrak D_\nu f]
			+N_{\mathscr D_\nu^{\rm ker}\mathsf K}[g,f]\notag\\
			&+N_{\mathsf K}[\mathfrak p_\nu g,f]
			+N_{\mathsf K}[g,\mathfrak p_\nu f]
			-\mathfrak p_\nu N_{\mathsf K}[g,f].
			\label{eq:complete-conjugated-kernel-deformation}
		\end{align}
		\item\label{it:iterated-kernels}
		Set $\mathsf K_\varnothing:=\phi$ and define recursively
		\begin{equation}
			\label{eq:iterated-coulomb-deformation-kernels}
			\mathsf K_{\nu I}:=\mathscr D_\nu^{\rm ker}\mathsf K_I.
		\end{equation}
		Then every $\mathsf K_I$, $|I|\le\mathrm N$, is homogeneous of degree $-1$ in $z$, smooth on the
		unit sphere, and satisfies
		\begin{equation}
			\label{eq:iterated-coulomb-kernel-bound}
			|\partial_x^\alpha\partial_z^\beta \mathsf K_I(x,z)|
			\le C|z|^{-1-|\beta|}
		\end{equation}
		whenever $|\alpha|+|\beta|+|I|\le\mathrm N+1$.  In particular every $\mathsf K_I$ with
		$|I|\le\mathrm N-1$ satisfies \eqref{eq:admissible-deformation-kernel}.
	\end{enumerate}
\end{prop}

\begin{proof}
	\emph{(\ref{it:kernel-deformation-identity}).}
	Fix \(x\) and a matrix \(M\).  Set
	\[
	V_M:=(Mv)\cdot\nabla_v,\quad A_s:=\mathrm{Id}+sM,\quad G_s(v):=G(A_sv),\quad F_s(v):=F(A_sv).
	\]
	The change of variables $w=A_sv$, $w'=A_sv'$ and the chain rule give
	\[
	\mathscr Q_{\mathsf K}(G_s,F_s)(v)
	=\mathscr Q_{\mathsf K^{(s)}}(G,F)(A_sv),
	\qquad
	\mathsf K^{(s)}(z)=(\det A_s)^{-1}A_s\mathsf K(A_s^{-1}z)A_s^T.
	\]
	The bounds \eqref{eq:admissible-deformation-kernel} justify differentiation at $s=0$ on the smooth core: $\partial_s\mathsf K^{(s)}$ and its first $z$-derivatives are bounded locally by $C|z|^{-1}$ and $C|z|^{-2}$, respectively, both locally integrable in three dimensions. Since
	\[
	\frac{\dd}{\dd s}\Big|_{0}(\det A_s)^{-1}=-\operatorname{tr}M,
	\qquad
	\frac{\dd}{\dd s}\Big|_{0}A_s^{-1}=-M,
	\]
	we obtain
	\[
	V_M\mathscr Q_{\mathsf K}(G,F)-\mathscr Q_{\mathsf K}(V_MG,F)-\mathscr Q_{\mathsf K}(G,V_MF)
	=\mathscr Q_{\mathfrak D_M^{\rm ker}\mathsf K}(G,F).
	\]
	At $M=M_\nu$, the product rule for $\tau_\nu\cdot\nabla_x$ adds the $x$-derivative of the kernel, proving \eqref{eq:complete-kernel-deformation-identity}. Moreover
	\[
	\mathfrak D_\nu\mu^{-1/2}=-\mathfrak p_\nu\mu^{-1/2},
	\qquad
	\mathfrak D_\nu(\sqrt\mu\,g)
	=\sqrt\mu\,(\mathfrak D_\nu g+\mathfrak p_\nu g),
	\]
	which gives \eqref{eq:complete-conjugated-kernel-deformation}.
	
	\emph{(\ref{it:iterated-kernels}).}
	Write $z=r\omega$ with $r>0$ and $\omega\in\mathbb S^2$, and set
	$\mathsf K_I(x,r\omega)=r^{-1}\Omega_I(x,\omega)$.
	Let $M=\nabla\tau_\nu$ and $P_\omega=\mathrm{Id}-\omega\otimes\omega$.
	The polar decomposition of $\nabla_z$ gives
	\[
	(Mz)\cdot\nabla_z(r^{-1}\Omega_I)
	=r^{-1}\bigl((P_\omega M\omega)\cdot\nabla_{\mathbb S^2}\Omega_I
	-(\omega\cdot M\omega)\Omega_I\bigr).
	\]
	Thus \eqref{eq:iterated-coulomb-deformation-kernels} is equivalent to
	\[
	\begin{aligned}
		\Omega_{\nu I}={}&\tau_\nu\cdot\nabla_x\Omega_I
		+(\operatorname{tr}M-\omega\cdot M\omega)\Omega_I
		-M\Omega_I-\Omega_I M^T\\
		&+(P_\omega M\omega)\cdot\nabla_{\mathbb S^2}\Omega_I.
	\end{aligned}
	\]
	This operator is first order on $\overline\cD\times\mathbb S^2$.
	A derivative of combined order $|\alpha|+b$ on its right-hand side
	involves derivatives of $\Omega_I$ of combined order at most
	$|\alpha|+b+1$, and coefficients controlled by derivatives of
	$\tau_\nu$ of order at most $|\alpha|+1$.
	Since $\Omega_\varnothing(\omega)=\mathrm{Id}-\omega\otimes\omega$
	is smooth, induction gives
	\[
	\sup_{x\in\overline\cD}
	\|\partial_x^\alpha\Omega_I(x,\cdot)\|_{C^b(\mathbb S^2)}
	\le C_{I,\alpha,b},
	\qquad |\alpha|+b+|I|\le\mathrm N+1.
	\]
	The constants depend only on the fixed differentiation order
	and the geometric fields in $C^{\mathrm N+1}$.
	Cartesian derivatives of $r^{-1}\Omega_I(x,\omega)$ therefore satisfy
	\[
	|\partial_x^\alpha\partial_z^\beta\mathsf K_I(x,z)|
	\le C_{I,\alpha,\beta}|z|^{-1-|\beta|},
	\qquad |\alpha|+|\beta|+|I|\le\mathrm N+1.
	\]
	This proves \eqref{eq:iterated-coulomb-kernel-bound}.  In particular,
	\eqref{eq:admissible-deformation-kernel} holds for $|I|\le\mathrm N-1$.
\end{proof}

\begin{lemma}[Polynomial degree bounds]
	\label{lem:deformation-degrees}
	Let $I$ be an ordered word over $\mathfrak A$, $k\ge0$ and $m:=|I|+k\le\mathrm N$.
	Repeated use of \eqref{eq:complete-conjugated-kernel-deformation} and the Leibniz rule expands $\mathfrak D_{\nu_1}\cdots\mathfrak D_{\nu_k}N_{\mathsf K_I}[g,f]$, for $\nu_1,\ldots,\nu_k\in\mathfrak A$, into terms of the form
	\[
	p\,N_{\mathsf K_J}[p_1a,p_2b].
	\]
	Here $a$ and $b$ are complete derivatives of order at most $k$ of $g$ and $f$,
	$\mathsf K_J$ is a kernel of \eqref{eq:iterated-coulomb-deformation-kernels}, and
	$p,p_1,p_2$ are polynomials in $v$ whose coefficients are polynomials in the
	derivatives of order at most $k$ of the fields $\tau_\nu$.
	Each term satisfies
	\begin{align}
		|J|+\tfrac12\bigl(\deg_vp+\deg_vp_1+\deg_vp_2\bigr)&\le m,
		\label{eq:deformation-combined-grading}\\
		\deg_v p+\deg_v p_1+\deg_v p_2&\le2m.
		\label{eq:deformation-degree-bound}
	\end{align}
	If $g=\sqrt\mu$ and every $\mathfrak D_\nu\sqrt\mu$ is rewritten as $\mathfrak p_\nu\sqrt\mu$, the same holds with $a=\sqrt\mu$.
\end{lemma}

\begin{proof}
	For $\Phi=pN_{\mathsf K_J}[p_1a,p_2b]$ set
	\[
	\operatorname{gr}\Phi:=|J|+\tfrac12(\deg_vp+\deg_vp_1+\deg_vp_2).
	\]
	Initially $\Phi=N_{\mathsf K_{I}}[g,f]$, $p=p_1=p_2=1$ and
	$\operatorname{gr}\Phi=|I|$.
	
	By \eqref{eq:complete-conjugated-kernel-deformation} and the Leibniz rule,
	\begin{align}\label{260916430}
		\begin{aligned}
			\mathfrak D_\nu\Phi
			={}&(\mathfrak D_\nu p)N_{\mathsf K_J}[p_1a,p_2b]
			+pN_{\mathsf K_J}[(\mathfrak D_\nu p_1)a,p_2b]
			+pN_{\mathsf K_J}[p_1a,(\mathfrak D_\nu p_2)b]\\
			&+pN_{\mathsf K_J}[p_1\mathfrak D_\nu a,p_2b]
			+pN_{\mathsf K_J}[p_1a,p_2\mathfrak D_\nu b]
			+pN_{\mathsf K_{\nu J}}[p_1a,p_2b]\\
			&+pN_{\mathsf K_J}[\mathfrak p_\nu p_1a,p_2b]
			+pN_{\mathsf K_J}[p_1a,\mathfrak p_\nu p_2b]
			-\mathfrak p_\nu pN_{\mathsf K_J}[p_1a,p_2b],
		\end{aligned}
	\end{align}
	with $\mathfrak p_\nu=-v^T(\nabla\tau_\nu)v$ by \eqref{eq:complete-kernel-deformation-operator}.
	
	By \eqref{eq:complete-phase-space-derivative},
	$\mathfrak D_\nu=\tau_\nu\cdot\nabla_x+((\nabla\tau_\nu)v)\cdot\nabla_v$; its velocity
	part preserves $\deg_v$.  Both $\mathfrak D_\nu$ and multiplication by $\mathfrak p_\nu$
	raise the order of the field derivatives in the coefficients by at most one.
	
	In \eqref{260916430}, differentiation of a polynomial factor or either argument does not increase the grading. Kernel deformation increases $|J|$ by one, while a Maxwellian-conjugation factor $\mathfrak p_\nu$ increases the total polynomial degree by at most two. Thus each term has grading at most $\operatorname{gr}\Phi+1$.
	
	Each letter also increases the derivative order of the coefficient fields by at most one. Induction proves \eqref{eq:deformation-combined-grading} and the asserted coefficient bounds; \eqref{eq:deformation-degree-bound} follows.

	If $\mathfrak D_\nu$ is applied to $\Phi$, then $\Phi$ arises from $N_{\mathsf K_I}[g,f]$ through at most $k-1$ letters. Since the initial value of $\operatorname{gr}$ is $|I|=m-k$,
	\[
	|J|\le\operatorname{gr}\Phi\le(m-k)+(k-1)=m-1\le\mathrm N-1.
	\]  
	So
	$\mathsf K_J$ satisfies \eqref{eq:admissible-deformation-kernel} by
	Proposition~\ref{prop:complete-lift-coulomb-deformation}\,(\ref{it:iterated-kernels}), and
	\eqref{eq:complete-conjugated-kernel-deformation} applies to $\Phi$.  Induction on the number of letters therefore justifies the expansion of $\mathfrak D_\nu\Phi$ above.
	
	When $a=\sqrt\mu$, the identity $\mathfrak D_\nu a=\mathfrak p_\nu a$ absorbs every derivative of the first argument into its polynomial factor. Thus $a=\sqrt\mu$ throughout the expansion.
\end{proof}

\begin{prop}[Deformed collision bounds]
	\label{prop:deformation-algebra}
	Let $1\le m\le\mathrm N$, let $\mathsf K_I$ be any kernel generated from the Coulomb
	tensor by \eqref{eq:iterated-coulomb-deformation-kernels} with $|I|\le m$, and
	let $p,p_1,p_2$ be polynomials in $v$ with
	\begin{equation}
		\label{eq:deformation-algebra-general-polynomials}
		\deg_vp+\deg_vp_1+\deg_vp_2\le2m,
		\qquad
		\text{coefficients in }L^\infty(\cD).
	\end{equation}
	Then, for
	$\frac14\le\varpi<\frac12$ and every $\ell\ge0$,
	\begin{equation}
		\label{eq:deformation-weighted-output}
		\|w^{\ell}pN_{\mathsf K_I}[p_1g,p_2f]\|_{H_\sigma^{-1}}
		\le C\|\mu^\varpi g\|_2\|w^{\ell+2m+2}f\|_{H_\sigma^1}.
	\end{equation}
	If $2\le q,r\le\infty$ and $q^{-1}+r^{-1}\le1/2$, then
	\begin{equation}
		\label{eq:deformation-weighted-output-spatial}
		\|w^{\ell}pN_{\mathsf K_I}[p_1g,p_2f]\|_{L_x^2H_\sigma^{-1}}
		\le C
		\|\mu^\varpi g\|_{L_x^qL_v^2}
		\|w^{\ell+2m+2}f\|_{L_x^rH_\sigma^1}.
	\end{equation}
	The constants may depend on $m$, $\ell$, $\varpi$, the geometric fields,
	and the $L^\infty(\cD)$ norms of the polynomial coefficients.
\end{prop}

\begin{proof}
	We first work at fixed $x$.  Write 
	\begin{align*}
		\begin{gathered}
			d:=\deg_vp+\deg_vp_1+\deg_vp_2\le2m,\\
			D^\pm:=\nabla_v\pm v,\quad H:=\mu^\varpi g,\quad F_1:=\sqrt\mu\,p_1g,\quad a:=\ell+2m+2,		
		\end{gathered}
	\end{align*}
	and
	\[
	A_1:=\mathsf K_I*F_1,
	\qquad
	B_{1,i}:=(\partial_{z_j}\mathsf K_{I,ij})*F_1 .
	\]
	By \eqref{eq:generalized-conjugated-landau-form}, with the $v'$-derivative moved onto the
	kernel, and one integration by parts in $v$, for $\psi\in C_c^\infty(\bR^3)$,
	\[
	\langle w^\ell pN_{\mathsf K_I}[p_1g,p_2f],\psi\rangle_v
	=-\int A_1D^-(p_2f)\cdot D^+(w^\ell p\psi)\,\dd v
	+\int B_1\,p_2f\cdot D^+(w^\ell p\psi)\,\dd v .
	\]
	
	Fix $\varpi<\frac12$ and let $c=(\frac12-\varpi)/4>0$.
	Since $p_1$ has bounded degree and coefficients bounded uniformly in $x$,
	$F_1=\mu^{1/2-\varpi}p_1H$ satisfies
	$|F_1|\le C_{\varpi,p_1}\ee^{-2c|v|^2}|H|$.  Split the
	convolutions at $|v_*|=|v|/2$, and write $B_1=B_1'+B_1''$, with $B_1''$ the part over
	$|v-v_*|\le1$.  On $|v_*|\le|v|/2$, we have $|v-v_*|\ge|v|/2$. On $|v_*|>|v|/2$, we have
	$\ee^{-2c|v_*|^2}\le\ee^{-c|v|^2/4}\ee^{-c|v_*|^2}$. On $|v-v_*|\le1$, we have
	$\ee^{-2c|v_*|^2}\le C\ee^{-c|v|^2}$. Moreover, $|z|^{-1}\mathbf 1_{|z|\le1}\in L^2_z$, which also gives $|A_1|\le C\|H\|_2$ on $|v|\le2$.  Hence
	\eqref{eq:iterated-coulomb-kernel-bound} at $|\alpha|=|\beta|=0$ gives, for $|v|\ge2$,
	\[
	\begin{aligned}
		|A_1(v)|\le{}&\frac{C}{|v|}\int_{|v_*|\le|v|/2}\ee^{-2c|v_*|^2}|H(v_*)|\,\dd v_*
		+C\ee^{-c|v|^2/4}\int_{|v_*|>|v|/2}\frac{\ee^{-c|v_*|^2}|H(v_*)|}{|v-v_*|}\,\dd v_*.
	\end{aligned}
	\]
	In addition, \eqref{eq:iterated-coulomb-kernel-bound} at $|\alpha|=0$ and $|\beta|=1$ gives, since $|v-v_*|>1$ on the domain of integration of $B_1'$, 
	\begin{align*}
		\begin{alignedat}{2}
			|B_1'(v)|\,&\le C\int_{\bR^3}\ee^{-2c|v_*|^2}|H(v_*)|\,\dd v_*\qquad &&\text{for }|v|\le2,\\
			|B_1'(v)|\,&\le\frac{C}{|v|^2}\int_{|v_*|\le|v|/2}\ee^{-2c|v_*|^2}|H(v_*)|\,\dd v_*&&\\
			&\qquad +C\ee^{-c|v|^2/4}\int_{|v_*|>|v|/2}\ee^{-c|v_*|^2}|H(v_*)|\,\dd v_*\qquad &&\text{for }|v|\ge2,\\
			|B_1''(v)|\,&\le C\int_{|v-v_*|\le1}|v-v_*|^{-2}\ee^{-2c|v_*|^2}|H(v_*)|\,\dd v_*\qquad &&\text{for every }v.		
		\end{alignedat}
	\end{align*}
	With the Cauchy--Schwarz inequality, these give
	\[
	|A_1|\le Cw^{-1}\|H\|_2,
	\qquad
	|B_1'|\le Cw^{-2}\|H\|_2,
	\qquad
	|B_1''|\le C\ee^{-c|v|^2}h,
	\]
	where $h:=(|z|^{-2}\mathbf 1_{|z|\le1})*|H|$ satisfies $\|h\|_{L^6}\le C\|H\|_2$ by the
	Hardy--Littlewood--Sobolev inequality.
	
	By Lemma~\ref{lem:landau-diffusion-eigenvalues} and
	\eqref{eq:coulomb-maxwellian-anisotropy},
	$\sigma\simeq w^{-3}P_v+w^{-1}(\mathrm{Id}-P_v)$ with
	$P_v$ as in \eqref{eq:radial-tangential-projections}, so $|\sigma^{-1/2}\xi|\le Cw^{3/2}|\xi|$ and
	$|\sigma^{1/2}\xi|\le Cw^{-1/2}|\xi|$ for $\xi\in\bR^3$.  Since \[ D^+(w^\ell p\psi)=w^\ell pD^+\psi+\nabla_v(w^\ell p)\psi, \quad D^-(p_2f)=p_2D^-f+(\nabla_vp_2)f, \quad |\nabla_v(w^\ell p)|\le Cw^{\ell+\deg_vp-1}, \] we have
	\[
	|\sigma^{1/2}D^+(w^\ell p\psi)|\le Cw^{\ell+\deg_vp}\Psi,
	\qquad
	|\sigma^{1/2}D^-(p_2f)|\le Cw^{\deg_vp_2}\Phi,
	\]
	where 
	$$
	\Psi:=|\sigma^{1/2}D^+\psi|+w^{-1/2}|\psi|,\qquad \Phi:=|\sigma^{1/2}D^-f|+w^{-1/2}|f|.
	$$
	By the bounds on $A_1$, $B_1'$ and $\sigma^{-1/2}$ above, 
	$$
	|\sigma^{-1/2}A_1\sigma^{-1/2}|\le Cw^2\|H\|_2,\qquad |\sigma^{-1/2}B_1'|\le Cw^{-1/2}\|H\|_2.
	$$
	Writing
	$$
	A_1=\sigma^{1/2}(\sigma^{-1/2}A_1\sigma^{-1/2})\sigma^{1/2},\qquad 
	B_1'=\sigma^{1/2}(\sigma^{-1/2}B_1'),
	$$
	we obtain, pointwise in $v$,
	\[
	\begin{aligned}
		&|A_1D^-(p_2f)\cdot D^+(w^\ell p\psi)|+|B_1'\,p_2f\cdot D^+(w^\ell p\psi)|\\
		&\qquad\le C\|H\|_2\bigl(w^2\,w^{\deg_vp_2}\Phi+w^{-1/2}w^{\deg_vp_2}|f|\bigr)w^{\ell+\deg_vp}\Psi
		\le C\|H\|_2\,w^{\ell+d+2}\Phi\Psi ,
	\end{aligned}
	\]
	because of $\deg_vp+\deg_vp_2\le d$ and $w^{-1/2}|f|\le\Phi$.
	The Cauchy--Schwarz inequality in $v$ gives
	\[
	\Bigl|\int A_1D^-(p_2f)\cdot D^+(w^\ell p\psi)\,\dd v\Bigr|
	+\Bigl|\int B_1'\,p_2f\cdot D^+(w^\ell p\psi)\,\dd v\Bigr|
	\le C\|H\|_2\|w^{\ell+d+2}\Phi\|_2\|\Psi\|_2 .
	\]
	By the embedding $H^1(\bR^3)\subset L^3(\bR^3)$, $|\nabla_v\log w^a|\le a$, and \eqref{eq:local-H1-H_sigma} applied to $w^af$,
	\[
	\|\ee^{-c|v|^2/4}f\|_{L^3}
	\le C\|\ee^{-c|v|^2/4}f\|_{H^1}
	\le C\|w^{-3/2}\nabla_v(w^af)\|_2+C\|w^{-1/2}w^af\|_2
	\le C\|w^af\|_{H_\sigma^1}.
	\]  After the polynomial factors are absorbed
	into the Gaussian, H\"older's inequality with exponents $6$, $3$ and $2$, and \eqref{eq:local-H1-H_sigma} for $\psi$, give
	\begin{align*}
		\Bigl|\int B_1''\,p_2f\cdot D^+(w^\ell p\psi)\,\dd v\Bigr|
		\,&\le C\|h\|_{L^6}\|\ee^{-c|v|^2/4}f\|_{L^3}
		\bigl\|\ee^{-c|v|^2/2}(|\nabla_v\psi|+|\psi|)\bigr\|_{L^2}\\
		&\le C\|H\|_2\|w^af\|_{H_\sigma^1}\|\psi\|_{H_\sigma^1}.
	\end{align*}
	Finally, $w^aD^-f=D^-(w^af)-(\nabla_vw^a)f$ and
	$|\sigma^{1/2}\nabla_vw^a|\le Cw^{a-1/2}$. Since $\ell+d+2\le a$,
	\eqref{eq:velalg-Dpm} and \eqref{eq:local-H1-H_sigma} give
	\[
	\|w^{\ell+d+2}\Phi\|_2
	\le\|w^a\sigma^{1/2}D^-f\|_2+\|w^{a-1/2}f\|_2
	\le C\|w^af\|_{H_\sigma^1},
	\qquad
	\|\Psi\|_2\le C\|\psi\|_{H_\sigma^1}.
	\]
	Hence
	\[
	|\langle w^\ell pN_{\mathsf K_I}[p_1g,p_2f],\psi\rangle_v|
	\le C\|\mu^\varpi g\|_2\|w^af\|_{H_\sigma^1}\|\psi\|_{H_\sigma^1},
	\]
	and duality gives
	\eqref{eq:deformation-weighted-output}; \eqref{eq:deformation-weighted-output-spatial} follows from it by
	H\"older's inequality in $x$ on the bounded set $\cD$.
\end{proof}

\begin{lemma}[Local complete commutators]
	\label{lem:nested-local-complete-commutators}
	For $X=\tau_\nu$, $\nu\in\mathfrak A$, on the smooth core,
	\begin{align}
		[\mathfrak D_X,\cA_\mu]f
		={}&\mathfrak p_XN_\phi[\sqrt{\mu},f]
		-2N_\phi[\mathfrak p_X\sqrt{\mu},f]
		-N_\phi[\sqrt{\mu},\mathfrak p_Xf]
		-N_{\mathfrak D_{M_X}^{\rm ker}\phi}[\sqrt{\mu},f].
		\label{eq:Amu-one-complete-exact}
	\end{align}
	For an ordered word $I=(\nu_1,\ldots,\nu_m)$ over $\mathfrak A$,
	$1\le m\le\mathrm N$, and $X_j:=\tau_{\nu_j}$, define
	\begin{equation}
		\label{eq:nested-local-complete-definition}
		\mathfrak L_{X_m}:=[\mathfrak D_{X_m},\cA_\mu],
		\qquad
		\mathfrak L_{X_1\mid\cdots\mid X_m}
		:=[\mathfrak D_{X_1},\mathfrak L_{X_2\mid\cdots\mid X_m}].
	\end{equation}
	\begin{enumerate}
		\item\label{it:local-nested-commutators}
		On the smooth core, $\mathfrak L_{X_1\mid\cdots\mid X_m}f$ is a finite sum of terms
		\begin{equation}
			\label{eq:nested-local-complete-expansion}
			c\,p\,N_{\mathsf K_A}\bigl[p_1\sqrt\mu,p_2f\bigr],
			\qquad |A|\le m,
			\qquad \deg_vp+\deg_vp_1+\deg_vp_2\le2m,
		\end{equation}
		with $c\in\bR$ and the coefficients of the polynomials $p,p_1,p_2$ in $v$ bounded on
		$\overline\cD$. Consequently, for every $\ell\ge0$,
		\begin{align}
			\|w^\ell\mathfrak L_{X_1\mid\cdots\mid X_m}f\|_{H_\sigma^{-1}}
			&\le C\|w^{\ell+2m+2}f\|_{H_\sigma^1},
			\label{eq:nested-local-complete-bound}\\
			\|w^\ell\mathfrak L_{X_1\mid\cdots\mid X_m}f\|_{L_x^2H_\sigma^{-1}}
			&\le C\|w^{\ell+2m+2}f\|_{L_x^2H_\sigma^1}.
			\label{eq:nested-local-complete-bound-spatial}
		\end{align}
		\item\label{it:local-ordered-commutators}
		Let $I$ be an ordered word over $\mathfrak A$ with $1\le|I|\le\mathrm N$.  On the
		smooth core,
		\begin{equation}\label{eq:local-L0-complete-expansion}
			[\mathfrak D_I,\cA_\mu]f
			=\sum_{|J|<|I|}\sum_{\lambda}
			c_{I,J,\lambda}(x)\, p_{I,J,\lambda}(x,v) N_{\mathsf K_{I,J,\lambda}} [m_{I,J,\lambda}(x,v),p'_{I,J,\lambda}(x,v)\mathfrak D_Jf],
		\end{equation}
		a finite sum.  Here $c_{I,J,\lambda}$ is bounded on $\overline\cD$, $p_{I,J,\lambda}$ and $p'_{I,J,\lambda}$ are
		polynomials in $v$ with coefficients bounded on $\overline\cD$, $m_{I,J,\lambda}$ is a
		Gaussian--polynomial function, and $\mathsf K_{I,J,\lambda}=\mathsf K_A$ is a kernel of
		\eqref{eq:iterated-coulomb-deformation-kernels} with $|A|\le|I|$.
		For every $\ell\ge0$,
		\begin{equation}\label{eq:local-L0-complete-bound}
			\|w^\ell[\mathfrak D_I,\cA_\mu]f\|_{H^{-1}_\sigma}
			\le C\sum_{|J|<|I|}
			\|w^{\ell+2|I|+2}\mathfrak D_Jf\|_{H^1_\sigma}.
		\end{equation}
	\end{enumerate}
\end{lemma}

\begin{proof}
	By the definition of $\cA_\mu$ before \eqref{eq:compact-integral-is-inflow-K} and
	$\mathscr Q_\phi=\cN$,
	\[
	\cA_\mu f=-\mu^{-1/2}\mathscr Q_\phi(\mu,\sqrt\mu f)=-N_\phi[\sqrt\mu,f].
	\]
	By \eqref{eq:complete-kernel-deformation-operator}, $\mathfrak D_X\sqrt{\mu}=\mathfrak p_X\sqrt{\mu}$ and
	$\mathscr D_X^{\rm ker}\phi=\mathfrak D_{M_X}^{\rm ker}\phi$. Hence
	\eqref{eq:complete-conjugated-kernel-deformation} at $\mathsf K=\phi$ gives \eqref{eq:Amu-one-complete-exact}.
	
	\emph{(\ref{it:local-nested-commutators}).}
	The case $m=1$ follows from \eqref{eq:Amu-one-complete-exact}.
	
	Let $m<\mathrm N$ and assume the expansion at order $m$.  By
	\eqref{eq:nested-local-complete-definition}, with $X$ the new outermost letter,
	\[
	\mathfrak L_{X\mid X_1\mid\cdots\mid X_m}
	=[\mathfrak D_X,\mathfrak L_{X_1\mid\cdots\mid X_m}]
	=\sum_B[\mathfrak D_X,B],
	\]
	the sum ranging over the terms of the order-$m$ expansion \eqref{eq:nested-local-complete-expansion}.
	Let
	\[
	Bf=c\,p(x,v)N_{\mathsf K}[p_1(x,v)\sqrt\mu,p_2(x,v)f]
	\]
	be such a term, with $\mathsf K=\mathsf K_A$ and $|A|\le m\le\mathrm N-1$. By Proposition~\ref{prop:complete-lift-coulomb-deformation}\,(\ref{it:iterated-kernels}), this kernel satisfies \eqref{eq:admissible-deformation-kernel}.
	By \eqref{eq:complete-conjugated-kernel-deformation} and the product rules
	\[
	\mathfrak D_X(p_2f)=(\mathfrak D_Xp_2)f+p_2\mathfrak D_Xf, \qquad \mathfrak D_X(p_1\sqrt\mu)=(\mathfrak D_Xp_1+\mathfrak p_Xp_1)\sqrt\mu,
	\] we have
	\begin{align*}
		[\mathfrak D_X,B]f
		={}&c(\mathfrak D_Xp)N_{\mathsf K}[p_1\sqrt\mu,p_2f]
		+cpN_{\mathsf K}[(\mathfrak D_Xp_1)\sqrt\mu,p_2f]\\
		&
		+cpN_{\mathsf K}[p_1\sqrt\mu,(\mathfrak D_Xp_2)f]
		+cpN_{\mathscr D_X^{\rm ker}\mathsf K}[p_1\sqrt\mu,p_2f]\\
		&+2cpN_{\mathsf K}[(\mathfrak p_Xp_1)\sqrt\mu,p_2f]
		+cpN_{\mathsf K}[p_1\sqrt\mu,\mathfrak p_Xp_2f]
		-c\mathfrak p_XpN_{\mathsf K}[p_1\sqrt\mu,p_2f].
	\end{align*}
	The term $cpN_{\mathsf K}[p_1\sqrt\mu,p_2\mathfrak D_Xf]$ cancels with $B(\mathfrak D_Xf)$; thus no complete derivative of $f$ remains.
	Only the kernel deformation changes $\mathsf K_A$: by \eqref{eq:iterated-coulomb-deformation-kernels}, $\mathscr D_X^{\rm ker}\mathsf K_A=\mathsf K_{\nu A}$ for $X=\tau_\nu$. Thus $|A|$ increases by at most one.
	Since 
	$$
	\deg_v(\mathfrak D_Xq)\le\deg_vq\quad(q=p,p_1,p_2)\quad \text{and}\quad \deg_v\mathfrak p_X\le2,
	$$
	$\deg_v p+\deg_v p_1+\deg_v p_2$ increases by at most two.
	Hence \eqref{eq:nested-local-complete-expansion} holds at order $m+1$.
	
	Applying \eqref{eq:deformation-weighted-output} termwise with $g=\sqrt\mu$ gives \eqref{eq:nested-local-complete-bound}; integration over $x$ gives \eqref{eq:nested-local-complete-bound-spatial}.  The bounds are termwise, so they hold
	for every $f$ with $w^{\ell+2m+2}f \in H_\sigma^1$ and $w^{\ell+2m+2}f\in L_x^2H_\sigma^1$, respectively.
	
	\emph{(\ref{it:local-ordered-commutators}).}
	By \eqref{eq:ad-word} at $\mathrm A=\cA_\mu$ and
	\eqref{eq:nested-local-complete-definition},
	\[
	[\mathfrak D_I,\cA_\mu]=\sum_{\varnothing\ne S\subset[|I|]}\mathfrak L_{I_S}\mathfrak D_{I_{S^c}},
	\]
	where $\mathfrak L_{I_S}:=\operatorname{ad}_{I_S}(\cA_\mu)=\mathfrak L_{X_{s_1}\mid\cdots\mid X_{s_k}}$ for $S=\{s_1<\cdots<s_k\}$.  By
	(\ref{it:local-nested-commutators}) of this lemma at $m=|S|$, each summand has the form
	\eqref{eq:local-L0-complete-expansion} with $J=I_{S^c}$ and $|J|<|I|$; in the notation of \eqref{eq:nested-local-complete-expansion},
	\[
	|A|\le|S|\le|I|,
	\qquad
	\deg_vp+\deg_vp_1+\deg_vp_2\le2|S|\le2|I|,
	\]
	and the factors $p_1\sqrt\mu$ and $p_2$ give $m_{I,J,\lambda}$ and $p'_{I,J,\lambda}$.  Then
	\eqref{eq:deformation-weighted-output} at $m=|I|$, with $g=\sqrt\mu$, gives
	\eqref{eq:local-L0-complete-bound}.  The bound is termwise, so it holds for each term of
	\eqref{eq:local-L0-complete-expansion} with $\mathfrak D_Jf$ replaced by any $q_J$ with
	$w^{\ell+2|I|+2}q_J\in H^1_\sigma$.
\end{proof}

\begin{lemma}[Maxwellian commutators]
	\label{lem:compact-maxwellian-complete-commutators}
	Use the decomposition \eqref{eq:compact-integral-is-inflow-K} and write
	\[
	\cK_0f=m_0(v)f+\cK_{k_0}f,
	\qquad
	(\cK_kf)(x,v):=\int_{\bR^3}k(x;v,v')f(x,v')\,\dd v',
	\]
	where $m_0=-8\pi\mu$ and $k_0=k_0^{\rm col}$ is given by
	\eqref{eq:compact-integral-kernel}.  For $X=\tau_\nu$, $\nu\in\mathfrak A$, set
	\[
	m_{\varnothing}:=m_0,
	\qquad
	m_{XI}:=\mathfrak D_Xm_I,
	\qquad
	k_{\varnothing}:=k_0,
	\]
	and define
	\begin{equation}
		\label{eq:iterated-compact-kernel-recursion}
		k_{XI}
		:=Xk_I+(M_Xv)\cdot\nabla_vk_I
		+(M_Xv')\cdot\nabla_{v'}k_I
		+(\operatorname{tr}M_X)k_I.
	\end{equation}
	Then, on the smooth core, with $m_I$ acting by multiplication,
	\[
	[\mathfrak D_X,m_I]=m_{XI},
	\qquad
	[\mathfrak D_X,\cK_{k_I}]=\cK_{k_{XI}}.
	\]
	Thus $\operatorname{ad}_K(\cK_0)=m_K+\cK_{k_K}$ for every ordered word
	$K\ne\varnothing$ over $\mathfrak A$, and, for every ordered word $I$ over $\mathfrak A$
	with $1\le|I|\le\mathrm N$,
	\begin{equation}
		\label{eq:compact-word-commutator-expansion}
		[\mathfrak D_I,\cK_0]f
		=\sum_{\varnothing\ne S\subset[|I|]}
		\operatorname{ad}_{I_S}(\cK_0)\,\mathfrak D_{I_{S^c}}f .
	\end{equation}
	
	For each ordered subword $J$ of $I$ with $|J|<|I|$, set
	\[
	\mathscr K_{I,J}
	:=\sum_{\substack{\varnothing\ne S\subset[|I|]\\ I_{S^c}=J}}\operatorname{ad}_{I_S}(\cK_0).
	\]  
	Each $\mathscr K_{I,J}$, initially defined on the smooth
	core, extends uniquely to a bounded map from $L^2_{x,v}$ to
	$L_x^2H_\sigma^{-1}$.  For every $\ell\ge0$ and every
	$f\in L^2_{x,v}$,
	\begin{equation}
		\label{eq:iterated-compact-full-operator-bound}
		\|w^\ell \mathscr K_{I,J}f\|_{L_x^2H^{-1}_\sigma}
		\le C\|f\|_{L_{x,v}^2},
		\qquad |J|<|I|\le\mathrm N.
	\end{equation}
\end{lemma}

\begin{proof}
	$[\mathfrak D_X,m_I]f=(\mathfrak D_Xm_I)f$.
	By \eqref{eq:complete-phase-space-derivative},
	\[
	[\mathfrak D_X,\cK_{k_I}]f
	=\int_{\bR^3}\bigl(Xk_I+(M_Xv)\cdot\nabla_vk_I\bigr)f\,\dd v'
	-\int_{\bR^3}k_I\,(M_Xv')\cdot\nabla_{v'}f\,\dd v'.
	\]
	Integrating the last term by parts in $v'$, using $\operatorname{div}_{v'}(M_Xv')=\operatorname{tr}M_X$, gives
	$[\mathfrak D_X,\cK_{k_I}]=\cK_{k_{XI}}$ with $k_{XI}$ of
	\eqref{eq:iterated-compact-kernel-recursion}.  Iteration gives $\operatorname{ad}_K(\cK_0)=m_K+\cK_{k_K}$; applying \eqref{eq:ad-word} with $\mathrm A=\cK_0$ yields \eqref{eq:compact-word-commutator-expansion}.
	
	Every $m_K$ is a Gaussian--polynomial profile with $x$-coefficients bounded on
	$\overline\cD$.  By \eqref{eq:local-H1-H_sigma}, for $\psi\in H_\sigma^1$,
	\[
	|\langle w^\ell m_Kf,\psi\rangle_v|
	\le\|w^{\ell+1/2}m_K\|_{L_v^\infty}\|f\|_{L_v^2}\|w^{-1/2}\psi\|_{L_v^2}
	\le C\|f\|_{L_v^2}\|\psi\|_{H_\sigma^1},
	\]
	which yields $\|w^\ell m_Kf\|_{H_\sigma^{-1}}\le C\|f\|_{L_v^2}$.
	
	By \eqref{eq:compact-integral-kernel}, $k_0$ is a finite sum of products of
	Gaussian--polynomial functions of $v$ and $v'$ with a function of $z=v-v'$ homogeneous
	of degree $-1$ and smooth on the unit sphere.  On functions of $z$,
	\[
	(M_Xv)\cdot\nabla_v+(M_Xv')\cdot\nabla_{v'}=(M_Xz)\cdot\nabla_z,
	\]
	so the singular factor retains degree $-1$.  The remaining derivatives preserve the Gaussian--polynomial form, while $X$ differentiates only its $x$-coefficients.  Consequently, for $|K|\le\mathrm N$, the coefficients remain bounded on $\overline\cD$ and, after reducing the Gaussian exponent,
	\[
	w(v)^{\ell+1/2}|k_K(x;v,v')|\le C\ee^{-c(|v|^2+|v'|^2)}|v-v'|^{-1}.
	\]
	
	The right side has uniformly bounded integrals in either velocity variable.  Schur's test and \eqref{eq:Hsigma-dual-weighted} give
	\[
	\|w^\ell\cK_{k_K}f\|_{H_\sigma^{-1}}
	\le C\|w^{\ell+1/2}\cK_{k_K}f\|_{L_v^2}
	\le C\|f\|_{L_v^2}.
	\]
	
	Both estimates are uniform in $x$.  Summation over $S$ with $I_{S^c}=J$ and integration in $x$ give \eqref{eq:iterated-compact-full-operator-bound} on the smooth core; density gives the asserted extension to $L^2_{x,v}$.
\end{proof}

\subsubsection{The collision expansion}\label{ss:diff-interior}

\begin{prop}[Ordered expansion of $N$]
	\label{prop:differentiated-interior}
	Let $I=(\nu_1,\ldots,\nu_m)\in\mathfrak A_\bullet^m$, $1\le m\le\mathrm N$, with
	position subwords defined by \eqref{eq:ordered-position-subwords}.  Set
	$\mathfrak C_{\mathrm t}:=0$ and, for $\nu\in\mathfrak A$,
	\[
	\mathfrak C_\nu[g,f]
	:=N_{\mathscr D_\nu^{\rm ker}\phi}[g,f]
	+N[\mathfrak p_\nu g,f]+N[g,\mathfrak p_\nu f]-\mathfrak p_\nu N[g,f].
	\]
	Then, on the smooth core,
	\begin{equation}
		\label{eq:collision-remainder-expansion}
		\mathfrak D_I N[g,f]
		=\sum_{S\subset[m]}
		N[\mathfrak D_{I_S}g,\mathfrak D_{I_{S^c}}f]
		+\mathfrak R_I[g,f],
	\end{equation}
	where
	\begin{align}
		\mathfrak R_\nu[g,f]
		&:=\mathfrak C_\nu[g,f],
		\label{eq:deformation-remainder-first}\\
		\mathfrak R_{\nu I}[g,f]
		&:=\mathfrak D_\nu\mathfrak R_I[g,f]
		+\sum_{S\subset[m]}
		\mathfrak C_\nu[\mathfrak D_{I_S}g,\mathfrak D_{I_{S^c}}f].
		\label{eq:deformation-remainder-recurrence}
	\end{align}
	Every term of $\mathfrak R_I[g,f]$ has the form
	\begin{equation}
		\label{eq:deformation-remainder-structure}
		pN_{\mathsf K_A}[p_1\mathfrak D_{J_1}g,p_2\mathfrak D_{J_2}f],
		\qquad
		|A|\le m,
		\quad
		|J_1|+|J_2|\le m-1,
		\quad
		\deg_vp+\deg_vp_1+\deg_vp_2\le2m,
	\end{equation}
	where $\mathsf K_A$ is defined by \eqref{eq:iterated-coulomb-deformation-kernels}
	and $p,p_1,p_2$ are polynomials in $v$, as in Lemma~\ref{lem:deformation-degrees}.  Each such term also satisfies
	\eqref{eq:deformation-combined-grading} with $J$ replaced by $A$:
	\[
	|A|+\tfrac12(\deg_vp+\deg_vp_1+\deg_vp_2)\le m .
	\]  Since $\mathfrak C_{\mathrm t}=0$,
	every term of $\mathfrak R_I[g,f]$ is a term of $\mathfrak R_{I'}[\partial_t^{\,i}g,\partial_t^{\,j(I)-i}f]$ for some $0\le i\le j(I)$, and
	$|A|\le|I'|$ in \eqref{eq:deformation-remainder-structure}.  Each term of
	\eqref{eq:deformation-remainder-structure} satisfies
	\eqref{eq:deformation-weighted-output} with $\mathfrak D_{J_1}g$ and $\mathfrak D_{J_2}f$ in place of $g$ and $f$.
\end{prop}

\begin{proof}
	For one letter, \eqref{eq:complete-conjugated-kernel-deformation} at $\mathsf K=\phi$,
	\eqref{eq:generalized-conjugated-landau-form} and
	\eqref{eq:complete-kernel-deformation-operator} give
	\eqref{eq:collision-remainder-expansion} and \eqref{eq:deformation-remainder-first}.
	At $\nu=\mathrm t$ the one-letter identity is the Leibniz rule in $t$.
	
	Assume \eqref{eq:collision-remainder-expansion} for $I$.  By the one-letter identity,
	for each $S\subset[m]$,
	\begin{align*}
		\mathfrak D_\nu
		N[\mathfrak D_{I_S}g,\mathfrak D_{I_{S^c}}f]
		={}&N[\mathfrak D_{\nu I_S}g,\mathfrak D_{I_{S^c}}f]
		+N[\mathfrak D_{I_S}g,\mathfrak D_{\nu I_{S^c}}f]+\mathfrak C_\nu[\mathfrak D_{I_S}g,\mathfrak D_{I_{S^c}}f].
	\end{align*}
	Summation over $S\subset[m]$ produces the two classes of subsets of $[m+1]$, according as the new position belongs to $S'$ or its complement.  The remaining terms are exactly \eqref{eq:deformation-remainder-recurrence}, proving \eqref{eq:collision-remainder-expansion} for $\nu I$.
	
	By \eqref{eq:deformation-remainder-first}--\eqref{eq:deformation-remainder-recurrence}, each remainder contains a kernel deformation or a Maxwellian-conjugation factor, and its two differentiated arguments have total order at most $m-1$.  By \eqref{eq:iterated-coulomb-deformation-kernels}, each letter raises $|A|$ by at most one, so $|A|\le m$.  By
	\eqref{eq:complete-kernel-deformation-operator}, each new Maxwellian factor has velocity
	degree two; differentiating an existing polynomial does not increase its degree.  So
	the total degree is at most $2m$, which is \eqref{eq:deformation-remainder-structure}.
	
	Every remainder term is produced from $N[a,b]=N_{\mathsf K_\varnothing}[a,b]$, with $a$
	and $b$ time derivatives of $g$ and $f$, by the $|I'|\le m$ spatial letters of $I$. Temporal derivatives act on the two arguments without changing the kernel: $\mathfrak C_{\mathrm t}=0$, and the complete fields, kernels, and polynomial coefficients are time independent.  Applying Lemma~\ref{lem:deformation-degrees} to $N_{\mathsf K_\varnothing}[a,b]$ with $k=|I'|$ gives
	\[
	|A|+\tfrac12(\deg_vp+\deg_vp_1+\deg_vp_2)\le|I'|\le m,
	\]
	so \eqref{eq:deformation-combined-grading} holds with $J$ replaced by $A$, and $|A|\le|I'|$.
	
	Each term has the form required by Proposition~\ref{prop:deformation-algebra}; applying \eqref{eq:deformation-weighted-output} with $\mathfrak D_{J_1}g$ and $\mathfrak D_{J_2}f$ in the two slots gives the asserted bound.
\end{proof}

\subsubsection{Temporal letters}\label{ss:time-letters}

For a distribution $u$ of order at most $K_\cD+1-|J|$, $\mathfrak D_Ju$ is defined by
$\langle\mathfrak D_Ju,\varphi\rangle:=\langle u,\mathfrak D_J^\dagger\varphi\rangle$,
$\varphi\in C_c^\infty$, with $\mathfrak D_J^\dagger$ of \eqref{eq:ordered-complete-formal-adjoint};
the complete pullbacks of \S\ref{app:diffuse-lifts} act in the same way, through their
adjoints, on distributions of order at most $K_\cD$.  This covers every distribution of order
at most $\mathrm N$ and every ordered word $J$ with $|J|\le\mathrm N$.
On $\partial\cD$ and $\gamma^\pm$, which are of class
$C^{K_\cD+2}$, and on their products with $(0,T)$, the same holds for $\tau_J$ and $\mathfrak D_J^\pm$, through their formal
adjoints there, with $C_c^{K_\cD+1}$ test functions in place of $C_c^\infty$.

\begin{lemma}
	\label{lem:dyn-time-commutation}
	Let $0<T<\infty$.  For every ordered word $I$ over $\mathfrak A$ with
	$|I|\le\mathrm N$,
	\begin{equation}
		\label{eq:dyn-time-commutation}
		[\partial_t,\mathfrak D_I]=0
		\quad\text{in }\cD'\bigl((0,T)\times\cD\times\bR^3\bigr).
	\end{equation}
	Consequently every mixed word may be normal ordered: for
	$I\in\mathfrak A_\bullet^m$ with spatial subword $I'$ and temporal order
	$j=j(I)$,
	\begin{equation}
		\label{eq:dyn-normal-ordering}
		\mathfrak D_I=\mathfrak D_{I'}\,\partial_t^{\,j}
		=\partial_t^{\,j}\,\mathfrak D_{I'}
		\quad\text{in }\cD' .
	\end{equation}
	Moreover,
	\begin{equation}
		\label{eq:dyn-spacetime-leibniz}
		\mathfrak D_{\mathrm t}N[a,b]
		=N[\mathfrak D_{\mathrm t}a,b]+N[a,\mathfrak D_{\mathrm t}b],
		\qquad
		[\mathfrak D_{\mathrm t},v\cdot\nabla_x]=0,
		\qquad
		[\mathfrak D_{\mathrm t},L]=0
		\quad\text{in }\cD' ,
	\end{equation}
	where the first identity holds whenever, for some $\tfrac14\le\varpi<\tfrac12$, $\mu^\varpi a,\mu^\varpi\partial_ta\in L^2_{\rm loc}((0,T)\times\cD;L^2_v)$ and $b,\partial_tb\in L^2_{\rm loc}((0,T)\times\cD;H^1_\sigma)$.
\end{lemma}

\begin{proof}
	Let $\nu\in\mathfrak A$, let $u$ be a distribution on $(0,T)\times\cD\times\bR^3$ of order at most $K_\cD-1$, and let
	$\Phi\in C_c^\infty((0,T)\times\cD\times\bR^3)$.  By
	\eqref{eq:complete-field-formal-adjoint},
	\begin{align*}
		\langle\partial_t\mathfrak D_\nu u,\Phi\rangle
		&=-\langle\mathfrak D_\nu u,\partial_t\Phi\rangle
		=\langle u,(\mathfrak D_\nu+d_\nu)\partial_t\Phi\rangle\\
		&=\langle u,\partial_t(\mathfrak D_\nu+d_\nu)\Phi\rangle
		=\langle\mathfrak D_\nu\partial_tu,\Phi\rangle ,
	\end{align*}
	where the third equality holds since $\mathfrak f_\nu$ and $d_\nu$ in
	\eqref{eq:complete-phase-space-derivative} and \eqref{eq:phase-divergence} are
	independent of $t$.
	So $[\partial_t,\mathfrak D_\nu]u=0$.  Since $\mathfrak D_\nu$ raises the order by at most one,
	iteration over the letters of $I$ gives \eqref{eq:dyn-time-commutation} on distributions of
	order at most $K_\cD-|I|$, in particular of order at most $\mathrm N$.
	
	Transposing the letters of a
	mixed word until all temporal letters stand to the right gives
	\eqref{eq:dyn-normal-ordering}.
	
	For the first identity of \eqref{eq:dyn-spacetime-leibniz}, let $J\Subset(0,T)$ and
	$K\Subset\cD$.  By \eqref{eq:velalg} at $\ell=0$, applied at each $(t,x)$,
	\[
	N:\ L^2(K;\mu^{-\varpi}L^2_v)\times L^2(K;H^1_\sigma)\to L^1(K;H^{-1}_\sigma)
	\]
	is bounded bilinear, and likewise with $J\times K$ in place of $K$. Here
	$\mu^{-\varpi}L^2_v$ is normed by $\|\mu^\varpi\cdot\|_{L^2_v}$.  The mollifications
	$a_\eta,b_\eta$ of $a,b$ in $t$ are $C^1$ on $J$ with values in these spaces.  The
	product rule for a bounded bilinear map gives the identity at $(a_\eta,b_\eta)$ at each
	$t\in J$.  As $\eta\to0$, $a_\eta$, $\partial_ta_\eta$, $b_\eta$ and $\partial_tb_\eta$
	converge to $a$, $\partial_ta$, $b$ and $\partial_tb$ in these spaces over $J\times K$.
	So $N[a_\eta,b_\eta]$ and the right side at $(a_\eta,b_\eta)$ converge in
	$L^1(J\times K;H^{-1}_\sigma)$, and the identity passes to the limit in $\cD'$.
	
	The other two identities of \eqref{eq:dyn-spacetime-leibniz} hold since
	$v\cdot\nabla_x$ and $L$ are independent of $t$.
\end{proof}

\subsubsection{The kinetic operator}\label{ss:complete-operators}

\begin{lemma}[Complete operator commutators]
	\label{lem:nested-complete-commutators}
	For an ordered word $I=(\nu_1,\ldots,\nu_r)$ over $\mathfrak A$,
	$1\le r\le\mathrm N$, and $X_j:=\tau_{\nu_j}$, set
	\[
	\mathfrak K_{X_r}:=[\mathfrak D_{X_r},(v\cdot\nabla_x+L_g)],\qquad
	\mathfrak K_{X_1\mid\cdots\mid X_r}
	:=[\mathfrak D_{X_1},\mathfrak K_{X_2\mid\cdots\mid X_r}].
	\]
	Then
	\begin{equation}
		\label{eq:nested-commutator-time-transparency}
		[\mathfrak D_I,\partial_t+v\cdot\nabla_x+L_g]
		=[\mathfrak D_I,v\cdot\nabla_x+L_g],
		\qquad
		\operatorname{ad}_I(\partial_t+v\cdot\nabla_x+L_g)
		=\operatorname{ad}_I(v\cdot\nabla_x+L_g).
	\end{equation}
	Moreover, let $g,u\in C^\infty(\bR_t\times\cD\times\bR^3)$ with
	$\operatorname{supp}g\cup\operatorname{supp}u\subset\bR_t\times K$ for a compact
	$K\Subset\cD\times\bR^3$.  Then
	\begin{equation}
		\label{eq:nested-commutator-decomposition}
		\mathfrak K_{X_1\mid\cdots\mid X_r}u
		=-N[\mathfrak D_{X_1}\cdots\mathfrak D_{X_r}g,u]
		+\mathfrak T_Iu+\mathfrak L_Iu+\mathfrak C_Iu
		+\mathfrak R_I[g,u].
	\end{equation}
	Here \(\mathfrak T_I\) is a finite sum of the vertical first-order operators
	in Proposition~\ref{prop:ordered-transport-weight}; \(\mathfrak L_I\) is the
	nested local commutator defined by
	\eqref{eq:nested-local-complete-definition};
	\(\mathfrak C_I=\mathscr K_{I,\varnothing}\), with \(\mathscr K_{I,J}\) as in
	Lemma~\ref{lem:compact-maxwellian-complete-commutators}; and every term of
	\(\mathfrak R_I[g,u]\), which differs from the remainder of \eqref{eq:collision-remainder-expansion}, has the form
	\begin{equation}
		\label{eq:nested-deformation-form}
		a(x)\,p\,N_{\mathsf K_J}
		\bigl[p_1\mathfrak D_{J_1}g,p_2u\bigr],
		\quad |J|\le r,
		\quad |J_1|\le r-1,
		\quad \deg_v p+\deg_v p_1+\deg_v p_2\le2r,
	\end{equation}
	with $a$ and the coefficients of $p,p_1,p_2$ bounded on $\overline\cD$.  The term
	$-N[\mathfrak D_{X_1}\cdots\mathfrak D_{X_r}g,u]$ is the only one that contains a complete derivative of order \(r\) of $g$.
	The terms \(\mathfrak T_Iu\), \(\mathfrak L_Iu\),
	\(\mathfrak C_Iu\) and those of \(\mathfrak R_I[g,u]\) satisfy
	\eqref{eq:ordered-transport-weight},
	\eqref{eq:nested-local-complete-bound},
	\eqref{eq:iterated-compact-full-operator-bound} and
	\eqref{eq:deformation-weighted-output}, respectively.  For other $g$ and $u$,
	$\mathfrak K_{X_1\mid\cdots\mid X_r}u$ denotes the right side of
	\eqref{eq:nested-commutator-decomposition} whenever $u\in L^2_{x,v}$ and, for almost every $x$, $w^{2r+2}u\in H_\sigma^1$
	and $\mu^\varpi\mathfrak D_Jg\in L_v^2$ for $|J|\le r$, for some $\tfrac14\le\varpi<\tfrac12$.  By these bounds at $\ell=0$, with
	\eqref{eq:real-weight-Hsigma-monotonicity}, and by \eqref{eq:deformation-weighted-output} for $N[\mathfrak D_{X_1}\cdots\mathfrak D_{X_r}g,u]$, every term of the
	right side of \eqref{eq:nested-commutator-decomposition} then lies in $H_\sigma^{-1}$ at almost every $x$.  For $g$ and $u$ depending on $t$, these conditions are imposed, and $\mathfrak K_{X_1\mid\cdots\mid X_r}u$ is defined, at almost every $t$.
\end{lemma}

\begin{proof}
	By \eqref{eq:dyn-time-commutation}, $[\mathfrak D_\nu,\partial_t]=0$ for every letter
	$\nu$.  So $[\mathfrak D_I,\partial_t]=0$ and $\operatorname{ad}_I(\partial_t)=0$, which
	give \eqref{eq:nested-commutator-time-transparency}.
	
	For \eqref{eq:nested-commutator-decomposition}, we argue by induction on $r$. For $r=1$, Proposition~\ref{prop:ordered-transport-weight},
	\eqref{eq:local-L0-complete-expansion} and
	\eqref{eq:compact-word-commutator-expansion} give the transport, local and compact
	terms (recall $L_g=L-N[g,\cdot]$ and $L=\cA_\mu+\cK_0$).  For the Coulomb part, by \eqref{eq:complete-conjugated-kernel-deformation},
	\[
	[\mathfrak D_X,N[g,\cdot]]u
	=N[\mathfrak D_Xg,u]+\mathfrak C_\nu[g,u],
	\qquad X=\tau_\nu,
	\]
	where each term of \(\mathfrak C_\nu[g,u]\) has the form
	\eqref{eq:nested-deformation-form}; the two occurrences of
	\(N[g,\mathfrak D_Xu]\) cancel.
	
	Let $r<\mathrm N$, assume \eqref{eq:nested-commutator-decomposition} at order $r$, and
	commute with \(\mathfrak D_X\).  The original Coulomb term gives
	\[
	-[\mathfrak D_X,N[\mathfrak D_Ig,\cdot]]u
	=-N[\mathfrak D_X\mathfrak D_Ig,u]
	-\mathfrak C_\nu[\mathfrak D_Ig,u];
	\]
	the two occurrences of $N[\mathfrak D_Ig,\mathfrak D_Xu]$ cancel.  On a term
	\eqref{eq:nested-deformation-form}, \(\mathfrak D_X\) falls on \(\mathfrak D_{J_1}g\),
	on \(\mathsf K_J\), on one of \(p,p_1,p_2\), or on \(a\).  By
	\eqref{eq:complete-conjugated-kernel-deformation} and the Leibniz rule, as in
	Lemma~\ref{lem:deformation-degrees}, the results have the form
	\eqref{eq:nested-deformation-form}, with $r+1$ in place of $r$.  The transport terms remain among the operators of
	Proposition~\ref{prop:ordered-transport-weight}.  By the definition of $\mathscr K_{I,J}$, $\mathfrak C_I=\operatorname{ad}_I(\cK_0)$, so
	$[\mathfrak D_X,\mathfrak C_I]=\operatorname{ad}_{XI}(\cK_0)=\mathfrak C_{XI}$.
	The local term is
	given by Lemma~\ref{lem:nested-local-complete-commutators}\,(\ref{it:local-nested-commutators}).
	The bounds follow from \eqref{eq:ordered-transport-weight},
	\eqref{eq:nested-local-complete-bound},
	\eqref{eq:iterated-compact-full-operator-bound} and
	\eqref{eq:deformation-weighted-output}.
\end{proof}

\subsection{Boundary algebra}\label{ss:complete-boundary-calculus}

We derive weak trace and boundary moment identities, spectral-cutoff commutator bounds, and differentiated diffuse boundary laws. The ordered wall commutators admit a finite-rank representation, while the spectral commutators satisfy bounds uniform in the cutoff. These estimates yield the diffuse law for the traces of differentiated solutions.

\begin{lemma}[Boundary Leibniz rule]
	\label{lem:boundary-word-Leibniz}
	Let $I$ be an ordered word over $\mathfrak A$ with $|I|=m\le\mathrm N$, let
	$a\in C^m(\gamma^-)$, and let $b\in L^1_{\rm loc}(\partial\cD)$, regarded on $\gamma^-$ as a
	function independent of $v$, satisfy $\tau_{I_S}b\in L^1_{\rm loc}(\partial\cD)$ for every
	$S\subset[m]$.  Then
	\begin{equation}
		\label{eq:boundary-word-Leibniz}
		\mathfrak D_I^-(ab)
		=\sum_{I_1\sqcup I_2=I}
		\bigl(\mathfrak D_{I_1}^-a\bigr)\bigl(\tau_{I_2}b\bigr)
		\qquad\text{in }\cD'(\gamma^-),
	\end{equation}
	with the two-coloring convention \eqref{eq:ordered-two-coloring-convention}.  On $(0,T)$, with
	$\mathfrak D_{\mathrm t}^\pm:=\partial_t$ and $\tau_{\mathrm t}:=\partial_t$, let $I$ be an ordered word
	over $\mathfrak A_\bullet$ with $|I'|\le\mathrm N$, let $a$ on $(0,T)\times\gamma^-$ have continuous derivatives of order at
	most $j(I)$ in $t$ and at most $|I'|$ in the variables of $\gamma^-$, and let $b\in L^1_{\rm loc}((0,T)\times\partial\cD)$ satisfy
	$\tau_{I_S}b\in L^1_{\rm loc}((0,T)\times\partial\cD)$ for every $S\subset[|I|]$; then \eqref{eq:boundary-word-Leibniz} holds in
	$\cD'((0,T)\times\gamma^-)$.
\end{lemma}

\begin{proof}
	By \eqref{eq:tangential-field-properties} and
	\eqref{eq:complete-derivative-preserves-grazing}, $\mathfrak f_\nu$ is tangent to
	$\gamma^-$, and its velocity part annihilates $b$; so $\mathfrak D_\nu^-b=\tau_\nu b$.  For
	$a'\in C^1(\gamma^-)$ and $g\in L^1_{\rm loc}(\gamma^-)$ with
	$\mathfrak D_\nu^-g\in L^1_{\rm loc}(\gamma^-)$, testing against $a'\psi$,
	$\psi\in C_c^1(\gamma^-)$, gives
	\begin{align}\label{260912439}	
		\mathfrak D_\nu^-(a'g)=(\mathfrak D_\nu^-a')g+a'\mathfrak D_\nu^-g\qquad \text{in }\cD'(\gamma^-).
	\end{align}
	Write $I=\nu J$ and differentiate the identity for $J$. For each term, \eqref{260912439} with $a'=\mathfrak D^-_{J_1}a$ and $g=\tau_{J_2}b$, for which $\mathfrak D_\nu^-g=\tau_{\nu J_2}b$, gives the terms
	$(\mathfrak D^-_{\nu J_1}a)(\tau_{J_2}b)$ and $(\mathfrak D^-_{J_1}a)(\tau_{\nu J_2}b)$ of \eqref{eq:boundary-word-Leibniz} for $I$, in which the first position of $I$ belongs to $S$ and to $S^c$, respectively.  On $(0,T)$ the same induction applies, since $\mathfrak D^-_{J_1}a$ is $C^1$ in the variables of the letter $\nu$.
	The temporal letter obeys the same product rule, with $\partial_tb=\tau_{\mathrm t}b$.
\end{proof}

\subsubsection{Boundary commutators}\label{ss:differentiated-wall-law}

For a boundary profile $a$ with the envelope
\eqref{eq:boundary-gaussian-polynomial-envelope} and for a profile
$p=p(x,v)$, set
\begin{align*}
	\cA_aq(x)&:=\int_{n_x\cdot v>0}q(x,v)\,a(x,v)\,n_x\cdot v\dd v,
	\\
	\cA_{p,\mu}q(x)&:=\int_{n_x\cdot v>0}q(x,v)\,p(x,v)\sqrt{\mu(v)}\,n_x\cdot v\dd v,
\end{align*}
so that $\cA_{1,\mu}=\cA_0$ of \eqref{eq:outgoing-diffuse-moment}.  For $\nu\in\mathfrak A$ set
\begin{align}
	r_\nu(x,v)
	&:=\operatorname{div}_x\tau_\nu(x)
	+((\nabla\tau_\nu(x))v)\cdot\nabla_v\log\sqrt\mu(v),
	\label{eq:rnu-profile}\\
	\mathscr E_\nu a
	&:=\mathfrak D_\nu a+\bigl(\operatorname{div}_x\tau_\nu\bigr)a ,
	\label{eq:outgoing-profile-generator}\\
	\mathscr E_{\nu,\mu}p
	&:=\mathfrak D_\nu p+r_\nu p ,\nonumber
\end{align}
so that $\mathscr E_\nu(p\sqrt\mu)=(\mathscr E_{\nu,\mu}p)\sqrt\mu$.  For the temporal letter,  we set $\mathscr E_{\mathrm t}:=0$ on $a$ and $p$ independent of $t$.
On $\partial\cD$, $\tau_\nu\cdot n=0$ and, by \eqref{eq:signed-distance-tangent-extension},
$n\cdot(\nabla\tau_\nu)n=0$; hence
$\operatorname{div}_{\partial\cD}\tau_\nu=\operatorname{div}_x\tau_\nu$. Since $\operatorname{div}_v((\nabla\tau_\nu)v)=\operatorname{div}_x\tau_\nu$ (see \eqref{eq:phase-divergence}) and $d_\nu=2\operatorname{div}_x\tau_\nu$, we obtain, on
$\partial\cD\times\bR^3$,
\begin{equation}
	\label{eq:boundary-phase-divergence}
	\operatorname{div}_{\partial\cD}\tau_\nu
	+\operatorname{div}_v((\nabla\tau_\nu)v)=d_\nu.
\end{equation}
For an operator $B:L^2(\gamma^+)\to L^2(\gamma^-)$, a letter $\nu\in\mathfrak A$ and an ordered
word $I=(\nu_1,\ldots,\nu_m)$ over $\mathfrak A$, set on the smooth trace core
\[
\operatorname{ad}^{\gamma}_{\nu}(B):=\mathfrak D_\nu^-B-B\mathfrak D_\nu^+,
\qquad
\operatorname{ad}^{\gamma}_I:=\operatorname{ad}^{\gamma}_{\nu_1}\cdots
\operatorname{ad}^{\gamma}_{\nu_m},
\qquad
\operatorname{ad}^{\gamma}_\varnothing:=\operatorname{Id},
\]
the boundary form of $\operatorname{ad}_I$ in \eqref{eq:ad-word}, and
\begin{equation}
	\label{eq:ordered-boundary-commutator-definition}
	P_I:=\operatorname{ad}^{\gamma}_I(P_w),
	\qquad\text{that is,}\qquad
	P_\varnothing=P_w,
	\quad
	P_{\nu I}=\mathfrak D_\nu^-P_I-P_I\mathfrak D_\nu^+ .
\end{equation}
For the truncated wall operator \eqref{eq:boundary-spectral-diffuse-operator},
\begin{equation*}
	P_{I,N}:=\operatorname{ad}^{\gamma}_I(P_{w,N}),
	\qquad |I|\le\mathrm N .
\end{equation*}
For the time-dependent operator $P_{w,N}(t)$, with $\mathfrak D_{\mathrm t}^\pm=\partial_t$,
\begin{equation}
	\label{eq:dyn-mixed-boundary-commutator}
	P_{I,N}:=\operatorname{ad}^{\gamma}_I(P_{w,N}),
	\qquad
	I\text{ an ordered word over }\mathfrak A_\bullet,\ |I|\le\mathrm N .
\end{equation}

\begin{lemma}[Moment and profile bounds]
	\label{lem:boundary-profile-bounds}
	\label{lem:moment-commutator}
	Let \(\nu\in\mathfrak A\).  For \(p\) and compactly velocity-supported \(q\), both of class \(C^1\), such that \(p\sqrt\mu\) satisfies \eqref{eq:boundary-gaussian-polynomial-envelope},
	\begin{equation}
		\label{eq:moment-commutator}
		\tau_\nu\cdot\nabla_x\cA_{p,\mu}q
		-\cA_{p,\mu}(\mathfrak D_\nu q)
		=\cA_{\mathscr E_{\nu,\mu}p,\mu}q.
	\end{equation}
	For an ordered word \(I=(\nu_1,\ldots,\nu_m)\) over \(\mathfrak A\), \(m\le\mathrm N\), set
	\begin{equation}
		\label{eq:moment-profile-recursion}
		p_\varnothing:=1,
		\qquad
		p_I:=\mathscr E_{\nu_1,\mu}\cdots
		\mathscr E_{\nu_m,\mu}1.
	\end{equation}
	Then \(p_I\sqrt\mu\) satisfies
	\eqref{eq:boundary-gaussian-polynomial-envelope}.
	If a boundary profile \(a\) satisfies
	\eqref{eq:boundary-gaussian-polynomial-envelope}, then for every finite
	\(L\ge0\) and every ordered word \(K\) over \(\mathfrak A\) with \(|K|\le\mathrm N\),
	\begin{equation}
		\label{eq:boundary-profile-flux-bound}
		\sup_{x\in\partial\cD}
		\int_{\{v:\,\pm n_x\cdot v>0\}}
		w^{2L}|\mathfrak D_K^{\pm}a(x,v)|^2
		|n_x\cdot v|\,\dd v
		\le C,
	\end{equation}
	where $C$ depends on $a$ only through the constants $C_a,c_a$ of
	\eqref{eq:boundary-gaussian-polynomial-envelope}.
\end{lemma}

\begin{proof}
	Set \(\alpha_+:=\max\{\bar n\cdot v,0\}\) on the signed-distance collar.  By
	\eqref{eq:complete-derivative-preserves-grazing} and the Lipschitz chain rule,
	\(\mathfrak D_\nu\alpha_+=0\) a.e.  Put \(F:=qp\sqrt\mu\,\alpha_+\), compactly supported
	in \(v\).  Integration by parts in \(v\), with
	\(\operatorname{div}_v((\nabla\tau_\nu)v)=\operatorname{div}_x\tau_\nu\), gives
	\[
	\tau_\nu\cdot\nabla_x\int_{\bR^3}F\,\dd v
	=\int_{\bR^3}\bigl(\mathfrak D_\nu F+(\operatorname{div}_x\tau_\nu)F\bigr)\dd v .
	\]
	Since
	\begin{gather*}
		\mathfrak D_\nu F=(\mathfrak D_\nu q)p\sqrt\mu\,\alpha_+
		+q\,\mathfrak D_\nu(p\sqrt\mu)\,\alpha_+,\\
		\mathfrak D_\nu(p\sqrt\mu)+(\operatorname{div}_x\tau_\nu)p\sqrt\mu
		=\mathscr E_\nu(p\sqrt\mu)=(\mathscr E_{\nu,\mu}p)\sqrt\mu,
	\end{gather*}
	this is \eqref{eq:moment-commutator} on $\partial\cD$, where $\bar n=n$.
	
	By \eqref{eq:boundary-profile-chart-representation} and
	\eqref{eq:rotating-frame-complete-lift}, each complete derivative in an
	adapted boundary chart has the form
	\[
	b(y')\cdot\nabla_{y'}+(B(y')\xi)\cdot\nabla_\xi,
	\]
	with bounded \(b,B\) of class \(C^{K_\cD}\).  Since \(K_\cD\ge2\mathrm N\), all coefficients
	created by at most \(\mathrm N\) complete derivatives and \(\mathrm N\) further
	chart derivatives are bounded on the fixed atlas.  Hence, for every ordered word \(K\) over \(\mathfrak A\) with \(|K|\le\mathrm N\),
	\(\mathfrak D_K\) is a finite sum of chart derivatives of order at most \(|K|\), with coefficients polynomial in \(\xi\).  By
	\eqref{def:reference-maxwellian},
	\(\nabla_v\log\sqrt\mu=-v\).  Thus \(r_\nu\) in
	\eqref{eq:rnu-profile} is polynomial in \(v\), with boundary
	coefficients of class \(C^{K_\cD}\).  Induction in \eqref{eq:moment-profile-recursion} shows that
	each \(p_I\) is polynomial in \(v\), with coefficients of class \(C^{K_\cD+1-|I|}\).  Every chart derivative through order
	\(\mathrm N\) of \(p_I\sqrt\mu\) is therefore a polynomial in \(\xi\)
	times \(\ee^{-|\xi|^2/2}\), hence
	\[
	\max_{\kappa\in\mathfrak K_{\p}}\max_{|\alpha|+|\beta|\le\mathrm N}|\partial_\xi^\alpha\partial_{y'}^\beta(p_I\sqrt\mu)_\kappa(y',\xi)|\le C\ee^{-|\xi|^2/4},
	\]
	which proves
	\eqref{eq:boundary-gaussian-polynomial-envelope} for \(p_I\sqrt\mu\).
	
	The same chart representation and
	\eqref{eq:boundary-gaussian-polynomial-envelope} give
	\[
	|w^L\mathfrak D_K^{\pm}a(x,v)|\le C\ee^{-c_a|v|^2/2},
	\qquad |K|\le\mathrm N .
	\]
	Integration in \(v\) proves
	\eqref{eq:boundary-profile-flux-bound}.
\end{proof}

A \emph{compact-velocity $C^2$ phase test} is a function $\Psi\in C^2(\overline\cD\times\bR^3)$ with compact support in $\overline\cD\times\bR^3$.  For an open interval $\mathcal T$, a \emph{compact-time, compact-velocity $C^2$ phase test} is a function $\Psi\in C^2(\mathcal T\times\overline\cD\times\bR^3)$ with compact support in $\mathcal T\times\overline\cD\times\bR^3$.

\begin{prop}[Weak complete-field boundary generator]
	\label{prop:weak-boundary-generator}\,
	
	\begin{enumerate}
		\item\label{it:weak-generator-stationary}
		Let $\nu\in\mathfrak A$ and let $q_J,q_{\nu J}\in\cH_{\kin,\sigma}$ satisfy
		\[
		\mathfrak D_\nu q_J=q_{\nu J}\qquad\text{in }\cD'(\cD\times\bR^3).
		\]
		Then, for every compact-velocity $C^2$ phase test $\Psi$ up to the physical
		boundary,
		\begin{equation}
			\label{eq:boundary-generator-identity}
			\int_\gamma (q_{\nu J})_\gamma\Psi\,(n\cdot v)\,\dd S\dd v
			=-\int_\gamma(q_J)_\gamma
			(\mathfrak D_\nu+d_\nu)\Psi\,(n\cdot v)\,\dd S\dd v.
		\end{equation}
		The identity \eqref{eq:boundary-generator-identity} holds also on $\gamma^+$ and on
		$\gamma^-$ separately.
		
		\item\label{it:weak-generator-evolutionary}
		Let $\nu\in\mathfrak A_\bullet$ and let $z_J,z_{\nu J}\in\cH_{\kin,\sigma}(\cD_T\times\bR_v^3)$ satisfy
		\[
		\mathfrak D_\nu z_J=z_{\nu J}
		\qquad\text{in }\cD'((0,T)\times\cD\times\bR^3).
		\]
		Then, for every compact-time, compact-velocity $C^2$ phase test $\Psi$ up to the
		physical boundary,
		\begin{equation}
			\label{eq:dynamic-boundary-generator}
			\int_{\gamma_T}(z_{\nu J})_\gamma\Psi\,
			\operatorname{sgn}(n\cdot v)\dd\gamma\dd t
			=-\int_{\gamma_T}(z_J)_\gamma
			(\mathfrak D_\nu+d_\nu)\Psi\,
			\operatorname{sgn}(n\cdot v)\dd\gamma\dd t.
		\end{equation}
		The identity \eqref{eq:dynamic-boundary-generator} holds also on $\gamma_T^+$ and on
		$\gamma_T^-$ separately.  At $\nu=\mathrm t$, where $d_{\mathrm t}=0$, it reads
		\begin{equation}
			\label{eq:temporal-boundary-generator}
			\partial_t(z_J)_{\gamma_T^\pm}=(z_{\mathrm tJ})_{\gamma_T^\pm}
			\qquad\text{in }\cD'\bigl((0,T);L^2(\gamma^\pm)\bigr).
		\end{equation}
		Consequently, if $1\le m\le\mathrm N$ and
		$f_0,\ldots,f_m\in\cH_{\kin,\sigma}(\cD_T\times\bR_v^3)$ satisfy $\partial_tf_j=f_{j+1}$ in
		$\cD'((0,T)\times\cD\times\bR^3)$ for $0\le j\le m-1$, then
		\begin{align}
			\label{eq:dynamic-trace-time-chain-boundary}
			\partial_t(f_j)_{\gamma^\pm}
			&=(f_{j+1})_{\gamma^\pm}
			\quad\text{in }\cD'\bigl((0,T);L^2(\gamma^\pm)\bigr),
			\qquad 0\le j\le m-1,\\
			\partial_t^j\cA_0(f_0)_{\gamma^+}
			&=\cA_0(f_j)_{\gamma^+}
			\quad\text{in }\cD'\bigl((0,T);L^2(\partial\cD)\bigr),
			\qquad 0\le j\le m,\notag
		\end{align}
		and
		\begin{equation}
			\label{eq:dynamic-outgoing-moment-Ht3}
			\|\cA_0(f_0)_{\gamma^+}\|_{H^m((0,T);L^2(\partial\cD))}
			\lesssim\sum_{j=0}^{m}\|f_j\|_{\cH_{\kin,\sigma}(\cD_T\times\bR_v^3)}.
		\end{equation}
	\end{enumerate}
\end{prop}

\begin{proof}
	\emph{(\ref{it:weak-generator-stationary}).}
	Write $W=q_J$, $q=q_{\nu J}$, $\mathfrak D=\mathfrak D_\nu$,
	$d=d_\nu$ and
	$\cR=[v\cdot\nabla_x,\mathfrak D]$.  For $u\in\cH_{\kin,\sigma}$ set
	\[
	\mathfrak G_u(\Psi):=
	\langle v\cdot\nabla_xu,\Psi\rangle_{x;H_\sigma^{-1},H_\sigma^1}
	+\iint_{\cD\times\bR^3}u\,v\cdot\nabla_x\Psi\,\dd v\dd x.
	\]
	Applying the Green identity \eqref{ibp_st}, which holds on $\cH_{\kin,\sigma}$ by Proposition~\ref{2607151218}, to $u$ and $\Psi$, we obtain
	\begin{equation}
		\label{eq:boundary-generator-green-functional}
		\mathfrak G_u(\Psi)=\int_\gamma u_\gamma\Psi\,(n\cdot v)\,\dd S\dd v.
	\end{equation}
	A direct commutator calculation gives
	\[
	\cR
	=\bigl(((v\cdot\nabla_x)M_\nu)v\bigr)\cdot\nabla_v,
	\qquad
	\operatorname{div}_{x,v}\cR
	=2v\cdot\nabla_x(\operatorname{div}_x\tau_\nu)
	=v\cdot\nabla_xd.
	\]
	In particular $\cR$ is a pure velocity field. The identities
	$$
	q=\mathfrak DW\quad \text{and}\quad v\cdot\nabla_xq=\mathfrak D(v\cdot\nabla_xW)+\cR W
	$$
	hold in $\cD'(\cD\times\bR^3)$.
	Fix $0<r_0\le r_\cD$, with $r_\cD$ the radius of the
	signed-distance collar $\cN_\cD$, and choose
	$\zeta\in C^\infty(\bR)$ with $0\le\zeta\le1$,
	$\zeta=0$ on $(-\infty,1]$, and $\zeta=1$ on $[2,\infty)$.  For
	$0<\varepsilon<r_0/8$, define
	\[
	\chi_\varepsilon(x)=\zeta(-\rho_\cD(x)/\varepsilon)
	\quad\text{when }-r_0/2<\rho_\cD(x)<0,
	\]
	and set $\chi_\varepsilon=1$ where $\rho_\cD\le-r_0/2$.  Then
	$\chi_\varepsilon\in C^{K_\cD+2}(\overline\cD)$,
	\[
	\chi_\varepsilon=0\ \text{if }-\rho_\cD\le\varepsilon,
	\qquad
	\chi_\varepsilon=1\ \text{if }-\rho_\cD\ge2\varepsilon,
	\qquad
	\operatorname{supp}(1-\chi_\varepsilon)\subset
	E_\varepsilon:=\{-2\varepsilon<\rho_\cD<0\},
	\]
	and $|E_\varepsilon|\le C\varepsilon$.  Since $\tau_\nu\cdot\nabla\rho_\cD=0$ on $\cN_\cD$
	and $\cR$ is a pure velocity field, $\mathfrak D\chi_\varepsilon=\cR\chi_\varepsilon=0$.
	Thus no derivative of $\chi_\varepsilon$ is created when $\mathfrak D$ or $\cR$ is
	transposed. Moreover, $\chi_\varepsilon\Psi$ has compact support in $\cD\times\bR^3$.  For every fixed smooth compact-velocity factor $\Phi$,
	\begin{equation*}
		\|\mathbf1_{E_\varepsilon}\Phi\|_{L_x^2H_\sigma^1}\longrightarrow0.
	\end{equation*}
	Hence
	\[
	\begin{aligned}
		\big|\langle v\cdot\nabla_xW,(1-\chi_\varepsilon)\Phi\rangle
		_{x;H_\sigma^{-1},H_\sigma^1}\big|
		&\le
		\|v\cdot\nabla_xW\|_{L_x^2H_\sigma^{-1}}
		\|\mathbf1_{E_\varepsilon}\Phi\|_{L_x^2H_\sigma^1}
		\longrightarrow0,
		\\
		\big|\langle v\cdot\nabla_xq,(1-\chi_\varepsilon)\Phi\rangle
		_{x;H_\sigma^{-1},H_\sigma^1}\big|
		&\le
		\|v\cdot\nabla_xq\|_{L_x^2H_\sigma^{-1}}
		\|\mathbf1_{E_\varepsilon}\Phi\|_{L_x^2H_\sigma^1}
		\longrightarrow0.
	\end{aligned}
	\]
	By the Cauchy--Schwarz inequality, the pointwise bound $|(1-\chi_\varepsilon)\Phi|\le\mathbf1_{E_\varepsilon}|\Phi|$ and the bound $|E_\varepsilon|\le C\varepsilon$, we have
	\[
	\Bigl|\iint W\,(1-\chi_\varepsilon)\Phi\Bigr|+\Bigl|\iint q\,(1-\chi_\varepsilon)\Phi\Bigr|
	\le C\varepsilon^{1/2}\bigl(\|W\|_{L^2_{x,v}}+\|q\|_{L^2_{x,v}}\bigr)
	\longrightarrow0,
	\]
	where $C$ depends on $\Phi$; $\|W\|_{L^2_{x,v}}+\|q\|_{L^2_{x,v}}<\infty$ by \eqref{260902137}.
	Applying this exhaustion to $v\cdot\nabla_xq=\mathfrak D(v\cdot\nabla_xW)+\cR W$ with test
	$\Psi$ and to $q=\mathfrak DW$ with test $v\cdot\nabla_x\Psi$ justifies the two
	transpositions below.
	
	Since $\operatorname{div}_{x,v}\mathfrak f_\nu=2\operatorname{div}_x\tau_\nu=d$,
	transposing
	$\mathfrak D$ and $\cR$ in these identities gives
	\begin{align*}
		\mathfrak G_q(\Psi)
		={}&-\langle v\cdot\nabla_xW,(\mathfrak D+d)\Psi\rangle
		-\iint W\,\bigl(\cR+(v\cdot\nabla_xd)\bigr)\Psi
		-\iint W\,(\mathfrak D+d)(v\cdot\nabla_x\Psi).
	\end{align*}
	Since
	\[
	\cR\Psi+\mathfrak D(v\cdot\nabla_x\Psi)=v\cdot\nabla_x(\mathfrak D\Psi)\quad\text{and}\quad(v\cdot\nabla_xd)\Psi+d\,(v\cdot\nabla_x\Psi)=v\cdot\nabla_x(d\Psi) ,
	\] we obtain
	\begin{equation}
		\label{eq:boundary-generator-green-commutation}
		\mathfrak G_q(\Psi)=-\mathfrak G_W((\mathfrak D+d)\Psi).
	\end{equation}
	Combining \eqref{eq:boundary-generator-green-functional} and
	\eqref{eq:boundary-generator-green-commutation} proves
	\eqref{eq:boundary-generator-identity} on the full signed boundary.
	
	The identity extends to compact-velocity $C^2$ phase tests: every term of $\mathfrak G_q(\Psi)$ and of $\mathfrak G_W((\mathfrak D+d)\Psi)$ in \eqref{eq:boundary-generator-green-commutation} is continuous in its
	$H_\sigma^{-1}$--$H_\sigma^1$ or $L^2$ pairing, and every boundary term is continuous in its pairing against the two	$L^2(\gamma^\pm)$ traces.
	
	By \eqref{eq:signed-distance},
	$\nabla\rho_\cD\in C^{K_\cD+1}(\overline{\cD})$.  Moreover
	$\nabla\rho_\cD=\bar n$ on the signed-distance collar, and
	$\nabla\rho_\cD=n$ on $\partial\cD$.
	For $\theta\in C^\infty(\bR)$ with
	$0\le\theta\le1$, $\theta=0$ on $(-\infty,1/2]$, and
	$\theta=1$ on $[1,\infty)$, put
	\begin{equation}
		\label{eq:boundary-generator-sign-cutoff}
		\theta_{k,\pm}(x,v)=\theta(\pm k\,\nabla\rho_\cD(x)\cdot v).
	\end{equation}
	Because $K_\cD\ge2\mathrm N$, $\theta_{k,\pm}\Psi$ is a compact-velocity $C^2$ phase test.  On the physical boundary,
	\eqref{eq:complete-derivative-preserves-grazing} gives
	$\mathfrak D_\nu\theta_{k,\pm}=0$.  Hence the boundary restriction of
	$(\mathfrak D_\nu+d_\nu)(\theta_{k,\pm}\Psi)$ equals
	$\theta_{k,\pm}(\mathfrak D_\nu+d_\nu)\Psi$.  Since $0\le\theta_{k,\pm}\le1$ and
	$\theta_{k,\pm}\to\mathbf1_{\gamma^\pm}$ almost everywhere with respect to
	$\dd\gamma$ of \eqref{def:dgamma}, dominated convergence and the Cauchy--Schwarz inequality against the
	two $L^2(\gamma^\pm)$ traces allow $k\to\infty$ in
	\eqref{eq:boundary-generator-identity} applied to $\theta_{k,\pm}\Psi$, which proves
	\eqref{eq:boundary-generator-identity} on $\gamma^+$ and on $\gamma^-$.
	
	\emph{(\ref{it:weak-generator-evolutionary}).}
	The proof of (\ref{it:weak-generator-stationary}) applies on $(0,T)\times\cD$ with
	$\partial_t+v\cdot\nabla_x$ in place of $v\cdot\nabla_x$, compact-time tests, an outer $L_t^2$
	on every norm, and \eqref{eq:dynamic-Green-identity} in place of \eqref{ibp_st}. For
	$U\in\cH_{\kin,\sigma}(\cD_T\times\bR_v^3)$ and a compact-time test $\Psi$, \eqref{eq:dynamic-Green-identity} gives
	\[
	\begin{aligned}
		&\int_0^T\!\!\int_\cD
		\bigl(\langle(\partial_t+v\cdot\nabla_x)U,\Psi\rangle_v
		+\langle U,(\partial_t+v\cdot\nabla_x)\Psi\rangle_v\bigr)\dd x\dd t=\int_{\gamma_T}U_\gamma\Psi\operatorname{sgn}(n\cdot v)\dd\gamma\dd t
	\end{aligned}
	\]
	in place of \eqref{eq:boundary-generator-green-functional}.  Since $\mathfrak D$, $d$,
	$\chi_\varepsilon$ and $\theta_{k,\pm}$ do not depend on $t$,
	$[\partial_t+v\cdot\nabla_x,\mathfrak D]=\cR$ and
	$\operatorname{div}_{x,v}\cR=(\partial_t+v\cdot\nabla_x)d$. Hence the two transpositions and the limits $\varepsilon\to0$ and $k\to\infty$ in the proof of (\ref{it:weak-generator-stationary}) are unchanged, and we obtain \eqref{eq:dynamic-boundary-generator} on $\gamma_T^+$ and on
	$\gamma_T^-$.  For $\nu=\mathrm t$, $\mathfrak D=\partial_t$, $d=0$ and $\cR=0$, and the same
	argument gives \eqref{eq:temporal-boundary-generator}.  Applied to the pairs
	$(f_j,f_{j+1})$, \eqref{eq:temporal-boundary-generator} gives the first line of
	\eqref{eq:dynamic-trace-time-chain-boundary}, and applying $\cA_0$ to it, with induction on $j$, gives
	the second.  By
	\eqref{eq:outgoing-moment-flux-bound} and \eqref{260902137},
	\[
	\|\cA_0(f_j)_{\gamma^+}\|_{L^2((0,T)\times\partial\cD)}
	\le C\|f_j\|_{\cH_{\kin,\sigma}(\cD_T\times\bR_v^3)}.
	\]
	With the second line of \eqref{eq:dynamic-trace-time-chain-boundary}, this gives
	\eqref{eq:dynamic-outgoing-moment-Ht3}.
\end{proof}

\begin{lemma}[Time slices of the boundary generator]
	\label{lem:dynamic-generator-slices}
	Fix a sign $\pm$.  Let $\mathcal T\subset\bR$ be an open interval, $\nu\in\mathfrak A$, and
	$r,r^{[\nu]}\in L^2(\mathcal T\times\gamma^\pm)$.
	The following are equivalent.
	\begin{enumerate}
		\item For every compact-time, compact-velocity $C^2$ phase test $\Psi$,
		\begin{equation}
			\label{eq:dynamic-generator-slice-spacetime}
			\int_{\mathcal T}\int_{\gamma^\pm}r^{[\nu]}\,\Psi\,\dd\gamma\dd t
			=-\int_{\mathcal T}\int_{\gamma^\pm}r\,(\mathfrak D_\nu+d_\nu)\Psi\,
			\dd\gamma\dd t .
		\end{equation}
		\item There is a null set $\mathcal N\subset\mathcal T$ such that, for
		$t\in\mathcal T\setminus\mathcal N$, $r(t),r^{[\nu]}(t)\in L^2(\gamma^\pm)$ and
		\begin{equation}
			\label{eq:dynamic-generator-slice-identity}
			\int_{\gamma^\pm}r^{[\nu]}(t)\,\Psi\,\dd\gamma
			=-\int_{\gamma^\pm}r(t)\,(\mathfrak D_\nu+d_\nu)\Psi\,\dd\gamma
		\end{equation}
		for every compact-velocity $C^2$ phase test $\Psi$.
	\end{enumerate}
	For finitely many pairs $(r,r^{[\nu]})$, the null set $\mathcal N$ in (2) can be chosen common to all of them.  In particular, if $z_J,z_{\nu J}$ are as in
	Proposition~\ref{prop:weak-boundary-generator}\,(\ref{it:weak-generator-evolutionary}), then at almost every
	$t\in(0,T)$ the traces $(z_J)_{\gamma^\pm}(t)$, $(z_{\nu J})_{\gamma^\pm}(t)$
	satisfy \eqref{eq:boundary-generator-identity} on $\gamma^+$ and on $\gamma^-$
	separately.
\end{lemma}

\begin{proof}
	The identity \eqref{eq:dynamic-generator-slice-spacetime} is the $\gamma_T^\pm$
	part of \eqref{eq:dynamic-boundary-generator}, and
	\eqref{eq:dynamic-generator-slice-identity} is the $\gamma^\pm$ part of
	\eqref{eq:boundary-generator-identity}.
	
	Fix $R\in\bN$ and let $\mathcal C_R$ be the set of $C^2$ phase tests with
	velocity support in $\{|v|\le R\}$.  The field $\mathfrak D_\nu$ is tangent to
	$\partial\cD\times\bR^3$, and its coefficients and $d_\nu$ are bounded on
	$\partial\cD\times\{|v|\le R\}$, so
	\[
	|(\mathfrak D_\nu+d_\nu)\Psi|\le C\|\Psi\|_{C^1(\partial\cD\times\{|v|\le R\})}
	\quad\text{on }\partial\cD\times\bR^3,\qquad\Psi\in\mathcal C_R .
	\]
	The Cauchy--Schwarz inequality gives,
	for $\Psi\in\mathcal C_R$ and $q,q'\in L^2(\gamma^\pm)$,
	\begin{equation}
		\label{eq:dynamic-generator-slice-continuity}
		\Bigl|\int_{\gamma^\pm}\bigl(q'\Psi+q(\mathfrak D_\nu+d_\nu)\Psi\bigr)\dd\gamma\Bigr|
		\le C\,\|\Psi\|_{C^1(\partial\cD\times\{|v|\le R\})}
		\bigl(\|q\|_{L^2(\gamma^\pm)}+\|q'\|_{L^2(\gamma^\pm)}\bigr),
	\end{equation}
	where $C=C(R)$.
	
	(2)$\Rightarrow$(1).  Let $\Psi$ be a compact-time test with velocity
	support in $\{|v|\le R\}$.  Then $\Psi(t,\cdot)\in\mathcal C_R$ for every $t$,
	so \eqref{eq:dynamic-generator-slice-identity} holds with $\Psi(t,\cdot)$ for $t\notin\mathcal N$.  By
	\eqref{eq:dynamic-generator-slice-continuity} both sides are integrable in $t$,
	and integration over $\mathcal T$ gives
	\eqref{eq:dynamic-generator-slice-spacetime}.
	
	(1)$\Rightarrow$(2).  For $\Psi\in\mathcal C_R$ set
	\[
	G_\Psi(t):=\int_{\gamma^\pm}\bigl(r^{[\nu]}(t)\Psi+r(t)(\mathfrak D_\nu+d_\nu)\Psi\bigr)\dd\gamma .
	\]
	By \eqref{eq:dynamic-generator-slice-continuity}, $G_\Psi\in L^2(\mathcal T)$.
	For $\alpha\in C_c^\infty(\mathcal T)$, the product $\alpha(t)\Psi(x,v)$ is a test in
	(1), and Fubini's theorem gives $\int_{\mathcal T}\alpha\,G_\Psi\dd t=0$.  Hence
	$G_\Psi=0$ outside a null set $\mathcal N_\Psi$.  The restrictions of
	$\mathcal C_R$ to $\partial\cD\times\{|v|\le R\}$ form a subset of the
	separable space $C^1(\partial\cD\times\{|v|\le R\})$, so there is a
	countable family $\{\Psi_{R,n}\}_{n}\subset\mathcal C_R$ that is dense in
	$\mathcal C_R$ for the $C^1$ norm.  Let $\mathcal N$ be the union of the
	$\mathcal N_{\Psi_{R,n}}$, $R,n\in\bN$, and of the null set where
	$r(t)\notin L^2(\gamma^\pm)$ or $r^{[\nu]}(t)\notin L^2(\gamma^\pm)$.  For
	$t\notin\mathcal N$ and $\Psi\in\mathcal C_R$,
	choose $n_k$ with $\Psi_{R,n_k}\to\Psi$ in $C^1(\partial\cD\times\{|v|\le R\})$; then
	\eqref{eq:dynamic-generator-slice-continuity} applied to $\Psi-\Psi_{R,n_k}$
	gives
	\[
	G_\Psi(t)=\lim_{k\to\infty}G_{\Psi_{R,n_k}}(t)=0,
	\]
	which proves (2).
	
	For the last assertion, apply (1)$\Rightarrow$(2) on $\mathcal T=(0,T)$ with $r=(z_J)_{\gamma^\pm}$ and $r^{[\nu]}=(z_{\nu J})_{\gamma^\pm}$, separately for the two signs. These lie in $L^2(\gamma_T^\pm)$ by \eqref{260902137}, and for them (1) is \eqref{eq:dynamic-boundary-generator} on $\gamma_T^\pm$.
\end{proof}

With $J_N$ and $P_{w,N}=\mathfrak m_{\Theta_w}J_N\cA_0$ of
\eqref{eq:boundary-spectral-diffuse-operator}, set, for $0\le j\le\mathrm N$ and a function $f$ on $(0,T)\times\cD\times\bR^3$,
\begin{equation}
	\label{eq:dyn-theta-triangular-notation}
	\mathfrak t_{j,N}[f]
	:=\sum_{k=1}^{j}\binom jk
	\bigl(\partial_t^kd_\delta\bigr)\,J_N\cA_0(\partial_t^{j-k}f)_{\gamma^+} ,
\end{equation}
the lower-order temporal contribution to the truncated wall law.  By \eqref{eq:wall-response-profiles}, 
$$
P_{w,N}(t)=\cP_{0,N}^-+d_\delta(t)J_N\cA_0,\qquad \cP_{0,N}^-:=\mathfrak m_\circ J_N\cA_0.
$$

\begin{lemma}[Boundary commutators]
	\label{lem:ordered-boundary-finite-rank}\,
	
	\begin{enumerate}
		\item\label{it:finite-rank-stationary}
		Every $P_I$, $0\le|I|\le\mathrm N$, extends from the smooth trace core to a
		finite-rank operator in the velocity fiber.  On every boundary chart, for
		$q\in L^2(\gamma^+)$,
		\begin{equation}
			\label{eq:ordered-boundary-finite-rank}
			P_Iq(x,v)=\sum_{j=1}^{N_I}m_{I,j}(x,v)\cA_{p_{I,j},\mu}q(x),
		\end{equation}
		where $m_{I,j}$ and $p_{I,j}\sqrt\mu$ satisfy
		\eqref{eq:boundary-gaussian-polynomial-envelope}.  Consequently, for every finite $\ell\ge0$,
		\begin{equation}
			\label{eq:ordered-boundary-finite-rank-bound}
			\|w^{\ell}P_Iq\|_{L^2(\gamma^-)}
			\le C\|q\|_{L^2(\gamma^+)},
		\end{equation}
		where $C$ is independent of $\delta$.  On every boundary chart, with $m_j:=m_{I,j}$ and $p_j:=p_{I,j}$ of \eqref{eq:ordered-boundary-finite-rank},
		\begin{equation}
			\label{eq:ordered-boundary-commutator-induction}
			P_{\nu I}q
			=\sum_{j=1}^{N_I}\left\{
			(\mathfrak D_\nu^- m_j)\cA_{p_j,\mu}q
			+m_j\cA_{\mathscr E_{\nu,\mu}p_j,\mu}q\right\}.
		\end{equation}
		
		\item\label{it:finite-rank-evolutionary}
		Let $0<T<\infty$.  The temporal commutators \eqref{eq:dyn-mixed-boundary-commutator} of the
		truncated wall are
		\begin{equation}
			\label{eq:dyn-mixed-temporal-commutator}
			P_{\mathrm t^k,N}=\bigl(\partial_t^kd_\delta\bigr)J_N\cA_0,
			\qquad 1\le k\le\mathrm N.
		\end{equation}
		For $0\le j\le\mathrm N$ and every finite $\ell\ge0$,
		\begin{equation}
			\label{eq:dyn-theta-triangular-bound}
			\|w^{\ell}\mathfrak t_{j,N}[f]\|_{L^2(\gamma_T^-)}
			\le C\,\delta\sum_{i<j}\|(\partial_t^if)_{\gamma^+}\|_{L^2(\gamma_T^+)},
		\end{equation}
		where $C$ is independent of $T$, $N$ and $\delta$.
	\end{enumerate}
\end{lemma}

\begin{proof}
	\emph{(\ref{it:finite-rank-stationary}).}
	Suppose that on one boundary chart $P_I$ has the form
	\eqref{eq:ordered-boundary-finite-rank}, and write $m_j:=m_{I,j}$,
	$p_j:=p_{I,j}$.  By \eqref{eq:ordered-boundary-commutator-definition} and Lemma~\ref{lem:boundary-word-Leibniz} for one letter, with $a=m_j$ and $b=\cA_{p_j,\mu}q$,
	\[
	P_{\nu I}q=\sum_{j=1}^{N_I}\bigl\{(\mathfrak D_\nu^-m_j)\cA_{p_j,\mu}q
	+m_j\bigl(\tau_\nu\cA_{p_j,\mu}q-\cA_{p_j,\mu}(\mathfrak D_\nu^+q)\bigr)\bigr\}.
	\]
	By \eqref{eq:moment-commutator}, the last bracket is $\cA_{\mathscr E_{\nu,\mu}p_j,\mu}q$, and we obtain
	\eqref{eq:ordered-boundary-commutator-induction}.  Since
	$P_\varnothing=P_w=\mathfrak m_{\Theta_w}\cA_0$ and $\cA_0=\cA_{1,\mu}$, induction proves
	\eqref{eq:ordered-boundary-finite-rank} through order $\mathrm N$.  By
	\eqref{eq:ordered-boundary-commutator-induction}, every $p_j$ is one
	of the ordered profiles $p_K$ generated by
	\eqref{eq:moment-profile-recursion}, with $|K|\le\mathrm N$. Hence
	Lemma~\ref{lem:boundary-profile-bounds} gives
	\eqref{eq:boundary-gaussian-polynomial-envelope} for $p_j\sqrt\mu$.
	Every $m_j$ is an ordered complete derivative of
	$\mathfrak m_{\Theta_w}=\mathfrak m_\circ+d_\delta$ of order at most $\mathrm N$, so
	Lemma~\ref{lem:wall-profile-derivatives} and the fixed Gaussian profile
	$\mathfrak m_\circ$ give \eqref{eq:boundary-gaussian-polynomial-envelope} for $m_j$.  For
	$\ell\ge0$ and a rank-one term $m\cA_{p,\mu}q$, the fiberwise Cauchy--Schwarz inequality and integration
	over $x\in\partial\cD$ give
	\begin{equation}
		\label{eq:ordered-boundary-rank-one-bound}
		\|w^{\ell}m\cA_{p,\mu}q\|_{L^2(\gamma^-)}^2
		\le
		\left(\sup_{x\in\partial\cD}\int_{n_x\cdot v<0}
		w^{2\ell}|m|^2|n_x\cdot v|\,\dd v\right)
		\left(\sup_{x\in\partial\cD}\int_{n_x\cdot v>0}
		|p|^2\mu\,n_x\cdot v\,\dd v\right)
		\|q\|_{L^2(\gamma^+)}^2.
	\end{equation}
	By \eqref{eq:boundary-profile-flux-bound} with $K=\varnothing$, at $a=m$ with $L=\ell$ and at $a=p\sqrt\mu$ with $L=0$, both suprema are finite for
	$(m,p)=(m_j,p_j)$ and bounded independently of $\delta$, since neither the $p_j$ nor the constants of Lemma~\ref{lem:wall-profile-derivatives} depend on $\delta$.
	Summing \eqref{eq:ordered-boundary-rank-one-bound} over
	\eqref{eq:ordered-boundary-finite-rank} proves
	\eqref{eq:ordered-boundary-finite-rank-bound} and defines the bounded
	extension of $P_I$ on $L^2(\gamma^+)$.
	
	\emph{(\ref{it:finite-rank-evolutionary}).}
	The decomposition $P_{w,N}(t)=\cP_{0,N}^-+d_\delta(t)J_N\cA_0$, with time-independent $\cP_{0,N}^-$ and $J_N\cA_0$, gives
	\eqref{eq:dyn-mixed-temporal-commutator}.  By \eqref{eq:dyn-theta-defect-time-bound}, integration in $v$ gives
	\[
	\|w^{\ell}(\partial_t^kd_\delta)h\|_{L^2(\gamma^-)}\le C\delta\|h\|_{L^2(\partial\cD)},
	\qquad h\in L^2(\partial\cD),
	\]
	uniformly in $t$. With $h=J_N\cA_0(\partial_t^{j-k}f)_{\gamma^+}$, the bound
	$\|J_N\|\le1$ on $L^2(\partial\cD)$ and \eqref{eq:outgoing-moment-flux-bound} give
	\eqref{eq:dyn-theta-triangular-bound}.
\end{proof}

Assume that
$$
f_{\gamma^-}=P_{w,N}(t)f_{\gamma^+}+b\quad \text{on }(0,T)\qquad\text{and}\qquad (\partial_t^if)_{\gamma^\pm}=\partial_t^if_{\gamma^\pm}\quad \text{in }\cD'\bigl((0,T);L^2(\gamma^\pm)\bigr)\quad \text{for }i\le j.
$$
Then the Leibniz rule for the product $d_\delta(t)J_N\cA_0f_{\gamma^+}$,
\[
\partial_t^j\bigl(d_\delta J_N\cA_0f_{\gamma^+}\bigr)
=d_\delta J_N\cA_0(\partial_t^jf)_{\gamma^+}+\mathfrak t_{j,N}[f],
\]
and the $t$-independence of $\cP_{0,N}^-$ give the truncated temporal law
\begin{equation}
	\label{eq:dyn-temporal-boundary-law}
	(\partial_t^jf)_{\gamma^-}=P_{w,N}(t)(\partial_t^jf)_{\gamma^+}+\mathfrak t_{j,N}[f]+\partial_t^jb .
\end{equation}
We suppress the time argument in $P_w(t)$ and $P_{w,N}(t)$ henceforth.

\subsubsection{Outgoing moments}\label{ss:common-boundary-moments}

For $H\in L^2(\partial\cD)$ the tangential derivative $\tau_\nu H$
is the boundary distribution
\begin{equation}
	\label{eq:surface-tangential-adjoint}
	\langle\tau_\nu H,\eta\rangle
	:=-\int_{\partial\cD}H\bigl(\tau_\nu\eta
	+(\operatorname{div}_{\partial\cD}\tau_\nu)\eta\bigr)\,\dd S,
	\qquad\eta\in C^{K_\cD+1}(\partial\cD).
\end{equation}

\begin{lemma}[Outgoing Gaussian moments]
	\label{lem:rough-outgoing-moments}
	Let $a$ be a boundary profile with the Gaussian envelope
	\eqref{eq:boundary-gaussian-polynomial-envelope}.  Then, for every
	$q\in L^2(\gamma^+)$,
	\begin{equation}
		\label{eq:word-profile-moment-bound}
		\|\cA_aq\|_{L^2(\partial\cD)}
		\le\Bigl(\sup_{x\in\partial\cD}
		\int_{n_x\cdot v>0}|a(x,v)|^2(n_x\cdot v)\,\dd v\Bigr)^{1/2}
		\|q\|_{L^2(\gamma^+)}
		\lesssim\|q\|_{L^2(\gamma^+)}.
	\end{equation}
	\begin{enumerate}
		\item\label{it:rough-moment-stationary}
		Let $\nu\in\mathfrak A$ and let $q,q^{[\nu]}\in L^2(\gamma^+)$ satisfy
		\begin{equation}
			\label{eq:rough-moment-hypothesis}
			\int_{\gamma^+}q^{[\nu]}\Psi\,(n\cdot v)\,\dd S\dd v
			=-\int_{\gamma^+}q\,(\mathfrak D_\nu+d_\nu)\Psi\,(n\cdot v)\,\dd S\dd v
		\end{equation}
		for every compact-velocity $C^2$ phase test $\Psi$ up to the physical
		boundary.  Then, with $\tau_\nu$ acting as in
		\eqref{eq:surface-tangential-adjoint},
		\begin{equation}
			\label{eq:word-rough-moment-one}
			\tau_\nu\cA_aq
			=\cA_aq^{[\nu]}+\cA_{\mathscr E_\nu a}q
			\qquad\text{in }\cD'(\partial\cD).
		\end{equation}
		Both terms on the right belong to $L^2(\partial\cD)$.
		\item\label{it:rough-moment-evolutionary}
		Let $0<T<\infty$, let $\nu\in\mathfrak A_\bullet$, and let
		$z,z^{[\nu]}\in\cH_{\kin,\sigma}(\cD_T\times\bR_v^3)$ satisfy
		$\mathfrak D_\nu z=z^{[\nu]}$ in $\cD'((0,T)\times\cD\times\bR^3)$.  Then
		\begin{equation}
			\label{eq:dynamic-one-letter-moment}
			\tau_\nu\cA_az_{\gamma^+}
			=\cA_az^{[\nu]}_{\gamma^+}
			+\cA_{\mathscr E_\nu a}z_{\gamma^+}
			\quad\text{in }\cD'((0,T)\times\partial\cD).
		\end{equation}
	\end{enumerate}
\end{lemma}

\begin{proof}
	The moment bound follows from fiberwise Cauchy--Schwarz and \eqref{eq:boundary-profile-flux-bound} with $K=\varnothing$ and $L=0$.
	
	\emph{(\ref{it:rough-moment-stationary}).}
	Fix $\eta\in C^{K_\cD+1}(\partial\cD)$ and choose a $C^{K_\cD+1}$ collar extension of
	$\eta$ and a $C^{\mathrm N}$ collar extension of $a$.  Let $\chi_R$ be radial, equal to one on $|v|\le R$, zero on
	$|v|\ge2R$, and satisfy $|\nabla_v\chi_R|\le C/R$.  Since
	\eqref{eq:complete-phase-space-derivative} gives
	\[
	\mathfrak D_\nu\chi_R
	=((\nabla\tau_\nu)v)\cdot\nabla_v\chi_R,
	\qquad
	|\mathfrak D_\nu\chi_R|
	\le C\mathbf 1_{\{R\le|v|\le2R\}},
	\]
	and, with \eqref{eq:boundary-gaussian-polynomial-envelope}, $|a|+|\mathfrak D_\nu a+d_\nu a|\le C\langle v\rangle\ee^{-c_a|v|^2}$, we have
	\begin{equation}\label{261016903}
		\begin{aligned}
			&\sup_{x\in\partial\cD}\int_{n_x\cdot v>0}
			\Bigl((1-\chi_R)^2\bigl(|a|^2+|\mathfrak D_\nu a+d_\nu a|^2\bigr)
			+|a|^2|\mathfrak D_\nu\chi_R|^2\Bigr)
			(n_x\cdot v)\,\dd v
			\\
			&\qquad \le C\int_R^\infty\langle s\rangle^2\ee^{-cs^2}s^3\,\dd s
			\longrightarrow0.
		\end{aligned}
	\end{equation}
	We insert the compact-velocity phase test $\Psi_R=\eta\chi_Ra$ into
	\eqref{eq:rough-moment-hypothesis}. The first inequality in
	\eqref{eq:word-profile-moment-bound} does not use \eqref{eq:boundary-gaussian-polynomial-envelope}, so it holds with any measurable $b$ on $\gamma^+$ in place of $a$. By this inequality and \eqref{261016903}, for $b=(1-\chi_R)a$, $(1-\chi_R)(\mathfrak D_\nu a+d_\nu a)$, and $a\,\mathfrak D_\nu\chi_R$,
	\begin{align*}
		&\|\cA_bq\|_{L^2(\partial\cD)}+\|\cA_bq^{[\nu]}\|_{L^2(\partial\cD)}
		\\
		&\qquad\le\Bigl(\sup_{x\in\partial\cD}\int_{n_x\cdot v>0}|b(x,v)|^2(n_x\cdot v)\,\dd v\Bigr)^{1/2}
		\bigl(\|q\|_{L^2(\gamma^+)}+\|q^{[\nu]}\|_{L^2(\gamma^+)}\bigr)
		\longrightarrow0
		\qquad(R\to\infty).
	\end{align*}
	Hence we may pass to the limit
	$R\to\infty$.  Since
	\[
	(\mathfrak D_\nu+d_\nu)(\eta\chi_Ra)
	=(\tau_\nu\eta)\chi_Ra
	+\eta\chi_R(\mathfrak D_\nu a+d_\nu a)
	+\eta a\mathfrak D_\nu\chi_R,
	\]
	we obtain
	\[
	\int_{\partial\cD}\eta\,\cA_aq^{[\nu]}\,\dd S
	=-\int_{\partial\cD}(\tau_\nu\eta)\cA_aq\,\dd S
	-\int_{\partial\cD}\eta\,\cA_{\mathfrak D_\nu a+d_\nu a}q\,\dd S.
	\]
	Therefore, by \eqref{eq:surface-tangential-adjoint},
	\begin{align*}
		\langle\tau_\nu\cA_aq,\eta\rangle
		={}&\int_{\partial\cD}\eta\,\cA_aq^{[\nu]}\,\dd S+\int_{\partial\cD}\eta\,
		\cA_{\mathfrak D_\nu a+
			(d_\nu-\operatorname{div}_{\partial\cD}\tau_\nu)a}q\,\dd S.
	\end{align*}
	Since \eqref{eq:boundary-phase-divergence} gives
	$d_\nu-\operatorname{div}_{\partial\cD}\tau_\nu
	=\operatorname{div}_x\tau_\nu$, the definition
	\eqref{eq:outgoing-profile-generator} yields
	\eqref{eq:word-rough-moment-one}.  By
	\eqref{eq:boundary-gaussian-polynomial-envelope} and
	\eqref{eq:outgoing-profile-generator},
	\[
	|\mathscr E_\nu a(x,v)|
	\le C w\,\ee^{-c_a|v|^2}.
	\]
	The fiberwise Cauchy--Schwarz estimate used in
	\eqref{eq:word-profile-moment-bound} therefore shows that both terms on the right of \eqref{eq:word-rough-moment-one} belong to $L^2(\partial\cD)$.
	
	\emph{(\ref{it:rough-moment-evolutionary}).}
	For $\nu\in\mathfrak A$, the proof of (\ref{it:rough-moment-stationary}) applies with
	$\eta\in C_c^{K_\cD+1}((0,T)\times\partial\cD)$ in place of $\eta\in C^{K_\cD+1}(\partial\cD)$ and with the $\gamma_T^+$ part of \eqref{eq:dynamic-boundary-generator}, which holds by
	Proposition~\ref{prop:weak-boundary-generator}\,(\ref{it:weak-generator-evolutionary}), in place of \eqref{eq:rough-moment-hypothesis}.  This gives \eqref{eq:dynamic-one-letter-moment}.  At $\nu=\mathrm t$, where $\tau_{\mathrm t}=\partial_t$ and $\mathscr E_{\mathrm t}=0$,
	\eqref{eq:dynamic-one-letter-moment} reads $\partial_t\cA_az_{\gamma^+}=\cA_az^{[\mathrm t]}_{\gamma^+}$. This follows by applying $\cA_a$ to \eqref{eq:temporal-boundary-generator} for the pair $(z,z^{[\mathrm t]})$.
\end{proof}

For $f\in\cH_{\kin,\sigma}(\cD_T\times\bR_v^3)$, write
\begin{equation}
	\label{eq:dynamic-outgoing-moment-definition}
	\phi[f](t,x):=\cA_0f_{\gamma^+}(t,x).
\end{equation}

\begin{lemma}[Tangential Sobolev regularity]
	\label{lem:tangential-word-sobolev}
	Let $0\le r\le\mathrm N$ and $H\in L^2(\partial\cD)$, and suppose that the distribution $\tau_IH$, with
	$\tau_I$ of \eqref{eq:ordered-complete-derivatives}, belongs to $L^2(\partial\cD)$
	for every ordered word $I$ over $\mathfrak A$ with $1\le|I|\le r$.  Then
	\begin{equation}
		\label{eq:tangential-word-sobolev}
		\|H\|_{H^r(\partial\cD)}
		\lesssim \|H\|_{L^2(\partial\cD)}
		+\sum_{1\le |I|\le r}\|\tau_IH\|_{L^2(\partial\cD)}.
	\end{equation}
	On $(0,T)\times\partial\cD$, \eqref{eq:tangential-word-sobolev} holds for ordered words over
	$\mathfrak A_\bullet$, with $\tau_{\mathrm t}=\partial_t$ and with $H^r((0,T)\times\partial\cD)$ and
	$L^2((0,T)\times\partial\cD)$ in place of $H^r(\partial\cD)$ and $L^2(\partial\cD)$; the
	constant is independent of $T$.
\end{lemma}

\begin{proof}
	By \eqref{eq:fixed-atlas-core-cover} and
	\eqref{eq:interior-affine-chart},
	$\{\widetilde U_\kappa\cap\partial\cD\}_{\kappa\in\mathfrak K_{\p}}$
	covers $\partial\cD$.  We choose
	$\varrho_\kappa\in C_c^{K_\cD+1}(\widetilde U_\kappa\cap\partial\cD)$ with $\sum_{\kappa\in\mathfrak K_{\p}}\varrho_\kappa=1$ on $\partial\cD$ and set
	\[
	\iota_\kappa(y'):=\phi_\kappa(y',0),
	\qquad
	H_\kappa:=H\circ\iota_\kappa,
	\qquad
	\varrho_\kappa^\flat:=\varrho_\kappa\circ\iota_\kappa.
	\]
	By \eqref{eq:coordinate-tangent-fields} and
	\eqref{eq:signed-distance-tangent-extension},
	\[
	\partial_{a_1}\cdots\partial_{a_j}H_\kappa
	=\bigl(\tau_{\kappa,a_1}\cdots\tau_{\kappa,a_j}H\bigr)
	\circ\iota_\kappa
	\quad\text{in }
	\cD'\!\left(\iota_\kappa^{-1}
	(\widetilde U_\kappa\cap\partial\cD)\right),
	\qquad 1\le j\le r.
	\]
	The right side of this identity is $(\tau_IH)\circ\iota_\kappa$ for the word $I=((\kappa,a_1),\ldots,(\kappa,a_j))$ over $\mathfrak A$.  Extend
	$\varrho_\kappa^\flat H_\kappa$ by zero to $\bR^2$.
	Since the surface Jacobian of $\iota_\kappa$ and its inverse are bounded on $\operatorname{supp}\varrho_\kappa^\flat$ by
	\eqref{eq:flattening-matrix-formulas}, the preceding identity and the distributional Leibniz rule give
	\begin{align*}
		\|\varrho_\kappa^\flat H_\kappa\|_{H^r(\bR^2)}
		&\lesssim
		\|H\|_{L^2(\widetilde U_\kappa\cap\partial\cD)}
		+\sum_{1\le |I|\le r}
		\|\tau_IH\|_{L^2(\widetilde U_\kappa\cap\partial\cD)}.
	\end{align*}
	Summing the chart estimates with the partition of unity yields
	\eqref{eq:tangential-word-sobolev}.  On $(0,T)\times\partial\cD$ the same argument applies on
	$(0,T)\times\bR^2$, with $\partial_t=\tau_{\mathrm t}$ added to the coordinate derivatives and
	the partition of unity acting in $x$ alone, so no constant depends on $T$.
\end{proof}

\begin{lemma}[Ordered outgoing moments]
	\label{lem:dynamic-boundary-moments}
	\label{cor:outgoing-moment-H3}
	Let $1\le m\le\mathrm N$.  For an ordered word $K=(\nu_1,\ldots,\nu_r)$ over $\mathfrak A_\bullet$ set
	\begin{equation}
		\label{eq:ordered-profile-words}
		\mathscr E_K:=\mathscr E_{\nu_1}\cdots\mathscr E_{\nu_r},
		\qquad
		\mathscr E_\varnothing:=\operatorname{Id},
	\end{equation}
	and let
	\begin{equation}
		\label{eq:outgoing-moment-static-definition}
		\phi:=\cA_0(f_\varnothing)_{\gamma^+}=\cA_{\sqrt\mu}(f_\varnothing)_{\gamma^+} ,
	\end{equation}
	which on $(0,T)\times\partial\cD$ is $\phi[f_\varnothing]$ of
	\eqref{eq:dynamic-outgoing-moment-definition}; here $f_\varnothing$ is the member of the family $(f_J)$ below with $J=\varnothing$.
	\begin{enumerate}
		\item\label{it:outgoing-moment-stationary}\label{it:outgoing-moment-H-stationary}
		Let $f_J\in\cH_{\kin,\sigma}(\cD\times\bR_v^3)$,
		$J\in\bigcup_{k\le m}\mathfrak A^k$, satisfy
		$\mathfrak D_\nu f_J=f_{\nu J}$ in $\cD'(\cD\times\bR^3)$ for
		$|\nu J|\le m$.  Then, for every $I=(\nu_1,\ldots,\nu_k)\in\mathfrak A^k$, $1\le k\le m$,
		\begin{equation}
			\label{eq:outgoing-moment-word}
			\tau_I\phi
			=\sum_{S\subset[k]}\cA_{\mathscr E_{I_S}\sqrt\mu}
			(f_{I_{S^c}})_{\gamma^+}
			\qquad\text{in }\cD'(\partial\cD),
		\end{equation}
		\begin{equation}
			\label{eq:outgoing-moment-word-L2}
			\|\tau_I\phi\|_{L^2(\partial\cD)}
			\lesssim\sum_{|J|\le k}\|(f_J)_{\gamma^+}\|_{L^2(\gamma^+)}
			\lesssim\sum_{|J|\le k}\|f_J\|_{\cH_{\kin,\sigma}(\cD\times\bR_v^3)},
		\end{equation}
		and consequently
		\begin{equation}
			\label{eq:stationary-outgoing-moment-H3}
			\|\phi\|_{H^{m}(\partial\cD)}
			\lesssim\sum_{|J|\le m}
			\|f_J\|_{\cH_{\kin,\sigma}(\cD\times\bR_v^3)}.
		\end{equation}
		\item\label{it:outgoing-moment-evolutionary}\label{it:outgoing-moment-H-evolutionary}
		Let $0<T<\infty$ and let $f_J\in\cH_{\kin,\sigma}(\cD_T\times\bR_v^3)$,
		$J\in\bigcup_{k\le m}\mathfrak A_\bullet^k$, satisfy
		$\mathfrak D_\nu f_J=f_{\nu J}$ in $\cD'((0,T)\times\cD\times\bR^3)$ for
		$|\nu J|\le m$.  Then \eqref{eq:outgoing-moment-word} and
		\eqref{eq:outgoing-moment-word-L2} hold for words $I$ over $\mathfrak A_\bullet$, in
		$\cD'((0,T)\times\partial\cD)$ and with $L^2((0,T)\times\partial\cD)$, $L^2(\gamma_T^+)$ and
		$\cH_{\kin,\sigma}(\cD_T\times\bR_v^3)$ in place of $L^2(\partial\cD)$, $L^2(\gamma^+)$ and
		$\cH_{\kin,\sigma}(\cD\times\bR_v^3)$, and
		\begin{equation}
			\label{eq:dynamic-outgoing-moment-H3}
			\|\phi\|_{H^{m}((0,T)\times\partial\cD)}
			\lesssim\sum_{|J|\le m}
			\|f_J\|_{\cH_{\kin,\sigma}(\cD_T\times\bR_v^3)}.
		\end{equation}
	\end{enumerate}
	The constants are independent of $T$.  If $\mathfrak D_Jf\in\cH_{\kin,\sigma}(\cD\times\bR_v^3)$ for $J\in\bigcup_{k\le m}\mathfrak A^k$, then the functions $f_J:=\mathfrak D_Jf$ satisfy the hypotheses of (\ref{it:outgoing-moment-stationary}) and $\phi=\cA_0f_{\gamma^+}$; the same holds for (\ref{it:outgoing-moment-evolutionary}) with $\mathfrak A_\bullet$ and $\cH_{\kin,\sigma}(\cD_T\times\bR_v^3)$ in place of $\mathfrak A$ and $\cH_{\kin,\sigma}(\cD\times\bR_v^3)$.
\end{lemma}

\begin{proof}
	\emph{(\ref{it:outgoing-moment-stationary}).}
	By $\mathscr E_\nu(p\sqrt\mu)=(\mathscr E_{\nu,\mu}p)\sqrt\mu$ and
	\eqref{eq:moment-profile-recursion}, $\mathscr E_{I_S}\sqrt\mu=p_{I_S}\sqrt\mu$, which
	satisfies \eqref{eq:boundary-gaussian-polynomial-envelope} by
	Lemma~\ref{lem:boundary-profile-bounds}.   For $|\nu J|\le m$,  the hypothesis of Lemma~\ref{lem:rough-outgoing-moments}\,(\ref{it:rough-moment-stationary}) holds for $q=(f_J)_{\gamma^+}$, $q^{[\nu]}=(f_{\nu J})_{\gamma^+}$ by Proposition~\ref{prop:weak-boundary-generator}\,(\ref{it:weak-generator-stationary}), since $\mathfrak D_\nu f_J=f_{\nu J}$. Hence, for every profile $a$ satisfying \eqref{eq:boundary-gaussian-polynomial-envelope}, Lemma~\ref{lem:rough-outgoing-moments} gives
	\[
	\tau_\nu\cA_a(f_J)_{\gamma^+}
	=\cA_a(f_{\nu J})_{\gamma^+}
	+\cA_{\mathscr E_\nu a}(f_J)_{\gamma^+}
	\quad\text{in }\cD'(\partial\cD).
	\]
	Taking $a=\sqrt\mu$ and $J=\varnothing$ gives the case $k=1$. For $k\ge2$, apply the same identity to the induction hypothesis: the preceding identity with $\nu=\nu_1$, applied to each term $\cA_{\mathscr E_{I'_{S'}}\sqrt\mu}(f_{I'_{S'^c}})_{\gamma^+}$ of \eqref{eq:outgoing-moment-word} for $I':=(\nu_2,\ldots,\nu_k)$, gives the terms with $1\notin S$ and with $1\in S$; this is \eqref{eq:outgoing-moment-word}.
	Applying \eqref{eq:word-profile-moment-bound}  to the $2^k$ terms and using the stationary line of \eqref{260902137}, we obtain \eqref{eq:outgoing-moment-word-L2}.  By \eqref{eq:outgoing-moment-word-L2}, Lemma~\ref{lem:tangential-word-sobolev} applies to $H=\phi$ with $r=m$; with
	\eqref{eq:outgoing-moment-flux-bound} and the stationary line of \eqref{260902137} for $\|\phi\|_{L^2(\partial\cD)}$, it gives \eqref{eq:stationary-outgoing-moment-H3}.
	
	\emph{(\ref{it:outgoing-moment-evolutionary}).}
	For every letter $\nu\in\mathfrak A_\bullet$ and every jet $f_J$, Proposition~\ref{prop:weak-boundary-generator} and Lemma~\ref{lem:rough-outgoing-moments} give, in distributions on $(0,T)\times\partial\cD$,
	\[
	\tau_\nu\cA_a(f_J)_{\gamma^+}
	=\cA_a(f_{\nu J})_{\gamma^+}
	+\cA_{\mathscr E_\nu a}(f_J)_{\gamma^+}.
	\]
	Here $a$ ranges over the profiles $\mathscr E_K\sqrt\mu$; their Gaussian envelopes are uniform on the fixed atlas. For the temporal letter, $\tau_{\mathrm t}=\partial_t$ and $\mathscr E_{\mathrm t}=0$. Iterating this distributional identity gives \eqref{eq:outgoing-moment-word}. The Gaussian moment bound and the outgoing graph trace then yield
	\[
	\|\tau_I\phi\|_{L^2((0,T)\times\partial\cD)}
	\le C_m\sum_{|J|\le |I|}\|f_J\|_{\cH_{\kin,\sigma}(\cD_T\times\bR_v^3)}.
	\]
	The chartwise characterization of $H^m((0,T)\times\partial\cD)$ in Lemma~\ref{lem:tangential-word-sobolev}, including the zeroth-order moment bound, proves \eqref{eq:dynamic-outgoing-moment-H3}. All identities are distributional; no pointwise boundary differentiability of the jets is assumed. The trace, Gaussian-profile and chart constants are independent of $T$.
\end{proof}

\subsubsection{The spectral cutoff}\label{app:dyn-boundary-interpolation}

Let $J_N$ be as in \eqref{eq:boundary-spectral-cutoff-definition}, and let $\tau_\nu$ act on boundary
scalars as in \eqref{eq:surface-tangential-adjoint}.  For ordered words $K$ over $\mathfrak A$ and $\nu\in\mathfrak A$ set
\begin{equation}
	\label{eq:boundary-spectral-word-operators}
	\mathfrak J_{\varnothing,N}:=J_N,
	\qquad
	\mathfrak J_{\nu K,N}:=[\tau_\nu,\mathfrak J_{K,N}].
\end{equation}

\begin{lemma}[Uniform commutators of the spectral cutoff]
	\label{lem:spectral-commutators}
	Let $K$ be an ordered word over $\mathfrak A$ with $1\le|K|\le\mathrm N$, and let
	$|K|-1\le s\le\mathrm N-1$.  Then, for every $N\in\bN$ and
	$\psi\in H^s(\partial\cD)$,
	\begin{equation}
		\label{eq:nested-spectral-commutator}
		\|\mathfrak J_{K,N}\psi\|_{H^{s-|K|+1}(\partial\cD)}
		\le C\|\psi\|_{H^s(\partial\cD)},
	\end{equation}
	with $C$ independent of $N$.
\end{lemma}

\begin{proof}
	Fix $\nu\in\mathfrak A$ and let $X:=\tau_\nu$, a tangent field on $\partial\cD$ of class
	$C^{K_\cD+1}\subset C^{\mathrm N+1}$.  Set $h:=N^{-1}$ and $P_h:=h^2A_\partial$, so that
	$J_N=\chi(P_h)$.  Use the $H^s$ norm of \eqref{eq:word-boundary-spectral-Hs}, equivalent to the chart norm for $0\le s\le2\mathrm N$ by Lemma~\ref{lem:word-boundary-spectral-chart}.
	The constants $C$ below are independent of $h$.
	
	For $z\notin\bR$ let $R_h(z)=(z-P_h)^{-1}$.  The spectral theorem gives,
	for $0\le q\le\mathrm N-1$,
	\begin{equation}
		\label{eq:spectral-resolvent-base}
		\|R_h(z)u\|_{H^q}\le |\Im z|^{-1}\|u\|_{H^q}.
	\end{equation}
	For $z$ in a fixed compact set with $0<|\Im z|\le2$ and $0<h\le1$,
	$h^2(1+\lambda)|z-h^2\lambda|^{-1}\le C|\Im z|^{-1}$ for $\lambda\ge0$; hence
	\begin{equation}
		\label{eq:spectral-resolvent-gain}
		\|R_h(z)u\|_{H^{q+2}}
		\le C h^{-2}|\Im z|^{-1}\|u\|_{H^q},
		\qquad 0\le q\le\mathrm N-1.
	\end{equation}
	
	The second-order operator $[X,A_\partial]$ has $C^{K_\cD-1}$ coefficients. The fixed atlas and \eqref{eq:order-and-domain-regularity} yield
	\begin{equation}
		\label{eq:spectral-XA-bound}
		\|[X,A_\partial]u\|_{H^q}
		\le C\|u\|_{H^{q+2}},
		\qquad 0\le q\le\mathrm N-1,
	\end{equation}
	and therefore
	$\|[X,P_h]u\|_{H^q}\le Ch^2\|u\|_{H^{q+2}}$.
	
	Extend $\chi$ to $C_c^\infty(\bR)$ without changing it on $[0,\infty)$, and let
	$\widetilde\chi\in C_c^\infty(\{|\Im z|<2\})$ be an almost-analytic extension with $|\bar\partial\widetilde\chi(z)|\le C|\Im z|^2$.
	The Helffer--Sj\"ostrand formula (see, \textit{e.g.}, \cite[Chapter~2]{Davies1995}), together with $[X,R_h(z)]=R_h(z)[X,P_h]R_h(z)$, gives
	\[
	[X,J_N]
	=-\frac1\pi\int_{\mathbb C}\bar\partial\widetilde\chi(z)
	R_h(z)[X,P_h]R_h(z)\dd A(z).
	\]
	By \eqref{eq:spectral-resolvent-base}--\eqref{eq:spectral-XA-bound},
	\begin{align*}
		\|R_h(z)[X,P_h]R_h(z)u\|_{H^q}
		\,&\le|\Im z|^{-1}\|[X,P_h]R_h(z)u\|_{H^q}\\
		&\le Ch^2|\Im z|^{-1}\|R_h(z)u\|_{H^{q+2}}
		\le C|\Im z|^{-2}\|u\|_{H^q},
		\qquad 0\le q\le\mathrm N-1,
	\end{align*}
	uniformly in $0<h\le1$.  Since $|\bar\partial\widetilde\chi(z)|\,|\Im z|^{-2}$
	is bounded on the compact support of $\widetilde\chi$, the integral converges
	absolutely in $H^q$, and hence
	\[
	\|[X,J_N]u\|_{H^q}\le C\|u\|_{H^q},
	\qquad 0\le q\le\mathrm N-1.
	\]
	The preceding estimates first justify the resolvent commutator identity on smooth functions. Since
	$R_h(z):H^q\to H^{q+2}$ and $[X,P_h]:H^{q+2}\to H^q$ have the stated bounds, the integral defines a bounded operator on $H^q$; density identifies it with $[X,J_N]$. Interpolation gives, for $0\le s\le\mathrm N-1$,
	\begin{equation*}
		\|[X,J_N]\psi\|_{H^s(\partial\cD)}\le C\|\psi\|_{H^s(\partial\cD)},
	\end{equation*}
	with $C$ independent of the cutoff $N$.
	
	Let $K=\nu K'$ with $k:=|K|\ge2$, and suppose the assertion holds for $K'$. For $k-1\le s\le\mathrm N-1$, the commutator identity, initially on smooth functions, gives
	\[
	\begin{aligned}
		\|\mathfrak J_{K,N}\psi\|_{H^{s-k+1}}
		&\le C\|\mathfrak J_{K',N}\psi\|_{H^{s-k+2}}
		+C\|\mathfrak J_{K',N}(\tau_\nu\psi)\|_{H^{s-k+1}}\\
		&\le C\bigl(\|\psi\|_{H^s}+\|\tau_\nu\psi\|_{H^{s-1}}\bigr)
		\le C\|\psi\|_{H^s}.
	\end{aligned}
	\]
	The induction hypothesis applies at both $s$ and $s-1$ because $s\ge k-1$. The result follows by density. The bound depends on the fixed word length and the boundary geometry, but not on $N$.
	
\end{proof}

\begin{lemma}[Sobolev bounds of the spectral cutoff]
	\label{lem:fixedN-range}
	Let $N\in\bN$ and $0\le r\le\mathrm N+2$.  Then, for every
	$\psi\in L^2(\partial\cD)$,
	\begin{equation}
		\label{eq:fixedN-Hr}
		\|J_N\psi\|_{H^r(\partial\cD)}\le C\,(1+N)^r\|\psi\|_{L^2(\partial\cD)},
	\end{equation}
	where $C$ is independent of $N$.
\end{lemma}

\begin{proof}
	Let $\{e_k\}$ and $\{\lambda_k\}$ be the eigenbasis and eigenvalues of $A_\partial$ in the
	proof of Lemma~\ref{lem:boundary-spectral-cutoffs}.  If
	$\operatorname{supp}\chi\subset[0,R_\chi]$, then $\chi(\lambda_k/N^2)\ne0$ implies
	$1+\lambda_k\le(1+R_\chi)(1+N)^2$.  Since $r\le\mathrm N+2\le2\mathrm N$,
	\eqref{eq:word-boundary-spectral-Hs} and
	Lemma~\ref{lem:word-boundary-spectral-chart} give
	\[
	\|J_N\psi\|_{H^r}^2
	\le C\sum_k(1+\lambda_k)^r\chi(\lambda_k/N^2)^2|\langle\psi,e_k\rangle|^2
	\le C\,(1+N)^{2r}\|\psi\|_2^2 ,
	\]
	as claimed.
\end{proof}

\begin{lemma}[Three-factor spectral calculus]
	\label{lem:spectral-three-factor-coloring}
	Let $a$ be a boundary profile with the envelope
	\eqref{eq:boundary-gaussian-polynomial-envelope}, with constants $C_a,c_a$.  For an
	ordered word $K=(\nu_1,\ldots,\nu_k)$ over $\mathfrak A$, $k\le\mathrm N$, set
	\[
	P^a_{K,N}:=\operatorname{ad}^{\gamma}_K\bigl(aJ_N\cA_0\bigr),
	\]
	so that $P^a_{K,N}=P_{K,N}$ at $a=\mathfrak m_{\Theta_w}$.  The constants $C$ below are
	independent of $N$ and depend on $a$ only through $c_a$.
	\begin{enumerate}
		\item\label{it:three-factor-stationary}
		For $\omega\in\{\mathrm L,\mathrm S,\mathrm R\}^{k}$ and
		$c\in\{\mathrm L,\mathrm S,\mathrm R\}$ let $K_{\omega=c}$ be the subword
		$(\nu_i)_{\omega_i=c}$ with $i$ increasing.  Then, on the smooth trace core,
		\begin{equation}
			\label{eq:three-color-expansion}
			P^a_{K,N}
			=\sum_{\omega\in\{\mathrm L,\mathrm S,\mathrm R\}^{k}}
			\bigl(\mathfrak D^-_{K_{\omega=\mathrm L}}a\bigr)\,
			\mathfrak J_{K_{\omega=\mathrm S},N}\,
			\cA_{p_{K_{\omega=\mathrm R}},\mu},
		\end{equation}
		with $p_K$ of \eqref{eq:moment-profile-recursion}, and, at fixed $N$, \eqref{eq:three-color-expansion} defines the bounded extension of $P^a_{K,N}$ to
		$L^2(\gamma^+)$.
		Let $0\le n\le\mathrm N$ and let $\{r_J:|J|\le n\}\subset L^2(\gamma^+)$ be an
		ordered boundary trace chain: for every $|\nu J|\le n$,
		\eqref{eq:rough-moment-hypothesis} holds with
		$(q,q^{[\nu]})=(r_J,r_{\nu J})$.  If $k+|J|\le n$, then, for every finite $\ell$,
		\begin{equation}
			\label{eq:uniform-spectral-trace-commutator}
			\|w^\ell P^a_{K,N}r_J\|_{L^2(\gamma^-)}
			\le C\,C_a
			\sum_{|H|\le|J|+\max\{k-1,0\}}\|r_H\|_{L^2(\gamma^+)} .
		\end{equation}
		If $k+|J|+1\le n$, then, for every $\nu\in\mathfrak A$,
		\begin{equation}
			\label{eq:rough-spectral-product-rule}
			\mathfrak D_\nu^-\bigl(P^a_{K,N}r_J\bigr)
			=P^a_{K,N}r_{\nu J}+P^a_{\nu K,N}r_J
			\quad\text{in }\cD'(\gamma^-),
		\end{equation}
		and, for every finite $\ell$,
		\[
		w^\ell P^a_{K,N}r_{\nu J},\,w^\ell P^a_{\nu K,N}r_J\in L^2(\gamma^-).
		\]
		
		\item\label{it:three-factor-spacetime}
		Let $\mathcal T\subset\bR$ be a bounded open interval, and let $a=a(t)$,
		$t\in\mathcal T$, have the envelope with constants $C_a,c_a$ independent of $t$, with
		$\mathfrak D_{K'}a\in C(\mathcal T\times\gamma^-)$ for $|K'|\le\mathrm N$.  Let
		$0\le n\le\mathrm N$ and let $\{r_J:|J|\le n\}\subset L^2(\mathcal T\times\gamma^+)$ be a
		spacetime trace chain: for every $|\nu J|\le n$ the pair $(r_J,r_{\nu J})$ satisfies
		\eqref{eq:dynamic-generator-slice-spacetime} on $\gamma^+$.  If $k+|J|\le n$, then,
		for every finite $\ell$,
		\begin{equation}
			\label{eq:uniform-spectral-trace-commutator-dynamic}
			\|w^\ell P^{a}_{K,N}r_J\|_{L^2(\mathcal T\times\gamma^-)}
			\le C\,C_a
			\sum_{|H|\le|J|+\max\{k-1,0\}}\|r_H\|_{L^2(\mathcal T\times\gamma^+)},
		\end{equation}
		with $C$ independent of $\mathcal T$.  If $k+|J|+1\le n$, then, for every
		$\nu\in\mathfrak A$,
		\begin{equation}
			\label{eq:rough-spectral-product-rule-dynamic}
			\mathfrak D_\nu^-\bigl(P^{a}_{K,N}r_J\bigr)
			=P^{a}_{K,N}r_{\nu J}+P^{a}_{\nu K,N}r_J
			\quad\text{in }\cD'(\mathcal T\times\gamma^-),
		\end{equation}
		in the sense of \eqref{eq:dynamic-generator-slice-spacetime} on $\gamma^-$, and, for every finite $\ell$,
		\[
		w^\ell P^{a}_{K,N}r_{\nu J},\,w^\ell P^{a}_{\nu K,N}r_J\in L^2(\mathcal T\times\gamma^-).
		\]
		The hypotheses on $a$ hold
		for $a=\mathfrak m_\circ+d_\delta(t)$, for which $P^a_{K,N}=P_{K,N}(t)$, with
		$C_a\le C$, and for $a=\partial_t^id_\delta(t)$, $1\le i\le\mathrm N$, with
		$C_a\le C\delta$; here $C$ is independent of $\delta$ and $t$.
	\end{enumerate}
\end{lemma}

\begin{proof}
	\emph{(\ref{it:three-factor-stationary}).}
	Let $\nu\in\mathfrak A$, let $b$ be a $C^1$ profile on $\gamma^-$, and let $K_1,K_2$ be ordered words over
	$\mathfrak A$.  On the smooth trace core, \eqref{eq:moment-commutator} and
	\eqref{eq:moment-profile-recursion} give
	\[
	\tau_\nu\cA_{p_{K_2},\mu}-\cA_{p_{K_2},\mu}\mathfrak D_\nu^+=\cA_{p_{\nu K_2},\mu}.
	\]
	Together with the one-letter boundary Leibniz rule and \eqref{eq:boundary-spectral-word-operators}, this gives
	\[
	\operatorname{ad}^{\gamma}_\nu\bigl(b\,\mathfrak J_{K_1,N}\cA_{p_{K_2},\mu}\bigr)
	=(\mathfrak D_\nu^-b)\,\mathfrak J_{K_1,N}\cA_{p_{K_2},\mu}
	+b\,\mathfrak J_{\nu K_1,N}\cA_{p_{K_2},\mu}
	+b\,\mathfrak J_{K_1,N}\cA_{p_{\nu K_2},\mu}.
	\]
	Since $aJ_N\cA_0=a\,\mathfrak J_{\varnothing,N}\cA_{p_\varnothing,\mu}$, applying this
	with $\nu=\nu_k,\nu_{k-1},\ldots,\nu_1$ proves \eqref{eq:three-color-expansion}; each
	letter enters exactly one of the three factors.  At fixed $N$, $J_N$ is self-adjoint and has finite rank
	with range in $H^{\mathrm N+2}(\partial\cD)\subset C^{\mathrm N}(\partial\cD)$
	(Lemmas~\ref{lem:boundary-spectral-cutoffs} and~\ref{lem:fixedN-range}).  Defining $J_N\tau_I\psi$ by duality through \eqref{eq:surface-tangential-adjoint}, we obtain
	\[
	\|J_N\tau_I\psi\|_{H^{\mathrm N+2}(\partial\cD)}\le C_N\|\psi\|_{L^2(\partial\cD)},
	\qquad \psi\in L^2(\partial\cD),\quad |I|\le\mathrm N.
	\]
	Since every
	$\mathfrak J_{K,N}$, $|K|\le\mathrm N$, is a finite sum of terms $\pm\tau_{I_1}J_N\tau_{I_2}$, it is bounded on $L^2(\partial\cD)$. Hence every term
	of \eqref{eq:three-color-expansion} is bounded from $L^2(\gamma^+)$ to $L^2(\gamma^-)$ by
	\eqref{eq:word-profile-moment-bound} and the envelope of $\mathfrak D^-_{K_{\omega=\mathrm L}}a$.
	
	Fix $\omega$, write $K_{\rm L},K_{\rm S},K_{\rm R}$ for
	$K_{\omega=\mathrm L},K_{\omega=\mathrm S},K_{\omega=\mathrm R}$ and
	$k_{\rm L},k_{\rm S},k_{\rm R}$ for their lengths, and set $j:=|J|$.  By
	\eqref{eq:word-rough-moment-one} with $p\sqrt\mu$ in place of $a$ and
	$\mathscr E_\nu(p\sqrt\mu)=(\mathscr E_{\nu,\mu}p)\sqrt\mu$,
	\[
	\tau_\nu\cA_{p,\mu}r_H=\cA_{p,\mu}r_{\nu H}+\cA_{\mathscr E_{\nu,\mu}p,\mu}r_H
	\quad\text{in }\cD'(\partial\cD),
	\qquad |\nu H|\le n .
	\]
	Iterating this identity along words $I$ over $\mathfrak A$ with $|I|\le s$, we write
	$\tau_I\cA_{p_{K_{\rm R}},\mu}r_J$ as a finite sum of terms $\cA_{p_{K'},\mu}r_H$ with
	$|H|\le j+s$ and $|K'|\le k_{\rm R}+s$.  For $0\le s\le k_{\rm S}$ with $j+s\le n$, the words
	$K'$ satisfy $|K'|\le k_{\rm R}+s\le k_{\rm R}+k_{\rm S}\le k\le\mathrm N$, so
	Lemma~\ref{lem:boundary-profile-bounds}, \eqref{eq:word-profile-moment-bound} and
	Lemma~\ref{lem:tangential-word-sobolev} give
	\begin{equation}
		\label{eq:r11-q2-moment-Hs}
		\|\cA_{p_{K_{\rm R}},\mu}r_J\|_{H^{s}(\partial\cD)}
		\le C\sum_{|H|\le j+s}\|r_H\|_{L^2(\gamma^+)} .
	\end{equation}
	By \eqref{eq:boundary-gaussian-polynomial-envelope}, as in the proof of
	Lemma~\ref{lem:boundary-profile-bounds},
	$|w^\ell\mathfrak D^-_{K_{\rm L}}a|\le C\,C_a\ee^{-c_a|v|^2/2}$, and  integration in $v$ against $|n_x\cdot v|\dd v$ gives
	\[
	\|w^\ell(\mathfrak D^-_{K_{\rm L}}a)\psi\|_{L^2(\gamma^-)}\le C\,C_a\|\psi\|_{L^2(\partial\cD)}.
	\]
	If $k_{\rm S}=0$, then $\mathfrak J_{K_{\rm S},N}=J_N$, and the bound for $\|w^\ell(\mathfrak D^-_{K_{\rm L}}a)\psi\|_{L^2(\gamma^-)}$ above with $\psi=J_N\cA_{p_{K_{\rm R}},\mu}r_J$, the bound $\|J_N\|\le1$ on $L^2(\partial\cD)$, and \eqref{eq:r11-q2-moment-Hs} at $s=0$ give, for the $\omega$-term of $P^a_{K,N}r_J$,
	\[
	\|w^\ell(\mathfrak D^-_{K_{\rm L}}a)\,J_N\cA_{p_{K_{\rm R}},\mu}r_J\|_{L^2(\gamma^-)}
	\le C\,C_a\sum_{|H|\le j}\|r_H\|_{L^2(\gamma^+)} .
	\]
	If $k_{\rm S}\ge1$,
	the bound for $\|w^\ell(\mathfrak D^-_{K_{\rm L}}a)\psi\|_{L^2(\gamma^-)}$ with $\psi=\mathfrak J_{K_{\rm S},N}\cA_{p_{K_{\rm R}},\mu}r_J$, and then \eqref{eq:nested-spectral-commutator} and \eqref{eq:r11-q2-moment-Hs}, both at
	$s=k_{\rm S}-1$, give
	\[
	\begin{aligned}
		&\|w^\ell(\mathfrak D^-_{K_{\rm L}}a)\,\mathfrak J_{K_{\rm S},N}\cA_{p_{K_{\rm R}},\mu}r_J\|_{L^2(\gamma^-)}
		\le C\,C_a\|\mathfrak J_{K_{\rm S},N}\cA_{p_{K_{\rm R}},\mu}r_J\|_{L^2(\partial\cD)}\\
		&\qquad \le C\,C_a\|\cA_{p_{K_{\rm R}},\mu}r_J\|_{H^{k_{\rm S}-1}(\partial\cD)}
		\le C\,C_a\sum_{|H|\le j+k_{\rm S}-1}\|r_H\|_{L^2(\gamma^+)} .
	\end{aligned}
	\]
	Summing over the three-colorings gives \eqref{eq:uniform-spectral-trace-commutator}.
	
	For \eqref{eq:rough-spectral-product-rule}, apply $\mathfrak D_\nu^-$ in $\cD'(\gamma^-)$
	to each term of \eqref{eq:three-color-expansion} at $r_J$.
	Lemma~\ref{lem:boundary-word-Leibniz} for one letter differentiates the product of
	$\mathfrak D^-_{K_{\rm L}}a$ and the scalar
	$\psi:=\mathfrak J_{K_{\rm S},N}\cA_{p_{K_{\rm R}},\mu}r_J$; the new factor $\mathfrak D^-_{\nu K_{\rm L}}a$ has order $k_{\rm L}+1\le\mathrm N$, since $k_{\rm L}\le k$ and $k+1\le n\le\mathrm N$.
	By \eqref{eq:boundary-spectral-word-operators},
	\[
	\tau_\nu\psi=\mathfrak J_{K_{\rm S},N}\tau_\nu\cA_{p_{K_{\rm R}},\mu}r_J
	+\mathfrak J_{\nu K_{\rm S},N}\cA_{p_{K_{\rm R}},\mu}r_J,
	\]
	and the identity for $\tau_\nu\cA_{p,\mu}r_H$ above, at $p=p_{K_{\rm R}}$ and $H=J$, differentiates the moment; here $\mathscr E_{\nu,\mu}p_{K_{\rm R}}=p_{\nu K_{\rm R}}$ by \eqref{eq:moment-profile-recursion}.  Since $k+j+1\le n$, \eqref{eq:r11-q2-moment-Hs} also holds for $\cA_{p_{K_{\rm R}},\mu}r_{\nu J}$ and $\cA_{p_{\nu K_{\rm R}},\mu}r_J$. These bounds, \eqref{eq:nested-spectral-commutator} and $\|J_N\|\le1$ on $L^2(\partial\cD)$ show that every term of $\tau_\nu\psi$ belongs to $L^2(\partial\cD)$.  The new letter either differentiates one of the three factors, producing $P^a_{\nu K,N}r_J$, or acts on the trace chain, producing $P^a_{K,N}r_{\nu J}$. This proves \eqref{eq:rough-spectral-product-rule}.  Both terms on
	the right satisfy \eqref{eq:uniform-spectral-trace-commutator}, at $(K,\nu J)$ and at
	$(\nu K,J)$.
	
	\emph{(\ref{it:three-factor-spacetime}).}  By Lemma~\ref{lem:dynamic-generator-slices},
	applied to the finitely many pairs of the chain, there is a null set
	$\mathcal N\subset\mathcal T$ such that, for $t\notin\mathcal N$,
	$\{r_J(t)\}_{|J|\le n}$ is an ordered boundary trace chain in $L^2(\gamma^+)$.  At such
	$t$, (\ref{it:three-factor-stationary}) applies with the profile $a(t)$ and constants
	independent of $t$ (the constants $C$ there depend on $a$ only through $c_a$, and $C_a$, $c_a$ do not depend on $t$); the terms of \eqref{eq:three-color-expansion} at $r_J(t)$ are measurable in $t$ since $\mathfrak D_{K'}a$,
	$|K'|\le\mathrm N$, are continuous in $t$.  Squaring
	\eqref{eq:uniform-spectral-trace-commutator} and integrating over $\mathcal T$ gives
	\eqref{eq:uniform-spectral-trace-commutator-dynamic}.  At each $t\notin\mathcal N$, we extend
	\eqref{eq:rough-spectral-product-rule} to tests reaching $\gamma^0$ by the sign cutoffs
	\eqref{eq:boundary-generator-sign-cutoff} as in the proof of
	Proposition~\ref{prop:weak-boundary-generator}. The extended identity is \eqref{eq:dynamic-generator-slice-identity} on
	$\gamma^-$ for
	\[
	r=P^a_{K,N}r_J(t),
	\qquad
	r^{[\nu]}=P^a_{K,N}r_{\nu J}(t)+P^a_{\nu K,N}r_J(t),
	\]
	which lie in
	$L^2(\mathcal T\times\gamma^-)$ by \eqref{eq:uniform-spectral-trace-commutator-dynamic}. Hence
	Lemma~\ref{lem:dynamic-generator-slices} gives
	\eqref{eq:rough-spectral-product-rule-dynamic}.  For $a=\mathfrak m_\circ+d_\delta(t)$
	and $a=\partial_t^id_\delta(t)$, $1\le i\le\mathrm N$,
	\eqref{eq:dyn-theta-defect-time-bound} gives,
	as in the paragraph after Definition~\ref{def:gaussian-polynomial-envelope}, the
	envelope uniformly in $t$ with the stated $C_a$. Since $\mathfrak m_\circ$ does not depend on $t$ and $d_\delta=c_\mu^{-1}h_w$, we have $\partial_t\mathfrak D_{K'}a=c_\mu^{-1}\mathfrak D_{K'}\partial_t^{i+1}h_w$ for $|K'|\le\mathrm N$, with $i=0$ for $a=\mathfrak m_\circ+d_\delta(t)$.  Hence \eqref{eq:dyn-theta-defect-complete-bound} at $L=0$ and $j=i+1\le\mathrm N+1$ gives
	\[
	|\partial_t\mathfrak D_{K'}a(t,x,v)|\le C\delta\ee^{-c|v|^2},
	\qquad t\in\mathcal T,\ (x,v)\in\gamma^-,\ |K'|\le\mathrm N,
	\]
	with $C$ and $c$ independent of $\delta$ and $t$.  Since $\mathfrak D_{K'}a(t)\in C(\gamma^-)$ for each $t$, this bound gives the continuity of $\mathfrak D_{K'}a$ on $\mathcal T\times\gamma^-$.
\end{proof}

\subsubsection{The differentiated diffuse law}\label{ss:dyn-differentiated-equations}

\begin{prop}[Ordered differentiated diffuse law]
	\label{prop:ordered-differentiated-diffuse-law}
	\begin{enumerate}
		\item\label{it:diffuse-law-physical}
		Let $f\in\cH_{\kin,\sigma}(\cD\times\bR_v^3)$ satisfy
		\begin{equation}
			\label{eq:stationary-boundary-hypothesis}
			f_{\gamma^-}=P_wf_{\gamma^+}+h_w
			\qquad\text{in }L^2(\gamma^-)
		\end{equation}
		for a datum $h_w$ from Lemma~\ref{lem:wall-profile-derivatives}.
		Fix $1\le m\le\mathrm N$ and assume that the ordered jets $f_J$,
		$|J|\le m$, with $f_\varnothing=f$, belong to
		$\cH_{\kin,\sigma}(\cD\times\bR_v^3)$, satisfy
		\begin{equation}
			\label{eq:ordered-boundary-rough-chain}
			\mathfrak D_\nu f_J=f_{\nu J}
			\quad\text{in }\cD'(\cD\times\bR^3),
			\qquad |\nu J|\le m,
		\end{equation}
		and carry the completion-defined traces of
		Definition~\ref{def:H_kin_sigma}.  Then, for every
		$I=(\nu_1,\ldots,\nu_k)\in\mathfrak A^k$ with $1\le k\le m$,
		\begin{equation}
			\label{eq:ordered-differentiated-diffuse-law}
			(f_I)_{\gamma^-}=P_w(f_I)_{\gamma^+}+b_I,
		\end{equation}
		\begin{equation}
			\label{eq:ordered-boundary-b}
			b_I=\sum_{\varnothing\ne S\subset[k]}
			P_{I_S}\bigl((f_{I_{S^c}})_{\gamma^+}\bigr)
			+\mathfrak D_I^-h_w,
		\end{equation}
		in $L^2(\gamma^-)$, where the sum is over sets $S$ of positions in $I$, as in \eqref{eq:ordered-two-coloring-convention}.  For every finite $\ell\ge0$,
		\begin{equation}
			\label{eq:ordered-boundary-triangular-solution}
			\bigl\|w^\ell(b_I-\mathfrak D_I^-h_w)\bigr\|_{L^2(\gamma^-)}
			\le C\sum_{|J|<|I|}
			\bigl\|(f_J)_{\gamma^+}\bigr\|_{L^2(\gamma^+)}
		\end{equation}
		and
		\begin{equation}
			\label{eq:ordered-boundary-triangular-full}
			\bigl\|w^{\ell}b_I\bigr\|_{L^2(\gamma^-)}
			\le C\Bigl(\delta+
			\sum_{|J|<|I|}\bigl\|(f_J)_{\gamma^+}\bigr\|_{L^2(\gamma^+)}\Bigr),
		\end{equation}
		where $C$ is independent of $\delta$.
		\item\label{it:diffuse-law-fixedN}
		Fix $N\in\bN$ and let $f$ satisfy the hypotheses of
		(\ref{it:diffuse-law-physical}) with
		\eqref{eq:stationary-boundary-hypothesis} replaced by
		\begin{equation}
			\label{eq:spectral-triangular-hypothesis}
			f_{\gamma^-}=P_{w,N}f_{\gamma^+}+h_w .
		\end{equation}
		Then, for $I$ as in (\ref{it:diffuse-law-physical}),
		\begin{equation}
			\label{eq:spectral-triangular-law}
			(f_I)_{\gamma^-}
			=P_{w,N}(f_I)_{\gamma^+} +b_{I,N},
			\qquad
			b_{I,N}
			:=\sum_{\varnothing\ne S\subset[k]}
			P_{I_S,N}\bigl((f_{I_{S^c}})_{\gamma^+}\bigr)
			+\mathfrak D_I^-h_w,
		\end{equation}
		in $L^2(\gamma^-)$, and, for every finite $\ell\ge0$,
		\begin{equation}
			\label{eq:spectral-triangular-bound}
			\|w^\ell b_{I,N}\|_{L^2(\gamma^-)}
			\le C\Bigl(\delta+\sum_{|J|<|I|}
			\|(f_J)_{\gamma^+}\|_{L^2(\gamma^+)}\Bigr),
		\end{equation}
		where $C$ is independent of $N$ and $\delta$.
	\end{enumerate}
\end{prop}

\begin{proof}
	\emph{(\ref{it:diffuse-law-physical}).}
	Let $|\nu J|\le m$, $q:=(f_J)_{\gamma^+}$ and $q^{[\nu]}:=(f_{\nu J})_{\gamma^+}$.
	Proposition~\ref{prop:weak-boundary-generator}\,(\ref{it:weak-generator-stationary})
	gives \eqref{eq:rough-moment-hypothesis} for $(q,q^{[\nu]})$ and, on $\gamma^-$,
	$\mathfrak D_\nu^-(f_J)_{\gamma^-}=(f_{\nu J})_{\gamma^-}$ in $\cD'(\gamma^-)$.  Let
	$|I|\le\mathrm N-1$.  On each boundary chart, Lemma~\ref{lem:boundary-word-Leibniz} for
	one letter differentiates each term $m_{I,j}\cA_{p_{I,j},\mu}q$ of
	\eqref{eq:ordered-boundary-finite-rank}, and \eqref{eq:word-rough-moment-one} at
	$a=p_{I,j}\sqrt\mu$ gives
	\[
	\tau_\nu\cA_{p_{I,j},\mu}q=\cA_{p_{I,j},\mu}q^{[\nu]}+\cA_{\mathscr E_{\nu,\mu}p_{I,j},\mu}q.
	\]
	Then \eqref{eq:ordered-boundary-commutator-induction} collects the terms:
	\begin{equation}
		\label{eq:rough-PI-derivative}
		\mathfrak D_\nu^-(P_Iq)=P_Iq^{[\nu]}+P_{\nu I}q
		\quad\text{in }\cD'(\gamma^-),
	\end{equation}
	and both terms on the right belong to $L^2(\gamma^-)$ at every finite weight by
	\eqref{eq:ordered-boundary-finite-rank-bound}.
	
	At $I=\varnothing$, \eqref{eq:ordered-differentiated-diffuse-law} with $b_\varnothing=h_w$ is
	\eqref{eq:stationary-boundary-hypothesis}.  Let $I=\nu J$ with $|I|\le m$, assume
	\eqref{eq:ordered-differentiated-diffuse-law}--\eqref{eq:ordered-boundary-b} at $J$, and
	apply $\mathfrak D_\nu^-$ in $\cD'(\gamma^-)$.  The identity
	$\mathfrak D_\nu^-(f_J)_{\gamma^-}=(f_{\nu J})_{\gamma^-}$ differentiates the left side, and
	\eqref{eq:rough-PI-derivative} differentiates $P_w(f_J)_{\gamma^+}=P_\varnothing(f_J)_{\gamma^+}$ and every $P_{J_S}\bigl((f_{J_{S^c}})_{\gamma^+}\bigr)$ on the right; $\mathfrak D_\nu^-\mathfrak D_J^-h_w=\mathfrak D_I^-h_w$
	classically, since $h_w$ is of class $C^{K_\cD}$.  In
	\[
	\mathfrak D_\nu^-\bigl(P_{J_S}(f_{J_{S^c}})_{\gamma^+}\bigr)
	=P_{J_S}(f_{\nu J_{S^c}})_{\gamma^+}+P_{\nu J_S}(f_{J_{S^c}})_{\gamma^+} ,
	\]
	the new position
	is added to $S^c$ or to $S$, and
	\[
	\mathfrak D_\nu^-\bigl(P_w(f_J)_{\gamma^+}\bigr)=P_w(f_I)_{\gamma^+}+P_\nu(f_J)_{\gamma^+}
	\]
	supplies the summand at $S=\{1\}$.  This is
	\eqref{eq:ordered-differentiated-diffuse-law}--\eqref{eq:ordered-boundary-b} at $I$,
	and all its terms lie in $L^2(\gamma^-)$.  Each of
	the $2^{|I|}-1$ terms of $b_I-\mathfrak D_I^-h_w$ acts on a trace $(f_{I_{S^c}})_{\gamma^+}$ with $|I_{S^c}|\le|I|-1$.  By
	\eqref{eq:ordered-boundary-finite-rank-bound},
	\[
	\|w^\ell P_{I_S}((f_{I_{S^c}})_{\gamma^+})\|_{L^2(\gamma^-)}\le C\|(f_{I_{S^c}})_{\gamma^+}\|_{L^2(\gamma^+)},
	\]
	and summation over $S$ gives
	\eqref{eq:ordered-boundary-triangular-solution}; with
	\eqref{eq:wall-profile-complete-bound} for $\mathfrak D_I^-h_w$ this gives
	\eqref{eq:ordered-boundary-triangular-full}.
	
	\emph{(\ref{it:diffuse-law-fixedN}).}
	The proof of (\ref{it:diffuse-law-physical}) applies with $P_{w,N}$ and $P_{I,N}$ in
	place of $P_w$ and $P_I$.  By
	Proposition~\ref{prop:weak-boundary-generator}\,(\ref{it:weak-generator-stationary}),
	the traces $r_J:=(f_J)_{\gamma^+}$, $|J|\le m$, form
	an ordered boundary trace chain. Then
	Lemma~\ref{lem:spectral-three-factor-coloring}\,(\ref{it:three-factor-stationary}) applies at $a=\mathfrak m_{\Theta_w}$, for which $P^a_{K,N}=P_{K,N}$ and  the envelope has $C_a\le C$. It supplies
	\eqref{eq:rough-spectral-product-rule} in place of \eqref{eq:rough-PI-derivative} and
	\eqref{eq:uniform-spectral-trace-commutator} in place of
	\eqref{eq:ordered-boundary-finite-rank-bound}, with $n=m$. The conditions of \eqref{eq:rough-spectral-product-rule} and \eqref{eq:uniform-spectral-trace-commutator} hold at $P_{J_S,N}\bigl((f_{J_{S^c}})_{\gamma^+}\bigr)$ in the induction and at $P_{I_S,N}\bigl((f_{I_{S^c}})_{\gamma^+}\bigr)$ in the bound, since 
	$$
	|J_S|+|J_{S^c}|+1=|I_S|+|I_{S^c}|=|I|\le m.
	$$
	The constant is independent of $N$ and $\delta$, by the uniform Maxwellian envelope bound.
	This gives \eqref{eq:spectral-triangular-law} and \eqref{eq:spectral-triangular-bound}.
\end{proof}

\begin{prop}[Ordered differentiated diffuse law: the evolutionary case]
	\label{prop:spectral-triangular-evolutionary}
	Fix $N\in\bN$, $0<T<\infty$ and $1\le m_\ast\le\mathrm N$, and let
	$f\in\cH_{\kin,\sigma}(\cD_T\times\bR_v^3)$ satisfy
	\begin{equation}
		\label{eq:dyn-temporal-boundary-hypothesis-N}
		f_{\gamma^-}=P_{w,N}(t)f_{\gamma^+}+b
		\qquad\text{in }L^2(\gamma_T^-)
	\end{equation}
	for a datum $b$ with the following regularity: for $0\le i\le m_\ast$ and every
	ordered word $K$ over $\mathfrak A$ with $|K|+i\le m_\ast$ there is
	$\mathfrak D_K^-\partial_t^ib\in L^2(\gamma_T^-)$, with
	$\mathfrak D_\varnothing^-\partial_t^ib=\partial_t^ib$, such that, for every
	$\nu\in\mathfrak A$ with $|\nu K|+i\le m_\ast$ and every compact-time,
	compact-velocity $C^2$ phase test $\Psi$,
	\begin{equation}
		\label{eq:dyn-datum-spatial-regularity}
		\int_{\gamma_T^-}\bigl(\mathfrak D_{\nu K}^-\partial_t^ib\bigr)\Psi\,\dd\gamma\dd t
		=-\int_{\gamma_T^-}\bigl(\mathfrak D_K^-\partial_t^ib\bigr)
		(\mathfrak D_\nu+d_\nu)\Psi\,\dd\gamma\dd t .
	\end{equation}
	Assume that the mixed jets $f_J$, $J$ an ordered word over $\mathfrak A_\bullet$ with
	$|J|\le m_\ast$ and $f_\varnothing=f$, belong to
	$\cH_{\kin,\sigma}(\cD_T\times\bR_v^3)$, carry the completion-defined traces of
	Definition~\ref{def:H_kin_sigma}, and satisfy
	\begin{equation}
		\label{eq:dyn-jet-chain}
		\mathfrak D_\nu f_J=f_{\nu J}
		\quad\text{in }\cD'\bigl((0,T)\times\cD\times\bR^3\bigr),
		\qquad
		\nu\in\mathfrak A_\bullet,\quad|\nu J|\le m_\ast .
	\end{equation}
	Let $I$ be an ordered word over $\mathfrak A_\bullet$ with $1\le|I|\le m_\ast$, spatial
	subword $I'$, $m':=|I'|$, and temporal order $j:=j(I)$.  Then
	\begin{equation}
		\label{eq:dyn-mixed-boundary-law-N}
		\begin{aligned}
			(f_I)_{\gamma^-}\,&=P_{w,N}(f_I)_{\gamma^+}
			+\mathfrak b_{I,N}[f_{\gamma^+}]
			+\mathfrak D_{I'}^-
			\bigl(\mathfrak t_{j,N}[f]+\partial_t^{\,j}b\bigr),\\
			\mathfrak b_{I,N}[f_{\gamma^+}]
			\,&:=\sum_{\varnothing\ne S\subset[m']}
			P_{I'_S,N}\bigl((\mathfrak D_{I'_{S^c}}\partial_t^{\,j}f)_{\gamma^+}\bigr),
		\end{aligned}
	\end{equation}
	in $L^2(\gamma_T^-)$, where the sum is over sets $S$ of positions in the spatial subword $I'$ only.
	Here $\mathfrak D_{I'}^-\partial_t^jb$ is the derivative in
	\eqref{eq:dyn-datum-spatial-regularity}, and
	$\mathfrak D_{I'}^-\mathfrak t_{j,N}[f]$, the weak derivative of
	\eqref{eq:dyn-theta-triangular-notation} along $I'$, is the finite sum
	\begin{equation}
		\label{eq:dyn-temporal-part-spatial-expansion-N}
		\mathfrak D_{I'}^-\mathfrak t_{j,N}[f]
		=\sum_{k=1}^{j}\binom jk\sum_{S\subset[m']}
		P_{I'_S\mathrm t^k,N}\bigl((f_{I'_{S^c}\mathrm t^{j-k}})_{\gamma^+}\bigr).
	\end{equation}
	For every finite $\ell\ge0$,
	\begin{equation}
		\label{eq:dyn-temporal-part-spatial-bound}
		\bigl\|w^{\ell}\mathfrak D_{I'}^-\mathfrak t_{j,N}[f]\bigr\|_{L^2(\gamma_T^-)}
		\le C\,\delta
		\sum_{|J|<|I|,\ j(J)<j}
		\bigl\|(f_J)_{\gamma^+}\bigr\|_{L^2(\gamma_T^+)}
	\end{equation}
	and
	\begin{equation}
		\label{eq:dyn-mixed-boundary-bound}
		\bigl\|w^{\ell}\mathfrak b_{I,N}[f_{\gamma^+}]\bigr\|_{L^2(\gamma_T^-)}
		\le C
		\sum_{|J|<|I|,\ j(J)=j}
		\bigl\|(f_J)_{\gamma^+}\bigr\|_{L^2(\gamma_T^+)},
	\end{equation}
	where $C$ is independent of $T$, $N$ and $\delta$.
\end{prop}

\begin{proof}
	Since $f_{\mathrm t^i}=\partial_t^if$ by \eqref{eq:dyn-jet-chain}, the time chain
	\eqref{eq:dynamic-trace-time-chain-boundary} gives
	$(\partial_t^if)_{\gamma^\pm}=\partial_t^if_{\gamma^\pm}$ for $i\le j$.  Hence the
	truncated temporal law \eqref{eq:dyn-temporal-boundary-law} holds. All its terms lie
	in $L^2(\gamma_T^-)$. On the right side, $P_{w,N}(t)(\partial_t^jf)_{\gamma^+}$ does so by \eqref{eq:wall-J-weighted-operator-bound} at $\cJ=J_N$, applied at each $t$ with $\Theta_w(t)$ in place of $\Theta_w$, with a bound uniform in $t$ by \eqref{eq:dyn-theta-defect-time-bound}; $\mathfrak t_{j,N}[f]$ does so by \eqref{eq:dyn-theta-triangular-bound}, and $\partial_t^jb$ by the hypothesis on $b$.
	
	By \eqref{eq:dyn-jet-chain} and induction on the letters, $f_J=\mathfrak D_Jf$
	in $\cD'$ for $|J|\le m_\ast$, and \eqref{eq:dyn-normal-ordering} gives
	$\mathfrak D_Jf=\mathfrak D_{J'\mathrm t^{j(J)}}f$.  Hence
	$f_J=f_{J'\mathrm t^{j(J)}}$ in $\cH_{\kin,\sigma}(\cD_T\times\bR_v^3)$, with the same
	traces, and we may take $I=I'\mathrm t^j$, $I'=(\nu_1,\ldots,\nu_{m'})$.  For $K$
	spatial, $\nu\in\mathfrak A$ and $|\nu K|+i\le m_\ast$,
	Proposition~\ref{prop:weak-boundary-generator}\,(\ref{it:weak-generator-evolutionary}),
	applied to the pair $(f_{K\mathrm t^i},f_{\nu K\mathrm t^i})$, gives
	\eqref{eq:dynamic-boundary-generator} on $\gamma_T^+$ and on $\gamma_T^-$.  On
	$\gamma_T^+$, where $\operatorname{sgn}(n\cdot v)=1$, this is
	\eqref{eq:dynamic-generator-slice-spacetime}; so, for every $0\le i\le m_\ast$, the traces
	$r^{(i)}_K:=(f_{K\mathrm t^i})_{\gamma^+}$, $K$ spatial with $|K|\le m_\ast-i$, form a
	spacetime trace chain on $\mathcal T=(0,T)$.
	
	Fix $1\le k\le j$ and let $a:=\partial_t^kd_\delta(t)$.  By
	\eqref{eq:dyn-mixed-temporal-commutator}, $P_{\mathrm t^k,N}(t)=aJ_N\cA_0$. Hence, by \eqref{eq:dyn-mixed-boundary-commutator},
	\[
	P_{K\mathrm t^k,N}(t)=\operatorname{ad}^{\gamma}_K(P_{\mathrm t^k,N}(t))=P^a_{K,N}
	\]
	for spatial $K$ with $|K|+k\le\mathrm N$, in particular for $|K|\le m'$, since
	$m'+j\le m_\ast\le\mathrm N$.
	Lemma~\ref{lem:spectral-three-factor-coloring}\,(\ref{it:three-factor-spacetime})
	applies on $\mathcal T=(0,T)$ to this $a$, with $C_a\le C\delta$, and to the chain
	$\{r^{(j-k)}_K\}$ defined above.  We prove, for every suffix $K$ of $I'$,
	\begin{equation}
		\label{eq:r11-q2-suffix-temporal}
		\mathfrak D_K^-\bigl(P_{\mathrm t^k,N}(f_{\mathrm t^{j-k}})_{\gamma^+}\bigr)
		=\sum_{S\subset[|K|]}P_{K_S\mathrm t^k,N}
		\bigl((f_{K_{S^c}\mathrm t^{j-k}})_{\gamma^+}\bigr)
		\quad\text{in }L^2(\gamma_T^-) .
	\end{equation}
	At $K=\varnothing$ it is trivial.  If it holds at $K$ and $\nu K$ is a suffix of $I'$,
	\eqref{eq:rough-spectral-product-rule-dynamic}, applied to each summand with $K_S$ as
	operator word and $K_{S^c}$ as chain index, gives it at $\nu K$, and the new position is added to
	$S^c$ or to $S$. The hypothesis for this chain, with $n=m_\ast-(j-k)$, follows from
	\[
	|K_S|+|K_{S^c}|+1=|\nu K|\le m'\le m_\ast-(j-k).
	\]  At $K=I'$, multiplication by $\binom jk$ and
	summation over $k$ give \eqref{eq:dyn-temporal-part-spatial-expansion-N}, and each
	successive weak derivative belongs to $L^2(\gamma_T^-)$.  By
	\eqref{eq:uniform-spectral-trace-commutator-dynamic} with $C_a\le C\delta$, each summand satisfies
	\[
	\bigl\|w^\ell P_{I'_S\mathrm t^k,N}\bigl((f_{I'_{S^c}\mathrm t^{j-k}})_{\gamma^+}\bigr)\bigr\|_{L^2(\gamma_T^-)}
	\le C\delta\sum_{|H|\le m'}\bigl\|(f_{H\mathrm t^{j-k}})_{\gamma^+}\bigr\|_{L^2(\gamma_T^+)}.
	\]
	Since $k\ge1$, every $H$ there satisfies $|H\mathrm t^{j-k}|\le m'+j-k<|I|$ and $j(H\mathrm t^{j-k})=j-k<j$, which gives
	\eqref{eq:dyn-temporal-part-spatial-bound}.
	
	Set $h:=\mathfrak t_{j,N}[f]+\partial_t^jb$ and, for a suffix $K$ of $I'$,
	$\mathfrak D_K^-h:=\mathfrak D_K^-\mathfrak t_{j,N}[f]+\mathfrak D_K^-\partial_t^jb$, with the first term given by \eqref{eq:r11-q2-suffix-temporal}, summed over $k$ as in \eqref{eq:dyn-temporal-part-spatial-expansion-N}, and the second by \eqref{eq:dyn-datum-spatial-regularity}.  We
	prove \eqref{eq:dyn-mixed-boundary-law-N}, with $I'$ replaced by every suffix $K$ of
	$I'$ and $\mathfrak D_{I'}^-(\cdots)$ by $\mathfrak D_K^-h$, by the induction in the proof of
	Proposition~\ref{prop:ordered-differentiated-diffuse-law}\,(\ref{it:diffuse-law-physical})
	in $\cD'(\gamma_T^-)$; at $K=\varnothing$ it is \eqref{eq:dyn-temporal-boundary-law}.  In
	the passage from $K$ to $\nu K$, \eqref{eq:dynamic-boundary-generator} on $\gamma_T^-$ for the pair $(f_{K\mathrm t^j},f_{\nu K\mathrm t^j})$ gives $(f_{\nu K\mathrm t^j})_{\gamma^-}$ as the
	weak $\mathfrak D_\nu^-$ derivative of $(f_{K\mathrm t^j})_{\gamma^-}$.
	Applied at $a=\mathfrak m_\circ+d_\delta(t)$ and to the chain $\{r^{(j)}_K\}$, \eqref{eq:rough-spectral-product-rule-dynamic} gives the weak derivative of each $P_{K_S,N}\bigl((f_{K_{S^c}\mathrm t^j})_{\gamma^+}\bigr)$, $S\subset[|K|]$, where $P_{\varnothing,N}=P_{w,N}(t)$. The condition of \eqref{eq:rough-spectral-product-rule-dynamic} for this chain holds since $|\nu K|\le m'\le m_\ast-j$. The new position is added to $S^c$ or to $S$. Finally, \eqref{eq:r11-q2-suffix-temporal} with
	\eqref{eq:dyn-datum-spatial-regularity} gives $\mathfrak D_{\nu K}^-h$.  All terms lie in
	$L^2(\gamma_T^-)$.  At $K=I'$ this is \eqref{eq:dyn-mixed-boundary-law-N}, by the identifications $f_J=\mathfrak D_Jf$ and $f_J=f_{J'\mathrm t^{j(J)}}$ above.
	
	Applying \eqref{eq:uniform-spectral-trace-commutator-dynamic} at $a=\mathfrak m_\circ+d_\delta(t)$ to each summand of $\mathfrak b_{I,N}$, with $I'_S$ as operator word and $I'_{S^c}$ as chain index, we obtain
	\[
	\bigl\|w^\ell P_{I'_S,N}\bigl((f_{I'_{S^c}\mathrm t^j})_{\gamma^+}\bigr)\bigr\|_{L^2(\gamma_T^-)}
	\le C\sum_{|H|\le m'-1}\bigl\|(f_{H\mathrm t^j})_{\gamma^+}\bigr\|_{L^2(\gamma_T^+)},
	\qquad \varnothing\ne S\subset[m'],
	\]
	and hence \eqref{eq:dyn-mixed-boundary-bound}, since
	$|H\mathrm t^j|<|I|$ and $j(H\mathrm t^j)=j$.  The constants of
	Lemma~\ref{lem:spectral-three-factor-coloring} depend on neither $N$ nor $\mathcal T$, and depend on $a$ only through $c_a$, which does not depend on $\delta$ by \eqref{eq:dyn-theta-defect-time-bound}. Lemma~\ref{lem:spectral-three-factor-coloring}\,(\ref{it:three-factor-spacetime}) gives $C_a\le C$ for $a=\mathfrak m_\circ+d_\delta(t)$ and $C_a\le C\delta$ for $a=\partial_t^kd_\delta(t)$, with $C$ independent of $\delta$. Thus the constants are independent of $T$, $N$ and $\delta$.
\end{proof}

\subsection{The complete flow}\label{app:diffuse-lifts}

For a tangential field $X$ we call $\widehat\Phi_s^X$ of \eqref{eq:complete-phase-flow} the \emph{complete flow} of $X$.  We prove global existence of the complete flows, preservation of $\gamma^\pm$ and $\gamma^0$, weighted $L_x^pL_v^2$, $H_\sigma^1$, $H_\sigma^{-1}$ and trace bounds for their pullbacks, and compatibility of the pulled-back traces.

\subsubsection{The flow and its pullbacks}\label{ss:flow-uniform-window}

\begin{lemma}
	\label{lem:complete-flow}
	Let $\nu\in\mathfrak A$, and let $\Phi_s^\nu$ and $\widehat\Phi_s^\nu$ be the
	flows \eqref{eq:complete-phase-flow} of $\tau_\nu$.  Then $\Phi_s^\nu$ and
	$\widehat\Phi_s^\nu$ are defined for every $s\in\bR$.  For every $s\in\bR$
	and $u\in C^1(\overline\cD\times\bR^3)$,
	\begin{equation}
		\label{eq:complete-flow-level-set}
		\rho_\cD(\Phi_s^\nu(x))=\rho_\cD(x),
		\qquad
		\frac{\dd}{\dd s}
		u(\widehat\Phi_s^\nu(x,v))
		=\mathfrak D_\nu u(\widehat\Phi_s^\nu(x,v)),
	\end{equation}
	where $x\in\bR^3$ in the first identity and
	$(x,v)\in\overline\cD\times\bR^3$ in the second.
	If $x\in\cN_\cD$, then $\Phi_s^\nu(x)\in\cN_\cD$ for every $s\in\bR$, and,
	for every $v\in\bR^3$,
	\begin{equation}
		\label{eq:complete-flow-normal-velocity-preservation}
		\bar n(\Phi_s^\nu(x))\cdot D\Phi_s^\nu(x)v
		=\bar n(x)\cdot v.
	\end{equation}
	In particular, $n_{\Phi_s^\nu(x)}\cdot D\Phi_s^\nu(x)v=n_x\cdot v$ for
	$x\in\p\cD$ and $v\in\bR^3$.
	With $J_s^\nu:=\det D\Phi_s^\nu$,
	\begin{equation}
		\label{eq:complete-flow-jacobians}
		|\det D\widehat\Phi_s^\nu|=|J_s^\nu|^2,
		\qquad
		|n_y\cdot w|\,\dd S_y\dd w
		=|J_s^\nu(x)|^2|n_x\cdot v|\,\dd S_x\dd v,
	\end{equation}
	where, in the second identity, $x\in\p\cD$, $v\in\bR^3$ and
	$(y,w):=\widehat\Phi_s^\nu(x,v)$.
	Moreover, there are $s_0>0$ and $C\ge1$, independent of $\nu\in\mathfrak A$,
	such that
	\begin{equation}
		\label{eq:complete-flow-uniform-bounds}
		\sup_{|s|\le2s_0,\;x\in\overline\cD}
		\Bigl(\|D\Phi_s^\nu(x)\|_{\rm op}+\|(D\Phi_s^\nu(x))^{-1}\|_{\rm op}
		+|J_s^\nu(x)|+|J_s^\nu(x)|^{-1}\Bigr)\le C.
	\end{equation}
\end{lemma}

\begin{proof}
	The field $\tau_\nu$ is compactly supported and $C^1$ by
	\eqref{eq:tangential-field-family}--\eqref{eq:tangential-field-properties}; hence both flows exist for all $s\in\bR$.  Moreover,
	$\tau_\nu\cdot\nabla\rho_\cD=0$ on $\cN_\cD$.  Since $\tau_\nu$ vanishes outside $\cN_\cD$,
	$\tau_\nu\cdot\nabla\rho_\cD=0$ on $\bR^3$, and the first identity in
	\eqref{eq:complete-flow-level-set} follows from the flow equation.  In
	particular $|\rho_\cD(\Phi_s^\nu(x))|=|\rho_\cD(x)|<r_\cD$ for
	$x\in\cN_\cD$, so $\Phi_s^\nu(x)\in\cN_\cD$ for $x\in\cN_\cD$ and
	$s\in\bR$.
	By \eqref{eq:signed-distance}, $\overline\cD=\{\rho_\cD\le0\}$, so the first identity also shows that $\Phi_s^\nu$ maps $\overline\cD$ into itself.  Differentiating the tangent map in $s$ shows that $\widehat\Phi_s^\nu$ is the flow of the field $\mathfrak f_{\tau_\nu}$ of \eqref{eq:complete-phase-space-derivative}. Since $\mathfrak D_\nu=\mathfrak f_{\tau_\nu}\cdot\nabla_{x,v}$ and $\widehat\Phi_s^\nu$ maps $\overline\cD\times\bR^3$ and $\cN_\cD\times\bR^3$ into themselves, the chain rule gives the second identity in
	\eqref{eq:complete-flow-level-set}, and the same identity for $u\in C^1(\cN_\cD\times\bR^3)$ and $(x,v)\in\cN_\cD\times\bR^3$.  Furthermore,
	$\mathfrak D_\nu(\bar n\cdot v)=0$ on $\cN_\cD\times\bR^3$ by
	\eqref{eq:complete-derivative-preserves-grazing}, so the identity on $\cN_\cD\times\bR^3$ with $u=\bar n\cdot v$ shows that $\bar n(\Phi_s^\nu(x))\cdot D\Phi_s^\nu(x)v$ does not depend on $s$.  This proves
	\eqref{eq:complete-flow-normal-velocity-preservation}.
	The derivative of $\widehat\Phi_s^\nu$ is block triangular with diagonal
	blocks $D\Phi_s^\nu$, so its phase Jacobian is $|J_s^\nu|^2$.  Since
	$\rho_\cD\circ\Phi_s^\nu=\rho_\cD$, differentiation gives, for $x\in\p\cD$,
	\[
	(D\Phi_s^\nu(x))^Tn_{\Phi_s^\nu(x)}=n_x.
	\]
	Let $v\in\bR^3$ and $(y,w)=\widehat\Phi_s^\nu(x,v)$.  By the change of variables formula for surface measure and the above identity,
	\[
	\dd S_y=|J_s^\nu(x)|\left|(D\Phi_s^\nu(x))^{-T}n_x\right|\dd S_x
	=|J_s^\nu(x)|\,\dd S_x .
	\]
	Moreover, $\dd w=|J_s^\nu(x)|\,\dd v$ since $w=D\Phi_s^\nu(x)v$, and $|n_y\cdot w|=|n_x\cdot v|$ by \eqref{eq:complete-flow-normal-velocity-preservation}.  Multiplying these
	three identities proves \eqref{eq:complete-flow-jacobians}.
	
	The variational equation and Liouville's formula give
	\[
	\frac{\dd}{\dd s}D\Phi_s^\nu=(\nabla\tau_\nu)(\Phi_s^\nu)\,D\Phi_s^\nu,
	\qquad
	J_s^\nu=\exp\int_0^s(\operatorname{div}\tau_\nu)(\Phi_r^\nu)\dd r.
	\]
	Gronwall's inequality for the forward and inverse flows yields
	\eqref{eq:complete-flow-uniform-bounds}, uniformly over the finite family $\mathfrak A$.
\end{proof}

For a tangential field $X$ write $A_s^X(x):=D\Phi_s^X(x)$ and
$J_s^X:=\det A_s^X$, and set $A_s^\nu:=A_s^{\tau_\nu}$; then \eqref{eq:complete-flow-uniform-bounds} and \eqref{eq:diffuse-part-weight} give
\begin{equation}
	\label{eq:flow-weight-comparison}
	C^{-1}w(v)\le w(A_s^\nu(x)v)\le Cw(v),
	\qquad |s|\le 2s_0 ,\quad x\in\overline\cD,
\end{equation}
with the constant $C$ of \eqref{eq:complete-flow-uniform-bounds}.
By \eqref{eq:complete-flow-level-set} and
\eqref{eq:complete-flow-normal-velocity-preservation},
$\widehat\Phi_s^\nu(\gamma^\pm)=\gamma^\pm$ and
$\widehat\Phi_s^\nu(\gamma^0)=\gamma^0$.  For a function $H$ on $\gamma^\pm$ define the boundary pullback
\begin{equation}
	\label{eq:boundary-flow-pullback}
	T_s^{\nu,\gamma_\pm}H:=H\circ\widehat\Phi_s^\nu\big|_{\gamma^\pm},
\end{equation}
and write $T_s^{X,\gamma_\pm}:=T_s^{\nu,\gamma_\pm}$ for $X=\tau_\nu$.

By \eqref{eq:complete-flow-jacobians}, for $X=\tau_\nu$ with $\nu\in\mathfrak A$, the $L^2(\dd x\dd v)$ adjoint of
\eqref{eq:bulk-flow-pullback} and the distributional pullback are
\begin{equation}
	\label{eq:bulk-flow-adjoint-quotient}
	(T_s^X)^\dagger\psi
	=|J_{-s}^X|^2T_{-s}^X\psi,
	\qquad
	\langle T_s^XQ,\psi\rangle
	:=\langle Q,(T_s^X)^\dagger\psi\rangle.
\end{equation}

For $X=\tau_\nu$ and $u\in C^1(\overline\cD\times\bR^3)$, the chain rule gives $\partial_sT_s^Xu=T_s^X\mathfrak D_Xu$. The pullback remains $C^1$, and the flow group law yields 
\[
\mathfrak D_XT_s^Xu=\partial_r|_{r=0}T_{s+r}^Xu=\partial_sT_s^Xu .
\] Hence
\begin{equation}
	\label{eq:r11-sd-flow-commutation}
	\mathfrak D_XT_s^Xu=\partial_sT_s^Xu=T_s^X\mathfrak D_Xu.
\end{equation}

By \eqref{eq:real-weight-Hsigma-equivalence},
\eqref{eq:coulomb-sigmavv-comparison} and \eqref{eq:local-H1-H_sigma},
for every $\mathrm L\in\bR$,
\begin{equation}
	\label{eq:flow-weighted-coercive-equivalence}
	\|w^{\mathrm L}q\|_{H_\sigma^1}^2
	\simeq_{\mathrm L}
	\int_{\bR^3}w^{2\mathrm L}
	\left(\sigma^{ij}\partial_i q\,\partial_j q+w^{-1}|q|^2\right)\dd v.
\end{equation}

\begin{lemma}
	\label{lem:complete-flow-weighted-bounds}
	Let $\nu\in\mathfrak A$, let $T_s:=T_s^\nu$ be the complete pullback
	\eqref{eq:bulk-flow-pullback}, and let $s_0$ be as in
	\eqref{eq:complete-flow-uniform-bounds}.
	\begin{enumerate}
		\item\label{it:complete-flow-fiber}
		Let $d\ge1$ and $\lambda\in\bR$, and let
		$H_\lambda:=L_v^2(\bR^3;\bR^d;w^{2\lambda}\dd v)$.  For every
		$1\le p\le\infty$ and every $|s|\le s_0$,
		\begin{equation}
			\label{eq:complete-flow-fiber-bound}
			\|T_sG\|_{L_x^pH_\lambda}\le C\|G\|_{L_x^pH_\lambda}
			\qquad(G\in L_x^pH_\lambda),
		\end{equation}
		with $C$ independent of $s$ and $p$; for $G$ on $(0,T)\times\cD\times\bR^3$ and $(T_sG)(t):=T_s\bigl(G(t)\bigr)$, the same bound holds with an
		outer $L_t^\rho$, $1\le\rho\le\infty$.  If $1\le p<\infty$, then
		\begin{equation}\label{eq:lifted-pullback-strong-continuity}
			\lim_{s\to0}\|T_sG-G\|_{L_x^pH_\lambda}=0
			\qquad(G\in L_x^pH_\lambda).
		\end{equation}
		In particular, for every $L\in\bR$ and every $q$ with
		$w^Lq\in L^2(\cD\times\bR_v^3)$,
		\begin{equation}
			\label{eq:complete-flow-L2-bound}
			\|w^LT_sq\|_{L^2(\cD\times\bR_v^3)}
			\le C\|w^Lq\|_{L^2(\cD\times\bR_v^3)}
			\qquad(|s|\le s_0),
		\end{equation}
		where $C$ is independent of $s$, and
		\begin{equation}
			\label{eq:complete-flow-strong-continuity}
			\lim_{s\to0}\|w^L(T_sq-q)\|_{L^2(\cD\times\bR_v^3)}=0.
		\end{equation}
		\item\label{it:complete-flow-landau}
		On $L_x^2H_\sigma^{-1}$,
		$T_s$ is the distributional pullback \eqref{eq:bulk-flow-adjoint-quotient}.
		For every $L\ge0$, every $1\le p\le\infty$ and every $|s|\le s_0$,
		\begin{equation}
			\label{eq:complete-flow-H1-bound}
			\|w^LT_sq\|_{L_x^pH_\sigma^1}
			\le C\|w^{L+1}q\|_{L_x^pH_\sigma^1},
		\end{equation}
		and
		\begin{equation}
			\label{eq:complete-flow-Hminus1-bound}
			\|w^LT_sQ\|_{L_x^2H_\sigma^{-1}}
			\le C\|w^{L+1}Q\|_{L_x^2H_\sigma^{-1}},
		\end{equation}
		for all $q$ and $Q$ with $w^{L+1}q\in L_x^pH_\sigma^1$ and $w^{L+1}Q\in L_x^2H_\sigma^{-1}$, where $C$ is
		independent of $s$.
		\item\label{it:complete-flow-time}
		On spacetime distributions of order at most $\mathrm N$,
		\begin{equation}
			[\partial_t,T_s]=0.
			\label{eq:complete-flow-time-commutation}
		\end{equation}
	\end{enumerate}
\end{lemma}

\begin{proof}
	Abbreviate $\Phi_s=\Phi_s^\nu$, $A_s=A_s^\nu$ and
	$\widehat\Phi_s=\widehat\Phi_s^\nu$.
	
	\emph{(\ref{it:complete-flow-fiber}).}
	By the first identity in \eqref{eq:complete-flow-level-set} and \eqref{eq:signed-distance}, $\Phi_s$ maps $\cD=\{\rho_\cD<0\}$ diffeomorphically onto
	itself.  Fix $x\in\cD$.  The change of variables
	$z=A_s(x)v$ in $v$, \eqref{eq:flow-weight-comparison} raised to the power $\lambda$, and
	\eqref{eq:complete-flow-uniform-bounds} for $|\det A_s(x)|^{-1}$ give
	\[
	\bigl\|w^\lambda(v)\,G\bigl(\Phi_s(x),A_s(x)v\bigr)\bigr\|_{L_v^2}
	\le C_1\,\Gamma\bigl(\Phi_s(x)\bigr),
	\qquad
	\Gamma(y):=\|w^\lambda G(y,\cdot)\|_{L_v^2}.
	\]
	For $p<\infty$, the change of variables $y=\Phi_s(x)$ and
	\eqref{eq:complete-flow-uniform-bounds} give
	\[
	\|\Gamma\circ\Phi_s\|_{L_x^p}\le C_0^{1/p}\|\Gamma\|_{L_x^p},
	\]
	where $C_0$ is the constant in \eqref{eq:complete-flow-uniform-bounds}.  At
	$p=\infty$, $\|\Gamma\circ\Phi_s\|_{L_x^\infty}=\|\Gamma\|_{L_x^\infty}$.  This
	gives \eqref{eq:complete-flow-fiber-bound}.
	The spacetime bound follows by applying the same estimate at each time.
	
	If $G\in C_c(\cD\times\bR^3)$, then
	$\bigcup_{|s|\le s_0}\widehat\Phi_{-s}(\operatorname{supp}G)$ is compact. Outside this set, $T_sG-G=0$ for $|s|\le s_0$. On it, $\Phi_sx\to x$ and
	$A_s(x)v\to v$ uniformly.  Moreover, $|T_sG-G|\le2\|G\|_{L^\infty(\cD\times\bR^3)}$ for $|s|\le s_0$, and $w^{2\lambda}$ is bounded on this compact set by a constant depending on $\lambda$ and $\operatorname{supp}G$.
	Dominated convergence gives \eqref{eq:lifted-pullback-strong-continuity}
	for such $G$.  Density of $C_c(\cD\times\bR^3)$ in $L_x^pH_\lambda$, for $p<\infty$, and the uniform bound \eqref{eq:complete-flow-fiber-bound} extend this convergence to every $G\in L_x^pH_\lambda$.
	
	With $p=2$, $d=1$ and $\lambda=L$, \eqref{eq:complete-flow-fiber-bound} and
	\eqref{eq:lifted-pullback-strong-continuity} give
	\eqref{eq:complete-flow-L2-bound} and \eqref{eq:complete-flow-strong-continuity}.
	
	\emph{(\ref{it:complete-flow-landau}).}
	For $x\in\overline\cD$, $|s|\le s_0$ and $v,\xi\in\bR^3$, let $z=A_s(x)v$.  The decomposition
	\eqref{eq:coulomb-eigen-decomposition} and the bounds
	\eqref{eq:eigenvalue-size} imply
	\begin{equation}
		\label{eq:flow-sigma-ellipticity}
		\sigma(v)\le Cw(v)^{-1}I,
		\qquad
		\sigma(z)\ge c w(z)^{-3}I.
	\end{equation}
	Since \eqref{eq:complete-flow-uniform-bounds} gives
	$|A_s^T\xi|\le C|\xi|$, \eqref{eq:flow-weight-comparison} and
	\eqref{eq:flow-sigma-ellipticity} yield
	\begin{equation}
		\label{eq:complete-flow-sigma-comparison}
		(A_s^T\xi)^T\sigma(v)(A_s^T\xi)
		\le Cw(z)^2\,\xi^T\sigma(z)\xi.
	\end{equation}
	Fix $x\in\cD$ and write $A:=A_s(x)$.  Let $f\in C_c^\infty(\bR^3)$.  By \eqref{eq:flow-weighted-coercive-equivalence}, \eqref{eq:flow-weight-comparison}, \eqref{eq:complete-flow-sigma-comparison} and the change of variables $z=Av$, for which $|\det A|^{-1}\le C$ by
	\eqref{eq:complete-flow-uniform-bounds},
	\[
	\begin{aligned}
		\|w^Lf(A\,\cdot)\|_{H_\sigma^1}^2
		&\lesssim\int_{\bR^3}w(v)^{2L}\Bigl((\nabla f)(Av)^TA\,\sigma(v)A^T(\nabla f)(Av)+w(v)^{-1}|f(Av)|^2\Bigr)\dd v\\
		&\lesssim\int_{\bR^3}w(z)^{2L}\Bigl(w(z)^2(\nabla f)(z)^T\sigma(z)(\nabla f)(z)+w(z)^{-1}|f(z)|^2\Bigr)\dd z
		\lesssim\|w^{L+1}f\|_{H_\sigma^1}^2 .
	\end{aligned}
	\]
	Hence
	\begin{equation}
		\label{eq:flow-Hsigma-pointwise}
		\|w^Lf(A\,\cdot)\|_{H_\sigma^1}
		\le C\|w^{L+1}f\|_{H_\sigma^1},
	\end{equation}
	where $C$ is independent of $x$ and $s$.  By density (Lemma~\ref{lem:H_sigma-Hilbert}\,(1)) and \eqref{eq:local-H1-H_sigma}, the bound \eqref{eq:flow-Hsigma-pointwise} extends to every $f$ with $w^{L+1}f\in H_\sigma^1$.
	
	Let $w^{L+1}q\in L_x^pH_\sigma^1$.  Since $\Phi_s$ is a diffeomorphism of $\cD$,
	$w^{L+1}q(\Phi_s(x),\cdot)\in H_\sigma^1$ for almost every $x$, and
	\eqref{eq:flow-Hsigma-pointwise} with $f=q(\Phi_s(x),\cdot)$ gives
	\[
	\|w^LT_sq(x,\cdot)\|_{H_\sigma^1}
	\le C\|w^{L+1}q(\Phi_s(x),\cdot)\|_{H_\sigma^1}
	\qquad\text{for almost every }x\in\cD.
	\]
	The $L_x^p$ norm and the change of variables $y=\Phi_s(x)$ with
	\eqref{eq:complete-flow-uniform-bounds} give \eqref{eq:complete-flow-H1-bound}.
	
	By Remark~\ref{rem:Riesz-H_sigma} and the density of
	$C_c^\infty(\cD\times\bR^3)$ in $L_x^2H_\sigma^1$ (Lemma~\ref{lem:H_sigma-Hilbert}\,(1)), $C_c^\infty(\cD\times\bR^3)$
	is dense in $L_x^2H_\sigma^{-1}$.
	
	For $R\in C_c^\infty(\cD\times\bR^3)$ set $Q=w^{-L-1}R$.  For
	$\varphi\in C_c^\infty(\cD\times\bR^3)$, the change of variables $(y,z)=\widehat\Phi_s(x,v)$ and \eqref{eq:complete-flow-jacobians} give
	\begin{equation}
		\label{eq:flow-negative-duality}
		\langle w^LT_sQ,\varphi\rangle
		=\langle R,r_{s,L}T_{-s}\varphi\rangle,
		\qquad
		r_{s,L}(y,z):=|\det A_{-s}(y)|^2
		w(z)^{-L-1}w(A_{-s}(y)z)^L.
	\end{equation}
	Since
	\[
	\nabla_zr_{s,L}(y,z)
	=r_{s,L}(y,z)\Bigl(-(L+1)w(z)^{-2}z+L\,w(A_{-s}(y)z)^{-2}A_{-s}(y)^TA_{-s}(y)z\Bigr),
	\]
	\eqref{eq:complete-flow-uniform-bounds} and \eqref{eq:flow-weight-comparison} give
	\begin{equation}
		\label{eq:flow-negative-multiplier}
		|r_{s,L}(y,z)|\le Cw(z)^{-1},
		\qquad
		|\nabla_zr_{s,L}(y,z)|\le Cw(z)^{-2}.
	\end{equation}
	Let $\widetilde z=A_{-s}(y)z$ and
	$\xi=(\nabla_v\varphi)(\widehat\Phi_{-s}(y,z))$.  By the product rule,
	\[
	\nabla_z(r_{s,L}T_{-s}\varphi)=(\nabla_zr_{s,L})T_{-s}\varphi+r_{s,L}A_{-s}^T\xi.
	\]
	Hence, by
	\eqref{eq:flow-sigma-ellipticity},
	\eqref{eq:flow-negative-multiplier},
	\eqref{eq:complete-flow-sigma-comparison} with $(-s,y,z,\widetilde z)$ in place of $(s,x,v,z)$, and
	\eqref{eq:flow-weight-comparison} for $w(\widetilde z)\lesssim w(z)$,
	\begin{align}
		r_{s,L}^2(A_{-s}^T\xi)^T\sigma(z)(A_{-s}^T\xi)
		&\lesssim w(z)^{-2}w(\widetilde z)^2\,\xi^T\sigma(\widetilde z)\xi
		\lesssim \xi^T\sigma(\widetilde z)\xi,
		\label{261016601}\\
		(\nabla_zr_{s,L})^T\sigma(z)(\nabla_zr_{s,L})
		&\lesssim w(z)^{-5}\le w(z)^{-1},
		\label{261016602}
	\end{align}
	and the zeroth term satisfies
	\begin{equation}
		\label{261016603}
		w(z)^{-1}r_{s,L}^2\lesssim w(z)^{-3}\le w(z)^{-1}.
	\end{equation}
	We use the change of variables
	\begin{equation}
		\label{261016701}
		(x,v)=\widehat\Phi_{-s}(y,z),
	\end{equation}
	under which $\widetilde z=A_{-s}(y)z=v$ and $\xi=(\nabla_v\varphi)(x,v)$ by \eqref{eq:complete-phase-flow}.  The first inequality below follows from \eqref{eq:flow-weighted-coercive-equivalence} at $\mathrm L=0$, the product rule and \eqref{261016601}--\eqref{261016603}; the second follows from the change of variables \eqref{261016701}, with \eqref{eq:complete-flow-jacobians}, \eqref{eq:complete-flow-uniform-bounds} and \eqref{eq:flow-weight-comparison} for its Jacobian and for $w(z)^{-1}\lesssim w(v)^{-1}$:
	\[
	\begin{aligned}
		\|r_{s,L}T_{-s}\varphi\|_{L_x^2H_\sigma^1}^2
		&\lesssim\int_\cD\int_{\bR^3}\Bigl(\xi^T\sigma(\widetilde z)\xi+w(z)^{-1}|T_{-s}\varphi|^2\Bigr)\dd z\dd y\\
		&\lesssim\int_\cD\int_{\bR^3}\Bigl((\nabla_v\varphi)^T\sigma(v)\nabla_v\varphi+w(v)^{-1}|\varphi|^2\Bigr)\dd v\dd x .
	\end{aligned}
	\]
	With \eqref{eq:flow-weighted-coercive-equivalence} at $\mathrm L=0$ on the right, this gives
	\begin{equation}
		\label{eq:flow-negative-test-bound}
		\|r_{s,L}T_{-s}\varphi\|_{L_x^2H_\sigma^1}
		\le C\|\varphi\|_{L_x^2H_\sigma^1}.
	\end{equation}
	Hence
	\eqref{eq:flow-negative-duality} gives, for every $R\in C_c^\infty(\cD\times\bR^3)$,
	\[
	\|w^LT_s(w^{-L-1}R)\|_{L_x^2H_\sigma^{-1}}
	\le C\|R\|_{L_x^2H_\sigma^{-1}}.
	\]
	Density extends this estimate to every $R\in L_x^2H_\sigma^{-1}$. The extension agrees with the distributional pullback \eqref{eq:bulk-flow-adjoint-quotient}, since \eqref{eq:flow-negative-duality} and \eqref{eq:flow-negative-test-bound} identify its action on test functions.  Taking $R=w^{L+1}Q$ proves
	\eqref{eq:complete-flow-Hminus1-bound}.
	
	\emph{(\ref{it:complete-flow-time}).}
	Since $\widehat\Phi_s$ is independent of time,
	direct pullback on test functions gives
	\eqref{eq:complete-flow-time-commutation}.
\end{proof}

\begin{lemma}[Flow slicing]
	\label{lem:distributional-flow-slicing}
	Let $X=\tau_\nu$ with $\nu\in\mathfrak A$, let $T_s:=T_s^X$, and let $s_0$
	be as in \eqref{eq:complete-flow-uniform-bounds}.  Let $d\ge1$ and
	$\lambda\in\bR$, let $H_\lambda$ be as in
	Lemma~\ref{lem:complete-flow-weighted-bounds}\,(\ref{it:complete-flow-fiber}), let
	$1\le m\le\mathrm N$, and let $G_0,\ldots,G_m\in L_x^1H_\lambda$ satisfy
	\[
	\mathfrak D_XG_k=G_{k+1}
	\quad\text{in }\cD'(\cD\times\bR^3),
	\qquad 0\le k<m,
	\]
	componentwise. Then, for almost every $x\in\cD$, the curve
	$F_x:(-s_0,s_0)\to H_\lambda$,
	\[
	F_x(s):=T_sG_0(x,\cdot),
	\]
	has a representative in $W^{m,1}((-s_0,s_0);H_\lambda)$, and
	\begin{equation}\label{eq:slicing-derivatives}
		F_x^{(k)}(s)=T_sG_k(x,\cdot)
		\quad\text{for a.e. }s,
		\qquad 0\le k\le m.
	\end{equation}
	For almost every $x\in\cD$ and $0\le k\le m-1$, the continuous
	representative of $F_x^{(k)}$ satisfies
	\begin{equation}\label{eq:slicing-zero-identification}
		F_x^{(k)}(0)=G_k(x,\cdot)
		\quad\text{in }H_\lambda.
	\end{equation}
\end{lemma}

\begin{proof}
	The phase vector field $\mathfrak f_X$ has flow $\widehat\Phi_s$ and divergence given by \eqref{eq:phase-divergence}. Since $G_k,G_{k+1}\in L^1_{\rm loc}$ and $\mathfrak D_XG_k=G_{k+1}$ distributionally, mollification gives, for $0\le k<m$,
	\begin{equation}\label{eq:C1-test-extension}
		\int G_{k+1}\psi\dd z=-\int G_k\,\operatorname{div}(\mathfrak f_X\psi)\dd z
		\qquad\text{for every }\psi\in C_c^1(\cD\times\bR^3).
	\end{equation}
	
	For a scalar
	component $G$ and $\varphi\in C_c^\infty(\cD\times\bR^3)$, \eqref{eq:bulk-flow-adjoint-quotient} and $\det D\widehat\Phi_{-s}=(J_{-s}^X)^2$ give
	\begin{equation}\label{eq:pullback-duality}
		\langle T_sG,\varphi\rangle
		=\int G(z)\varphi(\widehat\Phi_{-s}z)\,
		\bigl(\det D\widehat\Phi_{-s}\bigr)(z)\dd z.
	\end{equation}
	The factor $K_s$ multiplying $G$ in \eqref{eq:pullback-duality} satisfies the continuity equation
	\begin{equation}\label{eq:adjoint-flow-factor}
		\partial_sK_s+\operatorname{div}_z(\mathfrak f_XK_s)=0
		\quad\text{in }\cD'(\cD\times\bR^3),
		\qquad
		K_s(z):=\varphi(\widehat\Phi_{-s}z)\bigl(\det D\widehat\Phi_{-s}\bigr)(z).
	\end{equation}
	To see this, change variables $z=\widehat\Phi_sy$ against $\psi\in C_c^\infty(\cD\times\bR^3)$:
	\[
	\int K_s(z)\psi(z)\dd z
	=\int\varphi(y)\psi(\widehat\Phi_sy)\dd y,
	\]
	and differentiation in $s$ gives
	\[
	\frac{\dd}{\dd s}\int K_s\psi
	=\int\varphi(y)(\mathfrak f_X\cdot\nabla\psi)(\widehat\Phi_sy)\dd y
	=\int K_s(z)\,\mathfrak f_X\cdot\nabla\psi(z)\dd z,
	\]
	which is \eqref{eq:adjoint-flow-factor}.  Since $K_s=(T_s^X)^\dagger\varphi$, it lies in $C_c^1(\cD\times\bR^3)$ by the mapping property of the adjoints of the complete pullbacks stated at the beginning of \S\ref{sec:diffuse-complete-equations}, and
	\eqref{eq:C1-test-extension} applies with $\psi=K_s$.  For
	$\alpha\in C_c^\infty((-s_0,s_0))$ and $\varphi\in C_c^\infty(\cD\times\bR^3)$,
	Fubini's theorem, \eqref{eq:pullback-duality}, \eqref{eq:adjoint-flow-factor} and
	\eqref{eq:C1-test-extension} give
	\begin{align*}
		&-\int\alpha'(s)\langle T_sG_k,\varphi\rangle\,\dd s
		= -\int\alpha'(s)\langle G_k,K_s\rangle\,\dd s\\
		=\,& \int\alpha(s)\langle G_k,\partial_sK_s\rangle\,\dd s= -\int\alpha(s)\langle G_k,
		\operatorname{div}(\mathfrak f_XK_s)\rangle\,\dd s\\
		=\,& \int\alpha(s)\langle G_{k+1},K_s\rangle\,\dd s= \int\alpha(s)\langle T_sG_{k+1},\varphi\rangle\,\dd s,
	\end{align*}
	where all pairings are absolutely integrable because
	$\bigcup_{|s|\le s_0}\operatorname{supp}K_s$ is a compact subset of
	$\cD\times\bR^3$.  Hence
	\begin{equation}\label{eq:global-weak-flow-chain}
		\partial_s(T_sG_k)=T_sG_{k+1}
		\quad\text{in }\cD'((-s_0,s_0)\times\cD\times\bR^3).
	\end{equation}
	
	Set $I:=(-s_0,s_0)$ and $F_k(s,x):=T_sG_k(x,\cdot)$ for $0\le k\le m$.  By
	Lemma~\ref{lem:complete-flow-weighted-bounds}\,(\ref{it:complete-flow-fiber})
	with $p=1$, $F_k\in L^1(I\times\cD;H_\lambda)$.  Testing
	\eqref{eq:global-weak-flow-chain} with $\alpha(s)\chi(x)w^{2\lambda}(v)\zeta(v)$,
	where $\alpha\in C_c^\infty(I)$, $\chi\in C_c^\infty(\cD)$ and
	$\zeta\in C_c^\infty(\bR^3;\bR^d)$, gives
	$\partial_sF_k=F_{k+1}$ in $\cD'(I\times\cD;H_\lambda)$ for $0\le k<m$.  Since
	$F_k,\partial_sF_k\in L^1(I\times\cD;H_\lambda)$, for almost every $x$ the slice
	$F_k(\cdot,x)$ belongs to $W^{1,1}(I;H_\lambda)$ with derivative
	$F_{k+1}(\cdot,x)$.  Iterating over $k$ proves \eqref{eq:slicing-derivatives} and $F_x\in W^{m,1}(I;H_\lambda)$.
	
	For $0\le k\le m-1$, the continuous representative of $F_k(\cdot,x)$ takes at
	$s=0$ the value $\lim_{r\downarrow0}\fint_{-r}^rT_sG_k(x,\cdot)\dd s$ in $H_\lambda$.  Since
	\[
	\Bigl\|\fint_{-r}^rT_sG_k\dd s-G_k\Bigr\|_{L_x^1H_\lambda}
	\le\sup_{|s|\le r}\|T_sG_k-G_k\|_{L_x^1H_\lambda}
	\]
	and the right-hand side tends to zero as $r\downarrow0$ by
	\eqref{eq:lifted-pullback-strong-continuity} with $p=1$, these averages converge
	to $G_k$ in $L_x^1H_\lambda$. Hence, along a sequence $r_n\downarrow0$,
	they converge to $G_k(x,\cdot)$ in $H_\lambda$ for almost every $x$.  This proves
	\eqref{eq:slicing-zero-identification}.
\end{proof}

\begin{lemma}[Complete-flow difference quotient]
	\label{lem:strong-complete-flow-generator-L2}
	Let $\nu\in\mathfrak A$, let $s_0$ be as in
	\eqref{eq:complete-flow-uniform-bounds}, and let
	$q,z\in L^2(\cD\times\bR^3)$ satisfy
	\[
	\mathfrak D_\nu q=z\qquad\text{in }\cD'(\cD\times\bR^3).
	\]
	Then, for every $s\in(-s_0,s_0)$,
	\begin{equation}
		\label{eq:strong-complete-flow-generator-FTC}
		T_s^\nu q-q
		=\int_0^sT_t^\nu z\dd t
		\qquad\text{in }L^2(\cD\times\bR_v^3),
	\end{equation}
	and
	\begin{equation}
		\label{eq:strong-complete-flow-generator-limit}
		\frac{T_s^\nu q-q}{s}\longrightarrow z
		\qquad\text{strongly in }L^2(\cD\times\bR_v^3)\quad(s\to0).
	\end{equation}
	Let $\ell\in\bR$.  If $q,z\in L^2(\cD\times\bR_v^3,w^{2\ell}\dd x\dd v)$
	and $\mathfrak D_\nu q=z$ in $\cD'(\cD\times\bR^3)$, then
	\eqref{eq:strong-complete-flow-generator-FTC} and
	\eqref{eq:strong-complete-flow-generator-limit} hold with
	$L^2(\cD\times\bR_v^3)$ replaced by
	$L^2(\cD\times\bR_v^3,w^{2\ell}\dd x\dd v)$.
\end{lemma}

\begin{proof}
	Lemma~\ref{lem:distributional-flow-slicing} applies with 
	$$
	X=\tau_\nu,\quad d=1,\quad \lambda=0,\quad m=1,\quad G_0=q,\quad G_1=z.
	$$
	By \eqref{eq:slicing-derivatives}
	and \eqref{eq:slicing-zero-identification}, for almost every $x$,
	\[
	T_s^\nu q(x,\cdot)-q(x,\cdot)
	=\int_0^sT_t^\nu z(x,\cdot)\dd t
	\quad\text{in }L_v^2
	\]
	for almost every $s$.  By \eqref{eq:complete-flow-L2-bound},
	\begin{equation}
		\label{eq:complete-flow-generator-integrability}
		\int_{I_s}\|T_t^\nu z\|_{L^2_{x,v}}\dd t
		\le C|s|\|z\|_{L^2_{x,v}},
		\qquad I_s=[\min\{0,s\},\max\{0,s\}].
	\end{equation}
	By \eqref{eq:complete-flow-generator-integrability} and Fubini's theorem,
	\eqref{eq:strong-complete-flow-generator-FTC} holds for almost every $s$.  For
	$s,s'\in(-s_0,s_0)$, $\widehat\Phi_{s-s'}^\nu\circ\widehat\Phi_{s'}^\nu=\widehat\Phi_s^\nu$
	gives $T_s^\nu q-T_{s'}^\nu q=T_{s'}^\nu(T_{s-s'}^\nu q-q)$. Hence
	\eqref{eq:complete-flow-L2-bound} and
	\eqref{eq:complete-flow-strong-continuity} at $L=0$ make the map
	$s\mapsto T_s^\nu q-q$ continuous from $(-s_0,s_0)$ into $L^2(\cD\times\bR_v^3)$. By
	\eqref{eq:complete-flow-generator-integrability}, so is $s\mapsto\int_0^sT_t^\nu z\dd t$.
	Hence \eqref{eq:strong-complete-flow-generator-FTC} holds for
	every $s\in(-s_0,s_0)$.
	
	For $s\ne0$, \eqref{eq:strong-complete-flow-generator-FTC} and
	\eqref{eq:complete-flow-strong-continuity} with $L=0$ give
	\[
	\left\|\frac{T_s^\nu q-q}{s}-z\right\|_{L^2}
	\le \frac1{|s|}\int_{I_s}\|T_t^\nu z-z\|_{L^2}\dd t
	\le \sup_{t\in I_s}\|T_t^\nu z-z\|_{L^2}
	\longrightarrow0.
	\]
	This is \eqref{eq:strong-complete-flow-generator-limit}.
	
	For the weighted assertion, the same argument applies in
	$L^2(\cD\times\bR_v^3,w^{2\ell}\dd x\dd v)$, with
	Lemma~\ref{lem:distributional-flow-slicing} at $\lambda=\ell$ and with
	\eqref{eq:complete-flow-L2-bound} and
	\eqref{eq:complete-flow-strong-continuity} at $L=\ell$ in place of $L=0$.
\end{proof}

\begin{lemma}[Coefficient continuity along complete flows]
	\label{lem:complete-flow-coefficient-continuity}
	Let $\nu\in\mathfrak A$, $X=\tau_\nu$, and let $T_s^X$ be the complete
	pullback \eqref{eq:bulk-flow-pullback}.  Let $B\ge0$ and
	$0<\alpha<1$.  If
	\begin{equation}
		\label{eq:complete-flow-coefficient-continuity-data}
		w^Bq\in L^\infty_{x,v},
		\qquad
		\mu^{1/4}q\in C^\alpha_{\kin}(\overline\cD\times\bR^3),
	\end{equation}
	then
	\begin{equation}
		\label{eq:complete-flow-coefficient-continuity}
		\omega_{X,q}(r):=
		\sup_{|s|\le r}
		\|\mu^{1/4}(T_s^Xq-q)\|_{L^\infty_{x,v}}
		\longrightarrow0
		\qquad(r\downarrow0).
	\end{equation}
	Moreover, for every $K<\infty$,
	\[
	\sup\Bigl\{\omega_{\tau_\nu,q}(r):\ \nu\in\mathfrak A,\
	\|w^Bq\|_{L^\infty}+\|\mu^{1/4}q\|_{C^\alpha_{\kin}}\le K\Bigr\}
	\longrightarrow0
	\qquad(r\downarrow0).
	\]
\end{lemma}

\begin{proof}
	Write $F:=\mu^{1/4}q$, and fix $x\in\cD$; the constants below do not depend on $x$.  Since $A_s^X=A_s^\nu$,
	\eqref{eq:flow-weight-comparison} gives, for every $R\ge1$,
	\begin{align}
		&\sup_{|s|\le 2s_0}\sup_{|v|\ge R}
		\mu^{1/4}(v)
		\bigl(|q(\Phi_s^X(x),A_s^X(x)v)|+|q(x,v)|\bigr)
		\le
		C\|w^Bq\|_{L^\infty}
		\sup_{|v|\ge R}\mu^{1/4}(v)w(v)^{-B}.
		\label{eq:complete-flow-coefficient-tail}
	\end{align}
	
	For a fixed $R\ge1$, with $C=C(R)$ from here on, the flow equation and \eqref{eq:complete-flow-uniform-bounds} give, for $|s|\le2s_0$ and $|v|\le R$,
	\[
	|\Phi_s^X(x)-x|\le C|s|,
	\qquad
	|A_s^X(x)v-v|\le C|s|.
	\]
	Bounding the minimum in \eqref{kineticdistance} by its value at the velocity $v$ gives
	\[
	\mathbf d\bigl((0,x,v),(0,\Phi_s^X(x),A_s^X(x)v)\bigr)
	\le\max\bigl(|\Phi_s^X(x)-x|^{1/3},|A_s^X(x)v-v|\bigr)
	\le\max\bigl(C|s|^{1/3},C|s|\bigr).
	\]
	Since $|s|\le2s_0$, it follows that
	\begin{equation}
		\label{eq:complete-flow-coefficient-distance}
		\mathbf d\bigl((0,x,v),(0,\Phi_s^X(x),A_s^X(x)v)\bigr)
		\le C|s|^{1/3},
		\qquad |s|\le 2s_0,\ |v|\le R.
	\end{equation}
	Moreover, the bound on $|A_s^X(x)v-v|$ above gives
	\begin{equation}
		\label{eq:complete-flow-coefficient-Maxwellian-ratio}
		\sup_{|v|\le R}
		\left|
		\frac{\mu^{1/4}(v)}{\mu^{1/4}(A_s^X(x)v)}-1
		\right|
		\le C|s|.
	\end{equation}
	We have
	\begin{align*}
		&\mu^{1/4}(v)\bigl(q(\Phi_s^X(x),A_s^X(x)v)-q(x,v)\bigr)\\
		&\qquad=\Bigl(\frac{\mu^{1/4}(v)}{\mu^{1/4}(A_s^X(x)v)}-1\Bigr)F(\Phi_s^X(x),A_s^X(x)v)
		+F(\Phi_s^X(x),A_s^X(x)v)-F(x,v).
	\end{align*}
	Hence, by \eqref{eq:complete-flow-coefficient-distance},
	\eqref{eq:complete-flow-coefficient-Maxwellian-ratio} and
	\eqref{eq:complete-flow-coefficient-continuity-data},
	\begin{align}
		&\sup_{|v|\le R}
		\mu^{1/4}(v)
		|q(\Phi_s^X(x),A_s^X(x)v)-q(x,v)|
		\le
		C\|F\|_{L^\infty}|s|
		+C[F]_{C^\alpha_{\kin}}|s|^{\alpha/3}.
		\label{eq:complete-flow-coefficient-local}
	\end{align}
	Given $\varepsilon>0$, choose $R$ so that
	\[
	C\|w^Bq\|_{L^\infty}\sup_{|v|\ge R}\mu^{1/4}(v)w(v)^{-B}\le\varepsilon,
	\]
	with $C$ the constant of \eqref{eq:complete-flow-coefficient-tail}.  Then
	\eqref{eq:complete-flow-coefficient-tail} and
	\eqref{eq:complete-flow-coefficient-local} give
	$\limsup_{r\downarrow0}\omega_{X,q}(r)\le\varepsilon$.  This proves
	\eqref{eq:complete-flow-coefficient-continuity}.  Both estimates depend on $q$ only through $\|w^Bq\|_{L^\infty}+\|F\|_{C^\alpha_{\kin}}$, which also gives the asserted uniformity.
\end{proof}

Recall $\cX_{\rm c}=\{\tau_\nu:\nu\in\mathfrak A\}$.  Unless stated otherwise, $X\in\cX_{\rm c}$ is fixed, and
\[
T_s:=T_s^X,
\qquad
\Phi_s:=\Phi_s^X,
\qquad
A_s:=A_s^X,
\qquad
\widehat\Phi_s(x,v)=(\Phi_s(x),A_s(x)v),
\]
with $T_s^X$ acting componentwise on a vector-valued function.

By \eqref{eq:flow-weighted-coercive-equivalence} at $\mathrm L=0$,
\begin{equation}
	\label{eq:auxiliary-coercive-equivalence}
	\|q\|_{H_\sigma^1}^2\simeq
	\int_{\bR^3}\bigl(\sigma^{ij}\partial_iq\,\partial_jq
	+w^{-1}|q|^2\bigr)\,\dd v.
\end{equation}

\subsubsection{Boundary traces}\label{ss:flow-boundary-traces}

\begin{lemma}[Boundary-flow continuity]
	\label{lem:boundary-flow-continuity}
	Let $\nu\in\mathfrak A$, let \(T_s^{\nu,\gamma_\pm}\) be given by
	\eqref{eq:boundary-flow-pullback}, and let $s_0$ be as in
	\eqref{eq:complete-flow-uniform-bounds}.  For every \(L\ge0\), every
	\(|s|\le s_0\) and every $H$ with $w^LH\in L^2(\gamma^\pm)$,
	\begin{equation}
		\label{eq:boundary-flow-bound}
		\|w^LT_s^{\nu,\gamma_\pm}H\|_{L^2(\gamma^\pm)}
		\le C\|w^LH\|_{L^2(\gamma^\pm)},
	\end{equation}
	with $C$ independent of $s$.
	Moreover,
	\[
	T_s^{\nu,\gamma_\pm}T_t^{\nu,\gamma_\pm}=T_{s+t}^{\nu,\gamma_\pm}
	\quad(|s|,|t|,|s+t|\le s_0),
	\qquad
	\|w^L(T_s^{\nu,\gamma_\pm}H-H)\|_{L^2(\gamma^\pm)}\longrightarrow0
	\quad(s\to0).
	\]
\end{lemma}

\begin{proof}
	For \(|s|\le s_0\), \eqref{eq:complete-flow-uniform-bounds}, the change of
	variables \eqref{eq:complete-flow-jacobians} on $\gamma$ and
	\eqref{eq:flow-weight-comparison} prove \eqref{eq:boundary-flow-bound}.
	Preservation of $\gamma^\pm$ and the flow group law give the identity for boundary pullbacks.
	If $H$ is the restriction to $\gamma^\pm$ of some
	$\varphi\in C_c^\infty(\bR_x^3\times\bR_v^3)$ with
	$\operatorname{supp}\varphi\cap(\partial\cD\times\bR^3)$ a compact subset of
	$\gamma^\pm$, then
	\[
	T_s^{\nu,\gamma_\pm}H(x,v)
	=H(\Phi_s^\nu(x),A_s^\nu(x)v)\longrightarrow H(x,v)
	\]
	pointwise, and the supports remain in one compact subset of $\gamma^\pm$ for
	small $|s|$. Hence dominated convergence gives the strong continuity for such $H$.
	Since $\gamma^0$ is $\dd\gamma$-null, these restrictions are dense in
	$L^2(\gamma^\pm,w^{2L}\dd\gamma)$. Thus \eqref{eq:boundary-flow-bound} extends
	the strong continuity to every $H$ with $w^LH\in L^2(\gamma^\pm)$.
\end{proof}

\begin{lemma}[Complete-flow traces]
	\label{lem:complete-flow-actual-traces}
	Let $\nu\in\mathfrak A$ and $X=\tau_\nu$, and let $s_0$ be as in
	\eqref{eq:complete-flow-uniform-bounds}.  For $|s|\le s_0$ set
	\begin{equation}
		\label{eq:flow-conjugated-transport-field}
		\beta_s(x,v):=A_s(x)^{-1}\bigl((v\cdot\nabla_x)A_s\bigr)(x)\,v .
	\end{equation}
	Then
	\begin{equation}
		\label{eq:flow-conjugated-field-bound}
		|\beta_s(x,v)|\le C|v|^2\le Cw(v)^2,
		\qquad
		|\nabla_v\beta_s(x,v)|\le C|v|,
		\qquad |s|\le s_0 .
	\end{equation}
	\begin{enumerate}
		\item\label{it:flow-traces-stationary}
		Let $q\in\cH_{\kin,\sigma}(\cD\times\bR_v^3)$ and $|s|\le s_0$.  Then
		\begin{equation}
			\label{eq:flow-conjugated-transport}
			v\cdot\nabla_x(T_s^Xq)
			=T_s^X(v\cdot\nabla_xq)+\beta_s\cdot\nabla_v(T_s^Xq)
			\qquad\text{in }\cD'(\cD\times\bR^3).
		\end{equation}
		If
		\begin{equation}
			\label{eq:flow-graph-coordinates}
			T_s^X q\in L^2(\cD;H^1_\sigma),
			\qquad
			v\cdot\nabla_x(T_s^X q)\in L^2(\cD;H^{-1}_\sigma),
		\end{equation}
		then
		\begin{equation}
			\label{eq:flow-actual-traces}
			T_s^X q\in\cH_{\kin,\sigma}(\cD\times\bR_v^3),
			\qquad
			(T_s^X q)_{\gamma^\pm}
			=T_s^{X,\gamma_\pm}q_{\gamma^\pm}.
		\end{equation}
		If $w^5q\in\cH_{\kin,\sigma}(\cD\times\bR_v^3)$, then
		\eqref{eq:flow-graph-coordinates} holds.
		\item\label{it:flow-traces-evolutionary}
		Let $0<T<\infty$,
		$q\in\cH_{\kin,\sigma}(\cD_T\times\bR_v^3)$, and $|s|\le s_0$.  Define
		\[
		(T_s^Xq)(t,x,v):=q\bigl(t,\widehat\Phi_s^X(x,v)\bigr),
		\qquad
		(T_s^{X,\gamma_\pm}H)(t,x,v):=H\bigl(t,\widehat\Phi_s^X(x,v)\bigr)
		\quad\text{on }\gamma_T^\pm .
		\]
		Then
		\begin{equation}
			\label{eq:flow-transport-evolutionary}
			(\partial_t+v\cdot\nabla_x)(T_s^Xq)
			=T_s^X\bigl((\partial_t+v\cdot\nabla_x)q\bigr)
			+\beta_s\cdot\nabla_v(T_s^Xq)
			\qquad\text{in }\cD'((0,T)\times\cD\times\bR^3).
		\end{equation}
		If
		\begin{equation}
			\label{eq:flow-graph-coordinates-evolutionary}
			T_s^Xq\in L^2(\cD_T;H^1_\sigma),
			\qquad
			(\partial_t+v\cdot\nabla_x)(T_s^Xq)\in L^2(\cD_T;H^{-1}_\sigma),
		\end{equation}
		then
		\begin{equation}
			\label{eq:flow-actual-traces-evolutionary}
			\begin{gathered}
				T_s^Xq\in\cH_{\kin,\sigma}(\cD_T\times\bR_v^3),
				\qquad
				(T_s^Xq)_{\gamma_T^\pm}=T_s^{X,\gamma_\pm}q_{\gamma_T^\pm},
				\\
				(T_s^Xq)(t)=T_s^X\bigl(q(t)\bigr)
				\quad(0\le t\le T),
			\end{gathered}
		\end{equation}
		where $(T_s^Xq)(t)$ and $q(t)$ are the values of the representatives in
		$C([0,T];L^2(\cD\times\bR_v^3))$ of Proposition~\ref{2607151218}.
		If $w^5q\in\cH_{\kin,\sigma}(\cD_T\times\bR_v^3)$, then
		\eqref{eq:flow-graph-coordinates-evolutionary} holds.
	\end{enumerate}
\end{lemma}

\begin{proof}
	Differentiating
	$\partial_sA_s=(\nabla X)(\Phi_s)A_s$ in $x$ gives
	\[
	\partial_sD_xA_s=(\nabla^2X)(\Phi_s)[A_s\,\cdot,A_s]+(\nabla X)(\Phi_s)D_xA_s,
	\]
	and Gronwall's inequality with \eqref{eq:complete-flow-uniform-bounds} gives
	\[
	\sup_{|s|\le s_0,\,x\in\overline\cD}|D_xA_s(x)|
	=\sup_{|s|\le s_0,\,x\in\overline\cD}|D^2\Phi_s(x)|\le C.
	\]
	Hence \eqref{eq:flow-conjugated-field-bound} holds, since $\beta_s$ is quadratic in $v$ and $\|A_s(x)^{-1}\|_{\rm op}\le C$ by \eqref{eq:complete-flow-uniform-bounds}.
	
	We prove (\ref{it:flow-traces-evolutionary}). For (\ref{it:flow-traces-stationary}), omit Steps~1 and~5 and replace \eqref{260730346} by \eqref{260902137} in Step~2 and Proposition~\ref{H_kin=weak_landau}\,(1) by Proposition~\ref{H_kin=weak_landau}\,(2) in Step~4.
	
	Take $q_k\in C_c^\infty([0,T]\times\overline\cD\times\bR^3)$ converging to $q$ in $\cH_{\kin,\sigma}(\cD_T\times\bR_v^3)$, as in Definition~\ref{def:H_kin_sigma}.
	
	\emph{Step 1.}  By Proposition~\ref{2607151218},
	$q\in C([0,T];L^2(\cD\times\bR_v^3))$.  By
	\eqref{eq:complete-flow-L2-bound} with $L=0$,
	\[
	\|T_s^X(q(t))-T_s^X(q(t'))\|_{L^2(\cD\times\bR_v^3)}\le C\|q(t)-q(t')\|_{L^2(\cD\times\bR_v^3)},
	\qquad t,t'\in[0,T],
	\]
	so $t\mapsto T_s^X(q(t))$ belongs to
	$C([0,T];L^2(\cD\times\bR_v^3))$. The
	slice of $G:=T_s^Xq$ at $t$ is $T_s^X(q(t))$ for almost every $t$.
	
	\emph{Step 2.}  By \eqref{260730346},
	\[
	\sup_{0\le t\le T}\|q_k(t)-q(t)\|_{L^2}+\|(q_k)_\gamma-q_\gamma\|_{L^2(\gamma_T)}
	\to0.
	\]  Set $g_k:=T_s^Xq_k$ and $Q_k:=(\partial_t+v\cdot\nabla_x)q_k$.  By the mapping property of the complete pullbacks on $C_c^1(\overline\cD\times\bR^3)$, stated at the beginning of \S\ref{sec:diffuse-complete-equations}, $g_k$ is $C^1$ on
	$[0,T]\times\overline\cD\times\bR^3$ with compact velocity support, and
	$g_k(t)=T_s^X(q_k(t))$.  The phase flow $\widehat\Phi_s$ preserves $\gamma^+$, $\gamma^-$ and $\gamma^0$
	(Lemma~\ref{lem:complete-flow}), so
	$g_k|_{\gamma_T^\pm}=T_s^{X,\gamma_\pm}(q_k|_{\gamma_T^\pm})$.
	For $u\in C^1(\overline\cD\times\bR^3)$, the chain rule gives
	\[
	\nabla_x(T_s^Xu)=A_s^T(\nabla_xu)\circ\widehat\Phi_s
	+\bigl(D_x(A_sv)\bigr)^T(\nabla_vu)\circ\widehat\Phi_s,
	\qquad
	\nabla_v(T_s^Xu)=A_s^T(\nabla_vu)\circ\widehat\Phi_s.
	\]
	By \eqref{eq:flow-conjugated-transport-field},
	\[
	\bigl(D_x(A_sv)\bigr)v=\bigl((v\cdot\nabla_x)A_s\bigr)v=A_s\beta_s.
	\]
	Since the velocity component of $\widehat\Phi_s(x,v)$ is $A_sv$, the formulas for $\nabla_x(T_s^Xu)$ and $\nabla_v(T_s^Xu)$ give
	\begin{equation}
		\label{eq:flow-conjugated-transport-smooth}
		v\cdot\nabla_x(T_s^Xu)
		=(A_sv)\cdot(\nabla_xu)\circ\widehat\Phi_s
		+(A_s\beta_s)\cdot(\nabla_vu)\circ\widehat\Phi_s
		=T_s^X(v\cdot\nabla_xu)+\beta_s\cdot\nabla_v(T_s^Xu).
	\end{equation}
	By \eqref{eq:complete-flow-time-commutation} and
	\eqref{eq:flow-conjugated-transport-smooth} at each $t$,
	\[
	(\partial_t+v\cdot\nabla_x)g_k=T_s^XQ_k+\beta_s\cdot\nabla_vg_k .
	\]  Let
	$\phi\in C_c^\infty([0,T]\times\overline\cD\times\bR^3)$ with
	$\operatorname{supp}\phi\subset[0,T]\times\overline\cD\times B_R$.  The divergence
	theorem in $(t,x)$ at fixed $v$, integrated in $v$, gives \eqref{2607151226} for
	$g_k$. By \eqref{eq:flow-conjugated-transport-smooth}, \eqref{eq:bulk-flow-adjoint-quotient} and an integration by parts in $v$ at each $t$, the first term on the left-hand side of \eqref{2607151226} for $g_k$ is
	\[
	\iint_{\cD_T}\bigl(\la Q_k,(T_s^X)^\dagger\phi\ra_v
	-\la g_k,\nabla_v\cdot(\beta_s\phi)\ra_v\bigr)\dd x\dd t,
	\]
	where
	$(T_s^X)^\dagger\phi\in L^2(\cD_T;H^1_\sigma)$ by the same mapping property for the adjoints and \eqref{eq:Hv-into-Hsigma}, and $\nabla_v\cdot(\beta_s\phi)$ and
	$(\partial_t+v\cdot\nabla_x)\phi$ belong to $L^\infty((0,T)\times\cD\times\bR^3)$ and have support in
	$[0,T]\times\overline\cD\times B_R$.  Every term passes to the limit:
	$$
	Q_k\to Q:=(\partial_t+v\cdot\nabla_x)q\quad \text{in }L^2(\cD_T;H^{-1}_\sigma),\qquad g_k\to G\quad \text{in }L^2(\cD_T\times\bR^3_v),
	$$
	and $g_k(t)\to T_s^X(q(t))$ in $L^2$ for
	$t=0,T$, by \eqref{eq:complete-flow-L2-bound} with $L=0$; and
	$T_s^{X,\gamma_\pm}(q_k)_{\gamma_T^\pm}\to T_s^{X,\gamma_\pm}q_{\gamma_T^\pm}$ in
	$L^2(\gamma_T^\pm)$, by \eqref{eq:boundary-flow-bound} with $L=0$ at each $t$; since
	$\gamma^0$ is $\dd\gamma$-null, this gives the convergence of the boundary term of \eqref{2607151226}.  Hence, for such $\phi$,
	\begin{equation}
		\label{eq:flow-green-limit-evolutionary}
		\begin{aligned}
			&\iint_{\cD_T}\Bigl(\la Q,(T_s^X)^\dagger\phi\ra_v
			-\la G,\nabla_v\cdot(\beta_s\phi)\ra_v
			+\la G,(\partial_t+v\cdot\nabla_x)\phi\ra_v\Bigr)\dd x\dd t\\
			&\qquad=\Bigl[\iint_{\cD\times\bR^3}T_s^X\bigl(q(t)\bigr)\phi(t)
			\dd x\dd v\Bigr]_{t=0}^{t=T}
			+\iint_{\gamma_T}\bigl(T_s^{X,\gamma}q_\gamma\bigr)\phi
			\operatorname{sgn}(n_x\cdot v)\dd\gamma\dd t ,
		\end{aligned}
	\end{equation}
	where $T_s^{X,\gamma}q_\gamma$ equals $T_s^{X,\gamma_\pm}q_{\gamma_T^\pm}$ on
	$\gamma_T^\pm$ and belongs to $L^2(\gamma_T)$.
	
	\emph{Step 3.}  For
	$\phi\in C_c^\infty((0,T)\times\cD\times\bR^3)$ the right-hand side of
	\eqref{eq:flow-green-limit-evolutionary} vanishes.  By the definition of the
	distributional pullback \eqref{eq:bulk-flow-adjoint-quotient} and of
	$\beta_s\cdot\nabla_vG$, this is
	\eqref{eq:flow-transport-evolutionary}.
	
	\emph{Step 4.}  Assume
	\eqref{eq:flow-graph-coordinates-evolutionary} and fix $R>0$.  For $\phi$ as in
	Step~2 with velocity support in $B_R$, the functionals
	\begin{align*}
		\ell_1(\phi)&:=\iint_{\cD_T}\la(\partial_t+v\cdot\nabla_x)G,\phi\ra_v\dd x\dd t,\\
		\ell_2(\phi)&:=\iint_{\cD_T}\bigl(\la Q,(T_s^X)^\dagger\phi\ra_v
		-\la G,\nabla_v\cdot(\beta_s\phi)\ra_v\bigr)\dd x\dd t
	\end{align*}
	satisfy
	\[
	|\ell_1(\phi)|+|\ell_2(\phi)|\le C\|\phi\|_{L^2(\cD_T;H^1_v)}\qquad\text{with } C=C(R) ,
	\]
	by
	\eqref{eq:flow-graph-coordinates-evolutionary}, \eqref{eq:Hv-into-Hsigma}, the change
	of variables $\widehat\Phi_{-s}$ and \eqref{eq:flow-conjugated-field-bound}. Moreover,
	$\ell_1=\ell_2$ on $C_c^\infty((0,T)\times\cD\times B_R)$ by Step~3.  Let $0\le\chi_n\le1$ be compactly supported cutoffs in $(0,T)\times\cD$ converging pointwise to $1$. Dominated convergence gives $\chi_n\phi\to\phi$ in $L^2(\cD_T;H^1_v)$; the preceding bound then extends $\ell_1=\ell_2$ from interior tests to $\phi$.
	Substituting this into \eqref{eq:flow-green-limit-evolutionary} gives
	\eqref{2607151226} for $G$ with
	\[
	G(0,\cdot)=T_s^X\bigl(q(0)\bigr),
	\qquad
	G(T,\cdot)=T_s^X\bigl(q(T)\bigr),
	\qquad
	G_\gamma=T_s^{X,\gamma}q_\gamma ,
	\]
	which belong to $L^2(\cD\times\bR^3_v)$ and $L^2(\gamma_T)$.  Together with
	\eqref{eq:flow-graph-coordinates-evolutionary},
	Proposition~\ref{H_kin=weak_landau}\,(1) gives
	$G\in\cH_{\kin,\sigma}(\cD_T\times\bR_v^3)$, and its uniqueness statement
	identifies the traces of $G$ with $T_s^{X,\gamma_\pm}q_{\gamma_T^\pm}$.
	
	\emph{Step 5.}  Both the continuous representative of $G$ and $t\mapsto T_s^X(q(t))$ belong to $C([0,T];L^2(\cD\times\bR_v^3))$. Their almost-everywhere equality extends to every $t\in[0,T]$.  This proves \eqref{eq:flow-actual-traces-evolutionary}.
	
	\emph{Step 6.}  Let
	$w^5q\in\cH_{\kin,\sigma}(\cD_T\times\bR_v^3)$.  Then $w^5q\in L^2(\cD_T;H^1_\sigma)$ and
	\[
	w^5Q=(\partial_t+v\cdot\nabla_x)(w^5q)\in L^2(\cD_T;H^{-1}_\sigma).
	\]
	At
	almost every $t$, \eqref{eq:complete-flow-H1-bound} with $(L,p)=(4,2)$ and
	\eqref{eq:complete-flow-Hminus1-bound} with $L=0$ give, respectively,
	\[
	\|w^4G(t)\|_{L_x^2H_\sigma^1}\le C\|w^5q(t)\|_{L_x^2H_\sigma^1},
	\qquad
	\|T_s^XQ(t)\|_{L_x^2H_\sigma^{-1}}\le C\|wQ(t)\|_{L_x^2H_\sigma^{-1}}.
	\]
	By \eqref{eq:Hsigma-multiplier} with $\zeta=w^{-4}$, applied to $w^4G(t)$ and to $w^5Q(t)$,
	$G\in L^2(\cD_T;H^1_\sigma)$ and $wQ\in L^2(\cD_T;H^{-1}_\sigma)$; hence $T_s^XQ\in L^2(\cD_T;H^{-1}_\sigma)$.
	For the last term of \eqref{eq:flow-transport-evolutionary},
	\eqref{eq:flow-conjugated-field-bound},
	\eqref{eq:coulomb-eigen-decomposition} and \eqref{eq:eigenvalue-size} give
	\[
	|\beta_s\cdot\nabla_vG|\le Cw^2|\nabla_vG|\le Cw^{7/2}|\sigma^{1/2}\nabla_vG|.
	\]
	Thus, pointwise in $(t,x)$ and for $\varphi\in H^1_\sigma$,
	\[
	\Bigl|\int_{\bR^3}(\beta_s\cdot\nabla_vG)\varphi\dd v\Bigr|
	\le C\|w^4\sigma^{1/2}\nabla_vG\|_{L_v^2}\|w^{-1/2}\varphi\|_{L_v^2}
	\le C\|w^4G\|_{H_\sigma^1}\|\varphi\|_{H_\sigma^1}
	\]
	by \eqref{eq:flow-weighted-coercive-equivalence} at $\mathrm L=4$ and
	\eqref{eq:auxiliary-coercive-equivalence}.  Hence
	$\beta_s\cdot\nabla_vG\in L^2(\cD_T;H^{-1}_\sigma)$, and
	\eqref{eq:flow-transport-evolutionary} gives
	\eqref{eq:flow-graph-coordinates-evolutionary}.
\end{proof}

\subsection{Finite differences}\label{app:flow-identities}

Complete-flow differences realize the differentiated kinetic and collision identities in the graph framework.

For $h\ne0$, define the bulk, boundary and adjoint difference quotients
\begin{equation}
	\label{eq:bulk-flow-difference}
	\Delta_h^Xq:=\frac{T_h^Xq-q}{h},
	\qquad
	\Delta_h^{X,\gamma_\pm}H:=\frac{T_h^{X,\gamma_\pm}H-H}{h},
	\qquad
	\Delta_h^{X,\dagger}:=\frac{(T_h^X)^\dagger-I}{h},
\end{equation}
with $T_h^{X,\gamma_\pm}=T_h^{\nu,\gamma_\pm}$ the boundary pullback
\eqref{eq:boundary-flow-pullback} of $X=\tau_\nu$ and $(T_h^X)^\dagger$ the adjoint
\eqref{eq:bulk-flow-adjoint-quotient}.  Unless stated otherwise, fix $0<h_0\le s_0/\mathrm N$, with $s_0$ as in
\eqref{eq:complete-flow-uniform-bounds}, and take $0<|h|<h_0$. All assertions are uniform in this choice of $h_0$.  An iterated difference quotient of
order $k$ is written as an integral over the unit cube $[0,1]^k$ in the variables
$\theta_1,\ldots,\theta_k$ through the forward-difference Peano formula
\[
\frac1{h^k}\sum_{j=0}^k(-1)^{k-j}\binom kjF(c+jh)
=\int_{[0,1]^k}F^{(k)}\bigl(c+(\theta_1+\cdots+\theta_k)h\bigr)\,\dd\theta_1\cdots\dd\theta_k
\]
for $F$ in $W^{k,1}$ on an interval containing $c$ and $c+kh$, with values in a Banach space, the left-hand side being evaluated at the continuous representative of $F$.  Auxiliary flow parameters $\bar h$ and $\bar h'$ are restricted explicitly wherever they occur.

Let $1\le m\le\mathrm N$ and $q\in\cH_{\kin,\sigma}$, and assume
\eqref{eq:flow-actual-traces} at $s=jh$ for $1\le j\le m$; on $(0,T)$, assume
instead $T_s^Xq\in\cH_{\kin,\sigma}(\cD_T\times\bR_v^3)$ and the trace identity of
\eqref{eq:flow-actual-traces-evolutionary}
at the same $s$.  The linearity of the Green trace maps gives
\begin{equation}
	\label{eq:direct-same-flow-actual-traces}
	\bigl((\Delta_h^X)^mq\bigr)_{\gamma^\pm}
	= (\Delta_h^{X,\gamma_\pm})^m q_{\gamma^\pm}.
\end{equation}

For $1\le r\le\mathrm N$ and an operator $\mathrm A$, define
\begin{equation}
	\label{eq:flow-operator-difference}
	\mathfrak Q_{r,h}^X(\mathrm A)
	:=\frac1{h^r}\sum_{j=0}^r(-1)^{r-j}\binom rj
	T_{jh}^X\mathrm AT_{-jh}^X.
\end{equation}

Let $\mathscr T_L:=v\cdot\nabla_x+L$, with $L$ as in \eqref{def:LAK}.  By \eqref{eq:ad-word},
$\operatorname{ad}_{\mathfrak D_X}^r=\operatorname{ad}_{X^r}$ for $X\in\cX_{\rm c}$ and $r\ge1$.

\begin{lemma}[Weighted phase-space density]
	\label{lem:weighted-phase-density}
	Let $M\ge1$ be an integer, and let $L_i\ge0$ and $1<p_i<\infty$,
	$1\le i\le M$.  If
	\[
	q\in L^2_{x,v},
	\qquad
	w^{L_i}q\in L_x^{p_i}H_\sigma^1,
	\qquad 1\le i\le M,
	\]
	then there are $q_k\in C_c^\infty(\cD\times\bR^3)$ such that
	\begin{equation}
		\label{eq:weighted-phase-density-limit}
		\|q_k-q\|_{L^2_{x,v}}
		+\sum_{i=1}^M
		\|w^{L_i}(q_k-q)\|_{L_x^{p_i}H_\sigma^1}
		\longrightarrow0
		\qquad(k\to\infty).
	\end{equation}
\end{lemma}

\begin{proof}
	Let $\chi_R$ be the cutoff in \eqref{eq:radial-cutoff}, choose
	$\eta_j\in C_c^\infty(\cD)$ with $0\le\eta_j\le1$ and $\eta_j\to1$ pointwise, and
	set $q_{j,R}:=\eta_j\chi_Rq$.  Dominated convergence, \eqref{eq:radial-cutoff-convergence}, and
	\eqref{eq:Hsigma-multiplier} show that $q_{j,R}\to q$ in all the stated norms,
	first as $R\to\infty$ and then as $j\to\infty$.  The velocity support of $q_{j,R}$ lies in
	$B_{2R}$, so \eqref{eq:local-coulomb-sobolev-comparison} gives
	$q_{j,R}\in\bigcap_{i=1}^ML_x^{p_i}H_v^1$.  For a standard phase-space mollifier
	$\rho_\varepsilon$ and small $\varepsilon$,
	$\rho_\varepsilon*q_{j,R}\in C_c^\infty(\cD\times\bR^3)$ has velocity support in
	one fixed ball and converges to $q_{j,R}$ in $L^2_{x,v}$ and in every
	$L_x^{p_i}H_v^1$, hence, by \eqref{eq:Hv-into-Hsigma}, in the norms of
	\eqref{eq:weighted-phase-density-limit}.  A diagonal sequence gives \eqref{eq:weighted-phase-density-limit}.
\end{proof}

For $X=\tau_\nu$, $\nu\in\mathfrak A$, and $F\in C_c(\cD\times\bR^3)$, the change of
variables under the phase flow \eqref{eq:complete-phase-flow} gives
\[
(T_s^X)^\dagger F(y)
=F(\widehat\Phi_{-s}^X(y))
|\det D\widehat\Phi_{-s}^X(y)|.
\]
By \eqref{eq:complete-phase-flow}, \eqref{eq:complete-flow-jacobians} and
\eqref{eq:complete-flow-uniform-bounds}, there are a compact set $K_F$ depending only on
$\operatorname{supp}F$ and $C=C(\|F\|_{L^\infty})$ such that
\begin{equation}
	\label{eq:complete-flow-adjoint-domination}
	\operatorname{supp}(T_s^X)^\dagger F\subset K_F\quad(|s|\le s_0),
	\qquad
	\sup_{|s|\le s_0}\|(T_s^X)^\dagger F\|_{L^\infty}
	\le C.
\end{equation}

\begin{lemma}
	\label{lem:repeated-same-flow-FTC}
	Let $X=\tau_\nu$ with $\nu\in\mathfrak A$.
	\begin{enumerate}
		\item\label{it:repeated-FTC-identified}
		Let $1\le m\le\mathrm N$, $L\ge0$,
		$2\le p<\infty$ and $1\le\rho\le\infty$.  Let $u$ be a function of $(x,v)$ or of $(t,x,v)$; in the first case $L_t^\rho$ is omitted from the norms below, and in the second case every identity below is understood
		to hold at almost every $t$.  Let $u_j:=\mathfrak D_X^ju$, $0\le j\le m$, in the
		sense of distributions in $(x,v)$, and assume
		\[
		w^{L+1}u_j\in L_t^\rho L_x^pH_\sigma^1,
		\qquad 0\le j\le m.
		\]
		Then, for $1\le k\le m$ and
		$0<|h|<h_0$,
		\begin{equation}
			\label{eq:repeated-same-flow-FTC}
			(\Delta_h^X)^ku
			=\int_{[0,1]^k}
			T_{(\theta_1+\cdots+\theta_k)h}^Xu_k\,\dd\theta_1\cdots\dd\theta_k
			\quad\text{in }\cD'.
		\end{equation}
		Moreover, for $1\le k\le m$,
		\begin{equation}
			\label{eq:repeated-same-flow-weighted-bound}
			\sup_{0<|h|<h_0}
			\|w^L(\Delta_h^X)^ku\|_{L_t^\rho L_x^pH_\sigma^1}
			\le C
			\|w^{L+1}u_k\|_{L_t^\rho L_x^pH_\sigma^1},
		\end{equation}
		where $C$ is independent of $h$, $\rho$ and $u$.
		\item\label{it:repeated-FTC-fiber}
		Let $\lambda\in\bR$ and
		$1\le k\le\mathrm N$, let $H_\lambda$ be as in
		Lemma~\ref{lem:complete-flow-weighted-bounds}\,(\ref{it:complete-flow-fiber}) with $d=1$, and let $G_0,\ldots,G_k\in L_x^1H_\lambda$ satisfy
		$\mathfrak D_XG_j=G_{j+1}$ in $\cD'(\cD\times\bR^3)$, $0\le j<k$.  If $h\ne0$ and $|\bar h|+k|h|<s_0$, then, for
		almost every $x\in\cD$, in $H_\lambda$,
		\begin{equation}
			\label{eq:repeated-same-flow-FTC-fiber}
			T_{\bar h}^X(\Delta_h^X)^kG_0
			=\int_{[0,1]^k}
			T_{\bar h+(\theta_1+\cdots+\theta_k)h}^XG_k\,\dd\theta_1\cdots\dd\theta_k .
		\end{equation}
		If, in addition, $w^\lambda G_k\in L_x^pL_v^2$ for some $1\le p\le\infty$, then
		\begin{equation}
			\label{eq:repeated-same-flow-fiber-bound}
			\|w^\lambda T_{\bar h}^X(\Delta_h^X)^kG_0\|_{L_x^pL_v^2}
			\le C\|w^\lambda G_k\|_{L_x^pL_v^2},
		\end{equation}
		with $C$ of \eqref{eq:complete-flow-fiber-bound}.  On $(0,T)$, if these
		hypotheses hold at almost every $t$, then \eqref{eq:repeated-same-flow-FTC-fiber} and
		\eqref{eq:repeated-same-flow-fiber-bound} hold at almost every $t$.
	\end{enumerate}
\end{lemma}

\begin{proof}
	\emph{(\ref{it:repeated-FTC-identified}).}  For time-dependent $u$, work at almost every $t$.  Since
	$p\ge2$ and $\cD$ is bounded, \eqref{eq:local-H1-H_sigma} and
	\eqref{eq:real-weight-Hsigma-monotonicity} give $u_j\in L^2(\cD\times\bR_v^3)$,
	$0\le j\le m$, so \eqref{eq:strong-complete-flow-generator-FTC} applies.
	Every flow parameter $s$ below satisfies $|s|\le k|h|<\mathrm Nh_0\le s_0$.
	The fundamental theorem of calculus \eqref{eq:strong-complete-flow-generator-FTC}
	and the group law give $T_h^Xu-u=\int_0^hT_s^Xu_1\dd s$.
	Iterating this identity, using the commutation of $\Delta_h^X$ with $T_s^X$, yields
	\eqref{eq:repeated-same-flow-FTC}.  The right-hand side of
	\eqref{eq:repeated-same-flow-FTC} belongs to
	$\{z:w^Lz\in L_x^pH_\sigma^1\}$ by
	\eqref{eq:complete-flow-H1-bound}, uniformly in
	$(\theta_1,\ldots,\theta_k)\in[0,1]^k$.  The distributional identity therefore holds in the stated weighted space.  Minkowski's inequality in $(\theta_1,\ldots,\theta_k)$,
	\eqref{eq:complete-flow-H1-bound} and integration in $t$ give
	\eqref{eq:repeated-same-flow-weighted-bound}.
	
	\emph{(\ref{it:repeated-FTC-fiber}).}  Set $F_x(s):=T_s^XG_0(x,\cdot)$. By Lemma~\ref{lem:distributional-flow-slicing},
	for almost every $x$ this curve has a representative in
	$W^{k,1}((-s_0,s_0);H_\lambda)$ satisfying
	$F_x^{(k)}(s)=T_s^XG_k(x,\cdot)$ for almost every $s$; below, $F_x$ denotes this
	continuous representative.  The group law, \eqref{eq:complete-flow-fiber-bound} at $p=1$, and
	\eqref{eq:lifted-pullback-strong-continuity} make $s\mapsto T_s^XG_0$
	continuous in $L_x^1H_\lambda$.  For $0\le j\le k$,
	$|\bar h+jh|\le|\bar h|+k|h|<s_0$. Hence, by this continuity, $F_x(\bar h+jh)=T_{\bar h+jh}^XG_0(x,\cdot)$ in $H_\lambda$ for almost every $x$, $0\le j\le k$, and thus
	\[
	T_{\bar h}^X(\Delta_h^X)^kG_0(x,\cdot)
	=\frac1{h^k}\sum_{j=0}^k(-1)^{k-j}\binom kj F_x(\bar h+jh).
	\]
	The $W^{k,1}$ Peano formula gives \eqref{eq:repeated-same-flow-FTC-fiber};
	Minkowski's inequality and \eqref{eq:complete-flow-fiber-bound}, applied for
	$|\bar h+(\theta_1+\cdots+\theta_k)h|\le|\bar h|+k|h|<s_0$, give
	\eqref{eq:repeated-same-flow-fiber-bound}.
\end{proof}

\begin{lemma}[Operator differences]
	\label{lem:direct-same-flow}
	Let $X=\tau_\nu$ with $\nu\in\mathfrak A$, and recall
	$\mathfrak Q_{r,h}^X(\mathrm A)$ from \eqref{eq:flow-operator-difference}.
	\begin{enumerate}
		\item\label{it:direct-same-flow-leibniz}
		For every linear operator $\mathrm A$, every $u$ for which the terms are defined, and
		$1\le m\le\mathrm N$,
		\begin{equation}
			\label{eq:forward-discrete-Leibniz}
			(\Delta_h^X)^m(\mathrm Au)
			=\mathrm A(\Delta_h^X)^mu
			+\sum_{r=1}^m\binom mr
			\mathfrak Q_{r,h}^X(\mathrm A)
			T_{rh}^X(\Delta_h^X)^{m-r}u.
		\end{equation}
		\item\label{it:direct-same-flow-equation}
		If $v\cdot\nabla_x q+L_gq=F$ in $\cD'$ and $1\le m\le\mathrm N$, then, whenever the terms are defined,
		\begin{equation}
			\label{eq:direct-same-flow-identity}
			(v\cdot\nabla_x+L_g)(\Delta_h^X)^mq
			=(\Delta_h^X)^mF
			-\sum_{r=1}^m\binom mr
			\mathfrak Q_{r,h}^X(v\cdot\nabla_x+L_g)
			T_{rh}^X(\Delta_h^X)^{m-r}q
			\quad\text{in }\cD'.
		\end{equation}
		\item\label{it:direct-same-flow-operator}
		Let $1\le r\le\mathrm N$, and let $\mathrm A=\mathscr T_L$ or
		$\mathrm A=\operatorname{ad}_J(\mathscr T_L)$, with $\operatorname{ad}_J$ of \eqref{eq:ad-word} and
		$J$ an ordered complete word with $|J|\ge1$ and $|J|+r\le\mathrm N$.  Then, as operators on
		$C_c^\infty(\cD\times\bR^3)$,
		\begin{equation}
			\label{eq:direct-operator-cube}
			\mathfrak Q_{r,h}^X(\mathrm A)
			=\int_{[0,1]^r}
			T_{(\theta_1+\cdots+\theta_r)h}^X
			\operatorname{ad}_{\mathfrak D_X}^r(\mathrm A)
			T_{-(\theta_1+\cdots+\theta_r)h}^X\,
			\dd\theta_1\cdots\dd\theta_r.
		\end{equation}
		Set $B_r:=\operatorname{ad}_{\mathfrak D_X}^r(\mathrm A)$.  Let $(L_i,p_i)$,
		$1\le i\le M$, satisfy $L_i\ge0$ and $1<p_i<\infty$, and assume that, for every
		$z\in C_c^\infty(\cD\times\bR^3)$,
		\begin{align}
			&\sup_{|\bar h|\le r h_0}
			\|T_{\bar h}^XB_rT_{-\bar h}^Xz\|_{L_x^2H_\sigma^{-1}}
			\notag\\*
			&\qquad\le C\left(
			\|z\|_{L^2_{x,v}}+\|z\|_{L_x^2H_\sigma^1}
			+\sum_{i=1}^M\|w^{L_i}z\|_{L_x^{p_i}H_\sigma^1}
			\right),
			\label{eq:operator-orbit-bound}
		\end{align}
		for some constant $C$ independent of $z$, and that, when $\mathrm A=\operatorname{ad}_J(\mathscr T_L)$,
		\eqref{eq:operator-orbit-bound} holds with $B_r$ replaced by $\mathrm A$.  Then
		\eqref{eq:direct-operator-cube} holds in $\cD'$ at every $u$ for which
		\[
		\|u\|_{L^2_{x,v}}+\|u\|_{L_x^2H_\sigma^1}+\sum_{i=1}^M\|w^{L_i}u\|_{L_x^{p_i}H_\sigma^1}<\infty.
		\]
	\end{enumerate}
\end{lemma}

\begin{proof}
	\emph{(\ref{it:direct-same-flow-leibniz}).}  Set $\mathrm A_a:=T_{ah}^X\mathrm AT_{-ah}^X$ and
	$u_c:=T_{ch}^Xu$.  After multiplying \eqref{eq:forward-discrete-Leibniz} by $h^m$ and including
	$\mathrm A(\Delta_h^X)^mu$ as the $r=0$ term, the coefficient of
	$\mathrm A_au_c$ on the right is
	\begin{align*}
		&\sum_{r=a}^c
		\binom mr(-1)^{r-a}\binom ra
		(-1)^{m-c}\binom{m-r}{c-r}\\
		&\qquad=(-1)^{m-c}\binom mc\binom ca
		\sum_{r=a}^c(-1)^{r-a}\binom{c-a}{r-a}
		=(-1)^{m-c}\binom mc\mathbf 1_{\{a=c\}}.
	\end{align*}
	Only the terms $a=c$ survive, with the coefficient of $T_{ch}^X(\mathrm Au)$, which proves \eqref{eq:forward-discrete-Leibniz}.
	
	\emph{(\ref{it:direct-same-flow-equation}).}  This is \eqref{eq:forward-discrete-Leibniz}
	with $\mathrm A=v\cdot\nabla_x+L_g$, $u=q$ and $\mathrm Aq=F$.
	
	\emph{(\ref{it:direct-same-flow-operator}).}  Set
	$B_j:=\operatorname{ad}_{\mathfrak D_X}^j(\mathrm A)$, $0\le j\le r$.  For $u,\phi\in C_c^\infty(\cD\times\bR^3)$ set
	$I:=[-\mathrm N|h|,\mathrm N|h|]\subset(-s_0,s_0)$, since
	$\mathrm N|h|<\mathrm Nh_0\le s_0$.  For $0\le j\le r$ set
	\[
	G_j(s):=
	\left\langle B_jT_{-s}^Xu,(T_s^X)^\dagger\phi\right\rangle.
	\]
	If $\mathrm A=\mathscr T_L$, \eqref{boundedness of L} and $L^2$ duality for the transport
	term give $G_0\in L^1(I)$.  If $\mathrm A=\operatorname{ad}_J(\mathscr T_L)$, the decomposition
	\eqref{eq:nested-commutator-decomposition} at $g=0$ reads $\mathrm A=\mathfrak T_J+\mathfrak L_J+\mathfrak C_J$, the sum of the transport, local and compact commutators. These three terms satisfy
	\eqref{eq:ordered-transport-weight-integrated},
	\eqref{eq:nested-local-complete-bound-spatial} and
	\eqref{eq:iterated-compact-full-operator-bound} at $\ell=0$, respectively, and these bounds, with \eqref{eq:complete-flow-H1-bound} at $p=2$ and \eqref{eq:complete-flow-L2-bound} at $L=0$ for $T_{-s}^Xu$, and with \eqref{eq:complete-flow-adjoint-domination}, \eqref{eq:complete-flow-uniform-bounds} and \eqref{eq:Hv-into-Hsigma} for $(T_s^X)^\dagger\phi$, give $G_0\in L^1(I)$.
	For $j\ge1$, the operators
	$\operatorname{ad}_{\mathfrak D_X}^j(\mathscr T_L)$ and
	$\operatorname{ad}_{\mathfrak D_X}^j\operatorname{ad}_J(\mathscr T_L)=\operatorname{ad}_{X^jJ}(\mathscr T_L)$ are finite sums of
	transport, local and compact commutators of \eqref{eq:nested-commutator-decomposition} at $g=0$, of total complete order at most
	$\mathrm N$, since $j\le r\le\mathrm N$ and $|J|+r\le\mathrm N$.  Hence \eqref{eq:ordered-transport-weight-integrated},
	\eqref{eq:nested-local-complete-bound-spatial},
	\eqref{eq:iterated-compact-full-operator-bound}, and
	\eqref{eq:complete-flow-H1-bound}--\eqref{eq:complete-flow-Hminus1-bound} give
	$G_j\in L^1(I)$ for $1\le j\le r$.  By
	\eqref{eq:strong-complete-flow-generator-FTC},
	\eqref{eq:complete-field-formal-adjoint} and \eqref{eq:r11-sd-flow-commutation}, in $\cD'(I)$,
	\begin{align}
		G_j'(s)
		&=-\left\langle B_j\mathfrak D_XT_{-s}^Xu,
		(T_s^X)^\dagger\phi\right\rangle
		+\left\langle \mathfrak D_XB_jT_{-s}^Xu,
		(T_s^X)^\dagger\phi\right\rangle\notag\\
		&=\left\langle [\mathfrak D_X,B_j]T_{-s}^Xu,
		(T_s^X)^\dagger\phi\right\rangle
		=\left\langle T_s^XB_{j+1}T_{-s}^Xu,\phi\right\rangle,
		\label{eq:direct-conjugated-pairing-derivative}
	\end{align}
	where $[\mathfrak D_X,B_j]=B_{j+1}$; the right-hand side belongs to $L^1(I)$ by the
	same bounds.  Hence
	$G_j\in W^{1,1}(I)$ for $0\le j<r$.  Iterating
	\eqref{eq:direct-conjugated-pairing-derivative}, we obtain
	\[
	\frac{\dd^r}{\dd s^r}
	\bigl(T_s^X\mathrm AT_{-s}^Xu\bigr)
	=T_s^XB_rT_{-s}^Xu
	\qquad\text{in }\cD'(I\times\cD\times\bR^3).
	\]
	The $W^{r,1}$ Peano formula, paired with $\phi$, now gives
	\eqref{eq:direct-operator-cube} for smooth $u$.
	
	Assume that
	\[
	\|u\|_{L^2_{x,v}}+\|u\|_{L_x^2H_\sigma^1}+\sum_{i=1}^M\|w^{L_i}u\|_{L_x^{p_i}H_\sigma^1}<\infty.
	\]
	Lemma~\ref{lem:weighted-phase-density}, applied with the pairs $(L_i,p_i)$,
	$1\le i\le M$, and $(0,2)$, gives $u_k\in C_c^\infty(\cD\times\bR^3)$ converging to
	$u$ in the norms on the right-hand side of \eqref{eq:operator-orbit-bound}.  By \eqref{eq:operator-orbit-bound}, the integrands in
	\eqref{eq:direct-operator-cube} at $u_k$ are Cauchy in $L_x^2H_\sigma^{-1}$,
	uniformly on $[0,1]^r$, since $|(\theta_1+\cdots+\theta_r)h|\le rh_0$. The limit is the integrand at $u$, since $B_r$ is a finite sum of the transport,
	local, and compact commutators in \eqref{eq:nested-commutator-decomposition}
	at $g=0$, each continuous into $\cD'$ along this approximation.  The left-hand side also converges in $\cD'$. For $\mathrm A=\mathscr T_L$, use
	$L^2$ duality against the tests $(T_{jh}^X)^\dagger\phi$. By the mapping property of the adjoints of the complete pullbacks, stated at the beginning of \S\ref{sec:diffuse-complete-equations}, these tests lie in $C_c^{K_\cD}(\cD\times\bR^3)$, so that $L(T_{jh}^X)^\dagger\phi\in L^2_{x,v}$. Since $L$ is symmetric by \eqref{eq:L-symmetry}, we have
	\[
	\langle LT_{-jh}^X(u_k-u),(T_{jh}^X)^\dagger\phi\rangle
	=\langle T_{-jh}^X(u_k-u),L(T_{jh}^X)^\dagger\phi\rangle.
	\]
	For
	$\mathrm A=\operatorname{ad}_J(\mathscr T_L)$, the convergence follows from \eqref{eq:operator-orbit-bound} with $B_r$
	replaced by $\mathrm A$.  This extends \eqref{eq:direct-operator-cube} to $u$.
\end{proof}

\begin{corollary}[Conjugated operator orbits]
	\label{cor:A0-orbits}
	Let $X=\tau_\nu$ with $\nu\in\mathfrak A$, let $1\le r\le\mathrm N$ and $\ell\ge0$.  For every ordered
	complete word $I$
	with $1\le|I|\le\mathrm N$, every $z$ with $w^{\ell+2|I|+4}z\in L_x^2H_\sigma^1$,  and every $|\bar h|\le s_0$, with $\operatorname{ad}_I(\mathscr T_L)$ given by \eqref{eq:nested-commutator-decomposition} at $g=0$,
	\begin{equation}
		\label{eq:A0-orbit-bound}
		\bigl\|w^\ell T_{\bar h}^X\operatorname{ad}_I(\mathscr T_L)
		T_{-\bar h}^Xz\bigr\|_{L_x^2H_\sigma^{-1}}
		\le C\|w^{\ell+2|I|+4}z\|_{L_x^2H_\sigma^1},
	\end{equation}
	with $C$ independent of $\bar h$ and $z$.
	In particular, at $\ell=0$, \eqref{eq:operator-orbit-bound} holds for $\mathrm A=\mathscr T_L$
	with the single pair $(L_1,p_1)=(2r+4,2)$.  For an ordered complete word $J$ with
	$|J|\ge1$ and $|J|+r\le\mathrm N$, \eqref{eq:operator-orbit-bound} and its form with
	$B_r$ replaced by $\mathrm A$ hold for $\mathrm A=\operatorname{ad}_J(\mathscr T_L)$ with the single pair
	$(L_1,p_1)=(2(r+|J|)+4,2)$.
\end{corollary}

\begin{proof}
	At $g=0$, \eqref{eq:nested-commutator-decomposition} reads
	$\operatorname{ad}_I(\mathscr T_L)=\mathfrak T_I+\mathfrak L_I+\mathfrak C_I$.
	By \eqref{eq:ordered-transport-weight-integrated},
	\eqref{eq:nested-local-complete-bound-spatial} at $m=|I|$ and
	\eqref{eq:iterated-compact-full-operator-bound}, for every $\ell'\ge0$ and every $u$ with
	\[
	\|w^{\ell'+4}u\|_{L_x^2H_\sigma^1}+\|w^{\ell'+2|I|+2}u\|_{L_x^2H_\sigma^1}+\|u\|_{L^2_{x,v}}<\infty,
	\]
	we have
	\begin{gather*}
		\|w^{\ell'}\mathfrak T_Iu\|_{L_x^2H_\sigma^{-1}}\le C\|w^{\ell'+4}u\|_{L_x^2H_\sigma^1},\\
		\|w^{\ell'}\mathfrak L_Iu\|_{L_x^2H_\sigma^{-1}}\le C\|w^{\ell'+2|I|+2}u\|_{L_x^2H_\sigma^1},\qquad \|w^{\ell'}\mathfrak C_Iu\|_{L_x^2H_\sigma^{-1}}\le C\|u\|_{L^2_{x,v}}.
	\end{gather*}
	The outer pullback costs one weight by \eqref{eq:complete-flow-Hminus1-bound}:
	\[
	\|w^\ell T_{\bar h}^X\operatorname{ad}_I(\mathscr T_L)T_{-\bar h}^Xz\|_{L_x^2H_\sigma^{-1}}
	\le C\|w^{\ell+1}\operatorname{ad}_I(\mathscr T_L)T_{-\bar h}^Xz\|_{L_x^2H_\sigma^{-1}}.
	\]
	Apply the estimates for $\mathfrak T_I$, $\mathfrak L_I$, and $\mathfrak C_I$
	with $\ell'=\ell+1$. The inner pullback costs one further weight, by
	\eqref{eq:complete-flow-H1-bound}:
	\[
	\|w^{\ell'+4}T_{-\bar h}^Xz\|_{L_x^2H_\sigma^1}\le C\|w^{\ell+6}z\|_{L_x^2H_\sigma^1},
	\qquad
	\|w^{\ell'+2|I|+2}T_{-\bar h}^Xz\|_{L_x^2H_\sigma^1}\le C\|w^{\ell+2|I|+4}z\|_{L_x^2H_\sigma^1},
	\]
	and $\|T_{-\bar h}^Xz\|_{L^2_{x,v}}\le C\|z\|_{L^2_{x,v}}$ by \eqref{eq:complete-flow-L2-bound}.  Hence
	\[
	\|w^\ell T_{\bar h}^X\operatorname{ad}_I(\mathscr T_L)T_{-\bar h}^Xz\|_{L_x^2H_\sigma^{-1}}
	\le C\bigl(\|w^{\ell+6}z\|_{L_x^2H_\sigma^1}+\|w^{\ell+2|I|+4}z\|_{L_x^2H_\sigma^1}+\|z\|_{L^2_{x,v}}\bigr).
	\]  Since $|I|\ge1$, the two weighted terms on the right are bounded by
	$C\|w^{\ell+2|I|+4}z\|_{L_x^2H_\sigma^1}$.  The $L^2$ term is controlled by
	\eqref{eq:local-H1-H_sigma}, and weight monotonicity
	\eqref{eq:real-weight-Hsigma-monotonicity} gives \eqref{eq:A0-orbit-bound}.  For $\mathrm A=\mathscr T_L$,
	\[
	B_r=\operatorname{ad}_{X^r}(\mathscr T_L)
	\]
	by \eqref{eq:ad-word}, so \eqref{eq:A0-orbit-bound} at $I=X^r$ and $\ell=0$ gives the pair $(2r+4,2)$. For
	$\mathrm A=\operatorname{ad}_J(\mathscr T_L)$, 
	\[
	B_r=\operatorname{ad}_{\mathfrak D_X}^r\operatorname{ad}_J(\mathscr T_L)
	=\operatorname{ad}_{X^rJ}(\mathscr T_L)\qquad\text{with } |X^rJ|=r+|J|\le\mathrm N ,
	\]
	so
	\eqref{eq:A0-orbit-bound} at $\ell=0$, at $I=X^rJ$ for $B_r$ and at $I=J$ for $\mathrm A$, with
	\eqref{eq:real-weight-Hsigma-monotonicity}, gives the pair $(2(r+|J|)+4,2)$.  Since $rh_0\le\mathrm Nh_0\le s_0$, \eqref{eq:A0-orbit-bound} covers the range $|\bar h|\le rh_0$ of
	\eqref{eq:operator-orbit-bound}.
\end{proof}

\begin{prop}[Complete-flow collision differences]
	\label{prop:complete-flow-collision-differences}
	Let $X=\tau_\nu$ with $\nu\in\mathfrak A$.
	Let $\mathsf K=\mathsf K_{K_0}$ be generated from the Coulomb tensor by
	\eqref{eq:iterated-coulomb-deformation-kernels}, let $N_{\mathsf K}$ be given by
	\eqref{eq:generalized-conjugated-landau-form}, and write $m_0:=|K_0|$; at $K_0=\varnothing$,
	$N_{\mathsf K_\varnothing}=N$.  Define on the smooth core
	\begin{align}
		\operatorname{Def}_{1,X}[a,b]
		&:={\mathfrak D}_XN_{\mathsf K}[a,b]
		-N_{\mathsf K}[{\mathfrak D}_Xa,b]-N_{\mathsf K}[a,{\mathfrak D}_Xb],
		\label{eq:collision-flow-defect-generator-one}\\
		\operatorname{Def}_{r+1,X}[a,b]
		&:={\mathfrak D}_X\operatorname{Def}_{r,X}[a,b]
		-\operatorname{Def}_{r,X}[{\mathfrak D}_Xa,b]
		-\operatorname{Def}_{r,X}[a,{\mathfrak D}_Xb],
		\qquad m_0+r\le\mathrm N-1.
		\label{eq:collision-flow-defect-generator-recursion}
	\end{align}
	At $K_0=\varnothing$, $\operatorname{Def}_{1,X}[a,b]=\mathfrak R_\nu[a,b]$ by \eqref{eq:collision-remainder-expansion} with $I=(\nu)$, since $\mathfrak D_X=\mathfrak D_\nu$.
	For $1\le r\le\mathrm N-m_0$, set
	\begin{equation}
		\label{eq:collision-flow-defect-difference}
		\mathscr C_{r,h}^X[a,b]
		:=\frac1{h^r}\sum_{j=0}^r(-1)^{r-j}\binom rj
		\mathscr N_{jh}^X[T_{rh}^Xa,T_{rh}^Xb],
		\qquad
		\mathscr N_s^X[a,b]:=T_s^XN_{\mathsf K}[T_{-s}^Xa,T_{-s}^Xb].
	\end{equation}
	\begin{enumerate}
		\item\label{it:collision-differences-Peano}
		For $1\le r\le\mathrm N-m_0$,
		\begin{equation}
			\label{eq:collision-flow-defect-Peano}
			\mathscr C_{r,h}^X[a,b]
			=\int_{[0,1]^r}
			T_{(\theta_1+\cdots+\theta_r)h}^X\operatorname{Def}_{r,X}
			[T_{(r-\theta_1-\cdots-\theta_r)h}^Xa,T_{(r-\theta_1-\cdots-\theta_r)h}^Xb]\,
			\dd\theta_1\cdots\dd\theta_r.
		\end{equation}
		Equation \eqref{eq:collision-flow-defect-Peano} holds in distributions for all $a$ and
		$b$ with
		\begin{equation}
			\label{eq:collision-flow-defect-Peano-class}
			w^{-\lambda}a\in L^2_{x,v}\ \text{ for some }\lambda\ge0,
			\qquad
			w^{2(m_0+r)+3}b\in L_x^2H_\sigma^1 .
		\end{equation}
		\item\label{it:collision-differences-pointwise}
		For $1\le r\le\mathrm N-m_0$, $1/4\le\varpi<1/2$ and every finite $L\ge0$,
		\begin{equation}
			\label{eq:collision-flow-defect-bound}
			\|w^L\operatorname{Def}_{r,X}[a,b]\|_{H_\sigma^{-1}}
			\le C
			\|\mu^\varpi a\|_{L_v^2}
			\|w^{L+2(m_0+r)+2}b\|_{H_\sigma^1}.
		\end{equation}
		By density (Lemma~\ref{lem:H_sigma-Hilbert}\,(1)), the bound
		\eqref{eq:collision-flow-defect-bound} extends, in distributions, to all $a,b$ with
		$\mu^\varpi a\in L_v^2$ and $w^{L+2(m_0+r)+2}b\in H_\sigma^1$.
		\item\label{it:collision-differences-Leibniz}
		Set $\mathscr C_{0,h}^X:=N_{\mathsf K}$.  For $1\le m\le\mathrm N-m_0$,
		\begin{equation}
			\label{eq:collision-mth-endpoint}
			(\Delta_h^X)^mN_{\mathsf K}[a,b]
			=\sum_{i+j+r=m}\frac{m!}{i!\,j!\,r!}\,
			\mathscr C_{r,h}^X\bigl[
			T_{(m-r-i)h}^X(\Delta_h^X)^ia,\ (\Delta_h^X)^jb\bigr],
		\end{equation}
		where the sum is taken over all $i,j,r\ge0$ with $i+j+r=m$.  Equation \eqref{eq:collision-mth-endpoint} holds
		for all $a$ and $b$ for which its terms are defined.
		\item\label{it:collision-differences-weighted}
		Let $1\le r\le\mathrm N-m_0$, $\ell\ge0$ and $1/4\le\varpi<1/2$, and let
		$2<q\le\infty$ and $2\le p<\infty$ satisfy $1/q+1/p\le1/2$.  Then, for
		all $a$ and $b$ satisfying \eqref{eq:collision-flow-defect-Peano-class} for which
		\[
		\sup_{\theta_1,\ldots,\theta_r\in[0,1]}\|\mu^\varpi T_{(r-\theta_1-\cdots-\theta_r)h}^Xa\|_{L_x^qL_v^2}\,\|w^{\ell+2(m_0+r)+4}b\|_{L_x^pH_\sigma^1}<\infty,
		\]
		we have
		\begin{equation}
			\label{eq:collision-defect-weighted}
			\|w^\ell\mathscr C_{r,h}^X[a,b]\|_{L_x^2H_\sigma^{-1}}
			\le C\sup_{\theta_1,\ldots,\theta_r\in[0,1]}
			\|\mu^\varpi T_{(r-\theta_1-\cdots-\theta_r)h}^Xa\|_{L_x^qL_v^2}\,
			\|w^{\ell+2(m_0+r)+4}b\|_{L_x^pH_\sigma^1},
		\end{equation}
		where $C$ is independent of $h$, $a$ and $b$.  On $(0,T)$, let $1/\rho+1/\rho'\le1/2$, and let
		\eqref{eq:collision-flow-defect-Peano-class} hold at almost every $t$.  Then
		\[
		\|w^\ell\mathscr C_{r,h}^X[a,b]\|_{L_t^2L_x^2H_\sigma^{-1}}
		\le C\sup_{\theta_1,\ldots,\theta_r\in[0,1]}
		\|\mu^\varpi T_{(r-\theta_1-\cdots-\theta_r)h}^Xa\|_{L_t^\rho L_x^qL_v^2}\,
		\|w^{\ell+2(m_0+r)+4}b\|_{L_t^{\rho'}L_x^pH_\sigma^1}.
		\]
		In particular, let $\lambda\ge0$.  For all $a$ and $b$ with $w^{-\lambda}a\in L_x^qL_v^2$ and
		$w^{\ell+2(m_0+r)+4}b\in L_x^pH_\sigma^1$, \eqref{eq:collision-flow-defect-Peano-class} holds and
		\begin{equation}
			\label{261016801}
			\|w^\ell\mathscr C_{r,h}^X[a,b]\|_{L_x^2H_\sigma^{-1}}
			\le C\|w^{-\lambda}a\|_{L_x^qL_v^2}\,
			\|w^{\ell+2(m_0+r)+4}b\|_{L_x^pH_\sigma^1},
		\end{equation}
		where $C$ depends on $\lambda$ and is independent of $h$, $a$ and $b$.  On $(0,T)$, let $1/\rho+1/\rho'\le1/2$.  For all $a$ and $b$ with
		$w^{-\lambda}a\in L_t^\rho L_x^qL_v^2$ and $w^{\ell+2(m_0+r)+4}b\in L_t^{\rho'}L_x^pH_\sigma^1$, \eqref{eq:collision-flow-defect-Peano-class} holds at almost every $t$,
		and \eqref{261016801} holds with $L_t^\rho L_x^qL_v^2$, $L_t^{\rho'}L_x^pH_\sigma^1$ and $L_t^2L_x^2H_\sigma^{-1}$ in place of
		$L_x^qL_v^2$, $L_x^pH_\sigma^1$ and $L_x^2H_\sigma^{-1}$.
	\end{enumerate}
	
\end{prop}

\begin{proof}
	On the smooth core, conjugation \eqref{eq:complete-conjugated-kernel-deformation},
	the defining recursion \eqref{eq:collision-flow-defect-generator-one}--\eqref{eq:collision-flow-defect-generator-recursion},
	and \eqref{eq:r11-sd-flow-commutation} give inductively
	\begin{equation}
		\label{eq:collision-flow-defect-orbit-derivative}
		\frac{\dd^r}{\dd s^r}\mathscr N_s^X[a,b]
		=T_s^X\operatorname{Def}_{r,X}
		[T_{-s}^Xa,T_{-s}^Xb],
		\qquad 1\le r\le\mathrm N-m_0.
	\end{equation}
	
	\emph{(\ref{it:collision-differences-Peano}).}  On the smooth core, the forward-difference Peano formula for
	$s\mapsto\mathscr N_s^X[T_{rh}^Xa,T_{rh}^Xb]$, with
	\eqref{eq:collision-flow-defect-orbit-derivative}, gives
	\eqref{eq:collision-flow-defect-Peano}.  Let $a,b$ satisfy
	\eqref{eq:collision-flow-defect-Peano-class}, fix
	$\phi\in C_c^\infty(\cD\times\bR^3)$ and $\varpi\in[\tfrac14,\tfrac12)$, and set
	$L_r:=2(m_0+r)+3$.  By \eqref{eq:complete-flow-adjoint-domination} and
	\eqref{eq:complete-flow-uniform-bounds}, the tests $(T_{\bar h}^X)^\dagger\phi$,
	$|\bar h|\le r h_0$, are supported in one compact set and are bounded, together
	with their velocity gradients, uniformly in $\bar h$; for the velocity gradients we use that the Jacobian of
	$\widehat\Phi_s^X$ depends only on $x$.  Hence
	\eqref{eq:collision-flow-defect-bound} at $L=0$ and the Cauchy--Schwarz inequality in $x$ give,
	for $|\bar h|,|\bar h'|\le rh_0$,
	\[
	\bigl|\bigl\langle T_{\bar h}^X\operatorname{Def}_{r,X}
	[T_{\bar h'}^Xa,T_{\bar h'}^Xb],\phi\bigr\rangle\bigr|
	\le C\|\mu^\varpi T_{\bar h'}^Xa\|_{L^2_{x,v}}
	\|w^{L_r-1}T_{\bar h'}^Xb\|_{L_x^2H_\sigma^1}
	\le C\|w^{-\lambda}a\|_{L^2_{x,v}}
	\|w^{L_r}b\|_{L_x^2H_\sigma^1},
	\]
	where $C=C(\phi,\lambda)$ is independent of $\bar h$ and $\bar h'$, and the last step uses
	$\mu^\varpi\le Cw^{-\lambda}$,
	\eqref{eq:complete-flow-L2-bound} at the real weight $-\lambda$, and
	\eqref{eq:complete-flow-H1-bound} at $p=2$.  The same argument with
	\eqref{eq:deformation-weighted-output} at the order $\max\{m_0,1\}$ in place of
	\eqref{eq:collision-flow-defect-bound} gives
	\[
	|\langle\mathscr N_{jh}^X[T_{rh}^Xa,T_{rh}^Xb],\phi\rangle|
	\le C\|w^{-\lambda}a\|_{L^2_{x,v}}\|w^{L_r}b\|_{L_x^2H_\sigma^1},
	\qquad 0\le j\le r,\quad 0<|h|<h_0,
	\]
	where $C=C(\phi,\lambda)$ is independent of $h$.  Here \eqref{eq:deformation-weighted-output} at $\ell=0$ and \eqref{eq:complete-flow-H1-bound} require the weight $2\max\{m_0,1\}+3$ on $b$, which is at most $L_r$.
	For fixed $h$, the two sides of \eqref{eq:collision-flow-defect-Peano},
	paired with $\phi$, are bilinear in $(a,b)$.  By \eqref{eq:collision-flow-defect-difference} and the bound for $\langle\mathscr N_{jh}^X[T_{rh}^Xa,T_{rh}^Xb],\phi\rangle$, $0\le j\le r$,
	\[
	|\langle\mathscr C_{r,h}^X[a,b],\phi\rangle|
	\le2^rC|h|^{-r}\|w^{-\lambda}a\|_{L^2_{x,v}}\|w^{L_r}b\|_{L_x^2H_\sigma^1}.
	\]
	For $(\theta_1,\ldots,\theta_r)\in[0,1]^r$, the bound for $\langle T_{\bar h}^X\operatorname{Def}_{r,X}[T_{\bar h'}^Xa,T_{\bar h'}^Xb],\phi\rangle$ applies at $\bar h=(\theta_1+\cdots+\theta_r)h$ and $\bar h'=(r-\theta_1-\cdots-\theta_r)h$, since $|\bar h|,|\bar h'|\le r|h|<rh_0$; integration over $[0,1]^r$ gives
	\[
	\begin{aligned}
		&\Bigl|\int_{[0,1]^r}\bigl\langle T_{(\theta_1+\cdots+\theta_r)h}^X\operatorname{Def}_{r,X}
		[T_{(r-\theta_1-\cdots-\theta_r)h}^Xa,T_{(r-\theta_1-\cdots-\theta_r)h}^Xb],\phi\bigr\rangle\,
		\dd\theta_1\cdots\dd\theta_r\Bigr|\\
		&\qquad\le C\|w^{-\lambda}a\|_{L^2_{x,v}}\|w^{L_r}b\|_{L_x^2H_\sigma^1}.
	\end{aligned}
	\]
	Here $C=C(\phi,\lambda)$, except for the factor $|h|^{-r}$ in the first bound;
	this is harmless because $h$ is fixed during the density argument.  The space
	$C_c^\infty(\cD\times\bR^3)$ is dense in $L^2(w^{-2\lambda}\dd x\dd v)$, and by
	Lemma~\ref{lem:weighted-phase-density} with $(L_1,p_1)=(L_r,2)$ it is dense
	for $\|b\|_{L^2_{x,v}}+\|w^{L_r}b\|_{L_x^2H_\sigma^1}$; here
	$b\in L^2_{x,v}$ by \eqref{eq:local-H1-H_sigma}.  Hence
	\eqref{eq:collision-flow-defect-Peano} holds in $\cD'$ for all $a,b$ with
	\eqref{eq:collision-flow-defect-Peano-class}.
	
	\emph{(\ref{it:collision-differences-pointwise}).}  Iterating \eqref{eq:complete-conjugated-kernel-deformation}, every summand of
	$\operatorname{Def}_{r,X}$ has the form $pN_{\mathsf K_A}[p_1a,p_2b]$ of
	Lemma~\ref{lem:deformation-degrees} with $k=r$ in which neither slot is differentiated,
	and contains at least one kernel or Maxwellian-conjugation deformation.  By \eqref{eq:deformation-combined-grading} with $J$ replaced by $A$, we have $|A|\le m_0+r$ and
	\[
	\deg_vp+\deg_vp_1+\deg_vp_2\le2(m_0+r);
	\]
	the coefficients are bounded on
	$\overline\cD$.  Proposition~\ref{prop:deformation-algebra} at
	$m=m_0+r$ therefore gives
	\eqref{eq:collision-flow-defect-bound}.
	
	\emph{(\ref{it:collision-differences-Leibniz}).}  The definition \eqref{eq:collision-flow-defect-difference} at $r=1$ gives
	\begin{equation}
		\label{eq:collision-endpoint-first}
		\Delta_h^XN_{\mathsf K}[a,b]
		=N_{\mathsf K}[T_h^Xa,\Delta_h^Xb]+N_{\mathsf K}[\Delta_h^Xa,b]
		+\mathscr C_{1,h}^X[a,b],
	\end{equation}
	which is \eqref{eq:collision-mth-endpoint} at $m=1$.
	The same definition gives, for $0\le r\le\mathrm N-m_0-1$,
	\begin{equation}
		\label{eq:collision-defect-endpoint-recursion}
		\Delta_h^X\mathscr C_{r,h}^X[a,b]
		=\mathscr C_{r,h}^X[T_h^Xa,\Delta_h^Xb]
		+\mathscr C_{r,h}^X[\Delta_h^Xa,b]
		+\mathscr C_{r+1,h}^X[a,b].
	\end{equation}
	Apply $\Delta_h^X$ to \eqref{eq:collision-mth-endpoint} and use
	\eqref{eq:collision-defect-endpoint-recursion}. A term indexed by $(i,j,r)$
	contributes to $(i,j+1,r)$, $(i+1,j,r)$ and $(i,j,r+1)$. Their shifts are again those of
	\eqref{eq:collision-mth-endpoint} at $m+1$, by the group law and since $\Delta_h^X$ commutes with $T_s^X$. Collecting them, we find that the coefficient
	of $(i,j,r)$ at order $m+1$ is
	\[
	\frac{m!}{i!\,(j-1)!\,r!}
	+\frac{m!}{(i-1)!\,j!\,r!}
	+\frac{m!}{i!\,j!\,(r-1)!}
	=\frac{(m+1)!}{i!\,j!\,r!},
	\]
	with the convention that a term with a negative factorial index vanishes.
	This proves the formula at order $m+1$.
	
	Since \eqref{eq:collision-endpoint-first} and
	\eqref{eq:collision-defect-endpoint-recursion} are finite algebraic identities in the
	bilinear maps $\mathscr N_{jh}^X$, \eqref{eq:collision-mth-endpoint} holds whenever its
	terms are defined.
	
	\emph{(\ref{it:collision-differences-weighted}).}  By \eqref{eq:collision-flow-defect-Peano} and Minkowski's inequality in
	$(\theta_1,\ldots,\theta_r)$, it suffices to bound, uniformly in $(\theta_1,\ldots,\theta_r)\in[0,1]^r$,
	$w^\ell T_{\bar h}^X\operatorname{Def}_{r,X}[T_{rh-\bar h}^Xa,T_{rh-\bar h}^Xb]$, with
	$\bar h:=(\theta_1+\cdots+\theta_r)h$; for the flow
	parameters here, $|\bar h|$ and $|rh-\bar h|$ are at most
	$r|h|<s_0$.  Put $a':=T_{rh-\bar h}^Xa$ and $b':=T_{rh-\bar h}^Xb$. By \eqref{eq:complete-flow-Hminus1-bound}, \eqref{eq:collision-flow-defect-bound} at $L=\ell+1$, H\"older's inequality in $x$ with $\frac1q+\frac1p\le\frac12$ on the bounded set $\cD$, and \eqref{eq:complete-flow-H1-bound},
	\[
	\begin{aligned}
		\|w^\ell T_{\bar h}^X\operatorname{Def}_{r,X}[a',b']\|_{L_x^2H_\sigma^{-1}}
		&\le C\|w^{\ell+1}\operatorname{Def}_{r,X}[a',b']\|_{L_x^2H_\sigma^{-1}}\\
		&\le C\bigl\|\|\mu^\varpi a'\|_{L_v^2}\|w^{\ell+2(m_0+r)+3}b'\|_{H_\sigma^1}\bigr\|_{L_x^2}\\
		&\le C\|\mu^\varpi a'\|_{L_x^qL_v^2}\|w^{\ell+2(m_0+r)+4}b\|_{L_x^pH_\sigma^1}.
	\end{aligned}
	\]
	This gives \eqref{eq:collision-defect-weighted}.  On $(0,T)$,
	\eqref{eq:collision-flow-defect-Peano} holds at almost every $t$. We apply the bounds above and
	H\"older's inequality in $t$ on the bounded interval to the integrand at each
	fixed $(\theta_1,\ldots,\theta_r)$, and then Minkowski's inequality in these
	variables in $L_t^2L_x^2H_\sigma^{-1}$; this gives the bound on $(0,T)$.
	For \eqref{261016801}, H\"older's inequality in $x$ on the
	bounded set $\cD$, with $q>2$ and $p\ge2$, and \eqref{eq:real-weight-Hsigma-monotonicity} with $\ell\ge0$ give
	\[
	\|w^{-\lambda}a\|_{L^2_{x,v}}\le C\|w^{-\lambda}a\|_{L_x^qL_v^2},
	\qquad
	\|w^{2(m_0+r)+3}b\|_{L_x^2H_\sigma^1}\le C\|w^{\ell+2(m_0+r)+4}b\|_{L_x^pH_\sigma^1},
	\]
	and hence \eqref{eq:collision-flow-defect-Peano-class}; on $(0,T)$, this holds at almost every $t$.  For $|s|\le r|h|$, where $r|h|<\mathrm Nh_0\le s_0$, the bound $\mu^\varpi\le Cw^{-\lambda}$ and
	\eqref{eq:complete-flow-fiber-bound} with $d=1$, $p=q$ and the real weight $-\lambda$ give
	\[
	\|\mu^\varpi T_s^Xa\|_{L_x^qL_v^2}\le C\|w^{-\lambda}T_s^Xa\|_{L_x^qL_v^2}\le C\|w^{-\lambda}a\|_{L_x^qL_v^2},
	\]
	and the same bound holds with an outer $L_t^\rho$ norm on both sides.  Since $|(r-\theta_1-\cdots-\theta_r)h|\le r|h|$ for $\theta_1,\ldots,\theta_r\in[0,1]$, it follows that
	\[
	\sup_{\theta_1,\ldots,\theta_r\in[0,1]}\|\mu^\varpi T_{(r-\theta_1-\cdots-\theta_r)h}^Xa\|_{L_x^qL_v^2}\le C\|w^{-\lambda}a\|_{L_x^qL_v^2}<\infty ,
	\]
	and \eqref{eq:collision-defect-weighted} gives \eqref{261016801}.  On $(0,T)$, the same argument with the outer $L_t^\rho$ norm and the bound for $\|w^\ell\mathscr C_{r,h}^X[a,b]\|_{L_t^2L_x^2H_\sigma^{-1}}$ in part~(\ref{it:collision-differences-weighted}) gives \eqref{261016801} with $L_t^\rho L_x^qL_v^2$, $L_t^{\rho'}L_x^pH_\sigma^1$ and $L_t^2L_x^2H_\sigma^{-1}$ in place of $L_x^qL_v^2$, $L_x^pH_\sigma^1$ and $L_x^2H_\sigma^{-1}$.
\end{proof}

\begin{corollary}
	\label{cor:deformation-term-differences}
	Let $X=\tau_\nu$ with $\nu\in\mathfrak A$, let $\Phi$ be a term $p\,N_{\mathsf K_A}[p_1a,p_2b]$
	of \eqref{eq:deformation-remainder-structure} at word length $m$, and let
	$0\le r\le\mathrm N-m$.  Then
	\begin{equation}
		\label{eq:collision-deformation-difference-expansion}
		(\Delta_h^X)^r\Phi
		=\sum_{k+s+i+j=r}\frac{r!}{k!\,s!\,i!\,j!}\,
		\bigl((\Delta_h^X)^kp\bigr)\,
		T_{kh}^X\mathscr C_{s,h}^X\bigl[
		T_{(r-k-s-i)h}^X(\Delta_h^X)^i(p_1a),\
		(\Delta_h^X)^j(p_2b)\bigr],
	\end{equation}
	where $\mathscr C_{s,h}^X$ is formed from $\mathsf K_A$ by
	\eqref{eq:collision-flow-defect-difference} and
	$\mathscr C_{0,h}^X=N_{\mathsf K_A}$.  Consequently, for every finite $\ell\ge0$, every $\mathrm L_{\rm G}\in\bR$
	and all $2\le q,p\le\infty$ with $1/q+1/p\le1/2$,
	\begin{equation}
		\label{eq:collision-deformation-difference-bound}
		\|w^\ell(\Delta_h^X)^r\Phi\|_{L_x^2H_\sigma^{-1}}
		\le C\max_{i+j\le r}\ \sup_{|\bar h|\le r|h|}
		\|w^{\mathrm L_{\rm G}}T_{\bar h}^X(\Delta_h^X)^ia\|_{L_x^qL_v^2}\,
		\|w^{\ell+2(m+r)+4}(\Delta_h^X)^jb\|_{L_x^pH_\sigma^1},
	\end{equation}
	where $C$ is independent of $h$, $a$ and $b$.  On $(0,T)$, for exponents with
	$1/q+1/p\le1/2$ and $1/\rho+1/\rho'\le1/2$ chosen separately for each $(i,j)$,
	\[
	\|w^\ell(\Delta_h^X)^r\Phi\|_{L_t^2L_x^2H_\sigma^{-1}}
	\le C\max_{i+j\le r}\ \sup_{|\bar h|\le r|h|}
	\|w^{\mathrm L_{\rm G}}T_{\bar h}^X(\Delta_h^X)^ia\|_{L_t^\rho L_x^qL_v^2}\,
	\|w^{\ell+2(m+r)+4}(\Delta_h^X)^jb\|_{L_t^{\rho'}L_x^pH_\sigma^1}.
	\]
\end{corollary}

\begin{proof}
	The discrete product rule for the composition operator $T_s^X$ is
	\[
	\Delta_h^X(uz)=\bigl(\Delta_h^Xu\bigr)\bigl(T_h^Xz\bigr)+u\,\Delta_h^Xz .
	\]
	Iterating this in the factor $p$ and expanding the remaining difference of
	$N_{\mathsf K_A}$ by \eqref{eq:collision-mth-endpoint} at $m_0=|A|\le m$ gives
	\eqref{eq:collision-deformation-difference-expansion}; the same rule
	separates $p_1$ and $p_2$ from $a$ and $b$.  These operations do not raise the velocity degree of the polynomials;
	their coefficients remain uniformly bounded for $0<|h|<h_0$.  Indeed, for a polynomial $P$ in $v$
	of degree $d'$, we have
	\[
	(T_s^XP)(x,v)=P\bigl(\Phi_s^X(x),A_s^X(x)v\bigr),
	\qquad
	(\Delta_h^X)^kP=\int_{[0,1]^k}
	T^X_{(\theta_1+\cdots+\theta_k)h}\mathfrak D_X^kP\,\dd\theta_1\cdots\dd\theta_k ,
	\]
	where $A_s^X$ is bounded for $|s|\le2s_0$ by \eqref{eq:complete-flow-uniform-bounds}. Moreover, $\mathfrak D_X$ does not raise the $v$-degree, since its velocity part $(\nabla X\,v)\cdot\nabla_v$ is linear in $v$.
	For $k\le r\le\mathrm N-m$, the coefficients of $\mathfrak D_X^kP$ involve
	at most $k$ spatial derivatives of those of $P$ and $k+1$ derivatives of
	$\tau_\nu$. These are bounded by Lemma~\ref{lem:deformation-degrees}
	and \eqref{eq:order-and-domain-regularity}.
	
	Write $d:=\deg_vp$, $d_1:=\deg_vp_1$, $d_2:=\deg_vp_2$.  Since $\Phi$ is a term
	of \eqref{eq:deformation-remainder-structure} at word length $m$,
	\eqref{eq:deformation-combined-grading} with $J$ replaced by $A$ gives
	\begin{equation}
		\label{eq:deformation-difference-weight-count}
		2|A|+d+d_2\le2|A|+d+d_1+d_2\le2m .
	\end{equation}
	By \eqref{eq:Hsigma-multiplier} with $\zeta=w^{-d'}P$, a polynomial factor $P$ of
	degree at most $d'$ with bounded coefficients adds $d'$ to the weight:
	\[
	\|Pf\|_{H_\sigma^1}\le C\|w^{d'}f\|_{H_\sigma^1},
	\qquad
	\|PF\|_{H_\sigma^{-1}}\le C\|w^{d'}F\|_{H_\sigma^{-1}}.
	\]
	In the second slot,
	$\Delta_h^X$ and the pullbacks turn $p_2b$ into a sum of terms
	$P_2\,T_{\bar h'}^X(\Delta_h^X)^{j'}b$ with $\deg_vP_2\le d_2$, $j'\le j$ and one
	composed pullback $T_{\bar h'}^X$, $|\bar h'|\le r|h|$; the first slot is
	treated in the same way with $p_1$.
	
	Let $s\ge1$.  By \eqref{eq:collision-flow-defect-Peano}, the summand is an
	integral over $(\theta_1,\ldots,\theta_s)\in[0,1]^s$, with
	$\bar h:=(\theta_1+\cdots+\theta_s)h$, of
	$(\Delta_h^X)^kp\cdot T_{kh+\bar h}^X\operatorname{Def}_{s,X}[\cdot,\cdot]$, where
	$\operatorname{Def}_{s,X}$ is formed from $\mathsf K_A$, so $m_0=|A|$ in
	\eqref{eq:collision-flow-defect-bound}.  The weight losses, from the outer factor inward, are
	$d$ for $(\Delta_h^X)^kp$; one for $T_{kh+\bar h}^X$ by
	\eqref{eq:complete-flow-Hminus1-bound}; $2(|A|+s)+2$ for $\operatorname{Def}_{s,X}$ by
	\eqref{eq:collision-flow-defect-bound}; $d_2$ for $P_2$; one for
	$T_{\bar h'}^X$ by \eqref{eq:complete-flow-H1-bound}.  By
	\eqref{eq:deformation-difference-weight-count} and $s\le r$,
	\[
	\ell+d+1+2(|A|+s)+2+d_2+1
	\le\ell+2m+2s+4
	\le\ell+2(m+r)+4 .
	\]
	Let $s=0$.  Write the summand as
	$T_{kh}^X\bigl(\tilde p\,N_{\mathsf K_A}[\cdot,\cdot]\bigr)$ with
	$\tilde p:=T_{-kh}^X(\Delta_h^X)^kp$, of degree at most $d$.  The pullback
	$T_{kh}^X$ adds one to the weight.  The remaining term is
	$\tilde pN_{\mathsf K_A}[P_1a',P_2b']$, where $|A|\le m$, $a',b'$ are the pulled-back differences of $a,b$, and $\tilde p,P_1,P_2$ are polynomials of degrees at most $d,d_1,d_2$ whose coefficients are bounded as shown above.  By \eqref{eq:deformation-difference-weight-count},
	\[
	\deg\tilde p+\deg P_1+\deg P_2\le d+d_1+d_2\le2m .
	\]
	Hence $\tilde p,P_1,P_2$ satisfy \eqref{eq:deformation-algebra-general-polynomials}, and
	\eqref{eq:deformation-weighted-output-spatial} at the order $m$ applies.  This estimate
	adds $2m+2$ to the weight, and $T_{\bar h'}^X$ adds one more.  The total is
	$\ell+2m+4\le\ell+2(m+r)+4$.
	
	In both cases the first slot is measured in $\mu^\varpi L_v^2$, with
	$\varpi\in[\frac14,\frac12)$ fixed; moreover, $\mu^\varpi|P_1|\le Cw^{\mathrm L_{\rm G}}$.
	H\"older's inequality in $x$ with $1/q+1/p\le1/2$ and the boundedness of $\cD$ give
	\eqref{eq:collision-deformation-difference-bound}.  On $(0,T)$, the bound of
	each term above holds at almost every $t$ and at each fixed
	$(\theta_1,\ldots,\theta_s)$. We apply H\"older's inequality in $t$ on the bounded
	interval, at the exponents chosen for the difference orders of the two factors of the term, and then Minkowski's inequality in
	$(\theta_1,\ldots,\theta_s)$ in $L_t^2L_x^2H_\sigma^{-1}$; this gives the bound on $(0,T)$.
\end{proof}

\subsection{Weak closure}\label{app:weak-convergence-closure}

We identify weak limits of complete-flow differences and prove the Friedrichs commutator estimate needed for $L_x^2H_v^1$.

For the stationary problem let $\Omega=\{\ast\}$; for the evolutionary problem let $\Omega=(0,T)$, $0<T<\infty$. Set
\begin{equation}
	\label{eq:complete-jet-closure-base}
	\begin{aligned}
		\mathsf Q_\Omega&:=\Omega\times\cD\times\bR_v^3,
		&\mathsf H_\Omega&:=\cH_{\kin,\sigma}(\mathsf Q_\Omega),
		&(\mathsf T_\Omega,\gamma_\Omega^\pm)&:=
		\begin{cases}
			(v\cdot\nabla_x,\ \gamma^\pm), &\Omega=\{\ast\},\\
			(\p_t+v\cdot\nabla_x,\ \gamma_T^\pm), &\Omega=(0,T).
		\end{cases}
	\end{aligned}
\end{equation}
Thus $\mathsf H_{(0,T)}=\cH_{\kin,\sigma}(\cD_T\times\bR_v^3)$; when $\Omega=\{\ast\}$, norms over $\Omega\times\cD$ mean norms over $\cD$.

\begin{lemma}[Limits of complete-flow differences]
	\label{lem:adjoint-complete-flow-differences}
	Let $\Omega$ and $\mathsf Q_\Omega$ be as in \eqref{eq:complete-jet-closure-base}, let
	$X=\tau_\nu$ with $\nu\in\mathfrak A$, let $1\le m\le\mathrm N$, and let
	$q\in L^2(\mathsf Q_\Omega)$.  Then
	\begin{equation}
		\label{eq:complete-flow-graph-closure-distribution}
		(\Delta_h^X)^mq\longrightarrow\mathfrak D_X^mq
		\quad\text{in }\cD'(\mathsf Q_\Omega)
		\quad\text{as }h\to0.
	\end{equation}
\end{lemma}

\begin{proof}
	For $\psi\in C_c^\infty(\cD\times\bR^3)$, the adjoint of
	$\partial_sT_s^X=\mathfrak D_XT_s^X$ from \eqref{eq:r11-sd-flow-commutation} is
	$\partial_s(T_s^X)^\dagger\psi=(T_s^X)^\dagger\mathfrak D_X^\dagger\psi$.
	Integration in $s$ gives
	\[
	\Delta_h^{X,\dagger}\psi
	=\frac1h\int_0^h (T_s^X)^\dagger
	\mathfrak D_X^\dagger\psi\dd s.
	\]
	Applying the same identity to $(\mathfrak D_X^\dagger)^k\psi$, $0\le k<m$, and using the group law yields
	\begin{equation}
		\label{eq:direct-same-flow-adjoint-cube}
		(\Delta_h^{X,\dagger})^m\psi
		=\int_{[0,1]^m}
		(T_{(\theta_1+\cdots+\theta_m)h}^X)^\dagger
		(\mathfrak D_X^\dagger)^m\psi\,\dd\theta_1\cdots\dd\theta_m.
	\end{equation}
	For $F=(\mathfrak D_X^\dagger)^m\psi$, the adjoint change-of-variables formula preceding \eqref{eq:complete-flow-adjoint-domination} reads
	\[
	(T_s^X)^\dagger F(y)=F(\widehat\Phi_{-s}^X(y))|\det D\widehat\Phi_{-s}^X(y)| .
	\]
	As $s\to0$, the flow and its Jacobian converge locally uniformly to the identity and one.
	The domination in \eqref{eq:complete-flow-adjoint-domination} gives
	$(T_s^X)^\dagger F\to F$ in $L^2$. Hence \eqref{eq:direct-same-flow-adjoint-cube} and Minkowski's inequality imply
	\begin{equation}
		\label{eq:complete-flow-adjoint-difference-limit}
		(\Delta_h^{X,\dagger})^m\psi
		\longrightarrow (\mathfrak D_X^\dagger)^m\psi
		\qquad\text{in }L^2(\cD\times\bR^3)
		\quad\text{as }h\to0.
	\end{equation}
	The duality identity
	$\langle(\Delta_h^X)^mq,\psi\rangle=\langle q,(\Delta_h^{X,\dagger})^m\psi\rangle$
	from \eqref{eq:bulk-flow-adjoint-quotient} and \eqref{eq:bulk-flow-difference} now proves
	\eqref{eq:complete-flow-graph-closure-distribution} in the stationary case. When
	$\Omega=(0,T)$ and $\psi\in C_c^\infty(\mathsf Q_\Omega)$, the same pairing holds at
	almost every $t$. By \eqref{eq:direct-same-flow-adjoint-cube} for $\psi(t)$ and \eqref{eq:complete-flow-adjoint-domination} for $F=(\mathfrak D_X^\dagger)^m\psi(t)$,
	\[
	|\langle q(t),(\Delta_h^{X,\dagger})^m\psi(t)\rangle|\le C\|q(t)\|_{L^2(\cD\times\bR^3)},
	\qquad 0<|h|<h_0,
	\]
	at almost every $t$, with $C$ independent of $h$ and $t$. The limit \eqref{eq:complete-flow-adjoint-difference-limit} at each $t$ and dominated convergence then prove \eqref{eq:complete-flow-graph-closure-distribution} in the evolutionary case.
\end{proof}

\begin{lemma}[Weak closure]
	\label{lem:complete-flow-graph-closure}
	Let $\Omega$, $\mathsf Q_\Omega$, $\mathsf H_\Omega$,
	$\mathsf T_\Omega$ and $\gamma_\Omega^\pm$ be as in
	\eqref{eq:complete-jet-closure-base}.
	Let $X=\tau_\nu$ with $\nu\in\mathfrak A$, let $1\le m\le\mathrm N$, and let
	$q\in\mathsf H_\Omega$.  For $0<|h|<h_0$ set $Z_h:=(\Delta_h^X)^m q$ and assume
	\begin{equation}
		\label{eq:complete-flow-graph-closure-bound}
		\sup_{0<|h|<h_0}\|Z_h\|_{\mathsf H_\Omega}<\infty.
	\end{equation}
	Then $\mathfrak D_X^m q\in\mathsf H_\Omega$ and, as $h\to0$,
	\begin{alignat}{2}
		Z_h&\rightharpoonup \mathfrak D_X^m q
		\qquad&&\text{in }L^2(\Omega\times\cD;H_\sigma^1),
		\label{eq:complete-flow-graph-closure-bulk}\\
		\mathsf T_\Omega Z_h&\rightharpoonup
		\mathsf T_\Omega\mathfrak D_X^m q
		\qquad&&\text{in }L^2(\Omega\times\cD;H_\sigma^{-1}),
		\label{eq:complete-flow-graph-closure-transport}\\
		(Z_h)_{\gamma_\Omega^\pm}&\rightharpoonup
		(\mathfrak D_X^m q)_{\gamma_\Omega^\pm}
		\qquad&&\text{in }L^2(\gamma_\Omega^\pm),
		\label{eq:complete-flow-graph-closure-trace}
	\end{alignat}
	where the traces in \eqref{eq:complete-flow-graph-closure-trace} are the
	completion-defined Green traces.  If $\Omega=(0,T)$, moreover
	\begin{equation}
		\label{eq:complete-flow-graph-closure-endpoints}
		Z_h(t)\rightharpoonup(\mathfrak D_X^m q)(t)
		\quad\text{in }L^2(\cD\times\bR_v^3),
		\qquad 0\le t\le T.
	\end{equation}
	In particular, $Z_h\rightharpoonup\mathfrak D_X^mq$ in $\mathsf H_\Omega$.
\end{lemma}

\begin{proof}Fix a sequence $h_j\to0$. By the graph bound, a subsequence of $Z_{h_j}$ converges weakly in the Hilbert space $\mathsf H_\Omega$ to an element $F$. Its bulk coordinate is $\mathfrak D_X^mq$: indeed, the continuous embedding into $L^2(\mathsf Q_\Omega)$ and Lemma~\ref{lem:adjoint-complete-flow-differences} identify this coordinate in distributions. Thus $\mathfrak D_X^mq$ belongs to $\mathsf H_\Omega$.
	
	The bulk coordinate uniquely determines all the other graph coordinates. Its transport derivative is fixed distributionally. For completeness, uniqueness of the traces follows by taking the difference of two Green identities with the same bulk function. Tests supported in arbitrarily small non-grazing patches of $\gamma_\Omega^+$ and $\gamma_\Omega^-$ identify the corresponding $L^2(|n\cdot v|\,dS\,dv)$ traces. The grazing set carries zero flux measure. In the evolutionary case, tests supported near either time endpoint identify the two $L^2$ time traces. This argument uses only the Green identity and the continuity of its trace maps, not a differentiability assumption on $q$.
	
	Consequently every weakly convergent subsequence has the same limit. Weak convergence in $\mathsf H_\Omega$ follows for the full family, and continuity of the transport, boundary and time-evaluation maps gives \eqref{eq:complete-flow-graph-closure-bulk}--\eqref{eq:complete-flow-graph-closure-endpoints}.
\end{proof}

\begin{lemma}[Identification of flow differences]
	\label{lem:complete-flow-weak-identification}
	Let $\Omega$ and $\mathsf Q_\Omega$ be as in
	\eqref{eq:complete-jet-closure-base}, and let $1<r<\infty$; when $\Omega=\{\ast\}$, the outer $L^r(\Omega;\cdot)$ is absent.
	Let $X=\tau_\nu$ with $\nu\in\mathfrak A$, let $1\le m\le\mathrm N$, let
	$1<p<\infty$, and let $q\in L^2(\mathsf Q_\Omega)$.  Assume
	\begin{equation}
		\label{eq:complete-flow-weak-identification-bound}
		\sup_{0<|h|<h_0}
		\|(\Delta_h^X)^m q\|_{L^r(\Omega;L_x^pH_\sigma^1)}<\infty.
	\end{equation}
	Then $\mathfrak D_X^mq\in L^r(\Omega;L_x^pH_\sigma^1)$ and
	\begin{equation}
		\label{eq:complete-flow-weak-identification-limit}
		(\Delta_h^X)^m q\rightharpoonup\mathfrak D_X^mq
		\qquad\text{in }L^r(\Omega;L_x^pH_\sigma^1)
		\quad\text{as }h\to0.
	\end{equation}
	In particular,
	\begin{equation}
		\label{eq:complete-flow-weak-identification-lsc}
		\|\mathfrak D_X^mq\|_{L^r(\Omega;L_x^pH_\sigma^1)}
		\le
		\liminf_{h\to0}
		\|(\Delta_h^X)^m q\|_{L^r(\Omega;L_x^pH_\sigma^1)}.
	\end{equation}
\end{lemma}

\begin{proof}
	Fix $h_j\to0$. Reflexivity of $L^r(\Omega;L_x^pH_\sigma^1)$ and
	\eqref{eq:complete-flow-weak-identification-bound} give a weakly convergent subsequence.
	Its limit in distributions is $\mathfrak D_X^mq$ by Lemma~\ref{lem:adjoint-complete-flow-differences}.
	Since this identifies the limit of every weakly convergent subsequence, the full family converges as in
	\eqref{eq:complete-flow-weak-identification-limit}. Weak lower semicontinuity gives
	\eqref{eq:complete-flow-weak-identification-lsc}.
\end{proof}

\begin{lemma}[Friedrichs regularization]
	\label{lem:friedrichs-complete-field}
	Let \(\nu\in\mathfrak A\), let \(\mathfrak D_\nu\) be as in
	\eqref{eq:complete-phase-space-derivative}, and let
	\(K\Subset K'\Subset\cD\times\bR^3\).  Let \(\rho\in C_c^\infty(\bR^6)\),
	\(\rho\ge0\), \(\operatorname{supp}\rho\subset\{|z|\le1\}\), \(\int\rho=1\),
	and \(\rho_\varepsilon(z):=\varepsilon^{-6}\rho(z/\varepsilon)\) for
	\(0<\varepsilon<\operatorname{dist}(K,\bR^6\setminus K')\).
	If \(u,\mathfrak D_\nu u\in L^2(K')\), then, as \(\varepsilon\to0\),
	\begin{equation}
		\label{eq:friedrichs-L2}
		\mathfrak D_\nu(\rho_\varepsilon*u)
		-\rho_\varepsilon*(\mathfrak D_\nu u)
		\longrightarrow0
		\quad\text{in }L^2(K).
	\end{equation}
	If, in addition,
	\(u,\mathfrak D_\nu u\in L_x^2H_v^1(K')\), then the convergence
	\eqref{eq:friedrichs-L2} holds in \(L_x^2H_v^1(K)\).
\end{lemma}

\begin{proof}
	By \eqref{eq:complete-phase-space-derivative}, \(\mathfrak D_\nu=\mathfrak f_\nu\cdot\nabla_{x,v}\) has bounded Lipschitz coefficients and bounded divergence on \(K'\). The Friedrichs commutator lemma~\cite[Lemma~II.1]{DiPernaLions1989} gives \eqref{eq:friedrichs-L2}. To obtain the \(H_v^1\) assertion, recall \(M_\nu=\nabla\tau_\nu\). Then
	\[
	\partial_{v_k}(\mathfrak D_\nu u)=\mathfrak D_\nu(\partial_{v_k}u)
	+(M_\nu e_k)\cdot\nabla_vu .
	\] Hence
	\(\partial_{v_k}u,\mathfrak D_\nu(\partial_{v_k}u)\in L^2(K')\), and
	\[
	\begin{aligned}
		\partial_{v_k}\bigl(\mathfrak D_\nu(\rho_\varepsilon*u)
		-\rho_\varepsilon*(\mathfrak D_\nu u)\bigr)
		={}&\mathfrak D_\nu(\rho_\varepsilon*\partial_{v_k}u)
		-\rho_\varepsilon*(\mathfrak D_\nu\partial_{v_k}u)\\
		&+(M_\nu e_k)\cdot\nabla_v(\rho_\varepsilon*u)
		-\rho_\varepsilon*\bigl((M_\nu e_k)\cdot\nabla_vu\bigr).
	\end{aligned}
	\]
	The field \((0,M_\nu e_k)\) is bounded and Lipschitz on \(K'\), and has zero divergence. Since \(M_\nu e_k\cdot\nabla_vu\in L^2(K')\), the Friedrichs lemma applies to both terms. Summing over \(k\) proves the claim.
\end{proof}

\subsection{Interpolation}\label{app:discrete-and-rungs}

The continuous, temporal, and discrete interpolation estimates below combine one-dimensional inequalities along complete flows with maximal-function bounds \cite{MazyaShaposhnikova1999,Kalamajska1994}.

For $m\ge1$ and exponents $a,p$, define
\begin{equation}
	\label{eq:GN-output-exponents}
	\frac1{\mathsf e(m;a,p)}
	:=\frac1m\cdot\frac1a+\Bigl(1-\frac1m\Bigr)\frac1p ,
	\qquad
	\mathsf e(m,p):=\mathsf e(m;2,p),
\end{equation}
with $m/0:=\infty$ and $\mathsf e(\infty;a,p):=p$.

For $2\le m\le\mathrm N$, $\mathrm L\in\bR$, $1<a\le p\le\infty$, and $j\ge0$, $k\ge1$ with $j+k\le m-1$, set $\tau=k/(m-j)$. Applying Proposition~\ref{prop:polynomial-weight-interpolation} with weights $(\mathrm L,\mathrm L+m-j)$ and exponents $(\mathsf e(m/j;a,p),a)$ gives
\begin{equation}\label{eq:output-weight-p}
	\|w^{\mathrm L+k}q\|_{L_x^{\mathsf e(m/(j+k);\,a,p)}H^1_\sigma}
	\le C\|w^{\mathrm L}q\|_{L_x^{\mathsf e(m/j;\,a,p)}H^1_\sigma}^{1-\tau}
	\|w^{\mathrm L+m-j}q\|_{L_x^{a}H^1_\sigma}^{\tau},
\end{equation}
provided that both norms on the right are finite. The interpolated weight is $\mathrm L+k$, and its spatial exponent is $p_\tau=\mathsf e(m/(j+k);a,p)$, since
\begin{equation}
	\label{eq:GN-rung-exponent-composition}
	\frac\tau a+\frac{1-\tau}{\mathsf e(m/j;a,p)}=\frac1{\mathsf e(m/(j+k);a,p)} .
\end{equation}
Indeed,
\(\tau+(1-\tau)j/m=(j+k)/m\) and \((1-\tau)(1-j/m)=1-(j+k)/m\).
The finite flow interval produces a lower-order term, which the following inequalities preserve under interpolation. For $b,G\ge0$,
\begin{equation}
	\label{eq:rung-exponent-monotonicity}
	b^{1-\alpha}G^{\alpha}\le b^{1-\beta}G^{\beta}+b ,
	\qquad 0\le\alpha\le\beta\le1.
\end{equation}
In particular, for $0\le j<i\le m$, $\tau=(i-j)/(m-j)$ and $A,B\ge0$,
\begin{equation}
	\label{eq:rung-composition-bound}
	\bigl(A^{1-j/m}B^{j/m}+A\bigr)^{1-\tau}B^{\tau}\le2A^{1-i/m}B^{i/m}+A .
\end{equation}
Indeed, by
\[
(x+y)^{1-\tau}\le x^{1-\tau}+y^{1-\tau},\quad \left(1-\frac j m\right)(1-\tau)=1-\frac i m,\quad \frac j m (1-\tau)+\tau=\frac i m,
\]
we have
\[
\bigl(A^{1-j/m}B^{j/m}+A\bigr)^{1-\tau}B^{\tau}\le A^{1-i/m}B^{i/m}+A^{1-\tau}B^{\tau},
\]
and, by $\tau\le i/m$ and \eqref{eq:rung-exponent-monotonicity},
\[
A^{1-\tau}B^{\tau}\le A^{1-i/m}B^{i/m}+A .
\]

\begin{lemma}[Temporal slicing]
	\label{lem:temporal-slicing}
	Let $0<T<\infty$, let $\cI$ be a family of ordered words of length at most $\mathrm N$ over
	$\mathfrak A$ such that $\varnothing\in\cI$ and $J\in\cI$ whenever $\nu J\in\cI$ and
	$\nu\in\mathfrak A$, and let
	$z\in\cD'((0,T)\times\cD\times\bR^3)$ satisfy
	\begin{equation}
		\label{eq:temporal-slicing-integrability}
		\mathfrak D_Iz\in L^1_{\rm loc}
		\bigl((0,T);L^2_{\rm loc}(\cD\times\bR^3)\bigr),
		\qquad I\in\cI,
	\end{equation}
	and
	\begin{equation}
		\label{eq:temporal-slicing-spacetime-chain}
		\mathfrak D_\nu(\mathfrak D_Jz)=\mathfrak D_{\nu J}z
		\quad\text{in }\cD'((0,T)\times\cD\times\bR^3),
		\qquad \nu J\in\cI .
	\end{equation}
	Then there is a set $E\subset(0,T)$ with $|(0,T)\setminus E|=0$ such that
	\begin{equation}
		\label{eq:temporal-slicing-chain}
		\mathfrak D_\nu(\mathfrak D_Jz(t))=\mathfrak D_{\nu J}z(t)
		\quad\text{in }\cD'(\cD\times\bR^3),
		\qquad \nu J\in\cI,\quad t\in E.
	\end{equation}
\end{lemma}

\begin{proof}
	Choose $L^2_{\rm loc}(\cD\times\bR^3)$-valued representatives of $\mathfrak D_Iz$ for all $I\in\cI$. For $\nu J\in\cI$,
	$\varphi\in C_c^\infty(\cD\times\bR^3)$ and $t\in(0,T)$ set
	\[
	g_\varphi(t):=\langle\mathfrak D_{\nu J}z(t),\varphi\rangle
	-\langle\mathfrak D_Jz(t),\mathfrak D_\nu^\dagger\varphi\rangle ,
	\]
	with $\mathfrak D_\nu^\dagger$ of \eqref{eq:complete-field-formal-adjoint}.  The coefficients of $\mathfrak D_\nu^\dagger$ are bounded on compact
	sets, so \eqref{eq:temporal-slicing-integrability} at $I=\nu J$ and at $I=J$ gives
	$g_\varphi\in L^1_{\rm loc}(0,T)$. Testing \eqref{eq:temporal-slicing-spacetime-chain} against $\alpha(t)\varphi(x,v)$ for $\alpha\in C_c^\infty(0,T)$ gives $\int_0^T\alpha g_\varphi\dd t=0$. Thus $g_\varphi$ vanishes outside a null set $Z_\varphi$.  Choose countable $C^1$-dense families of test functions supported in the members $K$ of a compact exhaustion of $\cD\times\bR^3$. Since $\cI$ is finite, the union of these null sets over the chosen test functions and all $\nu J\in\cI$ is still a null set $Z$. For $t\notin Z$, local $L^2$ integrability makes both pairings in $g_\varphi(t)$ continuous in the $C^1$ norm on each $C_c^\infty(K)$. They therefore agree on all test functions, proving the assertion with $E=(0,T)\setminus Z$.
\end{proof}

\subsubsection{Gagliardo--Nirenberg inequalities}\label{ss:flow-GN}

To retain the distinct radial and tangential Coulomb weights, we use the square root of the diffusion matrix. In the radial decomposition \eqref{eq:coulomb-eigen-decomposition}, set
\begin{equation}\label{eq:coulomb-square-root-symbol}
	b_\parallel:=\lambda_\parallel^{1/2},
	\qquad b_\perp:=\lambda_\perp^{1/2},
	\qquad B(v):=\sigma(v)^{1/2}
	=b_\parallel(|v|)P_v+b_\perp(|v|)Q_v.
\end{equation}
Then \eqref{eq:eigenvalue-size} and \eqref{eq:eigenvalue-log-derivatives} of
Lemma~\ref{lem:coulomb-symbol-bounds} imply
\begin{align}
	b_\parallel(r)\simeq\la r\ra^{-3/2},
	\qquad
	b_\perp(r)\simeq\la r\ra^{-1/2},
	\label{eq:sqrt-eigenvalue-size}
\end{align}
and
\begin{align}
	\abs{(r\p_r)^kb_\parallel(r)}\le Cb_\parallel(r),
	\qquad
	\abs{(r\p_r)^kb_\perp(r)}\le Cb_\perp(r),
	\qquad k\ge1.
	\label{eq:sqrt-eigenvalue-log}
\end{align}
Since $\sigma$ is smooth and positive definite by Lemma~\ref{lem:coulomb-symbol-bounds}, both $B$ and $B^{-1}$ extend smoothly across $v=0$.

Introduce the four-component Coulomb frame
\begin{equation}\label{eq:frame-def}
	\mathbb A q
	:=\begin{pmatrix}
		B(v)\nabla_vq\\
		w(v)^{-1/2}q
	\end{pmatrix}
	\in\bR^4.
\end{equation}
Then, by \eqref{eq:flow-weighted-coercive-equivalence} at $\mathrm L=0$,
\begin{equation}\label{eq:frame-exact-norm}
	\|\mathbb A q\|_{L_v^2(\bR^4)}^2
	=\int_{\bR^3}\bigl(\sigma_{ij}\partial_iq\partial_j q
	+w^{-1}|q|^2\bigr)\,\dd v
	\simeq \|q\|_{H^1_\sigma}^2.
\end{equation}

As $w\ge1$, the same equivalence gives, for real $\mathrm L'\le\mathrm L$,
\begin{equation}
	\|w^{\mathrm L'}\mathbb Aq\|_{L_v^2}
	=\left(\int w^{2\mathrm L'}\bigl(\sigma_{ij}\partial_iq\,\partial_j q
	+w^{-1}|q|^2\bigr)\dd v\right)^{1/2}
	\le C\|w^{\mathrm L}q\|_{H^1_\sigma}.
	\label{eq:weighted-frame-monotonicity}
\end{equation}

\begin{lemma}[Frame commutator]\label{lem:first-frame-commutator}
	Let $X\in\cX_{\rm c}$, and let $B$ and $\mathbb A$ be given by
	\eqref{eq:coulomb-square-root-symbol} and \eqref{eq:frame-def}.  There is a
	$4\times4$ matrix $R_1(x,v)$ such that, for every smooth $q$,
	\begin{equation}\label{eq:first-relative-commutator}
		[\mathfrak D_X,\mathbb A]q=R_1(x,v)\mathbb A q .
	\end{equation}
	Moreover, for $0\le k\le\mathrm N$,
	\begin{equation}\label{eq:R1-C-derivatives}
		\abs{\mathfrak D_X^kR_1(x,v)}_{\rm op}
		\le C\la v\ra.
	\end{equation}
\end{lemma}

\begin{proof}
	With $M=\nabla X$, differentiation in $v$ gives
	\[
	\mathfrak D_X(\nabla_vq)=\nabla_v(\mathfrak D_Xq)-M^T\nabla_vq .
	\]  Since $B$
	depends only on $v$,
	\[
	\mathfrak D_X(B\nabla_vq)
	=\bigl((Mv\cdot\nabla_v)B\bigr)\nabla_vq+B\,\mathfrak D_X(\nabla_vq)
	=\bigl((Mv\cdot\nabla_v)B-BM^T\bigr)\nabla_vq+B\nabla_v(\mathfrak D_Xq),
	\]
	that is,
	\begin{equation*}
		[\mathfrak D_X,B\nabla_v]q
		=\bigl((Mv\cdot\nabla_v)B-BM^T\bigr)\nabla_vq .
	\end{equation*}
	The derivative block therefore factors through $B\nabla_vq$ with coefficient
	\begin{equation}\label{eq:R1D}
		R_{1,D}:=\bigl((Mv\cdot\nabla_v)B-BM^T\bigr)B^{-1}.
	\end{equation}
	For the multiplication block $a_0(v)=w(v)^{-1/2}$,
	\[
	[\mathfrak D_X,a_0]q=(Mv\cdot\nabla_va_0)q
	=R_{1,0}a_0q,
	\qquad
	R_{1,0}:=Mv\cdot\nabla_v\log a_0=-\frac12\la v\ra^{-2}\,v\cdot Mv.
	\]
	Hence \eqref{eq:first-relative-commutator} holds with $R_1:=\operatorname{diag}(R_{1,D},R_{1,0})$.
	
	For $r=|v|\ge1$, put $D=r\partial_r$. The bounds \eqref{eq:sqrt-eigenvalue-size}--\eqref{eq:sqrt-eigenvalue-log}, the quotient rule and the identity $r^j\partial_r^j=D(D-1)\cdots(D-j+1)$ yield, for every $j\ge0$,
	\[
	|\partial_r^j(b_\parallel^{\pm1})|\le Cr^{\mp3/2-j},
	\qquad
	|\partial_r^j(b_\perp^{\pm1})|\le Cr^{\mp1/2-j}. \] Moreover, $P_v$ and $Q_v$ are homogeneous of degree zero in $v$, so that
	$|\partial_v^\beta P_v|+|\partial_v^\beta Q_v|\le C|v|^{-|\beta|}$ for $|v|\ge1$.
	Since $B=b_\parallel P_v+b_\perp Q_v$ and
	$B^{-1}=b_\parallel^{-1}P_v+b_\perp^{-1}Q_v$ are smooth, the chain and
	Leibniz rules give, for every multi-index $\beta$,
	\[
	|\partial_v^\beta B|_{\rm op}\le C\la v\ra^{-1/2-|\beta|},
	\qquad
	|\partial_v^\beta B^{-1}|_{\rm op}\le C\la v\ra^{3/2-|\beta|}.
	\]
	Since $|\partial_x^\alpha\partial_v^\beta(Mv)|\le C\la v\ra^{1-|\beta|}$ for
	$|\alpha|\le\mathrm N$, the Leibniz rule applied to \eqref{eq:R1D} and to
	$R_{1,0}$ gives, for $|\alpha|\le\mathrm N$ and
	$|\beta|\le\mathrm N$,
	\begin{equation}\label{eq:R1-symbol-class}
		|\partial_x^\alpha\partial_v^\beta R_1(x,v)|_{\rm op}
		\le C\la v\ra^{1-|\beta|},
		\qquad (x,v)\in\cD\times\bR^3 .
	\end{equation}
	
	The symbol class in \eqref{eq:R1-symbol-class} is stable under $\mathfrak D_X$ at the cost of one derivative. Indeed, the bounds on $Mv$ and $\partial_x^\alpha X$ give, for $1\le n\le\mathrm N$ and any matrix $F$,
	\[
	|\partial_x^\alpha\partial_v^\beta F|_{\rm op}\le C\la v\ra^{1-|\beta|}\ \ (|\alpha|+|\beta|\le n)
	\implies
	|\partial_x^\alpha\partial_v^\beta\mathfrak D_XF|_{\rm op}\le C'\la v\ra^{1-|\beta|}\ \ (|\alpha|+|\beta|\le n-1),
	\]
	with $C'$ depending only on $C$ and $X$.
	Applying this $k\le\mathrm N$ times to \eqref{eq:R1-symbol-class}, starting from $n=\mathrm N$,
	proves \eqref{eq:R1-C-derivatives}.
\end{proof}

Define matrices $R_k$ recursively by
\begin{equation}\label{eq:R-recursion}
	R_0:=I,
	\qquad
	R_{k+1}:=\mathfrak D_XR_k+R_kR_1,
	\qquad k\ge0 .
\end{equation}
Here $\mathfrak D_X$ acts entrywise and $R_1$ is given by Lemma~\ref{lem:first-frame-commutator}. The factor order in $R_kR_1$ is essential.

\begin{lemma}[Iterated commutator expansion]
	\label{lem:iterated-frame-commutator}
	Let $X\in\cX_{\rm c}$ and let $R_k$ be given by \eqref{eq:R-recursion}.  For
	$1\le k\le\mathrm N$,
	\begin{equation}\label{eq:ad-frame}
		\operatorname{ad}_{\mathfrak D_X}^k(\mathbb A)=R_k\mathbb A,
	\end{equation}
	and
	\begin{equation}\label{eq:Rk-bound}
		|R_k(x,v)|_{\rm op}\le C\la v\ra^k.
	\end{equation}
	Moreover, for $0\le k\le\mathrm N$ and every smooth $u$,
	\begin{equation}
		\mathfrak D_X^k(\mathbb Au)
		=\sum_{j=0}^{k}\binom kj\,R_j\,\mathbb A\,\mathfrak D_X^{k-j}u .
		\label{eq:app-frame-expansion}
	\end{equation}
	If $0\le k\le\mathrm N$ and $u\in L^1_{\rm loc}(\cD\times\bR^3)$, then
	\eqref{eq:app-frame-expansion} holds in $\cD'(\cD\times\bR^3;\bR^4)$, with
	$\mathfrak D_X^{k-j}u$ taken in $\cD'(\cD\times\bR^3)$.
\end{lemma}

\begin{proof}
	If $\operatorname{ad}_{\mathfrak D_X}^k(\mathbb A)=R_k\mathbb A$, then
	\[
	\begin{aligned}
		\operatorname{ad}_{\mathfrak D_X}^{k+1}(\mathbb A)
		&=[\mathfrak D_X,R_k\mathbb A]\\
		&=(\mathfrak D_XR_k)\mathbb A+R_k[\mathfrak D_X,\mathbb A]\\
		&=(\mathfrak D_XR_k+R_kR_1)\mathbb A.
	\end{aligned}
	\]
	The first identity follows by induction from \eqref{eq:first-relative-commutator}.
	
	To control the coefficient growth while retaining the derivatives needed in the induction, we prove the stronger bounds
	\begin{equation*}
		|\mathfrak D_X^aR_k|_{\rm op}\le C\la v\ra^k,
		\qquad a+k\le\mathrm N .
	\end{equation*}
	At $k=1$ this is \eqref{eq:R1-C-derivatives}.  If it holds at $k$, the
	recursion and the Leibniz rule give
	\[
	\mathfrak D_X^aR_{k+1}
	=\mathfrak D_X^{a+1}R_k
	+\sum_{i=0}^{a}\binom ai(\mathfrak D_X^iR_k)(\mathfrak D_X^{a-i}R_1),
	\]
	and for $a+k+1\le\mathrm N$,
	\[
	|\mathfrak D_X^aR_{k+1}|_{\rm op}\le|\mathfrak D_X^{a+1}R_k|_{\rm op}+\sum_{i=0}^{a}\binom ai|\mathfrak D_X^iR_k|_{\rm op}|\mathfrak D_X^{a-i}R_1|_{\rm op}\le C\la v\ra^{k}+C\la v\ra^{k}\la v\ra\le C\la v\ra^{k+1}.
	\]
	Taking $a=0$ yields \eqref{eq:Rk-bound}.
	
	Finally, \eqref{eq:app-frame-expansion} holds at $k=0$ because $R_0=I$.  If it holds at
	$k$, applying $\mathfrak D_X$ to each summand gives, by \eqref{eq:first-relative-commutator} and
	\eqref{eq:R-recursion},
	\[
	\mathfrak D_X\bigl(R_j\mathbb A\mathfrak D_X^{k-j}u\bigr)
	=R_{j+1}\mathbb A\mathfrak D_X^{k-j}u
	+R_j\mathbb A\mathfrak D_X^{k+1-j}u ,
	\]
	and Pascal's rule collects the two
	resulting sums into \eqref{eq:app-frame-expansion} at $k+1$.
	
	For the distributional assertion, the matrices $R_j$, $j\le\mathrm N$, belong to $C^{K_\cD+1-j}\subset C^{\mathrm N+1}$ by \eqref{eq:R-recursion},
	$\nabla X\in C^{K_\cD}$ and \eqref{eq:order-and-domain-regularity}, and so do the
	coefficients of $\mathfrak D_X$ and $\mathbb A$.  Hence
	\eqref{eq:first-relative-commutator} and the product rule
	$\mathfrak D_X(RF)=(\mathfrak D_XR)F+R\,\mathfrak D_XF$, which hold for smooth $q$
	and $F$, hold by duality for distributions $q$ of order at most $\mathrm N-1$ and
	$F$ of order at most $\mathrm N$.  In the step from $k<\mathrm N$ they are applied
	to $q=\mathfrak D_X^{k-j}u$ and $F=\mathbb Aq$, which have these orders since $u$ is
	locally integrable; the induction above can therefore be carried out in $\cD'$.
\end{proof}

Fix $s_0$ as in Lemma~\ref{lem:complete-flow}.  For a field
$X\in\cX_{\rm c}$ with spatial flow $\Phi_s=\Phi_s^X$ and a scalar function $\varphi$
on $\cD$, define
\[
\cM_X\varphi(x)
:=\sup_{0<r\le s_0}\frac1{2r}\int_{-r}^r|\varphi(\Phi_sx)|\dd s.
\]

\begin{lemma}[Flow maximal estimate]
	\label{lem:flow-maximal}
	Let $1<p\le\infty$.
	\begin{enumerate}
		\item\label{it:flow-maximal-continuous}
		Let $X\in\cX_{\rm c}$ and $\varphi\in L^p(\cD)$.  Then
		\begin{equation}\label{eq:flow-maximal-Lp}
			\|\cM_X\varphi\|_{L^p(\cD)}
			\le C\norm\varphi_{L^p(\cD)}.
		\end{equation}
		\item\label{it:flow-maximal-discrete}
		For $X\in\cX_{\rm c}$, $h\ne0$ and a scalar
		function $\varphi$ on $\cD$, set
		\[
		\cM_{X,h}^{\rm d}\varphi(x):=
		\sup_{1\le n\le\lfloor s_0/|h|\rfloor}
		\frac1{2n+1}\sum_{|j|\le n}|\varphi(\Phi_{jh}^Xx)|.
		\]
		If $0<|h|\le s_0$ and
		$\varphi\in L^p(\cD)$, then
		\begin{equation}
			\label{eq:discrete-maximal}
			\|\cM_{X,h}^{\rm d}\varphi\|_{L_x^p}
			\le C\|\varphi\|_{L_x^p},
		\end{equation}
		where $C$ is independent of $h$.
	\end{enumerate}
\end{lemma}

\begin{proof}
	Let $M_{\bR}$ and $M_{\mathbb Z}$ be the centered Hardy--Littlewood maximal
	operators on $\bR$ and on $\mathbb Z$; by the Hardy--Littlewood maximal theorem,
	$$
	\|M_{\bR}F\|_{L^p(\bR)}\le C_p\|F\|_{L^p(\bR)}\quad \text{and}\quad \|M_{\mathbb Z}F\|_{\ell^p(\mathbb Z)}\le C_p\|F\|_{\ell^p(\mathbb Z)}\qquad \text{for all }1<p\le\infty.
	$$
	
	\emph{(\ref{it:flow-maximal-continuous}).}  By the group law,
	\[
	\cM_X\varphi(\Phi_tx)\le
	M_{\bR}\bigl[\varphi(\Phi_\cdot x)\mathbf 1_{|\cdot|<2s_0}\bigr](t)\qquad\text{for } |t|\le s_0 .
	\]
	For $p<\infty$, averaging the $p$-th power over $|t|\le s_0$ and using the Jacobian bounds in \eqref{eq:complete-flow-uniform-bounds}, Fubini's theorem and the $L^p(\bR)$ bound for $M_{\bR}$, we obtain
	\[
	\|\cM_X\varphi\|_{L_x^p}^p
	\le C\int_\cD\int_{-s_0}^{s_0}|\cM_X\varphi(\Phi_tx)|^p\dd t\dd x
	\le C\int_\cD\int_{-2s_0}^{2s_0}|\varphi(\Phi_tx)|^p\dd t\dd x
	\le C\norm\varphi_{L_x^p}^p .
	\]
	For $p=\infty$, $\Phi_t$ maps null sets to null sets by
	\eqref{eq:complete-flow-uniform-bounds}, so, by Fubini's theorem, for almost every $x$,
	$|\varphi(\Phi_tx)|\le\|\varphi\|_{L^\infty}$ for almost every $|t|<2s_0$, and
	$\cM_X\varphi(x)\le\|\varphi\|_{L^\infty}$.
	
	\emph{(\ref{it:flow-maximal-discrete}).}  We may take $h>0$, since
	$\cM_{X,-h}^{\rm d}=\cM_{X,h}^{\rm d}$.  The argument of (\ref{it:flow-maximal-continuous}), with
	$M_{\mathbb Z}$ applied to
	$k\mapsto\varphi(\Phi_{kh+\tau}^Xx)\mathbf 1_{|kh+\tau|<2s_0}$ and integrated over
	$\tau\in[0,h)$, gives \eqref{eq:discrete-maximal}. The constant is independent of $h$, since summing over $k\in\mathbb Z$ the integrals over $\tau\in[0,h)$ gives integrals over $t=kh+\tau\in\bR$.
\end{proof}

For $m\ge2$ choose a smooth kernel
\begin{equation}
	\label{eq:Hilbert-interpolation-mollifier}
	\eta^{[m]}\in C_c^\infty((-1,1)),
	\qquad
	\int_{-1}^1\eta^{[m]}(s)\,\dd s=1,
	\qquad
	\int_{-1}^1s^n\eta^{[m]}(s)\,\dd s=0
	\quad(1\le n\le m-2),
\end{equation}
and write $\eta^{[m]}_r(s):=r^{-1}\eta^{[m]}(s/r)$.  These moment conditions can be imposed on a finite linear combination of translated smooth bumps.

We use the following form of the Landau--Kolmogorov interpolation inequality.

\begin{lemma}[Landau--Kolmogorov interpolation]
	\label{lem:Hilbert-one-dimensional}
	Let $H$ be a Hilbert space, $r>0$, $m\ge2$ and $1\le i\le m-1$.  If
	$U\in W^{m,1}((-r,r);H)$, then
	\begin{equation}\label{eq:Hilbert-general}
		\|U^{(i)}(0)\|_H
		\le C\left\{r^{-i}A_r(U)+r^{m-i}B_r^{(m)}(U)\right\},
	\end{equation}
	where
	\[
	A_r(U)=\fint_{-r}^r\|U(s)\|_H\dd s,
	\qquad
	B_r^{(m)}(U)=\fint_{-r}^r\|U^{(m)}(s)\|_H\dd s,
	\]
	and $U^{(i)}(0)$ is the value at $0$ of the continuous representative of $U^{(i)}$.
	Consequently, if $A,B\ge0$ and $U\in W^{m,1}((-s_0,s_0);H)$ satisfies
	$A_r(U)\le A$ and $B_r^{(m)}(U)\le B$ for $0<r\le s_0$, then
	\begin{equation}\label{eq:finite-window-optimized}
		\|U^{(i)}(0)\|_H\le C\{A^{1-i/m}B^{i/m}+A\}.
	\end{equation}
	The constant $C$ is independent of $r$, $A$, $B$, $U$ and $H$.
\end{lemma}

\begin{proof}
	Integration by parts $i$ times gives
	\[
	\int_{-r}^r\eta^{[m]}_r(s)U^{(i)}(s)\dd s
	=(-1)^i\int_{-r}^r\bigl(\eta^{[m]}_r\bigr)^{(i)}(s)U(s)\dd s .
	\]
	Since $\|(\eta^{[m]}_r)^{(i)}\|_{L^\infty}\le Cr^{-1-i}$,
	\begin{equation}\label{eq:smoothed-derivative-bound}
		\Bigl\|\int\eta^{[m]}_r(s)U^{(i)}(s)\dd s\Bigr\|_H
		\le Cr^{-1-i}\int_{-r}^r\|U(s)\|_H\dd s
		=2Cr^{-i}A_r(U).
	\end{equation}
	Taylor's formula in $H$ to order $m-i-1$ with integral remainder reads
	\begin{equation}\label{eq:Taylor-Uprime}
		U^{(i)}(s)
		=\sum_{n=0}^{m-i-1}\frac{s^n}{n!}U^{(i+n)}(0)
		+\int_0^s\frac{(s-t)^{m-i-1}}{(m-i-1)!}U^{(m)}(t)\dd t.
	\end{equation}
	Integrate \eqref{eq:Taylor-Uprime} against $\eta^{[m]}_r$: by
	\eqref{eq:Hilbert-interpolation-mollifier} the terms of degree
	$1\le n\le m-i-1\le m-2$ vanish and the constant term is reproduced, so
	\[
	\Bigl\|U^{(i)}(0)-\int\eta^{[m]}_r(s)U^{(i)}(s)\dd s\Bigr\|_H
	\le Cr^{m-i-1}\int_{-r}^r\|U^{(m)}(t)\|_H\dd t
	=2Cr^{m-i}B_r^{(m)}(U).
	\]
	Together with \eqref{eq:smoothed-derivative-bound}, this proves \eqref{eq:Hilbert-general}.
	
	For \eqref{eq:finite-window-optimized}, \eqref{eq:Hilbert-general} gives
	\[
	Y:=\|U^{(i)}(0)\|_H\le C(r^{-i}A+r^{m-i}B)\qquad\text{for } 0<r\le s_0 .
	\]
	If $A=0$, let
	$r\downarrow0$; then $Y=0$.  If $A>0$ and $B\le As_0^{-m}$, take $r=s_0$; then $Y\le2Cs_0^{-i}A$.
	Otherwise, put $r=(A/B)^{1/m}$, so that 
	$$
	r<s_0\quad \text{and}\quad r^{-i}A=r^{m-i}B=A^{1-i/m}B^{i/m}.
	$$
	Then $Y\le2CA^{1-i/m}B^{i/m}$.
\end{proof}

For $2\le m\le\mathrm N$, $\mathrm L\in\bR$ and exponents
\begin{equation}
	\label{eq:GN-slots}
	1<a\le p\le\infty ,
\end{equation}
set
\begin{equation}
	\label{eq:GN-lower-order}
	\cG^{\mathrm L}_{m,a}(u):=\sum_{|J|\le m}
	\bigl\|w^{\mathrm L+m-|J|}\mathfrak D_Ju\bigr\|_{L_x^{a}H^1_\sigma},
	\qquad \mathfrak D_\varnothing=\operatorname{Id},
\end{equation}
where the sum is taken over the ordered complete words $J$ with $|J|\le m$. In the plain fiber norm, set
\begin{equation}
	\label{eq:GN-lower-order-gaussian}
	\cG^{\sharp,\mathrm L}_{m,a}(u):=\sum_{|J|\le m}
	\bigl\|w^{\mathrm L+m-|J|}\mathfrak D_Ju\bigr\|_{L_x^{a}L_v^2}.
\end{equation}

\begin{thm}[Complete-flow Gagliardo--Nirenberg inequalities]
	\label{thm:complete-flow-GN}
	Let $X\in\cX_{\rm c}$, $2\le m\le\mathrm N$, $1\le i\le m-1$ and
	$\mathrm L\in\bR$, let $a,p$ satisfy \eqref{eq:GN-slots}, and let $\mathsf e$ be as in \eqref{eq:GN-output-exponents}.
	\begin{enumerate}
		\item\label{it:GN-coercive}
		Let $w^{\mathrm L}u\in L_x^pH^1_\sigma$ and $\cG^{\mathrm L}_{m,a}(u)<\infty$,
		with $\cG^{\mathrm L}_{m,a}$ of \eqref{eq:GN-lower-order}.  Then
		\begin{equation}
			\label{eq:complete-flow-GN}
			\|w^{\mathrm L}\mathfrak D_X^{\,i}u\|
			_{L_x^{\mathsf e(m/i;\,a,p)}H^1_\sigma}
			\le C\|w^{\mathrm L}u\|_{L_x^pH^1_\sigma}^{1-i/m}\,
			\cG^{\mathrm L}_{m,a}(u)^{i/m}
			+C\|w^{\mathrm L}u\|_{L_x^pH^1_\sigma}.
		\end{equation}
		\item\label{it:GN-fiber}
		Let $w^{\mathrm L}u\in L_x^pL_v^2$ and $\cG^{\sharp,\mathrm L}_{m,a}(u)<\infty$,
		with $\cG^{\sharp,\mathrm L}_{m,a}$ of \eqref{eq:GN-lower-order-gaussian}.
		Then
		\begin{equation}
			\label{eq:complete-flow-GN-gaussian}
			\|w^{\mathrm L}\mathfrak D_X^{\,i}u\|_{L_x^{\mathsf e(m/i;\,a,p)}L_v^2}
			\le C\|w^{\mathrm L}u\|_{L_x^pL_v^2}^{1-i/m}
			\cG^{\sharp,\mathrm L}_{m,a}(u)^{i/m}
			+C\|w^{\mathrm L}u\|_{L_x^pL_v^2}.
		\end{equation}
	\end{enumerate}
\end{thm}
\begin{proof}
	\emph{(\ref{it:GN-coercive}).}
	Let $\mathbb A$ be the frame of \eqref{eq:frame-def}, write $\widehat\Phi_s(x,v)=(\Phi_sx,A_s(x)v)$ for the complete phase flow, and set $u_{\mathrm L}=w^{\mathrm L}u$.  We proceed by induction on $i$. For $i=0$, the assertion follows from
	$\mathsf e(\infty;a,p)=p$. 
	
	We first commute the polynomial weight through complete derivatives. By Lemma~\ref{lem:weight-leibniz} at $\ell=\pm\mathrm L$, the functions
	$w^{\mp\mathrm L}\mathfrak D_X^{\,l}w^{\pm\mathrm L}$, $l\le\mathrm N$, are bounded in
	$C^1(\bR_v^3)$ uniformly in $x$. The Leibniz rule and \eqref{eq:Hsigma-multiplier} therefore give, for $j,k\ge0$,
	\begin{equation}
		\label{eq:GN-weight-commutation}
		\bigl\|\la v\ra^{k}\mathfrak D_X^{\,j}u_{\mathrm L}\bigr\|_{H^1_\sigma}
		\le C\sum_{n=0}^{j}
		\bigl\|w^{\mathrm L+k}\mathfrak D_X^{\,n}u\bigr\|_{H^1_\sigma},
	\end{equation}
	and
	\[
	\bigl\|w^{\mathrm L+k}\mathfrak D_X^{\,j}u\bigr\|_{H^1_\sigma}
	\le C\sum_{n=0}^{j}
	\bigl\|\la v\ra^{k}\mathfrak D_X^{\,n}u_{\mathrm L}\bigr\|_{H^1_\sigma}.
	\]
	
	For the moment suppose that $u$ is smooth. Along the complete phase flow consider the $L_v^2(\bR^4)$-valued curve
	\[
	V_x(s):=(\mathbb Au_{\mathrm L})\circ\widehat\Phi_s(x,\cdot).
	\]
	The complete-flow identity \eqref{eq:complete-flow-level-set} and Lemma~\ref{lem:iterated-frame-commutator} yield, for $0\le j\le m$,
	\[
	V_x^{(j)}=\Bigl(\sum_{k=0}^{j}\binom jk
	R_k\,\mathbb A\,\mathfrak D_X^{\,j-k}u_{\mathrm L}\Bigr)\circ\widehat\Phi_s ,
	\qquad |R_k(x,v)|_{\rm op}\le C\la v\ra^{k}.
	\]
	Put
	\[
	\varphi(x):=\|u_{\mathrm L}(x)\|_{H^1_\sigma},
	\qquad
	\psi(x):=\sum_{j=0}^{m}
	\bigl\|\la v\ra^{\,m-j}\mathfrak D_X^{\,j}u_{\mathrm L}(x)
	\bigr\|_{H^1_\sigma}.
	\]
	Changing variables by $z=A_s(x)v$ and using the Jacobian bound in \eqref{eq:complete-flow-uniform-bounds}, Lemma~\ref{lem:real-weight-Hsigma-equivalence}, and \eqref{eq:weighted-frame-monotonicity} with $(\mathrm L',\mathrm L)=(k,m-j+k)$ for each $R_k$-term in $V_x^{(j)}$, $j\ge1$, gives
	\[
	\|V_x(s)\|_{L_v^2}\le C\varphi(\Phi_sx),
	\qquad
	\|V_x^{(j)}(s)\|_{L_v^2}\le C\psi(\Phi_sx),
	\qquad 1\le j\le m .
	\]
	By \eqref{eq:GN-weight-commutation} and
	\eqref{eq:real-weight-Hsigma-monotonicity},
	\begin{equation}
		\label{eq:GN-b-bound}
		\|\psi\|_{L_x^{a}}\le C\,\cG^{\mathrm L}_{m,a}(u).
	\end{equation}
	
	Consequently, the flow maximal operator controls the one-dimensional averages: 
	$$
	A_r(V_x)\le C\cM_X\varphi(x)\quad \text{and}\quad B_r^{(m)}(V_x)\le C\cM_X\psi(x)\qquad \text{for }0<r\le s_0.
	$$
	Lemma~\ref{lem:Hilbert-one-dimensional}, applied at order $i$ to $V_x$, then gives, pointwise in $x$,
	\begin{equation}
		\label{eq:complete-flow-one-dimensional}
		\|V_x^{(i)}(0)\|_{L_v^2}
		\le C(\cM_X\varphi)^{1-i/m}(\cM_X\psi)^{i/m}+C\cM_X\varphi .
	\end{equation}
	
	To recover the desired derivative from the derivative of the transported frame, use the commutator expansion at $s=0$:
	\[
	\mathbb A\mathfrak D_X^{\,i}u_{\mathrm L}
	=V_x^{(i)}(0)-\sum_{k=1}^{i}\binom ik
	R_k\,\mathbb A\,\mathfrak D_X^{\,i-k}u_{\mathrm L} .
	\]
	Fix $1\le k\le i$ and put $j=i-k$, $\tau=k/(m-j)$. Since $j+k=i\le m-1$, the $k$ velocity weights lost in $R_k$ can be interpolated between the $j$th rung and its top-weight counterpart. Indeed, $|R_k|_{\rm op}\le C\la v\ra^k$, \eqref{eq:weighted-frame-monotonicity}, the weight-commutation estimates, and \eqref{eq:output-weight-p} applied to $q=\mathfrak D_X^{\,j}u$ give
	\begin{align*}
		\begin{aligned}
			\bigl\|R_k\mathbb A\mathfrak D_X^{\,j}u_{\mathrm L}\bigr\|
			_{L_x^{\mathsf e(m/i;\,a,p)}L_v^2}
			&\le C\bigl\|w^{\mathrm L}\mathfrak D_X^{\,j}u\bigr\|
			_{L_x^{\mathsf e(m/j;\,a,p)}H^1_\sigma}^{1-\tau}
			\bigl\|w^{\mathrm L+m-j}\mathfrak D_X^{\,j}u\bigr\|
			_{L_x^{a}H^1_\sigma}^{\tau}\\
			&\quad+C\sum_{n<j}\bigl\|w^{\mathrm L+k}\mathfrak D_X^{\,n}u\bigr\|
			_{L_x^{\mathsf e(m/i;\,a,p)}H^1_\sigma}.
		\end{aligned}
	\end{align*}
	For the remainder of the induction, abbreviate $A=\|w^{\mathrm L}u\|_{L_x^pH^1_\sigma}$ and $B=\cG^{\mathrm L}_{m,a}(u)$. Thus
	\[
	\|w^{\mathrm L+m-j}\mathfrak D_X^{\,j}u\|_{L_x^{a}H^1_\sigma}\le B\quad \text{and}\quad  \|w^{\mathrm L}\mathfrak D_X^{\,j}u\|_{L_x^{\mathsf e(m/j;\,a,p)}H^1_\sigma}\le C(A^{1-j/m}B^{j/m}+A).
	\]
	For $j=0$ the latter bound is the definition of $A$; for $1\le j<i$ it is the induction hypothesis. As $\tau=(i-j)/(m-j)$, the interpolation inequality \eqref{eq:rung-composition-bound} yields
	\[
	C\bigl\|w^{\mathrm L}\mathfrak D_X^{\,j}u\bigr\|_{L_x^{\mathsf e(m/j;\,a,p)}H^1_\sigma}^{1-\tau}\bigl\|w^{\mathrm L+m-j}\mathfrak D_X^{\,j}u\bigr\|_{L_x^{a}H^1_\sigma}^{\tau}
	\le C(A^{1-i/m}B^{i/m}+A).
	\]
	The weight-commutation remainder contains derivatives of order $n<j$. Applying the same interpolation with $(n+k,n)$ instead of $(i,j)$ yields
	\[
	\|w^{\mathrm L+k}\mathfrak D_X^{\,n}u\|_{L_x^{\mathsf e(m/(n+k);\,a,p)}H^1_\sigma}\le C(A^{1-(n+k)/m}B^{(n+k)/m}+A).
	\]
	Since $n+k\le i$ and $\cD$ is bounded, $L_x^{\mathsf e(m/(n+k);\,a,p)}\subset L_x^{\mathsf e(m/i;\,a,p)}$; \eqref{eq:rung-exponent-monotonicity} then gives
	\[
	C(A^{1-(n+k)/m}B^{(n+k)/m}+A)\le C(A^{1-i/m}B^{i/m}+A).
	\]
	Thus every term in the commutator correction satisfies
	\begin{equation}\label{eq:GN-Rk-term-bound}
		\|R_k\mathbb A\mathfrak D_X^{\,j}u_{\mathrm L}\|_{L_x^{\mathsf e(m/i;\,a,p)}L_v^2}\le C(A^{1-i/m}B^{i/m}+A).
	\end{equation}
	
	For the leading term, the reciprocal-exponent identity \eqref{eq:GN-output-exponents} and H\"older's inequality in $x$ give
	\begin{align*}
		\bigl\|(\cM_X\varphi)^{1-i/m}(\cM_X\psi)^{i/m}\bigr\|_{L_x^{\mathsf e(m/i;\,a,p)}}
		\le\|\cM_X\varphi\|_{L_x^p}^{1-i/m}\|\cM_X\psi\|_{L_x^{a}}^{i/m}
		\le CA^{1-i/m}B^{i/m}.
	\end{align*}
	Here Lemma~\ref{lem:flow-maximal}\,(\ref{it:flow-maximal-continuous}) applies because $a,p>1$, and the last inequality uses $\|\varphi\|_{L_x^p}=A$ and \eqref{eq:GN-b-bound}. Combining this with \eqref{eq:complete-flow-one-dimensional}, \eqref{eq:GN-Rk-term-bound}, and \eqref{eq:frame-exact-norm} yields
	\[
	\bigl\|\mathfrak D_X^{\,i}u_{\mathrm L}\bigr\|_{L_x^{\mathsf e(m/i;\,a,p)}H^1_\sigma}
	\le C\bigl(A^{1-i/m}B^{i/m}+A\bigr)+C\|\cM_X\varphi\|_{L_x^{\mathsf e(m/i;\,a,p)}} .
	\]
	Commuting $w^{\mathrm L}$ back through $\mathfrak D_X^{\,i}$ gives
	\[
	\|w^{\mathrm L}\mathfrak D_X^{\,i}u\|_{H^1_\sigma}
	\le C\|\mathfrak D_X^{\,i}u_{\mathrm L}\|_{H^1_\sigma}
	+C\sum_{n<i}\|w^{\mathrm L}\mathfrak D_X^{\,n}u\|_{H^1_\sigma}.
	\]  Because $\cD$ is bounded and $\mathsf e(m/i;a,p)\le p$,
	\[
	\|w^{\mathrm L}u\|_{L_x^{\mathsf e(m/i;\,a,p)}H^1_\sigma}+\|\cM_X\varphi\|_{L_x^{\mathsf e(m/i;\,a,p)}}\le C\|w^{\mathrm L}u\|_{L_x^pH^1_\sigma},
	\]
	where Lemma~\ref{lem:flow-maximal}\,(\ref{it:flow-maximal-continuous}) is used for the second term; for
	$1\le n<i$,
	\[
	\|w^{\mathrm L}\mathfrak D_X^{\,n}u\|_{L_x^{\mathsf e(m/i;\,a,p)}H^1_\sigma}\le C(A^{1-i/m}B^{i/m}+A)
	\]
	by the induction hypothesis at order $n$, the embedding $L_x^{\mathsf e(m/n;\,a,p)}\subset L_x^{\mathsf e(m/i;\,a,p)}$, and \eqref{eq:rung-exponent-monotonicity}. This closes the induction for smooth $u$.
	
	It remains to justify the argument under the stated distributional hypotheses. Define
	\[
	G_j:=\sum_{k=0}^{j}\binom jk
	R_k\,\mathbb A\,\mathfrak D_X^{\,j-k}u_{\mathrm L} ,
	\qquad 0\le j\le m .
	\]
	The two norm hypotheses, the weight-commutation estimates, \eqref{eq:weighted-frame-monotonicity}, \eqref{eq:Rk-bound}, and the boundedness of $\cD$ imply that $u_{\mathrm L}$ is locally integrable and $G_0,\ldots,G_m\in L_x^1L_v^2$. The distributional assertion of
	Lemma~\ref{lem:iterated-frame-commutator} applies to $u_{\mathrm L}$ and gives
	$G_j=\mathfrak D_X^{\,j}(\mathbb Au_{\mathrm L})$ in $\cD'$, hence $\mathfrak D_XG_j=G_{j+1}$ for $j<m$. Applying Lemma~\ref{lem:distributional-flow-slicing} with $d=4$ and $\lambda=0$, we obtain for almost every $x$ a representative $V_x$ of $s\mapsto T_s^XG_0(x,\cdot)$ in $W^{m,1}((-s_0,s_0);L_v^2(\bR^4))$. Its derivatives through order $m-1$ have continuous representatives at $s=0$, with $V_x^{(i)}(0)=G_i(x)$. The maximal-function and commutator estimates already proved therefore apply to these representatives and give \eqref{eq:complete-flow-GN}.
	
	\emph{(\ref{it:GN-fiber}).} In the plain $L_v^2$ norm there is no frame commutator. Use instead the curve
	\[
	U_x(s,v):=w^{\mathrm L}(v)\,u(\Phi_sx,A_s(x)v).
	\]
	For smooth $u$, by the second identity in \eqref{eq:complete-flow-level-set},
	\[
	\partial_s^{\,j}U_x(s,v)
	=w^{\mathrm L}(v)(\mathfrak D_X^{\,j}u)(\Phi_sx,A_s(x)v),
	\qquad 0\le j\le m.
	\]
	The fixed weight in this curve requires neither the previous weight-commutation argument nor the frame correction. Set
	$$
	\varphi(x):=\|w^{\mathrm L}u(x)\|_{L_v^2}\quad \text{and}\quad \psi(x):=\|w^{\mathrm L}\mathfrak D_X^{\,m}u(x)\|_{L_v^2}.
	$$
	The change of variables $z=A_s(x)v$, \eqref{eq:complete-flow-uniform-bounds} and the
	bound $w^{\mathrm L}(A_s^{-1}z)\le Cw^{\mathrm L}(z)$ from \eqref{eq:flow-weight-comparison}
	give
	$$
	\|U_x(s)\|_{L_v^2}\le C\varphi(\Phi_sx)\quad \text{and}\quad \|\partial_s^{\,m}U_x(s)\|_{L_v^2}\le C\psi(\Phi_sx).
	$$
	Moreover, by \eqref{eq:GN-lower-order-gaussian},
	$\|\psi\|_{L_x^{a}}\le\cG^{\sharp,\mathrm L}_{m,a}(u)$. The one-dimensional inequality \eqref{eq:finite-window-optimized} applies with $\partial_s^{\,i}U_x(0)=w^{\mathrm L}\mathfrak D_X^{\,i}u(x)$; H\"older's inequality and Lemma~\ref{lem:flow-maximal}\,(\ref{it:flow-maximal-continuous}) then yield
	\eqref{eq:complete-flow-GN-gaussian}.  For general $u$,
	Lemma~\ref{lem:distributional-flow-slicing} at $d=1$, $\lambda=\mathrm L$ and $G_k=\mathfrak D_X^{\,k}u$, $0\le k\le m$,
	which belong to $L_x^1H_{\mathrm L}$ by the top-norm hypothesis and the boundedness of $\cD$, gives for almost every $x$ a $W^{m,1}((-s_0,s_0);H_{\mathrm L})$ representative of $s\mapsto T_s^Xu(x,\cdot)$, with derivatives through order $m-1$ identified at $s=0$. Multiplication by $w^{\mathrm L}$ is an isometry from $H_{\mathrm L}$ to $L_v^2$, and the preceding argument applies to this representative.
\end{proof}

The time exponent is determined by the interpolation ratio $i/m$; the spatial estimate below uses the additional margin $\theta_\circ>0$.

\begin{lemma}[Temporal interpolation]
	\label{lem:temporal-letter-unequal}
	Let $I\subset\bR$ be an interval with $|I|\ge1$, let $H$ be a Hilbert
	space continuously embedded in $L^1_{\rm loc}(\bR^3)$, let $m\ge2$ and $1\le i\le m-1$, let
	\begin{equation}
		\label{eq:temporal-unequal-hypotheses}
		2\le\rho_{\mathrm{in}}\le\infty,
		\qquad
		1\le a\le p_0\le\infty,
		\qquad
		\theta_\circ>0,
		\qquad
		\gamma:=\frac im+\theta_\circ\le1 ,
	\end{equation}
	and let $v$ on $I\times\cD$ satisfy
	$\mathfrak D_{\mathrm t}^{\,k}v\in L^1(I;L_x^1(H))$ for $0\le k\le m$, the
	derivatives taken in $\cD'(I\times\cD\times\bR^3)$.  Write
	\begin{equation}
		\label{eq:temporal-unequal-quantities}
		\mathfrak b:=\|v\|_{L_t^{\rho_{\mathrm{in}}}L_x^{p_0}H},
		\qquad
		\mathfrak b_2:=\|v\|_{L_t^{2}L_x^{a}H},
		\qquad
		\mathfrak G:=\bigl\|\mathfrak D_{\mathrm t}^{\,m}v\bigr\|_{L_t^{2}L_x^{a}H},
	\end{equation}
	where the norms in $t$ are taken over $I$. Let the output exponents be
	\begin{equation}
		\label{eq:temporal-unequal-output}
		\rho_1:=\mathsf e(m/i;2,\rho_{\mathrm{in}}),
		\qquad
		\frac1{p_1}:=\frac{1-\gamma}{p_0}+\frac{\gamma}{a} .
	\end{equation}
	Then
	\begin{equation}
		\label{eq:temporal-unequal}
		\bigl\|\mathfrak D_{\mathrm t}^{\,i}v\bigr\|_{L_t^{\rho_1}L_x^{p_1}H}
		\le C\,\mathfrak b^{1-\gamma}\bigl(\mathfrak G+\mathfrak b_2\bigr)^{\gamma}
		+C\,\mathfrak b ,
	\end{equation}
	with $C$ independent of $H$, $v$ and $|I|$.
\end{lemma}

\begin{proof}
	Testing the distributional identities against $\alpha(t)\varphi(x,v)$, with $\alpha\in C_c^\infty(I)$ and $\varphi\in C_c^\infty(\cD\times\bR^3)$, identifies the time derivatives in $\cD'(I;L_x^1(H))$. Thus $v\in W^{m,1}(I;L_x^1(H))$. We may assume that $\mathfrak b$, $\mathfrak b_2$ and $\mathfrak G$ are finite.
	
	Let $\eta:=\eta^{[m+1]}$ be the kernel
	\eqref{eq:Hilbert-interpolation-mollifier}, whose moments satisfy $\int s^n\eta=0$ for $1\le n\le m-1$. Define
	\begin{equation*}
		\theta(s):=2\eta(2s)-\eta(s),
		\qquad
		r_n:=\frac{2^{-n-1}}{m+1},
		\qquad
		\theta_n(s):=r_n^{-1}\theta(s/r_n)=\eta_{r_{n+1}}(s)-\eta_{r_n}(s),
	\end{equation*}
	for $n\ge0$. The difference $\theta$ has vanishing moments through order $m-1$, including its zeroth moment. Consequently $\theta=\Theta^{(m)}$ for some $\Theta\in C_c^\infty((-1,1))$. The kernels above are supported in $(-1,1)$; setting $\Theta_n(s):=r_n^{m-1}\Theta(s/r_n)$ gives
	$\theta_n=\Theta_n^{(m)}$ and, for $1\le q_2\le\infty$,
	\begin{equation}
		\label{eq:dyadic-kernel-norms}
		\bigl\|\theta_n^{(i)}\bigr\|_{L^{q_2}(\bR)}
		=r_n^{-i-1+1/q_2}\bigl\|\theta^{(i)}\bigr\|_{L^{q_2}(\bR)},
		\qquad
		\bigl\|\Theta_n^{(i)}\bigr\|_{L^{q_2}(\bR)}
		=r_n^{\,m-1-i+1/q_2}\bigl\|\Theta^{(i)}\bigr\|_{L^{q_2}(\bR)} .
	\end{equation}
	For every Banach space $X$ and every $V$ with
	$V,\mathfrak D_{\mathrm t}^{\,m}V\in L^1_{\rm loc}(\bR;X)$,
	\begin{equation}
		\label{eq:dyadic-two-readings}
		\theta_n^{(i)}\ast V=\Theta_n^{(i)}\ast\mathfrak D_{\mathrm t}^{\,m}V ,
	\end{equation}
	and
	\begin{equation}
		\label{eq:banach-young}
		\|K\ast V\|_{L_t^{q_1}(X)}
		\le\|K\|_{L^{q_2}(\bR)}\,\|V\|_{L_t^{q_0}(X)},
		\qquad
		1+\frac1{q_1}=\frac1{q_2}+\frac1{q_0} ,
	\end{equation}
	for $1\le q_0,q_1,q_2\le\infty$ and scalar kernels $K$, by Young's inequality for convolution.

	The condition $|I|\ge1$ allows an order-$m$ reflection extension to $\widetilde I:=\{t:\operatorname{dist}(t,I)<\tfrac1{m+1}\}$ with constants depending only on $m$ (see, \textit{e.g.}, \cite[Chapter~4, Lemma~4.1]{TaylorPDE1}). For every Banach space $X$ and $1\le\rho\le\infty$, this extension is bounded from $L^\rho(I;X)$ to $L^\rho(\widetilde I;X)$, and its $m$th time derivative is obtained by an operator of the same kind acting on $\mathfrak D_{\mathrm t}^{\,m}v$.  The three norms in \eqref{eq:temporal-unequal-quantities}, taken over $\widetilde I$, are bounded by $C\mathfrak b$, $C\mathfrak b_2$ and $C\mathfrak G$; we retain the same notation for these bounds. Since $\eta_{r_n}$, $\theta_n$ and $\Theta_n$ are supported in $(-r_0,r_0)$, with $r_0<1/(m+1)$, their convolutions at $t\in I$ sample only the extension on $\widetilde I$. Integration by parts therefore gives \eqref{eq:dyadic-two-readings} and $\eta_{r_n}\ast\mathfrak D_{\mathrm t}^{\,i}v=\eta_{r_n}^{(i)}\ast v$ at such $t$, without boundary contributions. For estimates on $I$, these convolutions agree with those of $\mathbf 1_{\widetilde I}v$ and $\mathbf 1_{\widetilde I}\mathfrak D_{\mathrm t}^{\,m}v$, so \eqref{eq:banach-young} applies.
	
	The dyadic differences telescope:
	$$
	\eta_{r_0}+\sum_{n=0}^{k}\theta_n=\eta_{r_{k+1}},\qquad \text{and}\quad  \eta_{r_{k+1}}\ast\mathfrak D_{\mathrm t}^{\,i}v\to\mathfrak D_{\mathrm t}^{\,i}v\quad \text{in }L^1_{\rm loc}(I;L_x^1(H)),
	$$
	so that
	\begin{equation}
		\label{eq:temporal-unequal-telescope}
		\mathfrak D_{\mathrm t}^{\,i}v
		=\eta_{r_0}^{(i)}\ast v+\sum_{n\ge0}\theta_n^{(i)}\ast v ,
		\qquad r_0=\tfrac1{2(m+1)} .
	\end{equation}
	
	Young's inequality and the fixed scale $r_0=1/[2(m+1)]$ give
	\[
	\bigl\|\eta_{r_0}^{(i)}\ast v\bigr\|_{L_t^{\rho_1}L_x^{p_1}H}
	\le C\|v\|_{L_t^{\rho_1}L_x^{p_1}H} .
	\]
	Since $a\le p_0$, the exponent $p_\ast:=\mathsf e(m/i;a,p_0)$ satisfies
	\[
	1/p_1=1/p_\ast+\theta_\circ(1/a-1/p_0)\ge1/p_\ast,
	\]
	and hence $L_x^{p_\ast}\subset L_x^{p_1}$ on the bounded domain $\cD$. The ordinary interpolation exponents satisfy
	\[
	\frac1{\rho_1}=\frac{1-i/m}{\rho_{\mathrm{in}}}+\frac{i/m}{2},\qquad
	\frac1{p_\ast}=\frac{1-i/m}{p_0}+\frac{i/m}{a},
	\]
	so H\"older's inequality gives
	\begin{equation}
		\label{eq:temporal-unequal-coarse}
		\|v\|_{L_t^{\rho_1}L_x^{p_\ast}H}
		\le\mathfrak b^{1-i/m}\mathfrak b_2^{\,i/m} .
	\end{equation}
	
	Set
	\begin{equation}
		\label{eq:temporal-unequal-b2}
		\frac{\gamma}{b}:=\frac1{\rho_1}-\frac{1-\gamma}{\rho_{\mathrm{in}}}
		=\frac i{2m}+\frac{\theta_\circ}{\rho_{\mathrm{in}}},
		\qquad
		\frac1{q_2}:=\frac12+\frac1b .
	\end{equation}
	Because $1/\rho_{\mathrm{in}}\le1/2$, these definitions imply $0<1/b\le1/2$ and $1\le q_2<2$, including the admissible endpoint $q_2=1$. Applying \eqref{eq:banach-young} to the two representations in \eqref{eq:dyadic-two-readings}, with exponent triples $(1,\rho_{\mathrm{in}},\rho_{\mathrm{in}})$ and $(q_2,2,b)$, and using \eqref{eq:dyadic-kernel-norms}, we obtain, for $n\ge0$,
	\[
	\bigl\|\theta_n^{(i)}\ast v\bigr\|_{L_t^{\rho_{\mathrm{in}}}L_x^{p_0}H}
	\le C\,r_n^{-i}\,\mathfrak b ,
	\qquad
	\bigl\|\theta_n^{(i)}\ast v\bigr\|_{L_t^{b}L_x^{a}H}
	\le C\,r_n^{\,m-1-i+1/q_2}\,\mathfrak G .
	\]  By \eqref{eq:temporal-unequal-output} and \eqref{eq:temporal-unequal-b2},
	\[
	\frac1{\rho_1}=\frac{1-\gamma}{\rho_{\mathrm{in}}}+\frac\gamma b,
	\qquad
	\frac1{p_1}=\frac{1-\gamma}{p_0}+\frac\gamma a,
	\]
	so H\"older's inequality gives
	\begin{equation}
		\label{eq:temporal-unequal-block}
		\bigl\|\theta_n^{(i)}\ast v\bigr\|_{L_t^{\rho_1}L_x^{p_1}H}
		\le C\,r_n^{\,\varepsilon}\,\mathfrak b^{1-\gamma}\mathfrak G^{\gamma},
		\qquad
		\varepsilon:=-(1-\gamma)i+\gamma\Bigl(m-1-i+\frac1{q_2}\Bigr).
	\end{equation}
	
	Using $1/q_2=\tfrac12+1/b$ and \eqref{eq:temporal-unequal-b2}, we obtain
	\[
	\varepsilon
	=-i+\gamma\Bigl(m-\frac12\Bigr)+\frac i{2m}+\frac{\theta_\circ}{\rho_{\mathrm{in}}},
	\]
	while $\gamma=i/m+\theta_\circ$ gives $\gamma(m-\tfrac12)=i-\tfrac i{2m}+\theta_\circ(m-\tfrac12)$. Consequently
	\begin{equation*}
		\varepsilon=\theta_\circ\Bigl(m-\frac12+\frac1{\rho_{\mathrm{in}}}\Bigr)>0 .
	\end{equation*}
	Since $\varepsilon>0$ and $r_n=2^{-n-1}/(m+1)$, \eqref{eq:temporal-unequal-block} is summable in $L_t^{\rho_1}L_x^{p_1}H$. By the telescoping identity \eqref{eq:temporal-unequal-telescope}, this sum agrees in distributions with $\mathfrak D_{\mathrm t}^{\,i}v-\eta_{r_0}^{(i)}\ast v$.   The fixed-scale estimate, the bounded-domain inclusion $L_x^{p_\ast}\subset L_x^{p_1}$ and \eqref{eq:temporal-unequal-coarse} bound the remaining low-frequency term by $C\mathfrak b^{1-i/m}\mathfrak b_2^{\,i/m}$. Combining the two estimates gives
	\[
	\bigl\|\mathfrak D_{\mathrm t}^{\,i}v\bigr\|_{L_t^{\rho_1}L_x^{p_1}H}
	\le C\,\mathfrak b^{1-i/m}\mathfrak b_2^{\,i/m}+C\,\mathfrak b^{1-\gamma}\mathfrak G^{\gamma}
	\le C\,\mathfrak b^{1-\gamma}\bigl(\mathfrak G+\mathfrak b_2\bigr)^{\gamma}+C\,\mathfrak b .
	\]
	The final inequality follows from \eqref{eq:rung-exponent-monotonicity}, since $i/m\le\gamma$.
\end{proof}

\subsubsection{Discrete interpolation}\label{ss:flow-discrete-interpolation}

\begin{lemma}[Finite-difference interpolation inequalities]
	\label{lem:discrete-Peano}
	Let $H$ be a Hilbert space, let $m\ge2$ and $1\le i\le m-1$, write
	$da_j=a_{j+1}-a_j$, and set
	\[
	M_nb:=\frac1{2n+1}\sum_{|j|\le n}\|b_j\|_H.
	\]
	There is an integer $n_0=n_0(m)\ge1$ such that, for every $n\ge n_0$ and every
	$(a_j)_{|j|\le n+m}\subset H$,
	\begin{equation}
		\label{eq:Peano-general}
		\|d^{\,i}a_0\|_H
		\le C\bigl\{n^{-i}M_{n+m}a+n^{\,m-i}M_n(d^ma)\bigr\},
	\end{equation}
	with $C$ independent of $n$ and $H$.
\end{lemma}

\begin{proof}
	Choose $\Gamma^{[k]}\in C_c^\infty((-1,1))$, $0\le k\le m-1$, satisfying
	\[
	\int_{-1}^1\Gamma^{[k]}(s)s^l\,\dd s=k!\,\delta_{lk},
	\qquad 0\le l\le m-1,
	\]
	which exist by the argument given after
	\eqref{eq:Hilbert-interpolation-mollifier}.  Put
	$\gamma_j^{[k]}:=n^{-1-k}\Gamma^{[k]}(j/n)$ and
	\[
	A_{lk}:=n^{\,k-l}\sum_{|j|\le n}\gamma^{[k]}_j\binom jl ,
	\qquad 0\le l,k\le m-1 .
	\]
	The leading coefficient of $\binom jl$ is $1/l!$, and the Riemann-sum
	error for each fixed smooth compactly supported integrand is $O(n^{-1})$. Thus
	\[
	A_{lk}=\frac{k!}{l!}\delta_{lk}+O(n^{-1})=\delta_{lk}+O(n^{-1}),
	\]
	so $A$ is invertible with $|A^{-1}|\le C$ once $n\ge n_0(m)$.  Set
	\begin{equation*}
		\gamma_j:=n^{-1-i}\sum_{k=0}^{m-1}(A^{-1}e_i)_k\,\Gamma^{[k]}(j/n),
		\qquad |j|\le n .
	\end{equation*}
	Then
	\[
	\sum_{|j|\le n}\gamma_j\binom jl=n^{\,l-i}(AA^{-1}e_i)_l=\delta_{li}\qquad\text{for } 0\le l\le m-1 ,
	\]
	while
	\begin{equation}
		\label{eq:Peano-weight-size}
		\max_{|j|\le n}|\gamma_j|\le Cn^{-1-i},
		\qquad
		\sum_{|j|\le n}|\gamma_j|\le Cn^{-i}.
	\end{equation}
	
	Newton's formula holds for both positive and negative lattice indices:
	\[
	a_j=\sum_{l=0}^{m-1}\binom jl\,d^la_0+R_j,
	\qquad
	R_j=\operatorname{sgn}(j)\sum_k\binom{j-k-1}{m-1}d^ma_k ,
	\]
	the sum being taken over $0\le k<j$ if $j>0$ and over $j\le k<0$ if $j<0$, so that
	$|k|\le n$ whenever $|j|\le n$; moreover
	$\bigl|\binom{j-k-1}{m-1}\bigr|\le Cn^{m-1}$ there.  Multiplying by $\gamma_j$ and summing over $|j|\le n$, the moment identities give
	\[
	d^{\,i}a_0=\sum_{|j|\le n}\gamma_ja_j-\sum_{|k|\le n}K_{n,k}\,d^ma_k,
	\qquad
	K_{n,k}:=\sum_j\operatorname{sgn}(j)\gamma_j\binom{j-k-1}{m-1}.
	\]
	Here the last sum ranges over $k<j\le n$ when $k\ge0$, and over $-n\le j\le k$ when $k<0$.
	By \eqref{eq:Peano-weight-size}, the first sum satisfies
	\[
	\Bigl\|\sum_{|j|\le n}\gamma_ja_j\Bigr\|_H
	\le\max_{|j|\le n}|\gamma_j|\sum_{|j|\le n}\|a_j\|_H
	\le C n^{-1-i}(2n+1)M_na\le Cn^{-i}M_{n+m}a.
	\]
	The bound on the binomial coefficients above and \eqref{eq:Peano-weight-size} give
	\[
	\max_k|K_{n,k}|\le Cn^{m-1}\sum_j|\gamma_j|\le Cn^{\,m-1-i},
	\]
	and hence the second sum satisfies
	\[
	\Bigl\|\sum_{|k|\le n}K_{n,k}\,d^ma_k\Bigr\|_H
	\le Cn^{\,m-1-i}(2n+1)M_n(d^ma)\le Cn^{\,m-i}M_n(d^ma).
	\]
	The required estimate follows.
\end{proof}

\begin{lemma}[Discrete interpolation along a flow]
	\label{lem:discrete-interpolation}
	Let $h_0>0$ satisfy $2(n_0(m)+m)h_0\le s_0$ for $2\le m\le\mathrm N$, with $n_0(m)$ of
	Lemma~\ref{lem:discrete-Peano} and $s_0$ of \eqref{eq:complete-flow-uniform-bounds}.
	Let $2\le m\le\mathrm N$ and $1\le i\le m-1$, let $X\in\cX_{\rm c}$ and
	$0<|h|<h_0$, and let $\mathsf e$ be as in \eqref{eq:GN-output-exponents}.
	\begin{enumerate}
		\item\label{it:discrete-interpolation-stationary}
		Let $\mathrm L\in\bR$ and $2\le p\le\infty$.  For a function $u$ set
		\begin{equation}
			\label{eq:discrete-lower-density}
			\cR_{m,h}^{X,\mathrm L}(u):=
			\sum_{j=0}^{m-1}\|w^{\mathrm L+m-j}(\Delta_h^X)^ju\|_{L_x^2H_\sigma^1}.
		\end{equation}
		Then, whenever
		\[
		\|w^{\mathrm L}u\|_{L_x^pH_\sigma^1}+\|w^{\mathrm L}(\Delta_h^X)^mu\|_{L_x^2H_\sigma^1}+\cR_{m,h}^{X,\mathrm L}(u)<\infty,
		\]
		we have
		\begin{equation}
			\label{eq:discrete-GN}
			\|w^{\mathrm L}(\Delta_h^X)^iu\|_{L_x^{\mathsf e(m/i,\,p)}H_\sigma^1}
			\le C\|w^{\mathrm L}u\|_{L_x^pH_\sigma^1}^{1-i/m}
			\left(\|w^{\mathrm L}(\Delta_h^X)^mu\|_{L_x^2H_\sigma^1}
			+\cR_{m,h}^{X,\mathrm L}(u)\right)^{i/m}
			+C\|w^{\mathrm L}u\|_{L_x^pH_\sigma^1}.
		\end{equation}
		The constant $C$ is independent of $h$ and $u$.
		\item\label{it:discrete-interpolation-evolutionary}
		Let $0<T<\infty$, let $2\le p_x\le\infty$, and let $u$ be a function on
		$(0,T)\times\cD\times\bR^3$.  With $\cR_{m,h}^{X,0}$ of
		\eqref{eq:discrete-lower-density}, set
		\begin{equation}
			\label{eq:dynamic-discrete-top-density}
			\mathfrak e_{m,h}^X[u](t)
			:=
			\|(\Delta_h^X)^m u(t)\|_{L_x^2H_\sigma^1}
			+\cR_{m,h}^{X,0}(u(t)).
		\end{equation}
		Then, for almost every $t$,
		\begin{equation}
			\label{eq:dynamic-discrete-pointwise}
			\|(\Delta_h^X)^iu(t)\|_{L_x^{\mathsf e(m/i,\,p_x)}H_\sigma^1}
			\le C\|u(t)\|_{L_x^{p_x}H_\sigma^1}^{1-i/m}
			\mathfrak e_{m,h}^X[u](t)^{i/m}
			+C\|u(t)\|_{L_x^{p_x}H_\sigma^1}.
		\end{equation}
		Let moreover $1\le p_t\le p_t^\flat\le\infty$ and
		\begin{equation}
			\label{eq:dynamic-discrete-time-exponent}
			r:=\mathsf e(m/i;p_t,p_t^\flat) ,
		\end{equation}
		so that $r\le p_t^\flat$.  If
		$\|u\|_{L_x^{p_x}H_\sigma^1}\in L^{p_t^\flat}(0,T)$ and
		$\mathfrak e_{m,h}^X[u]\in L^{p_t}(0,T)$, then
		\begin{equation}
			\label{eq:dynamic-discrete-time}
			\|(\Delta_h^X)^iu\|_{L_t^{r}L_x^{\mathsf e(m/i,\,p_x)}H_\sigma^1}
			\le C\|u\|_{L_t^{p_t^\flat}L_x^{p_x}H_\sigma^1}^{1-i/m}
			\|\mathfrak e_{m,h}^X[u]\|_{L_t^{p_t}}^{i/m}
			+CT^{1/r-1/p_t^\flat}
			\|u\|_{L_t^{p_t^\flat}L_x^{p_x}H_\sigma^1}.
		\end{equation}
		The constant is uniform in the time interval and the step size.
	\end{enumerate}
\end{lemma}

\begin{proof}
	\emph{(\ref{it:discrete-interpolation-stationary}).}
	Abbreviate 
	$$
	u_{\mathrm L}:=w^{\mathrm L}u,\quad A:=\|w^{\mathrm L}u\|_{L_x^pH_\sigma^1},\quad \cR:=\cR_{m,h}^{X,\mathrm L}(u),\quad B:=\|w^{\mathrm L}(\Delta_h^X)^mu\|_{L_x^2H_\sigma^1}+\cR.
	$$
	We argue by induction on $i$. The case $i=0$ follows from
	$\mathsf e(\infty,p)=p$. 
	
	We first control the frame under the lifted flow. For $|s|\le2s_0$,
	let $Q_s(x,v)$ be the $4\times4$ matrix determined by
	\begin{equation}
		\label{eq:discrete-frame-transport}
		T_s^X\mathbb A(w^{\mathrm L}q)=Q_s\,\mathbb A(w^{\mathrm L}T_s^Xq)
	\end{equation}
	for every $q$.  Write $v':=A_s(x)v$ and $\omega:=w^{\mathrm L}(v')/w^{\mathrm L}(v)$.  By \eqref{eq:frame-def} and $T_s^Xq=q\circ\widehat\Phi_s$, the matrix $Q_s$ is block upper
	triangular, with diagonal blocks $\omega B(v')A_s^{-T}B(v)^{-1}$ and
	$(w(v')/w(v))^{\mathrm L-1/2}$, and with off-diagonal block
	\[
	w(v)^{1/2-\mathrm L}B(v')\bigl((\nabla w^{\mathrm L})(v')-\omega A_s^{-T}(\nabla w^{\mathrm L})(v)\bigr),
	\]
	where $B(\cdot)$ is the matrix \eqref{eq:coulomb-square-root-symbol}.
	The group law gives $Q_{s_2}^{-1}\circ\widehat\Phi_{s_1}=Q_{s_1}Q_{s_1+s_2}^{-1}$,
	in which $B(v)$ cancels.  Let $F$ be a function on $\bR^3$ with
	$|\partial_v^\beta F|\le C\la v\ra^{a-|\beta|}$. 
	Fix $x\in\overline\cD$ and $v\in\bR^3$.  Since $\partial_sA_s=(\nabla X)(\Phi_s)A_s$, the chain rule gives
	\begin{align}\label{260915747}
		\partial_s\bigl[F(A_sv)\bigr]
		=(\nabla F)(A_sv)\,(\nabla X)(\Phi_s)A_sv
		=\sum_{j,l=1}^3(\partial_jX_l)(\Phi_s)\,(A_sv)_j\,(\partial_lF)(A_sv) .
	\end{align}
	Since moreover $\partial_s\Phi_s=X(\Phi_s)$, for $1\le l\le3$ and every multi-index $\alpha$ with $|\alpha|\le\mathrm N-1$,
	\begin{align}\label{260915748}
		\partial_s(A_sv)_l=\sum_{j=1}^3(\partial_jX_l)(\Phi_s)\,(A_sv)_j ,
		\qquad
		\partial_s\bigl[(\partial^\alpha X_l)(\Phi_s)\bigr]=\sum_{j=1}^3X_j(\Phi_s)\,(\partial_j\partial^\alpha X_l)(\Phi_s) .
	\end{align}
	Combining \eqref{260915747} and \eqref{260915748}, we obtain by induction on $k$ that, for $0\le k\le\mathrm N$, $\partial_s^k[F(A_sv)]$ is a finite sum, whose number of terms depends only on $k$, of products of $k$ factors $(\partial^\alpha X_l)(\Phi_s)$ with $|\alpha|\le k$, of $|\beta|$ components of $A_sv$ and of $(\partial_v^\beta F)(A_sv)$, where $|\beta|\le k$.  Since $X\in C^{K_\cD+1}$ has compact support and $\mathrm N\le K_\cD$, we have $|(\partial^\alpha X_l)(\Phi_s)|\le C$ for $|\alpha|\le\mathrm N$.  By \eqref{eq:flow-weight-comparison} and \eqref{eq:diffuse-part-weight}, $C^{-1}\la v\ra\le\la A_sv\ra\le C\la v\ra$ for $|s|\le2s_0$, and the hypothesis on $F$ gives
	\[
	|A_sv|^{|\beta|}\,\bigl|(\partial_v^\beta F)(A_sv)\bigr|
	\le C\la A_sv\ra^{a}\le C\la v\ra^{a} .
	\]
	Therefore
	\[
	\bigl|\partial_s^k[F(A_sv)]\bigr|\le C\la v\ra^{a},
	\qquad 0\le k\le\mathrm N,\quad |s|\le2s_0 ,
	\]
	with $C$ independent of $x\in\overline\cD$ and $v\in\bR^3$.  The bounds for $\partial_v^\beta B^{\pm1}$ from Lemma~\ref{lem:first-frame-commutator},
	together with \eqref{eq:complete-flow-uniform-bounds}, therefore yield
	\[
	|\partial_{s_1}^{k_1}\partial_{s_2}^{k_2}(Q_{s_1}Q_{s_1+s_2}^{-1})|_{\rm op}\le C\la v\ra,
	\qquad k_1+k_2\le\mathrm N,\quad |s_1|,|s_2|\le s_0 .
	\]
	Put
	\[
	E_h:=h^{-1}(I-Q_h^{-1})=-h^{-1}\int_0^h\partial_sQ_s^{-1}\,\dd s .
	\]
	By the group law,
	\[
	E_h\circ\widehat\Phi_s=-h^{-1}\int_0^h\partial_{s'}(Q_sQ_{s+s'}^{-1})\,\dd s'\qquad\text{with } |s'|<h_0\le s_0 ,
	\]
	and hence $|\partial_s^k(E_h\circ\widehat\Phi_s)|_{\rm op}\le C\la v\ra$ for
	$0\le k\le\mathrm N-1$ and $|s|\le s_0$, uniformly in $x$ and in $0<|h|<h_0$.
	Consequently, \eqref{eq:discrete-frame-transport} at $s=h$ gives
	\begin{equation}
		\label{eq:discrete-frame-product-rule}
		\Delta_h^X\mathbb A(w^{\mathrm L}q)
		=\mathbb A(w^{\mathrm L}\Delta_h^Xq)+E_hT_h^X\mathbb A(w^{\mathrm L}q) .
	\end{equation}
	
	For fixed $x$, set $V_x(jh):=T_{jh}^X(\mathbb Au_{\mathrm L})(x)$ and apply
	Lemma~\ref{lem:discrete-Peano} with $H=L_v^2(\bR^4)$, at the
	order $i$, to $a_j:=V_x(jh)$.  Since the flow is autonomous, $T_{jh}^X=(T_h^X)^j$, so that
	$$
	d^{\,i}a_0=h^{\,i}(\Delta_h^X)^iV_x(0)\quad \text{and}\quad d^ma_j=h^mT_{jh}^X(\Delta_h^X)^m(\mathbb Au_{\mathrm L})(x).
	$$
	Dividing \eqref{eq:Peano-general} by $|h|^i$, we obtain
	\begin{equation}
		\label{eq:discrete-peano-window}
		\|(\Delta_h^X)^iV_x(0)\|_{L_v^2}
		\le C\{(n|h|)^{-i}\cM_{X,h}^{\rm d}\varphi
		+(n|h|)^{\,m-i}\cM_{X,h}^{\rm d}\psi_{m,h}\},
	\end{equation}
	where
	\[
	\varphi(x):=\|u_{\mathrm L}(x)\|_{H_\sigma^1}\quad\text{and}\quad \psi_{m,h}(x):=\|(\Delta_h^X)^m(\mathbb Au_{\mathrm L})(x)\|_{L_v^2} .
	\]
	Indeed, by \eqref{eq:complete-flow-uniform-bounds} and \eqref{eq:frame-exact-norm}, 
	\[
	M_{n+m}a\le C\cM_{X,h}^{\rm d}\varphi(x)\quad\text{and}\quad |h|^{-m}M_n(d^ma)\le C\cM_{X,h}^{\rm d}\psi_{m,h}(x)\qquad\text{whenever } n+m\le s_0/|h| .
	\]
	
	By the choice of $h_0$, we have $n_0|h|<s_0/2$, $4|h|<s_0$ and $m|h|<s_0/2$. Thus
	$(n+m)|h|<s_0$ whenever $n|h|\le s_0/2$.  Hence \eqref{eq:discrete-peano-window} holds at
	every integer $n$ with $n_0\le n\le s_0/(2|h|)$.
	
	At the fixed point $x$, write 
	$$
	\alpha=\cM_{X,h}^{\rm d}\varphi(x),\quad \beta=\cM_{X,h}^{\rm d}\psi_{m,h}(x),\quad r_*:=(\alpha/\beta)^{1/m},
	$$
	with
	$r_*=\infty$ if $\beta=0$.  If $r_*>s_0/2$, take $n=\lfloor s_0/(2|h|)\rfloor$;
	then $n\ge n_0$, $n|h|\ge s_0/4$ and $\beta\le\alpha(2/s_0)^m$, so \eqref{eq:discrete-peano-window} gives
	\[
	\|(\Delta_h^X)^iV_x(0)\|_{L_v^2}\le C\bigl((n|h|)^{-i}\alpha+(n|h|)^{\,m-i}\beta\bigr)\le C\alpha.
	\]
	If
	$n_0|h|\le r_*\le s_0/2$, take $n=\lfloor r_*/|h|\rfloor$; then $n\ge n_0\ge1$
	gives 
	$r_*/2<n|h|\le r_*$.  Since $m-i\ge1$ and
	$r_*^{-i}\alpha=r_*^{m-i}\beta=\alpha^{1-i/m}\beta^{i/m}$, \eqref{eq:discrete-peano-window} gives
	\[
	\|(\Delta_h^X)^iV_x(0)\|_{L_v^2}
	\le C\bigl((n|h|)^{-i}\alpha+(n|h|)^{\,m-i}\beta\bigr)
	\le C\bigl(2^ir_*^{-i}\alpha+r_*^{m-i}\beta\bigr)
	=C(2^i+1)\alpha^{1-i/m}\beta^{i/m}.
	\]
	If
	$r_*<n_0|h|$, then
	\[
	\|(\Delta_h^X)^iV_x(0)\|_{L_v^2}
	\le|h|^{-i}\sum_{k=0}^{i}\binom ik\|V_x(kh)\|_{L_v^2}
	\le C|h|^{-i}\sum_{k=0}^{i}\varphi(\Phi_{kh}^Xx)
	\le C|h|^{-i}\alpha .
	\]
	Here $\|V_x(kh)\|_{L_v^2}\le C\varphi(\Phi_{kh}^Xx)$ follows from \eqref{eq:complete-flow-uniform-bounds} and \eqref{eq:frame-exact-norm}; the last inequality uses the average at $n=i$ in $\cM_{X,h}^{\rm d}$, which is admissible since $i<m\le s_0/|h|$. Since $r_*<n_0|h|$, we have $\alpha<(n_0|h|)^m\beta$, which gives
	$|h|^{-i}\alpha\le n_0^{\,i}\alpha^{1-i/m}\beta^{i/m}$.  Consequently
	\begin{equation}
		\label{eq:discrete-peano-optimized}
		\|(\Delta_h^X)^iV_x(0)\|_{L_v^2}
		\le C(\cM_{X,h}^{\rm d}\varphi)^{1-i/m}
		(\cM_{X,h}^{\rm d}\psi_{m,h})^{i/m}
		+C\cM_{X,h}^{\rm d}\varphi .
	\end{equation}
	
	The order-$m$ expansion of \eqref{eq:discrete-frame-product-rule} has leading term $\mathbb A(w^{\mathrm L}(\Delta_h^X)^mu)$. Its other terms contain flow translates of $\mathbb A(w^{\mathrm L}(\Delta_h^X)^ju)$, $j<m$, with coefficients bounded by $C\langle v\rangle^{m-j}$. This follows from the discrete product rule and $\Delta_h^XT_h^X=T_h^X\Delta_h^X$, using difference quotients of $E_h\circ\widehat\Phi_s$ only up to order $m-1\le\mathrm N-1$. The weighted frame inequality, \eqref{eq:flow-weight-comparison}, and the flow Jacobian bounds then give $\|\psi_{m,h}\|_{L_x^2}\le CB$, uniformly in $h$.

	Taking the $L_x^{\mathsf e(m/i,\,p)}$ norm in \eqref{eq:discrete-peano-optimized}, Lemma~\ref{lem:flow-maximal}\,(\ref{it:flow-maximal-discrete}) at the exponents $p$ and $2$, followed
	by H\"older's inequality in $x$ with \eqref{eq:GN-output-exponents} at the ratio $i/m$
	and $a=2$, gives
	\begin{equation}
		\label{eq:discrete-frame-difference}
		\bigl\|(\Delta_h^X)^i(\mathbb Au_{\mathrm L})\bigr\|
		_{L_x^{\mathsf e(m/i,\,p)}L_v^2}
		\le CA^{1-i/m}B^{i/m}+CA ,
	\end{equation}
	since $\|\varphi\|_{L_x^p}=A$ and $\|\psi_{m,h}\|_{L_x^2}\le CB$, as proved above.
	
	It remains to commute the differences through the weighted frame. By
	\eqref{eq:frame-exact-norm},
	\[
	\|\mathbb A(w^{\mathrm L}(\Delta_h^X)^iu)\|_{L_v^2}\simeq\|w^{\mathrm L}(\Delta_h^X)^iu\|_{H_\sigma^1}.
	\]
	Iterating \eqref{eq:discrete-frame-product-rule} using the product rule
	\[
	\Delta_h^X(GF)=(G\circ\widehat\Phi_h)\Delta_h^XF+h^{-1}(G\circ\widehat\Phi_h-G)F\quad\text{for a matrix } G
	\]
	and the commutation $\Delta_h^XT_h^X=T_h^X\Delta_h^X$,
	gives
	\[
	(\Delta_h^X)^i(\mathbb Au_{\mathrm L})
	=\mathbb A\bigl(w^{\mathrm L}(\Delta_h^X)^iu\bigr)
	+\sum_{j=0}^{i-1}\sum_{l=1}^{i}G_{j,l}\,
	T_{lh}^X\mathbb A\bigl(w^{\mathrm L}(\Delta_h^X)^ju\bigr),
	\]
	where each $G_{j,l}$ is a sum of products of at most $i-j$ matrices, each a difference
	quotient of $s\mapsto E_h\circ\widehat\Phi_s$ of order at most $i-1$ with nodes in $[0,ih]$; by the bounds on $E_h$ above,
	$|G_{j,l}|_{\rm op}\le C\la v\ra^{i-j}$.  By \eqref{eq:complete-flow-uniform-bounds} and
	\eqref{eq:flow-weight-comparison},
	\[
	\|w^aT_{lh}^XF(x)\|_{L_v^2}\le C\|w^aF(\Phi_{lh}x)\|_{L_v^2} .
	\] Moreover, \eqref{eq:weighted-frame-monotonicity} with both weight exponents equal to $a$,
	applied to $w^{\mathrm L}q$ in place of $q$, gives $\|w^a\mathbb A(w^{\mathrm L}q)\|_{L_v^2}\le C\|w^{\mathrm L+a}q\|_{H_\sigma^1}$.  Hence, by these two facts at $a=i-j$, pointwise in $x$,
	\[
	\Bigl\|\sum_{l=1}^{i}G_{j,l}\,T_{lh}^X\mathbb A\bigl(w^{\mathrm L}(\Delta_h^X)^ju\bigr)(x)\Bigr\|_{L_v^2}
	\le C\sum_{l=1}^{i}T_{lh}^X\bigl(\|w^{\mathrm L+i-j}(\Delta_h^X)^ju\|_{H_\sigma^1}\bigr)(x) .
	\]
	For $0\le j<i$, set $\tau:=(i-j)/(m-j)\in(0,1]$.  By \eqref{eq:output-weight-p} with the top exponent $2$, $k=i-j$ and $q=(\Delta_h^X)^ju$,
	\[
	\|w^{\mathrm L+i-j}(\Delta_h^X)^ju\|_{L_x^{\mathsf e(m/i,\,p)}H_\sigma^1}
	\le C\|w^{\mathrm L}(\Delta_h^X)^ju\|_{L_x^{\mathsf e(m/j,\,p)}H_\sigma^1}^{1-\tau}
	\|w^{\mathrm L+m-j}(\Delta_h^X)^ju\|_{L_x^2H_\sigma^1}^{\tau}.
	\]
	Here
	$\|w^{\mathrm L+m-j}(\Delta_h^X)^ju\|_{L_x^2H_\sigma^1}\le\cR\le B$, since it is the $j$th term of $\cR$, and
	\[
	\|w^{\mathrm L}(\Delta_h^X)^ju\|_{L_x^{\mathsf e(m/j,\,p)}H_\sigma^1}\le C(A^{1-j/m}B^{j/m}+A)
	\]
	by the induction hypothesis at the order $j<i$.  The interpolation inequality \eqref{eq:rung-composition-bound} now gives
	\[
	\|w^{\mathrm L+i-j}(\Delta_h^X)^ju\|_{L_x^{\mathsf e(m/i,\,p)}H_\sigma^1}
	\le C(A^{1-i/m}B^{i/m}+A).
	\]
	Finally $T_{lh}^X$, $l|h|<s_0$, acts here on a scalar
	function of $x$, so the Jacobian bounds in \eqref{eq:complete-flow-uniform-bounds} give
	\[
	\bigl\|T_{lh}^X\bigl(\|w^{\mathrm L+i-j}(\Delta_h^X)^ju\|_{H_\sigma^1}\bigr)\bigr\|_{L_x^{\mathsf e(m/i,\,p)}}
	\le C\|w^{\mathrm L+i-j}(\Delta_h^X)^ju\|_{L_x^{\mathsf e(m/i,\,p)}H_\sigma^1}.
	\]
	Combining these bounds with \eqref{eq:discrete-frame-difference} yields \eqref{eq:discrete-GN}.
	
	These bounds hold through order $m$ because only differences of $E_h\circ\widehat\Phi_s$ of order at most $m-1$ occur.
	
	\emph{(\ref{it:discrete-interpolation-evolutionary}).}
	Applying (\ref{it:discrete-interpolation-stationary}) with $\mathrm L=0$ and $p=p_x$
	to $u(t)$ for almost every $t$ such that
	$\|u(t)\|_{L_x^{p_x}H_\sigma^1}+\mathfrak e_{m,h}^X[u](t)<\infty$ gives
	\eqref{eq:dynamic-discrete-pointwise}, since the bracket in \eqref{eq:discrete-GN} is
	then $\mathfrak e_{m,h}^X[u](t)$.
	
	To obtain \eqref{eq:dynamic-discrete-time}, take the $L_t^{r}$ norm in
	\eqref{eq:dynamic-discrete-pointwise}.  By
	\eqref{eq:dynamic-discrete-time-exponent} and \eqref{eq:GN-output-exponents},
	$\frac1r=\frac{1-i/m}{p_t^\flat}+\frac{i/m}{p_t}$, so H\"older's inequality in $t$ with the
	exponents $\bigl(\tfrac{p_t^\flat}{1-i/m},\tfrac{p_t}{i/m}\bigr)$ gives
	\[
	\bigl\|\|u\|_{L_x^{p_x}H_\sigma^1}^{1-i/m}\,\mathfrak e_{m,h}^X[u]^{i/m}\bigr\|_{L^r(0,T)}
	\le\|u\|_{L_t^{p_t^\flat}L_x^{p_x}H_\sigma^1}^{1-i/m}\|\mathfrak e_{m,h}^X[u]\|_{L_t^{p_t}}^{i/m},
	\]
	and, since $r\le p_t^\flat$, H\"older's inequality on $(0,T)$ gives
	\[
	\bigl\|\|u\|_{L_x^{p_x}H_\sigma^1}\bigr\|_{L^r(0,T)}
	\le T^{1/r-1/p_t^\flat}\|u\|_{L_t^{p_t^\flat}L_x^{p_x}H_\sigma^1}.
	\]
\end{proof}

\subsubsection{Bounds along ordered words}\label{ss:flow-ordered-words}

The preceding inequalities apply to ordered compositions of complete derivatives. The order of the letters is retained because these derivatives need not commute.

\begin{corollary}[Gagliardo--Nirenberg bounds along ordered words]
	\label{cor:ordered-word-gaussian-rungs}
	Let $\mathrm L_{\rm G}\in\bR$, $2\le a_{\rm G}\le p_{\rm G}\le\infty$, $2\le p_{\rm c}\le\infty$ and
	$2\le M\le\mathrm N$, let $\mathsf e$ be as in
	\eqref{eq:GN-output-exponents}, and let $\cG^{0}_{m,2}$ and
	$\cG^{\sharp,\mathrm L_{\rm G}}_{m,a_{\rm G}}$ be \eqref{eq:GN-lower-order} at
	$(\mathrm L,a)=(0,2)$ and \eqref{eq:GN-lower-order-gaussian} at
	$(\mathrm L,a)=(\mathrm L_{\rm G},a_{\rm G})$.
	\begin{enumerate}
		\item\label{it:rungs-stationary}
		Let $f$ be a function on $\cD\times\bR^3$ such that, for all ordered complete
		words $J$ and $\nu J$,
		\begin{equation}
			\mathfrak D_\nu(\mathfrak D_Jf)=\mathfrak D_{\nu J}f
			\quad\text{in }\cD'(\cD\times\bR^3),
			\qquad |\nu J|\le\mathrm N,
			\label{eq:pointwise-slice-chain-hypothesis}
		\end{equation}
		and
		\begin{equation}
			\|w^{\mathrm L_{\rm G}}f\|_{L_x^{p_{\rm G}}L_v^2}
			+\|f\|_{L_x^{p_{\rm c}}H_\sigma^1}
			+\cG^{0}_{\mathrm N,2}(f)+\cG^{\sharp,\mathrm L_{\rm G}}_{\mathrm N,a_{\rm G}}(f)<\infty .
			\label{eq:dynamic-pointwise-finite-norms}
		\end{equation}
		Then every ordered complete word $I$ with
		$1\le|I|=k\le M-1$ satisfies
		\begin{equation}
			\label{eq:ordered-word-gaussian-rungs}
			\|w^{\mathrm L_{\rm G}}\mathfrak D_If\|
			_{L_x^{\mathsf e(M/k;\,a_{\rm G},p_{\rm G})}L_v^2}
			\le C\|w^{\mathrm L_{\rm G}}f\|_{L_x^{p_{\rm G}}L_v^2}^{1-k/M}
			\cG^{\sharp,\mathrm L_{\rm G}}_{M,a_{\rm G}}(f)^{k/M}
			+C\|w^{\mathrm L_{\rm G}}f\|_{L_x^{p_{\rm G}}L_v^2},
		\end{equation}
		and
		\[
		\|w^{\mathrm L_{\rm G}}\mathfrak D_If\|_{L_x^{\mathsf e(M/k;\,a_{\rm G},p_{\rm G})}L_v^2}+\|w^{\mathrm L_{\rm G}}f\|_{L_x^{p_{\rm G}}L_v^2}+\cG^{\sharp,\mathrm L_{\rm G}}_{M,a_{\rm G}}(f)<\infty.
		\]
		
		Let $\mathrm L\ge0$, and let $f$ be a function on $\cD\times\bR^3$ that
		satisfies \eqref{eq:pointwise-slice-chain-hypothesis} and
		$$
		\|w^{\mathrm L}f\|_{L_x^{p_{\rm c}}H_\sigma^1}+\cG^{\mathrm L}_{M,2}(f)<\infty,
		$$
		with $\cG^{\mathrm L}_{M,2}$ of \eqref{eq:GN-lower-order}.  Then every ordered
		complete word $I$ with $1\le|I|=k\le M-1$ satisfies
		\begin{equation}
			\label{eq:ordered-word-coercive-rungs-weighted}
			\|w^{\mathrm L}\mathfrak D_If\|_{L_x^{\mathsf e(M/k,\,p_{\rm c})}H_\sigma^1}
			\le C\|w^{\mathrm L}f\|_{L_x^{p_{\rm c}}H_\sigma^1}^{1-k/M}
			\cG^{\mathrm L}_{M,2}(f)^{k/M}
			+C\|w^{\mathrm L}f\|_{L_x^{p_{\rm c}}H_\sigma^1}.
		\end{equation}
		\item\label{it:rungs-evolutionary}
		Let $t_0<t_1$ be real, and let $u$ be a
		function on $(t_0,t_1)\times\cD\times\bR^3$ satisfying
		\eqref{eq:temporal-slicing-integrability}--\eqref{eq:temporal-slicing-spacetime-chain}
		with $z=u$, $(t_0,t_1)$ in place of $(0,T)$ and $\cI$ the set of ordered
		complete words of length at most $\mathrm N$, and satisfying
		\eqref{eq:dynamic-pointwise-finite-norms} at $f=u(t)$ for almost every
		$t\in(t_0,t_1)$.  Let $2\le r_0\le\infty$ be the
		temporal base exponent and $2\le r_1\le r_0$ the temporal top exponent.  With the norms in
		$t$ taken over $(t_0,t_1)$, write
		\begin{equation}
			\label{eq:slab-rung-quantities}
			\begin{gathered}
				\mathfrak b^\sharp_\rho:=\|w^{\mathrm L_{\rm G}}u\|_{L_t^{\rho}L_x^{p_{\rm G}}L_v^2},
				\qquad
				\mathfrak G^\sharp_{r_1}
				:=\bigl\|\cG^{\sharp,\mathrm L_{\rm G}}_{M,a_{\rm G}}(u)\bigr\|
				_{L_t^{r_1}}.
			\end{gathered}
		\end{equation}
		Then every ordered complete word $J$ with $1\le|J|=k\le M-1$ satisfies
		\begin{equation}
			\|w^{\mathrm L_{\rm G}}\mathfrak D_Ju\|
			_{L_t^{\mathsf e(M/k;\,r_1,r_0)}
				L_x^{\mathsf e(M/k;\,a_{\rm G},p_{\rm G})}L_v^2}
			\le C\,\bigl(\mathfrak b^\sharp_{r_0}\bigr)^{1-k/M}
			\bigl(\mathfrak G^\sharp_{r_1}
			+\mathfrak b^\sharp_{r_1}\bigr)^{k/M},
			\label{eq:slab-word-gaussian-rungs}
		\end{equation}
		where $C$ is independent of $u$ and $t_1-t_0$.
	\end{enumerate}
	Both (\ref{it:rungs-stationary}) and (\ref{it:rungs-evolutionary}) hold with $\mathrm N$ replaced by any integer $M'$ with
	$M\le M'\le\mathrm N$.
\end{corollary}

\begin{proof}
	\emph{(\ref{it:rungs-stationary}).}
	The two estimates follow from the same induction. We take
	$(a,p,\lambda)=(a_{\rm G},p_{\rm G},\mathrm L_{\rm G})$ and
	Theorem~\ref{thm:complete-flow-GN}\,(\ref{it:GN-fiber}) for the former and
	$(a,p,\lambda)=(2,p_{\rm c},\mathrm L)$ and
	Theorem~\ref{thm:complete-flow-GN}\,(\ref{it:GN-coercive}) for the latter.
	Put
	\begin{align*}
		(b,G)\,&:=\bigl(\|w^{\mathrm L_{\rm G}}f\|_{L_x^{p_{\rm G}}L_v^2},\,\cG^{\sharp,\mathrm L_{\rm G}}_{M,a_{\rm G}}(f)\bigr)
		\quad\text{for \eqref{eq:ordered-word-gaussian-rungs}},\\
		(b,G)\,&:=\bigl(\|w^{\mathrm L}f\|_{L_x^{p_{\rm c}}H_\sigma^1},\,\cG^{\mathrm L}_{M,2}(f)\bigr)
		\quad\text{for \eqref{eq:ordered-word-coercive-rungs-weighted}},
	\end{align*}
	and write $P_k$ for the supremum of $\|w^{\lambda}\mathfrak D_If\|_{L_x^{\mathsf e(M/k;a,p)}L_v^2}$, respectively $\|w^{\lambda}\mathfrak D_If\|_{L_x^{\mathsf e(M/k;a,p)}H_\sigma^1}$, over the
	ordered complete words $I$ of length $k$; thus $P_0=b$. The estimates reduce to
	$$
	P_k\le C(b^{1-k/M}G^{k/M}+b)\qquad \text{for }0\le k\le M-1.
	$$  In the Gaussian case $b,G<\infty$ by \eqref{eq:dynamic-pointwise-finite-norms}; in the coercive case both are finite by hypothesis.
	
	Fix $1\le k\le M-1$ and write $I=\nu J$, $|J|=k-1$.  Apply Theorem~\ref{thm:complete-flow-GN} to
	$\mathfrak D_Jf$ at $X=\tau_\nu$, $(m,i)=(M-k+1,1)$, the weight $\lambda$ and the
	base exponent $\mathsf e(M/(k-1);a,p)$ of $P_{k-1}$.  The required inequalities are
	$$
	m\ge2,\quad m\le M\le\mathrm N,\quad a\le\mathsf e(M/(k-1);a,p)\le p,
	$$
	and $\|w^{\lambda}\mathfrak D_Jf\|_{L_x^{\mathsf e(M/(k-1);a,p)}L_v^2}$, respectively $\|w^{\lambda}\mathfrak D_Jf\|_{L_x^{\mathsf e(M/(k-1);a,p)}H_\sigma^1}$, is at most
	$P_{k-1}<\infty$ by the induction, and $\cG^{\sharp,\lambda}_{m,a}(\mathfrak D_Jf)$, respectively $\cG^{\lambda}_{m,a}(\mathfrak D_Jf)$, is at most $G<\infty$. Indeed, \eqref{eq:pointwise-slice-chain-hypothesis} gives
	$\mathfrak D_K\mathfrak D_Jf=\mathfrak D_{KJ}f$ with $|KJ|\le M$, and the term
	of $\mathfrak D_K\mathfrak D_Jf$ in this sum carries the weight
	\[
	\lambda+m-|K|=\lambda+M-|KJ| ,
	\]
	which is the weight of $\mathfrak D_{KJ}f$ in $G$.  By
	\eqref{eq:GN-rung-exponent-composition} with $(m,j,k)$ replaced by
	$(M,k-1,1)$, the output exponent is $\mathsf e(M/k;a,p)$.  Hence
	\[
	P_k\le CP_{k-1}^{1-1/m}G^{1/m}+CP_{k-1},
	\qquad m=M-k+1 .
	\]
	By induction,
	\[
	P_{k-1}\le C(b^{1-(k-1)/M}G^{(k-1)/M}+b) .
	\]
	By \eqref{eq:rung-composition-bound} with $(i,j,m,A,B)=(k,k-1,M,b,G)$ and $\tau=1/(M-k+1)$,
	\[
	CP_{k-1}^{1-1/m}G^{1/m}
	\le C\bigl(b^{1-(k-1)/M}G^{(k-1)/M}+b\bigr)^{1-1/m}G^{1/m}
	\le C\bigl(b^{1-k/M}G^{k/M}+b\bigr),
	\]
	and \eqref{eq:rung-exponent-monotonicity} at $(\alpha,\beta)=((k-1)/M,k/M)$ gives
	$CP_{k-1}\le C(b^{1-k/M}G^{k/M}+b)$.
	
	\emph{(\ref{it:rungs-evolutionary}).}
	By Lemma~\ref{lem:temporal-slicing} on $(t_0,t_1)$, the distributional chain \eqref{eq:pointwise-slice-chain-hypothesis} holds for $u(t)$ at almost every time. Hence the hypotheses of
	(\ref{it:rungs-stationary}) hold at $u(t)$ for almost every $t\in(t_0,t_1)$, at
	the same $\mathrm L_{\rm G}$, $a_{\rm G}$, $M$ and $(p_{\rm G},p_{\rm c})$.  For almost every $t$, \eqref{eq:ordered-word-gaussian-rungs} gives
	\[
	P(t)\le C\,\mathfrak a(t)^{1-k/M}\mathfrak g(t)^{k/M}+C\,\mathfrak a(t),
	\]
	with
	\[
	\begin{gathered}
		P(t):=\|w^{\mathrm L_{\rm G}}\mathfrak D_Ju(t)\|_{L_x^{\mathsf e(M/k;\,a_{\rm G},p_{\rm G})}L_v^2},\\
		\mathfrak a(t):=\|w^{\mathrm L_{\rm G}}u(t)\|_{L_x^{p_{\rm G}}L_v^2},
		\qquad
		\mathfrak g(t):=\cG^{\sharp,\mathrm L_{\rm G}}_{M,a_{\rm G}}(u(t)),
	\end{gathered}
	\]
	all measurable in $t$.  Write
	\begin{equation}
		\label{eq:slab-rung-time-reciprocal}
		\rho:=\mathsf e(M/k;r_1,r_0),
		\qquad
		\frac1\rho=\frac kM\cdot\frac1{r_1}
		+\Bigl(1-\frac kM\Bigr)\frac1{r_0},
	\end{equation}
	so that $r_1\le\rho\le r_0$.  H\"older's inequality after taking the $L_t^\rho$ norm in $t$ with the exponents
	$\bigl(\tfrac{r_0}{1-k/M},\tfrac{r_1}{k/M}\bigr)$, admissible by
	\eqref{eq:slab-rung-time-reciprocal}, bounds the first term. Moreover, H\"older's inequality with the same exponents, applied to
	$\mathfrak a=\mathfrak a^{1-k/M}\mathfrak a^{k/M}$, bounds the
	second:
	\[
	\bigl\|\mathfrak a^{1-k/M}\mathfrak g^{k/M}\bigr\|_{L_t^{\rho}}
	\le\|\mathfrak a\|_{L_t^{r_0}}^{1-k/M}
	\|\mathfrak g\|_{L_t^{r_1}}^{k/M},
	\qquad
	\|\mathfrak a\|_{L_t^{\rho}}
	\le\|\mathfrak a\|_{L_t^{r_0}}^{1-k/M}
	\|\mathfrak a\|_{L_t^{r_1}}^{k/M}.
	\]
	These bounds yield \eqref{eq:slab-word-gaussian-rungs}.
	
	The argument also applies with $M'$, $M\le M'\le\mathrm N$, in place of $\mathrm N$: it uses \eqref{eq:pointwise-slice-chain-hypothesis} only at
	$|\nu J|\le M$, the top-order norms in \eqref{eq:dynamic-pointwise-finite-norms} only through
	$\cG^{\sharp,\mathrm L_{\rm G}}_{M,a_{\rm G}}(f)\le\cG^{\sharp,\mathrm L_{\rm G}}_{M',a_{\rm G}}(f)$,
	and Lemma~\ref{lem:temporal-slicing} on the ordered complete words of length at
	most $M'$.
\end{proof}

\begin{corollary}[Mixed-word Gagliardo--Nirenberg estimates]
	\label{cor:ordered-word-rungs-mixed}
	Let $\mathrm L_{\rm G}\in\bR$, let $2\le a_{\rm G}\le p_{\rm G}\le\infty$ and $2\le p_{\rm c}\le\infty$, where $p_{\rm G}$ and $p_{\rm c}$ are spatial base exponents. Let
	$t_0<t_1$ be real with $t_1-t_0\ge1$, let $M'$ be an integer with
	$2\le M'\le\mathrm N$, and let $u$ be a function on
	$(t_0,t_1)\times\cD\times\bR^3$ such that, for the ordered words $J$ and
	$\nu J$ of the mixed alphabet \eqref{eq:dyn-mixed-alphabet},
	\begin{equation}
		\label{eq:mixed-rung-chain}
		\mathfrak D_\nu(\mathfrak D_Ju)=\mathfrak D_{\nu J}u
		\quad\text{in }\cD'\bigl((t_0,t_1)\times\cD\times\bR^3\bigr),
		\qquad|\nu J|\le M',
	\end{equation}
	and
	\begin{equation}
		\label{eq:mixed-rung-local-integrability}
		\mathfrak D_Ju\in L^1_{\rm loc}\bigl((t_0,t_1);L^2_{\rm loc}(\cD\times\bR^3)\bigr),
		\qquad|J|\le M'.
	\end{equation}
	Let $2\le M\le M'$, let $2\le r_0\le\infty$ be the temporal base exponent,
	let $0<\theta_\circ\le1/(2M)$,
	let $\mathrm L\ge0$, let $\mathfrak b^\sharp_\rho$ be as in
	\eqref{eq:slab-rung-quantities}, and set $\mathfrak b^0_\rho:=\|w^{\mathrm L}u\|_{L_t^{\rho}L_x^{p_{\rm c}}H_\sigma^1}$.  With $\cG^{\mathrm L}_{M,2}$ and
	$\cG^{\sharp,\mathrm L_{\rm G}}_{M,a_{\rm G}}$ given by
	\eqref{eq:GN-lower-order}--\eqref{eq:GN-lower-order-gaussian} with the sums
	taken over the ordered words of $\mathfrak A_\bullet$ in place of the ordered
	complete words, write
	\begin{equation}
		\label{eq:mixed-rung-tops}
		\mathfrak G^\sharp_\bullet
		:=\bigl\|\cG^{\sharp,\mathrm L_{\rm G}}_{M,a_{\rm G}}(u)\bigr\|_{L_t^2},
		\qquad
		\mathfrak G^0_\bullet:=\bigl\|\cG^{\mathrm L}_{M,2}(u)\bigr\|_{L_t^2} .
	\end{equation}
	For $1\le k\le M-1$ and
	$(a,p)\in\{(a_{\rm G},p_{\rm G}),(2,p_{\rm c})\}$,
	define the exponents by
	\begin{equation}
		\label{eq:mixed-rung-exponent}
		\frac1{\mathsf p_k(a)}
		:=\frac1{\mathsf e(M/k;a,p)}
		+2k\theta_\circ\Bigl(1-\frac kM\Bigr)\Bigl(\frac1a-\frac1p\Bigr) .
	\end{equation}
	Here $\mathsf p_k(a)$ is understood relative to the pair $(a,p)$: $p=p_{\rm G}$ for the Gaussian norm and $p=p_{\rm c}$ for the coercive norm. The exponent families may differ even when $a_{\rm G}=2$.
	There are $\Gamma_k\in[k/M,1]$ such that every ordered
	word $J\in\mathfrak A_\bullet^{k}$ satisfies
	\begin{align}
		\|w^{\mathrm L_{\rm G}}\mathfrak D_Ju\|
		_{L_t^{\mathsf e(M/k,\,r_0)}L_x^{\mathsf p_k(a_{\rm G})}L_v^2}
		&\le C\,\bigl(\mathfrak b^\sharp_{r_0}\bigr)^{1-\Gamma_k}
		\bigl(\mathfrak G^\sharp_\bullet+\mathfrak b^\sharp_2\bigr)^{\Gamma_k}
		+C\,\mathfrak b^\sharp_{r_0},
		\label{eq:mixed-word-gaussian-rungs}\\
		\|w^{\mathrm L}\mathfrak D_Ju\|
		_{L_t^{\mathsf e(M/k,\,r_0)}L_x^{\mathsf p_k(2)}H_\sigma^1}
		&\le C\,\bigl(\mathfrak b^0_{r_0}\bigr)^{1-\Gamma_k}
		\bigl(\mathfrak G^0_\bullet+\mathfrak b^0_2\bigr)^{\Gamma_k}
		+C\,\mathfrak b^0_{r_0},
		\label{eq:mixed-word-coercive-rungs}
	\end{align}
	where $C$ is independent of $u$ and $t_1-t_0$. 
	Moreover,
	\[
	\begin{gathered}
		\|w^{\mathrm L_{\rm G}}\mathfrak D_Ju\|_{L_t^{\mathsf e(M/k,\,r_0)}L_x^{\mathsf p_k(a_{\rm G})}L_v^2}<\infty
		\quad\text{whenever}\quad
		\mathfrak b^\sharp_{r_0}+\mathfrak G^\sharp_\bullet+\mathfrak b^\sharp_2<\infty,\\
		\|w^{\mathrm L}\mathfrak D_Ju\|_{L_t^{\mathsf e(M/k,\,r_0)}L_x^{\mathsf p_k(2)}H_\sigma^1}<\infty
		\quad\text{whenever}\quad
		\mathfrak b^0_{r_0}+\mathfrak G^0_\bullet+\mathfrak b^0_2<\infty.
	\end{gathered}
	\]
\end{corollary}

\begin{proof}
	We prove \eqref{eq:mixed-word-gaussian-rungs}, at
	$(a,p)=(a_{\rm G},p_{\rm G})$ and with $\mathfrak b_\rho:=\mathfrak b^\sharp_\rho$,
	$\mathfrak G:=\mathfrak G^\sharp_\bullet$; \eqref{eq:mixed-word-coercive-rungs}
	is the same at $(a,p)=(2,p_{\rm c})$,
	with Theorem~\ref{thm:complete-flow-GN}\,(\ref{it:GN-coercive}) in place of
	Theorem~\ref{thm:complete-flow-GN}\,(\ref{it:GN-fiber}), $\cG^{\mathrm L}_{M,2}$ in place of
	$\cG^{\sharp,\mathrm L_{\rm G}}_{M,a_{\rm G}}$ and $\mathrm L$ in place of $\mathrm L_{\rm G}$.  For $2\le r\le\infty$, let $P_k^{r}$ denote the supremum of the left
	side of \eqref{eq:mixed-word-gaussian-rungs}, with $r$ in place of $r_0$, over the ordered words of length
	$k$, with $\mathsf p_0(a):=p$, the value of \eqref{eq:mixed-rung-exponent} at $k=0$. Thus $P_0^{r_0}=\mathfrak b_{r_0}$. Set
	$\mathfrak B_{\mathrm{tot}}:=\mathfrak G+\mathfrak b_2$.  We may assume $\mathfrak b_{r_0}+\mathfrak B_{\mathrm{tot}}<\infty$; otherwise the assertion is immediate.  Then
	every word $J$ of $\mathfrak A_\bullet$ with $|J|\le M$ satisfies
	$\|w^{\mathrm L_{\rm G}}\mathfrak D_Ju\|_{L_t^2L_x^{a}L_v^2}\le\mathfrak B_{\mathrm{tot}}$. Hence
	$w^{\mathrm L_{\rm G}}\mathfrak D_Ju\in L^1((t_0,t_1);L_x^1L_v^2)$, which is the
	integrability required by Lemma~\ref{lem:temporal-letter-unequal} above.
	
	Fix $1\le k\le M-1$, write $J'=\nu J$ with $|J|=k-1$ and $\nu\in\mathfrak A_\bullet$,
	and put $m':=M-k+1\ge2$. When $\nu$ is fixed below, $P_k^{r_0}$ denotes the supremum over words $J'$ of the indicated type.  By $\theta_\circ\le1/(2M)$ and $k-1\le M$,
	\begin{equation}
		\label{261016920}
		2(k-1)\theta_\circ\le1,
		\qquad\text{hence}\qquad
		\frac1{\mathsf p_{k-1}(a)}
		=\frac1a-\Bigl(1-\frac{k-1}M\Bigr)\bigl(1-2(k-1)\theta_\circ\bigr)
		\Bigl(\frac1a-\frac1p\Bigr)
		\le\frac1a ,
	\end{equation}
	so $a\le\mathsf p_{k-1}(a)$, and $\theta_\circ\le\tfrac1{2M}\le\tfrac1{m'}$.
	
	If $\nu=\mathrm t$, apply Lemma~\ref{lem:temporal-letter-unequal} on $I=(t_0,t_1)$, $|I|=t_1-t_0\ge1$, to
	$v=w^{\mathrm L_{\rm G}}\mathfrak D_Ju$, respectively $v=w^{\mathrm L}\mathfrak D_Ju$, at
	$(m,i)=(m',1)$, at the base exponents
	\[
	(\rho_{\mathrm{in}},p_0)=(\mathsf e(M/(k-1),r_0),\mathsf p_{k-1}(a))
	\]
	and at the top
	exponent $a$.  Its hypotheses \eqref{eq:temporal-unequal-hypotheses} hold:
	$\mathsf e(M/(k-1),r_0)\ge2$, $\theta_\circ>0$ and, by the previous paragraph,
	$$
	1\le a\le\mathsf p_{k-1}(a)\qquad \text{and}\qquad \gamma=1/m'+\theta_\circ\le2/m'\le1.
	$$
	We take $H=L_v^2$ for $v=w^{\mathrm L_{\rm G}}\mathfrak D_Ju$ and 
	$H=H_\sigma^1$ for $v=w^{\mathrm L}\mathfrak D_Ju$.  For $0\le i\le m'-1$, the word
	$\mathrm t^{i+1}J$ has length $|\mathrm t^{i+1}J|\le(k-1)+m'=M\le M'$, so
	\eqref{eq:mixed-rung-chain} at $\nu=\mathrm t$, with $\mathrm t^iJ$ in place of $J$, gives
	$\mathfrak D_{\mathrm t}(\mathfrak D_{\mathrm t^iJ}u)=\mathfrak D_{\mathrm t^{i+1}J}u$.
	Since $\mathfrak D_{\mathrm t}=\partial_t$ commutes with multiplication by
	$w^{\mathrm L_{\rm G}}$ and $w^{\mathrm L}$, these identities for $0\le i\le m'-1$ give
	\[
	\mathfrak D_{\mathrm t}^{\,m'}\bigl(w^{\mathrm L_{\rm G}}\mathfrak D_Ju\bigr)=w^{\mathrm L_{\rm G}}\mathfrak D_{\mathrm t^{m'}J}u,
	\qquad
	\mathfrak D_{\mathrm t}^{\,m'}\bigl(w^{\mathrm L}\mathfrak D_Ju\bigr)=w^{\mathrm L}\mathfrak D_{\mathrm t^{m'}J}u .
	\]
	For a word $K$ of $\mathfrak A_\bullet$ with $|K|\le M$, the bound $w\ge1$, respectively
	\eqref{eq:real-weight-Hsigma-monotonicity} with $(\mathrm L',\mathrm L)$ replaced by
	$(\mathrm L,\mathrm L+M-|K|)$, gives, at almost every $t$,
	\[
	\|w^{\mathrm L_{\rm G}}\mathfrak D_Ku\|_{L_x^{a_{\rm G}}L_v^2}\le\cG^{\sharp,\mathrm L_{\rm G}}_{M,a_{\rm G}}(u),
	\qquad
	\|w^{\mathrm L}\mathfrak D_Ku\|_{L_x^2H_\sigma^1}\le C\,\cG^{\mathrm L}_{M,2}(u) .
	\]
	We apply these bounds at $K=J$ and $K=\mathrm t^{m'}J$, of lengths $k-1$ and $M$.  By the
	identities above, the last two quantities in \eqref{eq:temporal-unequal-quantities} at
	$m=m'$ satisfy
	\[
	\|v\|_{L^2_tL^{a}_xH}+\|\mathfrak D_{\mathrm t}^{\,m'}v\|_{L^2_tL^{a}_xH}\le C\mathfrak B_{\mathrm{tot}} .
	\]
	The temporal output of the lemma is
	\[
	\mathsf e(m',\mathsf e(M/(k-1),r_0))=\mathsf e(M/k,r_0) ,
	\] by \eqref{eq:GN-rung-exponent-composition} with $(m,j,k)$ replaced by $(M,k-1,1)$ and $(a,p)$ by $(2,r_0)$. Its spatial
	output $p_1$ is given by
	\[
	\frac1{p_1}
	=\frac1{m'}\cdot\frac1a+\Bigl(1-\frac1{m'}\Bigr)\frac1{\mathsf p_{k-1}(a)}
	+\theta_\circ\Bigl(\frac1a-\frac1{\mathsf p_{k-1}(a)}\Bigr)
	\le\frac1{\mathsf p_k(a)} ,
	\]
	and the inequality follows from $M-k\ge1$, $a\le p$ and
	\begin{equation}
		\label{261016301}
		M\Bigl(\frac1{p_1}-\frac1{\mathsf p_k(a)}\Bigr)
		=-\bigl(\theta_\circ(M-k-1)+2(k-1)(M-k+1)\theta_\circ^2\bigr)
		\Bigl(\frac1a-\frac1p\Bigr).
	\end{equation}
	By \eqref{eq:temporal-unequal-output} at $(m,i)=(m',1)$ and $p_0=\mathsf p_{k-1}(a)$, the formula \eqref{261016920} for $1/\mathsf p_{k-1}(a)$ and $1-\frac{k-1}M=\frac{m'}M$,
	\begin{align}
		\frac1{p_1}&=\frac1a-\Bigl(1-\frac1{m'}-\theta_\circ\Bigr)\Bigl(\frac1a-\frac1{\mathsf p_{k-1}(a)}\Bigr)\notag\\
		&=\frac1a-\Bigl(1-\frac kM-\frac{(M-k+1)\theta_\circ}M\Bigr)\bigl(1-2(k-1)\theta_\circ\bigr)\Bigl(\frac1a-\frac1p\Bigr).\label{261016302}
	\end{align}
	By \eqref{eq:mixed-rung-exponent} and \eqref{eq:GN-output-exponents},
	\begin{equation}
		\label{261016303}
		\frac1{\mathsf p_k(a)}=\frac1a-\Bigl(1-\frac kM\Bigr)(1-2k\theta_\circ)\Bigl(\frac1a-\frac1p\Bigr).
	\end{equation}
	Subtracting \eqref{261016303} from \eqref{261016302} and multiplying by $M$ gives \eqref{261016301}.
	Hence, since $\cD$ is bounded, \eqref{eq:temporal-unequal} gives
	\begin{equation}
		\label{eq:mixed-rung-step-temporal}
		P_k^{r_0}\le C\bigl(P_{k-1}^{r_0}\bigr)^{1-\gamma}\mathfrak B_{\mathrm{tot}}^{\gamma}
		+C\,P_{k-1}^{r_0},
		\qquad \gamma=\frac1{m'}+\theta_\circ .
	\end{equation}
	
	If $\nu\in\mathfrak A$, write, with
	$J_j$ the final subword of $J'$ of length $j$ and
	$\cG^{\sharp,\mathrm L_{\rm G}}_{M,a_{\rm G}}$ of \eqref{eq:mixed-rung-tops},
	\[
	\mathfrak P_j(t):=\|w^{\mathrm L_{\rm G}}\mathfrak D_{J_j}u(t)\|
	_{L_x^{\mathsf p_j(a)}L_v^2},
	\qquad
	\mathfrak g:=\cG^{\sharp,\mathrm L_{\rm G}}_{M,a_{\rm G}}(u) .
	\]
	By \eqref{eq:mixed-rung-chain} and \eqref{eq:mixed-rung-local-integrability},
	Lemma~\ref{lem:temporal-slicing} on $(t_0,t_1)$ applies at $z=\mathfrak D_Ju$, with $\cI$ the
	ordered complete words of length at most $m'$, and gives
	$\mathfrak D_K(\mathfrak D_Ju(t))=\mathfrak D_{KJ}u(t)$ for $|K|\le m'$ at almost every $t$.
	It suffices to consider times at which $\mathfrak P_{k-1}(t)<\infty$; otherwise the asserted bound is immediate.
	Then Theorem~\ref{thm:complete-flow-GN}\,(\ref{it:GN-fiber}) at $X=\tau_\nu$ applies at
	almost every such $t$ to $\mathfrak D_Ju(t)$ at $(m,i)=(m',1)$, the weight $\mathrm L_{\rm G}$ and the base
	exponent $\mathsf p_{k-1}(a)$. This is admissible since $m'\ge2$, $(k-1)+m'=M\le M'\le\mathrm N$,
	$a\le\mathsf p_{k-1}(a)$, the top $\cG^{\sharp,\mathrm L_{\rm G}}_{m',a}(\mathfrak D_Ju(t))$ of the theorem is at most $\mathfrak g(t)$,
	and $\mathfrak g(t)<\infty$ at almost every $t$ since $\mathfrak G<\infty$. The bound by $\mathfrak g(t)$ holds because the term of
	$\mathfrak D_K(\mathfrak D_Ju(t))$ in that top carries the weight
	\[
	\mathrm L_{\rm G}+m'-|K|=\mathrm L_{\rm G}+M-|KJ| ,
	\] which is the weight of
	$\mathfrak D_{KJ}u(t)$ in $\mathfrak g(t)$.  The theorem has
	spatial output $\mathsf e(m';a,\mathsf p_{k-1}(a))$.  By
	\eqref{eq:GN-output-exponents} and $a\le\mathsf p_{k-1}(a)$, with $p_1$ of the
	case $\nu=\mathrm t$,
	\[
	\frac1{\mathsf e(m';a,\mathsf p_{k-1}(a))}
	=\frac1{p_1}-\theta_\circ\Bigl(\frac1a-\frac1{\mathsf p_{k-1}(a)}\Bigr)
	\le\frac1{p_1}\le\frac1{\mathsf p_k(a)} ,
	\]
	so, since $\cD$ is bounded, the theorem yields
	\[
	\mathfrak P_k(t)\le C\mathfrak P_{k-1}(t)^{1-1/m'}\mathfrak g(t)^{1/m'}
	+C\mathfrak P_{k-1}(t).
	\]
	
	First take $r_0=2$, when every temporal output $\mathsf e(M/k,2)$ equals $2$. Hence, in the case $\nu\in\mathfrak A$ of the preceding estimate, both factors and the output carry the exponent $2$ in $t$, and
	H\"older's inequality in $t$ gives
	\[
	P_k^{2}\le C(P_{k-1}^{2})^{1-1/m'}\mathfrak B_{\mathrm{tot}}^{1/m'}+CP_{k-1}^{2} .
	\]
	Together with \eqref{eq:mixed-rung-step-temporal}, and with
	\eqref{eq:rung-exponent-monotonicity} applied in the case $\nu\in\mathfrak A$ at
	\[
	(b,G,\alpha,\beta)=(P_{k-1}^{2},\mathfrak B_{\mathrm{tot}},1/m',1/m'+\theta_\circ),
	\]
	both cases give
	\[
	P_k^{2}\le C\bigl(P_{k-1}^{2}\bigr)^{1-\gamma_k}\mathfrak B_{\mathrm{tot}}^{\gamma_k}
	+C\,P_{k-1}^{2},
	\qquad \gamma_k:=\frac1{m'}+\theta_\circ .
	\]
	Set $\Gamma_0:=0$ and $1-\Gamma_k:=(1-\gamma_k)(1-\Gamma_{k-1})$. Since
	\[
	\prod_{1\le j\le k}(1-1/(M-j+1))=1-k/M\quad\text{and}\quad 1/(M-j+1)\le\gamma_j\le1\quad\text{(as noted above)} ,
	\]
	the exponents satisfy $k/M\le\Gamma_k\le1$. Since both factors in the definition of $1-\Gamma_k$ lie in $[0,1]$, also $\Gamma_{k-1}\le\Gamma_k$ and $\gamma_k\le\Gamma_k$.  We argue by induction on $k$, starting from $P_0^{2}=\mathfrak b_2$.  By
	$$
	(x+y)^{1-\gamma_k}\le x^{1-\gamma_k}+y^{1-\gamma_k}\qquad \text{and}\qquad (1-\Gamma_{k-1})(1-\gamma_k)=1-\Gamma_k,
	$$
	we have
	\[
	\bigl(\mathfrak b_2^{1-\Gamma_{k-1}}\mathfrak B_{\mathrm{tot}}^{\Gamma_{k-1}}+\mathfrak b_2\bigr)^{1-\gamma_k}\mathfrak B_{\mathrm{tot}}^{\gamma_k}
	\le\mathfrak b_2^{1-\Gamma_k}\mathfrak B_{\mathrm{tot}}^{\Gamma_k}+\mathfrak b_2^{1-\gamma_k}\mathfrak B_{\mathrm{tot}}^{\gamma_k} .
	\]
	Inserting the induction hypothesis $P_{k-1}^{2}\le C\bigl(\mathfrak b_2^{1-\Gamma_{k-1}}\mathfrak B_{\mathrm{tot}}^{\Gamma_{k-1}}+\mathfrak b_2\bigr)$ into the recursion above and applying
	\eqref{eq:rung-exponent-monotonicity} at $(\alpha,\beta)=(\gamma_k,\Gamma_k)$ and $(\Gamma_{k-1},\Gamma_k)$ gives
	\begin{equation}
		\label{eq:mixed-rung-base-two}
		P_k^{2}\le C\bigl(\mathfrak b_2^{1-\Gamma_k}\mathfrak B_{\mathrm{tot}}^{\Gamma_k}
		+\mathfrak b_2\bigr)\le C\,\mathfrak B_{\mathrm{tot}} ,
	\end{equation}
	where the last inequality follows from $\mathfrak b_2\le\mathfrak B_{\mathrm{tot}}$ and $\Gamma_k\le1$. 
	
	For arbitrary $r_0$, put $\rho_k:=\mathsf e(M/k,r_0)$.  In the
	case $\nu\in\mathfrak A$ of the spatial one-letter estimate, take the $L_t^{\rho_k}$ norm and apply H\"older's inequality with the
	exponents $\bigl(\tfrac{\rho_{k-1}}{1-1/m'},\tfrac{2}{1/m'}\bigr)$, admissible
	because $\tfrac1{\rho_k}=\tfrac{1-1/m'}{\rho_{k-1}}+\tfrac1{2m'}$.  For the lower-order term, H\"older's inequality with the same exponents, applied to
	$\mathfrak P_{k-1}=\mathfrak P_{k-1}^{1-1/m'}\mathfrak P_{k-1}^{1/m'}$, gives
	\[
	\|\mathfrak P_{k-1}\|_{L_t^{\rho_k}}
	\le\bigl(P_{k-1}^{r_0}\bigr)^{1-1/m'}\bigl(P_{k-1}^{2}\bigr)^{1/m'}
	\le C\bigl(P_{k-1}^{r_0}\bigr)^{1-1/m'}\mathfrak B_{\mathrm{tot}}^{1/m'}
	\]
	by \eqref{eq:mixed-rung-base-two}. Thus the case $\nu\in\mathfrak A$ satisfies
	\eqref{eq:mixed-rung-step-temporal} with $\gamma=1/m'$, and both cases give
		the same recursion as for $r_0=2$, now with $P^{r_0}$ in place of $P^{2}$: $P_k^{r_0}\le C\bigl(P_{k-1}^{r_0}\bigr)^{1-\gamma_k}\mathfrak B_{\mathrm{tot}}^{\gamma_k}+C\,P_{k-1}^{r_0}$.
	
	The same induction, starting from $P^{r_0}_0=\mathfrak b_{r_0}$, gives
	$P^{r_0}_k\le C(\mathfrak b_{r_0}^{1-\Gamma_k}\mathfrak B_{\mathrm{tot}}^{\Gamma_k}
	+\mathfrak b_{r_0})$, which is \eqref{eq:mixed-word-gaussian-rungs}.
\end{proof}

\section{\texorpdfstring{Construction of the Stationary State under Diffuse Reflection}{Construction of the Stationary State under Diffuse Reflection}}\label{sec:diffuse-common}

We first establish the weighted stationary estimates and remove the spectral cutoff. The nonlinear construction and proof of Theorem~\ref{thm:stationary-diffuse} follow in \S\ref{ss:diffuse-fixed-point}.  The \emph{structural data} are $\cD$, the fixed fields $\cX_{\rm c}$ of Proposition~\ref{prop:complete-fields}, the order $\mathrm N$, and the target weight $\ell_\sharp$ of \eqref{eq:finite-hierarchy-weight-choice}. The exponents, weights, and $\alpha_{\rm bs}$ of \eqref{eq:holder-exponent} depend only on these data.  Unless stated otherwise, $C$ and the constants implicit in $\lesssim$ and $\simeq$ depend only on the structural data and may change from line to line.  A subscript records an additional dependence, as in $C_k$, $C_\ell$ or $C_L$.

In this section the wall is stationary: $\Theta_w$,
$\mathfrak m_{\Theta_w}$, $d_\delta$, $P_w$ and $h_w$ denote
$\Theta_w^\infty$, $\mathfrak m_{\Theta_w^\infty}$, $d_\delta^\infty$,
$P_w^\infty$ and $h_w^\infty$ of \eqref{eq:limiting-wall-profiles}, so that the
wall data below are those of the boundary law of
\eqref{eq:stationary-full-diffuse-problem}.

\subsection{Exponents, weights, and norms}\label{ss:weight-interpolation}

The following exponents and weights enter the reduced norm $\widetilde{\mathfrak Q}_{\mathrm N}$ defined in \eqref{eq:reduced-Q3}.  For $2\le p\le\infty$, the interpolation exponents $\mathsf e(m,p)$ of \eqref{eq:GN-output-exponents} lie between $2$ and $p$ and are nondecreasing in $m$. Fix
\begin{equation}
	\label{eq:dyn-own-order-endpoints}
	q_\flat:=\tfrac{30}{13},
	\qquad
	p_\flat:=18,
	\qquad\text{so that}\qquad
	\delta_\flat:=\frac12-\frac1{q_\flat}-\frac1{p_\flat}=\frac1{90}>0 .
\end{equation}
The strict inequality $1/q_\flat+1/p_\flat<1/2$ leaves room in spatial H\"older estimates.  On the bounded domain $\cD$, H\"older's inequality gives $\|u_1u_2\|_{L_x^2}\le|\cD|^{\delta_\flat}\|u_1\|_{L_x^{q_\flat}}\|u_2\|_{L_x^{p_\flat}}$; we use the same estimate for the sub-conjugate pairs below. Also,
\begin{equation*}
	2<q_\flat<\tfrac73 .
\end{equation*}

The first argument of the nonlinear collision operator is measured in $w^{-2}L_v^2$; this corresponds to $\mathrm L=-2$ in \eqref{eq:velalg-poly}. For the deformation terms we use \eqref{eq:deformation-weighted-output-spatial} and $\mu^\varpi\le Cw^{-2}$.

For jets of order $j$, the first factor in a bilinear collision term is measured in $L_x^{\mathsf g(j)}$, while a differentiated solution factor is measured in $L_x^{\mathsf s(j)}$, where
\begin{equation}
	\label{eq:dyn-own-order-exponents}
	\mathsf g(j):=\mathsf e(\mathrm N/j;q_\flat,\infty),
	\qquad
	\mathsf s(j):=\mathsf e(\mathrm N/j,p_\flat),
	\qquad 0\le j\le\mathrm N,
\end{equation}
with the convention of \eqref{eq:GN-output-exponents} at $j=0$.  These are the exponents in \eqref{eq:ordered-word-gaussian-rungs} and \eqref{eq:ordered-word-coercive-rungs-weighted}, respectively, with $M=\mathrm N$, $a_{\rm G}=q_\flat$, $(p_{\rm G},p_{\rm c})=(\infty,p_\flat)$ as in \eqref{eq:dynamic-pointwise-finite-norms}, and word length $j$.
Since 
$$
\frac1{\mathsf g(k)}=\frac k{\mathrm N}\cdot\frac1{q_\flat},\qquad \frac1{\mathsf s(j)}=\frac1{p_\flat}+\frac j{\mathrm N}\bigl(\frac12-\frac1{p_\flat}\bigr),\qquad \frac1{q_\flat}=\frac12-\frac1{p_\flat}-\delta_\flat
$$
by \eqref{eq:dyn-own-order-exponents} and \eqref{eq:dyn-own-order-endpoints}, we have, for $k+j=m\le\mathrm N$,
\begin{equation}
	\label{eq:dyn-own-order-pairing}
	\frac1{\mathsf g(k)}+\frac1{\mathsf s(j)}
	=\frac m{\mathrm N}\Bigl(\frac12-\frac1{p_\flat}\Bigr)+\frac1{p_\flat}
	{}-\frac k{\mathrm N}\delta_\flat .
\end{equation}
For a graph row $U_m$ or $W_m$ of order $m$, the solution jet of order $j\le m$ instead uses $L_x^{\mathsf s(j\mathrm N/m)}$, where $\mathsf s(j\mathrm N/m)=\mathsf e(m/j,p_\flat)$, with the value $p_\flat$ at $j=0$. This is \eqref{eq:ordered-word-coercive-rungs-weighted} with top order $M=m$ and base exponent $p_{\rm c}=p_\flat$. Since $\frac1{\mathsf s(j\mathrm N/m)}=\frac1{p_\flat}+\frac jm\bigl(\frac12-\frac1{p_\flat}\bigr)$, with $j=m-k$ and the expression for $1/q_\flat$ above, the computation leading to \eqref{eq:dyn-own-order-pairing} gives
\begin{align}
	\label{eq:dyn-row-pairing}
	\begin{aligned}
		\frac1{\mathsf g(k)}+\frac1{\mathsf s(j\mathrm N/m)}
		=\frac12-k\Bigl(\frac12-\frac1{p_\flat}\Bigr)
		\Bigl(\frac1m-\frac1{\mathrm N}\Bigr)
		{}-\frac k{\mathrm N}\delta_\flat
		\ <\ \frac12 \qquad\qquad\qquad &\\
		\text{whenever}
		\quad k+j=m\le\mathrm N,\quad k\ge1.&
	\end{aligned}
\end{align}

Fix an integer target weight $\ell_\sharp\ge0$. The largest weight loss from the deformation terms is
\begin{equation}
	\label{eq:Adef}
	A_\diamond:=2\mathrm N+2 ,
\end{equation}
the value of $2m+2$ in \eqref{eq:deformation-weighted-output} at $m=\mathrm N$.  This weight covers every source in the differentiated system: the original-Coulomb and compact terms in \eqref{eq:velalg} and \eqref{eq:iterated-compact-full-operator-bound} incur no additional loss; the deformation terms in \eqref{eq:deformation-weighted-output} and \eqref{eq:local-L0-complete-bound} lose at most $2m+2\le A_\diamond$; and $\cV_{I,J}$ in \eqref{eq:ordered-transport-weight} loses at most $4\le A_\diamond$. Set
\begin{equation}
	\label{eq:finite-hierarchy-weight-choice}
	\begin{gathered}
		\ell_{\mathrm N}:=\ell_\sharp,
		\qquad
		\ell_m:=\ell_{m+1}+A_\diamond\quad(0\le m\le\mathrm N-1),
		\qquad
		A_\flat:=\ell_0+3 .
	\end{gathered}
\end{equation}

These nonnegative integer weights satisfy
\begin{equation}
	\label{eq:finite-hierarchy-strict-ladder}
	\ell_0>\ell_1>\cdots>\ell_{\mathrm N-1}>\ell_{\mathrm N}=\ell_\sharp,
	\qquad
	\ell_{\mathrm N-1}\ge A_\diamond.
\end{equation}

Let $\alpha_{\rm bs}$ be the kinetic H\"older exponent given by Theorem~\ref{260526316} for $\cD$:
\begin{equation}
	\label{eq:holder-exponent}
	\alpha_{\rm bs}:=\alpha(\cD)\in(0,1);
\end{equation}
it is the exponent denoted $\alpha_0$ in Lemma~\ref{260527559}.

\begin{prop}[Undifferentiated prescribed-inflow regularity]
	\label{prop:base-regularity}
	Let $\ell\ge0$ be fixed, and let $B^{\st}_{\rm in}:=B^{\st}_1(\alpha_{\rm bs})$ be the
	gradient weight loss \eqref{eq:stationary-B1} of Lemma~\ref{260527559} at $\alpha=\alpha_{\rm bs}$.
	Set $\varepsilon_{\rm bs}:=\varepsilon_{\rm in}$, the coefficient threshold
	of \eqref{eq:weighted-inflow-estimate}.
	Assume
	\begin{equation}
		\label{eq:KS-coefficient}
		\|\mu^{1/4}g\|_{L^\infty}\le\varepsilon_{\rm bs},
		\qquad
		\|\mu^{1/4}g\|_{C^{\alpha_{\rm bs}}_{\kin}}\le M<\infty.
	\end{equation}
	If $U\in\cH_{\kin,\sigma}$ is the solution of
	\[
	v\cdot\nabla_x U+L_gU=0,\qquad U_{\gamma^-}=H,
	\]
	with $L_g$ of \eqref{def:L_g}, and the incoming datum satisfies
	$\|w^{\ell+B^{\st}_{\rm in}}H\|_{(C_x^{0,1}\cap C_v^{1,1})(\gamma^-)}<\infty$, then for any $1<p<\infty$,
	\begin{equation}
		\label{eq:base-regularity}
		\|w^\ell U\|_{\cH_{\kin,\sigma}}
		+\|w^\ell U\|_{L^\infty}
		+[w^\ell U]_{C^{\alpha_{\rm bs}}_{\kin}}
		+\|w^\ell \nabla_vU\|_{L_x^pL_v^2}
		\le C_{\ell,M,p}\|w^{\ell+B^{\st}_{\rm in}}H\|_{(C_x^{0,1}\cap C_v^{1,1})(\gamma^-)}.
	\end{equation}
\end{prop}

\begin{proof}
	The weighted prescribed-inflow estimate
	\eqref{eq:weighted-inflow-estimate}, at output weight $\ell$, gives
	\[
	\|w^\ell U\|_{\cH_{\kin,\sigma}}
	\le C_\ell\|w^\ell H\|_{L^2(\gamma^-)}.
	\]
	Since $B^{\st}_{\rm in}\ge B^{\st}_0=37$ by \eqref{eq:stationary-B1} and \eqref{eq:stationary-B0}, the supremum
	bound on $w^{\ell+B^{\st}_{\rm in}}H$ and
	$\int_{\bR^3}w^{-2B^{\st}_{\rm in}}|n_x\cdot v|\dd v<\infty$
	give
	\[
	\|w^\ell H\|_{L^2(\gamma^-)}\le C_\ell \|w^{\ell+B^{\st}_{\rm in}}H\|_{(C_x^{0,1}\cap C_v^{1,1})(\gamma^-)}.
	\]
	The weighted inflow estimate at the threshold $\varepsilon_{\rm bs}=\varepsilon_{\rm in}$ gives the solution in $\cH_{\kin,\sigma}$ and its completion-defined traces.  Lemma~\ref{260527559} then yields the weighted $L^\infty$ and kinetic H\"older
	bounds, and the mixed velocity-gradient bound.
	Uniqueness in that lemma identifies its solution with $U$, since both satisfy the same equation and incoming condition.  The first two regularity
	conclusions of that lemma follow from $\|\mu^{1/4}g\|_\infty\le\varepsilon_{\rm bs}
	=\varepsilon_0(\cD)$, while the
	$D_vL_x^p L_v^2$ conclusion uses the undifferentiated kinetic H\"older bound
	in \eqref{eq:KS-coefficient}.
	The H\"older exponent of that lemma is
	$\alpha_0=\alpha_{\rm bs}$ by \eqref{eq:holder-exponent}.  This gives
	\eqref{eq:base-regularity}.
\end{proof}

All kinetic H\"older seminorms below have exponent $\alpha_{\rm bs}$ from \eqref{eq:holder-exponent}. Define
\begin{alignat}{3}
	&U_m(u)&&:=\sum_{I\in\mathfrak A^m}
	\|\mathfrak D_Iu\|_{\cH_{\kin,\sigma}},
	\qquad &&0\le m\le\mathrm N,
	\label{eq:order-two-Um}\\
	&W_{m}(u)&&:=\sum_{I\in\mathfrak A^m}
	\|w^{\ell_m}\mathfrak D_Iu\|_{\cH_{\kin,\sigma}},
	\qquad &&0\le m\le\mathrm N,
	\label{eq:order-two-Wm}\\
	&R_m(u)&&:=\sum_{I\in\mathfrak A^m}
	\|w^{\,\mathrm N-m-2}\mathfrak D_Iu\|_{L_x^{q_\flat}L_v^2},
	\qquad &&1\le m\le\mathrm N,
	\label{eq:order-two-Rm}\\
	&\cB_{\rm bs}(u)&&:=
	\|w^{A_\flat}u\|_{L^\infty}
	+[w^{A_\flat}u]_{C^{\alpha_{\rm bs}}_{\kin}}
	+\|w^{A_\flat}\nabla_vu\|_{L_x^{p_\flat}L_v^2},
	\label{eq:base-strong}\\
	&\widetilde{\mathfrak Q}_{\mathrm N}(u)&&:=\cB_{\rm bs}(u)
	+\sum_{m=0}^{\mathrm N}W_{m}(u),
	\label{eq:reduced-Q3}
\end{alignat}
for measurable $u$ on $\cD\times\bR_v^3$. The quantities $U_m,W_m$ are the unweighted and weighted graph rows, $R_m$ the mixed rows, $\cB_{\rm bs}$ the strong norm, and $\widetilde{\mathfrak Q}_{\mathrm N}$ the reduced norm. We set $\mathfrak A^0=\{\varnothing\}$ and $\mathfrak D_\varnothing={\rm Id}$. All these quantities are extended-valued: a row is finite only if its distributional jets belong to the indicated spaces, including their completion-defined Green traces; otherwise it equals $+\infty$.

Since $w^{-2-A_\flat}\in L_v^2$, H\"older's inequality in $v$ gives, for every measurable $u$ on $\cD\times\bR_v^3$,
\begin{equation}
	\label{eq:base-fiber-Lsharp}
	\|w^{-2}u\|_{L_x^\infty L_v^2}
	\le\|w^{-2-A_\flat}\|_{L_v^2}\|w^{A_\flat}u\|_{L^\infty_{x,v}}
	\le C\cB_{\rm bs}(u)
	\le C\widetilde{\mathfrak Q}_{\mathrm N}(u).
\end{equation}

\begin{lemma}[Mixed integrability from the graph norm]
	\label{lem:dynamic-gaussian-mixed-gain}
	Let $k\ge0$.  Every $h\in\cH_{\kin,\sigma}(\cD\times\bR_v^3)$ satisfies
	\begin{equation}
		\label{eq:graph-to-mixed-integrability}
		\|w^{-2-k}h\|_{L_x^{q_\flat}L_v^2}
		\le C_k\norm h_{\cH_{\kin,\sigma}(\cD\times\bR_v^3)},
	\end{equation}
	and, for every $0<T<\infty$, every
	$h\in\cH_{\kin,\sigma}(\cD_T\times\bR_v^3)$ satisfies
	\begin{equation}
		\label{eq:graph-to-mixed-integrability-evolutionary}
		\|w^{-2-k}h\|_{L_t^{\varrho}L_x^{q_\flat}L_v^2}
		\le C_k\|h\|_{\cH_{\kin,\sigma}(\cD_T\times\bR_v^3)},
		\qquad 2\le\varrho\le q_\flat ,
	\end{equation}
	with $C_k$ depending only on the structural data and on $k$, and not on $T$ or on $\varrho$.
\end{lemma}

\begin{proof}
	Write $w^{-2-k}=m_1m_2$ with $m_1=w^{-3/2}$ and $m_2=w^{-1/2-k}$, and choose $r$ by $1/2=1/r+1/q_\flat$. Since $q_\flat<7/3$, we have $r>14$ and $r(1/2+k)>3$, hence $m_2\in L_v^r$.
	
	We use the elementary bounds for $m_1$.  Here $m_1^2=\la v\ra^{-3}$ and
	$|\nabla_vm_1|^2\lesssim\la v\ra^{-5}\le\la v\ra^{-1}$.  Then
	\eqref{eq:coulomb-maxwellian-anisotropy} and
	\eqref{eq:local-H1-H_sigma} give
	\[
	\|m_1h\|_{L_v^2}^2+\|\nabla_v(m_1h)\|_{L_v^2}^2
	\lesssim
	\int\!\left(\la v\ra^{-3}|\nabla_vh|^2
	+\la v\ra^{-1}|h|^2\right)\dd v
	\lesssim\|h\|_{H_\sigma^1}^2,
	\]
	and therefore
	$ \|m_1h\|_{H_v^1}\le C\|h\|_{H_\sigma^1}$.
	By the multiplier bound \eqref{eq:Hsigma-multiplier} at $\zeta=m_1$ followed by the embedding $H_\sigma^{-1}\subset H_v^{-1}$ of \eqref{eq:Hv-into-Hsigma}, by $m_1=m_1(v)$ and by $\|m_1\|_\infty<\infty$, respectively,
	\begin{equation}
		\label{eq:mixed-multiplier-dual-transport-trace}
		\begin{gathered}
			\|m_1F\|_{H_v^{-1}}\le C\|F\|_{H_\sigma^{-1}},\qquad
			v\cdot\nabla_x(m_1h)=m_1(v\cdot\nabla_xh),\\
			\|(m_1h)_{\gamma^-}\|_{L^2(\gamma^-)}\le
			\|m_1\|_\infty\|h_{\gamma^-}\|_{L^2(\gamma^-)}.
		\end{gathered}
	\end{equation}

	For the evolutionary bound, the $H_v^1$ estimate and the first inequality in \eqref{eq:mixed-multiplier-dual-transport-trace}, applied pointwise in $(t,x)$, give
	\[
	\|m_1h\|_{L_t^2L_x^2H_v^1}
	\le C\|h\|_{L_t^2L_x^2H_\sigma^1},
	\qquad
	\|m_1F\|_{L_t^2L_x^2H_v^{-1}}
	\le C\|F\|_{L_t^2L_x^2H_\sigma^{-1}}.
	\]
	Moreover the initial values and the incoming traces defined by \eqref{lim_transp_sigma} satisfy
	\begin{equation}
		\label{eq:dynamic-mixed-multiplier-incoming}
		\begin{aligned}
			\|(m_1h)(0)\|_{L^2(\cD\times\bR_v^3)}
			&\le \|m_1\|_\infty\|h(0)\|_{L^2(\cD\times\bR_v^3)},
			\\
			\|(m_1h)_{\gamma_T^-}\|_{L^2(\gamma_T^-)}
			&\le \|m_1\|_\infty\|h_{\gamma_T^-}\|_{L^2(\gamma_T^-)}.
		\end{aligned}
	\end{equation}
	Using $m_1=m_1(v)$ and applying the spacetime multiplier bounds above and
	\eqref{eq:dynamic-mixed-multiplier-incoming} to a defining sequence in
	Definition~\ref{def:H_kin_sigma}, we obtain
	\[
	\|m_1h\|_{\cH_\kin(\cD_T\times\bR_v^3)}
	\le C\|h\|_{\cH_{\kin,\sigma}(\cD_T\times\bR_v^3)},
	\]
	with $C$ depending only on $\cD$.  Hence
	\begin{align*}
		\|w^{-2-k}h\|_{L_t^{\varrho}L_x^{q_\flat}L_v^2}
		\,&=\|m_2\,(m_1h)\|_{L_t^{\varrho}L_x^{q_\flat}L_v^2}
		\\
		&\le C_k\|m_1h\|_{\cH_\kin(\cD_T\times\bR_v^3)}
		\le C_k\|h\|_{\cH_{\kin,\sigma}(\cD_T\times\bR_v^3)},
		\qquad 2\le\varrho\le q_\flat,
	\end{align*}
	where the first inequality follows from H\"older's inequality in $v$ and \eqref{eq:dynamic-global-mixed-integrability} with $\varrho$ in place of $r$ and $q=q_\flat$, both applied to
	$u=m_1h$.
	This is \eqref{eq:graph-to-mixed-integrability-evolutionary}.
	
	The stationary estimate \eqref{eq:graph-to-mixed-integrability} follows without the time variable: the multiplier bounds hold in the stationary graph space, and \eqref{eq:stationary-strong-subcritical-sharp} at $q_\flat<12/5$ replaces \eqref{eq:dynamic-global-mixed-integrability} to yield \eqref{eq:graph-to-mixed-integrability}.
\end{proof}

\begin{lemma}[Weighted velocity control]
	\label{lem:base-slot-control}
	For every $\ell\ge0$, every $1\le r\le\infty$ and every measurable $u$ on $\cD\times\bR_v^3$,
	\begin{equation}
		\label{eq:base-slot-pointwise}
		\|w^{\ell}u\|_{L^r_xH_\sigma^1}
		\le C_\ell\|w^{\ell}\nabla_vu\|_{L^r_xL^2_v}
		+C_\ell|\cD|^{1/r}\|w^{\ell+2}u\|_{L^\infty} .
	\end{equation}
	In particular, for every $0\le L\le A_\flat-2$,
	\begin{equation}
		\label{eq:base-slot-control}
		\|w^Lu\|_{L_x^{p_\flat}H_\sigma^1}
		\le C_L\cB_{\rm bs}(u).
	\end{equation}
\end{lemma}

\begin{proof}
	By \eqref{eq:H_sigma_norm}, \eqref{eq:coulomb-maxwellian-anisotropy} and \eqref{eq:coulomb-sigmavv-comparison}, we have, for almost every $x\in\cD$,
	\[
	\|w^\ell u(x,\cdot)\|_{H_\sigma^1}
	\le C_\ell\|w^\ell\nabla_vu(x,\cdot)\|_{L^2_v}
	+C_\ell\|w^{\ell-1/2}u(x,\cdot)\|_{L^2_v} ,
	\]
	and, by H\"older's inequality,
	\[
	\|w^{\ell-1/2}u(x,\cdot)\|_{L^2_v}\le\|w^{-5/2}\|_{L^2_v}\|w^{\ell+2}u(x,\cdot)\|_{L^\infty_v}.
	\]
	Taking the $L^r_x$ norm over $\cD$, we obtain
	\eqref{eq:base-slot-pointwise}.  For \eqref{eq:base-slot-control}, apply \eqref{eq:base-slot-pointwise} at $\ell=L$ and
	$r=p_\flat$; since $L+2\le A_\flat$, \eqref{eq:base-strong} gives
	\[
	C_L\|w^L\nabla_vu\|_{L_x^{p_\flat}L_v^2}+C_L|\cD|^{1/p_\flat}\|w^{L+2}u\|_{L^\infty}\le C_L\cB_{\rm bs}(u).
	\]
	Combining this with \eqref{eq:base-slot-pointwise}, we obtain \eqref{eq:base-slot-control}.
\end{proof}

\begin{lemma}[Graph and mixed rows]
	\label{lem:mixed-rows-dominated}
	Every measurable $u$ on $\cD\times\bR_v^3$ satisfies, for
	$1\le m\le\mathrm N$,
	\[
	R_m(u)\le C\,\widetilde{\mathfrak Q}_{\mathrm N}(u).
	\]
\end{lemma}

\begin{proof}
	For $1\le m\le\mathrm N$, put
	$k:=\ell_m-(\mathrm N-m)$, so that $k\ge0$: at $m=\mathrm N$, $k=\ell_\sharp$, and for $m\le\mathrm N-1$, by \eqref{eq:finite-hierarchy-strict-ladder}, $\ell_m\ge\ell_{\mathrm N-1}\ge A_\diamond>\mathrm N$. Then
	\eqref{eq:graph-to-mixed-integrability} at that $k$, applied to
	$w^{\ell_m}\mathfrak D_Iu$, reads
	\[
	\|w^{\,\mathrm N-m-2}\mathfrak D_Iu\|_{L_x^{q_\flat}L_v^2}
	\le C_k\|w^{\ell_m}\mathfrak D_Iu\|_{\cH_{\kin,\sigma}},
	\]
	and summing over $I\in\mathfrak A^m$ gives $R_m(u)\le CW_m(u)\le C\widetilde{\mathfrak Q}_{\mathrm N}(u)$.
\end{proof}

The reduced norm also omits $U_m$, since weighted graph rows dominate them: for $0\le m\le\mathrm N$ and every measurable $u$
on $\cD\times\bR_v^3$, Lemma~\ref{lem:compact-velocity-support}\,(1) on
$\cH_{\kin,\sigma}(\cD\times\bR_v^3)$ with the $C^1_v$ multiplier $w^{-\ell_m}$, applied to
$w^{\ell_m}\mathfrak D_Iu$ for $I\in\mathfrak A^m$, gives
\begin{equation}
	\label{eq:unweighted-row-dominated}
	U_m(u)\le CW_m(u).
\end{equation}

By \eqref{eq:order-two-Rm} and \eqref{eq:GN-lower-order-gaussian}, $\cG^{\sharp,-2}_{\mathrm N,q_\flat}(u)$ is the sum of $\|w^{\mathrm N-2}u\|_{L_x^{q_\flat}L_v^2}$ and $R_1(u),\ldots,R_{\mathrm N}(u)$.  Since $A_\flat\ge\mathrm N$ by \eqref{eq:finite-hierarchy-weight-choice} and \eqref{eq:finite-hierarchy-strict-ladder}, the first term satisfies
\[
\|w^{\mathrm N-2}u\|_{L_x^{q_\flat}L_v^2}
\le|\cD|^{1/q_\flat}\|w^{\mathrm N-2-A_\flat}\|_{L_v^2}
\|w^{A_\flat}u\|_{L^\infty}
\le C\cB_{\rm bs}(u).
\]
Hence
\begin{equation}
	\label{eq:reduced-Q-gaussian-density}
	\cG^{\sharp,-2}_{\mathrm N,q_\flat}(u)
	\le C\Bigl(\cB_{\rm bs}(u)+\sum_{m=1}^{\mathrm N}R_m(u)\Bigr)
	\le C\,\widetilde{\mathfrak Q}_{\mathrm N}(u),
\end{equation}
where the last inequality follows from Lemma~\ref{lem:mixed-rows-dominated}.

The reduced norm also controls the coefficient bounds required for the spectral problem in \S\ref{ss:diffuse-hierarchy-block}.

\begin{lemma}[Coefficient bounds from the reduced norm]
	\label{lem:coefficient-transmission}
	There is a constant $C$, depending only on the structural data, such that every measurable $u$ on $\cD\times\bR_v^3$ satisfies
	\begin{equation}
		\label{eq:coefficient-transmission}
		\|\mu^{1/4}u\|_{L^\infty}
		+\|\mu^{1/4}u\|_{C^{\alpha_{\rm bs}}_{\kin}}
		\le C\cB_{\rm bs}(u).
	\end{equation}
	Consequently there is $\eta_1>0$, depending only on the structural data, such that
	every measurable $g$ on $\cD\times\bR_v^3$ with $\widetilde{\mathfrak Q}_{\mathrm N}(g)\le\eta_1$ satisfies
	\begin{equation}
		\label{eq:construction-coefficient-regularity}
		\|\mu^{1/4}g\|_{L^\infty}
		\le\min\{\varepsilon_{\rm bs},\varepsilon_{\rm dif}\},
		\qquad
		\|\mu^{1/4}g\|_{C^{\alpha_{\rm bs}}_{\kin}}\le1,
	\end{equation}
	where $\varepsilon_{\rm bs}$ is the coefficient threshold of Proposition~\ref{prop:base-regularity} and $\varepsilon_{\rm dif}$ is the constant of \eqref{eq:diffuse-graph-smallness}.
\end{lemma}

\begin{proof}
	Put
	$M(v):=\mu^{1/4}(v)w^{-A_\flat}(v)$, so that $\mu^{1/4}u=Mw^{A_\flat}u$.  Since $M$ is a fixed smooth
	Gaussian--polynomial multiplier,
	$\|M\|_{L^\infty}+[M]_{C^{\alpha_{\rm bs}}_{\kin}}<\infty$.  By \eqref{eq:base-strong} and the Leibniz bound for the kinetic H\"older seminorm, this gives \eqref{eq:coefficient-transmission}.  The second assertion follows from \eqref{eq:coefficient-transmission} and \eqref{eq:reduced-Q3} by taking $\eta_1$ so small that $C\eta_1\le\min\{\varepsilon_{\rm bs},\varepsilon_{\rm dif},1\}$.
\end{proof}

\begin{lemma}[Interpolation bounds for the background]
	\label{lem:background-finite-p}
	Let $\eta>0$ and let $g$ be measurable on $\cD\times\bR_v^3$ with $\widetilde{\mathfrak Q}_{\mathrm N}(g)\le\eta$.  Then, at the rungs
	\eqref{eq:dyn-own-order-exponents},
	\begin{equation}
		\label{eq:Dg-rungs}
		\sum_{k=1}^{\mathrm N-1}\ \sum_{|I|=k}
		\|w^{-2}\mathfrak D_Ig\|_{L_x^{\mathsf g(k)}L_v^2}
		+\sum_{|I|=\mathrm N}\|w^{-2}\mathfrak D_Ig\|_{L_x^{q_\flat}L_v^2}
		\le C\eta ,
	\end{equation}
	where the sums run over ordered complete words and $C$ depends only on the structural data.
\end{lemma}

\begin{proof}
	The last sum is $R_{\mathrm N}(g)$, because
	the weight $w^{\mathrm N-m-2}$ of \eqref{eq:order-two-Rm} equals $w^{-2}$ at $m=\mathrm N$. By Lemma~\ref{lem:mixed-rows-dominated},
	$R_{\mathrm N}(g)\le C\widetilde{\mathfrak Q}_{\mathrm N}(g)\le C\eta$.
	For $1\le k\le\mathrm N-1$ we apply \eqref{eq:ordered-word-gaussian-rungs} to
	$g$ at $(\mathrm L_{\rm G},a_{\rm G})=(-2,q_\flat)$, $M=\mathrm N$ and $p_{\rm G}=\infty$.  Its hypothesis
	\eqref{eq:pointwise-slice-chain-hypothesis} holds by \eqref{eq:ordered-complete-derivatives}, since the words act on $g$ in the sense of distributions.  Its hypothesis
	\eqref{eq:dynamic-pointwise-finite-norms}, at the base exponents
	$(p_{\rm G},p_{\rm c})=(\infty,p_\flat)$, holds: its four norms are bounded as follows.
	\begin{itemize}
		\item $\|w^{-2}g\|_{L_x^\infty L_v^2}\le C\cB_{\rm bs}(g)$, by \eqref{eq:base-fiber-Lsharp}.
		\item $\|g\|_{L_x^{p_\flat}H_\sigma^1}\le C\cB_{\rm bs}(g)$, by
		\eqref{eq:base-slot-control} at $L=0$.
		\item $\cG^{\sharp,-2}_{\mathrm N,q_\flat}(g)\le C\widetilde{\mathfrak Q}_{\mathrm N}(g)$,
		by \eqref{eq:reduced-Q-gaussian-density}.
		\item $\cG^{0}_{\mathrm N,2}(g)\le C\widetilde{\mathfrak Q}_{\mathrm N}(g)$, by \eqref{eq:real-weight-Hsigma-monotonicity} and \eqref{innerprod_st_sigma}, since
		$\ell_m\ge\mathrm N-m$ for $0\le m\le\mathrm N$ by
		\eqref{eq:finite-hierarchy-weight-choice} and \eqref{eq:finite-hierarchy-strict-ladder}.
	\end{itemize}
	By these four bounds and \eqref{eq:reduced-Q3},
	\[
	\|w^{-2}g\|_{L_x^\infty L_v^2}+\|g\|_{L_x^{p_\flat}H_\sigma^1}+\cG^{\sharp,-2}_{\mathrm N,q_\flat}(g)+\cG^{0}_{\mathrm N,2}(g)\le C\eta .
	\]
	Since $\mathsf e(\mathrm N/k;q_\flat,\infty)=\mathsf g(k)$ by \eqref{eq:dyn-own-order-exponents}, every ordered complete word $I$ with $1\le|I|=k\le\mathrm N-1$ satisfies, by \eqref{eq:ordered-word-gaussian-rungs} and the inequality above,
	\[
	\|w^{-2}\mathfrak D_Ig\|_{L_x^{\mathsf g(k)}L_v^2}
	\le C\|w^{-2}g\|_{L_x^\infty L_v^2}^{1-k/\mathrm N}\,\cG^{\sharp,-2}_{\mathrm N,q_\flat}(g)^{k/\mathrm N}+C\|w^{-2}g\|_{L_x^\infty L_v^2}
	\le C\eta .
	\]
	Summing over the words $I$ with $1\le|I|=k\le\mathrm N-1$, we obtain
	\[
	\sum_{k=1}^{\mathrm N-1}\sum_{|I|=k}\|w^{-2}\mathfrak D_Ig\|_{L_x^{\mathsf g(k)}L_v^2}\le C\eta.
	\]
	Combining this with $R_{\mathrm N}(g)\le C\eta$, we obtain \eqref{eq:Dg-rungs}.
\end{proof}

\subsection{The differentiated system and the spectral problem}\label{ss:diffuse-hierarchy-block}

The ordered interior identities extend to $\cH_{\kin,\sigma}$ jet chains and yield the differentiated equations and mass hierarchy for the fixed-cutoff problem.  A \emph{jet chain} in $\cH_{\kin,\sigma}$ is a family $(q_J)_{|J|\le m}\subset\cH_{\kin,\sigma}$, indexed by ordered words and satisfying $\mathfrak D_\nu q_J=q_{\nu J}$ in distributions for $|\nu J|\le m$. The same definition applies in $\cH_{\kin,\sigma}(\cD_T\times\bR_v^3)$.  The differentiated equations below are used to realize the spectral jets in Proposition~\ref{prop:fixedN-order-m-realization} and the evolutionary jets in Proposition~\ref{prop:dyn-differentiated-equation}. The mass identities enter the stationary and evolutionary hierarchy estimates in Propositions~\ref{prop:uniform-hierarchy} and~\ref{prop:dyn-hierarchy}.

\begin{prop}[Ordered interior algebra]
	\label{prop:graph-ordered-interior-algebra}
	Fix $1\le m\le\mathrm N$.  Let $q_I\in\cH_{\kin,\sigma}(\cD\times\bR^3_v)$ for
	ordered words $|I|\le m$, with $q_\varnothing=f$, and let
	$g\in\cD'(\cD\times\bR^3)$.  Assume
	\[
	\mathfrak D_\nu q_J=q_{\nu J}
	\quad\text{in }\cD'(\cD\times\bR^3),
	\qquad |\nu J|\le m,
	\]
	\[
	\sum_{|J|<m}\|w^{A_\diamond}q_J\|_{L_x^2H_\sigma^1}<\infty,
	\qquad
	\mu^{1/4}\mathfrak D_Jg\in L^2(\cD\times\bR_v^3),
	\quad |J|\le m.
	\]
	Then \eqref{eq:ordered-transport-commutator},
	\eqref{eq:local-L0-complete-expansion},
	\eqref{eq:compact-word-commutator-expansion}, and
	\eqref{eq:collision-remainder-expansion}, with every $\mathfrak D_Jf$
	replaced by $q_J$, hold in $\cD'(\cD\times\bR^3)$ for $1\le|I|\le m$.  The right sides of \eqref{eq:ordered-transport-commutator}, \eqref{eq:local-L0-complete-expansion} and \eqref{eq:compact-word-commutator-expansion} belong to $L_x^2H_\sigma^{-1}$, and the terms $N[\mathfrak D_{I_S}g,q_{I_{S^c}}]$ and $\mathfrak R_I[g,f]$ of \eqref{eq:collision-remainder-expansion} to $L_x^1H_\sigma^{-1}$.
\end{prop}

\begin{proof}
	The transport formula \eqref{eq:ordered-transport-commutator} follows from the jet-chain identity; Proposition~\ref{prop:ordered-transport-weight} gives its $L_x^2H_\sigma^{-1}$ bound.
	
	To justify the collision identity, let $\chi_R(v)=\chi(v/R)$, $R\ge1$, with $\chi$ as in \eqref{eq:radial-cutoff}.  For every nonempty ordered word $J$, $|J|\le\mathrm N$,
	\begin{equation}
		\label{eq:radial-cutoff-complete-word}
		|\mathfrak D_J\chi_R|
		+R|\nabla_v\mathfrak D_J\chi_R|
		\lesssim_J\mathbf1_{\{R\le|v|\le2R\}}.
	\end{equation}
	Moreover,
	\begin{equation}
		\label{eq:radial-cutoff-coercive}
		(\nabla_v\chi_R)^T\sigma(v)\nabla_v\chi_R
		\lesssim R^{-4}w(v)^{-1}\mathbf1_{\{R\le|v|\le2R\}}.
	\end{equation}
	On the shell $\{R\le|v|\le2R\}$, $|\nabla_v^k\chi_R|\lesssim_k R^{-k}$ and $|(\nabla\tau_\nu)v|\lesssim R$. Since the velocity part of $\mathfrak D_\nu$ is $((\nabla\tau_\nu)v)\cdot\nabla_v$, induction gives \eqref{eq:radial-cutoff-complete-word}.  Since
	$\nabla_v\chi_R$ is radial, \eqref{eq:coulomb-sigmavv-comparison} and
	$w\simeq R$ on the shell give \eqref{eq:radial-cutoff-coercive}.
	
	Insert $\chi_R$ in both slots of $N[\cdot,\cdot]$. The following shell estimates permit passage to $R\to\infty$.  If a complete derivative
	falls on the cutoff in the first slot, then
	\eqref{eq:radial-cutoff-complete-word} and absolute continuity give
	\[
	\|\mathbf1_{\{R\le|v|\le2R\}}
	\mu^{1/4}\mathfrak D_Jg\|_{L_{x,v}^2}\longrightarrow0.
	\]
	For a second-slot jet $q_J$ with $|J|<|I|$, write $h=w^{A_\diamond}q_J$.  By
	\eqref{eq:radial-cutoff-complete-word}
	and \eqref{eq:coulomb-eigen-decomposition}--\eqref{eq:eigenvalue-size},
	\[
	\|(\mathfrak D_{K}\chi_R)h\|_{H_\sigma^1}^2
	\lesssim_{K}
	\int_{R\le|v|\le2R}
	\bigl(\sigma_{ij}\partial_i h\,\partial_j h
	+w^{-1}|h|^2\bigr)\,\dd v,
	\qquad 1\le|K|\le\mathrm N.
	\]
	Its integral in $x$ tends to zero by absolute continuity. Indeed, \eqref{eq:auxiliary-coercive-equivalence} and $h\in L_x^2H_\sigma^1$ yield
	\[
	\int_{\cD}\int_{\bR^3}\bigl(\sigma_{ij}\partial_ih\,\partial_jh+w^{-1}|h|^2\bigr)\,\dd v\,\dd x
	\le C\|h\|_{L_x^2H_\sigma^1}^2<\infty.
	\]
	For the current jet,
	\eqref{eq:radial-cutoff-convergence} gives $\chi_Rq_I(x)\to q_I(x)$ in
	$H_\sigma^1$ for almost every $x$, and \eqref{eq:radial-cutoff-coercive} gives
	$\|\chi_Rq_I(x)\|_{H_\sigma^1}\lesssim\|q_I(x)\|_{H_\sigma^1}$ uniformly in
	$R\ge1$.  Dominated convergence in $x$ gives
	$\chi_Rq_I\to q_I$ in $L_x^2H_\sigma^1$.
	
	Fix $R$ and $U\Subset\cD$. Test \eqref{eq:complete-conjugated-kernel-deformation} against functions supported in $U\times\bR^3$. Both input slots are then evaluated only on the compact set $\overline U\times\operatorname{supp}\chi_R$.  Since $\operatorname{supp}\chi_R\subset\{|v|\le2R\}$, \eqref{eq:local-H1-H_sigma} gives, for every $u\in L_x^2H_\sigma^1$,
	\[
	\|\mathbf1_{\{|v|\le2R\}}u\|_{L^2_{x,v}}+\|\mathbf1_{\{|v|\le2R\}}\nabla_vu\|_{L^2_{x,v}}
	\le C(1+R)^{3/2}\|u\|_{L_x^2H_\sigma^1}.
	\]  The hypothesis
	$\mu^{1/4}\mathfrak D_Jg\in L^2$, $|J|\le m$, gives the local $L^2$
	hypotheses of Lemma~\ref{lem:friedrichs-complete-field} for the first slot
	and its next complete derivative; $\mathfrak D_\nu q_J=q_{\nu J}$ and
	$q_J\in L_x^2H_\sigma^1$ give the local $H_v^1$ hypotheses of that lemma for
	the second slot.
	Lemma~\ref{lem:friedrichs-complete-field} gives $L^2$ regularization of the first slot and $H_v^1$ regularization of the second. Apply \eqref{eq:complete-conjugated-kernel-deformation} to these approximations.
	As the regularization is removed, the Friedrichs commutators vanish by \eqref{eq:friedrichs-L2}, while
	Lemma~\ref{lem:velalg} and Proposition~\ref{prop:deformation-algebra} give the convergence of the collision terms.  Sending $R\to\infty$ by the shell
	estimates gives the one-letter identity \eqref{eq:complete-conjugated-kernel-deformation}.
	
	Iterating the one-letter identity gives the differentiated equation once the nonlinear terms are identified in $L_x^1H_\sigma^{-1}$ and the fixed linear part is handled.  The recurrence \eqref{eq:deformation-remainder-recurrence} gives \eqref{eq:collision-remainder-expansion} distributionally by ordered iteration.  The term $N[g,q_I]$, with $S=\varnothing$, needs only the unweighted graph norm of $q_I$; every other term has in its second slot a jet $q_J$ with $|J|<|I|$, and $\|w^{A_\diamond}q_J\|_{L_x^2H_\sigma^1}<\infty$ by the hypothesis
	$\sum_{|J|<m}\|w^{A_\diamond}q_J\|_{L_x^2H_\sigma^1}<\infty$.
	Lemma~\ref{lem:velalg} and Proposition~\ref{prop:deformation-algebra},
	followed by Cauchy--Schwarz in $x$, bound all nonlinear terms in
	$L_x^1H_\sigma^{-1}$.
	
	For $L$, use the decomposition $L=\cA_\mu+\cK_0$ preceding \eqref{eq:compact-integral-is-inflow-K}.  The local identity
	\eqref{eq:Amu-one-complete-exact} consists of collision-type terms of the form treated above, with first slot $\sqrt\mu$ or $\mathfrak p_X\sqrt\mu$ in place of a jet of $g$.  The preceding cutoff and Friedrichs passage,
	followed by ordered induction, gives
	\eqref{eq:local-L0-complete-expansion} distributionally, and \eqref{eq:local-L0-complete-bound} at $\ell=0$, \eqref{eq:real-weight-Hsigma-monotonicity} and the hypothesis on $q_J$ give, since $2|I|+2\le A_\diamond$ by \eqref{eq:Adef},
	\[
	\Bigl\|\sum_{|J|<|I|}\sum_{\lambda}c_{I,J,\lambda}\,p_{I,J,\lambda}N_{\mathsf K_{I,J,\lambda}}[m_{I,J,\lambda},p'_{I,J,\lambda}q_J]\Bigr\|_{L_x^2H_\sigma^{-1}}
	\le C\sum_{|J|<|I|}\|w^{2|I|+2}q_J\|_{L_x^2H_\sigma^1}<\infty .
	\]
	For $\cK_0$, \eqref{eq:compact-integral-is-inflow-K}, \eqref{eq:iterated-compact-kernel-recursion}, and \eqref{eq:iterated-compact-full-operator-bound} justify the $L^2$ Friedrichs passage locally in $\cD\times\bR^3$.  The
	Gaussian factors of $m_0$ and $k_0$ in Lemma~\ref{lem:compact-maxwellian-complete-commutators} remove the velocity cutoff on the current jet $q_I$.  The chain $\mathfrak D_\nu q_J=q_{\nu J}$ identifies each derivative of a jet with a jet of the next length, and ordered iteration gives
	\eqref{eq:compact-word-commutator-expansion}.  By \eqref{eq:iterated-compact-full-operator-bound}, the right side of
	\eqref{eq:compact-word-commutator-expansion} belongs to $L_x^2H_\sigma^{-1}$.
\end{proof}

\begin{thm}[Weak differentiated equations through complete order $\mathrm N$]
	\label{thm:weak-differentiated-system}
	Fix $1\le m\le\mathrm N$ and let $g$ be measurable on
	$\cD\times\bR_v^3$, its ordered derivatives $\mathfrak D_Jg$ being taken in
	$\cD'(\cD\times\bR^3)$.  For each ordered word $I$ with $|I|\le m$, let
	$f_I\in\cH_{\kin,\sigma}$, with $f_\varnothing=f$, and assume
	\begin{equation}
		\label{eq:distributional-jet-chain}
		\mathfrak D_\nu f_J=f_{\nu J}
		\quad\text{in }\cD'(\cD\times\bR^3),
		\qquad \nu\in\mathfrak A,
		\qquad |\nu J|\le m.
	\end{equation}
	Assume further that
	\begin{equation}
		\label{eq:lower-weight-assumption}
		\sum_{|J|<m}\|w^{A_\diamond}f_J\|_{L_x^2H_\sigma^1}<\infty
	\end{equation}
	and
	\begin{equation}
		\label{eq:background-jet-assumption}
		\mu^{1/4}\mathfrak D_Jg\in L^2(\cD\times\bR_v^3),
		\qquad |J|\le m.
	\end{equation}
	Then, for every nonempty $I$ with $|I|\le m$,
	\begin{align}
		[\mathfrak D_I,v\cdot\nabla_x]f
		&=\sum_{|J|<|I|}\cV_{I,J}f_J,
		\label{eq:weak-transport-word}\\
		\mathfrak D_I(Lf)
		&=Lf_I+[\mathfrak D_I,L]f,
		\label{eq:weak-L0-word}\\
		\mathfrak D_I N[g,f]
		&=\sum_{S\subset[|I|]}
		N[\mathfrak D_{I_S}g,f_{I_{S^c}}]
		+\mathfrak R_I[g,f]
		\label{eq:weak-collision-word}
	\end{align}
	in distributions, where the local and compact parts of
	$[\mathfrak D_I,L]f$ are \eqref{eq:local-L0-complete-expansion} and
	\eqref{eq:compact-word-commutator-expansion}, with each lower
	$\mathfrak D_Jf$ replaced by $f_J$; in $\mathfrak R_I[g,f]$, whose terms have the form \eqref{eq:deformation-remainder-structure}, $\mathfrak D_{J_2}f$ is replaced by $f_{J_2}$.
	Moreover,
	\begin{align}
		[\mathfrak D_I,(v\cdot\nabla_x+L_g)]f
		={}&\sum_{|J|<|I|}\cV_{I,J}f_J
		+[\mathfrak D_I,L]f\notag\\
		&-\sum_{\varnothing\ne S\subset[|I|]}
		N[\mathfrak D_{I_S}g,f_{I_{S^c}}]
		-\mathfrak R_I[g,f].
		\label{eq:weak-triangular-commutator}
	\end{align}
	The right side of \eqref{eq:weak-transport-word} and the two parts
	of $[\mathfrak D_I,L]f$ in \eqref{eq:weak-L0-word} belong to
	$L_x^2H_\sigma^{-1}$; the terms $N[\mathfrak D_{I_S}g,f_{I_{S^c}}]$ and
	$\mathfrak R_I[g,f]$ belong to $L_x^1H_\sigma^{-1}$.  If $v\cdot\nabla_x f+L_gf=0$ in distributions, then
	\begin{equation}
		\label{eq:weak-differentiated-equation}
		v\cdot\nabla_x f_I+L_gf_I=-[\mathfrak D_I,(v\cdot\nabla_x+L_g)]f
		\qquad\text{in }\cD'(\cD\times\bR^3).
	\end{equation}
	
\end{thm}

\begin{proof}
	Proposition~\ref{prop:graph-ordered-interior-algebra}, applied to $q_J=f_J$ under \eqref{eq:distributional-jet-chain}--\eqref{eq:background-jet-assumption}, yields \eqref{eq:weak-transport-word}--\eqref{eq:weak-collision-word} with the asserted source spaces.  Combining \eqref{def:L_g}, \eqref{eq:weak-transport-word} and
	\eqref{eq:weak-collision-word}, we obtain \eqref{eq:weak-triangular-commutator}, since the summand $N[g,f_I]$ with $S=\varnothing$ in \eqref{eq:weak-collision-word} cancels the term $-N[g,f_I]$ of $L_gf_I$.  Applying $\mathfrak D_I$ to
	$v\cdot\nabla_x f+L_gf=0$ in distributions, we obtain
	\eqref{eq:weak-differentiated-equation}.
\end{proof}

The bulk mass of a complete derivative is determined by the preceding jet.
For \(\nu\in\mathfrak A\), write $T_s^\nu$ for the complete pullback
\eqref{eq:bulk-flow-pullback} along $\tau_\nu$.  For
\(f\) with \(w^{-1/2}f\in L^2(\cD\times\bR^3)\), define the \emph{mass selectors} \(\cM_\nu\),
\(\nu\in\mathfrak A\),
\begin{equation}
	\label{eq:first-mass-selector}
	\cM_\nu[f]:=-\int_{\cD}\int_{\bR^3}
	f\sqrt\mu\big(2\operatorname{div}\tau_\nu
	+((\nabla \tau_\nu)v)\cdot\nabla_v\log\sqrt\mu\big)\,\dd v\dd x.
\end{equation}

\begin{prop}[Complete-derivative bulk mass hierarchy]
	\label{prop:complete-derivative-mass-hierarchy}
	Let $f\in\cD'(\cD\times\bR^3)$, let $I=(\nu_1,\ldots,\nu_m)\in\mathfrak A^m$ with $1\le m\le\mathrm N$, and write
	$I_\ast:=(\nu_2,\ldots,\nu_m)$.
	On the smooth core,
	\begin{equation}
		\label{eq:complete-derivative-mass-identity}
		\iint_{\cD\times\bR^3}\mathfrak D_If\sqrt\mu\dd v\dd x
		=\cM_{\nu_1}[\mathfrak D_{I_\ast}f].
	\end{equation}
	Identity \eqref{eq:complete-derivative-mass-identity}
	also holds for $f$ outside the smooth core provided
	$w^{-1/2}\mathfrak D_{I_\ast}f\in L^2(\cD\times\bR_v^3)$ and there is a sequence
	$s_j\to0$ such that
	\begin{equation}
		\label{eq:complete-derivative-mass-formula}
		\frac{T_{s_j}^{\nu_1}\mathfrak D_{I_\ast}f-\mathfrak D_{I_\ast}f}{s_j}
		\rightharpoonup \mathfrak D_If
		\quad\text{weakly in }L^2(\cD\times\bR_v^3,w^{2\ell}\dd x\dd v)
		\text{ for some }\ell\in\bR.
	\end{equation}
	Condition \eqref{eq:complete-derivative-mass-formula} holds with $\ell=0$
	whenever $\mathfrak D_{I_\ast}f,\mathfrak D_If\in L^2(\cD\times\bR_v^3)$; by \eqref{260902137}, this covers every jet chain of $\cH_{\kin,\sigma}$ and a.e.\ time
	slice of a jet chain of $\cH_{\kin,\sigma}(\cD_T\times\bR_v^3)$.
	In this case,
	\begin{equation}
		\label{eq:complete-derivative-mass-sharp}
		\left|\iint\mathfrak D_If\sqrt\mu\right|
		\le C_I\|\mathfrak D_{I_\ast}f\|_{L^2(\cD\times\bR_v^3)}.
	\end{equation}
\end{prop}

\begin{proof}
	Condition \eqref{eq:complete-derivative-mass-formula} holds with $\ell=0$ under the stated $L^2$ hypothesis by \eqref{eq:strong-complete-flow-generator-limit} of
	Lemma~\ref{lem:strong-complete-flow-generator-L2}.
	
	On the smooth core, set
	\begin{equation}
		\label{eq:mass-selector-density}
		B_\nu(x,v):=
		2\operatorname{div}\tau_\nu(x)
		+((\nabla\tau_\nu(x))v)\cdot\nabla_v\log\sqrt\mu(v),
	\end{equation}
	so that $\cM_\nu[f]=-\iint f\sqrt\mu\,B_\nu$ by \eqref{eq:first-mass-selector}. By \eqref{eq:phase-divergence} and Proposition~\ref{prop:complete-fields}\,(\ref{it:cf-tangential}), we have $\operatorname{div}_{x,v}\bigl(\tau_\nu,(\nabla\tau_\nu)v\bigr)=2\operatorname{div}\tau_\nu$ and $\tau_\nu\cdot n=0$ on $\partial\cD$. Since $\sqrt\mu$ depends only on $v$, the first of these identities and \eqref{eq:mass-selector-density} give
	\[
	\operatorname{div}_{x,v}\bigl(\sqrt\mu\,(\tau_\nu,(\nabla\tau_\nu)v)\bigr)
	=\sqrt\mu\,\bigl(2\operatorname{div}\tau_\nu+((\nabla\tau_\nu)v)\cdot\nabla_v\log\sqrt\mu\bigr)
	=\sqrt\mu\,B_\nu .
	\]
	Integrate by parts, using this identity and $\tau_\nu\cdot n=0$ on $\partial\cD$, to obtain
	\[
	\iint \mathfrak D_{\nu_1}\mathfrak D_{I_\ast}f\sqrt\mu
	=-\iint \mathfrak D_{I_\ast}f\sqrt\mu\,B_{\nu_1}
	=\cM_{\nu_1}[\mathfrak D_{I_\ast}f],
	\]
	which is \eqref{eq:complete-derivative-mass-identity}.
	
	For the nonsmooth case, assume \eqref{eq:complete-derivative-mass-formula} and define
	\begin{equation}
		\label{eq:mass-flow-adjoint-weight}
		\Psi_s^{\nu_1}(y,z)
		:=|J_{-s}^{\nu_1}(y)|^2
		\sqrt{\mu(A_{-s}^{\nu_1}(y)z)}.
	\end{equation}
	By \eqref{eq:bulk-flow-adjoint-quotient},
	$\Psi_s^{\nu_1}=(T_s^{\nu_1})^\dagger\sqrt\mu$, and hence
	\begin{equation}
		\label{eq:mass-flow-pairing}
		\iint \frac{T_s^{\nu_1}\mathfrak D_{I_\ast}f-\mathfrak D_{I_\ast}f}{s}\sqrt\mu
		=\iint \mathfrak D_{I_\ast}f(y,z)\,
		\frac{\Psi_s^{\nu_1}(y,z)-\Psi_0^{\nu_1}(y,z)}{s}\,\dd z\dd y.
	\end{equation}
	Since $A_{-s}^{\nu_1}=D\Phi_{-s}^{\nu_1}$ and $J_{-s}^{\nu_1}=\det A_{-s}^{\nu_1}$, the variational equation and Liouville's formula give
	\begin{equation}
		\label{eq:mass-flow-ode}
		\partial_sA_{-s}^{\nu_1}(y)
		=-(\nabla\tau_{\nu_1})(\Phi_{-s}^{\nu_1}(y))A_{-s}^{\nu_1}(y),
		\qquad
		\partial_sJ_{-s}^{\nu_1}(y)
		=-(\operatorname{div}\tau_{\nu_1})(\Phi_{-s}^{\nu_1}(y))J_{-s}^{\nu_1}(y).
	\end{equation}
	At $s=0$ we have $A_0^{\nu_1}=\mathbf I$ and $J_0^{\nu_1}=1$; differentiating \eqref{eq:mass-flow-adjoint-weight} using \eqref{eq:mass-flow-ode} gives
	\[
	\left.\partial_s\Psi_s^{\nu_1}(y,z)\right|_{s=0}
	=-2(\operatorname{div}\tau_{\nu_1})(y)\sqrt{\mu(z)}
	-(\nabla_v\sqrt\mu)(z)\cdot(\nabla\tau_{\nu_1}(y))z .
	\]
	By this identity and \eqref{eq:mass-selector-density}, and by \eqref{eq:mass-flow-adjoint-weight},
	\eqref{eq:mass-flow-ode} and \eqref{eq:complete-flow-uniform-bounds} (with its $s_0$), there is $c>0$ such that
	\begin{equation}
		\label{eq:mass-flow-weight-derivative}
		\left.\partial_s\Psi_s^{\nu_1}\right|_{s=0}
		=-\sqrt\mu\,B_{\nu_1},
		\qquad
		\sup_{|s|\le s_0}|\partial_s\Psi_s^{\nu_1}(y,z)|
		\lesssim_{\nu_1}(1+|z|^2)e^{-c|z|^2}.
	\end{equation}
	The Gaussian bound in \eqref{eq:mass-flow-weight-derivative} and dominated convergence yield
	\begin{equation}
		\label{eq:mass-flow-weight-limit}
		\frac{\Psi_s^{\nu_1}-\Psi_0^{\nu_1}}{s}
		\longrightarrow-\sqrt\mu\,B_{\nu_1}
		\qquad\text{strongly in }w^{a}L^2_{y,z}
		\ \text{ for every finite }a,
	\end{equation}
	since $w^{|a|}(1+|z|^2)e^{-c|z|^2}\in L^2_{y,z}$.  Put $a=\tfrac12$.  Then
	$w^{-a}\mathfrak D_{I_\ast}f\in L_{x,v}^2$ by hypothesis.
	Since $\sqrt\mu\,w^{-\ell}\in L^2_{x,v}$, \eqref{eq:complete-derivative-mass-formula} may be tested against $\sqrt\mu$. Letting $s=s_j\to0$ in \eqref{eq:mass-flow-pairing}, with \eqref{eq:complete-derivative-mass-formula} on the left and \eqref{eq:mass-flow-weight-limit} on the right, we obtain
	$$
	\iint\mathfrak D_If\sqrt\mu\dd v\dd x=-\iint\mathfrak D_{I_\ast}f\sqrt\mu\,B_{\nu_1}\dd v\dd x,
	$$
	which is \eqref{eq:complete-derivative-mass-identity}.
	
	Since $\sqrt\mu\,B_{\nu_1}\in L^2_{x,v}$, Cauchy--Schwarz in \eqref{eq:complete-derivative-mass-identity} gives
	\eqref{eq:complete-derivative-mass-sharp}.
\end{proof}

We define the fixed-cutoff problem.  The first bound in \eqref{eq:construction-coefficient-regularity}, required by Proposition~\ref{prop:inhomogeneous-full-diffuse-graph} and Theorem~\ref{260322442}\,(2), follows from Lemma~\ref{lem:coefficient-transmission} when $\widetilde{\mathfrak Q}_{\mathrm N}(g)\le\eta_1$. From here to the end of \S\ref{app:fixedN-realization-induction}
we assume $\widetilde{\mathfrak Q}_{\mathrm N}(g)<\infty$, the first bound in \eqref{eq:construction-coefficient-regularity}, and
the condition $0<\delta\le\delta_{\rm dif}$ of \eqref{eq:diffuse-graph-smallness}.

Fix $N$ in the boundary cutoff $J_N$ of \eqref{eq:boundary-spectral-cutoff-definition}; it is distinct from the derivative order $\mathrm N$. The data $F=0$, $b=h_w$, and $m=0$ satisfy the compatibility condition of Corollary~\ref{cor:fixedN-spectral-solvability} for every $N$. Hence there is a unique
$f_N\in\cH_{\kin,\sigma}(\cD\times\bR_v^3)$ with
\begin{equation}
	\label{eq:spectral-problem}
	v\cdot\nabla_x f_N+L_gf_N=0,
	\qquad
	(f_N)_{\gamma^-}=P_{w,N}(f_N)_{\gamma^+}+h_w,
	\qquad
	\iint f_N\sqrt\mu=0.
\end{equation}
Define $\phi_N$ and $H_N$ by
\begin{equation}
	\label{eq:spectral-boundary-data}
	\phi_N:=\cA_0(f_N)_{\gamma^+},
	\qquad
	H_N:=P_{w,N}(f_N)_{\gamma^+}+h_w
	= (\mathfrak m_\circ+d_\delta)J_N\phi_N+h_w,
	\qquad
	(f_N)_{\gamma^-}=H_N,
\end{equation}
where the second expression for $H_N$ follows from \eqref{eq:boundary-spectral-diffuse-operator} and $\mathfrak m_{\Theta_w}=\mathfrak m_\circ+d_\delta$ (see \eqref{eq:limiting-wall-profiles}).

\begin{lemma}[Incoming diffuse data]
	\label{lem:incoming-datum-bound}
	For every finite $L\ge0$ there is $C_L$, depending only on the structural data and on $L$, such that every $\psi\in H^3(\partial\cD)$ satisfies
	\begin{equation}
		\label{eq:incoming-datum-bound}
		\bigl\|w^L(\mathfrak m_{\Theta_w}\psi+h_w)\bigr\|_{(C_x^{0,1}\cap C_v^{1,1})(\gamma^-)}
		+\bigl\|w^L(\mathfrak m_{\Theta_w}\psi+h_w)\bigr\|_{L^2(\gamma^-)}
		\le C_L\bigl(\|\psi\|_{H^3(\partial\cD)}+\delta\bigr).
	\end{equation}
\end{lemma}

\begin{proof}
	Sobolev embedding on the two-dimensional boundary \cite[Chapter~4, Proposition~3.3]{TaylorPDE1} gives 
	$$
	\|\psi\|_{C^{0,1}(\partial\cD)}\le C\|\psi\|_{H^3(\partial\cD)}.
	$$
	Put $H:=\mathfrak m_{\Theta_w}\psi+h_w$.  By \eqref{eq:wall-profile-chart-bound} with $K_x=1$ and $K_v=2$, on every boundary chart,
	\[
	\sup_{|\alpha|\le2,\ |\beta|\le1}
	\bigl(|w^L\partial_{y'}^\beta\partial_v^\alpha d_\delta|
	+|w^L\partial_{y'}^\beta\partial_v^\alpha h_w|\bigr)
	\le C_L\delta .
	\]  Since $\mathfrak m_\circ$ is a fixed Gaussian profile, $\mathfrak m_{\Theta_w}=\mathfrak m_\circ+d_\delta$ satisfies 
	$$
	\sup_{|\alpha|\le2,\ |\beta|\le1}|w^L\partial_{y'}^\beta\partial_v^\alpha\mathfrak m_{\Theta_w}|\le C_L
	$$ on every boundary chart.  Therefore we have
	\[
	\|w^LH\|_{(C_x^{0,1}\cap C_v^{1,1})(\gamma^-)}\le C_L\bigl(\|\psi\|_{C^{0,1}(\partial\cD)}+\delta\bigr).
	\]
	For the incoming $L^2$ norm,
	\[
	\|w^LH\|_{L^2(\gamma^-)}\le C\|w^{L+3}H\|_{L^\infty(\gamma^-)}
	\le C\|w^{L+3}H\|_{(C_x^{0,1}\cap C_v^{1,1})(\gamma^-)}
	\le C_L\bigl(\|\psi\|_{C^{0,1}(\partial\cD)}+\delta\bigr),
	\] where the last inequality is the preceding bound with $L+3$ in place of $L$.
\end{proof}

\subsection{Realization of the complete jets}\label{app:fixedN-realization-induction}

We first realize all ordered complete jets of the mass-normalized spectral solution $f_N$ at a fixed cutoff $N$. Constants in this subsection may depend on $N$, the weight, the background, and the wall perturbation; cutoff-independent bounds are proved in $\S$\ref{ss:diffuse-uniform-hierarchy}.  

\begin{lemma}[Complete jets at fixed spectral cutoff]
	\label{lem:fixedN-complete-jets}
	Fix $N$
	and let $f_N$ be the mass-normalized spectral solution
	\eqref{eq:spectral-problem}.  Then, for every ordered complete word $I$ with
	$|I|\le\mathrm N$,
	\begin{equation}
		\label{eq:fixedN-jet-fullgraph}
		q_{I,N}:=\mathfrak D_If_N\in\cH_{\kin,\sigma}.
	\end{equation}
	Its two flux traces are the completion-defined Green traces, and
	\[
	\|w^Lq_{I,N}\|_{\cH_{\kin,\sigma}}\le C
	\]
	at every finite weight $L$.  Moreover, at $f=f_N$ and with $f_I:=q_{I,N}$,
	\eqref{eq:distributional-jet-chain} and \eqref{eq:weak-differentiated-equation} hold for $|I|\le\mathrm N$.
\end{lemma}

We prove Lemma~\ref{lem:fixedN-complete-jets} by induction on the word length, using the assertions $\mathsf P_m$ of Definition~\ref{def:fixedN-package}:
Lemma~\ref{lem:fixedN-base-estimates} gives $\mathsf P_0$, and
Proposition~\ref{prop:fixedN-order-m-realization} gives $\mathsf P_m$ from
$\mathsf P_{m-1}$.  The requisite bounds on the incoming data $\mathfrak D_I^-H_N$ are established in Lemma~\ref{lem:fixedN-inflow-separation}.

\begin{lemma}[Incoming flow bounds]
	\label{lem:fixedN-inflow-separation}
	Fix $N$, let $f_N$ be the mass-normalized spectral solution \eqref{eq:spectral-problem},
	let $\phi_N,H_N$ be as in \eqref{eq:spectral-boundary-data}, and let $L\ge0$
	be finite.  Then, for every ordered word $I$ with $|I|\le\mathrm N$,
	\begin{equation}
		\label{eq:fixedN-incoming-word-bound}
		\|w^L\mathfrak D_I^-H_N\|_{L^2(\gamma^-)}
		<\infty .
	\end{equation}
	There is $h_0>0$, depending only on the structural data, such that the following holds.  If $m\ge1$ and
	$|I|+m\le\mathrm N$, then, for every $X=\tau_\xi$ with $\xi\in\mathfrak A$,
	\begin{align}
		\label{eq:fixedN-inflow-separation}
		&\sup_{0<|h|<h_0}
		\|(\Delta_h^{X,\gamma_-})^m\mathfrak D_I^-H_N\|_{L^2(\gamma^-)}
		<\infty,\\
		\label{eq:fixedN-incoming-difference-limit}
		&(\Delta_h^{X,\gamma_-})^m\mathfrak D_I^-H_N
		\longrightarrow
		(\mathfrak D_X^-)^m\mathfrak D_I^-H_N
		\quad\text{in }L^2(\gamma^-)
		\qquad(h\to0).
	\end{align}
\end{lemma}

\begin{proof}
	Set $a:=\mathfrak m_\circ+d_\delta$ and $b:=J_N\phi_N$. Then
	$H_N=ab+h_w$ by \eqref{eq:spectral-boundary-data}.  By
	\eqref{eq:boundary-word-Leibniz}, for every ordered word $K$ with
	$|K|\le\mathrm N$,
	\[
	\mathfrak D_K^-H_N
	=\sum_{K_1\sqcup K_2=K}\bigl(\mathfrak D_{K_1}^-a\bigr)\bigl(\tau_{K_2}b\bigr)
	+\mathfrak D_K^-h_w .
	\]
	Lemma~\ref{lem:wall-profile-derivatives} and
	\eqref{eq:wall-response-profiles} imply that $a$ satisfies
	\eqref{eq:boundary-gaussian-polynomial-envelope}, so that
	\eqref{eq:boundary-profile-flux-bound} gives
	\[
	\sup_{x\in\partial\cD}
	\int_{n_x\cdot v<0}w^{2L}|\mathfrak D_{K_1}^-a(x,v)|^2|n_x\cdot v|\,\dd v \le C.
	\]
	By \eqref{eq:fixedN-Hr} at $r=\mathrm N$ and
	\eqref{eq:outgoing-moment-flux-bound},
	\[
	\|\tau_{K_2}b\|_{L^2(\partial\cD)}
	\le C\|\phi_N\|_{L^2(\partial\cD)} \le C\|(f_N)_{\gamma^+}\|_{L^2(\gamma^+)},
	\]
	and \eqref{eq:wall-profile-complete-bound} gives
	$\|w^L\mathfrak D_K^-h_w\|_{L^2(\gamma^-)}\le C\delta$.  Since $b$ depends only on $x\in\partial\cD$, integrating first in $v$ and using the bound on $\mathfrak D_{K_1}^-a$, we have
	\[
	\|w^L(\mathfrak D_{K_1}^-a)(\tau_{K_2}b)\|_{L^2(\gamma^-)}^2
	\le C\|\tau_{K_2}b\|_{L^2(\partial\cD)}^2 .
	\]
	The bounds on $\tau_{K_2}b$ and $\mathfrak D_K^-h_w$, together with the Leibniz expansion at $K=I$, give \eqref{eq:fixedN-incoming-word-bound}.
	
	For the difference bounds, let $s_0$ be as in \eqref{eq:complete-flow-uniform-bounds}, and take $0<h_0\le s_0/\mathrm N$.
	By \eqref{eq:fixedN-Hr} at $r=\mathrm N+2$ and the two-dimensional Sobolev
	embedding, $b\in C^{\mathrm N,\alpha}(\partial\cD)$ for every $0<\alpha<1$.  Moreover,
	$a,h_w\in C^{\mathrm N}(\partial\cD\times\bR^3)$ by
	\eqref{eq:boundary-gaussian-polynomial-envelope} and
	\eqref{eq:wall-profile-chart-bound}, so that $H_N=ab+h_w\in C^{\mathrm N}(\partial\cD\times\bR^3)$.  Set $G:=\mathfrak D_I^-H_N$.  Since
	$|I|+m\le\mathrm N$, the functions $(\mathfrak D_X^-)^jG$, $0\le j\le m-1$,
	are $C^1$ on $\partial\cD\times\bR^3$, being derivatives of $H_N$ of order $|I|+j\le\mathrm N-1$.  Hence
	\eqref{eq:complete-flow-level-set}, applied to a $C^1$ extension of $(\mathfrak D_X^-)^jG$ to $\overline\cD\times\bR^3$ and restricted to $\gamma^-$, which the flow preserves by \eqref{eq:complete-flow-normal-velocity-preservation}, gives
	\[
	\frac{\dd}{\dd s}T_s^{X,\gamma_-}(\mathfrak D_X^-)^jG
	=T_s^{X,\gamma_-}(\mathfrak D_X^-)^{j+1}G,
	\qquad 0\le j\le m-1,\quad |s|\le s_0 .
	\]
	Integrating this identity $m$ times, with the group law
	$T_s^{X,\gamma_-}T_t^{X,\gamma_-}=T_{s+t}^{X,\gamma_-}$ of
	Lemma~\ref{lem:boundary-flow-continuity}, we obtain
	\begin{equation}
		\label{eq:fixedN-inflow-repeated-FTC}
		(\Delta_h^{X,\gamma_-})^mG
		=\int_{[0,1]^m}
		T_{(\theta_1+\cdots+\theta_m)h}^{X,\gamma_-}
		(\mathfrak D_X^-)^mG\,
		\dd\theta_1\cdots\dd\theta_m,
	\end{equation}
	where
	\[
	|(\theta_1+\cdots+\theta_m)h|\le m|h|<\mathrm Nh_0\le s_0 .
	\]
	Hence \eqref{eq:boundary-flow-bound} and
	\eqref{eq:fixedN-incoming-word-bound} at $L=0$ and the word $\xi^mI$, for which $(\mathfrak D_X^-)^mG=\mathfrak D_{\xi^mI}^-H_N$ and $|\xi^mI|=|I|+m\le\mathrm N$, give
	\eqref{eq:fixedN-inflow-separation}.
	By
	\eqref{eq:fixedN-inflow-repeated-FTC},
	\[
	\bigl\|(\Delta_h^{X,\gamma_-})^mG-(\mathfrak D_X^-)^mG
	\bigr\|_{L^2(\gamma^-)}
	\le\sup_{|s|\le m|h|}
	\bigl\|\bigl(T_s^{X,\gamma_-}(\mathfrak D_X^-)^mG
	-(\mathfrak D_X^-)^mG\bigr)\bigr\|_{L^2(\gamma^-)},
	\]
	and the right-hand side tends to zero as $h\to0$ by the strong continuity in Lemma~\ref{lem:boundary-flow-continuity}, which proves
	\eqref{eq:fixedN-incoming-difference-limit}.
\end{proof}

\begin{defn}[Induction assertions] \label{def:fixedN-package}
	Let $0\le m\le\mathrm N$, let $I$ range over the ordered complete words and let
	$L$ range over the finite weights.  We denote by $\mathsf P_m$ the conjunction of the
	following statements.
	\begin{enumerate}
		\item\label{it:package-graph}
		For $|I|\le m$, $q_{I,N}:=\mathfrak D_If_N\in\cH_{\kin,\sigma}$, its two
		flux traces are the completion-defined Green traces, and
		\[
		\|w^Lq_{I,N}\|_{\cH_{\kin,\sigma}}\le C,
		\qquad
		(q_{I,N})_{\gamma^-}=\mathfrak D_I^-H_N .
		\]
		\item\label{it:package-equation}
		For $|I|\le m$, at $f=f_N$ and with $f_I:=q_{I,N}$, \eqref{eq:weak-differentiated-equation} holds.
		\item\label{it:package-forcing}
		For $|I|\le m$,
		\[
		\|w^L(v\cdot\nabla_x q_{I,N}+L_gq_{I,N})\|_{L_x^2H_\sigma^{-1}}
		\le C .
		\]
	\end{enumerate}
\end{defn}

\begin{lemma}[Fixed-cutoff base estimates]
	\label{lem:fixedN-base-estimates}
	Fix $N$
	and let $f_N$ be the mass-normalized spectral solution
	\eqref{eq:spectral-problem}.  Set
	$\cB_N:=\cB_{\rm bs}(f_N)$.  Then, for every finite $L\ge0$,
	\begin{equation}
		\label{eq:fixedN-preliminary-weighted-base}
		\cB_N+\|w^Lf_N\|_{\cH_{\kin,\sigma}} \le C.
	\end{equation}
	Moreover,
	\begin{equation}
		\label{eq:fixedN-weighted-base-pflat}
		\|w^Lf_N\|_{L_x^{p_\flat}H_\sigma^1}\le C.
	\end{equation}
	
\end{lemma}

\begin{proof}
	We first prove \eqref{eq:fixedN-preliminary-weighted-base}.  We have $\|f_N\|_{\cH_{\kin,\sigma}}<\infty$ by
	Corollary~\ref{cor:fixedN-spectral-solvability}.  Since $\|J_N\phi_N\|_{H^3(\partial\cD)}\le C\|f_N\|_{\cH_{\kin,\sigma}}$ by Lemma~\ref{lem:fixedN-range} at $r=3$ and \eqref{eq:outgoing-moment-flux-bound}, Lemma~\ref{lem:incoming-datum-bound} at $\psi=J_N\phi_N$ gives, for every finite $L$,
	\[
	\|w^{L+B^{\st}_{\rm in}}H_N\|_{(C_x^{0,1}\cap C_v^{1,1})(\gamma^-)}
	+\|w^LH_N\|_{L^2(\gamma^-)}\le C.
	\]
	The hypothesis \eqref{eq:KS-coefficient} of Proposition~\ref{prop:base-regularity} holds: its first bound is the first bound in \eqref{eq:construction-coefficient-regularity}, assumed before \eqref{eq:spectral-problem}, and its second holds with $M=C\widetilde{\mathfrak Q}_{\mathrm N}(g)<\infty$ by \eqref{eq:coefficient-transmission} and \eqref{eq:reduced-Q3}.
	Applying \eqref{eq:base-regularity} with $\ell=A_\flat$ and
	\eqref{eq:weighted-inflow-estimate} at the weight $L$ to
	\eqref{eq:spectral-problem}, viewed as the prescribed-inflow problem with
	incoming datum $H_N$, we obtain
	\eqref{eq:fixedN-preliminary-weighted-base}.  For \eqref{eq:fixedN-weighted-base-pflat}, we apply
	\eqref{eq:base-regularity} with $\ell=L+3$ to obtain
	\[
	\|w^{L+3}f_N\|_{L^\infty}+\|w^{L+3}\nabla_vf_N\|_{L_x^{p_\flat}L_v^2}\le C ;
	\] then
	\eqref{eq:base-slot-pointwise} at $\ell=L$ and $r=p_\flat$ gives \eqref{eq:fixedN-weighted-base-pflat}.
\end{proof}

\noindent
At order zero, $q_{\varnothing,N}=f_N$. The bound \eqref{eq:fixedN-preliminary-weighted-base} and the boundary identity in \eqref{eq:spectral-boundary-data} give (\ref{it:package-graph}) of $\mathsf P_0$. The spectral equation $v\cdot\nabla_xf_N+L_gf_N=0$ of \eqref{eq:spectral-problem} gives (\ref{it:package-equation}), since the commutator in \eqref{eq:weak-differentiated-equation} vanishes for the empty word; the same equation makes (\ref{it:package-forcing}) immediate.

\begin{prop}[Order-$m$ jet realization]
	\label{prop:fixedN-order-m-realization}
	Fix $N$, let $f_N$ be the mass-normalized spectral solution \eqref{eq:spectral-problem}, and let $1\le m\le\mathrm N$.  If
	$\mathsf P_{m-1}$ of Definition~\ref{def:fixedN-package} holds, then $\mathsf P_m$ holds.
\end{prop}

\begin{proof}
	We construct the order-$m$ jets by induction on $(r,J)$, then establish $\mathsf P_m$.  Each step of this induction is divided into
	parts~(i)--(ix).
	
	\medskip
	\noindent\textbf{The induction over $(r,J)$.} Every ordered complete word of length $m$ is $\nu^rJ$ with $\nu$ a letter, $r\ge1$, and $|J|=m-r$; with $X=\tau_\nu$ we write it $X^rJ$.  For each letter $\nu$, the pairs $(r,J)$ with $r+|J|=m$ are treated in decreasing order of $r$.  Throughout the step, $X=\tau_\nu$ and $J$ are
	fixed, and we put $j_0:=|J|=m-r$; with $\mathscr T_L$ of \S\ref{app:flow-identities}, we have
	$v\cdot\nabla_x+L_g=\mathscr T_L-N[g,\cdot]$.  We also put
	\[
	Z_h:=(\Delta_h^X)^rq_{J,N},
	\qquad
	E_h:=\|Z_h\|_{\cH_{\kin,\sigma}} .
	\]
	We take $h_0<s_0/(2\mathrm N)$, with $s_0$ of \eqref{eq:complete-flow-uniform-bounds}; moreover, $h_0$ is at most the constant $h_0$ of Lemma~\ref{lem:fixedN-inflow-separation} and satisfies the hypothesis on $h_0$ of Lemma~\ref{lem:discrete-interpolation}.  Then
	every flow parameter of absolute value at most $2\mathrm Nh_0$ lies in the range $|s|\le s_0$ of Lemma~\ref{lem:complete-flow-weighted-bounds}.  Part~(viii) decreases $h_0$ further.  By the
	convention stated after \eqref{eq:bulk-flow-difference}, cube variables are written
	$\theta_1,\theta_2,\ldots$ and $\bar h$ denotes a flow parameter that arises from difference
	quotients and is not written out.  In the step, a jet carrying $j$ derivatives in a solution slot is measured in $L_x^{\mathsf e(m/j,p_\flat)}H_\sigma^1$, as for a graph row of order $m$.  The step assumes the unweighted bound
	\begin{equation}
		\label{eq:fixedN-step-input}
		\|q_{J',N}\|_{L_x^{\mathsf e(m/j_0,p_\flat)}H_\sigma^1}
		\le C \qquad\text{for every }|J'|=j_0,
	\end{equation}
	with the convention of \eqref{eq:GN-output-exponents} at $j_0=0$, and proves $\|q_{XJ,N}\|_{L_x^{\mathsf e(m/(j_0+1),p_\flat)}H_\sigma^1}\le C$.  Part~(viii) also gives difference bounds, used in part~(v) of the later steps.  At $r=m$, where $j_0=0$ and the only word of length $0$ gives $q_{\varnothing,N}=f_N$,
	\eqref{eq:fixedN-step-input} is \eqref{eq:fixedN-weighted-base-pflat} at
	$L=0$; at $r<m$ it is given by the steps at $r+1$.
	
	\medskip
	\noindent\emph{(i) Exponents.}
	We match the exponents of discrete interpolation to those of the order-$m$ solution derivatives and verify the corresponding H\"older pairing with the background.  Since \[
	1/\mathsf e(m/j,p_\flat)=\frac jm\cdot\frac12+\bigl(1-\frac jm\bigr)\frac1{p_\flat} ,
	\] which is affine in $j$ and equals $\frac12$ at $j=m$, and since, by $m-j_0=r$,
	\[
	\Bigl(1-\frac ir\Bigr)j_0+\frac ir\,m=j_0+\frac ir\,(m-j_0)=j_0+i,
	\]
	the exponents \eqref{eq:GN-output-exponents} satisfy,
	for $0\le i\le r-1$,
	\begin{equation}
		\label{eq:fixedN-step-exponent}
		\Bigl(1-\frac ir\Bigr)\frac1{\mathsf e(m/j_0,p_\flat)}+\frac i{2r}
		=\frac1{\mathsf e(m/(j_0+i),p_\flat)}=\frac 1{\mathsf e(r/i,\mathsf e(m/j_0,p_\flat))}.
	\end{equation}
	Hence,
	\eqref{eq:discrete-GN}, applied at the difference order $r$ and the base exponent
	$\mathsf e(m/j_0,p_\flat)$, yields the rung $\mathsf e(m/(j_0+i),p_\flat)$ of order
	$j_0+i$.  Moreover
	$k:=r-i\ge1$ satisfies $k+(j_0+i)=m$, so
	$(\mathsf g(k),\mathsf e(m/(j_0+i),p_\flat))$ is an admissible H\"older pair
	by \eqref{eq:dyn-row-pairing}.
	
	\medskip
	\noindent\emph{(ii) The lower differences.}
	The lower-order graph estimates control the difference quotients at every finite weight.
	For a subword $J'$ of $J$ and an integer $i\ge0$ with $|J'|+i\le m-1$, (\ref{it:package-graph}) of $\mathsf P_{m-1}$ gives
	$\|w^Lq_{X^iJ',N}\|_{\cH_{\kin,\sigma}}\le C$ for every finite $L\ge0$.  By \eqref{eq:repeated-same-flow-weighted-bound} at $u=q_{J',N}$ and $u_j=q_{X^jJ',N}$, these jets control the
	differences $(\Delta_h^X)^iq_{J',N}$, and hence the quantity of \eqref{eq:discrete-lower-density},
	\[
	\cR^{X,0}_{r,h}(q_{J,N})=\sum_{i=0}^{r-1}\|w^{r-i}(\Delta_h^X)^iq_{J,N}\|_{L_x^2H_\sigma^1},
	\]
	so that
	\begin{equation}
		\label{eq:fixedN-step-discrete-lower}
		\sup_{0<|h|<h_0}
		\Bigl(\cR_{r,h}^{X,0}(q_{J,N})
		+\sum_{J'}\sum_{i=0}^{m-1-|J'|}\|w^L(\Delta_h^X)^iq_{J',N}\|_{L_x^2H_\sigma^1}\Bigr)
		\le C \end{equation}
	for every finite $L\ge0$, where the outer sum is taken over the subwords $J'$ of $J$.
	
	\medskip
	\noindent\emph{(iii) Operator orbits and traces.}
	The operator difference identity \eqref{eq:direct-operator-cube} must hold distributionally at the arguments below; we also need $Z_h\in\cH_{\kin,\sigma}$.
	By Corollary~\ref{cor:A0-orbits} at $\ell=0$ and at every order $1\le s\le r$,
	\eqref{eq:operator-orbit-bound} holds for $\mathrm A=\mathscr T_L$ at the order $s$ with
	$(L_1,p_1)=(2s+4,2)$.  At the arguments
	$z=T_{sh}^X(\Delta_h^X)^{r-s}q_{J,N}$, $1\le s\le r$, of the difference equation of part~(iv),
	\[
	\|z\|_{L^2_{x,v}}+\|z\|_{L_x^2H_\sigma^1}+\|w^{2s+4}z\|_{L_x^2H_\sigma^1}<\infty
	\]
	by \eqref{eq:complete-flow-H1-bound}, \eqref{eq:real-weight-Hsigma-monotonicity}, \eqref{eq:local-H1-H_sigma}, and \eqref{eq:fixedN-step-discrete-lower} at $L=A_\diamond+3$, since the pullback in \eqref{eq:complete-flow-H1-bound} raises the weight $2s+4$ to $2s+5\le2\mathrm N+5=A_\diamond+3$.
	Thus \eqref{eq:direct-operator-cube} holds in $\cD'$ for $\mathrm A=\mathscr T_L$ at these arguments.
	
	By (\ref{it:package-graph}) of $\mathsf P_{m-1}$,
	$w^5q_{J,N}\in\cH_{\kin,\sigma}$, since $|J|\le m-1$.  Lemma~\ref{lem:complete-flow-actual-traces}\,(\ref{it:flow-traces-stationary}), applied
	to $q=q_{J,N}$ at the flow parameters $h,2h,\ldots,rh$, therefore gives
	\begin{equation}
		\label{eq:fixedN-step-flow-membership}
		T_{sh}^Xq_{J,N}\in\cH_{\kin,\sigma},
		\qquad
		(T_{sh}^Xq_{J,N})_{\gamma^\pm}
		=T_{sh}^{X,\gamma_\pm}(q_{J,N})_{\gamma^\pm},
		\qquad 1\le s\le r .
	\end{equation}
	In particular, since $Z_h$ is a linear combination of the $T_{sh}^Xq_{J,N}$, $0\le s\le r$, by \eqref{eq:bulk-flow-difference}, $Z_h\in\cH_{\kin,\sigma}$ and $E_h<\infty$ for each fixed $h\neq0$, and \eqref{eq:direct-same-flow-actual-traces} gives $(Z_h)_{\gamma^\pm}=(\Delta_h^{X,\gamma_\pm})^r(q_{J,N})_{\gamma^\pm}$.
	
	\medskip
	\noindent\emph{(iv) The difference equation.}
	Apply the difference operator to the weak jet equation.
	By (\ref{it:package-equation}) of $\mathsf P_{m-1}$ and
	\eqref{eq:weak-triangular-commutator},
	\begin{equation}
		\label{eq:fixedN-step-source-split}
		\mathscr T_Lq_{J,N}
		=N[g,q_{J,N}]
		+\sum_{\varnothing\ne S\subset[j_0]}
		N[\mathfrak D_{J_S}g,q_{J_{S^c},N}]
		+\mathfrak R_J[g,f_N]
		-\mathfrak G_Jf_N ,
	\end{equation}
	where $\mathfrak G_J:=[\mathfrak D_J,\mathscr T_L]$ collects the transport,
	local and compact commutators of
	\eqref{eq:weak-triangular-commutator}.  By the last identity of \eqref{eq:ad-word} at $I=J$ and $\mathrm A=\mathscr T_L$, and by induction on
	\eqref{eq:deformation-remainder-recurrence},
	the jet in every solution slot of $\mathfrak G_Jf_N$ and of $\mathfrak R_J[g,f_N]$ is $q_{J',N}$ for a proper subword $J'$ of $J$.  In particular,
	\[
	\mathfrak G_Jf_N=\sum_{\varnothing\ne S\subset[j_0]}\operatorname{ad}_{J_S}(\mathscr T_L)\,q_{J_{S^c},N}
	\qquad\text{in }\cD',
	\]
	since both sides expand $[\mathfrak D_J,\mathscr T_L]$ with each lower complete derivative
	replaced by its jet.  Each $\operatorname{ad}_{J_S}(\mathscr T_L)$ is, by
	\eqref{eq:nested-commutator-decomposition} at $g=0$, a sum of transport, local and compact
	commutators that contains no complete derivative.  By Corollary~\ref{cor:A0-orbits}, \eqref{eq:operator-orbit-bound} and its form
	with $B_s$ replaced by $\mathrm A$ hold for $\mathrm A=\operatorname{ad}_{J_S}(\mathscr T_L)$ at the order $s$ with
	$(L_1,p_1)=(2(s+|S|)+4,2)$.
	At the arguments $z=T_{sh}^X(\Delta_h^X)^{r-s}q_{J_{S^c},N}$,
	\[
	\|z\|_{L^2_{x,v}}+\|z\|_{L_x^2H_\sigma^1}+\|w^{2(s+|S|)+4}z\|_{L_x^2H_\sigma^1}<\infty
	\]
	by \eqref{eq:complete-flow-H1-bound}, \eqref{eq:real-weight-Hsigma-monotonicity}, \eqref{eq:local-H1-H_sigma} and
	\eqref{eq:fixedN-step-discrete-lower} at $L=A_\diamond+3$, since the pullback raises the weight $2(s+|S|)+4$ to $2(s+|S|)+5\le2\mathrm N+5=A_\diamond+3$ by $s+|S|\le m\le\mathrm N$; so, as for $\mathscr T_L$ in part~(iii), \eqref{eq:direct-operator-cube} holds in $\cD'$ for $\operatorname{ad}_{J_S}(\mathscr T_L)$ at these arguments.
	
	Apply $(\Delta_h^X)^r$ to
	\eqref{eq:fixedN-step-source-split}, expand $\mathscr T_L$ at $u=q_{J,N}$ and each
	$\operatorname{ad}_{J_S}(\mathscr T_L)$ at $u=q_{J_{S^c},N}$ by
	\eqref{eq:forward-discrete-Leibniz}, expand each original-Coulomb product by
	\eqref{eq:collision-mth-endpoint}, and subtract $N[g,Z_h]$.  The term of $(\Delta_h^X)^rN[g,q_{J,N}]$ carrying the full $r$ differences on the second slot and $-N[g,Z_h]$ combine into
	$N[T_{rh}^Xg-g,Z_h]$.  Indeed, in \eqref{eq:collision-mth-endpoint} at the difference order $r$ and $[a,b]=[g,q_{J,N}]$, the term with all $r$ differences on $b$ is $N[T_{rh}^Xg,Z_h]$, since, with $K_0=\varnothing$ there, $\mathscr C_{0,h}^X=N_{\mathsf K_\varnothing}=N$.  The resulting equation reads
	\begin{equation}
		\label{eq:fixedN-step-endpoint-equation}
		v\cdot\nabla_x Z_h+L_gZ_h
		=\mathscr R^{\mathscr C_0}_{r,h}
		+\mathscr R^{\mathscr C_+}_{r,h}
		+\mathscr R^{\mathscr T_L}_{r,h}
		+N[T_{rh}^Xg-g,Z_h]
		\quad\text{in }\cD' ,
	\end{equation}
	where $\mathscr R^{\mathscr C_0}_{r,h}$ collects the original-Coulomb terms of
	\eqref{eq:collision-mth-endpoint} with $S\ne\varnothing$ or with at least
	one difference on the first slot, $\mathscr R^{\mathscr C_+}_{r,h}$ the defect
	terms $\mathscr C^X_{s,h}$, $s\ge1$, of \eqref{eq:collision-flow-defect-difference}
	together with the differences of $\mathfrak R_J[g,f_N]$, and
	\begin{equation*}
		\begin{aligned}
			\mathscr R^{\mathscr T_L}_{r,h}
			:={}&-\sum_{s=1}^{r}\binom rs
			\mathfrak Q_{s,h}^X(\mathscr T_L)
			T_{sh}^X(\Delta_h^X)^{r-s}q_{J,N}\\
			&-\sum_{\varnothing\ne S\subset[j_0]}\Bigl(
			\operatorname{ad}_{J_S}(\mathscr T_L)(\Delta_h^X)^rq_{J_{S^c},N}
			+\sum_{s=1}^{r}\binom rs
			\mathfrak Q_{s,h}^X\bigl(\operatorname{ad}_{J_S}(\mathscr T_L)\bigr)
			T_{sh}^X(\Delta_h^X)^{r-s}q_{J_{S^c},N}\Bigr).
		\end{aligned}
	\end{equation*}
	Every solution factor in $\mathscr R^{\mathscr T_L}_{r,h}$ therefore has complete order at most $m-1$: in the first sum $|J|+r-s=m-s$, and in the second $|J_{S^c}|+i\le |J|-1+r=m-1$. The flow translations do not change these orders.
	
	\medskip
	\noindent\emph{(v) The solution factors.}
	The solution factors in $\mathscr R^{\mathscr C_0}_{r,h}$ and $\mathscr R^{\mathscr C_+}_{r,h}$ are estimated at their corresponding spatial exponents.  We have
	\[
	\|(\Delta_h^X)^rq_{J,N}\|_{L_x^2H_\sigma^1}=\|Z_h\|_{L_x^2H_\sigma^1}\le E_h .
	\]
	With this, \eqref{eq:discrete-GN} at the difference order $r$ and the base exponent
	$\mathsf e(m/j_0,p_\flat)$, together with \eqref{eq:fixedN-step-exponent},
	\eqref{eq:fixedN-step-input} and \eqref{eq:fixedN-step-discrete-lower}, gives
	\begin{equation}
		\label{eq:fixedN-step-GN}
		\|(\Delta_h^X)^iq_{J,N}\|
		_{L_x^{\mathsf e(m/(j_0+i),p_\flat)}H_\sigma^1}
		\lesssim 1+E_h^{\,i/r},
		\qquad 1\le i\le r-1 .
	\end{equation}
	Every solution factor $(\Delta_h^X)^iq_{J',N}$ with $J'$ a subword of $J$ and
	$|J'|+i\le m-1$ is measured in $L_x^{\mathsf e(m/(|J'|+i),p_\flat)}H_\sigma^1$, and, with $k:=m-|J'|-i$, $(\mathsf g(k),\mathsf e(m/(|J'|+i),p_\flat))$ is an admissible H\"older pair by \eqref{eq:dyn-row-pairing}, as in part~(i).  If $i=0$, the bound is \eqref{eq:fixedN-step-input} with $|J'|$ in place of $j_0$.  If
	$i\ge1$ and $J'=J$, it is \eqref{eq:fixedN-step-GN}.  If $i\ge1$ and $|J'|<j_0$, then the
	step at $(r',J')$, $r':=m-|J'|>r$, with the same letter $\nu$, precedes the present one, and, since $i\le r'-1$, the bound is \eqref{eq:fixedN-step-GN} of that step, with $E_h$ replaced by $C$.
	
	\medskip
	\noindent\emph{(vi) The background factors.}
	The background differences are controlled by the weighted hierarchy; the remaining coefficient difference is handled by a continuity modulus.
	The complete $X$-jets of $g$ satisfy $\mathfrak D_X\mathfrak D_{X^jK}g=\mathfrak D_{X^{j+1}K}g$ in $\cD'$ for $j+1+|K|\le\mathrm N$, and
	$w^{-2}\mathfrak D_Kg\in L_x^1L_v^2$ for $|K|\le\mathrm N$ by
	\eqref{eq:base-fiber-Lsharp} if $K=\varnothing$ and by \eqref{eq:Dg-rungs} if $|K|\ge1$.  Hence, for every ordered complete word
	$K$ and every $1\le k_1\le\mathrm N-|K|$, the hypotheses of Lemma~\ref{lem:repeated-same-flow-FTC}\,(\ref{it:repeated-FTC-fiber}) hold for $G_j:=\mathfrak D_{X^jK}g$,
	$0\le j\le k_1$, at $\lambda=-2$, and
	\eqref{eq:repeated-same-flow-FTC-fiber} applies for every $\bar h$ with
	$|\bar h|\le2\mathrm Nh_0-k_1|h|$, since, for such $\bar h$ and all
	$\theta_1,\ldots,\theta_{k_1}\in[0,1]$,
	\[
	|\bar h+(\theta_1+\cdots+\theta_{k_1})h|\le|\bar h|+k_1|h|\le2\mathrm Nh_0<s_0 .
	\]
	By
	\eqref{eq:base-fiber-Lsharp} and
	\eqref{eq:complete-flow-fiber-bound},
	\begin{equation}
		\label{eq:fixedN-step-background-zero-shifts}
		\sup_{|\bar h|\le2\mathrm Nh_0}
		\|w^{-2}T_{\bar h}^Xg\|_{L_x^\infty L_v^2}
		\le C\widetilde{\mathfrak Q}_{\mathrm N}(g).
	\end{equation}
	By \eqref{eq:repeated-same-flow-fiber-bound}, \eqref{eq:complete-flow-fiber-bound},
	\eqref{eq:Dg-rungs}, and \eqref{eq:reduced-Q3}, for every ordered complete word $K$ and every integer $k_1\ge0$,
	\begin{align}
		\label{eq:fixedN-step-background-differences}
		\begin{gathered}
			\|w^{-2}T_{\bar h}^X(\Delta_h^X)^{k_1}
			\mathfrak D_Kg\|_{L_x^{\mathsf g(k)}L_v^2}\le C,\\
			\qquad 1\le k=k_1+|K|\le\mathrm N,
			\quad 0<|h|<h_0,\quad |\bar h|\le2\mathrm Nh_0-k_1|h| ,
		\end{gathered}
	\end{align}
	with $\mathsf g(\mathrm N)=q_\flat$.
	For the coefficient difference $T_{rh}^Xg-g$ in the last term of
	\eqref{eq:fixedN-step-endpoint-equation}, we use the modulus
	\[
	\omega_g(\rho):=
	\max_{\nu\in\mathfrak A}\ \sup_{|s|\le\rho}
	\|\mu^{1/4}(T_s^{\tau_\nu}g-g)\|_{L^\infty_{x,v}} ;
	\]
	then \eqref{eq:reduced-Q3}, Lemma~\ref{lem:coefficient-transmission}, and
	Lemma~\ref{lem:complete-flow-coefficient-continuity} give
	$\omega_g(\rho)\to0$ as $\rho\downarrow0$.
	
	\medskip
	\noindent\emph{(vii) The right-hand side.}
	It remains to estimate the four terms in \eqref{eq:fixedN-step-endpoint-equation} in $L_x^2H_\sigma^{-1}$.
	
	Every summand of $\mathscr R^{\mathscr C_0}_{r,h}$ is a product of a background factor of the form \eqref{eq:fixedN-step-background-differences}, of some
	complete order $k\ge1$, with a solution factor
	$(\Delta_h^X)^iq_{J',N}$, where
	$J'$ is a subword of $J$ with $k+|J'|+i=m$ and $|J'|+i\le m-1$.
	By part~(v), $\|(\Delta_h^X)^iq_{J',N}\|_{L_x^{\mathsf e(m/(|J'|+i),p_\flat)}H_\sigma^1}\le C$ unless $J'=J$ and $1\le i\le r-1$.  Hence, by \eqref{eq:velalg} in the form \eqref{eq:velalg-poly},
	\eqref{eq:fixedN-step-background-differences}, and the bounds of part~(v),
	\begin{equation}
		\label{eq:fixedN-step-material}
		\|\mathscr R^{\mathscr C_0}_{r,h}\|_{L_x^2H_\sigma^{-1}}
		\le C+C\sum_{i=1}^{r-1} \|(\Delta_h^X)^iq_{J,N}\|
		_{L_x^{\mathsf e(m/(j_0+i),p_\flat)}H_\sigma^1}.
	\end{equation}
	
	For $\mathscr R^{\mathscr C_+}_{r,h}$, let $1/q+1/p\le1/2$ and $p<\infty$.  Proposition~\ref{prop:complete-flow-collision-differences}
	at the order $s$, $\ell=0$, $m_0=0$ and $\varpi=1/4$, together with
	$\mu^{1/4}\le Cw^{-2}$ and $2s+4\le A_\diamond+2$, gives
	\begin{equation}\label{261016940}
		\|\mathscr C_{s,h}^X[a,b]\|_{L_x^2H_\sigma^{-1}}
		\lesssim
		\sup_{\theta_1,\ldots,\theta_s\in[0,1]}
		\|w^{-2}T_{(s-\theta_1-\cdots-\theta_s)h}^Xa\|_{L_x^qL_v^2}
		\|w^{A_\diamond+2}b\|_{L_x^pH_\sigma^1},
		\qquad 1\le s\le r .
	\end{equation}
	The hypotheses of \eqref{eq:collision-flow-defect-Peano-class} follow from the right-hand side of \eqref{261016940}. Indeed, setting all $\theta_j=1$ gives $w^{-2}a\in L_x^qL_v^2\subset L^2_{x,v}$, while $p\ge2$ and $2s+3\le A_\diamond+2$ give $w^{2s+3}b\in L_x^2H_\sigma^1$ by \eqref{eq:real-weight-Hsigma-monotonicity}.
	Each defect term of
	$\mathscr R^{\mathscr C_+}_{r,h}$ arises from \eqref{eq:collision-mth-endpoint} at the difference order $r$, applied to $N[\mathfrak D_{J_S}g,q_{J_{S^c},N}]$ with $S\subset[j_0]$, and is $\frac{r!}{k_1!\,i!\,s!}\,\mathscr C^X_{s,h}[a,b]$ with
	\[
	a=T_{(r-s-k_1)h}^X(\Delta_h^X)^{k_1}\mathfrak D_{J_S}g,
	\qquad
	b=(\Delta_h^X)^iq_{J_{S^c},N},
	\qquad s\ge1,\quad k_1+i+s=r .
	\]
	Thus $a$ is a background factor of complete order $k:=k_1+|S|$, and $b$ is a solution factor of complete order
	\[
	j':=|J_{S^c}|+i=m-k-s\le m-k-1.
	\]
	Moreover, for $\theta_1,\ldots,\theta_s\in[0,1]$, 
	$$
	T_{(s-\theta_1-\cdots-\theta_s)h}^Xa=T_{\bar h}^X(\Delta_h^X)^{k_1}\mathfrak D_{J_S}g\quad \text{and}\quad |\bar h|\le(r-k_1)|h|\le2\mathrm Nh_0-k_1|h|.
	$$
	At $k=0$, where $a=T_{(r-s)h}^Xg$ and $b=(\Delta_h^X)^{r-s}q_{J,N}$, we take $(q,p)=(\infty,2)$: $\|w^{A_\diamond+2}b\|_{L_x^2H_\sigma^1}\le C$, uniformly in $0<|h|<h_0$, by
	\eqref{eq:fixedN-step-discrete-lower} at $L=A_\diamond+2$, since $r-s\le r-1$, and, by
	\eqref{eq:fixedN-step-background-zero-shifts},
	\[
	\sup_{\theta_1,\ldots,\theta_s\in[0,1]}
	\|w^{-2}T_{(r-\theta_1-\cdots-\theta_s)h}^Xg\|_{L_x^\infty L_v^2}
	\le C\widetilde{\mathfrak Q}_{\mathrm N}(g) .
	\]
	At $k\ge1$ we apply Lemma~\ref{lem:dyn-hier-weight-trade} to the solution factor $b$ at $\ell=A_\diamond+2$, at the ratio $\theta=1/m$ and at $p=\mathsf e(m/j',p_\flat)$.
	The exponent $\widetilde p$ of
	\eqref{eq:weight-trade-output-exponents} at $p=\mathsf e(m/j',p_\flat)$ satisfies, by
	\eqref{eq:GN-output-exponents}, \eqref{eq:dyn-own-order-exponents} and
	\eqref{eq:dyn-own-order-endpoints},
	\begin{equation}
		\label{261016501}
		\frac12-\frac1{\widetilde p}
		=\Bigl(1-\frac1m\Bigr)\frac{m-j'}m\Bigl(\frac12-\frac1{p_\flat}\Bigr)
		\ge\frac km\Bigl(\frac12-\frac1{p_\flat}\Bigr)
		>\frac k{\mathrm N\,q_\flat}=\frac1{\mathsf g(k)} ,
	\end{equation}
	since $m-j'\ge k+1$, $(1-\frac1m)(k+1)\ge k$ for $k\le m-1$, $m\le\mathrm N$ and $\frac1{q_\flat}<\frac12-\frac1{p_\flat}$.
	Hence $(q,p)=(\mathsf g(k),\widetilde p)$ is admissible in the application of \eqref{eq:collision-flow-defect-Peano-class} above, and Lemma~\ref{lem:dyn-hier-weight-trade} gives
	\[
	\|w^{A_\diamond+2}b\|_{L_x^{\widetilde p}H_\sigma^1}
	\le C\|b\|_{L_x^{\mathsf e(m/j',p_\flat)}H_\sigma^1}^{1-\theta}
	\|w^{m(A_\diamond+2)}b\|_{L_x^2H_\sigma^1}^{\theta}.
	\]
	For the $L_x^2$ endpoint, \eqref{eq:fixedN-step-discrete-lower} at $L=m(A_\diamond+2)$ gives
	$\|w^{m(A_\diamond+2)}b\|_{L_x^2H_\sigma^1}\le C$, uniformly in $0<|h|<h_0$.  The unweighted endpoint at
	$\mathsf e(m/j',p_\flat)$ is the bound of part~(v) for that factor $(\Delta_h^X)^iq_{J',N}$, where $J'$ is a subword of $J$ and $|J'|+i=j'$.  Hence
	\eqref{261016940}, with \eqref{eq:fixedN-step-background-zero-shifts} and \eqref{eq:fixedN-step-discrete-lower} at $k=0$, and with \eqref{eq:fixedN-step-background-differences} at $K=J_S$ and the bound of Lemma~\ref{lem:dyn-hier-weight-trade} above, whose two endpoints are bounded as just described, at $k\ge1$, together with \eqref{eq:collision-deformation-difference-bound} for the differences of $\mathfrak R_J[g,f_N]$, gives, using
	$(1+E_h^{\,i/r})^{1-\theta}\le1+E_h^{\,i/r}$,
	\begin{equation}
		\label{eq:fixedN-step-defect-bound}
		\|\mathscr R^{\mathscr C_+}_{r,h}\|_{L_x^2H_\sigma^{-1}}
		\le C+C\sum_{i=1}^{r-1}E_h^{\,i/r}, \qquad 0<|h|<h_0 .
	\end{equation}
	
	By \eqref{eq:direct-operator-cube}, each
	operator quotient in $\mathscr R^{\mathscr T_L}_{r,h}$ is an average of
	$T_{\bar h}^X\operatorname{ad}_{X^sI}(\mathscr T_L)T_{-\bar h}^X$ with $I=\varnothing$ or $I=J_S$,
	the word $X^sI$ having length at most $m$.  Since, for every operator $B$,
	\[
	T_{\bar h}^XBT_{-\bar h}^XT_{sh}^X=T_{sh}^XT_{\bar h-sh}^XBT_{-(\bar h-sh)}^X,
	\]
	each operator quotient composed with $T_{sh}^X$ is bounded from $w^{2m+5}L_x^2H_\sigma^1$ into
	$L_x^2H_\sigma^{-1}$ by \eqref{eq:complete-flow-Hminus1-bound} and
	\eqref{eq:A0-orbit-bound} at $\ell=1$.  The operators $\operatorname{ad}_{J_S}(\mathscr T_L)$ are
	bounded from $w^{2m+4}L_x^2H_\sigma^1$ into $L_x^2H_\sigma^{-1}$ by
	\eqref{eq:A0-orbit-bound} at $\bar h=0$.  Since the jets in the solution slots of $\mathscr R^{\mathscr T_L}_{r,h}$ are the differences listed at the end of part~(iv), they are bounded by
	\eqref{eq:fixedN-step-discrete-lower} at the weight $2m+5$, and we obtain
	\begin{equation}
		\label{eq:fixedN-step-linear}
		\sup_{0<|h|<h_0}
		\|\mathscr R^{\mathscr T_L}_{r,h}\|_{L_x^2H_\sigma^{-1}}\le C .
	\end{equation}
	
	Since $r|h|\le\mathrm Nh_0$, the collision bound \eqref{260911445} and the definition of $\omega_g$ give
		\[
		\|N[T_{rh}^Xg-g,Z_h]\|_{L_x^2H_\sigma^{-1}}
		\le C\omega_g(\mathrm Nh_0)E_h .
		\]
		Here $C$ is independent of the spectral cutoff.
	
	\medskip
	\noindent\emph{(viii) The bound uniform in $h$.}
	We first bound the incoming trace of $Z_h$.  By \eqref{eq:fixedN-step-flow-membership},
	\eqref{eq:direct-same-flow-actual-traces}, and (\ref{it:package-graph}) of
	$\mathsf P_{m-1}$,
	\[
	(Z_h)_{\gamma^-}=(\Delta_h^{X,\gamma_-})^r\mathfrak D_J^-H_N ;
	\]
	hence \eqref{eq:fixedN-inflow-separation}, applied at the word $J$ and the
	difference order $r$, which is possible since $|J|+r=m\le\mathrm N$, gives
	\[
	\sup_{0<|h|<h_0}\|(Z_h)_{\gamma^-}\|_{L^2(\gamma^-)}\le C .
	\]
	The constant in the graph estimate \eqref{eq:perturbed-landau-stationary-estimate} and the coefficient of $\omega_g$ above are independent of $N$. Since $\omega_g(\rho)\to0$ as $\rho\downarrow0$, we may reduce $h_0$ so that, after applying the graph estimate, the last term of \eqref{eq:fixedN-step-endpoint-equation} contributes at most $E_h/4$. The bounds \eqref{eq:fixedN-step-material}--\eqref{eq:fixedN-step-linear} and the incoming trace estimate then give 
	$$
	E_h\le C+C\sum_{i=1}^{r-1}E_h^{\,i/r}+\frac14E_h.
	$$
	Young's inequality absorbs the terms $E_h^{i/r}$, $1\le i<r$. Since $E_h<\infty$ for each nonzero $h$ by part~(iii), we obtain
	\begin{equation}
		\label{eq:fixedN-step-graph}
		\sup_{0<|h|<h_0}E_h\le C.
	\end{equation}
	The bound in \eqref{eq:fixedN-step-GN} follows as well. Finite-rank regularity is used only to establish $E_h<\infty$ at fixed $N$. The coefficient absorbed above is independent of the spectral cutoff, but the inhomogeneous bound at this stage may depend on $N$ through the incoming datum. Cutoff-independent row estimates are obtained in Proposition~\ref{prop:uniform-hierarchy}.
	
	\medskip
	\noindent\emph{(ix) Identification of the limit.}
	The uniform graph bound identifies the limiting jet through graph closure. A separate weak-identification argument provides the spatial exponent needed for the next step.
	By \eqref{eq:fixedN-step-graph} and
	Lemma~\ref{lem:complete-flow-graph-closure} at the order $r$, $\mathfrak D_X^rq_{J,N}\in\cH_{\kin,\sigma}$, and, as $h\to0$, $Z_h$ converges weakly to $\mathfrak D_X^rq_{J,N}$ in $\cH_{\kin,\sigma}$, its flux traces converging weakly to those of $\mathfrak D_X^rq_{J,N}$.  Thus
	\begin{equation}
		\label{eq:fixedN-step-limit}
		q_{X^rJ,N}=\mathfrak D_X^rq_{J,N}\in\cH_{\kin,\sigma},
		\qquad
		(q_{X^rJ,N})_{\gamma^\pm}
		=\operatorname*{w-lim}_{h\to0}(Z_h)_{\gamma^\pm} .
	\end{equation}
	By \eqref{eq:fixedN-step-limit}, the identity for $(Z_h)_{\gamma^-}$ in
	part~(viii), and \eqref{eq:fixedN-incoming-difference-limit} at the word $J$ and
	the difference order $r$,
	\[
	(q_{X^rJ,N})_{\gamma^-}=(\mathfrak D_X^-)^r\mathfrak D_J^-H_N
	=\mathfrak D_{X^rJ}^-H_N .
	\]
	For $r\ge2$, Lemma~\ref{lem:complete-flow-weak-identification} and
	\eqref{eq:complete-flow-weak-identification-lsc}, applied to $q=q_{J,N}$ at the
	difference order $1$ and at $p=\mathsf e(m/(j_0+1),p_\flat)$, turn \eqref{eq:fixedN-step-GN} at $i=1$, with $E_h$ replaced by $C$, into $\|q_{XJ,N}\|_{L_x^{\mathsf e(m/(j_0+1),p_\flat)}H_\sigma^1}\le C$.  At $r=1$, where $j_0+1=m$ and hence $\mathsf e(m/(j_0+1),p_\flat)=\mathsf e(1,p_\flat)=2$, the bound $\|q_{XJ,N}\|_{L_x^2H_\sigma^1}\le C$ follows from \eqref{eq:fixedN-step-graph} and \eqref{eq:fixedN-step-limit}.
	
	\medskip
	Taking $r=1$ for each letter $\nu$ and every word $J$ of length $m-1$ realizes all jets of order $m$ and establishes
	\eqref{eq:distributional-jet-chain} through order $m$.
	
	\medskip
	\noindent\textbf{Verification of $\mathsf P_m$.}
	For words of length $m$, we first establish the differentiated equation, then the forcing bound; the weighted graph bound follows.
	
	\emph{(\ref{it:package-equation}).}
	The hypotheses of Theorem~\ref{thm:weak-differentiated-system} hold at the order $m$.
	\begin{itemize}
		\item The jets of length at most $m$ lie in $\cH_{\kin,\sigma}$ and satisfy the chain \eqref{eq:distributional-jet-chain}, by
		(\ref{it:package-graph}) of $\mathsf P_{m-1}$, \eqref{eq:fixedN-step-limit} at $r=1$ and \eqref{eq:ordered-complete-derivatives}.
		\item The lower-order bound \eqref{eq:lower-weight-assumption} holds by
		(\ref{it:package-graph}) of $\mathsf P_{m-1}$ at $L=A_\diamond$.
		\item The background condition \eqref{eq:background-jet-assumption} holds by \eqref{eq:reduced-Q3}.
	\end{itemize}
	Theorem~\ref{thm:weak-differentiated-system} therefore gives \eqref{eq:weak-differentiated-equation} for every
	ordered word of length $m$.  This is (\ref{it:package-equation}) of $\mathsf P_m$.
	
	\emph{Reduction to the forcing.}
	Fix a finite weight $L\ge0$ and an ordered complete word $I$ of length $m$.  By
	part~(ix), $q_{I,N}\in\cH_{\kin,\sigma}$ and
	$(q_{I,N})_{\gamma^-}=\mathfrak D_I^-H_N$, and
	$\|w^L\mathfrak D_I^-H_N\|_{L^2(\gamma^-)}\le C$ by \eqref{eq:fixedN-incoming-word-bound}.  Hence
	Theorem~\ref{thm:weighted-steady}\,(2), in the form
	\eqref{eq:weighted-steady-estimate}, gives
	\begin{equation}
		\label{eq:fixedN-order-m-graph-reduction}
		\|w^Lq_{I,N}\|_{\cH_{\kin,\sigma}}
		\le C\bigl(\|w^L(v\cdot\nabla_xq_{I,N}+L_gq_{I,N})\|_{L_x^2H_\sigma^{-1}} +1\bigr)
	\end{equation}
	whenever
	\[
	\|w^L(v\cdot\nabla_xq_{I,N}+L_gq_{I,N})\|_{L_x^2H_\sigma^{-1}}<\infty.
	\]
	By (\ref{it:package-equation}) of $\mathsf P_m$,
	verified above, and \eqref{eq:fixedN-step-source-split} with $I$ in place of $J$ (recall that $v\cdot\nabla_x+L_g=\mathscr T_L-N[g,\cdot]$),
	\begin{equation}
		\label{eq:fixedN-order-m-forcing}
		v\cdot\nabla_xq_{I,N}+L_gq_{I,N}
		=\sum_{\varnothing\ne S\subset[m]}
		N[\mathfrak D_{I_S}g,q_{I_{S^c},N}]
		+\mathfrak R_I[g,f_N]
		-\mathfrak G_If_N ,
	\end{equation}
	the jet in every solution slot on the right being that of a word of length at most $m-1$.  It
	remains to bound the three terms on the right in $w^LL_x^2H_\sigma^{-1}$ by $C$, which is (\ref{it:package-forcing}) of $\mathsf P_m$.
	
	\emph{The commutators $\mathfrak G_If_N$.}
	The transport, local and compact commutators in $\mathfrak G_If_N$ involve solution jets measured in $L_x^2H_\sigma^1$ or $L_{x,v}^2$.  By (\ref{it:package-graph}) of $\mathsf P_{m-1}$ at the weights $L+4$, $L+2m+2$ and $0$, and by \eqref{260902137},
	\[
	\|w^{L+4}q_{J,N}\|_{L_x^2H_\sigma^1}+\|w^{L+2m+2}q_{J,N}\|_{L_x^2H_\sigma^1}+\|q_{J,N}\|_{L^2_{x,v}}\le C,
	\qquad |J|\le m-1;
	\]
	hence \eqref{eq:ordered-transport-weight-integrated}, \eqref{eq:local-L0-complete-bound} and \eqref{eq:iterated-compact-full-operator-bound} give $\|w^L\mathfrak G_If_N\|_{L_x^2H_\sigma^{-1}}\le C$.
	
	\emph{The original-Coulomb terms.} Let $|S|=k\ge1$ and $j:=m-k$.  By \eqref{eq:velalg}, in the form \eqref{eq:velalg-poly} at $\ell=L$, and H\"older's inequality in $x$,
	$\|w^LN[\mathfrak D_{I_S}g,q_{I_{S^c},N}]\|_{L_x^2H_\sigma^{-1}}<\infty$ provided that
	$w^{-2}\mathfrak D_{I_S}g\in L_x^qL_v^2$ and $w^Lq_{I_{S^c},N}\in L_x^pH_\sigma^1$ with $\frac1q+\frac1p\le\frac12$.  By \eqref{eq:Dg-rungs}, $w^{-2}\mathfrak D_{I_S}g\in L_x^{\mathsf g(k)}L_v^2$.
	At $j=0$ the jet in the solution slot is $f_N$, and $w^Lf_N\in L_x^{p_\flat}H_\sigma^1$ by
	\eqref{eq:fixedN-weighted-base-pflat}, and the pair $(\mathsf g(m),p_\flat)$ is admissible by
	\eqref{eq:dyn-row-pairing} at $j=0$.  For
	$1\le j\le m-1$ we apply Lemma~\ref{lem:dyn-hier-weight-trade} to $q_{I_{S^c},N}$ at a ratio $\theta$, between \eqref{eq:fixedN-step-input} with $j$ in place of $j_0$ and (\ref{it:package-graph}) of
	$\mathsf P_{m-1}$ at the weight $L/\theta$, which gives
	\[
	\|w^Lq_{I_{S^c},N}\|_{L_x^{\widetilde p}H_\sigma^1}\le C,
	\]
	with $\widetilde p$ of \eqref{eq:weight-trade-output-exponents} at
	$p=\mathsf e(m/j,p_\flat)$.
	The ratio depends on whether $m\le\mathrm N-1$.
	
	If $m\le\mathrm N-1$, take $\theta=1/(m+1)$. Then, by \eqref{eq:weight-trade-output-exponents} and \eqref{eq:GN-output-exponents},
	\[
	\frac1{\widetilde p}
	=\frac1{2(m+1)}+\frac m{m+1}\Bigl(\frac j{2m}+\frac{m-j}m\cdot\frac1{p_\flat}\Bigr)
	=\frac{j+1}{m+1}\cdot\frac12+\Bigl(1-\frac{j+1}{m+1}\Bigr)\frac1{p_\flat},
	\]
	that is, $\widetilde p=\mathsf e((m+1)/(j+1),p_\flat)$.  Since $k+(j+1)=m+1\le\mathrm N$, the pair $(\mathsf g(k),\widetilde p)$ is admissible by \eqref{eq:dyn-row-pairing} with $(m+1,j+1)$ in place of $(m,j)$.
	
	If $m=\mathrm N$, take $\theta=1/80$, so that $\theta\bigl(\frac12-\frac1{p_\flat}\bigr)=\frac1{80}\cdot\frac49=\frac12\delta_\flat$ by \eqref{eq:dyn-own-order-endpoints}.  By
	\eqref{eq:GN-output-exponents}, \eqref{eq:dyn-own-order-exponents},
	\eqref{eq:dyn-own-order-endpoints} and $k+j=\mathrm N$, we have
	\[
	\frac1{\mathsf g(k)}+\frac1{\widetilde p}-\frac12
	=\frac k{\mathrm N}\Bigl(\theta\Bigl(\frac12-\frac1{p_\flat}\Bigr)-\delta_\flat\Bigr)
	=-\frac k{2\mathrm N}\,\delta_\flat .
	\]
	In either case, \eqref{eq:velalg-poly} at $\ell=L$, H\"older's inequality in $x$, and the bounds \eqref{eq:Dg-rungs} and $\|w^Lq_{I_{S^c},N}\|_{L_x^{\widetilde p}H_\sigma^1}\le C$ give
	\[
	\|w^LN[\mathfrak D_{I_S}g,q_{I_{S^c},N}]\|_{L_x^2H_\sigma^{-1}}
	\le C\|w^{-2}\mathfrak D_{I_S}g\|_{L_x^{\mathsf g(k)}L_v^2}\|w^Lq_{I_{S^c},N}\|_{L_x^{\widetilde p}H_\sigma^1}
	\le C .
	\]
	
	\emph{The deformation term.}
	Every term of $\mathfrak R_I[g,f_N]$ has the form
	\eqref{eq:deformation-remainder-structure}, with a background factor $\mathfrak D_{J_1}g$ of order $k=|J_1|$
	and a jet $q_{J_2,N}$ of order $j=|J_2|$ in the solution slot, $k+j\le m-1$, and is bounded by \eqref{eq:deformation-weighted-output-spatial}, which requires a bound on this jet at a weight at most $L+A_\diamond$.
	At $k=0$ the pair is $(\infty,2)$, as in part~(vii): $\|w^{-2}g\|_{L_x^\infty L_v^2}\le C\widetilde{\mathfrak Q}_{\mathrm N}(g)$ by \eqref{eq:base-fiber-Lsharp}, and $\|w^{L+A_\diamond}q_{J_2,N}\|_{L_x^2H_\sigma^1}\le C$ by (\ref{it:package-graph}) of $\mathsf P_{m-1}$ at the weight $L+A_\diamond$.
	At $k\ge1$, we use \eqref{eq:Dg-rungs} for the background factor and, for the jet, argue as for the solution factors of the defect terms in part~(vii): since $m-j\ge k+1$, we apply Lemma~\ref{lem:dyn-hier-weight-trade} at the ratio $\theta=1/m$ between \eqref{eq:fixedN-step-input} with $j$ in place of $j_0$ and (\ref{it:package-graph}) of
	$\mathsf P_{m-1}$ at the weight $m(L+A_\diamond)$; the pair $(\mathsf g(k),\widetilde p)$ is admissible by \eqref{261016501} with $j$ in place of $j'$.  Thus, for a term of $\mathfrak R_I[g,f_N]$ of the form \eqref{eq:deformation-remainder-structure} with $|J_1|=k\ge1$ and $|J_2|=j$, the estimate \eqref{eq:deformation-weighted-output-spatial} at $\ell=L$ and $\varpi=\frac14$, together with $\mu^{1/4}\le Cw^{-2}$ and \eqref{eq:real-weight-Hsigma-monotonicity}, gives the first inequality below, and \eqref{eq:Dg-rungs} and Lemma~\ref{lem:dyn-hier-weight-trade} give the second:
	\[
	\|w^LpN_{\mathsf K_A}[p_1\mathfrak D_{J_1}g,p_2q_{J_2,N}]\|_{L_x^2H_\sigma^{-1}}
	\le C\|w^{-2}\mathfrak D_{J_1}g\|_{L_x^{\mathsf g(k)}L_v^2}\|w^{L+A_\diamond}q_{J_2,N}\|_{L_x^{\widetilde p}H_\sigma^1}
	\le C .
	\]
	
	\emph{Conclusion.}
	By \eqref{eq:fixedN-order-m-forcing} and the bounds above,
	\[
	\|w^L(v\cdot\nabla_xq_{I,N}+L_gq_{I,N})\|_{L_x^2H_\sigma^{-1}}\le C,
	\]
	which is (\ref{it:package-forcing}) of $\mathsf P_m$.  Consequently, \eqref{eq:fixedN-order-m-graph-reduction} gives $\|w^Lq_{I,N}\|_{\cH_{\kin,\sigma}}\le C$ at every finite weight $L$, which is the bound in (\ref{it:package-graph}) of $\mathsf P_m$, and
	part~(ix) identifies the flux traces required in (\ref{it:package-graph}).  The induction is complete.
\end{proof}

\begin{proof}[\textbf{Proof of Lemma~\ref{lem:fixedN-complete-jets}}]
	Lemma~\ref{lem:fixedN-base-estimates} gives $\mathsf P_0$.
	Proposition~\ref{prop:fixedN-order-m-realization}, applied successively at
	$m=1,\ldots,\mathrm N$, gives $\mathsf P_{\mathrm N}$.  Then
	(\ref{it:package-graph}) of $\mathsf P_{\mathrm N}$ gives \eqref{eq:fixedN-jet-fullgraph}, the
	identification of the two flux traces and the weighted bound; with \eqref{eq:ordered-complete-derivatives}, it also gives \eqref{eq:distributional-jet-chain}, and
	(\ref{it:package-equation}) of $\mathsf P_{\mathrm N}$ gives \eqref{eq:weak-differentiated-equation}.
\end{proof}

\subsection{The hierarchy and the removal of the cutoff}\label{ss:diffuse-uniform-hierarchy}

\begin{thm}[Fixed-background estimate]
	\label{thm:fixed-background-construction}
	There exist $\eta_\st,\delta_\st\in(0,1]$, depending only on the structural data, such that, for every measurable $g$ on $\cD\times\bR_v^3$, if
	\begin{equation}
		\label{eq:fixed-background-smallness}
		\widetilde{\mathfrak Q}_{\mathrm N}(g)\le\eta_\st,
		\qquad 0<\delta\le\delta_\st,
	\end{equation}
	then the unique mass-normalized solution
	$f\in\cH_{\kin,\sigma}(\cD\times\bR_v^3)$ of
	\eqref{eq:full-diffuse-linear-problem} at the background $g$ with $F=0$, $b=h_w$ and $m=0$, furnished by
	Proposition~\ref{prop:full-diffuse-graph-solvability}, satisfies
	\begin{equation}
		\label{eq:fixed-background-construction}
		\widetilde{\mathfrak Q}_{\mathrm N}(f)
		+\|\cA_0f_{\gamma^+}\|_{H^{\mathrm N}(\partial\cD)}
		\le C\delta,
	\end{equation}
	with $C$ independent of $g$ and $\delta$.
\end{thm}

We first obtain bounds independent of the spectral cutoff $N$ and then pass to the limit at the fixed background $g$.

\begin{lemma}[Ladder closure]
	\label{lem:ladder-closure}
	Fix $M\in\bN$ and constants $C_0,\ldots,C_M\ge1$. There are
	$\theta_*>0$ and $C>0$, depending only on these constants and $M$,
	with the following property. Let $\cR_0,\ldots,\cR_M\in[0,\infty]$,
	$\cN_m:=\sum_{k=0}^m\cR_k$, $\cN_{-1}:=0$, and $\cQ\in[0,\infty)$.
	If $\cN_M<\infty$ and, for some $0\le\theta\le\theta_*$,
	\begin{equation}
		\label{eq:ladder-step}
		\cR_m\le C_m\bigl(\cQ+\cN_{m-1}\bigr)
		+C_m\theta\cR_m,\qquad 0\le m\le M,
	\end{equation}
	then $\cN_M\le C\cQ$.
\end{lemma}

\begin{proof}
	Choose $\theta_*$ so that $C_m\theta_*\le1/2$ for $0\le m\le M$. Thus
	$\cR_m\le 2C_m(\cQ+\cN_{m-1})$.
	Induction on $m$ gives $\cN_M\le C\cQ$, with $C$ depending only
	on $M$ and $C_0,\ldots,C_M$.
\end{proof}

\begin{prop}[Uniform reduced hierarchy]
	\label{prop:uniform-hierarchy}
	Let $g$ be measurable on $\cD\times\bR_v^3$, and set $\eta_g:=\widetilde{\mathfrak Q}_{\mathrm N}(g)$.
	There is $\eta_0>0$, depending only on the structural data,
	such that if
	$\eta_g\le\eta_0$ and $0<\delta\le\delta_{\rm dif}$, then the
	mass-normalized spectral solutions $f_N$ of \eqref{eq:spectral-problem} and
	the boundary scalars $\phi_N$ of \eqref{eq:spectral-boundary-data} satisfy
	\begin{equation}
		\label{eq:spectral-bound}
		\sup_N\left[
		\widetilde{\mathfrak Q}_{\mathrm N}(f_N)
		+\|\phi_N\|_{H^{\mathrm N}(\partial\cD)}
		\right] \le C\delta . \end{equation}
\end{prop}

\begin{proof}
	Let $\cB_N:=\cB_{\rm bs}(f_N)$ be the strong coordinate \eqref{eq:base-strong} of $f_N$, as in Lemma~\ref{lem:fixedN-base-estimates}. The hypotheses of Corollary~\ref{cor:spectral-uniform-estimate} and \S\ref{app:fixedN-realization-induction} hold at each cutoff $N$.
	Decrease $\eta_0$, if necessary, so that $\eta_0\le\min\{1,\eta_1\}$, with $\eta_1$ of Lemma~\ref{lem:coefficient-transmission}; then
	Lemma~\ref{lem:coefficient-transmission} gives
	\eqref{eq:construction-coefficient-regularity} whenever $\eta_g\le\eta_0$.
	This implies the coefficient smallness in
	\eqref{eq:diffuse-graph-smallness} and \eqref{eq:KS-coefficient} with $M=1$.  Since
	$0<\delta\le\delta_{\rm dif}$, the remaining smallness hypothesis of
	Corollary~\ref{cor:spectral-uniform-estimate} is verified.  Hence, at every
	$N$, Lemma~\ref{lem:fixedN-complete-jets} gives
	\eqref{eq:fixedN-jet-fullgraph} for $|I|\le\mathrm N$ and
	Lemma~\ref{lem:fixedN-base-estimates} gives $\cB_N<\infty$.  For the background factors, Lemma~\ref{lem:background-finite-p} at $\eta=\eta_g$ gives \eqref{eq:Dg-rungs}.

	Throughout, $\mathrm N$ is the fixed differentiation order and $N$ is the spectral cutoff. We use the exponents of \eqref{eq:dyn-own-order-exponents}. Although the fixed-cutoff construction gives finite complete-word rows, their bounds must be uniform in $N$. We obtain these bounds by Lemma~\ref{lem:ladder-closure} before estimating the boundary scalar and $\cB_N$.
	
	\medskip
	\noindent\emph{Step 1. The row of order zero.}
	Corollary~\ref{cor:spectral-uniform-estimate}, through
	\eqref{eq:spectral-full-diffuse-estimate}, gives
	\begin{equation}
		\label{eq:uniform-U0}
		U_0(f_N)\le C\delta.
	\end{equation}
	For each fixed finite $L$, the incoming datum $H_N=P_{w,N}(f_N)_{\gamma^+}+h_w$ of \eqref{eq:spectral-boundary-data} satisfies
	\[
	\|w^LH_N\|_{L^2(\gamma^-)}
	\lesssim_L \|(f_N)_{\gamma^+}\|_{L^2(\gamma^+)}+\delta
	\lesssim_L \delta.
	\]
	The first bound uses \eqref{eq:wall-J-weighted-operator-bound} and \eqref{eq:boundary-gaussian-polynomial-envelope}; the second follows from \eqref{eq:uniform-U0}. The bulk source vanishes, so the weighted inflow estimate \eqref{eq:weighted-inflow-estimate} at $L=\ell_0$ gives $W_0(f_N)\le C\delta$. Thus
	\begin{equation}
		\label{eq:uniform-base}
		U_0(f_N)+W_0(f_N)\le C\delta .
	\end{equation}
	
	\medskip
	\noindent\emph{Step 2. The row estimate.}
	For $0\le m\le\mathrm N$, put
	\begin{equation}
		\label{eq:hierarchy-rows}
		\cR_m:=W_m(f_N),
		\qquad
		\cN_m:=\sum_{k=0}^{m}\cR_k,
		\qquad
		\cQ:=\delta+\eta_g\cB_N.
	\end{equation}
	By \eqref{eq:uniform-base}, $\cR_0\le C\cQ$. Lemma~\ref{lem:fixedN-complete-jets} gives $\cN_{\mathrm N}<\infty$ for each fixed cutoff, while \eqref{eq:unweighted-row-dominated} gives $U_k(f_N)\le C\cR_k$. For $1\le m\le\mathrm N$ we shall prove
	\begin{equation}
		\label{eq:hierarchy-row}
		\cR_m\le C\cQ+C\cN_{m-1}+C\eta_g\cR_m .
	\end{equation}
	Throughout, $L$ denotes the weight of the row being estimated: $L=0$
	for $U_m$ and $L=\ell_m$ for $W_m$.
	
	\emph{(i) Solution derivatives.}
	As in \S\ref{ss:weight-interpolation}, in a row of order $m$ a jet of order $j<m$ in the solution slot is measured in $L_x^{\mathsf s(j\mathrm N/m)}H_\sigma^1$, where $\mathsf s(j\mathrm N/m)=\mathsf e(m/j,p_\flat)$ is the exponent of \eqref{eq:dyn-row-pairing}. In (iii), the jet in the solution slot has order $j\le m-1$ and is measured in $L_x^{\mathsf e((m-1)/j,p_\flat)}H_\sigma^1$, as in a row of order $m-1$.  
	Since $L\le\ell_m$ and $\ell_{i-1}-\ell_i{}=A_\diamond\ge1$ for $1\le i\le\mathrm N$ by \eqref{eq:finite-hierarchy-weight-choice}, the weights of the gradings $\cG^{L}_{m,2}(f_N)$ and
	$\cG^{L+A_\diamond}_{m-1,2}(f_N)$ of \eqref{eq:GN-lower-order} satisfy
	\[
	L+m-|K|\le\ell_{|K|},
	\qquad
	L+A_\diamond+m-1-|K|\le\ell_{|K|},
	\qquad 0\le|K|\le m-1,
	\]
	and the terms of $\cG^{L}_{m,2}(f_N)$ with $|K|=m$ carry the weight $L$.  With
	\eqref{eq:real-weight-Hsigma-monotonicity} this gives
	\begin{equation}
		\label{eq:hierarchy-tops}
		\cG^{L}_{m,2}(f_N)\le C\bigl(\cR_m+\cN_{m-1}\bigr),
		\qquad
		\cG^{L+A_\diamond}_{m-1,2}(f_N)\le C\cN_{m-1}.
	\end{equation}
	Every weight used below is at most $L+A_\diamond$, and
	\[
	L+A_\diamond\le\ell_{m-1}\le\ell_0=A_\flat-3,
	\]
	where the first inequality follows from \eqref{eq:finite-hierarchy-weight-choice} if $L=\ell_m$ and from \eqref{eq:finite-hierarchy-strict-ladder} if $L=0$, the second holds since the weights $\ell_i$ decrease in $i$, and $\ell_0=A_\flat-3$ by \eqref{eq:finite-hierarchy-weight-choice}. Hence
	\eqref{eq:base-slot-control} gives
	$\|w^{L'}f_N\|_{L_x^{p_\flat}H_\sigma^1}\le C\cB_N$ for $0\le L'\le L+A_\diamond$.
	Lemma~\ref{lem:fixedN-complete-jets} gives
	\eqref{eq:pointwise-slice-chain-hypothesis}.  Therefore
	\eqref{eq:ordered-word-coercive-rungs-weighted} applies at $f=f_N$ and
	$p_{\rm c}=p_\flat$, at $(M,\mathrm L)=(m,L)$ if $m\ge2$ and at
	$(M,\mathrm L)=(m-1,L+A_\diamond)$ if $m\ge3$. We claim that, for $0\le j\le m-1$, the following inequalities hold. At $(m,j)=(1,0)$ the exponent in \eqref{eq:hierarchy-slot-deformation} is understood to be $2$:
	\begin{align}
		\label{eq:hierarchy-slot}
		\sum_{|J|=j}\|w^{L}\mathfrak D_Jf_N\|
		_{L_x^{\mathsf e(m/j,p_\flat)}H_\sigma^1}
		&\le C\bigl(\cB_N+\cR_m+\cN_{m-1}\bigr),\\
		\label{eq:hierarchy-slot-deformation}
		\sum_{|J|=j}\|w^{L+A_\diamond}\mathfrak D_Jf_N\|
		_{L_x^{\mathsf e((m-1)/j,p_\flat)}H_\sigma^1}
		&\le C\bigl(\cB_N+\cN_{m-1}\bigr).
	\end{align}
	For $1\le j\le m-1$ in \eqref{eq:hierarchy-slot} we apply \eqref{eq:ordered-word-coercive-rungs-weighted} at $(M,\mathrm L)=(m,L)$, bound $\|w^Lf_N\|_{L_x^{p_\flat}H_\sigma^1}$ by $C\cB_N$ as above and $\cG^{L}_{m,2}(f_N)$ by the first bound of \eqref{eq:hierarchy-tops}, and use Young's inequality: for $|J|=j$,
	\begin{align*}
		\|w^L\mathfrak D_Jf_N\|_{L_x^{\mathsf e(m/j,p_\flat)}H_\sigma^1}
		&\le C\|w^Lf_N\|_{L_x^{p_\flat}H_\sigma^1}^{1-j/m}\cG^{L}_{m,2}(f_N)^{j/m}+C\|w^Lf_N\|_{L_x^{p_\flat}H_\sigma^1}\\
		&\le C\cB_N^{1-j/m}\bigl(\cR_m+\cN_{m-1}\bigr)^{j/m}+C\cB_N
		\le C\bigl(\cB_N+\cR_m+\cN_{m-1}\bigr).
	\end{align*}
	Summing over the words $J$ with $|J|=j$ gives \eqref{eq:hierarchy-slot}.  For $1\le j\le m-2$, the same argument at $(M,\mathrm L)=(m-1,L+A_\diamond)$, with the bound $\|w^{L+A_\diamond}f_N\|_{L_x^{p_\flat}H_\sigma^1}\le C\cB_N$ above and the second bound of \eqref{eq:hierarchy-tops}, gives \eqref{eq:hierarchy-slot-deformation}.  At $j=0$ the claim is
	\eqref{eq:base-slot-control}, since the exponent there is $p_\flat$ (recall $\mathsf s(0)=p_\flat$), except at $(m,j)=(1,0)$ in \eqref{eq:hierarchy-slot-deformation}, which falls under the next case. At $j=m-1$ in \eqref{eq:hierarchy-slot-deformation} the solution exponent is $2$. The second inequality of \eqref{eq:hierarchy-tops} applies when $m\ge2$; for $m=1$ this is the weighted order-zero estimate $W_0(f_N)$, since $L+A_\diamond\le\ell_0$.
	
	\emph{(ii) Differentiated-background collision terms.}
	By \eqref{eq:weak-collision-word} the differentiated-background terms of the
	order-$m$ equation are $N[\mathfrak D_{I_S}g,\mathfrak D_{I_{S^c}}f_N]$ with
	$1\le|S|=k\le m$; write $j:=m-k$, so that $0\le j\le m-1$.  By \eqref{eq:velalg}, which loses no weight, and H\"older's inequality in $x$ with the exponents of \eqref{eq:dyn-row-pairing}, each of them satisfies
	\[
	\bigl\|w^{L}N[\mathfrak D_{I_S}g,\mathfrak D_{I_{S^c}}f_N]\bigr\|_{L_x^2H_\sigma^{-1}}
	\le C\|w^{-2}\mathfrak D_{I_S}g\|_{L_x^{\mathsf g(k)}L_v^2}\,
	\bigl\|w^{L}\mathfrak D_{I_{S^c}}f_N\bigr\|
	_{L_x^{\mathsf s(j\mathrm N/m)}H_\sigma^1} .
	\]
	The background factor is bounded by $C\eta_g$ in \eqref{eq:Dg-rungs}, while \eqref{eq:hierarchy-slot} gives
	\[
	\|w^{L}\mathfrak D_{I_{S^c}}f_N\|_{L_x^{\mathsf s(j\mathrm N/m)}H_\sigma^1}\le C(\cB_N+\cR_m+\cN_{m-1}),
	\]
	since $\mathsf s(j\mathrm N/m)=\mathsf e(m/j,p_\flat)$ by (i).  Hence
	\begin{equation}
		\label{eq:hierarchy-coulomb}
		\sup_{|I|=m}\ \sum_{\varnothing\ne S\subset[m]}
		\bigl\|w^{L}N[\mathfrak D_{I_S}g,\mathfrak D_{I_{S^c}}f_N]
		\bigr\|_{L_x^2H_\sigma^{-1}}
		\le C\eta_g\bigl(\cB_N+\cR_m+\cN_{m-1}\bigr).
	\end{equation}
	
	\emph{(iii) Deformation remainder.}
	Every term of $\mathfrak R_I[g,f_N]$ has the form
	\eqref{eq:deformation-remainder-structure} with $|J_1|=k$, $|J_2|=j$ and
	$k+j\le m-1$. The estimate \eqref{eq:deformation-weighted-output-spatial} at $\ell=L$ bounds such a term in $w^LL_x^2H_\sigma^{-1}$, with the jet $\mathfrak D_{J_2}f_N$ in the solution slot at the weight $L+2m+2\le L+A_\diamond$ (recall $A_\diamond=2\mathrm N+2$ and $m\le\mathrm N$), and with
	any admissible exponent pair.  By \eqref{eq:real-weight-Hsigma-monotonicity} and \eqref{eq:hierarchy-slot-deformation},
	\[
	\|w^{L+2m+2}\mathfrak D_{J_2}f_N\|_{L_x^{\mathsf e((m-1)/j,p_\flat)}H_\sigma^1}\le C(\cB_N+\cN_{m-1}).
	\]
	For $k\ge1$, the pairing is strictly sub-conjugate. Indeed, $k+j\le m-1$ and \eqref{eq:dyn-row-pairing} give
	\[
	\frac1{\mathsf g(k)}+\frac1{\mathsf e((m-1)/j,p_\flat)}
	\le \frac1{\mathsf g(k)}+\frac1{\mathsf e((k+j)/j,p_\flat)}
	\le\frac12-\frac{k}{90\mathrm N} <\frac12.
	\]
	Here $j=0$ means $\mathsf e((m-1)/0,p_\flat)=\mathsf e(k/0,p_\flat)=p_\flat$. When $k=0$, the background is in $L_x^\infty$ and the solution exponent is at least $2$, so H\"older's inequality still applies; at $(m,k,j)=(1,0,0)$ the solution exponent is $2$.  By
	\eqref{eq:Dg-rungs}, $\|w^{-2}\mathfrak D_{J_1}g\|_{L_x^{\mathsf g(k)}L_v^2}\le C\eta_g$ when $k\ge1$, and, by
	\eqref{eq:base-fiber-Lsharp}, $\|w^{-2}g\|_{L_x^\infty L_v^2}\le C\eta_g$
	when $k=0$.  Thus \eqref{eq:deformation-weighted-output-spatial} at $\ell=L$ and $\varpi=\frac14$, together with $\mu^{1/4}\le Cw^{-2}$, gives the first inequality below, and the bounds on the background factor and on the jet stated above give the second:
	\begin{align*}
		\bigl\|w^{L}pN_{\mathsf K_A}[p_1\mathfrak D_{J_1}g,p_2\mathfrak D_{J_2}f_N]\bigr\|_{L_x^2H_\sigma^{-1}}
		&\le C\|w^{-2}\mathfrak D_{J_1}g\|_{L_x^{\mathsf g(k)}L_v^2}\|w^{L+2m+2}\mathfrak D_{J_2}f_N\|_{L_x^{\mathsf e((m-1)/j,p_\flat)}H_\sigma^1}\\
		&\le C\eta_g\bigl(\cB_N+\cN_{m-1}\bigr).
	\end{align*}
	The number of subword decompositions and the polynomial coefficients depend only on the fixed differentiation order and the complete fields, not on the spectral cutoff. Consequently,
	\begin{equation}
		\label{eq:hierarchy-deformation}
		\sup_{|I|=m}\bigl\|w^{L}\mathfrak R_I[g,f_N]\bigr\|_{L_x^2H_\sigma^{-1}}
		\le C\eta_g\bigl(\cB_N+\cN_{m-1}\bigr).
	\end{equation}
	
	\emph{(iv) Lower-order commutators and boundary data.}
	The transport, local and compact terms of
	\eqref{eq:weak-triangular-commutator} carry jets of complete order strictly below $m$ in their solution slots, and each raises the weight by at most
	$A_\diamond$, so that the raised weight is at most $L+A_\diamond\le\ell_{m-1}$ by (i).  Hence
	\eqref{eq:ordered-transport-weight-integrated},
	\eqref{eq:local-L0-complete-bound}, and
	\eqref{eq:iterated-compact-full-operator-bound} give, for $|I|=m$,
	\[
	\|w^L[\mathfrak D_I,\cA_\mu]f_N\|_{L_x^2H_\sigma^{-1}}
	+\sum_{|J|<m}\bigl(\|w^L\cV_{I,J}q_{J,N}\|_{L_x^2H_\sigma^{-1}}
	+\|w^L\mathscr K_{I,J}q_{J,N}\|_{L_x^2H_\sigma^{-1}}\bigr)
	\le C\cN_{m-1},
	\]
	where $\mathscr K_{I,J}$ is taken to be zero when $J$ is not an ordered subword of $I$.  For $|I|=m$, the terms $b_{I,N}$ and $P_{w,N}(q_{I,N})_{\gamma^+}$ of \eqref{eq:spectral-triangular-law} at $f=f_N$ and the bulk mass of $q_{I,N}$ satisfy, by \eqref{eq:spectral-triangular-bound}, \eqref{eq:wall-J-weighted-operator-bound} at $\cJ=J_N$ and $L=\ell_m$, and \eqref{eq:complete-derivative-mass-sharp} at $f=f_N$,
	\[
	\|w^{\ell_m}b_{I,N}\|_{L^2(\gamma^-)}\le C(\delta+\cN_{m-1}),
	\qquad
	\|w^{\ell_m}P_{w,N}(q_{I,N})_{\gamma^+}\|_{L^2(\gamma^-)}\le CU_m(f_N),
	\qquad
	\Bigl|\iint q_{I,N}\sqrt\mu\Bigr|\le C\cN_{m-1}.
	\]
	Here $I_\ast$ is as in Proposition~\ref{prop:complete-derivative-mass-hierarchy}, which applies at $f=f_N$ since $q_{I,N},q_{I_\ast,N}\in\cH_{\kin,\sigma}\subset L^2(\cD\times\bR_v^3)$ by \eqref{eq:fixedN-jet-fullgraph} and \eqref{260902137}, and the last bound uses
	\[
	\|q_{I_\ast,N}\|_{L^2}\le C\|q_{I_\ast,N}\|_{\cH_{\kin,\sigma}}\le CU_{m-1}(f_N),
	\]
	again by \eqref{260902137}.
	
	\emph{(v) The order-$m$ estimate.}
	Apply \eqref{eq:spectral-full-diffuse-estimate} to the order-$m$ equation
	\eqref{eq:weak-differentiated-equation}, with \eqref{eq:hierarchy-coulomb}
	and \eqref{eq:hierarchy-deformation} at $L=0$ and the bounds of
	(iv), to obtain
	\begin{align}\label{261015818}
		U_m(f_N)\le C\delta+C\cN_{m-1}
		+C\eta_g\bigl(\cB_N+\cR_m+\cN_{m-1}\bigr).
	\end{align}
	The boundary law \eqref{eq:spectral-triangular-law}
	with the bounds on $\|w^{\ell_m}b_{I,N}\|_{L^2(\gamma^-)}$ and $\|w^{\ell_m}P_{w,N}(q_{I,N})_{\gamma^+}\|_{L^2(\gamma^-)}$ in (iv) gives
	\[
	\sup_{|I|=m}\|w^{\ell_m}(q_{I,N})_{\gamma^-}\|_{L^2(\gamma^-)}
	\le C\delta+CU_m(f_N)+C\cN_{m-1}.
	\]
	Applying the weighted inflow estimate \eqref{eq:weighted-inflow-estimate} at $L=\ell_m$, with the source bounds \eqref{eq:hierarchy-coulomb}--\eqref{eq:hierarchy-deformation} and the incoming trace bound above, gives
	\[
	W_m(f_N)\le C\delta+CU_m(f_N)+C\cN_{m-1}
	+C\eta_g\bigl(\cB_N+\cR_m+\cN_{m-1}\bigr) .
	\]
	Together with \eqref{261015818}, this gives
	\[
	\cR_m=W_m(f_N)\le C\delta+C\cN_{m-1}+C\eta_g\bigl(\cB_N+\cR_m+\cN_{m-1}\bigr) .
	\]
	As $\eta_g\le1$ and $\cQ=\delta+\eta_g\cB_N$, this is \eqref{eq:hierarchy-row}.
	
	\medskip
	\noindent\emph{Step 3. Closing the induction.}
	The constants in \eqref{eq:hierarchy-row} and in the strong-coordinate estimate \eqref{eq:base-absorption-row} below are independent of $\eta_0$. Indeed, apart from their displayed factors $\eta_g$, the estimates use only $\eta_g\le\min\{1,\eta_1\}$, through \eqref{eq:construction-coefficient-regularity} and $\eta_g\le1$.
	Choose $\eta_0$ sufficiently small, depending only on the structural
	data, to apply Lemma~\ref{lem:ladder-closure} to
	\eqref{eq:hierarchy-rows} with $M=\mathrm N$ and $\theta=\eta_g$.
	The zeroth- and higher-order estimates, together with $\cN_{\mathrm N}<\infty$, give its
	hypotheses; here $\cQ=\delta+\eta_g\cB_N<\infty$ because $\cB_N<\infty$.
	Consequently,
	\begin{equation}
		\label{eq:uniform-high}
		\cN_{\mathrm N}\le C\bigl(\delta+\eta_g\cB_N\bigr).
	\end{equation}
	
	\medskip
	\noindent\emph{Step 4. The boundary scalar and $\cB_N$.}
	By \eqref{eq:fixedN-jet-fullgraph},
	Lemma~\ref{cor:outgoing-moment-H3}\,(\ref{it:outgoing-moment-H-stationary})
	applies at $m=\mathrm N$ and yields
	\begin{equation}
		\label{eq:spectral-scalar-row}
		\|\phi_N\|_{H^{\mathrm N}(\partial\cD)}
		\le C\sum_{m=0}^{\mathrm N}U_m(f_N)
		\le C\delta+C\eta_g\cB_N,
	\end{equation}
	where the last inequality follows from \eqref{eq:uniform-base}, \eqref{eq:unweighted-row-dominated} and \eqref{eq:uniform-high}. Moreover, Lemma~\ref{lem:word-boundary-spectral-chart} gives
	\[
	\|J_N\phi_N\|_{H^{\mathrm N}(\partial\cD)}\le C\|\phi_N\|_{H^{\mathrm N}(\partial\cD)} .
	\]
	
	To estimate $\cB_N$, apply Lemma~\ref{lem:incoming-datum-bound} with $\psi=J_N\phi_N$ and $L=A_\flat+B^{\st}_{\rm in}$; here $H^{\mathrm N}(\partial\cD)\subset H^3(\partial\cD)$. This gives
	\begin{equation}
		\label{eq:spectral-incoming-datum}
		\|w^{A_\flat+B^{\st}_{\rm in}}H_N\|_{
			(C_x^{0,1}\cap C_v^{1,1})(\gamma^-)}
		\le C\bigl(\delta+\|J_N\phi_N\|_{H^{\mathrm N}(\partial\cD)}\bigr).
	\end{equation}  By \eqref{eq:base-regularity} at $U=f_N$ and $H=H_N$, with $\ell=A_\flat$,
	$M=1$ and $p=p_\flat$, together with \eqref{eq:spectral-incoming-datum}, the bound of Lemma~\ref{lem:word-boundary-spectral-chart} and \eqref{eq:spectral-scalar-row}, we obtain
	\begin{equation}
		\label{eq:base-absorption-row}
		\cB_N
		\le C\delta+C\eta_g\cB_N.
	\end{equation}
	Decreasing $\eta_0$ once more, we absorb the last term in
	\eqref{eq:base-absorption-row} and obtain $\cB_N\le C\delta$.
	Together with \eqref{eq:uniform-high} and
	\eqref{eq:unweighted-row-dominated}, this yields 
	$$
	U_m(f_N)+W_m(f_N)\le C\delta\qquad \text{for }1\le m\le\mathrm N,
	$$ 
	since $\eta_g\le1$.  Recall from \eqref{eq:reduced-Q3} that \[
	\widetilde{\mathfrak Q}_{\mathrm N}(f_N)=\cB_N+\sum_{m=0}^{\mathrm N}W_m(f_N) .
	\] Combining
	\eqref{eq:uniform-base}, \eqref{eq:uniform-high},
	\eqref{eq:spectral-scalar-row} and \eqref{eq:base-absorption-row}, we obtain
	\eqref{eq:spectral-bound}.
\end{proof}

We pass to the limit $N\to\infty$ in \eqref{eq:spectral-bound}. Lower semicontinuity will control the reduced norm, while a separate strong-convergence argument identifies the limiting solution at the fixed background $g$. The following lower-semicontinuity lemma also covers the evolutionary case.

\begin{lemma}[Lower semicontinuity of the strong coordinate]
	\label{lem:strong-coordinate-lsc}
	Let $1<p<\infty$, $A\ge0$ and $\alpha\in(0,1)$.
	\begin{enumerate}
		\item\label{it:strong-coordinate-lsc-st}
		If $u_n\to u$ in $L^2_{\rm loc}(\cD\times\bR_v^3)$, then
		\begin{equation}
			\label{eq:strong-coordinate-lsc}
			\begin{aligned}
				\|w^Au\|_{L^\infty}&\le\liminf_{n\to\infty}\|w^Au_n\|_{L^\infty},\\
				[w^Au]_{C^\alpha_{\kin}}&\le\liminf_{n\to\infty}[w^Au_n]_{C^\alpha_{\kin}},\\
				\|w^A\nabla_vu\|_{L_x^pL_v^2}&\le\liminf_{n\to\infty}\|w^A\nabla_vu_n\|_{L_x^pL_v^2},
			\end{aligned}
		\end{equation}
		where all norms are taken over $\cD\times\bR_v^3$ and the kinetic H\"older seminorm is taken with respect to
		$\mathbf d_{\st}$ of \eqref{eq:stationary-kinetic-distance-main}.
		\item\label{it:strong-coordinate-lsc-dyn}
		Let $0<T<\infty$.  If $u_n\to u$ in $L^2_{\rm loc}(\cD_T\times\bR_v^3)$, then
		\begin{equation}
			\label{eq:strong-coordinate-lsc-dyn}
			\begin{aligned}
				\|w^Au\|_{L^\infty}&\le\liminf_{n\to\infty}\|w^Au_n\|_{L^\infty},\\
				[w^Au]_{C^\alpha_{\kin}}&\le\liminf_{n\to\infty}[w^Au_n]_{C^\alpha_{\kin}},\\
				\|w^A\nabla_vu\|_{L_t^\infty L_x^pL_v^2}&\le\liminf_{n\to\infty}\|w^A\nabla_vu_n\|_{L_t^\infty L_x^pL_v^2},
			\end{aligned}
		\end{equation}
		where all norms are taken over $\cD_T\times\bR_v^3$ and the kinetic H\"older seminorm is taken with respect to
		$\mathbf d$ of \eqref{kineticdistance}.
	\end{enumerate}
	In either case, the sum of the three quantities is lower semicontinuous along $u_n\to u$.
\end{lemma}

\begin{proof}
	We prove (\ref{it:strong-coordinate-lsc-dyn}); the stationary case follows by omitting the time variable. For each inequality, pass to a subsequence realizing the finite lower limit. Local $L^2$ convergence permits a further subsequence with $u_n\to u$ almost everywhere.  The first inequality
	then follows from the almost everywhere convergence, and the third from the weak-$*$ lower semicontinuity of the
	norm of $L_t^\infty L_x^pL_v^2=(L_t^1L_x^{p'}L_v^2)^*$, since $w^A\nabla_vu_n\to w^A\nabla_vu$ in distributions.
	
	For the second inequality, let $K:=\liminf_{n\to\infty}[w^Au_n]_{C^\alpha_{\kin}}$.  Letting $n\to\infty$ in the H\"older bound of $w^Au_n$
	gives \[
	|w^Au(\mathsf Z)-w^Au(\mathsf Z')|\le K\,\mathbf d(\mathsf Z,\mathsf Z')^\alpha
	\quad\text{for all }\mathsf Z,\mathsf Z'\text{ outside a null set} .
	\]  Since $\mathbf d$ induces the Euclidean topology, $w^Au$ has a continuous representative
	satisfying the same bound for all $\mathsf Z,\mathsf Z'$, that is, $[w^Au]_{C^\alpha_{\kin}}\le K$.
\end{proof}

\begin{corollary}[Lower semicontinuity of $\cB_{\rm bs}$]
	\label{lem:base-lsc}
	Let $u_n\to u$ in $L^2_{\rm loc}(\cD\times\bR^3)$ and assume
	$\liminf_n\cB_{\rm bs}(u_n)<\infty$.  Then
	\begin{equation}
		\cB_{\rm bs}(u)\le\liminf_{n\to\infty}\cB_{\rm bs}(u_n).
		\label{eq:base-lsc}
	\end{equation}
\end{corollary}

\begin{proof}
	Apply Lemma~\ref{lem:strong-coordinate-lsc}\,(\ref{it:strong-coordinate-lsc-st}) with $A=A_\flat$, $p=p_\flat$ and $\alpha=\alpha_{\rm bs}$. Its three quantities are precisely the terms of
	\eqref{eq:base-strong}, so the assertion of Lemma~\ref{lem:strong-coordinate-lsc} on their sum gives \eqref{eq:base-lsc}.
\end{proof}

\begin{lemma}[Closure of the reduced norm]
	\label{lem:reduced-norm-closure}
	Let $u_n\to u$ in $L^2_{\rm loc}(\cD\times\bR^3)$ and assume
	$\liminf_n\widetilde{\mathfrak Q}_{\mathrm N}(u_n)<\infty$.  Then
	\begin{align}
		\label{eq:reduced-norm-closure}
		\begin{gathered}
			\mathfrak D_Iu\in\cH_{\kin,\sigma}\ \ (|I|\le\mathrm N),
			\qquad
			w^{\ell_m}\mathfrak D_Iu\in\cH_{\kin,\sigma}\ \ (I\in\mathfrak A^m,\ m\le\mathrm N),
			\\
			\widetilde{\mathfrak Q}_{\mathrm N}(u)
			\le\liminf_{n\to\infty}\widetilde{\mathfrak Q}_{\mathrm N}(u_n).
		\end{gathered}
	\end{align}
\end{lemma}

\begin{proof}
	Choose a subsequence $u_{n_k}$ with
	$\widetilde{\mathfrak Q}_{\mathrm N}(u_{n_k})<\infty$ for every $k$ and
	\[
	\widetilde{\mathfrak Q}_{\mathrm N}(u_{n_k})
	\longrightarrow
	\liminf_{n\to\infty}\widetilde{\mathfrak Q}_{\mathrm N}(u_n).
	\]
	By \eqref{eq:unweighted-row-dominated}
	and \eqref{eq:reduced-Q3},
	\begin{equation}
		\label{eq:reduced-closure-rows}
		\sup_k\Bigl(\sum_{m=0}^{\mathrm N}U_m(u_{n_k})
		+\sum_{m=0}^{\mathrm N}W_m(u_{n_k})\Bigr)
		\le C\sup_k\widetilde{\mathfrak Q}_{\mathrm N}(u_{n_k})<\infty.
	\end{equation}
	For an ordered complete word $I$ with $|I|\le\mathrm N$, put $q_{I,k}:=\mathfrak D_Iu_{n_k}$, with $\mathfrak D_\varnothing={\rm Id}$.  For $\psi\in C_c^\infty(\cD\times\bR^3)$, the ordered adjoint $\mathfrak D_I^\dagger\psi$ of \eqref{eq:ordered-complete-formal-adjoint} satisfies $\mathfrak D_I^\dagger\psi\in L^\infty(\cD\times\bR^3)$ and is compactly supported in $\cD\times\bR^3$, so the local $L^2$ convergence gives
	\[
	\langle q_{I,k},\psi\rangle=\langle u_{n_k},\mathfrak D_I^\dagger\psi\rangle
	\longrightarrow\langle u,\mathfrak D_I^\dagger\psi\rangle=\langle\mathfrak D_Iu,\psi\rangle .
	\]
	Hence $q_{I,k}\to\mathfrak D_Iu$ and $w^{\ell_{|I|}}q_{I,k}\to w^{\ell_{|I|}}\mathfrak D_Iu$ in $\cD'(\cD\times\bR^3)$, and, by \eqref{eq:reduced-closure-rows},
	\[
	\sup_k\bigl(\|q_{I,k}\|_{\cH_{\kin,\sigma}}+\|w^{\ell_{|I|}}q_{I,k}\|_{\cH_{\kin,\sigma}}\bigr)<\infty .
	\]Lemma~\ref{lem:dynamic-weak-graph-bulk-identification}\,(1) therefore gives $\mathfrak D_Iu,\,w^{\ell_{|I|}}\mathfrak D_Iu\in\cH_{\kin,\sigma}$, as asserted in \eqref{eq:reduced-norm-closure}, and
	\[
	q_{I,k}\rightharpoonup\mathfrak D_Iu,
	\qquad
	w^{\ell_{|I|}}q_{I,k}\rightharpoonup w^{\ell_{|I|}}\mathfrak D_Iu
	\qquad\text{in }\cH_{\kin,\sigma}.
	\]
	By \eqref{eq:base-lsc} for the strong coordinate $\cB_{\rm bs}$, and by the weak lower semicontinuity of the norm of $\cH_{\kin,\sigma}$ for the weighted rows $W_m$, we have
	\[
	\cB_{\rm bs}(u)\le\liminf_{k\to\infty}\cB_{\rm bs}(u_{n_k}),
	\qquad
	\sum_{m=0}^{\mathrm N}W_m(u)\le\liminf_{k\to\infty}\sum_{m=0}^{\mathrm N}W_m(u_{n_k}) .
	\]
	Adding the two bounds and
	using \eqref{eq:reduced-Q3}, we obtain
	\[
	\widetilde{\mathfrak Q}_{\mathrm N}(u)
	\le\lim_{k\to\infty}\widetilde{\mathfrak Q}_{\mathrm N}(u_{n_k})
	=\liminf_{n\to\infty}\widetilde{\mathfrak Q}_{\mathrm N}(u_n),
	\]
	which proves the last assertion of \eqref{eq:reduced-norm-closure}.
\end{proof}

From here to the end of this subsection, $g$ and $\delta$ satisfy
\eqref{eq:fixed-background-smallness}, with $\eta_\st$ and $\delta_\st$ chosen
so that $\eta_\st\le\eta_0$ and $\delta_\st\le\delta_{\rm dif}$,
$\eta_0$ being that of Proposition~\ref{prop:uniform-hierarchy}.  Since $\eta_0\le\eta_1$ by the proof of Proposition~\ref{prop:uniform-hierarchy}, Lemma~\ref{lem:coefficient-transmission} gives \eqref{eq:construction-coefficient-regularity}, so that \eqref{eq:diffuse-graph-smallness} holds at $g$.  Hence Proposition~\ref{prop:full-diffuse-graph-solvability}, at the compatible data $F=0$, $b=h_w$ and $m=0$, furnishes the unique mass-normalized solution $f$ of Theorem~\ref{thm:fixed-background-construction}, and \eqref{eq:inhomogeneous-full-diffuse-graph-estimate} applies at $g$.  Recall that $f_N$ is the mass-normalized spectral solution of \eqref{eq:spectral-problem} and $\phi_N=\cA_0(f_N)_{\gamma^+}$ (see \eqref{eq:spectral-boundary-data}), and put
$\phi:=\cA_0f_{\gamma^+}$.

\begin{prop}[Strong convergence]
	\label{prop:fixedN-limits-base}
	As $N\to\infty$,
	\begin{equation}
		\label{eq:strong-base-limit}
		\|f_N-f\|_{\cH_{\kin,\sigma}}
		+\|\phi_N-\phi\|_{L^2(\partial\cD)}
		\longrightarrow0 .
	\end{equation}
\end{prop}

\begin{proof}
	The spectral tail $(\mathbf I-J_N)\phi_N$ tends to zero in $L^2(\partial\cD)$. Choose $c_\chi>0$ such that $\chi=1$ on
	$[0,c_\chi]$, with $\chi$ the cutoff of \eqref{eq:boundary-spectral-cutoff-definition}.  The operator $\mathbf I-J_N$ acts on each eigenfunction of $-\Delta_{\partial\cD}$ with eigenvalue $\lambda$ as multiplication by $1-\chi(\lambda/N^2)$, which vanishes for $\lambda\le c_\chi N^2$, so that
	\begin{equation}
		\label{eq:spectral-cutoff-tail-multiplier}
		|1-\chi(\lambda/N^2)|
		\lesssim_{\chi,k} N^{-k}(1+\lambda)^{k/2},
		\qquad \lambda\ge0,
		\qquad k\ge0 .
	\end{equation}
	Hence Lemma~\ref{lem:word-boundary-spectral-chart},
	\eqref{eq:spectral-cutoff-tail-multiplier} at $k=\mathrm N$ and
	\eqref{eq:spectral-bound} give
	\[
	\|(\mathbf I-J_N)\phi_N\|_{L^2(\partial\cD)}
	\lesssim N^{-\mathrm N}\|\phi_N\|_{H^{\mathrm N}(\partial\cD)}
	\longrightarrow0 .
	\]
	
	Set $r_N:=f_N-f$.  Both $f_N$ and $f$ satisfy the bulk equation of
	\eqref{eq:full-diffuse-linear-problem} at
	the common background $g$ with $F=0$ and are mass-normalized, so
	$v\cdot\nabla_xr_N+L_gr_N=0$ in $L^2(\cD;H_\sigma^{-1})$ and $\iint r_N\sqrt\mu=0$.  Subtracting the boundary law of \eqref{eq:full-diffuse-linear-problem} with $b=h_w$ from that of \eqref{eq:spectral-problem} gives
	\[
	(r_N)_{\gamma^-}
	=P_w(r_N)_{\gamma^+}
	+(\mathfrak m_\circ+d_\delta)(J_N-\mathbf I)\phi_N .
	\]%
	Moreover, the Gaussian
	multiplier bound in \eqref{eq:boundary-gaussian-polynomial-envelope} gives
	\[
	\| (\mathfrak m_\circ+d_\delta)(J_N-\mathbf I)\phi_N\|_{L^2(\gamma^-)}
	\lesssim\|(J_N-\mathbf I)\phi_N\|_{L^2(\partial\cD)}\longrightarrow0 .
	\]
	Applying \eqref{eq:inhomogeneous-full-diffuse-graph-estimate} to $r_N$, with zero bulk source and zero mass, yields
	\[
	\|r_N\|_{\cH_{\kin,\sigma}}\le C\|(\mathbf I-J_N)\phi_N\|_{L^2(\partial\cD)}\le C\delta N^{-\mathrm N}.
	\] Since $\phi_N-\phi=\cA_0(r_N)_{\gamma^+}$ by \eqref{eq:spectral-boundary-data} and the definition of $\phi$, \eqref{eq:outgoing-moment-flux-bound} and \eqref{260902137} give
	\[
	\|\phi_N-\phi\|_{L^2(\partial\cD)}
	\lesssim\|(r_N)_{\gamma^+}\|_{L^2(\gamma^+)}
	\lesssim\|r_N\|_{\cH_{\kin,\sigma}}.
	\]
	Thus the cutoff error satisfies
	\[
	\|f_N-f\|_{\cH_{\kin,\sigma}}+\|\phi_N-\phi\|_{L^2(\partial\cD)}\le C\delta N^{-\mathrm N},
	\]
	which proves \eqref{eq:strong-base-limit}.
\end{proof}

\begin{proof}[\textbf{Proof of Theorem~\ref{thm:fixed-background-construction}}]
	By \eqref{eq:strong-base-limit} and \eqref{eq:local-H1-H_sigma},
	$f_N\to f$ in $L^2_{\rm loc}(\cD\times\bR^3)$, and \eqref{eq:spectral-bound} gives
	$\sup_N\widetilde{\mathfrak Q}_{\mathrm N}(f_N)\le C\delta$.  Hence
	\eqref{eq:reduced-norm-closure} gives
	$\widetilde{\mathfrak Q}_{\mathrm N}(f)\le C\delta$.
	For the boundary scalar, \eqref{eq:spectral-bound} and \eqref{eq:strong-base-limit} give $\sup_N\|\phi_N\|_{H^{\mathrm N}(\partial\cD)}\le C\delta$ and $\phi_N\to\phi$ in $L^2(\partial\cD)$.  Hence $\phi_N\rightharpoonup\phi$ in $H^{\mathrm N}(\partial\cD)$, and the weak lower semicontinuity of the norm gives
	\[
	\phi\in H^{\mathrm N}(\partial\cD),
	\qquad
	\|\phi\|_{H^{\mathrm N}(\partial\cD)}
	\le\liminf_{N\to\infty}\|\phi_N\|_{H^{\mathrm N}(\partial\cD)}
	\le C\delta .
	\]
	Since $\phi=\cA_0f_{\gamma^+}$, these bounds give \eqref{eq:fixed-background-construction}. The constant is independent of $g$ and $\delta$ by \eqref{eq:spectral-bound}.
\end{proof}

\subsection{The fixed point and the steady state}\label{ss:diffuse-fixed-point}

The estimate \eqref{eq:fixed-background-construction} yields an invariant ball controlled by the reduced norm. The fixed-background map is a contraction on this ball in the kinetic graph metric \eqref{eq:base-graph-metric}.  

\begin{defn}[Boundary-constrained invariant ball]
	\label{def:invariant-ball}
	For $R>0$ set
	\begin{equation}
		\label{eq:invariant-ball}
		\widetilde{\mathbb B}_R(\delta):=
		\left\{u\in\cH_{\kin,\sigma}(\cD\times\bR_v^3):\ \widetilde{\mathfrak Q}_{\mathrm N}(u)\le R\delta,\quad
		u_{\gamma^-}=P_wu_{\gamma^+}+h_w,\quad \iint u\sqrt\mu=0\right\},
	\end{equation}
	and equip it with the metric
	\begin{equation}
		\label{eq:base-graph-metric}
		d_{\cH}(u_1,u_2)
		:=\|u_1-u_2\|_{\cH_{\kin,\sigma}}.
	\end{equation}
\end{defn}

\begin{lemma}[Completeness of the invariant ball]
	\label{lem:ball-closed}
	For fixed $R>0$ and $\delta>0$, $\widetilde{\mathbb B}_R(\delta)$ is complete under the metric
	\eqref{eq:base-graph-metric}.
\end{lemma}

\begin{proof}
	Let $u_n$ be a $d_{\cH}$-Cauchy sequence in $\widetilde{\mathbb B}_R(\delta)$.  Since
	$\cH_{\kin,\sigma}$ is complete,
	$u_n\to u$ in $\cH_{\kin,\sigma}$ for some $u$.  Hence $u_n\to u$ in $L^2_{\rm loc}(\cD\times\bR^3)$ by \eqref{eq:local-H1-H_sigma}.
	Each constraint in \eqref{eq:invariant-ball} passes to the limit.
	For the bound, since $\widetilde{\mathfrak Q}_{\mathrm N}(u_n)\le R\delta$, Lemma~\ref{lem:reduced-norm-closure} gives
	$\widetilde{\mathfrak Q}_{\mathrm N}(u)\le R\delta$.
	
	For the boundary law, \eqref{260902137} gives $(u_n)_{\gamma^\pm}\to u_{\gamma^\pm}$ in
	$L^2(\gamma^\pm)$, and the boundedness
	\eqref{eq:wall-J-operator-bound} of $P_w$ gives
	\[
	h_w=(u_n)_{\gamma^-}-P_w(u_n)_{\gamma^+}
	\longrightarrow u_{\gamma^-}-P_wu_{\gamma^+}
	\qquad\text{in }L^2(\gamma^-),
	\]
	so $u_{\gamma^-}=P_wu_{\gamma^+}+h_w$.
	
	For the mass, $u_n\to u$ in $L^2(\cD\times\bR_v^3)$ by \eqref{260902137}, so
	$\iint u\sqrt\mu=\lim_n\iint u_n\sqrt\mu=0$.
	Thus $u\in\widetilde{\mathbb B}_R(\delta)$ and $d_\cH(u_n,u)\to0$.
\end{proof}

\begin{defn}[Reduced fixed-background map]
	\label{def:map}
	Whenever \eqref{eq:fixed-background-smallness} holds, let
	$\cS_\st(g)\in\cH_{\kin,\sigma}(\cD\times\bR_v^3)$ be the
	unique solution of
	\begin{equation}
		\label{eq:map}
		v\cdot\nabla_x \cS_\st(g)+L_g\cS_\st(g)=0,
		\qquad
		\cS_\st(g)_{\gamma^-}=P_w\cS_\st(g)_{\gamma^+}+h_w,
		\qquad
		\iint\cS_\st(g)\sqrt\mu=0.
	\end{equation}
\end{defn}

\begin{proof}[Well-definedness]
	By the paragraph before Proposition~\ref{prop:fixedN-limits-base}, \eqref{eq:fixed-background-smallness} implies \eqref{eq:diffuse-graph-smallness}, and Proposition~\ref{prop:full-diffuse-graph-solvability} gives the unique mass-normalized solution $\cS_\st(g)\in\cH_{\kin,\sigma}(\cD\times\bR^3_v)$ of \eqref{eq:map}.  This is the solution $f$ of Theorem~\ref{thm:fixed-background-construction}, which then gives \eqref{eq:fixed-background-construction}.
\end{proof}

\begin{lemma}[Self-map on the invariant ball]
	\label{lem:self-map}
	Let $R_\st$ be the constant $C$ of \eqref{eq:fixed-background-construction}, which is independent of $g$ and $\delta$.  Decrease $\delta_\st$, if
	necessary, so that
	\begin{equation}
		\label{eq:self-map-radius}
		R_\st\delta_\st\le\eta_\st.
	\end{equation}
	Then, for $0<\delta\le\delta_\st$,
	\(
	\cS_\st:\widetilde{\mathbb B}_{R_\st}(\delta)
	\longrightarrow
	\widetilde{\mathbb B}_{R_\st}(\delta)
	\).
	Moreover,
	\(
	\cS_\st(0)\in\widetilde{\mathbb B}_{R_\st}(\delta)\).
\end{lemma}

\begin{proof}
	For $g\in\widetilde{\mathbb B}_{R_\st}(\delta)$ we have
	$\widetilde{\mathfrak Q}_{\mathrm N}(g)\le R_\st\delta\le\eta_\st$.
	Theorem~\ref{thm:fixed-background-construction} and
	\eqref{eq:map} therefore give
	\[
	\widetilde{\mathfrak Q}_{\mathrm N}(\cS_\st(g))
	+\|\cA_0\cS_\st(g)_{\gamma^+}\|_{H^{\mathrm N}}
	\le R_\st\delta
	\]
	together with the required boundary law and mass constraint.
	Thus the ball is invariant.  The same estimate applies at $g=0$;
	its image belongs to the ball, which is therefore nonempty.
\end{proof}

\begin{prop}[Strict contraction in the invariant ball]
	\label{prop:strict-contraction}
	After decreasing $\delta_\st$, if necessary, whenever
	\[
	0<\delta\le\delta_\st,
	\qquad
	g_1,g_2\in\widetilde{\mathbb B}_{R_\st}(\delta),
	\]
	we have
	\begin{equation}
		\label{eq:strict-contraction}
		d_{\cH}(\cS_\st(g_1),\cS_\st(g_2))
		\le\frac12d_{\cH}(g_1,g_2).
	\end{equation}
\end{prop}

\begin{proof}
	Set $u_i=\cS_\st(g_i)$, $d=g_1-g_2$ and $z=u_1-u_2$.
	Since $L_g=L-N[g,\cdot]$, bilinearity gives
	\begin{equation*}
		v\cdot\nabla_xz+L_{g_1}z=N[d,u_2],
		\qquad z_{\gamma^-}=P_wz_{\gamma^+},
		\qquad \iint z\sqrt\mu=0.
	\end{equation*}
	Since $g_1$ satisfies \eqref{eq:fixed-background-smallness}, the homogeneous estimate \eqref{eq:homogeneous-diffuse-coercivity} applies.  Its source belongs to $L_x^2H_\sigma^{-1}$ by
	Lemma~\ref{260715201} and Lemma~\ref{lem:velalg}.
	The bilinear bound \eqref{eq:velalg-poly}, H\"older on the bounded
	domain $\cD$ with the strictly sub-conjugate pair
	\eqref{eq:dyn-own-order-endpoints}, and the mixed gain
	\eqref{eq:graph-to-mixed-integrability} yield
	\[
	\begin{aligned}
		\|z\|_{\cH_{\kin,\sigma}}
		&\lesssim \|N[d,u_2]\|_{L_x^2H_\sigma^{-1}}\\
		&\lesssim \|w^{-2}d\|_{L_x^{q_\flat}L_v^2}
		\|u_2\|_{L_x^{p_\flat}H_\sigma^1}
		\lesssim \delta\|d\|_{\cH_{\kin,\sigma}}.
	\end{aligned}
	\]
	In the last inequality we used \eqref{eq:base-slot-control},
	\eqref{eq:reduced-Q3} and \eqref{eq:fixed-background-construction}.
	The constants are independent of $g_1,g_2$ and $\delta$. Choosing $\delta_\st$ smaller gives \eqref{eq:strict-contraction}.
\end{proof}

\begin{thm}[Reduced-norm nonlinear solution]
	\label{thm:reduced-fixed-point}
	For $0<\delta\le\delta_\st$, with $\delta_\st$ as decreased in
	Proposition~\ref{prop:strict-contraction}, the nonlinear problem
	\eqref{eq:stationary-full-diffuse-problem} has a mass-normalized solution
	$f_\st$ satisfying
	\begin{equation}
		\label{eq:fixed-point-bound}
		\widetilde{\mathfrak Q}_{\mathrm N}(f_\st)
		+\|\cA_0(f_\st)_{\gamma^+}\|_{H^{\mathrm N}(\partial\cD)}
		\le C\delta
	\end{equation}
	and
	\begin{equation}
		\label{eq:fixed-point-invariant-bound}
		\widetilde{\mathfrak Q}_{\mathrm N}(f_\st)
		\le\eta_\st ,
	\end{equation}
	with $\eta_\st$ of
	Theorem~\ref{thm:fixed-background-construction}.
\end{thm}

\begin{proof}
	By Lemma~\ref{lem:ball-closed}, $(\widetilde{\mathbb B}_{R_\st}(\delta),d_\cH)$ is a
	complete metric space.  By Lemma~\ref{lem:self-map}, it contains $\cS_\st(0)$
	and $\cS_\st$ maps it into itself.  By \eqref{eq:strict-contraction} and the
	Banach fixed point theorem, there is
	$f_\st\in\widetilde{\mathbb B}_{R_\st}(\delta)$ with $\cS_\st(f_\st)=f_\st$.  By \eqref{eq:map} with $g=f_\st$, $f_\st$ solves
	\eqref{eq:stationary-full-diffuse-problem} and $\iint f_\st\sqrt\mu=0$.  By
	\eqref{eq:invariant-ball} and \eqref{eq:self-map-radius}, $f_\st$ satisfies
	\eqref{eq:fixed-point-invariant-bound}.
	Hence \eqref{eq:fixed-background-smallness} holds with $g=f_\st$, and
	Theorem~\ref{thm:fixed-background-construction}, applied to $\cS_\st(f_\st)=f_\st$, gives
	\eqref{eq:fixed-point-bound}.
\end{proof}

Since \eqref{eq:fixed-point-bound} gives $\|\cA_0(f_\st)_{\gamma^+}\|_{H^{\mathrm N}(\partial\cD)}\le C\delta$, Lemma~\ref{lem:incoming-datum-bound}, with $H^{\mathrm N}(\partial\cD)\subset H^3(\partial\cD)$, shows that the incoming datum of $f_\st$ satisfies, for every finite $L\ge0$, $\|w^L(f_\st)_{\gamma^-}\|_{(C_x^{0,1}\cap C_v^{1,1})(\gamma^-)}\le C_L\delta$, and Proposition~\ref{prop:base-regularity} then bounds $f_\st$ at every finite weight.

\begin{prop}[Finite-weight regularity of the stationary state]
	\label{prop:stationary-finite-weight-regularity}
	Let $0<\delta\le\delta_\st$, let $f_\st$ be the mass-normalized stationary solution
	furnished by Theorem~\ref{thm:reduced-fixed-point}, and set
	\begin{equation}
		\phi_s:=\cA_0(f_\st)_{\gamma^+},
		\qquad
		H_s:=(f_\st)_{\gamma^-}=\mathfrak m_{\Theta_w}\phi_s+h_w.
		\label{eq:stationary-incoming-datum}
	\end{equation}
	For every fixed $\ell\ge0$ and $1<p<\infty$,
	\begin{align}
		\|w^\ell f_\st\|_{\cH_{\kin,\sigma}}
		+\|w^\ell \nabla_vf_\st\|_{L_x^{p}L_v^2}
		+\|w^\ell f_\st\|_{L_x^{p}H_\sigma^1}
		&\lesssim_{\ell,p} \delta,
		\label{eq:stationary-finite-weight-graph}\\
		\|w^\ell f_\st\|_{L^\infty_{x,v}}
		+[w^\ell f_\st]_{C^{\alpha_{\rm bs}}_{\kin}}
		&\lesssim_\ell \delta,
		\label{eq:stationary-finite-weight-pointwise}\\
		\|w^\ell H_s\|_{(C_x^{0,1}\cap C_v^{1,1})(\gamma^-)}
		&\lesssim_\ell \delta.
		\label{eq:stationary-incoming-regularity}
	\end{align}
\end{prop}

\begin{proof}
	The bound \eqref{eq:stationary-incoming-regularity} follows from the outgoing moment estimate.
	Lemma~\ref{lem:incoming-datum-bound} at $\psi=\phi_s$, with $H^{\mathrm N}(\partial\cD)\subset H^3(\partial\cD)$ ($\mathrm N\ge4$) and $\|\phi_s\|_{H^{\mathrm N}(\partial\cD)}\le C\delta$ by \eqref{eq:fixed-point-bound}, gives, for every finite $L\ge0$,
	\begin{equation}
		\|w^L H_s\|_{(C_x^{0,1}\cap C_v^{1,1})(\gamma^-)}
		\lesssim_L\delta ,
		\label{eq:stat-inflow-reg}
	\end{equation}
	which is \eqref{eq:stationary-incoming-regularity}.
	
	For the remaining estimates, we use the prescribed-inflow equation
	\begin{equation}
		v\cdot\nabla_x f_\st+L_{f_\st}f_\st=0,
		\qquad (f_\st)_{\gamma^-}=H_s.
		\label{eq:stationary-prescribed-inflow}
	\end{equation}
	By \eqref{eq:fixed-point-invariant-bound}, $\eta_\st\le\eta_0\le\eta_1$ and
	\eqref{eq:construction-coefficient-regularity},
	\eqref{eq:KS-coefficient} holds with $g=f_\st$ and $M=1$.
	We apply Proposition~\ref{prop:base-regularity}, at the given $p$, to \eqref{eq:stationary-prescribed-inflow} at the weights $\ell$ and $\ell+2$.  The two required incoming-data bounds follow from
	\eqref{eq:stat-inflow-reg} with $L=\ell+B^{\st}_{\rm in}$ and
	$L=\ell+2+B^{\st}_{\rm in}$.  We therefore obtain, with the last term on the left coming from the weight $\ell+2$,
	\begin{equation}
		\label{eq:r12-stat-base}
		\|w^\ell f_\st\|_{\cH_{\kin,\sigma}}
		+\|w^\ell f_\st\|_{L^\infty}
		+[w^\ell f_\st]_{C^{\alpha_{\rm bs}}_{\kin}}
		+\|w^\ell\nabla_vf_\st\|_{L_x^{p}L_v^2}
		+\|w^{\ell+2}f_\st\|_{L^\infty}
		\lesssim_{\ell,p}\delta.
	\end{equation}
	
	For the $L_x^pH_\sigma^1$ estimate, \eqref{eq:base-slot-pointwise} at $r=p$ and \eqref{eq:r12-stat-base} give
	\[
	\|w^\ell f_\st\|_{L_x^pH_\sigma^1}
	\le C_\ell\|w^\ell\nabla_vf_\st\|_{L_x^pL_v^2}+C_\ell|\cD|^{1/p}\|w^{\ell+2}f_\st\|_{L^\infty}
	\lesssim_{\ell,p}\delta .
	\]
	This proves \eqref{eq:stationary-finite-weight-graph}.  Moreover, the pointwise terms of \eqref{eq:r12-stat-base} at $p=p_\flat$ give \eqref{eq:stationary-finite-weight-pointwise}; since $p_\flat$ is fixed by \eqref{eq:dyn-own-order-endpoints}, the constants in \eqref{eq:stationary-finite-weight-pointwise} do not depend on $p$.
\end{proof}

\begin{thm}[Stationary solution and regularity]
	\label{thm:stationary-diffuse}
	After decreasing $\delta_\st$ of Theorem~\ref{thm:reduced-fixed-point}, if necessary, there is $C>0$ such that the following hold for $0<\delta\le\delta_\st$; the new $\delta_\st$ and $C$ depend only on the structural data.
	\begin{enumerate}
		\item\label{it:stat-existence} (\emph{Existence.})
		The problem
		\eqref{eq:stationary-full-diffuse-problem} has a solution
		$f_\st$ with $\int_{\cD\times \bR^3_v}f_\st \sqrt\mu \dd v\dd x=0$ such that
		\begin{align}\label{2601011505}
			\begin{aligned}
				\|w^{A_\flat}f_\st\|_{L^\infty}
				+[w^{A_\flat}f_\st]_{C^{\alpha_{\rm bs}}_{\kin}}
				+\|w^{A_\flat}\nabla_vf_\st\|_{L_x^{p_\flat}L_v^2}+\sum_{m=0}^{\mathrm N}\sum_{I\in\mathfrak A^m}
				\|w^{\ell_m}\mathfrak D_If_\st\|_{\cH_{\kin,\sigma}} \le C\delta ,
			\end{aligned}
		\end{align}
		and
		\begin{equation}
			\label{eq:boundary-moment-H3}
			\|\cA_0(f_\st)_{\gamma^+}\|_{H^{\mathrm N}(\partial\cD)}\le C\delta .
		\end{equation}
		Consequently, for every $0<\beta<1$, there is $C_\beta>0$ such that
		\begin{equation}
			\label{eq:boundary-moment-holder}
			\|\cA_0(f_\st)_{\gamma^+}\|_{C^{1,\beta}(\partial\cD)}
			\le C_\beta\delta.
		\end{equation}
		
		\item\label{it:stat-uniqueness} (\emph{Uniqueness.})
		The solution $f_\st$ in (\ref{it:stat-existence}) is unique among mass-normalized solutions
		$f\in\cH_{\kin,\sigma}(\cD\times \bR^3_v)$ of \eqref{eq:stationary-full-diffuse-problem} satisfying
		\begin{equation}
			\label{eq:uniqueness-ball}
			\|\mu^{1/4}f\|_{L^\infty(\cD\times \bR^3_v)}\leq \varepsilon_{\rm dif},
		\end{equation}
		where $\varepsilon_{\rm dif}$ is the constant in \eqref{eq:diffuse-graph-smallness}.

		\item\label{it:stat-regularity} (\emph{Finite-weight regularity.})
		For every $\ell\ge0$ and $1<p<\infty$, there are $C_\ell,C_{\ell,p}>0$ such that
		\begin{equation}
			\label{eq:main-full-regularity-norm}
			\begin{aligned}
				\|w^\ell f_\st\|_{\cH_{\kin,\sigma}}
				+\|w^\ell\nabla_v f_\st\|_{L_x^{p}L_v^2}
				+\|w^\ell f_\st\|_{L_x^{p}H_\sigma^1}
				&\le C_{\ell,p}\delta,\\
				\|w^\ell f_\st\|_{L^\infty_{x,v}}
				+[w^\ell f_\st]_{C^{\alpha_{\rm bs}}_{\kin}}
				&\le C_\ell\delta,\\
				\|w^\ell(f_\st)_{\gamma^-}\|_{(C_x^{0,1}\cap C_v^{1,1})(\gamma^-)}
				&\le C_\ell\delta.
			\end{aligned}
		\end{equation}
		
		\item\label{it:stat-nonnegativity}
		(\emph{Physical nonnegativity.})  The associated physical stationary distribution
		\[
		F_\st:=\mu+\sqrt\mu\,f_\st
		\]
		is nonnegative almost everywhere in $\cD\times\bR_v^3$.
	\end{enumerate}
\end{thm}

Recall from \eqref{eq:reduced-Q3} that the left-hand side of \eqref{2601011505} coincides with $\widetilde{\mathfrak Q}_{\mathrm N}(f_\st)$, so we have
\begin{equation}
	\label{eq:solution-bound}
	\widetilde{\mathfrak Q}_{\mathrm N}(f_\st)\le C\delta.
\end{equation}

\begin{proof}[\textbf{Proof of parts (\ref{it:stat-existence})--(\ref{it:stat-regularity}) of Theorem~\ref{thm:stationary-diffuse}}]
	Parts (\ref{it:stat-existence}) and (\ref{it:stat-regularity}) follow from
	Theorem~\ref{thm:reduced-fixed-point} and
	Proposition~\ref{prop:stationary-finite-weight-regularity} at every
	$0<\delta\le\delta_\st$; $\delta_\st$ is decreased only in the proof of (\ref{it:stat-nonnegativity}).  For (\ref{it:stat-existence}),
	Theorem~\ref{thm:reduced-fixed-point} gives the mass-normalized solution
	$f_\st$ and the bound \eqref{eq:fixed-point-bound}.  The first term of \eqref{eq:fixed-point-bound} gives
	\eqref{eq:solution-bound}, which is \eqref{2601011505}, and the second gives \eqref{eq:boundary-moment-H3}.  Since $H^{\mathrm N}(\partial\cD)\subset H^3(\partial\cD)\subset C^{1,\beta}(\partial\cD)$ for $\mathrm N\ge4$, \eqref{eq:boundary-moment-H3} gives \eqref{eq:boundary-moment-holder}.  This proves (\ref{it:stat-existence}).
	
	For (\ref{it:stat-uniqueness}), let $0<\delta\le\delta_\st$ (recall from \S\ref{ss:diffuse-uniform-hierarchy} that $\delta_\st\le\delta_{\rm dif}$), and
	let $f\in \cH_{\kin,\sigma}(\cD\times \bR^3_v)$ be a mass-normalized solution of \eqref{eq:stationary-full-diffuse-problem} satisfying \eqref{eq:uniqueness-ball}.
	Put $\widetilde f:=f_\st - f$.
	By bilinearity, \[
	N[f_\st,f_\st]-N[f,f]=N[\widetilde f,f_\st]+N[f,\widetilde f] ,
	\] and $h_w$ cancels in the difference of the boundary laws; hence $\widetilde f$ satisfies
	\begin{equation*}
		\begin{cases}
			v\cdot\nabla_x\widetilde f+L\widetilde f-N[f,\widetilde f]=N[\widetilde f, f_\st]
			&\text{in }L^2(\cD;H_\sigma^{-1}),\\
			\widetilde f_{\gamma^-}=P_w \widetilde f_{\gamma^+}
			&\text{in }L^2(\gamma^-,\dd\gamma),
		\end{cases}
		\qquad \int_{\cD}\int_{\bR^3_v}\widetilde f \sqrt \mu \dd v\dd x=0.
	\end{equation*}
	By \eqref{eq:uniqueness-ball} and $\delta\le\delta_{\rm dif}$,
	\eqref{eq:diffuse-graph-smallness} holds at the background $f$, and \eqref{eq:homogeneous-diffuse-coercivity} gives
	\[
	\|\widetilde f\|_{\cH_{\kin,\sigma}}
	\le C\|N[\widetilde f,f_\st]\|_{L^2(\cD;H_\sigma^{-1})} .
	\]
	As in the proof of Proposition~\ref{prop:strict-contraction}, with $(\widetilde f,f_\st)$ in place of $(d,u_2)$,
	\[
	C\|N[\widetilde f,f_\st]\|_{L^2(\cD;H_\sigma^{-1})}\le C\delta\|\widetilde f\|_{\cH_{\kin,\sigma}}\le\tfrac12\|\widetilde f\|_{\cH_{\kin,\sigma}},
	\]
	with $C$ as in that proof and by the choice of $\delta_\st$ there; here $f_\st=\cS_\st(f_\st)$ with $f_\st\in\widetilde{\mathbb B}_{R_\st}(\delta)$ (see the proof of Theorem~\ref{thm:reduced-fixed-point}) plays the role of $u_2=\cS_\st(g_2)$.  Hence $\widetilde f=0$, that is, $f=f_\st$.
	Since $f_\st$ satisfies \eqref{eq:uniqueness-ball} by \eqref{eq:fixed-point-invariant-bound}, $\eta_\st\le\eta_1$ and \eqref{eq:construction-coefficient-regularity}, this proves (\ref{it:stat-uniqueness}).
	
	For (\ref{it:stat-regularity}), Proposition~\ref{prop:stationary-finite-weight-regularity}
	at the given $\ell$ and $p$ gives \eqref{eq:main-full-regularity-norm}.
\end{proof}

\begin{proof}[\textbf{Proof of part (\ref{it:stat-nonnegativity}) of Theorem~\ref{thm:stationary-diffuse}}]
	After decreasing $\delta_\st$ below the coefficient threshold in
	Lemma~\ref{lem:perturbed-landau-ellipticity} (this is possible since $\|\mu^{1/4}f_\st\|_{L^\infty}\le\|f_\st\|_{L^\infty}\lesssim\delta$ by \eqref{eq:stationary-finite-weight-pointwise} at $\ell=0$), set
	\[
	F:=F_\st=\mu+\sqrt\mu\,f_\st,
	\qquad
	a^{ij}:=\phi^{ij}*F,
	\qquad
	b^i:=\partial_{v_j}a^{ij},
	\qquad
	J^i:=a^{ij}\partial_{v_j}F-b^iF .
	\]
	The perturbative equation for $f_\st$ is equivalent to the physical stationary Landau equation
	$v\cdot\nabla_xF=\partial_{v_i}J^i$
	in distributions.  Moreover, by \eqref{eq:double-divergence-sigma},
	$ \partial_{v_i}b^i=-8\pi F$,
	and Lemma~\ref{lem:perturbed-landau-ellipticity} gives
	\[
	\frac12\sigma^{ij}(v)\xi_i\xi_j
	\le a^{ij}(x,v)\xi_i\xi_j
	\le2\sigma^{ij}(v)\xi_i\xi_j .
	\]
	
	\smallskip
	\noindent\emph{Step 1. An $L^1$ negative-part balance.}
	For $\varepsilon>0$ define
	\[
	\beta_\varepsilon(r):=
	\begin{cases}
		-r-\varepsilon\arctan((-r)/\varepsilon),&r<0,\\
		0,&r\ge0.
	\end{cases}
	\]
	For $r<0$, writing $s=-r>0$,
	\[
	\beta_\varepsilon'(r)=-\frac{(s/\varepsilon)^2}{1+(s/\varepsilon)^2},
	\qquad
	\beta_\varepsilon''(r)=\frac{2s/\varepsilon^2}{(1+(s/\varepsilon)^2)^2}.
	\]
	Hence $\beta_\varepsilon\in C^2(\bR)$ is convex,
	$\beta_\varepsilon(0)=\beta_\varepsilon'(0)=0$, $-1<\beta_\varepsilon'\le0$,
	$\|\beta_\varepsilon'\|_\infty+\|\beta_\varepsilon''\|_\infty<\infty$ for fixed $\varepsilon$, and
	$0\le\beta_\varepsilon(r)\le r_-$,
	$\beta_\varepsilon(r)\to r_-$.
	Put
	\[
	H_\varepsilon(r):=r\beta_\varepsilon'(r)-\beta_\varepsilon(r).
	\]
	For $r<0$, with $s=-r$ as above,
	\[
	H_\varepsilon(r)
	=\varepsilon\left(\arctan\frac{s}{\varepsilon}
	-\frac{s/\varepsilon}{1+(s/\varepsilon)^2}\right)\ge0,
	\qquad
	rH_\varepsilon(r)\le0,
	\]
	and
	\[
	|H_\varepsilon(r)|\le r_-,
	\qquad
	|r\beta_\varepsilon''(r)|
	=\frac{2(s/\varepsilon)^2}{(1+(s/\varepsilon)^2)^2}\le\frac12 .
	\]
	
	The velocity cutoffs are justified by the finite-weight bounds.  Proposition~\ref{prop:stationary-finite-weight-regularity} gives
	$f_\st\in\cH_{\kin,\sigma}$ and its finite-weight bounds.  Lemma~\ref{lem:local-embedding-H_sigma} gives
	\[
	\|w^{-1/2}f_\st\|_{L^2_v}
	+\|w^{-3/2}\nabla_vf_\st\|_{L^2_v}
	\lesssim \|f_\st\|_{H^1_\sigma}.
	\]
	Since $\sqrt\mu$ absorbs any polynomial weight, this gives, for one fixed $m>1/2$,
	\[
	w^mF,\quad w^m\nabla_vF,\quad w^mJ,\quad w^mbF
	\in L^2(\cD\times\bR^3_v).
	\]
	Here the last two statements use the coefficient bounds in
	Lemma~\ref{lem:perturbed-landau-ellipticity} and the bounds \eqref{eq:coulomb-maxwellian-anisotropy}--\eqref{eq:coulomb-sigmavv-comparison} of Lemma~\ref{lem:landau-diffusion-eigenvalues} for the Maxwellian coefficient $\sigma$.
	Moreover, $\sqrt\mu\,f_\st\in\cH_{\kin,\sigma}$ by the multiplier property in
	Lemma~\ref{lem:compact-velocity-support}, while $\mu\in\cH_{\kin,\sigma}$.  Hence
	$F=\mu+\sqrt\mu\,f_\st\in\cH_{\kin,\sigma}$.
	
	Choose radial $\chi_R,\zeta_R\in C_c^\infty(\bR^3)$ with
	$0\le\chi_R,\zeta_R\le1$,
	\[
	\chi_R=1\ \hbox{on }B_R,
	\qquad \supp\chi_R\subset B_{2R},
	\qquad |\nabla_v\chi_R|\le C/R,\qquad \zeta_R=1\ \hbox{on }B_{3R},
	\qquad \supp\zeta_R\subset B_{4R}.
	\]
	Lemma~\ref{lem:compact-velocity-support} gives
	$\zeta_RF\in\cH_{\kin}$, and
	\[
	v\cdot\nabla_x(\zeta_RF)
	=-J\cdot\nabla_v\zeta_R+\partial_{v_i}(\zeta_RJ^i)
	\]
	in $L^2_xH_v^{-1}$.  Thus, in Proposition~\ref{prop:renormalized-identity}, one may take
	$g^0=-J\cdot\nabla_v\zeta_R$ and $g^i=\zeta_RJ^i$.  Because $\chi_R$ is supported in $B_{2R}$ whereas
	$\zeta_R\equiv1$ on $B_{3R}$, the $g^0$ term vanishes on the support of $\chi_R$, and $g^i=J^i$ there.
	Consequently the stationary identity in Proposition~\ref{prop:renormalized-identity}, applied to $\zeta_RF$ with
	$\beta_\varepsilon$ and test $\chi_R$, gives
	\[
	\int_\gamma \beta_\varepsilon(F_\gamma)\chi_R\sgn(n_x\cdot v)\,\dd\gamma
	=-\int \beta_\varepsilon''(F)\,\partial_{v_i}F\,J^i\chi_R
	-\int \beta_\varepsilon'(F)J\cdot\nabla_v\chi_R .
	\]
	All unlabeled bulk integrals in this proof are over $\cD\times\bR_v^3$.
	Since $J=a\nabla_vF-bF$, and since $H_\varepsilon'(r)=r\beta_\varepsilon''(r)$ gives $\beta_\varepsilon''(F)F\nabla_vF=\nabla_vH_\varepsilon(F)$, we have
	\[
	-\beta_\varepsilon''(F)\nabla_vF\cdot J
	=-\beta_\varepsilon''(F)a^{ij}\partial_{v_i}F\partial_{v_j}F
	+b\cdot\nabla_vH_\varepsilon(F).
	\]
	Here
	$H_\varepsilon'(F)=F\beta_\varepsilon''(F)$ is bounded (by $|r\beta_\varepsilon''(r)|\le\frac12$), so
	$H_\varepsilon(F)\chi_R\in H_0^1(\bR^3_v)$ for a.e. $x$.  Approximating it by smooth compactly supported functions and using
	$\nabla_v\cdot b=-8\pi F$, we obtain
	\[
	\int b\cdot\nabla_vH_\varepsilon(F)\chi_R
	=8\pi\int FH_\varepsilon(F)\chi_R
	-\int bH_\varepsilon(F)\cdot\nabla_v\chi_R .
	\]
	We have $-\beta_\varepsilon''(F)a^{ij}\partial_{v_i}F\partial_{v_j}F\le0$ and $8\pi FH_\varepsilon(F)\le0$ (since $\beta_\varepsilon''\ge0$, $(a^{ij})$ is positive definite by $a^{ij}\xi_i\xi_j\ge\frac12\sigma^{ij}(v)\xi_i\xi_j$, and $rH_\varepsilon(r)\le0$).  Hence
	\[
	\int_{\gamma^+}\beta_\varepsilon(F_{\gamma^+})\chi_R\,\dd\gamma
	-\int_{\gamma^-}\beta_\varepsilon(F_{\gamma^-})\chi_R\,\dd\gamma
	\le \mathcal E_{R,\varepsilon},
	\]
	where
	\[
	\mathcal E_{R,\varepsilon}:=-\int\beta_\varepsilon'(F)J\cdot\nabla_v\chi_R-\int bH_\varepsilon(F)\cdot\nabla_v\chi_R .
	\]
	Since $|\beta_\varepsilon'|\le1$, $|H_\varepsilon(F)|\le|F|$ and $|\nabla_v\chi_R|\le C/R$,
	\[
	|\mathcal E_{R,\varepsilon}|
	\le \frac{C}{R}
	\int_{\cD\times\{R<|v|<2R\}}
	\bigl(|J|+|bF|\bigr)\,\dd x\dd v .
	\]
	By Cauchy--Schwarz (note that $\|w^{-m}\|_{L^2(\{R<|v|<2R\})}\lesssim R^{3/2-m}$ and $\cD$ is bounded) and $w^mJ,\ w^mbF\in L^2(\cD\times\bR^3_v)$,
	\[
	|\mathcal E_{R,\varepsilon}|
	\lesssim R^{1/2-m}
	\bigl(\|w^mJ\|_2+\|w^mbF\|_2\bigr)\longrightarrow0
	\quad\text{uniformly in }\varepsilon .
	\]
	Furthermore,
	\[
	F_-\le\sqrt\mu\,|f_\st|\in L^1(\cD\times\bR^3_v),
	\qquad
	(F_\gamma)_-\le\sqrt\mu\,|(f_\st)_\gamma|\in L^1(\gamma),
	\]
	where the trace estimate follows from Proposition~\ref{2607151218} and Cauchy--Schwarz.
	Thus dominated convergence, first as $R\to\infty$ and then as $\varepsilon\downarrow0$, yields
	\begin{equation}\label{261016980}
		\int_{\gamma^+}(F_{\gamma^+})_-\,\dd\gamma
		-\int_{\gamma^-}(F_{\gamma^-})_-\,\dd\gamma
		\le0.
	\end{equation}
	
	\smallskip
	\noindent\emph{Step 2. Nonnegativity of the boundary traces.}
	For a.e. $x\in\partial\cD$ set
	\[
	n(x):=\int_{n_x\cdot v>0}(F_{\gamma^+})_-(x,v)n_x\cdot v\,\dd v,
	\qquad
	p(x):=\int_{n_x\cdot v>0}(F_{\gamma^+})_+(x,v)n_x\cdot v\,\dd v.
	\]
	Then, by \eqref{eq:outgoing-diffuse-moment},
	\[
	\int_{n_x\cdot v>0}F_{\gamma^+}\,n_x\cdot v\,\dd v=c_\mu+\cA_0(f_\st)_{\gamma^+}=p(x)-n(x).
	\]
	Multiplying the perturbative diffuse law by $\sqrt\mu$ and using
	$\sqrt\mu\,\mathfrak m_{\Theta_w}=M_{\Theta_w}$ and $\sqrt\mu\,h_w=c_\mu M_{\Theta_w}-\mu$ (see \eqref{eq:wall-response-profiles}--\eqref{eq:nonisothermal-diffuse-response}) gives
	\[
	\begin{aligned}
		F_{\gamma^-}(x,v)
		&=\mu(v)+\sqrt{\mu(v)}\,\bigl(\mathfrak m_{\Theta_w}\cA_0(f_\st)_{\gamma^+}+h_w\bigr)(x,v)\\
		&=M_{\Theta_w(x)}(v)\bigl(c_\mu+\cA_0(f_\st)_{\gamma^+}(x)\bigr)
		=M_{\Theta_w(x)}(v)\,[p(x)-n(x)].
	\end{aligned}
	\]
	Because $M_{\Theta_w}>0$ and has unit incoming flux (by the normalization in \eqref{new:eq:intro-wall-maxwellian}),
	\[
	\int_{n_x\cdot v<0}(F_{\gamma^-})_-(x,v)|n_x\cdot v|\,\dd v
	=(n(x)-p(x))_+.
	\]
	Consequently, for a.e. $x\in\partial\cD$,
	\[
	\int_{n_x\cdot v>0}(F_{\gamma^+})_-n_x\cdot v\,\dd v-\int_{n_x\cdot v<0}(F_{\gamma^-})_-|n_x\cdot v|\,\dd v
	=n(x)-(n(x)-p(x))_+=\min\{n(x),p(x)\}\ge0.
	\]
	Hence, by \eqref{261016980} of Step~1,
	\[
	\int_{\partial\cD}\min\{n(x),p(x)\}\,\dd S_x=0,
	\qquad
	\min\{n(x),p(x)\}=0
	\quad\hbox{for a.e. }x.
	\]
	
	The smallness of the perturbation excludes $p(x)=0$.  By
	Proposition~\ref{prop:stationary-finite-weight-regularity} (see \eqref{eq:stationary-finite-weight-pointwise} at $\ell=0$),
	$\|f_\st\|_{L^\infty_{x,v}}\le {\rm C}_2\delta$,
	with ${\rm C}_2$ depending only on the structural data.
	For each $x\in\partial\cD$ put
	\[
	V_x:=B_{1/8}\!\left(\frac12n_x\right).
	\]
	Then $V_x\subset B_{5/8}$ and $n_x\cdot v\ge3/8$ on $V_x$.  Set
	\[
	m_*:=\min_{|v|\le5/8}\sqrt\mu(v)>0
	\]
	and decrease $\delta_\st$ so that ${\rm C}_2\delta_\st\le m_*/2$.
	To transfer the bulk bound to the kinetic trace, choose
	$\zeta\in C_c^\infty(\bR^3)$ with $\zeta=1$ on $B_{3/4}$.  By
	Lemma~\ref{lem:compact-velocity-support}, $\zeta f_\st\in\cH_{\kin}$ and its trace is
	$\zeta(f_\st)_\gamma$.  The set
	$\widetilde{\cO}:=\overline\cD\times B_{3/4}$ is relatively open in
	$\overline\cD\times\bR_v^3$.  Since $\zeta=1$ on $B_{3/4}$ and
	$\|f_\st\|_{L^\infty_{x,v}}\le {\rm C}_2\delta$, we have
	\[
	-{\rm C}_2\delta
	\le \operatorname*{ess\,inf}_{\widetilde{\cO}}(\zeta f_\st)
	\le \operatorname*{ess\,sup}_{\widetilde{\cO}}(\zeta f_\st)
	\le {\rm C}_2\delta.
	\]
	Applying the upper and lower trace bounds in
	\eqref{trace-int-st} gives
	\[
	-{\rm C}_2\delta\le \zeta(f_\st)_\gamma\le {\rm C}_2\delta
	\quad \dd\gamma\text{-a.e. on }\gamma\cap\widetilde{\cO}.
	\]
	Since $\zeta=1$ on $B_{3/4}$, in particular
	$|(f_\st)_\gamma(x,v)|\le {\rm C}_2\delta$ for $\dd\gamma\text{-a.e. }(x,v)$ with $|v|\le5/8$.
	Therefore
	\[
	F_\gamma(x,v)
	=\sqrt\mu(v)\bigl(\sqrt\mu(v)+(f_\st)_\gamma(x,v)\bigr)
	\ge\frac12m_*^2
	\]
	for such velocities.  Integrating over $V_x$, we obtain, for a.e. $x$,
	\[
	p(x)\ge\int_{V_x}F_{\gamma^+}(x,v)n_x\cdot v\,\dd v
	\ge \frac12m_*^2\frac38|B_{1/8}|>0.
	\]
	Thus, as $\min\{n(x),p(x)\}=0$, $n(x)=0$ for a.e. $x$, and then the diffuse law reads $F_{\gamma^-}=M_{\Theta_w}\,p$.  Consequently
	$F_{\gamma^+}\ge0$ and $F_{\gamma^-}\ge0$, $\dd\gamma$-almost everywhere.
	
	\smallskip
	\noindent\emph{Step 3. Nonnegativity in the bulk.}
	Let $\beta_*:=\beta_1$ and $H_*:=H_1$ be the functions of Step~1 at $\varepsilon=1$, so that $\beta_*(r)=r_- -\arctan(r_-)$. By Step~1, $\beta_*$ is convex and satisfies the assumptions of Proposition~\ref{prop:renormalized-identity}, $|\beta_*'|\le1$ and $\sup_r|r\beta_*''(r)|\le\frac12$. For $r=-s<0$, $H_*(r)=\arctan s-s/(1+s^2)>0$, since this function of $s$ vanishes at $s=0$ and has derivative $2s^2/(1+s^2)^2>0$; hence $rH_*(r)<0$ for $r<0$. Moreover, $0\le H_*(r)\le r_-$ by Step~1, so that $|bH_*(F)|\le|bF|$ pointwise.
	Apply the same localized stationary renormalized identity as in Step~1 with $\beta_*$.
	The boundary term is now exactly zero because $F_\gamma\ge0$ and $\beta_*$ vanishes on $[0,\infty)$.  The bound $|bH_*(F)|\le|bF|$, together with $|\beta_*'|\le1$, shows, as in Step~1, that $\mathcal E_{R,1}$ satisfies
	\[
	|\mathcal E_{R,1}|\lesssim R^{1/2-m}(\|w^mJ\|_2+\|w^mbF\|_2)
	\]
	and therefore tends to zero as $R\to\infty$.  Also
	$0\le H_*(r)\le\pi/2$ and
	$|FH_*(F)|\le(\pi/2)F_-\in L^1$.
	Since $0\le\beta_*''\le C$ and $0\le a\le2\sigma$, the diffusion density satisfies
	\[
	0\le\beta_*''(F)a^{ij}\partial_{v_i}F\partial_{v_j}F\chi_R\le2C\sigma^{ij}\partial_{v_i}F\partial_{v_j}F,
	\]
	and
	$\sigma^{ij}\partial_{v_i}F\partial_{v_j}F\in L^1(\cD\times\bR_v^3)$ since $F\in\cH_{\kin,\sigma}\subset L_x^2H_\sigma^1$.  Therefore dominated convergence, as $R\to\infty$, gives
	\[
	0
	=-\int \beta_*''(F)a^{ij}\partial_{v_i}F\partial_{v_j}F
	+8\pi\int FH_*(F).
	\]
	Here $\int\beta_*''(F)a^{ij}\partial_{v_i}F\partial_{v_j}F\ge0$ and $\int FH_*(F)\le0$, and $FH_*(F)<0$ almost everywhere on $\{F<0\}$. Hence
	\[
	\int_{\cD\times\bR_v^3}(-F)H_*(F)\,\dd x\dd v=0,
	\]
	and consequently $F\ge0$ almost everywhere. This proof is independent of the evolutionary stability theorem.
\end{proof}

\section{Global Diffuse Dynamics and Exponential Relaxation}\label{sec:diffuse-evolution}

We study the evolution about the stationary state $f_\st$ of Theorem~\ref{thm:stationary-diffuse}. Uniform estimates for the fixed-background problem lead to the nonlinear solution by contraction. The regularized construction underlying those estimates and the passage to the limit are proved in the next section. The highest differentiation order is
\[
\mathrm M:=\mathrm N-1 ,
\]
one below the stationary order; thus $\mathrm M\ge3$.
In this section and the next, the \emph{structural data} are the domain $\cD$ with the fields $\cX_{\rm c}$ and the order $\mathrm N$ of
\S\ref{sec:diffuse-common}, the evolutionary target weight $\ell^\dy_\sharp$
fixed below, and the rate $\lambda_w$ of \eqref{eq:dyn-theta-relaxation}; the exponents,
the weights and the time cutoff fixed in this section depend only on them.  The decay rate $\lambda$
is used before it is fixed in \eqref{eq:dyn-decay-rate-fixed-point}; it depends only on
the structural data.  Unless otherwise stated, constants depend only on the structural data, never on $T$; additional dependencies are indicated by subscripts.

\subsection{The evolutionary problem and its centering}\label{ss:diffuse-evolution-problem}

The physical problem is
\begin{equation}
	\label{eq:dyn-physical-problem}
	\begin{cases}
		(\partial_t+v\cdot\nabla_x)f_\pt+Lf_\pt=N[f_\pt,f_\pt]
		&\text{in }(0,T)\times\cD\times\bR_v^3,\\
		(f_\pt)_{\gamma^-}=P_w(t)(f_\pt)_{\gamma^+}+h_w(t)
		&\text{on }\gamma_T^-,\\
		f_\pt(0)=f_0
		&\text{in }\cD\times\bR_v^3.
	\end{cases}
\end{equation}
Here $P_w(t)$ and $h_w(t)$ are the diffuse operator and boundary datum of \eqref{eq:nonisothermal-diffuse-response}, with wall temperature $\Theta_w(t)$ satisfying \eqref{eq:dyn-theta-amplitude}.  The uniform bounds \eqref{eq:dyn-theta-defect-time-bound} for $h_w$, through $\mathrm N=\mathrm M+1$ time derivatives, imply for each finite $T>0$ that
\begin{equation}
	\label{eq:dyn-wall-datum-class}
	h_w\in H^{\mathrm M+1}\bigl((0,T);L^2(\gamma^-)\bigr).
\end{equation}
Since $h_w=c_\mu d_\delta$, the unit-flux identity \eqref{eq:wall-defect-flux-orthogonality} also gives
\begin{equation}
	\label{eq:dyn-wall-datum-flux}
	\beta_{h_w}=0 ,
\end{equation}
with $\beta$ of \eqref{eq:dynamic-boundary-mass-scalars}.

Fix a nonnegative integer evolutionary target weight $\ell^\dy_\sharp$.
The stationary target weight is
$$
\ell_\sharp:=\ell^\dy_\sharp+\mathrm M A_\diamond,
$$
with $A_\diamond$ of \eqref{eq:Adef}; $\ell_m$ and $\widetilde{\mathfrak Q}_{\mathrm N}$ denote
those of \eqref{eq:finite-hierarchy-weight-choice} and \eqref{eq:reduced-Q3} at this target.
Fix $\delta_\star>0$ sufficiently small for
Theorem~\ref{thm:stationary-diffuse} at $\ell_\sharp$ and for the
extended-wall estimates of Lemma~\ref{lem:dyn-wall-extended} and
Proposition~\ref{prop:dynamic-alpha-uniform-graph}.
We assume $0<\delta\le\delta_\star$ and decrease this structural
threshold when needed.

At the target weight $\ell_\sharp$, Theorem~\ref{thm:stationary-diffuse} applies to the limiting wall datum $h_w^\infty$ of \eqref{eq:limiting-wall-profiles}, since $\delta\le\delta_\star\le\delta_\st$. The resulting stationary solution $f_\st$ satisfies
\begin{equation}
	\label{eq:dyn-theta-background}
	(f_\st)_{\gamma^-}=P_w^\infty(f_\st)_{\gamma^+}+h_w^\infty,
	\qquad
	\widetilde{\mathfrak Q}_{\mathrm N}(f_\st)\le C\delta .
\end{equation}
The boundary identity is its stationary diffuse law, while the bound follows from \eqref{eq:solution-bound}; the constant depends only on the structural data because $\ell_\sharp$ is fixed by $\mathrm N$ and $\ell^\dy_\sharp$.  

Center the physical solution and the background $g$ at $f_\st$ by setting
\[
f_\dy:=f_\pt-f_\st,
\qquad
g_\dy:=g-f_\st .
\]
The stationary term belongs to the finite-weight graph and trace classes of Proposition~\ref{prop:stationary-finite-weight-regularity} for every finite $\ell\ge0$ and $1<p<\infty$.

Subtracting the stationary diffuse law in \eqref{eq:dyn-theta-background} from the time-dependent boundary condition gives
\begin{equation}
	\label{eq:dyn-theta-residual}
	(f_\dy)_{\gamma^-}=P_w(t)(f_\dy)_{\gamma^+}+b_w(t),
	\qquad
	b_w(t):=\bigl(d_\delta(t)-d_\delta^\infty\bigr)\cA_0(f_\st)_{\gamma^+}
	+\bigl(h_w(t)-h_w^\infty\bigr).
\end{equation}
Here $d_\delta$ and $d_\delta^\infty$ are the wall profiles of \eqref{eq:wall-response-profiles} and \eqref{eq:limiting-wall-profiles}. Since $h_w=c_\mu d_\delta$ and $h_w^\infty=c_\mu d_\delta^\infty$, we may write \[
b_w(t)=\bigl(d_\delta(t)-d_\delta^\infty\bigr)\bigl(\cA_0(f_\st)_{\gamma^+}+c_\mu\bigr) .
\]  The decaying Gaussian bound \eqref{eq:dyn-theta-defect-relaxation-bound} and the stationary outgoing-moment estimate \eqref{eq:boundary-moment-H3} therefore yield
\begin{equation}
	\label{eq:dyn-theta-residual-bound}
	\|b_w(t)\|_{L^2(\gamma^-)}\le C\,\delta\,\ee^{-\lambda_wt},
	\qquad t\ge0 .
\end{equation}
Consequently,
\begin{equation*}
	\int_s^\infty\|b_w(r)\|_{L^2(\gamma^-)}^2\dd r\le C\delta^2\ee^{-\lambda_ws},
	\qquad s\ge0 .
\end{equation*}
The second factor in the formula for $b_w$ is independent of $v$. Applying \eqref{eq:wall-defect-flux-orthogonality} at the current and limiting wall temperatures gives
\begin{equation}
	\label{eq:dyn-theta-residual-flux}
	\beta_{b_w}=0 .
\end{equation}  For any background $g$, bilinearity gives
\[
N[g,f_\pt]-N[f_\st,f_\st]=N[g,f_\dy]+N[g_\dy,f_\st] .
\]
Thus, at $g=f_\pt$, the perturbation satisfies the following fixed-background problem, with $L_gf=Lf-N[g,f]$ as in \eqref{def:L_g}:
\begin{equation}
	\label{eq:dyn-perturbation-problem}
	\begin{cases}
		(\partial_t+v\cdot\nabla_x)f_\dy+L_gf_\dy
		=\mathfrak S:=N[g_\dy,f_\st]
		&\text{in }(0,T)\times\cD\times\bR_v^3,\\
		(f_\dy)_{\gamma^-}=P_w(t)(f_\dy)_{\gamma^+}+b_w(t)
		&\text{on }\gamma_T^-,\\
		f_\dy(0)=f_{\dy,0}:=f_0-f_\st
		&\text{in }\cD\times\bR_v^3 .
	\end{cases}
\end{equation}

We use the mixed alphabet \eqref{eq:dyn-mixed-alphabet} for complete words $|I|\le\mathrm M$, with $j(I)$ the number of temporal letters. As the stationary state is independent of time and the spatial complete fields are time-independent,
\begin{equation}
	\label{eq:dyn-fs-time-independent}
	\mathfrak D_If_\st=0\quad\text{if }j(I)\ge1 .
\end{equation}

\subsection{Exponents, weights, and norms}\label{ss:dyn-rungs}

Set
\begin{gather}
	\ell^\dy_{\mathrm M}:=\ell^\dy_\sharp,\qquad\ell^\dy_m:=\ell^\dy_{m+1}+A_\diamond\quad(0\le m\le\mathrm M-1),
	\label{eq:dyn-weight-choice}\\
	A^\dy_\flat:=\ell^\dy_0+3 .
	\label{eq:dyn-data-weights}
\end{gather}
These weights satisfy
\[
\ell^\dy_0>\ell^\dy_1>\cdots>\ell^\dy_{\mathrm M-1}>\ell^\dy_{\mathrm M}=\ell^\dy_\sharp,\qquad\ell^\dy_{\mathrm M-1}\ge A_\diamond;
\]
moreover, by the definition of $\ell_\sharp$, $\ell^\dy_0=\ell^\dy_\sharp+\mathrm M A_\diamond=\ell_\sharp$, so that \eqref{eq:dyn-stationary-slots} below holds at every weight $\ell\le\ell^\dy_0$.

For $1\le m\le\mathrm M$, let $K_m$ be the integer \eqref{eq:dyn-Km}, which covers the weight losses at order $m$ in \S\ref{ss:dyn-ordered-words}. It depends only on $m$ and the fixed exponents.  Starting from $\Lambda_{\mathrm M}=\ell^\dy_\sharp$,
choose the remaining weights backwards so that
\begin{equation}\label{eq:dyn-realization-weights}
	\Lambda_{m-1}\ge K_m(\Lambda_m+A_\diamond+1)
	\quad(1\le m\le\mathrm M),
	\qquad \Lambda_0\ge A^\dy_\flat.
\end{equation}
In particular,
$\Lambda_m\ge\ell^\dy_m$ and $\Lambda_{m-1}\ge\Lambda_m$.
The order-$m$ jets will be constructed for $L\le\Lambda_m$, and the base estimates for $L\le\Lambda_0$. The weights are fixed independently of $T$, the spectral cutoff, and the regularization parameter.

Let $B^{\dy}_{\rm in}:=B^{\dy}_1(\cD,\alpha_{\rm bs})$, where $\alpha_{\rm bs}$ is given by \eqref{eq:holder-exponent}. Lemma~\ref{260527559_dy} at $\ell=L$ requires data at weight $L+B^{\dy}_{\rm in}+1$, and \eqref{eq:dynamic-B1} gives $B^{\dy}_{\rm in}\ge B^{\dy}_0$ of \eqref{eq:dynamic-B0}. Fix an integer $\check A$ sufficiently large in terms of $\Lambda_0$, $B^{\dy}_{\rm in}$, $\mathrm M$, and $A_\diamond$.

The spatial base exponents $(p_{\rm G},p_{\rm c})$ in Corollary~\ref{cor:ordered-word-rungs-mixed} are $(P_\natural,q)$ for the evolutionary estimates, rather than the stationary pair $(\infty,p_\flat)$ of \eqref{eq:dyn-own-order-exponents}. Choosing a finite $P_\natural$ allows the $L_t^2L_x^{P_\natural}$ estimate of Lemma~\ref{lem:dyn-natural-base}. For the estimates \eqref{eq:dyn-slot-gaussian-L2} and \eqref{eq:dyn-slot-coercive-L2}, used when one collision factor is stationary, the base pair is $(P_\natural,p_\flat)$. All these estimates are uniform in $T$.

We define the exponents $p_\sharp$ and $q$ by
\begin{equation}
	\label{eq:dyn-base-choice}
	p_\sharp:=2p_\flat,
	\qquad
	\frac1q:=\frac{1-\theta}{p_\sharp}+\frac\theta2 ,
\end{equation}
where $0<\theta<1$, $P_\natural<\infty$ and $\theta_\sharp>0$ are fixed, depending only on
$\mathrm N$ and $\ell^\dy_\sharp$, such that
\begin{equation}
	\label{eq:dyn-finite-base}
	\ell^\dy_0\le(1-\theta)(A^\dy_\flat-2),
	\qquad
	(p_\flat-1)\,\theta_\sharp\le\frac1{2\mathrm M},
	\qquad
	\frac{\mathrm M-1}{P_\natural}+\frac1q
	+2\theta_\sharp(\mathrm M-1)\Bigl(\frac1{q_\flat}+\frac12\Bigr)\le\frac1{p_\flat},
\end{equation}
and
\begin{equation}
	\label{eq:dyn-polarized-margin}
	\frac1{P_\natural}+\frac{2\theta_\sharp(\mathrm M-1)}{q_\flat}
	\le\frac1{q_\flat\,\mathrm N} .
\end{equation}
These conditions give
\[
p_\flat<q<p_\sharp,
\qquad
\theta_\sharp\le\frac1{2\mathrm M},
\qquad
P_\natural>q_\flat\mathrm N,
\qquad
2\theta_\sharp(\mathrm M-1)\le\frac1{\mathrm N} :
\]
the upper bound on $q$ follows from $\theta>0$ and $p_\sharp>2$, the lower one from the last condition of \eqref{eq:dyn-finite-base}, whose other two terms are positive, the bound on $\theta_\sharp$ from the second condition, since $p_\flat-1\ge1$, and the last two from \eqref{eq:dyn-polarized-margin}.
These parameters can be chosen: the first condition of \eqref{eq:dyn-finite-base} is $\theta(\ell^\dy_0+1)\le1$, while the remaining conditions hold for sufficiently small positive $\theta$, $1/P_\natural$, and $\theta_\sharp$, since $1/p_\sharp=1/(2p_\flat)<1/p_\flat$. The last condition of \eqref{eq:dyn-finite-base} and \eqref{eq:dyn-polarized-margin} are used in Lemma~\ref{lem:dyn-pairings}; the first condition of \eqref{eq:dyn-finite-base} enters Proposition~\ref{prop:dyn-hierarchy}.

Fix $1\le m\le\mathrm M$. For $1\le j\le m-1$, the spatial exponents $\mathsf g^\natural(j;m)$ and $\mathsf s^\natural(j;m)$ are those of \eqref{eq:mixed-rung-exponent} at top order $m$, with $\theta_\circ=\theta_\sharp$, $a_{\rm G}=q_\flat$ and $(p_{\rm G},p_{\rm c})=(P_\natural,q)$:
\begin{equation}
	\label{eq:dyn-slot-exponents}
	\begin{aligned}
		\frac1{\mathsf g^\natural(j;m)}
		&=\frac j{m\,q_\flat}+\frac{m-j}{m\,P_\natural}
		+2j\theta_\sharp\Bigl(1-\frac jm\Bigr)\Bigl(\frac1{q_\flat}-\frac1{P_\natural}\Bigr),\\
		\frac1{\mathsf s^\natural(j;m)}
		&=\frac j{2m}+\frac{m-j}{m\,q}
		+2j\theta_\sharp\Bigl(1-\frac jm\Bigr)\Bigl(\frac12-\frac1q\Bigr),
	\end{aligned}
\end{equation}
and at the endpoints we set
\begin{equation}
	\label{eq:dyn-slot-endpoints}
	\mathsf g^\natural(0;m):=\infty,
	\quad
	\mathsf g^\natural(m;m):=q_\flat,
	\qquad
	\mathsf s^\natural(0;m):=p_\sharp,
	\quad
	\mathsf s^\natural(m;m):=2 .
\end{equation}

In the bound \eqref{eq:dyn-bilinear-allocation} of a collision term of order $m$, a jet of order $j$ is measured, by Lemma~\ref{lem:dyn-slots}, in
\[
L_t^{2m/j}L_x^{\mathsf g^\natural(j;m)}L_v^2\quad\text{in the first slot},
\qquad
L_t^{2m/j}L_x^{\mathsf s^\natural(j;m)}H_\sigma^1\quad\text{in the second slot},
\]
with $2m/0:=\infty$.

For a spatial base exponent $p$, $2\le a\le p$ and $0\le j\le m$, put
\begin{equation}
	\label{eq:dyn-pj}
	\frac1{\mathsf p_j(a)}:=\frac1{\mathsf e(m/j;a,p)}+2j\theta_\sharp\Bigl(1-\frac jm\Bigr)\Bigl(\frac1a-\frac1p\Bigr).
\end{equation}
For $1\le j\le m-1$, this is \eqref{eq:mixed-rung-exponent} at top order $m$ and $\theta_\circ=\theta_\sharp$, while $\mathsf p_0(a)=p$ and $\mathsf p_m(a)=a$. The dependence of $\mathsf p_j(a)$ on $p$ is implicit and will always be specified when the notation is used. The next lemma establishes the spatial H\"older pairings for two evolutionary jets and for each ordering of an evolutionary jet with a stationary one. The former are used in Lemma~\ref{lem:dyn-collision-terms}, and the latter in Lemma~\ref{cor:dyn-forcing-estimate}.

\begin{lemma}
	\label{lem:dyn-pairings}
	Let $1\le m\le\mathrm M$ and $0\le k\le m$.  Here $\mathsf p_j(a)$ is \eqref{eq:dyn-pj} at this $m$, and $\mathsf g$ and $\mathsf s$ are the rungs \eqref{eq:dyn-own-order-exponents} of $f_\st$,
	\[
	\frac1{\mathsf g(j)}=\frac j{\mathrm N\,q_\flat},
	\qquad
	\frac1{\mathsf s(j)}=\frac1{p_\flat}+\frac j{\mathrm N}\Bigl(\frac12-\frac1{p_\flat}\Bigr),
	\qquad 0\le j\le\mathrm N.
	\]
	These stationary exponents are defined at order $\mathrm N$, independently of the evolutionary order $m$.
	\begin{enumerate}
		\item\label{it:dyn-pairing-slots}
		The slot exponents satisfy
		\begin{equation}
			\label{eq:dyn-allocation-space}
			\frac1{\mathsf g^\natural(k;m)}+\frac1{\mathsf s^\natural(m-k;m)}\le\frac12 .
		\end{equation}
		\item\label{it:dyn-pairing-first}
		With $\mathsf p_j$ taken at $p=p_\flat$,
		\begin{equation}
			\label{eq:dyn-polarized-first}
			\frac1{\mathsf g(k)}+\frac1{\mathsf p_{m-k}(2)}\le\frac12 .
		\end{equation}
		\item\label{it:dyn-pairing-second}
		With $\mathsf p_j$ taken at $p=P_\natural$,
		\begin{equation}
			\label{eq:dyn-polarized-second}
			\frac1{\mathsf p_{k}(q_\flat)}+\frac1{\mathsf s(m-k)}\le\frac12 .
		\end{equation}
	\end{enumerate}
\end{lemma}

\begin{proof}
	(\ref{it:dyn-pairing-slots}) For $k=0$ and $k=m$, \eqref{eq:dyn-slot-endpoints}, \eqref{eq:dyn-base-choice} and \eqref{eq:dyn-own-order-endpoints} give the values $\frac12$ and $\frac1{q_\flat}+\frac1{p_\sharp}=\frac12-\delta_\flat-\frac1{2p_\flat}$ of the left side, respectively.  Let $1\le k\le m-1$.  The coefficients of $\theta_\sharp$ in \eqref{eq:dyn-slot-exponents} at $j=k$ and at $j=m-k$ coincide,
	\[
	2k\Bigl(1-\frac km\Bigr)=2(m-k)\Bigl(1-\frac{m-k}m\Bigr)=\frac{2k(m-k)}m ,
	\]
	and \eqref{eq:dyn-slot-exponents} and \eqref{eq:dyn-own-order-endpoints} give
	\[
	\frac1{\mathsf g^\natural(k;m)}+\frac1{\mathsf s^\natural(m-k;m)}
	=\frac12-\frac k{m\,p_\flat}{}-\frac km\delta_\flat
	+\frac{m-k}{m\,P_\natural}+\frac k{m\,q}
	+2\theta_\sharp\,\frac{k(m-k)}m
	\Bigl(\frac1{q_\flat}+\frac12-\frac1{P_\natural}-\frac1q\Bigr) .
	\]
	By $k\ge1$, $m-k\le\mathrm M-1$ and the last condition of \eqref{eq:dyn-finite-base},
	\[
	\frac{m-k}{m\,P_\natural}+\frac k{m\,q}+2\theta_\sharp\,\frac{k(m-k)}m\Bigl(\frac1{q_\flat}+\frac12-\frac1{P_\natural}-\frac1q\Bigr)
	\le\frac km\Bigl(\frac{\mathrm M-1}{P_\natural}+\frac1q
	+2\theta_\sharp(\mathrm M-1)\Bigl(\frac1{q_\flat}+\frac12\Bigr)\Bigr)
	\le\frac k{m\,p_\flat} ,
	\]
	so the right side is at most $\frac12-\frac km\delta_\flat\le\frac12$.
	
	(\ref{it:dyn-pairing-first}) By \eqref{eq:dyn-pj} at $p=p_\flat$,
	\[
	\frac1{\mathsf p_{m-k}(2)}=\frac12-\frac km\Bigl(\frac12-\frac1{p_\flat}\Bigr)+\frac{2\theta_\sharp k(m-k)}m\Bigl(\frac12-\frac1{p_\flat}\Bigr) ,
	\]
	and this, together with
	\eqref{eq:dyn-own-order-endpoints} and
	$\tfrac1{\mathsf g(k)}=\tfrac k{q_\flat\mathrm N}$ (see \eqref{eq:dyn-own-order-exponents}), gives
	\[
	\frac1{\mathsf g(k)}+\frac1{\mathsf p_{m-k}(2)}-\frac12
	=k\Bigl(\frac12-\frac1{p_\flat}\Bigr)
	\Bigl(\frac1{\mathrm N}-\frac1m
	+\frac{2\theta_\sharp(m-k)}m\Bigr)
	{}-\frac k{\mathrm N}\delta_\flat ,
	\]
	which vanishes at $k=0$. For $k\ge1$, it is nonpositive, since $\delta_\flat>0$ and the second factor in parentheses is nonpositive: $m-k\le\mathrm M-1$ and, by
	\eqref{eq:dyn-polarized-margin} and $m\le\mathrm M=\mathrm N-1$,
	\[
	2\theta_\sharp(\mathrm M-1)\ \le\ \frac1{\mathrm N}
	\ \le\ 1-\frac m{\mathrm N}
	\ =\ m\Bigl(\frac1m-\frac1{\mathrm N}\Bigr) .
	\]
	
	(\ref{it:dyn-pairing-second}) By \eqref{eq:dyn-own-order-exponents} and \eqref{eq:dyn-pj} at $p=P_\natural$,
	\[
	\frac1{\mathsf s(m-k)}=\frac12-\frac{\mathrm N-m+k}{\mathrm N}\Bigl(\frac12-\frac1{p_\flat}\Bigr),
	\qquad
	\frac1{\mathsf p_k(q_\flat)}=\frac k{m\,q_\flat}+\frac{m-k}{m\,P_\natural}+\frac{2\theta_\sharp k(m-k)}m\Bigl(\frac1{q_\flat}-\frac1{P_\natural}\Bigr) .
	\]
	Hence, by \eqref{eq:dyn-own-order-endpoints} and $\tfrac km-\tfrac{\mathrm N-m+k}{\mathrm N}=-\tfrac{m-k}m\cdot\tfrac{\mathrm N-m}{\mathrm N}$,
	\[
	\frac1{\mathsf p_k(q_\flat)}+\frac1{\mathsf s(m-k)}-\frac12
	=\frac{m-k}m\Bigl(-\frac{\mathrm N-m}{\mathrm N}
	\Bigl(\frac12-\frac1{p_\flat}\Bigr)
	+\frac1{P_\natural}
	+2\theta_\sharp k\Bigl(\frac1{q_\flat}-\frac1{P_\natural}\Bigr)\Bigr)
	{}-\frac km\delta_\flat ,
	\]
	which equals $-\delta_\flat$ at $k=m$. For $k\le m-1$, it is nonpositive, since $k\le\mathrm M-1$, \eqref{eq:dyn-polarized-margin}, $\tfrac1{q_\flat}\le\tfrac12-\tfrac1{p_\flat}$ and $\mathrm N-m\ge1$ give
	\[
	\frac1{P_\natural}+2\theta_\sharp k\Bigl(\frac1{q_\flat}-\frac1{P_\natural}\Bigr)
	\le\frac1{P_\natural}+\frac{2\theta_\sharp(\mathrm M-1)}{q_\flat}
	\le\frac1{q_\flat\mathrm N}
	\le\frac{\mathrm N-m}{\mathrm N}\Bigl(\frac12-\frac1{p_\flat}\Bigr) .
	\]
\end{proof}

For $0\le t<\infty$, $\cH_{\kin,\sigma}(\cD_{(t,\infty)}\times\bR_v^3)$ consists of the $f$ on $\cD_{(t,\infty)}\times\bR_v^3$ whose restrictions to $\cD_{(t,T')}\times\bR_v^3$ lie in $\cH_{\kin,\sigma}(\cD_{(t,T')}\times\bR_v^3)$ with norms bounded uniformly in $T'\in(t,\infty)$; its norm is the supremum of these norms.  We write $\cD_\infty:=\cD_{(0,\infty)}$.  For $0<T\le\infty$, every quantity over $(0,T)$ below, including those of Definition~\ref{def:dyn-coordinates}, is nondecreasing in $T$ and at $T=\infty$ equals its limit as $T\uparrow\infty$; hence an estimate between such quantities with a constant independent of $T$ holds at $T=\infty$ if it holds at every finite $T$.  In this section $T\ge1$ means $1\le T<\infty$, and $T=\infty$ is admitted only where it is written.

\begin{defn}[The coordinates and the rows]
	\label{def:dyn-coordinates}
	Let $f$ be defined on $\cD_T\times\bR_v^3$.  The \emph{strong}, the \emph{decay} and the \emph{base coordinates} of $f$ are
	\begin{alignat}{2}
		\cA_T[f]&:=\|w^{A^\dy_\flat}f\|_{L^\infty}
		+\bigl[w^{A^\dy_\flat}f\bigr]_{C^{\alpha_{\rm bs}}_{\kin}}
		+\|w^{A^\dy_\flat}\nabla_vf\|_{L_t^\infty L_x^{p_\sharp}L_v^2},
		\label{eq:dyn-coordinate-A}\\
		\cB_T[f]&:=\sup_{0\le t\le T}\ee^{\lambda t}
		\|f\|_{\cH_{\kin,\sigma}(\cD_{(t,T)}\times\bR_v^3)}^2,
		\label{eq:dyn-coordinate-B}\\
		\cS_T^{\ell}[f]&:=\|w^{\ell}f\|_{L_t^2L_x^{q}H_\sigma^1},
		\qquad
		\cS_T[f]:=\cS_T^{0}[f],
		\qquad \ell\ge0,
		\label{eq:dyn-base-row}
	\end{alignat}
	where $\lambda$ is the decay rate and $q$ is that of \eqref{eq:dyn-base-choice}.  In the first slot, $w^{-2}f$ is also measured in $L_t^2$, at the spatial exponent $P_\natural$ of \eqref{eq:dyn-finite-base}:
	\begin{equation}
		\label{eq:dyn-natural-base}
		\mathsf G_\natural(T)[f]:=\bigl\|w^{-2}f\bigr\|_{L_t^2L_x^{P_\natural}L_v^2}.
	\end{equation}
	The \emph{rows} of the evolutionary hierarchy are the rows \eqref{eq:order-two-Um}--\eqref{eq:order-two-Wm} of \S\ref{ss:weight-interpolation} over the mixed alphabet, at the top order $\mathrm M$ and the weights \eqref{eq:dyn-weight-choice}, and $\cN_m$ is the sum of the weighted rows of order at most $m$: for $0\le m\le\mathrm M$,
	\begin{alignat}{2}
		\mathsf U_m(T)[f]&:=\sum_{I\in\mathfrak A_\bullet^m}
		\|\mathfrak D_If\|_{\cH_{\kin,\sigma}(\cD_T\times\bR_v^3)},
		\notag\\
		\mathsf W_m(T)[f]&:=\sum_{I\in\mathfrak A_\bullet^m}
		\|w^{\ell^\dy_m}\mathfrak D_If\|_{\cH_{\kin,\sigma}(\cD_T\times\bR_v^3)},
		\label{eq:dyn-hierarchy-rows}\\
		\cN_m(T)[f]&:=\sum_{j\le m}\mathsf W_j(T)[f].
		\label{eq:dyn-hier-rows-aggregate}
	\end{alignat}
	For $g$ on $\cD_T\times\bR_v^3$ with $1\le T\le\infty$, the coordinates and the rows are collected in the \emph{reduced quantity}, the counterpart of $\widetilde{\mathfrak Q}_{\mathrm N}$ of \eqref{eq:reduced-Q3}:
	\begin{equation}
		\label{eq:dyn-reduced-Q}
		\widetilde{\mathfrak Q}^\dy_T(g):=\cA_T[g]+\cB_T[g_\dy]^{1/2}+\cS_T[g_\dy]
		+\cN_{\mathrm M}(T)[g_\dy].
	\end{equation}
\end{defn}

The strong coordinate $\cA_T$, the evolutionary counterpart of $\cB_{\rm bs}$ in \eqref{eq:base-strong}, is used for $f_\pt$, $f_\dy$, and $g$. The decay and base coordinates and the rows are used for the centered functions $f_\dy$ and $g_\dy$. We suppress $[f]$ when $f=f_\dy$. The rows may take the value $+\infty$, as in \S\ref{ss:weight-interpolation}.  The graph norm in the rows controls both traces: the incoming trace is one of its terms by Definition~\ref{def:H_kin_sigma}, and, by \eqref{260902137}, the outgoing trace satisfies
\[
\|f_{\gamma^+}\|_{L^2(\gamma^+_T)}\le C\|f\|_{\cH_{\kin,\sigma}(\cD_T\times\bR_v^3)} .
\]

Since $w^{-2-A^\dy_\flat}\in L^2_v$, H\"older's inequality in $v$ gives, for every
measurable $f$ on $\cD\times\bR_v^3$,
\begin{equation}
	\label{eq:dyn-strong-fiber-mixed}
	\|w^{-2}f\|_{L^2_v}\le C\|w^{A^\dy_\flat}f\|_{L^\infty_v} .
\end{equation}
Hence, by \eqref{eq:base-slot-pointwise} at $r=p_\sharp$, every measurable $f$ on $\cD_T\times\bR_v^3$ satisfies, for $0\le\ell\le A^\dy_\flat-2$,
\begin{equation}
	\label{eq:dyn-coordinate-A-reading}
	\|w^{-2}f\|_{L_t^\infty L_x^\infty L_v^2}
	+\|w^{\ell}f\|_{L_t^\infty L_x^{p_\sharp}H_\sigma^1}
	\le C_\ell\,\cA_T[f] .
\end{equation}

\begin{lemma}[Control of jets]
	Every measurable $f$ on $\cD_T\times\bR_v^3$ satisfies, for $0\le m\le\mathrm M$ and $0\le\ell\le\ell^\dy_m$,
	\begin{equation}
		\label{eq:dyn-slot-tops-rows}
		\sum_{|K|\le m}\Bigl(\bigl\|w^{\,m-|K|-2}\mathfrak D_Kf\bigr\|_{L_t^2L_x^{q_\flat}L_v^2}
		+\bigl\|w^{\,\ell+m-|K|}\mathfrak D_Kf\bigr\|_{L_t^2L_x^2H_\sigma^1}\Bigr)
		\le C_\ell\,\cN_m(T)[f],
	\end{equation}
	the sum being taken over the ordered words $K$ of $\mathfrak A_\bullet$, with $C_\ell$ depending only on the structural data and $\ell$.
\end{lemma}

\begin{proof}
	The two weights in \eqref{eq:dyn-slot-tops-rows} are $m-|K|-2$ and $\ell+m-|K|$, respectively.
	For $|K|<m$, by $\ell\le\ell^\dy_m$, $A_\diamond\ge1$ and \eqref{eq:dyn-weight-choice},
	\[
	m-|K|\le\ell+m-|K|\le\ell^\dy_m+(m-|K|)A_\diamond=\ell^\dy_{|K|} ,
	\]
	and at $|K|=m$ the weights $0$ and $\ell$ are at most $\ell^\dy_m$.  For $|K|=j$, put $k:=\ell^\dy_j-(m-j)\ge0$; then \eqref{eq:graph-to-mixed-integrability-evolutionary} at $\varrho=2$ and this $k$, applied to $w^{\ell^\dy_j}\mathfrak D_Kf$, gives
	\[
	\|w^{m-j-2}\mathfrak D_Kf\|_{L_t^2L_x^{q_\flat}L_v^2}\le C\|w^{\ell^\dy_j}\mathfrak D_Kf\|_{\cH_{\kin,\sigma}(\cD_T\times\bR_v^3)},
	\]
	and \eqref{eq:real-weight-Hsigma-monotonicity} gives
	\[
	\|w^{\ell+m-j}\mathfrak D_Kf\|_{L_t^2L_x^2H_\sigma^1}\le C_\ell\|w^{\ell^\dy_j}\mathfrak D_Kf\|_{\cH_{\kin,\sigma}(\cD_T\times\bR_v^3)}.
	\]
	Summing over $|K|=j$ gives
	\[
	\sum_{|K|=j}\Bigl(\|w^{m-j-2}\mathfrak D_Kf\|_{L_t^2L_x^{q_\flat}L_v^2}+\|w^{\ell+m-j}\mathfrak D_Kf\|_{L_t^2L_x^2H_\sigma^1}\Bigr)\le C_\ell\mathsf W_j(T)[f] .
	\]
	Summing over $0\le j\le m$, by \eqref{eq:dyn-hier-rows-aggregate}, gives \eqref{eq:dyn-slot-tops-rows}.
\end{proof}

The unweighted rows $\mathsf U_m$ are not listed in \eqref{eq:dyn-reduced-Q}, since the weighted rows dominate them: Lemma~\ref{lem:compact-velocity-support}\,(1) with the $C^1_v$ multiplier $w^{-\ell^\dy_m}$, applied to $w^{\ell^\dy_m}\mathfrak D_If$, gives
\begin{equation}
	\label{eq:dyn-unweighted-row-dominated}
	\mathsf U_m(T)[f]\le C\,\mathsf W_m(T)[f],
	\qquad 0\le m\le\mathrm M .
\end{equation}

At $j=0$, the endpoints \eqref{eq:dyn-slot-endpoints} are supplied by \eqref{eq:dyn-coordinate-A-reading}; at $j=m$, they follow from \eqref{eq:dyn-slot-tops-rows}. The next lemma interpolates between these bounds. Its chain and local-integrability hypotheses follow, in particular, from $\cN_{\mathrm M}(T)[f]<\infty$.

\begin{lemma}
	\label{lem:dyn-slots}
	Let $T\ge1$, let $1\le m\le\mathrm M$, and let $f$ satisfy, when $m\ge2$, \eqref{eq:mixed-rung-chain} and \eqref{eq:mixed-rung-local-integrability} on $(0,T)\times\cD\times\bR_v^3$ at the order $M'=m$ of Corollary~\ref{cor:ordered-word-rungs-mixed}.  Then, with $2m/0:=\infty$ and with $\mathsf p_k$ of \eqref{eq:dyn-pj} taken at $p=P_\natural$ in \eqref{eq:dyn-slot-gaussian-L2} and at $p=p_\flat$ in \eqref{eq:dyn-slot-coercive-L2}, every ordered word $I$ of the mixed alphabet
	\eqref{eq:dyn-mixed-alphabet} with $|I|=k$ and $0\le k\le m$ satisfies
	\begin{align}
		\bigl\|w^{-2}\mathfrak D_If\bigr\|
		_{L_t^{2m/k}
			L_x^{\mathsf g^\natural(k;m)}L_v^2}
		&\le C\bigl(\|w^{-2}f\|_{L_t^\infty L_x^\infty L_v^2}+\cN_m(T)[f]
		+\mathsf G_\natural(T)[f]\bigr),
		\label{eq:dyn-slot-gaussian}\\
		\bigl\|w^{\ell}\mathfrak D_If\bigr\|
		_{L_t^{2m/k}
			L_x^{\mathsf s^\natural(k;m)}H_\sigma^1}
		&\le C_\ell\bigl(\|w^{\ell}f\|_{L_t^\infty L_x^{p_\sharp}H_\sigma^1}+\cN_m(T)[f]
		+\cS_T^{\ell}[f]\bigr),
		\label{eq:dyn-slot-coercive}\\
		\bigl\|w^{-2}\mathfrak D_If\bigr\|_{L_t^2L_x^{\mathsf p_k(q_\flat)}L_v^2}
		&\le C\bigl(\mathsf G_\natural(T)[f]+\cN_m(T)[f]\bigr),
		\label{eq:dyn-slot-gaussian-L2}\\
		\bigl\|w^{\ell}\mathfrak D_If\bigr\|_{L_t^2L_x^{\mathsf p_k(2)}H_\sigma^1}
		&\le C_\ell\bigl(\cS_T^{\ell}[f]+\cN_m(T)[f]\bigr),
		\label{eq:dyn-slot-coercive-L2}
	\end{align}
	for every $0\le\ell\le\ell^\dy_m$, with $C$ depending only on the structural data and $C_\ell$ only on the structural data and $\ell$.
\end{lemma}

\begin{proof}
	At $k=0$ the left sides of \eqref{eq:dyn-slot-gaussian} and \eqref{eq:dyn-slot-coercive} are the norms
	$\|w^{-2}f\|_{L_t^\infty L_x^\infty L_v^2}$ and $\|w^{\ell}f\|_{L_t^\infty L_x^{p_\sharp}H_\sigma^1}$,
	by \eqref{eq:dyn-slot-endpoints} and the convention $2m/0:=\infty$; at $k=m$, by \eqref{eq:dyn-slot-endpoints} and $2m/m=2$,
	they are summands of the left side of \eqref{eq:dyn-slot-tops-rows}.  Let $1\le k\le m-1$.
	
	Apply Corollary~\ref{cor:ordered-word-rungs-mixed} on $(0,T)$; this is admissible since $T\ge1$ and since $m\ge2$ in this range of $k$. We take $(\mathrm L_{\rm G},a_{\rm G})=(-2,q_\flat)$, the weight $\mathrm L=\ell$, the order $M'=m$, the top order $m$, the temporal base exponent $r_0=\infty$, the spatial base exponents
	$(p_{\rm G},p_{\rm c})=(P_\natural,q)$ and the parameter $\theta_\circ=\theta_\sharp$ of \eqref{eq:dyn-finite-base}, which satisfies
	$\theta_\sharp\le1/(2\mathrm M)\le1/(2m)$.
	The temporal exponent of \eqref{eq:mixed-word-gaussian-rungs} and \eqref{eq:mixed-word-coercive-rungs} is $\mathsf e(m/k,\infty)=2m/k$ and their spatial exponents are $\mathsf g^\natural(k;m)$ and $\mathsf s^\natural(k;m)$, by
	\eqref{eq:dyn-slot-exponents}.  The quantities $\mathfrak b^\sharp_\rho$ of \eqref{eq:slab-rung-quantities} and $\mathfrak b^0_\rho$ of Corollary~\ref{cor:ordered-word-rungs-mixed} at $\rho\in\{\infty,2\}$ are
	\[
	\mathfrak b^\sharp_{\infty}=\|w^{-2}f\|_{L_t^\infty L_x^{P_\natural}L_v^2}
	\le C\,\|w^{-2}f\|_{L_t^\infty L_x^\infty L_v^2},
	\qquad
	\mathfrak b^\sharp_{2}=\mathsf G_\natural(T)[f],
	\]
	by the boundedness of $\cD$ and \eqref{eq:dyn-natural-base}, and
	\[
	\mathfrak b^0_{\infty}=\|w^{\ell}f\|_{L_t^\infty L_x^{q}H_\sigma^1}
	\le C\,\|w^{\ell}f\|_{L_t^\infty L_x^{p_\sharp}H_\sigma^1},
	\qquad
	\mathfrak b^0_{2}=\cS_T^{\ell}[f] ,
	\]
	by $q<p_\sharp$, the boundedness of $\cD$ and \eqref{eq:dyn-base-row}.  For $\mathfrak G^\sharp_\bullet$ and $\mathfrak G^0_\bullet$ of \eqref{eq:mixed-rung-tops}, Minkowski's inequality in $t$ and \eqref{eq:dyn-slot-tops-rows}, at $\ell=0$ for $\mathfrak G^\sharp_\bullet$, give $\mathfrak G^\sharp_\bullet\le C\cN_m(T)[f]$ and $\mathfrak G^0_\bullet\le C_\ell\cN_m(T)[f]$.  The bounds \eqref{eq:mixed-word-gaussian-rungs} and \eqref{eq:mixed-word-coercive-rungs}, together with Young's inequality, give \eqref{eq:dyn-slot-gaussian} and
	\eqref{eq:dyn-slot-coercive}.
	
	Since $\mathsf p_0(a)=p$ and $\mathsf p_m(a)=a$, at $k=0$ the left side of \eqref{eq:dyn-slot-gaussian-L2} is $\mathsf G_\natural(T)[f]$, by \eqref{eq:dyn-natural-base}, and that of \eqref{eq:dyn-slot-coercive-L2} is at most $C\cS_T^{\ell}[f]$, by $p_\flat<q$, the boundedness of $\cD$ and \eqref{eq:dyn-base-row}; at $k=m$ both left sides are summands of the left side of \eqref{eq:dyn-slot-tops-rows}.  For $1\le k\le m-1$ we apply Corollary~\ref{cor:ordered-word-rungs-mixed} with the choices above, except that the temporal base exponent is $r_0=2$ and $p_{\rm c}=p_\flat$.  The temporal exponent of \eqref{eq:mixed-word-gaussian-rungs} and \eqref{eq:mixed-word-coercive-rungs} is then $\mathsf e(m/k,2)=2$, their spatial exponents are $\mathsf p_k(q_\flat)$ and $\mathsf p_k(2)$ of \eqref{eq:dyn-pj}, at $p=P_\natural$ and at $p=p_\flat$, respectively, and
	\[
	\mathfrak b^\sharp_2=\mathsf G_\natural(T)[f],
	\qquad
	\mathfrak b^0_2=\|w^{\ell}f\|_{L_t^2L_x^{p_\flat}H_\sigma^1}\le C\,\cS_T^{\ell}[f],
	\]
	by \eqref{eq:dyn-natural-base}, and by $p_\flat<q$, the boundedness of $\cD$ and \eqref{eq:dyn-base-row}.  At $r_0=2$ the right sides of \eqref{eq:mixed-word-gaussian-rungs} and \eqref{eq:mixed-word-coercive-rungs} are at most $C(\mathfrak G^\sharp_\bullet+\mathfrak b^\sharp_2)$ and $C(\mathfrak G^0_\bullet+\mathfrak b^0_2)$, and the bounds on $\mathfrak G^\sharp_\bullet$ and $\mathfrak G^0_\bullet$ above give \eqref{eq:dyn-slot-gaussian-L2} and \eqref{eq:dyn-slot-coercive-L2}.
\end{proof}

The next lemma is used in Proposition~\ref{prop:dyn-hierarchy}.

\begin{lemma}
	\label{lem:dyn-natural-base}
	Let $\vartheta_\natural:=2/P_\natural\in(0,1]$. The time exponent $2\vartheta_\natural$ below may be less than one; the estimate uses only the corresponding integral. Let $f$ on $\cD_T\times\bR_v^3$ have
	$\|w^{-2}f\|_{L_t^\infty L_x^\infty L_v^2}<\infty$ and $f\in L_t^{2\vartheta_\natural}L_x^2H_\sigma^1$.  Then
	\begin{equation}
		\label{eq:dyn-natural-base-bound}
		\mathsf G_\natural(T)[f]
		\le C\,\|w^{-2}f\|_{L_t^\infty L_x^\infty L_v^2}^{1-\vartheta_\natural}\,
		\|f\|_{L_t^{2\vartheta_\natural}L_x^2H_\sigma^1}^{\vartheta_\natural}.
	\end{equation}
\end{lemma}

\begin{proof}
	Put $\varphi:=\|w^{-2}f\|_{L_v^2}$.  Since $P_\natural\ge2$ and $\frac1{P_\natural}=\frac{\vartheta_\natural}2$, H\"older's inequality gives, at each $t$,
	\[
	\|\varphi(t)\|_{L_x^{P_\natural}}\le\|\varphi(t)\|_{L_x^\infty}^{1-\vartheta_\natural}
	\|\varphi(t)\|_{L_x^2}^{\vartheta_\natural} .
	\]
	Moreover, $\|w^{-2}f\|_{L_v^2}\le C\|f\|_{H_\sigma^1}$ by
	\eqref{eq:local-H1-H_sigma}. Taking the $L_t^2$ norm, we obtain
	\[
	\begin{aligned}
		\mathsf G_\natural(T)[f]
		&\le\|w^{-2}f\|_{L_t^\infty L_x^\infty L_v^2}^{1-\vartheta_\natural}\,\bigl\|\,\|\varphi\|_{L_x^2}^{\vartheta_\natural}\bigr\|_{L_t^2}
		=\|w^{-2}f\|_{L_t^\infty L_x^\infty L_v^2}^{1-\vartheta_\natural}\,\|\varphi\|_{L_t^{2\vartheta_\natural}L_x^2}^{\vartheta_\natural}\\
		&\le C\,\|w^{-2}f\|_{L_t^\infty L_x^\infty L_v^2}^{1-\vartheta_\natural}\,
		\|f\|_{L_t^{2\vartheta_\natural}L_x^2H_\sigma^1}^{\vartheta_\natural},
	\end{aligned}
	\]
	which proves \eqref{eq:dyn-natural-base-bound}.
\end{proof}

Let $0\le m\le\mathrm M$.  By \eqref{eq:collision-remainder-expansion}, for an ordered word $I$ of the mixed alphabet with $|I|=m$, $\mathfrak D_IN[a,b]$ is the sum of the remainder $\mathfrak R_I[a,b]$ and of terms $N[\mathfrak D_{I_1}a,\mathfrak D_{I_2}b]$, where $I_1$ and $I_2$ are ordered words of the mixed alphabet with $|I_1|=k$, $|I_2|=m-k$ and $0\le k\le m$.  By \eqref{eq:velalg-poly}, for $\ell\ge0$ and pointwise in $(t,x)$,
\begin{equation}
	\label{eq:dyn-bilinear-pointwise}
	\bigl\|w^{\ell}N[\mathfrak D_{I_1}a,\mathfrak D_{I_2}b](t,x)\bigr\|_{H^{-1}_\sigma}
	\le C_\ell\,\bigl\|w^{-2}\mathfrak D_{I_1}a(t,x)\bigr\|_{L_v^2}\,
	\bigl\|w^{\ell}\mathfrak D_{I_2}b(t,x)\bigr\|_{H_\sigma^1} .
\end{equation}
The next lemma, used in \S\ref{ss:dyn-differentiated}, bounds these terms and the terms \eqref{eq:deformation-remainder-structure} of $\mathfrak R_I[a,b]$; \eqref{eq:dyn-bilinear-allocation} and \eqref{eq:dyn-bilinear-allocation-deformation} are the counterparts of \eqref{eq:hierarchy-coulomb} and \eqref{eq:hierarchy-deformation}.

\begin{lemma}
	\label{lem:dyn-collision-terms}
	Let $T\ge1$, $1\le m\le\mathrm M$ and $0\le\ell\le\ell^\dy_m$, and let $a$ and $b$ satisfy, when $m\ge2$, \eqref{eq:mixed-rung-chain}--\eqref{eq:mixed-rung-local-integrability} on $(0,T)\times\cD\times\bR_v^3$ at the order $M'=m$.  Then the following hold, with $C_\ell$ depending only on the structural data and $\ell$.
	\begin{enumerate}
		\item\label{it:dyn-collision-main}
		For every ordered word $I$ of the mixed alphabet with $|I|=m$, every term $N[\mathfrak D_{I_1}a,\mathfrak D_{I_2}b]$ of the expansion \eqref{eq:collision-remainder-expansion} of $\mathfrak D_IN[a,b]$ satisfies
		\begin{align}
			\label{eq:dyn-bilinear-allocation}
			\begin{aligned}
				&\bigl\|w^{\ell}N[\mathfrak D_{I_1}a,\mathfrak D_{I_2}b]\bigr\|
				_{L_t^2L_x^2H^{-1}_\sigma}\\
				&\qquad \le C_\ell\,
				\Bigl(\|w^{-2}a\|_{L_t^\infty L_x^\infty L_v^2}+\cN_m(T)[a]
				{}+\mathsf G_\natural(T)[a]\Bigr)
				\Bigl(\|w^{\ell}b\|_{L_t^\infty L_x^{p_\sharp}H_\sigma^1}+\cN_m(T)[b]
				{}+\cS_T^{\ell}[b]\Bigr) .
			\end{aligned}
		\end{align}
		\item\label{it:dyn-collision-deformation}
		For every ordered word $I$ of the mixed alphabet with $|I|=m$, every term \eqref{eq:deformation-remainder-structure} of $\mathfrak R_I[a,b]$ satisfies
		\begin{align}
			\label{eq:dyn-bilinear-allocation-deformation}
			&\begin{aligned}&\bigl\|w^{\ell}\,p\,N_{\mathsf K_A}[p_1\mathfrak D_{J_1}a,p_2\mathfrak D_{J_2}b]\bigr\|
				_{L_t^2L_x^2H^{-1}_\sigma}\\
				&\qquad \le C_\ell\,
				\Bigl(\|w^{-2}a\|_{L_t^\infty L_x^\infty L_v^2}+\cN_{m-1}(T)[a]
				{}+\mathsf G_\natural(T)[a]\Bigr)
				\Bigl(\|w^{\ell+A_\diamond}b\|_{L_t^\infty L_x^{p_\sharp}H_\sigma^1}
				+\cN_{m-1}(T)[b]
				{}+\cS_T^{\ell+A_\diamond}[b]\Bigr) .
			\end{aligned}
		\end{align}
	\end{enumerate}
\end{lemma}

\begin{proof}
	(\ref{it:dyn-collision-main}) For $0\le k\le m$ put
	\[
	\Bigl(\frac1{r_1},\frac1{q_1}\Bigr)
	:=\Bigl(\frac k{2m},\ \frac1{\mathsf g^\natural(k;m)}\Bigr),
	\qquad
	\Bigl(\frac1{r_2},\frac1{q_2}\Bigr)
	:=\Bigl(\frac{m-k}{2m},\ \frac1{\mathsf s^\natural(m-k;m)}\Bigr) .
	\]
	By \eqref{eq:dyn-slot-gaussian}, the first factor on the right side of \eqref{eq:dyn-bilinear-pointwise} is bounded in $L_t^{r_1}L_x^{q_1}$, and, by \eqref{eq:dyn-slot-coercive} with $m-k$ in place of $k$, the second in $L_t^{r_2}L_x^{q_2}$.  The temporal reciprocals add to $1/2$, so H\"older's inequality in time is uniform in $T$. This requires the exponents at the actual differentiated order $m$, rather than those at order $\mathrm N$ in \eqref{eq:Dg-rungs}.  In space, $\frac1{q_1}+\frac1{q_2}\le\frac12$ by Lemma~\ref{lem:dyn-pairings}\,(\ref{it:dyn-pairing-slots}).  Hence \eqref{eq:dyn-bilinear-pointwise}, H\"older's inequality with the exponents $(r_1,q_1)$ and $(r_2,q_2)$, and Lemma~\ref{lem:dyn-slots} give \eqref{eq:dyn-bilinear-allocation}.
	
	(\ref{it:dyn-collision-deformation}) Set
	$j=|J_1|+|J_2|\le m-1$.  By
	\eqref{eq:deformation-weighted-output-spatial} with $\varpi=\frac14$ and
	$\mu^{1/4}\lesssim w^{-2}$, each remainder is bounded by the product
	of the first-slot norm of $w^{-2}\mathfrak D_{J_1}a$ and the
	second-slot norm of $w^{\ell+A_\diamond}\mathfrak D_{J_2}b$.
	Here $2m+2\le A_\diamond$, and the replacement of the velocity
	weight follows from \eqref{eq:real-weight-Hsigma-monotonicity}.
	If $j\ge1$, write $k=|J_1|$ and $h=|J_2|$, so that $k+h=j$.
	The two factors are measured, respectively, in
	\[
	L_t^{2j/k}L_x^{\mathsf g^\natural(k;j)}L_v^2,
	\qquad
	L_t^{2j/h}L_x^{\mathsf s^\natural(h;j)}H_\sigma^1,
	\]
	with $2j/0=\infty$, and the second factor carries the weight
	$\ell+A_\diamond$.  Lemma~\ref{lem:dyn-slots} applies at order $j$,
	since $\ell+A_\diamond\le\ell^\dy_{m-1}\le\ell^\dy_j$.
	The time reciprocals add to $1/2$, and the spatial reciprocals
	add to at most $1/2$ by
	Lemma~\ref{lem:dyn-pairings}\,(\ref{it:dyn-pairing-slots}) at order $j$.
	H\"older's inequality and the pointwise deformation estimate thus
	give the claimed product bound with $\cN_j(T)\le\cN_{m-1}(T)$.
	None of the embedding constants depends on $T$, $N$, or $\varepsilon$.
	
	If $j=0$, independently of $m$, both factors have order zero.
	Using the exponent pairs $(\infty,\infty)$ and $(2,2)$ in time and
	space, respectively, gives
	\[
	\begin{aligned}
		&\|w^\ell pN_{\mathsf K_A}[p_1a,p_2b]\|_{L_t^2L_x^2H_\sigma^{-1}}\\
		&\quad\lesssim
		\|w^{-2}a\|_{L_t^\infty L_x^\infty L_v^2}
		\|w^{\ell+A_\diamond}b\|_{L_t^2L_x^2H_\sigma^1}
		\lesssim
		\|w^{-2}a\|_{L_t^\infty L_x^\infty L_v^2}\,\cN_0(T)[b].
	\end{aligned}
	\]
	The last inequality follows from \eqref{eq:dyn-slot-tops-rows}
	at order zero, since
	$\ell+A_\diamond\le\ell^\dy_{m-1}\le\ell^\dy_0$.
	This includes the undifferentiated remainder terms at every order
	$1\le m\le\mathrm M$.
\end{proof}

The centered collision terms also contain derivatives of $f_\st$, both through $g=f_\st+g_\dy$ and through the source $\mathfrak S$ in \eqref{eq:dyn-perturbation-problem}.

\begin{lemma}
	The steady state $f_\st$ of \eqref{eq:dyn-theta-background} satisfies, for every $0\le\ell\le\ell_\sharp$ and every
	$0\le j\le\mathrm M$,
	\begin{equation}
		\label{eq:dyn-stationary-slots}
		\sum_{|I|=j}\bigl\|w^{-2}\mathfrak D_If_\st\bigr\|
		_{L_x^{\mathsf g(j)}L_v^2}
		+\sum_{|I|=j}\bigl\|w^{\ell}\mathfrak D_If_\st\bigr\|
		_{L_x^{\mathsf s(j)}H_\sigma^1}
		\le C_\ell\,\delta ,
	\end{equation}
	with $C_\ell$ depending only on the structural data and $\ell$.
\end{lemma}

\begin{proof}
	By \eqref{eq:dyn-fs-time-independent}, only the words without temporal letters contribute.  For the first sum with $1\le j\le\mathrm M$, Lemma~\ref{lem:background-finite-p}, applied to $f_\st$, gives
	\[
	\sum_{|I|=j}\|w^{-2}\mathfrak D_If_\st\|_{L_x^{\mathsf g(j)}L_v^2}\le C\delta ;
	\]
	its hypothesis holds by \eqref{eq:dyn-theta-background}, with $\eta=C\delta$. At $j=0$ the sum is $\|w^{-2}f_\st\|_{L_x^\infty L_v^2}$, since $\mathsf g(0)=\infty$, and, by \eqref{eq:base-fiber-Lsharp} and \eqref{eq:dyn-theta-background},
	\[
	\|w^{-2}f_\st\|_{L_x^\infty L_v^2}\le C\widetilde{\mathfrak Q}_{\mathrm N}(f_\st)\le C\delta .
	\]
	For the second sum with $1\le j\le\mathrm M$, recall that $\mathsf s(j)=\mathsf e(\mathrm N/j,p_\flat)$ by \eqref{eq:dyn-own-order-exponents}.  Hence \eqref{eq:ordered-word-coercive-rungs-weighted} of
	Corollary~\ref{cor:ordered-word-gaussian-rungs}, with $\mathrm N$ in place of $M$, $\mathrm L=\ell$ and
	$p_{\rm c}=p_\flat$, and Young's inequality give
	\[
	\sum_{|I|=j}\bigl\|w^{\ell}\mathfrak D_If_\st\bigr\|_{L_x^{\mathsf s(j)}H_\sigma^1}
	\le C_\ell\bigl(\|w^{\ell}f_\st\|_{L_x^{p_\flat}H_\sigma^1}+\cG^{\ell}_{\mathrm N,2}(f_\st)\bigr);
	\]
	its hypothesis \eqref{eq:pointwise-slice-chain-hypothesis} holds for $f_\st$, as verified in the proof of Lemma~\ref{lem:background-finite-p}, and its finiteness hypothesis follows from the bounds below. At $j=0$ the sum is $\|w^{\ell}f_\st\|_{L_x^{p_\flat}H_\sigma^1}$, since $\mathsf s(0)=p_\flat$. It remains to bound $\|w^{\ell}f_\st\|_{L_x^{p_\flat}H_\sigma^1}$ and $\cG^{\ell}_{\mathrm N,2}(f_\st)$.  The weights of
	$\cG^{\ell}_{\mathrm N,2}(f_\st)$ in \eqref{eq:GN-lower-order} satisfy
	$\ell+\mathrm N-|K|\le\ell_{|K|}$ for $0\le|K|\le\mathrm N$,
	because $\ell\le\ell_\sharp=\ell_{\mathrm N}$ and
	$\ell_m\ge\ell_{m+1}+1$ by \eqref{eq:finite-hierarchy-weight-choice}.  Hence \eqref{eq:real-weight-Hsigma-monotonicity} and \eqref{innerprod_st_sigma} give the first bound below, and \eqref{eq:base-slot-control} at $L=\ell\le\ell_0\le A_\flat-2$ gives the second,
	\[
	\cG^{\ell}_{\mathrm N,2}(f_\st)\le C_\ell\sum_{m=0}^{\mathrm N}W_m(f_\st),
	\qquad
	\|w^{\ell}f_\st\|_{L_x^{p_\flat}H_\sigma^1}\le C_\ell\cB_{\rm bs}(f_\st),
	\]
	and, by \eqref{eq:reduced-Q3} and \eqref{eq:dyn-theta-background},
	\[
	\sum_{m=0}^{\mathrm N}W_m(f_\st)+\cB_{\rm bs}(f_\st)=\widetilde{\mathfrak Q}_{\mathrm N}(f_\st)\le C\delta.
	\]
	Combining the bounds for the two sums proves \eqref{eq:dyn-stationary-slots} for $0\le j\le\mathrm M$.
\end{proof}

\subsection{The differentiated system}\label{ss:dyn-differentiated}

We differentiate \eqref{eq:dyn-perturbation-problem} along mixed complete words and estimate the resulting forcing uniformly in $T$.

\begin{prop}[Differentiated system]
	\label{prop:dyn-differentiated-equation}
	Let $1\le m\le\mathrm M$ and let $f_\dy$ solve the interior equation of
	\eqref{eq:dyn-perturbation-problem} in $\cD'(\cD_T\times\bR_v^3)$.  Assume that the mixed jets form a chain,
	\begin{equation}
		\label{eq:dyn-differentiated-chain}
		\mathfrak D_Jf_\dy\in\cH_{\kin,\sigma}(\cD_T\times\bR_v^3)
		\quad(|J|\le m),
		\qquad
		\mathfrak D_\nu(\mathfrak D_Jf_\dy)=\mathfrak D_{\nu J}f_\dy
		\quad\text{in }\cD'(\cD_T\times\bR_v^3)\ \ (|\nu J|\le m);
	\end{equation}
	that the jets of length less than $m$ have finite weighted norms,
	\begin{equation}
		\label{eq:dyn-differentiated-lower-weight}
		\sum_{|J|<m}
		\bigl\|w^{A_\diamond}\mathfrak D_Jf_\dy\bigr\|_{L^2_tL^2_xH^1_\sigma}<\infty ;
	\end{equation}
	and that the jets of $g$ and $g_\dy$ satisfy
	\begin{equation}
		\label{eq:dyn-differentiated-first-slots}
		\mu^{1/4}\mathfrak D_Jg,\ \mu^{1/4}\mathfrak D_Jg_\dy\in L^2_tL^2_{x,v},
		\qquad |J|\le m .
	\end{equation}
	Then, for every ordered word
	$I$ of $\mathfrak A_\bullet$ with $|I|\le m$,
	\begin{equation}
		\label{eq:dyn-differentiated-equation}
		(\partial_t+v\cdot\nabla_x)\mathfrak D_If_\dy+L_g\mathfrak D_If_\dy
		=\mathfrak F_I[f_\dy]
		\qquad\text{in }\cD'(\cD_T\times\bR_v^3),
	\end{equation}
	with $\mathfrak F_I$ the sum of
	\begin{align}
		&-\sum_{|J|<|I|}\cV_{I,J}\mathfrak D_Jf_\dy-[\mathfrak D_I,L]f_\dy,
		\label{eq:dyn-forcing-linear}\\
		&\sum_{I_1\sqcup I_2=I,\ I_1\neq\varnothing}
		N\bigl[\mathfrak D_{I_1}g,\ \mathfrak D_{I_2}f_\dy\bigr],
		\label{eq:dyn-forcing-coefficient}\\
		&\sum_{I_1\sqcup I_2=I}
		N\bigl[\mathfrak D_{I_1}g_\dy,\ \mathfrak D_{I_2}f_\st\bigr],
		\label{eq:dyn-forcing-polarized}\\
		&\mathfrak R_I[g,f_\dy]+\mathfrak R_I[g_\dy,f_\st],
		\label{eq:dyn-forcing-deformation}
	\end{align}
	where the last line consists of the deformation remainders in \eqref{eq:collision-remainder-expansion}. For a word $I$ with spatial subword $I'$ and $j=j(I)$, the normal ordering \eqref{eq:dyn-normal-ordering} and the time Leibniz rule give
	\[
	\mathfrak R_I[a,b]=\sum_{i=0}^{j}\binom ji\,\mathfrak R_{I'}\bigl[\partial_t^{\,i}a,\partial_t^{\,j-i}b\bigr].
	\]
	Here $\mathfrak R_{I'}$ is defined by \eqref{eq:collision-remainder-expansion}. Each summand has the structure \eqref{eq:deformation-remainder-structure} at order $|I|$, with $J_1,J_2$ ordered words of the mixed alphabet. Moreover,
	$\mathfrak F_\varnothing[f_\dy]=\mathfrak S$.  If $\mathfrak S$ is replaced by $\mathfrak S+F$ in
	\eqref{eq:dyn-perturbation-problem}, then $\mathfrak D_IF$ is added to $\mathfrak F_I[f_\dy]$.
\end{prop}

\begin{proof}
	For the spatial subword we use the distributional identities \eqref{eq:weak-transport-word}--\eqref{eq:weak-collision-word} of Theorem~\ref{thm:weak-differentiated-system} at order $m\le\mathrm M<\mathrm N$. Their proof follows that of
	Proposition~\ref{prop:graph-ordered-interior-algebra}, which applies without change with test functions in $(t,x,v)$ and with $L^2_tL^2_xH^1_\sigma$ and
	$L^2_tL^2_{x,v}$ in place of $L^2_xH^1_\sigma$ and $L^2_{x,v}$, because the spatial letters have time-independent coefficients by \eqref{eq:complete-phase-space-derivative}. The required hypotheses, including those for $f_\st$ in the second collision slot, hold as follows.
	\begin{itemize}
		\item With $f_I:=\mathfrak D_If_\dy$, the condition $f_I\in\cH_{\kin,\sigma}$ and \eqref{eq:distributional-jet-chain} are
		\eqref{eq:dyn-differentiated-chain}.
		\item The lower-weight bound \eqref{eq:lower-weight-assumption} for $f_\dy$ is
		\eqref{eq:dyn-differentiated-lower-weight}.
		\item The background condition \eqref{eq:background-jet-assumption} is
		\eqref{eq:dyn-differentiated-first-slots}.
		\item The function $f_\st$ in the second slot of \eqref{eq:dyn-forcing-polarized} and
		\eqref{eq:dyn-forcing-deformation} has $\mathfrak D_Jf_\st\in\cH_{\kin,\sigma}$, $|J|\le\mathrm N$,
		by \eqref{eq:dyn-theta-background} and \eqref{eq:unweighted-row-dominated}; these jets are
		distributional derivatives, so they satisfy \eqref{eq:distributional-jet-chain}.  The
		lower-weight bound \eqref{eq:lower-weight-assumption} for $f_\st$ holds by
		\eqref{eq:dyn-theta-background} and \eqref{eq:real-weight-Hsigma-monotonicity}, since $\ell_k\ge A_\diamond$ for $k\le\mathrm N-1$ in
		\eqref{eq:finite-hierarchy-weight-choice} and the jets $\mathfrak D_Jf_\st$ whose words $J$ contain a temporal letter vanish by \eqref{eq:dyn-fs-time-independent}.
	\end{itemize}
	For the temporal letters we use \eqref{eq:dyn-spacetime-leibniz}. The time derivative commutes with every letter by \eqref{eq:dyn-time-commutation}.  The hypotheses of the first identity of \eqref{eq:dyn-spacetime-leibniz} are \eqref{eq:dyn-differentiated-chain},
	\eqref{eq:dyn-differentiated-first-slots} and, for $f_\st$ in the second slot, the properties of the jets $\mathfrak D_Jf_\st$ verified above,
	with $\partial_tf_\st=0$.  In the normal ordering \eqref{eq:dyn-normal-ordering}, temporal differentiation acts before the spatial letters. Thus \eqref{eq:dyn-spacetime-leibniz} is applied to $N$, without differentiating the operators $N_{\mathsf K_A}$ in \eqref{eq:deformation-remainder-structure}.  Applying \eqref{eq:dyn-spacetime-leibniz} to the
	temporal letters and then \eqref{eq:weak-transport-word}--\eqref{eq:weak-collision-word}
	to the spatial subword, with \eqref{eq:mixed-transport-commutator} for the transport coefficients, we obtain
	\eqref{eq:dyn-forcing-linear}--\eqref{eq:dyn-forcing-deformation}; the $I_1=\varnothing$ term $N[g,\mathfrak D_If_\dy]$ belongs to $L_g\mathfrak D_If_\dy$ by \eqref{def:L_g} and is therefore absent from \eqref{eq:dyn-forcing-coefficient}. By
	\eqref{eq:dyn-fs-time-independent}, the summands of \eqref{eq:dyn-forcing-polarized} whose word $I_2$ contains a temporal letter vanish.  This completes the proof.
\end{proof}

\begin{lemma}
	\label{cor:dyn-forcing-estimate}
	Let $T\ge1$, $1\le m\le\mathrm M$ and $0\le\ell\le\ell^\dy_m$, let $I$ be an ordered word
	of $\mathfrak A_\bullet$ with $|I|=m$, and let $f_\dy$ on $\cD_T\times\bR_v^3$ and the background $g$ satisfy
	\eqref{eq:dyn-differentiated-chain}--\eqref{eq:dyn-differentiated-first-slots}.
	Then the forcing $\mathfrak F_I[f_\dy]$ of Proposition~\ref{prop:dyn-differentiated-equation} satisfies
	\begin{equation}
		\label{eq:dyn-forcing-estimate}
		\begin{aligned}
			\bigl\|w^{\ell}\mathfrak F_I[f_\dy]\bigr\|_{L^2_tL^2_xH^{-1}_\sigma}
			\le{}&C_\ell\sum_{|J|<m}\bigl\|w^{\ell+A_\diamond}\mathfrak D_Jf_\dy\bigr\|
			_{L^2_tL^2_xH^1_\sigma}\\
			&+C_\ell\Bigl(\|w^{-2}g_\dy\|_{L_t^\infty L_x^\infty L_v^2}+\cN_m(T)[g_\dy]
			{}+\mathsf G_\natural(T)[g_\dy]\Bigr)\\
			&\qquad\times\Bigl(\|w^{\ell+A_\diamond}f_\dy\|_{L_t^\infty L_x^{p_\sharp}H_\sigma^1}
			+\cN_m(T)[f_\dy]
			{}+\cS_T^{\ell+A_\diamond}[f_\dy]\Bigr)\\
			&+C_\ell\,\delta\,\Bigl(\cS_T^{\ell+A_\diamond}[f_\dy]+\cN_m(T)[f_\dy]\Bigr)\\
			&+C_\ell\,\delta\,\Bigl(\mathsf G_\natural(T)[g_\dy]
			+\cN_m(T)[g_\dy]\Bigr) ,
		\end{aligned}
	\end{equation}
	with $C_\ell$ depending only on the structural data and $\ell$ (in particular independent of $T$).
\end{lemma}

\begin{proof}
	In this proof the constants may depend on $\ell$. Since $g=g_\dy+f_\st$, the terms \eqref{eq:dyn-forcing-coefficient} and \eqref{eq:dyn-forcing-deformation} split, and $\mathfrak F_I[f_\dy]$ is the sum of \eqref{eq:dyn-forcing-linear} and of
	\begin{align*}
		&\sum N\bigl[\mathfrak D_{I_1}g_\dy,\mathfrak D_{I_2}f_\dy\bigr]+\mathfrak R_I[g_\dy,f_\dy],\\
		&\sum N\bigl[\mathfrak D_{I_1}f_\st,\mathfrak D_{I_2}f_\dy\bigr]+\mathfrak R_I[f_\st,f_\dy],\\
		&\sum_{I_1\sqcup I_2=I}N\bigl[\mathfrak D_{I_1}g_\dy,\mathfrak D_{I_2}f_\st\bigr]+\mathfrak R_I[g_\dy,f_\st],
	\end{align*}
	the sums in the first two lines being taken over $I_1\sqcup I_2=I$ with $I_1\neq\varnothing$.  We estimate the linear commutators and the three displayed groups separately. Since $f_\st$ is stationary, the factors involving its derivatives are in $L_t^\infty$; hence their products with the $L_t^2$ factors are controlled independently of $T$.
	
	\noindent\emph{Linear commutators.}
	We show that
	\[
	\Bigl\|w^{\ell}\Bigl(\sum_{|J|<m}\cV_{I,J}\mathfrak D_Jf_\dy+[\mathfrak D_I,L]f_\dy\Bigr)\Bigr\|_{L^2_tL^2_xH^{-1}_\sigma}
	\le C\sum_{|J|<m}\bigl\|w^{\ell+A_\diamond}\mathfrak D_Jf_\dy\bigr\|_{L^2_tL^2_xH^1_\sigma} ,
	\]
	which accounts for the first line of \eqref{eq:dyn-forcing-estimate}. The coefficients $\cV_{I,J}$ are given by \eqref{eq:mixed-transport-commutator}, and $[\mathfrak D_I,L]f_\dy=[\mathfrak D_{I'},L]\partial_t^{\,j(I)}f_\dy$ by \eqref{eq:dyn-normal-ordering} and \eqref{eq:dyn-spacetime-leibniz}, with $I'$ the spatial subword of $I$ and $[\mathfrak D_{I'},L]$ that of \eqref{eq:weak-triangular-commutator}.  The $x$-coefficients of $\cV_{I,J}$ and of $[\mathfrak D_{I'},L]$ are bounded on $\overline\cD$.  Since $L=\cA_\mu+\cK_0$, $[\mathfrak D_{I'},L]=[\mathfrak D_{I'},\cA_\mu]+[\mathfrak D_{I'},\cK_0]$.  At each fixed $t$, \eqref{eq:ordered-transport-weight-integrated}, \eqref{eq:local-L0-complete-bound}
	integrated in $x$, and \eqref{eq:iterated-compact-full-operator-bound} with
	\eqref{eq:weighted-frame-monotonicity} at $\mathrm L'=\tfrac12$ and $\mathrm L=\ell+A_\diamond$ give, respectively,
	\[
	\begin{aligned}
		\|w^\ell\cV_{I,J}\mathfrak D_Jf_\dy\|_{L^2_xH^{-1}_\sigma}
		&\le C\|w^{\ell+4}\mathfrak D_Jf_\dy\|_{L^2_xH^1_\sigma},
		\qquad |J|<m,\\
		\|w^\ell[\mathfrak D_{I'},\cA_\mu]\partial_t^{\,j(I)}f_\dy\|_{L^2_xH^{-1}_\sigma}
		&\le C\sum_{|J|<m}\|w^{\ell+2m+2}\mathfrak D_Jf_\dy\|_{L^2_xH^1_\sigma},\\
		\|w^\ell[\mathfrak D_{I'},\cK_0]\partial_t^{\,j(I)}f_\dy\|_{L^2_xH^{-1}_\sigma}
		&\le C\sum_{|J|<m}\|\mathfrak D_Jf_\dy\|_{L^2_{x,v}}
		\le C\sum_{|J|<m}\|w^{\ell+A_\diamond}\mathfrak D_Jf_\dy\|_{L^2_xH^1_\sigma}.
	\end{aligned}
	\]
	Since $4\le2m+2\le2\mathrm M+2<A_\diamond$, \eqref{eq:real-weight-Hsigma-monotonicity} bounds both weighted terms by the $\ell+A_\diamond$ norm. Summation over $|J|<m$ and integration in time give the asserted bound, independently of $T$.
	
	\noindent\emph{The $g_\dy$--$f_\dy$ terms.}
	For $(a,b)=(g_\dy,f_\dy)$, the hypotheses of Lemma~\ref{lem:dyn-collision-terms} follow from \eqref{eq:dyn-differentiated-chain} and \eqref{eq:dyn-differentiated-first-slots}; the jets of $g_\dy$ are understood distributionally.  Lemma~\ref{lem:dyn-collision-terms}, with $\cN_{m-1}\le\cN_m$ and with the weight $\ell$ in the first and last terms of the second factor of \eqref{eq:dyn-bilinear-allocation} raised to $\ell+A_\diamond$ by \eqref{eq:real-weight-Hsigma-monotonicity}, gives
	\[
	\begin{aligned}
		&\Bigl\|w^{\ell}\Bigl(\sum N\bigl[\mathfrak D_{I_1}g_\dy,\mathfrak D_{I_2}f_\dy\bigr]+\mathfrak R_I[g_\dy,f_\dy]\Bigr)\Bigr\|_{L^2_tL^2_xH^{-1}_\sigma}\\
		&\qquad\le C\Bigl(\|w^{-2}g_\dy\|_{L_t^\infty L_x^\infty L_v^2}+\cN_m(T)[g_\dy]+\mathsf G_\natural(T)[g_\dy]\Bigr)\\
		&\qquad\qquad\times\Bigl(\|w^{\ell+A_\diamond}f_\dy\|_{L_t^\infty L_x^{p_\sharp}H_\sigma^1}+\cN_m(T)[f_\dy]+\cS_T^{\ell+A_\diamond}[f_\dy]\Bigr) .
	\end{aligned}
	\]
	
	\noindent\emph{The $f_\st$--$f_\dy$ terms.}
	Let $I_1\sqcup I_2=I$ with $I_1\neq\varnothing$, and put $k:=|I_1|$.  The chain \eqref{eq:dyn-differentiated-chain} and finite coordinate norms permit Lemma~\ref{lem:dyn-slots} to be applied to $f_\dy$. Using \eqref{eq:dyn-bilinear-pointwise}, the first pairing of Lemma~\ref{lem:dyn-pairings}, \eqref{eq:dyn-stationary-slots}, and \eqref{eq:dyn-slot-coercive-L2} at order $m-k$, we obtain
	\[
	\begin{aligned}
		\bigl\|w^{\ell}N[\mathfrak D_{I_1}f_\st,\mathfrak D_{I_2}f_\dy]\bigr\|_{L^2_tL^2_xH^{-1}_\sigma}
		&\le C\bigl\|w^{-2}\mathfrak D_{I_1}f_\st\bigr\|_{L_x^{\mathsf g(k)}L_v^2}\,\bigl\|w^{\ell}\mathfrak D_{I_2}f_\dy\bigr\|_{L^2_tL_x^{\mathsf p_{m-k}(2)}H^1_\sigma}\\
		&\le C\delta\,\bigl(\cS_T^{\ell}[f_\dy]+\cN_m(T)[f_\dy]\bigr) .
	\end{aligned}
	\]
	For a term \eqref{eq:deformation-remainder-structure} of $\mathfrak R_I[f_\st,f_\dy]$, put $j:=|J_1|+|J_2|\le m-1$.  Then \eqref{eq:deformation-weighted-output-spatial} at $\varpi=\frac14$ replaces \eqref{eq:dyn-bilinear-pointwise}, with the weight $\ell+A_\diamond$ in the second factor, as in the proof of Lemma~\ref{lem:dyn-collision-terms}\,(\ref{it:dyn-collision-deformation}).  If $j\ge1$, the bound above holds with $(J_1,J_2)$ in place of $(I_1,I_2)$, with $j$ in place of $m$ in Lemma~\ref{lem:dyn-pairings}\,(\ref{it:dyn-pairing-first}) and in \eqref{eq:dyn-slot-coercive-L2}, and with the weight $\ell+A_\diamond$ in place of $\ell$, which is admissible in \eqref{eq:dyn-slot-coercive-L2} since $\ell+A_\diamond\le\ell^\dy_{m-1}\le\ell^\dy_j$; the right side of the resulting bound is then $C\delta\bigl(\cS_T^{\ell+A_\diamond}[f_\dy]+\cN_j(T)[f_\dy]\bigr)$, and $\cN_j\le\cN_m$.  If $j=0$, we apply \eqref{eq:deformation-weighted-output-spatial} with $(\infty,2)$ in place of its exponents $(q,r)$ and measure $w^{-2}f_\st$ in $L_x^\infty L_v^2$ and $w^{\ell+A_\diamond}f_\dy$ in $L^2_tL^2_xH^1_\sigma$, where
	\[
	\|w^{-2}f_\st\|_{L_x^\infty L_v^2}\le C\delta,
	\qquad
	\|w^{\ell+A_\diamond}f_\dy\|_{L^2_tL^2_xH^1_\sigma}\le C\cN_0(T)[f_\dy],
	\]
	by \eqref{eq:dyn-stationary-slots} and by \eqref{eq:dyn-slot-tops-rows} at $m=0$, since $\ell+A_\diamond\le\ell^\dy_0$.  Since $\cS_T^{\ell}\le C\cS_T^{\ell+A_\diamond}$ by \eqref{eq:real-weight-Hsigma-monotonicity}, this gives
	\[
	\Bigl\|w^{\ell}\Bigl(\sum N\bigl[\mathfrak D_{I_1}f_\st,\mathfrak D_{I_2}f_\dy\bigr]+\mathfrak R_I[f_\st,f_\dy]\Bigr)\Bigr\|_{L^2_tL^2_xH^{-1}_\sigma}
	\le C\delta\bigl(\cS_T^{\ell+A_\diamond}[f_\dy]+\cN_m(T)[f_\dy]\bigr) .
	\]
	
	\noindent\emph{The $g_\dy$--$f_\st$ terms.}
	Let $I_1\sqcup I_2=I$ and $k:=|I_1|$.  Lemma~\ref{lem:dyn-slots} applies to the distributional jets of $g_\dy$ under \eqref{eq:dyn-differentiated-first-slots}. The second pairing of Lemma~\ref{lem:dyn-pairings}, together with \eqref{eq:dyn-bilinear-pointwise}, \eqref{eq:dyn-slot-gaussian-L2}, and \eqref{eq:dyn-stationary-slots} (recall $\ell\le\ell^\dy_0$), gives
	\[
	\begin{aligned}
		\bigl\|w^{\ell}N[\mathfrak D_{I_1}g_\dy,\mathfrak D_{I_2}f_\st]\bigr\|_{L^2_tL^2_xH^{-1}_\sigma}
		&\le C\bigl\|w^{-2}\mathfrak D_{I_1}g_\dy\bigr\|_{L^2_tL_x^{\mathsf p_k(q_\flat)}L_v^2}\,\bigl\|w^{\ell}\mathfrak D_{I_2}f_\st\bigr\|_{L_x^{\mathsf s(m-k)}H^1_\sigma}\\
		&\le C\delta\bigl(\mathsf G_\natural(T)[g_\dy]+\cN_m(T)[g_\dy]\bigr) .
	\end{aligned}
	\]
	Every term of $\mathfrak R_I[g_\dy,f_\st]$ has the form \eqref{eq:deformation-remainder-structure}; for such a term put $j:=|J_1|+|J_2|\le m-1$.  As in the proof of Lemma~\ref{lem:dyn-collision-terms}\,(\ref{it:dyn-collision-deformation}), \eqref{eq:deformation-weighted-output-spatial} at $\varpi=\frac14$ replaces \eqref{eq:dyn-bilinear-pointwise}, with the weight $\ell+A_\diamond$ in the second factor.
	
	If $j\ge1$, the bound above holds with $(J_1,J_2)$ in place of $(I_1,I_2)$.  Here Lemma~\ref{lem:dyn-pairings}\,(\ref{it:dyn-pairing-second}) and \eqref{eq:dyn-slot-gaussian-L2} are used at the order $j$ in place of $m$.  In the second factor the weight $\ell+A_\diamond$ takes the place of $\ell$; it is admissible in \eqref{eq:dyn-stationary-slots}, since $\ell+A_\diamond\le\ell^\dy_{m-1}\le\ell^\dy_0=\ell_\sharp$.  The right side of the resulting bound is $C\delta\bigl(\mathsf G_\natural(T)[g_\dy]+\cN_j(T)[g_\dy]\bigr)$, and $\cN_j\le\cN_m$.
	
	For $j=0$ we use \eqref{eq:deformation-weighted-output-spatial} with $(q,r)=(P_\natural,p_\flat)$. The spatial H\"older condition $1/P_\natural+1/p_\flat\le1/2$ follows from Lemma~\ref{lem:dyn-pairings}\,(\ref{it:dyn-pairing-second}) at $k=0$, using $\mathsf s(0)=p_\flat\ge\mathsf s(m)$.  We measure $w^{-2}g_\dy$ in $L^2_tL_x^{P_\natural}L_v^2$, that is, by $\mathsf G_\natural(T)[g_\dy]$ of \eqref{eq:dyn-natural-base}.  We measure $w^{\ell+A_\diamond}f_\st$ in $L_x^{p_\flat}H^1_\sigma$, where $\|w^{\ell+A_\diamond}f_\st\|_{L_x^{p_\flat}H^1_\sigma}\le C\delta$ by \eqref{eq:dyn-stationary-slots} at $j=0$.
	
	Together with the bound for $N[\mathfrak D_{I_1}g_\dy,\mathfrak D_{I_2}f_\st]$, the cases $j=0$ and $j\ge1$ yield
	\[
	\Bigl\|w^{\ell}\Bigl(\sum_{I_1\sqcup I_2=I}N\bigl[\mathfrak D_{I_1}g_\dy,\mathfrak D_{I_2}f_\st\bigr]+\mathfrak R_I[g_\dy,f_\st]\Bigr)\Bigr\|_{L^2_tL^2_xH^{-1}_\sigma}
	\le C\delta\bigl(\mathsf G_\natural(T)[g_\dy]+\cN_m(T)[g_\dy]\bigr) .
	\]
\end{proof}

\subsection{Exponential decay and the initial data}\label{ss:dyn-decay-data}

We first prove exponential decay for the fixed-background problem, uniformly in the time horizon, and then construct the initial jets for the regularized evolution.  For $0\le t_0<t_1\le\infty$ and $f$ on $\cD_{(t_0,t_1)}\times\bR_v^3$, the dissipation rate is
\begin{equation}
	\label{eq:dyn-decay-density}
	\mathfrak d_f(r)
	:=\|f(r)\|_{L^2_xH^1_\sigma}^2
	+\|(\partial_r+v\cdot\nabla_x)f(r)\|_{L^2_xH^{-1}_\sigma}^2
	+\|(\mathbf I-\cP_0^+)f_{\gamma^+}(r)\|_{L^2(\gamma^+)}^2 ,
\end{equation}
and the dissipation is
\[
\cD_f(s,t):=\int_s^t\mathfrak d_f(r)\dd r ,
\qquad t_0\le s<t\le t_1 .
\]

For $z\in\cH_{\kin,\sigma}(\cD_{(t_0,t_1)}\times\bR_v^3)$,
Proposition~\ref{prop:full-diffuse-trace-recovery}\,(\ref{it:recovery-evolutionary})
gives
\begin{equation}
	\label{eq:dyn-decay-signed-trace-by-dissipation}
	\|z_{\gamma^+}\|_{L^2(\gamma^+_{(t_0,t_1)})}^2
	\le C\bigl(\mathsf J_z(t_0)-\mathsf J_z(t_1)+\cD_z(t_0,t_1)\bigr),
	\qquad 0\le t_0<t_1<\infty,
\end{equation}
where $\mathsf J_z:=C^{-1}\mathfrak e_z$, $\mathfrak e_z$ is the function of Proposition~\ref{prop:full-diffuse-trace-recovery}\,(\ref{it:recovery-evolutionary}), and $C\ge1$ is the constant of that proposition, enlarged if necessary; $C$ does not depend on $(t_0,t_1)$.  The proof of that proposition expresses $\mathfrak e_z$ as a fixed multiple of a finite sum of outgoing endpoint functionals \eqref{eq:nongrazing-endpoint-functional}. Their cutoffs $\psi$ and $\omega$ are independent of time and of $(t_0,t_1)$, so $\mathsf J_z(t)$ depends only on $z(t)$.  By
\eqref{eq:recovery-endpoint-bound}, $0\le\mathsf J_z(t)\le\|z(t)\|_{L^2_{x,v}}^2$.

The steady state satisfies
\begin{equation}
	\label{eq:dynamic-global-stationary-second-slot-small}
	\|w^{A^\dy_\flat}f_\st\|_{L^\infty}
	+\|f_\st\|_{L^{p_\flat}_xH^1_\sigma}
	\le C\delta ,
\end{equation}
by Proposition~\ref{prop:stationary-finite-weight-regularity} at $p=p_\flat$ and
$\ell\in\{0,A^\dy_\flat\}$.

\begin{lemma}[Forcing and graph estimates]
	\label{lem:dyn-decay-graph-forcing}
	There is $C>0$, depending only on the structural data (in particular
	independent of $\tau$, $t$ and $N$), such that, for $0\le\tau<t<\infty$, the following hold.
	\begin{enumerate}
		\item
		For every $z\in\cH_{\kin,\sigma}(\cD_{(\tau,t)}\times\bR_v^3)$ and every measurable
		$f$ on $\cD_{(\tau,t)}\times\bR_v^3$,
		\begin{equation}
			\label{eq:dyn-decay-forcing}
			\|N[z,f]\|_{L^2(\cD_{(\tau,t)};H^{-1}_\sigma)}
			\le C\,\|f\|_{L^\infty_tL^{p_\flat}_xH^1_\sigma(\cD_{(\tau,t)}\times\bR_v^3)}\,
			\|z\|_{\cH_{\kin,\sigma}(\cD_{(\tau,t)}\times\bR_v^3)} .
		\end{equation}
		\item
		If $b\in L^2(\gamma^-_{(\tau,t)})$ and
		$z\in\cH_{\kin,\sigma}(\cD_{(\tau,t)}\times\bR_v^3)$ satisfies the wall law
		\begin{equation}
			\label{eq:dyn-decay-diffuse-law}
			z_{\gamma^-}(r)=P(r)z_{\gamma^+}(r)+b(r)
			\qquad\text{on }\gamma^-,\quad \tau<r<t ,\qquad P\in\{P_w,P_{w,N}\},\ N\in\bN,
		\end{equation}
		then
		\begin{equation}
			\label{eq:dyn-decay-graph}
			\|z\|_{\cH_{\kin,\sigma}(\cD_{(\tau,t)}\times\bR_v^3)}
			\le C\Bigl(\|z(\tau)\|_{L^2_{x,v}}+\cD_z(\tau,t)^{1/2}
			+\|b\|_{L^2(\gamma^-_{(\tau,t)})}\Bigr).
		\end{equation}
	\end{enumerate}
\end{lemma}

\begin{proof}
	For \eqref{eq:dyn-decay-forcing}, apply \eqref{eq:dyn-bilinear-pointwise} with $\ell=0$, $|I_1|=|I_2|=0$, taking $z$ as the first argument of the collision operator. H\"older's inequality with $(q_\flat,p_\flat)$ from \eqref{eq:dyn-own-order-endpoints} and the $L_t^\infty$ bound on $f$ reduce the estimate to
	\[
	\|w^{-2}z\|_{L^2_tL^{q_\flat}_xL^2_v(\cD_{(\tau,t)}\times\bR_v^3)}
	\le C\|z\|_{\cH_{\kin,\sigma}(\cD_{(\tau,t)}\times\bR_v^3)} ,
	\]
	which is \eqref{eq:graph-to-mixed-integrability-evolutionary} at $\varrho=2$ and $k=0$.
	
	To prove \eqref{eq:dyn-decay-graph}, the graph norm bound \eqref{innerprod_sigma} on $(\tau,t)$ gives
	\[
	\|z\|_{\cH_{\kin,\sigma}(\cD_{(\tau,t)}\times\bR_v^3)}^2
	\le\|z(\tau)\|_{L^2_{x,v}}^2+\cD_z(\tau,t)+\|z_{\gamma^-}\|_{L^2(\gamma^-_{(\tau,t)})}^2 .
	\]
	By \eqref{eq:dyn-decay-diffuse-law} and \eqref{eq:wall-J-operator-bound},
	\[
	\|z_{\gamma^-}\|_{L^2(\gamma^-_{(\tau,t)})}
	\le C\|z_{\gamma^+}\|_{L^2(\gamma^+_{(\tau,t)})}+\|b\|_{L^2(\gamma^-_{(\tau,t)})} ,
	\]
	where the constant of \eqref{eq:wall-J-operator-bound} depends only on $\cD$ by Lemma~\ref{lem:wall-boundary-form}, hence is uniform in $r$ and the same for $P_w$ and $P_{w,N}$, $N\in\bN$.  By \eqref{eq:dyn-decay-signed-trace-by-dissipation} at $(t_0,t_1)=(\tau,t)$ and $0\le\mathsf J_z\le\|z\|_{L^2_{x,v}}^2$,
	\[
	\|z_{\gamma^+}\|_{L^2(\gamma^+_{(\tau,t)})}^2
	\le C\bigl(\|z(\tau)\|_{L^2_{x,v}}^2+\cD_z(\tau,t)\bigr) .
	\]
	These bounds give \eqref{eq:dyn-decay-graph}.
\end{proof}

\begin{prop}[Exponential stability]
	\label{prop:dynamic-global-exponential-stability}
	There are $\varepsilon_\star>0$ and $C$, depending only on the
	structural data, and $\lambda_{\rm exp}>0$, depending only on the structural data other than $\lambda_w$, such that, after decreasing $\delta_\star$, if necessary, the following holds.  Let $T\ge1$, $0<\lambda\le\min\{\lambda_{\rm exp},\lambda_w\}$,
	$0\le\varepsilon_\dy\le\varepsilon_\star$ and $\mathsf G\ge0$.  Let $g=f_\st+g_\dy$ be a frozen background with
	\begin{equation}
		\label{eq:dyn-decay-background-hypotheses}
		\|w^{A^\dy_\flat}g_\dy\|_{L^\infty}\le\varepsilon_\dy,
		\qquad
		\sup_{0\le t\le T}\ee^{\lambda t}
		\|g_\dy\|_{\cH_{\kin,\sigma}(\cD_{(t,T)}\times\bR_v^3)}^2\le\mathsf G,
	\end{equation}
	and let
	\begin{equation}
		\label{eq:dyn-decay-solution-hypotheses}
		f_\dy\in\cH_{\kin,\sigma}(\cD_T\times\bR_v^3)\ \text{solve }
		\eqref{eq:dyn-perturbation-problem},
		\qquad
		\mathfrak m_{f_\dy}(0)=0 ,
	\end{equation}
	where the wall law of \eqref{eq:dyn-perturbation-problem} may be replaced by
	\begin{equation}
		\label{eq:dyn-decay-general-wall}
		(f_\dy)_{\gamma^-}=P_\cJ(t)(f_\dy)_{\gamma^+}+b_w+c
		\qquad\text{on }\gamma_T^-.
	\end{equation}
	Here $P_\cJ(t)$ is the operator \eqref{eq:wall-contraction-operator} at $\Theta_w(t)$, with $\cJ\in\{\mathbf I,J_N\}$, $N\in\bN$, so that $P_\cJ\in\{P_w,P_{w,N}\}$, and
	$c\in L^2(\gamma_T^-)$ satisfies, for some $K_c\ge0$,
	\begin{equation}
		\label{eq:dyn-decay-extra-datum}
		\int_{\partial\cD}\beta_c(t,x)\dd S_x=0\quad\text{for a.e. }t,
		\qquad
		\int_s^T\|c(r)\|_{L^2(\gamma^-)}^2\dd r\le K_c\,\ee^{-\lambda s},
		\quad 0\le s\le T,
	\end{equation}
	with $\beta_c$ as in \eqref{eq:dynamic-boundary-mass-scalars}.  The wall law of
	\eqref{eq:dyn-perturbation-problem} is the case $\cJ=\mathbf I$, $c=0$, $K_c=0$.  Then
	\begin{equation}
		\label{eq:dynamic-global-exponential-stability-finite-horizon}
		\|f_\dy(t)\|_{L^2_{x,v}}^2+\cD_{f_\dy}(t,T)
		+\|(f_\dy)_{\gamma^+}\|_{L^2(\gamma^+_{(t,T)})}^2
		\le C\Bigl(\|f_{\dy,0}\|_{L^2_{x,v}}^2+\delta^2\bigl(1+\mathsf G\bigr)+K_c\Bigr)\ee^{-\lambda t}
	\end{equation}
	for $0\le t\le T$.  If \eqref{eq:dyn-decay-background-hypotheses}, \eqref{eq:dyn-decay-solution-hypotheses} and, in the case \eqref{eq:dyn-decay-general-wall}, \eqref{eq:dyn-decay-extra-datum} hold at $T=\infty$, then
	\eqref{eq:dynamic-global-exponential-stability-finite-horizon} holds for every
	$t\ge0$ with $\cD_{f_\dy}(t,\infty)$ and
	$\|(f_\dy)_{\gamma^+}\|_{L^2(\gamma^+_{(t,\infty)})}^2$ on the left, and with
	the same $C$.  In either case, if $b_w\equiv0$, the bound holds with $\delta^2\mathsf G$ in place of $\delta^2(1+\mathsf G)$.
\end{prop}

\begin{proof}
	Put $f=f_\dy$, $F=N[g,f]+N[g_\dy,f_\st]$, and $b=b_w+c$.
	Collision invariance, \eqref{eq:dyn-theta-residual-flux}, and the condition on
	$\beta_c$ preserve the zero-mass constraint.  Thus $\mathfrak m_f(t)=0$
	by Proposition~\ref{prop:dynamic-alpha-uniform-graph}.
	For $0\le s<t\le T$, write
	\[
	\begin{gathered}
		E(t)=\|f(t)\|_2^2,\qquad D=\cD_f(s,t),\qquad
		B=\|f_{\gamma^+}\|_{L^2(\gamma^+_{(s,t)})}^2,\\
		R=\|b\|_{L^2(\gamma^-_{(s,t)})}^2,\qquad
		S=\|F\|_{L^2((s,t)\times\cD;H^{-1}_\sigma)}^2,\qquad
		J(t)=\mathsf J_f(t).
	\end{gathered}
	\]
	Constants are uniform in $s,t,T$ and the spectral cutoff; the coercivity constants are independent of $\lambda_w$.  The data satisfy
	\begin{equation}\label{eq:dyn-decay-nu-b}
		R\le C(\delta^2+K_c)\ee^{-\lambda s};
	\end{equation}
	when $b_w=0$, the $\delta^2$ term is absent.
	
	Set
	\[
	H=\|(\mathbf I-\mathbf P)f\|_{L^2((s,t)\times\cD;H^1_\sigma)}^2
	+\|(\mathbf I-\cP_0^+)f_{\gamma^+}\|_{L^2(\gamma^+_{(s,t)})}^2
	+\mathfrak d_{\cJ,(s,t)}^2.
	\]
	The boundary Green identity \eqref{eq:dynamic-quadratic-Green-identity},
	coercivity \eqref{eq:L-coercivity}, and the wall decomposition
	$P_\cJ f_{\gamma^+}=p_\cJ+r_\cJ$ give, for a structural $a>0$,
	\begin{equation}\label{eq:dyn-decay-energy}
		E(t)+aH\le E(s)+2\int_s^t|\langle F,f\rangle|\dd r
		+C(\delta^2B+R)+2\int_s^t|\langle p_\cJ,c\rangle_{\gamma^-}|\dd r.
	\end{equation}
	Here $\langle p_\cJ,b_w\rangle_{\gamma^-}=0$ by
	\eqref{eq:dyn-theta-residual-flux}, while
	$\|r_\cJ\|_{L^2(\gamma^-)}\lesssim\delta\|f_{\gamma^+}\|_{L^2(\gamma^+)}$.
	
	The macroscopic estimate of Lemma~\ref{lem:dynamic-alpha-macroscopic-estimate}
	and the transport-tail and trace estimates
	\eqref{eq:transport-tail-L2-evolutionary},
	\eqref{eq:dyn-decay-signed-trace-by-dissipation} yield continuous functionals
	$\eta,\bar\eta$, depending only on $f(t)$, with
	$|\eta(t)|+|\bar\eta(t)|\lesssim E(t)$ and
	\begin{align}
		\int_s^t\|\mathbf Pf\|_2^2\dd r
		&\lesssim H+S+\delta^2B+R+\eta(s)-\eta(t),
		\label{eq:dyn-decay-macro}\\
		\int_s^t E(r)\dd r
		&\lesssim J(s)-J(t)+D+R+\bar\eta(t)-\bar\eta(s),
		\label{eq:dyn-decay-mass-input}\\
		B&\lesssim J(s)-J(t)+D,\qquad 0\le J(r)\le E(r).
		\label{eq:dyn-decay-trace-input}
	\end{align}
	The endpoint functionals are independent of $(s,t)$.  Moreover, the
	boundedness of $L:H^1_\sigma\to H^{-1}_\sigma$ and the equivalence of norms
	on the macroscopic space imply, from the equation for $f$,
	\begin{equation*}
		D\lesssim H+\int_s^t\|\mathbf Pf\|_2^2\dd r+S.
	\end{equation*}
	Finally, \eqref{260911445} for $N[g,f]$,
	Lemma~\ref{lem:dyn-decay-graph-forcing} for $N[g_\dy,f_\st]$,
	\eqref{eq:dyn-decay-background-hypotheses} and
	\eqref{eq:dynamic-global-stationary-second-slot-small} give
	\begin{equation*}
		S\lesssim(\varepsilon_\dy+\delta)^2D
		+\delta^2\mathsf G\ee^{-\lambda s}.
	\end{equation*}
	
	Choose a small structural $\rho>0$ and add $\rho$ times
	\eqref{eq:dyn-decay-macro} to \eqref{eq:dyn-decay-energy}.
	The graph reduction and Young's inequality give
	\[
	E(t)+\rho\eta(t)+aD
	\le E(s)+\rho\eta(s)+C(S+\delta^2B+R)
	+2\int_s^t|\langle p_\cJ,c\rangle_{\gamma^-}|\dd r
	\]
	for some smaller structural $a>0$.  Since
	$\|p_\cJ\|_{L^2(\gamma^-)}\le\|f_{\gamma^+}\|_{L^2(\gamma^+)}$,
	the last term is bounded by $\theta B+C_\theta\|c\|_{L^2(\gamma^-_{(s,t)})}^2$.
	Choose $\theta$ first, then $\delta$ and $\varepsilon_\dy$ sufficiently small to absorb the terms involving $S$ and $B$ into the dissipation.  Finally, add a small structural
	multiple $h>0$ of \eqref{eq:dyn-decay-mass-input}, choosing $h$ so
	that its contribution to $D$ is absorbed as well.  Thus for some
	fixed $\zeta,a>0$,
	\[
	\mathscr E(t)+aD+a\int_s^t E(r)\dd r
	\le \mathscr E(s)+C\bigl(\delta^2(1+\mathsf G)+K_c\bigr)\ee^{-\lambda s},
	\qquad
	\mathscr E=E+\rho\eta+\zeta J-h\bar\eta.
	\]
	Here $\zeta$ is the sum of the outgoing-trace coefficient in the
	first absorption and the coefficient from the $h$-multiple of
	\eqref{eq:dyn-decay-mass-input}.  Taking $\rho$ and $h$ sufficiently
	small gives
	$\tfrac12E\le\mathscr E\le CE$, because $0\le J\le E$ and
	$|\eta|+|\bar\eta|\lesssim E$.  Consequently,
	\begin{equation}\label{eq:dyn-decay-gronwall}
		\mathscr E(t)+aD+a\int_s^t\mathscr E(r)\dd r
		\le\mathscr E(s)
		+C\bigl(\delta^2(1+\mathsf G)+K_c\bigr)\ee^{-\lambda s},
	\end{equation}
	with structural constants independent of $\lambda_w$.
	
	Choose $\lambda_{\rm exp}>0$ sufficiently small in terms of the
	dissipation coefficient in \eqref{eq:dyn-decay-gronwall}.
	For $0<\lambda\le\min\{\lambda_{\rm exp},\lambda_w\}$,
	Lemma~\ref{lem:delayed-dissipation-decay}, with its forcing rate
	set equal to $\lambda$, yields
	$\mathscr E(t)\lesssim\bigl(\mathscr E(0)+\delta^2(1+\mathsf G)+K_c\bigr)\ee^{-\lambda t}$.
	For finite $T$, extend $\mathscr E$ by zero beyond $T$. Inequality \eqref{eq:dyn-decay-gronwall} remains valid for this extension since $\mathscr E(T)\ge0$.  Applying the same inequality on $(t,T)$
	controls $D$, and \eqref{eq:dyn-decay-trace-input} controls the
	outgoing trace.  Since $\mathscr E(0)\lesssim\|f_{\dy,0}\|_2^2$,
	this proves \eqref{eq:dynamic-global-exponential-stability-finite-horizon}.
	The constants are uniform in $T$, so the estimate passes to $T=\infty$.
	When $b_w=0$, \eqref{eq:dyn-decay-nu-b} contains no separate
	$\delta^2$ term, and the last assertion follows.
\end{proof}

The general wall law \eqref{eq:dyn-decay-general-wall} is needed in Corollary~\ref{cor:dyn-decay-B}, where $\cJ=J_N$.

From now on we fix the rate
\begin{equation}
	\label{eq:dyn-decay-rate-fixed-point}
	\lambda:=\min\{\lambda_{\rm exp},\lambda_w\},
\end{equation}
where $\lambda_{\rm exp}$ is fixed by Proposition~\ref{prop:dynamic-global-exponential-stability}, independently of $\lambda_w$. The chosen $\lambda$ also depends on the wall decay rate $\lambda_w$.

The equation \eqref{eq:dyn-physical-problem} determines the time derivatives of $f_\pt$ at $t=0$ from the initial datum. We realize them by a lift $\check E$ on $[-1,0]$, vanishing near $t=-1$, so that the regularized centered problem can be solved on $(-1,T)$ with zero initial value at $t=-1$. 

Let
\[
I_T^-:=(-1,T),
\qquad
\cD_{T,\ast}:=I_T^-\times\cD,
\qquad
\gamma_{T,\ast}^\pm:=I_T^-\times\gamma^\pm ,
\]
and let
\begin{equation}
	\label{eq:dynamic-common-cutoff}
	\varkappa\in C^\infty(\bR),
	\qquad
	0\le\varkappa\le1,
	\qquad
	\varkappa=0\ \text{ on }[-1,-\tfrac34]\cup[\mathrm M+1,\infty),
	\qquad
	\varkappa=1\ \text{ on }[-\tfrac14,\mathrm M-1];
\end{equation}
write $\chi_-:=\varkappa|_{[-1,0]}$ and $\chi_+:=\varkappa|_{[0,\infty)}$.

On $[-1,0]$ we extend the wall temperature by its Taylor
polynomial,
\begin{equation}
	\label{eq:dyn-theta-past-extension}
	\vartheta(t,x):=\sum_{k=0}^{\mathrm N}\frac{t^k}{k!}\,\partial_t^k\vartheta(0^+,x),
	\qquad -1\le t\le0 .
\end{equation}
By \eqref{eq:dyn-theta-amplitude}, the bounded derivative $\partial_t^{k+1}\vartheta$ ensures the existence of $\partial_t^k\vartheta(0^+)$ in $C^{K_\cD}(\partial\cD)$ for $k\le\mathrm N$. The extension therefore matches derivatives through order $\mathrm N$ at $t=0$, while $\partial_t^{\mathrm N+1}\vartheta=0$ on $[-1,0)$.

\begin{lemma}
	\label{lem:dyn-wall-extended}
	Let $\vartheta$ be extended to $[-1,0]$ by \eqref{eq:dyn-theta-past-extension}.
	\begin{enumerate}
		\item\label{it:dyn-wall-ext-bounds}
		The extension satisfies
		\begin{align}\label{26101591}
			\max_{0\le j\le\mathrm N+1}\ \sup_{-1\le t<0}\
			\|\partial_t^j\vartheta(t,\cdot)\|_{C^{K_\cD}(\partial\cD)}
			\le\sum_{i=0}^{\mathrm N}\frac1{i!}\le\ee .
		\end{align}
		Hence the pair $(\ee\delta,\ee^{-1}\vartheta)$ satisfies \eqref{eq:dyn-theta-amplitude} with $t>-1$, $t\ne0$, in place of $t>0$.
		\item\label{it:dyn-wall-ext-results}
		For $-1<t<T$, $t\ne0$, the bounds \eqref{eq:wall-profile-chart-bound} and \eqref{eq:wall-profile-complete-bound} of Lemma~\ref{lem:wall-profile-derivatives} and the bounds \eqref{eq:dyn-theta-defect-time-bound} and the first bound in \eqref{eq:dyn-theta-defect-complete-bound} of Lemma~\ref{lem:dyn-theta-profile-time-derivatives} hold with $\ee C$ and $\ee C_{I,L}$ in place of $C$ and $C_{I,L}$, and $\partial_t^kd_\delta$ is continuous across $t=0$ for $k\le\mathrm N$.  Propositions~\ref{prop:dynamic-alpha-uniform-graph} and~\ref{prop:dynamic-inhomogeneous-physical-solvability} hold with $I_T^-$ in place of $(0,T)$.
		\item\label{it:dyn-wall-ext-difference}
		For $-1<t<T$,
		\begin{equation}
			\label{eq:dynamic-spectral-Pw-difference-exact}
			\bigl(P_w(t)-P_{w,N}(t)\bigr)h_{\gamma^+}
			=\bigl(\mathfrak m_\circ+d_\delta(t)\bigr)(\mathbf I-J_N)\cA_0h_{\gamma^+} .
		\end{equation}
		\item\label{it:dyn-wall-ext-profile}
		For $\psi\in L^2(\partial\cD)$, every finite $L\ge0$ and every integer $0\le k\le\mathrm M+1$,
		\begin{equation}
			\label{eq:dynamic-spectral-profile-L2-isometry-bound}
			\bigl\|w^L\bigl(\mathfrak m_\circ+d_\delta(t)\bigr)\psi\bigr\|_{L^2(\gamma^-)}
			\le C_L\|\psi\|_{L^2(\partial\cD)},
			\qquad
			\bigl\|w^L\bigl(\partial_t^kd_\delta\bigr)\psi\bigr\|_{L^2(\gamma^-)}
			\le C_L\,\delta\,\|\psi\|_{L^2(\partial\cD)},
		\end{equation}
		with $C_L$ independent of $t\in(-1,T)$.
	\end{enumerate}
\end{lemma}

\begin{proof}
	By the choice of $\delta_\star$, the extended wall amplitude satisfies
	$\ee\delta\le\min\{\tfrac12,\delta_{\rm dyn}\}$.
	
	(\ref{it:dyn-wall-ext-bounds}) For $k\le\mathrm N$, $\|\partial_t^k\vartheta(0^+)\|_{C^{K_\cD}(\partial\cD)}\le1$ by \eqref{eq:dyn-theta-amplitude}, and $|t|\le1$; this gives \eqref{26101591}.  By \eqref{26101591} for $-1<t<0$ and by \eqref{eq:dyn-theta-amplitude} for $t>0$,
	\[
	\max_{0\le j\le\mathrm N+1}\ \sup_{t>-1,\,t\ne0}\
	\|\partial_t^j(\ee^{-1}\vartheta)(t,\cdot)\|_{C^{K_\cD}(\partial\cD)}\le1 .
	\]
	
	(\ref{it:dyn-wall-ext-results}) Since $\Theta_w=\tfrac12\bigl(1+(\ee\delta)(\ee^{-1}\vartheta)\bigr)$, Lemmas~\ref{lem:wall-profile-derivatives} and~\ref{lem:dyn-theta-profile-time-derivatives} and Propositions~\ref{prop:dynamic-alpha-uniform-graph} and~\ref{prop:dynamic-inhomogeneous-physical-solvability} apply, by (\ref{it:dyn-wall-ext-bounds}) and $\ee\delta\le\delta_{\rm dyn}$, with $(\ee\delta,\ee^{-1}\vartheta)$ in place of $(\delta,\vartheta)$; since the bounds of Lemmas~\ref{lem:wall-profile-derivatives} and~\ref{lem:dyn-theta-profile-time-derivatives} listed in (\ref{it:dyn-wall-ext-results}) are linear in $\delta$, their constants $C$ and $C_{I,L}$ become $\ee C$ and $\ee C_{I,L}$.  Moreover, $d_\delta(t)=\mathfrak m_{\Theta_w(t)}-\mathfrak m_\circ$, the map $\Theta\mapsto\mathfrak m_\Theta$ is smooth on $[\tfrac14,\tfrac34]$, and $\partial_t^k\vartheta$ is continuous across $t=0$ for $k\le\mathrm N$; hence so is $\partial_t^kd_\delta$.
	
	(\ref{it:dyn-wall-ext-difference}) This follows from \eqref{eq:boundary-spectral-diffuse-operator} and \eqref{eq:dyn-theta-wall-split}.
	
	(\ref{it:dyn-wall-ext-profile}) As $\psi$ is independent of $v$, \eqref{eq:wall-profile-chart-bound} and \eqref{eq:dyn-theta-defect-time-bound} at $K_v=K_x=0$, together with $\mathfrak m_\circ=c_\mu^{-1}\sqrt\mu$, give the bounds for $t\ne0$. Continuity from (\ref{it:dyn-wall-ext-results}) covers $t=0$, since $k\le\mathrm M+1=\mathrm N$.
\end{proof}

The extension defines $d_\delta$, $P_w(t)$, $P_{w,N}(t)$ and the triangular terms $\mathfrak t_{j,N}[f]$ on $(-1,T)$. The prescribed data $h_w$ and $b_w$ are used only for $t>0$; on $(-1,0)$ the incoming datum is determined by the lift $\check E$ through \eqref{eq:dynamic-negative-time-spectral-inflow}.

With
\begin{equation}
	\label{eq:dynamic-initial-p-u}
	z_0:=f_0,
	\qquad\text{so that}\qquad
	f_{\dy,0}=z_0-f_\st ,
\end{equation}
define $z_1,\ldots,z_{\mathrm M}$ recursively by
\begin{equation}
	\label{eq:dynamic-generated-q0}
	z_j:=\sum_{a=0}^{j-1}\binom{j-1}{a}
	N[z_a,z_{j-1-a}]
	-\bigl(v\cdot\nabla_x+L\bigr)z_{j-1}.
\end{equation}
For a smooth solution of \eqref{eq:dyn-physical-problem} with initial value $z_0$, the recursion gives $z_j=\partial_t^jf_\pt(0)$. Inductively, $z_j$ is a sum of velocity-nonlocal multilinear expressions of degrees $1$ through $j+1$ in derivatives of $z_0$: the linear term in \eqref{eq:dynamic-generated-q0} has degree at most $j$, while the quadratic terms have degree at most $(a+1)+(j-a)=j+1$.  Each step of \eqref{eq:dynamic-generated-q0} adds at most one $x$-derivative, through $v\cdot\nabla_x$, or two $v$-derivatives, through $L$ and $N$.  Thus
\begin{equation}
	\label{eq:dynamic-generated-r0}
	z_j=\cP_j\bigl[(\partial_x^\alpha\partial_v^\beta z_0)_{|\alpha|\le j,\,|\beta|\le2j}\bigr],
	\qquad
	\deg\cP_j\le j+1 .
\end{equation}
We write
\begin{equation*}
	z:=\bigl(z_0,z_1,\ldots,z_{\mathrm M}\bigr)
\end{equation*}
and call $z_0,\ldots,z_{\mathrm M}$ the initial jets.
The compatibility conditions are
\begin{equation}
	\label{eq:dynamic-initial-u-compatibility}
	(z_j)_{\gamma^-}=P_w(0)(z_j)_{\gamma^+}
	+\mathfrak t_j[z]+\partial_t^{\,j}h_w(0),
	\qquad
	\mathfrak t_j[z]:=\sum_{k=1}^{j}\binom jk
	\bigl(\partial_t^{\,k}d_\delta\bigr)(0)\,\cA_0(z_{j-k})_{\gamma^+},
\end{equation}
for $0\le j\le\mathrm M$, where $\partial_t^{\,j}h_w(0)\in L^2(\gamma^-)$ is well defined, since
$h_w\in H^{\mathrm M+1}((0,T);L^2(\gamma^-))$ by \eqref{eq:dyn-wall-datum-class} and $j\le\mathrm M$. These are precisely the $j$th time derivatives of the incoming wall law at $t=0$: the derivatives falling on $P_w(t)=(\mathfrak m_\circ+d_\delta(t))\cA_0$ give $\mathfrak t_j[z]$. The
case $j=0$ reads
\begin{equation*}
	(z_0)_{\gamma^-}=P_w(0)(z_0)_{\gamma^+}+h_w(0) .
\end{equation*}

The zero-mass condition on the data is
\begin{equation}
	\label{eq:dynamic-initial-zero-perturbation-mass}
	\int_\cD\int_{\bR^3}f_{\dy,0}\sqrt\mu\dd v\dd x=0 .
\end{equation}
We impose this condition only in results that require it.

The negative-time lift $\check{E}$ and its source $\check{F}$ are given by
\begin{align}
	\check{\cP}(t)&:=\sum_{j=0}^{\mathrm M}\frac{t^j}{j!}\,z_j,
	\label{eq:dynamic-negative-time-Taylor-polynomial}\\
	\check{E}(t)&:=\chi_-(t)\check{\cP}(t),\qquad -1\le t\le0,
	\label{eq:dynamic-negative-time-lift}\\
	\check{F}&:=(\partial_t+v\cdot\nabla_x)\check{E}+L\check{E}-N[\check{E},\check{E}].
	\label{eq:dynamic-negative-time-source}
\end{align}
Write $\bar g_0$ for the extension of the background to $\cD_{T,\ast}$,
\begin{equation}
	\label{eq:dynamic-base-background}
	\bar g_0(t):=
	\begin{cases}
		\check{E}(t),&-1<t<0,\\
		g(t),&0<t<T,
	\end{cases}
	\qquad
	\bar g_0\big|_{(-1,-3/4]}=0 ,
\end{equation}
where the second identity holds because $\chi_-=0$ on $[-1,-3/4]$.  For the coefficient bounds of Lemma~\ref{lem:dyn-strong-coefficient-class}, which are taken over $(-\infty,T)$, we
extend $\bar g_0$ by zero to $t\le-1$:
\begin{equation}
	\label{eq:dyn-strong-past-extension}
	\hat g:=g\,\mathbf 1_{(0,T)}+\check{E}\,\mathbf 1_{(-1,0]},
\end{equation}
so that $\hat g=\bar g_0$ on $I_T^-$ and $\hat g=0$ for $t\le-3/4$.  There is no jump at
$t=-1$, nor at $t=0$ when $g(0)=z_0$, since $\check{E}(0)=z_0$ by $\chi_-(0)=1$.  For an object defined on $I_T^-$ and vanishing on $(-1,-3/4]$, such as $g^\varepsilon$ and $f_\pt^{\varepsilon,N}$ below, the hat denotes its extension by zero to $t\le-1$.

The following gluing lemma joins graph-space functions across a common $L^2$ time trace. We use it at $b=0$.

\begin{lemma}
	\label{lem:dyn-time-gluing}
	Let $a<b<c$, and let $u$ be a function on $(a,c)\times\cD\times\bR_v^3$ such that
	\[
	u_-:=u|_{(a,b)}\in\cH_{\kin,\sigma}((a,b)\times\cD\times\bR_v^3),
	\qquad
	u_+:=u|_{(b,c)}\in\cH_{\kin,\sigma}((b,c)\times\cD\times\bR_v^3)
	\]
	and $u_-(b)=u_+(b)$ in $L^2(\cD\times\bR_v^3)$.  Then $u\in\cH_{\kin,\sigma}((a,c)\times\cD\times\bR_v^3)$ and
	\begin{equation}
		\label{eq:dyn-time-gluing}
		\|u\|_{\cH_{\kin,\sigma}((a,c)\times\cD\times\bR_v^3)}^2
		\le\|u_-\|_{\cH_{\kin,\sigma}((a,b)\times\cD\times\bR_v^3)}^2
		+\|u_+\|_{\cH_{\kin,\sigma}((b,c)\times\cD\times\bR_v^3)}^2 .
	\end{equation}
\end{lemma}

\begin{proof}
	Adding the Green identities \eqref{eq:dynamic-Green-identity} on $(a,b)$ and $(b,c)$ cancels the common endpoint trace, so $u$ satisfies \eqref{2607151226} on $(a,c)$ with the piecewise transport derivative and boundary traces. Proposition~\ref{H_kin=weak_landau}\,(1) gives the graph-space membership, and \eqref{innerprod_sigma} yields \eqref{eq:dyn-time-gluing}.
\end{proof}

The size of the data is
\begin{equation}
	\label{eq:dyn-data-size}
	\mathsf I:=\|f_{\dy,0}\|_{L^2_{x,v}}+\sum_{|\alpha|\le2\mathrm M+2}\ \sum_{|\beta|\le3\mathrm M+2}
	\|w^{\check{A}}\partial_x^\alpha\partial_v^\beta z_0\|_{L^\infty(\cD\times\bR^3)} .
\end{equation}
Henceforth, through \S\ref{sec:dyn-realization-section}, assume the compatibility conditions \eqref{eq:dynamic-initial-u-compatibility} and $\mathsf I\le1$.  The boundary scalars of the jets are
$\varphi_j^0:=\cA_0(z_j)_{\gamma^+}$, $0\le j\le\mathrm M$.  

\begin{lemma}
	\label{lem:dynamic-zero-initial-source-inflow-lift}
	Let $\check{E}$ and $\check{F}$ be defined by
	\eqref{eq:dynamic-negative-time-Taylor-polynomial}--%
	\eqref{eq:dynamic-negative-time-source} from $z$.
	\begin{enumerate}
		\item\label{it:dyn-lift-weights}
		For every $L\ge0$ with $L+\mathrm MA_\diamond\le \check{A}$,
		\begin{equation}
			\label{eq:dyn-breg-lift-weights}
			\bigl\|w^{L}\partial_x^\alpha\partial_v^\beta z_j\bigr\|
			_{L^\infty(\cD\times\bR^3)}
			\le C\,\mathsf I,
			\qquad |\alpha|+|\beta|\le\mathrm M+2,
			\quad 0\le j\le\mathrm M .
		\end{equation}
		\item\label{it:dyn-data-jets}
		The jets satisfy
		\begin{equation}
			\label{eq:dyn-data-jets}
			\sum_{j\le\mathrm M}\sum_{|I|\le\mathrm M}
			\Bigl(\|w^{\ell^\dy_0}\mathfrak D_Iz_j\|_{\cH_{\kin,\sigma}}
			+\|w^{A^\dy_\flat}\mathfrak D_Iz_j\|_{L^\infty}\Bigr)
			\le C\,\mathsf I
		\end{equation}
		and
		\begin{equation}
			\label{eq:dyn-data-graph}
			\sum_{j\le\mathrm M}\|z_j\|_{\cH_{\kin,\sigma}(\cD\times\bR_v^3)}\le C\,\mathsf I ;
		\end{equation}
		in particular, by \eqref{260902137},
		\begin{equation}
			\label{eq:dyn-lift-size-relations}
			\sum_{j\le\mathrm M}\Bigl(\|z_j\|_{L^2_{x,v}}+\|z_j\|_{L_x^2H_\sigma^1}
			+\|v\cdot\nabla_xz_j\|_{L_x^2H_\sigma^{-1}}+\|(z_j)_{\gamma^\pm}\|_{L^2(\gamma^\pm)}\Bigr)
			\le C\,\mathsf I ,
		\end{equation}
		and the boundary scalars satisfy
		\begin{equation}
			\label{eq:dyn-data-scalar-H3}
			\varphi_j^0\in H^{\mathrm M}(\partial\cD),
			\qquad 0\le j\le\mathrm M .
		\end{equation}
		\item\label{it:dyn-lift-basic}
		The lift satisfies
		\begin{align}
			\check{E}&=0
			&&\text{on }[-1,-3/4],
			\label{eq:dynamic-negative-time-zero-near-minus-one}\\
			\partial_t^j\check{E}(0)&=z_j,
			&&0\le j\le\mathrm M,
			\label{eq:dynamic-negative-time-jet-matching}
		\end{align}
		and $\check{E}\in\cH_{\kin,\sigma}((-1,0)\times\cD\times\bR_v^3)$ with
		\begin{equation}
			\label{eq:dyn-lift-graph}
			\|\check{E}\|_{\cH_{\kin,\sigma}((-1,0)\times\cD\times\bR_v^3)}\le C\mathsf I .
		\end{equation}
		Moreover
		\begin{equation}
			\partial_t^j\check{F}(0)=0,
			\qquad 0\le j\le\mathrm M-1.
			\label{eq:dynamic-negative-time-source-matching}
		\end{equation}
		\item\label{it:dyn-breg-lift-bounds}
		For $0\le L\le\Lambda_0+B^{\dy}_{\rm in}+6$,
		\begin{equation*}
			\begin{aligned}
				&\|w^L\check{F}\|_{L^\infty((-1,0)\times\cD\times\bR^3)}
				+\sum_{j\le\mathrm M}\|w^L\partial_t^{\,j}\check{E}\|_{L^\infty((-1,0)\times\cD\times\bR^3)}
				+\|w^L(\check{E})_{\gamma^-}\|_{(C^{0,1}_{t,x}\cap C^{1,1}_v)((-1,0)\times\gamma^-)}\\
				&\qquad+\sum_{|I|+m\le\mathrm M}
				\|w^L\mathfrak D_I^-\partial_t^m(\check{E})_{\gamma^-}\|_{L^2((-1,0)\times\gamma^-)}
				+\sum_{|I|\le\mathrm M}\sum_{j\le\mathrm M}
				\|w^L\mathfrak D_I\partial_t^j\check{F}\|_{L^2((-1,0)\times\cD\times\bR^3)}
				\le C\mathsf I .
			\end{aligned}
		\end{equation*}
		\item\label{it:dyn-breg-source-bounds}
		For every $L\ge0$,
		\begin{equation}
			\label{eq:dyn-breg-source-bounds}
			\|w^L\check{F}\|_{L^2((-1,0)\times\cD;H^{-1}_\sigma)}
			\le C\|w^{L+3}\check{F}\|_{L^\infty((-1,0)\times\cD\times\bR^3)} .
		\end{equation}
	\end{enumerate}
\end{lemma}

\begin{proof}
	(\ref{it:dyn-lift-weights}) By \eqref{eq:dynamic-generated-r0}, a derivative $\partial_x^\alpha\partial_v^\beta z_j$ with
	$|\alpha|+|\beta|\le\mathrm M+2$ and $j\le\mathrm M$ uses only
	$\partial_x^{\alpha'}\partial_v^{\beta'}z_0$ with
	$$
	|\alpha'|\le|\alpha|+j\le2\mathrm M+2\quad\text{and}\quad |\beta'|\le|\beta|+2j\le3\mathrm M+2,
	$$
	which occur in $\mathsf I$.  For $\ell\ge0$ and functions $u,u_1,u_2$ on $\cD\times\bR^3$ for which the right sides are finite,
	\begin{align}
		\|w^\ell\partial_x^\alpha\partial_v^\beta(v\cdot\nabla_x+L)u\|_{L^\infty}
		&\le C\sum_{\substack{|\alpha'|\le|\alpha|+1,\ |\beta'|\le|\beta|+2\\|\alpha'|+|\beta'|\le|\alpha|+|\beta|+2}}\|w^{\ell+1}\partial_x^{\alpha'}\partial_v^{\beta'}u\|_{L^\infty},
		\label{eq:dyn-jet-transport-bound}\\
		\|w^\ell\partial_x^\alpha\partial_v^\beta N[u_1,u_2]\|_{L^\infty}
		&\le C\sum_{\substack{|\alpha'|,|\alpha''|\le|\alpha|\\|\beta'|\le|\beta|+1,\ |\beta''|\le|\beta|+2}}\|\partial_x^{\alpha'}\partial_v^{\beta'}u_1\|_{L^\infty}\|w^{\ell+1}\partial_x^{\alpha''}\partial_v^{\beta''}u_2\|_{L^\infty},
		\label{eq:dyn-jet-collision-bound}
	\end{align}
	with norms over $\cD\times\bR^3$ and $C$ depending only on $\ell$, $|\alpha|$ and $|\beta|$.  These follow from \eqref{def_Ag}--\eqref{def:L_g}, \eqref{eq:N-divergence-expansion}, \eqref{eq:perturbation-sigma-bound} and \eqref{26720259}:
	\[
	Lu=-\partial_{v_i}(\sigma^{ij}\partial_{v_j}u)+(v_iv_j\sigma^{ij}-12\pi\mu)u-\cK u,
	\qquad
	\cK u=\sqrt\mu(4v_iv_j-2\delta_{ij})\sigma^{ij}_{\sqrt\mu u},
	\qquad
	\partial_v^\beta\sigma^{ij}_G=\sigma^{ij}_{\partial_v^\beta G},
	\]
	and the coefficients of $L+\cK$ and of $N[u_1,\cdot]$ are sums of products of polynomials of degree at most $2$ in $v$ with $\mu$, $\sqrt\mu u_1$, $\sigma^{ij}_G$ or $\partial_v\sigma^{ij}_G$, $G\in\{\mu,\sqrt\mu u_1\}$.  We apply \eqref{eq:dyn-jet-transport-bound} and \eqref{eq:dyn-jet-collision-bound} $j$ times to \eqref{eq:dynamic-generated-q0}, at weights $\ell\le L+j-1$.  Then every factor $\partial_x^{\alpha'}\partial_v^{\beta'}z_0$ is measured in $L^\infty(\cD\times\bR^3)$ at a weight at most $L+j\le L+jA_\diamond\le \check{A}$, and \eqref{eq:dyn-data-size} gives
	\[
	\|w^{L+jA_\diamond}\partial_x^{\alpha'}\partial_v^{\beta'}z_0\|_{L^\infty(\cD\times\bR^3)}\le\mathsf I .
	\]
	Since $\cP_j$ of \eqref{eq:dynamic-generated-r0} has degree at most $j+1$ and no constant term, and $\mathsf I\le1$, this gives
	\eqref{eq:dyn-breg-lift-weights}.
	
	(\ref{it:dyn-data-jets}) A temporal letter of $\mathfrak D_I$ annihilates the $t$-independent $z_j$, so we may take $I$ spatial.  For a spatial word $I$, $\mathfrak D_I$ is a composition of $|I|$ fields $\mathfrak D_\nu=\mathfrak D_{\tau_\nu}$ of
	\eqref{eq:complete-phase-space-derivative}, where $\tau_\nu$ is compactly supported and of class $C^{K_\cD+1}$ by \eqref{eq:tangential-field-family}, and $|I|\le\mathrm M<K_\cD$ by \eqref{eq:order-and-domain-regularity}.  Hence $\mathfrak D_I$ is a sum of the $\partial_x^\alpha\partial_v^\beta$, $|\alpha|+|\beta|\le|I|$, with coefficients bounded by $Cw^{|I|}$ together with their first derivatives, so that, for $L\ge0$ and with norms over $\cD\times\bR^3$,
	\begin{align}
		\label{eq:dyn-field-derivative-bound}
		\begin{gathered}
			\|w^L\mathfrak D_Iu\|_{L^\infty}\le C\sum_{|\alpha|+|\beta|\le|I|}\|w^{L+|I|}\partial_x^\alpha\partial_v^\beta u\|_{L^\infty},\\
			\|w^L\nabla_{x,v}\mathfrak D_Iu\|_{L^\infty}\le C\sum_{|\alpha|+|\beta|\le|I|+1}\|w^{L+|I|}\partial_x^\alpha\partial_v^\beta u\|_{L^\infty} .
		\end{gathered}
	\end{align}
	By (\ref{it:dyn-lift-weights}) at $L=A^\dy_\flat+\mathrm M+4$,
	which is admissible by the choice of $\check A$, we have
	\[
	\|w^{A^\dy_\flat+\mathrm M+4}\partial_x^\alpha\partial_v^\beta z_j\|_{L^\infty(\cD\times\bR^3)}\le C\mathsf I ,
	\qquad |\alpha|+|\beta|\le\mathrm M+2,\quad 0\le j\le\mathrm M .
	\]
	Hence $u:=\mathfrak D_Iz_j$ satisfies
	\[
	\|w^{A^\dy_\flat+3}u\|_{L^\infty(\cD\times\bR^3)}
	+\|w^{A^\dy_\flat+3}\nabla_xu\|_{L^\infty(\cD\times\bR^3)}
	+\|w^{A^\dy_\flat+3}\nabla_vu\|_{L^\infty(\cD\times\bR^3)}
	\le C\mathsf I ,
	\]
	and $u$, $\nabla_xu$ and $\nabla_vu$ are continuous up to
	$\partial\cD$, since the derivatives of order at most $\mathrm M+1$ of $z_j$ are
	Lipschitz in $(x,v)$ on bounded sets.  By \eqref{eq:auxiliary-coercive-equivalence} with $\sigma$ bounded, \eqref{eq:Hsigma-dual-weighted},
	$\|w^{-2}\|_{L^2_v}+\||v|^{1/2}w^{-3}\|_{L^2_v}<\infty$ and $|\cD|+|\partial\cD|<\infty$,
	\[
	\|w^{\ell^\dy_0}u\|_{L^2(\cD\times\bR^3_v)}+\|w^{\ell^\dy_0}u\|_{L^2(\cD;H^1_\sigma)}+\|v\cdot\nabla_x(w^{\ell^\dy_0}u)\|_{L^2(\cD;H^{-1}_\sigma)}+\|(w^{\ell^\dy_0}u)_{\gamma}\|_{L^2(\gamma)}\le C\mathsf I .
	\]
	The Green formula on $\cD$
	gives \eqref{eq:H_kin_sigma-weak-stationary} for $w^{\ell^\dy_0}u$.  So
	Proposition~\ref{H_kin=weak_landau}\,(2) gives $w^{\ell^\dy_0}u\in\cH_{\kin,\sigma}$, and
	\eqref{eq:dyn-data-jets} follows.
	
	By Lemma~\ref{lem:compact-velocity-support}\,(1) on $\cH_{\kin,\sigma}(\cD\times\bR_v^3)$ at $\zeta=w^{-\ell^\dy_0}$ and \eqref{eq:dyn-data-jets} at $I=\varnothing$, $\|z_j\|_{\cH_{\kin,\sigma}}\le C\|w^{\ell^\dy_0}z_j\|_{\cH_{\kin,\sigma}}\le C\mathsf I$, which gives \eqref{eq:dyn-data-graph}, and \eqref{eq:dyn-lift-size-relations} follows from \eqref{eq:dyn-data-graph} and \eqref{260902137}.
	It remains to prove \eqref{eq:dyn-data-scalar-H3}. As above,
	\[
	\|w^{A^\dy_\flat}\partial_x^\alpha\partial_v^\beta z_j\|_{L^\infty(\cD\times\bR^3)}\le C\mathsf I ,
	\qquad |\alpha|+|\beta|\le\mathrm M+1,
	\]
	and these derivatives of $z_j$ are continuous up to $\partial\cD$.  On $\partial\cD\cap U_\kappa$, with the frame $R_\kappa$ of Proposition~\ref{prop:complete-fields}\,(\ref{it:cf-frame}), \[
	\varphi_j^0(x)=\int_{\xi_3>0}z_j(x,R_\kappa(x)\xi)\sqrt{\mu(\xi)}\,\xi_3\dd\xi ,
	\] and $R_\kappa$ is of class $C^{K_\cD+1}$ with $K_\cD\ge2\mathrm N$ by \eqref{eq:order-and-domain-regularity}.  So a tangential derivative of order at most $\mathrm M$ of $\varphi_j^0$ is a finite sum of integrals of the same form, in which $z_j$ is replaced by a derivative of $z_j$ in $(x,v)$ of order at most $\mathrm M$ and $\sqrt{\mu(\xi)}\,\xi_3$ is multiplied by a polynomial in $\xi$ with bounded coefficients. Hence the boundary scalar $\varphi_j^0$ has bounded tangential derivatives of order at most $\mathrm M$ on $\partial\cD$, which gives \eqref{eq:dyn-data-scalar-H3}.
	
	(\ref{it:dyn-lift-basic}) Equations \eqref{eq:dynamic-common-cutoff},
	\eqref{eq:dynamic-negative-time-Taylor-polynomial} and
	\eqref{eq:dynamic-negative-time-lift} give
	\eqref{eq:dynamic-negative-time-zero-near-minus-one} and
	\eqref{eq:dynamic-negative-time-jet-matching}, since $\chi_-=0$ on $[-1,-3/4]$, $\chi_-=1$ on $[-1/4,0]$ and $\partial_t^j\check{\cP}(0)=z_j$.
	
	Next, we prove \eqref{eq:dyn-lift-graph}.  By \eqref{eq:dyn-data-graph}, each $z_i$ belongs to $\cH_{\kin,\sigma}(\cD\times\bR_v^3)$.  Let $a\in C^1([-1,0])$ and $\phi\in C_c^\infty([-1,0]\times\overline\cD\times\bR_v^3)$. Applying
	\eqref{eq:H_kin_sigma-weak-stationary} for $z_i$ with $\phi(t,\cdot)$ in place of $\phi$, multiplying by $a(t)$, integrating in $t$
	and integrating by parts in $t$, we obtain \eqref{2607151226} for $az_i$ on $(-1,0)$,
	with \[
	(\partial_t+v\cdot\nabla_x)(az_i)=a'z_i+a\,v\cdot\nabla_xz_i .
	\]  So
	Proposition~\ref{H_kin=weak_landau}\,(1) gives
	$az_i\in\cH_{\kin,\sigma}((-1,0)\times\cD\times\bR_v^3)$.  Moreover, by \eqref{eq:Hsigma-dual-weighted}, \eqref{eq:auxiliary-coercive-equivalence}, $\ell^\dy_0\ge1$ and \eqref{eq:dyn-data-jets} at $I=\varnothing$,
	\[
	\|z_i\|_{L^2_xH^{-1}_\sigma}\le C\|w^{1/2}z_i\|_{L^2_{x,v}}\le C\|w^{\ell^\dy_0}z_i\|_{L^2_xH^1_\sigma}\le C\mathsf I .
	\]
	By \eqref{innerprod_sigma}, \eqref{eq:dyn-lift-size-relations} and this bound,
	\[
	\|az_i\|_{\cH_{\kin,\sigma}((-1,0)\times\cD\times\bR_v^3)}\le C\|a\|_{C^1}\mathsf I .
	\]
	Since $\check{E}=\sum_{i\le\mathrm M}a_iz_i$ with
	$a_i:=\chi_-t^i/i!$, by \eqref{eq:dynamic-negative-time-Taylor-polynomial} and
	\eqref{eq:dynamic-negative-time-lift}, we obtain
	$\check{E}\in\cH_{\kin,\sigma}((-1,0)\times\cD\times\bR_v^3)$ and \eqref{eq:dyn-lift-graph}.
	
	It remains to prove \eqref{eq:dynamic-negative-time-source-matching}.  We differentiate \eqref{eq:dynamic-negative-time-source} $j$ times, evaluate
	at $t=0$ and use \eqref{eq:dynamic-negative-time-jet-matching} and the
	Leibniz rule; this gives, for $0\le j\le\mathrm M-1$,
	\begin{align*}
		\partial_t^j\check{F}(0)
		&={}
		z_{j+1}+v\cdot\nabla_xz_j+Lz_j
		-\sum_{a=0}^{j}\binom ja N[z_a,z_{j-a}]
		=0,
	\end{align*}
	where the last equality follows from \eqref{eq:dynamic-generated-q0} at $j+1\le\mathrm M$.  This proves
	\eqref{eq:dynamic-negative-time-source-matching}.

	(\ref{it:dyn-breg-lift-bounds}) By (\ref{it:dyn-lift-weights}), the derivatives of order at most $\mathrm M+2$ of each $z_j$ are bounded, so those of order at most $\mathrm M+1$ are Lipschitz in $(x,v)$ on bounded sets and continuous up to $\partial\cD$; hence $\mathfrak D_I^-\partial_t^m(\check{E})_{\gamma^-}$, $|I|+m\le\mathrm M$, are defined pointwise.  By \eqref{eq:dynamic-negative-time-lift} and \eqref{eq:dynamic-negative-time-source}, the source
	\[
	\check{F}=\chi_-'\check{\cP}
	+\chi_-(\partial_t+v\cdot\nabla_x+L)\check{\cP}-\chi_-^2N[\check{\cP},\check{\cP}]
	\]
	contains at most one $x$-derivative or two $v$-derivatives of the $z_j$ and is of degree at
	most $2\mathrm M+2$ in the derivatives of $z_0$.  Since $\chi_-$ depends only on $t$, and $\check{\cP}$ and $\partial_t\check{\cP}$ are combinations of the $z_k$ with coefficients depending only on $t$, $\check F$ is such a combination of the $z_k$, the $(v\cdot\nabla_x+L)z_k$ and the $N[z_a,z_b]$; \eqref{eq:dyn-jet-transport-bound} and \eqref{eq:dyn-jet-collision-bound} apply at $\ell=L$ to the last two, with the weight $L+1\le L+A_\diamond$ on the right side.  By the Leibniz rule,
	\[
	\begin{aligned}
		\partial_t^{\,j}\check{F}
		={}&\sum_{k\le j}\binom jk\Bigl(\chi_-^{(k+1)}\partial_t^{\,j-k}\check{\cP}
		+\chi_-^{(k)}(\partial_t+v\cdot\nabla_x+L)\partial_t^{\,j-k}\check{\cP}\Bigr)\\
		&-\sum_{k+a+b=j}\frac{j!}{k!\,a!\,b!}\,(\chi_-^2)^{(k)}
		N\bigl[\partial_t^{\,a}\check{\cP},\partial_t^{\,b}\check{\cP}\bigr],
	\end{aligned}
	\]
	and $\partial_t^{\,i}\check{\cP}=\sum_{k=i}^{\mathrm M}\frac{t^{k-i}}{(k-i)!}\,z_k$.  Hence
	$\partial_t^{\,j}\check{F}$ is a combination of the same form.  A temporal letter of $\mathfrak D_I$ is $\partial_t$ and commutes with the spatial letters, and the first bound of \eqref{eq:dyn-field-derivative-bound} applies to the spatial subword of $I$.  Since $\|w^{-2}\|_{L^2_v}+\||v|^{1/2}w^{-3}\|_{L^2_v}<\infty$ and
	$|\cD|+|\partial\cD|<\infty$, for $u$ supported in
	$[-1,0]\times\overline\cD\times\bR^3$,
	\[
	\|u\|_{L^2}\le C\|w^2u\|_{L^\infty},
	\qquad
	\|u_{\gamma^\pm}\|_{L^2(\gamma^\pm)}\le C\|w^3u\|_{L^\infty} .
	\]
	The largest weight required by the five terms in
	(\ref{it:dyn-breg-lift-bounds}) is $L+\mathrm M+2+A_\diamond$:
	the complete derivatives of $\check F$ account for this loss.
	For $L\le\Lambda_0+B^{\dy}_{\rm in}+6$,
	\eqref{eq:dyn-breg-lift-weights} applies at that weight by our choice
	of $\check A$.  Thus
	\[
	\|w^{L+2+\mathrm M+A_\diamond}\partial_x^\alpha\partial_v^\beta z_j\|_{L^\infty(\cD\times\bR^3)}\le C\mathsf I ,
	\qquad |\alpha|+|\beta|\le\mathrm M+2,\quad 0\le j\le\mathrm M .
	\]
	This bound, the form of $\partial_t^{\,j}\check{F}$ described above, the two inequalities above for $u$ supported in $[-1,0]\times\overline\cD\times\bR^3$, and $\mathsf I\le1$ give the bound in (\ref{it:dyn-breg-lift-bounds}).
	
	(\ref{it:dyn-breg-source-bounds}) By \eqref{eq:Hsigma-dual-weighted} with $\|w^{-5/2}\|_{L^2_v}<\infty$,
	\[
	\|w^L\check{F}\|_{L^2((-1,0)\times\cD;H^{-1}_\sigma)}
	\le C\|w^{L+1/2}\check{F}\|_{L^2}
	\le C|\cD|^{1/2}\|w^{L+3}\check{F}\|_{L^\infty((-1,0)\times\cD\times\bR^3)} .
	\]
	This completes the proof.
\end{proof}

\subsection{The coefficient class and the seed}\label{ss:dyn-strong}

We verify the coefficient hypotheses for the negative-time extension \eqref{eq:dyn-strong-past-extension} when $g(0)=z_0$ and $\widetilde{\mathfrak Q}^\dy_T(g)+\mathsf I$ is small, and construct a seed realizing the prescribed initial time jets.

\begin{lemma}
	\label{lem:dyn-strong-coefficient-class}
	There is a constant $C$, depending only on the structural data, such that every $g$ on
	$\cD_T\times\bR^3_v$ satisfies
	\begin{equation}
		\label{eq:dyn-strong-coefficient-forward}
		\|\mu^{1/4}g\|_{L^\infty(\cD_T\times\bR^3)}\le C\,\cA_T[g] ,
	\end{equation}
	the lift satisfies
	\begin{equation}
		\label{eq:dyn-strong-coefficient-lift-size}
		\bigl\|w\,\mu^{1/4}\check{E}\bigr\|
		_{(C^{0,1}_{t,x}\cap C^{1,1}_v)((-1,0]\times\cD\times\bR^3)}
		\le C\,\mathsf I ,
	\end{equation}
	and, if $g(0)=z_0$, the past extension $\hat g$ of $g$ by
	\eqref{eq:dyn-strong-past-extension}, with $\mu^{1/4}\hat g$ extended to $\bR^7$ by Lemma~\ref{lem:norm-preserving-holder-extension}, satisfies
	\begin{align}
		\|\mu^{1/4}\hat g\|_{L^\infty(\bR^7)}
		&\le C\bigl(\cA_T[g]+\mathsf I\bigr) ,
		\label{eq:dyn-strong-coefficient-smallness}\\
		\bigl[\mu^{1/4}\hat g\bigr]
		_{C^{\alpha_{\rm bs}}_{\kin}((-\infty,T)\times\overline\cD\times\bR^3)}
		&\le C\bigl(\cA_T[g]+\mathsf I\bigr) .
		\label{eq:dyn-strong-coefficient-holder}
	\end{align}
	Consequently there is $\eta_{\star,1}\in(0,1]$, depending only on the structural data,
	such that, if $g(0)=z_0$ and
	\begin{equation}
		\label{eq:dyn-strong-coefficient-choice}
		\widetilde{\mathfrak Q}^\dy_T(g)+\mathsf I\le\eta_{\star,1} ,
	\end{equation}
	then
	\begin{equation}
		\label{eq:dyn-strong-coefficient-output}
		\|\mu^{1/4}\hat g\|_{L^\infty(\bR^7)}\le\tfrac12\min\{\varepsilon_{\rm dyn},\varepsilon_{\mathrm{ell}},\varepsilon_0(\cD)\},
		\qquad
		\bigl\|\mu^{1/4}\hat g\bigr\|_{C^{\alpha_{\rm bs}}_{\kin}}\le1 ,
	\end{equation}
	where $\varepsilon_{\rm dyn}$, $\varepsilon_{\mathrm{ell}}$ and $\varepsilon_0(\cD)$ are the constants of Proposition~\ref{prop:dynamic-alpha-uniform-graph}, Lemma~\ref{lem:perturbed-landau-ellipticity} and Theorem~\ref{260322442}.
\end{lemma}

\begin{proof}
	\noindent\emph{The negative-time lift.} Put $K:=w\mu^{1/4}\check{E}$ on $(-1,0]\times\cD\times\bR^3$. By \eqref{eq:dynamic-negative-time-Taylor-polynomial}--\eqref{eq:dynamic-negative-time-lift} and \eqref{eq:dyn-breg-lift-weights} at $L=0$,
	\[
	\|K\|_{(C^{0,1}_{t,x}\cap C^{1,1}_v)((-1,0]\times\cD\times\bR^3)}\le C\mathsf I .
	\]
	Indeed, $\check E$ is a finite linear combination of $z_0,\ldots,z_{\mathrm M}$ with smooth time coefficients, while $w\mu^{1/4}$ and its first two velocity derivatives are bounded.  This is \eqref{eq:dyn-strong-coefficient-lift-size}; it does not involve $g$.
	
	\noindent\emph{Bounds on the two half-lines.} We first estimate $\mu^{1/4}\hat g$ separately for positive and negative times. By Lemma~\ref{lem:dynamic-weighted-Lipschitz-to-kinetic-holder}, the multiplier $\mu^{1/4}w^{-A^\dy_\flat}$ is bounded with bounded kinetic
	H\"older seminorm. On $t\ge0$, the first two terms of \eqref{eq:dyn-coordinate-A} and the product inequality bound the $L^\infty$ norm and kinetic H\"older seminorm of $\mu^{1/4}g$ by $C\cA_T[g]$.  The
	first of these bounds is \eqref{eq:dyn-strong-coefficient-forward}.
	On $(-\infty,0]$, $\hat g=\check{E}\mathbf 1_{(-1,0]}$ vanishes for $t\le-3/4$.  The lift bound and Lemma~\ref{lem:dynamic-weighted-Lipschitz-to-kinetic-holder}, applied with $m=0$, give
	\[
	\|\mu^{1/4}\hat g\|_{L^\infty((-\infty,0]\times\overline\cD\times\bR^3)}
	+\bigl[\mu^{1/4}\hat g\bigr]_{C^{\alpha_{\rm bs}}_{\kin}((-\infty,0]\times\overline\cD\times\bR^3)}
	\le C\mathsf I .
	\]
	
	\noindent\emph{Across $t=0$.} Assume $g(0)=z_0$. Write $G:=\mu^{1/4}\hat g$ and $\alpha:=\alpha_{\rm bs}$. The identity $\check E(0)=z_0=g(0)$ makes $G$ continuous across $t=0$.  Let $Z=(t,x,v)$ and
	$Z'=(t',x',v')$ with $t<0\le t'$, $x,x'\in\overline\cD$, $d:=\mathbf d(Z,Z')$, and let
	$\hat v$ be the minimizer in \eqref{kineticdistance}.  Since $|t|,t'\le t'-t$, the definition of $\mathbf d$ gives $|t|,t'\le d^2$,
	$|v-\hat v|,|v'-\hat v|\le d$ and $|x-x'-(t-t')\hat v|\le d^3$.  
	
	\emph{Geometric estimate.} As $\partial\cD$ is compact and $C^2$ by the assumption in \S\ref{sec:diffuse-setting}, there are
	$\bar{r}_0\in(0,1]$ and $C_\cD$ such that
	\begin{equation}
		\label{eq:dyn-coefficient-chord}
		x_1,x_2\in\overline\cD,\ |x_1-x_2|\le \bar{r}_0
		\ \Longrightarrow\
		\operatorname{dist}\bigl((1-\theta)x_1+\theta x_2,\overline\cD\bigr)
		\le C_\cD|x_1-x_2|^2,
		\quad 0\le\theta\le1 .
	\end{equation}
	\emph{Proof of the claim.} If $\operatorname{dist}(x_1,\partial\cD)>\bar{r}_0$, the segment $[x_1,x_2]$ is contained in $\cD$. Otherwise,
	it is contained in one boundary chart, where $\cD=\{y_3>\phi(y')\}$ with
	$\|\nabla^2\phi\|_{L^\infty}\le C$. Put $\psi(\theta):=y_3-\phi(y')$ at $y=(1-\theta)x_1+\theta x_2$. Then $\psi(0),\psi(1)\ge0$ and $|\psi''|\le C|x_1-x_2|^2$ on $[0,1]$, hence $\psi\ge-C|x_1-x_2|^2/8$ on $[0,1]$.
	
	(a) For $d$ bounded away from zero, the $L^\infty$ bound for $G$ suffices. We may therefore assume $d$ is small enough,
	depending only on $\cD$, that $2d^{3/2}\le\bar r_0$ and
	$d\le\tfrac12d^{-1/2}$.
	
	(b) If $|\hat v|\le d^{-1/2}$, then
	\[
	|x-x'|\le d^2|\hat v|+d^3\le2d^{3/2}\le \bar{r}_0 .
	\]
	Put $\theta=-t/(t'-t)$ and $P=(1-\theta)x+\theta x'$.
	The kinetic distance and \eqref{eq:dyn-coefficient-chord} give
	\[
	|P-(x-t\hat v)|\le d^3,
	\qquad \operatorname{dist}(P,\overline\cD)\lesssim_\cD d^3.
	\]
	Choose $y\in\overline\cD$ with $|y-P|\lesssim_\cD d^3$ and set
	$\widetilde Q=(0,y,\hat v)$.  Then
	\[
	|y-(x'-t'\hat v)|\lesssim_\cD d^3,
	\qquad \mathbf d(\widetilde Q,Z')\lesssim_\cD d,
	\]
	so the positive-time H\"older estimate yields
	$|G(\widetilde Q)-G(Z')|\lesssim_\cD\cA_T[g]d^\alpha$.
	On the negative-time side $G=w^{-1}K$, where $K$ satisfies the lift estimate above.
	Since
	\[
	|x-y|\le |t|\,|\hat v|+C_\cD d^3,
	\qquad |t|\le d^2,
	\qquad |v-\hat v|\le d,
	\]
	and $w(\hat v)\ge\langle\hat v\rangle$, we obtain
	\[
	|G(Z)-G(\widetilde Q)|\lesssim_\cD
	\mathsf I\bigl(|t|+w(\hat v)^{-1}|x-y|+|v-\hat v|\bigr)
	\lesssim_\cD\mathsf I d
	\lesssim_\cD\mathsf I d^\alpha.
	\]
	
	(c) If $|\hat v|>d^{-1/2}$, then $|v|,|v'|\ge\tfrac12d^{-1/2}$, since $|v-\hat v|,|v'-\hat v|\le d\le\tfrac12d^{-1/2}$.  The time-Lipschitz
	bound at fixed $(x,v)$ gives \[
	|G(Z)-G(0,x,v)|\le C\mathsf I\,|t|\le C\mathsf I\,d^2 ,
	\]
	and $(0,x,v)$, $Z'$ belong to the half $[0,T)\times\overline\cD\times\bR^3$, where
	\[
	|G|\le\mu^{1/4}(v)\,\cA_T[g]\le C\cA_T[g]\,w(v)^{-2}
	\] by the first summand of
	\eqref{eq:dyn-coordinate-A}.  Hence
	\[
	\begin{aligned}
		|G(Z)-G(Z')|\,&\le|G(Z)-G(0,x,v)|+|G(0,x,v)|+|G(Z')|\\
		&\le C\mathsf I\,d^2+C\cA_T[g]\,d\le C(\cA_T[g]+\mathsf I)d^\alpha ,
	\end{aligned}
	\]
	where the second inequality follows from $w(v)^{-2},w(v')^{-2}\le4d$.
	
	\noindent\emph{Extension.} The preceding cases, together with the estimates on the two time half-lines, yield \[
	[G]_{C^\alpha_\kin((-\infty,T)\times\overline\cD\times\bR^3)}\le C(\cA_T[g]+\mathsf I) .
	\]  We extend
	$\mu^{1/4}\hat g$ off $\overline\cD$ in $x$, and past $T$ in $t$, by
	Lemma~\ref{lem:norm-preserving-holder-extension} on $(\bR^7,\mathbf d)$; only its values on $\cD$ enter the equations in which $\hat g$ is a coefficient.  By \eqref{eq:norm-preserving-holder-extension},
	\[
	\begin{aligned}
		\|\mu^{1/4}\hat g\|_{L^\infty(\bR^7)}+\bigl[\mu^{1/4}\hat g\bigr]_{C^{\alpha_{\rm bs}}_{\kin}((-\infty,T)\times\overline\cD\times\bR^3)}\le C\bigl(\cA_T[g]+\mathsf I\bigr) .
	\end{aligned}
	\]
	This proves
	\eqref{eq:dyn-strong-coefficient-smallness} and \eqref{eq:dyn-strong-coefficient-holder}.
	Since $\cA_T[g]\le\widetilde{\mathfrak Q}^\dy_T(g)$, choosing
	$\eta_{\star,1}$ sufficiently small gives
	\eqref{eq:dyn-strong-coefficient-output}.
\end{proof}

Both \eqref{eq:dyn-strong-coefficient-output} and the weaker condition
\begin{equation}
	\label{eq:dyn-strong-coefficient-output-reg}
	\|\mu^{1/4}\hat g\|_{L^\infty(\bR^7)}\le\min\{\varepsilon_{\rm dyn},\varepsilon_{\mathrm{ell}},\varepsilon_0(\cD)\},
	\qquad
	\bigl\|\mu^{1/4}\hat g\bigr\|_{C^{\alpha_{\rm bs}}_{\kin}}\le2,
\end{equation}
give the coefficient hypothesis of Theorems~\ref{260504554}, \ref{260526316} and
\ref{thm:weighted-steady}\,(1) and of Proposition~\ref{prop:dynamic-alpha-uniform-graph}, and, with $K=2$, that of Lemma~\ref{260527559_dy}; the constants of those results do not depend on $\hat g$. The relaxed form is used for the regularized backgrounds $g^\varepsilon$ introduced below, with $\hat g^\varepsilon$ in place of $\hat g$.

Set $\zeta_0:=f_{\dy,0}$ and $\zeta_j:=z_j$ for $1\le j\le\mathrm M$.

A background $g$ on $\cD_T\times\bR_v^3$, with $g_\dy:=g-f_\st$, matches the initial jets if
\begin{equation}
	\label{eq:dyn-ball-jets}
	\partial_t^{\,j}g_\dy(0)=\zeta_j,
	\qquad 0\le j\le\mathrm M .
\end{equation}
If $\widetilde{\mathfrak Q}^\dy_T(g)\le r$, then, by the definitions \eqref{eq:dyn-reduced-Q} and \eqref{eq:dyn-hier-rows-aggregate},
\begin{align}
	\cA_T[g]&\le r,
	\label{eq:dyn-ball-A}\\
	\cB_T[g_\dy]&\le r^2,
	\label{eq:dyn-ball-B}\\
	\cS_T[g_\dy]+\cN_{\mathrm M}(T)[g_\dy]
	&\le r .
	\label{eq:dyn-ball-rows}
\end{align}

For the sum of the two rows in \eqref{eq:dyn-ball-rows}, write
\begin{equation*}
	\Sigma_\tau[z]:=\cS_\tau[z]+\cN_{\mathrm M}(\tau)[z],
	\qquad 0<\tau\le T\le\infty ,
\end{equation*}
and write $\Sigma_J$ for the same sum with all norms taken over an interval $J$ in place
of $(0,\tau)$.

A time-localized Taylor polynomial provides a background with the required initial jets and reduced-norm bound.

\begin{lemma}
	\label{lem:dyn-fp-seed}
	Let $\varkappa$ be the cutoff \eqref{eq:dynamic-common-cutoff}, and let
	\begin{equation}
		\label{eq:dyn-fp-seed}
		g^{(0)}:=f_\st+\varkappa\,\check{\pi},
		\qquad
		\check{\pi}(t):=\sum_{j\le\mathrm M}\frac{t^j}{j!}\zeta_j .
	\end{equation}
	The following estimates hold for every $1\le T\le\infty$, with $C\ge1$ depending only on the structural data.
	\begin{enumerate}
		\item $\cA_T[f_\st]\le C\,\delta$, $g^{(0)}$ matches the initial jets
		\eqref{eq:dyn-ball-jets}, and
		\begin{equation}
			\label{eq:dyn-fp-seed-bound}
			\cA_T[g^{(0)}]\le C\bigl(\delta+\mathsf I\bigr),
			\qquad
			\cB_T\bigl[g^{(0)}-f_\st\bigr]
			\le C\bigl(\mathsf I^2+\delta^2\bigr),
			\qquad
			\Sigma_T\bigl[g^{(0)}-f_\st\bigr]\le C\bigl(\delta+\mathsf I\bigr);
		\end{equation}
		in particular
		$\widetilde{\mathfrak Q}^\dy_T(g^{(0)})\le C(\delta+\mathsf I)$.
		\item For every $-1/4\le a<0$,
		\begin{equation}
			\label{eq:dyn-fp-seed-past}
			\begin{gathered}
				\|w^{A^\dy_\flat}\check{\pi}\|_{L^\infty((a,0)\times\cD\times\bR^3)}+\|w^{A^\dy_\flat}\nabla_v\check{\pi}\|_{L_t^\infty(a,0)L_x^{p_\sharp}L_v^2}\le C\bigl(\delta+\mathsf I\bigr),\\
				\|\check{\pi}\|^2_{\cH_{\kin,\sigma}((a,0)\times\cD\times\bR_v^3)}
				\le C\bigl(\mathsf I^2+\delta^2\bigr),
				\qquad
				\Sigma_{(a,0)}[\check{\pi}]\le C\bigl(\delta+\mathsf I\bigr).
			\end{gathered}
		\end{equation}
		\item The past extension $\hat g^{(0)}$ of $g^{(0)}$ by
		\eqref{eq:dyn-strong-past-extension} is smooth in $t$ on $(-\infty,T)$, vanishes for
		$t\le-3/4$, and satisfies
		\begin{equation}
			\label{eq:dyn-fp-seed-time-derivative}
			\begin{aligned}
				&\|w^{A^\dy_\flat}\partial_t\hat g^{(0)}\|_{L^\infty((-1,T)\times\cD\times\bR^3)}+\|w^{A^\dy_\flat}\nabla_v\partial_t\hat g^{(0)}\|_{L_t^\infty(-1,T)L_x^{p_\sharp}L_v^2}\\
				&\quad+\Sigma_{(-1,T)}\bigl[\partial_t\hat g^{(0)}\bigr]
				+\bigl\|\mu^{1/4}\partial_t\hat g^{(0)}\bigr\|
				_{C^{\alpha_{\rm bs}}_{\kin}((-\infty,T)\times\overline\cD\times\bR^3)}
				\le C\bigl(\delta+\mathsf I\bigr) .
			\end{aligned}
		\end{equation}
		On $I_T^-$, each $\partial_t^j\hat g^{(0)}$, $j\ge0$, is a finite combination of
		$z_0,\ldots,z_{\mathrm M}$ and $f_\st$ with coefficients smooth and bounded in $t$; for
		$T<\infty$ it belongs to $\cH_{\kin,\sigma}(\cD_{T,\ast}\times\bR_v^3)$.
	\end{enumerate}
\end{lemma}

\begin{proof}
	\emph{(1).} We first bound $\cA_T[f_\st]$. Proposition~\ref{prop:stationary-finite-weight-regularity} at $\ell=A^\dy_\flat$ and $p=p_\sharp$
	gives
	\[
	\|w^{A^\dy_\flat}f_\st\|_{L^\infty}+[w^{A^\dy_\flat}f_\st]_{C^{\alpha_{\rm bs}}_{\kin}}+\|w^{A^\dy_\flat}\nabla_vf_\st\|_{L_t^\infty L_x^{p_\sharp}L_v^2}\le C\delta ;
	\]
	for the H\"older summand this uses that $f_\st$ is time independent: by \eqref{2607271135}, for all $t,t'\in\bR$,
	\begin{equation}
		\label{261016901}
		\mathbf d_\st((x,v),(y,v'))\le\mathbf d((t,x,v),(t',y,v')),
	\end{equation}
	so that
	\[
	[w^{A^\dy_\flat}f_\st]_{C^{\alpha_{\rm bs}}_\kin(\cD_T\times\bR_v^3)}
	\le[w^{A^\dy_\flat}f_\st]_{C^{\alpha_{\rm bs}}_\kin(\cD\times\bR_v^3)} .
	\]
	So
	\begin{equation}
		\label{261016902}
		\cA_T[f_\st]\le C\delta .
	\end{equation}
	
	Since $g^{(0)}-f_\st=\varkappa\check\pi$ and $\varkappa=1$ near $t=0$ by \eqref{eq:dynamic-common-cutoff}, differentiation at $t=0$ gives the prescribed jets \eqref{eq:dyn-ball-jets}.  To prove the first and third estimates in \eqref{eq:dyn-fp-seed-bound}, write $\varkappa\check\pi=\sum_{j\le\mathrm M}\varkappa(t)t^j\zeta_j/j!$. The time coefficients are smooth and vanish for $t\ge\mathrm M+1$, so, uniformly in $T$,
	\begin{align}\label{26115915}
		\begin{aligned}
			&\cA_T[\varkappa\check{\pi}]+\Sigma_T[\varkappa\check{\pi}]\\
			\le \,&C\sum_{j\le\mathrm M}\Bigl(\|w^{A^\dy_\flat}\zeta_j\|_{C^{\alpha_{\rm bs}}_{\kin}(\cD\times\bR_v^3)}
			+\|w^{A^\dy_\flat}\nabla_v\zeta_j\|_{L_x^{p_\sharp}L_v^2}
			+\|\zeta_j\|_{L_x^qH_\sigma^1}\\
			&\qquad\qquad +\sum_{|I|\le\mathrm M}\bigl(\|w^{\ell^\dy_{|I|}}\mathfrak D_I\zeta_j\|_{\cH_{\kin,\sigma}(\cD\times\bR_v^3)}
			+\|w^{\ell^\dy_{|I|}+1/2}\mathfrak D_I\zeta_j\|_{L^2_{x,v}}\bigr)\Bigr),
		\end{aligned}
	\end{align}
	where the inner sum is taken over the ordered words $I$ over $\mathfrak A$.
	
	To verify this inequality, let $J$ be a word over $\mathfrak A_\bullet$ of complete order $m\le\mathrm M$ containing $k$ temporal letters, and let $I$ be its spatial subword, so that $|I|=m-k$.  By \eqref{eq:dyn-normal-ordering}, the temporal letters fall on the coefficients, and \eqref{innerprod_sigma} and \eqref{260902137} give
	\[
	\begin{aligned}
		\Bigl\|w^{\ell^\dy_m}\mathfrak D_J\Bigl(\varkappa(t)\frac{t^j}{j!}\zeta_j\Bigr)\Bigr\|_{\cH_{\kin,\sigma}(\cD_T\times\bR_v^3)}
		&=\Bigl\|\partial_t^k\Bigl(\varkappa(t)\frac{t^j}{j!}\Bigr)w^{\ell^\dy_m}\mathfrak D_I\zeta_j\Bigr\|_{\cH_{\kin,\sigma}(\cD_T\times\bR_v^3)}\\
		&\le C\bigl(\|w^{\ell^\dy_{|I|}}\mathfrak D_I\zeta_j\|_{\cH_{\kin,\sigma}(\cD\times\bR_v^3)}+\|w^{\ell^\dy_{|I|}+1/2}\mathfrak D_I\zeta_j\|_{L^2_{x,v}}\bigr).
	\end{aligned}
	\]
	Here we use $\|\partial_t^i(\varkappa(t)t^j/j!)\|_{L^2(0,T)}\le C$ for $i\le\mathrm M+1$, since $\varkappa$ vanishes for $t\ge\mathrm M+1$, and $|\partial_t^k(\varkappa(t)t^j/j!)|_{t=0}|\le1$, since $\varkappa=1$ near $t=0$.  We also use Lemma~\ref{lem:compact-velocity-support}\,(1) with the $C^1_v$ multiplier $w^{\ell^\dy_m-\ell^\dy_{|I|}}$, as before \eqref{eq:dyn-unweighted-row-dominated}, since $\ell^\dy_m\le\ell^\dy_{|I|}$, and \eqref{eq:Hsigma-dual-weighted} with $\la v\ra\le Cw$.  For the H\"older summand,
	\[
	\Bigl[w^{A^\dy_\flat}\varkappa(t)\frac{t^j}{j!}\zeta_j\Bigr]_{C^{\alpha_{\rm bs}}_\kin(\cD_T\times\bR_v^3)}
	\le C\|w^{A^\dy_\flat}\zeta_j\|_{C^{\alpha_{\rm bs}}_\kin(\cD\times\bR_v^3)},
	\]
	by \eqref{261016901} and since the coefficients are bounded and Lipschitz in $t$, with $|t-t'|\le\mathbf d((t,x,v),(t',y,v'))^2$.
	Since $\zeta_j=z_j$ for $j\ge1$ and $\zeta_0=z_0-f_\st$, it suffices to bound the right-hand side of \eqref{26115915} with $\zeta_j$ replaced by $z_j$ and by $f_\st$.
	
	For the $z_j$ and every word $I$ over $\mathfrak A$ with $|I|\le\mathrm M$, \eqref{eq:dyn-data-jets} gives
	\[
	\begin{aligned}
		\|w^{\ell^\dy_{|I|}}\mathfrak D_Iz_j\|_{\cH_{\kin,\sigma}(\cD\times\bR_v^3)}&\le C\mathsf I,\qquad 
		\|w^{\ell^\dy_{|I|}+1/2}\mathfrak D_Iz_j\|_{L^2_{x,v}}\le C\|w^{A^\dy_\flat}\mathfrak D_Iz_j\|_{L^\infty}\le C\mathsf I ;
	\end{aligned}
	\]
	here we use $\ell^\dy_{|I|}\le\ell^\dy_0$ in the first inequality and, in the second, $A^\dy_\flat=\ell^\dy_0+3\ge\ell^\dy_{|I|}+3$ by \eqref{eq:dyn-data-weights} and the boundedness of $\cD$.
	By Lemma~\ref{lem:dynamic-weighted-Lipschitz-to-kinetic-holder} at $m=A^\dy_\flat$ for the H\"older norm, and by Lemma~\ref{lem:dynamic-zero-initial-source-inflow-lift}\,(\ref{it:dyn-lift-weights}) at $L=A^\dy_\flat+2$,
	\[
	\begin{aligned}
		&\|w^{A^\dy_\flat}z_j\|_{C^{\alpha_{\rm bs}}_{\kin}(\cD\times\bR_v^3)}
		+\|w^{A^\dy_\flat}\nabla_vz_j\|_{L_x^{p_\sharp}L_v^2}
		+\|z_j\|_{L_x^qH_\sigma^1}\\
		&\qquad\le C\sum_{|\alpha|+|\beta|\le2}\|w^{A^\dy_\flat+2}\partial_x^\alpha\partial_v^\beta z_j\|_{L^\infty(\cD\times\bR^3)}
		\le C\mathsf I .
	\end{aligned}
	\]
	
	For $f_\st$, let $I$ be a spatial word of length $m$. Since $\ell^\dy_m\le\ell^\dy_0=\ell_{\mathrm N}\le\ell_m$, Lemma~\ref{lem:compact-velocity-support}\,(1) with the $C^1_v$ multiplier $w^{\ell^\dy_m-\ell_m}$ bounds the graph norm of $w^{\ell^\dy_m}\mathfrak D_If_\st$ by that of $w^{\ell_m}\mathfrak D_If_\st$; since $\ell^\dy_m+1\le\ell_{\mathrm N}+A_\diamond=\ell_{\mathrm N-1}\le\ell_m$, \eqref{eq:auxiliary-coercive-equivalence} and \eqref{eq:real-weight-Hsigma-monotonicity} bound the $L^2_{x,v}$ norm of $w^{\ell^\dy_m+1/2}\mathfrak D_If_\st$:
	\[
	\begin{aligned}
		\|w^{\ell^\dy_m}\mathfrak D_If_\st\|_{\cH_{\kin,\sigma}}&\le C\,W_m(f_\st),\\
		\|w^{\ell^\dy_m+1/2}\mathfrak D_If_\st\|_{L^2_{x,v}}&\le C\|w^{\ell^\dy_m+1}\mathfrak D_If_\st\|_{L^2_xH^1_\sigma}\le C\,W_m(f_\st),
	\end{aligned}
	\]
	with $W_m$ of \eqref{eq:order-two-Wm}. Here $W_m(f_\st)\le\widetilde{\mathfrak Q}_{\mathrm N}(f_\st)\le C\delta$ by \eqref{eq:reduced-Q3} and \eqref{eq:dyn-theta-background}. In addition, $\|f_\st\|_{L^q_xH^1_\sigma}\le C\delta$ by
	\eqref{eq:stationary-finite-weight-graph} at $\ell=0$ and $p=q$. Combining these bounds with \eqref{261016902}, we obtain the first and third bounds of
	\eqref{eq:dyn-fp-seed-bound}, with constants independent of $T$, since the interval $[0,\mathrm M+1]$ does not depend on $T$.  For the second bound,
	$\varkappa\check{\pi}=0$ for $t\ge\mathrm M+1$, so the supremum in
	\eqref{eq:dyn-coordinate-B} may be restricted to $0\le t\le\mathrm M+1$, where
	$\ee^{\lambda t}\le\ee^{\lambda_{\rm exp}(\mathrm M+1)}$, since $\lambda\le\lambda_{\rm exp}$ by \eqref{eq:dyn-decay-rate-fixed-point}. Hence
	\[
	\cB_T[\varkappa\check{\pi}]\le C\ee^{\lambda_{\rm exp}(\mathrm M+1)}\|\varkappa\check{\pi}\|^2_{\cH_{\kin,\sigma}((0,\mathrm M+1)\times\cD\times\bR_v^3)}\le C(\mathsf I^2+\delta^2),
	\]
	where the first inequality uses \eqref{260730346} for the value at $t$, and the last inequality follows from the argument for the first and third bounds of \eqref{eq:dyn-fp-seed-bound}, carried out on $(0,\mathrm M+1)$.
	
	\emph{(2).} Applying the same estimates on the fixed interval $[-1/4,0]$ proves \eqref{eq:dyn-fp-seed-past}.
	
	\emph{(3).} On $I_T^-$, by \eqref{eq:dynamic-negative-time-lift} with $z_0=f_\st+\zeta_0$ and $z_j=\zeta_j$ ($j\ge1$), the extension \eqref{eq:dyn-strong-past-extension} reads
	\[
	\hat g^{(0)}=\varkappa\,\check{\pi}+\psi\,f_\st ,
	\qquad
	\psi:=\varkappa\ \text{on }(-1,0],
	\quad
	\psi:=1\ \text{on }[0,T),
	\]
	and, since $\varkappa=1$ on $[-1/4,0]$ by \eqref{eq:dynamic-common-cutoff}, $\psi$ is smooth. Both summands vanish on $(-1,-3/4]$.  Hence
	$\partial_t\hat g^{(0)}=\partial_t(\varkappa\check{\pi})+\psi'f_\st$ is supported in
	$[-1,\mathrm M+1]$ and is a finite combination of $\zeta_0,\ldots,\zeta_{\mathrm M}$ and
	$f_\st$ with coefficients smooth in $t$.  The argument for \eqref{eq:dyn-fp-seed-bound}, carried out on $[-1,\mathrm M+1]$, gives
	\[
	\|w^{A^\dy_\flat}\partial_t\hat g^{(0)}\|_{L^\infty((-1,T)\times\cD\times\bR^3)}+\|w^{A^\dy_\flat}\nabla_v\partial_t\hat g^{(0)}\|_{L_t^\infty(-1,T)L_x^{p_\sharp}L_v^2}+\Sigma_{(-1,T)}\bigl[\partial_t\hat g^{(0)}\bigr]\le C(\delta+\mathsf I) .
	\]
	For the H\"older term, regroup with
	$\check{\cP}=f_\st+\check{\pi}$ of \eqref{eq:dynamic-negative-time-Taylor-polynomial}: on
	$I_T^-$,
	\[
	\partial_t\hat g^{(0)}=\partial_t\bigl(\varkappa\,\check{\cP}\bigr)
	-\varkappa'\mathbf 1_{(0,\infty)}\,f_\st ,
	\]
	where the coefficient
	$\varkappa'\mathbf 1_{(0,\infty)}$ is smooth, since $\varkappa'=0$ on
	$[-1/4,\mathrm M-1]$.  The first part is a combination of $z_0,\ldots,z_{\mathrm M}$
	with coefficients smooth in $t$, so, as in \eqref{eq:dyn-strong-coefficient-lift-size},
	\[
	\|w\mu^{1/4}\partial_t(\varkappa\check{\cP})\|_{(C^{0,1}_{t,x}\cap C^{1,1}_v)(I_T^-\times\cD\times\bR^3)}\le C\mathsf I,
	\]
	and, by the inequality above and
	Lemma~\ref{lem:dynamic-weighted-Lipschitz-to-kinetic-holder} at $m=0$ on $(-\infty,T)$, since $\mathbf 1_{I_T^-}\mu^{1/4}\partial_t(\varkappa\check{\cP})$ vanishes for $t\le-3/4$,
	\[
	\|\mathbf 1_{I_T^-}\mu^{1/4}\partial_t(\varkappa\check{\cP})\|_{L^\infty((-\infty,T)\times\overline\cD\times\bR^3)}
	+[\mathbf 1_{I_T^-}\mu^{1/4}\partial_t(\varkappa\check{\cP})]_{C^{\alpha_{\rm bs}}_\kin((-\infty,T)\times\overline\cD\times\bR^3)}
	\le C\mathsf I .
	\]
	For the second part,
	\[
	\|\mu^{1/4}f_\st\|_{L^\infty}+[\mu^{1/4}f_\st]_{C^{\alpha_{\rm bs}}_\kin}\le C\delta
	\] by
	\eqref{eq:stationary-finite-weight-pointwise} at $\ell=0$, \eqref{2607271135} and the product
	inequality. Since the coefficient
	$\varkappa'\mathbf 1_{(0,\infty)}$ is bounded and Lipschitz in $t$, this gives
	\[
	\|\mu^{1/4}\varkappa'\mathbf 1_{(0,\infty)}f_\st\|_{L^\infty((-\infty,T)\times\overline\cD\times\bR^3)}+\bigl[\mu^{1/4}\varkappa'\mathbf 1_{(0,\infty)}f_\st\bigr]_{C^{\alpha_{\rm bs}}_\kin((-\infty,T)\times\overline\cD\times\bR^3)}\le C\delta .
	\]
	For the last assertion of (3), note that, by $\hat g^{(0)}=\varkappa\,\check{\pi}+\psi f_\st$ and $\zeta_0=z_0-f_\st$, every
	$\partial_t^j\hat g^{(0)}$, $j\ge0$, is on $I_T^-$ a finite combination of $z_0,\ldots,z_{\mathrm M}$ and
	$f_\st$ with coefficients smooth in $t$ and bounded on $I_T^-$.  For $T<\infty$, each product of such a
	coefficient with $z_i$ or $f_\st$ lies in $\cH_{\kin,\sigma}(\cD_{T,\ast}\times\bR_v^3)$, as
	$\check{E}$ does in the proof of
	Lemma~\ref{lem:dynamic-zero-initial-source-inflow-lift}\,(\ref{it:dyn-lift-basic}), by
	\eqref{eq:dyn-data-jets} and \eqref{eq:stationary-finite-weight-graph}.
	
	The definition \eqref{eq:dyn-reduced-Q} and \eqref{eq:dyn-fp-seed-bound} give the asserted reduced-norm estimate, after increasing the structural constant $C$.
\end{proof}

\subsection{The regularized problem}\label{ss:dyn-realization}

For a fixed background $g$, the regularized problem is posed on $I_T^-$ with zero initial data at $t=-1$. The incoming diffuse operator is truncated at spectral cutoff $N$, the negative-time lift supplies the extended data, and a past-directed mollifier acts only on $g-g^{(0)}$. Both regularizations are removed in the next section.

For $-1<t<T$, let $\mathfrak k_N$ denote the residual of $f_\st$ under the truncated wall operator $P_{w,N}(t)$, and let $\mathfrak d_N$ denote its spectral truncation term:
\begin{equation}
	\label{eq:dyn-limits-steady-residual}
	\mathfrak k_N:=(f_\st)_{\gamma^-}-P_{w,N}(t)(f_\st)_{\gamma^+},
	\qquad
	\mathfrak d_N:=\bigl(\mathfrak m_\circ+d_\delta(t)\bigr)(\mathbf I-J_N)\cA_0(f_\st)_{\gamma^+},
\end{equation}
and $\mathfrak k_N=(h_w-b_w)+\mathfrak d_N$ for $t>0$, by
\eqref{eq:dyn-theta-background}, \eqref{eq:dyn-theta-residual} and
\eqref{eq:dynamic-spectral-Pw-difference-exact}.

For fixed $N$, define the spectral corner corrector and the corrected inflow
residual by
\begin{align}
	c_N^\pt(t,x,v)
	:={}&\chi_+(t)\bigl(\mathfrak m_\circ+d_\delta(t)\bigr)
	(\mathbf I-J_N)
	\Bigl(\sum_{j=0}^{\mathrm M}\frac{t^j}{j!}\,\varphi_j^0\Bigr)(x),
	\label{eq:dynamic-spectral-corner-correction}\\
	b_N&:=h_w+c_N^\pt+(1-\chi_+)\,\mathfrak d_N,
	\label{eq:dynamic-spectral-corrected-inflow-residual}
\end{align}
with $d_\delta$ of \eqref{eq:wall-response-profiles}, $h_w$ of
\eqref{eq:nonisothermal-diffuse-response} and $\mathfrak d_N$ of
\eqref{eq:dyn-limits-steady-residual}.  By \eqref{eq:dynamic-common-cutoff},
$(1-\chi_+)\mathfrak d_N=0$ on $[0,\mathrm M-1]$.  The initial compatibility conditions \eqref{eq:dynamic-initial-u-compatibility} concern the full operator $P_w$. Near $t=0$, the corrector $c_N^\pt$ accounts for the difference between $P_w$ and $P_{w,N}$; the resulting matching of time jets is proved in Proposition~\ref{prop:dyn-fixedN-data}\,(\ref{it:dyn-real-corner}).

Define the negative-time spectral inflow lift by
\begin{equation}
	\label{eq:dynamic-negative-time-spectral-inflow}
	\check{b}_N(t)
	:=(\check{E})_{\gamma^-}(t)-P_{w,N}(t)(\check{E})_{\gamma^+}(t),
	\qquad -1\le t\le0,
\end{equation}
with $\check{E}$ of \eqref{eq:dynamic-negative-time-lift}.

For $0\le j\le\mathrm M$, set
\begin{equation}
	\label{eq:dynamic-zero-initial-extended-sources}
	\bar F_j(t):=
	\begin{cases}
		\partial_t^j\check{F}(t),&-1<t<0,\\
		0,&0<t<T.
	\end{cases}
\end{equation}
The extended spectral inflow is
\begin{equation}
	\label{eq:dynamic-zero-initial-extended-spectral-inflow}
	\bar b_N(t)
	:=\begin{cases}
		\check{b}_N(t),&-1<t<0,\\
		b_N(t),&0<t<T.
	\end{cases}
\end{equation}

Subtracting $f_\st$ from the truncated wall law $f_{\gamma^-}=P_{w,N}(t)f_{\gamma^+}+\bar b_N$ replaces the datum by $\bar b_N-\mathfrak k_N$. The following lemma identifies the resulting datum on $(0,T)$ and gives bounds independent of $T$ and $N$.

\begin{lemma}
	\label{lem:dyn-limits-wall}
	Let $T\ge1$ and $N\in\bN$.  The constants below depend only on the structural data and
	$L$ (in particular they are independent of $T$ and $N$).
	\begin{enumerate}
		\item\label{it:dyn-limits-wall-form}
		On $(0,T)$, $\bar b_N-\mathfrak k_N=b_w+c_N$, where
		\begin{equation}
			\label{eq:dyn-limits-corner}
			c_N(t):=\chi_+(t)\,(\mathfrak m_\circ+d_\delta)(\mathbf I-J_N)\Bigl(\sum_{k\le\mathrm M}\frac{t^k}{k!}\widetilde\varphi_k^0\Bigr),
			\qquad t\ge0,
		\end{equation}
		with $\widetilde\varphi_k^0:=\cA_0(\zeta_k)_{\gamma^+}$, $\zeta_k$ being the centered initial jets defined before \eqref{eq:dyn-ball-jets}, and $c_N$ is supported in $[0,\mathrm M+1]$.  On $(0,T)$, the corrector $c_N^\pt$ of \eqref{eq:dynamic-spectral-corner-correction} and the datum $\bar b_N$ satisfy
		\begin{equation}
			\label{eq:dyn-limits-corner-physical}
			c_N^\pt=c_N+\chi_+\,\mathfrak d_N,
			\qquad
			\bar b_N=h_w+\mathfrak d_N+c_N .
		\end{equation}
		The term $\mathfrak d_N$ of
		\eqref{eq:dyn-limits-steady-residual} satisfies, for every finite $L$,
		\begin{equation}
			\label{eq:dyn-limits-steady-residual-size}
			\sup_{t>-1}\sum_{j+|I'|\le\mathrm M}
			\bigl\|w^{L}\mathfrak D^-_{I'}\partial_t^{\,j}\mathfrak d_N(t)\bigr\|_{L^2(\gamma^-)}
			+\bigl\|w^{L}\mathfrak d_N\bigr\|_{(C^{0,1}_{t,x}\cap C^{1,1}_v)(\gamma^-_T)}
			\le C_L\,\delta .
		\end{equation}
		\item\label{it:dyn-limits-wall-sign}
		For $h\in L^2(\gamma^+)$ and $-1<t<T$,
		\begin{equation}
			\label{eq:dyn-limits-wall-flux}
			\int_{\gamma^-}\bigl(P_{w,N}(t)h\bigr)\sqrt\mu\dd\gamma
			=\int_{\gamma^+}h\sqrt\mu\dd\gamma ,
		\end{equation}
		and, for $0<t<T$,
		\begin{equation}
			\label{eq:dyn-limits-corner-massfree}
			\int_{\gamma^-}\bar b_N\sqrt\mu\dd\gamma=0 .
		\end{equation}
		\item\label{it:dyn-limits-wall-size}
		Let
		\begin{equation}
			\label{eq:dyn-limits-corner-size}
			\varsigma_N:=\sum_{k\le\mathrm M}
			\bigl\|(\mathbf I-J_N)\widetilde\varphi_k^0\bigr\|_{H^{\mathrm M}(\partial\cD)} .
		\end{equation}
		Then $\varsigma_N\to0$ as $N\to\infty$, and, for every finite $L$,
		\begin{equation}
			\label{eq:dyn-limits-corner-bound}
			\sum_{j+|I'|\le\mathrm M}
			\bigl\|w^{L}\mathfrak D^-_{I'}\partial_t^{\,j}c_N\bigr\|
			_{L^2(\gamma^-_T)}
			+\bigl\|w^{L}c_N\bigr\|_{(C^{0,1}_{t,x}\cap C^{1,1}_v)(\gamma^-_T)}
			\le C_L\,\varsigma_N .
		\end{equation}
		The wall residual satisfies
		\begin{equation}
			\label{eq:dyn-limits-wall-datum-size}
			\sum_{j+|I'|\le\mathrm M}
			\bigl\|w^{\ell^\dy_0}\mathfrak D^-_{I'}\partial_t^{\,j}b_w\bigr\|
			_{L^2(\gamma^-_T)}
			\le C\,\delta .
		\end{equation}
	\end{enumerate}
\end{lemma}

\begin{proof}
	(\ref{it:dyn-limits-wall-form}).  The first identity in \eqref{eq:dyn-limits-corner-physical} follows from the linearity of
	\eqref{eq:dynamic-spectral-corner-correction} in the scalars, since
	$\varphi_j^0-\widetilde\varphi_j^0=\delta_{j0}\cA_0(f_\st)_{\gamma^+}$ (recall $\varphi_j^0=\cA_0(z_j)_{\gamma^+}$, $\widetilde\varphi_j^0=\cA_0(\zeta_j)_{\gamma^+}$ and $\zeta_j=z_j-\delta_{j0}f_\st$); the second follows from \eqref{eq:dynamic-spectral-corrected-inflow-residual} and the first, since $\bar b_N=b_N$ on $(0,T)$.
	The corrector $c_N$ is supported in
	$[0,\mathrm M+1]$ by \eqref{eq:dynamic-common-cutoff}.  Then, since
	$\mathfrak k_N=(h_w-b_w)+\mathfrak d_N$ on $(0,T)$, the second identity of \eqref{eq:dyn-limits-corner-physical} gives
	\[
	\bar b_N-\mathfrak k_N=b_w+c_N .
	\]
	For
	\eqref{eq:dyn-limits-steady-residual-size}, split $\mathfrak D^-_{I'}\partial_t^j\mathfrak d_N$
	by the Leibniz rule \eqref{eq:boundary-word-Leibniz} with $a=\mathfrak m_\circ+d_\delta(t)$ and $b=(\mathbf I-J_N)\cA_0(f_\st)_{\gamma^+}$.  Since $b$ does not depend on $t$, the factors carrying $a$ are
	$\partial_t^{\,j_1}\mathfrak D^-_{I_1'}(\mathfrak m_\circ+d_\delta)$, $j_1+|I_1'|\le\mathrm M$.  For $j_1=0$ they satisfy
	\eqref{eq:boundary-profile-flux-bound} with $a=\mathfrak m_\circ+d_\delta(t)$ and $K=I_1'$ at each $t$, with a constant independent of $t>-1$, since by
	Lemma~\ref{lem:wall-profile-derivatives}, applied for $-1<t<0$ as in Lemma~\ref{lem:dyn-wall-extended}\,(\ref{it:dyn-wall-ext-results}), the constants $C_a,c_a$ of \eqref{eq:boundary-gaussian-polynomial-envelope} for $\mathfrak m_\circ+d_\delta(t)$ do not depend on $t$.  For $j_1\ge1$ they equal
	$c_\mu^{-1}\partial_t^{\,j_1}\mathfrak D^-_{I_1'}h_w$ (since $\mathfrak m_\circ$ does not depend on $t$ and $h_w=c_\mu d_\delta$), and, since $\mathfrak D^-_{I_1'}h_w$ is the restriction of $\mathfrak D_{I_1'}h_w$ to $\gamma^-$ by \eqref{eq:boundary-complete-derivative}, the first bound of
	\eqref{eq:dyn-theta-defect-complete-bound} with $(j,I)=(j_1,I_1')$ gives, for $t>0$,
	\[
	\bigl|w^L\partial_t^{\,j_1}\mathfrak D^-_{I_1'}(\mathfrak m_\circ+d_\delta)(t,x,v)\bigr|\le C_L\,\delta\ee^{-c|v|^2},
	\qquad (x,v)\in\gamma^-,
	\]
	which holds for $-1<t<0$ as well by
	Lemma~\ref{lem:dyn-wall-extended}\,(\ref{it:dyn-wall-ext-results}).  The other factors are $\tau_{I_2'}b$, $|I_2'|\le\mathrm M$.  By $0\le1-\chi\le1$ for the spectral cutoff $\chi$ of $J_N$ and Lemma~\ref{lem:word-boundary-spectral-chart}, and by \eqref{eq:boundary-moment-H3}, respectively,
	\[
	\|(\mathbf I-J_N)\psi\|_{H^{\mathrm M}(\partial\cD)}\le C\|\psi\|_{H^{\mathrm M}(\partial\cD)},
	\qquad
	\|\cA_0(f_\st)_{\gamma^+}\|_{H^{\mathrm M}(\partial\cD)}\le C\delta .
	\]  Multiplying the bounds of the two factors, we obtain \eqref{eq:dyn-limits-steady-residual-size}. The constant is independent of $N$, since the bound of $\mathbf I-J_N$ on $H^{\mathrm M}(\partial\cD)$ is, and of $T$, since the bounds of the factors carrying $a$ are uniform in $t>-1$.  For the Lipschitz bound in \eqref{eq:dyn-limits-steady-residual-size}, we have $\|w^L(\mathfrak m_\circ+d_\delta)\|_{(C^{0,1}_{t,x}\cap C^{1,1}_v)(\gamma^-_T)}\le C_L$, with $C_L$ independent of $T$: in $x$ and $v$ by \eqref{eq:boundary-gaussian-polynomial-envelope} and Lemma~\ref{lem:wall-profile-derivatives}, and in $t$ by the bound on $\partial_td_\delta$ of Lemma~\ref{lem:dyn-theta-profile-time-derivatives}.  The scalar $(\mathbf I-J_N)\cA_0(f_\st)_{\gamma^+}$ does not depend on $t$ and $v$, and its $C^{0,1}(\partial\cD)$ norm is at most $C\delta$ by the two bounds above and the embedding $H^{\mathrm M}(\partial\cD)\hookrightarrow C^{0,1}(\partial\cD)$, $\mathrm M\ge3$, on the two-dimensional boundary.  Multiplying the two bounds gives the Lipschitz bound in \eqref{eq:dyn-limits-steady-residual-size}.
	
	(\ref{it:dyn-limits-wall-sign}).  At each $t$, \eqref{eq:dyn-limits-wall-flux} follows by integrating over $\partial\cD$ the pointwise flux identity stated before \eqref{eq:truncated-wall-flux-integrated}, with $\Theta_w(t)$ in place of $\Theta_w$, and applying \eqref{eq:truncated-wall-flux-integrated}.  To prove \eqref{eq:dyn-limits-corner-massfree}, integrate the boundary mass scalars \eqref{eq:dynamic-boundary-mass-scalars} of the three terms in $\bar b_N=h_w+c_N^\pt+(1-\chi_+)\mathfrak d_N$ over $\partial\cD$.  First, $\beta_{h_w}\equiv0$ by
	\eqref{eq:dyn-wall-datum-flux}.  Since $\int_{n_x\cdot v<0}\mathfrak m_{\Theta_w(t)}\sqrt\mu\,|n_x\cdot v|\dd v=1$, the scalars
	$\beta_{c_N^\pt}$ and $\beta_{(1-\chi_+)\mathfrak d_N}=(1-\chi_+)\beta_{\mathfrak d_N}$ are of the form $(\mathbf I-J_N)\phi$, and
	$\int_{\partial\cD}(\mathbf I-J_N)\phi\dd S_x=0$ by \eqref{eq:truncated-wall-flux-integrated}.  Summing over the three summands gives
	\eqref{eq:dyn-limits-corner-massfree}.
	
	(\ref{it:dyn-limits-wall-size}).  To prove $\varsigma_N\to0$, note that
	$\widetilde\varphi^0_k\in H^{\mathrm M}(\partial\cD)$ by \eqref{eq:dyn-data-scalar-H3} and
	\eqref{eq:boundary-moment-H3} (note that $\widetilde\varphi^0_0=\varphi^0_0-\cA_0(f_\st)_{\gamma^+}$ and $\widetilde\varphi^0_k=\varphi^0_k$ for $k\ge1$). Since $J_N$ commutes with $(1+A_\partial)^{\mathrm M/2}$, \eqref{eq:boundary-spectral-strong-convergence} applied to $(1+A_\partial)^{\mathrm M/2}\widetilde\varphi^0_k$ gives $(\mathbf I-J_N)\widetilde\varphi^0_k\to0$ in the spectral norm \eqref{eq:word-boundary-spectral-Hs} with $s=\mathrm M$, which is equivalent to the norm of $H^{\mathrm M}(\partial\cD)$ by Lemma~\ref{lem:word-boundary-spectral-chart}.  Hence $\varsigma_N\to0$.  By \eqref{eq:dyn-limits-corner},
	$$
	c_N=\chi_+\,(\mathfrak m_\circ+d_\delta)\,\psi_N,\qquad \psi_N:=(\mathbf I-J_N)\sum_{k\le\mathrm M}\frac{t^k}{k!}\widetilde\varphi_k^0.
	$$
	On every boundary chart and for $j,|\beta|,|\alpha|\le\mathrm M$,
	\[
	\bigl|w^L\partial_t^j\partial_{y'}^\beta\partial_v^\alpha(\mathfrak m_\circ+d_\delta)\bigl(t,\phi_\kappa(y',0),v\bigr)\bigr|\le C_L\ee^{-c|v|^2},
	\]
	uniformly in $t>0$, by Lemma~\ref{lem:wall-profile-derivatives}
	and Lemma~\ref{lem:dyn-theta-profile-time-derivatives}; since $L$ is arbitrary and $|v|\le w$, the inequality above also holds with $w^L$ replaced by $|v|^iw^L$, $i\ge0$, and $C_L$ and $c$ replaced by constants depending on $L$ and $i$.  Since $\chi_+$ and the coefficients $t^k/k!$ are bounded on $[0,\mathrm M+1]$ with their derivatives and the $H^{\mathrm M}(\partial\cD)$ norms of the scalars $(\mathbf I-J_N)\widetilde\varphi^0_k$ add up to $\varsigma_N$,
	\[
	\sup_{0<t<\mathrm M+1}\sum_{j\le\mathrm M}
	\|\partial_t^j(\chi_+\psi_N)(t)\|_{H^{\mathrm M}(\partial\cD)}
	\le C\varsigma_N .
	\]  Moreover,
	\[
	\|\chi_+\psi_N\|_{W^{1,\infty}((0,\mathrm M+1)\times\partial\cD)}\le C\varsigma_N
	\] by the
	embedding $H^{\mathrm M}(\partial\cD)\hookrightarrow C^{0,1}(\partial\cD)$, $\mathrm M\ge3$, on the
	two-dimensional boundary.  With \eqref{eq:boundary-word-Leibniz}, these give \eqref{eq:dyn-limits-corner-bound}.  Since $c_N$ is supported in $[0,\mathrm M+1]$, these bounds do not depend on $T$. 
	
	For \eqref{eq:dyn-limits-wall-datum-size}, consider the first summand
	$(d_\delta(t)-d_\delta^\infty)\cA_0(f_\st)_{\gamma^+}$ of $b_w$ in \eqref{eq:dyn-theta-residual}.  By
	\eqref{eq:boundary-word-Leibniz}, each letter of $\mathfrak D^-_{I'}\partial_t^j$ falls on the
	profile $d_\delta(t)-d_\delta^\infty$ or on the scalar; since the scalar does not depend on $t$, the profile factors are $\partial_t^{\,j}\mathfrak D^-_{I_1'}(d_\delta(t)-d_\delta^\infty)$, $I_1'$ a subword of $I'$.  Since $\mathfrak D^-_{I_1'}$ is the restriction of $\mathfrak D_{I_1'}$ to $\gamma^-$ by \eqref{eq:boundary-complete-derivative}, the second bound of \eqref{eq:dyn-theta-defect-complete-bound} with $L=\ell^\dy_0$ gives, for $t>0$,
	\[
	\bigl|w^{\ell^\dy_0}\partial_t^{\,j}\mathfrak D^-_{I_1'}\bigl(d_\delta(t)-d_\delta^\infty\bigr)(x,v)\bigr|\le C\delta\ee^{-\lambda_wt}\ee^{-c|v|^2},
	\qquad (x,v)\in\gamma^-,
	\]
	since
	$d_\delta(t)-d_\delta^\infty=c_\mu^{-1}(h_w(t)-h_w^\infty)$; that bound covers the words with $|I'|\le\mathrm M\le\mathrm N$. The scalar satisfies
	$\|\cA_0(f_\st)_{\gamma^+}\|_{H^{\mathrm M}(\partial\cD)}\le C\delta$ by
	\eqref{eq:boundary-moment-H3}. For the second summand $h_w(t)-h_w^\infty$, the second bound of
	\eqref{eq:dyn-theta-defect-complete-bound} with $L=\ell^\dy_0$ gives, for $t>0$,
	\[
	\bigl|w^{\ell^\dy_0}\mathfrak D^-_{I'}\partial_t^j\bigl(h_w(t)-h_w^\infty\bigr)(x,v)\bigr|\le C\delta\ee^{-\lambda_wt}\ee^{-c|v|^2},
	\qquad (x,v)\in\gamma^-.
	\]
	Combining these bounds, we obtain, for $t>0$,
	\[
	\|w^{\ell^\dy_0}\mathfrak D^-_{I'}\partial_t^jb_w(t)\|_{L^2(\gamma^-)}\le C\delta\ee^{-\lambda_wt} ,
	\]
	and integration over $(0,T)$ gives \eqref{eq:dyn-limits-wall-datum-size} uniformly in $T$.
\end{proof}

\begin{prop}[Glued data]
	\label{prop:dyn-fixedN-data}
	Let $T\ge1$ and $N\in\bN$.
	\begin{enumerate}
		\item\label{it:dyn-real-corner}
		(\emph{Gluing at $t=0$.})  The inflow data satisfy
		\begin{equation}
			\label{eq:dynamic-spectral-inflow-gluing}
			\partial_t^j\check{b}_N(0)=\partial_t^jb_N(0),
			\qquad 0\le j\le\mathrm M.
		\end{equation}
		\item\label{it:dyn-real-extended}
		(\emph{The extended data.})  Let $g\in\cH_{\kin,\sigma}(\cD_T\times\bR_v^3)$ with
		$g(0)=z_0$, and let $\bar g_0$ be defined from $g$ by \eqref{eq:dynamic-base-background}.  Then
		\begin{align}
			\partial_t\bar F_j&=\bar F_{j+1}
			&&\text{in }\cD'(I_T^-;L_x^2H_\sigma^{-1}),
			&&0\le j\le\mathrm M-1,
			\label{eq:dynamic-zero-initial-extended-source-chain}\\
			\bar b_N&\in H^{\mathrm M}(I_T^-;L^2(\gamma^-)),
			\label{eq:dynamic-zero-initial-extended-inflow-H3}\\
			\bar g_0&\in\cH_{\kin,\sigma}(\cD_{T,\ast}\times\bR_v^3),
			\label{eq:dynamic-zero-initial-extended-coefficient-graphs}\\
			\bar F_j&\in L^2(I_T^-;L_x^2H_\sigma^{-1}),
			&&0\le j\le\mathrm M.
			\label{eq:dynamic-zero-initial-extended-source-spaces}
		\end{align}
		For $0\le j\le\mathrm M$, $\bar g_0=0$, $\bar F_j=0$ and $\partial_t^j\bar b_N=0$ on
		$(-1,-3/4]$.
	\end{enumerate}
\end{prop}

\begin{proof}
	(1) On $[-1/4,0]$, $\chi_-=1$ by \eqref{eq:dynamic-common-cutoff}, so
	$\check{E}=\check{\cP}$ by \eqref{eq:dynamic-negative-time-lift}.  Inserting this into
	\eqref{eq:dynamic-negative-time-spectral-inflow}, with
	\eqref{eq:dynamic-spectral-Pw-difference-exact} and the scalars 
	$\varphi_j^0=\cA_0(z_j)_{\gamma^+}$, gives, for $-1/4\le t\le0$,
	\begin{align*}
		\check{b}_N(t)
		&={}
		\beta(t)
		+\bigl(\mathfrak m_\circ+d_\delta(t)\bigr)(\mathbf I-J_N)
		\sum_{j=0}^{\mathrm M}\frac{t^j}{j!}\,\varphi_j^0 ,
		\qquad
		\beta(t):=\bigl(\check{\cP}(t)\bigr)_{\gamma^-}
		-P_w(t)\bigl(\check{\cP}(t)\bigr)_{\gamma^+} .
	\end{align*}
	Since $\partial_t^i\check{\cP}(0)=z_i$,
	$P_w(t)=(\mathfrak m_\circ+d_\delta(t))\cA_0$ and $\partial_t^kd_\delta$ is continuous across $t=0$ for $k\le\mathrm N$ (so the derivatives of $d_\delta$ at $t=0$ from the left are the values $(\partial_t^kd_\delta)(0)$ in $\mathfrak t_j[z]$), the Leibniz rule and
	\eqref{eq:dynamic-initial-u-compatibility} give
	\[
	\partial_t^j\beta(0)
	=(z_j)_{\gamma^-}-P_w(0)(z_j)_{\gamma^+}-\mathfrak t_j[z]
	=\partial_t^jh_w(0),
	\qquad 0\le j\le\mathrm M .
	\]
	The second summand in the expression for $\check{b}_N(t)$ is \eqref{eq:dynamic-spectral-corner-correction} without
	the factor $\chi_+$, and $\chi_+=1$ on $[0,\mathrm M-1]$ by
	\eqref{eq:dynamic-common-cutoff}.  So its derivatives of order
	$\le\mathrm M$ at $t=0$ are those of $c_N^\pt$.  Since also
	$(1-\chi_+)\mathfrak d_N=0$ on $[0,\mathrm M-1]$,
	\eqref{eq:dynamic-spectral-corrected-inflow-residual} gives
	\[
	\partial_t^j\check{b}_N(0)
	=\partial_t^jh_w(0)+\partial_t^jc_N^\pt(0)
	=\partial_t^jb_N(0),
	\qquad 0\le j\le\mathrm M,
	\]
	which is \eqref{eq:dynamic-spectral-inflow-gluing}.
	
	(2) The data have no jump at $t=0$. Indeed, $g(0)=z_0$ and \eqref{eq:dynamic-negative-time-jet-matching} at $j=0$ give $\bar g_0(0^-)=z_0=\bar g_0(0^+)$. Likewise, \eqref{eq:dynamic-negative-time-source-matching} gives $\bar F_j(0^-)=0=\bar F_j(0^+)$ for $0\le j\le\mathrm M-1$, and $\partial_t\bar F_j=\bar F_{j+1}$ holds on $(-1,0)$ and on $(0,T)$ by \eqref{eq:dynamic-zero-initial-extended-sources}, so no jump term appears at $t=0$. This proves \eqref{eq:dynamic-zero-initial-extended-source-chain}.
	For \eqref{eq:dynamic-zero-initial-extended-inflow-H3}, on $(-1,0)$, $\check{b}_N=(\check{E})_{\gamma^-}-(\mathfrak m_\circ+d_\delta)J_N\cA_0(\check{E})_{\gamma^+}$ by \eqref{eq:dynamic-negative-time-spectral-inflow} and \eqref{eq:boundary-spectral-diffuse-operator}.  For $0\le m\le\mathrm M$, $\partial_t^m(\check{E})_{\gamma^-}\in L^2((-1,0)\times\gamma^-)$ by Lemma~\ref{lem:dynamic-zero-initial-source-inflow-lift}\,(\ref{it:dyn-breg-lift-bounds}) at $I=\varnothing$, and, by the Leibniz rule in $t$, \eqref{eq:dynamic-spectral-profile-L2-isometry-bound} and $0\le J_N\le\mathbf I$, for $-1<t<0$,
	\[
	\begin{aligned}
		&\bigl\|\partial_t^m\bigl((\mathfrak m_\circ+d_\delta)J_N\cA_0(\check{E})_{\gamma^+}\bigr)(t)\bigr\|_{L^2(\gamma^-)}\\
		&\qquad\le C\sum_{k\le m}\|\cA_0(\partial_t^k\check{E})_{\gamma^+}(t)\|_{L^2(\partial\cD)}
		\le C\sum_{k\le m}\|w^3\partial_t^k\check{E}\|_{L^\infty((-1,0)\times\cD\times\bR^3)}
		\le C\mathsf I ,
	\end{aligned}
	\]
	by \eqref{eq:outgoing-moment-flux-bound}, $|\partial\cD|<\infty$, $\||v|^{1/2}w^{-3}\|_{L^2_v}<\infty$ and Lemma~\ref{lem:dynamic-zero-initial-source-inflow-lift}\,(\ref{it:dyn-breg-lift-bounds}).  Hence $\check{b}_N\in H^{\mathrm M}((-1,0);L^2(\gamma^-))$.
	On $(0,T)$, $\bar b_N=h_w+\mathfrak d_N+c_N$ by \eqref{eq:dyn-limits-corner-physical}, and
	\[
	\begin{aligned}
		h_w&\in H^{\mathrm M}\bigl((0,T);L^2(\gamma^-)\bigr)
		&&\text{by \eqref{eq:dyn-wall-datum-class}},\\
		\mathfrak d_N&\in H^{\mathrm M}\bigl((0,T);L^2(\gamma^-)\bigr)
		&&\text{by \eqref{eq:dyn-limits-steady-residual-size} at }L=0\text{, as }T<\infty,\\
		c_N&\in H^{\mathrm M}\bigl((0,T);L^2(\gamma^-)\bigr)
		&&\text{by \eqref{eq:dyn-limits-corner-bound} at }L=0 .
	\end{aligned}
	\]
	By \eqref{eq:dynamic-spectral-inflow-gluing},
	$\partial_t^j\check{b}_N(0^-)=\partial_t^jb_N(0^+)$ for
	$0\le j\le\mathrm M$, so $\bar b_N$ has no jump in
	$\partial_t^j$, $0\le j\le\mathrm M-1$.  Hence
	\eqref{eq:dynamic-zero-initial-extended-inflow-H3} follows from
	the bound on $(-1,0)$ and the three statements
	above.
	
	Next, \eqref{eq:dynamic-zero-initial-extended-coefficient-graphs} follows from \eqref{eq:dyn-time-gluing} at $(a,b,c)=(-1,0,T)$, since $\check{E}\in\cH_{\kin,\sigma}((-1,0)\times\cD\times\bR_v^3)$ by Lemma~\ref{lem:dynamic-zero-initial-source-inflow-lift}\,(\ref{it:dyn-lift-basic}), $g\in\cH_{\kin,\sigma}(\cD_T\times\bR_v^3)$, and $\bar g_0(0^-)=\bar g_0(0^+)$.
	For \eqref{eq:dynamic-zero-initial-extended-source-spaces}, it suffices by \eqref{eq:dynamic-zero-initial-extended-sources} to bound $\partial_t^j\check{F}$ in $L^2((-1,0);L_x^2H_\sigma^{-1})$. Since $\|F\|_{H_\sigma^{-1}}\le C\|\la v\ra^{1/2}F\|_{L^2_v}$ by \eqref{eq:Hsigma-dual-weighted}, this follows from
	Lemma~\ref{lem:dynamic-zero-initial-source-inflow-lift}\,(\ref{it:dyn-breg-lift-bounds}) at
	$I=\varnothing$ and $L=\tfrac12$.
	Finally, on $(-1,-3/4]$, $\check{E}=0$ by \eqref{eq:dynamic-negative-time-zero-near-minus-one};
	hence $\check{F}=0$ by \eqref{eq:dynamic-negative-time-source} and $\check{b}_N=0$ by
	\eqref{eq:dynamic-negative-time-spectral-inflow}, and the vanishing of $\bar g_0$, $\bar F_j$
	and $\partial_t^j\bar b_N$ there follows from \eqref{eq:dynamic-base-background} and
	\eqref{eq:dynamic-zero-initial-extended-sources}--%
	\eqref{eq:dynamic-zero-initial-extended-spectral-inflow}.
\end{proof}

We regularize the background by a past-directed time mollifier, so that $(\varrho_\varepsilon*z)(t)$ uses only values of $z$ on $(t-\varepsilon,t)$:
\begin{align}
	\label{eq:dyn-real-mollifier}
	\begin{gathered}
		\varrho\in C^\infty_c((0,1)),
		\quad
		\varrho\ge0,\quad 
		\int_0^1\varrho(s)\dd s=1,\\
		\varrho_\varepsilon(s)=\varepsilon^{-1}\varrho(s/\varepsilon),
		\quad
		(\varrho_\varepsilon*z)(t)=\int_0^\varepsilon\varrho_\varepsilon(s)z(t-s)\dd s .
	\end{gathered}
\end{align}
Let $\varrho_\varepsilon^{(j)}$ be the $j$-th derivative of $\varrho_\varepsilon$, $\hat g$ the
extension \eqref{eq:dyn-strong-past-extension}, and $\hat g^{(0)}$ the past extension of the
seed \eqref{eq:dyn-fp-seed}.  The regularized background is
\begin{equation}
	\label{eq:dyn-real-mollified-objects}
	g^\varepsilon:=\hat g^{(0)}+\varrho_\varepsilon*\bigl(\hat g-\hat g^{(0)}\bigr),
	\qquad
	g^\varepsilon_j:=\partial_t^j\hat g^{(0)}
	+\varrho_\varepsilon^{(j)}*\bigl(\hat g-\hat g^{(0)}\bigr),
	\qquad 0\le j\le\mathrm M,
\end{equation}
on $I_T^-$, and we put $g^\varepsilon_\dy:=g^\varepsilon-f_\st$.  The source $\bar F_0$ in \eqref{eq:dynamic-zero-initial-extended-sources} remains unchanged. Mollifying only $\hat g-\hat g^{(0)}$ preserves the initial time jets: the next lemma gives $\partial_t^jg^\varepsilon(0)=\partial_t^j\hat g^{(0)}(0)$ for $0\le j\le\mathrm M$.

For $0<T<\infty$ put
$Q:=(-\infty,T)\times\overline\cD\times\bR^3$.  For $s\ge0$, the map $(t,x,v)\mapsto(t-s,x,v)$ sends $Q$ into $Q$ and preserves the kinetic distance \eqref{kineticdistance}, and $\mu^{1/4}$ depends only on $v$, so
\[
\|\mu^{1/4}u(\cdot-s)\|_{C^{\alpha_{\rm bs}}_{\kin}(Q)}\le\|\mu^{1/4}u\|_{C^{\alpha_{\rm bs}}_{\kin}(Q)} .
\]
Consequently, whenever $\|\mu^{1/4}u\|_{C^{\alpha_{\rm bs}}_{\kin}(Q)}<\infty$, past-directed convolution satisfies
\begin{equation}
	\label{eq:dyn-strong-mollification-holder}
	\|\mu^{1/4}\varrho_\varepsilon*u\|_{C^{\alpha_{\rm bs}}_{\kin}(Q)}
	\le\|\mu^{1/4}u\|_{C^{\alpha_{\rm bs}}_{\kin}(Q)} .
\end{equation}

\begin{lemma}[Regularized background]
	\label{lem:dyn-real-increment}
	Let $T\ge1$, $0<\varepsilon<1/4$ and $g\in\cH_{\kin,\sigma}(\cD_T\times\bR_v^3)$, and let $\bar g_0$ and $g^\varepsilon$ be defined from $g$ by
	\eqref{eq:dynamic-base-background} and \eqref{eq:dyn-real-mollified-objects}.
	Then, on $I_T^-$, with $\bar g_0$ extended by zero to $t\le-1$, and for
	$0\le j\le\mathrm M$,
	\begin{equation}
		\label{eq:dyn-real-increment-split}
		g^\varepsilon=\bar g_0\ \text{ on }(-1,0],
		\qquad
		g^\varepsilon_j=\varrho_\varepsilon^{(j)}*\bar g_0
		+\bigl(\partial_t^j\hat g^{(0)}-\varrho_\varepsilon*\partial_t^j\hat g^{(0)}\bigr),
		\qquad
		\partial_t^jg^\varepsilon(0)=\delta_{j0}f_\st+\zeta_j .
	\end{equation}
	There is $C$, depending only on the structural data, such
	that the increment $\hat g^{(0)}-\varrho_\varepsilon*\hat g^{(0)}$ satisfies
	\begin{equation}
		\label{eq:dyn-real-increment-bound}
		\begin{aligned}
			&\|w^{A^\dy_\flat}(\hat g^{(0)}-\varrho_\varepsilon*\hat g^{(0)})\|_{L^\infty((-1,T)\times\cD\times\bR^3)}+\|w^{A^\dy_\flat}\nabla_v(\hat g^{(0)}-\varrho_\varepsilon*\hat g^{(0)})\|_{L_t^\infty(-1,T)L_x^{p_\sharp}L_v^2}\\
			&+\Sigma_{(-1,T)}\bigl[\hat g^{(0)}-\varrho_\varepsilon*\hat g^{(0)}\bigr]
			+\Sigma_T\bigl[\hat g^{(0)}-\varrho_\varepsilon*\hat g^{(0)}\bigr]+\bigl\|\mu^{1/4}(\hat g^{(0)}-\varrho_\varepsilon*\hat g^{(0)})\bigr\|
			_{C^{\alpha_{\rm bs}}_{\kin}(Q)}\\
			&\quad \le C\,\varepsilon\,\bigl(\delta+\mathsf I\bigr),\\
			&\cB_T\bigl[\hat g^{(0)}-\varrho_\varepsilon*\hat g^{(0)}\bigr]
			\le\bigl(C\,\varepsilon\,(\delta+\mathsf I)\bigr)^2 ,
		\end{aligned}
	\end{equation}
	with $Q$ of \eqref{eq:dyn-strong-mollification-holder}, and $\Sigma_T$ and $\cB_T$ taken
	over $(0,T)$.
\end{lemma}

\begin{proof}
	By \eqref{eq:dyn-strong-past-extension}, $\hat g-\hat g^{(0)}=(g-g^{(0)})\mathbf 1_{(0,T)}$
	vanishes for $t\le0$. Since the convolution is past directed,
	$\varrho_\varepsilon^{(j)}*(\hat g-\hat g^{(0)})=0$ on $(-\infty,0]$.  With
	$\hat g=\bar g_0$ on $I_T^-$ and $\hat g^{(0)}=\check{E}=\hat g$ on $(-1,0]$ by \eqref{eq:dyn-strong-past-extension}, this gives the first identity of
	\eqref{eq:dyn-real-increment-split}.  The second follows from
	$\varrho_\varepsilon^{(j)}*\hat g^{(0)}=\varrho_\varepsilon*\partial_t^j\hat g^{(0)}$, obtained by
	integration by parts in $s$, since $\hat g^{(0)}$ is smooth in $t$ and vanishes for
	$t\le-3/4$ by Lemma~\ref{lem:dyn-fp-seed}.  The term $\varrho_\varepsilon*(\hat g-\hat g^{(0)})$ of
	\eqref{eq:dyn-real-mollified-objects} is smooth in $t$ and vanishes, with all its
	$t$-derivatives, at $t=0$.  Hence
	\[
	\partial_t^jg^\varepsilon(0)=\partial_t^j\hat g^{(0)}(0)=\partial_t^j(f_\st+\check{\pi})(0){}=\delta_{j0}f_\st+\zeta_j,
	\]
	since $\varkappa=1$ on $[-1/4,\mathrm M-1]$ and $\partial_t^j\check{\pi}(0)=\zeta_j$ by \eqref{eq:dyn-fp-seed}.  This is the third identity.
	
	For the first inequality of \eqref{eq:dyn-real-increment-bound}, we write the increment as
	\[
	\bigl(\hat g^{(0)}-\varrho_\varepsilon*\hat g^{(0)}\bigr)(t)
	=\int_0^\varepsilon\varrho_\varepsilon(s)\int_0^s
	\partial_t\hat g^{(0)}(t-r)\dd r\dd s ,
	\]
	a superposition of the translates $\partial_t\hat g^{(0)}(\cdot-r)$, $0<r<\varepsilon$,
	against a nonnegative measure of mass $\int_0^\varepsilon s\varrho_\varepsilon(s)\dd s\le\varepsilon$.
	Each norm on the left of the first inequality of \eqref{eq:dyn-real-increment-bound}, other than the H\"older one and $\Sigma_T$, is a sum of norms $\|\mathsf T z\|_{L^a(J;\mathcal Y)}$ over
	intervals $J\subset I_T^-$ and of $L^2_{x,v}$ norms of values at $t=-1$; here
	$\mathsf T$ acts in $(x,v)$ only, or is a temporal letter $\partial_t^{\,k}$ or
	$\partial_t+v\cdot\nabla_x$.  Each such $\mathsf T$ commutes with translation in $t$, and the
	bound
	\[
	\|\mu^{1/4}\partial_t\hat g^{(0)}(\cdot-r)\|_{C^{\alpha_{\rm bs}}_{\kin}(Q)}\le\|\mu^{1/4}\partial_t\hat g^{(0)}\|_{C^{\alpha_{\rm bs}}_{\kin}(Q)}
	\]
	holds by the argument before
	\eqref{eq:dyn-strong-mollification-holder}.  Since $\partial_t\hat g^{(0)}=0$ for
	$t\le-3/4$, the translates, restricted to $I_T^-$, involve $\partial_t\hat g^{(0)}$ only on $(-1,T)$, and their values at $t=-1$ vanish.  Hence Minkowski's integral inequality bounds each of these norms of the increment by $\varepsilon$ times
	the same norm of $\partial_t\hat g^{(0)}$ over $(-1,T)$, and the H\"older norm on $Q$ of $\mu^{1/4}$ times the increment by $\varepsilon$ times that of $\mu^{1/4}\partial_t\hat g^{(0)}$:
	\[
	\begin{aligned}
		\|\mathsf T(\hat g^{(0)}-\varrho_\varepsilon*\hat g^{(0)})\|_{L^a(J;\mathcal Y)}&\le\varepsilon\|\mathsf T\partial_t\hat g^{(0)}\|_{L^a((-1,T);\mathcal Y)},\\
		\|\mu^{1/4}(\hat g^{(0)}-\varrho_\varepsilon*\hat g^{(0)})\|_{C^{\alpha_{\rm bs}}_{\kin}(Q)}&\le\varepsilon\|\mu^{1/4}\partial_t\hat g^{(0)}\|_{C^{\alpha_{\rm bs}}_{\kin}(Q)} .
	\end{aligned}
	\]
	These norms of
	$\partial_t\hat g^{(0)}$, and the H\"older norm of $\mu^{1/4}\partial_t\hat g^{(0)}$, satisfy
	\[
	\begin{aligned}
		&\|w^{A^\dy_\flat}\partial_t\hat g^{(0)}\|_{L^\infty((-1,T)\times\cD\times\bR^3)}+\|w^{A^\dy_\flat}\nabla_v\partial_t\hat g^{(0)}\|_{L_t^\infty(-1,T)L_x^{p_\sharp}L_v^2}\\
		&\quad+\Sigma_{(-1,T)}\bigl[\partial_t\hat g^{(0)}\bigr]+\bigl\|\mu^{1/4}\partial_t\hat g^{(0)}\bigr\|_{C^{\alpha_{\rm bs}}_{\kin}(Q)}\le C(\delta+\mathsf I)
	\end{aligned}
	\]
	by \eqref{eq:dyn-fp-seed-time-derivative}.  It remains to bound $\Sigma_T$.  Each term of $\Sigma_T$ other than a graph norm is at most the corresponding term of $\Sigma_{(-1,T)}$.  For a graph norm, \eqref{innerprod_sigma} and \eqref{260730346} give, for $u\in\cH_{\kin,\sigma}(\cD_{T,\ast}\times\bR_v^3)$,
	\[
	\|u\|_{\cH_{\kin,\sigma}(\cD_T\times\bR_v^3)}
	\le\|u\|_{\cH_{\kin,\sigma}(\cD_{T,\ast}\times\bR_v^3)}+\|u(0)\|_{L^2_{x,v}}
	\le C\|u\|_{\cH_{\kin,\sigma}(\cD_{T,\ast}\times\bR_v^3)}
	\]
	with a universal constant $C$, which gives the bound on $\Sigma_T$.
	
	For $\cB_T$, the increment vanishes on $(0,T)$ for $t\ge\mathrm M+2$, because
	$\hat g^{(0)}=f_\st$ on $(t-\varepsilon,t)$ there ($\varkappa=0$ on $[\mathrm M+1,\infty)$ and $\varepsilon<1/4$).  So, in \eqref{eq:dyn-coordinate-B},
	the graph norm over $(t,T)$ vanishes for $t\ge\mathrm M+2$.  For $0\le t<\mathrm M+2$ we have
	$\ee^{\lambda t}\le\ee^{\lambda_{\rm exp}(\mathrm M+2)}$, since $\lambda\le\lambda_{\rm exp}$, and, by \eqref{260730346} for the value at $t$ and by \eqref{eq:dyn-unweighted-row-dominated},
	\[
	\bigl\|\hat g^{(0)}-\varrho_\varepsilon*\hat g^{(0)}\bigr\|_{\cH_{\kin,\sigma}(\cD_{(t,T)}\times\bR_v^3)}\le C'\,\Sigma_T\bigl[\hat g^{(0)}-\varrho_\varepsilon*\hat g^{(0)}\bigr],
	\]
	with $C'$ depending only on the structural data, and $\Sigma_T[\hat g^{(0)}-\varrho_\varepsilon*\hat g^{(0)}]\le C\varepsilon(\delta+\mathsf I)$ by the first inequality of
	\eqref{eq:dyn-real-increment-bound}.  Hence
	\[
	\begin{aligned}
		\cB_T\bigl[\hat g^{(0)}-\varrho_\varepsilon*\hat g^{(0)}\bigr]
		&\le\ee^{\lambda_{\rm exp}(\mathrm M+2)}\Bigl(C'\,\Sigma_T\bigl[\hat g^{(0)}-\varrho_\varepsilon*\hat g^{(0)}\bigr]\Bigr)^2\\
		&\le\bigl(C'\ee^{\lambda_{\rm exp}(\mathrm M+2)/2}\,C\,\varepsilon\,(\delta+\mathsf I)\bigr)^2 ,
	\end{aligned}
	\]
	which is the bound on $\cB_T$ in
	\eqref{eq:dyn-real-increment-bound}, after enlarging $C$ by the factor $C'\ee^{\lambda_{\rm exp}(\mathrm M+2)/2}$.
\end{proof}

Let $C$ be the constant of \eqref{eq:dyn-real-increment-bound}.  Fix
$0<\varepsilon_\lambda\le1/8$, depending only on the structural data, so
small that
\begin{equation}
	\label{eq:dyn-ball-eps-restriction}
	\ee^{\lambda\varepsilon_\lambda}\le\tfrac65,
	\qquad
	4C\varepsilon_\lambda\le\min\{1,\varepsilon_{\rm dyn},\varepsilon_{\mathrm{ell}},\varepsilon_0(\cD)\} .
\end{equation}
The mollification parameter is restricted to $0<\varepsilon\le\varepsilon_\lambda$.

The regularized background is controlled in the form
\begin{equation}
	\label{eq:dyn-ball-eps-A}
	\|w^{A^\dy_\flat}g_\dy\|_{L^\infty(\cD_T\times\bR^3)}
	+\|w^{A^\dy_\flat}g\|_{L^\infty(\cD_T\times\bR^3)}+\|w^{A^\dy_\flat}\nabla_vg\|_{L_t^\infty L_x^{p_\sharp}L_v^2}
	+\cB_T[g_\dy]^{1/2}+\cS_T[g_\dy]+\cN_{\mathrm M}(T)[g_\dy]\le\rho
\end{equation}
at a radius $\rho\ge0$, with $g_\dy=g-f_\st$.  Compared with $\widetilde{\mathfrak Q}^\dy_T(g)\le\rho$, \eqref{eq:dyn-ball-eps-A} omits the kinetic H\"older seminorm of $w^{A^\dy_\flat}g$ contained in $\cA_T[g]$, and adds the $L^\infty$ norm of $w^{A^\dy_\flat}g_\dy$; the H\"older control of $g^\varepsilon$ is provided instead by \eqref{eq:dyn-strong-coefficient-output-reg} for $\hat g^\varepsilon$, the extension of $g^\varepsilon$ by zero to $t\le-1$ (see after \eqref{eq:dyn-strong-past-extension}).

\begin{lemma}
	\label{lem:dyn-regbg-ball}
	Let $T\ge1$ and $0<\varepsilon\le\varepsilon_\lambda$, let $g$ satisfy \eqref{eq:dyn-ball-jets} and $\widetilde{\mathfrak Q}^\dy_T(g)\le r$ for some $r\ge0$, and assume \eqref{eq:dyn-strong-coefficient-choice}.  Let $\bar g_0$ and $g^\varepsilon$ be defined from $g$ by \eqref{eq:dynamic-base-background} and \eqref{eq:dyn-real-mollified-objects}.  There is $C$, depending only on the structural data, such that the following hold.
	\begin{enumerate}
		\item\label{it:dyn-regbg-ball}
		The regularized background $g^\varepsilon$, restricted to $(0,T)$, matches the initial jets and satisfies \eqref{eq:dyn-ball-eps-A} at $\rho=2r+C(\delta+\mathsf I)$.
		\item\label{it:dyn-regbg-coefficient}
		The extension $\hat g^\varepsilon$, with $\mu^{1/4}\hat g^\varepsilon$ extended to $\bR^7$ by Lemma~\ref{lem:norm-preserving-holder-extension}, satisfies \eqref{eq:dyn-strong-coefficient-output-reg}; in particular,
		\begin{equation}
			\label{eq:dyn-real-coefficient-hypotheses}
			\|\mu^{1/4}g^\varepsilon\|_{L^\infty(I_T^-\times\cD\times\bR^3)}\le\varepsilon_{\rm dyn} .
		\end{equation}
		\item\label{it:dyn-regbg-graph}
		Its deviation from $f_\st$ satisfies
		\begin{equation}
			\label{eq:dyn-regbg-graph}
			\|\hat g^\varepsilon-f_\st\|_{\cH_{\kin,\sigma}(\cD_{T,\ast}\times\bR^3_v)}\le C\bigl(\widetilde{\mathfrak Q}^\dy_T(g)+\delta+\mathsf I\bigr)\le C .
		\end{equation}
	\end{enumerate}
\end{lemma}

\begin{proof}
	By \eqref{eq:dyn-ball-jets} at $j=0$ and \eqref{eq:dynamic-initial-p-u}, $g(0)=z_0$, and $g\in\cH_{\kin,\sigma}(\cD_T\times\bR_v^3)$, since $g_\dy\in\cH_{\kin,\sigma}(\cD_T\times\bR_v^3)$ by \eqref{eq:dyn-ball-rows} and \eqref{eq:dyn-unweighted-row-dominated}, and $f_\st\in\cH_{\kin,\sigma}(\cD_T\times\bR_v^3)$ by \eqref{eq:dyn-theta-background}, as $T<\infty$.  Hence Lemma~\ref{lem:dyn-real-increment} applies.
	
	\emph{(\ref{it:dyn-regbg-ball}).}  The restriction of $g^\varepsilon$ to $(0,T)$ matches the initial jets
	by the third identity of \eqref{eq:dyn-real-increment-split}.  Since $\bar g_0=\hat g$ on $I_T^-$,
	the second identity at $j=0$ gives
	$g^\varepsilon=\varrho_\varepsilon*\hat g+(\hat g^{(0)}-\varrho_\varepsilon*\hat g^{(0)})$.  For
	$0<t<T$, $(\varrho_\varepsilon*\hat g)(t)$ is an average of $\hat g$ over $(t-\varepsilon,t)$,
	and on $(-\varepsilon,0)$ we have $\hat g=\check{\cP}=f_\st+\check{\pi}$, with
	$\check{\pi}$ of \eqref{eq:dyn-fp-seed}, since $\chi_-=1$ on $[-1/4,0]$ and
	$\varepsilon\le\varepsilon_\lambda\le1/8$.  So the lift $f_\st+\check{\pi}$ enters the average only for
	$0<t<\varepsilon$.  The increment
	$\hat g^{(0)}-\varrho_\varepsilon*\hat g^{(0)}$ satisfies
	\eqref{eq:dyn-real-increment-bound} with $C\varepsilon(\delta+\mathsf I)$ and $(C\varepsilon(\delta+\mathsf I))^2$ replaced by $C(\delta+\mathsf I)$ and $(C(\delta+\mathsf I))^2$, since $\varepsilon\le1$.  We bound the summands of \eqref{eq:dyn-ball-eps-A} at
	$g^\varepsilon$ in turn.
	\begin{itemize}
		\item By Minkowski's inequality, for each of the two norms $\|w^{A^\dy_\flat}\varphi\|_{L^\infty(\cD\times\bR^3)}$ and $\|w^{A^\dy_\flat}\nabla_v\varphi\|_{L_x^{p_\sharp}L_v^2}$ of $\varphi=(\varrho_\varepsilon*\hat g)(t)$, the essential supremum over $0<t<T$ is at most the larger of the essential suprema of that norm of $g(t)$ over $0<t<T$ and of $f_\st+\check{\pi}(t)$ over $-\varepsilon<t<0$.
		Summing the two bounds and applying the triangle inequality to $g^\varepsilon=\varrho_\varepsilon*\hat g+(\hat g^{(0)}-\varrho_\varepsilon*\hat g^{(0)})$,
		\[
		\|w^{A^\dy_\flat}g^\varepsilon\|_{L^\infty(\cD_T\times\bR^3)}+\|w^{A^\dy_\flat}\nabla_vg^\varepsilon\|_{L_t^\infty L_x^{p_\sharp}L_v^2}
		\le\cA_T[g]+C(\delta+\mathsf I),
		\]
		where, by Lemma~\ref{lem:dyn-fp-seed} and \eqref{eq:dyn-real-increment-bound},
		\[
		\begin{gathered}
			\|w^{A^\dy_\flat}(f_\st+\check{\pi})\|_{L^\infty((-\varepsilon,0)\times\cD\times\bR^3)}+\|w^{A^\dy_\flat}\nabla_v(f_\st+\check{\pi})\|_{L_t^\infty(-\varepsilon,0)L_x^{p_\sharp}L_v^2}\le C(\delta+\mathsf I),\\
			\|w^{A^\dy_\flat}(\hat g^{(0)}-\varrho_\varepsilon*\hat g^{(0)})\|_{L^\infty((-1,T)\times\cD\times\bR^3)}\\
			\quad+\|w^{A^\dy_\flat}\nabla_v(\hat g^{(0)}-\varrho_\varepsilon*\hat g^{(0)})\|_{L_t^\infty(-1,T)L_x^{p_\sharp}L_v^2}\le C(\delta+\mathsf I) .
		\end{gathered}
		\]
		For the first summand $\|w^{A^\dy_\flat}g^\varepsilon_\dy\|_{L^\infty}$ of \eqref{eq:dyn-ball-eps-A}, which is the $L^\infty$ summand of \eqref{eq:dyn-coordinate-A} evaluated at $g^\varepsilon_\dy$ rather than at $g^\varepsilon$, the contribution of the values $\hat g-f_\st=g-f_\st$ on $(0,T)$ to the average $\varrho_\varepsilon*(\hat g-f_\st)$ is bounded by $\|w^{A^\dy_\flat}(g-f_\st)\|_{L^\infty(\cD_T\times\bR^3)}$, which is at most
		\[
		\|w^{A^\dy_\flat}g\|_{L^\infty}+\|w^{A^\dy_\flat}f_\st\|_{L^\infty}\le\cA_T[g]+C\delta,
		\]
		and, with \eqref{eq:dyn-fp-seed-past} for the values $\hat g-f_\st=\check{\pi}$ on $(-\varepsilon,0)$ and \eqref{eq:dyn-real-increment-bound} for the increment,
		\[
		\|w^{A^\dy_\flat}g^\varepsilon_\dy\|_{L^\infty(\cD_T\times\bR^3)}\le\cA_T[g]+C(\delta+\mathsf I).
		\]
		\item For the rows \eqref{eq:dyn-hierarchy-rows}, with $I\in\mathfrak A_\bullet^m$, $0\le m\le\mathrm M$,
		\[
		\mathfrak D_Ig^\varepsilon_\dy=\varrho_\varepsilon*\mathfrak D_I(\hat g-f_\st)
		+\mathfrak D_I(\hat g^{(0)}-\varrho_\varepsilon*\hat g^{(0)}) .
		\]  By Minkowski's inequality, and since the translation in $t$ preserves the graph norm,
		\[
		\|w^{\ell^\dy_m}\varrho_\varepsilon*\mathfrak D_I(\hat g-f_\st)\|_{\cH_{\kin,\sigma}(\cD_T\times\bR_v^3)}
		\le\sup_{0<s<\varepsilon}\|w^{\ell^\dy_m}\mathfrak D_I(\hat g-f_\st)\|_{\cH_{\kin,\sigma}((-s,T-s)\times\cD\times\bR_v^3)} ,
		\]
		and the same holds without the weight $w^{\ell^\dy_m}$ and for the
		$L^2_tL^q_xH^1_\sigma$ norm of $\cS_T$.  Since $\hat g-f_\st$ has no jump at $t=0$ by
		\eqref{eq:dyn-ball-jets} and \eqref{eq:dynamic-negative-time-jet-matching}, it follows from \eqref{eq:dyn-time-gluing} at $(a,b,c)=(-s,0,T-s)$ that
		\[
		\begin{aligned}
			\|w^{\ell^\dy_m}\mathfrak D_I(\hat g-f_\st)\|^2_{\cH_{\kin,\sigma}((-s,T-s)\times\cD\times\bR_v^3)}
			\le{}&\|w^{\ell^\dy_m}\mathfrak D_I(\hat g-f_\st)\|^2_{\cH_{\kin,\sigma}((-s,0)\times\cD\times\bR_v^3)}\\
			&+\|w^{\ell^\dy_m}\mathfrak D_I(\hat g-f_\st)\|^2_{\cH_{\kin,\sigma}(\cD_T\times\bR_v^3)} ,
		\end{aligned}
		\]
		and the same holds without the weight $w^{\ell^\dy_m}$; moreover, since $(-s,T-s)\subset(-s,0]\cup(0,T)$,
		\[
		\|\hat g-f_\st\|^2_{L^2((-s,T-s);L^q_xH^1_\sigma)}
		\le\|\hat g-f_\st\|^2_{L^2((-s,0);L^q_xH^1_\sigma)}
		+\|\hat g-f_\st\|^2_{L^2((0,T);L^q_xH^1_\sigma)} .
		\]
		With \eqref{eq:dyn-fp-seed-past} for the past part and
		\eqref{eq:dyn-real-increment-bound} for the increment,
		\[
		\cS_T[g^\varepsilon_\dy]+\cN_{\mathrm M}(T)[g^\varepsilon_\dy]\le\Sigma_T[g_\dy]+C(\delta+\mathsf I) .
		\]
		\item For the decay coordinate, use the same splitting, in which the interval $(t-s,0)$
		occurs only when $t<\varepsilon$, together with
		$\ee^{\lambda t}\le\ee^{\lambda\varepsilon}\ee^{\lambda(t-s)}$.  With
		\eqref{eq:dyn-fp-seed-past} for the past part and $\ee^{\lambda\varepsilon}\le\tfrac65$ by
		\eqref{eq:dyn-ball-eps-restriction}, this gives
		\[
		\cB_T\bigl[\varrho_\varepsilon*(\hat g-f_\st)\bigr]
		\le\tfrac65\bigl(\cB_T[g_\dy]+C(\mathsf I^2+\delta^2)\bigr) .
		\]
		The square root of $\cB_T$ is subadditive, and
		\[
		\cB_T[\hat g^{(0)}-\varrho_\varepsilon*\hat g^{(0)}]^{1/2}\le C\varepsilon(\delta+\mathsf I)\le C(\delta+\mathsf I)
		\]
		by \eqref{eq:dyn-real-increment-bound}, so
		\[
		\begin{aligned}
			\cB_T[g^\varepsilon_\dy]^{1/2}
			&\le\bigl(\tfrac65\bigr)^{1/2}\bigl(\cB_T[g_\dy]+C(\mathsf I^2+\delta^2)\bigr)^{1/2}+C(\delta+\mathsf I)\\
			&\le2\cB_T[g_\dy]^{1/2}+C(\delta+\mathsf I) .
		\end{aligned}
		\]
	\end{itemize}
	Adding the four bounds, we see that, with $g^\varepsilon$ in place of $g$ in \eqref{eq:dyn-ball-eps-A},
	\[
	\begin{aligned}
		&\|w^{A^\dy_\flat}g^\varepsilon_\dy\|_{L^\infty(\cD_T\times\bR^3)}
		+\|w^{A^\dy_\flat}g^\varepsilon\|_{L^\infty(\cD_T\times\bR^3)}+\|w^{A^\dy_\flat}\nabla_vg^\varepsilon\|_{L_t^\infty L_x^{p_\sharp}L_v^2}
		+\cB_T[g^\varepsilon_\dy]^{1/2}+\cS_T[g^\varepsilon_\dy]+\cN_{\mathrm M}(T)[g^\varepsilon_\dy]\\
		&\qquad\le2\widetilde{\mathfrak Q}^\dy_T(g)+C(\delta+\mathsf I)\le2r+C(\delta+\mathsf I),
	\end{aligned}
	\]
	where the last inequality follows from the hypothesis $\widetilde{\mathfrak Q}^\dy_T(g)\le r$.  This is the bound in (\ref{it:dyn-regbg-ball}).
	
	\emph{(\ref{it:dyn-regbg-coefficient}).}  The bound \eqref{eq:dyn-strong-coefficient-output-reg} for $\hat g^\varepsilon$ follows from \eqref{eq:dyn-strong-coefficient-output} for $\hat g$, \eqref{eq:dyn-strong-mollification-holder} and the increment
	bound.  Since $g(0)=z_0$, Lemma~\ref{lem:dyn-strong-coefficient-class} applies to $g$ under
	\eqref{eq:dyn-strong-coefficient-choice}, and gives \eqref{eq:dyn-strong-coefficient-output} for
	$\hat g$.  It suffices to bound the two norms of \eqref{eq:dyn-strong-coefficient-output-reg} for $\hat g^\varepsilon$ on $Q$, since
	the extension of Lemma~\ref{lem:norm-preserving-holder-extension} off $\overline\cD$ and past
	$T$, applied to $\mu^{1/4}\hat g^\varepsilon$, preserves both the supremum and the seminorm.  On $Q$,
	\eqref{eq:dyn-strong-coefficient-output} holds with $\varrho_\varepsilon*\hat g$ in place of $\hat g$, by
	\eqref{eq:dyn-strong-mollification-holder} and, for the supremum, Minkowski's inequality in $t$.  The increment satisfies
	\[
	\|\mu^{1/4}(\hat g^{(0)}-\varrho_\varepsilon*\hat g^{(0)})\|_{L^\infty(Q)}\le\|\mu^{1/4}(\hat g^{(0)}-\varrho_\varepsilon*\hat g^{(0)})\|_{C^{\alpha_{\rm bs}}_{\kin}(Q)}\le C\varepsilon(\delta+\mathsf I),
	\]
	with $C$ the constant of
	\eqref{eq:dyn-real-increment-bound}.  By $\mathsf I\le1$, $\delta\le1$ and
	\eqref{eq:dyn-ball-eps-restriction},
	\[
	C\varepsilon(\delta+\mathsf I)\le2C\varepsilon_\lambda\le\tfrac12\min\{1,\varepsilon_{\rm dyn},
	\varepsilon_{\mathrm{ell}},\varepsilon_0(\cD)\} .
	\]
	So
	$\hat g^\varepsilon$ satisfies \eqref{eq:dyn-strong-coefficient-output-reg}.  Since $\hat g^\varepsilon=g^\varepsilon$ on $I_T^-$, the first bound of \eqref{eq:dyn-strong-coefficient-output-reg} gives \eqref{eq:dyn-real-coefficient-hypotheses}.
	
	\emph{(\ref{it:dyn-regbg-graph}).}  We apply (\ref{it:dyn-regbg-ball}) with $r=\widetilde{\mathfrak Q}^\dy_T(g)$.  The term $t=0$ of \eqref{eq:dyn-coordinate-B} then gives
	\[
	\|g^\varepsilon-f_\st\|_{\cH_{\kin,\sigma}(\cD_T\times\bR^3_v)}\le\cB_T[g^\varepsilon_\dy]^{1/2}\le2\widetilde{\mathfrak Q}^\dy_T(g)+C(\delta+\mathsf I) .
	\]
	On $(-1,0]$, $g^\varepsilon=\bar g_0=\check{E}$ by the first identity of \eqref{eq:dyn-real-increment-split} and by \eqref{eq:dynamic-base-background}, and
	\[
	\|\check{E}-f_\st\|_{\cH_{\kin,\sigma}((-1,0)\times\cD\times\bR^3_v)}\le C(\delta+\mathsf I)
	\]
	by \eqref{eq:dyn-lift-graph} for $\check{E}$ and \eqref{eq:dyn-theta-background} for $f_\st$.  Since $g^\varepsilon(0)=f_\st+\zeta_0=z_0=\check{E}(0)$ by the third identity of \eqref{eq:dyn-real-increment-split} at $j=0$ and by \eqref{eq:dynamic-negative-time-jet-matching}, \eqref{eq:dyn-time-gluing} at $(a,b,c)=(-1,0,T)$ applies to $\hat g^\varepsilon-f_\st$ on $\cD_{T,\ast}$; with the two bounds above, it gives the first inequality of \eqref{eq:dyn-regbg-graph}.  The second follows from \eqref{eq:dyn-strong-coefficient-choice}, which gives $\widetilde{\mathfrak Q}^\dy_T(g)+\mathsf I\le\eta_{\star,1}\le1$, and from $\delta\le1$. This completes the proof.
\end{proof}

The regularized problem on $\cD_{T,\ast}$ has zero data at $t=-1$. Its temporal jets satisfy the differentiated equations and boundary laws below. Write
\begin{equation}
	\label{eq:dyn-real-centered-unknown}
	f_\dy^{\varepsilon,N}:=f_\pt^{\varepsilon,N}-f_\st ,
	\qquad
	\mathfrak D_{\mathrm t}^{\,j}f_\dy^{\varepsilon,N}:=\mathfrak D_{\mathrm t}^{\,j}f_\pt^{\varepsilon,N}
	\quad(1\le j\le\mathrm M) .
\end{equation}

\begin{prop}[Regularized problem]
	\label{prop:dyn-fixedN-realization}
	Let $T\ge1$, $N\in\bN$ and $0<\varepsilon\le\varepsilon_\lambda$, let $g$ satisfy \eqref{eq:dyn-ball-jets} and
	$\widetilde{\mathfrak Q}^\dy_T(g)\le r$ for some $r\ge0$, and assume
	\eqref{eq:dyn-strong-coefficient-choice}.  Let $\bar g_0$ and $g^\varepsilon$ be defined from $g$ by
	\eqref{eq:dynamic-base-background} and \eqref{eq:dyn-real-mollified-objects}.
	\begin{enumerate}
		\item\label{it:dyn-real-solution}
		There is a unique $f_\pt^{\varepsilon,N}\in\cH_{\kin,\sigma}(\cD_{T,\ast}\times\bR_v^3)$ with
		\begin{equation}
			\label{eq:dyn-real-mollified-problem}
			(\partial_t+v\cdot\nabla_x)f_\pt^{\varepsilon,N}
			+L_{g^\varepsilon}f_\pt^{\varepsilon,N}=\bar F_0,
			\quad
			f_\pt^{\varepsilon,N}(-1)=0,
			\quad
			\bigl(f_\pt^{\varepsilon,N}\bigr)_{\gamma^-}
			=P_{w,N}\bigl(f_\pt^{\varepsilon,N}\bigr)_{\gamma^+}+\bar b_N .
		\end{equation}
		It satisfies
		\begin{equation}
			\label{eq:dyn-real-uepsN-graph}
			\bigl\|f_\pt^{\varepsilon,N}\bigr\|_{\cH_{\kin,\sigma}(\cD_{T,\ast}\times\bR_v^3)}
			\le C_T ,
		\end{equation}
		with $C_T$ independent of $N$ and $\varepsilon$.  Moreover, $f_\pt^{\varepsilon,N}$ vanishes on $(-1,-3/4]$.
		\item\label{it:dyn-real-prop-chain}
		The temporal jets $\mathfrak D_{\mathrm t}^{\,j}f_\pt^{\varepsilon,N}$, $0\le j\le\mathrm M$, exist in
		$\cH_{\kin,\sigma}(\cD_{T,\ast}\times\bR_v^3)$, and
		\begin{equation}
			\label{eq:dyn-real-prop-chain}
			\partial_tg^\varepsilon_j=g^\varepsilon_{j+1},
			\qquad
			\partial_t\mathfrak D_{\mathrm t}^{\,j}f_\pt^{\varepsilon,N}
			=\mathfrak D_{\mathrm t}^{\,j+1}f_\pt^{\varepsilon,N},
			\qquad 0\le j\le\mathrm M-1 .
		\end{equation}
		
		\item\label{it:dyn-real-prop-equations}
		For $1\le j\le\mathrm M$,
		\begin{equation}
			\label{eq:dyn-real-prop-equations}
			(\partial_t+v\cdot\nabla_x)\mathfrak D_{\mathrm t}^{\,j}f_\pt^{\varepsilon,N}
			+L_{g^\varepsilon}\mathfrak D_{\mathrm t}^{\,j}f_\pt^{\varepsilon,N}
			=\sum_{a=1}^{j}\binom ja
			N\bigl[g^\varepsilon_a,\ \mathfrak D_{\mathrm t}^{\,j-a}f_\pt^{\varepsilon,N}\bigr]
			+\bar F_j .
		\end{equation}
		\item\label{it:dyn-real-prop-boundary}
		For $0\le j\le\mathrm M$,
		\begin{equation}
			\label{eq:dyn-real-prop-boundary}
			\bigl(\mathfrak D_{\mathrm t}^{\,j}f_\pt^{\varepsilon,N}\bigr)_{\gamma^-}
			=P_{w,N}\bigl(\mathfrak D_{\mathrm t}^{\,j}f_\pt^{\varepsilon,N}\bigr)_{\gamma^+}
			+\mathfrak t_{j,N}\bigl[f_\pt^{\varepsilon,N}\bigr]+\partial_t^{\,j}\bar b_N .
		\end{equation}
		\item\label{it:dyn-real-prop-initial}
		For $0\le j\le\mathrm M$,
		\begin{equation}
			\label{eq:dyn-real-prop-initial}
			\mathfrak D_{\mathrm t}^{\,j}f_\pt^{\varepsilon,N}=\partial_t^j\check{E}
			\quad\text{on }(-1,0),
			\qquad
			\mathfrak D_{\mathrm t}^{\,j}f_\pt^{\varepsilon,N}(0)=z_j .
		\end{equation}
		\item\label{it:dyn-real-prop-spatial}
		For every ordered word $I$ of $\mathfrak A_\bullet$ with $|I|\le\mathrm M$,
		$\mathfrak D_If_\pt^{\varepsilon,N}\in\cH_{\kin,\sigma}(\cD_{T,\ast}\times\bR_v^3)$ with
		completion-defined traces, and the rows of
		\eqref{eq:dyn-hierarchy-rows} satisfy $\mathsf U_m(T)[f_\dy^{\varepsilon,N}]+\mathsf W_m(T)[f_\dy^{\varepsilon,N}]<\infty$ for $0\le m\le\mathrm M$.
		\item\label{it:dyn-real-restriction}
		The jets of $f_\pt^{\varepsilon,N}$ form a chain on $I_T^-$, and hence on $(0,T)$: $\mathfrak D_\nu(\mathfrak D_Jf_\pt^{\varepsilon,N})=\mathfrak D_{\nu J}f_\pt^{\varepsilon,N}$ in $\cD'$ for $\nu\in\mathfrak A_\bullet$ and every word $\nu J$ of the word set
		\begin{equation}
			\label{eq:dyn-real-spatial-words-enlarged}
			\bigl\{I\in\mathfrak A_\bullet^m:\ m\le\mathrm M\bigr\} .
		\end{equation}
		
		\item\label{it:dyn-real-mass}
		If in addition \eqref{eq:dynamic-initial-zero-perturbation-mass} holds, then $f_\dy^{\varepsilon,N}$ of \eqref{eq:dyn-real-centered-unknown}
		satisfies
		\begin{equation}
			\label{eq:dyn-real-prop-mass}
			\mathfrak m_{f_\dy^{\varepsilon,N}}(t)=0,
			\qquad 0\le t\le T,
		\end{equation}
		with $\mathfrak m$ the mass scalar \eqref{eq:dynamic-mass-scalar}.
		\item\label{it:dyn-real-strong}
		$\cA_T[f_\dy^{\varepsilon,N}]<\infty$ and, with $\varsigma_N$ of
		\eqref{eq:dyn-limits-corner-size} and with a constant $C$ depending only on the
		structural data (in particular independent of $T$, $N$, $\varepsilon$ and $r$),
		\begin{equation}
			\label{eq:dyn-strong-A-linear-defect}
			\cA_T\bigl[f_\dy^{\varepsilon,N}\bigr]\le C\bigl(\delta+\mathsf I+\varsigma_N
			+\cN_{\mathrm M}(T)\bigl[f_\dy^{\varepsilon,N}\bigr]\bigr).
		\end{equation}
	\end{enumerate}
\end{prop}

The proof in \S\ref{app:dyn-doubly-regularized} uses only the preceding results and the preliminary facts of \S\ref{ss:dyn-past-mollification}.

For the estimates in the next subsection we rewrite the system in terms of $f_\dy^{\varepsilon,N}$. Under the hypotheses of Proposition~\ref{prop:dyn-fixedN-realization}, with
$g^\varepsilon_\dy=g^\varepsilon-f_\st$ as after \eqref{eq:dyn-real-mollified-objects}, on $\cD_{T,\ast}$, $f_\dy^{\varepsilon,N}$ satisfies
\begin{equation}
	\label{eq:dyn-limits-fixed-system}
	\begin{cases}
		(\partial_t+v\cdot\nabla_x)f_\dy^{\varepsilon,N}+L_{g^\varepsilon}f_\dy^{\varepsilon,N}
		=N[g^\varepsilon_\dy,f_\st]
		&\text{on }(0,T),\\[2pt]
		(\partial_t+v\cdot\nabla_x)f_\dy^{\varepsilon,N}+L_{g^\varepsilon}f_\dy^{\varepsilon,N}
		=N[g^\varepsilon_\dy,f_\st]+\bar F_0
		&\text{on }(-1,0],\\[2pt]
		\bigl(\mathfrak D_{\mathrm t}^{\,j}f_\dy^{\varepsilon,N}\bigr)_{\gamma^-}
		=P_{w,N}\bigl(\mathfrak D_{\mathrm t}^{\,j}f_\dy^{\varepsilon,N}\bigr)_{\gamma^+}
		+\mathfrak t_{j,N}\bigl[f_\dy^{\varepsilon,N}\bigr]+\partial_t^{\,j}\bigl(\bar b_N-\mathfrak k_N\bigr)
		\quad\text{on }\gamma^-_{T,\ast},
		&0\le j\le\mathrm M,\\[2pt]
		f_\dy^{\varepsilon,N}(-1)=-f_\st,
		\quad
		\mathfrak D_{\mathrm t}^{\,j}f_\dy^{\varepsilon,N}(0)=\zeta_j,
		&0\le j\le\mathrm M ,
	\end{cases}
\end{equation}
with
\begin{equation}
	\label{eq:dyn-limits-fixed-datum}
	\bar b_N-\mathfrak k_N=
	\begin{cases}
		\check{b}_N-\mathfrak k_N,&-1<t<0,\\
		b_w+c_N,&0<t<T,
	\end{cases}
\end{equation}
and $c_N$ of \eqref{eq:dyn-limits-corner}.  The first two
lines are \eqref{eq:dyn-real-mollified-problem} minus the first line of \eqref{eq:stationary-full-diffuse-problem}, which reads \[
v\cdot\nabla_xf_\st+L_{g^\varepsilon}f_\st=-N[g^\varepsilon_\dy,f_\st]
\] since $L_{g^\varepsilon}=L-N[g^\varepsilon,\cdot\,]$ by \eqref{def:L_g}; the source vanishes on $(0,T)$, since
$\bar F_0=0$ on $(0,T)$ by \eqref{eq:dynamic-zero-initial-extended-sources}.  The third is
\eqref{eq:dyn-real-prop-boundary} minus the definition of $\mathfrak k_N$ in \eqref{eq:dyn-limits-steady-residual}; here, for $j\ge1$, since $f_\st$ does not depend on $t$ and $P_{w,N}(t)=(\mathfrak m_\circ+d_\delta(t))J_N\cA_0$ by \eqref{eq:boundary-spectral-diffuse-operator} and \eqref{eq:wall-response-profiles}, the definitions \eqref{eq:dyn-theta-triangular-notation} and \eqref{eq:dyn-limits-steady-residual} give
\[
\mathfrak t_{j,N}[f_\pt^{\varepsilon,N}]=\mathfrak t_{j,N}[f_\dy^{\varepsilon,N}]
+(\partial_t^jd_\delta)J_N\cA_0(f_\st)_{\gamma^+},
\qquad
(\partial_t^jd_\delta)J_N\cA_0(f_\st)_{\gamma^+}=-\partial_t^j\mathfrak k_N .
\]
The fourth is $f_\pt^{\varepsilon,N}(-1)=0$ and \eqref{eq:dyn-real-prop-initial}, since
$\mathfrak D_{\mathrm t}^{\,j}f_\dy^{\varepsilon,N}-\zeta_j
=\mathfrak D_{\mathrm t}^{\,j}f_\pt^{\varepsilon,N}-z_j$ by
\eqref{eq:dyn-real-centered-unknown} and $\zeta_j=z_j-\delta_{j0}f_\st$.
Since $f_\pt^{\varepsilon,N}=0$ on $(-1,-3/4]$ by
Proposition~\ref{prop:dyn-fixedN-realization}\,(\ref{it:dyn-real-solution}),
\eqref{eq:dyn-real-prop-chain} gives
$\mathfrak D_{\mathrm t}^{\,j}f_\dy^{\varepsilon,N}(-1)
=\mathfrak D_{\mathrm t}^{\,j}f_\pt^{\varepsilon,N}(-1)=0$ for $1\le j\le\mathrm M$.

\begin{lemma}
	\label{lem:dyn-limits-wall-law}
	Let $T\ge1$, $N\in\bN$ and $0<\varepsilon\le\varepsilon_\lambda$, let $g$ satisfy \eqref{eq:dyn-ball-jets} and
	$\widetilde{\mathfrak Q}^\dy_T(g)\le r$ for some $r\ge0$, and assume
	\eqref{eq:dyn-strong-coefficient-choice}.  Let $f_\dy^{\varepsilon,N}$ be given by \eqref{eq:dyn-real-centered-unknown}, and let $I$ be an
	ordered word of $\mathfrak A_\bullet$ with $1\le|I|\le\mathrm M$,
	spatial subword $I'$ of length $m'$ and temporal order $j=j(I)$.  Then
	\begin{equation}
		\label{eq:dyn-limits-temporal-triangular}
		\bigl(\mathfrak D_If_\dy^{\varepsilon,N}\bigr)_{\gamma^-}
		=P_{w,N}\bigl(\mathfrak D_If_\dy^{\varepsilon,N}\bigr)_{\gamma^+}
		+\mathfrak t_{I,N}[f_\dy^{\varepsilon,N}]+\mathfrak r_{I,N}
		\qquad\text{in }L^2(\gamma^-_T),
	\end{equation}
	with
	\begin{equation}
		\label{eq:dyn-limits-wall-mixed-terms}
		\mathfrak t_{I,N}[f_\dy^{\varepsilon,N}]:=\mathfrak D_{I'}^-\mathfrak t_{j,N}[f_\dy^{\varepsilon,N}],
		\qquad
		\mathfrak r_{I,N}:=\mathfrak b_{I,N}\bigl[(f_\dy^{\varepsilon,N})_{\gamma^+}\bigr]
		+\mathfrak D_{I'}^-\partial_t^{\,j}\bigl(\bar b_N-\mathfrak k_N\bigr) .
	\end{equation}
	Here $\mathfrak D_{I'}^-\mathfrak t_{j,N}[f_\dy^{\varepsilon,N}]$ is the finite sum
	\eqref{eq:dyn-temporal-part-spatial-expansion-N}, the bracket
	$\mathfrak b_{I,N}[\,\cdot\,]$ is the coloring of
	\eqref{eq:dyn-mixed-boundary-law-N}, and
	$\mathfrak D_{I'}^-\partial_t^j(\bar b_N-\mathfrak k_N)$ is the weak derivative of
	\eqref{eq:dyn-datum-spatial-regularity}.  For every $\ell\ge0$,
	\begin{equation}
		\label{eq:dyn-limits-temporal-triangular-bound}
		\bigl\|w^{\ell}\mathfrak t_{I,N}[f_\dy^{\varepsilon,N}]\bigr\|_{L^2(\gamma^-_T)}
		\le C_{I,\ell}\,\delta\sum_{|J|\le|I|,\ j(J)<j(I)}
		\bigl\|(\mathfrak D_Jf_\dy^{\varepsilon,N})_{\gamma^+}\bigr\|_{L^2(\gamma^+_T)} ,
	\end{equation}
	\begin{equation}
		\label{eq:dyn-limits-wall-coloring-bound}
		\bigl\|w^{\ell}\mathfrak b_{I,N}\bigl[(f_\dy^{\varepsilon,N})_{\gamma^+}\bigr]\bigr\|
		_{L^2(\gamma^-_T)}
		\le C_{I,\ell}\sum_{|J|<|I|,\ j(J)=j(I)}
		\bigl\|(\mathfrak D_Jf_\dy^{\varepsilon,N})_{\gamma^+}\bigr\|_{L^2(\gamma^+_T)} ,
	\end{equation}
	with $C_{I,\ell}$ depending only on the structural data, $I$ and $\ell$ (in particular independent of $T$, $N$ and $\varepsilon$).
\end{lemma}

\begin{proof}
	Write $f_\dy=f_\dy^{\varepsilon,N}|_{(0,T)}$. For a purely temporal word, $m'=0$, $\mathfrak b_{I,N}[(f_\dy)_{\gamma^+}]=0$,
	$\mathfrak t_{I,N}=\mathfrak t_{j,N}$ and $\mathfrak r_{I,N}=\partial_t^j(\bar b_N-\mathfrak k_N)$, so
	\eqref{eq:dyn-limits-temporal-triangular} is the third line of
	\eqref{eq:dyn-limits-fixed-system}.  Moreover, \eqref{eq:dyn-limits-temporal-triangular-bound}
	follows from \eqref{eq:dyn-theta-triangular-bound} of
	Lemma~\ref{lem:ordered-boundary-finite-rank}\,(\ref{it:finite-rank-evolutionary})
	for $\mathfrak t_{j,N}$.  At a word with $m'\ge1$ we use Proposition~\ref{prop:spectral-triangular-evolutionary} at
	$m_\ast=|I|\le\mathrm M$, with $f=f_\dy$ and the datum $b=b_w+c_N$ of
	\eqref{eq:dyn-limits-fixed-datum} on $(0,T)$; its constants depend on the wall only through
	Lemma~\ref{lem:dyn-theta-profile-time-derivatives}.  The
	hypotheses of that proposition hold.
	\begin{itemize}
		\item The wall law \eqref{eq:dyn-temporal-boundary-hypothesis-N} is the third line of
		\eqref{eq:dyn-limits-fixed-system} at $j=0$.
		\item Since $f_\dy=f_\pt^{\varepsilon,N}-f_\st$ on $(0,T)$, the jets
		$\mathfrak D_Jf_\dy$, $|J|\le|I|$, are in $\cH_{\kin,\sigma}(\cD_T\times\bR_v^3)$, satisfy the
		jet chain \eqref{eq:dyn-jet-chain} and carry completion-defined traces: for $f_\pt^{\varepsilon,N}$ by
		Proposition~\ref{prop:dyn-fixedN-realization}\,(\ref{it:dyn-real-prop-spatial}),~(\ref{it:dyn-real-restriction}),
		and for $f_\st$ by \eqref{eq:dyn-theta-background}, \eqref{eq:unweighted-row-dominated} and
		\eqref{eq:dyn-fs-time-independent}, which give $\mathfrak D_Kf_\st\in\cH_{\kin,\sigma}$ for $|K|\le\mathrm N$.
		\item The datum satisfies \eqref{eq:dyn-datum-spatial-regularity} for $|K|+i\le|I|$: it is
		a finite sum of products of $d_\delta(t)-d_\delta^\infty$ or $\mathfrak m_\circ+d_\delta(t)$ with scalars in
		$H^{\mathrm M}(\partial\cD)$; its weak derivatives are given by
		\eqref{eq:boundary-word-Leibniz} and are in $L^2(\gamma^-_T)$ by
		\eqref{eq:dyn-limits-corner-bound} and \eqref{eq:dyn-limits-wall-datum-size}.
	\end{itemize}
	Then \eqref{eq:dyn-mixed-boundary-law-N} is \eqref{eq:dyn-limits-temporal-triangular} with
	\eqref{eq:dyn-limits-wall-mixed-terms}, and
	\eqref{eq:dyn-temporal-part-spatial-bound} and \eqref{eq:dyn-mixed-boundary-bound} of
	Proposition~\ref{prop:spectral-triangular-evolutionary} give \eqref{eq:dyn-limits-temporal-triangular-bound} and
	\eqref{eq:dyn-limits-wall-coloring-bound}.
\end{proof}

\subsection{The hierarchy and the fixed-background solution}\label{ss:dyn-coordinates}

For backgrounds in the prescribed small class, the fixed-background problem \eqref{eq:dyn-perturbation-problem} admits a reduced estimate uniform in the time horizon and regularization parameters. This estimate is used in the nonlinear fixed point.

\begin{thm}[The fixed-background solution]
	\label{thm:dyn-physical-solution}
	There is $\eta_\star\in(0,1]$, depending only on the structural data, such that, after
	decreasing $\delta_\star$, if necessary, the following holds.  Let $T\ge1$, let the datum satisfy
	\eqref{eq:dynamic-initial-zero-perturbation-mass} and $\delta+\mathsf I\le\delta_\star$, and
	let $g$ satisfy \eqref{eq:dyn-ball-jets} and $\widetilde{\mathfrak Q}^\dy_T(g)\le\eta_\star$.
	Then:
	\begin{enumerate}
		\item\label{it:dyn-phys-existence} there is exactly one
		$f_\dy\in\cH_{\kin,\sigma}(\cD_T\times\bR_v^3)$ solving \eqref{eq:dyn-perturbation-problem} with
		$\partial_t^{\,j}f_\dy(0)=\zeta_j$ for $1\le j\le\mathrm M$;
		\item\label{it:dyn-phys-mass} $\mathfrak m_{f_\dy}(t)=0$ for $0\le t\le T$;
		\item\label{it:dyn-phys-coords} with $C$ depending only on the structural data (in particular independent of $g$, $\delta$ and $\mathsf I$),
		\begin{equation}
			\label{eq:dyn-physical-coordinates}
			\widetilde{\mathfrak Q}^\dy_T\bigl(f_\st+f_\dy\bigr)\le C\bigl(\delta+\mathsf I\bigr);
		\end{equation}
		\item\label{it:dyn-phys-halfline} if $g$ is defined on $(0,\infty)$ and
		$\widetilde{\mathfrak Q}^\dy_\infty(g)\le\eta_\star$, there is exactly one
		$f_\dy\in\bigcap_{T\ge1}\cH_{\kin,\sigma}(\cD_T\times\bR_v^3)$ solving
		the problem of (\ref{it:dyn-phys-existence}) at every $T\ge1$, and it satisfies
		\eqref{eq:dyn-physical-coordinates} at every $T$ with the same constant.
	\end{enumerate}
\end{thm}

Parts (\ref{it:dyn-phys-existence})--(\ref{it:dyn-phys-coords}) of Theorem~\ref{thm:dyn-physical-solution} follow in \S\ref{app:dyn-limit-removal} from the fixed-parameter construction and the uniform estimates proved below.

\begin{proof}[\textbf{Proof of Theorem~\ref{thm:dyn-physical-solution}\,(\ref{it:dyn-phys-halfline})}]
	By monotonicity in $T$ (see Definition~\ref{def:dyn-coordinates}), $\widetilde{\mathfrak Q}^\dy_T(g)\le\widetilde{\mathfrak Q}^\dy_\infty(g)\le\eta_\star$. Parts (\ref{it:dyn-phys-existence})--(\ref{it:dyn-phys-coords}) therefore apply on every finite interval. Uniqueness makes these solutions consistent under restriction, defining a solution on $(0,\infty)$ that is unique on each finite interval. The constant in \eqref{eq:dyn-physical-coordinates} is independent of $T$, as required.
\end{proof}

We establish bounds for $f_\dy^{\varepsilon,N}$ in the coordinates and weighted rows, uniformly in $T$, $\varepsilon$ and $N$. The only regularization error in these bounds is $\varsigma_N$ from \eqref{eq:dyn-limits-corner-size}, which tends to zero by Lemma~\ref{lem:dyn-limits-wall}\,(\ref{it:dyn-limits-wall-size}). The following corollary applies both to $\cJ=J_N$ and, upon removal of the cutoff in \S\ref{app:dyn-limit-removal}, to $\cJ=\mathbf I$. We fix the smallness thresholds at the end of the subsection.

\begin{corollary}[The decay coordinate]
	\label{cor:dyn-decay-B}
	Let $T\ge1$ and let $g$ satisfy either
	$\widetilde{\mathfrak Q}^\dy_T(g)\le\rho$ or \eqref{eq:dyn-ball-eps-A},
	where
	\begin{equation}
		\label{eq:dyn-decay-B-smallness}
		0\le\rho\le1,\qquad \rho+\delta\le\varepsilon_{\rm dec} .
	\end{equation}
	Here $\varepsilon_{\rm dec}>0$ is a sufficiently small constant depending only on the structural data.
	Let either $\cJ=\mathbf I$ and $b=b_w$, or $\cJ=J_N$ and $b=b_w+c_N$ for some
	$N\in\bN$, with $c_N$ of \eqref{eq:dyn-limits-corner}; in the first case, replace
	$\varsigma_N$ by $0$ below.  Let $u\in\cH_{\kin,\sigma}(\cD_T\times\bR_v^3)$ solve \eqref{eq:dyn-perturbation-problem} with $u$, $P_\cJ(t)$ and $b$ in place
	of $f_\dy$, $P_w(t)$ and $b_w(t)$, and assume
	$\mathfrak m_u(t)=0$ for $0\le t\le T$.  Then
	\begin{equation}
		\label{eq:dyn-decay-B-closed}
		\cB_T[u]\le 			 C\bigl(\mathsf I^2+\delta^2+\varsigma_N^2\bigr)
	\end{equation}
	and
	\begin{equation}
		\label{eq:dyn-base-input-B-mass-full}
		\|u(t)\|_{L^2_{x,v}}^2+\cD_u(t,T)
		\le C\bigl(\mathsf I^2+\delta^2+\varsigma_N^2\bigr)\ee^{-\lambda t},
		\qquad 0\le t\le T ,
	\end{equation}
	with $C$ depending only on the structural data (in particular independent of $T$, $N$,
	$\varepsilon$, $\delta$, $\mathsf I$ and $\rho$).
\end{corollary}

\begin{proof}
	The stationary estimate of Lemma~\ref{lem:dyn-fp-seed} gives
	$\|w^{A^\dy_\flat}f_\st\|_{L^\infty}\lesssim\delta$.
	Consequently, both background hypotheses in the statement imply
	\[
	\|w^{A^\dy_\flat}g_\dy\|_{L^\infty}\lesssim\rho+\delta,
	\qquad \cB_T[g_\dy]\le\rho^2.
	\]
	Fix the structural constant in the first inequality and take
	$\varepsilon_{\rm dec}$ small enough that
	\eqref{eq:dyn-decay-background-hypotheses} holds with
	$\varepsilon_\dy\le\varepsilon_\star$ and $\mathsf G=\rho^2$.
	Proposition~\ref{prop:dynamic-global-exponential-stability}
	therefore applies at the rate in \eqref{eq:dyn-decay-rate-fixed-point}.

	For $\cJ=\mathbf I$, the additional boundary datum vanishes and $K_c=0$. For $\cJ=J_N$, the boundary law is
	\eqref{eq:dyn-decay-general-wall} at $c=c_N$, and we show that \eqref{eq:dyn-decay-extra-datum} holds with
	$K_c=C\varsigma_N^2$, $C$ depending only on the structural data.
	For the first condition of \eqref{eq:dyn-decay-extra-datum}, we apply \eqref{eq:dynamic-physical-mass-balance} at $\alpha=1$,
	$\cJ=J_N$, $S=N[g,u]+\mathfrak S$ by \eqref{def:L_g}, and $b=b_w+c_N$.
	Since
	$$
	\mathfrak m_u\equiv0,\quad \langle S,\sqrt\mu\rangle_v=0,\quad \beta_{b_w}\equiv0
	$$ (the last by
	\eqref{eq:dyn-theta-residual-flux}), the mass balance \eqref{eq:dynamic-physical-mass-balance} gives
	$\int_{t_0}^t\int_{\partial\cD}\beta_{c_N}\dd S_x\dd r=0$ for all $0\le t_0<t\le T$; hence
	$\int_{\partial\cD}\beta_{c_N}\dd S_x=0$ for a.e.\ $t$.
	For the second, $c_N=0$ for $t>\mathrm M+1$ by \eqref{eq:dynamic-common-cutoff} and
	$\|c_N\|_{L^2(\gamma^-_T)}\le C\varsigma_N$ by \eqref{eq:dyn-limits-corner-bound}. Since $\lambda\le\lambda_{\rm exp}$, we have $\ee^{-\lambda s}\ge\ee^{-(\mathrm M+1)\lambda_{\rm exp}}$ for $0\le s\le\mathrm M+1$, and hence
	\[
	\int_s^T\|c_N\|_{L^2(\gamma^-)}^2\le C\,\ee^{(\mathrm M+1)\lambda_{\rm exp}}\varsigma_N^2\,\ee^{-\lambda s},
	\qquad 0\le s\le T,
	\]
	and $\int_s^T\|c_N\|_{L^2(\gamma^-)}^2=0$ for $s\ge\mathrm M+1$; that is,
	$\sup_{0\le t\le T}\ee^{\lambda t}\|c_N\|^2_{L^2(\gamma^-_{(t,T)})}\le C\varsigma_N^2$.
	
	Now \eqref{eq:dynamic-global-exponential-stability-finite-horizon} applies in both cases, and, with
	$\|f_{\dy,0}\|_{L^2_{x,v}}\le\mathsf I$ by \eqref{eq:dyn-data-size}, $\mathsf G=\rho^2\le1$ by $\rho\le1$, and $K_c\le C\varsigma_N^2$, it gives
	\[
	\|u(t)\|_{L^2_{x,v}}^2+\cD_u(t,T)
	\le C\bigl(\mathsf I^2+\delta^2(1+\mathsf G)+K_c\bigr)\ee^{-\lambda t}
	\le C\bigl(\mathsf I^2+\delta^2+\varsigma_N^2\bigr)\ee^{-\lambda t},
	\]
	which is \eqref{eq:dyn-base-input-B-mass-full}.
	
	It remains to prove \eqref{eq:dyn-decay-B-closed}. By \eqref{eq:dyn-theta-residual-bound} and $\lambda\le\lambda_w$,
	\[
	\sup_{0\le t\le T}\ee^{\lambda t}\|b_w\|^2_{L^2(\gamma^-_{(t,T)})}\le C\delta^2 .
	\]
	Since $u$ satisfies \eqref{eq:dyn-decay-diffuse-law} with the datum $b$ of the corollary, we square \eqref{eq:dyn-decay-graph} on each $(t,T)$, multiply by $\ee^{\lambda t}$ and take the supremum over $0\le t\le T$ (see \eqref{eq:dyn-coordinate-B}). This gives
	\[
	\cB_T[u]\le C\sup_{0\le t\le T}\ee^{\lambda t}\Bigl(\|u(t)\|_{L^2_{x,v}}^2
	+\cD_u(t,T)+\|b\|^2_{L^2(\gamma^-_{(t,T)})}\Bigr),
	\]
	and the boundary term satisfies
	\[
	\sup_{0\le t\le T}\ee^{\lambda t}\|b\|^2_{L^2(\gamma^-_{(t,T)})}\le C(\delta^2+\varsigma_N^2)
	\]
	by the bound on $b_w$ above, $\sup_{0\le t\le T}\ee^{\lambda t}\|c_N\|^2_{L^2(\gamma^-_{(t,T)})}\le C\varsigma_N^2$ and $b=b_w+c_N$.  Combining this with
	\eqref{eq:dyn-base-input-B-mass-full}, we obtain \eqref{eq:dyn-decay-B-closed}.
\end{proof}

\begin{lemma}
	\label{lem:dyn-base-decay-to-mixed}
	Let $T\ge1$, let $\mathsf K\ge0$, and let $z$ satisfy
	\begin{equation}
		\label{eq:dyn-base-decay-hypothesis}
		\|z(t)\|_{L^2_{x,v}}^2+\cD_{z}(t,T)\le\mathsf K^2\ee^{-\lambda t},
		\qquad 0\le t\le T .
	\end{equation}
	Then for every $\tau\in(0,1]$
	\begin{equation}
		\label{eq:dyn-base-decay-mixed}
		\|z\|_{L_t^{2\tau}L_x^2H_\sigma^1}
		\le C\,\mathsf K,
	\end{equation}
	with $C$ depending only on the structural data and $\tau$ (in particular independent of $T$).
\end{lemma}

\begin{proof}
	Partition $(0,T)$ into $I_k=[k,k+1)\cap[0,T)$, $k\ge0$. On each interval,
	\[
	\int_{I_k}\|z\|_{L^2_xH^1_\sigma}^{2\tau}
	\le\Bigl(\int_{I_k}\|z\|^2_{L^2_xH^1_\sigma}\Bigr)^{\tau}
	{}\le\cD_{z}(k,T)^{\tau}\le\mathsf K^{2\tau}\ee^{-\lambda\tau k} ,
	\]
	by H\"older's inequality,
	\eqref{eq:dyn-decay-density} and \eqref{eq:dyn-base-decay-hypothesis}.  Summing the geometric series gives
	$\|z\|_{L_t^{2\tau}L_x^2H_\sigma^1}\le C_{\lambda,\tau}\mathsf K$,
	where $C_{\lambda,\tau}$ depends on $\lambda$ and $\tau$, but not on $T$.
	This proves \eqref{eq:dyn-base-decay-mixed}.
\end{proof}

The base coordinate is obtained by interpolating, at each time, between the weighted $L^{p_\sharp}_x$ norm controlled by the strong coordinate and the unweighted $L^2_x$ norm controlled by Lemma~\ref{lem:dyn-base-decay-to-mixed}.

\begin{prop}[Base estimate]
	\label{prop:dyn-base}
	Let $T\ge1$, let $\theta$, $q$ be those of \eqref{eq:dyn-base-choice}, and let $z$
	satisfy \eqref{eq:dyn-base-decay-hypothesis} for some $\mathsf K\ge0$.  Then, for
	$0\le\ell\le(1-\theta)(A^\dy_\flat-2)$,
	\begin{equation}
		\label{eq:dyn-base-product}
		\cS_T^{\ell}[z]\le C\,\|w^{A^\dy_\flat-2}z\|_{L_t^\infty L_x^{p_\sharp}H_\sigma^1}^{1-\theta}\,\mathsf K^{\theta},
	\end{equation}
	with $C$ depending only on the structural data and $\ell$ (in particular independent of $T$).
	If $u$ satisfies the hypotheses of Corollary~\ref{cor:dyn-decay-B}, then
	\begin{equation}
		\label{eq:dyn-base-unconditional-mass}
		\cS^{\ell}_T[u]
		\le C\bigl(\cA_T[u]+\delta+\mathsf I
		+\varsigma_N\bigr) ,
		\qquad 0\le\ell\le(1-\theta)(A^\dy_\flat-2),
	\end{equation}
	with $C$ depending only on the structural data and $\ell$ (in particular independent of $T$, $N$, $\varepsilon$ and $\rho$).
\end{prop}

\begin{proof}
	At each $t$, Proposition~\ref{prop:polynomial-weight-interpolation} with the pair
	$(p_0,p_1)=(p_\sharp,2)$, the ratio $\theta$ and the weights
	$(\mathrm L_0,\mathrm L_1)=(A^\dy_\flat-2,0)$ gives
	\[
	\|w^{\ell}z(t)\|_{L^{q}_xH^1_\sigma}
	\le C_\ell\|w^{A^\dy_\flat-2}z(t)\|_{L^{p_\sharp}_xH^1_\sigma}^{1-\theta}
	\|z(t)\|_{L^2_xH^1_\sigma}^{\theta},
	\qquad \ell\le(1-\theta)(A^\dy_\flat-2) ,
	\]
	where, by \eqref{eq:polynomial-weight-slots} and \eqref{eq:dyn-base-choice}, the exponent and the weight given by that proposition are
	\[
	\frac1{p_\theta}=\frac\theta2+\frac{1-\theta}{p_\sharp}=\frac1q,
	\qquad
	\mathrm L_\theta=\theta\cdot0+(1-\theta)(A^\dy_\flat-2)=(1-\theta)(A^\dy_\flat-2),
	\]
	and $\mathrm L_\theta$ is lowered to $\ell$ by
	\eqref{eq:real-weight-Hsigma-monotonicity}.
	
	\medskip
	Taking the $L^2$ norm in time and using H\"older's inequality and the identity
	\[
	\Bigl(\int_0^T\|z(t)\|_{L^2_xH^1_\sigma}^{2\theta}\dd t\Bigr)^{1/2}
	=\|z\|_{L_t^{2\theta}L_x^2H_\sigma^1}^{\theta},
	\]
	we obtain
	\begin{equation}
		\label{eq:dyn-base-interpolation}
		\|w^{\ell}z\|_{L_t^2L_x^{q}H_\sigma^1}
		\le C_\ell\,
		\|w^{A^\dy_\flat-2}z\|_{L_t^\infty L_x^{p_\sharp}H_\sigma^1}^{1-\theta}\,
		\|z\|_{L_t^{2\theta}L_x^2H_\sigma^1}^{\theta},
		\qquad
		0\le\ell\le(1-\theta)(A^\dy_\flat-2) .
	\end{equation}
	
	\medskip
	Because $0<\theta<1$, Lemma~\ref{lem:dyn-base-decay-to-mixed} at $\tau=\theta$ gives
	$\|z\|_{L_t^{2\theta}L_x^2H_\sigma^1}\le C\mathsf K$. Inserting this into
	\eqref{eq:dyn-base-interpolation}, we obtain \eqref{eq:dyn-base-product}.
	
	\medskip
	For the second assertion, let $u$ satisfy the hypotheses of Corollary~\ref{cor:dyn-decay-B}. Then
	\eqref{eq:dyn-base-input-B-mass-full} gives \eqref{eq:dyn-base-decay-hypothesis} at $z=u$ and
	$\mathsf K=C^{1/2}\bigl(\delta+\mathsf I+\varsigma_N\bigr)$, with $C$ of that corollary. Moreover,
	$\|w^{A^\dy_\flat-2}u\|_{L_t^\infty L_x^{p_\sharp}H_\sigma^1}\le C\cA_T[u]$ by
	\eqref{eq:dyn-coordinate-A-reading} at $\ell=A^\dy_\flat-2$. Hence \eqref{eq:dyn-base-product}, and then Young's inequality, give
	\[
	\cS^{\ell}_T[u]\le C\,\cA_T[u]^{1-\theta}\bigl(\delta+\mathsf I+\varsigma_N\bigr)^{\theta}
	\le C\bigl(\cA_T[u]+\delta+\mathsf I+\varsigma_N\bigr) ,
	\]
	which is \eqref{eq:dyn-base-unconditional-mass}. Its constant is independent of $T$, $N$, $\varepsilon$ and $\rho$, since so are the constants of Corollary~\ref{cor:dyn-decay-B} and of \eqref{eq:dyn-coordinate-A-reading}.
\end{proof}

The evolutionary hierarchy couples the unweighted and weighted rows, while \eqref{eq:dyn-strong-A-linear-defect} recovers the strong coordinate. Its closure is the analogue of Proposition~\ref{prop:uniform-hierarchy}.

\begin{prop}[Hierarchy estimate]
	\label{prop:dyn-hierarchy}
	Let $T\ge1$, $N\in\bN$ and $0<\varepsilon\le\varepsilon_\lambda$, let the datum satisfy
	\eqref{eq:dynamic-initial-zero-perturbation-mass}, let $g$ satisfy \eqref{eq:dyn-ball-jets},
	and set $\eta^\dy_g:=\widetilde{\mathfrak Q}^\dy_T(g)$.  There is
	$\eta_{\star,0}\in(0,\eta_{\star,1}]$, depending only on the structural data,
	such that if $\eta^\dy_g+\delta+\mathsf I\le\eta_{\star,0}$, then $f_\dy^{\varepsilon,N}$ of \eqref{eq:dyn-real-centered-unknown}, constructed from $g$, satisfies
	\begin{equation}
		\label{eq:dyn-strong-closed-A}
		\cA_T[f_\dy^{\varepsilon,N}]\le C\bigl(\delta+\mathsf I+\varsigma_N\bigr),
	\end{equation}
	\begin{equation}
		\label{eq:dyn-hierarchy-closed}
		\cS_T[f_\dy^{\varepsilon,N}]+\cN_{\mathrm M}(T)[f_\dy^{\varepsilon,N}]\le C\bigl(\delta+\mathsf I+\varsigma_N\bigr),
	\end{equation}
	with $C$ depending only on the structural data (in particular independent of $T$, $N$ and
	$\varepsilon$). Here $\eta_{\star,1}$ is that of Lemma~\ref{lem:dyn-strong-coefficient-class}.
\end{prop}

\begin{proof}
	We apply Lemma~\ref{lem:ladder-closure} with $M=\mathrm M$ and $\cR_m=\mathsf W_m(T)$. The reflected boundary trace is controlled first by the unweighted graph estimate; only then is the weighted row estimated. After the rows have been closed, \eqref{eq:dyn-strong-A-linear-defect} controls the strong coordinate.
	
	\medskip
	Write $f_\dy$ for the restriction of $f_\dy^{\varepsilon,N}$ to $(0,T)$.   Since $\eta^\dy_g+\mathsf I\le\eta_{\star,0}\le\eta_{\star,1}$, \eqref{eq:dyn-strong-coefficient-choice} holds.  Hence Proposition~\ref{prop:dyn-fixedN-realization}\,(\ref{it:dyn-real-prop-spatial}),
	(\ref{it:dyn-real-mass}), (\ref{it:dyn-real-strong}) and Lemma~\ref{lem:dyn-regbg-ball} apply at $r=\eta^\dy_g$ and give the following: the mixed jets of
	$f_\dy$ through complete order $\mathrm M$ lie in $\cH_{\kin,\sigma}$ and $\sum_{m=0}^{\mathrm M}(\mathsf U_m(T)+\mathsf W_m(T))<\infty$,
	$\mathfrak m_{f_\dy}\equiv0$ on $[0,T]$, $\cA_T[f_\dy]<\infty$ and
	\eqref{eq:dyn-strong-A-linear-defect} holds. Moreover, the regularized background $g^\varepsilon$ matches
	the initial jets and satisfies \eqref{eq:dyn-ball-eps-A} at
	$\rho:=2\eta^\dy_g+C(\delta+\mathsf I)$, and $\hat g^\varepsilon$ satisfies
	\eqref{eq:dyn-strong-coefficient-output-reg}.  By \eqref{eq:dyn-limits-fixed-system} at $j=0$ and \eqref{eq:dyn-limits-fixed-datum}, and since $\mathfrak t_{0,N}=0$ by \eqref{eq:dyn-theta-triangular-notation}, $f_\dy$ solves the problem \eqref{eq:dyn-perturbation-problem} with background $g^\varepsilon$ and the truncated wall law,
	\[
	\begin{cases}
		(\partial_t+v\cdot\nabla_x)f_\dy+L_{g^\varepsilon}f_\dy=N[g^\varepsilon_\dy,f_\st]
		&\text{in }(0,T)\times\cD\times\bR_v^3,\\
		(f_\dy)_{\gamma^-}=P_{w,N}(f_\dy)_{\gamma^+}+b_w+c_N
		&\text{on }\gamma_T^-,\\
		f_\dy(0)=\zeta_0=f_{\dy,0}
		&\text{in }\cD\times\bR_v^3 .
	\end{cases}
	\]
	
	Set $\kappa:=\eta^\dy_g+\delta+\mathsf I$. Since
	$\rho=2\eta^\dy_g+C(\delta+\mathsf I)$, we have
	$\rho\lesssim\kappa$, uniformly in $T,N,\varepsilon$. Choose $\eta_{\star,0}\le\eta_{\star,1}$ so small that
	$\rho\le1$ and $\rho+\delta\le\varepsilon_{\rm dec}$ whenever
	$\kappa\le\eta_{\star,0}$.  We may also assume
	$\eta_{\star,0}\le\varepsilon_{\rm dec}$.  Then
	Corollary~\ref{cor:dyn-decay-B} with $\cJ=J_N$, and hence
	\eqref{eq:dyn-base-unconditional-mass}, applies to $f_\dy$. Put
	\[
	\cQ:=\bigl(\delta+\mathsf I\bigr)+\varsigma_N
	+\kappa\,\cA_T[f_\dy] .
	\]
	The background factors of \eqref{eq:dyn-forcing-estimate}, evaluated at $g^\varepsilon_\dy$, satisfy 
	\begin{equation}
		\label{eq:dyn-hier-background-rungs}
		\|w^{-2}g^\varepsilon_\dy\|_{L_t^\infty L_x^\infty L_v^2}+\cN_{\mathrm M}(T)[g^\varepsilon_\dy]
		{}+\mathsf G_\natural(T)[g^\varepsilon_\dy]
		\le C\,\kappa .
	\end{equation}
	Indeed, by \eqref{eq:dyn-strong-fiber-mixed} and \eqref{eq:dyn-ball-eps-A},
	\[
	\|w^{-2}g^\varepsilon_\dy\|_{L_t^\infty L_x^\infty L_v^2}\le C\|w^{A^\dy_\flat}g^\varepsilon_\dy\|_{L^\infty}\le C\rho .
	\]
	Moreover, the bound $\cB_T[g^\varepsilon_\dy]\le\rho^2$ of \eqref{eq:dyn-ball-eps-A} and
	\eqref{260902137} give \eqref{eq:dyn-base-decay-hypothesis} for $g^\varepsilon_\dy$ at $\mathsf K=C\rho$. Hence
	Lemma~\ref{lem:dyn-base-decay-to-mixed} at $\tau=\vartheta_\natural$ (recall $\vartheta_\natural=2/P_\natural\in(0,1]$ from Lemma~\ref{lem:dyn-natural-base}) gives
	$\|g^\varepsilon_\dy\|_{L_t^{2\vartheta_\natural}L_x^2H_\sigma^1}\le C\rho$, and
	Lemma~\ref{lem:dyn-natural-base}, with the bound above on $\|w^{-2}g^\varepsilon_\dy\|_{L_t^\infty L_x^\infty L_v^2}$, gives
	$\mathsf G_\natural(T)[g^\varepsilon_\dy]\le C\rho$.  With the rows in \eqref{eq:dyn-ball-eps-A}, and since $\rho\le(2+C)\kappa$, these bounds give \eqref{eq:dyn-hier-background-rungs}.

	Let $0\le m\le\mathrm M$ and let $I$ be an ordered word of $\mathfrak A_\bullet$ with
	$|I|=m$.  If $m\ge1$, $\mathfrak D_If_\dy$ solves \eqref{eq:dyn-differentiated-equation} by Proposition~\ref{prop:dyn-differentiated-equation}, applied with $g^\varepsilon$ and $g^\varepsilon_\dy$ in place of $g$ and $g_\dy$, whose hypotheses are verified below for each order; at $I=\varnothing$, \eqref{eq:dyn-differentiated-equation} with $\mathfrak F_\varnothing[f_\dy]=N[g^\varepsilon_\dy,f_\st]$ is the first line of the problem above.  We view \eqref{eq:dyn-differentiated-equation} as the prescribed-inflow problem of
	\eqref{eq:weighted-graph-convention} for $\mathfrak D_If_\dy$, with source
	$\mathfrak F_I[f_\dy]$, initial datum $\mathfrak D_If_\dy(0)$ and incoming datum
	$(\mathfrak D_If_\dy)_{\gamma^-}$.  By Theorem~\ref{thm:weighted-steady}\,(1), applied to
	$\mathfrak D_If_\dy$ in the form \eqref{eq:weighted-graph-convention} at
	$L=\ell^\dy_m$, if
	\[
	\bigl\|w^{\ell^\dy_m}\mathfrak F_I[f_\dy]\bigr\|_{L^2(\cD_T;H^{-1}_\sigma)}+\bigl\|w^{\ell^\dy_m}\mathfrak D_If_\dy(0)\bigr\|_{L^2_{x,v}}+\bigl\|w^{\ell^\dy_m}(\mathfrak D_If_\dy)_{\gamma^-}\bigr\|_{L^2(\gamma_T^-)}<\infty,
	\]
	then $w^{\ell^\dy_m}\mathfrak D_If_\dy\in\cH_{\kin,\sigma}(\cD_T\times\bR_v^3)$ and
	\begin{align}
		\label{eq:dyn-weighted-graph}
		\begin{aligned}
			&\bigl\|w^{\ell^\dy_m}\mathfrak D_If_\dy\bigr\|_{\cH_{\kin,\sigma}(\cD_T\times\bR_v^3)}
			\\
			\le\,&C\Bigl(
			\bigl\|w^{\ell^\dy_m}\mathfrak F_I[f_\dy]\bigr\|_{L^2(\cD_T;H^{-1}_\sigma)}
			+\bigl\|w^{\ell^\dy_m}\mathfrak D_If_\dy(0)\bigr\|_{L^2_{x,v}}
			+\bigl\|w^{\ell^\dy_m}(\mathfrak D_If_\dy)_{\gamma^-}\bigr\|_{L^2(\gamma_T^-)}
			\Bigr),
		\end{aligned}
	\end{align}
	with $C$ depending only on the structural data (in particular independent of $T$, $N$ and $\varepsilon$).  The estimate \eqref{eq:dyn-hier-mass-graph} of Proposition~\ref{prop:dynamic-alpha-uniform-graph} at $\alpha=1$ is used below. The coefficient hypothesis \eqref{eq:weighted-dynamic-coefficient} of Theorem~\ref{thm:weighted-steady} and the hypothesis \eqref{eq:dynamic-alpha-graph-smallness} of Proposition~\ref{prop:dynamic-alpha-uniform-graph} hold: indeed, $\delta\le\delta_\star\le\delta_{\rm dyn}$, and \eqref{eq:dyn-strong-coefficient-output-reg} for
	$\hat g^\varepsilon$, whose restriction to $\cD_T$ is $g^\varepsilon$, gives
	\[
	\|\mu^{1/4}g^\varepsilon\|_{L^\infty(\cD_T\times\bR^3)}\le\min\{\varepsilon_{\rm dyn},\varepsilon_0(\cD)\} .
	\]
	
	For $0\le m\le\mathrm M$ and $|I|=m$, the reflected trace satisfies
	\begin{equation}
		\label{eq:dyn-hier-reflection}
		\bigl\|w^{\ell^\dy_m}P_{w,N}(\mathfrak D_If_\dy)_{\gamma^+}\bigr\|_{L^2(\gamma_T^-)}
		\le C\bigl\|(\mathfrak D_If_\dy)_{\gamma^+}\bigr\|_{L^2(\gamma_T^+)}\le C\,\mathsf U_m(T) ,
	\end{equation}
	where the first inequality follows from \eqref{eq:wall-J-weighted-operator-bound} at $\cJ=J_N$, applied at a.e.\ $t$ with $\Theta_w(t)$ in place of $\Theta_w$ and integrated in $t$, and the second from \eqref{260902137}.  The constant of \eqref{eq:wall-J-weighted-operator-bound} is uniform in $0<t<T$, since its proof uses $\Theta_w(t)$ only through \eqref{eq:wall-profile-chart-bound} for $d_\delta(t)$, which holds with one constant for $t>0$ by \eqref{eq:dyn-theta-amplitude}.
	
	\medskip
	At order $m=0$ the mass datum of
	\eqref{eq:dyn-hier-mass-graph} vanishes since $\mathfrak m_{f_\dy}\equiv0$.  Hence
	\eqref{eq:dyn-hier-mass-graph} at $f=f_\dy$ gives
	\[
	\mathsf U_0(T)=\|f_\dy\|_{\cH_{\kin,\sigma}(\cD_T\times\bR_v^3)}
	\le C\Bigl(\|f_{\dy,0}\|_{L^2_{x,v}}+\|N[g^\varepsilon_\dy,f_\st]\|_{L^2_tL^2_xH^{-1}_\sigma}+\|b_w+c_N\|_{L^2(\gamma_T^-)}\Bigr),
	\]
	with $C$ depending only on $\cD$; we now estimate the three terms on the right.  The initial datum
	$f_\dy(0)=\zeta_0=f_{\dy,0}$ satisfies $\|f_{\dy,0}\|_{L^2_{x,v}}\le\mathsf I$ by \eqref{eq:dyn-limits-fixed-system}
	and \eqref{eq:dyn-data-size}.  For the source, \eqref{eq:dyn-bilinear-pointwise} at $I_1=I_2=\varnothing$ and H\"older's inequality in $x$ at the exponent pair $(q_\flat,p_\flat)$ of \eqref{eq:dyn-own-order-endpoints} give, as in the proof of \eqref{eq:dyn-decay-forcing}, for $\ell\in\{0,\ell^\dy_0\}$,
	\begin{equation}
		\label{eq:dyn-hier-base-source}
		\bigl\|w^{\ell}N[g^\varepsilon_\dy,f_\st]\bigr\|_{L^2_tL^2_xH^{-1}_\sigma}
		\le C\bigl\|w^{-2}g^\varepsilon_\dy\bigr\|_{L^2_tL^{q_\flat}_xL^2_v}\bigl\|w^{\ell}f_\st\bigr\|_{L^{p_\flat}_xH^1_\sigma}
		\le C\delta\kappa\le C\delta ,
	\end{equation}
	since the first factor is at most $C\cN_0(T)[g^\varepsilon_\dy]\le C\kappa$ by \eqref{eq:dyn-slot-tops-rows} at $m=0$ and \eqref{eq:dyn-hier-background-rungs}, the second is at most $C\delta$ by \eqref{eq:dyn-stationary-slots} at $j=0$, as $\ell^\dy_0=\ell_\sharp$, and $\kappa\le\tfrac12$.  The incoming datum $b_w+c_N$ satisfies $\|b_w+c_N\|_{L^2(\gamma^-_T)}\le C\varsigma_N+C\delta$ by \eqref{eq:dyn-limits-corner-bound} and \eqref{eq:dyn-limits-wall-datum-size}.  Hence
	\begin{equation}
		\label{eq:dyn-hier-base-U0}
		\mathsf U_0(T)\le C\bigl(\delta+\mathsf I\bigr)+C\,\varsigma_N .
	\end{equation}
	For the weighted row, the initial datum is bounded by
	$C(\delta+\mathsf I)$ through \eqref{eq:dyn-data-jets},
	\eqref{eq:dyn-theta-background} and \eqref{260902137}.
	The source is bounded by $C\delta\kappa$ by
	\eqref{eq:dyn-hier-base-source}. The prescribed part of the incoming
	trace satisfies
	\[
	\|w^{\ell^\dy_0}(b_w+c_N)\|_{L^2(\gamma_T^-)}
	\le C(\delta+\varsigma_N)
	\]
	by \eqref{eq:dyn-limits-corner-bound} and
	\eqref{eq:dyn-limits-wall-datum-size}. For its reflected part,
	\eqref{eq:dyn-hier-reflection} and \eqref{eq:dyn-hier-base-U0}
	give the same bound, up to the term $C\mathsf I$.
	Applying \eqref{eq:dyn-weighted-graph} therefore yields
	\[
	\mathsf W_0(T)\le C(\delta+\mathsf I+\varsigma_N)
	\le C_0\cQ.
	\]
	This is \eqref{eq:ladder-step} for $m=0$.
	
	\medskip
	For $1\le m\le\mathrm M$ and $|I|=m$,  Proposition~\ref{prop:dyn-differentiated-equation} is applied with $g^\varepsilon$ and $g^\varepsilon_\dy$ in place of $g$ and $g_\dy$.  Its hypotheses hold
	for $m\ge1$:
	\begin{itemize}
		\item \eqref{eq:dyn-differentiated-chain}: the jets of $f_\dy$ lie in $\cH_{\kin,\sigma}$ and are its
		distributional derivatives, by Proposition~\ref{prop:dyn-fixedN-realization}\,(\ref{it:dyn-real-restriction});
		\item \eqref{eq:dyn-differentiated-lower-weight}: by $\cN_{m-1}(T)<\infty$ and \eqref{eq:real-weight-Hsigma-monotonicity}, since
		$\ell^\dy_k\ge\ell^\dy_m+A_\diamond\ge A_\diamond$ for $k<m$ by \eqref{eq:dyn-weight-choice};
		
		\item \eqref{eq:dyn-differentiated-first-slots}: by \eqref{eq:dyn-ball-eps-A} for
		$g^\varepsilon$, \eqref{eq:dyn-theta-background},
		$g^\varepsilon=g^\varepsilon_\dy+f_\st$ and $T<\infty$, since, for $|J|\le m$, \eqref{260902137} and $\mu^{1/4}\le1\le w$ give
		\[
		\begin{aligned}
			\|\mu^{1/4}\mathfrak D_Jg^\varepsilon_\dy\|_{L^2_tL^2_{x,v}}&\le C\|w^{\ell^\dy_{|J|}}\mathfrak D_Jg^\varepsilon_\dy\|_{\cH_{\kin,\sigma}(\cD_T\times\bR_v^3)},\\
			\|\mu^{1/4}\mathfrak D_Jf_\st\|_{L^2_tL^2_{x,v}}&\le C\,T^{1/2}\|w^{\ell_{|J|}}\mathfrak D_Jf_\st\|_{\cH_{\kin,\sigma}(\cD\times\bR_v^3)} .
		\end{aligned}
		\]
	\end{itemize}
	Write $\mathsf b$ for the weight in the graph estimate: $\mathsf b=0$ for $\mathsf U_m$ and $\mathsf b=\ell^\dy_m$ for $\mathsf W_m$.  We apply \eqref{eq:dyn-hier-mass-graph} at $\mathsf b=0$,
	respectively \eqref{eq:dyn-weighted-graph} at $\mathsf b=\ell^\dy_m$, to
	\eqref{eq:dyn-differentiated-equation}, and estimate the data in (i)--(iv); the mass, (iii) below, is a datum of
	\eqref{eq:dyn-hier-mass-graph} only.
	
	\emph{(i) The initial datum.}
	The datum is
	$$
	\mathfrak D_If_\dy(0)=\mathfrak D_If_\pt^{\varepsilon,N}(0)-\mathfrak D_If_\st. 
	$$
	We first identify $\mathfrak D_If_\pt^{\varepsilon,N}(0)$. Let
	$\mathfrak D_I=\mathfrak D_{I'}\partial_t^{\,j}$ be the normal ordering
	\eqref{eq:dyn-normal-ordering}.  By
	Proposition~\ref{prop:dyn-fixedN-realization}\,(\ref{it:dyn-real-prop-spatial}),
	$\mathfrak D_If_\pt^{\varepsilon,N}\in\cH_{\kin,\sigma}(\cD_{T,\ast}\times\bR^3)$, which
	embeds into $C([-1,T];L^2_{x,v})$ by Proposition~\ref{2607151218}, so its value at $t=0$ is
	its endpoint value on $(-1,0)$.  On $(-1,0)$, by
	Proposition~\ref{prop:dyn-fixedN-realization}\,(\ref{it:dyn-real-restriction}),
	\eqref{eq:dyn-real-prop-chain} and \eqref{eq:dyn-real-prop-initial},
	\[
	\mathfrak D_If_\pt^{\varepsilon,N}=\mathfrak D_{I'}\mathfrak D_{\mathrm t}^{\,j}f_\pt^{\varepsilon,N}
	=\mathfrak D_{I'}\partial_t^{\,j}\check{E}
	\] in the sense of distributions.  By
	\eqref{eq:dynamic-negative-time-Taylor-polynomial} and \eqref{eq:dynamic-negative-time-lift},
	the right side is a finite sum of products $a(t)\mathfrak D_{I'}z_i$ with scalar coefficients $a\in C^\infty([-1,0])$, and
	$\mathfrak D_{I'}z_i\in L^2_{x,v}$ by \eqref{eq:dyn-data-jets}.  Hence $\mathfrak D_If_\pt^{\varepsilon,N}(0)=\mathfrak D_{I'}z_j$, by
	\eqref{eq:dynamic-negative-time-jet-matching}.
	
	It remains to bound $\mathfrak D_{I'}z_j$ and $\mathfrak D_If_\st$. Since
	$0\le\mathsf b\le\ell^\dy_0$,
	$\|w^{\mathsf b}\mathfrak D_{I'}z_j\|_{\cH_{\kin,\sigma}(\cD\times\bR_v^3)}\le C\mathsf I$ by \eqref{eq:dyn-data-jets}. The term $\mathfrak D_If_\st$ vanishes
	when $I$ has a temporal letter, by \eqref{eq:dyn-fs-time-independent}, and otherwise satisfies
	\[
	\|w^{\mathsf b}\mathfrak D_If_\st\|_{\cH_{\kin,\sigma}(\cD\times\bR_v^3)}\le C\,W_m(f_\st)\le C\,\widetilde{\mathfrak Q}_{\mathrm N}(f_\st)\le C\delta
	\]
	by \eqref{eq:order-two-Wm}, \eqref{eq:unweighted-row-dominated}, \eqref{eq:reduced-Q3} and \eqref{eq:dyn-theta-background}, since $\mathsf b\le\ell^\dy_0=\ell_{\mathrm N}\le\ell_m$.  With \eqref{260902137} we obtain
	$\bigl\|w^{\mathsf b}\mathfrak D_If_\dy(0)\bigr\|_{L^2_{x,v}}\le C(\delta+\mathsf I)$.
	
	\emph{(ii) The source.}
	Apply Lemma~\ref{cor:dyn-forcing-estimate} with the regularized background
	at weight $\mathsf b\in\{0,\ell^\dy_m\}$. The lower-order terms
	of \eqref{eq:dyn-forcing-estimate} are bounded by $C\cN_{m-1}(T)$,
	since $\mathsf b+A_\diamond\le\ell^\dy_k$ for $k<m$.
	By \eqref{eq:dyn-hier-background-rungs}, its background factor is
	at most $C\kappa$. Moreover,
	\eqref{eq:dyn-coordinate-A-reading} and
	\eqref{eq:dyn-base-unconditional-mass} give
	\[
	\|w^{\mathsf b+A_\diamond}f_\dy\|_{L_t^\infty L_x^{p_\sharp}H^1_\sigma}
	+\cS_T^{\mathsf b+A_\diamond}[f_\dy]
	\le C(\cA_T[f_\dy]+\delta+\mathsf I+\varsigma_N).
	\]
	Here $\mathsf b+A_\diamond\le\ell^\dy_{m-1}\le
	\ell^\dy_0=A^\dy_\flat-3$ by \eqref{eq:dyn-weight-choice};
	the same weight lies in the range of
	\eqref{eq:dyn-base-unconditional-mass} by
	\eqref{eq:dyn-finite-base}.
	The stationary factors in the remaining terms are $O(\delta)$,
	and $\delta\le\kappa$. Writing
	$\cN_m=\cN_{m-1}+\mathsf W_m$ and using the definition of
	$\cQ$, we conclude that
	\[
	\|w^{\mathsf b}\mathfrak F_I[f_\dy]\|_{L^2(\cD_T;H^{-1}_\sigma)}
	\le C\bigl(\cQ+\cN_{m-1}(T)+\kappa\mathsf W_m(T)\bigr).
	\]
	All constants are uniform in $T,N,\varepsilon$.
	
	\emph{(iii) The mass.}
	The mass is the fourth datum $\|\mathfrak m_{\mathfrak D_If_\dy}\|_{L^2(0,T)}$ of \eqref{eq:dyn-hier-mass-graph}.  Suppose first that
	$j(I)<|I|$, so that $|I'|\ge1$ in the normal ordering of~(i).  We apply Proposition~\ref{prop:complete-derivative-mass-hierarchy} for $I'$ and
	$\partial_t^{\,j(I)}f_\dy(t)$ at a.e.\ $t\in(0,T)$; this is allowed, since the jets of $f_\dy$ form a jet chain of $\cH_{\kin,\sigma}(\cD_T\times\bR_v^3)$ (the chain \eqref{eq:dyn-differentiated-chain}, verified above). By \eqref{eq:complete-derivative-mass-sharp}, we have
	\[
	|\mathfrak m_{\mathfrak D_If_\dy}(t)|
	\le C\bigl\|\mathfrak D_{I'_\ast}\partial_t^{\,j(I)}f_\dy(t)\bigr\|_{L^2_{x,v}} ,
	\]
	where $I'_\ast$ is $I'$ with its first letter removed, as in that proposition. Since
	$\mathfrak D_{I'_\ast}\partial_t^{\,j(I)}=\mathfrak D_J$ with $|J|=|I|-1$ by
	\eqref{eq:dyn-normal-ordering}, integration in $t$ and \eqref{260829706} give
	\[
	\bigl\|\mathfrak m_{\mathfrak D_If_\dy}\bigr\|_{L^2(0,T)}
	\le C\sum_{|J|<|I|}\bigl\|\mathfrak D_Jf_\dy\bigr\|_{\cH_{\kin,\sigma}(\cD_T\times\bR_v^3)},
	\qquad 1\le|I|\le\mathrm M,\quad j(I)<|I|,
	\]
	with $C$ independent of $T$, $N$ and $\varepsilon$. Hence, by \eqref{eq:dyn-unweighted-row-dominated}, $\|\mathfrak m_{\mathfrak D_If_\dy}\|_{L^2(0,T)}\le C\cN_{m-1}(T)$.  If $I$ is purely temporal, then no complete derivative acts, so
	$\mathfrak m_{\mathfrak D_{\mathrm t}^{\,j}f_\dy}=\partial_t^{\,j}\mathfrak m_{f_\dy}$ in $\cD'(0,T)$,
	and the mass vanishes, since $\mathfrak m_{f_\dy}\equiv0$.
	
	\emph{(iv) The incoming datum.}
	It is given by \eqref{eq:dyn-limits-temporal-triangular}, with $\mathfrak t_{I,N}[f_\dy]$ and
	$\mathfrak r_{I,N}$ of \eqref{eq:dyn-limits-wall-mixed-terms}.  By \eqref{eq:dyn-limits-temporal-triangular-bound}, \eqref{260902137} and \eqref{eq:dyn-unweighted-row-dominated},
	its triangular part $\mathfrak t_{I,N}[f_\dy]$ satisfies
	\[
	\begin{aligned}
		\bigl\|\mathfrak t_{I,N}[f_\dy]\bigr\|_{L^2(\gamma^-_T)}
		&\le C\delta\sum_{|J|\le m,\ j(J)<j(I)}\bigl\|(\mathfrak D_Jf_\dy)_{\gamma^+}\bigr\|_{L^2(\gamma^+_T)}\\
		&\le C\delta\bigl(\mathsf U_0(T)+\dots+\mathsf U_m(T)\bigr)\le C\delta\bigl(\cN_{m-1}(T)+\mathsf W_m(T)\bigr) .
	\end{aligned}
	\]
	The term $\mathfrak r_{I,N}$ is the sum of the coloring $\mathfrak b_{I,N}[(f_\dy)_{\gamma^+}]$, which is present when $I$ has a spatial letter and involves only words of length
	less than $m$, and the derivative $\mathfrak D_{I'}^-\partial_t^j(\bar b_N-\mathfrak k_N)$ of the datum $\bar b_N-\mathfrak k_N$.  We have
	\[
	\bigl\|\mathfrak b_{I,N}[(f_\dy)_{\gamma^+}]\bigr\|_{L^2(\gamma^-_T)}\le C\cN_{m-1}(T),
	\qquad
	\bigl\|\mathfrak D_{I'}^-\partial_t^j(\bar b_N-\mathfrak k_N)\bigr\|_{L^2(\gamma^-_T)}\le C\varsigma_N+C\delta,
	\]
	the first by \eqref{eq:dyn-limits-wall-coloring-bound}, \eqref{260902137} and \eqref{eq:dyn-unweighted-row-dominated}, and the second
	by Lemma~\ref{lem:dyn-limits-wall}\,(\ref{it:dyn-limits-wall-form}) and~(\ref{it:dyn-limits-wall-size}).  These three parts are the boundary datum of
	\eqref{eq:dyn-hier-mass-graph} at $\mathsf b=0$.  At $\mathsf b=\ell^\dy_m$ the triangular
	and coloring bounds are \eqref{eq:dyn-limits-temporal-triangular-bound} and
	\eqref{eq:dyn-limits-wall-coloring-bound} at $\ell=\ell^\dy_m$, and the bound on the datum holds with the weight $w^{\ell^\dy_m}$ by \eqref{eq:dyn-limits-corner-bound} and \eqref{eq:dyn-limits-wall-datum-size}, since $\ell^\dy_m\le\ell^\dy_0$ by \eqref{eq:dyn-weight-choice}.  The
	incoming datum of \eqref{eq:dyn-weighted-graph} is the whole trace $(\mathfrak D_If_\dy)_{\gamma^-}$, so at
	$\mathsf b=\ell^\dy_m$ it also contains the reflected part, which is at most $C\,\mathsf U_m(T)$ by \eqref{eq:dyn-hier-reflection}.  This
	term carries no small factor, and it gives the term $C'_m\mathsf U_m(T)$ in the bound for $\mathsf W_m(T)$ in (v).
	
	\emph{(v) Assembly.}
	Inserting the bounds (i)--(iv) into \eqref{eq:dyn-hier-mass-graph} and \eqref{eq:dyn-weighted-graph}, and summing over the words $I$ with $|I|=m$, we obtain, with a
	constant $C'_m\ge1$ depending only on the structural data,
	\[
	\begin{aligned}
		\mathsf U_m(T)&\le C'_m\bigl(\cQ+\cN_{m-1}(T)\bigr)
		+C'_m\kappa\,\mathsf W_m(T) ,\\
		\mathsf W_m(T)&\le C'_m\bigl(\cQ+\cN_{m-1}(T)\bigr)+C'_m\,\mathsf U_m(T)
		+C'_m\kappa\,\mathsf W_m(T) .
	\end{aligned}
	\]
	Inserting the first line into the second, we obtain
	\[
	\mathsf W_m(T)\le C_m\bigl(\cQ+\cN_{m-1}(T)\bigr)
	+C_m\kappa\,\mathsf W_m(T) ,
	\qquad 1\le m\le\mathrm M ,
	\]
	for a constant $C_m$ depending only on the structural data. This is
	\eqref{eq:ladder-step}.
	
	\medskip
	It remains to close the hierarchy. Set $D:=\delta+\mathsf I+\varsigma_N$.  The row inequalities
	satisfy \eqref{eq:ladder-step} with $\theta=\kappa$ and with
	structural constants independent of $T,N,\varepsilon$ and $\kappa$.
	Since the fixed-cutoff realization already gives
	$\cN_{\mathrm M}(T)<\infty$, Lemma~\ref{lem:ladder-closure}
	applies after reducing the structural threshold for $\kappa$.
	Together with \eqref{eq:dyn-strong-A-linear-defect}, it gives
	\[
	\cN_{\mathrm M}(T)\le C(D+\kappa\cA_T[f_\dy]),
	\qquad
	\cA_T[f_\dy]\le C(D+\cN_{\mathrm M}(T)).
	\]
	The constants are independent of $T,N,\varepsilon$ and $\kappa$.
	Substitution yields
	$\cN_{\mathrm M}(T)\le C D+C\kappa\cN_{\mathrm M}(T)$.
	Choosing $\eta_{\star,0}$ sufficiently small, after all structural
	constants have been fixed, absorbs the last term.  Hence
	$\cN_{\mathrm M}(T)+\cA_T[f_\dy]\lesssim D$.
	Finally, \eqref{eq:dyn-base-unconditional-mass} at $\ell=0$
	yields $\cS_T[f_\dy]\lesssim D$ and completes the proof of
	\eqref{eq:dyn-strong-closed-A}--\eqref{eq:dyn-hierarchy-closed}.
\end{proof}

The coordinate bounds are preserved under strong graph convergence, as needed both for the fixed point and for removal of the regularizations.

\begin{lemma}
	\label{lem:dyn-limits-coordinate-closure}
	Let $T\ge1$, let $z_n\to z$ in $\cH_{\kin,\sigma}(\cD_T\times\bR_v^3)$, and assume
	$\sup_n\cN_{\mathrm M}(T)[z_n]<\infty$.  Then, for every $I\in\mathfrak A_\bullet^m$ with
	$m\le\mathrm M$, the jets $\mathfrak D_Iz$ and $w^{\ell^\dy_m}\mathfrak D_Iz$ lie in
	$\cH_{\kin,\sigma}(\cD_T\times\bR_v^3)$, and
	\begin{equation}
		\label{eq:dyn-limits-closure-jets}
		\begin{gathered}
			\mathfrak D_Iz_n\rightharpoonup\mathfrak D_Iz,
			\qquad
			w^{\ell^\dy_m}\mathfrak D_Iz_n\rightharpoonup w^{\ell^\dy_m}\mathfrak D_Iz
			\qquad\text{in }\cH_{\kin,\sigma}(\cD_T\times\bR_v^3),
			\\
			(\mathfrak D_Iz_n)(0)\rightharpoonup(\mathfrak D_Iz)(0)
			\qquad\text{in }L^2(\cD\times\bR_v^3).
		\end{gathered}
	\end{equation}
	Moreover
	\begin{equation}
		\label{eq:dyn-limits-closure}
		\cA_T[z]\le\liminf_{n\to\infty}\cA_T[z_n],
		\qquad
		\cS_T[z]\le\liminf_{n\to\infty}\cS_T[z_n],
	\end{equation}
	\begin{equation}
		\label{eq:dyn-limits-closure-rows}
		\mathsf W_m(T)[z]\le\liminf_{n\to\infty}\mathsf W_m(T)[z_n],
	\end{equation}
	for $0\le m\le\mathrm M$.
\end{lemma}

\begin{proof}
	By the convention after \eqref{eq:reduced-Q3}, recalled after
	Definition~\ref{def:dyn-coordinates}, the bound on
	$\cN_{\mathrm M}(T)[z_n]$ and \eqref{eq:dyn-unweighted-row-dominated} show that the jets $\mathfrak D_Iz_n$ and $w^{\ell^\dy_m}\mathfrak D_Iz_n$, $I\in\mathfrak A_\bullet^m$, $m\le\mathrm M$, lie in $\cH_{\kin,\sigma}(\cD_T\times\bR_v^3)$ with
	\[
	\|\mathfrak D_Iz_n\|_{\cH_{\kin,\sigma}(\cD_T\times\bR_v^3)}+\|w^{\ell^\dy_m}\mathfrak D_Iz_n\|_{\cH_{\kin,\sigma}(\cD_T\times\bR_v^3)}\le C\,\cN_{\mathrm M}(T)[z_n] ,
	\]
	and the right side is bounded uniformly in $n$.
	Write
	$\mathfrak D_I=\mathfrak D_{I'}\partial_t^{\,j}$ in the normal ordering \eqref{eq:dyn-normal-ordering} and
	$\mathfrak D_I^\dagger=(-\partial_t)^{\,j}\mathfrak D_{I'}^\dagger$, with the ordered adjoint
	\eqref{eq:ordered-complete-formal-adjoint} of $I'$.  For $\psi\in C_c^\infty(\cD_T\times\bR^3)$,
	\[
	\langle w^{\ell^\dy_m}\mathfrak D_Iz_n,\psi\rangle
	=\langle z_n,\mathfrak D_I^\dagger(w^{\ell^\dy_m}\psi)\rangle
	\longrightarrow\langle z,\mathfrak D_I^\dagger(w^{\ell^\dy_m}\psi)\rangle
	=\langle w^{\ell^\dy_m}\mathfrak D_Iz,\psi\rangle ,
	\]
	and likewise without the weight.  So the two jets converge to $\mathfrak D_Iz$ and
	$w^{\ell^\dy_m}\mathfrak D_Iz$ in $\cD'(\cD_T\times\bR_v^3)$, and
	Lemma~\ref{lem:dynamic-weak-graph-bulk-identification}\,(\ref{it:closure-evolutionary}) gives
	$\mathfrak D_Iz,\,w^{\ell^\dy_m}\mathfrak D_Iz\in\cH_{\kin,\sigma}(\cD_T\times\bR_v^3)$ and \eqref{eq:dyn-limits-closure-jets}.
	
	For the coordinates, Lemma~\ref{lem:strong-coordinate-lsc}\,(\ref{it:strong-coordinate-lsc-dyn}) with
	$A=A^\dy_\flat$, $p=p_\sharp$ and $\alpha=\alpha_{\rm bs}$ gives \eqref{eq:dyn-limits-closure} for $\cA_T$, since the three quantities of \eqref{eq:strong-coordinate-lsc-dyn} are then the three summands of
	\eqref{eq:dyn-coordinate-A}.  The coordinate $\cS_T$ of
	\eqref{eq:dyn-base-row} is the norm of the reflexive space $L_t^2L_x^{q}H_\sigma^1$.
	We may assume that $\liminf_n\cS_T[z_n]<\infty$, and pass to a subsequence realizing it. By weak compactness, a further subsequence has a weak limit, identified with $z$ by testing
	against $C_c^\infty$. The weak lower semicontinuity of the norm then gives
	\eqref{eq:dyn-limits-closure} for $\cS_T$.  Finally, each summand of $\mathsf W_m$ in
	\eqref{eq:dyn-hierarchy-rows} is the norm in $\cH_{\kin,\sigma}(\cD_T\times\bR_v^3)$ of a
	weighted jet that converges weakly there by
	\eqref{eq:dyn-limits-closure-jets}, and weak lower semicontinuity gives
	\eqref{eq:dyn-limits-closure-rows}.
\end{proof}

Fix the structural threshold $\eta_\star$ of
Theorem~\ref{thm:dyn-physical-solution} sufficiently small, and then
reduce $\delta_\star$ so that
\begin{equation}
	\label{eq:dyn-threshold-choice}
	\eta_\star+\delta_\star\le\eta_{\star,0}.
\end{equation}
Thus $\eta_g^\dy+\delta+\mathsf I\le\eta_{\star,0}$ whenever
$\widetilde{\mathfrak Q}^\dy_T(g)\le\eta_\star$ and
$\delta+\mathsf I\le\delta_\star$.
This gives the hypotheses of Proposition~\ref{prop:dyn-hierarchy}
and Lemma~\ref{lem:dyn-strong-coefficient-class}.  It also gives
\eqref{eq:dyn-decay-B-smallness} for an unregularized background with
$\rho=\widetilde{\mathfrak Q}^\dy_T(g)$; for regularized backgrounds,
the required bound follows from $\rho\lesssim\kappa$ in the proof of
Proposition~\ref{prop:dyn-hierarchy}.  We retain $\delta_\star\le1$.

\subsection{The fixed point and the evolutionary solution}\label{ss:dyn-fixed-point}

The estimate \eqref{eq:dyn-physical-coordinates} gives an invariant reduced class, while the kinetic graph norm provides the contraction metric.
Throughout this
subsection the datum satisfies \eqref{eq:dynamic-initial-zero-perturbation-mass} and $\delta>0$.

\begin{defn}[The evolutionary class]
	\label{def:dyn-ball}
	For $1\le T\le\infty$ and $r\ge0$, let
	\[
	\begin{aligned}
		\widetilde{\mathbb B}^\dy_T(r):=\bigl\{g\text{ on }\cD_T\times\bR_v^3:{}&\ g\text{ satisfies \eqref{eq:dyn-ball-jets}},\ \widetilde{\mathfrak Q}^\dy_T(g)\le r,\\
		&\ (g_\dy)_{\gamma^-}=P_w(t)(g_\dy)_{\gamma^+}+b_w(t)\ \text{on }\gamma_T^-,\\
		&\ \mathfrak m_{g_\dy}(t)=0\ \text{for every finite }0\le t\le T\bigr\}.
	\end{aligned}
	\]
	We equip it with the metric
	\begin{equation}
		\label{eq:dyn-fp-metric-infty}
		d^\dy_\cH(g^1,g^2):=\|g^1-g^2\|_{\cH_{\kin,\sigma}(\cD_T\times\bR_v^3)} ,
	\end{equation}
	with $d^\dy_\cH(g^1,g^2)<\infty$ by the row $\mathsf W_0$ of \eqref{eq:dyn-ball-rows} and \eqref{eq:dyn-unweighted-row-dominated}; at $T=\infty$ it is the
	supremum over finite $T$.
\end{defn}

\begin{lemma}
	\label{lem:dyn-fp-closed}
	For $1\le T<\infty$ and $r\ge0$, the class $\widetilde{\mathbb B}^\dy_T(r)$ is complete under
	the metric \eqref{eq:dyn-fp-metric-infty}.
\end{lemma}

\begin{proof}
	Let $y^{(n)}\in\widetilde{\mathbb B}^\dy_T(r)$ be a Cauchy sequence.  Since $\cH_{\kin,\sigma}(\cD_T\times\bR_v^3)$ is complete, $y^{(n)}\to y$ in
	$\cH_{\kin,\sigma}(\cD_T\times\bR_v^3)$ for some $y$.  We pass each defining condition to the limit.
	\begin{itemize}
		\item By \eqref{eq:dyn-ball-rows}, $\cN_{\mathrm M}(T)[y^{(n)}-f_\st]\le r$ for every $n$, so
		Lemma~\ref{lem:dyn-limits-coordinate-closure} applies at $z_n=y^{(n)}-f_\st$ and
		$z=y-f_\st$.  By its bounds \eqref{eq:dyn-limits-closure}--\eqref{eq:dyn-limits-closure-rows},
		\[
		\cS_T[y-f_\st]\le\liminf_{n\to\infty}\cS_T[y^{(n)}-f_\st],
		\qquad
		\mathsf W_m(T)[y-f_\st]\le\liminf_{n\to\infty}\mathsf W_m(T)[y^{(n)}-f_\st],
		\qquad 0\le m\le\mathrm M .
		\]
		Its endpoint convergence \eqref{eq:dyn-limits-closure-jets} at the temporal words gives
		\eqref{eq:dyn-ball-jets} for $y$, since $(\partial_t^{\,j}z_n)(0)=\zeta_j$ for every $n$.
		\item As in the proof of Lemma~\ref{lem:dyn-limits-coordinate-closure},
		Lemma~\ref{lem:strong-coordinate-lsc}\,(\ref{it:strong-coordinate-lsc-dyn}) gives
		$\cA_T[y]\le\liminf_n\cA_T[y^{(n)}]$.
		\item The decay coordinate $\cB_T$ is continuous in the metric \eqref{eq:dyn-fp-metric-infty}.  Indeed, $\cB_T^{1/2}$ satisfies the triangle inequality, and, by \eqref{260902137} and \eqref{eq:dyn-coordinate-B}, $\cB_T[z]\le C\ee^{\lambda T}\|z\|_{\cH_{\kin,\sigma}(\cD_T\times\bR_v^3)}^2$.  Hence
		\[
		\cB_T[y-f_\st]=\lim_n\cB_T[y^{(n)}-f_\st] .
		\]
		\item By \eqref{260902137},
		\[
		\|(y^{(n)}-y)_\gamma\|_{L^2(\gamma_T)}\le C\|y^{(n)}-y\|_{\cH_{\kin,\sigma}(\cD_T\times\bR_v^3)} ,
		\]
		so the traces converge in
		$L^2(\gamma_T^\pm)$, and the wall law passes to the limit by
		\eqref{eq:wall-J-operator-bound}.
		\item The zero-mass constraint passes to the limit, since, by \eqref{260902137},
		\[
		\sup_{0\le t\le T}\|(y^{(n)}-y)(t)\|_{L^2(\cD\times\bR^3)}\le C\|y^{(n)}-y\|_{\cH_{\kin,\sigma}(\cD_T\times\bR_v^3)} .
		\]
	\end{itemize}
	The first three bounds cover the four terms of \eqref{eq:dyn-reduced-Q} and give
	\[
	\widetilde{\mathfrak Q}^\dy_T(y)\le\liminf_n\widetilde{\mathfrak Q}^\dy_T(y^{(n)})\le r .
	\]
	Hence
	$y\in\widetilde{\mathbb B}^\dy_T(r)$, and $y^{(n)}\to y$ in the metric \eqref{eq:dyn-fp-metric-infty}.  This completes the proof.
\end{proof}

\begin{defn}[The centered self-map]
	\label{def:dyn-selfmap}
	For $1\le T\le\infty$, whenever $\delta+\mathsf I\le\delta_\star$ and $g$ satisfies
	\eqref{eq:dyn-ball-jets} and $\widetilde{\mathfrak Q}^\dy_T(g)\le\eta_\star$, let $f_\dy[g]$ be the solution
	given by Theorem~\ref{thm:dyn-physical-solution}\,(\ref{it:dyn-phys-existence}), or by (\ref{it:dyn-phys-halfline}) of that theorem if $T=\infty$, and set
	\[
	\Phi_T(g):=f_\st+f_\dy[g] .
	\]
\end{defn}

For $\delta+\mathsf I\le\delta_\star$ and $r\le\eta_\star$, a function $g\in\widetilde{\mathbb B}^\dy_T(r)$ solves \eqref{eq:dyn-physical-problem} if and only if $\Phi_T(g)=g$: the centered equation \eqref{eq:dyn-perturbation-problem} at background $g$ is equivalent to the nonlinear equation at $g=f_\pt$. Its graph solution is unique by Theorem~\ref{thm:dyn-physical-solution}; the required graph membership follows from the row $\mathsf W_0$ and \eqref{eq:dyn-unweighted-row-dominated}, while the initial jets agree by definition of the class. Uniqueness on each finite interval also gives
\begin{equation}
	\label{eq:dyn-fp-selfmap-consistency}
	\Phi_\infty(g)\big|_{[0,T]}=\Phi_T\bigl(g|_{[0,T]}\bigr),
	\qquad T\ge1 .
\end{equation}

\begin{lemma}
	\label{prop:dyn-fp-selfmap}
	Fix a structural constant $R\ge1$ large enough that both
	\eqref{eq:dyn-physical-coordinates} and the seed estimate of
	Lemma~\ref{lem:dyn-fp-seed} hold with $R$ in place of their constants.
	Decrease $\delta_\star$, if necessary, so that
	\begin{equation}
		\label{eq:dyn-fp-radii}
		R\delta_\star\le\eta_\star.
	\end{equation}
	Then, for $\delta+\mathsf I\le\delta_\star$ and $1\le T\le\infty$,
	\[
	\Phi_T:\widetilde{\mathbb B}^\dy_T\bigl(R(\delta+\mathsf I)\bigr)\longrightarrow
	\widetilde{\mathbb B}^\dy_T\bigl(R(\delta+\mathsf I)\bigr).
	\]
	Moreover, with the seed $g^{(0)}$ of \eqref{eq:dyn-fp-seed},
	$\Phi_T(g^{(0)})\in\widetilde{\mathbb B}^\dy_T\bigl(R(\delta+\mathsf I)\bigr)$.
\end{lemma}

\begin{proof}
	The seed has the prescribed initial jets by Lemma~\ref{lem:dyn-fp-seed}, as do elements of the class. In either case,
	\[
	\widetilde{\mathfrak Q}^\dy_T(g)
	\le R(\delta+\mathsf I)\le\eta_\star.
	\]
	Theorem~\ref{thm:dyn-physical-solution} gives
	$\widetilde{\mathfrak Q}^\dy_T(\Phi_T(g))\le R(\delta+\mathsf I)$.
	Its existence and mass assertions give the initial jets, wall law and
	zero-mass constraint for $\Phi_T(g)-f_\st$.
	This proves the assertion for finite $T$.
	For $T=\infty$, use the half-line assertion of the same theorem and
	\eqref{eq:dyn-fp-selfmap-consistency}; the bounds are uniform in $T$.
\end{proof}

\begin{prop}[Contraction]
	\label{prop:dyn-fp-contraction}
	After decreasing $\delta_\star$, if necessary, whenever $\delta+\mathsf I\le\delta_\star$,
	$1\le T<\infty$ and $g^1,g^2\in\widetilde{\mathbb B}^\dy_T(R(\delta+\mathsf I))$, we have
	\begin{equation}
		\label{eq:dyn-fp-contraction}
		d^\dy_\cH\bigl(\Phi_T(g^1),\Phi_T(g^2)\bigr)
		\le\tfrac12\,d^\dy_\cH\bigl(g^1,g^2\bigr) .
	\end{equation}
\end{prop}

\begin{proof}
	Put $d=g^1-g^2$, $f^i_\dy=f_\dy[g^i]$ and $\omega=f^1_\dy-f^2_\dy$. Subtracting \eqref{eq:dyn-perturbation-problem} at the two backgrounds gives
	\[
	\begin{gathered}
		(\partial_t+v\cdot\nabla_x)\omega+L_{g^1}\omega
		=N[d,f_\st+f^2_\dy]=N[d,\Phi_T(g^2)],\\
		\omega(0)=0,\qquad
		\omega_{\gamma^-}=P_w(t)\omega_{\gamma^+}.
	\end{gathered}
	\]
	Collision invariance and the homogeneous diffuse law give $\mathfrak m_\omega\equiv0$.  Since $g^1(0)=z_0$ and
	\[
	\widetilde{\mathfrak Q}^\dy_T(g^1)+\mathsf I
	\le R(\delta+\mathsf I)+\mathsf I
	\le\eta_\star+\delta_\star\le\eta_{\star,0},
	\]
	Lemma~\ref{lem:dyn-strong-coefficient-class} and the choice
	$\delta_\star\le\delta_{\rm dyn}$ verify
	\eqref{eq:dynamic-alpha-graph-smallness} for $g^1$.
	The homogeneous graph estimate \eqref{eq:dyn-hier-mass-graph}
	at $\alpha=1$ and \eqref{eq:dyn-decay-forcing} therefore give
	\[
	\begin{aligned}
		\|\omega\|_{\cH_{\kin,\sigma}}
		&\lesssim \|N[d,\Phi_T(g^2)]\|_{L^2_tL^2_xH_\sigma^{-1}}\\
		&\lesssim \|\Phi_T(g^2)\|_{L_t^\infty L_x^{p_\flat}H_\sigma^1}
		\|d\|_{\cH_{\kin,\sigma}}\\
		&\lesssim \cA_T[\Phi_T(g^2)]\|d\|_{\cH_{\kin,\sigma}}
		\lesssim R(\delta+\mathsf I)\|d\|_{\cH_{\kin,\sigma}}.
	\end{aligned}
	\]
	The third inequality uses \eqref{eq:dyn-coordinate-A-reading},
	$p_\sharp=2p_\flat$ and the boundedness of $\cD$;
	the last uses Lemma~\ref{prop:dyn-fp-selfmap}.
	The constants are independent of $T$. Fixing $R$ before $\delta_\star$ and then decreasing $\delta_\star$ makes the coefficient at most $1/2$, proving
	\eqref{eq:dyn-fp-contraction}.
\end{proof}

\begin{prop}[Global solution]
	\label{prop:dyn-fp-global}
	For $\delta+\mathsf I\le\delta_\star$, with $\delta_\star$ as decreased in
	Proposition~\ref{prop:dyn-fp-contraction}, problem \eqref{eq:dyn-physical-problem} on
	$[0,\infty)$ has a solution $f_\pt$ in
	$\widetilde{\mathbb B}^\dy_\infty(R(\delta+\mathsf I))$, and
	\begin{equation}
		\label{eq:dyn-fp-global-bounds}
		\widetilde{\mathfrak Q}^\dy_\infty(f_\pt)\le R\bigl(\delta+\mathsf I\bigr)\le\eta_\star ,
	\end{equation}
	with $R$ of Lemma~\ref{prop:dyn-fp-selfmap} and $\eta_\star$ of
	Theorem~\ref{thm:dyn-physical-solution}.
\end{prop}

\begin{proof}
	To verify completeness of the half-line class, let $(g_n)$ be Cauchy in $d^\dy_\cH$ on $[0,\infty)$.  For every finite $T\ge1$, its restrictions
	are Cauchy in $\widetilde{\mathbb B}^\dy_T(R(\delta+\mathsf I))$.
	By Lemma~\ref{lem:dyn-fp-closed} they converge to some $g_T$.
	Uniqueness of the limit gives $g_{T'}|_{[0,T]}=g_T$ whenever $T'\ge T$;
	these functions therefore define $g$ on $[0,\infty)$.  Since
	$\widetilde{\mathfrak Q}^\dy_T(g)=\widetilde{\mathfrak Q}^\dy_T(g_T)
	\le R(\delta+\mathsf I)$ for every $T$, the monotone half-line
	definition before Definition~\ref{def:dyn-coordinates} gives
	$g\in\widetilde{\mathbb B}^\dy_\infty(R(\delta+\mathsf I))$.
	To obtain convergence in the half-line metric, fix $\epsilon>0$ and
	choose $n_0$ so that $d^\dy_\cH(g_n,g_m)\le\epsilon$ for
	$m,n\ge n_0$.  Passing $m\to\infty$ on each finite interval gives
	$\|g_n-g\|_{\cH_{\kin,\sigma}(\cD_T\times\bR_v^3)}
	\le\epsilon$ for every $T$ and every $n\ge n_0$.
	Taking the supremum in $T$ proves $d^\dy_\cH(g_n,g)\to0$.
	Thus the half-line class is complete.  By Lemma~\ref{prop:dyn-fp-selfmap}, it contains $\Phi_\infty(g^{(0)})$ and $\Phi_\infty$
	maps it into itself.  Let $g^1,g^2\in\widetilde{\mathbb B}^\dy_\infty(R(\delta+\mathsf I))$ and $1\le T<\infty$.  Since $\widetilde{\mathfrak Q}^\dy_T\le\widetilde{\mathfrak Q}^\dy_\infty$, the restrictions $g^i|_{[0,T]}$ lie in $\widetilde{\mathbb B}^\dy_T(R(\delta+\mathsf I))$, so \eqref{eq:dyn-fp-selfmap-consistency} and \eqref{eq:dyn-fp-contraction} give
	\[
	d^\dy_\cH\bigl(\Phi_\infty(g^1)|_{[0,T]},\Phi_\infty(g^2)|_{[0,T]}\bigr)
	\le\tfrac12\,d^\dy_\cH\bigl(g^1|_{[0,T]},g^2|_{[0,T]}\bigr) .
	\]
	Since $d^\dy_\cH$ at $T=\infty$ is the supremum over finite $T$, $\Phi_\infty$ is a $\tfrac12$-contraction.  By the Banach fixed point theorem there is
	$f_\pt\in\widetilde{\mathbb B}^\dy_\infty(R(\delta+\mathsf I))$ with $\Phi_\infty(f_\pt)=f_\pt$.  Since $R(\delta+\mathsf I)\le R\delta_\star\le\eta_\star$ by \eqref{eq:dyn-fp-radii}, $f_\pt$
	solves \eqref{eq:dyn-physical-problem} by the paragraph after
	Definition~\ref{def:dyn-selfmap}, and \eqref{eq:dyn-fp-global-bounds} follows from $f_\pt\in\widetilde{\mathbb B}^\dy_\infty(R(\delta+\mathsf I))$.
\end{proof}

The fixed point, decay estimate, boundary regularity and renormalization yield the global theorem. It suffices to state the weighted bounds for integer target weights: a real weight $\ell\ge0$ is covered by $\ell^\dy_\sharp=\lceil\ell\rceil$.

\begin{thm}[Global nonlinear evolution and exponential asymptotic stability]
	\label{thm:dynamic-global}
	Assume the wall hypotheses
	(\ref{it:dyn-theta-amplitude})--(\ref{it:dyn-theta-relaxation}).
	Fix a nonnegative integer target weight \(\ell^\dy_\sharp\) as in
	\eqref{eq:dyn-weight-choice}, and let \(f_\st\) be the stationary
	solution of Theorem~\ref{thm:stationary-diffuse} at
	\(\ell_\sharp=\ell_0^\dy\) for the limiting wall datum
	\(h_w^\infty\).  Let \(f_0\) have finite data size \(\mathsf I\)
	as in \eqref{eq:dyn-data-size}.  Suppose that its generated temporal
	jets \(z_0,\ldots,z_{\mathrm M}\), defined by
	\eqref{eq:dynamic-initial-p-u}--\eqref{eq:dynamic-generated-q0},
	satisfy \eqref{eq:dynamic-initial-u-compatibility} for every
	\(0\le j\le\mathrm M\), and that
	\eqref{eq:dynamic-initial-zero-perturbation-mass} holds.
	There are a structural threshold \(\delta_\star>0\) and a constant
	\(C\), depending only on the structural data, such that the following assertions hold whenever
	\(0<\delta\) and \(\delta+\mathsf I\le\delta_\star\).
	\begin{enumerate}
		\item\label{it:dyn-global-existence}
		(\emph{Existence.})  Problem \eqref{eq:dyn-physical-problem} has a global solution $f_\pt$ on $[0,\infty)$ such that
		\[
		\int_{\cD\times\bR^3_v}f_\pt(t)\sqrt\mu\dd v\dd x=0,\qquad t\ge0,
		\]
		and, for $f_\dy:=f_\pt-f_\st$,
		\begin{equation}
			\label{eq:dyn-global-regularity}
			\begin{aligned}
				&\|w^{A^\dy_\flat}f_\pt\|_{L^\infty}
				+\bigl[w^{A^\dy_\flat}f_\pt\bigr]_{C^{\alpha_{\rm bs}}_{\kin}}
				+\|w^{A^\dy_\flat}\nabla_vf_\pt\|_{L^\infty_tL^{p_\sharp}_xL^2_v}
				+\sup_{t\ge0}\ee^{\lambda t/2}\|f_\dy\|_{\cH_{\kin,\sigma}(\cD_{(t,\infty)}\times\bR_v^3)}\\
				&\qquad+\|f_\dy\|_{L_t^2L_x^{q}H_\sigma^1}
				+\sum_{m=0}^{\mathrm M}\sum_{I\in\mathfrak A_\bullet^m}
				\|w^{\ell^\dy_m}\mathfrak D_If_\dy\|_{\cH_{\kin,\sigma}}
				\le C\bigl(\delta+\mathsf I\bigr),
			\end{aligned}
		\end{equation}
		where $\lambda=\min\{\lambda_{\rm exp},\lambda_w\}>0$ is the rate of \eqref{eq:dyn-decay-rate-fixed-point} and the other norms are taken over $(0,\infty)\times\cD\times\bR_v^3$.  The outgoing moment of $f_\dy$ satisfies
		\begin{equation}
			\label{eq:dyn-global-boundary-moment-H}
			\|\cA_0(f_\dy)_{\gamma^+}\|_{H^{\mathrm M}((0,\infty)\times\partial\cD)}\le C\bigl(\delta+\mathsf I\bigr).
		\end{equation}
		Consequently, for every $0<\beta<\frac12$, there is $C_\beta>0$ such that
		\begin{equation}
			\label{eq:dyn-global-boundary-moment-holder}
			\|\cA_0(f_\dy)_{\gamma^+}\|_{C^{1,\beta}([0,\infty)\times\partial\cD)}
			\le C_\beta\bigl(\delta+\mathsf I\bigr).
		\end{equation}
		Moreover,
		\begin{equation}
			\label{eq:dyn-global-decay}
			\begin{aligned}
				&\|f_\dy(t)\|_{L^2_{x,v}}^2
				+\int_t^\infty\Bigl(\|f_\dy(s)\|_{L^2_xH^1_\sigma}^2
				+\|(\partial_s+v\cdot\nabla_x)f_\dy(s)\|_{L^2_xH^{-1}_\sigma}^2\Bigr)\dd s\\
				&\qquad+\|(f_\dy)_{\gamma^+}\|_{L^2((t,\infty)\times\gamma^+)}^2
				\le C\bigl(\|f_0-f_\st\|_{L^2_{x,v}}^2+\delta^2\bigr)\ee^{-\lambda t},
				\qquad t\ge0 .
			\end{aligned}
		\end{equation}
		If $\vartheta$ is independent of $t$, then, with $\lambda_w=\lambda_{\rm exp}$, \eqref{eq:dyn-global-decay} holds at the rate $\lambda=\lambda_{\rm exp}$ and without the term $\delta^2$.

		\item\label{it:dyn-global-uniqueness}
		(\emph{Uniqueness.})  The solution $f_\pt$ in (\ref{it:dyn-global-existence}) is unique among solutions $g_\pt$ of \eqref{eq:dyn-physical-problem}
		on $[0,\infty)$ such that
		$g_\pt\in\cH_{\kin,\sigma}(\cD_T\times\bR_v^3)$ for every $T\ge1$ and
		\begin{equation}
			\label{eq:dyn-global-uniqueness-class}
			\|\mu^{1/4}g_\pt\|_{L^\infty((0,\infty)\times\cD\times\bR^3_v)}\le\varepsilon_{\rm dyn},
		\end{equation}
		with $\varepsilon_{\rm dyn}$ of Proposition~\ref{prop:dynamic-alpha-uniform-graph}.
		\item\label{it:dyn-global-nonnegativity}
		(\emph{Physical nonnegativity.})  If the initial physical distribution
		\[
		F_0:=\mu+\sqrt\mu\,f_0
		\]
		is nonnegative almost everywhere, then
		$F(t):=\mu+\sqrt\mu\,f_\pt(t)\ge0$ almost everywhere in
		$\cD\times\bR_v^3$ for every $t\ge0$.
	\end{enumerate}
\end{thm}

By Definition~\ref{def:dyn-coordinates}, \eqref{eq:dyn-global-regularity} is equivalently
\begin{equation}
	\label{eq:dyn-global-Q}
	\widetilde{\mathfrak Q}^\dy_\infty(f_\pt)\le C(\delta+\mathsf I).
\end{equation}

\begin{proof}[\textbf{Proof of parts (\ref{it:dyn-global-existence})--(\ref{it:dyn-global-uniqueness}) of Theorem~\ref{thm:dynamic-global}}]
	Proposition~\ref{prop:dyn-fp-global} gives $f_\pt\in\widetilde{\mathbb B}^\dy_\infty(R(\delta+\mathsf I))$ solving \eqref{eq:dyn-physical-problem}, where $R$ is chosen in Lemma~\ref{prop:dyn-fp-selfmap}.  By Definition~\ref{def:dyn-ball}, $\mathfrak m_{f_\dy}(t)=0$ for every $t\ge0$; since $f_\st$ has zero mass by Theorem~\ref{thm:stationary-diffuse}\,(\ref{it:stat-existence}), this is the mass condition in (\ref{it:dyn-global-existence}).  Since $R$ depends only on the structural data, \eqref{eq:dyn-fp-global-bounds} gives \eqref{eq:dyn-global-Q}, which is \eqref{eq:dyn-global-regularity}.
	
	By \eqref{eq:dynamic-outgoing-moment-H3} of Lemma~\ref{lem:dynamic-boundary-moments} at $m=\mathrm M$, \eqref{eq:dyn-unweighted-row-dominated} and \eqref{eq:dyn-global-regularity}, $\|\cA_0(f_\dy)_{\gamma^+}\|_{H^{\mathrm M}((0,T)\times\partial\cD)}\le C(\delta+\mathsf I)$ with $C$ independent of $T$; letting $T\to\infty$ gives \eqref{eq:dyn-global-boundary-moment-H}.
	Since $\mathrm M\ge3$ and $(s,s+2)\times\partial\cD$ is three-dimensional, $H^{\mathrm M}((s,s+2)\times\partial\cD)\subset C^{1,\beta}([s,s+2]\times\partial\cD)$ for $0<\beta<\frac12$, with a constant independent of $s\ge0$.  Two points of $[0,\infty)\times\partial\cD$ at time distance less than $1$ lie in one such $[s,s+2]\times\partial\cD$, and for two points at time distance at least $1$ the H\"older quotient of the first derivatives is at most twice their supremum; hence \eqref{eq:dyn-global-boundary-moment-H} gives \eqref{eq:dyn-global-boundary-moment-holder}.
	
	For the decay estimate, apply
	Proposition~\ref{prop:dynamic-global-exponential-stability} at $T=\infty$ and at
	$\lambda$ of \eqref{eq:dyn-decay-rate-fixed-point}, with the background $g=f_\pt$,
	so that $g_\dy=f_\dy$.  Its hypotheses hold.
	\begin{itemize}
		\item By \eqref{eq:dyn-fp-global-bounds} and Lemma~\ref{lem:dyn-fp-seed},
		\[
		\|w^{A^\dy_\flat}f_\dy\|_{L^\infty}
		\le \cA_\infty[f_\pt]+\cA_\infty[f_\st]
		\le(R+C)(\delta+\mathsf I).
		\]
		Since $R$ is fixed before the final choice of $\delta_\star$, we may
		decrease $\delta_\star$ so that this coefficient norm is at most
		$\varepsilon_\star$.  This gives the first hypothesis of
		\eqref{eq:dyn-decay-background-hypotheses}.  Also
		\[
		\mathsf G=\cB_\infty[f_\dy]
		\le\widetilde{\mathfrak Q}^\dy_\infty(f_\pt)^2
		\le\eta_\star^2\le1
		\]
		by \eqref{eq:dyn-fp-global-bounds}.
		\item \eqref{eq:dyn-decay-solution-hypotheses} holds: since $\Phi_\infty(f_\pt)=f_\pt$,
		$f_\dy=f_\dy[f_\pt]$ solves
		\eqref{eq:dyn-perturbation-problem} at the background $f_\pt$ on every $(0,T)$; $f_\dy\in\cH_{\kin,\sigma}(\cD_\infty\times\bR_v^3)$ by the row $\mathsf W_0$ and \eqref{eq:dyn-unweighted-row-dominated}, and
		$\mathfrak m_{f_\dy}(0)=0$ by \eqref{eq:dynamic-initial-zero-perturbation-mass}.
	\end{itemize}
	Since $\delta^2\mathsf G\le\delta^2$,
	\eqref{eq:dynamic-global-exponential-stability-finite-horizon} at $T=\infty$, with no extra
	datum ($c=0$ and $K_c=0$ in \eqref{eq:dyn-decay-general-wall}--\eqref{eq:dyn-decay-extra-datum}), gives \eqref{eq:dyn-global-decay}; here $f_{\dy,0}=f_0-f_\st$, and, by \eqref{eq:dyn-decay-density},
	\[
	\begin{aligned}
		&\|f_\dy(t)\|_{L^2_{x,v}}^2+\int_t^\infty\Bigl(\|f_\dy(s)\|_{L^2_xH^1_\sigma}^2+\|(\partial_s+v\cdot\nabla_x)f_\dy(s)\|_{L^2_xH^{-1}_\sigma}^2\Bigr)\dd s+\|(f_\dy)_{\gamma^+}\|_{L^2((t,\infty)\times\gamma^+)}^2\\
		&\qquad\le\|f_\dy(t)\|_{L^2_{x,v}}^2+\cD_{f_\dy}(t,\infty)+\|(f_\dy)_{\gamma^+}\|_{L^2(\gamma^+_{(t,\infty)})}^2 .
	\end{aligned}
	\]
	
	If $\vartheta$ is independent of $t$, then \eqref{eq:dyn-theta-relaxation} holds with $\vartheta_\infty=\vartheta$ for every $\lambda_w>0$, so $b_w\equiv0$ by \eqref{eq:dyn-theta-residual}.
	By Proposition~\ref{prop:dynamic-global-exponential-stability}, $\lambda_{\rm exp}$ does not depend on $\lambda_w$; we take $\lambda_w=\lambda_{\rm exp}$ in the structural data from the beginning, so that $\lambda=\lambda_{\rm exp}$ and the bounds above, including $\mathsf G\le1$, hold at this $\lambda$.  
	Let $S$ be the supremum over $t\ge0$ of $\ee^{\lambda t}$ times the left side of \eqref{eq:dyn-global-decay}.  As above, the last assertion of Proposition~\ref{prop:dynamic-global-exponential-stability} gives $S\le C(\|f_0-f_\st\|_{L^2_{x,v}}^2+\delta^2\mathsf G)$; in particular $S<\infty$.  
	By (\ref{it:diffuse-P0-projection}) of Lemma~\ref{lem:diffuse-boundary-algebra}, $\|f_\dy(t)\|_{L^2_{x,v}}^2+\cD_{f_\dy}(t,\infty)\le\ee^{-\lambda t}S$ for $t\ge0$, so \eqref{eq:dyn-decay-graph} at $z=f_\dy$ on $(t,T)$ with $b=0$, and then $T\to\infty$, give $\mathsf G\le C_*S$, with $C_*$ depending only on the structural data.  
	Decreasing $\delta_\star$ once more, we absorb the term
	$C\delta^2\mathsf G$ using $\mathsf G\le C_*S$.
	Consequently $S\lesssim\|f_0-f_\st\|_{L^2_{x,v}}^2$,
	which proves \eqref{eq:dyn-global-decay} at $\lambda=\lambda_{\rm exp}$
	without the $\delta^2$ term and completes (\ref{it:dyn-global-existence}).
	
	To prove (\ref{it:dyn-global-uniqueness}), we first check that $f_\pt$ belongs to the stated uniqueness class.  By \eqref{eq:dyn-fp-global-bounds}, $\mathsf W_0(T)<\infty$ for every $T\ge1$; hence $f_\dy\in\cH_{\kin,\sigma}(\cD_T\times\bR_v^3)$ for every $T\ge1$ by \eqref{eq:dyn-unweighted-row-dominated}.  Since the $t$-independent $f_\st$ lies in $\cH_{\kin,\sigma}(\cD\times\bR_v^3)$ by Proposition~\ref{prop:stationary-finite-weight-regularity}, it follows that $f_\pt=f_\st+f_\dy\in\cH_{\kin,\sigma}(\cD_T\times\bR_v^3)$ for every $T\ge1$.  Moreover, $f_\pt$ satisfies \eqref{eq:dyn-global-uniqueness-class}.
	Indeed, $f_\pt(0)=z_0$ by \eqref{eq:dyn-ball-jets} at $j=0$, and, for every finite $T\ge1$, $\widetilde{\mathfrak Q}^\dy_T(f_\pt)\le\widetilde{\mathfrak Q}^\dy_\infty(f_\pt)\le\eta_\star$ by \eqref{eq:dyn-fp-global-bounds}, so \eqref{eq:dyn-strong-coefficient-choice} holds by the paragraph after \eqref{eq:dyn-threshold-choice}.  Hence
	Lemma~\ref{lem:dyn-strong-coefficient-class} gives
	$\|\mu^{1/4}f_\pt\|_{L^\infty(\cD_T\times\bR^3)}\le\tfrac12\varepsilon_{\rm dyn}$, and, since $T$ is arbitrary, \eqref{eq:dyn-global-uniqueness-class} follows.  Let $g_\pt$ be a
	solution as in (\ref{it:dyn-global-uniqueness}), let $T\ge1$, and put
	\[
	\omega:=g_\pt-f_\pt\in\cH_{\kin,\sigma}(\cD_T\times\bR_v^3) .
	\]  Subtracting
	the two problems and using \[
	N[g_\pt,g_\pt]-N[f_\pt,f_\pt]=N[g_\pt,\omega]+N[\omega,f_\pt] ,
	\] we obtain
	\[
	(\partial_t+v\cdot\nabla_x)\omega+L\omega-N[g_\pt,\omega]=N[\omega,f_\pt],
	\qquad
	\omega(0)=0,
	\qquad
	\omega_{\gamma^-}=P_w(t)\omega_{\gamma^+} .
	\]
	By \eqref{eq:dyn-global-uniqueness-class} and $\delta\le\delta_\star\le\delta_{\rm dyn}$, Proposition~\ref{prop:dynamic-alpha-uniform-graph} at $\alpha=1$
	applies with the background $g_\pt$.  The source has zero mass, since $\langle N[a,b],\sqrt\mu\rangle_v=0$, and the initial value and the wall datum vanish, so $\mathfrak m_\omega\equiv0$, and \eqref{eq:dyn-hier-mass-graph} gives
	\[
	\|\omega\|_{\cH_{\kin,\sigma}}\le C\|N[\omega,f_\pt]\|_{L^2_tL^2_xH^{-1}_\sigma} .
	\]  As in the
	proof of Proposition~\ref{prop:dyn-fp-contraction}, by \eqref{eq:dyn-decay-forcing}, by \eqref{eq:dyn-coordinate-A-reading} at $\ell=0$ with H\"older's inequality in $x$, and by $\cA_T[f_\pt]\le\widetilde{\mathfrak Q}^\dy_\infty(f_\pt)\le R(\delta+\mathsf I)$ (see \eqref{eq:dyn-fp-global-bounds}), we have
	\[
	\begin{aligned}
		C\|N[\omega,f_\pt]\|_{L^2_tL^2_xH^{-1}_\sigma}
		&\le C\|f_\pt\|_{L^\infty_tL^{p_\flat}_xH^1_\sigma}\|\omega\|_{\cH_{\kin,\sigma}}
		\le C\cA_T[f_\pt]\,\|\omega\|_{\cH_{\kin,\sigma}}\\
		&\le CR(\delta+\mathsf I)\|\omega\|_{\cH_{\kin,\sigma}}\le\tfrac12\|\omega\|_{\cH_{\kin,\sigma}} ,
	\end{aligned}
	\]
	with $C$ as in that proof; the last inequality follows from the choice of $\delta_\star$ there.  Hence, as $\|\omega\|_{\cH_{\kin,\sigma}}<\infty$,
	$\omega=0$ on $(0,T)$ for every $T$, and $g_\pt=f_\pt$.  This proves (\ref{it:dyn-global-uniqueness}).
\end{proof}

\begin{proof}[\textbf{Proof of part (\ref{it:dyn-global-nonnegativity}) of Theorem~\ref{thm:dynamic-global}}]
	Fix $T>0$ and set
	\[
	F:=\mu+\sqrt\mu\,f_\pt,
	\qquad
	a^{ij}:=\phi^{ij}*F,
	\qquad
	b^i:=\partial_{v_j}a^{ij},
	\qquad
	J^i:=a^{ij}\partial_{v_j}F-b^iF .
	\]
	Since $f_\pt$ solves \eqref{eq:dyn-physical-problem}, $F$ solves the physical Landau equation, which is the first identity below; the second is \eqref{eq:double-divergence-sigma}. Both hold in distributions on $(0,T)\times\cD\times\bR_v^3$,
	\[
	(\partial_t+v\cdot\nabla_x)F=\partial_{v_i}J^i,
	\qquad
	\partial_{v_i}b^i=-8\pi F.
	\]
	Moreover, $\|\mu^{1/4}f_\pt\|_{L^\infty((0,\infty)\times\cD\times\bR^3)}\le\tfrac12\varepsilon_{\mathrm{ell}}$ by \eqref{eq:dyn-strong-coefficient-output} for $g=f_\pt$ at every finite $T\ge1$, which applies since $f_\pt(0)=z_0$ and, by the proof of (\ref{it:dyn-global-uniqueness}), \eqref{eq:dyn-strong-coefficient-choice} holds for $g=f_\pt$.  Hence
	Lemma~\ref{lem:perturbed-landau-ellipticity}, applied at each $(t,x)$, gives
	\[
	\tfrac12\sigma^{ij}\xi_i\xi_j\le a^{ij}\xi_i\xi_j\le2\sigma^{ij}\xi_i\xi_j .
	\]
	As in the proof of (\ref{it:dyn-global-uniqueness}), $f_\pt\in\cH_{\kin,\sigma}((0,T)\times\cD\times\bR_v^3)$.
	By the multiplier property in Lemma~\ref{lem:compact-velocity-support},
	$\sqrt\mu\,f_\pt\in\cH_{\kin,\sigma}$; hence
	$F\in\cH_{\kin,\sigma}((0,T)\times\cD\times\bR_v^3)$.
	As in the proof of (\ref{it:stat-nonnegativity}) of Theorem~\ref{thm:stationary-diffuse}, Lemma~\ref{lem:local-embedding-H_sigma} and the Gaussian factor give, for one fixed $m>1/2$,
	\[
	w^mF,\quad w^m\nabla_vF,\quad w^mJ,\quad w^mbF
	\in L^2((0,T)\times\cD\times\bR_v^3).
	\]
	
	\medskip
	\noindent\emph{Negative-part balance.}
	For $\varepsilon>0$ let $\beta_\varepsilon$, $H_\varepsilon$ and the velocity cutoffs $\chi_R,\zeta_R$ be those of Step~1 of the proof of Theorem~\ref{thm:stationary-diffuse}\,(\ref{it:stat-nonnegativity}).
	Lemma~\ref{lem:compact-velocity-support} gives
	$\zeta_RF\in\cH_{\kin}$, and
	\[
	(\partial_t+v\cdot\nabla_x)(\zeta_RF)
	=-J\cdot\nabla_v\zeta_R+\partial_{v_i}(\zeta_RJ^i).
	\]
	On $\supp\chi_R\subset B_{2R}$, $g^0:=-J\cdot\nabla_v\zeta_R$ and $g^i:=\zeta_RJ^i$ satisfy $g^0=0$ and
	$g^i=J^i$, since $\zeta_R\equiv1$ on $B_{3R}$.  Proposition~\ref{prop:renormalized-identity}, applied on $(0,t)$ to $\zeta_RF$ with
	these $g^0$, $g^i$ and test $\chi_R$, gives, for $0\le t\le T$,
	\begin{align*}
		&\int\beta_\varepsilon(F(t))\chi_R
		-\int\beta_\varepsilon(F_0)\chi_R
		+\int_0^t\int_\gamma
		\beta_\varepsilon(F_\gamma)\chi_R\sgn(n_x\cdot v)\,\dd\gamma\dd s\\
		&\qquad
		=-\int_0^t\int
		\beta_\varepsilon''(F)\,\partial_{v_i}F\,J^i\chi_R
		-\int_0^t\int
		\beta_\varepsilon'(F)J\cdot\nabla_v\chi_R .
	\end{align*}
	By $J=a\nabla_vF-bF$, the integrand $-\beta_\varepsilon''(F)\,\partial_{v_i}F\,J^i$ of the first term on the right equals $-\beta_\varepsilon''(F)a^{ij}\partial_{v_i}F\partial_{v_j}F+b\cdot\nabla_vH_\varepsilon(F)$, the sum of a diffusion contribution and a drift term.  For the drift term, we use the distributional identity
	$\nabla_v\cdot b=-8\pi F$ exactly as in the stationary proof, and obtain
	\[
	\int b\cdot\nabla_vH_\varepsilon(F)\chi_R
	=8\pi\int FH_\varepsilon(F)\chi_R
	-\int bH_\varepsilon(F)\cdot\nabla_v\chi_R .
	\]
	We have $-\beta_\varepsilon''(F)a^{ij}\partial_{v_i}F\partial_{v_j}F\le0$ and $8\pi FH_\varepsilon(F)\le0$, since $\beta_\varepsilon''\ge0$, $(a^{ij})$ is positive definite by $a^{ij}\xi_i\xi_j\ge\tfrac12\sigma^{ij}\xi_i\xi_j$, and $rH_\varepsilon(r)\le0$.  Hence
	\begin{align*}
		&\int\beta_\varepsilon(F(t))\chi_R
		+\int_0^t\int_\gamma
		\beta_\varepsilon(F_\gamma)\chi_R\sgn(n_x\cdot v)\,\dd\gamma\dd s\\
		&\qquad\le
		\int\beta_\varepsilon(F_0)\chi_R+\mathcal E_{R,\varepsilon}(t),
	\end{align*}
	where $\mathcal E_{R,\varepsilon}(t)$ is the sum of the two terms containing $\nabla_v\chi_R$,
	\[
	\mathcal E_{R,\varepsilon}(t)=-\int_0^t\int\beta_\varepsilon'(F)J\cdot\nabla_v\chi_R-\int_0^t\int bH_\varepsilon(F)\cdot\nabla_v\chi_R .
	\]
	Since $|\beta_\varepsilon'|\le1$, $|H_\varepsilon(F)|\le|F|$ and $|\nabla_v\chi_R|\le C/R$, the Cauchy--Schwarz inequality and $w^mJ,\,w^mbF\in L^2((0,T)\times\cD\times\bR_v^3)$ give
	\[
	|\mathcal E_{R,\varepsilon}(t)|
	\le\frac{C}{R}
	\int_0^t\int_{\cD\times\{R<|v|<2R\}}
	\bigl(|J|+|bF|\bigr)
	\lesssim R^{1/2-m}\longrightarrow0
	\quad\text{uniformly in }\varepsilon .
	\]  We also have
	\[
	F_-\le\sqrt\mu\,|f_\pt|,
	\qquad
	(F_\gamma)_-\le\sqrt\mu\,|(f_\pt)_\gamma|,
	\]
	and these are $L^1$ on every finite time interval by Proposition~\ref{2607151218} and Cauchy--Schwarz.
	Dominated convergence in $R$ and then in $\varepsilon$ therefore yields
	\[
	\int F_-(t)
	+\int_0^t\left[
	\int_{\gamma^+}(F_{\gamma^+})_-
	-\int_{\gamma^-}(F_{\gamma^-})_-
	\right]\dd s
	\le\int(F_0)_-.
	\]
	
	\medskip
	\noindent\emph{Diffuse boundary flux.}
	For a.e. $(s,x)$ write
	\[
	n(s,x):=\int_{n_x\cdot v>0}(F_{\gamma^+})_-n_x\cdot v\,\dd v,
	\qquad
	p(s,x):=\int_{n_x\cdot v>0}(F_{\gamma^+})_+n_x\cdot v\,\dd v.
	\]
	The physical diffuse law, obtained from the wall law of \eqref{eq:dyn-physical-problem} as in Step~2 of the proof of Theorem~\ref{thm:stationary-diffuse}\,(\ref{it:stat-nonnegativity}), gives
	$F_{\gamma^-}=M_{\Theta_w(s,x)}(p-n)$.  Hence, by the unit incoming flux of $M_{\Theta_w(s,x)}$,
	\[
	\int_{n_x\cdot v<0}(F_{\gamma^-})_-|n_x\cdot v|\,\dd v
	=(n-p)_+\le n.
	\]
	Hence the boundary contribution is nonnegative: for a.e. $s$, $\int_{\gamma^+}(F_{\gamma^+})_--\int_{\gamma^-}(F_{\gamma^-})_-\ge0$, since it equals $\int_{\partial\cD}\bigl(n-(n-p)_+\bigr)\dd S_x$.  If $F_0\ge0$, then $\int(F_0)_-=0$, so
	$\int F_-(t)=0$.  Finally, $t\mapsto\int F_-(t)$ is continuous.  Indeed,
	\[
	|F_-(t)-F_-(s)|\le\sqrt\mu\,|f_\pt(t)-f_\pt(s)|,
	\]
	and $f_\pt\in C_tL^2_{x,v}$ by the $C_tL^2_{x,v}$ trace of $\cH_{\kin,\sigma}$, so the Cauchy--Schwarz inequality gives
	\[
	\|F_-(t)-F_-(s)\|_{L^1_{x,v}}\le C\|f_\pt(t)-f_\pt(s)\|_{L^2_{x,v}};
	\]
	thus $t\mapsto F_-(t)$ is continuous in $L^1_{x,v}$.  Hence, since $T>0$ is arbitrary, $F(t)\ge0$ almost everywhere for every $t\ge0$.
\end{proof}

\section{Regularized Evolutionary Construction and Passage to the Limit}
\label{sec:dyn-realization-section}

We prove Proposition~\ref{prop:dyn-fixedN-realization} and parts
(\ref{it:dyn-phys-existence})--(\ref{it:dyn-phys-coords}) of Theorem~\ref{thm:dyn-physical-solution}, using the complete differential calculus of \S\ref{sec:diffuse-complete-equations} and the stationary estimates of \S\ref{sec:diffuse-common}. The evolutionary construction also requires past mollification and realization of the temporal derivatives.

\subsection{The past mollification and the zero-initial temporal chain}\label{ss:dyn-past-mollification}

We establish the graph-space properties of past mollification and construct the temporal derivatives of zero-initial solutions by difference quotients. Both will be used in Proposition~\ref{prop:dyn-fixedN-realization}.

\subsubsection{The past mollification}\label{ss:dyn-past-moll-kernel}

Let $0<T<\infty$.  The kernel $\varrho$, the mollifier $\varrho_\varepsilon$ and its
derivatives $\varrho_\varepsilon^{(j)}$ are those of \eqref{eq:dyn-real-mollifier} and
\eqref{eq:dyn-real-mollified-objects}.

Since $z=0$ on $(-1,-3/4]$, extend it by zero to $t\le-1$. The past convolutions are then defined on all of $I_T^-$.

\begin{prop}[The past mollification]
	\label{prop:dyn-fixedN-mollification}
	Let $0<\varepsilon<1/4$, let $z\in\cH_{\kin,\sigma}(\cD_{T,\ast}\times\bR_v^3)$ satisfy
	$z=0$ on $(-1,-3/4]$, and put
	\[
	z^\varepsilon:=\varrho_\varepsilon*z,
	\qquad
	z_j^\varepsilon:=\varrho_\varepsilon^{(j)}*z
	\qquad\text{on }I_T^-.
	\]
	\begin{enumerate}
		\item\label{it:dyn-real-moll-graph}
		$z_j^\varepsilon\in\cH_{\kin,\sigma}(\cD_{T,\ast}\times\bR_v^3)$ for every $j\ge0$, with
		\begin{equation}
			\label{eq:dyn-real-moll-graph-bound}
			\|z^\varepsilon\|_{\cH_{\kin,\sigma}(\cD_{T,\ast}\times\bR_v^3)}
			\le\|z\|_{\cH_{\kin,\sigma}(\cD_{T,\ast}\times\bR_v^3)},
			\qquad
			\|z_j^\varepsilon\|_{\cH_{\kin,\sigma}(\cD_{T,\ast}\times\bR_v^3)}
			\le\|\varrho_\varepsilon^{(j)}\|_{L^1}
			\|z\|_{\cH_{\kin,\sigma}(\cD_{T,\ast}\times\bR_v^3)},
		\end{equation}
		and
		\[
		(\partial_t+v\cdot\nabla_x)z_j^\varepsilon
		=\varrho_\varepsilon^{(j)}*\bigl((\partial_t+v\cdot\nabla_x)z\bigr),
		\qquad
		(z_j^\varepsilon)_{\gamma^\pm}=\varrho_\varepsilon^{(j)}*z_{\gamma^\pm}.
		\]
		\item\label{it:dyn-real-moll-chain}
		$\partial_tz_j^\varepsilon=z_{j+1}^\varepsilon$ in
		$\cD'(I_T^-\times\cD\times\bR^3)$ for every $j\ge0$; in particular
		$z_j^\varepsilon=\partial_t^jz^\varepsilon$ and $z^\varepsilon$ is smooth
		in $t$.
	\end{enumerate}
\end{prop}

\begin{proof}
	For (\ref{it:dyn-real-moll-graph}), choose $\chi\in C^\infty(\bR)$ with $\chi=0$ on
	$(-\infty,-7/8]$ and $\chi=1$ on $[-3/4,\infty)$, and smooth
	$z_n\to z$ in $\cH_{\kin,\sigma}(\cD_{T,\ast}\times\bR_v^3)$.  Multiplication by $\chi$ preserves the graph space. The only additional transport term is $\chi'z$, and by \eqref{eq:Hsigma-dual-weighted}
	and the tail bound \eqref{260829706},
	\[
	\|\chi'z\|_{L^2(I_T^-;L_x^2H^{-1}_\sigma)}\le C\|\la v\ra^{1/2}z\|_{L^2(I_T^-\times\cD\times\bR_v^3)}
	\le C\|z\|_{\cH_{\kin,\sigma}(\cD_{T,\ast}\times\bR_v^3)} .
	\]  Since $\chi z=z$, smooth approximants $z_n$ may be replaced by $\chi z_n$ and thus taken to vanish for
	$t\le-7/8$. For $0<s<\varepsilon$, consider the zero-extended delays $z(\cdot-s)$.  Each $z_n(\cdot-s)$ is smooth on $\overline{I_T^-}$ because $z_n$ vanishes for $t\le-7/8$. Its transport term and traces are delayed likewise, and the
	terms of \eqref{innerprod_sigma} at $z_n(\cdot-s)$ are those of $z_n$ over $(-1,T-s)$, the
	endpoint term being $0$.  So $\|z_n(\cdot-s)\|\le\|z_n\|$ in
	$\cH_{\kin,\sigma}(\cD_{T,\ast}\times\bR_v^3)$, and $s\mapsto z_n(\cdot-s)$ is continuous
	into this space.  By density and the bound uniform in $s$, $z_n(\cdot-s)\to z(\cdot-s)$ in this space, the map $s\mapsto z(\cdot-s)$ is continuous into this space, and
	\[
	\|z(\cdot-s)\|_{\cH_{\kin,\sigma}(\cD_{T,\ast}\times\bR_v^3)}\le\|z\|_{\cH_{\kin,\sigma}(\cD_{T,\ast}\times\bR_v^3)},
	\qquad 0<s<\varepsilon .
	\]
	Consequently,
	$z_j^\varepsilon=\int_0^\varepsilon\varrho_\varepsilon^{(j)}(s)\,z(\cdot-s)\dd s$ is a
	Bochner integral in $\cH_{\kin,\sigma}(\cD_{T,\ast}\times\bR_v^3)$. The transport operator and both boundary traces, including the outgoing trace on $\gamma^+$, are bounded linear maps on the graph space by \eqref{260902137}, and therefore commute with this Bochner integral. This gives the identities
	in (\ref{it:dyn-real-moll-graph}).  Minkowski's inequality gives
	$\|z_j^\varepsilon\|\le\|\varrho_\varepsilon^{(j)}\|_{L^1(0,\varepsilon)}\|z\|$, which is
	\eqref{eq:dyn-real-moll-graph-bound}.
	
	For (\ref{it:dyn-real-moll-chain}), let $\varphi\in C_c^\infty(I_T^-\times\cD\times\bR^3)$. By Fubini's theorem and integration by parts in the convolution variable,
	\[
	\int z_j^\varepsilon\,\partial_t\varphi
	=-\int z_{j+1}^\varepsilon\,\varphi
	\]
	with no boundary terms because $\varrho_\varepsilon^{(j)}\in C_c^\infty((0,\varepsilon))$. Thus $z_j^\varepsilon=\partial_t^jz^\varepsilon$ for every $j$.
\end{proof}

\subsubsection{The zero-initial temporal chain}\label{app:dyn-zero-initial}

\begin{prop}[Temporal chain]
	\label{prop:dynamic-fixedN-zero-initial-temporal-realization}
	Let $N\in\bN$ and $0<T<\infty$.  Let
	$k_0,\ldots,k_{\mathrm M}$ be measurable on $I_T^-\times\cD\times\bR^3$,
	$S_0,\ldots,S_{\mathrm M}\in L^2(I_T^-;L_x^2H_\sigma^{-1})$ and
	$b\in H^{\mathrm M}(I_T^-;L^2(\gamma^-))$ satisfy
	\begin{align}\label{eq:dynamic-zero-initial-abstract-data}		
		\begin{gathered}
			\partial_tk_j=k_{j+1},
			\qquad
			\partial_tS_j=S_{j+1}
			\quad\text{in }\cD'\ \ (0\le j\le\mathrm M-1),\\
			w^{-2}k_j\in L^\infty(I_T^-;L_x^\infty L_v^2)\ \ (0\le j\le\mathrm M),
		\end{gathered}
	\end{align}
	and $\|\mu^{1/4}k_0\|_{L^\infty(I_T^-\times\cD\times\bR^3)}\le\varepsilon_{\rm dyn}$, with
	$\varepsilon_{\rm dyn}$ of Proposition~\ref{prop:dynamic-alpha-uniform-graph}.  Let $Z_0\in\cH_{\kin,\sigma}(\cD_{T,\ast}\times\bR_v^3)$ satisfy
	\[
	Z_0=0\ \text{on }(-1,-3/4],
	\qquad
	(\partial_t+v\cdot\nabla_x)Z_0+L_{k_0}Z_0=S_0,
	\qquad
	(Z_0)_{\gamma^-}=P_{w,N}(Z_0)_{\gamma^+}+b .
	\]
	Then there are
	$Z_1,\ldots,Z_{\mathrm M}\in\cH_{\kin,\sigma}(\cD_{T,\ast}\times\bR_v^3)$
	such that
	\begin{align}
		Z_j&=0
		&&\text{on }(-1,-3/4],\quad0\le j\le\mathrm M,
		\label{eq:dynamic-zero-initial-temporal-left-zero}\\
		\partial_tZ_j&=Z_{j+1},
		&&0\le j\le\mathrm M-1,
		\label{eq:dynamic-zero-initial-temporal-chain}\\
		(\partial_t+v\cdot\nabla_x)Z_j+L_{k_0}Z_j
		&=\sum_{a=1}^{j}\binom ja
		N[k_a,Z_{j-a}]+S_j,
		&&0\le j\le\mathrm M,
		\label{eq:dynamic-zero-initial-temporal-Z0}\\
		(Z_j)_{\gamma^-}
		&=P_{w,N}(Z_j)_{\gamma^+}+\mathfrak t_{j,N}[Z_0]
		+\partial_t^jb,
		&&0\le j\le\mathrm M.
		\label{eq:dynamic-zero-initial-temporal-boundary}
	\end{align}
	Here $\mathfrak t_{j,N}$ is that of \eqref{eq:dyn-theta-triangular-notation}.
\end{prop}

\begin{proof}
	\noindent\emph{Step 0. Notation and elementary bounds.}  For $0<h<1/4$ write
	\[
	I_{T,h}^-:=(-1,T-h),
	\qquad
	(\tau_hz)(t):=z(t+h),
	\qquad
	\Delta_hz:=\frac{\tau_hz-z}{h}
	\]
	for the shortened interval, the translation in time and the difference quotient.
	
	Let $X$ be a Banach space.  For $z\in H^1(I_T^-;X)$,
	$\Delta_hz(t)=h^{-1}\int_0^h\partial_tz(t+s)\dd s$, so Minkowski's inequality gives
	\begin{equation}
		\label{eq:dynamic-difference-quotient-bounds}
		\|\Delta_hz\|_{L^2(I_{T,h}^-;X)}\le\|\partial_tz\|_{L^2(I_T^-;X)},
		\qquad
		\|w^{-2}\Delta_hz\|_{L^\infty(I_{T,h}^-;L_x^\infty L_v^2)}
		\le\|w^{-2}\partial_tz\|_{L^\infty(I_T^-;L_x^\infty L_v^2)},
	\end{equation}
	the second following from the same integral identity in the indicated weighted mixed norm, whenever its right-hand side is finite.  The derivative identities and weighted $L^\infty$ bounds in \eqref{eq:dynamic-zero-initial-abstract-data} imply $k_a\in H^1(I_T^-;w^2L^2_{x,v})$ for $a\le\mathrm M-1$.
	The translation $t\mapsto t+h$ commutes with $\partial_t+v\cdot\nabla_x$ and with the traces, so that
	\[
	\begin{gathered}
		\|\tau_hz\|_{L^2(I_{T,h}^-;L_x^2H^1_\sigma)}\le\|z\|_{L^2(I_T^-;L_x^2H^1_\sigma)},
		\qquad
		\|(\tau_hz)_{\gamma^-}\|_{L^2(I_{T,h}^-;L^2(\gamma^-))}\le\|z_{\gamma^-}\|_{L^2(\gamma^-_{T,\ast})},\\
		\|(\partial_t+v\cdot\nabla_x)\tau_hz\|_{L^2(I_{T,h}^-;L_x^2H^{-1}_\sigma)}\le\|(\partial_t+v\cdot\nabla_x)z\|_{L^2(I_T^-;L_x^2H^{-1}_\sigma)},
	\end{gathered}
	\]
	while the endpoint term at $t=-1$ becomes $\|z(-1+h)\|_{L^2_{x,v}}$, which is at most $C\|z\|_{\cH_{\kin,\sigma}(\cD_{T,\ast}\times\bR_v^3)}$ by \eqref{260730346}.  By density, $\tau_h$ therefore maps $\cH_{\kin,\sigma}(\cD_{T,\ast}\times\bR_v^3)$ into $\cH_{\kin,\sigma}(I_{T,h}^-\times\cD\times\bR_v^3)$ and
	\[
	\|\tau_hz\|_{\cH_{\kin,\sigma}(I_{T,h}^-\times\cD\times\bR_v^3)}\le\|z\|_{\cH_{\kin,\sigma}(\cD_{T,\ast}\times\bR_v^3)}
	\]
	for $z$ vanishing on $(-1,-3/4]$ (for these, $z(-1+h)=0$ since $h<1/4$).
	
	Since $P_{w,N}(t)=\cP_{0,N}^-+d_\delta(t)J_N\cA_0$ (see after
	\eqref{eq:dyn-theta-triangular-notation}), the Leibniz rule for the difference quotient of
	$d_\delta(t)\,J_N\cA_0z_{\gamma^+}(t)$ gives, for $N\in\bN$,
	$z\in\cH_{\kin,\sigma}(\cD_{T,\ast}\times\bR_v^3)$, $0<h<1/4$ and finite $L\ge0$,
	\begin{equation}
		\label{eq:dynamic-difference-quotient-wall-defect}
		\Delta_h\bigl(P_{w,N}z_{\gamma^+}\bigr)
		=P_{w,N}\bigl(\Delta_hz\bigr)_{\gamma^+}
		+\bigl(\Delta_hd_\delta\bigr)J_N\cA_0(\tau_hz)_{\gamma^+} ,
		\qquad
		\bigl\|w^L\Delta_hd_\delta\bigr\|_{L^\infty}\le C\delta ,
	\end{equation}
	where the bound follows from \eqref{eq:dyn-theta-defect-time-bound} at $j=1$ and the mean-value theorem, since, by Lemma~\ref{lem:dyn-wall-extended}\,(\ref{it:dyn-wall-ext-results}), \eqref{eq:dyn-theta-defect-time-bound} holds for $-1<t<T$, $t\ne0$, and $d_\delta$ is continuous across $t=0$.
	
	By \eqref{eq:velalg-poly} at $(\ell,\mathrm L)=(0,-2)$, applied at each $(t,x)$,
	for every $0<h<1/4$,
	\begin{equation}
		\label{eq:dynamic-difference-quotient-source-slot}
		\|N[k,z]\|_{L^2(I_{T,h}^-;L_x^2H_\sigma^{-1})}
		\le C\|w^{-2}k\|_{L^\infty(I_{T,h}^-;L_x^\infty L_v^2)}
		\|z\|_{L^2(I_{T,h}^-;L_x^2H_\sigma^{1})}.
	\end{equation}
	
	We construct $Z_1,\ldots,Z_{\mathrm M}$ by induction on $j$.  At $j=0$, \eqref{eq:dynamic-zero-initial-temporal-left-zero},
	\eqref{eq:dynamic-zero-initial-temporal-Z0} and
	\eqref{eq:dynamic-zero-initial-temporal-boundary} are the hypotheses on $Z_0$, with
	$\mathfrak t_{0,N}=0$.
	
	\noindent\emph{Step 1. The difference-quotient problem.}  Let $0\le j\le\mathrm M-1$ and suppose
	$Z_0,\ldots,Z_j$ have been constructed so that
	\eqref{eq:dynamic-zero-initial-temporal-left-zero}--\eqref{eq:dynamic-zero-initial-temporal-boundary} hold for these functions, and
	$$
	\partial_t(Z_i)_{\gamma^\pm}=(Z_{i+1})_{\gamma^\pm}\qquad \text{in }\cD'(I_T^-;L^2(\gamma^\pm))\quad \text{for }i<j.
	$$
	Fix $0<h<1/4$ and set
	\[
	Y_h:=\Delta_hZ_j
	\quad\text{on }I_{T,h}^- .
	\]
	Apply $\Delta_h$ to \eqref{eq:dynamic-zero-initial-temporal-Z0} for $Z_j$, and use the bilinearity of $N$ in the form
	\begin{equation*}
		\Delta_h\bigl(N[k_a,Z_{j-a}]\bigr)
		=N[\tau_hk_a,\Delta_hZ_{j-a}]
		+N[\Delta_hk_a,Z_{j-a}]
	\end{equation*}
	for $0\le a\le j$; the case $a=0$ enters through $L_{k_0}=L-N[k_0,\cdot\,]$.  Since $\tau_h$ commutes with $\partial_t+v\cdot\nabla_x$ and with $L$, we obtain, on $I_{T,h}^-$,
	\begin{equation}
		\label{eq:dynamic-difference-quotient-equation}
		\begin{aligned}
			(\partial_t+v\cdot\nabla_x)Y_h+L_{\tau_hk_0}Y_h
			={}&\Delta_hS_j
			+N[\Delta_hk_0,Z_j]\\
			&+\sum_{a=1}^{j}\binom ja
			\Bigl(N[\tau_hk_a,\Delta_hZ_{j-a}]
			+N[\Delta_hk_a,Z_{j-a}]\Bigr),
		\end{aligned}
	\end{equation}
	together with
	\begin{align}
		\label{eq:dynamic-difference-quotient-data}
		\begin{aligned}
			Y_h(-1)\,&=0,\\
			(Y_h)_{\gamma^-}\,&=P_{w,N}(Y_h)_{\gamma^+}
			+\bigl(\Delta_hd_\delta\bigr)J_N\cA_0\bigl(\tau_hZ_j\bigr)_{\gamma^+}
			+\Delta_h\bigl(\mathfrak t_{j,N}[Z_0]\bigr)
			+\Delta_h\bigl(\partial_t^jb\bigr).
		\end{aligned}
	\end{align}
	The initial value in \eqref{eq:dynamic-difference-quotient-data} vanishes
	because $Z_j=0$ on $(-1,-3/4]$ by
	\eqref{eq:dynamic-zero-initial-temporal-left-zero} and $h<1/4$.  The boundary identity follows by taking $\Delta_h$ of \eqref{eq:dynamic-zero-initial-temporal-boundary} and using \eqref{eq:dynamic-difference-quotient-wall-defect} and commutation of time translation with the traces.
	
	\noindent\emph{Step 2. Bounds uniform in $h$.}  Since
	\[
	\|\mu^{1/4}\tau_hk_0\|_{L^\infty(I_{T,h}^-\times\cD\times\bR^3)}\le\|\mu^{1/4}k_0\|_{L^\infty(I_T^-\times\cD\times\bR^3)},
	\]
	the smallness of $k_0$ holds with $(\tau_hk_0,I_{T,h}^-)$ in place of $(k_0,I_T^-)$.
	Hence, by Lemma~\ref{lem:dyn-wall-extended}\,(\ref{it:dyn-wall-ext-results}),
	Proposition~\ref{prop:dynamic-inhomogeneous-physical-solvability}\,(\ref{it:solv-spectral}) applies to
	\eqref{eq:dynamic-difference-quotient-equation}--%
	\eqref{eq:dynamic-difference-quotient-data}, with a constant independent of
	$h$ and $N$.
	
	For the bulk terms, \eqref{eq:dynamic-difference-quotient-bounds} and the chain of the
	$k_a$ in \eqref{eq:dynamic-zero-initial-abstract-data} give
	\[
	\|w^{-2}\Delta_hk_a\|_{L^\infty(I_{T,h}^-;L_x^\infty L_v^2)}\le\|w^{-2}k_{a+1}\|_{L^\infty(I_T^-;L_x^\infty L_v^2)}
	\]
	for
	$0\le a\le\mathrm M-1$, which covers every $a\le j$ in \eqref{eq:dynamic-difference-quotient-equation}, since $j\le\mathrm M-1$.  With \eqref{eq:dynamic-zero-initial-temporal-chain}, the first bound of \eqref{eq:dynamic-difference-quotient-bounds} also gives
	\[
	\|\Delta_hZ_{j-a}\|_{L^2(I_{T,h}^-;L_x^2H_\sigma^1)}
	\le\|Z_{j-a+1}\|_{\cH_{\kin,\sigma}(\cD_{T,\ast}\times\bR_v^3)},
	\qquad 1\le a\le j.
	\]
	We also have
	\[
	\|w^{-2}\tau_hk_a\|_{L^\infty(I_{T,h}^-;L_x^\infty L_v^2)}\le\|w^{-2}k_a\|_{L^\infty(I_T^-;L_x^\infty L_v^2)} .
	\]
	These bounds and \eqref{eq:dynamic-difference-quotient-source-slot} give, for $1\le a\le j$,
	\[
	\begin{aligned}
		&\|N[\tau_hk_a,\Delta_hZ_{j-a}]\|_{L^2(I_{T,h}^-;L_x^2H_\sigma^{-1})}
		+\|N[\Delta_hk_a,Z_{j-a}]\|_{L^2(I_{T,h}^-;L_x^2H_\sigma^{-1})}\\
		&\qquad\le C\|w^{-2}k_a\|_{L^\infty(I_T^-;L_x^\infty L_v^2)}
		\|Z_{j-a+1}\|_{\cH_{\kin,\sigma}(\cD_{T,\ast}\times\bR_v^3)}
		+C\|w^{-2}k_{a+1}\|_{L^\infty(I_T^-;L_x^\infty L_v^2)}
		\|Z_{j-a}\|_{\cH_{\kin,\sigma}(\cD_{T,\ast}\times\bR_v^3)} ,
	\end{aligned}
	\]
	and, in the same way,
	\[
	\|N[\Delta_hk_0,Z_j]\|_{L^2(I_{T,h}^-;L_x^2H_\sigma^{-1})}\le C\|w^{-2}k_1\|_{L^\infty(I_T^-;L_x^\infty L_v^2)}\|Z_j\|_{\cH_{\kin,\sigma}(\cD_{T,\ast}\times\bR_v^3)}.
	\]
	Thus all collision terms in \eqref{eq:dynamic-difference-quotient-equation} are bounded independently of $h$. In addition, \eqref{eq:dynamic-difference-quotient-bounds} gives
	\[
	\|\Delta_hS_j\|_{L^2(I_{T,h}^-;L_x^2H_\sigma^{-1})}\le\|S_{j+1}\|_{L^2(I_T^-;L_x^2H_\sigma^{-1})} .
	\]
	
	For the boundary data, \eqref{eq:dynamic-difference-quotient-bounds} provides
	\[
	\|\Delta_h\partial_t^jb\|_{L^2(I_{T,h}^-;L^2(\gamma^-))}\le\|\partial_t^{j+1}b\|_{L^2(I_T^-;L^2(\gamma^-))},\qquad 
	\|\Delta_h\mathfrak t_{j,N}[Z_0]\|_{L^2(I_{T,h}^-;L^2(\gamma^-))}\le\|\partial_t\mathfrak t_{j,N}[Z_0]\|_{L^2(I_T^-;L^2(\gamma^-))} .
	\]
	By the
	Leibniz rule applied to \eqref{eq:dyn-theta-triangular-notation} and the trace chain for $i<j$, and then by Pascal's rule with $\binom j{j+1}=0$,
	\[
	\begin{aligned}
		\partial_t\mathfrak t_{j,N}[Z_0]
		&=\sum_{k=1}^{j}\binom jk\Bigl((\partial_t^{k+1}d_\delta)J_N\cA_0(Z_{j-k})_{\gamma^+}+(\partial_t^kd_\delta)J_N\cA_0(Z_{j+1-k})_{\gamma^+}\Bigr)\\
		&=\mathfrak t_{j+1,N}[Z_0]-(\partial_td_\delta)J_N\cA_0(Z_j)_{\gamma^+}.
	\end{aligned}
	\]
	Hence, by this identity, \eqref{eq:dyn-theta-triangular-bound} at $j+1$, \eqref{eq:dyn-theta-defect-time-bound} and \eqref{260902137},
	\[
	\begin{aligned}
		\|\partial_t\mathfrak t_{j,N}[Z_0]\|_{L^2(\gamma^-_{T,\ast})}
		&\le\|\mathfrak t_{j+1,N}[Z_0]\|_{L^2(\gamma^-_{T,\ast})}+\|(\partial_td_\delta)J_N\cA_0(Z_j)_{\gamma^+}\|_{L^2(\gamma^-_{T,\ast})}\\
		&\le C\delta\sum_{i\le j}\|Z_i\|_{\cH_{\kin,\sigma}(\cD_{T,\ast}\times\bR_v^3)} .
	\end{aligned}
	\]
	The bound \eqref{eq:dyn-theta-triangular-bound} extends to $I_T^-$: its proof uses only \eqref{eq:dyn-theta-defect-time-bound}, which holds there almost everywhere by Lemma~\ref{lem:dyn-wall-extended}\,(\ref{it:dyn-wall-ext-results}).  The term
	$(\Delta_hd_\delta)J_N\cA_0(\tau_hZ_j)_{\gamma^+}$ satisfies, by \eqref{eq:dynamic-difference-quotient-wall-defect}, $0\le J_N\le\mathbf I$, \eqref{eq:outgoing-moment-flux-bound} and \eqref{260902137},
	\[
	\|(\Delta_hd_\delta)J_N\cA_0(\tau_hZ_j)_{\gamma^+}\|_{L^2(I_{T,h}^-;L^2(\gamma^-))}\le C\delta\|Z_j\|_{\cH_{\kin,\sigma}(\cD_{T,\ast}\times\bR_v^3)} ,
	\]
	with $C$ independent of $h$ and $N$.
	Proposition~\ref{prop:dynamic-inhomogeneous-physical-solvability}\,(\ref{it:solv-spectral})
	then gives
	\begin{equation}
		\label{eq:dynamic-difference-quotient-uniform}
		\begin{aligned}
			&\sup_{0<h<1/4}\|Y_h\|_{\cH_{\kin,\sigma}(I_{T,h}^-\times\cD\times\bR_v^3)}\\
			\le{}&C\Bigl(
			\|S_{j+1}\|_{L^2(I_T^-;L_x^2H_\sigma^{-1})}
			+\|\partial_t^{j+1}b\|_{L^2(\gamma_{T,\ast}^-)}
			+\sum_{i=0}^{j}
			\|Z_i\|_{\cH_{\kin,\sigma}(\cD_{T,\ast}\times\bR_v^3)}\Bigr),
		\end{aligned}
	\end{equation}
	with $C$ depending on $T$, $\delta$ and $\|w^{-2}k_a\|_{L^\infty(I_T^-;L_x^\infty L_v^2)}$, $a\le j+1$, and not
	on $h$.
	
	\noindent\emph{Step 3. The limit $h\to0$.}  Fix $T'<T$ and set $Z_{j+1}:=\partial_tZ_j$ in $\cD'(I_T^-\times\cD\times\bR^3)$.  For
	$\varphi\in C_c^\infty((-1,T')\times\cD\times\bR^3)$ and $h$ small,
	\[
	\int Y_{h}\,\varphi
	=-\int Z_j\,\Delta_{-h}\varphi
	\longrightarrow-\int Z_j\,\partial_t\varphi,
	\]
	so $Y_h\to Z_{j+1}$ in $\cD'$.  The uniform bound \eqref{eq:dynamic-difference-quotient-uniform} and Lemma~\ref{lem:dynamic-weak-graph-bulk-identification}\,(\ref{it:closure-evolutionary}), applied on $(-1,T')$, give $Z_{j+1}\in\cH_{\kin,\sigma}((-1,T')\times\cD\times\bR_v^3)$ and $Y_h\rightharpoonup Z_{j+1}$ weakly in that space, with
	\[
	\|Z_{j+1}\|_{\cH_{\kin,\sigma}((-1,T')\times\cD\times\bR_v^3)}\le\sup_{0<h<1/4}\|Y_h\|_{\cH_{\kin,\sigma}(I_{T,h}^-\times\cD\times\bR_v^3)}
	\qquad\text{for every }T'<T .
	\]
	The slices then satisfy
	\[
	\|Z_{j+1}(T')\|_{L^2(\cD\times\bR_v^3)}\le C\sup_{0<h<1/4}\|Y_h\|_{\cH_{\kin,\sigma}(I_{T,h}^-\times\cD\times\bR_v^3)},
	\]
	with $C$ independent of $T'$, by
	\eqref{260902137}; let $Z_{j+1}(T)$ be a weak limit of $Z_{j+1}(T'_k)$, $T'_k\uparrow T$.  The
	identity \eqref{2607151226} on $(-1,T'_k)$ then passes to $(-1,T)$, and
	Proposition~\ref{H_kin=weak_landau}\,(1) gives
	$Z_{j+1}\in\cH_{\kin,\sigma}(\cD_{T,\ast}\times\bR_v^3)$.
	By \eqref{eq:temporal-boundary-generator} of Proposition~\ref{prop:weak-boundary-generator}\,(\ref{it:weak-generator-evolutionary}) for the pair $(Z_j,Z_{j+1})$, applied on $I_T^-$,
	$\partial_t\bigl((Z_j)_{\gamma^\pm}\bigr)=(Z_{j+1})_{\gamma^\pm}$ in $\cD'(I_T^-;L^2(\gamma^\pm))$.
	
	\noindent\emph{Step 4. The equation for $Z_{j+1}$.}  Differentiate \eqref{eq:dynamic-zero-initial-temporal-Z0} for $Z_j$ in
	$\cD'(I_T^-\times\cD\times\bR^3)$.  The Leibniz rule \eqref{eq:dyn-spacetime-leibniz}
	on $I_T^-$ at $\varpi=\tfrac14$ applies to $N[k_a,Z_{j-a}]$, since,
	by \eqref{eq:dynamic-zero-initial-abstract-data} and $\mu^{1/4}\le Cw^{-2}$, $\mu^{1/4}k_a$ and $\mu^{1/4}k_{a+1}$ lie in $L^2_{\rm loc}(I_T^-\times\cD;L^2_v)$ and
	$Z_i\in L^2(I_T^-;L_x^2H^1_\sigma)$ for $i\le j+1$.  With $L_{k_0}=L-N[k_0,\cdot\,]$,
	$[\partial_t,v\cdot\nabla_x]=[\partial_t,L]=0$ and the chains of the $k_a$, the $S_j$ and the
	$Z_i$, this gives
	\[
	(\partial_t+v\cdot\nabla_x)Z_{j+1}+L_{k_0}Z_{j+1}
	=N[k_1,Z_j]
	+\sum_{a=1}^{j}\binom ja\bigl(N[k_{a+1},Z_{j-a}]+N[k_a,Z_{j+1-a}]\bigr)
	+S_{j+1}.
	\]
	Pascal's rule gives \eqref{eq:dynamic-zero-initial-temporal-Z0} at order $j+1$: shifting the index in the $k_{a+1}$ terms and including $N[k_1,Z_j]$ yields
	\[
	\begin{aligned}
		&N[k_1,Z_j]+\sum_{a=1}^{j}\binom ja\bigl(N[k_{a+1},Z_{j-a}]+N[k_a,Z_{j+1-a}]\bigr)=\sum_{a=1}^{j+1}\binom{j+1}aN[k_a,Z_{j+1-a}] .
	\end{aligned}
	\]
	
	\noindent\emph{Step 5. The boundary condition for $Z_{j+1}$.}  The trace chain established above gives $\partial_t(Z_i)_{\gamma^\pm}=(Z_{i+1})_{\gamma^\pm}$ for $i\le j$. Moreover, Lemma~\ref{lem:dyn-wall-extended}\,(\ref{it:dyn-wall-ext-results}) provides continuity across $t=0$ and bounds on $I_T^-$ for $\partial_t^kd_\delta$, $k\le j+1$. Thus the Leibniz rule in the truncated wall law \eqref{eq:dyn-temporal-boundary-law}, at order $j+1$ with $f=Z_0$, gives \eqref{eq:dynamic-zero-initial-temporal-boundary} for $Z_{j+1}$ on $I_T^-$.
	Finally, $Z_{j+1}=\partial_tZ_j$ vanishes on $(-1,-3/4]$, which completes
	\eqref{eq:dynamic-zero-initial-temporal-left-zero}--%
	\eqref{eq:dynamic-zero-initial-temporal-chain} for $Z_0,\ldots,Z_{j+1}$.  Together with the trace chain at $i=j$, these identities close the induction through $j+1=\mathrm M$.
\end{proof}

\subsection{The doubly regularized solution}\label{app:dyn-doubly-regularized}

We prove Proposition~\ref{prop:dyn-fixedN-realization}. The boundary and regularity estimates below are formulated for general backgrounds.

\begin{corollary}[The mixed spacetime rungs]
	\label{cor:dyn-mixed-rungs}
	Let $(t_0,t_1)$ be an interval with $t_1-t_0\ge1$, and let $f$ satisfy
	\eqref{eq:mixed-rung-chain} and \eqref{eq:mixed-rung-local-integrability} over the
	ordered complete words on
	$(t_0,t_1)\times\cD\times\bR_v^3$ at the order
	$M'=\mathrm M$ of Corollary~\ref{cor:ordered-word-rungs-mixed}, and, for almost every $t\in(t_0,t_1)$, \eqref{eq:dynamic-pointwise-finite-norms} at $f(t)$, with $\mathrm M$ in place of $\mathrm N$,
	$(\mathrm L_{\rm G},a_{\rm G})=(-2,q_\flat)$ and the spatial base exponents $(p_{\rm G},p_{\rm c})=(\infty,p_\sharp)$.  Let all norms in $t$ below
	be taken over $(t_0,t_1)$.
	Then
	\begin{equation}
		\label{eq:dynamic-spatial-G-rung}
		\sum_{I\in\mathfrak A^k}\bigl\|w^{-2}\mathfrak D_If\bigr\|
		_{L_t^{q_\flat\mathrm M/k}
			L_x^{q_\flat\mathrm M/k}L_v^2}
		\le C\,\|w^{-2}f\|_{L_t^\infty L_x^\infty L_v^2}^{1-k/\mathrm M}\,\mathfrak T_{\mathrm M}[f]^{k/\mathrm M},
		\qquad 1\le k\le\mathrm M-1 ,
	\end{equation}
	where
	\begin{equation}
		\label{eq:dyn-spatial-gaussian-top}
		\mathfrak T_{\mathrm M}[f]:=\bigl\|\cG^{\sharp,-2}_{\mathrm M,q_\flat}(f)\bigr\|
		_{L_t^{q_\flat}}
		+\|w^{-2}f\|_{L_t^{q_\flat}L_x^\infty L_v^2} ,
	\end{equation}
	with $C$ depending only on the structural data (in particular independent of $t_0$ and $t_1$).
\end{corollary}

\begin{proof}
	Corollary~\ref{cor:ordered-word-gaussian-rungs}\,(\ref{it:rungs-evolutionary}) applies on $(t_0,t_1)$ with $\mathrm N$ replaced by $\mathrm M$. Its temporal-slicing hypotheses \eqref{eq:temporal-slicing-integrability}--\eqref{eq:temporal-slicing-spacetime-chain} are precisely \eqref{eq:mixed-rung-local-integrability} and \eqref{eq:mixed-rung-chain} for the ordered complete words; the pointwise bounds \eqref{eq:dynamic-pointwise-finite-norms} are assumed.  We take $(\mathrm L_{\rm G},a_{\rm G})=(-2,q_\flat)$, the temporal base exponent $r_0=\infty$, the temporal top exponent $q_\flat\le r_0$ and the Gaussian spatial base exponent $p_{\rm G}=\infty$. The temporal and spatial exponents then coincide:
	$\mathsf e(\mathrm M/k;q_\flat,\infty)=q_\flat\mathrm M/k$.  
\end{proof}

\begin{lemma}[Extended background]
	\label{lem:dyn-real-hypotheses-ball}
	\label{lem:dyn-breg-background}
	Let $1\le T<\infty$, let $g$ satisfy \eqref{eq:dyn-ball-jets} and $\widetilde{\mathfrak Q}^\dy_T(g)\le r$ for some
	$r\ge0$, let $0<\varepsilon\le\varepsilon_\lambda$, and let $\bar g_0$ and $g^\varepsilon$ be defined from $g$ by
	\eqref{eq:dynamic-base-background} and \eqref{eq:dyn-real-mollified-objects}.
	\begin{enumerate}
		\item\label{it:dyn-real-hyp-spatial}
		For every ordered complete spatial word $I$ with $|I|\le\mathrm M$, the background
		norms over $I_T^-$ satisfy
		\begin{equation}
			\label{eq:dyn-real-background-norms}
			\|w^{-2}\bar g_0\|_{L^\infty(I_T^-;L_x^\infty L_v^2)}
			+\|\bar g_0\|_{L^\infty(I_T^-;L_x^{p_\sharp}H_\sigma^1)}
			+\mathfrak T_{\mathrm M}[\bar g_0]<\infty,
			\qquad
			w^{\mathrm M-|I|}\mathfrak D_I\bar g_0\in\cH_{\kin,\sigma}(\cD_{T,\ast}\times\bR_v^3),
		\end{equation}
		with $\mathfrak T_{\mathrm M}$ of
		\eqref{eq:dyn-spatial-gaussian-top}.
		
		\item\label{it:dyn-real-hyp-basic}
		$g\in\cH_{\kin,\sigma}(\cD_T\times\bR_v^3)$ and $g(0)=z_0$.
		
		\item\label{it:dyn-breg-bg-norms}
		Let $0\le a\le\mathrm M$, and let $\mathsf T$ be a weight, a complete derivative
		$\mathfrak D_K$, a flow translate $T_s^X$ with $X=\tau_\nu$, $\nu\in\mathfrak A$, the transport operator, or a composition of
		these.  Then, for $1\le p\le\infty$ and every normed space $\mathcal Y$ of functions of $(x,v)$,
		\begin{equation}
			\label{eq:dyn-breg-background}
			\|\mathsf Tg^\varepsilon_a\|_{L^p(I_T^-;\mathcal Y)}
			\le\|\varrho_\varepsilon^{(a)}\|_{L^1}\|\mathsf T\bar g_0\|_{L^p(I_T^-;\mathcal Y)}
			+2\|\mathsf T\partial_t^a\hat g^{(0)}\|_{L^p(I_T^-;\mathcal Y)} ,
		\end{equation}
		and \eqref{eq:dyn-breg-background} also holds with the norm of
		$\cH_{\kin,\sigma}(\cD_{T,\ast}\times\bR_v^3)$ in place of that of $L^p(I_T^-;\mathcal Y)$
		when $\mathsf T=w^\ell\mathfrak D_K$, with $\ell\ge0$ and $K$ an ordered complete spatial
		word.  Moreover, \eqref{eq:dyn-real-background-norms} holds with
		$\partial_t^a\hat g^{(0)}$ in place of $\bar g_0$.  Finally, with $s_0$ of
		\eqref{eq:complete-flow-uniform-bounds} and the norms in $t$ taken over $I_T^-$, for every
		ordered complete spatial word $K$ with $|K|\le\mathrm M-1$ and every
		$X=\tau_\nu$ with $\nu\in\mathfrak A$,
		\begin{equation}
			\label{261016990}
			\sup_{|s|\le s_0}
			\bigl\|w^{-2}T_s^X\mathfrak D_K\partial_t^a\hat g^{(0)}\bigr\|
			_{L_t^{\mathsf g_\dy(|K|)}L_x^{\mathsf g_\dy(|K|)}L_v^2}<\infty ,
		\end{equation}
		where $\mathsf g_\dy(|K|):=q_\flat\mathrm M/|K|$, with $\mathsf g_\dy(0):=\infty$.
		\item\label{it:dyn-breg-bg-rungs}
		For every ordered complete spatial word $I$ with $|I|\le\mathrm M$, the functions
		$w^{\mathrm M-|I|}\mathfrak D_Ig^\varepsilon$ and $w^{\mathrm M-|I|}\mathfrak D_I\bar g_0$ belong to
		$\cH_{\kin,\sigma}(\cD_{T,\ast}\times\bR_v^3)$, so \eqref{eq:graph-to-mixed-integrability-evolutionary}
		applies to them on $I_T^-$, and $g^\varepsilon$ and $\bar g_0$ satisfy, on $I_T^-$, the rung bounds
		\eqref{eq:dynamic-spatial-G-rung}.
		
	\end{enumerate}
\end{lemma}

\begin{proof}

	\emph{(\ref{it:dyn-real-hyp-spatial}).}  On $(0,T)$, $\|w^{-2}g\|_{L_t^\infty L_x^\infty L_v^2}$ and
	$\|g\|_{L^\infty_tL_x^{p_\sharp}H_\sigma^1}$ are at most
	$C\cA_T[g]\le Cr$, by \eqref{eq:dyn-coordinate-A-reading} and \eqref{eq:dyn-ball-A}.  Next, the function
	\[
	w^{\mathrm M-|I|}\mathfrak D_Ig=w^{\mathrm M-|I|}(\mathfrak D_Ig_\dy+\mathfrak D_If_\st)
	\]
	is in $\cH_{\kin,\sigma}(\cD_T\times\bR_v^3)$.  For $g_\dy$, the row
	$\mathsf W_{|I|}(T)[g_\dy]\le r$ of \eqref{eq:dyn-ball-rows} carries the
	weight \[
	\ell^\dy_{|I|}=\ell^\dy_\sharp+(\mathrm M-|I|)A_\diamond\ge\mathrm M-|I|
	\] (by \eqref{eq:dyn-weight-choice}, since $A_\diamond\ge1$),
	and Lemma~\ref{lem:compact-velocity-support}\,(1) with the $C^1_v$ multiplier
	$w^{\mathrm M-|I|-\ell^\dy_{|I|}}$ lowers the weight to $\mathrm M-|I|$.  For $f_\st$,
	$w^{\ell_{|I|}}\mathfrak D_If_\st\in\cH_{\kin,\sigma}(\cD\times\bR_v^3)$ by
	\eqref{eq:dyn-theta-background} ($|I|\le\mathrm M\le\mathrm N-1$), with
	$\ell_{|I|}\ge\mathrm N-|I|>\mathrm M-|I|$ by \eqref{eq:finite-hierarchy-strict-ladder}. For finite $T$, the stationary function also lies in $\cH_{\kin,\sigma}(\cD_T\times\bR_v^3)$.  In $\mathfrak T_{\mathrm M}$, the top norm $\|\cG^{\sharp,-2}_{\mathrm M,q_\flat}(g_\dy)\|_{L_t^{q_\flat}}$ is at
	most $C\cN_{\mathrm M}(T)[g_\dy]$, by the proof of \eqref{eq:dyn-slot-tops-rows} with \eqref{eq:graph-to-mixed-integrability-evolutionary} at the time exponent $q_\flat$ in place of $2$; moreover, $\|\cG^{\sharp,-2}_{\mathrm M,q_\flat}(f_\st)\|_{L_t^{q_\flat}}<\infty$
	on $(0,T)$ by \eqref{eq:graph-to-mixed-integrability}; and
	\[
	\|w^{-2}g\|_{L_t^{q_\flat}L_x^\infty L_v^2}\le T^{1/q_\flat}\|w^{-2}g\|_{L_t^\infty L_x^\infty L_v^2}.
	\]
	
	On $(-1,0)$, $\bar g_0=\check{E}=\chi_-\check{\cP}$ is a cutoff times a polynomial in
	$t$ with the coefficients $z_j$, so
	\[
	\|w^{-2}\check{E}\|_{L^\infty((-1,0);L_x^\infty L_v^2)}+\|\check{E}\|_{L^\infty((-1,0);L_x^{p_\sharp}H_\sigma^1)}+\mathfrak T_{\mathrm M}[\check{E}]<\infty ,
	\]
	with $\mathfrak T_{\mathrm M}$ taken over $(-1,0)$, by $\mathsf I<\infty$,
	\eqref{eq:dyn-breg-lift-weights}, \eqref{eq:graph-to-mixed-integrability} and
	$w^{-2}\in L_v^2$.  A letter
	$\tau_\nu\cdot\nabla_x+((\nabla\tau_\nu)v)\cdot\nabla_v$ loses one power of $w$.  So
	$w^{\mathrm M-|I|}\mathfrak D_I\check{E}$ and its image under $\partial_t+v\cdot\nabla_x$ require
	\eqref{eq:dyn-breg-lift-weights} only at the weight $L=\mathrm M+4$ and for derivatives
	of order at most $\mathrm M+2$; as permitted by the choice of $\check A$,
	\eqref{eq:dyn-breg-lift-weights} gives
	\[
	\sum_{j\le\mathrm M}\sum_{|\alpha|+|\beta|\le\mathrm M+2}\bigl\|w^{\mathrm M+4}\partial_x^\alpha\partial_v^\beta z_j\bigr\|_{L^\infty(\cD\times\bR^3)}\le C\mathsf I .
	\]
	Moreover,
	\[
	w^{\mathrm M-|I|}\mathfrak D_I\check{E}=\chi_-\sum_{j\le\mathrm M}\tfrac{t^j}{j!}\,w^{\mathrm M-|I|}\mathfrak D_Iz_j,
	\]
	and this function lies in $\cH_{\kin,\sigma}((-1,0)\times\cD\times\bR_v^3)$ by the argument for $\check{E}$ in the proof of Lemma~\ref{lem:dynamic-zero-initial-source-inflow-lift}\,(\ref{it:dyn-lift-basic}), with $w^{\mathrm M-|I|}\mathfrak D_Iz_j$ in place of $z_j$: $w^{\mathrm M-|I|}\mathfrak D_Iz_j\in\cH_{\kin,\sigma}(\cD\times\bR_v^3)$ by \eqref{eq:dyn-data-jets} and Lemma~\ref{lem:compact-velocity-support}\,(1) with the $C^1_v$ multiplier $w^{\mathrm M-|I|-\ell^\dy_0}$, and $\|w^{\mathrm M-|I|}\mathfrak D_Iz_j\|_{L^2_xH^{-1}_\sigma}\le C\|w^{\ell^\dy_0}\mathfrak D_Iz_j\|_{L^2_xH^1_\sigma}$ by \eqref{eq:Hsigma-dual-weighted} and \eqref{eq:auxiliary-coercive-equivalence}, since $\mathrm M-|I|+1\le\ell^\dy_0$.
	
	At $t=0$, the background and the negative-time lift have the same zeroth time trace: $g(0)=z_0$ by \eqref{eq:dyn-ball-jets} at $j=0$ and
	\eqref{eq:dynamic-initial-p-u}, and $z_0=\check{E}(0)$ by \eqref{eq:dynamic-negative-time-jet-matching}.  The distributional identities for the time-independent complete fields, together with the graph-space time traces, identify both traces of $\mathfrak D_I\bar g_0$ at $t=0$ with $\mathfrak D_Iz_0$. Hence \eqref{eq:dyn-time-gluing} at $(a,b,c)=(-1,0,T)$ gives
	$w^{\mathrm M-|I|}\mathfrak D_I\bar g_0\in\cH_{\kin,\sigma}(\cD_{T,\ast}\times\bR_v^3)$.
	
	\emph{(\ref{it:dyn-real-hyp-basic}).}  We have
	$g\in\cH_{\kin,\sigma}(\cD_T\times\bR_v^3)$ by (\ref{it:dyn-real-hyp-spatial}) of this lemma at $I=\varnothing$, restricted to $(0,T)$, together with Lemma~\ref{lem:compact-velocity-support}\,(1).  Moreover, $g(0)=z_0$ by \eqref{eq:dyn-ball-jets} at
	$j=0$ and \eqref{eq:dynamic-initial-p-u}.

	\emph{(\ref{it:dyn-breg-bg-norms}).}  By (\ref{it:dyn-real-hyp-basic}), Lemma~\ref{lem:dyn-real-increment} applies.  The second identity of \eqref{eq:dyn-real-increment-split} gives \eqref{eq:dyn-breg-background} by commutation of $\mathsf T$ with time convolution and Minkowski's inequality, using $\|\varrho_\varepsilon\|_{L^1}=1$.  For the graph norm of $\cH_{\kin,\sigma}$, \eqref{eq:dyn-breg-background} is
	\eqref{eq:dyn-real-moll-graph-bound} at $\mathsf T\bar g_0$ and at
	$\mathsf T\partial_t^a\hat g^{(0)}$, which vanish on $(-1,-3/4]$ (as required there) by
	\eqref{eq:dynamic-base-background} and Lemma~\ref{lem:dyn-fp-seed}.  Next, we show that \eqref{eq:dyn-real-background-norms} holds with
	$\partial_t^a\hat g^{(0)}$ in place of $\bar g_0$.  By Lemma~\ref{lem:dyn-fp-seed},
	$\partial_t^a\hat g^{(0)}$ is on $I_T^-$ a finite combination of
	$z_0,\ldots,z_{\mathrm M}$ and $f_\st$ with coefficients smooth in $t$. Hence the argument used for
	$\check{E}$ and $f_\st$ in the proof of (\ref{it:dyn-real-hyp-spatial}) applies. It uses
	\eqref{eq:dyn-breg-lift-weights}, $\mathsf I<\infty$ and
	\eqref{eq:dyn-theta-background}, and, for the first two norms of \eqref{eq:dyn-real-background-norms} at $f_\st$ on $(0,T)$,
	\eqref{eq:dyn-coordinate-A-reading} at $\ell=0$ and the bound $\cA_T[f_\st]\le C\delta$ of Lemma~\ref{lem:dyn-fp-seed}.
	
	It remains to prove \eqref{261016990}.  Corollary~\ref{cor:dyn-mixed-rungs} applies at
	$\partial_t^a\hat g^{(0)}$ with $(t_0,t_1)=I_T^-$; its hypotheses follow from the norms of \eqref{eq:dyn-real-background-norms}
	as those at $\bar g_0$ do in the proof of (\ref{it:dyn-breg-bg-rungs}) below.  This gives \eqref{261016990} at $s=0$ for $|K|\ge1$, by \eqref{eq:dynamic-spatial-G-rung} at $k=|K|\le\mathrm M-1$; at $K=\varnothing$, where $\mathsf g_\dy(0)=\infty$, it is the bound
	$\|w^{-2}\partial_t^a\hat g^{(0)}\|_{L^\infty(I_T^-;L_x^\infty L_v^2)}<\infty$ proved above. Then
	\eqref{eq:complete-flow-fiber-bound}, taken with the outer $L_t^{\mathsf g_\dy(|K|)}$, gives it for $|s|\le s_0$.
	
	\emph{(\ref{it:dyn-breg-bg-rungs}).}  At $\bar g_0$ the functions $w^{\mathrm M-|I|}\mathfrak D_I\bar g_0$ lie in $\cH_{\kin,\sigma}(\cD_{T,\ast}\times\bR_v^3)$ by
	\eqref{eq:dyn-real-background-norms}.  We apply Corollary~\ref{cor:dyn-mixed-rungs} on $I_T^-$ to $\bar g_0$. The required hypotheses follow from \eqref{eq:dyn-real-background-norms}.  Indeed, the weight $w^{\mathrm M-|I|}$ of
	\eqref{eq:dyn-real-background-norms} is that of
	$\cG^{0}_{\mathrm M,2}$ in \eqref{eq:dynamic-pointwise-finite-norms} with $\mathrm M$ in place of $\mathrm N$; since $\cG^{0}_{\mathrm M,2}$ is a sum of $L_x^2H^1_\sigma$ norms (see \eqref{eq:GN-lower-order}) and $\cH_{\kin,\sigma}\subset L^2_tL^2_xH^1_\sigma$, this gives
	$\cG^{0}_{\mathrm M,2}(\bar g_0)\in L^2(I_T^-)$.  As $w^{-k}$, $k\ge0$, is a
	$C^1_v$ multiplier, \eqref{eq:Hsigma-multiplier} also gives
	$\mathfrak D_I\bar g_0\in\cH_{\kin,\sigma}(\cD_{T,\ast}\times\bR_v^3)$, hence
	\eqref{eq:mixed-rung-local-integrability}.  By $\mathfrak T_{\mathrm M}[\bar g_0]<\infty$,
	$\cG^{\sharp,-2}_{\mathrm M,q_\flat}(\bar g_0)\in L^{q_\flat}(I_T^-)$.  The jets are
	the distributional derivatives of $\bar g_0$, so \eqref{eq:mixed-rung-chain} holds over the ordered complete words; the same
	holds at $g^\varepsilon$.  Finally, since the weights act in $v$ only and $\bar g_0=0$ on
	$(-1,-3/4]$, (\ref{it:dyn-breg-bg-norms}) applies to the norms of \eqref{eq:dyn-real-background-norms}, the graph norms of $w^{\mathrm M-|I|}\mathfrak D_I\bar g_0$ included: by \eqref{eq:dyn-breg-background} at $a=0$, \eqref{eq:dyn-real-background-norms} holds with $g^\varepsilon$ in place of $\bar g_0$,
	since it holds for $\bar g_0$ by (\ref{it:dyn-real-hyp-spatial}) and, by (\ref{it:dyn-breg-bg-norms}), with $\hat g^{(0)}$ in place of $\bar g_0$.  Hence $w^{\mathrm M-|I|}\mathfrak D_Ig^\varepsilon\in\cH_{\kin,\sigma}(\cD_{T,\ast}\times\bR_v^3)$, and the rung bounds hold at $g^\varepsilon$ as at $\bar g_0$.
\end{proof}

By
Lemma~\ref{lem:dyn-real-hypotheses-ball}\,(\ref{it:dyn-real-hyp-spatial}), \eqref{eq:dyn-real-background-norms} holds, and by
Lemma~\ref{lem:dyn-real-hypotheses-ball}\,(\ref{it:dyn-real-hyp-basic}), $g\in\cH_{\kin,\sigma}(\cD_T\times\bR_v^3)$ and $g(0)=z_0$.  In particular, the hypotheses of
Proposition~\ref{prop:dyn-fixedN-data}\,(\ref{it:dyn-real-extended}) hold.

The incoming data satisfy the following bounds uniformly in $T$.

\begin{lemma}
	\label{lem:dyn-strong-wall-datum}
	Let $H_s:=(f_\st)_{\gamma^-}$.  On $(-1,0)$, with $H_s$ viewed as time independent, for $0\le L\le\Lambda_0+B^{\dy}_{\rm in}+3$,
	\[
	\bigl\|w^L(\check{E})_{\gamma^-}\bigr\|_{(C^{0,1}_{t,x}\cap C^{1,1}_v)((-1,0)\times\gamma^-)}
	+\bigl\|w^L(\check{E})_{\gamma^-}\bigr\|_{L^2((-1,0)\times\gamma^-)}
	+\bigl\|w^LH_s\bigr\|_{L^2((-1,0)\times\gamma^-)}
	\le C_L\bigl(\delta+\mathsf I\bigr) .
	\]
	Let moreover $T\ge1$, let $u\in\cH_{\kin,\sigma}(\cD_T\times\bR_v^3)$, let
	$\cJ\in\{\mathbf I,J_N\}$ with $N\in\bN$, and let
	\[
	H_u:=P_\cJ(t)u_{\gamma^+}+b_w
	\qquad\text{on }\gamma^-_T,
	\]
	with $P_\cJ(t)$ the operator \eqref{eq:wall-contraction-operator} at $\Theta_w(t)$ and $b_w$ of
	\eqref{eq:dyn-theta-residual}.  Then, for every finite $L$,
	\begin{equation}
		\label{eq:dyn-strong-wall-C11}
		\bigl\|w^{L}H_u\bigr\|_{(C^{0,1}_{t,x}\cap C^{1,1}_v)(\gamma^-_T)}
		\le C_L\bigl(\cN_{\mathrm M}(T)[u]+\delta\bigr),
	\end{equation}
	and
	\begin{equation}
		\label{eq:dyn-strong-wall-L2}
		\bigl\|w^{L}H_u\bigr\|_{L^2(\gamma^-_T)}
		\le C_L\bigl(\mathsf U_0(T)[u]+\delta\bigr),
	\end{equation}
	with $C_L$ depending only on the structural data and $L$ (in particular independent of $T$ and $N$).
\end{lemma}

\begin{proof}
	For the negative-time estimate, fix $0\le L\le\Lambda_0+B^{\dy}_{\rm in}+3$.  Then
	\[
	\|w^L(\check{E})_{\gamma^-}\|_{(C^{0,1}_{t,x}\cap C^{1,1}_v)((-1,0)\times\gamma^-)}+\|w^L(\check{E})_{\gamma^-}\|_{L^2((-1,0)\times\gamma^-)}\le C\mathsf I
	\]
	by
	Lemma~\ref{lem:dynamic-zero-initial-source-inflow-lift}\,(\ref{it:dyn-breg-lift-bounds}).  Since $\||v|^{1/2}w^{-3}\|_{L^2_v}<\infty$ and $|\partial\cD|<\infty$, \eqref{eq:stationary-incoming-regularity} at the weight $L+3$ gives
	\[
	\|w^LH_s\|_{L^2((-1,0)\times\gamma^-)}\le C\|w^{L+3}H_s\|_{L^\infty(\gamma^-)}\le C_L\delta .
	\]
	
	To prove \eqref{eq:dyn-strong-wall-C11}, write $\phi_u=\cA_0u_{\gamma^+}$. It suffices to consider $\cN_{\mathrm M}(T)[u]<\infty$; the jets then form the distributional chain required by Lemma~\ref{lem:dynamic-boundary-moments}, and that lemma, \eqref{eq:dyn-unweighted-row-dominated} and $\mathrm M\ge3$ give
	\[
	\|\phi_u\|_{H^3((0,T)\times\partial\cD)}
	\le C\cN_{\mathrm M}(T)[u].
	\]
	The operator $J_N$ acts only on the spatial boundary variable, commutes with $\partial_t$, and is contractive in every spectral $H^s(\partial\cD)$ norm, $0\le s\le3$. Applying this to $\partial_t^j\phi_u$ at $s=3-j$ and using the fixed-chart equivalence of spatial Sobolev norms shows that both $\phi_u$ and $J_N\phi_u$ satisfy the same spacetime $H^3$ bound, with a constant independent of $N$ and $T$. Cover $(0,T)$ by intervals of length between $1$ and $2$ with bounded overlap. The embedding $H^3(J\times\partial\cD)\hookrightarrow W^{1,\infty}(J\times\partial\cD)$ is uniform over these intervals, and hence
	\[
	\|\cJ\phi_u\|_{W^{1,\infty}((0,T)\times\partial\cD)}
	\le C\cN_{\mathrm M}(T)[u],
	\qquad \cJ\in\{\mathbf I,J_N\}.
	\]
	The bound for $\cJ\phi_u$, the identity $H_u=\mathfrak m_{\Theta_w(t)}\cJ\phi_u+b_w$, and the weighted wall-profile bounds of Lemmas~\ref{lem:wall-profile-derivatives} and~\ref{lem:dyn-theta-profile-time-derivatives}, together with \eqref{eq:dyn-theta-defect-relaxation-bound} and \eqref{eq:boundary-moment-holder}, give \eqref{eq:dyn-strong-wall-C11}. The constants depend on $L$ and the structural data, but not on $T$ or $N$.
	
	For
	\eqref{eq:dyn-strong-wall-L2}, the first term of $H_u$ satisfies
	\[
	\|w^LP_\cJ(t)u_{\gamma^+}\|_{L^2(\gamma^-_T)}\le C_L\|\phi_u\|_{L^2((0,T)\times\partial\cD)},
	\]
	by
	\eqref{eq:dynamic-spectral-profile-L2-isometry-bound} and $0\le\cJ\le\mathbf I$. By
	\eqref{eq:outgoing-moment-flux-bound} and \eqref{260902137},
	\[
	\|\phi_u\|_{L^2((0,T)\times\partial\cD)}\le C\|u_{\gamma^+}\|_{L^2(\gamma^+_T)}\le C\mathsf U_0(T)[u].
	\]
	By \eqref{eq:dyn-theta-defect-relaxation-bound} and \eqref{eq:boundary-moment-holder}, we have $\|w^Lb_w\|_{L^2(\gamma^-_T)}\le C_L\delta$.  Combining these bounds, we obtain \eqref{eq:dyn-strong-wall-L2}.
\end{proof}

For $g$ as in Lemma~\ref{lem:dyn-regbg-ball}, with $\hat g^\varepsilon$ the extension of $g^\varepsilon$ by zero to $t\le-1$, and
$\hat f_\pt$ on $(-\infty,T)\times\cD\times\bR_v^3$, the source $\mathsf F$ is defined by
\begin{equation}
	\label{eq:dyn-strong-quasilinear-regrouping}
	(\partial_t+v\cdot\nabla_x)\hat f_\pt+L_{\hat g^\varepsilon}\hat f_\pt=\mathsf F .
\end{equation}

Lemma~\ref{260527559_dy}, applied about $f_\st$, bounds the weighted graph norm of $f_\dy$ and the strong coordinates of $\hat f_\pt$ in terms of the source and incoming data.

\begin{prop}[Regularity]
	\label{prop:dyn-strong}
	Let $T\ge1$ and $0<\varepsilon\le\varepsilon_\lambda$, let $g$ satisfy \eqref{eq:dyn-ball-jets} and
	\eqref{eq:dyn-strong-coefficient-choice}, and let $\Lambda\ge0$, $\mathsf N\ge0$ and
	$\mathsf S\ge0$.  Let $\hat f_\pt\in\cH_{\kin,\sigma}(\cD_{T,\ast}\times\bR_v^3)$ vanish
	for $t\le-3/4$ and solve
	\eqref{eq:dyn-strong-quasilinear-regrouping} with a source $\mathsf F$ and
	$(\hat f_\pt)_{\gamma^-}=\hat H$, and assume that, with
	$H_s=(f_\st)_{\gamma^-}$, for $0\le L\le\Lambda+B^{\dy}_{\rm in}+1$,
	\begin{equation}
		\label{eq:dyn-strong-chain-datum}
		\bigl\|w^{L}\hat H\bigr\|_{(C^{0,1}_{t,x}\cap C^{1,1}_v)(\gamma^-_{T,\ast})}
		+\bigl\|w^{L}\bigl(\hat H-H_s\bigr)\bigr\|_{L^2(\gamma^-_{T,\ast})}\le\mathsf N ,
	\end{equation}
	\begin{equation}
		\label{eq:dyn-strong-chain-source}
		\bigl\|w^{L}\mathsf F\bigr\|_{L^\infty}
		+\bigl\|w^{L}\mathsf F\bigr\|_{L^2(\cD_{T,\ast};H^{-1}_\sigma)}\le\mathsf S .
	\end{equation}
	Then, with $f_\dy:=\hat f_\pt-f_\st$ on $I_T^-$ and all norms over $\cD_{T,\ast}\times\bR_v^3$,
	\begin{align}
		&\bigl\|w^{\Lambda}f_\dy\bigr\|_{\cH_{\kin,\sigma}(\cD_{T,\ast}\times\bR_v^3)}
		+\bigl\|w^{\Lambda}\hat f_\pt\bigr\|_{L^\infty}
		+\bigl[w^{\Lambda}\hat f_\pt\bigr]_{C^{\alpha_{\rm bs}}_{\kin}}
		+\bigl\|w^{\Lambda}\nabla_v\hat f_\pt\bigr\|_{L^\infty_tL^{p_\sharp}_xL^2_v}
		\notag\\
		&\qquad\le C\bigl(\delta+\mathsf N+\mathsf S\bigr),
		\label{eq:dyn-strong-chain-output}
	\end{align}
	with $C$ depending only on the structural data and $\Lambda$ (in particular independent of $T$).
\end{prop}

\begin{proof}
	Let $a=\hat g^\varepsilon$, $u=\hat f_\pt-f_\st$, and
	$F_\st=N[a-f_\st,f_\st]$.  The stationary equation
	\eqref{eq:stationary-prescribed-inflow} gives
	\[
	(\partial_t+v\cdot\nabla_x)u+L_a u=\mathsf F+F_\st,
	\qquad u(-1)=-f_\st,\qquad u_{\gamma^-}=\hat H-H_s.
	\]
	The initial and incoming data are compatible: $\hat f_\pt$ and its
	incoming trace vanish near $t=-1$, whereas $f_\st$ is stationary with
	incoming trace $H_s$.
	
	We verify the stationary terms needed in the inflow argument.  Put
	\[
	\mathcal G_T=\|\mu^{1/8}(a-f_\st)\|_{\cH_{\kin,\sigma}(\cD_{T,\ast}\times\bR_v^3)}.
	\]
	The multiplier estimate \eqref{eq:Hsigma-multiplier} and
	\eqref{eq:dyn-regbg-graph} give $\mathcal G_T\lesssim1$, uniformly in
	$T,N,\varepsilon$.  For every fixed weight $m\ge0$, the Coulomb
	bilinear estimate used in \eqref{eq:dynamic-stationary-source-bound} yields
	\[
	\begin{aligned}
		\|w^mF_\st\|_{L^2(\cD_{T,\ast};H_\sigma^{-1})}
		&\lesssim_m\mathcal G_T\bigl(
		\|w^{m+5}f_\st\|_{L^\infty_{x,v}}+
		\|w^{m+2}\nabla_v f_\st\|_{L_x^{30}L_v^2}\bigr)
		\\&\lesssim_m\delta\mathcal G_T.
	\end{aligned}
	\]
	Here the last inequality is precisely
	\eqref{eq:stationary-finite-weight-graph}--\eqref{eq:stationary-finite-weight-pointwise}.  Those bounds, together
	with \eqref{eq:stationary-incoming-regularity}, also give
	\[
	\|w^m f_\st\|_{L^2_{x,v}}+
	\|w^m H_s\|_{(C_x^{0,1}\cap C_v^{1,1})(\gamma^-)}
	\lesssim_m\delta.
	\]
	In particular, on every time interval of length at most four,
	\(\|w^m f_\st\|_{L^2_{t,x,v}}\lesssim_m\delta\), with a constant
	independent of the position of the interval.
	
	The coefficient bounds \eqref{eq:dyn-strong-coefficient-output-reg}
	ensure \(\|\mu^{1/4}a\|_{L^\infty}\le\varepsilon_0(\cD)\)
	and \(\|\mu^{1/4}a\|_{C^{\alpha_{\rm bs}}_{\kin}}\le2\).
	Thus the linear energy, local pointwise, kinetic H\"older and gradient
	theorems used in the proof of Lemma~\ref{260527559_dy} apply with
	constants uniform in $T,N,\varepsilon$.  The stationary input required here consists of the preceding estimates. In particular,
	Theorem~\ref{thm:weighted-steady}\,(1), applied
	to the displayed equation for $u$, controls its weighted graph norm
	by the source, its initial trace and $\hat H-H_s$.  For the pointwise
	and kinetic H\"older estimates, lift the zero initial datum and the
	incoming datum $\hat H$ as in
	\eqref{eq:dynamic-lift-derivative-bound}.  The local $L^2$ input on
	each translated time interval is bounded by the graph estimate for $u$,
	the lift, and the preceding stationary bound on intervals of length
	four.  Thus the local estimates of Theorems~\ref{260504554} and
	\ref{260526316} have constants independent of $T$.  Finally,
	Theorem~\ref{260527139} gives the mixed velocity-gradient estimate:
	its additional coefficient assumption holds with
	$\alpha=\alpha_{\rm bs}$ and $K=2$ by
	\eqref{eq:dyn-strong-coefficient-output-reg}.  The stationary inputs
	in each of these applications are exactly the bounds displayed above;
	no smallness hypothesis on the boundary datum from
	Theorem~\ref{Landau_st_exi} is invoked.
	
	At $m=\Lambda+B^{\dy}_{\rm in}$, the data requirements of these
	estimates are satisfied by \eqref{eq:dyn-strong-chain-datum} and
	\eqref{eq:dyn-strong-chain-source}, since they cover all weights up to
	$m+1$.  The stationary contribution is at most
	$C_\Lambda\delta(1+\mathcal G_T)\lesssim_\Lambda\delta$.
	The resulting weighted graph, pointwise, H\"older and velocity-gradient
	bounds therefore give \eqref{eq:dyn-strong-chain-output}.
\end{proof}

\subsubsection{The solution and its temporal jets}\label{ss:dyn-proof-base}

\begin{proof}[\textbf{Proof of Proposition~\ref{prop:dyn-fixedN-realization}\,(\ref{it:dyn-real-solution})}]
	By \eqref{eq:dyn-real-coefficient-hypotheses} and Lemma~\ref{lem:dyn-wall-extended}\,(\ref{it:dyn-wall-ext-results}), Proposition~\ref{prop:dynamic-inhomogeneous-physical-solvability}\,(\ref{it:solv-spectral}) applies on $I_T^-$ with $g^\varepsilon$ in place of $g$ and with the data
	\[
	F=\bar F_0,
	\qquad
	b=\bar b_N ,
	\]
	and zero initial value,
	where $\bar F_0\in L^2(I_T^-;L_x^2H_\sigma^{-1})$ by
	\eqref{eq:dynamic-zero-initial-extended-source-spaces} at $j=0$.  It gives existence and uniqueness of the solution $f_\pt^{\varepsilon,N}$ of \eqref{eq:dyn-real-mollified-problem}, and
	\eqref{eq:dynamic-spectral-solvability-estimate} gives
	\[
	\bigl\|f_\pt^{\varepsilon,N}\bigr\|_{\cH_{\kin,\sigma}(\cD_{T,\ast}\times\bR_v^3)}
	\le C_T\Bigl(\|\bar F_0\|_{L^2(I_T^-;L_x^2H_\sigma^{-1})}
	+\|\bar b_N\|_{L^2(\gamma_{T,\ast}^-)}\Bigr),
	\]
	with $C_T$ independent of $N$ and $\varepsilon$.  With the bound on
	$\sup_N\|\bar b_N\|_{L^2(\gamma^-_{T,\ast})}$ below, this gives
	\eqref{eq:dyn-real-uepsN-graph}.  On
	$(-1,-3/4]$ all the data of \eqref{eq:dyn-real-mollified-problem} vanish, by
	\eqref{eq:dyn-real-increment-split}, \eqref{eq:dynamic-base-background} and
	Proposition~\ref{prop:dyn-fixedN-data}\,(\ref{it:dyn-real-extended}).  The uniqueness clause of
	Proposition~\ref{prop:dynamic-inhomogeneous-physical-solvability}\,(\ref{it:solv-spectral}), applied on that
	subinterval with zero data, therefore gives $f_\pt^{\varepsilon,N}=0$ on $(-1,-3/4]$. The required bound $\sup_{N\in\bN}\|\bar b_N\|_{L^2(\gamma^-_{T,\ast})}<\infty$ follows from the two time intervals.
	On $(0,T)$, $\bar b_N=h_w+\mathfrak d_N+c_N$ by \eqref{eq:dyn-limits-corner-physical}, with $h_w\in L^2(\gamma^-_T)$.  Moreover, $\|\mathfrak d_N\|_{L^2(\gamma_T^-)}\le CT^{1/2}\delta$ by \eqref{eq:dyn-limits-steady-residual-size} at $L=0$ and $\|c_N\|_{L^2(\gamma_T^-)}\le C\varsigma_N$ by \eqref{eq:dyn-limits-corner-bound} at $L=j=0$, with $C$ independent of $N$, and $\varsigma_N$ is bounded in $N$, since $\varsigma_N\to0$ by Lemma~\ref{lem:dyn-limits-wall}\,(\ref{it:dyn-limits-wall-size}).  On $(-1,0)$, $\bar b_N$ is the combination
	\eqref{eq:dynamic-negative-time-spectral-inflow} of the traces of $\check{E}$, which satisfy
	$\|(\check{E})_{\gamma^\pm}\|_{L^2((-1,0)\times\gamma^\pm)}\le C\mathsf I$ by \eqref{eq:dyn-lift-size-relations}.  Moreover, $\|P_{w,N}(t)q\|_{L^2(\gamma^-)}\le C\|q\|_{L^2(\gamma^+)}$ for $-1<t<0$ and $q\in L^2(\gamma^+)$, by
	\eqref{eq:wall-J-operator-bound} at $\cJ=J_N$, with $C$ independent of $N$, as asserted.
\end{proof}

\begin{proof}[\textbf{Proof of Proposition~\ref{prop:dyn-fixedN-realization}\,(\ref{it:dyn-real-prop-chain})--(\ref{it:dyn-real-prop-initial})}]
	The first identity of \eqref{eq:dyn-real-prop-chain} is the chain of the $k_j$ verified below.  The rest of (\ref{it:dyn-real-prop-chain}), (\ref{it:dyn-real-prop-equations}) and
	(\ref{it:dyn-real-prop-boundary}) of this proposition is
	Proposition~\ref{prop:dynamic-fixedN-zero-initial-temporal-realization}, applied with
	\[
	(k_j,S_j,b,Z_0)=\bigl(g_j^\varepsilon,\bar F_j,\bar b_N,f_\pt^{\varepsilon,N}\bigr);
	\]
	its $Z_j$ is $\mathfrak D_{\mathrm t}^{\,j}f_\pt^{\varepsilon,N}$ by \eqref{eq:dynamic-zero-initial-temporal-chain}.  We verify its hypotheses.
	\begin{itemize}
		\item The chain of the $k_j$ in \eqref{eq:dynamic-zero-initial-abstract-data} follows from
		\eqref{eq:dyn-real-mollified-objects}, term by term.  In the first term, $\hat g^{(0)}$ is smooth
		in $t$ by Lemma~\ref{lem:dyn-fp-seed}.  To the second,
		Proposition~\ref{prop:dyn-fixedN-mollification}\,(\ref{it:dyn-real-moll-chain}) applies at
		$z=\hat g-\hat g^{(0)}=\bar g_0-\hat g^{(0)}$, which belongs to $\cH_{\kin,\sigma}(\cD_{T,\ast}\times\bR_v^3)$ by
		\eqref{eq:dynamic-zero-initial-extended-coefficient-graphs} and
		Lemma~\ref{lem:dyn-fp-seed}, and vanishes on $(-1,0]$ by
		\eqref{eq:dyn-strong-past-extension}.  The chain of the $S_j$ is
		\eqref{eq:dynamic-zero-initial-extended-source-chain}.
		\item The sources satisfy $\bar F_j\in L^2(I_T^-;L_x^2H_\sigma^{-1})$ by
		\eqref{eq:dynamic-zero-initial-extended-source-spaces}, and
		$\bar b_N\in H^{\mathrm M}(I_T^-;L^2(\gamma^-))$ by
		\eqref{eq:dynamic-zero-initial-extended-inflow-H3}.
		\item The coefficients satisfy $w^{-2}g^\varepsilon_j\in L^\infty(I_T^-;L_x^\infty L_v^2)$, $0\le j\le\mathrm M$:
		by the second identity of \eqref{eq:dyn-real-increment-split} and Minkowski's inequality in
		$t$,
		\[
		\|w^{-2}g^\varepsilon_j\|_{L^\infty(I_T^-;L_x^\infty L_v^2)}\le\|\varrho_\varepsilon^{(j)}\|_{L^1}
		\|w^{-2}\bar g_0\|_{L^\infty(I_T^-;L_x^\infty L_v^2)}+2\|w^{-2}\partial_t^j\hat g^{(0)}\|_{L^\infty(I_T^-;L_x^\infty L_v^2)} .
		\]
		Here $\|\varrho_\varepsilon^{(j)}\|_{L^1}<\infty$ and $\|w^{-2}\bar g_0\|_{L^\infty(I_T^-;L_x^\infty L_v^2)}<\infty$, the latter
		by \eqref{eq:dyn-real-background-norms}.  For the second term, by
		Lemma~\ref{lem:dyn-fp-seed}, $\partial_t^j\hat g^{(0)}$ is on $I_T^-$ a finite combination
		of $z_0,\ldots,z_{\mathrm M}$ and $f_\st$ with bounded coefficients, with
		$\|w^{-2}z_i\|_{L_x^\infty L_v^2}\le C\mathsf I$ by \eqref{eq:dyn-data-jets}, and
		$w^{-2}f_\st\in L_x^\infty L_v^2$ by \eqref{eq:base-fiber-Lsharp} and \eqref{eq:dyn-theta-background}.  Moreover
		$\|\mu^{1/4}g^\varepsilon_0\|_{L^\infty}\le\varepsilon_{\rm dyn}$ by
		\eqref{eq:dyn-real-coefficient-hypotheses}.
		\item The function $f_\pt^{\varepsilon,N}\in\cH_{\kin,\sigma}(\cD_{T,\ast}\times\bR_v^3)$ of
		(\ref{it:dyn-real-solution}) solves \eqref{eq:dyn-real-mollified-problem}, which is the case $j=0$ of \eqref{eq:dynamic-zero-initial-temporal-Z0} and
		\eqref{eq:dynamic-zero-initial-temporal-boundary}, and
		$f_\pt^{\varepsilon,N}=0$ on $(-1,-3/4]$ by (\ref{it:dyn-real-solution}).
	\end{itemize}
	
	For (\ref{it:dyn-real-prop-initial}), $\check{E}$ solves
	\eqref{eq:dyn-real-mollified-problem} on $(-1,0)$.  Indeed, on $(-1,0)$,
	$g^\varepsilon=\bar g_0=\check{E}$ by \eqref{eq:dyn-real-increment-split} and
	\eqref{eq:dynamic-base-background}.  Hence $L_{g^\varepsilon}\check{E}=L\check{E}-N[\check{E},\check{E}]$
	by \eqref{def:L_g}, and \eqref{eq:dynamic-negative-time-source} and
	\eqref{eq:dynamic-zero-initial-extended-sources} give
	\[
	(\partial_t+v\cdot\nabla_x)\check{E}+L_{g^\varepsilon}\check{E}=\check{F}=\bar F_0
	\qquad\text{on }(-1,0).
	\]
	The boundary condition
	$(\check{E})_{\gamma^-}=P_{w,N}(\check{E})_{\gamma^+}+\bar b_N$ on $(-1,0)$ is
	\eqref{eq:dynamic-negative-time-spectral-inflow} with
	\eqref{eq:dynamic-zero-initial-extended-spectral-inflow}, and the initial condition
	$\check{E}(-1)=0$ is \eqref{eq:dynamic-negative-time-zero-near-minus-one}.  Also
	$\check{E}\in\cH_{\kin,\sigma}((-1,0)\times\cD\times\bR_v^3)$ by
	Lemma~\ref{lem:dynamic-zero-initial-source-inflow-lift}\,(\ref{it:dyn-lift-basic}).  Since the
	restriction of $f_\pt^{\varepsilon,N}$ to $(-1,0)$ solves the same problem, the uniqueness
	clause of
	Proposition~\ref{prop:dynamic-inhomogeneous-physical-solvability}\,(\ref{it:solv-spectral}),
	applied on $(-1,0)$, gives
	$f_\pt^{\varepsilon,N}=\check{E}$ on $(-1,0)$.  By \eqref{eq:dyn-real-prop-chain},
	$\mathfrak D_{\mathrm t}^{\,j+1}f_\pt^{\varepsilon,N}=\partial_t\mathfrak D_{\mathrm t}^{\,j}f_\pt^{\varepsilon,N}$
	in $\cD'((-1,0)\times\cD\times\bR^3)$, and induction on $j$ gives the first identity of
	\eqref{eq:dyn-real-prop-initial}.  Each $\mathfrak D_{\mathrm t}^{\,j}f_\pt^{\varepsilon,N}$
	belongs to $\cH_{\kin,\sigma}(\cD_{T,\ast}\times\bR_v^3)$ by
	Proposition~\ref{prop:dynamic-fixedN-zero-initial-temporal-realization}, and
	$\partial_t^j\check{E}$ is continuous on $[-1,0]$ with values in $L^2_{x,v}$.  Hence both sides of the first identity of \eqref{eq:dyn-real-prop-initial} are continuous on $[-1,0]$ with values in $L^2_{x,v}$, the left side by the continuous embedding of Proposition~\ref{2607151218}, so they agree at $t=0$, and \eqref{eq:dynamic-negative-time-jet-matching} gives $\mathfrak D_{\mathrm t}^{\,j}f_\pt^{\varepsilon,N}(0)=\partial_t^j\check{E}(0)=z_j$.  This proves the second identity of \eqref{eq:dyn-real-prop-initial}, and completes the proof.
\end{proof}

\begin{proof}[\textbf{Proof of Proposition~\ref{prop:dyn-fixedN-realization}\,(\ref{it:dyn-real-mass})}]
	Write $f:=f_\pt^{\varepsilon,N}$ and $\mathfrak m(t):=\mathfrak m_f(t)$.  Move
	$N[g^\varepsilon,f]$ to the right side of \eqref{eq:dyn-real-mollified-problem}; the new
	source belongs to $L^2(I_T^-;L_x^2H_\sigma^{-1})$ by \eqref{eq:velalg-poly} at $(\ell,\mathrm L)=(0,-2)$, since
	$w^{-2}g^\varepsilon\in L^\infty(I_T^-;L_x^\infty L_v^2)$ by the proof of (\ref{it:dyn-real-prop-chain})--(\ref{it:dyn-real-prop-initial}) of this proposition at $j=0$.  The
	$\sqrt\mu$-moment of $N[g^\varepsilon,f]$ vanishes, $\langle N[g^\varepsilon,f],\sqrt\mu\rangle_v=0$, by \eqref{eq:N-mass-kernel}.  By \eqref{eq:dyn-limits-wall-flux}, the wall operator $P_{w,N}(t)$
	preserves the total flux:
	\[
	\int_{\gamma^-}(P_{w,N}(t)f_{\gamma^+})\sqrt\mu\dd\gamma=\int_{\gamma^+}f_{\gamma^+}\sqrt\mu\dd\gamma.
	\]
	So the computation of
	Lemma~\ref{lem:dynamic-alpha-mass-scalar} at $\alpha=1$ and $\cJ=J_N$, carried out at a.e.\ $t$,
	gives
	\eqref{eq:dynamic-physical-mass-balance} on $(0,t)$, that is,
	\[
	\mathfrak m(t)
	=\mathfrak m(0)+\int_0^{t}\bigl\langle\bar F_0(r),\sqrt\mu\bigr\rangle_{x,v}\dd r
	+\int_0^{t}\int_{\partial\cD}\beta_{\bar b_N}(r,x)\dd S_x\dd r ,
	\qquad 0\le t\le T .
	\]
	For $0<r<T$, $\bar F_0(r)=0$ by \eqref{eq:dynamic-zero-initial-extended-sources} and
	$\int_{\partial\cD}\beta_{\bar b_N}(r)\dd S_x=0$ by \eqref{eq:dyn-limits-corner-massfree}.
	Hence $\mathfrak m(t)=\mathfrak m(0)=\mathfrak m_{z_0}$ for $0\le t\le T$, the last by
	(\ref{it:dyn-real-prop-initial}) at $j=0$.  Finally,
	$\mathfrak m_{f_\dy^{\varepsilon,N}}=\mathfrak m-\mathfrak m_{f_\st}$ by \eqref{eq:dyn-real-centered-unknown}, and
	$\mathfrak m_{z_0}-\mathfrak m_{f_\st}=\mathfrak m_{f_{\dy,0}}=0$ by
	\eqref{eq:dynamic-initial-zero-perturbation-mass}.  Combining these identities, we obtain \eqref{eq:dyn-real-prop-mass}.
\end{proof}

\subsubsection{The base level of the ladder}\label{app:dyn-boundary-regularization}

Throughout \S\S\ref{app:dyn-boundary-regularization}--\ref{ss:dyn-fixedN-consequences}, $T$, $N$, $\varepsilon$ and $g$ are as in Proposition~\ref{prop:dyn-fixedN-realization}, and its hypotheses are in force.
We derive the weighted base estimates on $I_T^-=(-1,T)$ for the full perturbation variable $f_\pt^{\varepsilon,N}$.
The argument uses the evolutionary background bounds of Lemma~\ref{lem:dyn-breg-background} in place of the stationary estimates \eqref{eq:base-fiber-Lsharp} and \eqref{eq:Dg-rungs}. Set
\begin{equation}
	\label{eq:dyn-breg-unknown}
	Z_j:=\mathfrak D_{\mathrm t}^{\,j}f_\pt^{\varepsilon,N}\ \text{ on }I_T^-,
	\quad 0\le j\le\mathrm M,
	\qquad
	F:=\bar F_0 ,
\end{equation}
so that, by Proposition~\ref{prop:dyn-fixedN-realization}\,(\ref{it:dyn-real-solution}),~(\ref{it:dyn-real-prop-chain}), $Z_0=f_\pt^{\varepsilon,N}$ solves \eqref{eq:dyn-real-mollified-problem} and the $Z_j$ form the chain \eqref{eq:dyn-real-prop-chain}.
The truncated incoming datum $H_N$ is defined through $\varphi_N$ and $\Phi_N$ by \begin{equation} \label{eq:dynamic-fixedN-strong-incoming}
	\varphi_N:=\cA_0(Z_0)_{\gamma^+},
	\qquad
	\Phi_N:=J_N\varphi_N,
	\qquad
	H_N:=P_{w,N}(Z_0)_{\gamma^+}+\bar b_N
	=(\mathfrak m_\circ+d_\delta)\Phi_N+\bar b_N ,
\end{equation}
so that $(Z_0)_{\gamma^-}=H_N$ on $\gamma^-_{T,\ast}$.  The initial
value of a jet is its value at the left endpoint $t=-1$.

Throughout this subsection, fixed-cutoff regularity constants may depend on $(T,N,\varepsilon)$ and the fixed data. Parameter-independent bounds are identified explicitly.

By \eqref{eq:dynamic-zero-initial-extended-sources} and \eqref{eq:dynamic-zero-initial-extended-spectral-inflow},
$\bar F_j=\mathbf 1_{(-1,0)}\partial_t^{\,j}\check{F}$ and $\bar b_N=\check{b}_N$ on $(-1,0)$.

The incoming datum satisfies the following estimates.

\begin{lemma}[Incoming flow bounds]
	\label{lem:dynamic-fixedN-incoming-flow-bounds}
	Let $H_N$ be \eqref{eq:dynamic-fixedN-strong-incoming}.
	\begin{enumerate}
		\item\label{it:dyn-incoming-JN}
		The scalar $\Phi_N$ satisfies
		\begin{equation}
			\Phi_N\in H^{\mathrm M}\bigl(I_T^-;H^{\mathrm M+2}(\partial\cD)\bigr) .
			\label{eq:dynamic-fixedN-incoming-JN-H5}
		\end{equation}
		\item\label{it:dyn-incoming-words}
		Let $0\le L\le\Lambda_0+B^{\dy}_{\rm in}+3$.  Then
		\begin{equation}
			\label{eq:dynamic-fixedN-incoming-strong-bound}
			\|w^LH_N\|_{(C_{t,x}^{0,1}\cap C_v^{1,1})(\gamma^-_{T,\ast})}<\infty ,
		\end{equation}
		and, for every ordered word $I$ over $\mathfrak A_\bullet$ with
		$|I|\le\mathrm M$, where $\mathfrak D^-_{\mathrm t}=\partial_t$ as after
		\eqref{eq:boundary-word-Leibniz},
		\begin{equation}
			\|w^L \mathfrak D_I^-H_N\|_{L^2(\gamma^-_{T,\ast})}<\infty .
			\label{eq:dynamic-fixedN-incoming-word-bound}
		\end{equation}
		\item\label{it:dyn-incoming-differences}
		There is a window $h_0$, depending only on the structural data, with $0<h_0\le s_0/(2\mathrm M)$, where $s_0$ is that of
		\eqref{eq:complete-flow-uniform-bounds}, such that, for every ordered word $I$ over
		$\mathfrak A_\bullet$, every
		$X=\tau_\nu$ with $\nu\in\mathfrak A$, and every integer $m\ge1$ with
		$|I|+m\le\mathrm M$,
		\begin{align}
			\sup_{0<|h|<h_0}
			\|(\Delta_h^{X,\gamma_-})^m
			\mathfrak D_I^-H_N\|_{L^2(\gamma^-_{T,\ast})}
			&<\infty,
			\label{eq:dynamic-fixedN-incoming-difference-bound}\\
			(\Delta_h^{X,\gamma_-})^m
			\mathfrak D_I^-H_N
			&\longrightarrow
			(\mathfrak D_X^-)^m
			\mathfrak D_I^-H_N
			\quad\text{in }L^2(\gamma^-_{T,\ast})
			\qquad(h\to0).
			\label{eq:dynamic-fixedN-incoming-difference-limit}
		\end{align}
	\end{enumerate}
\end{lemma}

\begin{proof}
	\emph{(\ref{it:dyn-incoming-JN}).}  The chain \eqref{eq:dyn-real-prop-chain} and
	Proposition~\ref{prop:weak-boundary-generator}\,(\ref{it:weak-generator-evolutionary}), applied on $I_T^-$ at $m=\mathrm M$, give
	$\partial_t^{\,j}\varphi_N=\cA_0(Z_j)_{\gamma^+}$, with
	\[
	\|\partial_t^{\,j}\varphi_N\|_{L^2(I_T^-\times\partial\cD)}
	\le C\|Z_j\|_{\cH_{\kin,\sigma}(\cD_{T,\ast}\times\bR_v^3)} ,
	\qquad 0\le j\le\mathrm M .
	\]
	Then \eqref{eq:fixedN-Hr} with $r=\mathrm M+2$ gives \eqref{eq:dynamic-fixedN-incoming-JN-H5}.
	
	\emph{(\ref{it:dyn-incoming-words}).}  Since $\partial_t$ commutes with the spatial letters on $\gamma^-_{T,\ast}$, we move the temporal letters of each word over $\mathfrak A_\bullet$ to the right of its spatial letters.
	By \eqref{eq:dynamic-trace-time-chain-boundary} of Proposition~\ref{prop:weak-boundary-generator}\,(\ref{it:weak-generator-evolutionary}), the distributional derivatives of $H_N$ satisfy
	$\partial_t^jH_N=(Z_j)_{\gamma^-}\in L^2(\gamma^-_{T,\ast})$ for $0\le j\le\mathrm M$.
	Thus $\partial_t^jH_N\in H^1(I_T^-;L^2(\gamma^-))$ for $0\le j<\mathrm M$.
	In particular, the traces at $t=0$ of the two restrictions agree through order $\mathrm M-1$; no distributional delta term is introduced by gluing their time derivatives.  We may therefore estimate the boundary words separately on $(-1,0)$ and $(0,T)$.
	
	On $(-1,0)$, $Z_0=\check{E}$ by \eqref{eq:dyn-real-prop-initial} at $j=0$, so $H_N=(\check{E})_{\gamma^-}$ by \eqref{eq:dynamic-fixedN-strong-incoming} and \eqref{eq:dynamic-negative-time-spectral-inflow}.  Since $L\le\Lambda_0+B^{\dy}_{\rm in}+3\le\Lambda_0+B^{\dy}_{\rm in}+6$,
	Lemma~\ref{lem:dynamic-zero-initial-source-inflow-lift}\,(\ref{it:dyn-breg-lift-bounds}) gives
	\[
	\sum_{|I|+m\le\mathrm M}\|w^L\mathfrak D_I^-\partial_t^m(\check{E})_{\gamma^-}\|_{L^2((-1,0)\times\gamma^-)}\le C\mathsf I .
	\]
	It also gives
	\[
	\|w^L(\check{E})_{\gamma^-}\|_{(C^{0,1}_{t,x}\cap C^{1,1}_v)((-1,0)\times\gamma^-)}\le C\mathsf I .
	\]
	
	On $(0,T)$, $H_N=(\mathfrak m_\circ+d_\delta)\Phi_N+h_w+\mathfrak d_N+c_N$ by \eqref{eq:dynamic-fixedN-strong-incoming} and \eqref{eq:dyn-limits-corner-physical}.
	The derivative $\mathfrak D_I^-\bigl((\mathfrak m_\circ+d_\delta)\Phi_N\bigr)$ is split by
	\eqref{eq:boundary-word-Leibniz} into the terms $\bigl(\mathfrak D_{I_1}^-(\mathfrak m_\circ+d_\delta)\bigr)\tau_{I_2}\Phi_N$, $I_1\sqcup I_2=I$.
	By \eqref{eq:boundary-profile-flux-bound}, applied at each $t\in (0,T)$ to $a=\partial_t^{\,j(I_1)}(\mathfrak m_\circ+d_\delta)(t)$ and $K=I_1'$, with the bounds of Lemmas~\ref{lem:wall-profile-derivatives} and~\ref{lem:dyn-theta-profile-time-derivatives}, and by \eqref{eq:dynamic-fixedN-incoming-JN-H5},
	\[
	\bigl\|w^L\bigl(\mathfrak D_{I_1}^-(\mathfrak m_\circ+d_\delta)\bigr)\tau_{I_2}\Phi_N\bigr\|_{L^2(\gamma_T^-)}
	\le C\|\tau_{I_2}\Phi_N\|_{L^2((0,T)\times\partial\cD)}
	\le C\|\Phi_N\|_{H^{\mathrm M}(I_T^-;H^{\mathrm M+2}(\partial\cD))} .
	\]
	For a spatial word $I$, $\mathfrak D^-_Ih_w$ is the restriction of $\mathfrak D_Ih_w$ to $\gamma^-$ by \eqref{eq:boundary-complete-derivative}, so the first bound of \eqref{eq:dyn-theta-defect-complete-bound} gives
	\[
	\sum_{m+|I|\le\mathrm M}\|w^L\mathfrak D_I^-\partial_t^mh_w\|_{L^2(\gamma_T^-)}\le C\delta .
	\]
	The words of $\mathfrak d_N$ and of $c_N$ are bounded by \eqref{eq:dyn-limits-steady-residual-size}, integrated over $(0,T)$, and by \eqref{eq:dyn-limits-corner-bound}, with $\varsigma_N<\infty$ by Lemma~\ref{lem:dyn-limits-wall}\,(\ref{it:dyn-limits-wall-size}).  The bounds on $(-1,0)$ and on $(0,T)$ give \eqref{eq:dynamic-fixedN-incoming-word-bound}.
	
	For \eqref{eq:dynamic-fixedN-incoming-strong-bound}, the bound on $(-1,0)$ is the second one above.  On $(0,T)$, since $\mathrm M\ge3$, \eqref{eq:dynamic-fixedN-incoming-JN-H5} gives $\Phi_N\in H^3(I_T^-\times\partial\cD)$, so $\Phi_N\in W^{1,\infty}(I_T^-\times\partial\cD)$ by the Sobolev embedding on the three-dimensional manifold $I_T^-\times\partial\cD$.  Lemmas~\ref{lem:wall-profile-derivatives} and~\ref{lem:dyn-theta-profile-time-derivatives} give
	\[
	\|w^L(\mathfrak m_\circ+d_\delta)\|_{(C^{0,1}_{t,x}\cap C^{1,1}_v)(\gamma_T^-)}+\|w^Lh_w\|_{(C^{0,1}_{t,x}\cap C^{1,1}_v)(\gamma_T^-)}\le C_L ,
	\]
	and, as $\Phi_N$ does not depend on $v$, the weighted Lipschitz norm of $(\mathfrak m_\circ+d_\delta)\Phi_N$ on $\gamma_T^-$ is finite.  Those of $\mathfrak d_N$ and $c_N$ are finite by \eqref{eq:dyn-limits-steady-residual-size} and \eqref{eq:dyn-limits-corner-bound}.  	Finally, $H_N\in H^1(I_T^-;L^2(\gamma^-))$, so the two one-sided traces at $t=0$ agree in $L^2(\gamma^-)$.  Both pieces admit Lipschitz representatives in $(t,x,v)$ with the weighted bounds just obtained.  Their traces therefore agree pointwise on $\gamma^-$, and gluing the two representatives preserves the stated Lipschitz bound.  This proves \eqref{eq:dynamic-fixedN-incoming-strong-bound}.
	
	\emph{(\ref{it:dyn-incoming-differences}).}  Fix $X=\tau_\nu$ with $\nu\in\mathfrak A$.
	By Lemma~\ref{lem:boundary-flow-continuity}, \eqref{eq:boundary-flow-bound}
	holds on $|s|\le s_0$, and $T_s^{X,\gamma_-}$ is strongly continuous on $L^2(\gamma^-)$.
	For $0\le j\le\mathrm M$, let
	\[
	\varphi_{j,n}^0:=J_n\varphi_j^0.
	\]
	Since $J_n$ commutes with $(1+A_\partial)^{\mathrm M/2}$, \eqref{eq:boundary-spectral-strong-convergence} applied to $(1+A_\partial)^{\mathrm M/2}\psi$ gives, for every $\psi\in H^{\mathrm M}(\partial\cD)$, $J_n\psi\to\psi$ in the norm
	\eqref{eq:word-boundary-spectral-Hs} at $s=\mathrm M$. As this norm is equivalent to the chart norm by
	Lemma~\ref{lem:word-boundary-spectral-chart}, we obtain $\varphi_{j,n}^0\to\varphi_j^0$ strongly in $H^{\mathrm M}(\partial\cD)$, since $\varphi_j^0\in H^{\mathrm M}(\partial\cD)$ by \eqref{eq:dyn-data-scalar-H3}. Let $c^\pt_{N,n}$ be \eqref{eq:dynamic-spectral-corner-correction} with
	$\varphi_j^0$ replaced by $\varphi_{j,n}^0$, let $\mathfrak d_{N,n}$ be $\mathfrak d_N$ with
	its scalar $(\mathbf I-J_N)\cA_0(f_\st)_{\gamma^+}\in H^{\mathrm N}(\partial\cD)$ replaced by $J_n(\mathbf I-J_N)\cA_0(f_\st)_{\gamma^+}$, and set
	\[
	H_{N,n}:=(\mathfrak m_\circ+d_\delta)\Phi_N+\bar b_{N,n},
	\qquad
	\bar b_{N,n}:=
	\begin{cases}
		\check{b}_N,&-1<t<0,\\
		h_w+c^\pt_{N,n}+(1-\chi_+)\mathfrak d_{N,n},&0<t<T.
	\end{cases}
	\]
	For almost every $t\ne0$, $H_{N,n}$ has the required spatial boundary derivatives; its time derivatives are taken separately on $(-1,0)$ and $(0,T)$.  No time compatibility of the approximants at $t=0$ is required here: the fundamental theorem of calculus below is applied only along spatial boundary flows at almost every fixed time.  On $(0,T)$, this follows from Lemmas~\ref{lem:fixedN-range}, \ref{lem:wall-profile-derivatives}
	and~\ref{lem:dyn-theta-profile-time-derivatives} and from \eqref{eq:dynamic-fixedN-incoming-JN-H5}.  On $(-1,0)$, $H_{N,n}=H_N=(\check{E})_{\gamma^-}$ as in the proof of (\ref{it:dyn-incoming-words}), and $\check{E}$ is a combination, with coefficients depending only on $t$, of the $z_j$ by \eqref{eq:dynamic-negative-time-Taylor-polynomial}--\eqref{eq:dynamic-negative-time-lift}; the derivatives of order at most $\mathrm M+2$ of the $z_j$ are bounded by \eqref{eq:dyn-breg-lift-weights}.  Moreover, for every $\ell\ge0$, every
	ordered word $I$ over $\mathfrak A_\bullet$, and every $m\ge0$ with
	$|I|+m\le\mathrm M$,
	\begin{align}
		&\|w^\ell(\mathfrak D_X^-)^m
		\mathfrak D_I^-(H_{N,n}-H_N)\|_{L^2(\gamma^-_{T,\ast})}
		\notag\\
		&\qquad\le C\Bigl(
		\sum_{j=0}^{\mathrm M}\|\varphi_{j,n}^0-\varphi_j^0\|_{H^{\mathrm M}(\partial\cD)}
		+\|(J_n-\mathbf I)(\mathbf I-J_N)\cA_0(f_\st)_{\gamma^+}\|_{H^{\mathrm M}(\partial\cD)}\Bigr)
		\longrightarrow0
		\qquad(n\to\infty).
		\label{eq:dynamic-fixedN-boundary-smooth-approximation}
	\end{align}
	The inequality follows from \eqref{eq:dynamic-spectral-corner-correction}, \eqref{eq:dyn-limits-steady-residual} and the Leibniz rule
	\eqref{eq:boundary-word-Leibniz}, applied as in the proof of
	\eqref{eq:dyn-limits-corner-bound}, with
	Lemma~\ref{lem:word-boundary-spectral-chart}.  The convergence follows from
	that of $\varphi_{j,n}^0$ above and its analogue
	\[
	J_n(\mathbf I-J_N)\cA_0(f_\st)_{\gamma^+}\to(\mathbf I-J_N)\cA_0(f_\st)_{\gamma^+}
	\quad\text{in }H^{\mathrm M}(\partial\cD) ,
	\] proved in the same way.
	Now fix an ordered word $I$ over $\mathfrak A_\bullet$ and $m\ge1$ with $|I|+m\le\mathrm M$, and let $h_0:=s_0/(2\mathrm M)$.
	For the smooth approximants, the fundamental theorem of calculus, applied $m$ times
	along the flow $T_s^{X,\gamma_-}$ as in the proof of \eqref{eq:fixedN-inflow-repeated-FTC} (with \eqref{eq:complete-flow-level-set}, \eqref{eq:complete-flow-normal-velocity-preservation} and the group law of Lemma~\ref{lem:boundary-flow-continuity}), gives, for a.e. $t$ and $0<|h|<h_0$,
	\begin{equation}
		(\Delta_h^{X,\gamma_-})^m\mathfrak D_I^-H_{N,n}
		=\int_{[0,1]^m}
		T_{(\theta_1+\cdots+\theta_m)h}^{X,\gamma_-}
		(\mathfrak D_X^-)^m\mathfrak D_I^-H_{N,n}
		\,\dd\theta_1\cdots\dd\theta_m.
		\label{eq:dynamic-fixedN-smooth-boundary-FTC}
	\end{equation}
	The flow parameter in \eqref{eq:dynamic-fixedN-smooth-boundary-FTC} satisfies
	\[
	|(\theta_1+\cdots+\theta_m)h|\le m|h|<\mathrm Mh_0=s_0/2.
	\]
	By \eqref{eq:boundary-flow-bound}
	and \eqref{eq:dynamic-fixedN-boundary-smooth-approximation}, we may let
	$n\to\infty$ in \eqref{eq:dynamic-fixedN-smooth-boundary-FTC}.  The limit involves the temporal
	derivatives of $H_N$ taken on $(-1,0)$ and on $(0,T)$ separately; these are its distributional
	ones, since $\partial_t^{\,j}H_N=(Z_j)_{\gamma^-}\in L^2(\gamma^-_{T,\ast})$ by \eqref{eq:dynamic-trace-time-chain-boundary} of Proposition~\ref{prop:weak-boundary-generator}\,(\ref{it:weak-generator-evolutionary}), applied on $I_T^-$ to the chain $Z_0,\ldots,Z_{\mathrm M}$ of \eqref{eq:dyn-real-prop-chain}.  Put
	\[
	G:=(\mathfrak D_X^-)^m\mathfrak D_I^-H_N,
	\]
	so that the limit reads
	\begin{equation}
		(\Delta_h^{X,\gamma_-})^m\mathfrak D_I^-H_N
		=\int_{[0,1]^m}T_{(\theta_1+\cdots+\theta_m)h}^{X,\gamma_-}G
		\,\dd\theta_1\cdots\dd\theta_m.
		\label{eq:dynamic-fixedN-boundary-FTC-H3}
	\end{equation}
	Here the bounds for $G$ are at fixed $T$, $N$ and $\varepsilon$;
	the flow window $h_0=s_0/(2\mathrm M)$ depends only on the
	boundary geometry and the chosen tangential fields.  Consequently,
	by \eqref{eq:dynamic-fixedN-boundary-FTC-H3} and \eqref{eq:boundary-flow-bound},
	\[
	\sup_{0<|h|<h_0}
	\|(\Delta_h^{X,\gamma_-})^m\mathfrak D_I^-H_N\|_{L^2(\gamma^-_{T,\ast})}
	\le C\|G\|_{L^2(\gamma^-_{T,\ast})},
	\]
	and $\|G\|_{L^2(\gamma^-_{T,\ast})}<\infty$ by
	\eqref{eq:dynamic-fixedN-incoming-word-bound} since $m+|I|\le\mathrm M$; this proves \eqref{eq:dynamic-fixedN-incoming-difference-bound}. For \eqref{eq:dynamic-fixedN-incoming-difference-limit}, for almost every $t$ the norm
	$\|T_s^{X,\gamma_-}G(t)-G(t)\|_{L^2(\gamma^-)}$ tends to $0$ by
	Lemma~\ref{lem:boundary-flow-continuity} and is at most $C\|G(t)\|_{L^2(\gamma^-)}$ by
	\eqref{eq:boundary-flow-bound}.
	Dominated convergence in $t$ yields $\|T_s^{X,\gamma_-}G-G\|_{L^2(\gamma^-_{T,\ast})}\to0$ as $s\to0$. Combining this with \eqref{eq:dynamic-fixedN-boundary-FTC-H3}, we obtain
	\[
	\left\|
	(\Delta_h^{X,\gamma_-})^m\mathfrak D_I^-H_N-G
	\right\|_{L^2(\gamma^-_{T,\ast})}
	\le
	\sup_{|s|\le m|h|}
	\|T_s^{X,\gamma_-}G-G\|_{L^2(\gamma^-_{T,\ast})}
	\longrightarrow0,
	\]
	which proves \eqref{eq:dynamic-fixedN-incoming-difference-limit}.
\end{proof}

By Proposition~\ref{prop:dyn-fixedN-realization}\,(\ref{it:dyn-real-solution}),
$Z_0=f_\pt^{\varepsilon,N}=0$ on $(-1,-3/4]$. Since $\mathfrak D_I$ acts locally in $t$, every jet $\mathfrak D_IZ_0$
vanishes there, and its value at $t=-1$ is $0$ whenever
it belongs to $\cH_{\kin,\sigma}(\cD_{T,\ast}\times\bR_v^3)$.  This gives the initial condition $Z_I(-1)=0$ in the following definition.

\begin{defn}[The assertions $\mathsf P^\dy_m$ and $\mathsf P^\bullet_m$]
	\label{def:dyn-package}
	Let $0\le m\le\mathrm M$, let $Z_0$ and $F$ be as in \eqref{eq:dyn-breg-unknown}, let $I$ range over the ordered complete spatial words
	and let $L$ range over the integers $0\le L\le\Lambda_{|I|}$, with the ladder $\Lambda_0,\ldots,\Lambda_{\mathrm M}$ of
	\eqref{eq:dyn-realization-weights}.  Write
	\[
	\mathscr L_{g^\varepsilon}:=\partial_t+v\cdot\nabla_x+L_{g^\varepsilon},
	\qquad
	\mathfrak K_I:=[\mathfrak D_I,(v\cdot\nabla_x+L_{g^\varepsilon})].
	\]
	Then $\mathsf P^\dy_m$ is the conjunction of the following.
	\begin{enumerate}
		\item\label{it:dyn-package-graph}
		For $|I|\le m$, $Z_I:=\mathfrak D_IZ_0\in\cH_{\kin,\sigma}(\cD_{T,\ast}\times\bR_v^3)$,
		its two flux traces are the completion-defined Green traces, and
		\[
		\|w^LZ_I\|_{\cH_{\kin,\sigma}(\cD_{T,\ast}\times\bR_v^3)}\le C,
		\qquad
		Z_I(-1)=0,
		\qquad
		(Z_I)_{\gamma^-}=\mathfrak D_I^-H_N .
		\]
		\item\label{it:dyn-package-equation}
		For $|I|\le m$, $\mathscr L_{g^\varepsilon} Z_I=-\mathfrak K_IZ_0+\mathfrak D_IF$ in
		$\cD'(I_T^-\times\cD\times\bR^3)$.
		\item\label{it:dyn-package-forcing}
		For $|I|\le m$, $\|w^L\mathfrak K_IZ_0\|_{L_t^2L_x^2H_\sigma^{-1}}\le C$.
	\end{enumerate}
	For the ordered words $I$ of $\mathfrak A_\bullet$ with $|I|\le m$, normally ordered by
	\eqref{eq:dyn-normal-ordering}, the same items define $\mathsf P^\bullet_m$; in particular
	$\mathsf P^\bullet_0=\mathsf P^\dy_0$.
\end{defn}

The weighted graph estimate for $Z_0$ follows from Theorem~\ref{thm:weighted-steady}\,(1); Proposition~\ref{prop:dyn-strong} controls the remaining strong coordinates of the zero extension.

\begin{lemma}[Fixed-cutoff base estimates]
	\label{prop:dyn-real-fixedN-strong}
	For every integer $0\le\ell\le\Lambda_0$, with all norms over $I_T^-$,
	\begin{align}
		&\|w^{\ell}Z_0\|_{\cH_{\kin,\sigma}(\cD_{T,\ast}\times\bR_v^3)}
		+\|w^\ell Z_0\|_{L_t^\infty L_x^{p_\sharp}H_\sigma^1}
		\notag\\
		&\qquad+\|w^{\ell+2}Z_0\|_{L^\infty}
		+\bigl[w^{\ell+2}Z_0\bigr]_{C^{\alpha_{\rm bs}}_{\kin}}
		+\|w^{\ell+2}\nabla_vZ_0\|_{L_t^\infty L_x^{p_\sharp}L_v^2}
		\le C .
		\label{eq:dynamic-fixedN-weighted-base-estimates}
	\end{align}
\end{lemma}

\begin{proof}
	\noindent\emph{Step 1. The regularity chain.}  Apply Proposition~\ref{prop:dyn-strong} at $\Lambda=\ell+2$ with
	\[
	(\hat f_\pt,\hat H,\mathsf F)=(\hat Z_0,\hat H_N,\hat F),
	\]
	the hats on the right denoting extension by zero to $t\le-1$.  Let $\mathsf N$ and $\mathsf S$ be the
	maxima over $0\le L\le\ell+B^{\dy}_{\rm in}+3$ (the range $0\le L\le\Lambda+B^{\dy}_{\rm in}+1$ of Proposition~\ref{prop:dyn-strong}) of
	\[
	\|w^L\hat H_N\|_{(C^{0,1}_{t,x}\cap C^{1,1}_v)(\gamma^-_{T,\ast})}+\|w^L(\hat H_N-H_s)\|_{L^2(\gamma^-_{T,\ast})} 
	\]
	and of 
	$$
	\|w^L\mathsf F\|_{L^\infty}+\|w^L\mathsf F\|_{L^2(\cD_{T,\ast};H^{-1}_\sigma)},
	$$ 
	respectively, as in \eqref{eq:dyn-strong-chain-datum} and \eqref{eq:dyn-strong-chain-source}.
	Then $\mathsf N<\infty$ and $\mathsf S<\infty$. For $\mathsf N$, this follows from
	\eqref{eq:dynamic-fixedN-incoming-strong-bound} and \eqref{eq:dynamic-fixedN-incoming-word-bound} at $I=\varnothing$ of Lemma~\ref{lem:dynamic-fixedN-incoming-flow-bounds}\,(\ref{it:dyn-incoming-words}), which apply since $\ell\le\Lambda_0$, together with
	$\|w^LH_s\|_{L^2(\gamma^-_{T,\ast})}\le C(T+1)^{1/2}\delta$; the bound on $H_s$ follows from
	\eqref{eq:stationary-incoming-regularity} at the weight $L+3$, as $H_s=(f_\st)_{\gamma^-}$ does not depend on $t$. For $\mathsf S$, see the first item below.  The
	extensions by zero vanish for $t\le-3/4$, since
	\[
	Z_0=0,
	\qquad
	H_N=(Z_0)_{\gamma^-}=0,
	\qquad
	F=0
	\qquad\text{on }(-1,-\tfrac34],
	\]
	by Proposition~\ref{prop:dyn-fixedN-realization}\,(\ref{it:dyn-real-solution}) for $Z_0$ and $H_N$ and by Proposition~\ref{prop:dyn-fixedN-data}\,(\ref{it:dyn-real-extended}) for $F$.
	The remaining hypotheses of Proposition~\ref{prop:dyn-strong} hold as follows.
	\begin{itemize}
		\item $\hat Z_0\in\cH_{\kin,\sigma}(\cD_{T,\ast}\times\bR_v^3)$ by
		Proposition~\ref{prop:dyn-fixedN-realization}\,(\ref{it:dyn-real-solution}).  By
		\eqref{eq:dyn-real-mollified-problem}, $\hat Z_0$ solves
		\eqref{eq:dyn-strong-quasilinear-regrouping} with background $\hat g^\varepsilon$ and
		source $\mathsf F=\hat F$, where $F=\bar F_0=\mathbf 1_{(-1,0)}\check{F}$ by
		\eqref{eq:dynamic-zero-initial-extended-sources}.  By
		Lemma~\ref{lem:dynamic-zero-initial-source-inflow-lift}\,(\ref{it:dyn-breg-source-bounds}), for $0\le L\le\ell+B^{\dy}_{\rm in}+3$,
		\[
		\|w^L\mathsf F\|_{L^\infty}
		+\|w^L\mathsf F\|_{L^2(\cD_{T,\ast};H^{-1}_\sigma)}
		\le C\|w^{\ell+B^{\dy}_{\rm in}+6}\check{F}\|_{L^\infty((-1,0)\times\cD\times\bR^3)} .
		\]
		Here
		\[
		\|w^{\ell+B^{\dy}_{\rm in}+6}\check{F}\|_{L^\infty((-1,0)\times\cD\times\bR^3)}<\infty
		\]
		by
		Lemma~\ref{lem:dynamic-zero-initial-source-inflow-lift}\,(\ref{it:dyn-breg-lift-bounds}). This gives \eqref{eq:dyn-strong-chain-source}.
		\item The incoming condition $(\hat Z_0)_{\gamma^-}=\hat H_N$ holds by
		\eqref{eq:dynamic-fixedN-strong-incoming}.
	\end{itemize}
	
	\noindent\emph{Step 2. Return to $Z_0$.}  The weighted graph estimate supplies the first term of \eqref{eq:dynamic-fixedN-weighted-base-estimates}; the remaining terms follow from \eqref{eq:dyn-strong-chain-output} in Step~1. Since $Z_0$ solves \eqref{eq:dyn-real-mollified-problem} with $Z_0(-1)=0$ and, by \eqref{eq:dynamic-fixedN-strong-incoming}, $(Z_0)_{\gamma^-}=H_N$, Theorem~\ref{thm:weighted-steady}\,(1), in the form \eqref{eq:weighted-graph-convention} at $L=\ell$ on $I_T^-$, gives
	\[
	\|w^\ell Z_0\|_{\cH_{\kin,\sigma}(\cD_{T,\ast}\times\bR_v^3)}
	\le C\bigl(\|w^\ell\bar F_0\|_{L^2(\cD_{T,\ast};H_\sigma^{-1})}+\|w^\ell H_N\|_{L^2(\gamma^-_{T,\ast})}\bigr).
	\]
	Its hypotheses hold: $w^{\ell}\bar F_0\in L_t^2L_x^2H_\sigma^{-1}$ by the first item of Step~1, $w^{\ell}H_N\in L^2(\gamma^-_{T,\ast})$ by \eqref{eq:dynamic-fixedN-incoming-word-bound} at $I=\varnothing$, and the coefficient condition \eqref{eq:weighted-dynamic-coefficient} for $g^\varepsilon$ holds by \eqref{eq:dyn-real-coefficient-hypotheses} and \eqref{eq:dynamic-KS-perturbed-smallness-choice}.  The pointwise, H\"older and $\nabla_v$ terms of \eqref{eq:dyn-strong-chain-output}, with $\hat f_\pt=\hat Z_0=Z_0$ on $I_T^-$, are the last three terms of \eqref{eq:dynamic-fixedN-weighted-base-estimates}.  By
	\eqref{eq:base-slot-pointwise}, applied at weight $\ell$ and $r=p_\sharp$,
	\[
	\|w^\ell Z_0\|_{L_t^\infty L_x^{p_\sharp}H_\sigma^1}\le C\|w^{\ell+2}\nabla_vZ_0\|_{L_t^\infty L_x^{p_\sharp}L_v^2}+C\|w^{\ell+2}Z_0\|_{L^\infty},
	\]
	since $w\ge1$. This proves \eqref{eq:dynamic-fixedN-weighted-base-estimates}.
\end{proof}

\noindent
By Proposition~\ref{prop:dyn-fixedN-realization}\,(\ref{it:dyn-real-solution})
and Lemma~\ref{prop:dyn-real-fixedN-strong} at $0\le\ell\le\Lambda_0$,
$\mathsf P^\dy_0$ holds.  Indeed, at $I=\varnothing$ we have $\mathfrak K_\varnothing=0$, so (\ref{it:dyn-package-forcing}) of $\mathsf P^\dy_0$ is trivial and
(\ref{it:dyn-package-equation}) is $\mathscr L_{g^\varepsilon} Z_0=F$, that is
\eqref{eq:dyn-real-mollified-problem}.  In (\ref{it:dyn-package-graph}), the weighted bound is \eqref{eq:dynamic-fixedN-weighted-base-estimates} at $\ell=L\le\Lambda_0$, $Z_0(-1)=0$ by the paragraph before Definition~\ref{def:dyn-package}, and $(Z_0)_{\gamma^-}=H_N$ by \eqref{eq:dynamic-fixedN-strong-incoming}.

\subsubsection{The realization ladder}\label{ss:dyn-ordered-words}

Write $g=g^\varepsilon$. Starting from $\mathsf P^\dy_0$, we realize the complete spatial jets by the flow-difference induction of Proposition~\ref{prop:fixedN-order-m-realization}, now on $I_T^-$ and with source $F$. The purely temporal chain and the mixed-word induction then give $\mathsf P^\bullet_{\mathrm M}$, hence parts (\ref{it:dyn-real-prop-spatial}) and (\ref{it:dyn-real-restriction}) of Proposition~\ref{prop:dyn-fixedN-realization}. The temporal exponents and weight losses used in this induction are as follows.  

By
\eqref{eq:graph-to-mixed-integrability-evolutionary} a jet in the first slot may be measured in $L_t^\varrho$ for any $\varrho\in[2,q_\flat]$, so $q_\flat$ is also the temporal endpoint, and the
temporal counterpart of $\delta_\flat$ of \eqref{eq:dyn-own-order-endpoints} is
\begin{equation}
	\label{eq:dyn-temporal-endpoint}
	\delta_t:=\frac12-\frac1{q_\flat}=\frac12-\frac{13}{30}=\frac1{15}>0 .
\end{equation}
The rungs of the realization in this section are those of
\eqref{eq:dyn-own-order-exponents} at the top order $\mathrm M$, together with the
temporal rung produced by the same map,
\begin{equation}
	\label{eq:dyn-dynamic-rungs}
	\mathsf g_\dy(j)=\mathsf e(\mathrm M/j;q_\flat,\infty),
	\qquad
	\mathsf s_\dy(j):=\mathsf e(\mathrm M/j,p_\flat),
	\qquad
	\mathsf r_\dy(j):=\mathsf e(\mathrm M/j,\infty),
	\qquad 0\le j\le\mathrm M .
\end{equation}
Thus $\mathsf g_\dy(j)=q_\flat\mathrm M/j$, as in Lemma~\ref{lem:dyn-breg-background}\,(\ref{it:dyn-breg-bg-norms}), and $\mathsf r_\dy(j)=2\mathrm M/j$, with $\mathrm M/0:=\infty$; in particular, $\mathsf g_\dy(k)=q_\flat\mathrm M/k$ is the exponent of \eqref{eq:dynamic-spatial-G-rung}.
A jet carrying $k$ of the $\mathrm M$ derivatives is measured in $L_t^{\mathsf g_\dy(k)}L_x^{\mathsf g_\dy(k)}$ in the first slot, and one carrying $j$ of them in $L_t^{\mathsf r_\dy(j)}$ in the solution slot;
$\mathsf r_\dy(\mathrm M)=2$ is the exponent of the graph norm. By \eqref{eq:dyn-temporal-endpoint} and the endpoint pair
\eqref{eq:dyn-own-order-endpoints}, the rungs are sub-conjugate in both variables in a row of order $m\le\mathrm M$, as in \eqref{eq:dyn-row-pairing},
\begin{equation}
	\label{eq:dyn-dynamic-conjugacy}
	\begin{aligned}
		\frac1{\mathsf g_\dy(k)}+\frac1{\mathsf r_\dy(j\mathrm M/m)}
		\,&=\frac12-\frac k2\Bigl(\frac1m-\frac1{\mathrm M}\Bigr)-\frac k{\mathrm M}\delta_t
		\ <\ \frac12,\\
		\frac1{\mathsf g_\dy(k)}+\frac1{\mathsf s_\dy(j\mathrm M/m)}
		\,&=\frac12-k\Bigl(\frac12-\frac1{p_\flat}\Bigr)\Bigl(\frac1m-\frac1{\mathrm M}\Bigr)-\frac k{\mathrm M}\delta_\flat
		\ <\ \frac12,
	\end{aligned}
	\qquad k+j=m\le\mathrm M,\quad k\ge1 .
\end{equation}

By Lemma~\ref{lem:dyn-hier-weight-trade}, a polynomial weight on the jet in the solution slot replaces the exponents of that jet by those of \eqref{eq:weight-trade-output-exponents}; for the original-Coulomb terms in the rows of order $\mathrm M$ below we take the ratio $\theta$ of that lemma to be
\begin{equation}\label{eq:dyn-dynamic-weight-ratio}
	\theta_\flat:=\frac{\delta_\flat}{2(\frac12-\frac1{p_\flat})},
	\qquad \frac{\theta_\flat}{2}<\delta_t.
\end{equation}
The last inequality follows from the fixed endpoints
\eqref{eq:dyn-own-order-endpoints} and \eqref{eq:dyn-temporal-endpoint}.

The induction at order $1\le m\le\mathrm M$ and weight $L\le\Lambda_m$
uses lower-order estimates only at weights bounded by
\begin{equation}\label{2610151131}
	2m+6,\quad m(A_\diamond+2)+1,\quad L+2m+2,\quad
	m(L+A_\diamond),\quad (m+1)L,
	\quad \frac L{\theta_\flat}\quad(m=\mathrm M).
\end{equation}
Indeed, the first two bounds cover the fixed polynomial losses, the next
three cover the commutators and the lower-order interpolations, and the last
is needed only at the top order. We set
\begin{equation}\label{eq:dyn-Km}
	K_m:=\max\{m+1,\,1/\theta_\flat\},\qquad 1\le m\le\mathrm M,
\end{equation}
which depends only on $m$ and the exponents fixed in
\eqref{eq:dyn-own-order-endpoints}.
Since $m\le\mathrm M<A_\diamond$, every weight in \eqref{2610151131} is at most
$K_m(L+A_\diamond+1)$ for every $L\ge0$.  Hence, by
\eqref{eq:dyn-realization-weights}, all are at most $\Lambda_{m-1}$.

We first bound the differences of the source.

\begin{lemma}
	\label{lem:dyn-breg-source-coefficient}
	Let $\nu\in\mathfrak A$, $X=\tau_\nu$ and $0\le j\le\mathrm M$; let $J$ be an ordered
	complete spatial word and $r\ge1$ with $|J|+r\le\mathrm M$; let $h\neq0$ with $|h|$ less than the window
	of Lemma~\ref{lem:dynamic-fixedN-incoming-flow-bounds}\,(\ref{it:dyn-incoming-differences}).  Then
	\begin{equation}
		\label{eq:dyn-breg-source-differences}
		\|(\Delta_h^X)^r\mathfrak D_J\bar F_j\|_{L^2(\cD_{T,\ast};H_\sigma^{-1})}
		\le C\|w\,\mathfrak D_{X^rJ}\bar F_j\|_{L^2(\cD_{T,\ast}\times\bR^3)}
		<\infty,
	\end{equation}
	with $C$ independent of $h$.
\end{lemma}

\begin{proof}
	By Lemma~\ref{lem:dynamic-zero-initial-source-inflow-lift}\,(\ref{it:dyn-breg-lift-bounds}) at $L=1$,
	\[
	\sum_{|I|\le\mathrm M}\sum_{k\le\mathrm M}\|w\,\mathfrak D_I\partial_t^k\check{F}\|_{L^2((-1,0)\times\cD\times\bR^3)}\le C\mathsf I ,
	\]
	so
	$w\,\mathfrak D_X^{\,i}\mathfrak D_J\bar F_j\in L^2$ for $0\le i\le r$.
	Hence Lemma~\ref{lem:strong-complete-flow-generator-L2} in its weighted form, applied $r$
	times at almost every $t$ as in the proof of Lemma~\ref{lem:repeated-same-flow-FTC}\,(\ref{it:repeated-FTC-identified}),
	gives
	\[
	(\Delta_h^X)^r\mathfrak D_J\bar F_j
	=\int_{[0,1]^r}T^X_{(\theta_1+\cdots+\theta_r)h}\mathfrak D_{X^rJ}\bar F_j\,\dd\theta_1\cdots\dd\theta_r ,
	\]
	where $|(\theta_1+\cdots+\theta_r)h|\le r|h|<s_0/2$, since $r\le\mathrm M$ and the window of
	Lemma~\ref{lem:dynamic-fixedN-incoming-flow-bounds}\,(\ref{it:dyn-incoming-differences}) is at most $s_0/(2\mathrm M)$, with $s_0$ of
	\eqref{eq:complete-flow-uniform-bounds}.  Combining this identity with \eqref{eq:complete-flow-L2-bound} at $L=1$ and
	\eqref{eq:Hsigma-dual-weighted}, we obtain
	\eqref{eq:dyn-breg-source-differences}; the constant does not depend on $h$, since that of \eqref{eq:complete-flow-L2-bound} does not depend on $s$.
\end{proof}

\begin{prop}[Derivative realization]
	\label{prop:dyn-order-m-realization}
	Let $1\le m\le\mathrm M$.  If $\mathsf P^\dy_{m-1}$ of Definition~\ref{def:dyn-package}
	holds, then $\mathsf P^\dy_m$ holds.
\end{prop}

\begin{proof}
	The spatial jets of length $m$ are constructed by induction over $(r,J)$; their equations and forcing estimates are then identified.
	
	\medskip
	\noindent\textbf{The induction over $(r,J)$.}
	Every ordered complete spatial word of length $m$ is $\nu^rJ$ with $\nu$ a letter, $r\ge1$ and $|J|=m-r$; with $X=\tau_\nu$ we write it $X^rJ$.  There is one step for each letter $\nu$ and each pair $(r,J)$ with $r+|J|=m$, taken in decreasing order of $r$.  Throughout the step, $X=\tau_\nu$ and $J$ are
	fixed, and
	\[
	Z_h:=(\Delta_h^X)^rZ_J,
	\qquad
	E_h:=\|Z_h\|_{\cH_{\kin,\sigma}(\cD_{T,\ast}\times\bR_v^3)},
	\qquad 0<|h|<h_0,
	\]
	with the window $h_0$ chosen below. By the convention stated after \eqref{eq:bulk-flow-difference}, cube variables are written
	$\theta_1,\theta_2,\ldots$ and $\bar h$ denotes a flow parameter that arises from
	difference quotients and is not written out.  The window $h_0$ does not depend on the
	step.  We take for $h_0$ the window of Lemma~\ref{lem:dynamic-fixedN-incoming-flow-bounds}\,(\ref{it:dyn-incoming-differences}), decreased, if necessary, to at most the window of Lemma~\ref{lem:discrete-interpolation} (used in part~(v) through
	\eqref{eq:dynamic-discrete-time}) and to less than $s_0/(2\mathrm M)$, with $s_0$ of
	\eqref{eq:complete-flow-uniform-bounds}.  Then $h_0<s_0/(2\mathrm M)\le s_0/\mathrm N$, since $2\mathrm M=2\mathrm N-2\ge\mathrm N$, so $h_0$ is a window of \S\ref{app:flow-identities}, and every flow parameter below has absolute value at most
	$2\mathrm Mh_0<s_0$; hence it belongs to the windows of \eqref{eq:repeated-same-flow-FTC-fiber}, \eqref{eq:boundary-flow-bound}, \eqref{eq:complete-flow-fiber-bound} and
	Lemma~\ref{lem:complete-flow-actual-traces}\,(\ref{it:flow-traces-evolutionary}).  Part~(viii) decreases $h_0$ further.  In the step, a jet carrying $j$ derivatives in a solution slot is measured in $L_t^{\mathsf r_\dy(j\mathrm M/m)}L_x^{\mathsf s_\dy(j\mathrm M/m)}H_\sigma^1$, as for the row of order $m$.  The step assumes the unweighted bound
	\begin{equation}
		\label{eq:dyn-step-input}
		\|Z_{J'}\|_{L_t^{\mathsf r_\dy(|J|\mathrm M/m)}
			L_x^{\mathsf s_\dy(|J|\mathrm M/m)}H_\sigma^1}
		\le C
		\qquad\text{for every }|J'|=|J|,
	\end{equation}
	and proves
	\[
	\|Z_{XJ}\|_{L_t^{\mathsf r_\dy((|J|+1)\mathrm M/m)}L_x^{\mathsf s_\dy((|J|+1)\mathrm M/m)}H_\sigma^1}\le C.
	\]
	At $r=m$, where $J=\varnothing$, $Z_J=Z_0$ and the norm in \eqref{eq:dyn-step-input} is $L_t^\infty L_x^{p_\flat}H_\sigma^1$,
	\eqref{eq:dyn-step-input} follows from
	\eqref{eq:dynamic-fixedN-weighted-base-estimates} of Lemma~\ref{prop:dyn-real-fixedN-strong} at $\ell=0$, since $p_\flat<p_\sharp$ and $\cD$ is bounded; at $r<m$ it is given by the steps at $r+1$.
	
	\medskip
	\noindent\emph{(i) Exponents.}
	The reciprocals of the space and time rungs are affine in the
	number of derivatives.  Since $|J|+r=m$ and both final exponents
	are $2$, discrete interpolation gives, for $0\le i<r$,
	\begin{equation}\label{eq:dyn-step-exponent}
		\begin{gathered}
			\left(1-\frac ir\right)\frac1{\mathsf s_\dy(|J|\mathrm M/m)}+\frac i{2r}
			=\frac1{\mathsf s_\dy((|J|+i)\mathrm M/m)},\\
			\left(1-\frac ir\right)\frac1{\mathsf r_\dy(|J|\mathrm M/m)}+\frac i{2r}
			=\frac1{\mathsf r_\dy((|J|+i)\mathrm M/m)}.
		\end{gathered}
	\end{equation}
	Thus \eqref{eq:dynamic-discrete-time} and
	\eqref{eq:GN-output-exponents} give the rung of order $|J|+i$
	in both variables.  With $k=r-i\ge1$, the strict inequalities
	\eqref{eq:dyn-dynamic-conjugacy} give the H\"older pairing for a
	background factor of order $k$ and a solution factor of order
	$|J|+i$.
	
	\medskip
	\noindent\emph{(ii) The lower differences.}
	We bound the differences of complete order at most $m-1$.
	For a subword $J'$ of $J$ and an integer $i\ge0$ with $|J'|+i\le m-1$, the jet $Z_{X^iJ'}$ satisfies $\|w^LZ_{X^iJ'}\|_{\cH_{\kin,\sigma}(\cD_{T,\ast}\times\bR_v^3)}\le C$ for $0\le L\le\Lambda_{m-1}$, by
	(\ref{it:dyn-package-graph}) of $\mathsf P^\dy_{m-1}$.
	We apply Lemma~\ref{lem:temporal-slicing} to $z=Z_{J'}$ and the words $X^j$,
	$0\le j\le m-1-|J'|$; the distributional derivatives $\mathfrak D_{X^j}Z_{J'}=Z_{X^jJ'}$
	belong to $\cH_{\kin,\sigma}(\cD_{T,\ast}\times\bR_v^3)$.  The lemma gives
	$\mathfrak D_X(Z_{X^jJ'}(t))=Z_{X^{j+1}J'}(t)$ for $t$ outside a null set, which may be taken independent of $J'$.  At these times
	Lemma~\ref{lem:repeated-same-flow-FTC}\,(\ref{it:repeated-FTC-identified}) applies, and
	\eqref{eq:repeated-same-flow-weighted-bound} bounds, by the jets $Z_{X^jJ'}$, the differences
	$(\Delta_h^X)^iZ_{J'}$ and the quantity of \eqref{eq:discrete-lower-density},
	\[
	\cR^{X,0}_{r,h}(Z_J)=\sum_{i=0}^{r-1}\|w^{r-i}(\Delta_h^X)^iZ_J\|_{L_x^2H_\sigma^1},
	\]
	and these bounds give
	\begin{equation}
		\label{eq:dyn-step-discrete-lower}
		\sup_{0<|h|<h_0}
		\Bigl(\|\cR_{r,h}^{X,0}(Z_J)\|_{L_t^2}
		+\sum_{J'}\sum_{i=0}^{m-1-|J'|}
		\|w^L(\Delta_h^X)^iZ_{J'}\|_{L_t^2L_x^2H_\sigma^1}\Bigr)
		\le C
	\end{equation}
	for $0\le L\le\Lambda_{m-1}-1$ (the right side of \eqref{eq:repeated-same-flow-weighted-bound} carries the weight $w^{L+1}$), the outer sum being taken over the subwords $J'$ of $J$.
	
	\medskip
	\noindent\emph{(iii) Operator orbits and traces.}
	The operator identity \eqref{eq:direct-operator-cube} is valid in $\cD'$ for $\mathrm A=\mathscr T_L$ at the arguments used below. We also need $Z_h\in\cH_{\kin,\sigma}(\cD_{T,\ast}\times\bR_v^3)$.
	By Corollary~\ref{cor:A0-orbits} at $\ell=0$ and at every order $1\le s\le r$,
	\eqref{eq:operator-orbit-bound} holds for $\mathrm A=\mathscr T_L=v\cdot\nabla_x+L$ at the order $s$ with
	$(L_1,p_1)=(2s+4,2)$.  Since the flows $T^X$ commute with $\partial_t$, $\mathfrak Q_{s,h}^X(\partial_t+\mathscr T_L)=\mathfrak Q_{s,h}^X(\mathscr T_L)$ by \eqref{eq:complete-flow-time-commutation}.  At the arguments
	$T_{sh}^X(\Delta_h^X)^{r-s}Z_J$ of the difference equation of part~(iv), the bounds \eqref{eq:complete-flow-H1-bound}, \eqref{eq:local-H1-H_sigma}, \eqref{eq:real-weight-Hsigma-monotonicity} and \eqref{eq:dyn-step-discrete-lower} at $L=A_\diamond+1\le\Lambda_{m-1}-1$ give
	\[
	\begin{aligned}
		&\|T_{sh}^X(\Delta_h^X)^{r-s}Z_J\|_{L_t^2L_{x,v}^2}
		+\|T_{sh}^X(\Delta_h^X)^{r-s}Z_J\|_{L_t^2L_x^2H_\sigma^1}
		+\|w^{2s+4}T_{sh}^X(\Delta_h^X)^{r-s}Z_J\|_{L_t^2L_x^2H_\sigma^1}\\
		&\qquad\le C\|w^{A_\diamond+1}(\Delta_h^X)^{r-s}Z_J\|_{L_t^2L_x^2H_\sigma^1}<\infty ,
		\qquad 1\le s\le r ,
	\end{aligned}
	\]
	since the pullback in \eqref{eq:complete-flow-H1-bound} raises the weight $2s+4$ to $2s+5\le2\mathrm M+5=A_\diamond+1$.
	These are the norms of \eqref{eq:operator-orbit-bound} at $(L_1,p_1)=(2s+4,2)$, taken in $L_t^2$.  Thus, for almost every $t\in I_T^-$ and $1\le s\le r$, \eqref{eq:direct-operator-cube} holds in $\cD'$ for $\mathrm A=\mathscr T_L$ at $\bigl(T_{sh}^X(\Delta_h^X)^{r-s}Z_J\bigr)(t)$.
	
	By (\ref{it:dyn-package-graph}) of $\mathsf P^\dy_{m-1}$,
	$w^5Z_J\in\cH_{\kin,\sigma}(\cD_{T,\ast}\times\bR_v^3)$, since $|J|\le m-1$ and
	$5\le\Lambda_{m-1}$.  Lemma~\ref{lem:complete-flow-actual-traces}\,(\ref{it:flow-traces-evolutionary}), applied
	to $q=Z_J$ at the flow times $h,2h,\ldots,rh$, therefore gives
	\begin{equation}
		\label{eq:dyn-step-flow-membership}
		T_{sh}^XZ_J\in\cH_{\kin,\sigma}(\cD_{T,\ast}\times\bR_v^3),
		\qquad
		(T_{sh}^XZ_J)_{\gamma^\pm}=T_{sh}^{X,\gamma_\pm}(Z_J)_{\gamma^\pm},
		\qquad 1\le s\le r ,
	\end{equation}
	and
	\begin{equation}
		\label{eq:dyn-step-flow-slices}
		(T_{sh}^XZ_J)(t)=T_{sh}^X\bigl(Z_J(t)\bigr),
		\qquad -1\le t\le T,
		\qquad 1\le s\le r .
	\end{equation}
	In particular, since $Z_h$ is a linear combination of the $T_{sh}^XZ_J$, $0\le s\le r$, by \eqref{eq:bulk-flow-difference}, $Z_h\in\cH_{\kin,\sigma}(\cD_{T,\ast}\times\bR_v^3)$ and $E_h<\infty$ for each fixed $h\neq0$, and \eqref{eq:direct-same-flow-actual-traces} gives $(Z_h)_{\gamma^\pm}=(\Delta_h^{X,\gamma_\pm})^r(Z_J)_{\gamma^\pm}$.
	
	\medskip
	\noindent\emph{(iv) The difference equation.}
	We derive the equation satisfied by $Z_h$.
	By (\ref{it:dyn-package-equation}) of $\mathsf P^\dy_{m-1}$,
	$\mathscr L_g Z_J=-\mathfrak K_JZ_0+\mathfrak D_JF$.  At $|J|=0$, $\mathfrak K_\varnothing=0$.  At $|J|\ge1$,
	\eqref{eq:weak-triangular-commutator} expands $\mathfrak K_JZ_0$ in
	$\cD'(I_T^-\times\cD\times\bR^3)$, as in the proof of
	Proposition~\ref{prop:dyn-differentiated-equation}.  Indeed, the hypotheses
	\eqref{eq:dyn-differentiated-chain}--\eqref{eq:dyn-differentiated-first-slots} of that
	proposition, for ordered complete spatial words, hold at the order $|J|$ by
	(\ref{it:dyn-package-graph}) of $\mathsf P^\dy_{m-1}$ at $L=0$ and $L=A_\diamond$
	and by \eqref{eq:dyn-real-background-norms} and Lemma~\ref{lem:dyn-breg-background}\,(\ref{it:dyn-breg-bg-norms}), as at the order $m$ in
	the verification of (\ref{it:dyn-package-equation}) below.  Since $\mathscr L_g=\partial_t+\mathscr T_L-N[g,\cdot\,]$, the equation of $Z_J$
	becomes
	\begin{equation}
		\label{eq:dyn-step-source-split}
		(\partial_t+\mathscr T_L) Z_J
		=N[g,Z_J]
		+\sum_{\varnothing\neq S\subset[|J|]}
		N[\mathfrak D_{J_S}g,Z_{J_{S^c}}]
		+\mathfrak R_J[g,Z_0]
		-\mathfrak G_JZ_0
		+\mathfrak D_JF ,
	\end{equation}
	where $\mathfrak G_J:=[\mathfrak D_J,\mathscr T_L]$, as in the proof of
	Proposition~\ref{prop:fixedN-order-m-realization}, and
	$\mathfrak G_\varnothing=\mathfrak R_\varnothing=0$.
	
	By the last identity of \eqref{eq:ad-word} at $I=J$ and $\mathrm A=\mathscr T_L$, and by induction on
	\eqref{eq:deformation-remainder-recurrence},
	the jet in every solution slot of $\mathfrak G_JZ_0$ and of $\mathfrak R_J[g,Z_0]$ is $Z_{J'}$ for a proper subword $J'$ of $J$.  In particular,
	\[
	\mathfrak G_JZ_0=\sum_{\varnothing\ne S\subset[|J|]}\operatorname{ad}_{J_S}(\mathscr T_L)Z_{J_{S^c}}
	\qquad\text{in }\cD',
	\]
	since both sides expand $[\mathfrak D_J,\mathscr T_L]$ with each lower complete derivative
	replaced by its jet.  Each $\operatorname{ad}_{J_S}(\mathscr T_L)$ is, by
	\eqref{eq:nested-commutator-decomposition} at $g=0$, a sum of transport, local and compact
	commutators that contains no complete derivative.  By Corollary~\ref{cor:A0-orbits},
	\eqref{eq:operator-orbit-bound} and its form with $B_s$ replaced by $\mathrm A$ hold for $\mathrm A=\operatorname{ad}_{J_S}(\mathscr T_L)$ at the order $s$ with $(L_1,p_1)=(2(s+|S|)+4,2)$.  At the arguments $T_{sh}^X(\Delta_h^X)^{r-s}Z_{J_{S^c}}$,
	\eqref{eq:complete-flow-H1-bound}, \eqref{eq:local-H1-H_sigma}, \eqref{eq:real-weight-Hsigma-monotonicity} and \eqref{eq:dyn-step-discrete-lower} at $L=A_\diamond+1$ give
	\[
	\begin{aligned}
		&\|T_{sh}^X(\Delta_h^X)^{r-s}Z_{J_{S^c}}\|_{L_t^2L_{x,v}^2}\\
		&+\|T_{sh}^X(\Delta_h^X)^{r-s}Z_{J_{S^c}}\|_{L_t^2L_x^2H_\sigma^1}
		+\|w^{2(s+|S|)+4}T_{sh}^X(\Delta_h^X)^{r-s}Z_{J_{S^c}}\|_{L_t^2L_x^2H_\sigma^1}\\
		\le\,&C\|w^{A_\diamond+1}(\Delta_h^X)^{r-s}Z_{J_{S^c}}\|_{L_t^2L_x^2H_\sigma^1}<\infty ,
		\qquad 1\le s\le r ,
	\end{aligned}
	\]
	since the pullback raises the weight $2(s+|S|)+4$ to $2(s+|S|)+5$ and, by $s+|S|\le r+|J|=m$,
	\[
	2(s+|S|)+5\le2\mathrm M+5=A_\diamond+1 ;
	\]
	so, as for $\mathscr T_L$ in part~(iii), for almost every $t\in I_T^-$ and $1\le s\le r$, \eqref{eq:direct-operator-cube} holds in $\cD'$ for $\mathrm A=\operatorname{ad}_{J_S}(\mathscr T_L)$ at $\bigl(T_{sh}^X(\Delta_h^X)^{r-s}Z_{J_{S^c}}\bigr)(t)$.
	
	Apply $(\Delta_h^X)^r$ to \eqref{eq:dyn-step-source-split}, expand $\partial_t+\mathscr T_L$ at $u=Z_J$
	and each $\operatorname{ad}_{J_S}(\mathscr T_L)$ at $u=Z_{J_{S^c}}$ by
	\eqref{eq:forward-discrete-Leibniz}, expand each original-Coulomb product
	by \eqref{eq:collision-mth-endpoint}, and subtract
	$N[g,Z_h]$. The term of $(\Delta_h^X)^rN[g,Z_J]$ carrying all $r$ differences on the solution slot combines with $-N[g,Z_h]$ into $N[T_{rh}^Xg-g,Z_h]$: in \eqref{eq:collision-mth-endpoint} at $[a,b]=[g,Z_J]$, the term with all $r$ differences on $b$ is $N[T_{rh}^Xg,Z_h]$, since, with $K_0=\varnothing$ there, $\mathscr C_{0,h}^X=N_{\mathsf K_\varnothing}=N$.  By \eqref{eq:complete-flow-time-commutation}, the resulting equation reads
	\begin{equation}
		\label{eq:dyn-step-endpoint-equation}
		\mathscr L_g Z_h
		=\mathscr R_{r,h}^{\mathscr C_0}
		+\mathscr R_{r,h}^{\mathscr C_+}
		+\mathscr R_{r,h}^{\mathscr T_L}
		{}+N[T_{rh}^Xg-g,\,Z_h]
		+(\Delta_h^X)^r\mathfrak D_JF
		\quad\text{in }\cD'(I_T^-\times\cD\times\bR^3),
	\end{equation}
	where $\mathscr R^{\mathscr C_0}_{r,h}$ collects the original-Coulomb terms
	of \eqref{eq:collision-mth-endpoint} with $S\ne\varnothing$ or
	with at least one difference on the first slot,
	$\mathscr R^{\mathscr C_+}_{r,h}$ the defect terms $\mathscr C_{s,h}^X$, $s\ge1$, of
	\eqref{eq:collision-flow-defect-difference} together with the differences of
	$\mathfrak R_J[g,Z_0]$, and
	\[
	\begin{aligned}
		\mathscr R_{r,h}^{\mathscr T_L}
		:={}&-\sum_{s=1}^r\binom rs
		\mathfrak Q_{s,h}^X(\mathscr T_L)
		T_{sh}^X(\Delta_h^X)^{r-s}Z_J\\
		&-\sum_{\varnothing\ne S\subset[|J|]}\Bigl(
		\operatorname{ad}_{J_S}(\mathscr T_L)(\Delta_h^X)^rZ_{J_{S^c}}
		+\sum_{s=1}^r\binom rs
		\mathfrak Q_{s,h}^X\bigl(\operatorname{ad}_{J_S}(\mathscr T_L)\bigr)
		T_{sh}^X(\Delta_h^X)^{r-s}Z_{J_{S^c}}\Bigr).
	\end{aligned}
	\]
	Every jet in a solution slot of $\mathscr R^{\mathscr T_L}_{r,h}$ is thus, up to a flow translate, a difference
	$(\Delta_h^X)^iZ_{J'}$ with $J'$ a subword of $J$ and $|J'|+i\le m-1$ (in the first sum $J'=J$ and $i=r-s\le r-1$; in the second $J'=J_{S^c}$, so that $|J'|\le|J|-1$, and $i\le r$).
	
	\medskip
	\noindent\emph{(v) The solution factors.}
	The solution factors on the right of \eqref{eq:dyn-step-endpoint-equation} have complete order at most $m-1$ and are estimated at the corresponding interpolation exponents.
	The differences of $Z_J$ itself are bounded by discrete interpolation:
	\begin{equation}
		\label{eq:dyn-step-GN}
		\|(\Delta_h^X)^iZ_J\|
		_{L_t^{\mathsf r_\dy((|J|+i)\mathrm M/m)}
			L_x^{\mathsf s_\dy((|J|+i)\mathrm M/m)}H_\sigma^1}
		\le C\bigl(1+E_h^{\,i/r}\bigr),
		\qquad 1\le i\le r-1 .
	\end{equation}
	Let $\|Z_J\|$ denote the norm of \eqref{eq:dyn-step-input}.  By \eqref{eq:dynamic-discrete-top-density} and \eqref{eq:dyn-step-discrete-lower},
	$\|\mathfrak e_{r,h}^X[Z_J]\|_{L_t^2}\le C(1+E_h)$, since \[
	\|(\Delta_h^X)^rZ_J\|_{L_t^2L_x^2H_\sigma^1}=\|Z_h\|_{L_t^2L_x^2H_\sigma^1}\le E_h
	\] by Definition~\ref{def:H_kin_sigma}.  For $1\le i\le r-1$, \eqref{eq:dynamic-discrete-time}, at the exponents of part~(i) and with the output exponents of \eqref{eq:dyn-step-exponent}, gives the first inequality below; the power of the length of the time interval in its last term is absorbed in $C$.  Then \eqref{eq:dyn-step-input} and the bound on $\|\mathfrak e_{r,h}^X[Z_J]\|_{L_t^2}$ give the second:
	\[
	\begin{aligned}
		\|(\Delta_h^X)^iZ_J\|_{L_t^{\mathsf r_\dy((|J|+i)\mathrm M/m)}L_x^{\mathsf s_\dy((|J|+i)\mathrm M/m)}H_\sigma^1}
		&\le C\|Z_J\|^{1-i/r}\|\mathfrak e_{r,h}^X[Z_J]\|_{L_t^2}^{\,i/r}+C\|Z_J\|\\
		&\le C(1+E_h)^{i/r}+C ,
	\end{aligned}
	\]
	and $(1+E_h)^{i/r}\le1+E_h^{\,i/r}$ gives \eqref{eq:dyn-step-GN}.
	Every solution factor $(\Delta_h^X)^iZ_{J'}$ with $J'$ a subword of $J$ and
	$|J'|+i\le m-1$ is measured in $L_t^{\mathsf r_\dy((|J'|+i)\mathrm M/m)}L_x^{\mathsf s_\dy((|J'|+i)\mathrm M/m)}H_\sigma^1$, and, with $k:=m-|J'|-i$, \eqref{eq:dyn-dynamic-conjugacy} makes
	\[
	\bigl(\mathsf g_\dy(k),\mathsf g_\dy(k)\bigr),
	\qquad
	\bigl(\mathsf r_\dy((|J'|+i)\mathrm M/m),\mathsf s_\dy((|J'|+i)\mathrm M/m)\bigr)
	\]
	a H\"older pair.  Such a factor is bounded in this norm as follows.
	\begin{itemize}
		\item At $i=0$, by \eqref{eq:dyn-step-input} with $|J'|$ in place of $|J|$, assumed in this or in an earlier step.
		\item At $i\ge1$ and $J'=J$, by \eqref{eq:dyn-step-GN}, since $|J|+i\le m-1$ gives $i\le r-1$.
		\item At $i\ge1$ and $|J'|<|J|$, by \eqref{eq:dyn-step-GN} of the earlier step at
		$(r',J')$, $r':=m-|J'|>r$, with the same letter $\nu$ and with $E_h$ replaced by $C$, as allowed by \eqref{eq:dyn-step-graph}, proved in part~(viii) of that step; these bounds cover $1\le i\le r'-1$.
	\end{itemize}
	
	\medskip
	\noindent\emph{(vi) The background factors.}
	The background factors are controlled by the interpolation bounds, while the coefficient modulus of $g$ vanishes at small flow displacement.
	At $k=0$, where $\mathsf g_\dy(0)=\infty$,
	Lemma~\ref{lem:dyn-breg-background}\,(\ref{it:dyn-breg-bg-norms}) at $a=0$ gives
	\[
	\|w^{-2}g\|\le\|w^{-2}\bar g_0\|+2\|w^{-2}\hat g^{(0)}\|
	\quad\text{with all norms in }L^\infty(I_T^-;L_x^\infty L_v^2) ,
	\]
	and, by \eqref{eq:dyn-real-background-norms} and the same lemma,
	\[
	\|w^{-2}\bar g_0\|_{L^\infty(I_T^-;L_x^\infty L_v^2)}+2\|w^{-2}\hat g^{(0)}\|_{L^\infty(I_T^-;L_x^\infty L_v^2)}<\infty .
	\]%
	With \eqref{eq:complete-flow-fiber-bound}, taken with an outer $L_t^\infty$, this gives
	\begin{equation}
		\label{eq:dyn-step-background-zero-shifts}
		\sup_{|\bar h|\le2\mathrm Mh_0}
		\bigl\|w^{-2}T_{\bar h}^Xg\bigr\|_{L_t^\infty L_x^\infty L_v^2}
		\le C .
	\end{equation}
	Let $K$ be an ordered complete spatial word with $k:=|K|\ge1$.  By Lemma~\ref{lem:dyn-breg-background}\,(\ref{it:dyn-breg-bg-rungs}), $g$ satisfies on $I_T^-$ \eqref{eq:dynamic-spatial-G-rung} if $k\le\mathrm M-1$, and \eqref{eq:graph-to-mixed-integrability-evolutionary} with the weight $w^{-2}$ if $k=\mathrm M$.  These give
	\[
	w^{-2}\mathfrak D_Kg\in L_t^{\mathsf g_\dy(k)}L_x^{\mathsf g_\dy(k)}L_v^2.
	\]
	With \eqref{eq:dyn-step-background-zero-shifts}, they give $w^{-2}\mathfrak D_Kg\in L_x^1L_v^2$ at almost every $t$ for
	$|K|\le\mathrm M$.  For an ordered complete spatial word $K$ and $1\le k_1\le\mathrm M-|K|$, the chain
	identities $\mathfrak D_X\mathfrak D_{X^jK}g=\mathfrak D_{X^{j+1}K}g$, $0\le j<k_1$, hold by
	\eqref{eq:temporal-slicing-chain} of Lemma~\ref{lem:temporal-slicing} at $z=g$, for the
	family of the words $X^jK$ and the suffixes of $K$; the hypotheses of that lemma are
	those of Corollary~\ref{cor:dyn-mixed-rungs} used for \eqref{eq:dynamic-spatial-G-rung}.
	Thus, at almost every $t$, the hypotheses of Lemma~\ref{lem:repeated-same-flow-FTC}\,(\ref{it:repeated-FTC-fiber}) hold for the jets $\mathfrak D_{X^jK}g$, $0\le j\le k_1$, at
	$\lambda=-2$, and
	\eqref{eq:repeated-same-flow-FTC-fiber} gives
	\[
	T_{\bar h}^X(\Delta_h^X)^{k_1}\mathfrak D_Kg
	=\int_{[0,1]^{k_1}}T^X_{\bar h+(\theta_1+\cdots+\theta_{k_1})h}\mathfrak D_{X^{k_1}K}g\,\dd\theta_1\cdots\dd\theta_{k_1} .
	\]  The flow parameters
	satisfy
	\[
	|\bar h+(\theta_1+\cdots+\theta_{k_1})h|\le|\bar h|+k_1|h|\le2\mathrm Mh_0<s_0
	\]
	in the range of $\bar h$ below.  Hence, for every ordered complete spatial word $K$ and every $0\le k_1\le\mathrm M-|K|$ with
	$k:=k_1+|K|\ge1$,
	\begin{equation}
		\label{eq:dyn-step-background-differences}
		\sup_{|\bar h|\le2\mathrm Mh_0-k_1|h|}
		\bigl\|w^{-2}T_{\bar h}^X(\Delta_h^X)^{k_1}\mathfrak D_Kg
		\bigr\|_{L_t^{\mathsf g_\dy(k)}
			L_x^{\mathsf g_\dy(k)}L_v^2}
		\le C,
		\qquad 1\le k\le\mathrm M,
		\quad 0<|h|<h_0 ,
	\end{equation}
	with $\mathsf g_\dy(\mathrm M)=q_\flat$ in both variables at
	$k=\mathrm M$.  Here \eqref{eq:complete-flow-fiber-bound}, taken with an outer $L_t^{\mathsf g_\dy(k)}$, removes the flow translates, and
	\[
	w^{-2}\mathfrak D_{X^{k_1}K}g\in L_t^{\mathsf g_\dy(k)}L_x^{\mathsf g_\dy(k)}L_v^2,
	\]
	as above with $X^{k_1}K$ in place of $K$.
	For the coefficient difference $T_{rh}^Xg-g$ in the fourth term of
	\eqref{eq:dyn-step-endpoint-equation}, we use the coefficient modulus of $g$, with
	$s_0$ as above,
	\[
	\omega_g(\rho):=
	\max_{\nu\in\mathfrak A}\ \sup_{|s|\le\rho}
	\bigl\|\mu^{1/4}\bigl(T_s^{\tau_\nu}g-g\bigr)\bigr\|
	_{L^\infty(\cD_{T,\ast}\times\bR^3)},
	\qquad 0<\rho\le s_0 .
	\]
	Then $\omega_g(\rho)\to0$ as $\rho\downarrow0$.  We obtain this from Lemma~\ref{lem:complete-flow-coefficient-continuity} at $(B,\alpha)=(A^\dy_\flat,\alpha_{\rm bs})$ and $q=g(t)$, $-1<t<T$, with a weaker H\"older condition on $\mu^{1/4}q$.  In the proof of that lemma, the H\"older seminorm $[\mu^{1/4}q]_{C^\alpha_{\kin}}$ enters only in \eqref{eq:complete-flow-coefficient-local}, where it bounds $|\mu^{1/4}q(\Phi_s^X(x),A_s^X(x)v)-\mu^{1/4}q(x,v)|$ by its product with $\mathbf d((0,x,v),(0,\Phi_s^X(x),A_s^X(x)v))^{\alpha}$, and this distance is bounded by \eqref{eq:complete-flow-coefficient-distance}.  There the seminorm may be replaced by
	\[
	\sup_{(x,v)\ne(y,v')}
	\frac{|\mu^{1/4}q(x,v)-\mu^{1/4}q(y,v')|}{\mathbf d((0,x,v),(0,y,v'))^{\alpha}} ,
	\]
	where $(x,v)$ and $(y,v')$ range over $\overline\cD\times\bR^3$; this supremum is at most $[\mu^{1/4}q]_{C^\alpha_{\kin}}$, since $\mathbf d_{\st}((x,v),(y,v'))\le\mathbf d((0,x,v),(0,y,v'))$ by \eqref{eq:stationary-kinetic-distance-main}.  Hence \eqref{eq:complete-flow-coefficient-continuity} and the last assertion of Lemma~\ref{lem:complete-flow-coefficient-continuity} hold with $\|\mu^{1/4}q\|_{C^\alpha_{\kin}}$ replaced by the sum of $\|\mu^{1/4}q\|_{L^\infty_{x,v}}$ and this supremum.  Since $\mathbf d((t,x,v),(t,y,v'))=\mathbf d((0,x,v),(0,y,v'))$ by \eqref{kineticdistance}, it suffices that
	\[
	\sup_{-1<t<T}\Bigl(\|w^{A^\dy_\flat}g(t)\|_{L^\infty_{x,v}}
	+\|\mu^{1/4}g(t)\|_{L^\infty_{x,v}}
	+\sup_{(x,v)\ne(y,v')}
	\frac{|\mu^{1/4}g(t,x,v)-\mu^{1/4}g(t,y,v')|}{\mathbf d((t,x,v),(t,y,v'))^{\alpha_{\rm bs}}}\Bigr)<\infty .
	\]
	
	For the background $\hat g$ of Proposition~\ref{prop:dyn-fixedN-realization}, extended as in \eqref{eq:dyn-strong-past-extension}, so that $\hat g=\check{E}$ on $(-1,0]$, and for $Q=(-\infty,T)\times\overline\cD\times\bR^3$ as in \eqref{eq:dyn-strong-mollification-holder},
	\[
	\|\mu^{1/4}\hat g\|_{C^{\alpha_{\rm bs}}_{\kin}(Q)}<\infty,
	\qquad
	\sup_{-1<t<T}\|w^{A^\dy_\flat}\hat g(t)\|_{L^\infty_{x,v}}<\infty ;
	\]
	here the first bound follows from \eqref{eq:dyn-strong-coefficient-smallness}--\eqref{eq:dyn-strong-coefficient-holder}, and the second follows from \eqref{eq:dyn-ball-A} on $(0,T)$, with the radius of Proposition~\ref{prop:dyn-fixedN-realization} in place of $r$, and from
	Lemma~\ref{lem:dynamic-zero-initial-source-inflow-lift}\,(\ref{it:dyn-breg-lift-bounds}) at $L=A^\dy_\flat\le\Lambda_0$ on $(-1,0]$.
	Next, $g=g^\varepsilon$ with $g^\varepsilon=\varrho_\varepsilon*\hat g+(\hat g^{(0)}-\varrho_\varepsilon*\hat g^{(0)})$ by
	\eqref{eq:dyn-real-mollified-objects}. By \eqref{eq:dyn-strong-mollification-holder} at $u=\hat g$ and the H\"older term of \eqref{eq:dyn-real-increment-bound},
	\[
	\|\mu^{1/4}g^\varepsilon\|_{C^{\alpha_{\rm bs}}_{\kin}(I_T^-\times\overline\cD\times\bR^3)}
	\le\|\mu^{1/4}\hat g\|_{C^{\alpha_{\rm bs}}_{\kin}(Q)}+C\varepsilon(\delta+\mathsf I) .
	\]
	Since $\varrho_\varepsilon\ge0$ has integral one and $\hat g$ vanishes for $t\le-3/4$, averaging in $t$ gives the first inequality below, and the first norm in \eqref{eq:dyn-real-increment-bound} gives the second:
	\[
	\begin{aligned}
		\sup_{-1<t<T}\|w^{A^\dy_\flat}g^\varepsilon(t)\|_{L^\infty_{x,v}}
		&\le\sup_{-1<t<T}\|w^{A^\dy_\flat}\hat g(t)\|_{L^\infty_{x,v}}
		+\sup_{-1<t<T}\|w^{A^\dy_\flat}(\hat g^{(0)}-\varrho_\varepsilon*\hat g^{(0)})(t)\|_{L^\infty_{x,v}}\\
		&\le\sup_{-1<t<T}\|w^{A^\dy_\flat}\hat g(t)\|_{L^\infty_{x,v}}+C\varepsilon(\delta+\mathsf I) .
	\end{aligned}
	\]
	Finally, for every $-1<t<T$ the points $(t,x,v)$ and $(t,y,v')$ lie in $I_T^-\times\overline\cD\times\bR^3$, so the definition of the space-time H\"older norm gives
	\[
	\|\mu^{1/4}g^\varepsilon(t)\|_{L^\infty_{x,v}}
	+\sup_{(x,v)\ne(y,v')}
	\frac{|\mu^{1/4}g^\varepsilon(t,x,v)-\mu^{1/4}g^\varepsilon(t,y,v')|}{\mathbf d((t,x,v),(t,y,v'))^{\alpha_{\rm bs}}}
	\le\|\mu^{1/4}g^\varepsilon\|_{C^{\alpha_{\rm bs}}_{\kin}(I_T^-\times\overline\cD\times\bR^3)} .
	\]
	This inequality and the two bounds on $g^\varepsilon$ above show that the condition for $\omega_g(\rho)\to0$ holds.
	
	\medskip
	\noindent\emph{(vii) The right-hand side.}
	We bound the five terms on the right of \eqref{eq:dyn-step-endpoint-equation} in $L_t^2L_x^2H_\sigma^{-1}$, in the order $\mathscr R^{\mathscr C_0}_{r,h}$, $\mathscr R^{\mathscr C_+}_{r,h}$, $\mathscr R^{\mathscr T_L}_{r,h}$, $(\Delta_h^X)^r\mathfrak D_JF$ and $N[T_{rh}^Xg-g,Z_h]$.
	
	Each summand of $\mathscr R^{\mathscr C_0}_{r,h}$ is a product of a background
	factor of the form \eqref{eq:dyn-step-background-differences}, of some
	complete order $k\ge1$, with a solution factor
	$(\Delta_h^X)^iZ_{J'}$, where $J'$ is a subword of $J$ with
	$k+|J'|+i=m$ and $|J'|+i\le m-1$.
	By part~(v),
	\[
	\|(\Delta_h^X)^iZ_{J'}\|_{L_t^{\mathsf r_\dy((|J'|+i)\mathrm M/m)}L_x^{\mathsf s_\dy((|J'|+i)\mathrm M/m)}H_\sigma^1}\le C
	\]
	unless $J'=J$ and $1\le i\le r-1$.  Hence, by \eqref{eq:velalg} in the form
	\eqref{eq:velalg-poly}, \eqref{eq:dyn-step-background-differences}, and the bounds of part~(v),
	\begin{equation}
		\label{eq:dyn-step-material}
		\|\mathscr R^{\mathscr C_0}_{r,h}\|_{L_t^2L_x^2H_\sigma^{-1}}
		\le C
		+C\sum_{i=1}^{r-1}
		\|(\Delta_h^X)^iZ_J\|
		_{L_t^{\mathsf r_\dy((|J|+i)\mathrm M/m)}
			L_x^{\mathsf s_\dy((|J|+i)\mathrm M/m)}H_\sigma^1} .
	\end{equation}
	
	For $\mathscr R^{\mathscr C_+}_{r,h}$, let $1/q+1/p\le1/2$, $p<\infty$ and $1/\rho+1/\rho'\le1/2$.  Proposition~\ref{prop:complete-flow-collision-differences}
	at the order $s$, $\ell=0$, $m_0=0$ and $\varpi=1/4$, in its time-integrated form, with $I_T^-$ in place of $(0,T)$,
	together with $\mu^{1/4}\le Cw^{-2}$ and $2s+4\le A_\diamond+2$,
	gives
	\[
	\|\mathscr C_{s,h}^X[a,b]\|_{L_t^2L_x^2H_\sigma^{-1}}
	\lesssim
	\sup_{\theta_1,\ldots,\theta_s\in[0,1]}
	\|w^{-2}T_{(s-\theta_1-\cdots-\theta_s)h}^Xa\|_{L_t^\rho L_x^qL_v^2}
	\|w^{A_\diamond+2}b\|_{L_t^{\rho'}L_x^pH_\sigma^1},
	\qquad 1\le s\le r .
	\]
	Here \eqref{eq:collision-flow-defect-Peano-class} holds at almost
	every $t$ with $\lambda=2$ whenever
	\[
	\sup_{\theta_1,\ldots,\theta_s\in[0,1]}
	\|w^{-2}T_{(s-\theta_1-\cdots-\theta_s)h}^Xa\|_{L_t^\rho L_x^qL_v^2}<\infty
	\quad\text{and}\quad
	\|w^{A_\diamond+2}b\|_{L_t^{\rho'}L_x^pH_\sigma^1}<\infty .
	\]
	Indeed, taking $\theta_1=\cdots=\theta_s=1$ in the supremum gives
	$w^{-2}a\in L_x^qL_v^2\subset L^2_{x,v}$; moreover, $w^{2s+3}b\in L_x^2H_\sigma^1$ by
	\eqref{eq:real-weight-Hsigma-monotonicity}, since $p\ge2$ and $2s+3\le A_\diamond+2$.
	
	By \eqref{eq:collision-mth-endpoint} at $[a,b]=[\mathfrak D_{J_S}g,Z_{J_{S^c}}]$, $S\subset[|J|]$, each defect term of $\mathscr R^{\mathscr C_+}_{r,h}$ is a multiple of
	\[
	\mathscr C_{s,h}^X\bigl[T_{(r-s-k_1)h}^X(\Delta_h^X)^{k_1}\mathfrak D_{J_S}g,\ (\Delta_h^X)^iZ_{J_{S^c}}\bigr],
	\qquad s\ge1,\quad k_1+i+s=r ,
	\]
	with the background factor $(\Delta_h^X)^{k_1}\mathfrak D_{J_S}g$, of complete order $k:=|S|+k_1$, and the solution factor $(\Delta_h^X)^iZ_{J_{S^c}}$, of complete order
	\[
	j':=|J|-|S|+i=m-k-s\le m-k-1.
	\]
	In the bound of $\mathscr C_{s,h}^X$ above, by the group law, the first argument enters through the flow translates
	\[
	T_{(s-\theta_1-\cdots-\theta_s)h}^XT_{(r-s-k_1)h}^X(\Delta_h^X)^{k_1}\mathfrak D_{J_S}g
	=T_{(r-k_1-\theta_1-\cdots-\theta_s)h}^X(\Delta_h^X)^{k_1}\mathfrak D_{J_S}g,
	\]
	where the flow parameter satisfies
	\[
	|(r-k_1-\theta_1-\cdots-\theta_s)h|\le(r-k_1)|h|\le2\mathrm Mh_0-k_1|h|.
	\]
	The differences of $\mathfrak R_J[g,Z_0]$ below have factors of the same kind.
	
	At $k=0$, that is, at $S=\varnothing$ and $k_1=0$, we take
	$(\rho,q)=(\infty,\infty)$ and $(\rho',p)=(2,2)$: the solution factor $(\Delta_h^X)^{r-s}Z_J$ is bounded in $L_t^2L_x^2H_\sigma^1$ with the weight $w^{A_\diamond+2}$ by \eqref{eq:dyn-step-discrete-lower}, since $r-s\le r-1$, and the background factor $g$, through the translates $T_{(r-\theta_1-\cdots-\theta_s)h}^Xg$, is bounded in $L_t^\infty L_x^\infty L_v^2$ with the weight $w^{-2}$, uniformly in $\theta_1,\ldots,\theta_s\in[0,1]$, by
	\eqref{eq:dyn-step-background-zero-shifts}, since $r|h|\le2\mathrm Mh_0$.  At $k\ge1$ the solution factor is bounded at
	the weight $A_\diamond+2$ by Lemma~\ref{lem:dyn-hier-weight-trade} at the ratio
	$\theta=1/m$, between the bound of part~(v) for that factor and
	\eqref{eq:dyn-step-discrete-lower} at the weight $m(A_\diamond+2)$, which uses (\ref{it:dyn-package-graph}) of $\mathsf P^\dy_{m-1}$ at the weight
	\[
	m(A_\diamond+2)+1\le(m+1)A_\diamond\le\Lambda_{m-1},
	\]
	where the first inequality holds since $2m+1\le A_\diamond$ and the second follows from \eqref{eq:dyn-realization-weights}.  The exponents $\widetilde r$, $\widetilde p$ of \eqref{eq:weight-trade-output-exponents} satisfy, by \eqref{eq:dyn-dynamic-rungs} and
	\eqref{eq:GN-output-exponents},
	\begin{equation}
		\label{261016992}
		\begin{aligned}
			\frac12-\frac1{\widetilde r}
			&=\Bigl(1-\frac1m\Bigr)\Bigl(\frac12-\frac1{\mathsf r_\dy(j'\mathrm M/m)}\Bigr)
			=\Bigl(1-\frac1m\Bigr)\frac{m-j'}{2m}
			\ge\frac k{2m}>\frac k{\mathrm Mq_\flat}=\frac1{\mathsf g_\dy(k)},\\
			\frac12-\frac1{\widetilde p}
			&=\Bigl(1-\frac1m\Bigr)\Bigl(\frac12-\frac1{\mathsf s_\dy(j'\mathrm M/m)}\Bigr)
			=\Bigl(1-\frac1m\Bigr)\frac{m-j'}m\Bigl(\frac12-\frac1{p_\flat}\Bigr)\\
			&\ge\frac km\Bigl(\frac12-\frac1{p_\flat}\Bigr)>\frac k{\mathrm Mq_\flat}=\frac1{\mathsf g_\dy(k)},
		\end{aligned}
	\end{equation}
	since $m-j'\ge k+1$, $(1-\frac1m)(k+1)\ge k$ for $k\le m-1$, $m\le\mathrm M$ and
	$\frac1{q_\flat}<\frac12-\frac1{p_\flat}<\frac12$ by \eqref{eq:dyn-own-order-endpoints}.  At $k\ge1$ we therefore take
	$(\rho,q)=(\mathsf g_\dy(k),\mathsf g_\dy(k))$ and $(\rho',p)=(\widetilde r,\widetilde p)$.  The background factors, through the translates $T_{(r-k_1-\theta_1-\cdots-\theta_s)h}^X(\Delta_h^X)^{k_1}\mathfrak D_{J_S}g$, are bounded in $L_t^{\mathsf g_\dy(k)}L_x^{\mathsf g_\dy(k)}L_v^2$ with the weight $w^{-2}$, uniformly in $\theta_1,\ldots,\theta_s\in[0,1]$, by
	\eqref{eq:dyn-step-background-differences} at $K=J_S$ and $\bar h=(r-k_1-\theta_1-\cdots-\theta_s)h$, and the solution factors $(\Delta_h^X)^iZ_{J_{S^c}}$ are bounded in $L_t^{\widetilde r}L_x^{\widetilde p}H_\sigma^1$ with the weight $w^{A_\diamond+2}$ as just described. The bound of part~(v) enters \eqref{eq:dyn-hier-weight-trade} with the power $1-\theta$, and $(1+E_h^{\,i/r})^{1-\theta}\le1+E_h^{\,i/r}$.
	
	The differences of $\mathfrak R_J[g,Z_0]$ are bounded as follows.  By \eqref{eq:deformation-remainder-structure} at the word length $|J|$, every term of $\mathfrak R_J[g,Z_0]$ has the form
	$\Phi=pN_{\mathsf K_A}[p_1\mathfrak D_{J_1}g,p_2Z_{J_2}]$ with $|J_1|+|J_2|\le|J|-1$, where $J_2$ is a proper subword of $J$ by part~(iv), and
	\[
	\|(\Delta_h^X)^r\Phi\|_{L_t^2L_x^2H_\sigma^{-1}}
	\le C\max_{k_1+i\le r}\ \sup_{|\bar h|\le r|h|}
	\|w^{-2}T_{\bar h}^X(\Delta_h^X)^{k_1}\mathfrak D_{J_1}g\|_{L_t^\rho L_x^qL_v^2}\,
	\|w^{A_\diamond+2}(\Delta_h^X)^iZ_{J_2}\|_{L_t^{\rho'}L_x^pH_\sigma^1},
	\]
	with the exponents chosen separately for each $(k_1,i)$.  The background factor $(\Delta_h^X)^{k_1}\mathfrak D_{J_1}g$ has the complete order $k:=|J_1|+k_1$, and the solution factor $(\Delta_h^X)^iZ_{J_2}$ the complete order $|J_2|+i\le m-k-1$; since $|\bar h|+k_1|h|\le2r|h|\le2\mathrm Mh_0$, both are bounded as those of the defect terms, at the same pairs.  The inequality above is the bound at $\ell=0$ and $\mathrm L_{\rm G}=-2$ given by \eqref{eq:real-weight-Hsigma-monotonicity} and by
	Corollary~\ref{cor:deformation-term-differences} in its time-integrated form, with $I_T^-$ in place of $(0,T)$ and with the jet in the solution slot at the weight \[
	2(|J|+r)+4=2m+4\le A_\diamond+2 .
	\]  Hence, at these pairs, the bound of $\|\mathscr C_{s,h}^X[a,b]\|_{L_t^2L_x^2H_\sigma^{-1}}$ above, with the bounds just given for its two factors and for the differences of $\mathfrak R_J[g,Z_0]$, gives
	\begin{equation}
		\label{eq:dyn-step-deformation}
		\|\mathscr R^{\mathscr C_+}_{r,h}\|_{L_t^2L_x^2H_\sigma^{-1}}
		\le C+C\sum_{i=1}^{r-1}E_h^{\,i/r},
		\qquad 0<|h|<h_0 .
	\end{equation}
	
	We bound the summands of $\mathscr R^{\mathscr T_L}_{r,h}$ one at a time.  By \eqref{eq:direct-operator-cube}, each
	operator quotient in $\mathscr R^{\mathscr T_L}_{r,h}$ is an average of
	$T_{\bar h}^X\operatorname{ad}_{X^sI}(\mathscr T_L)T_{-\bar h}^X$ with $I=\varnothing$ or $I=J_S$,
	the word $X^sI$ having length at most $m$.  Since, for every operator $B$,
	\[
	T_{\bar h}^XBT_{-\bar h}^XT_{sh}^X=T_{sh}^XT_{\bar h-sh}^XBT_{-(\bar h-sh)}^X,
	\]
	each operator quotient composed with $T_{sh}^X$ is bounded at almost every $t$ from $w^{2m+5}L_x^2H_\sigma^1$ into
	$L_x^2H_\sigma^{-1}$ by \eqref{eq:complete-flow-Hminus1-bound} and
	\eqref{eq:A0-orbit-bound} at $\ell=1$.  The operators $\operatorname{ad}_{J_S}(\mathscr T_L)$ are
	bounded in the same way from $w^{2m+4}L_x^2H_\sigma^1$ by \eqref{eq:A0-orbit-bound} at
	$\bar h=0$.  The jets in the solution slots of $\mathscr R^{\mathscr T_L}_{r,h}$ are, up to flow translates, the differences $(\Delta_h^X)^iZ_{J'}$ with $J'$ a subword of $J$ and $|J'|+i\le m-1$ (see the end of part~(iv)); we bound them by
	\eqref{eq:dyn-step-discrete-lower} at the weight $2m+5$, and integrate in $t$ to obtain
	\begin{equation}
		\label{eq:dyn-step-linear}
		\sup_{0<|h|<h_0}
		\|\mathscr R^{\mathscr T_L}_{r,h}\|_{L_t^2L_x^2H_\sigma^{-1}}
		\le C .
	\end{equation}
	
	The source term is $(\Delta_h^X)^r\mathfrak D_JF$ with $F=\bar F_0$. Since $h_0$ is at most the window of Lemma~\ref{lem:dynamic-fixedN-incoming-flow-bounds}\,(\ref{it:dyn-incoming-differences}) and $|J|+r=m\le\mathrm M$, Lemma~\ref{lem:dyn-breg-source-coefficient} at $j=0$ gives
	\begin{equation}
		\label{eq:dyn-step-source}
		\sup_{0<|h|<h_0}
		\|(\Delta_h^X)^r\mathfrak D_JF\|_{L_t^2L_x^2H_\sigma^{-1}}
		\le C.
	\end{equation}
	
	Finally, since
	$r|h|\le\mathrm Mh_0$ and $\|Z_h\|_{L_t^2L_x^2H_\sigma^1}\le E_h$, the definition of $\omega_g$ and \eqref{260911445}, applied
	pointwise in $(t,x)$ and integrated, give, with a universal constant $C$,
	\begin{equation}
		\label{261016991}
		\bigl\|N[T_{rh}^Xg-g,\,Z_h]\bigr\|_{L_t^2L_x^2H_\sigma^{-1}}
		\le C\,\omega_g(\mathrm Mh_0)\,E_h .
	\end{equation}
	
	\medskip
	\noindent\emph{(viii) The bound uniform in $h$.}
	The initial trace and incoming boundary data of $Z_h$ are controlled separately.
	By \eqref{eq:dyn-step-flow-slices} at $t=-1$, $Z_h(-1)=(\Delta_h^X)^r\bigl(Z_J(-1)\bigr)=0$, since $Z_J(-1)=0$ by (\ref{it:dyn-package-graph}) of $\mathsf P^\dy_{m-1}$.  By \eqref{eq:dyn-step-flow-membership},
	\eqref{eq:direct-same-flow-actual-traces} and (\ref{it:dyn-package-graph}) of $\mathsf P^\dy_{m-1}$,
	\[
	(Z_h)_{\gamma^-}=(\Delta_h^{X,\gamma_-})^r(Z_J)_{\gamma^-}=(\Delta_h^{X,\gamma_-})^r\mathfrak D_J^-H_N ,
	\]
	and \eqref{eq:dynamic-fixedN-incoming-difference-bound} of
	Lemma~\ref{lem:dynamic-fixedN-incoming-flow-bounds}\,(\ref{it:dyn-incoming-differences}) at the word $J$ and the difference order $r$, since $|J|+r=m\le\mathrm M$, gives
	\[
	\sup_{0<|h|<h_0}
	\|(Z_h)_{\gamma^-}\|_{L^2(\gamma^-_{T,\ast})}\le C .
	\]
	Decrease $h_0$ further so that $C\,\omega_g(\mathrm Mh_0)\le\tfrac14$, where $C$ is the constant of Theorem~\ref{260322442}\,(1) times that of \eqref{261016991}.
	Applying Theorem~\ref{260322442}\,(1) with background $g$, which satisfies its coefficient hypothesis by \eqref{eq:dyn-real-coefficient-hypotheses} and \eqref{eq:dynamic-KS-perturbed-smallness-choice}, to
	\eqref{eq:dyn-step-endpoint-equation}, with $Z_h(-1)=0$, gives
	\[
	E_h\le C\bigl(
	\|\mathscr R^{\mathscr C_0}_{r,h}\|+\|\mathscr R^{\mathscr C_+}_{r,h}\|
	+\|\mathscr R^{\mathscr T_L}_{r,h}\|
	+\|N[T_{rh}^Xg-g,Z_h]\|
	+\|(\Delta_h^X)^r\mathfrak D_JF\|
	+\|(Z_h)_{\gamma^-}\|_{L^2(\gamma^-_{T,\ast})}\bigr),
	\]
	the norms without subscript being $L_t^2L_x^2H_\sigma^{-1}$ norms.  By
	\eqref{eq:dyn-step-material} with \eqref{eq:dyn-step-GN}, and by \eqref{eq:dyn-step-deformation},
	\[
	\|\mathscr R^{\mathscr C_0}_{r,h}\|+\|\mathscr R^{\mathscr C_+}_{r,h}\|\le C+C\sum_{i=1}^{r-1}E_h^{\,i/r}.
	\]
	Moreover,
	\[
	\|\mathscr R^{\mathscr T_L}_{r,h}\|+\|(\Delta_h^X)^r\mathfrak D_JF\|+\|(Z_h)_{\gamma^-}\|_{L^2(\gamma^-_{T,\ast})}\le C
	\]
	by
	\eqref{eq:dyn-step-linear}, \eqref{eq:dyn-step-source} and the trace bound above.  Finally, by \eqref{261016991} and this choice of $h_0$, the term $C\|N[T_{rh}^Xg-g,Z_h]\|$ of the bound on $E_h$ above is at most $\frac14E_h$.  Hence
	\begin{equation}
		\label{eq:dyn-step-preabsorption}
		E_h\le C+C\sum_{i=1}^{r-1}E_h^{\,i/r}
		{}+\tfrac14E_h .
	\end{equation}
	Since $i<r$ in each term of the sum and $r\le\mathrm M$ is fixed,
	Young's inequality absorbs the sublinear powers into the left side
	of \eqref{eq:dyn-step-preabsorption}.  The quantities $E_h$ are
	finite for fixed $h\ne0$ by part~(iii); hence
	\begin{equation}
		\label{eq:dyn-step-graph}
		\sup_{0<|h|<h_0}E_h\le C,
	\end{equation}
	and then \eqref{eq:dyn-step-GN} holds with $E_h$ replaced by $C$; these are the
	difference bounds of the step, used in part~(v) of the later steps.
	
	\medskip
	\noindent\emph{(ix) Identification of the limit.}
	The $h$-uniform estimate permits weak graph closure, which identifies the limit as $Z_{XJ}$; the interpolation bound stated after \eqref{eq:dyn-step-input} then follows.
	Lemma~\ref{lem:complete-flow-graph-closure}, in the evolutionary case $\Omega=(0,T+1)$ of
	\eqref{eq:complete-jet-closure-base},
	applies to $q=Z_J$ at the order $r$, its hypothesis being \eqref{eq:dyn-step-graph}.  It gives
	$\mathfrak D_X^rZ_J\in\cH_{\kin,\sigma}(\cD_{T,\ast}\times\bR_v^3)$.  Moreover, as $h\to0$, the whole family $Z_h$ converges weakly to $\mathfrak D_X^rZ_J$ in this space, the traces $(Z_h)_{\gamma^\pm}$ converge weakly to those of $\mathfrak D_X^rZ_J$, and, by
	\eqref{eq:complete-flow-graph-closure-endpoints}, the time slices $Z_h(t)$,
	$-1\le t\le T$, converge weakly to $(\mathfrak D_X^rZ_J)(t)$ in $L^2(\cD\times\bR_v^3)$.  Thus
	\begin{equation}
		\label{eq:dyn-step-limit}
		Z_{X^rJ}=\mathfrak D_X^rZ_J
		\in\cH_{\kin,\sigma}(\cD_{T,\ast}\times\bR_v^3),
		\qquad
		(Z_{X^rJ})_{\gamma^\pm}
		=\operatorname*{w-lim}_{h\to0}(Z_h)_{\gamma^\pm} .
	\end{equation}
	By part~(viii), $(Z_h)_{\gamma^-}=(\Delta_h^{X,\gamma_-})^r\mathfrak D_J^-H_N$, which
	converges to $\mathfrak D_{X^rJ}^-H_N$ by
	\eqref{eq:dynamic-fixedN-incoming-difference-limit} at the order $r$ and the word
	$J$; so $(Z_{X^rJ})_{\gamma^-}=\mathfrak D_{X^rJ}^-H_N$.  The slices at $t=-1$ give
	\[
	Z_{X^rJ}(-1)=\operatorname*{w-lim}_{h\to0}Z_h(-1)=0 ,
	\] by part~(viii).
	
	For $r\ge2$, we apply
	Lemma~\ref{lem:complete-flow-weak-identification} and
	\eqref{eq:complete-flow-weak-identification-lsc} in the same evolutionary
	case to $q=Z_J$, at the difference order $1$ and at the pair
	\[
	(\mathsf r_\dy((|J|+1)\mathrm M/m),\mathsf s_\dy((|J|+1)\mathrm M/m));
	\]
	they turn
	\eqref{eq:dyn-step-GN} at $i=1$, with $E_h$ replaced by $C$, into
	\[
	\|Z_{XJ}\|_{L_t^{\mathsf r_\dy((|J|+1)\mathrm M/m)}L_x^{\mathsf s_\dy((|J|+1)\mathrm M/m)}H_\sigma^1}\le C.
	\]
	At $r=1$, where $|J|+1=m$, this pair is $(\mathsf r_\dy(\mathrm M),\mathsf s_\dy(\mathrm M))=(2,2)$, and the bound $\|Z_{XJ}\|_{L_t^2L_x^2H_\sigma^1}\le C$ follows from \eqref{eq:dyn-step-graph} and \eqref{eq:dyn-step-limit}.
	
	\medskip
	The case $r=1$, for all letters $\nu$ and words $J$ of length $m-1$, gives every spatial jet $Z_I$ of order $m$, and hence
	(\ref{it:dyn-package-graph}) of $\mathsf P^\dy_m$ at the weight $0$.
	
	\medskip
	\noindent\textbf{Closure of $\mathsf P^\dy_m$.}
	For words of length $m$, we establish the differentiated equation and forcing bound; the graph estimate then follows.
	
	\emph{(\ref{it:dyn-package-equation}).}
	The chain \eqref{eq:dyn-differentiated-chain} follows from $\mathsf P^\dy_{m-1}$ and \eqref{eq:dyn-step-limit} at $r=1$. Its lower-order weighted bound \eqref{eq:dyn-differentiated-lower-weight} follows from the induction hypothesis at $L=A_\diamond\le\Lambda_{m-1}$, by \eqref{eq:dyn-realization-weights}. For the first slots, \eqref{eq:dyn-real-background-norms} and Lemma~\ref{lem:dyn-breg-background}\,(\ref{it:dyn-breg-bg-norms}), together with $\mu^{1/4}\le Cw^{-1/2}$ and \eqref{eq:weighted-frame-monotonicity} at $\mathrm L'=0$, give \eqref{eq:dyn-differentiated-first-slots} for $g=g^\varepsilon$.
	At a spatial word $I$, the proof of Proposition~\ref{prop:dyn-differentiated-equation} uses these hypotheses only at spatial words and relies on \eqref{eq:weak-transport-word}--\eqref{eq:weak-collision-word}.  It therefore applies at
	$(f_\st,g_\dy)=(0,g)$, with $f=Z_0$ in these identities, to $\mathscr L_g Z_0=F$ on $I_T^-$.  By the last clause of the proposition, the extra
	source $F$ contributes the single summand $\mathfrak D_IF$.  Hence, with $\mathfrak K_IZ_0$ expanded by \eqref{eq:weak-triangular-commutator},
	$\mathscr L_g Z_I=-\mathfrak K_IZ_0+\mathfrak D_IF$ for $|I|\le m$, which is
	(\ref{it:dyn-package-equation}) of $\mathsf P^\dy_m$.
	
	\emph{Reduction to the forcing.}
	Fix an integer $0\le L\le\Lambda_m$ and an ordered complete spatial word $I$ of
	length $m$.  By part~(ix), $Z_I\in\cH_{\kin,\sigma}(\cD_{T,\ast}\times\bR_v^3)$,
	$(Z_I)_{\gamma^-}=\mathfrak D_I^-H_N$ and $Z_I(-1)=0$.  By
	\eqref{eq:dynamic-fixedN-incoming-word-bound}, since $L\le\Lambda_m\le\Lambda_0$,
	$\|w^L(Z_I)_{\gamma^-}\|_{L^2(\gamma^-_{T,\ast})}<\infty$, and
	Lemma~\ref{lem:dynamic-zero-initial-source-inflow-lift}\,(\ref{it:dyn-breg-lift-bounds}) with
	\eqref{eq:Hsigma-dual-weighted} bounds $w^L\mathfrak D_IF$ in
	$L_t^2L_x^2H_\sigma^{-1}$.  Hence Theorem~\ref{thm:weighted-steady}\,(1), in the form \eqref{eq:weighted-graph-convention},
	applied to $\mathscr L_g Z_I=-\mathfrak K_IZ_0+\mathfrak D_IF$, gives
	\begin{equation}
		\label{eq:dyn-order-m-graph-reduction}
		\|w^LZ_I\|_{\cH_{\kin,\sigma}(\cD_{T,\ast}\times\bR_v^3)}
		\le C\bigl(\|w^L\mathfrak K_IZ_0\|_{L_t^2L_x^2H_\sigma^{-1}}+1\bigr)
	\end{equation}
	whenever $\|w^L\mathfrak K_IZ_0\|_{L_t^2L_x^2H_\sigma^{-1}}<\infty$.  By (\ref{it:dyn-package-equation}) of
	$\mathsf P^\dy_m$, verified above, and \eqref{eq:weak-triangular-commutator} at the word $I$,
	in the form \eqref{eq:dyn-step-source-split},
	\begin{equation}
		\label{eq:dyn-order-m-forcing}
		-\mathfrak K_IZ_0
		=\sum_{\varnothing\ne S\subset[m]}
		N[\mathfrak D_{I_S}g,Z_{I_{S^c}}]
		+\mathfrak R_I[g,Z_0]
		-\mathfrak G_IZ_0 ,
	\end{equation}
	every solution factor on the right being a jet of a word of length at most $m-1$.
	Each of the three groups on the right is estimated in $w^LL_t^2L_x^2H_\sigma^{-1}$.
	
	\emph{The commutators $\mathfrak G_IZ_0$.}
	The term $\mathfrak G_IZ_0$ collects the transport, local and compact
	commutators, in which the jets in the solution slots are measured in $L_t^2L_x^2H_\sigma^1$ or $L_t^2L_{x,v}^2$.  The three kinds are bounded at almost every $t$ by \eqref{eq:ordered-transport-weight-integrated}, by
	\eqref{eq:local-L0-complete-bound} integrated in $x$, and by
	\eqref{eq:iterated-compact-full-operator-bound}, respectively, with these jets at the weights $L+4$, $L+2m+2$ and $0$.  These weights are at most $\Lambda_{m-1}$ by \eqref{eq:dyn-realization-weights}, so
	(\ref{it:dyn-package-graph}) of $\mathsf P^\dy_{m-1}$ bounds these jets, and
	integration in $t$ gives
	\[
	\|w^L\mathfrak G_IZ_0\|_{L_t^2L_x^2H_\sigma^{-1}}
	\le C\sum_{|J|<m}\Bigl(\|w^{L+4}Z_J\|_{L_t^2L_x^2H_\sigma^1}+\|w^{L+2m+2}Z_J\|_{L_t^2L_x^2H_\sigma^1}+\|Z_J\|_{L_t^2L_{x,v}^2}\Bigr)\le C ,
	\]
	the sum being taken over the ordered complete spatial words $J$ with $|J|<m$.
	
	\emph{The original-Coulomb terms.}
	Let $|S|=k\ge1$, $K:=I_S$, $I':=I_{S^c}$ and $j:=|I'|=m-k$.  By \eqref{eq:velalg},
	in the form \eqref{eq:velalg-poly}, and H\"older's inequality in $t$ and in $x$,
	\[
	\|w^LN[\mathfrak D_Kg,Z_{I'}]\|_{L_t^2L_x^2H_\sigma^{-1}}
	\le C\|w^{-2}\mathfrak D_Kg\|_{L_t^aL_x^aL_v^2}
	\|w^LZ_{I'}\|_{L_t^{\rho'}L_x^pH_\sigma^1},
	\qquad
	\frac1a+\frac1{\rho'}\le\frac12,
	\quad
	\frac1a+\frac1p\le\frac12.
	\]
	The background factor is bounded at $\mathsf g_\dy(k)$ in $t$ and in $x$, by \eqref{eq:dynamic-spatial-G-rung} and, at
	$k=\mathrm M$, by \eqref{eq:graph-to-mixed-integrability-evolutionary} with the weight $w^{-2}$.  At $j=0$ the jet in the solution slot is $Z_0$, bounded at $(\rho',p)=(\infty,p_\flat)$ by \eqref{eq:dynamic-fixedN-weighted-base-estimates} at
	$\ell=L\le\Lambda_0$ (as $p_\flat<p_\sharp$ and $\cD$ is bounded), and the pair is sub-conjugate in both variables by \eqref{eq:dyn-dynamic-conjugacy} at $j=0$.  For
	$1\le j\le m-1$ we apply Lemma~\ref{lem:dyn-hier-weight-trade} to $U=Z_{I'}$ at $\ell=L$ and a ratio $\theta$; the two factors on the right of \eqref{eq:dyn-hier-weight-trade} are bounded by \eqref{eq:dyn-step-input} with $j$ in place of $|J|$ and by (\ref{it:dyn-package-graph}) of $\mathsf P^\dy_{m-1}$ at the weight $L/\theta$.
	
	If $m\le\mathrm M-1$, take $\theta=1/(m+1)$, so that the weight is
	$(m+1)L\le\Lambda_{m-1}$.  The exponents $\widetilde r$ and $\widetilde p$ of \eqref{eq:weight-trade-output-exponents} are
	$\mathsf r_\dy((j+1)\mathrm M/(m+1))$ and $\mathsf s_\dy((j+1)\mathrm M/(m+1))$ by
	\eqref{eq:GN-output-exponents}; in time, for instance, \[
	\frac1{2(m+1)}+\Bigl(1-\frac1{m+1}\Bigr)\frac j{2m}=\frac{j+1}{2(m+1)} .
	\]
	Hence Lemma~\ref{lem:dyn-hier-weight-trade}, with its two factors bounded as above, gives
	\begin{equation}
		\label{eq:dyn-step-weighted-slot}
		\|w^LZ_{I'}\|
		_{L_t^{\mathsf r_\dy((j+1)\mathrm M/(m+1))}
			L_x^{\mathsf s_\dy((j+1)\mathrm M/(m+1))}H_\sigma^1}
		\le C ,
		\qquad 1\le j\le m-1 .
	\end{equation}
	Since $k+(j+1)=m+1\le\mathrm M$, \eqref{eq:dyn-dynamic-conjugacy} with $(m+1,j+1)$ in place of $(m,j)$ gives
	\[
	\begin{gathered}
		\frac1{\mathsf g_\dy(k)}+\frac1{\mathsf r_\dy((j+1)\mathrm M/(m+1))}=\frac12-\frac k2\Bigl(\frac1{m+1}-\frac1{\mathrm M}\Bigr)-\frac k{\mathrm M}\delta_t<\frac12,
		\\
		\frac1{\mathsf g_\dy(k)}+\frac1{\mathsf s_\dy((j+1)\mathrm M/(m+1))}=\frac12-k\Bigl(\frac12-\frac1{p_\flat}\Bigr)\Bigl(\frac1{m+1}-\frac1{\mathrm M}\Bigr)-\frac k{\mathrm M}\delta_\flat<\frac12 .
	\end{gathered}
	\]
	
	If $m=\mathrm M$, take $\theta=\theta_\flat$ of \eqref{eq:dyn-dynamic-weight-ratio}, so that
	the weight is $L/\theta_\flat\le\Lambda_{\mathrm M-1}$, and let $\widetilde r$, $\widetilde p$ be the exponents of \eqref{eq:weight-trade-output-exponents}, which satisfy
	\[
	\frac1{\widetilde p}-\frac1{\mathsf s_\dy(j)}
	=\theta_\flat\Bigl(\frac12-\frac1{\mathsf s_\dy(j)}\Bigr)=\frac k{2\mathrm M}\delta_\flat,
	\qquad
	\frac1{\widetilde r}-\frac1{\mathsf r_\dy(j)}=\frac k{\mathrm M}\cdot\frac{\theta_\flat}2 .
	\]
	Then \eqref{eq:dyn-dynamic-conjugacy} gives
	\[
	\frac1{\mathsf g_\dy(k)}+\frac1{\widetilde p}=\frac12-\frac k{2\mathrm M}\delta_\flat,
	\qquad
	\frac1{\mathsf g_\dy(k)}+\frac1{\widetilde r}
	=\frac12-\frac k{\mathrm M}\Bigl(\delta_t-\frac{\theta_\flat}2\Bigr)<\frac12 .
	\]
	In both cases $\|w^LN[\mathfrak D_Kg,Z_{I'}]\|_{L_t^2L_x^2H_\sigma^{-1}}\le C$.
	
	\emph{The deformation term.}
	Every term $\Phi$ of $\mathfrak R_I[g,Z_0]$ has the form
	$\Phi=pN_{\mathsf K_A}[p_1\mathfrak D_{J_1}g,p_2Z_{J_2}]$ of \eqref{eq:deformation-remainder-structure}, with a background factor $\mathfrak D_{J_1}g$ of order $k:=|J_1|$
	and a jet $Z_{J_2}$ of order $j:=|J_2|$ in the solution slot, $k+j\le m-1$. By
	\eqref{eq:deformation-weighted-output-spatial} at the order $m$, $\ell=L$ and $\varpi=\frac14$, which measures $\mathfrak D_{J_1}g$ and $Z_{J_2}$ without the polynomials $p_1$ and $p_2$, by $\mu^{1/4}\le Cw^{-2}$, by \eqref{eq:real-weight-Hsigma-monotonicity} with $L+2m+2\le L+A_\diamond$, and by H\"older's inequality in $t$,
	\[
	\|w^L\Phi\|_{L_t^2L_x^2H_\sigma^{-1}}
	\le C\|w^{-2}\mathfrak D_{J_1}g\|_{L_t^\rho L_x^qL_v^2}
	\|w^{L+A_\diamond}Z_{J_2}\|_{L_t^{\rho'}L_x^pH_\sigma^1}
	\]
	whenever $\frac1\rho+\frac1{\rho'}\le\frac12$ and $\frac1q+\frac1p\le\frac12$.
	At $k=0$ we take $(\rho,q)=(\infty,\infty)$ and $(\rho',p)=(2,2)$: the background factor $g$ is bounded in $L_t^\infty L_x^\infty L_v^2$ with the weight $w^{-2}$ by
	\eqref{eq:dyn-step-background-zero-shifts} at $\bar h=0$, and the jet $Z_{J_2}$ is bounded in $L_t^2L_x^2H_\sigma^1$ with the weight $w^{L+A_\diamond}$ by (\ref{it:dyn-package-graph}) of $\mathsf P^\dy_{m-1}$ at the weight $L+A_\diamond$.
	At $k\ge1$ we take $(\rho,q)=(\mathsf g_\dy(k),\mathsf g_\dy(k))$ and $(\rho',p)=(\widetilde r,\widetilde p)$, the exponents of \eqref{eq:weight-trade-output-exponents} at $\theta=1/m$ and the input pair $(\mathsf r_\dy(j\mathrm M/m),\mathsf s_\dy(j\mathrm M/m))$.  The background factor $\mathfrak D_{J_1}g$ is bounded in $L_t^{\mathsf g_\dy(k)}L_x^{\mathsf g_\dy(k)}L_v^2$ with the weight $w^{-2}$ by \eqref{eq:dynamic-spatial-G-rung}, and the jet is bounded as for the solution factors of the defect terms in part~(vii): Lemma~\ref{lem:dyn-hier-weight-trade} at $U=Z_{J_2}$, $\ell=L+A_\diamond$ and the ratio $\theta=1/m$ gives
	\[
	\|w^{L+A_\diamond}Z_{J_2}\|_{L_t^{\widetilde r}L_x^{\widetilde p}H_\sigma^1}
	\le C\|Z_{J_2}\|_{L_t^{\mathsf r_\dy(j\mathrm M/m)}L_x^{\mathsf s_\dy(j\mathrm M/m)}H_\sigma^1}^{1-1/m}
	\|w^{m(L+A_\diamond)}Z_{J_2}\|_{L_t^2L_x^2H_\sigma^1}^{1/m},
	\]
	and the two norms on the right are bounded by \eqref{eq:dyn-step-input} with $j$ in place of $|J|$ and by (\ref{it:dyn-package-graph}) of $\mathsf P^\dy_{m-1}$ at the weight $m(L+A_\diamond)\le\Lambda_{m-1}$.  Since $m-j\ge k+1$, the inequalities \eqref{261016992} with $j$ in place of $j'$ show that the pair is strictly sub-conjugate in both variables.
	
	\emph{Conclusion.}
	The three groups of terms of \eqref{eq:dyn-order-m-forcing} are therefore bounded in
	$w^LL_t^2L_x^2H_\sigma^{-1}$ by $C$, which is (\ref{it:dyn-package-forcing}) of $\mathsf P^\dy_m$.  Consequently \eqref{eq:dyn-order-m-graph-reduction} gives the bound on the graph norm of $w^LZ_I$ in
	(\ref{it:dyn-package-graph}) of $\mathsf P^\dy_m$, and part~(ix) gives the traces and
	the initial value required in (\ref{it:dyn-package-graph}).  This completes the proof.
\end{proof}

The temporal jets $Z_j$ are given by Proposition~\ref{prop:dyn-fixedN-realization}\,(\ref{it:dyn-real-prop-chain}). We record their equations, weighted estimates and interpolation bounds.

\begin{lemma}
	\label{lem:dyn-fixedN-temporal-words}
	\begin{enumerate}
		\item\label{it:dyn-temporal-words-identities}
		For $1\le j\le\mathrm M$,
		\begin{equation}
			\label{eq:dyn-fixedN-temporal-identities}
			\mathscr L_g Z_j=\sum_{a=1}^{j}\binom ja N\bigl[g^\varepsilon_a,Z_{j-a}\bigr]+\bar F_j,
			\qquad
			(Z_j)_{\gamma^-}=\partial_t^{\,j}H_N .
		\end{equation}
		\item\label{it:dyn-temporal-words-graph}
		For $0\le j\le\mathrm M$ and $0\le L\le\Lambda_j$,
		\[
		\|w^LZ_j\|_{\cH_{\kin,\sigma}(\cD_{T,\ast}\times\bR_v^3)}\le C,
		\qquad
		\|w^LN[g^\varepsilon_a,Z_{j-a}]\|_{L_t^2L_x^2H_\sigma^{-1}}\le C,
		\quad 1\le a\le j .
		\]
		\item\label{it:dyn-temporal-words-rung}
		For $2\le m\le\mathrm M$ and $1\le j\le m-1$,
		\begin{equation}
			\label{eq:dyn-mixed-temporal-rung}
			\|Z_j\|_{L_t^{\mathsf r_\dy(j\mathrm M/m)}L_x^{\mathsf s_\dy(j\mathrm M/m)}H_\sigma^1}
			\le C .
		\end{equation}
	\end{enumerate}
\end{lemma}

\begin{proof}
	The identities in \eqref{eq:dyn-fixedN-temporal-identities} follow from
	\eqref{eq:dyn-real-prop-equations}, \eqref{eq:dyn-real-prop-boundary}
	and the weak boundary generator of
	Proposition~\ref{prop:weak-boundary-generator}\,(\ref{it:weak-generator-evolutionary}).
	
	We prove the weighted estimates by induction on $j$.  The case $j=0$
	is $\mathsf P^\dy_0$.  Suppose the estimates hold below $j$.  For
	$1\le a\le j$ and $L\le\Lambda_j$, the velocity bilinear estimate gives
	\begin{equation*}
		\|w^LN[g^\varepsilon_a,Z_{j-a}]\|_{L_t^2L_x^2H_\sigma^{-1}}
		\lesssim
		\|w^{-2}g^\varepsilon_a\|_{L_t^\infty L_x^\infty L_v^2}
		\|w^LZ_{j-a}\|_{\cH_{\kin,\sigma}}\le C.
	\end{equation*}
	Here we use Lemma~\ref{lem:dyn-breg-background}\,(\ref{it:dyn-breg-bg-norms})
	and $\Lambda_j\le\Lambda_{j-a}$.  The source $w^L\bar F_j$ belongs
	to $L_t^2L_x^2H_\sigma^{-1}$ by
	Lemma~\ref{lem:dynamic-zero-initial-source-inflow-lift}\,(\ref{it:dyn-breg-lift-bounds}),
	and $w^L\partial_t^jH_N$ belongs to $L^2(\gamma^-_{T,\ast})$
	by \eqref{eq:dynamic-fixedN-incoming-word-bound}.  Moreover,
	$Z_j(-1)=0$ by Proposition~\ref{prop:dyn-fixedN-realization}\,
	(\ref{it:dyn-real-solution})--(\ref{it:dyn-real-prop-chain}).
	The coefficient hypothesis is \eqref{eq:dyn-real-coefficient-hypotheses},
	with \eqref{eq:dynamic-KS-perturbed-smallness-choice}.
	The weighted graph estimate
	\eqref{eq:weighted-graph-convention} applied to
	\eqref{eq:dyn-fixedN-temporal-identities} now gives
	$\|w^LZ_j\|_{\cH_{\kin,\sigma}}\le C$.
	This closes the induction and proves
	(\ref{it:dyn-temporal-words-graph}).
	
	For the last assertion fix $2\le m\le\mathrm M$ and $1\le j<m$.
	Apply Lemma~\ref{lem:temporal-letter-unequal} to $Z_0$ on $I_T^-$
		with $(\rho_{\mathrm{in}},p_0)=(\infty,p_\sharp)$, $a=2$,
	$\theta_\circ=\theta_\sharp$ and $H=H_\sigma^1$.
	Its derivative chain is furnished by
	Proposition~\ref{prop:dyn-fixedN-realization}\,(\ref{it:dyn-real-prop-chain});
	its base and top norms are finite by
	Lemma~\ref{prop:dyn-real-fixedN-strong} and the weighted estimates just proved.
	The hypotheses on $\gamma=j/m+\theta_\sharp$ follow from
	$\theta_\sharp(p_\flat-1)\le 1/(2\mathrm M)$ and $j<m$,
	as fixed in \eqref{eq:dyn-finite-base}.
	The temporal output exponent is $2m/j=\mathsf r_\dy(j\mathrm M/m)$.
	Writing $p_1$ for the spatial output exponent, we have
	\[
	\frac1{p_1}=\frac{1-j/m-\theta_\sharp}{p_\sharp}
	+\frac{j/m+\theta_\sharp}{2},
	\qquad
	\frac1{\mathsf s_\dy(j\mathrm M/m)}-\frac1{p_1}
	\ge \frac{1/m-\theta_\sharp(p_\flat-1)}{2p_\flat}>0.
	\]
	Thus $p_1>\mathsf s_\dy(j\mathrm M/m)$, and the boundedness of
	$\cD$ gives \eqref{eq:dyn-mixed-temporal-rung}.
\end{proof}

\begin{lemma}
	\label{lem:dyn-fixedN-mixed-words}
	Let $1\le m\le\mathrm M$.  If $\mathsf P^\bullet_{m-1}$ of Definition~\ref{def:dyn-package} holds, then $\mathsf P^\bullet_m$ holds.
\end{lemma}

\begin{proof}
	By normal ordering, $Z_{I'\mathrm t^j}=\mathfrak D_{I'}Z_j$.
	For $I=\mathrm t^m$, Lemma~\ref{lem:dyn-fixedN-temporal-words}
	gives the equation, incoming trace and graph estimate.  Its forcing
	bound includes the commutator
	\[
	[\partial_t^m,L_g]Z_0
	=-\sum_{a=1}^m\binom ma N[g_a^\varepsilon,Z_{m-a}],
	\]
	while $\mathfrak D_{\mathrm t^m}F=\bar F_m$ by
	\eqref{eq:dynamic-zero-initial-extended-source-chain}.
	
	Consider a word containing a spatial letter, written $I=X^rJ$,
	where $X=\tau_\nu$, $r\ge1$ and $|J|=m-r$.  We use the decreasing-$r$
	induction of Proposition~\ref{prop:dyn-order-m-realization}.
	The initial rung is supplied by the preceding spatial step if $J$
	has a spatial letter, and by \eqref{eq:dyn-mixed-temporal-rung}
	otherwise.  The interpolation exponents and their strict H\"older
	gaps are those of \eqref{eq:dyn-step-exponent} and
	\eqref{eq:dyn-dynamic-conjugacy}, with the total length of the word.
	The equation of $Z_I$ is obtained as in the verification of
	(\ref{it:dyn-package-equation}) in Proposition~\ref{prop:dyn-order-m-realization}.
	For the temporal letters, the hypotheses of
	Proposition~\ref{prop:dyn-differentiated-equation} are supplied by
	Proposition~\ref{prop:dyn-fixedN-realization}\,(\ref{it:dyn-real-prop-chain}) and
	Lemmas~\ref{lem:dyn-fixedN-temporal-words}\,(\ref{it:dyn-temporal-words-graph})
	and~\ref{lem:dyn-breg-background}\,(\ref{it:dyn-breg-bg-norms}).
	The spatial commutators are unchanged: $\mathscr T_L$ is independent of $t$,
	and $\operatorname{ad}_{J_S}(\mathscr T_L)=0$ if $J_S$ contains a temporal
	letter. It remains to control the differentiated background in the
	collision and deformation terms.
	
	Such a factor is $\mathfrak D_Kg_a^\varepsilon$ with $a\ge1$ and
	$K$ spatial.  Put $k=k_1+|K|+a\le m$, where $k_1$ is the number of
	spatial flow differences falling on this factor.  The convolution
	estimate \eqref{eq:dyn-breg-background}, the complete-flow bound and
	\eqref{261016990} yield, for $k_1=0$,
	\begin{equation}\label{eq:dyn-mixed-background-rung}
		\sup_{|s|\le2\mathrm Mh_0}
		\|w^{-2}T_s^X\mathfrak D_Kg_a^\varepsilon\|
		_{L_t^{\mathsf g_\dy(k)}L_x^{\mathsf g_\dy(k)}L_v^2}
		\le C_{T,\varepsilon}\left(1+
		\sup_{|s|\le2\mathrm Mh_0}
		\|w^{-2}T_s^X\mathfrak D_K\bar g_0\|
		_{L_t^{\mathsf g_\dy(|K|)}L_x^{\mathsf g_\dy(|K|)}L_v^2}
		\right).
	\end{equation}
	The right-hand side is finite by
	\eqref{eq:dyn-step-background-differences} for $|K|\ge1$ and by
	\eqref{eq:dyn-step-background-zero-shifts} for $K=\varnothing$.
	Here $|K|\le m-1$.  The passage to the smaller time exponent
	$\mathsf g_\dy(k)$ uses the finite length of $I_T^-$; no bound
	uniform in $T$ or $\varepsilon$ is asserted.
	
	For $k_1\ge1$, put $K_b=X^bK$, $0\le b\le k_1$.  Since
	$|K_b|+a\le m\le\mathrm M$ and $a\ge1$, we have
	$|K_b|\le\mathrm M-1$.  Applied to each $K_b$, the preceding
	estimate gives $\mathfrak D_{K_b}g_a^\varepsilon\in
	L^1_{\mathrm{loc}}(I_T^-;L^2_{\mathrm{loc}})$, since
	$\mathsf g_\dy(|K_b|+a)\ge q_\flat>2$.
	To identify the derivatives, use
	\eqref{eq:dyn-real-mollified-objects}.  The distributional identities
	$\mathfrak D_X\mathfrak D_{K_b}\bar g_0=
	\mathfrak D_{XK_b}\bar g_0$ hold for the ordered spatial words.
	The corresponding identities hold distributionally for the finite-jet seed.
	The spatial fields are time independent, hence commute with convolution
	in $t$ in the sense of distributions.  Consequently,
	$\mathfrak D_X(\mathfrak D_{K_b}g_a^\varepsilon)=
	\mathfrak D_{XK_b}g_a^\varepsilon$ on $I_T^-\times\cD\times\bR^3$.
	The local integrability just proved permits temporal slicing of these
	identities.  The repeated fundamental theorem of calculus along the
	complete flow now yields, at almost every time,
	\[
	T_s^X(\Delta_h^X)^{k_1}\mathfrak D_Kg_a^\varepsilon
	=\int_{[0,1]^{k_1}}
	T_{s+(\theta_1+\cdots+\theta_{k_1})h}^X
	\mathfrak D_{X^{k_1}K}g_a^\varepsilon\dd\theta,
	\]
	for $|s|+k_1|h|\le2\mathrm Mh_0$.  Minkowski's inequality and
	\eqref{eq:dyn-mixed-background-rung}, with $K$ replaced by
	$X^{k_1}K$, yield the flow-difference estimate
	\eqref{eq:dyn-step-background-differences} for the mixed factors.
	
	The remaining estimates follow with the same pairings as in the
	spatial induction.  A background factor of complementary order
	$k\ge1$ is in $L_t^{\mathsf g_\dy(k)}L_x^{\mathsf g_\dy(k)}L_v^2$;
	the solution factor of order $j$ has the rung
	$L_t^{\mathsf r_\dy(j\mathrm M/m)}
	L_x^{\mathsf s_\dy(j\mathrm M/m)}H_\sigma^1$.
	Their exponents are strictly sub-conjugate by
	\eqref{eq:dyn-dynamic-conjugacy}.
	For the original-Coulomb term when $m<\mathrm M$, the
	$(m+1,j+1)$ pairing in \eqref{eq:dyn-dynamic-conjugacy} is used
	with the weighted solution bound \eqref{eq:dyn-step-weighted-slot}.
	When $m=\mathrm M$, its weighted factor is obtained from
	Lemma~\ref{lem:dyn-hier-weight-trade} with ratio $\theta_\flat$;
	the deformation remainder uses ratio $1/m$.
	The weights in each case are among those of \eqref{2610151131}
	and are bounded by $\Lambda_{m-1}$.
	
	The source identity
	$\mathfrak D_JF=\mathfrak D_{J''}\bar F_j$ and its flow bounds follow
	from Lemmas~\ref{lem:dynamic-zero-initial-source-inflow-lift}
	and \ref{lem:dyn-breg-source-coefficient}.
	The incoming trace bounds and strong $L^2(\gamma^-_{T,\ast})$ convergence
	of the boundary differences follow from
	Lemma~\ref{lem:dynamic-fixedN-incoming-flow-bounds}.
	The lower-order jets vanish at $t=-1$; since the spatial flows preserve
	time, their spatial difference quotients also have zero initial trace.
	The fixed-cutoff weighted graph estimate therefore applies at each
	induction step.  Its constants may depend on $T$, $N$ and
	$\varepsilon$ here.  With the collision and weight bounds above it
	closes $\mathsf P_m^\bullet$.
\end{proof}

\begin{proof}[\textbf{Proof of Proposition~\ref{prop:dyn-fixedN-realization}\,(\ref{it:dyn-real-prop-spatial})--(\ref{it:dyn-real-restriction})}]
	By Definition~\ref{def:dyn-package} and the paragraph after Lemma~\ref{prop:dyn-real-fixedN-strong}, $\mathsf P^\bullet_0=\mathsf P^\dy_0$ holds, and Lemma~\ref{lem:dyn-fixedN-mixed-words}, applied successively at $m=1,\ldots,\mathrm M$, gives $\mathsf P^\bullet_{\mathrm M}$.
	
	\emph{(\ref{it:dyn-real-prop-spatial}).}  By (\ref{it:dyn-package-graph}) of
	$\mathsf P^\bullet_{\mathrm M}$,
	$w^L\mathfrak D_If_\pt^{\varepsilon,N}\in\cH_{\kin,\sigma}(\cD_{T,\ast}\times\bR_v^3)$ for
	$|I|\le\mathrm M$ and $0\le L\le\Lambda_{|I|}$; since $\Lambda_m\ge\ell^\dy_m$ at $m=|I|$, Lemma~\ref{lem:compact-velocity-support}\,(1) with the $C^1_v$ multiplier $w^{\ell^\dy_m-\Lambda_m}$ and
	\eqref{260902137} turn this into 
	$$
	\mathsf U_m(T)[f_\pt^{\varepsilon,N}]+\mathsf W_m(T)[f_\pt^{\varepsilon,N}]<\infty,\qquad 0\le m\le\mathrm M,
	$$
	where $\mathsf U_m$ and $\mathsf W_m$ are the rows \eqref{eq:dyn-hierarchy-rows}.  Also $\mathsf U_m(T)[f_\st]+\mathsf W_m(T)[f_\st]<\infty$ by the bound $\widetilde{\mathfrak Q}_{\mathrm N}(f_\st)\le C\delta$ of \eqref{eq:dyn-theta-background}, since $\ell^\dy_m\le\ell^\dy_0=\ell_{\mathrm N}\le\ell_m$.  Hence $\mathsf U_m(T)[f_\dy^{\varepsilon,N}]+\mathsf W_m(T)[f_\dy^{\varepsilon,N}]<\infty$, where $f_\dy^{\varepsilon,N}=f_\pt^{\varepsilon,N}-f_\st$, and the traces of $\mathfrak D_If_\pt^{\varepsilon,N}$ are completion-defined by (\ref{it:dyn-package-graph}) of $\mathsf P^\bullet_{\mathrm M}$.
	
	\emph{(\ref{it:dyn-real-restriction}).}  The jets of $f_\pt^{\varepsilon,N}$ form a chain on $I_T^-$ over the word set \eqref{eq:dyn-real-spatial-words-enlarged}, by the identity $Z_I=\mathfrak D_IZ_0$ in (\ref{it:dyn-package-graph}) of $\mathsf P^\bullet_{\mathrm M}$, together with \eqref{eq:dyn-normal-ordering}.  The chain passes to the open subinterval $(0,T)\subset I_T^-$, since restriction maps $\cH_{\kin,\sigma}(\cD_{T,\ast}\times\bR_v^3)$ into $\cH_{\kin,\sigma}(\cD_T\times\bR_v^3)$: each term of \eqref{260902137} on $(0,T)$ is at most the same term on $I_T^-$.
\end{proof}

\subsubsection{The strong coordinate at fixed cutoff}\label{ss:dyn-fixedN-consequences}

The remaining strong-coordinate estimate follows from the fixed-cutoff regularity and Proposition~\ref{prop:dyn-strong}.

\begin{proof}[\textbf{Proof of Proposition~\ref{prop:dyn-fixedN-realization}\,(\ref{it:dyn-real-strong})}]
	The spatial realization has already been proved in
	Proposition~\ref{prop:dyn-fixedN-realization}\,(\ref{it:dyn-real-prop-spatial})--(\ref{it:dyn-real-restriction}); it gives
	$\cN_{\mathrm M}(T)[f_\dy^{\varepsilon,N}]<\infty$ for every fixed
	$(T,N,\varepsilon)$.  The last three terms of
	\eqref{eq:dynamic-fixedN-weighted-base-estimates}, taken at
	$\ell=A^\dy_\flat-2\le\Lambda_0$ and restricted to $(0,T)$, give
	$\cA_T[Z_0]<\infty$.  Since
	$f_\dy^{\varepsilon,N}=Z_0-f_\st$ and
	$\cA_T[f_\st]\lesssim\delta$ by Lemma~\ref{lem:dyn-fp-seed},
	\[
	\cN_{\mathrm M}(T)[f_\dy^{\varepsilon,N}]
	+\cA_T[f_\dy^{\varepsilon,N}]<\infty.
	\]
	This finiteness is established for fixed $(T,N,\varepsilon)$; no uniform hierarchy estimate is invoked.
	
	For \eqref{eq:dyn-strong-A-linear-defect}, we apply Proposition~\ref{prop:dyn-strong} at $\Lambda=A^\dy_\flat$ with $(\hat f_\pt,\hat H,\mathsf F)=(\hat Z_0,\hat H_N,\hat F)$, the hats denoting extension by zero to $t\le-1$, where $F=\bar F_0$ and $H_N$ is the incoming trace of $f_\pt^{\varepsilon,N}$ in \eqref{eq:dynamic-fixedN-strong-incoming}.  Its hypotheses on $\hat Z_0$ hold: $\hat Z_0\in\cH_{\kin,\sigma}(\cD_{T,\ast}\times\bR_v^3)$ vanishes for $t\le-3/4$ by Proposition~\ref{prop:dyn-fixedN-realization}\,(\ref{it:dyn-real-solution}), solves \eqref{eq:dyn-strong-quasilinear-regrouping} with $\mathsf F=\hat F$ by \eqref{eq:dyn-real-mollified-problem}, and $(\hat Z_0)_{\gamma^-}=\hat H_N$ by \eqref{eq:dynamic-fixedN-strong-incoming}.  	We verify \eqref{eq:dyn-strong-chain-datum} with
	\[
	\mathsf N=C(\cN_{\mathrm M}(T)[f_\dy^{\varepsilon,N}]
	+\delta+\mathsf I+\varsigma_N),
	\]
	where $C$ depends only on the fixed weights and structural data.
	Here the row norm is finite at fixed $(T,N,\varepsilon)$; its uniform bound is not yet assumed.
	\begin{itemize}
		\item On $(0,T)$, \[
		H_N-H_s=P_{w,N}(f_\dy^{\varepsilon,N})_{\gamma^+}+b_w+c_N ,
		\] by
		\eqref{eq:dyn-limits-steady-residual} and
		Lemma~\ref{lem:dyn-limits-wall}\,(\ref{it:dyn-limits-wall-form}).  The part
		$P_{w,N}(f_\dy^{\varepsilon,N})_{\gamma^+}+b_w$ satisfies, for every finite $L$,
		\[
		\begin{aligned}
			&\bigl\|w^L\bigl(P_{w,N}(f_\dy^{\varepsilon,N})_{\gamma^+}+b_w\bigr)\bigr\|_{(C^{0,1}_{t,x}\cap C^{1,1}_v)(\gamma^-_T)}\\
			&\qquad+\bigl\|w^L\bigl(P_{w,N}(f_\dy^{\varepsilon,N})_{\gamma^+}+b_w\bigr)\bigr\|_{L^2(\gamma^-_T)}
			\le C_L\bigl(\cN_{\mathrm M}(T)[f_\dy^{\varepsilon,N}]+\delta\bigr),
		\end{aligned}
		\]
		by \eqref{eq:dyn-strong-wall-C11}--\eqref{eq:dyn-strong-wall-L2} at $\cJ=J_N$ and $u=f_\dy^{\varepsilon,N}$, whose jets are its distributional derivatives on $(0,T)$ by (\ref{it:dyn-real-prop-chain}), (\ref{it:dyn-real-prop-spatial}) and (\ref{it:dyn-real-restriction}) together with \eqref{eq:dyn-fs-time-independent}, with \eqref{eq:dyn-unweighted-row-dominated}.  The corrector satisfies
		\[
		\|w^Lc_N\|_{(C^{0,1}_{t,x}\cap C^{1,1}_v)(\gamma^-_T)}+\|w^Lc_N\|_{L^2(\gamma^-_T)}\le C_L\varsigma_N
		\]
		by \eqref{eq:dyn-limits-corner-bound}.  For the Lipschitz norm of
		$\hat H=H_N$ itself, we also use $\|w^LH_s\|_{(C^{0,1}_{t,x}\cap C^{1,1}_v)(\gamma^-_T)}\le C_L\delta$, which holds by \eqref{eq:stationary-incoming-regularity}, as $H_s$ does not depend on $t$.
		Together, for every finite $L$,
		\[
		\|w^LH_N\|_{(C^{0,1}_{t,x}\cap C^{1,1}_v)(\gamma^-_T)}+\|w^L(H_N-H_s)\|_{L^2(\gamma^-_T)}
		\le C_L\bigl(\cN_{\mathrm M}(T)[f_\dy^{\varepsilon,N}]+\delta+\varsigma_N\bigr) .
		\]
		\item On $(-1,0)$, $f_\pt^{\varepsilon,N}=\check{E}$ by \eqref{eq:dyn-real-prop-initial}, so
		$H_N=(\check{E})_{\gamma^-}$ by \eqref{eq:dynamic-negative-time-spectral-inflow}.  Hence, for \[
		0\le L\le A^\dy_\flat+B^{\dy}_{\rm in}+1\le\Lambda_0+B^{\dy}_{\rm in}+3 ,
		\] 
		the two norms in \eqref{eq:dyn-strong-chain-datum} satisfy
		\[
		\|w^LH_N\|_{(C^{0,1}_{t,x}\cap C^{1,1}_v)((-1,0)\times\gamma^-)}+\|w^L(H_N-H_s)\|_{L^2((-1,0)\times\gamma^-)}\le C(\delta+\mathsf I),
		\]
		by the bound on $(-1,0)$ of Lemma~\ref{lem:dyn-strong-wall-datum}.
		\item The fixed-cutoff construction proves $H_N\in H^1(I_T^-;L^2(\gamma^-))$ before the present strong-coordinate estimate; see the proof of \eqref{eq:dynamic-fixedN-incoming-strong-bound}. Thus the two one-sided time traces agree in $L^2(\gamma^-)$. Each side admits the Lipschitz representative established above, so the representatives agree at $t=0$ and glue with a Lipschitz constant controlled by the larger of the two one-sided constants. The resulting bound is independent of $T,N,\varepsilon$; the fixed-cutoff $H^1$ bound is used only to identify the common trace, not to estimate this constant.
	\end{itemize}
	The source satisfies \eqref{eq:dyn-strong-chain-source} with $\mathsf S\le C\mathsf I$, $C$
	independent of $T$, $N$ and $\varepsilon$.  Since $\bar F_0=\check{F}$ on $(-1,0)$ and
	$\bar F_0=0$ for $t>0$, this follows from
	Lemma~\ref{lem:dynamic-zero-initial-source-inflow-lift}\,(\ref{it:dyn-breg-source-bounds}) and from
	Lemma~\ref{lem:dynamic-zero-initial-source-inflow-lift}\,(\ref{it:dyn-breg-lift-bounds});
	the latter applies at this weight by the choice of $\Lambda_0$.
	
	The outgoing-moment estimates for $H_N$ from
	Lemma~\ref{lem:dyn-strong-wall-datum}, the source bound supported on
	$(-1,0)$, and the coefficient-difference estimate
	\eqref{eq:dyn-regbg-graph} have constants independent of
	$T,N,\varepsilon$. Proposition~\ref{prop:dyn-strong} therefore
	applies at the fixed weight $\Lambda=A^\dy_\flat$ with the same
	uniformity.  Restricting its conclusion to $(0,T)$ and using
	$Z_0=f_\dy^{\varepsilon,N}+f_\st$ gives
	\[
	\begin{aligned}
		\cA_T[f_\dy^{\varepsilon,N}]
		&\le\cA_T[Z_0]+\cA_T[f_\st]\\
		&\le C(\delta+\mathsf N+\mathsf S)
		\le C\bigl(\delta+\mathsf I+\varsigma_N
		+\cN_{\mathrm M}(T)[f_\dy^{\varepsilon,N}]\bigr).
	\end{aligned}
	\]
	Hence \eqref{eq:dyn-strong-A-linear-defect} is obtained before invoking Proposition~\ref{prop:dyn-hierarchy}; the latter is applied only after realization of all rows.
\end{proof}

\subsection{The removal of the two regularizations}\label{app:dyn-limit-removal}

We prove parts (\ref{it:dyn-phys-existence})--(\ref{it:dyn-phys-coords}) of Theorem~\ref{thm:dyn-physical-solution}.
Fix $T\ge1$ and a background $g$ with
$r=\widetilde{\mathfrak Q}^\dy_T(g)\le\eta_\star$.
By \eqref{eq:dyn-threshold-choice},
Proposition~\ref{prop:dyn-fixedN-realization} supplies the fixed-parameter
realization of $f_\dy^{\varepsilon,N}$, while
Proposition~\ref{prop:dyn-hierarchy} gives bounds
with constants independent of $T,N,\varepsilon$.

Proposition~\ref{prop:dynamic-inhomogeneous-physical-solvability}
\,(\ref{it:solv-physical}) gives a unique $\cH_{\kin,\sigma}$ solution $f_\dy$ of
\eqref{eq:dyn-perturbation-problem}, with background $g$,
source $\mathfrak S$, incoming datum $b_w$, and initial datum
$f_{\dy,0}$. Its hypotheses hold:
Lemma~\ref{lem:dyn-strong-coefficient-class} and
\eqref{eq:dyn-strong-coefficient-choice} give
$\|\mu^{1/4}g\|_{L^\infty}\le\varepsilon_{\rm dyn}/2$;
also $\delta\le\delta_{\rm dyn}$.
The source belongs to $L^2(\cD_T;H^{-1}_\sigma)$ by
\eqref{eq:dyn-decay-forcing}, the bound on the background rows,
and \eqref{eq:dynamic-global-stationary-second-slot-small}.
The initial and boundary data belong to their stated $L^2$
spaces by \eqref{eq:dyn-data-size} and
\eqref{eq:dyn-limits-wall-datum-size}.

\begin{prop}[Strong convergence]
	\label{prop:dyn-limits-base}
	As $(\varepsilon,N)\to(0,\infty)$ with $0<\varepsilon\le\varepsilon_\lambda$,
	\begin{equation}
		\label{eq:dyn-limits-base-strong}
		\bigl\|f_\dy^{\varepsilon,N}-f_\dy\bigr\|_{\cH_{\kin,\sigma}(\cD_T\times\bR_v^3)}
		\longrightarrow0 .
	\end{equation}
\end{prop}

\begin{proof}
	Write $r^{\varepsilon,N}=f_\dy^{\varepsilon,N}-f_\dy$
	and $\phi^{\varepsilon,N}=\cA_0(f_\dy^{\varepsilon,N})_{\gamma^+}$.
	Subtracting the two equations gives
	\begin{equation}
		\label{eq:dyn-limits-defect}
		\begin{gathered}
			(\partial_t+v\cdot\nabla_x)r^{\varepsilon,N}
			+L_gr^{\varepsilon,N}
			=N[g^\varepsilon-g,f_\pt^{\varepsilon,N}],
			\qquad r^{\varepsilon,N}(0)=0,\\
			(r^{\varepsilon,N})_{\gamma^-}
			=P_w(r^{\varepsilon,N})_{\gamma^+}
			+\mathfrak e_{\varepsilon,N},\qquad
			\mathfrak e_{\varepsilon,N}
			=(\mathfrak m_\circ+d_\delta(t))(J_N-\mathbf I)
			\phi^{\varepsilon,N}+c_N.
		\end{gathered}
	\end{equation}
	Indeed, $L_{g^\varepsilon}=L_g-N[g^\varepsilon-g,\cdot]$
	by \eqref{def:L_g}, and
	$g^\varepsilon_\dy-g_\dy=g^\varepsilon-g$.
	Thus the interior coefficient differences combine into
	$N[g^\varepsilon-g,f_\dy^{\varepsilon,N}+f_\st]$.
	The initial jets agree at order zero by
	\eqref{eq:dyn-limits-fixed-system} and
	\eqref{eq:dyn-limits-fixed-datum}.
	The truncated wall law and
	\eqref{eq:dynamic-spectral-Pw-difference-exact}, together with
	$\mathfrak t_{0,N}=0$ from \eqref{eq:dyn-theta-triangular-notation},
	give the boundary error with the sign in \eqref{eq:dyn-limits-defect}.
	
	For the interior defect, uniformly in $0<\varepsilon\le\varepsilon_\lambda$ and $N\in\bN$,
	\begin{equation}
		\label{eq:dyn-limits-eps-interior}
		\bigl\|N\bigl[g^\varepsilon-g,\,f_\pt^{\varepsilon,N}\bigr]\bigr\|
		_{L^2(\cD_T;H_\sigma^{-1})}
		\le C\Bigl(\varepsilon^{\alpha_{\rm bs}/3}\,
		\bigl[\mu^{1/4}\hat g\bigr]_{C^{\alpha_{\rm bs}}_{\kin}}
		+\varepsilon\,(\delta+\mathsf I)\Bigr)
		\|f_\pt^{\varepsilon,N}\|_{L^2(\cD_T;H_\sigma^1)} .
	\end{equation}
	Here the seminorm is that of \eqref{eq:dyn-strong-coefficient-holder}, and the constant
	is independent of $T$, $N$ and $\varepsilon$.  On $(0,T)$, where $\hat g=g$ by
	\eqref{eq:dyn-strong-past-extension}, \eqref{eq:dyn-real-increment-split} gives
	\[
	g^\varepsilon-g
	=\int_0^\varepsilon\varrho_\varepsilon(s)\bigl(\hat g(\cdot-s)-\hat g\bigr)\dd s
	+\bigl(\hat g^{(0)}-\varrho_\varepsilon*\hat g^{(0)}\bigr) .
	\]
	For $0<s\le1$, choosing the same velocity at both points in \eqref{kineticdistance} gives
	\[
	\mathbf d\bigl((t-s,x,v),(t,x,v)\bigr)
	\le\max\bigl\{s^{1/2},(s|v|)^{1/3}\bigr\}
	\le s^{1/3}(1+|v|)^{1/3}.
	\]
	The kinetic H\"older seminorm controls the first difference, while
	\eqref{eq:dyn-real-increment-bound} controls the seed increment:
	\[
	\bigl\|\mu^{3/8}(g^\varepsilon-g)\bigr\|_{L^\infty((0,T)\times\cD;L_v^2)}
	\le\bigl\|\mu^{1/8}(1+|v|)^{\alpha_{\rm bs}/3}\bigr\|_{L_v^2}\,
	\bigl[\mu^{1/4}\hat g\bigr]_{C^{\alpha_{\rm bs}}_{\kin}}\,
	\varepsilon^{\alpha_{\rm bs}/3}
	+\bigl\|\mu^{1/8}\bigr\|_{L_v^2}\,C\,\varepsilon\,(\delta+\mathsf I) .
	\]
	Applying \eqref{eq:velalg} with $(\ell,\varpi)=(0,3/8)$ at each $(t,x)$ yields
	\eqref{eq:dyn-limits-eps-interior}.   The coefficient bound \eqref{eq:dyn-strong-coefficient-output} gives
	$\bigl[\mu^{1/4}\hat g\bigr]_{C^{\alpha_{\rm bs}}_{\kin}}\le1$.  The graph bound
	\eqref{eq:dyn-real-uepsN-graph} gives $\|f_\pt^{\varepsilon,N}\|_{L^2(\cD_T;H_\sigma^1)}\le C_T$, independent of $N$
	and $\varepsilon$. Hence the source of \eqref{eq:dyn-limits-defect} satisfies
	\[
	\|N[g^\varepsilon-g,\,f_\pt^{\varepsilon,N}]\|_{L^2(\cD_T;H_\sigma^{-1})}\le C\bigl(\varepsilon^{\alpha_{\rm bs}/3}+\varepsilon(\delta+\mathsf I)\bigr),
	\]
	with $C$ depending on $T$ but not on $\varepsilon$ or $N$.
	
	For the incoming error, use
	Lemma~\ref{cor:outgoing-moment-H3}
	\,(\ref{it:outgoing-moment-H-evolutionary}),
	\eqref{eq:dyn-hierarchy-closed}, and
	\eqref{eq:dyn-unweighted-row-dominated} to obtain
	\[
	\|\phi^{\varepsilon,N}\|_{L_t^2H^{\mathrm M}(\partial\cD)}
	\le C\sum_{|I|\le\mathrm M}
	\|\mathfrak D_If_\dy^{\varepsilon,N}\|_{\cH_{\kin,\sigma}}
	\le C,
	\]
	uniformly in $T,N,\varepsilon$.
	Here Proposition~\ref{prop:dyn-fixedN-realization}
	\,(\ref{it:dyn-real-prop-chain}),
	(\ref{it:dyn-real-prop-spatial}) and
	(\ref{it:dyn-real-restriction}) identifies the regularized
	jets with the distributional complete derivatives of
	$f_\dy^{\varepsilon,N}$; the stationary part is time-independent by
	\eqref{eq:dyn-fs-time-independent}.
	The spectral multiplier estimate
	\eqref{eq:spectral-cutoff-tail-multiplier} at $k=\mathrm M$
	and Lemma~\ref{lem:word-boundary-spectral-chart} yield
	\[
	\|(\mathbf I-J_N)\phi^{\varepsilon,N}\|_{L^2((0,T)\times\partial\cD)}
	\le C N^{-\mathrm M}.
	\]
	The profile isometry
	\eqref{eq:dynamic-spectral-profile-L2-isometry-bound}
	at $L=0$ and the corner bound
	\eqref{eq:dyn-limits-corner-bound} therefore give
	\[
	\|\mathfrak e_{\varepsilon,N}\|_{L^2(\gamma_T^-)}
	\le C(N^{-\mathrm M}+\varsigma_N).
	\]
	Here $\varsigma_N\to0$ by
	Lemma~\ref{lem:dyn-limits-wall}
	\,(\ref{it:dyn-limits-wall-size}).
	
	The zero-initial-value estimate of
	Proposition~\ref{prop:dynamic-inhomogeneous-physical-solvability}
	\,(\ref{it:solv-physical}) applies to \eqref{eq:dyn-limits-defect}
	with the fixed background $g$. Its constant $C_T$ is independent of
	$N$ and $\varepsilon$.
	The two defect bounds give
	\[
	\|r^{\varepsilon,N}\|_{\cH_{\kin,\sigma}(\cD_T\times\bR_v^3)}
	\le C_T\bigl(\varepsilon^{\alpha_{\rm bs}/3}
	+\varepsilon(\delta+\mathsf I)+N^{-\mathrm M}+\varsigma_N\bigr).
	\]
	The right-hand side
	therefore tends to zero along every sequence $\varepsilon\downarrow0$,
	$N\to\infty$.  By the graph trace estimate \eqref{260902137},
	this convergence also gives
	\[
	\sup_{0\le t\le T}\|r^{\varepsilon,N}(t)\|_{L^2_{x,v}}
	+\|(r^{\varepsilon,N})_{\gamma^+}\|_{L^2(\gamma_T^+)}
	+\|(r^{\varepsilon,N})_{\gamma^-}\|_{L^2(\gamma_T^-)}\longrightarrow0.
	\]
	This proves \eqref{eq:dyn-limits-base-strong}.
\end{proof}

\begin{proof}[\textbf{Proof of Theorem~\ref{thm:dyn-physical-solution}\,(\ref{it:dyn-phys-existence})--(\ref{it:dyn-phys-coords})}]
	\emph{(\ref{it:dyn-phys-existence}).}
	Take the graph solution $f_\dy$ constructed above. Let
	$\varepsilon_n\downarrow0$, $N_n\to\infty$, and put
	$z_n=f_\dy^{\varepsilon_n,N_n}$, $a_n=g^{\varepsilon_n}$.
	By \eqref{eq:dyn-limits-base-strong}, $z_n\to f_\dy$ strongly
	in $\cH_{\kin,\sigma}(\cD_T\times\bR_v^3)$.  Hence the
	incoming and outgoing traces and the time-endpoint values converge
	strongly in $L^2$, by \eqref{260902137}.  In particular the
	truncated law
	$(z_n)_{\gamma^-}=P_w(z_n)_{\gamma^+}+b_w+
	\mathfrak e_{\varepsilon_n,N_n}$ passes to the diffuse law:
	$P_w:L^2(\gamma_T^+)\to L^2(\gamma_T^-)$ is bounded by
	\eqref{eq:wall-J-operator-bound}, and the preceding proof gives
	$\|\mathfrak e_{\varepsilon_n,N_n}\|_{L^2(\gamma_T^-)}\to0$.
	The collision term converges in $L^2(\cD_T;H_\sigma^{-1})$: with
	$F_n=f_\st+z_n$ and $F=f_\st+f_\dy$, bilinearity gives
	\[
	N[a_n,F_n]-N[g,F]
	=N[a_n-g,F_n]+N[g,z_n-f_\dy].
	\]
	The first term vanishes by \eqref{eq:dyn-limits-eps-interior}
	and \eqref{eq:dyn-real-uepsN-graph}; the second vanishes by
	\eqref{eq:velalg}, since
	$\|\mu^{3/8}g\|_{L^\infty_{t,x}L_v^2}<\infty$ by the coefficient
	hypothesis and $z_n-f_\dy\to0$ in $L^2(\cD_T;H_\sigma^1)$.
	Hence the bulk collision term converges in $H_\sigma^{-1}$, while
	the diffuse wall law converges in the flux $L^2$ trace space.
	
	The uniform jet bounds \eqref{eq:dyn-hierarchy-closed} and the
	distributional chain identities of
	Proposition~\ref{prop:dyn-fixedN-realization} together verify the hypotheses of
	Lemma~\ref{lem:dyn-limits-coordinate-closure}.   At the temporal words $I=\mathrm t^j$,
	\eqref{eq:dyn-limits-closure-jets}, together with the fourth line of
	\eqref{eq:dyn-limits-fixed-system}, gives
	$\partial_t^jf_\dy(0)=\zeta_j$ for $1\le j\le\mathrm M$.
	This proves (\ref{it:dyn-phys-existence}).
	
	\emph{(\ref{it:dyn-phys-mass}).}   Proposition~\ref{prop:dynamic-alpha-uniform-graph} at $\alpha=1$ gives
	$\mathfrak m_{f_\dy}(t)=0$ for $0\le t\le T$.  The
	smallness hypothesis \eqref{eq:dynamic-alpha-graph-smallness} of that proposition is checked in the paragraph before Proposition~\ref{prop:dyn-limits-base}. The initial
	mass of $f_{\dy,0}$ vanishes by \eqref{eq:dynamic-initial-zero-perturbation-mass}, and
	\eqref{eq:dynamic-physical-zero-mass-data} holds at $F=\mathfrak S$ and $b=b_w$, since
	$\langle N[a,b],\sqrt\mu\rangle_v=0$ and $\beta_{b_w}=0$ by \eqref{eq:dyn-theta-residual-flux}.
	
	\emph{(\ref{it:dyn-phys-coords}).}
	The convergence $z_n\to f_\dy$ in the graph norm and the uniform
	bound \eqref{eq:dyn-hierarchy-closed} meet the hypotheses of
	Lemma~\ref{lem:dyn-limits-coordinate-closure}.  Thus all weighted
	jets converge weakly to the corresponding distributional derivatives
	of $f_\dy$.  Hence \eqref{eq:dyn-limits-closure}--\eqref{eq:dyn-limits-closure-rows},
	applied to the bounds \eqref{eq:dyn-strong-closed-A} and
	\eqref{eq:dyn-hierarchy-closed}, give
	\begin{equation}
		\label{eq:dyn-limits-physical-A}
		\cA_T[f_\dy]+\cS_T[f_\dy]+\cN_{\mathrm M}(T)[f_\dy]
		\le C(\delta+\mathsf I).
	\end{equation}
	The approximating bound contains only the vanishing correction
	$C\varsigma_{N_n}$; its constant is independent of $T,N_n$ and
	$\varepsilon_n$.  In contrast, the constant used earlier for
	strong graph convergence was allowed to depend on the fixed $T$.
	
	For the decay coordinate, apply Corollary~\ref{cor:dyn-decay-B} with $\cJ=\mathbf I$, $\rho=r$, and background $g$ to
	$f_\dy$: its smallness hypothesis holds since $\rho=r\le\eta_\star\le1$
	and $\rho+\delta\le\eta_\star+\delta_\star\le\eta_{\star,0}
	\le\varepsilon_{\rm dec}$ by \eqref{eq:dyn-threshold-choice},
	$f_\dy\in\cH_{\kin,\sigma}(\cD_T\times\bR_v^3)$ by (\ref{it:dyn-phys-existence}), $f_\dy$ solves
	\eqref{eq:dyn-perturbation-problem}, and $\mathfrak m_{f_\dy}\equiv0$
	on $[0,T]$ by (\ref{it:dyn-phys-mass}).  So \eqref{eq:dyn-decay-B-closed} at $\cJ=\mathbf I$ gives
	\[
	\cB_T[f_\dy]\le C(\mathsf I^2+\delta^2)\le C(\delta+\mathsf I)^2 .
	\]
	By \eqref{eq:dyn-reduced-Q}, this bound,
	\eqref{eq:dyn-limits-physical-A} and $\cA_T[f_\st]\le C\delta$ of
	Lemma~\ref{lem:dyn-fp-seed} give
	\[
	\widetilde{\mathfrak Q}^\dy_T(f_\st+f_\dy)\le\cA_T[f_\st]+\cA_T[f_\dy]+\cB_T[f_\dy]^{1/2}+\cS_T[f_\dy]+\cN_{\mathrm M}(T)[f_\dy]\le C(\delta+\mathsf I),
	\]
	with $C$ depending only on the structural data.
\end{proof}

	\section*{Acknowledgements}
	C.K. is partially supported by NSF CAREER Grant DMS-2047681. This material is partly based upon work supported by the National Science Foundation under Grant No. DMS-2424139, while C.K. was in residence at the Simons Laufer Mathematical Sciences Institute in Berkeley, California, during the Fall 2025 semester. 
	J.S. has been supported by a KIAS Individual Grant (MG095802) at Korea Institute for Advanced Study and by the National Research Foundation of Korea (NRF) grant funded by the Korea government (Ministry of Science and ICT) (No. RS-2026-25477288). J.S. thanks the Department of Mathematics at the University of Wisconsin--Madison for its kind hospitality from March to July 2025 and in August 2026, SLMath for its kind hospitality in December 2025, and KIAS for supporting these research visits. Part of this work was carried out during these visits. 
	We thank Luis Silvestre for bringing the work of Niebel and Valentini \cite{NiVa} to our attention.
	We thank Minhyun Kim for carefully reading the manuscript and for helpful comments. The authors are also working with him on kinetic boundary problems involving fractional-order collision operators.

    {ChatGPT plus was used to assist with language editing, and Claude with proofreading. The authors take full responsibility for the content of this manuscript.}

	\appendix

	\section{Auxiliary results from stochastic analysis}\label{sec:stochastic-appendix}
	
	We collect the general probabilistic facts used in
	Subsections~\ref{subsec:grazing} and~\ref{subsec:prob-solutions}.
	
	Throughout, $(\Omega,\cF,(\cF_t)_{t\ge0},\mathbb P)$ is a
	filtered probability space satisfying the usual conditions. Unless explicitly
	stated otherwise, all processes are adapted to $(\cF_t)_{t\ge0}$, and
	every Brownian motion mentioned below is Brownian with respect to this
	filtration. Independent Brownian motions are equipped with the augmented joint
	filtration.
	
	\begin{lemma}\label{lem:BC}
		Let $(E_n)_{n\ge1}\subset\cF$.
		\begin{enumerate}
			\item[(I)] If $\sum_{n\ge1}\mathbb P(E_n)<\infty$, then
			\[
			\mathbb P\bigl(E_n\ \text{infinitely often}\bigr)=0.
			\]
			\item[(II)] If the events $E_n$ are independent and
			$\sum_{n\ge1}\mathbb P(E_n)=\infty$, then
			\[
			\mathbb P\bigl(E_n\ \text{infinitely often}\bigr)=1.
			\]
		\end{enumerate}
	\end{lemma}
	See, for example, \cite[Section~2.3]{durrett2019probability}.
	
	\begin{lemma}\label{lem:hitting_time}
		Let $(\Xi_t)_{t\ge0}$ be a continuous adapted $\bR^N$-valued process, let
		$F\subset\bR^N$ be closed, and let $q\ge0$. With the convention
		$\inf\emptyset:=+\infty$, define
		\[
		\zeta_q:=\inf\{t\ge q:\ \Xi_t\in F\},
		\qquad
		\zeta_q^+:=\inf\{t>q:\ \Xi_t\in F\}.
		\]
		Then $\zeta_q$ and $\zeta_q^+$ are stopping times, and
		\[
		\Xi_{\zeta_q}\in F\quad\text{on }\{\zeta_q<\infty\}.
		\]
		If $\Xi_q\notin F$, then $\zeta_q=\zeta_q^+$. Moreover, if
		$\cO\subset\bR^N$ is open, $\Xi_q\in\cO$, and
		\[
		\tau_{\cO}:=\inf\{t>q:\ \Xi_t\notin\cO\},
		\]
		then $\tau_{\cO}$ is a stopping time, $\tau_{\cO}>q$, and
		\[
		\Xi_{\tau_{\cO}}\in\partial\cO
		\quad\text{on }\{\tau_{\cO}<\infty\}.
		\]
	\end{lemma}
	These are the standard properties for continuous adapted processes; see
	\cite[Chapter~I]{RevuzYor} or \cite[Chapter~1]{KaratzasShreve}.
	
	\begin{lemma}[It\^o's formula]\label{lem:ito}
		Let $(B_t)_{t\ge0}$ be an $m$-dimensional Brownian motion, and let
		$\Xi_t=(\Xi_t^1,\ldots,\Xi_t^N)$ be an It\^o process satisfying
		\[
		\dd\Xi_t^i=b_t^i\dd t+
		\sum_{k=1}^m\sigma_t^{ik}\dd B_t^k,
		\qquad i=1,\ldots,N.
		\]
		Assume that $b^i$ and $\sigma^{ik}$ are progressively measurable and that,
		for every $T>0$,
		\[
		\int_0^T\left(\sum_{i=1}^N|b_s^i|
		+\sum_{i=1}^N\sum_{k=1}^m|\sigma_s^{ik}|^2\right)\dd s<\infty
		\qquad\text{almost surely}.
		\]
		Set
		\[
		a_t^{ij}:=\sum_{k=1}^m\sigma_t^{ik}\sigma_t^{jk}.
		\]
		Then, for every $\varphi\in C^{1,2}([0,\infty)\times\bR^N)$,
		\begin{align*}
			\varphi(t,\Xi_t)
			={}&\varphi(0,\Xi_0)
			+\int_0^t\partial_s\varphi(s,\Xi_s)\dd s+\sum_{i=1}^N\int_0^t
			b_s^i\,\partial_i\varphi(s,\Xi_s)\dd s\\
			&\qquad +\frac12\sum_{i,j=1}^N\int_0^t
			a_s^{ij}\,\partial_{ij}\varphi(s,\Xi_s)\dd s+\sum_{i=1}^N\sum_{k=1}^m\int_0^t
			\sigma_s^{ik}\,\partial_i\varphi(s,\Xi_s)\dd B_s^k
		\end{align*}
		almost surely for all $t\ge0$. The last term is a continuous local
		martingale.
		The identity holds simultaneously for all $t\ge0$ on an event of
		probability one. It therefore remains valid with $t$ replaced by a bounded
		stopping time.
	\end{lemma}
	See \cite[Chapter~3, Theorem~3.6]{KaratzasShreve} or
	\cite[Chapter~IV]{RevuzYor}.
	
	\begin{lemma}\label{lem:stochastic_integral}
		Let $(W_t)_{t\ge0}$ be an $m$-dimensional Brownian motion, let $H$ be an
		$\bR^m$-valued progressively measurable process with
		$\int_0^t|H_s|^2\dd s<\infty$ almost surely for every $t\ge0$, and set
		$M_t:=\int_0^tH_s\cdot\dd W_s$.
		\begin{enumerate}
			\item[(i)]  Let $\sigma$ be a stopping
			time, let $A\in\cF_0$, and let $K$ be a second integrand as above.
			If
			\[
			\mathbf 1_A\int_0^\sigma |H_s-K_s|^2\dd s=0
			\qquad\text{almost surely},
			\]
			then, almost surely on $A$,
			\[
			\int_0^tH_s\cdot\dd W_s=\int_0^tK_s\cdot\dd W_s
			\qquad(0\le t\le\sigma).
			\]
			\item[(ii)] If $\tau$ is a stopping
			time and
			\[
			\mathbb E\int_0^\tau|H_s|^2\dd s<\infty,
			\]
			then
			$(M_{t\wedge\tau})_{t\ge0}$ is a square-integrable martingale; in
			particular,
			\[
			\mathbb E[M_\tau\mid\cF_0]=0,
			\qquad
			\mathbb E[\mathbf 1_A M_\tau]=0
			\quad(A\in\cF_0).
			\]
			\item[(iii)]  If
			$\mathbb E\int_0^T|H_s|^2\dd s<\infty$, then
			\[
			\mathbb E\Bigl[\sup_{0\le t\le T}|M_t|^2\Bigr]
			\le4\,\mathbb E\int_0^T|H_s|^2\dd s .
			\]
			\item[(iv)] $M$ is a continuous
			local martingale with
			\[
			\langle M\rangle_t=\int_0^t|H_s|^2\dd s .
			\]
		\end{enumerate}
	\end{lemma}
	All four assertions are standard consequences of the construction and
	localization of the It\^o integral, the It\^o isometry, and Doob's
	inequality; see \cite[Section~3.2]{KaratzasShreve} or
	\cite[Chapter~IV]{RevuzYor}.
	
	\begin{lemma}\label{lem:brownian_shift}
		Let $(B_t)_{t\ge0}$ be an $m$-dimensional Brownian motion and let $q\ge0$.
		Then
		\[
		\widehat B_t:=B_{q+t}-B_q,
		\qquad
		\widehat{\cF}_t:=\cF_{q+t},
		\qquad t\ge0,
		\]
		defines an $m$-dimensional
		$(\widehat{\cF}_t)$-Brownian motion, independent of
		$\cF_q$. Moreover, for every admissible progressively measurable
		integrand $H$,
		\[
		\int_q^{q+t}H_s\cdot\dd B_s
		=\int_0^tH_{q+s}\cdot\dd\widehat B_s,
		\qquad t\ge0,
		\]
		almost surely.
	\end{lemma}
	This is the stationarity and independence of Brownian increments together
	with the time-shift property of the It\^o integral; see
	\cite[Section~3.2]{KaratzasShreve} or \cite[Chapter~IV]{RevuzYor}.
	
	\begin{lemma}
		\label{lem:Levy_characterization}
		Let $(M_t)_{t\ge0}$ be a real-valued continuous local martingale with
		$M_0=0$ and $\langle M\rangle_t=t$ for all $t\ge0$ almost surely. Then
		$M$ is a standard one-dimensional $(\cF_t)$-Brownian motion.
	\end{lemma}
	See \cite[Chapter~IV, Theorem~3.6]{RevuzYor} or
	\cite[Chapter~3, Theorem~3.16]{KaratzasShreve}.
	
	\begin{lemma}
		\label{lem:feynman_kac}
		Let $\cO\subset\bR^N$ be a bounded open set, let
		$b:\bR^N\to\bR^N$ and $\sigma:\bR^N\to\bR^{N\times m}$ be globally
		Lipschitz, and let $\Xi(\xi)$ be the unique non-exploding strong solution
		of
		\[
		\Xi_t(\xi)=\xi+\int_0^t b(\Xi_s)\dd s
		+\sqrt2\int_0^t\sigma(\Xi_s)\dd B_s .
		\]
		For $1\le i,j\le N$, set
		\[
		a_{ij}(\xi):=\sum_{k=1}^m\sigma_{ik}(\xi)\sigma_{jk}(\xi),
		\]
		and define
		\[
		\cA\phi
		:=\sum_{i,j=1}^N a_{ij}\,\partial_{ij}\phi
		+\sum_{i=1}^N b_i\,\partial_i\phi.
		\]
		For $\xi\in\cO$, set
		$\tau:=\inf\{s>0:\Xi_s(\xi)\notin\cO\}$, with
		$\inf\emptyset:=+\infty$.
		\begin{enumerate}
			\item[(i)] Let $T\in(0,\infty)$,
			$\lambda\ge0$, $g\in C_b([0,T]\times\overline{\cO})$,
			$f_0\in C(\overline{\cO})$, and
			$\psi\in C([0,T]\times\partial\cO)$ satisfy
			$\psi(0,\cdot)=f_0$ on $\partial\cO$. If
			\[
			u\in C^{1,2}\bigl((0,T)\times\cO\bigr)
			\cap C_b\bigl([0,T]\times\overline{\cO}\bigr)
			\]
			solves
			\[
			\partial_tu=\cA u-\lambda u+g,\qquad
			u(0,\cdot)=f_0,\qquad
			u|_{(0,T)\times\partial\cO}=\psi,
			\]
			then, for $(t,\xi)\in(0,T)\times\cO$,
			\begin{align}\label{eq:FK_evolution}
				u(t,\xi)
				=\mathbb E\Big[&
				\int_0^{\tau\wedge t}\ee^{-\lambda s}
				g\bigl(t-s,\Xi_s(\xi)\bigr)\dd s \notag\\
				&+\ee^{-\lambda\tau}
				\psi\bigl(t-\tau,\Xi_\tau(\xi)\bigr)
				\mathbf 1_{\{\tau\le t\}}
				+\ee^{-\lambda t}f_0\bigl(\Xi_t(\xi)\bigr)
				\mathbf 1_{\{\tau>t\}}\Big].
			\end{align}
			\item[(ii)] Let $\lambda>0$,
			$g\in C(\overline{\cO})$, and
			$\psi\in C(\partial\cO)$. If
			$u\in C^2(\cO)\cap C_b(\overline{\cO})$ solves
			$\cA u-\lambda u+g=0$ in $\cO$ and
			$u=\psi$ on $\partial\cO$, then
			\begin{align}\label{eq:FK_stationary}
				u(\xi)=\mathbb E\Big[
				\int_0^\tau \ee^{-\lambda s}g\bigl(\Xi_s(\xi)\bigr)\dd s
				+\ee^{-\lambda\tau}\psi\bigl(\Xi_\tau(\xi)\bigr)
				\mathbf 1_{\{\tau<\infty\}}\Big].
			\end{align}
		\end{enumerate}
	\end{lemma}
	
	The Lipschitz SDE theory is standard; see \cite[Section~5.2]{KaratzasShreve}.
	For Feynman--Kac formulae, see \cite[Section~5.7]{KaratzasShreve} and \cite{DynkinMarkov}.
	For the stopped form above, apply It\^o's formula to $\ee^{-\lambda s}u(t-s,\Xi_s)$, or to $\ee^{-\lambda s}u(\Xi_s)$ in the stationary case, stop first on compact subdomains of $\cO$, and pass to the exit time by bounded convergence.

	\section{De Giorgi argument for boundary H\"older regularity}\label{proofofholder}
	
	The two main ingredients for the proof of Proposition~\ref{2601271146} are
	the boundary isoperimetric estimate, Proposition~\ref{prop:isop},
	and the unit-scale boundary oscillation decay, Proposition~\ref{cor:boundary-oscillation-decay}.
	These are proved in \S\ref{260701421} and \S\ref{260701422},
	respectively.

	\begin{prop}[Boundary isoperimetric estimate]
		\label{prop:isop}
		Let
		\[
		Z_0=(t_0,x_0,v_0)
		\in \bR\times\overline{\bR^d_+}\times\bR^d,
		\qquad
		\widehat Z_0:=(t_0-1,x_0-v_0,v_0).
		\]
		Assume that
		$a^{ij},\,b^i\in L^\infty(Q_2^+(Z_0))$, $g\in L^\infty(Q_2^+(Z_0))$,
		and that, for some \(0<\lambda\le\Lambda\) and \(M\ge0\),
		\[
		\lambda |\xi|^2
		\le
		a^{ij}(Z)\xi_i\xi_j,
		\qquad
		|a^{ij}(Z)|\le \Lambda,
		\qquad
		|b(Z)|\le M
		\]
		for a.e. \(Z\in Q_2^+(Z_0)\) and every \(\xi\in\bR^d\).  Let
		\(f\in L^\infty(Q_2^+(Z_0))\) be a weak solution of
		\[
		(\partial_t+v\cdot\nabla_x)f
		=
		D_i(a^{ij}D_j f)+b^iD_i f+g
		\qquad\text{in }Q_2^+(Z_0),
		\]
		with
		\[
		f=h\le0
		\qquad
		\text{on }(\bR\times\gamma^-)\cap Q_2(Z_0)
		\]
		in the sense of traces.
		
		Let
		\[
		0<\omega\le \frac14,
		\qquad
		\delta_1,\delta_2\in(0,1).
		\]
		Then there exist constants $\theta\in(0,1)$ and $\nu\in(0,1)$,
		depending only on $d,\,\lambda,\,\Lambda,\,M,\,\omega,\,\delta_1,\,\delta_2$,
		such that the following holds.  Assume in addition that
		\[
		f\le1
		\qquad\text{a.e. in }Q_2^+(Z_0),
		\qquad
		|g|\le1
		\qquad\text{a.e. in }Q_2^+(Z_0).
		\]
		If
		\begin{align*}
			\bigl|\{f\ge 1-\theta\}\cap Q_\omega^+(Z_0)\bigr|
			&\ge
			\delta_1\,|Q_\omega(Z_0)|, \\
			\bigl|
			\bigl(\{f\le 0\}\cap Q_\omega^+(\widehat Z_0)\bigr)
			\cup
			Q_\omega^-(\widehat Z_0)
			\bigr|
			&\ge
			\delta_2\,|Q_\omega(\widehat Z_0)|,
		\end{align*}
		then
		\[
		\bigl|
		\{0<f<1-\theta\}
		\cap
		\bigl(Q_1^+(Z_0)\cup Q_1^+(\widehat Z_0)\bigr)
		\bigr|
		\ge
		\nu\,|Q_1(Z_0)\cup Q_1(\widehat Z_0)|.
		\]
	\end{prop}
	
	\begin{prop}[Unit-scale boundary oscillation decay]
		\label{cor:boundary-oscillation-decay}
		Let
		$Z_0\in \bR\times\overline{\bR^d_+}\times\bR^d$.
		Assume that
		$a^{ij},\,b^i\in L^\infty(Q_2^+(Z_0))$, $g\in L^\infty(Q_2^+(Z_0))$,
		and that, for some \(0<\lambda\le\Lambda\) and \(M\ge0\),
		\[
		\lambda |\xi|^2
		\le
		a^{ij}(Z)\xi_i\xi_j,
		\qquad
		|a^{ij}(Z)|\le \Lambda,
		\qquad
		|b(Z)|\le M
		\]
		for a.e. \(Z\in Q_2^+(Z_0)\) and every \(\xi\in\bR^d\).  Let
		\(f\in L^\infty(Q_2^+(Z_0))\) be a weak solution of
		\[
		(\partial_t+v\cdot\nabla_x)f
		=
		D_i(a^{ij}D_j f)+b^iD_i f+g
		\qquad\text{in }Q_2^+(Z_0),
		\]
		with $f=0$ on $(\bR\times\gamma^-)\cap Q_2(Z_0)$ in the sense of traces.
		
		Then there exist constants
		$\lambda_0\in(0,1)$, $\omega\in(0,1)$, $C>0$, 
		depending only on $d,\,\lambda,\,\Lambda,\,M$,
		such that
		\[
		\operatorname*{osc}_{Q_{\omega/2}^+(Z_0)} f
		\le
		\left(1-\frac{\lambda_0}{2}\right)
		\operatorname*{osc}_{Q_2^+(Z_0)} f
		+
		C\|g\|_{L^\infty(Q_2^+(Z_0))}.
		\]
	\end{prop}

	\subsection{Hypoellipticity for measure sources}
	
	This subsection collects the two ingredients used in Step~4 of
	Lemma~\ref{lem:positive-part-defect}: the hypoelliptic regularity estimate of
	Lemma~\ref{lem:bouchut-hypo} and the Bessel-potential reduction of measure
	sources in Proposition~\ref{prop:velocity-average-compactness}. Together with
	compact Sobolev embedding, they yield the \(L^1_{\mathrm{loc}}\)-compactness
	used in Proposition~\ref{prop:isop}.
	
	\begin{lemma}[Bouchut-type hypoelliptic regularity]
		\label{lem:bouchut-hypo}
		Let \(1<p<\infty\), \(0\le r\le 1\), and let
		\[
		f\in L^p(\bR_t\times\bR_x^d\times\bR_v^d),
		\qquad
		\nabla_v f\in L^p(\bR_t\times\bR_x^d\times\bR_v^d).
		\]
		Assume that there exist functions
		\(g^0,g^1,\ldots,g^d\) on
		\(\bR_t\times\bR_x^d\times\bR_v^d\) such that
		\[
		\|(1+|v|)g^0\|_{L^p(\bR^{1+2d})}
		+
		\sum_{i=1}^d
		\|(1+|v|)^2g^i\|_{L^p(\bR^{1+2d})}
		=:K_g<\infty,
		\]
		and
		\[
		(\partial_t+v\cdot\nabla_x)f
		=
		(1-\Delta_{(t,x)})^{r/2}
		\bigl(g^0+\partial_{v^i}g^i\bigr)
		\]
		in the sense of distributions. Then, with
		\[
		s:=\frac{1-r}{3}\ge 0,
		\]
		we have
		\[
		f\in W_{t,x}^{s,p}(\bR_t\times\bR_x^d\times\bR_v^d),
		\]
		and
		\begin{align*}
			\|D_x^s f\|_{L^p}
			+
			\|D_t^s f\|_{L^p}
			\le {}&
			C(d,p,r)
			\Bigl(
			\|f\|_{L^p}
			+
			\|D_v f\|_{L^p}^{1/2}K_g^{1/2}
			+
			\|D_v f\|_{L^p}^{2/3}K_g^{1/3}
			\Bigr).
		\end{align*}
	\end{lemma}
	
	\begin{proof}
		Apply \cite[Theorem~1.3]{bouchut2002hypoelliptic} with
		\[
		N=d,\qquad \beta=1,\qquad m=1,\qquad \kappa=\Omega=1,
		\]
		and with the multi-index family consisting of \(g^0\) and the \(d\)
		first-order terms \(g^i\).  In that theorem the source multiplier is
		\((1-\Delta_x-\partial_t^2)^{r/2}=(1-\Delta_{(t,x)})^{r/2}\), and its
		velocity weights
		\((1+|v|^2)^{(|\alpha|+1)/2}\) are equivalent to the weights used here.
		The resulting exponent is
		\[
		s=(1-r)\frac{\beta}{m+1+\beta}=\frac{1-r}{3}.
		\]
		After bounding each of the finitely many source norms by \(K_g\), the
		two families of terms in estimate~(1.13) of that theorem are bounded by
		\(\|D_vf\|_{L^p}^{1/2}K_g^{1/2}\) and
		\(\|D_vf\|_{L^p}^{2/3}K_g^{1/3}\), respectively.  This gives the stated estimate, including the inhomogeneous term \(\|f\|_{L^p}\).
	\end{proof}
	
	\begin{lemma}
		\label{lem:bessel-kernel-behavior}
		Let \(N\ge 1\) and \(r>0\). Let \(J_r^N\) be the Bessel kernel on \(\bR^N\),
		defined by
		\[
		\cF_N(J_r^N)(\xi)=(1+|\xi|^2)^{-r/2},
		\]
		where \(\cF_N\) denotes the Fourier transform on \(\bR^N\). Then
		\(J_r^N\) is smooth and strictly positive on \(\bR^N\setminus\{0\}\), and the
		following estimates hold.
		
		\begin{enumerate}
			\item (\emph{Near the origin.}) There exists \(C=C(N,r)>0\) such that, for
			\(0<|z|\le 1\),
			\[
			J_r^N(z)\le C
			\begin{cases}
				|z|^{r-N}, & 0<r<N,\\[1mm]
				\log \dfrac{e}{|z|}, & r=N,\\[2mm]
				1, & r>N.
			\end{cases}
			\]
			
			\item (\emph{At infinity.}) There exists \(C=C(N,r)>0\) such that, for
			\(|z|\ge 1\),
			\[
			J_r^N(z)\le C \ee^{-|z|/2}.
			\]
		\end{enumerate}
		
		In particular, if \(1\le p\le \infty\) satisfies
		\[
		\frac1p>1-\frac rN,
		\]
		with the convention \(1/\infty=0\), then for any $k\in\bN_0$,
		\begin{align}\label{2606301107}
			(1+|z|)^kJ_r^N(z)\in L^p(\bR^N).
		\end{align}
	\end{lemma}
	
	\begin{proof}
		See \cite[Section~4, formulas~(4.1)--(4.6), pp.~416--417]{AronszajnSmith1961}.
	\end{proof}
	
	\begin{prop}
		\label{prop:velocity-average-compactness}
		Let
		\[
		H^0,\ldots,H^d\in L^2(\bR^{1+2d}),
		\qquad
		\mu\in \cM(\bR^{1+2d})
		\]
		and assume that there exists \(R>0\) such that
		\begin{align}\label{2606301039}
			\operatorname{supp}(H^0,\ldots,H^d,\mu)
			\subset B_R^{1+2d}(0).
		\end{align}
		Let \(r\in(0,1)\) and let \(p\in(1,\infty)\) satisfy
		\begin{align}\label{2606301029}
			1<p<\frac{d+1}{d+1-r}.
		\end{align}
		Then there exist functions
		\(\widehat H^0,\ldots,\widehat H^d\) on \(\bR^{1+2d}\) such that
		\begin{align}\label{2606301040}
			&\|(1+|v|)\widehat H^0\|_{L^p(\bR^{1+2d})}
			+\sum_{i=1}^d
			\|(1+|v|)^2\widehat H^i\|_{L^p(\bR^{1+2d})} \notag\\
			&\qquad\le
			C(d,R,p,r)
			\left(
			\sum_{i=0}^d\|H^i\|_{L^2(\bR^{1+2d})}
			+\|\mu\|_{\cM(\bR^{1+2d})}
			\right),
		\end{align}
		and
		\begin{align}\label{2606301032}
			H^0+\partial_{v^i}H^i+\mu
			=
			(1-\Delta_{(t,x)})^{r/2}
			\left(
			\widehat H^0+\partial_{v^i}\widehat H^i
			\right)
			\quad
			\text{in }\cD'(\bR^{1+2d}).
		\end{align}
	\end{prop}
	\begin{proof}
		Write
		\[
		z:=(t,x)\in\bR^{1+d}.
		\]
		Let \(J_r^{d+1}\) be the Bessel kernel on \(\bR_z^{1+d}\), let \(J_s^d\) be
		the Bessel kernel on \(\bR_v^d\), and set
		\[
		K_{r,s}(z,v):=J_r^{d+1}(z)J_s^d(v).
		\]
		
		Define
		\begin{align*}
			\widetilde H^0
			&:=
			(1-\Delta_z)^{-r/2}H^0
			+
			(1-\Delta_z)^{-r/2}(1-\Delta_v)^{-1}\mu,\\
			\widetilde H^i
			&:=
			(1-\Delta_z)^{-r/2}H^i
			+
			(1-\Delta_z)^{-r/2}
			\bigl(-\partial_{v^i}\bigr)(1-\Delta_v)^{-1}\mu,
			\qquad i=1,\ldots,d.
		\end{align*}
		Equivalently,
		\[
		\widetilde H^0
		=
		J_r^{d+1}*_z H^0+K_{r,2}*\mu,
		\qquad
		\widetilde H^i
		=
		J_r^{d+1}*_z H^i-\partial_{v^i}K_{r,2}*\mu .
		\]
		Since the operators involved commute in \(\cD'(\bR^{1+2d})\), we have
		\begin{align}\label{2606301181}
			\widetilde H^0+\partial_{v^i}\widetilde H^i
			=
			(1-\Delta_z)^{-r/2}
			\left(
			H^0+\partial_{v^i}H^i+\mu
			\right).
		\end{align}
		
		We first prove the unweighted \(L^p\)-estimate for \(\widetilde H^i\). Since
		\(J_r^{d+1}\in L^1(\bR^{1+d})\), Young's inequality and \eqref{2606301039} give
		\[
		\|J_r^{d+1}*_zH^i\|_{L^2}
		\le C(d,r)\|H^i\|_{L^2},\qquad
		\|J_r^{d+1}*_zH^i\|_{L^1}
		\le C(d,r)\|H^i\|_{L^1}
		\le C(d,R,r)\|H^i\|_{L^2}.
		\]
		As \eqref{2606301029} implies \(1<p<2\), interpolation yields
		\begin{align}\label{2606301182}
			\|J_r^{d+1}*_zH^i\|_{L^p(\bR^{1+2d})}
			\le
			C(d,R,p,r)\|H^i\|_{L^2(\bR^{1+2d})},
			\qquad i=0,\ldots,d.
		\end{align}
		
		We next estimate the terms involving \(\mu\).  Since
		\(K_{r,m}(z,v)=J^{d+1}_r(z)J^d_m(v)\), we have
		\(\|K_{r,m}\|_{L^p(\bR^{1+2d})}
		=\|J^{d+1}_r\|_{L^p(\bR^{1+d})}\|J^d_m\|_{L^p(\bR^d)}\).  By
		Lemma~\ref{lem:bessel-kernel-behavior} the first factor is finite exactly
		under \eqref{2606301029}.
		If \(d\ge2\), the second factor is finite for
		\(m\in\{1,2\}\) as soon as \(p<\frac d{d-1}\), and this follows from
		\eqref{2606301029}, because
		\((d+1)(d-1)=d^2-1<d(d+1-r)\) for \(r<1\). If \(d=1\),
		Lemma~\ref{lem:bessel-kernel-behavior} instead gives the same conclusion
		for every finite \(p\), which is all that is needed here.
		Hence
		\[
		K_{r,1},\,K_{r,2}\in L^p(\bR^{1+2d}).
		\]
		Hence Young's inequality for convolution with a finite measure gives
		\begin{align}\label{2606301183}
			\|K_{r,1}*\mu\|_{L^p}
			+
			\|K_{r,2}*\mu\|_{L^p}
			\le
			C(d,p,r)\|\mu\|_{\cM(\bR^{1+2d})}.
		\end{align}
		For the derivative term, set
		\[
		T_j:=\partial_{v^j}(1-\Delta_v)^{-1/2}.
		\]
		The operator \(T_j\) is a Mihlin multiplier in the \(v\)-variable, and
		\[
		\partial_{v^j}K_{r,2}
		=
		T_jK_{r,1}.
		\]
		Therefore, by the \(L^p\)-boundedness of \(T_j\) and \eqref{2606301183},
		\begin{align}\label{2606301184}
			\|\partial_{v^j}K_{r,2}*\mu\|_{L^p}
			=
			\|T_j(K_{r,1}*\mu)\|_{L^p}\le C(d,p)\|K_{r,1}*\mu\|_{L^p}
			\le
			C(d,p,r)\|\mu\|_{\cM(\bR^{1+2d})}.
		\end{align}
		Combining \eqref{2606301182}, \eqref{2606301183}, and \eqref{2606301184}, we get
		\begin{align}\label{2606301186}
			\|\widetilde H^0\|_{L^p}
			+
			\sum_{i=1}^d
			\|\widetilde H^i\|_{L^p}
			\le
			C(d,R,p,r)
			\left(
			\sum_{i=0}^d\|H^i\|_{L^2}
			+\|\mu\|_{\cM(\bR^{1+2d})}
			\right).
		\end{align}
		
		It remains to localize in the \(v\)-variable. Choose
		\[
		\chi\in C_c^\infty(\bR_v^d),
		\qquad
		0\le \chi\le 1,
		\qquad
		\chi\equiv 1 \text{ on } B_{3R/2}^d(0),
		\qquad
		\supp\chi\subset B_{2R}^d(0).
		\]
		Define
		\[
		\widehat H^0
		:=
		\chi\widetilde H^0-(\partial_{v^i}\chi)\widetilde H^i,
		\qquad
		\widehat H^i
		:=
		\chi\widetilde H^i,
		\quad i=1,\ldots,d.
		\]
		Then
		\[
		\widehat H^0+\partial_{v^i}\widehat H^i
		=
		\chi\left(\widetilde H^0+\partial_{v^i}\widetilde H^i\right).
		\]
		By \eqref{2606301181}, the factor in parentheses is supported in
		\(\bR_z^{1+d}\times B_R^d(0)\), because \((1-\Delta_z)^{-r/2}\) acts only in the
		\(z\)-variables. Since \(\chi\equiv1\) on \(B_R^d(0)\), we obtain
		\[
		\widehat H^0+\partial_{v^i}\widehat H^i
		=
		(1-\Delta_z)^{-r/2}
		\left(
		H^0+\partial_{v^i}H^i+\mu
		\right).
		\]
		Applying \((1-\Delta_z)^{r/2}\) gives \eqref{2606301032}.
		
		Finally,
		\[
		\supp_v\widehat H^0,\ \supp_v\widehat H^i\subset B_{2R}^d(0),
		\]
		so the weights \((1+|v|)\) and \((1+|v|)^2\) are bounded by \(C_R\) on their
		supports. Thus, by the definitions of \(\widehat H^i\), the boundedness of
		\(\chi\) and \(\nabla_v\chi\), and \eqref{2606301186},
		\[
		\|(1+|v|)\widehat H^0\|_{L^p}
		+
		\sum_{i=1}^d
		\|(1+|v|)^2\widehat H^i\|_{L^p}
		\le
		C(d,R,p,r)
		\left(
		\sum_{i=0}^d\|H^i\|_{L^2}
		+\|\mu\|_{\cM(\bR^{1+2d})}
		\right).
		\]
		This proves \eqref{2606301040} and completes the proof.
	\end{proof}

	\subsection{Auxiliary tools for Proposition~\ref{prop:isop}}
	
	\begin{lemma}
		\label{lem:H1-binary-constant}
		Let \(\Omega\subset \bR^n\) be a connected open set.
		If \(u\in H^1(\Omega)\) and
		\[
		u(x)\in\{0,1\}
		\qquad\text{for a.e. }x\in \Omega,
		\]
		then \(u\) is constant a.e. on \(\Omega\).
	\end{lemma}
	\begin{proof}
		The proof is elementary; since $u^2=u$ a.e., $\nabla u = \nabla(u^2)=2u\nabla u$ in the sense of distributions. Thus $(1-2u)\nabla u=0$ a.e.  But $1-2u\in\{\pm1\}$ almost everywhere, hence $\nabla u=0$ almost everywhere.
		Therefore $u$ is constant almost everywhere on each connected component of $\Omega$;
		since $\Omega$ is connected, $u$ is constant almost everywhere on $\Omega$. 
	\end{proof}

	\begin{lemma}
		\label{lem:no-upward-jumps}
		Let \(v_0\in\bR^d\) be fixed, and let
		\(\Omega\subset \bR_t\times\bR^d_x\) be open. Let
		\(U\in L^\infty(\Omega)\) satisfy
		\[
		U(t,x)\in\{0,1\}
		\qquad\text{for a.e. }(t,x)\in\Omega.
		\]
		Assume that there exists \(g\in L^1_{\mathrm{loc}}(\Omega)\) such that
		\[
		(\partial_t+v_0\cdot\nabla_x)U\le g
		\qquad\text{in }\cD'(\Omega),
		\]
		that is, for every nonnegative \(\varphi\in C_c^\infty(\Omega)\),
		\[
		\left\langle
		(\partial_t+v_0\cdot\nabla_x)U,\varphi
		\right\rangle
		\le
		\int_\Omega g\varphi\dd x\dd t.
		\]
		Then
		\[
		(\partial_t+v_0\cdot\nabla_x)U\le 0
		\qquad\text{in }\cD'(\Omega).
		\]
		Equivalently, in the straightened variables, for almost every characteristic
		line the one-dimensional restriction of \(U\) admits a nonincreasing
		representative on each connected component of its intersection with
		\(\Omega\). In particular, it may jump from \(1\) to \(0\), but not from
		\(0\) to \(1\).
	\end{lemma}

	\begin{proof}
		\noindent\emph{Step 1. A one-dimensional fact.}
		We first prove the following one-dimensional fact. Let \(I\subset\bR\)
		be an open interval, let \(u\in L^\infty(I)\) satisfy
		\[
		u(s)\in\{0,1\}
		\qquad\text{for a.e. }s\in I,
		\]
		and let \(h\in L^1_{\mathrm{loc}}(I)\). Suppose that
		\[
		\partial_s u\le h
		\qquad\text{in }\cD'(I).
		\]
		We claim that
		\[
		\partial_s u\le0
		\qquad\text{in }\cD'(I).
		\]
		
		Throughout this one-dimensional argument, \(D\) denotes the distributional
		derivative in the \(s\)-variable. Thus, if \(w\in BV_{\mathrm{loc}}(I)\), then
		\(Dw\) is the associated locally finite signed Radon measure, characterized by
		\[
		\langle Dw,\varphi\rangle
		=
		-\int_I w(s)\varphi'(s)\dd s,
		\qquad
		\varphi\in C_c^\infty(I).
		\]
		We also identify \(h\in L^1_{\mathrm{loc}}(I)\) with the absolutely continuous
		Radon measure \(h(s)\dd s\).
		
		Fix \(s_*\in I\) and define
		\[
		H(s):=\int_{s_*}^s h(\tau)\dd \tau,
		\qquad s\in I.
		\]
		Then \(H\in W^{1,1}_{\mathrm{loc}}(I)\cap C(I)\) and \(H'=h\) a.e. in \(I\).
		Hence
		\[
		\partial_s(u-H)\le0
		\qquad\text{in }\cD'(I).
		\]
		Equivalently, \(-\partial_s(u-H)\) is a positive distribution. By the Riesz
		representation theorem for positive distributions, there exists a
		nonnegative Radon measure \(\mu\) on \(I\) such that
		\[
		\partial_s(u-H)=-\mu
		\qquad\text{in }\cD'(I).
		\]
		Consequently,
		\[
		Du=h(s)\dd s-\mu
		\qquad\text{as Radon measures on }I.
		\]
		In particular, \(u\in BV_{\mathrm{loc}}(I)\).
		
		We next choose a right-continuous representative of \(u\). Define
		\[
		F(s):=
		\begin{cases}
			-\mu((s_*,s]),
			& s>s_*,\\[1mm]
			0,
			& s=s_*,\\[1mm]
			\mu((s,s_*]),
			& s<s_*.
		\end{cases}
		\]
		Then \(F\) is right-continuous, belongs to \(BV_{\mathrm{loc}}(I)\), and
		satisfies
		\[
		DF=-\mu
		\qquad\text{in }\cD'(I).
		\]
		Since
		\[
		D(u-H-F)=0
		\qquad\text{in }\cD'(I),
		\]
		there exists a constant \(c\in\bR\) such that, after changing \(u\) on a
		null set,
		\[
		u(s)=c+H(s)+F(s)
		\qquad\text{for every }s\in I.
		\]
		We henceforth use this representative. It is right-continuous and belongs to
		\(BV_{\mathrm{loc}}(I)\).
		Since this $u$ is valued in $\{0,1\}$ a.e. in $I$, the right-continuity of $u$ yields that $u(s)\in \{0,1\}$ for all $s\in I$.
		
		Let \(J=(a,b)\Subset I\). Since \(u\in BV(J)\), the pointwise total variation
		of the right-continuous representative on \(J\) is finite and equals
		\(|Du|(J)\). Since \(u\) takes only the values \(0\) and \(1\), every genuine
		change of value contributes exactly \(1\) to the total variation. Hence \(u\)
		can change its value only finitely many times on \(J\). Therefore there exist
		finitely many points
		\[
		a<s_1<\cdots<s_N<b
		\]
		such that \(u\) is constant on each connected component of
		\[
		J\setminus\{s_1,\dots,s_N\}.
		\]
		These points are precisely the jump points of \(u\) in \(J\).
		
		For each \(k\), set
		\[
		\alpha_k:=u(s_k+)-u(s_k-),
		\]
		where the one-sided limits exist because \(u\in BV(J)\). Since \(u\) is
		\(\{0,1\}\)-valued and \(s_k\) is a genuine jump point, we have
		\[
		\alpha_k\in\{-1,1\}.
		\]
		For every \(\varphi\in C_c^\infty(J)\), integration by parts on the connected
		components of \(J\setminus\{s_1,\dots,s_N\}\) gives
		\[
		\begin{aligned}
			\langle Du,\varphi\rangle
			&=
			-\int_J u(s)\varphi'(s)\dd s =
			\sum_{k=1}^N \alpha_k\varphi(s_k).
		\end{aligned}
		\]
		Thus
		\[
		Du=\sum_{k=1}^N \alpha_k\delta_{s_k}
		\qquad\text{as a measure on }J.
		\]
		
		On the other hand, since $Du=h(s)\dd s-\mu$
		and \(h(s)\dd s\) has no atoms, we have, for each \(k\),
		\[
		\alpha_k
		=
		Du(\{s_k\})
		=
		-\mu(\{s_k\})
		\le0.
		\]
		Therefore \(\alpha_k=-1\) for every \(k\). Hence
		\[
		Du\le0
		\qquad\text{as a measure on }J.
		\]
		Since \(J\Subset I\) was arbitrary, we conclude that
		\[
		\partial_s u\le0
		\qquad\text{in }\cD'(I).
		\]
		This proves the one-dimensional fact.
		
		\noindent\emph{Step 2. Straightening the characteristics.}
		We now return to the transport inequality. Set
		\[
		s:=t,\qquad y:=x-tv_0,
		\qquad x=y+sv_0.
		\]
		Let
		\[
		\widetilde\Omega
		:=
		\{(s,y)\in\bR\times\bR^d:(s,y+sv_0)\in\Omega\},
		\]
		and define
		\[
		\widetilde U(s,y):=U(s,y+sv_0),
		\qquad
		\widetilde g(s,y):=g(s,y+sv_0).
		\]
		Then \(\widetilde U\in L^\infty(\widetilde\Omega)\),
		\[
		\widetilde U(s,y)\in\{0,1\}
		\qquad\text{for a.e. }(s,y)\in\widetilde\Omega,
		\]
		and \(\widetilde g\in L^1_{\mathrm{loc}}(\widetilde\Omega)\). Moreover, by
		the above change of variables, the assumption becomes
		\[
		\partial_s\widetilde U\le \widetilde g
		\qquad\text{in }\cD'(\widetilde\Omega).
		\]
		
		For \(y\in\bR^d\), set
		\[
		\widetilde\Omega_y
		:=
		\{s\in\bR:(s,y)\in\widetilde\Omega\}.
		\]
		We claim that for a.e. \(y\in\bR^d\),
		\[
		\partial_s\widetilde U(\cdot,y)
		\le
		\widetilde g(\cdot,y)
		\qquad\text{in }\cD'(\widetilde\Omega_y).
		\]
		
		To prove this, let \(I\times B\Subset\widetilde\Omega\), where
		\(I\subset\bR\) is an open interval and \(B\subset\bR^d\) is an
		open ball. For \(\rho\in C_c^\infty(I)\), \(\rho\ge0\), define
		\[
		A_\rho(y)
		:=
		-\int_I \widetilde U(s,y)\rho'(s)\dd s
		-
		\int_I \widetilde g(s,y)\rho(s)\dd s,
		\qquad y\in B.
		\]
		By Fubini's theorem, \(A_\rho\in L^1(B)\). Testing
		\(\partial_s\widetilde U\le\widetilde g\) with
		\[
		\psi(s,y)=\rho(s)\eta(y),
		\qquad
		\eta\in C_c^\infty(B),\quad \eta\ge0,
		\]
		gives
		\[
		\int_B A_\rho(y)\eta(y)\dd y\le0.
		\]
		Therefore \(A_\rho(y)\le0\) for a.e. \(y\in B\).
		
		Choose a countable dense family of nonnegative functions in
		\(C_c^\infty(I)\), dense with respect to the \(C^1\)-topology on \(I\).
		Intersecting the corresponding full-measure sets and using the continuity of
		the map
		\[
		\rho\mapsto
		-\int_I \widetilde U(s,y)\rho'(s)\dd s
		-
		\int_I \widetilde g(s,y)\rho(s)\dd s,
		\]
		we obtain, for a.e. \(y\in B\),
		\[
		-\int_I \widetilde U(s,y)\rho'(s)\dd s
		\le
		\int_I \widetilde g(s,y)\rho(s)\dd s
		\]
		for every nonnegative \(\rho\in C_c^\infty(I)\). Hence
		\[
		\partial_s\widetilde U(\cdot,y)
		\le
		\widetilde g(\cdot,y)
		\qquad\text{in }\cD'(I)
		\]
		for a.e. \(y\in B\).
		
		Now take a countable family of cylinders \(I_m\times B_m\Subset
		\widetilde\Omega\) with rational centers and radii such that this family forms
		a basis of \(\widetilde\Omega\). Intersecting the corresponding full-measure
		sets in the \(y\)-variable, we obtain a full-measure set \(Y\subset\bR^d\)
		such that the preceding slice inequality holds on every such cylinder whose
		spatial ball contains \(y\in Y\). It follows, by a partition of unity in the
		\(s\)-variable, that for every \(y\in Y\),
		\[
		\partial_s\widetilde U(\cdot,y)
		\le
		\widetilde g(\cdot,y)
		\qquad\text{in }\cD'(\widetilde\Omega_y).
		\]
		
		For every \(y\in Y\), the one-dimensional fact applies on each connected
		component of \(\widetilde\Omega_y\). Therefore
		\[
		\partial_s\widetilde U(\cdot,y)\le0
		\qquad\text{in }\cD'(\widetilde\Omega_y)
		\]
		for a.e. \(y\in\bR^d\).
		
		Finally, let \(\psi\in C_c^\infty(\widetilde\Omega)\) be nonnegative. By
		Fubini's theorem and the slice inequality just obtained,
		\[
		\begin{aligned}
			\left\langle \partial_s\widetilde U,\psi\right\rangle
			&=
			-\int_{\bR^d}
			\int_{\widetilde\Omega_y}
			\widetilde U(s,y)\partial_s\psi(s,y)\dd s\dd y \le 0.
		\end{aligned}
		\]
		Thus
		\[
		\partial_s\widetilde U\le0
		\qquad\text{in }\cD'(\widetilde\Omega).
		\]
		Undoing the change of variables gives
		\[
		(\partial_t+v_0\cdot\nabla_x)U\le0
		\qquad\text{in }\cD'(\Omega).
		\]
		The monotonicity statement along almost every characteristic follows from the
		one-dimensional characterization of functions whose distributional derivative
		is a nonpositive measure.
	\end{proof}
	
	\begin{lemma}[Propagation of zero-density points along characteristics]
		\label{lem:propagate-zero-density}
		Let \(v_0\in \bR^d\), let		$\Omega\subset \bR_t\times \bR_x^d$
		be open, and let \(U\in L^\infty(\Omega)\) satisfy
		\begin{align}\label{260629710}
			0\le U\le 1
			\qquad\text{a.e. in }\Omega
		\end{align}
		and
		\begin{align}\label{260629709}
			(\partial_t+v_0\cdot\nabla_x)U\le 0
			\qquad\text{in }\cD'(\Omega).
		\end{align}
		Let \(Y_0=(t_0,x_0)\in\Omega\) be a zero-density point of \(U\), namely
		\begin{align}\label{260629717}
			\lim_{r\downarrow0}
			\fint_{B_r^{d+1}(Y_0)}U(t,x)\dd x\dd t=0.
		\end{align}
		Assume that \(s>0\) is such that the characteristic segment
		\[
		\Gamma([0,s])
		:=
		\{(t_0+\tau,x_0+\tau v_0):0\le \tau\le s\}
		\]
		is compactly contained in \(\Omega\). Then
		\[
		Y_s:=(t_0+s,x_0+s v_0)
		\]
		is also a zero-density point of \(U\), i.e.
		\begin{align}\label{260629718}
			\lim_{r\downarrow0}
			\fint_{B_r^{d+1}(Y_s)}U(t,x)\dd x\dd t=0.
		\end{align}
	\end{lemma}
	
	\begin{proof}
		Set \(N:=d+1\), write \(Y=(t,x)\), and put
		\[
		\cV:=(1,v_0).
		\]
		Thus $\Gamma(\tau)=Y_0+\tau\cV$.
		Since \(\Gamma([0,s])\Subset\Omega\), there exists \(\rho_*>0\) such that
		\[
		\cN_{2\rho_*}
		:=
		\{Y\in\bR^{1+d}:\operatorname{dist}(Y,\Gamma([0,s]))<2\rho_*\}
		\Subset\Omega.
		\]
		
		Choose a nonnegative \(\eta\in C_c^\infty(B_1^N(0))\) with
		$\int_{\bR^{1+d}}\eta(Y)\dd Y=1$,
		and set
		\[
		\eta_\varepsilon(Y):=\varepsilon^{-N}\eta(\varepsilon^{-1}Y).
		\]
		For \(0<\varepsilon<\rho_*\), define
		\[
		U^\varepsilon(Y)
		:=
		\int_\Omega U(\bar Y)\eta_\varepsilon(Y-\bar Y)\dd \bar Y,
		\qquad
		Y\in\cN_{\rho_*}.
		\]
		This is well-defined and smooth in \(\cN_{\rho_*}\), because
		\(\cN_{2\rho_*}\Subset\Omega\). Moreover, by \eqref{260629709},
		\begin{align}\label{260629719}
			(\partial_t+v_0\cdot\nabla_x)U^\varepsilon\le0
			\qquad\text{in }\cN_{\rho_*}.
		\end{align}
		Indeed, for each \(Y\in\cN_{\rho_*}\), the function
		\(\bar Y\mapsto\eta_\varepsilon(Y-\bar Y)\) is a nonnegative test function
		compactly supported in \(\Omega\), and
		\[
		(\partial_t+v_0\cdot\nabla_x)U^\varepsilon(Y)
		=
		\left\langle
		(\partial_t+v_0\cdot\nabla_x)U,
		\eta_\varepsilon(Y-\cdot)
		\right\rangle
		\le0.
		\]
		
		Let \(0<r<\rho_*\). If \(Y\in B_r^N(Y_0)\), then
		$Y+\tau\cV\in\cN_{\rho_*}$
		for all $0\le \tau\le s$.
		Hence, by \eqref{260629719},
		\[
		\frac{\dd}{\dd\tau}U^\varepsilon(Y+\tau\cV)
		=
		(\partial_t+v_0\cdot\nabla_x)U^\varepsilon(Y+\tau\cV)
		\le0.
		\]
		Therefore
		\[
		U^\varepsilon(Y+s\cV)\le U^\varepsilon(Y)
		\qquad\text{for every }Y\in B_r^N(Y_0).
		\]
		Integrating over \(B_r^N(Y_0)\) and using the translation
		\(Y\mapsto Y+s\cV\), we obtain
		\begin{align}\label{260629720}
			\int_{B_r^N(Y_s)}U^\varepsilon(Y)\dd Y
			=\int_{B_r^N(Y_0)}U^\varepsilon(Y+s\cV)\dd Y\le
			\int_{B_r^N(Y_0)}U^\varepsilon(Y)\dd Y.
		\end{align}
		
		Since \(B_r^N(Y_0)\cup B_r^N(Y_s)\Subset\Omega\), we have
		\[
		U^\varepsilon\to U
		\qquad\text{in }L^1\bigl(B_r^N(Y_0)\cup B_r^N(Y_s)\bigr).
		\]
		Letting \(\varepsilon\downarrow0\) in \eqref{260629720}, we get
		\[
		\int_{B_r^N(Y_s)}U(Y)\dd Y
		\le
		\int_{B_r^N(Y_0)}U(Y)\dd Y.
		\]
		Equivalently,
		\[
		\fint_{B_r^N(Y_s)}U(Y)\dd Y
		\le
		\fint_{B_r^N(Y_0)}U(Y)\dd Y
		\qquad
		\text{for every }0<r<\rho_*.
		\]
		Letting \(r\downarrow0\) and using \eqref{260629717}, we obtain
		\[
		\lim_{r\downarrow0}
		\fint_{B_r^N(Y_s)}U(Y)\dd Y=0.
		\]
		Thus \(Y_s\) is a zero-density point of \(U\).
	\end{proof}

	\subsection{Proof of Proposition~\ref{prop:isop}}\label{260701421}
	
	We argue by contradiction; the argument is modeled on the proof of the De Giorgi
	intermediate-value lemma in \cite[Lemma~4.1]{2019GIMV}, but we recast it in a
	boundary-adapted form for the compactness limit obtained in
	Lemma~\ref{lem:positive-part-defect}.
	
	Suppose that the conclusion fails. Then there exist a sequence
	\[
	\theta_k\downarrow 0,
	\]
	a sequence of points
	\[
	Z_k=(t_k,x_k,v_k)\in \bR\times \overline{\bR^d_+}\times \bR^d,
	\qquad
	\widehat Z_k:=(t_k-1,x_k-v_k,v_k),
	\]
	and a sequence of weak solutions \(f_k\), in the sense of
	Definition~\ref{def:L2sol}, of
	\[
	\partial_t f_k+v\cdot\nabla_x f_k
	=
	D_i\!\bigl(a_k^{ij}D_j f_k\bigr)+b_k^iD_i f_k+g_k
	\qquad \text{in } Q_2^+(Z_k),
	\]
	where \(a_k^{ij},b_k^i,g_k\in L^\infty(Q_2^+(Z_k))\) and \(a_k,b_k\)
	satisfy \eqref{con_ab} on \(Q_2^+(Z_k)\) with the same constants
	\(\lambda,\Lambda,M\) for every \(k\), such that
	\[
	f_k\le 0
	\qquad\text{on }(\bR\times\gamma^-)\cap Q_2(Z_k)
	\]
	in the sense of traces,
	\[
	f_k\le 1,
	\qquad
	|g_k|\le 1
	\qquad\text{a.e. in }Q_2^+(Z_k),
	\]
	and
	\begin{align}
		\bigl|\{f_k\ge 1-\theta_k\}\cap Q_\omega^+(Z_k)\bigr|
		&\ge \delta_1\,|Q_\omega(Z_k)|, \label{eq:isop-contr-1}\\
		\bigl|\bigl(\{f_k\le 0\}\cap Q_\omega^+(\widehat Z_k)\bigr)\cup Q_\omega^-(\widehat Z_k)\bigr|
		&\ge \delta_2\,|Q_\omega(\widehat Z_k)|, \notag\\
		\bigl|\{0<f_k<1-\theta_k\}\cap\bigl(Q_1^+(Z_k)\cup Q_1^+(\widehat Z_k)\bigr)\bigr|
		&\longrightarrow 0. \label{eq:isop-contr-3}
	\end{align}

	The following key lemma will be used below. Its proof is postponed to
	the end of the proof of the proposition.
	
	\begin{lemma}
		\label{lem:positive-part-defect}
		Let \(Z_k=(t_k,x_k,v_k)\in \bR\times \overline{\bR^d_+}\times \bR^d\), and let
		\(f_k\) be a weak solution of
		\[
		\partial_t f_k+v\cdot\nabla_x f_k
		=
		D_i\!\bigl(a_k^{ij}D_j f_k\bigr)+b_k^iD_i f_k+g_k
		\qquad\text{in }Q_2^+(Z_k),
		\]
		with
		\[
		f_k\le 0
		\qquad\text{on }(\bR\times \gamma^-)\cap Q_2(Z_k)
		\]
		in the sense of traces. Assume moreover that
		\[
		f_k\le 1,
		\qquad
		|g_k|\le 1
		\qquad\text{a.e. in }Q_2^+(Z_k),
		\]
		and that the coefficients satisfy \eqref{con_ab} uniformly in \(k\).
		
		Set
		\[
		U_R:=\{(t,x,v): |t|<R^2,\ |x|<R^3,\ |v|<R\}.
		\]
		Extending \(f_k^+1_{Q_2^+(Z_k)}\) by zero outside \(Q_2(Z_k)\), define
		\[
		\widehat f_k(Z):=\left(f_k^+\,1_{Q_2^+(Z_k)}\right)(Z_k\circ Z),
		\qquad Z\in U_2.
		\]
		Then the following hold.
		
		\begin{enumerate}
			\item There exist
			\[
			H_{0,k}\in L^2_{\mathrm{loc}}(U_2),
			\qquad
			H_{1,k}\in L^2_{\mathrm{loc}}(U_2;\bR^d),
			\]
			and a nonnegative Radon measure
			\[
			\mu_k\in \cM_{\mathrm{loc}}(U_2)
			\]
			such that
			\begin{equation}\label{eq:positive-part-defect-identity}
				(\partial_t+v\cdot\nabla_x)\widehat f_k
				=
				H_{0,k}+\nabla_v\cdot H_{1,k}-\mu_k
				\qquad\text{in }\cD'(U_2).
			\end{equation}
			
			\item For every \(0<R<2\),
			\begin{equation}\label{eq:positive-part-local-bounds}
				\sup_k\Big(
				\|\widehat f_k\|_{L^\infty(U_R)}
				+\|\nabla_v \widehat f_k\|_{L^2(U_R)}
				+\|H_{0,k}\|_{L^2(U_R)}
				+\|H_{1,k}\|_{L^2(U_R)}
				+\mu_k(U_R)
				\Big)<\infty .
			\end{equation}
			
			\item The sequence \((\widehat f_k)_{k\in\bN}\) is precompact in
			\(L^1_{\mathrm{loc}}(U_2)\).
		\end{enumerate}
	\end{lemma}
	
	\begin{proof}[\textbf{Proof of Proposition~\ref{prop:isop}}]
		
		\smallskip
		\noindent
		\emph{Step 1. Normalization and positive-part extension.}
		Set
		\[
		\widehat Z_k:=Z_k\circ(-1,0,0),
		\]
		and
		\[
		Q_r^\circ:=(-r^2,0)\times B_{r^3}^d(0)\times B_r^d(0).
		\]
		Define
		\[
		\widehat f_k(Z)
		:=
		\left(f_k^+\,1_{Q_2^+(Z_k)}\right)(Z_k\circ Z),
		\qquad Z\in U_2.
		\]
		Since \(f_k\le1\), we have
		\[
		0\le \widehat f_k\le1
		\qquad\text{a.e. in }U_2.
		\]
		
		\smallskip
		\noindent
		\emph{Step 2. Compactness and the limit equation.}
		By Lemma~\ref{lem:positive-part-defect}, there exist
		\[
		H_{0,k}\in L^2_{\mathrm{loc}}(U_2),
		\qquad
		H_{1,k}\in L^2_{\mathrm{loc}}(U_2;\bR^d),
		\qquad
		\mu_k\in\cM_{\mathrm{loc}}(U_2),
		\]
		with \(\mu_k\ge0\), such that
		\begin{equation}\label{eq:defect-id-translated}
			(\partial_t+v\cdot\nabla_x)\widehat f_k
			=
			H_{0,k}+\nabla_v\cdot H_{1,k}-\mu_k
			\qquad\text{in }\cD'(U_2).
		\end{equation}
		Moreover, for every \(0<R<2\),
		\begin{equation}\label{eq:defect-bounds-translated}
			\sup_k\Big(
			\|\widehat f_k\|_{L^\infty(U_R)}
			+\|\nabla_v\widehat f_k\|_{L^2(U_R)}
			+\|H_{0,k}\|_{L^2(U_R)}
			+\|H_{1,k}\|_{L^2(U_R)}
			+\mu_k(U_R)
			\Big)<\infty,
		\end{equation}
		and \((\widehat f_k)_{k\in\bN}\) is precompact in
		\(L^1_{\mathrm{loc}}(U_2)\).
		
		Since $\overline{Q_{3/2}^\circ}\Subset U_2$, 
		after passing to a subsequence, there exists
		\[
		F\in L^1(Q_{3/2}^\circ)
		\]
		such that
		\begin{equation}\label{eq:L1-strong-positive-part}
			\widehat f_k\to F
			\qquad\text{strongly in }L^1(Q_{3/2}^\circ).
		\end{equation}
		After extracting a further subsequence, we may also assume that
		\begin{equation*}
			\widehat f_k(Z)\to F(Z)
			\qquad\text{for a.e. }Z\in Q_{3/2}^\circ.
		\end{equation*}
		Since \(0\le \widehat f_k\le1\), it follows that
		\[
		0\le F\le1
		\qquad\text{a.e. in }Q_{3/2}^\circ.
		\]
		Moreover, \eqref{eq:L1-strong-positive-part} and the uniform bound
		\(0\le \widehat f_k,F\le1\) imply
		\begin{equation*}
			\widehat f_k\to F
			\qquad\text{strongly in }L^2(Q_{3/2}^\circ).
		\end{equation*}
		Indeed,
		\[
		\|\widehat f_k-F\|_{L^2(Q_{3/2}^\circ)}^2
		\le
		\|\widehat f_k-F\|_{L^1(Q_{3/2}^\circ)}.
		\]
		
		Choose \(R_*\in(3/2,2)\). Since $
		Q_{3/2}^\circ\Subset U_{R_*}$, 
		\eqref{eq:defect-bounds-translated} implies that
		\[
		(\nabla_v\widehat f_k)_k
		\quad\text{is bounded in }L^2(Q_{3/2}^\circ;\bR^d).
		\]
		Hence, after passing to a further subsequence,
		\[
		\nabla_v\widehat f_k\rightharpoonup W
		\qquad\text{weakly in }L^2(Q_{3/2}^\circ;\bR^d)
		\]
		for some \(W\in L^2(Q_{3/2}^\circ;\bR^d)\). Since
		\(\widehat f_k\to F\) strongly in \(L^2(Q_{3/2}^\circ)\), we have
		\[
		W=\nabla_vF
		\]
		in the sense of distributions. Therefore,
		\begin{equation}\label{eq:F-H1v}
			F\in
			L^2\bigl(
			(-(3/2)^2,0)\times B_{(3/2)^3}(0);
			H^1(B_{3/2}(0))
			\bigr).
		\end{equation}
		
		Likewise, by \eqref{eq:defect-bounds-translated}, after passing to
		subsequences if necessary, we may assume that
		\begin{align*}
			H_{0,k}\,&\rightharpoonup H_0
			\qquad\text{weakly in }L^2(Q_{3/2}^\circ),\\
			H_{1,k}\,&\rightharpoonup H_1
			\qquad\text{weakly in }L^2(Q_{3/2}^\circ;\bR^d),
		\end{align*}
		and
		\[
		\mu_k\overset{\ast}{\rightharpoonup}\mu
		\qquad\text{weakly-* in }\cM(Q_{3/2}^\circ),
		\]
		where \(\mu\) is a nonnegative Radon measure. Passing to the limit in
		\eqref{eq:defect-id-translated}, we obtain
		\begin{equation}\label{eq:F-limit-eq}
			(\partial_t+v\cdot\nabla_x)F
			=
			H_0+\nabla_v\cdot H_1-\mu
			\qquad\text{in }\cD'(Q_{3/2}^\circ).
		\end{equation}
		
		\smallskip
		\noindent
		\emph{Step 3. Measure bounds for the limit.}
		Set
		\[
		Q_\omega:=Q_\omega((0,0,0)),
		\qquad
		\widehat Q_\omega:=Q_\omega((-1,0,0)),
		\qquad
		Q_1:=Q_1((0,0,0)),
		\qquad
		\widehat Q_1:=Q_1((-1,0,0)).
		\]
		All cylinders used in this step are contained in \(Q_{3/2}^\circ\), up to time
		slices of measure zero. Hence \eqref{eq:L1-strong-positive-part} implies
		convergence in measure on each of them.
		
		By construction of \(\widehat f_k\), the assumptions
		\eqref{eq:isop-contr-1}--\eqref{eq:isop-contr-3} become
		\begin{align}
			\bigl|\{\widehat f_k\ge 1-\theta_k\}\cap Q_\omega\bigr|
			&\ge \delta_1 |Q_\omega|,
			\notag\\
			\bigl|\{\widehat f_k=0\}\cap \widehat Q_\omega\bigr|
			&\ge \delta_2 |\widehat Q_\omega|,
			\label{eq:limit-measure-2}
			\\
			\bigl|\{0<\widehat f_k<1-\theta_k\}\cap(Q_1\cup\widehat Q_1)\bigr|
			&\longrightarrow0 .
			\label{eq:limit-measure-3}
		\end{align}
		Here, in \eqref{eq:limit-measure-2} and \eqref{eq:limit-measure-3}, we used
		that the positive-part extension vanishes both where \(f_k\le0\) in the upper
		half-space and on the lower half-space.
		
		We now pass to the limit in these three measure relations. First, fix
		\(0<\varepsilon<1/2\). Since \(\theta_k\to0\), for all sufficiently large \(k\),
		\[
		\{\widehat f_k\ge1-\theta_k\}\subset\{\widehat f_k>1-\varepsilon\}.
		\]
		Using convergence in measure on \(Q_\omega\), we get
		\[
		\delta_1|Q_\omega|
		\le
		\liminf_{k\to\infty}
		|\{\widehat f_k>1-\varepsilon\}\cap Q_\omega|
		\le
		|\{F\ge 1-\varepsilon\}\cap Q_\omega|.
		\]
		Letting \(\varepsilon\downarrow0\) and using \(0\le F\le1\), we obtain
		\begin{equation}\label{eq:F-one}
			|\{F=1\}\cap Q_\omega|
			\ge \delta_1|Q_\omega|.
		\end{equation}
		
		Similarly, from \eqref{eq:limit-measure-2}, for every \(0<\varepsilon<1/2\),
		\[
		\delta_2|\widehat Q_\omega|
		\le
		\liminf_{k\to\infty}
		|\{\widehat f_k<\varepsilon\}\cap \widehat Q_\omega|
		\le
		|\{F\le\varepsilon\}\cap \widehat Q_\omega|.
		\]
		Letting \(\varepsilon\downarrow0\), we get
		\begin{equation}\label{eq:F-zero}
			|\{F=0\}\cap \widehat Q_\omega|
			\ge \delta_2|\widehat Q_\omega|.
		\end{equation}
		
		Finally, fix \(0<\varepsilon<1/4\). For all sufficiently large \(k\),
		\[
		\{\varepsilon\le \widehat f_k\le1-\varepsilon\}
		\subset
		\{0<\widehat f_k<1-\theta_k\}.
		\]
		Hence, by \eqref{eq:limit-measure-3},
		\[
		\bigl|\{\varepsilon\le \widehat f_k\le1-\varepsilon\}
		\cap(Q_1\cup\widehat Q_1)\bigr|
		\longrightarrow0.
		\]
		Using convergence in measure on \(Q_1\cup\widehat Q_1\), we obtain
		\[
		|\{2\varepsilon\le F\le1-2\varepsilon\}\cap(Q_1\cup\widehat Q_1)|=0.
		\]
		Letting \(\varepsilon\downarrow0\), we conclude that
		\begin{equation}\label{eq:F-binary}
			|\{0<F<1\}\cap(Q_1\cup\widehat Q_1)|=0.
		\end{equation}
		Since \(0\le F\le1\), this means
		\[
		F(Z)\in\{0,1\}
		\qquad\text{for a.e. }Z\in Q_1\cup\widehat Q_1.
		\]
		
		\smallskip
		\noindent
		\emph{Step 4. Elimination of the \(v\)-variable.}
		Let
		\[
		\cO:=(-2,0)\times B_1.
		\]
		By the definition of the normalized cylinders,
		\[
		Q_1=Q_1((0,0,0))=(-1,0]\times B_1\times B_1,
		\qquad
		\widehat Q_1=Q_1((-1,0,0))=(-2,-1]\times B_1\times B_1.
		\]
		In particular,
		\[
		Q_1\cup \widehat Q_1=(\cO\times B_1)\cup (\{0\}\times B_1\times B_1).
		\]
		Therefore, by \eqref{eq:F-binary},
		\[
		F(Z)\in\{0,1\}
		\qquad\text{for a.e. }Z\in \cO\times B_1.
		\]
		
		By \eqref{eq:F-H1v} and Fubini's theorem, for a.e. \((t,x)\in \cO\) the map
		$v\mapsto F(t,x,v)$
		belongs to \(H^1(B_1)\). Hence Lemma~\ref{lem:H1-binary-constant} implies that for
		a.e. \((t,x)\in \cO\), the function \(v\mapsto F(t,x,v)\) is a.e. constant on
		\(B_1\).
		
		Define
		\[
		U(t,x):=\fint_{B_1} F(t,x,v)\dd v
		:=
		\frac{1}{|B_1|}\int_{B_1}F(t,x,v)\dd v
		\qquad\text{for a.e. }(t,x)\in \cO.
		\]
		Then \(U\in L^\infty(\cO)\),
		\[
		U(t,x)\in\{0,1\}
		\qquad\text{for a.e. }(t,x)\in\cO,
		\]
		and
		\begin{equation}\label{eq:F-equals-U}
			F(t,x,v)=U(t,x)
			\qquad\text{for a.e. }(t,x,v)\in \cO\times B_1.
		\end{equation}
		
		Set
		\[
		A_+:=(-\omega^2,0]\times B_{\omega^3},
		\qquad
		A_-:=(-1-\omega^2,-1]\times B_{\omega^3}.
		\]
		Since \(B_\omega\subset B_1\), from \eqref{eq:F-equals-U} and
		\eqref{eq:F-one}--\eqref{eq:F-zero} we deduce
		\begin{equation}\label{eq:U-one}
			|\{U=1\}\cap A_+| \ge \delta_1 |A_+|,
		\end{equation}
		and
		\begin{equation}\label{eq:U-zero}
			|\{U=0\}\cap A_-| \ge \delta_2 |A_-|.
		\end{equation}
		
		\smallskip
		\noindent
		\emph{Step 5. One-sided transport inequalities for \(U\).}
		Fix \(v_0\in B_\omega\). Choose a nonnegative function
		\[
		\zeta\in C_c^\infty(B_{1/2}(0))\quad \text{satisfying }\int_{\bR^d}\zeta(v)\dd v=1,
		\quad
		\int_{\bR^d}v\zeta(v)\dd v=0.
		\]
		Since \(\omega\le \frac14\), we have
		\[
		\supp \zeta(\cdot-v_0)\subset B_{3/4}(0)\subset B_1(0).
		\]
		
		Let \(\phi\in C_c^\infty(\cO)\) be nonnegative, and set
		\[
		\Phi(t,x,v):=\phi(t,x)\zeta(v-v_0).
		\]
		Since $\supp\Phi\Subset \cO\times B_1\subset Q_{3/2}^\circ$, we have
		$\Phi\in C_c^\infty(Q_{3/2}^\circ)$.
		Testing \eqref{eq:F-limit-eq} against \(\Phi\), and using \(\mu\ge0\), gives
		\[
		\left\langle (\partial_t+v\cdot\nabla_x)F,\Phi\right\rangle
		\le
		\iiint_{\cO\times B_1} H_0\,\Phi\dd v\dd x\dd t
		-
		\iiint_{\cO\times B_1} H_1\cdot \nabla_v\Phi\dd v\dd x\dd t .
		\]
		
		We now compute the left-hand side. Since \(F=U\) a.e. on
		\(\cO\times B_1\), and since
		\[
		\int_{B_1}\zeta(v-v_0)\dd v=1,
		\qquad
		\int_{B_1}v\zeta(v-v_0)\dd v=v_0,
		\]
		we have
		\[
		\left\langle (\partial_t+v\cdot\nabla_x)F,\Phi\right\rangle
		=
		\left\langle
		(\partial_t+v_0\cdot\nabla_x)U,\phi
		\right\rangle .
		\]
		Therefore
		\[
		\left\langle
		(\partial_t+v_0\cdot\nabla_x)U,\phi
		\right\rangle
		\le
		\iint_{\cO} g_{v_0}(t,x)\phi(t,x)\dd x\dd t,
		\]
		where
		\[
		g_{v_0}(t,x)
		:=
		\int_{B_1} H_0(t,x,v)\zeta(v-v_0)\dd v
		-
		\int_{B_1} H_1(t,x,v)\cdot \nabla_v\zeta(v-v_0)\dd v.
		\]
		Thus
		\[
		(\partial_t+v_0\cdot\nabla_x)U\le g_{v_0}
		\qquad\text{in }\cD'(\cO).
		\]
		Since \(H_0,H_1\in L^2(Q_{3/2}^\circ)\), we have 
		\[
		g_{v_0}\in L^2(\cO)\subset L^1(\cO).
		\]
		Moreover,
		$U\in\{0,1\}$ almost everywhere in $\cO$.
		By Lemma~\ref{lem:no-upward-jumps}, it follows that
		\begin{equation}\label{eq:U-transport-nonincreasing}
			(\partial_t+v_0\cdot\nabla_x)U\le 0
			\qquad\text{in }\cD'(\cO)
			\qquad\text{for every }v_0\in B_\omega.
		\end{equation}
		
		\smallskip
		\noindent
		\emph{Step 6. Propagation of zero-density points.}
		Let
		\[
		A_-^\circ:=(-1-\omega^2,-1)\times B_{\omega^3},
		\qquad
		A_+^\circ:=(-\omega^2,0)\times B_{\omega^3}.
		\]
		Since \(\partial A_-\) has zero \((d+1)\)-dimensional measure,
		\eqref{eq:U-zero} gives
		\[
		|\{U=0\}\cap A_-^\circ|
		=
		|\{U=0\}\cap A_-|
		>0.
		\]
		Choose a Lebesgue density point
		\[
		Y_-:=(t_-,x_-)\in \{U=0\}\cap A_-^\circ
		\]
		of the set \(\{U=0\}\cap A_-^\circ\). Since \(A_-^\circ\) is open, \(Y_-\) is
		also a density point of \(\{U=0\}\). As \(U\in\{0,1\}\) a.e., this means that
		\(U\) has density \(0\) at \(Y_-\), namely
		\[
		\lim_{r\downarrow 0}
		\fint_{B_r^{d+1}(Y_-)} U(t,x)\dd x\dd t
		=0.
		\]
		
		Now fix an arbitrary point
		\[
		Y:=(t,x)\in A_+^\circ.
		\]
		Set
		\[
		s:=t-t_->0,
		\qquad
		v_0:=\frac{x-x_-}{s}.
		\]
		Since $t\in(-\omega^2,0)$ and $t_-\in(-1-\omega^2,-1)$,
		we have
		$1-\omega^2<s<1+\omega^2$.
		Moreover, $|x-x_-|\le |x|+|x_-|<2\omega^3$,
		and hence
		\[
		|v_0|
		\le \frac{2\omega^3}{1-\omega^2}
		\le \omega
		\qquad\text{because }\omega\le \frac14.
		\]
		Thus \(v_0\in B_\omega\).
		
		We claim that the characteristic segment
		\[
		\Gamma:=\{(t_-+\tau,x_-+\tau v_0):0\le \tau\le s\}
		\]
		is compactly contained in \(\cO=(-2,0)\times B_1\). Indeed, for
		\(0\le \tau\le s\), we have $t_-+\tau\in [t_-,t]\subset(-2,0)$
		and
		\[
		|x_-+\tau v_0|
		\le |x_-|+\tau|v_0|
		\le \omega^3+(1+\omega^2)\omega
		\le \omega+2\omega^3
		\le \frac14+\frac{2}{4^3}
		<1.
		\]
		Hence \(\Gamma\Subset\cO\).
		
		Since \(v_0\in B_\omega\), the inequality \eqref{eq:U-transport-nonincreasing}
		holds for this \(v_0\). Therefore Lemma~\ref{lem:propagate-zero-density}
		applies and yields that \(U\) has density \(0\) at \(Y=(t,x)\). Since
		\(Y\in A_+^\circ\) was arbitrary, \(U\) has density \(0\) at every point of
		\(A_+^\circ\).
		
		We now conclude that \(U=0\) a.e. on \(A_+^\circ\). Indeed, if
		\[
		|\{U=1\}\cap A_+^\circ|>0,
		\]
		then, by the Lebesgue differentiation theorem, there would exist a density
		point of \(\{U=1\}\cap A_+^\circ\). Since \(A_+^\circ\) is open, this point
		would also be a density point of \(\{U=1\}\), and hence the density of \(U\)
		there would be \(1\). This contradicts the fact that \(U\) has density \(0\) at
		every point of \(A_+^\circ\). Therefore
		$U=0$ almost everywhere in $A_+^\circ$.
		Since \(\partial A_+\) has zero \((d+1)\)-dimensional measure, it follows that
		$U=0$ almost everywhere in $A_+$.
		Hence
		\[
		|\{U=1\}\cap A_+|=0.
		\]
		This contradicts \eqref{eq:U-one}, which states that
		\[
		|\{U=1\}\cap A_+|
		\ge \delta_1 |A_+|>0.
		\]
		The contradiction proves the desired estimate and completes the proof of
		Proposition~\ref{prop:isop}.
	\end{proof}

	\begin{proof}[\textbf{Proof of Lemma~\ref{lem:positive-part-defect}}]
		We first prove the defect-measure formulation and the uniform local bounds, and
		then prove the compactness assertion.
		
		For \(R>0\), set
		\[
		U_R(Z_k):=Z_k\circ U_R
		=
		\left\{
		(t,x,v): |t-t_k|<R^2,\ 
		|x-x_k-(t-t_k)v_k|<R^3,\ 
		|v-v_k|<R
		\right\}.
		\]
		Define
		\[
		\bar f_k:=f_k^+\,1_{Q_2^+(Z_k)}
		\qquad\text{on }U_2(Z_k),
		\]
		where \(\bar f_k\) is understood as the zero extension of \(f_k^+\) from
		\(Q_2^+(Z_k)\) to \(U_2(Z_k)\). Then
		\[
		\widehat f_k(Z)=\bar f_k(Z_k\circ Z),
		\qquad Z\in U_2.
		\]
		
		\smallskip
		\noindent
		\emph{Step 1. Local \(L^2\)-bound for \(\nabla_v\bar f_k\).}
		Fix \(0<R<2\), and choose \(R'\in(R,2)\). Take
		\[
		\psi\in C_c^\infty(U_{R'})\qquad \text{satisfying }0\le \psi\le 1,
		\quad
		\psi\equiv 1\quad\text{on }U_R,
		\]
		and set
		\[
		\psi_k(Y):=\psi(Z_k^{-1}\circ Y).
		\]
		Then
		\[
		\psi_k\in C_c^\infty(U_{R'}(Z_k)),
		\qquad
		0\le\psi_k\le1,
		\qquad
		\psi_k\equiv1\quad\text{on }U_R(Z_k).
		\]
		
		Let \(\alpha_n\in C^\infty(\bR)\) satisfy
		\[
		\alpha_n(s)=0\quad\text{for }s\le0,
		\qquad
		0\le \alpha_n(s)\le s_+,
		\qquad
		\alpha_n'\equiv0\quad\text{on }[2,\infty),
		\]
		\[
		\alpha_n(s)\uparrow s_+
		\quad\text{and}\quad
		0\le\alpha_n'(s)\uparrow 1_{\{s>0\}}
		\qquad\text{for }s\le\tfrac32,
		\qquad
		\alpha_n''\ge0
		\qquad\text{for }s\le1 .
		\]
		Such a family exists: mollify \(s\mapsto(s-\tfrac1n)_+\) at scale
		\(\tfrac1{4n}\), and then modify the result on \([\tfrac32,2]\) so that
		\(\alpha_n'\equiv0\) on \([2,\infty)\); the modification leaves
		\(\alpha_n\) unchanged for \(s\le\tfrac32\).
		Set
		\[
		\beta_n(s):=\alpha_n(s)^2.
		\]
		Then \(\alpha_n\), \(\alpha_n'\), and \(\alpha_n''\) are bounded, being
		constant outside \([0,2]\); hence \(\beta_n\in C^2(\bR)\),
		\(\beta_n\ge0\), \(\beta_n(0)=\beta_n'(0)=0\) and
		\(\sup_{\bR}(|\beta_n'|+|\beta_n''|)<\infty\), as required in
		Lemma~\ref{260312141}.  Since \(f_k\le1\) a.e.\ in \(Q_2^+(Z_k)\),
		only the values \(s\le1\) are evaluated below, and there
		\[
		0\le \beta_n(f_k)\le 1,
		\qquad
		0\le \beta_n'(f_k)
		=
		2\alpha_n(f_k)\alpha_n'(f_k)
		\le 2\alpha_n'(f_k),
		\]
		and, for \(s\le1\),
		\[
		\beta_n''(s)
		=
		2(\alpha_n'(s))^2+2\alpha_n(s)\alpha_n''(s)
		\ge 2(\alpha_n'(s))^2\ge0 .
		\]
		
		Applying the one-sided renormalized inequality \eqref{260113330} to \(f_k\),
		with \(\beta=\beta_n\) and test function \(\psi_k^2\), the incoming boundary
		contribution vanishes because
		\[
		f_k\le0\quad\text{on }(\bR\times\gamma^-)\cap Q_2(Z_k),
		\qquad
		\beta_n=0\quad\text{on }(-\infty,0].
		\]
		Hence
		\begin{align*}
			&\iiint_{Q_2^+(Z_k)}
			\beta_n''(f_k)a_k^{ij}D_i f_kD_j f_k\,\psi_k^2\dd Y
			\\
			\le\;&
			\iiint_{Q_2^+(Z_k)}
			\beta_n'(f_k)g_k\psi_k^2\dd Y
			+
			\iiint_{Q_2^+(Z_k)}
			\beta_n(f_k)
			\bigl|(\partial_t+v\cdot\nabla_x)(\psi_k^2)\bigr|\dd Y
			\\
			&\quad
			-
			2\iiint_{Q_2^+(Z_k)}
			\beta_n'(f_k)a_k^{ij}D_j f_k\,\psi_kD_i\psi_k\dd Y
			+
			\iiint_{Q_2^+(Z_k)}
			\beta_n'(f_k)b_k^iD_i f_k\,\psi_k^2\dd Y .
		\end{align*}
		The first two terms on the right-hand side are bounded uniformly in \(k,n\),
		because
		\[
		0\le\beta_n(f_k)\le1,
		\qquad
		0\le\beta_n'(f_k)\le2,
		\qquad
		|g_k|\le1,
		\]
		and
		\[
		\nabla_v \psi_k(Z)=\big(\nabla_v \psi\big)(Z_k^{-1}\circ Z),\qquad (\p_t+v\cdot \nabla_x)\psi_k(Z)=\big((\p_t+v\cdot \nabla_x) \psi\big)(Z_k^{-1}\circ Z).
		\]
		
		For the third term, using \(\beta_n'(f_k)\le2\alpha_n'(f_k)\),
		\(\beta_n''\ge2(\alpha_n')^2\), the boundedness of \(a_k\), and Young's
		inequality, we get
		\begin{align*}
			&\left|
			2\iiint_{Q_2^+(Z_k)}
			\beta_n'(f_k)a_k^{ij}D_j f_k\,\psi_kD_i\psi_k\dd Y
			\right|
			\\
			\le\;&
			\frac{\lambda}{4}
			\iiint_{Q_2^+(Z_k)}
			\beta_n''(f_k)|\nabla_v f_k|^2\psi_k^2\dd Y
			+
			C(R,R').
		\end{align*}
		Similarly,
		\[
		\left|
		\iiint_{Q_2^+(Z_k)}
		\beta_n'(f_k)b_k^iD_i f_k\,\psi_k^2\dd Y
		\right|
		\le
		\frac{\lambda}{4}
		\iiint_{Q_2^+(Z_k)}
		\beta_n''(f_k)|\nabla_v f_k|^2\psi_k^2\dd Y
		+
		C(R,R').
		\]
		Absorbing these terms into the left-hand side and using ellipticity, we obtain
		\[
		\iiint_{Q_2^+(Z_k)}
		\beta_n''(f_k)|\nabla_v f_k|^2\psi_k^2\dd Y
		\le
		C(R,R').
		\]
		Therefore,
		\[
		\iiint_{Q_2^+(Z_k)}
		(\alpha_n'(f_k))^2|\nabla_v f_k|^2\psi_k^2\dd Y
		\le
		C(R,R').
		\]
		Letting \(n\to\infty\), we get
		\[
		\iiint_{U_R(Z_k)}
		|\nabla_v\bar f_k|^2\dd Y
		=
		\iiint_{U_R(Z_k)\cap Q_2^+(Z_k)\cap\{f_k>0\}}
		|\nabla_v f_k|^2\dd Y
		\le
		C(R,R').
		\]
		Since $\widehat f_k(Z)=\bar f_k(Z_k\circ Z)$
		and kinetic translations have unit Jacobian and preserve \(\nabla_v\), this
		gives
		\begin{align}\label{2607014501}
			\|\nabla_v\widehat f_k\|_{L^2(U_R)}
			\le C(R,R').
		\end{align}
		
		\smallskip
		\noindent
		\emph{Step 2. Defect-measure formulation on \(U_2(Z_k)\).}
		Applying the same one-sided renormalized inequality, this time with
		\(\beta=\alpha_n\), and then letting \(n\to\infty\), we obtain, for every
		nonnegative \(\varphi\in C_c^\infty(U_2(Z_k))\),
		\begin{align*}
			&\iiint_{U_2(Z_k)}
			\bar f_k(-\partial_t\varphi-v\cdot\nabla_x\varphi)\dd Y
			+
			\iiint_{U_2(Z_k)}
			1_{Q_2^+(Z_k)}a_k^{ij}D_j\bar f_kD_i\varphi\dd Y
			\\
			&\quad
			-
			\iiint_{U_2(Z_k)}
			1_{Q_2^+(Z_k)}b_k^iD_i\bar f_k\,\varphi\dd Y
			\\
			\le\;&
			\iiint_{U_2(Z_k)}
			1_{Q_2^+(Z_k)\cap\{f_k>0\}}g_k\,\varphi\dd Y .
		\end{align*}
		Define
		\[
		\bar H_{1,k}
		:=
		\left(
		1_{Q_2^+(Z_k)}a_k^{ij}D_j\bar f_k
		\right)_{i=1,\ldots,d},
		\qquad
		\bar H_{0,k}
		:=
		1_{Q_2^+(Z_k)\cap\{f_k>0\}}g_k
		+
		1_{Q_2^+(Z_k)}b_k^iD_i\bar f_k .
		\]
		Then the preceding inequality becomes
		\[
		\left\langle
		(\partial_t+v\cdot\nabla_x)\bar f_k
		-\bar H_{0,k}
		-\nabla_v\cdot\bar H_{1,k},
		\varphi
		\right\rangle
		\le0
		\]
		for every nonnegative \(\varphi\in C_c^\infty(U_2(Z_k))\). By the
		Riesz--Schwartz theorem, there exists a nonnegative Radon measure
		$\bar\mu_k\in\cM_{\mathrm{loc}}(U_2(Z_k))$
		such that
		\[
		(\partial_t+v\cdot\nabla_x)\bar f_k
		=
		\bar H_{0,k}+\nabla_v\cdot\bar H_{1,k}-\bar\mu_k
		\qquad\text{in }\cD'(U_2(Z_k)).
		\]
		
		Now pull this identity back to \(U_2\). Define
		\[
		H_{0,k}(Z):=\bar H_{0,k}(Z_k\circ Z),
		\qquad
		H_{1,k}(Z):=\bar H_{1,k}(Z_k\circ Z),
		\]
		and define \(\mu_k\in\cM_{\mathrm{loc}}(U_2)\) by
		\[
		\mu_k(E):=\bar\mu_k(Z_k\circ E),
		\qquad E\Subset U_2.
		\]
		Since the kinetic translation preserves
		$\partial_t+v\cdot\nabla_x$, $\nabla_v$, and $\dd Z$,
		we get
		\begin{align*}
			(\partial_t+v\cdot\nabla_x)\widehat f_k
			=
			H_{0,k}+\nabla_v\cdot H_{1,k}-\mu_k
			\qquad\text{in }\cD'(U_2).
		\end{align*}
		This proves \eqref{eq:positive-part-defect-identity}.
		
		\smallskip
		\noindent
		\emph{Step 3. Local bounds for \(H_{0,k}\), \(H_{1,k}\), and \(\mu_k\).}
		Since \(0\le\widehat f_k\le1\), we have
		$\|\widehat f_k\|_{L^\infty(U_R)}\le1$.
		By \eqref{con_ab} and \eqref{2607014501},
		\begin{align*}
			\begin{aligned}
				\|H_{1,k}\|_{L^2(U_R)}
				\le
				C\|\nabla_v\widehat f_k\|_{L^2(U_R)}
				\,&\le C(R),\\
				\|H_{0,k}\|_{L^2(U_R)}
				\le
				\|1_{\{0<\widehat f_k\}}\,g_k(Z_k\circ\cdot)\|_{L^2(U_R)}
				+
				C\|\nabla_v\widehat f_k\|_{L^2(U_R)}
				\,&\le C(R).
			\end{aligned}
		\end{align*}
		
		It remains to estimate \(\mu_k(U_R)\). Choose \(R<R''<2\), and take
		\[
		\chi\in C_c^\infty(U_{R''})\qquad \text{satisfying }0\le\chi\le1,
		\quad
		\chi\equiv1\quad\text{on }U_R.
		\]
		Testing \eqref{eq:positive-part-defect-identity} against \(\chi\), we have
		\begin{align*}
			\mu_k(U_R)
			&\le
			\int_{U_2}\chi\dd\mu_k
			=
			\iiint_{U_2}
			\left(
			H_{0,k}\chi
			-
			H_{1,k}\cdot\nabla_v\chi
			+
			\widehat f_k(\partial_t+v\cdot\nabla_x)\chi
			\right)\dd Z .
		\end{align*}
		Therefore,
		\begin{align*}
			\mu_k(U_R)
			\le\;&
			\|H_{0,k}\|_{L^2(U_{R''})}\|\chi\|_{L^2}
			+
			\|H_{1,k}\|_{L^2(U_{R''})}\|\nabla_v\chi\|_{L^2}
			\\
			&\quad
			+
			\|\widehat f_k\|_{L^\infty(U_{R''})}
			\|(\partial_t+v\cdot\nabla_x)\chi\|_{L^1}.
		\end{align*}
		The right-hand side is bounded uniformly in \(k\). Hence
		$\sup_k\mu_k(U_R)<\infty$.
		This proves \eqref{eq:positive-part-local-bounds}.
		
		\smallskip
		\noindent
		\emph{Step 4. Precompactness in \(L^1_{\mathrm{loc}}(U_2)\).}
		Fix \(0<R<2\), and choose
		\[
		\phi\in C_c^\infty(U_2)\qquad \text{satisfying }0\le \phi\le 1,
		\quad
		\phi\equiv1
		\quad\text{on a neighborhood of }U_R.
		\]
		Set
		\[
		F_k:=\phi\widehat f_k .
		\]
		Since \(\operatorname{supp}\phi\Subset U_2\), there exists \(R_\phi<2\) such
		that
		$\operatorname{supp}\phi\subset U_{R_\phi}$.
		From \eqref{eq:positive-part-defect-identity}, we have
		\[
		(\partial_t+v\cdot\nabla_x)F_k
		=
		G_{0,k}+\nabla_v\cdot G_{1,k}+\nu_k
		\qquad\text{in }\cD'(\bR^{1+2d}),
		\]
		where
		\[
		G_{0,k}
		:=
		\phi H_{0,k}
		-
		H_{1,k}\cdot\nabla_v\phi
		+
		\widehat f_k(\partial_t+v\cdot\nabla_x)\phi,
		\]
		\[
		G_{1,k}:=\phi H_{1,k},
		\qquad
		\nu_k:=-\phi\mu_k .
		\]
		Here \(F_k,G_{0,k},G_{1,k}\), and \(\nu_k\) are extended by zero outside
		\(\operatorname{supp}\phi\).
		
		By \eqref{eq:positive-part-local-bounds}, we have
		\begin{align}\label{eq:step4-uniform-bounds}
			\begin{aligned}
				\sup_k\Big(
				&\|F_k\|_{L^2(\bR^{1+2d})}
				+
				\|\nabla_vF_k\|_{L^2(\bR^{1+2d})}
				+
				\|G_{0,k}\|_{L^2(\bR^{1+2d})}
				\\
				&\qquad
				+
				\|G_{1,k}\|_{L^2(\bR^{1+2d})}
				+
				\|\nu_k\|_{\cM(\bR^{1+2d})}
				\Big)
				<\infty .
			\end{aligned}
		\end{align}
		Moreover, all these objects are supported in the fixed compact set
		\(\operatorname{supp}\phi\).
		
		We now apply Proposition~\ref{prop:velocity-average-compactness} to the
		right-hand side of the equation for \(F_k\).  Fix \(r:=\frac12\) and
		\(p\in\bigl(1,\frac{d+1}{d+1-r}\bigr)\); since
		\(\frac{d+1}{d+1-r}<\frac{d+1}{d}\le2\), we have \(1<p<2\), and the
		compact support of \(F_k\) gives \(L^2\subset L^p\) on
		\(\operatorname{supp}\phi\).  Together with
		Lemma~\ref{lem:bouchut-hypo}, applied with this \(r\), this gives,
		with \(s:=\frac{1-r}{3}=\frac16\) and independently of \(k\),
		\begin{equation}\label{eq:step4-hypo-bound}
			\sup_k
			\left(
			\|D_t^sF_k\|_{L^p(\bR^{1+2d})}
			+
			\|D_x^sF_k\|_{L^p(\bR^{1+2d})}
			\right)
			<\infty .
		\end{equation}
		On the other hand, \eqref{eq:step4-uniform-bounds} and the compact support of
		\(F_k\) imply
		\[
		\sup_k
		\left(
		\|F_k\|_{L^p(\bR^{1+2d})}
		+
		\|\nabla_vF_k\|_{L^p(\bR^{1+2d})}
		\right)
		<\infty .
		\]
		The quantities in \eqref{eq:step4-hypo-bound} are Bessel-potential
		(\(H^{s,p}_{t,x}\)) norms; since \(H^{s,p}\hookrightarrow W^{s',p}\)
		for every \(0<s'<s\) and every \(1<p<\infty\), they control the
		Sobolev--Slobodeckij norm of any strictly smaller order.  Fix
		\(s'\in(0,s)\).  Since \(\nabla_vF_k\in L^p\), the last estimate gives
		the corresponding \(W_v^{s',p}\)-bound, and we obtain
		\[
		\sup_k
		\|F_k\|_{W^{s',p}(\bR^{1+2d})}
		<\infty .
		\]
		
		Choose a ball \(B\subset\bR^{1+2d}\) such that
		$\operatorname{supp}\phi\subset B$.
		Since the supports of \(F_k\) are contained in \(B\), and since the embedding
		$W^{s',p}(B)\hookrightarrow L^p(B)$
		is compact, it follows that \((F_k)_{k\in\bN}\) is precompact in \(L^p(B)\).
		
		Because \(\phi\equiv1\) on \(U_R\), we have
		\[
		F_k=\widehat f_k
		\qquad\text{on }U_R.
		\]
		Thus \((\widehat f_k)_{k\in\bN}\) is precompact in \(L^p(U_R)\). Since
		\(|U_R|<\infty\), it is also precompact in \(L^1(U_R)\). As \(0<R<2\) is
		arbitrary, a diagonal argument yields that \((\widehat f_k)_{k\in\bN}\) is
		precompact in \(L^1_{\mathrm{loc}}(U_2)\).
	\end{proof}

	\subsection{A measure-to-pointwise estimate}
	
	\begin{lemma}
		\label{lem:boundary-measure-to-pointwise}
		Fix \(0<\omega\le \frac14\) and \(\delta_2\in(0,1)\). Then there exist
		constants
		\[
		\lambda_{\mathrm{bd}}
		\in(0,1),
		\qquad
		\beta_{\mathrm{bd}}
		>0,
		\]
		depending only on $d,\lambda,\Lambda,M,\omega,\delta_2$,
		such that the following holds.
		
		Let
		\[
		Z_0=(t_0,x_0,v_0)\in
		\bR\times \overline{\bR^d_+}\times \bR^d,
		\qquad
		\widehat Z_0:=(t_0-1,x_0-v_0,v_0),
		\]
		and suppose that \(a^{ij}\), \(b^i\), \(f\), and \(g\) satisfy, in
		\(Q_2^+(Z_0)\), the equation and structural assumptions appearing in
		Theorem~\ref{cor:scaled-boundary-holder-estimate} with \(r=1\).
		Assume moreover that
		\[
		f\le 1
		\qquad\text{and}\qquad
		|g|\le \beta_{\mathrm{bd}}
		\quad\text{a.e. in }Q_2^+(Z_0),
		\]
		and
		\[
		\left|
		\bigl(\{f\le 0\}\cap Q_\omega^+(\widehat Z_0)\bigr)
		\cup
		Q_\omega^-(\widehat Z_0)
		\right|
		\ge
		\delta_2\,|Q_\omega(\widehat Z_0)|.
		\]
		Then
		\[
		\operatorname*{ess\,sup}_{Q_{\omega/2}^+(Z_0)} f
		\le
		1-\lambda_{\mathrm{bd}}.
		\]
	\end{lemma}
	
	\begin{proof}
		By Corollary~\ref{cor:scaled-smallness-implies-half}, applied with
		\(r=\omega\) and \(K=1\), there exists
		\[
		\kappa_\omega=\kappa_\omega(d,\lambda,\Lambda,M,\omega)>0
		\]
		such that the following implication holds: If \(u\) is a weak solution in
		\(Q_\omega^+(Z_0)\), with source term \(G\), satisfying
		\[
		u\le 0
		\qquad\text{on }(\bR\times\gamma^-)\cap Q_\omega(Z_0)
		\]
		in the sense of traces, and
		\[
		\|u_+\|_{L^2(Q_\omega^+(Z_0))}\le \kappa_\omega,
		\qquad
		\|G\|_{L^\infty(Q_\omega^+(Z_0))}\le 1,
		\]
		then
		\begin{equation}\label{eq:bd-mtp-half}
			u\le \frac12
			\qquad\text{a.e. in }Q_{\omega/2}^+(Z_0).
		\end{equation}
		
		Set
		\[
		\delta_1
		:=
		\min\left\{
		\frac12,\frac{\kappa_\omega^2}{|Q_\omega(0)|}
		\right\}\in(0,1).
		\]
		Let $\theta\in(0,1)$ and $\nu\in(0,1)$
		be the constants furnished by Proposition~\ref{prop:isop}. Choose
		\(k_0\in\bN\) such that
		\[
		k_0\nu>1.
		\]
		Define
		\[
		\beta_{\mathrm{bd}}:=\theta^{k_0},
		\qquad
		\lambda_{\mathrm{bd}}:=\frac12\,\theta^{k_0}.
		\]
		
		For \(k=0,1,\ldots,k_0\), set
		\[
		f_k:=\theta^{-k}\bigl(f-(1-\theta^k)\bigr).
		\]
		Then \(f_0=f\), and
		\[
		f_{k+1}=\frac1\theta\bigl(f_k-(1-\theta)\bigr).
		\]
		Since \(f\le1\), we have
		\[
		f_k\le1
		\qquad\text{a.e. in }Q_2^+(Z_0)
		\]
		for every \(k\). Moreover, \(f_k\) satisfies
		\[
		(\partial_t+v\cdot\nabla_x)f_k
		=
		D_i(a^{ij}D_j f_k)+b^iD_i f_k+g_k
		\qquad\text{in }Q_2^+(Z_0),
		\]
		where
		\[
		g_k:=\theta^{-k}g.
		\]
		Since \(|g|\le\beta_{\mathrm{bd}}=\theta^{k_0}\), we have
		\[
		|g_k|\le1
		\qquad\text{a.e. in }Q_2^+(Z_0),
		\qquad 0\le k\le k_0.
		\]
		Also, since \(f=0\) on \((\bR\times\gamma^-)\cap Q_2(Z_0)\), we have
		\[
		f_k=1-\theta^{-k}\le0
		\qquad\text{on }(\bR\times\gamma^-)\cap Q_2(Z_0)
		\]
		in the sense of traces.
		
		The sets \(\{f_k\le0\}\) are increasing in \(k\), and
		$\{f\le0\}\subset\{f_k\le0\}$.
		Therefore the lower-phase assumption implies that, for every \(0\le k\le k_0\),
		\begin{equation}\label{eq:bd-mtp-lower-iterates}
			\left|
			\bigl(\{f_k\le0\}\cap Q_\omega^+(\widehat Z_0)\bigr)
			\cup Q_\omega^-(\widehat Z_0)
			\right|
			\ge
			\delta_2\,|Q_\omega(\widehat Z_0)|.
		\end{equation}
		
		We claim that there exists \(k_1\in\{1,\ldots,k_0\}\) such that
		\begin{equation}\label{eq:bd-mtp-small-positive}
			|\{f_{k_1}>0\}\cap Q_\omega^+(Z_0)|
			<
			\delta_1\,|Q_\omega(Z_0)|.
		\end{equation}
		Suppose not. Then, for every \(1\le k\le k_0\),
		\[
		|\{f_k>0\}\cap Q_\omega^+(Z_0)|
		\ge
		\delta_1\,|Q_\omega(Z_0)|.
		\]
		Since
		\[
		\{f_{k+1}>0\}
		=
		\{f_k>1-\theta\}
		\subset
		\{f_k\ge1-\theta\},
		\]
		we get, for every \(0\le k\le k_0-1\),
		\[
		|\{f_k\ge1-\theta\}\cap Q_\omega^+(Z_0)|
		\ge
		\delta_1\,|Q_\omega(Z_0)|.
		\]
		Combining this with \eqref{eq:bd-mtp-lower-iterates}, Proposition~\ref{prop:isop}
		applied to \(f_k\) gives
		\[
		|\{0<f_k<1-\theta\}
		\cap
		\bigl(Q_1^+(Z_0)\cup Q_1^+(\widehat Z_0)\bigr)|
		\ge
		\nu\,|Q_1(Z_0)\cup Q_1(\widehat Z_0)|
		\]
		for every \(0\le k\le k_0-1\).
		
		Set
		\[
		\cQ^+:=Q_1^+(Z_0)\cup Q_1^+(\widehat Z_0),
		\qquad
		\cQ:=Q_1(Z_0)\cup Q_1(\widehat Z_0).
		\]
		Since
		\[
		\{f_k\le0\}\cup\{0<f_k<1-\theta\}
		\subset
		\{f_{k+1}\le0\},
		\]
		we have
		\[
		|\{f_{k+1}\le0\}\cap\cQ^+|
		\ge
		|\{f_k\le0\}\cap\cQ^+|
		+
		\nu\,|\cQ|
		\]
		for every \(0\le k\le k_0-1\). Iterating yields
		\[
		|\{f_{k_0}\le0\}\cap\cQ^+|
		\ge
		k_0\nu\,|\cQ|.
		\]
		However,
		\[
		|\{f_{k_0}\le0\}\cap\cQ^+|
		\le
		|\cQ^+|
		\le
		|\cQ|,
		\]
		which contradicts \(k_0\nu>1\). This proves
		\eqref{eq:bd-mtp-small-positive}.
		
		Fix \(k_1\in\{1,\ldots,k_0\}\) satisfying
		\eqref{eq:bd-mtp-small-positive}. Since \(f_{k_1}\le1\),
		\[
		\|(f_{k_1})_+\|_{L^2(Q_\omega^+(Z_0))}^2
		\le
		|\{f_{k_1}>0\}\cap Q_\omega^+(Z_0)|
		<
		\delta_1\,|Q_\omega(Z_0)|
		=
		\delta_1\,|Q_\omega(0)|
		\le
		\kappa_\omega^2.
		\]
		Together with \(|g_{k_1}|\le1\), \eqref{eq:bd-mtp-half} applied to
		\(u=f_{k_1}\) gives
		$f_{k_1}\le\frac12$ almost everywhere in $Q_{\omega/2}^+(Z_0)$.
		Returning to \(f\), we obtain
		\[
		f
		=
		\theta^{k_1}f_{k_1}+(1-\theta^{k_1})
		\le
		1-\frac12\,\theta^{k_1}
		\le
		1-\frac12\,\theta^{k_0}
		=
		1-\lambda_{\mathrm{bd}}
		\]
		a.e. in \(Q_{\omega/2}^+(Z_0)\). This completes the proof.
	\end{proof}

	\begin{corollary}
		\label{cor:scaled-boundary-measure-to-pointwise}
		Fix \(0<\omega\le \frac14\) and \(\delta\in(0,1)\), and let $\lambda_{\mathrm{bd}}\in(0,1)$ and $\beta_{\mathrm{bd}}>0$
		be the constants from Lemma~\ref{lem:boundary-measure-to-pointwise} with \(\delta_2:=\delta\).
		Let \(0<r\le1\), and let
		\[
		Z_0=(t_0,x_0,v_0)\in
		\bR\times \overline{\bR^d_+}\times \bR^d,
		\qquad
		\widehat Z_{0,r}:=(t_0-r^2,\;x_0-r^2v_0,\;v_0),
		\]
		and suppose that \(a^{ij}\), \(b^i\), \(f\), and \(g\) satisfy the equation and structural assumptions appearing in
		Theorem~\ref{cor:scaled-boundary-holder-estimate}.
		Assume moreover that
		\[
		f\le 1
		\qquad\text{and}\qquad
		|g|\le \beta_{\mathrm{bd}}
		\quad\text{a.e. in }Q_{2r}^+(Z_0),
		\]
		and
		\[
		\left|
		\bigl(\{f\le 0\}\cap Q_{\omega r}^+(\widehat Z_{0,r})\bigr)
		\cup
		Q_{\omega r}^-(\widehat Z_{0,r})
		\right|
		\ge
		\delta\,|Q_{\omega r}(\widehat Z_{0,r})|.
		\]
		Then
		\[
		\operatorname*{ess\,sup}_{Q_{\omega r/2}^+(Z_0)} f
		\le
		1-\lambda_{\mathrm{bd}}.
		\]
	\end{corollary}
	
	\begin{proof}
		Set \(K:=r^{-1}\) and, in the notation of Lemma~\ref{lem:kinetic-scaling},
		\[
		\tilde f:=f_K,\qquad
		\tilde a^{ij}:=a^{ij}_K,\qquad
		\tilde b^i:=\frac{b^i_K}{K},\qquad
		\tilde g:=\frac{g_K}{K^2}.
		\]
		The kinetic dilation \(Z\mapsto Z_K\) preserves \(\{x^d>0\}\),
		\(\gamma^-\), \(\gamma^+\) and \(\gamma^0\), and maps
		\(Q_{2r}^\pm(Z_0)\) onto \(Q_2^\pm(Z_{0,K})\). Since
		\[
		\bigl(\widehat Z_{0,r}\bigr)_K
		=\bigl(r^{-2}t_0-1,\ r^{-3}x_0-r^{-1}v_0,\ r^{-1}v_0\bigr)
		=\widehat{\bigl(Z_{0,K}\bigr)},
		\]
		it maps \(Q_{\omega r}^\pm(\widehat Z_{0,r})\) onto
		\(Q_\omega^\pm\bigl(\widehat{Z_{0,K}}\bigr)\).  In particular, the
		measure hypothesis and the zero-inflow trace condition are carried over
		verbatim. Moreover, \(Z_{0,K}\in\bR\times\overline{\bR^d_+}\times\bR^d\).
		By Lemma~\ref{lem:kinetic-scaling}, \(\tilde f\) is a weak solution of
		\[
		(\partial_t+v\cdot\nabla_x)\tilde f
		=
		D_i(\tilde a^{ij}D_j\tilde f)
		+\tilde b^iD_i\tilde f
		+\tilde g
		\qquad\text{in }Q_2^+(Z_{0,K}),
		\]
		with the same ellipticity constants \(\lambda,\Lambda\) and, since
		\(0<r\le1\),
		\[
		|\tilde b|\le rM\le M,
		\qquad
		|\tilde g|\le r^2\beta_{\mathrm{bd}}\le\beta_{\mathrm{bd}}.
		\]
		Therefore Lemma~\ref{lem:boundary-measure-to-pointwise}, applied at the
		base point \(Z_{0,K}\), gives
		\[
		\operatorname*{ess\,sup}_{Q_{\omega/2}^+(Z_{0,K})}\tilde f
		\le
		1-\lambda_{\mathrm{bd}},
		\]
		and undoing the dilation yields
		\[
		\operatorname*{ess\,sup}_{Q_{\omega r/2}^+(Z_0)} f
		\le
		1-\lambda_{\mathrm{bd}}.
		\]
		This completes the proof.
	\end{proof}
	
	\begin{corollary}
		\label{cor:full-cylinder-measure-to-pointwise}
		Fix \(0<\omega\le \frac14\) and \(\delta\in(0,1)\), and let $\lambda_{\mathrm{bd}}\in(0,1)$ and $
		\beta_{\mathrm{bd}}
		>0
		$
		be the constants from Lemma~\ref{lem:boundary-measure-to-pointwise} with \(\delta_2:=\delta\).
		Let \(0<r\le1\), and let
		\[
		Z_0=(t_0,x_0,v_0)\in\bR^{1+2d},
		\qquad
		\widehat Z_{0,r}:=(t_0-r^2,\;x_0-r^2v_0,\;v_0).
		\]
		Assume that \(a^{ij}\) and \(b^i\) satisfy \eqref{con_ab} in
		\(Q_{2r}(Z_0)\). Let \(f\) be a weak solution of
		\[
		(\partial_t+v\cdot\nabla_x)f
		=
		D_i(a^{ij}D_jf)+b^iD_if+g
		\qquad\text{in }Q_{2r}(Z_0).
		\]
		Assume moreover that
		\[
		f\le 1
		\qquad\text{and}\qquad
		|g|\le \beta_{\mathrm{bd}}
		\quad\text{a.e. in }Q_{2r}(Z_0),
		\]
		and
		\[
		\left|
		\{f\le 0\}\cap Q_{\omega r}(\widehat Z_{0,r})
		\right|
		\ge
		\delta\,|Q_{\omega r}(\widehat Z_{0,r})|.
		\]
		Then
		\begin{align}\label{260630446}
			\operatorname*{ess\,sup}_{Q_{\omega r/2}(Z_0)} f
			\le
			1-\lambda_{\mathrm{bd}}.
		\end{align}
	\end{corollary}

	\begin{proof}
		Choose \(N>0\) so large that, with
		\[
		Z_0^N:=(t_0,x_0+Ne_d,v_0),
		\qquad
		\widehat Z_{0,r}^N:=(t_0-r^2,x_0+Ne_d-r^2v_0,v_0),
		\]
		we have
		$Q_{2r}(Z_0^N)\subset\{x^d>0\}$.
		Then
		\[
		Q_{2r}^+(Z_0^N)=Q_{2r}(Z_0^N),
		\qquad
		Q_{\omega r}^+(\widehat Z_{0,r}^N)=Q_{\omega r}(\widehat Z_{0,r}^N),
		\qquad
		Q_{\omega r}^-(\widehat Z_{0,r}^N)=\emptyset.
		\]
		
		Translate the solution and coefficients by setting
		\[
		F(t,x,v):=f(t,x-Ne_d,v),
		\]
		and defining the translated coefficients and source term analogously. Then
		\(F\) satisfies the hypotheses of
		Corollary~\ref{cor:scaled-boundary-measure-to-pointwise} in
		\(Q_{2r}^+(Z_0^N)=Q_{2r}(Z_0^N)\); the zero-inflow condition is vacuous, and the measure assumption is exactly the translated version of the one assumed for \(f\).
		Therefore
		\[
		\operatorname*{ess\,sup}_{Q_{\omega r/2}^+(Z_0^N)}F
		\le
		1-\lambda_{\mathrm{bd}}.
		\]
		Since \(Q_{\omega r/2}^+(Z_0^N)=Q_{\omega r/2}(Z_0^N)\), translating back gives \eqref{260630446}.
	\end{proof}

	\subsection{Proposition~\ref{cor:boundary-oscillation-decay} and conclusion}\label{260701422}

	\begin{lemma}
		\label{prop:boundary-oscillation-improvement}
		Let \(Z_0\in \bR\times\overline{\bR^d_+}\times\bR^d\), and suppose that \(a^{ij}\), \(b^i\), \(f\), and \(g\) satisfy, in
		\(Q_2^+(Z_0)\), the equation and structural assumptions appearing in
		Theorem~\ref{cor:scaled-boundary-holder-estimate} with \(r=1\).
		
		Then there exist constants 
		\[ 
		\lambda_0\in(0,1), \qquad \omega\in\left(0,\frac12\right), \qquad \beta>0, 
		\] 
		depending only on \(d\), \(\lambda\), \(\Lambda\), and \(M\), such that if 
		\[
		\operatorname*{osc}_{Q_2^+(Z_0)} f\le 2, \qquad |g|\le \beta \quad\text{a.e. in }Q_2^+(Z_0),
		\] 
		then
		\[
		\operatorname*{osc}_{Q_{\omega/2}^+(Z_0)} f \le 2-\lambda_0.
		\]
	\end{lemma}
	
	\begin{proof}
		Fix
		\[
		\rho:=\frac14.
		\]
		Let \(\varepsilon_->0\) be the constant in
		Lemma~\ref{lem:grazing-shifted-minus-volume} corresponding to
		the value \(\rho\) of its parameter \(\omega\). Thus, if
		\[
		Q_{\rho^3/(8d)}(Z_0)\cap(\bR\times\gamma^0)\neq\emptyset,
		\]
		then
		\[
		|Q_{\rho}^-(\widehat Z_0)|
		\ge
		\varepsilon_-\,|Q_{\rho}(\widehat Z_0)|,
		\qquad
		\widehat Z_0:=(t_0-1,x_0-v_0,v_0).
		\]
		Set
		\begin{align*}
			\delta_*:=\min\left\{\frac12,\varepsilon_-\right\}.
		\end{align*}
		Let $(\lambda_0,\beta_0)$ be the common constants furnished by Corollaries~\ref{cor:full-cylinder-measure-to-pointwise} and
		\ref{cor:scaled-boundary-measure-to-pointwise}, with their parameters $\omega$ and $\delta$ equal to $\rho$ and $\delta_*$.
		Define
		\[
		r_0:=\frac{\rho^3}{16d},
		\qquad
		\omega:=r_0\rho,
		\qquad
		\beta:=\beta_0.
		\]
		Then
		\[
		\lambda_0\in(0,1),
		\qquad
		\omega\in\left(0,\frac12\right),
		\qquad
		\beta>0.
		\]
		
		Assume that \(f\) satisfies the hypotheses of the lemma. We prove
		\[
		\operatorname*{osc}_{Q_{\omega/2}^+(Z_0)}f
		\le
		2-\lambda_0.
		\]
		We distinguish three cases according to the local boundary geometry of
		\(Q_{2r_0}(Z_0)\).
		
		\smallskip
		\noindent
		\textbf{Case 1.}
		$Q_{2r_0}(Z_0)\cap(\bR\times\gamma)
		\subset
		\bR\times\gamma^-$.
		
		Extend the coefficients to the full cylinder by
		\[
		\bar a^{ij}:=a^{ij}\mathbf 1_{Q_{2r_0}^+(Z_0)}
		+\lambda\delta^{ij}\mathbf 1_{Q_{2r_0}(Z_0)\setminus Q_{2r_0}^+(Z_0)},
		\qquad
		\bar b^i:=b^i\mathbf 1_{Q_{2r_0}^+(Z_0)},
		\qquad
		\bar g:=g\mathbf 1_{Q_{2r_0}^+(Z_0)} ;
		\]
		these satisfy \eqref{con_ab} on \(Q_{2r_0}(Z_0)\) with the same
		constants \(\lambda,\Lambda,M\), and they agree with \(a,b\) on
		\(Q_{2r_0}^+(Z_0)\). All the statements below are applied with them.
		Let \(\tilde f\) be the zero extension of \(f\) from \(Q_{2r_0}^+(Z_0)\) to
		\(Q_{2r_0}(Z_0)\). By Lemma~\ref{lem:inflow-zero-extension}, \(\tilde f\) is an
		\(L^2\)-solution in \(Q_{2r_0}(Z_0)\), with source
		\[
		\tilde g:=g\,1_{Q_{2r_0}^+(Z_0)}.
		\]
		Moreover,
		\[
		\operatorname*{osc}_{Q_{2r_0}(Z_0)}\tilde f
		=
		\operatorname*{osc}_{Q_{2r_0}^+(Z_0)}f
		\le 2.
		\]
		
		Set
		\[
		\tilde m:=
		\frac12\left(
		\operatorname*{ess\,sup}_{Q_{2r_0}(Z_0)}\tilde f
		+
		\operatorname*{ess\,inf}_{Q_{2r_0}(Z_0)}\tilde f
		\right),
		\qquad
		\tilde u:=\tilde f-\tilde m.
		\]
		Then \(\tilde u\) solves the same equation with source \(\tilde g\), and
		\[
		\|\tilde u\|_{L^\infty(Q_{2r_0}(Z_0))}
		=
		\frac12\operatorname*{osc}_{Q_{2r_0}(Z_0)}\tilde f
		\le1.
		\]
		
		Define
		\[
		\widehat Z_{0,r_0}:=(t_0-r_0^2,x_0-r_0^2v_0,v_0),
		\qquad
		Q_*:=Q_{\rho r_0}(\widehat Z_{0,r_0}).
		\]
		Since $
		\{\tilde u\le0\}\cup\{-\tilde u\le0\}=Q_*$
		up to a null set, at least one of the two sets
		\[
		\{\tilde u\le0\}\cap Q_*,
		\qquad
		\{-\tilde u\le0\}\cap Q_*
		\]
		has measure at least $\frac12|Q_*|\ge \delta_*|Q_*|$.
		Applying Corollary~\ref{cor:full-cylinder-measure-to-pointwise} at scale
		\(r_0\) to either \(\tilde u\) or \(-\tilde u\), and using
		$|\tilde g|\le |g|\le\beta=\beta_0$,
		we obtain either
		\[
		\operatorname*{ess\,sup}_{Q_{\rho r_0/2}(Z_0)}\tilde u
		\le
		1-\lambda_0\qquad \text{or}\qquad 
		\operatorname*{ess\,sup}_{Q_{\rho r_0/2}(Z_0)}(-\tilde u)
		\le
		1-\lambda_0.
		\]
		Since \(|\tilde u|\le1\), either alternative gives
		\[
		\operatorname*{osc}_{Q_{\rho r_0/2}(Z_0)}\tilde f
		=
		\operatorname*{osc}_{Q_{\rho r_0/2}(Z_0)}\tilde u
		\le
		2-\lambda_0.
		\]
		Since \(\tilde f=f\) on \(Q_{2r_0}^+(Z_0)\) and
		\[
		Q_{\omega/2}^+(Z_0)
		=
		Q_{\rho r_0/2}^+(Z_0)
		\subset
		Q_{\rho r_0/2}(Z_0),
		\]
		we conclude that
		\[
		\operatorname*{osc}_{Q_{\omega/2}^+(Z_0)}f
		\le
		2-\lambda_0.
		\]
		
		\smallskip
		\noindent
		\textbf{Case 2.} $
		Q_{2r_0}(Z_0)\cap(\bR\times\gamma)
		\subset
		\bR\times\gamma^+$.
		
		In this case, $(\bR\times\gamma^-)\cap Q_{2r_0}(Z_0)=\emptyset$.
		Set
		\[
		m:=
		\frac12\left(
		\operatorname*{ess\,sup}_{Q_{2r_0}^+(Z_0)}f
		+
		\operatorname*{ess\,inf}_{Q_{2r_0}^+(Z_0)}f
		\right),
		\qquad
		u:=f-m.
		\]
		Then \(u\) solves the same equation with source \(g\) in \(Q_{2r_0}^+(Z_0)\),
		and the zero-inflow boundary condition in
		Corollary~\ref{cor:scaled-boundary-measure-to-pointwise} is vacuous. Moreover,
		\[
		\|u\|_{L^\infty(Q_{2r_0}^+(Z_0))}
		=
		\frac12\operatorname*{osc}_{Q_{2r_0}^+(Z_0)}f
		\le1.
		\]
		
		Define
		\[
		\widehat Z_{0,r_0}:=(t_0-r_0^2,x_0-r_0^2v_0,v_0),
		\]
		and write
		\[
		Q_*^+:=Q_{\rho r_0}^+(\widehat Z_{0,r_0}),
		\qquad
		Q_*^-:=Q_{\rho r_0}^-(\widehat Z_{0,r_0}),
		\qquad
		Q_*:=Q_{\rho r_0}(\widehat Z_{0,r_0}).
		\]
		Since $\{u\le0\}\cup\{-u\le0\}=Q_*^+$ 
		up to a null set, we have
		\[
		\begin{aligned}
			&\left|
			(\{u\le0\}\cap Q_*^+)\cup Q_*^-
			\right|
			+
			\left|
			(\{-u\le0\}\cap Q_*^+)\cup Q_*^-
			\right|
			\ge 2|Q_*^-|+|Q_*^+|
			\ge
			|Q_*|.
		\end{aligned}
		\]
		Hence at least one of the two quantities is at least $\frac12|Q_*|\ge \delta_*|Q_*|$.
		Applying Corollary~\ref{cor:scaled-boundary-measure-to-pointwise} at scale
		\(r_0\) to either \(u\) or \(-u\), and using 
		$|g|\le\beta=\beta_0$,
		we obtain either 
		\[
		\operatorname*{ess\,sup}_{Q_{\rho r_0/2}^+(Z_0)}u
		\le
		1-\lambda_0\quad \text{or}\quad 
		\operatorname*{ess\,sup}_{Q_{\rho r_0/2}^+(Z_0)}(-u)
		\le
		1-\lambda_0.
		\]
		Since \(|u|\le1\), either alternative gives
		\[
		\operatorname*{osc}_{Q_{\rho r_0/2}^+(Z_0)}f
		=\operatorname*{osc}_{Q_{\rho r_0/2}^+(Z_0)}u
		\le
		2-\lambda_0.
		\]
		Since 
		$Q_{\omega/2}^+(Z_0)
		=
		Q_{\rho r_0/2}^+(Z_0)$, 
		we conclude that
		\[
		\operatorname*{osc}_{Q_{\omega/2}^+(Z_0)}f
		\le
		2-\lambda_0.
		\]
		
		\smallskip
		\noindent
		\textbf{Case 3.} $Q_{2r_0}(Z_0)\cap(\bR\times\gamma^0)\neq\emptyset$.
		
		Since
		$2r_0=\frac{\rho^3}{8d}$, 
		Lemma~\ref{lem:grazing-shifted-minus-volume} gives
		\[
		|Q_{\rho}^-(\widehat Z_0)|
		\ge
		\varepsilon_-\,|Q_{\rho}(\widehat Z_0)|
		\ge
		\delta_*\,|Q_{\rho}(\widehat Z_0)|.
		\]
		Moreover, since \(Q_{2r_0}(Z_0)\) meets \(\bR\times\gamma^0\), its
		\(v\)-section is the ball \(B_{2r_0}(v_0)\), which contains a point
		with \(v^d=0\) and hence an open set of points with \(v^d>0\); therefore
		\[
		(\bR\times\gamma^-)\cap Q_{2r_0}(Z_0)
		\]
		has positive \(\dd\gamma\dd t\)-measure, and the boundary condition gives
		\(f_\gamma=0\) there. Hence, by
		Corollary~\ref{cor:trace-bound-scaled-cylinder} applied with
		\(R=2r_0\),
		\[
		\operatorname*{ess\,inf}_{Q_{2r_0}^+(Z_0)} f
		\le
		\operatorname*{ess\,inf}_{(\bR\times\gamma)\cap Q_{2r_0}(Z_0)}f_\gamma
		\le 0
		\le
		\operatorname*{ess\,sup}_{(\bR\times\gamma)\cap Q_{2r_0}(Z_0)}f_\gamma
		\le
		\operatorname*{ess\,sup}_{Q_{2r_0}^+(Z_0)} f.
		\]
		In particular,
		\[
		\operatorname*{ess\,inf}_{Q_2^+(Z_0)} f
		\le 0
		\le
		\operatorname*{ess\,sup}_{Q_2^+(Z_0)} f.
		\]
		
		Set
		\[
		A:=\operatorname*{ess\,sup}_{Q_2^+(Z_0)}f,
		\qquad
		B:=-\operatorname*{ess\,inf}_{Q_2^+(Z_0)}f,
		\qquad
		L:=\max\{A,B,1\}.
		\]
		Choose \(\sigma\in\{\pm1\}\) so that
		\[
		A_\sigma:=
		\operatorname*{ess\,sup}_{Q_2^+(Z_0)}(\sigma f)
		=
		\max\{A,B\}.
		\]
		Define
		\[
		u:=\frac{\sigma f}{L}.
		\]
		Then
		\[
		u\le1
		\quad\text{a.e. in }Q_2^+(Z_0),\qquad u=0
		\quad\text{on }(\bR\times\gamma^-)\cap Q_2(Z_0)
		\]
		in the sense of traces, and
		\[
		\left|\frac{\sigma g}{L}\right|
		\le |g|
		\le \beta
		=
		\beta_0
		\qquad\text{a.e. in }Q_2^+(Z_0).
		\]
		The lower-phase condition in
		Corollary~\ref{cor:scaled-boundary-measure-to-pointwise}, with scale \(1\), is
		satisfied because \(Q_{\rho}^-(\widehat Z_0)\) is contained in
		\[
		\bigl(\{u\le0\}\cap Q_{\rho}^+(\widehat Z_0)\bigr)
		\cup
		Q_{\rho}^-(\widehat Z_0).
		\]
		Therefore
		\[
		\operatorname*{ess\,sup}_{Q_{\rho/2}^+(Z_0)}u
		\le
		1-\lambda_0.
		\]
		Equivalently,
		\[
		\operatorname*{ess\,sup}_{Q_{\rho/2}^+(Z_0)}(\sigma f)
		\le
		L(1-\lambda_0).
		\]
		
		Let
		\[
		B_\sigma:=
		-\operatorname*{ess\,inf}_{Q_2^+(Z_0)}(\sigma f).
		\]
		Since \(A,B\ge0\) and \(A_\sigma=\max\{A,B\}\), we have
		\[
		0\le B_\sigma\le A_\sigma\le L.
		\]
		Thus
		\[
		\operatorname*{osc}_{Q_{\rho/2}^+(Z_0)}f
		=
		\operatorname*{osc}_{Q_{\rho/2}^+(Z_0)}(\sigma f)
		\le
		L(1-\lambda_0)+B_\sigma.
		\]
		If \(L=A_\sigma\), then \(A_\sigma\ge1\), and
		\[
		L(1-\lambda_0)+B_\sigma
		=
		A_\sigma+B_\sigma-\lambda_0A_\sigma
		\le
		\operatorname*{osc}_{Q_2^+(Z_0)}f-\lambda_0
		\le
		2-\lambda_0.
		\]
		If \(L=1>A_\sigma\), then \(B_\sigma\le A_\sigma<1\), and hence
		\[
		L(1-\lambda_0)+B_\sigma
		<
		2-\lambda_0.
		\]
		Therefore
		\[
		\operatorname*{osc}_{Q_{\rho/2}^+(Z_0)}f
		\le
		2-\lambda_0.
		\]
		Since
		\[
		Q_{\omega/2}^+(Z_0)
		=
		Q_{\rho r_0/2}^+(Z_0)
		\subset
		Q_{\rho/2}^+(Z_0),
		\]
		we conclude that
		\[
		\operatorname*{osc}_{Q_{\omega/2}^+(Z_0)}f
		\le
		2-\lambda_0.
		\]
		
		The three cases exhaust all possibilities.  Indeed,
		\(Q_{2r_0}(Z_0)=E\times B_{2r_0}(v_0)\) with
		\(E:=\{(t,x):t\in(t_0-4r_0^2,t_0],\ |x-x_0-(t-t_0)v_0|<8r_0^3\}\).
		Hence, if \(Q_{2r_0}(Z_0)\) meets \(\bR\times\gamma\), every
		\(v\in B_{2r_0}(v_0)\) occurs in \(Q_{2r_0}(Z_0)\cap(\bR\times\gamma)\).
		If some occurring \(v\) has \(v^d>0\) and another has \(v^d<0\), then,
		by connectedness of \(B_{2r_0}(v_0)\), some occurring \(v\) has
		\(v^d=0\), which is Case~3.  If \(Q_{2r_0}(Z_0)\) does not meet
		\(\bR\times\gamma\), Case~1 holds vacuously.  This completes the proof.
	\end{proof}

	\begin{proof}[\textbf{Proof of Proposition~\ref{cor:boundary-oscillation-decay}}]
		Let
		\[
		\lambda_0\in(0,1),
		\qquad
		\omega\in\left(0,\frac12\right),
		\qquad
		\beta>0
		\]
		be the constants from Lemma~\ref{prop:boundary-oscillation-improvement}.
		Set
		\[
		N_0
		:=
		\frac12\operatorname*{osc}_{Q_2^+(Z_0)}f
		+
		\beta^{-1}\|g\|_{L^\infty(Q_2^+(Z_0))}.
		\]
		If \(N_0=0\) or \(N_0=\infty\), then there is nothing to prove. Hence we may
		assume that \(0<N_0<\infty\).
		
		Define
		\[
		\tilde f:=N_0^{-1}f,
		\qquad
		\tilde g:=N_0^{-1}g.
		\]
		Then \(\tilde f\) satisfies the same equation with source term \(\tilde g\),
		and the zero-inflow boundary condition is preserved. Moreover,
		\[
		\operatorname*{osc}_{Q_2^+(Z_0)}\tilde f
		\le 2,
		\qquad
		\|\tilde g\|_{L^\infty(Q_2^+(Z_0))}
		\le \beta.
		\]
		Therefore Lemma~\ref{prop:boundary-oscillation-improvement} applied to
		\(\tilde f\) gives
		\[
		\operatorname*{osc}_{Q_{\omega/2}^+(Z_0)}\tilde f
		\le
		2-\lambda_0.
		\]
		Scaling back, we obtain
		\[
		\operatorname*{osc}_{Q_{\omega/2}^+(Z_0)}f
		\le
		(2-\lambda_0)N_0.
		\]
		By the definition of \(N_0\),
		\[
		(2-\lambda_0)N_0
		=
		\left(1-\frac{\lambda_0}{2}\right)
		\operatorname*{osc}_{Q_2^+(Z_0)}f
		+
		(2-\lambda_0)\beta^{-1}
		\|g\|_{L^\infty(Q_2^+(Z_0))}.
		\]
		Thus
		\[
		\operatorname*{osc}_{Q_{\omega/2}^+(Z_0)}f
		\le
		\left(1-\frac{\lambda_0}{2}\right)
		\operatorname*{osc}_{Q_2^+(Z_0)}f
		+
		2\beta^{-1}
		\|g\|_{L^\infty(Q_2^+(Z_0))}.
		\]
		The desired estimate follows with \(C:=2\beta^{-1}\).
	\end{proof}

	\begin{proof}[\textbf{Proof of Proposition~\ref{2601271146}}]
		\noindent\emph{Step 1. The case \(M\le 1\).}
		Let
		\[
		\lambda_0\in(0,1),
		\qquad
		\omega\in\left(0,\frac12\right),
		\qquad
		C_0>0
		\]
		be the constants in Proposition~\ref{cor:boundary-oscillation-decay}, with
		\(M=1\). Set
		\[
		\rho_0:=\frac{\omega}{4}.
		\]
		Since \(\rho_0\le\frac18<\frac12\le1-\frac{\lambda_0}{2}\), the number
		\(\alpha:=\log(1-\lambda_0/2)/\log\rho_0\) lies in \((0,1)\).  With this
		choice,
		\[
		1-\frac{\lambda_0}{2}\le \rho_0^\alpha.
		\]
		
		We first prove that, for every \(Z\in Q_1^+(Z_0)\),
		\begin{align}\label{260401533}
			\sup_{n\ge1}\left[
			\rho_0^{-\alpha n}
			\operatorname*{osc}_{Q_{\rho_0^n}^+(Z)}f
			\right]\le
			C\left(
			\|f\|_{L^\infty(Q_2^+(Z_0))}
			+
			\|g\|_{L^\infty(Q_2^+(Z_0))}
			\right).
		\end{align}
		Indeed, \(Z\in Q_1(Z_0)\) and \(\rho_0^n\le\rho_0\le1\), so
		Lemma~\ref{lem:cylinder-inclusion} gives
		\(Q_{\rho_0^n}(Z)\subset Q_{1+\rho_0^n}(Z_0)\subset Q_2(Z_0)\); hence
		the cylinders appearing below are contained in \(Q_2^+(Z_0)\).
		
		Fix \(n\ge1\) and set
		\[
		s_n:=\frac{\rho_0^n}{2}.
		\]
		Applying Lemma~\ref{lem:kinetic-scaling} with \(K=s_n^{-1}\), \(c=0\),
		\(g^i=0\), and \(g^0=g\), Proposition~\ref{cor:boundary-oscillation-decay}
		gives the scaled estimate
		\[
		\operatorname*{osc}_{Q_{\omega s_n/2}^+(Z)}f
		\le
		\left(1-\frac{\lambda_0}{2}\right)
		\operatorname*{osc}_{Q_{2s_n}^+(Z)}f
		+
		C_0s_n^2
		\|g\|_{L^\infty(Q_{2s_n}^+(Z))}.
		\]
		Since
		\[
		2s_n=\rho_0^n,
		\qquad
		\frac{\omega s_n}{2}
		=
		\frac{\omega}{4}\rho_0^n
		=
		\rho_0^{n+1},
		\]
		we obtain
		\[
		\operatorname*{osc}_{Q_{\rho_0^{n+1}}^+(Z)}f
		\le
		\left(1-\frac{\lambda_0}{2}\right)
		\operatorname*{osc}_{Q_{\rho_0^n}^+(Z)}f
		+
		C\rho_0^{2n}
		\|g\|_{L^\infty(Q_2^+(Z_0))}.
		\]
		Therefore, by the choice of \(\alpha\),
		\[
		\rho_0^{-\alpha(n+1)}
		\operatorname*{osc}_{Q_{\rho_0^{n+1}}^+(Z)}f
		\le
		\rho_0^{-\alpha n}
		\operatorname*{osc}_{Q_{\rho_0^n}^+(Z)}f
		+
		C\rho_0^{n(2-\alpha)}
		\|g\|_{L^\infty(Q_2^+(Z_0))}.
		\]
		Iterating yields
		\[
		\rho_0^{-\alpha n}
		\operatorname*{osc}_{Q_{\rho_0^n}^+(Z)}f
		\le
		C\left(
		\|f\|_{L^\infty(Q_2^+(Z_0))}
		+
		\|g\|_{L^\infty(Q_2^+(Z_0))}
		\right),
		\]
		which proves \eqref{260401533}.
		
		For \(Z=(t,x,v)\in Q_1^+(Z_0)\) one has \(Q_\rho(Z)=Q_\rho^+(Z)\)
		whenever \(\rho^2|v^d|+\rho^3<x^d\), so the cylinders
		\(Q_{\rho_0^n}^+(Z)\) form an admissible differentiation basis at
		\(Z\).  Since \(\mathbf d\) induces the Euclidean topology and the
		measure is doubling for \(\mathbf d\), Lebesgue's differentiation
		theorem and \eqref{260401533} produce a representative of \(f\), still
		denoted by \(f\), defined at every \(Z\in Q_1^+(Z_0)\) and satisfying
		\(|f(Z)-f(W)|\le\operatorname*{osc}_{Q_{\rho_0^n}^+(Z)}f\) for every
		\(W\in Q_{\rho_0^n}^+(Z)\).
		
		Let now \(Z\ne W\) in \(Q_1^+(Z_0)\), labeled so that the time
		coordinate of \(W\) does not exceed that of \(Z\), and put
		\(\sigma:=\mathbf d(Z,W)\).  If \(2\sigma<\rho_0\), choose
		\(n\ge1\) with \(\rho_0^{n+1}\le2\sigma<\rho_0^n\); by
		Lemma~\ref{260127226}, \(W\in Q_{2\sigma}(Z)\cap Q_1^+(Z_0)\subset
		Q_{\rho_0^n}^+(Z)\), so \eqref{260401533} gives
		\(|f(Z)-f(W)|\le C\rho_0^{\alpha n}(\cdots)\le
		C(2/\rho_0)^\alpha\sigma^\alpha(\cdots)\).  If
		\(2\sigma\ge\rho_0\), then
		\(|f(Z)-f(W)|\le2\|f\|_{L^\infty(Q_2^+(Z_0))}
		\le2(2/\rho_0)^\alpha\sigma^\alpha\|f\|_{L^\infty(Q_2^+(Z_0))}\).
		In both cases
		\begin{align}\label{2603181216}
			\frac{|f(Z)-f(W)|}{\mathbf d(Z,W)^\alpha}
			\le
			C\left(
			\|f\|_{L^\infty(Q_2^+(Z_0))}
			+
			\|g\|_{L^\infty(Q_2^+(Z_0))}
			\right)
		\end{align}
		for all distinct \(Z,W\in Q_1^+(Z_0)\).  Since
		\(1\le(M+1)^\alpha\le2\) and \((M+1)^{-2}\ge\frac14\) when
		\(M\le1\), this proves the desired estimate in that case.
		
		\noindent\emph{Step 2. Reduction to the case \(M\le1\).}
		We now consider the general case \(M>1\). Put
		\[
		R:=2(M+1).
		\]
		For \(Z=(t,x,v)\), write
		\[
		Z_R:=(R^2t,R^3x,Rv).
		\]
		Define
		\[
		\widetilde f(Z):=f(Z_{R^{-1}}),
		\qquad
		\widetilde a^{ij}(Z):=a^{ij}(Z_{R^{-1}}),
		\]
		\[
		\widetilde b^i(Z):=b^i(Z_{R^{-1}}),
		\qquad
		\widetilde g(Z):=g(Z_{R^{-1}}).
		\]
		By Lemma~\ref{lem:kinetic-scaling}, applied with \(K=R\), \(c=0\),
		\(g^i=0\), and \(g^0=g\), the function \(\widetilde f\) is an \(L^2\)-solution of
		\[
		(\partial_t+v\cdot\nabla_x)\widetilde f
		=
		D_i(\widetilde a^{ij}D_j\widetilde f)
		+
		\frac{\widetilde b^i}{R}D_i\widetilde f
		+
		\frac{\widetilde g}{R^2}
		\qquad\text{in }Q_{2R}^+(Z_{0,R}),
		\]
		with zero-inflow boundary condition on
		$(\bR\times\gamma^-)\cap Q_{2R}(Z_{0,R})$.
		Moreover,
		\[
		\left\|\frac{\widetilde b}{R}\right\|_{L^\infty(Q_{2R}^+(Z_{0,R}))}
		\le
		\frac{M}{2(M+1)}
		\le 1.
		\]
		Thus Step 1 can be applied locally to \(\widetilde f\).
		
		Let \(Z,W\in Q_1^+(Z_0)\), \(Z\neq W\), and set
		\[
		\widetilde Z:=Z_R,
		\qquad
		\widetilde W:=W_R.
		\]
		First assume that
		\[
		\mathbf d(Z,W)<\frac{1}{4R}.
		\]
		Then \(\mathbf d(\widetilde Z,\widetilde W)=R\,\mathbf d(Z,W)<\frac14\),
		so Lemma~\ref{260127226}(2) gives \(\widetilde W\in Q_{1/2}(\widetilde Z)\),
		up to interchanging \(Z\) and \(W\).  Moreover,
		\(\widetilde Z\in Q_R(Z_{0,R})\) and \(R\ge2\), so
		Lemma~\ref{lem:cylinder-inclusion} gives
		\[
		Q_2^+(\widetilde Z)\subset Q_{R+2}^+(Z_{0,R})\subset Q_{2R}^+(Z_{0,R}).
		\]
		Applying \eqref{2603181216} to \(\widetilde f\) on the unit cylinder centered
		at \(\widetilde Z\), and using
		$\mathbf d(\widetilde Z,\widetilde W)
		=
		R\,\mathbf d(Z,W)$,
		we get
		\[
		\frac{|f(Z)-f(W)|}{\mathbf d(Z,W)^\alpha}
		=
		R^\alpha
		\frac{|\widetilde f(\widetilde Z)-\widetilde f(\widetilde W)|}
		{\mathbf d(\widetilde Z,\widetilde W)^\alpha}
		\]
		\[
		\le
		C R^\alpha
		\left(
		\|\widetilde f\|_{L^\infty(Q_2^+(\widetilde Z))}
		+
		R^{-2}
		\|\widetilde g\|_{L^\infty(Q_2^+(\widetilde Z))}
		\right).
		\]
		Since \(Q_2^+(\widetilde Z)\subset Q_{2R}^+(Z_{0,R})\), this implies
		\[
		\frac{|f(Z)-f(W)|}{\mathbf d(Z,W)^\alpha}
		\le
		C R^\alpha
		\left(
		\|f\|_{L^\infty(Q_2^+(Z_0))}
		+
		R^{-2}\|g\|_{L^\infty(Q_2^+(Z_0))}
		\right).
		\]
		
		On the other hand, if $\mathbf d(Z,W)\ge \frac{1}{4R}$, 
		then the trivial estimate gives
		\[
		\frac{|f(Z)-f(W)|}{\mathbf d(Z,W)^\alpha}
		\le
		2(4R)^\alpha
		\|f\|_{L^\infty(Q_2^+(Z_0))}.
		\]
		Combining the two cases and using \(R=2(M+1)\), we obtain
		\[
		\frac{|f(Z)-f(W)|}{\mathbf d(Z,W)^\alpha}
		\le
		C(M+1)^\alpha
		\left(
		\|f\|_{L^\infty(Q_2^+(Z_0))}
		+
		(M+1)^{-2}\|g\|_{L^\infty(Q_2^+(Z_0))}
		\right).
		\]
		This proves Theorem~\ref{cor:scaled-boundary-holder-estimate} when \(r=1\).
	\end{proof}

	\section{One-dimensional barrier at the grazing boundary}\label{sec:half-line-barrier}
	
	This appendix records the constant-coefficient barrier used in the proof of
	Lemma~\ref{260605518}.
	
	Let
	\begin{align*}
		\Phi(r):=
		\begin{cases}
			(-r)\, U\left(\frac{1}{6},\frac{4}{3},\frac{-r^3}{9}\right),
			& r<0,\\[2mm]
			9^{1/3}\dfrac{\Gamma(1/3)}{\Gamma(1/6)},
			& r=0,\\[2mm]
			\dfrac{r}{6}\exp\left(-\frac {r^3} 9\right)
			U\left(\frac{7}{6},\frac{4}{3},\frac{r^3}{9}\right),
			& r>0.
		\end{cases}
	\end{align*}
	Here \(U(a,b,z)\) denotes Tricomi's confluent hypergeometric function,
	defined by
	\[
	U(a,b,z)
	=
	\frac{1}{\Gamma(a)}
	\int_0^\infty \ee^{-zt}t^{a-1}(1+t)^{b-a-1}\dd t
	\]
	for \(a>0\) and \(z>0\). 
	Define
	\[
	H(x,v):=x^{1/6}\Phi\bigl(x^{-1/3}v\bigr),
	\qquad x>0,\ v\in\bR,
	\]
	which is a classical explicit solution of the stationary
	Kolmogorov equation in the half-space with absorbing boundary condition
	\cite{goldman1971first,groeneboom1999integrated}.
	
	\begin{lemma}\label{260603514}
		\begin{enumerate}
			\item \(H>0\) on \((0,\infty)\times\bR\), and
			\begin{align}\label{260605412}
				H(x,v)\simeq
				\begin{cases}
					(-v)^{1/2},
					& \text{if } v<-x^{1/3},\\[1mm]
					x^{1/6},
					& \text{if } |v|<x^{1/3},\\[1mm]
					xv^{-5/2}\exp\left(-\dfrac{v^3}{9x}\right),
					& \text{if } x^{1/3}<v.
				\end{cases}
			\end{align}
			Moreover,
			\begin{align}\label{260602546}
				\lim_{(x,v)\to(0,v_0)}H(x,v)
				=
				9^{1/6}|v_0|^{1/2}\mathbf 1_{\{v_0<0\}}
				\qquad\text{for every }v_0\in\bR .
			\end{align}
			In particular, combining the three regimes in \eqref{260605412},
			\begin{align}\label{260605413}
				H(x,v)\gtrsim|v|^{1/2}\mathbf 1_{\{v<0\}},
				\qquad
				H(x,v)\lesssim x^{1/6}+|v|^{1/2},
				\qquad x>0,\ v\in\bR .
			\end{align}
			
			\item We extend \(H\) to \([0,\infty)\times\bR\) by the boundary values in
			\eqref{260602546}. Then
			\[
			H\in C([0,\infty)\times\bR)
			\cap
			C^\infty\left(\bigl([0,\infty)\times\bR\bigr)\setminus\{(0,0)\}\right).
			\]
			Moreover,
			\begin{align}\label{260602541}
				v\partial_x H-\partial_{vv}H=0
				\qquad\text{in }(0,\infty)\times\bR .
			\end{align}
			
			\item Define
			\[
			\eta(r):=1+|r|^3 1_{\{r>0\}}.
			\]
			Then
			\begin{align}\label{260603518}
				\frac{|\p_v H(x,v)|}{H(x,v)}
				&\simeq
				\frac{\eta(x^{-1/3}v)}{x^{1/3}+|v|},\qquad
				\frac{|\p_{vv}H(x,v)|}{H(x,v)}
				\lesssim
				\left(
				\frac{\eta(x^{-1/3}v)}{x^{1/3}+|v|}
				\right)^2 .
			\end{align}
		\end{enumerate}
	\end{lemma}
	\begin{proof}
		(1) For \(a>0\) and \(z>0\), the integral representation of \(U(a,b,z)\)
		implies \(U(a,b,z)>0\). Moreover, the standard asymptotic formulas for
		Tricomi's function \(U\), see \textit{e.g.},
		\cite[Section~4.1.2]{slater1960confluent}, give
		\begin{align}\label{260602542}
			\lim_{z\to0+}
			\frac{U(a,b,z)}{z^{1-b}}
			=
			\frac{\Gamma(b-1)}{\Gamma(a)},
			\qquad
			\lim_{z\to+\infty}
			\frac{U(a,b,z)}{z^{-a}}
			=
			1,
			\qquad a>0,\ b>1.
		\end{align}
		Applying these estimates to the definition of \(\Phi\), we obtain
		\begin{align*}
			\Phi(r)
			&=
			9^{1/6}|r|^{1/2}
			\bigl(1+O(|r|^{-3})\bigr)
			\quad \text{as } r\to-\infty,
			\\
			\Phi(r)
			&=
			\frac{9^{7/6}}{6}
			r^{-5/2}
			\exp\left(-\frac{r^3}{9}\right)
			\bigl(1+O(r^{-3})\bigr)
			\quad \text{as } r\to+\infty.
		\end{align*}
		Together with the positivity of \(\Phi\), these asymptotics imply the estimates in \eqref{260605412}.
		
		It remains to prove the boundary limit.
		If \(v_0<0\), then \(x^{-1/3}v\to-\infty\), and therefore \[ H(x,v) = x^{1/6}\Phi(x^{-1/3}v) \to 9^{1/6}|v_0|^{1/2}. \]
		If \(v_0>0\), then \(x^{-1/3}v\to+\infty\), and the exponential decay of \(\Phi(r)\) as \(r\to+\infty\) gives $H(x,v)\to0$.
		Finally, if \(v_0=0\), the conclusion follows by combining the three regimes in \eqref{260605412}.
		Hence \[ \lim_{(x,v)\to(0,v_0)}H(x,v) = 9^{1/6}|v_0|^{1/2}\mathbf 1_{\{v_0<0\}} \qquad\text{for every }v_0\in\bR .
		\]
		
		\medskip
		\noindent(2)
		We first recall that
		\begin{align}\label{26060262}
			\partial_z U(a,b,z)=-aU(a+1,b+1,z).
		\end{align}
		On \((-\infty,0)\) and on \((0,\infty)\) the two branches of \(\Phi\)
		solve one and the same linear equation
		$\Phi''+\frac{r^2}{3}\Phi'-\frac r6\Phi=0$
		with polynomial coefficients.
		To verify the gluing data explicitly, use the
		two-term expansion
		\[
		U\left(a,\frac43,z\right)
		=
		\frac{\Gamma(1/3)}{\Gamma(a)}z^{-1/3}
		+
		\frac{\Gamma(-1/3)}{\Gamma(a-1/3)}
		+O(z^{2/3})
		\qquad(z\downarrow0).
		\]
		For both branches this gives
		\[
		\Phi(r)
		=
		9^{1/3}\frac{\Gamma(1/3)}{\Gamma(1/6)}
		-
		\frac{\Gamma(2/3)}{2\Gamma(5/6)}\,r
		+O(|r|^3)
		\qquad(r\to0^\pm).
		\]
		Thus the two branches have the same value and first derivative at \(0\).
		Since \(r=0\) is an ordinary point of the displayed ODE, uniqueness for its
		Cauchy problem shows that the two branches are restrictions of a single
		entire solution, and $\Phi\in C^\infty(\bR)$.
		Moreover, since \(U(a,b,z)\) satisfies
		\[
		\left(z\partial_{zz}+(b-z)\partial_z-a\right)U(a,b,z)=0,
		\]
		a direct computation gives
		\[
		\left(
		\partial_{rr}
		+\frac{r^2}{3}\partial_r
		-\frac r6
		\right)\Phi(r)=0
		\qquad \text{for } r\neq0.
		\]
		Since \(\Phi\in C^\infty(\bR)\), this identity holds on all of \(\bR\).
		
		With \(r:=x^{-1/3}v\), we obtain
		\[
		v\partial_xH-\partial_{vv}H
		=
		x^{-1/2}
		\left(
		\frac r6\Phi(r)-\frac{r^2}{3}\Phi'(r)-\Phi''(r)
		\right)
		=0,
		\]
		which proves \eqref{260602541}.
		
		It remains to discuss the smooth extendability of \(H\) near
		\(\{x=0,\ v\neq0\}\). For \(r\le0\), the integral representation of \(U\)
		gives
		\begin{align}\label{260705611}
			\Phi(r)
			=
			C_- \int_0^\infty
			\ee^{-u^3/9}u^{-1/2}\bigl((-r)^3+u^3\bigr)^{1/6}\dd u
			\qquad \text{for all }r\le 0,
		\end{align}
		where \(C_->0\). Similarly, for \(r\ge0\),
		\begin{align}\label{260705612}
			\Phi(r)
			=
			C_+ \int_r^\infty
			\ee^{-u^3/9}u^{-1/2}(u^3-r^3)^{1/6}\dd u
			\qquad \text{for all }r\ge 0,
		\end{align}
		where \(C_+>0\).
		
		Let first \(v_0<0\). From \eqref{260705611}, for \(x\ge0\) and \(v<0\),
		we have
		\[
		H(x,v)
		=
		C_-(-v)^{1/2}
		F\left(\frac{x}{(-v)^3}\right),
		\]
		where
		\[
		F(\xi):=
		\int_0^\infty
		\ee^{-u^3/9}u^{-1/2}(1+\xi u^3)^{1/6}\dd u,
		\qquad \xi\ge0.
		\]
		Differentiating under the integral sign shows that \(F\) is smooth on
		\((0,\infty)\), and that all derivatives of \(F\) extend continuously to
		\([0,\infty)\). Indeed, the exponential factor controls all polynomial powers
		of \(u\). By Seeley's extension theorem \cite{seeley1964extension}, there
		exists
		$\widetilde F\in C^\infty(\bR)$
		such that
		$\widetilde F(\xi)=F(\xi)$
		for $\xi\ge0$.
		Define
		\[
		\widetilde H(x,v)
		:=
		C_-(-v)^{1/2}
		\widetilde F\left(\frac{x}{(-v)^3}\right),
		\qquad (x,v)\in\bR\times(-\infty,0).
		\]
		Then
		$\widetilde H\in C^\infty\bigl(\bR\times(-\infty,0)\bigr)$
		and
		$\widetilde H(x,v)=H(x,v)$ for $x\geq 0$ and $v<0$.
		Thus \(H\) admits a smooth local extension across \(x=0\) near every point
		\((0,v_0)\) with \(v_0<0\). In particular,
		\[
		H(0,v)
		=
		C_-(-v)^{1/2}F(0)
		=
		9^{1/6}(-v)^{1/2},
		\qquad v<0.
		\]
		
		Next let \(v_0>0\), and choose a compact interval
		\(I\Subset(0,\infty)\) containing \(v_0\). For \(x>0\) and \(v\in I\), using
		\eqref{260705612} and making the change of variables
		\(u^3=r^3+s^3\), where \(r=x^{-1/3}v\), we obtain
		\[
		H(x,v)
		=
		C_+ x\,\exp\left(-\frac{v^3}{9x}\right)
		\int_0^\infty
		\ee^{-s^3/9}s^{5/2}
		\bigl(v^3+xs^3\bigr)^{-5/6}\dd s .
		\]
		Differentiating this expression, every \(\partial_x^k\partial_v^m H\) with
		\(0<x\le1\) and \(v\in I\) is bounded by
		\(C_{k,m}(I)x^{-N_{k,m}}\exp(-c_I/x)\) for some \(N_{k,m}\) and some
		\(c_I>0\), because \(v^3\ge c_I'>0\) on \(I\).  Hence all derivatives of
		\(H\) tend to zero as \(x\downarrow0\), locally uniformly for
		\(v\in I\). Therefore \(H\) admits a smooth local extension across \(x=0\)
		near every point \((0,v_0)\) with \(v_0>0\), for instance by setting it equal
		to \(0\) for \(x\le0\) near \((0,v_0)\).
		
		We now combine the two extensions. Let \(\widetilde F\in C^\infty(\bR)\) be
		the Seeley extension of \(F\) constructed above, and set
		\[
		\cU:=\bR^2\setminus\bigl((-\infty,0]\times \{0\}\bigr).
		\]
		Define \(\widetilde H:\cU\to\bR\) by
		\[
		\widetilde H(x,v):=
		\begin{cases}
			C_-(-v)^{1/2}\widetilde F\left(\dfrac{x}{(-v)^3}\right),
			& v<0,\\[3mm]
			x^{1/6}\Phi(x^{-1/3}v),
			& x>0,\ v\ge0,\\[2mm]
			0,
			& x\le0,\ v>0.
		\end{cases}
		\]
		The three cases are disjoint and the pieces glue smoothly across the two
		interfaces. Indeed, for \(x\ge0\) and
		\(v<0\), the first line agrees with \(H\), since
		\(\widetilde F=F\) on \([0,\infty)\). For \(x>0\), the first two lines match
		smoothly across \(v=0\), because \(H(x,v)=x^{1/6}\Phi(x^{-1/3}v)\) and
		\(\Phi\in C^\infty(\bR)\). Finally, across \(\{x=0,\ v>0\}\), the second line
		is flat in \(x\), as shown above, and therefore matches smoothly with the
		third line. Hence
		$\widetilde H\in C^\infty(\cU)$.
		Moreover,
		\[
		\widetilde H=H
		\qquad\text{on }\bigl([0,\infty)\times\bR\bigr)\setminus\{(0,0)\}.
		\]
		In the convention of this paper, this proves
		\[
		H\in
		C^\infty\left(\bigl([0,\infty)\times\bR\bigr)\setminus\{(0,0)\}\right).
		\]
		
		\smallskip
		\noindent
		(3) For \(r:=x^{-1/3}v\), we have
		\begin{align}\label{26021134}
			\begin{aligned}
				\frac{\partial_v H(x,v)}{H(x,v)}
				&=
				x^{-1/3}\frac{\Phi'(r)}{\Phi(r)}
				=
				\frac{1+|r|}{x^{1/3}+|v|}
				\frac{\Phi'(r)}{\Phi(r)},
				\\
				\frac{\partial_{vv}H(x,v)}{H(x,v)}
				&=
				x^{-2/3}\frac{\Phi''(r)}{\Phi(r)}
				=
				\frac{(1+|r|)^2}{(x^{1/3}+|v|)^2}
				\frac{\Phi''(r)}{\Phi(r)}.
			\end{aligned}
		\end{align}
		
		\smallskip
		\noindent
		\emph{Step 1. Monotonicity of \(\Phi\).}
		Differentiating the integral representations of $\Phi$ in \eqref{260705611} and \eqref{260705612}
		shows that
		$\Phi'(r)<0$ for $r\neq0$.
		Thus \(\Phi\) is strictly decreasing on \(\bR\).
		
		\smallskip
		\noindent
		\emph{Step 2. Logarithmic derivative estimates.}
		By \eqref{26060262}, we have
		\[
		\frac{-\Phi'(r)}{\Phi(r)}
		=
		\begin{cases}
			\dfrac{1}{-r}
			-\dfrac{r^2}{18}
			\dfrac{U\left(7/6,7/3,(-r^3)/9\right)}
			{U\left(1/6,4/3,(-r^3)/9\right)},
			& r<0,\\[3mm]
			\dfrac{\Gamma(2/3)\Gamma(1/6)}
			{2\cdot 9^{1/3}\Gamma(5/6)\Gamma(1/3)},
			& r=0,\\[3mm]
			-\dfrac{1}{r}
			+\dfrac{r^2}{3}
			+\dfrac{7r^2}{18}
			\dfrac{U\left(13/6,7/3,r^3/9\right)}
			{U\left(7/6,4/3,r^3/9\right)},
			& r>0.
		\end{cases}
		\]
		Using \eqref{260602542}, we obtain
		\[
		\lim_{r\to-\infty}
		\frac{\bigl(-\Phi'/\Phi\bigr)(r)}{(-r)^{-1}}
		=
		\frac12,
		\qquad
		\lim_{r\to+\infty}
		\frac{\bigl(-\Phi'/\Phi\bigr)(r)}{r^2}
		=
		\frac13.
		\]
		Together with the positivity and continuity of \(-\Phi'/\Phi\), this implies
		$\frac{-\Phi'(r)}{\Phi(r)}
		\simeq
		\frac{\eta(r)}{1+|r|}$.
		
		Next, using the ODE for \(\Phi\),
		\[
		\frac{\Phi''(r)}{\Phi(r)}
		=
		-\frac{r^2}{3}\frac{\Phi'(r)}{\Phi(r)}
		+\frac r6.
		\]
		As \(r\to+\infty\), the previous estimate gives
		\begin{align}\label{260602845}
			\left|\frac{\Phi''(r)}{\Phi(r)}\right|
			\le
			\frac{r^2}{3}\left|\frac{\Phi'(r)}{\Phi(r)}\right|
			+\frac r6
			\lesssim r^4.
		\end{align}
		On the other hand, using the sharper asymptotic formula
		\[
		U(a,b,z)
		=
		z^{-a}
		\left(
		1-a(1+a-b)z^{-1}+O(z^{-2})
		\right)
		\qquad \text{as } z\to+\infty,
		\]
		we obtain, as \(r\to-\infty\), with \(z:=(-r^3)/9\),
		\begin{align}\label{260602846}
			\frac{\Phi''(r)}{\Phi(r)}
			&=
			-\frac{r^2}{3}\frac{\Phi'(r)}{\Phi(r)}
			+\frac r6 \notag\\
			&=
			\frac r6
			\left(
			\frac{1+7/(36z)+O(z^{-2})}
			{1+1/(36z)+O(z^{-2})}
			-1
			\right)=
			-\frac{1}{4r^2}
			+
			O(|r|^{-5}).
		\end{align}
		Therefore, \eqref{260602845}, \eqref{260602846}, and the boundedness of
		\(\Phi''/\Phi\) on compact subsets of \(\bR\) imply
		\[
		\left|\frac{\Phi''(r)}{\Phi(r)}\right|
		\lesssim
		\left(
		\frac{\eta(r)}{1+|r|}
		\right)^2.
		\]
		Combining these estimates with \eqref{26021134}, we get
		\[
		\frac{|\partial_vH(x,v)|}{H(x,v)}
		\simeq
		\frac{\eta(x^{-1/3}v)}{x^{1/3}+|v|}
		\]
		and
		\[
		\frac{|\partial_{vv}H(x,v)|}{H(x,v)}
		\lesssim
		\left(
		\frac{\eta(x^{-1/3}v)}{x^{1/3}+|v|}
		\right)^2.
		\]
		This proves \eqref{260603518} and completes the proof.
	\end{proof}
	
	\section{Kinetic mollification and a classical comparison principle}\label{app:kinetic-mollification}
	We record the mollification and comparison tools used to transfer the grazing
	barrier estimate to weak solutions in the proof of Theorem~\ref{260129331}.
	
	\begin{lemma}[Kinetic mollification]\label{2606061211}
		Let \(\zeta\in C_c^\infty(\bR^{1+2d})\) satisfy
		\[
		\zeta\ge 0,
		\qquad
		\iiint_{\bR^{1+2d}}\zeta(t,x,v)\dd v\dd x\dd t=1.
		\]
		For \(\varepsilon>0\), define
		\[
		\zeta_{\varepsilon,\star}(t,x,v)
		=
		\varepsilon^{-(2+4d)}
		\zeta(\varepsilon^{-2}t,\varepsilon^{-3}x,\varepsilon^{-1}v).
		\]
		For \(Z=(t,x,v)\) and \(Y=(s,y,w)\), set
		\[
		Y^{-1}\circ Z
		:=
		(t-s,x-y-(t-s)w,v-w).
		\]
		We define the kinetic mollification of \(f\) by
		\begin{align}\label{2606061226}
			f^{(\varepsilon),\star}(Z)
			:=
			\iiint_{\bR^{1+2d}}
			\zeta_{\varepsilon,\star}(Y)
			f(Y^{-1}\circ Z)\dd Y .
		\end{align}
		Then the following hold:
		\begin{enumerate}
			\item For any \(f\in L^1_{\mathrm{loc}}(\bR^{1+2d})\),
			\[
			f^{(\varepsilon),\star}\in C^\infty(\bR^{1+2d}),
			\]
			and
			\[
			f^{(\varepsilon),\star}(Z)\to f(Z)
			\quad\text{as }\varepsilon\downarrow0,
			\text{ for a.e. }Z\in\bR^{1+2d}.
			\]
			
			\item If \(f\in L^p(\bR^{1+2d})\), \(1\le p< \infty\), then
			\[
			\|f^{(\varepsilon),\star}\|_{L^p(\bR^{1+2d})}
			\le
			\|f\|_{L^p(\bR^{1+2d})},
			\qquad
			\lim_{\varepsilon\to0}
			\|f-f^{(\varepsilon),\star}\|_{L^p(\bR^{1+2d})}
			=0.
			\]
			
			\item If \(f\in L^\infty(\bR^{1+2d})\), then
			\[
			\|f^{(\varepsilon),\star}\|_{L^\infty(\bR^{1+2d})}
			\le
			\|f\|_{L^\infty(\bR^{1+2d})}.
			\]
			If, in addition, \(f\in C_{\mathrm{loc}}(\bR^{1+2d})\), then
			\[
			f^{(\varepsilon),\star}\to f
			\quad \text{locally uniformly in } \bR^{1+2d}.
			\]
			
			\item In the sense of distributions, the mollification commutes with
			the velocity derivatives and with the kinetic transport operator:
			\[
			D_{v_i}f^{(\varepsilon),\star}
			=
			(D_{v_i}f)^{(\varepsilon),\star},
			\qquad
			\cT f^{(\varepsilon),\star}
			=
			(\cT f)^{(\varepsilon),\star},
			\quad i=1,\dots,d,
			\]
			where \(\cT=\partial_t+v\cdot\nabla_x\).
		\end{enumerate}
	\end{lemma}
	
	\begin{proof}
		The kinetic group with dilations
		\[
		\mathscr S_\varepsilon(t,x,v)
		=(\varepsilon^2t,\varepsilon^3x,\varepsilon v)
		\]
		is a homogeneous group of homogeneous dimension \(2+4d\), and
		\(\{\zeta_{\varepsilon,\star}\}_{\varepsilon>0}\) is an approximate
		identity on this group. After the measure-preserving change of variables
		\(W=Y^{-1}\circ Z\), differentiation falls on the smooth kernel; hence
		\(f^{(\varepsilon),\star}\in C^\infty\).
		The almost-everywhere convergence in (1) follows from
		\cite[Theorem~2.6]{folland1982hardy}, while the contraction, \(L^p\)
		convergence, and locally uniform convergence in (2)--(3) are the standard
		approximate-identity properties; see
		\cite[Propositions~1.18 and~1.20]{folland1982hardy}.
		
		Finally, (4) follows by direct differentiation for smooth \(f\), and the distributional identities follow by duality.
	\end{proof}
	
	\begin{lemma}[Classical comparison principle in a kinetic half-cylinder]
		\label{lem:classical-kinetic-comparison}
		Let \(r>0\), and assume that
		\[
		a^{ij},\,\bar b^i\in C(\overline{Q_r^+(0)}),
		\qquad i,j=1,\dots,d,
		\]
		and
		\[
		a^{ij}(Z)\xi_i\xi_j\ge 0\qquad \text{for all }Z\in \overline{Q_r^+(0)}
		\,,\,\,\xi\in\bR^d.
		\]
		Here \(\gamma^-\), \(\gamma^+\) and \(\gamma^0\) are the incoming,
		outgoing and grazing parts of \(\{x^d=0\}\), namely
		\(\{x^d=0,\ v^d>0\}\), \(\{x^d=0,\ v^d<0\}\) and
		\(\{x^d=0,\ v^d=0\}\).
		Let
		\[
		w\in C(\overline{Q_r^+(0)})
		\cap
		C^\infty\Big(
		Q_r^+(0)
		\cup
		\big((\bR\times\gamma^+)\cap Q_r(0)\big)
		\Big),
		\]
		where \(C^\infty(E)\) means that the function has a \(C^\infty\) extension
		to some open set containing \(E\). Suppose that
		\[
		(\partial_t+v\cdot\nabla_x)w
		-
		a^{ij}D_{ij}w
		-
		\bar b^iD_iw
		\ge 0
		\quad\text{in }Q_r^+(0).
		\]
		Assume also that
		\[
		w\ge0
		\quad\text{on }
		\cA_r
		\cup
		\big((\bR\times(\gamma^-\cup\gamma^0))\cap Q_r(0)\big),
		\]
		where
		\[
		\cA_r
		:=
		\overline{Q_r^+(0)}
		\cap
		\left(
		\{t=-r^2\}\cup\{|x|=r^3\}\cup\{|v|=r\}
		\right).
		\]
		Then
		\[
		w\ge0
		\quad\text{in }\overline{Q_r^+(0)}.
		\]
	\end{lemma}
	
	\begin{proof}
		Set
		\[
		\cL w
		:=
		(\partial_t+v\cdot\nabla_x)w
		-
		a^{ij}D_{ij}w
		-
		\bar b^iD_iw .
		\]
		For \(\eta>0\), define
		\[
		w_\eta(t,x,v):=w(t,x,v)+\eta(t+r^2).
		\]
		Since \(t+r^2\ge0\) in \(\overline{Q_r^+(0)}\), we have
		\[
		w_\eta\ge0
		\quad\text{on }
		\cA_r
		\cup
		\big((\bR\times(\gamma^-\cup\gamma^0))\cap Q_r(0)\big).
		\]
		Moreover,
		\[
		\cL w_\eta
		=
		\cL w+\eta
		\ge \eta
		\quad\text{in }Q_r^+(0).
		\]
		
		We claim that
		\[
		w_\eta\ge0
		\quad\text{in }Q_r^+(0).
		\]
		Suppose, to the contrary, that \(w_\eta\) is negative somewhere in
		\(Q_r^+(0)\). Since \(w_\eta\in C(\overline{Q_r^+(0)})\), there exists
		\[
		Z_*=(t_*,x_*,v_*)\in \overline{Q_r^+(0)}
		\]
		such that
		\[
		w_\eta(Z_*)
		=
		\min_{\overline{Q_r^+(0)}} w_\eta
		<0.
		\]
		By the nonnegativity of \(w_\eta\) on the comparison boundary,
		\[
		Z_*\notin
		\cA_r
		\cup
		\big((\bR\times(\gamma^-\cup\gamma^0))\cap Q_r(0)\big).
		\]
		Hence
		\[
		Z_*\in
		Q_r^+(0)
		\cup
		\big((\bR\times\gamma^+)\cap Q_r(0)\big),
		\]
		which is precisely the set on which \(w\) is assumed to be smooth.
		
		We first record the consequences of the minimality in the \(v\)-variables.
		Since \(Z_*\notin\cA_r\), we have
		\[
		|v_*|<r.
		\]
		Thus \(Z_*\) is an interior point with respect to the \(v\)-variables. Hence,
		for every \(\xi\in\bR^d\), the curve
		\[
		s\mapsto (t_*,x_*,v_*+s\xi)
		\]
		stays in
		\[
		Q_r^+(0)
		\cup
		\big((\bR\times\gamma^+)\cap Q_r(0)\big)
		\]
		for all sufficiently small \(|s|\). Since \(Z_*\) is a minimum point of
		\(w_\eta\), the function
		\[
		s\mapsto w_\eta(t_*,x_*,v_*+s\xi)
		\]
		has a local minimum at \(s=0\). Therefore
		\[
		D_vw_\eta(Z_*)\cdot\xi=0,
		\qquad
		\xi_i\xi_jD_{ij}w_\eta(Z_*)\ge0
		\]
		for every \(\xi\in\bR^d\).
		
		Next we consider the transport direction. We claim that
		\[
		(\partial_t+v\cdot\nabla_x)w_\eta(Z_*)\le0.
		\]
		Indeed, the backward characteristic
		\[
		s\mapsto (t_*-s,x_*-s v_*,v_*),
		\qquad s\ge0,
		\]
		remains in \(\overline{Q_r^+(0)}\) for all sufficiently small \(s>0\):
		since \(Z_*\notin\cA_r\) we have \(t_*>-r^2\), \(|x_*|<r^3\)
		and \(|v_*|<r\), so the first three constraints persist for small
		\(s>0\); and if \(x_*^d=0\), then \(Z_*\in\gamma^+\) forces
		\(v_*^d<0\), whence \(x_*^d-sv_*^d>0\).
		Since \(Z_*\) is a minimum point of \(w_\eta\) on
		\(\overline{Q_r^+(0)}\), we have
		\[
		w_\eta(t_*-s,x_*-s v_*,v_*)
		\ge
		w_\eta(t_*,x_*,v_*).
		\]
		Dividing by \(s>0\) and letting \(s\downarrow0\), we obtain
		\[
		-(\partial_t+v\cdot\nabla_x)w_\eta(Z_*)\ge0,
		\]
		which proves the claim.
		
		Using the nonnegative definiteness of \(a\), together with the
		nonnegativity of \(D_v^2w_\eta(Z_*)\), we get
		\[
		a^{ij}(Z_*)D_{ij}w_\eta(Z_*)\ge0.
		\]
		Also, since \(D_vw_\eta(Z_*)=0\),
		\[
		\bar b^i(Z_*)D_iw_\eta(Z_*)=0.
		\]
		Therefore
		\[
		\begin{aligned}
			\cL w_\eta(Z_*)
			&=
			(\partial_t+v\cdot\nabla_x)w_\eta(Z_*)
			-
			a^{ij}(Z_*)D_{ij}w_\eta(Z_*)
			-
			\bar b^i(Z_*)D_iw_\eta(Z_*)  \\
			&\le 0.
		\end{aligned}
		\]
		On the other hand, \(\cL w_\eta\ge\eta\) in \(Q_r^+(0)\). By the
		assumed regularity, the left-hand side extends continuously to
		\[
		Q_r^+(0)
		\cup
		\big((\bR\times\gamma^+)\cap Q_r(0)\big),
		\]
		and hence
		\[
		\cL w_\eta(Z_*)\ge\eta>0.
		\]
		This is a contradiction. Thus
		\[
		w_\eta\ge0
		\quad\text{in }Q_r^+(0).
		\]
		Letting \(\eta\downarrow0\), we conclude that \(w\ge0\) in
		\(Q_r^+(0)\), and hence, by continuity, in
		\(\overline{Q_r^+(0)}\).
	\end{proof}

	\bibliographystyle{plain}

\end{document}